\documentclass[reqno,11pt]{amsart}
\usepackage[letterpaper,hmargin=1.25in,vmargin=0.9in]{geometry}
\usepackage{fontspec,microtype,amsmath,amssymb,amsthm,mathtools,unicode-math}
\usepackage{longtable,booktabs,array,enumitem,xurl,xcolor}
\usepackage[colorlinks=true,linkcolor=black,citecolor=blue!35!black,urlcolor=blue!35!black]{hyperref}
\newtheorem{theorem}{Theorem}[section]
\newtheorem{lemma}[theorem]{Lemma}
\newtheorem{proposition}[theorem]{Proposition}
\newtheorem{corollary}[theorem]{Corollary}
\theoremstyle{definition}
\newtheorem{definition}[theorem]{Definition}
\newtheorem{remark}[theorem]{Remark}

\newcommand{\R}{\mathbb R}\newcommand{\Z}{\mathbb Z}
\newcommand{\tr}{\operatorname{tr}}\newcommand{\diag}{\operatorname{diag}}
\newcommand{\cost}{\operatorname{cost}}
\newcommand{\Cstar}{\frac{180483}{100000}}
\newcommand{\Hcap}{\frac{29025933167887355473205403}{39478417604337498218323684}}
\newcommand{\Hzero}{\frac{1832290554767847573595563}{2491868845732503577653889}}
\providecommand{\tightlist}{\setlength{\itemsep}{0pt}\setlength{\parskip}{0pt}}
\setlist{itemsep=3pt,topsep=5pt}
\numberwithin{equation}{section}
\makeatletter
\def\@settitle{\begin{center}\baselineskip14\p@\relax\bfseries\@title\end{center}}
\makeatother

\newenvironment{proofappendix}{\begingroup\let\proofsubsection\subsection\let\proofsubsubsection\subsubsection\let\proofparagraph\paragraph\let\section\proofsubsection\let\subsection\proofsubsubsection\let\subsubsection\proofparagraph}{\endgroup}
\providecommand{\doi}[1]{DOI: \href{https://doi.org/#1}{#1}}
\space
\providecommand{\Qcap}{\frac{520}{33}}
\providecommand{\Nmat}{\mathsf N}
\providecommand{\Bmat}{\mathsf B}
\providecommand{\Gram}{\mathsf C}
\providecommand{\budget}{C_*}
\providecommand{\recip}{\mathfrak d}

\title[The Reinhardt conjecture]{A proof of the Reinhardt conjecture}
\author{Fei Xue}\date{\today}
\hypersetup{pdftitle={A proof of the Reinhardt conjecture},pdfauthor={Fei Xue},pdfsubject={All-length classification and exclusion of the long minimizing types}}
\begin{document}
\begin{abstract}
We prove the Reinhardt conjecture: among centrally symmetric convex disks in the plane, the smoothed regular octagon has the least optimal lattice packing density, and equality holds precisely for its nonsingular affine images. Local cost inequalities reduce arbitrarily long switching cycles of a minimizer to six types. We exclude those types by geometric comparisons and rigorous finite verification, and identify the remaining short minimizer.
\end{abstract}
\maketitle
\section{Introduction}\label{sec:introduction}
For a centrally symmetric convex disk $D\subset\mathbb R^2$ with nonempty interior, let
\[
 \delta_L(D)=\sup_{\mathcal L}\frac{|D|}{\det\mathcal L},
\]
where the supremum is over lattices whose translates of $\operatorname{int}D$ are pairwise disjoint. The Reinhardt problem asks for the minimum of $\delta_L(D)$. Nonsingular affine transformations preserve this density.

In 1934, Reinhardt proposed the smoothed regular octagon as a minimizer \cite[p.~230]{Reinhardt1934}. At each corner of a regular octagon, one replaces the vertex by an arc of the hyperbola tangent to the two sides meeting there and asymptotic to the lines extending the next two sides. We prove Reinhardt's conjecture, including the equality statement.
\begin{theorem}[The Reinhardt conjecture]\label{ten:full-reduction}\label{five:main}
Let $D_{\rm oct}$ be the smoothed regular octagon. Every centrally symmetric convex disk $D\subset\mathbb R^2$ with nonempty interior satisfies
\[
 \delta_L(D)\ge\delta_L(D_{\rm oct})
 =\frac{8-4\sqrt2-\log2}{2\sqrt2-1}.
\]
Equality holds precisely when $D$ is a nonsingular affine image of $D_{\rm oct}$.
\end{theorem}

\subsection{Outline of the proof}
After normalizing the area of a least-area circumscribed symmetric
hexagon, we use a harmonic bound on the matrix costate to show that
every global minimizer has a finite cycle of control arcs
(Theorem~\ref{fw:finiteness}). The arcs have three possible control
labels. At each switch
the label changes by $+1$ or $-1$ modulo three, with a prescribed
label shift at the end of a period (Section~\ref{fd:section}). Let
$B$ count the essential arcs and $c_i$ denote their costs. The normalized
hexagon has area $2\sqrt3$ and the disk has area $\sqrt3 C$, where $C$
is the total cost. Thus $\delta_L(D)=C/2$, and the smoothed octagon,
of cost $C_{\rm oct}$, gives the competitor bound
\[
 C=\sum_{i=0}^{B-1}c_i\le C_{\rm oct}<C_*,
 \qquad C_*=\frac{180483}{100000}.
\]
For long cycles the argument proceeds as follows.
\[
 \text{minimizing cycle with } B\ge8
 \;\longrightarrow\;21\text{ types}
 \;\longrightarrow\;\mathrm{R01}\text{--}\mathrm{R06}
 \;\longrightarrow\;\text{none}.
\]

\smallskip
\noindent\emph{Local cost bounds.}
Section~\ref{run:catalogue} records inequalities for short windows of
successive arcs, including the arcs needed to establish their
hypotheses. An inequality may use a whole window and charge only
some of its arcs. This distinction matters when two windows overlap.

\smallskip
\noindent\emph{From local costs to six types.}
Section~\ref{as:section} adds local bounds with total weight at most
one on each arc. A sum of at least $C_*$ rules out the word.
Section~\ref{run:structure} bounds the lengths of maximal runs of equal
increments. A mixed candidate with $B\ge8$ has at most three positive runs
and at most one negative run longer than two. For each possible number of positive
runs we combine compatible cost bounds and impose the label shift.
This leaves six mixed types with one positive run, thirteen with two and one with
three. The all-minus eight-arc word gives the twenty-first type.
Section~\ref{run:twenty-one-to-six} groups the mixed words by their
negative run lengths. Endpoint, return and cost comparisons leave five
of them, which together with the all-minus word form R01--R06.

\smallskip
\noindent\emph{The six exclusions.}
Full signed frame closure and finite covers exclude the three eight-arc
types. For R04 the reduced cost exceeds $C_{\rm oct}$ at every point of
its four-dimensional open domain. For R05 and R06 we combine the full
frame and source conditions and exclude the candidates on closed covers.
This gives $B\le7$. Among the short cycles only \texttt{++++} can
minimize (Section~\ref{five:completion}). Its four control arcs give eight
hyperbolic boundary arcs in opposite pairs. Area comparison in this
family singles out the smoothed regular octagon up to affine equivalence.
Appendix~\ref{app:verification} gives the inputs and procedures for the
finite checks.

\subsection{Historical background}
Minkowski studied lattice packings through critical lattices \cite{Minkowski1904,Minkowski1907}. For a disk centered at the origin, a critical lattice has least determinant among lattices with no nonzero point in the interior of the disk. This least value is the \emph{critical determinant}. Minkowski related critical lattices in the plane to centrally symmetric hexagons inscribed in the disk. The question of the least possible packing density appears in Blaschke's 1923 book, where it is attributed to Courant, who proposed the circular disk \cite[\S27, Problem~17]{Blaschke1923}. The circle has density $\pi/\sqrt{12}=0.906899\ldots$, which exceeds that of the smoothed octagon.

Reinhardt proved the existence of a minimizer and established geometric properties that underlie the variational approach \cite{Reinhardt1934}. For a centrally symmetric disk its best lattice-packing density is the ratio of its area to that of a least-area centrally symmetric circumscribed hexagon. A minimizing disk has a unique tangent at every boundary point, and each boundary point is the midpoint of a side of such a least-area hexagon. Thus the hexagons of optimal packings can be used to follow the boundary.

Mahler independently studied the problem through critical lattices in a series of papers in 1947 \cite{MahlerIrreducible1947,Mahler1947,MahlerHexagons1947}. He proved that a minimizing disk is irreducible. In particular, shrinking it to any proper centrally symmetric convex subdisk with the same center strictly decreases its critical determinant. Mahler's variational analysis led him to suggest that a minimizing boundary consists of straight segments and hyperbolic arcs \cite[p.~117]{Mahler1947}. The assertion that finitely many such pieces suffice is known as Mahler's first conjecture. Mahler also showed that the least density among centrally symmetric octagons is $4(3-\sqrt2)/7$, attained only by nonsingular affine images of the regular octagon. He then rounded its corners, computed the lower density of the resulting disk, and proposed that it might be the global minimum \cite[Theorem~10 and \S12]{MahlerHexagons1947}. Ledermann and Mahler determined the least density among centrally symmetric decagons. They called the rounded candidate the \emph{smoothed octagon} and showed that the analogous smoothed decagon packs slightly more densely \cite{LedermannMahler1949}.

In 1950 Fejes~T\'oth used a theorem of Dowker on circumscribed polygons \cite[Theorem~3]{Dowker1944} to prove that, for every centrally symmetric convex disk, the greatest packing density over all congruent copies equals $\delta_L$ \cite[pp.~62--63]{FejesToth1950}. Fejes~T\'oth also recorded that Mahler had announced this conclusion in 1949. Rogers proved that for any planar convex disk the greatest density over packings by translates equals its lattice-packing density \cite{Rogers1951}. These three packing densities therefore agree for centrally symmetric convex disks.

Mahler obtained the universal lower bound $\sqrt3/2$ in 1946 \cite[Lemma~2]{Mahler1946}. Fejes~T\'oth also obtained this bound in 1948 \cite{FejesToth1948}. Ennola improved it to $0.88$ in 1961, and Tammela raised it to $0.892656$ in 1970 \cite{Ennola1961,Tammela1970}. Nazarov proved in 1986 that the smoothed octagon is a local minimizer in the Hausdorff metric \cite{Nazarov1988}. For surveys of these results and the surrounding planar packing problems, see Brass, Moser and Pach \cite[\S1.2]{BrassMoserPach2005} and Zong \cite{Zong2014}.

Tammela's planar estimate also entered Smith's three-dimensional packing lower bound of $0.46421\ldots$ for centrally symmetric convex bodies \cite{Smith1999}. In the same dimension, Kallus showed that the Euclidean ball locally minimizes the best lattice-packing density among centrally symmetric convex bodies \cite{Kallus2014}.

Hales made the variational problem explicit in 2011 by expressing the area as an integral along the motion of six contact points on the boundary \cite{Hales2011}. The Euler--Lagrange equations single out the circle up to affine equivalence, and the second variation at the circle gives a direction of decreasing area. Hales also proved that a piecewise analytic minimizer is a smoothed polygon. In 2017 he formulated the motion as an optimal-control problem on $\mathrm{SL}_2(\mathbb R)$, with area as the cost \cite{Hales2017}. A constant control at a vertex of the control triangle traces one arc of a smoothed polygon. Hales showed that the smoothed regular $(6k+2)$-gons, $k\ge1$, satisfy the necessary conditions of the Pontryagin maximum principle, and established a strict local comparison near the octagon's lifted trajectory \cite[Theorems~4.4.1 and~4.5.1]{Hales2017}.

Vajjha developed the Hamiltonian formulation and analyzed the dynamics near the singular locus, where the conditions for choosing a control can degenerate \cite{Vajjha2022}. Hales and Vajjha confined the state of a global minimizer to a compact domain and treated Mahler's first conjecture through the dynamics of accumulating switches \cite[Theorem~5.0.2 and \S16.3]{HV2024}. We also give a proof of finite switching for a global minimizer in Theorem~\ref{fw:finiteness}.

\space
\section*{Acknowledgments}
The author proposed the argument giving an explicit finite upper bound for the number of control arcs. ChatGPT dot carried out the sustained exploration of local cost lower bounds. ChatGPT 6.1 Sol carried out the proofs of the six remaining types R01--R06 and the organization of the manuscript. The author retains responsibility for the mathematical claims and their presentation.

The author thanks Chuanming Zong and Thomas Hales for the seminar discussions around 2015 through which they introduced him to this problem. He also thanks Zong for sharing historical references. For more than a decade, leaving this conjecture unsolved made him restless. The purpose of this paper is to show that the conjecture can be resolved. The present proof is long, and the author hopes that further work will lead to a shorter and more beautiful proof.

The-Anh Vu-Le \cite{VuLe2026} and Lucas Schmutz \cite{Schmutz2026} independently completed proofs of the Reinhardt conjecture, including the equality case, during the same period as the present work. Vu-Le's proof is formalized in Lean~4, and Schmutz's proof is computer-assisted. Both works state the result for packings by congruent copies.
In the same period, Ho \cite{Ho2026} formalized in Lean the strict
lower bound $0.8950112$ for the optimal density of packings by
congruent copies of any centrally symmetric convex disk. These
statements apply to lattice packing by Fejes~T\'oth's theorem.

The accompanying reader archive (version~1118R01) contains all the algorithms used in the computer-assisted steps, together with the fixed inputs and instructions for reproducing the computations. The computations can be reproduced on little-endian x86-64 Linux with the dependencies specified in the archive. The archive provides a Dockerfile based on Ubuntu~24.04 for this platform.
\space
\clearpage\setcounter{tocdepth}{1}\tableofcontents\clearpage

\section{The geometric problem and its normalization}\label{fd:section}
The finite reduction is applied to a global minimizer of the original
geometric problem.  We first specify the trajectory, the full lift, and the
changes of coordinates under which the local estimates will be used.  The foundational interfaces are stated here; their detailed proofs
appear after the finite assembly, and are referenced explicitly below.

\subsection{The original minimizing class}
For a centrally symmetric convex disk $D$ with nonempty interior, write
\[
 \delta_L(D)=\sup_{\mathcal L}\frac{|D|}{\det\mathcal L},
\]
where the supremum is over lattices whose translates of $D$ have pairwise
disjoint interiors.  An \emph{original global minimizer} is a disk attaining
the infimum of this quantity over all centrally symmetric convex disks.
The density is unchanged by nonsingular affine transformations.  The
normalization used here fixes the area of a minimum-area circumscribed
centrally symmetric hexagon to $2\sqrt3$ and represents its moving frame
in $\mathrm{SL}_2(\mathbb R)$.

The geometric formulation, existence and reconstruction results are the
inputs of \cite[Sections 2.2--2.3]{Hales2017} and
\cite[Chapters 2--3]{HV2024}.  In particular, the regular representation
of an arbitrary original minimizer uses
\cite[Lemma 2.2.5, Lemmas 2.4.2--2.4.3, and Proposition 3.1.9]{HV2024}.
The statements below concern that minimizing disk and its actual
trajectory, rather than merely an equal-cost representative of a relaxed
problem.

Set
\[
 D_* =\operatorname{diag}(3^{1/4},3^{-1/4}),\qquad
 J=\begin{pmatrix}0&-1\\1&0\end{pmatrix},\qquad
 R_* =\frac12\begin{pmatrix}1&-3\\1&1\end{pmatrix}
      =D_*e^{\pi J/3}D_*^{-1}.
\]
For $W=(u,v)$, put
\begin{equation}\label{fd:domain}
 \Omega=\{(u,v):-1<u<1,\ v>0,\ u^2+v^2>1\},\qquad
 S=u^2+v^2-1.
\end{equation}
The three vertex fields and the frame velocity are
\begin{align}\label{fd:fields}
 F_0&=\left(v,\frac{v^2}{1+u}\right),&
 F_1&=\left(-\frac{v(1-u^2+v^2)}{S},\frac{2uv^2}{S}\right),&
 F_2&=\left(v,\frac{v^2}{u-1}\right),\\[-2pt]
 X(W)&=\frac1v\begin{pmatrix}u&-u^2-v^2\\1&-u\end{pmatrix},&
 A_*&=\begin{pmatrix}0&3/2\\-1/2&0\end{pmatrix}.\nonumber
\end{align}
The admissible controls are all measurable maps into the full triangle
$\Delta=\{\rho\in\mathbb R^3:\rho_j\ge0,\ \sum_j\rho_j=1\}$.
The equations and physical cost are
\begin{gather}\label{fd:dynamics}
 W'=\sum_{j=0}^2\rho_jF_j(W),\qquad G'=GX(W),\qquad
 G(t)\in\mathrm{SL}_2(\mathbb R),\\
 \ell(W)=\operatorname{tr}(A_*X(W))
        =\frac{3+u^2+v^2}{2v},\qquad
 C=\int_0^T\ell(W(t))\,dt,\qquad |D|=\sqrt3 C.
 \label{fd:cost}
\end{gather}
Here $t$ is physical time; the identity $X^2=-I$ fixes its scale.
The six curves
\[
 G(t)q_j,\qquad q_j=D_*(\cos(j\pi/3),\sin(j\pi/3))^t,
 \qquad 0\le j\le5,
\]
reconstruct the normalized boundary.

The group and state endpoints are
\begin{equation}\label{fd:endpoints}
 G(0)=I,\qquad G(T)=R_*,\qquad
 W(T)=R_*^{-1}\!\cdot W(0),
\end{equation}
with the usual fractional linear action on $u+iv$. The extended frame
satisfies $G_e(t+T)=G_e(t)R_*$ by \cite[Section~3.5]{HV2024}.
Differentiation gives $X(T)=R_*^{-1}X(0)R_*$ and hence the state endpoint.

The original boundary also determines the real winding of the frame.
Put $\widetilde G=D_*^{-1}GD_*$ and $U=e^{\pi J/3}$. The deformation
retraction onto $\mathrm{SO}(2)$ gives a homotopy with fixed endpoints
from $\widetilde G$ to $e^{(\pi/3+2\pi n)Jt/T}$ for some integer $n$.
Apply this homotopy to each of the six boundary arcs
$D_*^{-1}G(t)q_j=\widetilde G(t)U^j e_1$, where $e_1=(1,0)^t$.
Their endpoints stay fixed and the arcs avoid the origin throughout
the homotopy. The six rotation arcs at the end have total angle
$6(\pi/3+2\pi n)$, so the boundary has winding number $1+6n$ about the
origin. The original convex boundary is traversed once counterclockwise
about this interior point. Thus $n=0$, and $\widetilde G$ is homotopic
with fixed endpoints to $e^{\pi Jt/(3T)}$.

For clarity, the weights in \eqref{fd:dynamics} preserve the original
measurable-control class.  If $\mu$ denotes the original barycentric
control, then
\begin{equation}\label{fd:controls}
 (d_0,d_1,d_2)=\left(\frac{1+u}{v},\frac{S}{2v},\frac{1-u}{v}\right),
 \qquad
 \rho_j=\frac{\mu_jd_j}{\sum_k\mu_kd_k},\qquad
 \mu_j=\frac{\rho_j/d_j}{\sum_k\rho_k/d_k}.
\end{equation}
All $d_j$ are positive on $\Omega$, so these maps are measurable inverses
and preserve the vertices.  The half-plane coordinates of \cite{HV2024}
are $(u,v)/\sqrt3$; its control $\alpha$ corresponds to
$\mu=(\alpha_2,\alpha_0,\alpha_1)$.  The matrix conjugation and the
commutator identity establishing this convention are given later in
equation~\eqref{fw:reciprocals}.  No change of physical time
is involved.

\subsection{Compactness and admissible replacements}
\begin{theorem}[Geometric compactness]\label{fd:compact}
Every normalized original global minimizer has its state trajectory in the
compact subset $K$ of $\Omega$ specified by the following six inequalities:
\begin{equation}\label{fd:compact-cuts}
 \begin{gathered}
 v\le15(1+u),\qquad v\le15(1-u),\qquad S\ge2v/15,
 \qquad v\le9\sqrt3/2,\\
 (u+1)^2+v^2\ge\frac{8v}{9\sqrt3},\qquad
 (u-1)^2+v^2\ge\frac{8v}{9\sqrt3}.
 \end{gathered}
\end{equation}
On this set, the rates
$a=v/(1+u)$, $b=v/(1-u)$ and $c=2v/S$ obey
\begin{equation}\label{fd:compact-bounds}
 0<a,b,c\le15,\qquad a+b,a+c,b+c\le\frac{750}{37},\qquad
 \frac{99}{200}\le v\le\frac{39}{5}.
\end{equation}
\end{theorem}
The compactness assertion is exactly
\cite[Definition 5.0.1 and Theorem 5.0.2]{HV2024}, written in the present
coordinates; the elementary consequences \eqref{fd:compact-bounds}
are proved in Lemma~\ref{fw:compact-bounds}.
These inequalities constrain the original minimizing trajectory.
They are not additional restrictions on a local competitor: the
admissible replacement problem remains the full control triangle in
$\Omega$.

\begin{lemma}[Five-endpoint replacement]\label{fd:replacement}\label{fw:replacement}
Let $[a,b]$ be a window of an admissible trajectory.  A sufficiently small
family of solutions of \eqref{fd:dynamics}, with states in $\Omega$ and
controls in $\Delta$, can replace this window if it has the same entrance
and exit values of both $W$ and $G$ and converges uniformly to the original
window after orientation-preserving reparametrization.  Its duration may
vary and tend to $b-a$.  The resulting curve is an admissible convex
competitor, and its area change is $\sqrt3$ times the window-cost change.
\end{lemma}
The five endpoint coordinates are the three group coordinates and the
two state coordinates. The proof in Section~\ref{sec:framework} checks
curvature, real winding, and variable-duration attachment.
Lemma~\ref{fw:balanced-reconstruction} gives the hexagon normalization
of the reconstructed boundary.
The five-constraint viewpoint also appears in
\cite[Section 7.6]{Hales2011}.  Durations are varied about positive arcs
with positive caps; signed analytic extensions used for differentiation
are not asserted to be admissible competitors.

\subsection{A full maximizing lift and its normalization}
\begin{definition}\label{fd:full-lift}
A \emph{full maximizing lift} on a genuine minimizing interval consists
of the constant cost multiplier $\eta\ge0$, the state covector
$p=(p_u,p_v)^t$, and the entire traceless matrix adjoint $\Lambda$,
not all identically zero, supplied by the free-terminal-time maximum
principle for that interval.  Put
\[
 Z=\Lambda-\eta A_*,\qquad U=\operatorname{tr}(ZX),\qquad
 H_j=p^tF_j+U.
\]
The equations and maximizing conditions used below are
\begin{equation}\label{fd:PMP}
 \begin{gathered}
 p'=-\sum_j\rho_j(DF_j)^tp-\nabla_WU,\qquad
 \Lambda'=[\Lambda,X],\\
 \max_j H_j=0,\qquad \rho_j>0\ \Longrightarrow\ H_j=0
 \quad\hbox{almost everywhere}.
 \end{gathered}
\end{equation}
The lift is \emph{normal} when $\eta>0$.  Its normalization means dividing
$\eta$, both components of $p$, and every entry of $\Lambda$ by that
same positive number, so that $\eta=1$.
\end{definition}
The source is \cite[Sections 6.2--6.8]{HV2024}, with matrix pairing
$\operatorname{tr}(\Lambda G^{-1}\delta G)$.  There is no extra endpoint
condition fixing $T$.  In particular an isolated binary arc does not
supply a freely selectable substitute for the full lift.

It is useful to fix the switching convention.  With $g_0=1-u^2$, set
\begin{align}\label{fd:switching-rows}
 \alpha_0&=(0,S/v),\nonumber\\
 \alpha_1&=\left(-g_0,\frac{(1-u)((1+u)^2-v^2)}{2v}\right),\\
 \alpha_2&=\left(g_0,\frac{(1+u)((1-u)^2-v^2)}{2v}\right),
 \qquad \sigma_i=\alpha_i p.\nonumber
\end{align}
Their sum vanishes, and with $\beta=g_0S/(2v^3)>0$,
\begin{equation}\label{fd:ties}
 H_0-H_2=\sigma_0/\beta,\qquad
 H_1-H_0=\sigma_1/\beta,\qquad
 H_2-H_1=\sigma_2/\beta.
\end{equation}
For measurable maximizing controls the $\sigma_i$ are $C^1$ with
absolutely continuous first derivatives
(Lemma~\ref{fw:regularity}).  Thus the contact derivatives
in later arguments belong to the same continuous switching response on
the two sides of an actual switch.

\begin{proposition}[Switching restrictions and normality]\label{fd:normality}
The following statements hold for the actual positive maximal control arcs.
\begin{enumerate}
 \item The words $2020$, $0101$, and $1212$ cannot occur on an original
 global minimizer, without a normality assumption.
 \item On an interval containing six consecutive positive maximal arcs
 and avoiding those three words, every selected nonzero full free-time
 maximizing lift is normal.
 \item The words $202101$, $010212$, and $121020$ cannot occur on a full
 maximizing lift with positive arcs and caps.  Their durations may be any
 positive numbers.
\end{enumerate}
\end{proposition}
These are Lemma~\ref{sw:four},
Theorem~\ref{sw:normality}, and Lemma~\ref{sw:six}, respectively.
The first proof uses an independent edge problem; it does not presume
that its edge multiplier is a projection of a chosen full lift.
The second proves normality on the very interval later used, without a
compactness or numerical cost hypothesis.  In all applications with
$B\ge8$ below, the genuine containing fundamental interval contains six
complete positive arcs, so a selected full lift can be jointly normalized
before any normal identity is used.

The determinant
\begin{equation}\label{fd:casimir}
 \mathcal C=\det\Lambda=-\tfrac12\operatorname{tr}(\Lambda^2)
\end{equation}
is constant along one lift.  It is distinct from the physical cost $C$.
Consider a normal full lift on a genuine fundamental interval, jointly
normalized to $\eta=1$. If $\mathcal C>0$ and
$\operatorname{tr}(\Lambda X)>0$ at one actual point, then
Theorem~\ref{cal:harmonic} gives
\begin{equation}\label{fd:harmonic}
 \mathcal C\le\frac34\left(\frac{2|D|}{\pi}-1\right)^2.
\end{equation}
For a global minimizer the comparison with the smoothed octagon in
\eqref{fd:budget} then gives $\mathcal C<H_*:=\Hcap$.
Normality, positive determinant and positive trace must all concern the
same selected lift.  Neither positive determinant nor the future-sheet
condition is presumed for an arbitrary normal extremal.  The proof uses
only the lawful three-copy transport of the matrix adjoint and never
assumes a periodic state covector.

\subsection{Whole symmetries and genuine changes of cut}
\begin{proposition}[Lawful transport]\label{fd:lawful-symmetry}
Cyclic relabeling, and time reversal combined with $j\mapsto k-j$, preserve
the problem when applied to the state, frame and whole lift together.
An odd permutation of the labels alone is not a permitted transformation.
\end{proposition}
Here are explicit formulas fixing this convention.  Put
$\Phi(w)=(w+3)/(1-w)=R_*^{-1}\!\cdot w$.
A cyclic phase acts by
\[
 W_{\rm new}=\Phi W,\quad G_{\rm new}=R_*^{-1}GR_*,\quad
 p_{\rm new}=(D\Phi)^{-t}p,\quad
 \Lambda_{\rm new}=R_*^{-1}\Lambda R_*,\quad j_{\rm new}=j+1.
\]
It preserves $\eta$, $\ell$, the Hamiltonians and physical time, since
$D\Phi F_j=F_{j+1}\circ\Phi$ and
$X(\Phi W)=R_*^{-1}X(W)R_*$.
For a full reflection and time reversal let
$\psi(u,v)=(-u,v)$ and $Q=\operatorname{diag}(1,-1)$.
With every old quantity on the right evaluated at $T-t$, set
\begin{equation}\label{fd:reflection}
 \begin{gathered}
 W_{\rm new}(t)=\psi W(T-t),\qquad
 G_{\rm new}(t)=QG(T)^{-1}G(T-t)Q,\\
 p_{\rm new}(t)=-(D\psi)^{-t}p(T-t),\qquad
 \Lambda_{\rm new}(t)=-Q\Lambda(T-t)Q,\qquad
 j_{\rm new}(t)=2-j(T-t).
 \end{gathered}
\end{equation}
These formulas retain $\eta$ and physical cost.  The identities
$X(\psi W)=-QX(W)Q$, $QR_*Q=R_*^{-1}$, and
$\psi\Phi\psi=\Phi^{-1}$ verify the equations and endpoints.
Reflection and reversal also preserve the prescribed real homotopy
class.  Composing with a cyclic phase gives each $j\mapsto k-j$.
The same formulas, with the appropriate entrance normalization of the
frame, apply to a local window.  Each cap, return, contact and charged
arc is transported along with the complete source.

The original trajectory has the twisted extension
\begin{equation}\label{fd:extension}
 \begin{gathered}
 G_e(t+kT)=G(t)R_*^k,\qquad W_e(t+kT)=\Phi^kW(t),\\
 \rho_{j+k,e}(t+kT)=\rho_j(t),\qquad 0\le t<T,\quad k\in\mathbb Z.
 \end{gathered}
\end{equation}
No state costate is extended by this definition.
\begin{lemma}[Genuine reanchoring]\label{fd:reanchor}\label{fw:reanchor}
For any $s\in\mathbb R$, the trajectory
\[
 G_s(t)=G_e(s)^{-1}G_e(s+t),\qquad W_s(t)=W_e(s+t),
 \qquad 0\le t\le T,
\]
is an original global minimizer with the same physical cost and the
original endpoints and real winding.  It therefore has its own nonzero
full maximizing lift.
\end{lemma}
The proof in Section~\ref{sec:framework} identifies the same boundary
under an area-preserving transformation and verifies the real homotopy.
One must select the lift on the new genuine interval; independently
selected state covectors are not identified across a change of cut.

\subsection{Physical arc counts and the global budget}
Theorem~\ref{fw:finiteness} gives a finite vertex-control representation of
every normalized original global minimizer. Its proof combines the
harmonic bound with the maximum principle on a fixed edge of $\Delta$.
Merge adjacent equal controls and discard zero-duration pieces.
Let $B$ be the number of essential switches on the twisted physical
circle.  It is also the number of complete physical arcs in that circle.
Writing their successive labels in the extension as $j_i$, we have
\begin{equation}\label{fd:twist}
 j_{i+B}=j_i+1\pmod3.
\end{equation}
If $N$ counts the positive control arcs after a prescribed fundamental
cut, that cut splits at most one complete physical arc.  Consequently
\begin{equation}\label{fd:cut-count}
 N\le B+1.
\end{equation}
The upper argument first bounds $B$; it then uses this inequality to
bound $N$, keeping the prescribed cut.

Let $c_i>0$ be the cost of complete physical arc $i$, with indices modulo
$B$.  Cost invariance under $\Phi$ gives $c_{i+B}=c_i$, even though the
labels have the twist \eqref{fd:twist}.  Therefore
\begin{equation}\label{fd:physical-cost}
 C=\sum_{i=0}^{B-1}c_i.
\end{equation}
The smoothed-octagon competitor of \cite[Section 7.5]{HV2024} gives
\begin{equation}\label{fd:budget}
 C\le C_{\rm oct}:=\frac{2(8-4\sqrt2-\log2)}{2\sqrt2-1}
   < C_*:=\frac{180483}{100000}.
\end{equation}
The exact rational comparison, including elementary error bounds for
$\sqrt2$, $\sqrt3$, $\log2$ and $\pi$, is supplied with
\eqref{fw:budget} below.  Positivity of every $c_i$ and
the strict upper budget are the global ingredients in the assembly
argument below.

Finally put $\varepsilon_i=1$ or $-1$ according as
$j_{i+1}-j_i\equiv1$ or $-1\pmod3$.  These increments are $B$-periodic,
and their integer walk $J(0)=0$, $J(i+1)=J(i)+\varepsilon_i$ has
$j_i\equiv j_0+J(i)\pmod3$ and satisfies
\begin{equation}\label{fd:increment-twist}
 J(B)\equiv1\pmod3,\qquad
 J(qB+r)=qJ(B)+J(r),\quad 0\le r<B.
\end{equation}
Cyclic rotation changes the chosen first switch.  Our type convention
uses the least rotation with $-<+$, without additionally identifying a
word with its reversal or sign complement.  In a mixed minimizing cycle with
$B\ge8$, a proposed occurrence of $\texttt{+-+}$ or $\texttt{+---+}$
uses four or six distinct arcs, respectively. Cutting inside an unused
arc places the whole source in one fundamental interval, even across
the original seam. The first and third restrictions of
Proposition~\ref{fd:normality} then exclude the two patterns.
Source words are always reconstructed from the twisted walk,
never from a periodic continuation of the control labels.
\space
\begin{table}[htbp]
\centering\small
\begin{tabular}{@{}p{.18\textwidth}p{.73\textwidth}@{}}
\toprule
Symbol & Meaning\\
\midrule
$B$ & Number of essential switches in a twisted cycle.\\
$N$ & Number of positive arcs after a prescribed fundamental cut; $N\le B+1$.\\
$\varepsilon_i$ & Increment $+1$ or $-1$ between successive control labels.\\
$J(i),j_i$ & Integer walk with $J(0)=0$ and labels $j_i\equiv j_0+J(i)\pmod3$.\\
$W=(u,v)$ & State in the normalized upper-half-plane coordinates.\\
$G$ & Moving critical-hexagon frame in $\mathrm{SL}_2(\mathbb R)$.\\
$\eta,p,\Lambda$ & Cost multiplier, state covector, and matrix adjoint of a full lift.\\
$\mathcal C$ & The conserved determinant of the matrix adjoint, not a physical cost.\\
$c_i,C$ & Complete physical arc cost and total physical cost $C=\sum_i c_i$.\\
$C_*$ & The rational comparison budget $180483/100000$.\\
\bottomrule
\end{tabular}
\caption{Global notation. Local comparison coordinates are defined at their first use.}
\label{tab:notation}
\end{table}
\space

\section{Local cost bounds and their supporting windows}\label{run:catalogue}
\subsection{Supporting windows and physical charges}\label{as:supporting-windows}
Fix the twisted period of an original global minimizer and its complete
physical costs $c_0,\ldots,c_{B-1}>0$.  A local \emph{support} is an
ordered interval of actual positive maximal arcs in the unwrapped
trajectory, together with any specified endpoints of that interval.
A \emph{charge} specifies which of those arc costs, and with which
nonnegative coefficients, occur on the left side of the local inequality.
Support and charge need not coincide.

Uncharged caps and returns remain hypotheses: deleting them may destroy
an actual switch sign or circle condition.

\begin{lemma}[Placement of a proper source]\label{as:placement}
A consecutive source of $L<B$ complete physical arcs can be placed inside
a genuine minimizing fundamental interval while retaining every source
arc complete and positive.  A full lift can then be chosen on that whole
interval.  If $B\ge8$, it is normal and can be jointly normalized.
\end{lemma}
\begin{proof}
Place the fundamental cut strictly inside the positive arc immediately
preceding the source.  Since $L<B$, the interval ends after the source,
before returning to that cut.  Lemma~\ref{fd:reanchor} supplies the genuine
minimizing trajectory with its original endpoints and real winding.
There are $B-1\ge7$ complete arcs between the two pieces of the cut arc.
The first assertion of Proposition~\ref{fd:normality} applies to this
minimizer, and its second assertion proves normality of any selected
nonzero full lift on the containing interval.  Only then normalize it.
\end{proof}
The strict inequality $L<B$ is a sufficient placement condition, not a
substitute for the reanchoring argument.  If a theorem uses a larger
union of supports on one common lift, the containing interval must retain
that whole union.  Checking properness separately for its constituent
pieces does not establish the common-lift hypothesis.

A qualified application records its complete unwrapped source (including
positive caps, returns and auxiliary arcs), its whole lawful phase or
reflected-time transformation, a genuine containing interval and the
required selected jointly normalized full lift, and every theorem-specific
root, contact, sign, mean or circle hypothesis. If identities link several
pieces, their common-lift requirement is part of these data. Finally it
records the charged offsets, coefficients, lower bound and strictness.
A word match alone supplies neither an endpoint sign nor an adjoint.

If local positions $k$ have coefficients $q_{rk}\ge0$ and map to physical
arc $\pi_r(k)\in\mathbb Z/B\mathbb Z$, their physical row is
\begin{equation}\label{as:physical-row}
 a_{ri}=\sum_{\pi_r(k)=i}q_{rk},\qquad
 a_r\cdot c=\sum_{i=0}^{B-1}a_{ri}c_i.
\end{equation}
Thus repeated charged occurrences of the same physical arc count with
multiplicity.  One cannot replace that sum by a set union unless a
separate union theorem proves the once-counted charge.  On a proper
consecutive source the positions themselves have distinct physical
images, but different sources can still charge the same arc.

\space
\subsection{Local inequalities used to bound the runs}
The following costs restrict the number and lengths of the runs. They
are statements about original global minimizers. A supporting source
contains all positive caps, returns and sign-producing arcs used in its
proof. Its charge is the complete physical cost of the displayed local
positions. Positions start at zero, and $r:s$ denotes
$\{r,r+1,\ldots,s\}$. For the run classification, each supporting
window in this table is shorter than one physical period. A lawful change
of chart transports the whole source, its full lift and its charge.
Put $P=C_*/5=180483/500000$.

\begin{longtable}{@{}l c c r@{}}
\caption{Complete-source costs for the all-length classification.
Two payments are weak, as indicated.}\label{run:catalogue-table}\\
\toprule Complete source&Charge&Lower bound&Proof\\\midrule\endfirsthead
\toprule Complete source&Charge&Lower bound&Proof\\\midrule\endhead
\bottomrule\endfoot
$\mathtt{102021}$&$1:4$&$>389/1000$&\ref{loc:A389}\\
$\mathtt{1020212}$&$1:4$&$>473/1000$&\ref{loc:E473}\\
$\mathtt{01020212}$&$2:5$&$>74/125$&\ref{lc:two-return}\\
$\mathtt{01020212}$&$3:4$&$>9/25$&\ref{new:own}\\
$\mathtt{01020212}$, qualified anchors&$3:4$&$>3/7$&\ref{new:own}\\
$\mathtt{02010}$&$1:3$&$>49/200$&\ref{f245:theorem}\\
$\mathtt{010201020}$&$3:5$&$>1/2$&\ref{run:crossed-nine}\\
$\mathtt{101202}$&$1:4$&$>41/100$&\ref{fee:four}\\
$\mathtt{10120121}$&$1:6$&$>609/1000$&\ref{loc:OwnSix609}\\
\midrule
$\mathtt{012021}$&$2:4$&$>61/200$&\ref{ek:short-sources}\\
$\mathtt{0120212}$&$2:4$&$>379/1000$&\ref{lc:returned}\\
$\mathtt{10201020}$&$1:3$&$>61/200$&\ref{ua:seven-sign}\\
$\mathtt{012021201}$&$2:6$&$\ge2P$&\ref{six:lem:H3}, \eqref{six:eq:union-input}\\
$\mathtt{2021201020}$&$2:7$&$>807/1000$&\ref{late:r24:W807}\\
\midrule
$\mathtt{10210}$&$1:3$&$\ge3/5$&\ref{wk:five600}\\
$\mathtt{2021020}$&$3$&$>299/250$&\ref{lc:direct-center}\\
$\mathtt{202102101}$&$3$&$>13/20$&\ref{pf:nine-singletons}\\
$\mathtt{20210212}$&$3:4$&$>13/20$&\ref{opt:pair650}\\
$\mathtt{102012021}$&$4$&$>63/500$&\ref{fee:opposite126}\\
$\mathtt{0102021202}$&$6$&$>3/40$&\ref{late:r18:cap}\\
$\mathtt{012012}$&$2:3$&$>P$&\ref{fp:final}\\
\end{longtable}

The first group pays for positive own arcs, possibly with their ports.
The second pays terminal triples or joins the two sides of a short
negative separator. The last supplies negative interiors, central arcs
and internal forward pairs. The table gives sufficient values; several
of the cited proofs give stronger constants.

\begin{lemma}[Complete-source qualifications]\label{run:source-qualification}\label{as:source-qualification}
Each row of Table~\ref{run:catalogue-table} has its stated physical
payment on its complete proper source. In the $3/7$ row, each of local
anchor positions $2,5$ must also be the owner at position $3$ of a
complete proper $201020$ source in its own whole lawful chart.
\end{lemma}
\begin{proof}
Lemma~\ref{as:placement} gives a genuine containing fundamental
interval, retaining the original state, signed frame endpoints and real
winding. Select a nonzero full maximizing lift there, prove normality,
and then normalize its objective multiplier, matrix adjoint and both
state covectors jointly. A cost theorem retains this one lift wherever
it couples means or switching responses. Independently proved scalar
costs may be added without identifying their selected covectors.

The capped-core, one-return and two-return estimates retain all exterior
arcs. The $74/125$ row uses both returns of $01020212$ and charges
four arcs; the $9/25$ row charges only its own pair. For $3/7$, the
two proper six-sources give physical anchor-density bounds $x_2,x_5<1/6$
by Theorem~\ref{new:sixdensity}. These state facts apply to the selected
own-pair lift even if the six-sources used different containing intervals.
Neither guard is charged. Lemma~\ref{run:crossed-nine} retains both
crossing contexts on one nine-source lift before giving $1/2$; a bare
$02010$ core does not qualify that row. The four-own-arc theorem gives
$103/250>41/100$, and the six-own-arc theorem retains its entire
eighth supporting arc.

In $012021$, Lemma~\ref{cmp:forward-prefix-forcing} gives a strictly
negative entrance for the suffix $12021$. Corollary~\ref{ek:short-sources}
applies Theorem~\ref{loc:ST305} on that same lift and charges positions
$2:4$. With the full return $0120212$, it gives $379/1000$ on the same
positions. For $10201020$, Lemma~\ref{ua:seven-sign} retains positions
$1:7$ for the sign $Ea>1$ and the reversed prefix $(4,3,2,1,0)$ for
the $12021$ tail with entrance $1-Ea<0$. Its payment charges positions
$1:3$.

The nine-source $012021201$ keeps both returned contexts, internal
forward entrance witnesses, outgoing leaves and ties of
Lemma~\ref{six:lem:H3}. The scalar contradiction then gives the weak
payment \eqref{six:eq:union-input}. The ten-source
$2021201020$ keeps the common switching jet and both actual returns
of Proposition~\ref{late:r24:W807}; the two entrance and terminal
circles are retained on its selected full lift. The returns are part of
the source and are not charged again.

The seven-source $2021020$ contains the rooted initial $202$, both
actual contact jets and the full terminal $020$. The nine-source
$202102101$ retains the strict rooted response, regenerated prefix
and complete terminal used in Lemma~\ref{pf:nine-singletons}. That
lemma proves separate singleton inequalities, not two simultaneous
consequences of a maximum. The reverse-pair source $20210212$ has
the complete caps of Theorem~\ref{opt:pair650}.

For $102012021$, its contexts $2012$ at $2:5$ and $0120$ at $3:6$
supply both outward negative entrances by Lemma~\ref{fin:forward-force}.
The two capped triples, crossing and endpoint jets remain on one lift,
so the opposite-center theorem applies without a HIGH predicate.
The complete ten-source $0102021202$ retains the returned four-core
and internal H7 suffix of Lemma~\ref{late:r18:cap}; they establish
the cap endpoint qualifications without a later exterior cap.
The forward-pair and weak-five rows keep their uncharged complete caps.
\end{proof}

We next combine these costs for arbitrary run lengths. The classification
will then follow from length bounds, the twist congruence and explicit
source placements. Finite calculations in some cost proofs concern
those local inequalities rather than an enumeration of increment words.
\space
\subsection{Costs assigned to positive and negative runs}\label{run:profiles}
Consider a mixed minimizing cycle with $B\ge8$. Its maximal positive
and negative runs alternate. Write their lengths in cyclic order as
\[
 (a_1,b_1),\ldots,(a_M,b_M),\qquad
 a_i\ge1,\quad b_i\ge2,\quad b_i\ne3.
\]
The two restrictions on $b_i$ are exactly the forbidden increment
factors $\texttt{+-+}$ and $\texttt{+---+}$. The lengths count
increments, while the costs below belong to complete physical arcs.
No quotient by reversing this cyclic list is taken.

If a positive run starts at increment $s$ and has length $a$, assign
its own arcs $s,s+1,\ldots,s+a$ to that run. A negative run of length
$b$, between positive runs ending at $t$ and starting at $t+b$,
has the interior arcs $t+1,\ldots,t+b-1$. These own sets and
interiors partition the physical period. The two arcs facing the
positive run at $s$ are
\[
 s-1,\qquad s+a+1;
\]
we call them its left and right \emph{ports}. The two ports belonging
to opposite ends of a negative separator coincide when $b=2$ and
are distinct when $b\ge4$.

A local cost profile records a proved lower payment together with
its coefficients on the own arcs and its two ports. Supports retain
all complete caps, returns and other hypotheses of the source theorem.
The port coefficients describe how neighboring profiles can be added;
they are not additional costs. Whenever formal ports represent the
same physical arc, their coefficients are added, as in
\eqref{as:physical-row}. All indices are unfolded before taking
physical positions modulo $B$.

\subsubsection{Short positive runs without port charges}
The own-pair payment has a useful consequence for this assignment:
a singleton positive run can pay for its own arcs while leaving its
ports available to its neighbors.

\begin{lemma}[A singleton between short negative runs]\label{run:own-singleton}
Suppose $B>8$, $a_i=1$ and $b_{i-1}=b_i=2$. The cost of the two
own arcs of this positive run is strictly greater than $9/25$.
If also $a_{i-1},a_{i+1}\ge2$, it is strictly greater than $3/7$.
Both assertions have zero port loads.
\end{lemma}
\begin{proof}
Let $s$ be the starting increment and use a whole cyclic chart so
that $j_s=0$. The complete source $[s-3,s+4]$, with whole label
map $j\mapsto j+2$, is $01020212$. It has charge $\{s,s+1\}$
and anchor owners $s-1,s+2$. Its length eight is less than $B$,
so Lemma~\ref{as:placement} retains the whole source on a genuine
containing interval. Theorem~\ref{new:own} gives the $9/25$ bound.
No cost on either anchor is charged.

For the stronger conclusion the actual sources and maps are
\[
\begin{array}{c|c|c|c}
\text{ordered support}&\text{whole map}&\text{word}&\text{density owner}\\\hline
[s-4,s+1]&j\mapsto j+2&201020&s-1\\
(s+5,s+4,\ldots,s)&j\mapsto-j&201020&s+2.
\end{array}
\]
The ascending witness uses the increment factor $\texttt{++--+}$;
the other comes from the chronological factor $\texttt{+--++}$
by whole reflection and time reversal.
The hypotheses on the two neighboring positive lengths retain the
additional positive increments. Each six-source is proper and
Theorem~\ref{new:sixdensity} gives the corresponding physical
anchor square less than $1/6$. Apply these two state facts to the
single lift selected for the eight-source, and then apply
Theorem~\ref{new:own} again.

The six-sources need not have a common containing interval or
identical multipliers. Only their physical state conclusions are
used. In particular, when $M=2$ the two neighboring positive runs
are the same run, but the two density owners are distinct. The two
six-sources may have a union of length ten even at $B=9$; this does
not identify their covectors or charge their shared support twice.
\end{proof}

\begin{lemma}[A complete nine-source for a crossed core]\label{run:crossed-nine}
A complete proper source $010201020$ in an original global minimizer
satisfies
\[
 c_3+c_4+c_5>\frac12,
\]
where the indices are local to this nine-source.
\end{lemma}
\begin{proof}
Choose one full maximizing lift on a genuine containing interval and
normalize it jointly after normality, as in Lemma~\ref{as:placement}.
Its middle five-source $[2,6]$ is $02010$. The two actual
seven-sources
\[
 [2,8]=0201020,\qquad
 (6,5,\ldots,0)\ \xrightarrow{\ j\mapsto-j\ }\ 0201020
\]
are retained on this same lift, with the whole reflected transport
in the second case. For the central own reciprocals $(p,A,y,B_0,q)$,
use the same coordinates as in Section~\ref{new:conditional}:
\[
 a=\frac{2y}{y+A},\qquad b=\frac{2y}{y+B_0},\qquad
 E=\frac{(y+A)(y+B_0)}{1+y^2}.
\]
Here $a,b$ are local comparison coordinates, not run lengths.
Lemma~\ref{fin:H150} gives
\[
 Ea>\max(1,b),\qquad Eb>\max(1,a).
\]
Whole reversal exchanges the two analytical sides and does not sort
these coordinates. The positive central duration and actual state
identities give $E>1$, $a,b>0$ and $Ea,Eb<2$.
The two positive caps, actual leaf means and all maximizing equations
are those of the selected full lift. Thus every hypothesis of
Definition~\ref{late:r24:crossed-source} holds, and
Theorem~\ref{late:r24:F500} gives the assertion.
This implication uses the complete nine-source, not an unqualified
five-letter occurrence of $02010$.
\end{proof}

\begin{corollary}[A positive run of length two]\label{run:two-positive}
If $B>9$, $a_i=2$ and $b_{i-1}=b_i=2$, its three own arcs have
cost greater than $1/2$, with zero port loads.
\end{corollary}
\begin{proof}
For a run beginning at $s$, set $j_s=0$. The complete source
$[s-3,s+5]$, under $j\mapsto j+2$, is $010201020$.
Its local positions $3,4,5$ are exactly $s,s+1,s+2$.
Apply Lemma~\ref{run:crossed-nine}.
\end{proof}

\subsubsection{Terminal payments and available port capacity}
Put $t_0=379/1000$, $j_0=63/500$ and
$P=180483/500000$. For a positive run of length $a\ge3$,
with both adjacent negative runs of length two, retain the two
returned sources
\[
\begin{array}{c|c|c|c}
\text{ordered support}&\text{whole map}&\text{word}&\text{charge}\\\hline
(s+3,\ldots,s-3)&j\mapsto-j&0120212&\{s-1,s,s+1\}\\[0pt]
[s+a-3,s+a+3]&j\mapsto j-(a-3)&0120212&
 \{s+a-1,s+a,s+a+1\}.
\end{array}
\]
As above $j_s=0$. Each forward prefix supplies its actual nonpositive
entrance and the last arc is its actual return. Theorem~\ref{lc:returned}
gives more than $t_0$ for each charged triple, without an additional
entrance-circle or LOW-mean premise.

The table below records usable positive-run profiles. Every own arc
has coefficient at most one. The pair in the last column gives the
left and right formal port coefficients; all supports are proper
under the displayed restrictions.
\begin{center}
\begin{tabular}{clc}\toprule
Run and neighbors&Strict lower payment&Port loads\\\midrule
$a=1$, adjacent $b=2$, $B>8$&$9/25$&$(0,0)$\\
same, neighboring positive lengths $\ge2$&$3/7$&$(0,0)$\\
$a=2$, adjacent $b=2$, $B>9$&$1/2$&$(0,0)$\\
$a=3$, adjacent $b=2$, $B>9$&$2t_0$&$(1,1)$\\
$a=4$, adjacent $b=2$, $B>9$&$2t_0+j_0$&$(1,1)$\\
$a\ge5$, adjacent $b=2$, $B>9$&
$2t_0+P\lfloor(a-3)/2\rfloor$&$(1,1)$\\\bottomrule
\end{tabular}
\end{center}
The first three rows were just proved. For $a=3$ the two terminal
triples have disjoint own charges. At $a=4$, the complete source
$[s-2,s+6]+2$ is $102012021$ and charges its middle position,
namely $s+2$, above $j_0$ by Theorem~\ref{pf:opposite126}.
Its internal forward contexts $[2,5]=2012$ and $[3,6]=0120$
supply both actual outward negative entrances, as required there.
For $a\ge5$, insert the disjoint internal pairs
\[
 \{s+2+2r,s+3+2r\},\qquad
 0\le r<\lfloor(a-3)/2\rfloor.
\]
Their complete six-sources $[s+2r,s+5+2r]$ lie in the positive
run and phase to $012012$; Proposition~\ref{fp:final} pays more
than $P$ for each. These pairs avoid both terminal charges.
All supports used in these constructions have at most nine arcs.

There is also a half-port choice: the positive and negative profiles
of Lemmas~\ref{comb:half-strips} and~\ref{comb:negative-strips}
have load at most $1/2$ on each incident formal port. Their complete
source proof applies for $B\ge10$. The choices $(0,0)$, the half-port
choices, and $(1,1)$ are alternatives; their lower values cannot all be
added on the same own arcs. Their benefit is that a zero-port choice
on one run can leave full capacity for a terminal payment next to it.

For a separator of length two the two facing ports are a single
physical arc. If a negative profile has coefficient $n$ there and
the two incident positive profiles have coefficients $r,l$, their
compatibility condition is $n+r+l\le1$. At a longer separator the
two ports are different, so each has its own capacity condition.
This applies also to repeated incidences when $M=1$ or $2$.
For example, if the two formal ports of a $(1,1)$ profile coincide,
its physical row has coefficient two there and cannot be selected
with weight one in a capacity-one assembly. Theorem~\ref{as:capacity}
adds precisely the simultaneously compatible profiles.

\subsubsection{Two uniform three-run exclusions}
The following applications show how short-run payments release the
ports needed by a neighboring longer run. They exclude whole
families without first assigning individual names to their words.

\begin{proposition}[Positive lengths $1,2,a$]\label{run:three-pattern}
Suppose a mixed increment cycle has exactly three negative runs,
all of length two, and positive lengths $(1,2,a)$ or $(1,a,2)$ in
cyclic order, where $a\ge4$. It cannot be the switching type of an
original global minimizer.
\end{proposition}
\begin{proof}
Here $B=a+9\ge13$, so every source in the profiles above is proper.
The length-one run has neighboring positive lengths $2,a\ge2$,
and Lemma~\ref{run:own-singleton} pays more than $3/7$ on its own
two arcs. Corollary~\ref{run:two-positive} pays more than $1/2$
on the own three arcs of the length-two run. Neither uses a port.
At $a=4$ the remaining positive run uses its two terminal payments
and middle-center payment. At $a\ge5$ it uses the two terminal
payments and at least one complete forward pair. Since $P>j_0$,
its payment in either case is greater than $2t_0+j_0$.
The three own sets are disjoint, and only the long run charges its
two incident ports. The third separator need not be charged.
Consequently
\[
 C>\frac37+\frac12+2\frac{379}{1000}+\frac{63}{500}
   =\frac{1586}{875}=C_*+\frac{5419}{700000}>C_*.
\]
The twist requires $a-3\equiv1\pmod3$, hence $a\equiv1\pmod3$;
no larger allowed value changes the source argument or weakens its
bound. Both cyclic orders are covered by their whole lawful source
placements. At $a=4$ the two orders give the words
$\mathtt{--+--++--++++}$ and $\mathtt{--+--++++--++}$.
\end{proof}

\begin{corollary}[Positive lengths $1,1,a$]\label{run:three-singletons}
A cycle with exactly three negative runs of length two and positive
lengths $(1,1,a)$, where $a\ge5$, cannot occur in an original global
minimizer.
\end{corollary}
\begin{proof}
Its two singleton positive runs are different physical runs. Each
pays more than $9/25$ on its own arcs and leaves both ports unused.
The remaining run pays more than $2t_0+P$. The charges are compatible,
so
\[
 C>2\frac9{25}+2\frac{379}{1000}+P
  =\frac{919483}{500000}>C_*.
\]
The twist here requires $a\equiv2\pmod3$. The first case $a=5$ gives
\[
 \mathtt{--+--+--+++++}.
\]
Distinctness of the two singleton runs is essential; if
$M=2$ one physical run cannot be charged twice.
\end{proof}
\space

\section{Combining local costs}\label{as:section}
The inequalities of Section~\ref{run:catalogue} constrain the same
physical costs $c_i$, even when their supporting windows overlap. We
give the rule for adding them and illustrate it on two actual words.

\subsection{The capacity theorem, including equality}
\begin{theorem}[Physical-cost assembly]\label{as:capacity}
Let $c_i\ge0$, $0\le i<B$, and $C=\sum_i c_i$.  Suppose finitely many valid
local conclusions on this trajectory have the form
\begin{equation}\label{as:rows}
 a_r\cdot c\ge b_r,\qquad a_r\in\mathbb R^B,
\end{equation}
with specified rows known to be strict. The algebra also allows real row
coefficients; physical charge rows have $a_{ri}\ge0$ as in
\eqref{as:physical-row}.  If numbers $\lambda_r\ge0$
satisfy
\begin{equation}\label{as:capacity-constraint}
 \gamma_i:=\sum_r\lambda_ra_{ri}\le1\qquad(0\le i<B),
\end{equation}
then
\begin{equation}\label{as:lower-sum}
 C\ge L:=\sum_r\lambda_rb_r.
\end{equation}
The conclusion is strict if a strict row has positive weight, or if
$\gamma_i<1$ and $c_i>0$ for some $i$. Equality $C=L$ holds precisely
when $(1-\gamma_i)c_i=0$ for every $i$ and every positively weighted
row is an equality. In particular, positive physical costs force all
capacities to be one in the equality case.
\end{theorem}
\begin{proof}
There is an exact identity
\begin{equation}\label{as:slack-identity}
 C-L
 =\sum_{i=0}^{B-1}(1-\gamma_i)c_i
  +\sum_r\lambda_r(a_r\cdot c-b_r).
\end{equation}
Every summand is nonnegative.  A strict source row with positive weight
makes a summand in the second sum positive.  An unsaturated capacity at a positive cost makes a summand in the
first sum positive.  This
proves the inequality and both strictness claims.  Equality is equivalent to every term in the two sums vanishing.
\end{proof}

\begin{corollary}[Exclusion at the competitor budget]\label{as:budget-exclusion}
For an original global minimizer, a valid capacity combination with
$L\ge C_*$ is impossible.  If only the weaker comparison $C\le C_*$ is
used, the same conclusion holds when $L>C_*$, or when $L=C_*$ and either
strictness condition in Theorem~\ref{as:capacity} holds.
\end{corollary}
\begin{proof}
The actual budget is $C<C_*$ by \eqref{fd:budget}, whereas
Theorem~\ref{as:capacity} gives $C\ge L$.  The assertions with a weak
budget follow from its strict version.
\end{proof}
A weak source row remains weak; strictness is supplied by the stated
positive-weight or residual-capacity condition. Under the actual strict
budget, a weak lower sum exactly $C_*$ already excludes the type.

Only physical costs occur in \eqref{as:slack-identity}.  Supporting windows
may overlap arbitrarily, and independently established scalar inequalities
may use different genuine containing intervals.  Their full lifts need
not be identified.  By contrast, a relation between two means, two
switching derivatives, or a common determinant is not a scalar cost row
until its stipulated common-lift argument has been completed.  It cannot
be assembled by choosing the most favorable component from independent
lifts.

\subsection{Exhaustive alternatives and conditional rows}
For example,
\begin{equation}\label{as:example-max}
 \max(c_i,c_j)>h
\end{equation}
is the alternative $c_i>h$ or $c_j>h$, not two simultaneous rows.

\begin{theorem}[Branchwise physical-cost assembly]\label{as:branches}
Fix a twisted increment word $w$ and let $\mathcal M(w)$ denote the
original minimizing trajectories with that word, equipped with any full
lifts required by the source theorems.  Suppose finitely many predicates
$P_\nu$, $\nu\in V$, cover $\mathcal M(w)$: every such trajectory, with
the stipulated choices of lifts, satisfies at least one $P_\nu$.
For each $\nu$, suppose that under $P_\nu$ the complete source theorems
prove rows
\[
 a_{\nu r}\cdot c\ge b_{\nu r},\qquad a_{\nu ri}\ge0,
\]
with their stated strictness. If for every $\nu\in V$ there exist nonnegative weights
$\lambda_{\nu r}$ satisfying
\[
 \sum_r\lambda_{\nu r}a_{\nu ri}\le1\quad\hbox{for every }i,
 \qquad
 L_\nu:=\sum_r\lambda_{\nu r}b_{\nu r}\ge C_*,
\]
then $\mathcal M(w)$ is empty.  With a weak global budget $C\le C_*$,
replace the last condition by $L_\nu>C_*$, or by $L_\nu=C_*$ together
with a strictness condition from Theorem~\ref{as:capacity}, in every branch.
\end{theorem}
\begin{proof}
Suppose a trajectory in $\mathcal M(w)$ exists.  Choose its lifts as
required by the source hypotheses and select a predicate $P_\nu$ that
it satisfies.  All rows in that branch now hold on its one physical cost
vector.  Corollary~\ref{as:budget-exclusion} gives a contradiction.
The proof uses no row from a different branch.  The predicates need not
be disjoint, since choosing any one satisfied branch suffices.
\end{proof}

Predicates involving means or adjoints must use the same lifts as their
rows; independently proved scalar costs may be included in every branch.
A finite branch tree yields the predicates by conjunction along paths.
Each node's children must cover its entire condition, including equality
boundaries unless separately excluded. The predicates may overlap; a
maximum such as \eqref{as:example-max} still requires a sufficient
complementary payment in each branch.

\subsection{An explicit ten-switch example}
We illustrate the source and charge bookkeeping with two actual
exclusions used in the upper reduction.  This example uses
three local estimates proved later in this paper:
\begin{enumerate}
 \item The complete seven-arc source $2021020$ charges its central
 position $3$. We use the sufficient consequence that its cost is greater
 than $461/500$ of Theorem~\ref{lc:direct-center}.
 \item The right one-return source $1020212$ charges positions $1,2,3,4$;
 the left source $1210102$ charges positions $2,3,4,5$.  Each cost is
 greater than
 \[
 E_*:=\frac{4154032208439}{8759490110000}
 \]
 (Theorem~\ref{lc:one-return}).
 \item The complete six-arc sources $101202$, $212010$, and $020121$
 charge their four interior arcs with cost greater than $41/100$
 (the four-own-arc estimate~\ref{lc:own-four}).
\end{enumerate}
These statements include their whole lawful transforms and require all
source arcs to be positive and complete on a genuine containing interval.
In the last statement the two actual returns supply the endpoint signs;
a four-letter substring alone would not supply them.

\begin{proposition}\label{as:two-ten}
The cyclic increment words
$\texttt{----+--+++}$ and $\texttt{----+++--+}$ cannot occur on an original
global minimizer.
\end{proposition}
\begin{proof}
For each word reconstruct $j_i=J(i)\pmod3$ using
\eqref{fd:increment-twist}, starting with $j_0=0$.
In both cases the source $[9,15]$ is $0102101$, whose whole cyclic phase
$+2$ is $2021020$.  Its centre is arc $12$, or physical arc $2$.
The actual incoming and outgoing switches and all three terminal arcs
are retained.  A cut inside arc $8$ puts this entire source in a genuine
containing interval, and the first source theorem gives
\[
 c_2>461/500.
\]

For $\texttt{----+--+++}$, the source $[2,8]$ is $1020212$.
It charges $\{3,4,5,6\}$, with caps $2,7$ and returned arc $8$.
The source $[6,11]$ is $212010$ and charges physical
$\{7,8,9,0\}$; its two exterior arcs remain part of the source.
Thus the three inequalities charge disjoint sets, leaving physical arc
$1$ unpaid.

For $\texttt{----+++--+}$, the source $[6,12]$ is $1210102$.
It charges physical $\{8,9,0,1\}$, with caps $7,12$ and returned arc $6$.
The source $[3,8]$ is $020121$ and charges $\{4,5,6,7\}$.
Together with the central charge $\{2\}$ these sets are disjoint,
leaving physical arc $3$ unpaid.

Each source contains at most seven complete arcs, fewer than the period
of ten.  Lemma~\ref{as:placement} supplies a genuine interval for each
one.  The local theorem's internal contact and return arguments use one
selected normal lift on that interval.  Only their scalar costs are added
between intervals.  Giving all three rows weight one satisfies every
physical capacity constraint.  In either case their exact sum is
\begin{align*}
 C&>E_*+\frac{41}{100}+\frac{461}{500}\\
  &=\frac{15821673034959}{8759490110000}
   =C_*+\frac{122824997277}{87594901100000}>C_*.
\end{align*}
This contradicts \eqref{fd:budget}.
\end{proof}
The supports in this example overlap, including at arcs that are charged
by another row.  This is harmless: capacity counts coefficients in the
charges, not the number of supporting windows through an arc.

\subsection{What a finite certificate proves}
A finite certificate consists of the qualified applications above, their
physical rows \eqref{as:physical-row}, exhaustive branch predicates and
nonnegative rational weights satisfying \eqref{as:capacity-constraint}
with the required lower sum. Discovery may use optimization; verification
uses only exact arithmetic after source qualification. A complete reduction
must also cover every excluded word. The differential theorems establish
geometric eligibility; the finite check establishes source maps, branches,
weights and capacities. Neither approximate optima nor reported counts
replace these obligations. A feasible surviving scalar vector is not a
trajectory without state and frame evolution, original endpoints, real
winding and a full maximizing lift.
\space
\section{From local costs to twenty-one types}\label{run:structure}
\label{analytic:run-structure}
The local costs now bound both the number and the lengths of the
increment runs. We use the assignment of own arcs and ports in
Section~\ref{run:profiles}. For a mixed cycle write
\[
 (a_1,b_1),\ldots,(a_M,b_M),\qquad
 a_i\ge1,\quad b_i\ge2,\quad b_i\ne3,
\]
and put
\[
 B=\sum_i(a_i+b_i),\qquad L=\#\{i:b_i\ge4\}.
\]
Thus $L$ counts the negative runs longer than two. The original
twist~\eqref{fd:increment-twist} gives
\begin{equation}\label{analytic:twist}
 \sum_i(a_i-b_i)\equiv1\pmod3.
\end{equation}
Throughout, $C$ is the total physical cost and
\[
 P=\frac{180483}{500000},\qquad C_* =5P=\frac{180483}{100000},
 \qquad C<C_*.
\]

\begin{theorem}[All-length run structure]
\label{analytic:structural-theorem}
For an original global minimizer with a mixed increment cycle and
$B\ge8$,
\[
 M\le3,\qquad a_i\le8,\qquad b_i\in\{2,5,7,8,9\}.
\]
At most one negative run has length greater than two.
\end{theorem}

We first bound $M$ by a uniform half-port payment, then bound the
negative lengths and their multiplicity. The few remaining lengths
four and six are excluded by explicit compatible charges. The last
step bounds the positive lengths. A source is always retained with
all its complete supporting arcs. Separately proved scalar costs
are added by Theorem~\ref{as:capacity}; their selected full lifts
need not be identified.

\subsection{A uniform payment for every positive run}
\begin{lemma}\label{analytic:three-positive-payments}
Suppose $B\ge10$. If $d\in\{0,1,2\}$ counts the incident negative
runs of length two, every positive run has a strict charged profile
of value at least $p_d$, where
\[
 p_0=\frac{389}{2000},\qquad
 p_1=\frac{473}{2000},\qquad
 p_2=\frac{119}{250}.
\]
The load is at most one on each own arc and at most one half on
each formal port.
\end{lemma}
\begin{proof}
For a positive run of length one, half the complete $102021$ or
$1020212$ four-arc profile gives $p_0$ or $p_1$. If $d=2$, the
complete eight-source $01020212$ pays more than $74/125$ on its
four-arc core at local positions $2,\ldots,5$ and more than $9/25$
on its two own arcs at positions $3,4$. Giving
each inequality weight one half yields
\[
 \frac12\left(\frac{74}{125}+\frac9{25}\right)
 =\frac{119}{250}.
\]
Each own arc has load one and each port load one half. The
four-arc and two-arc bounds are Theorem~\ref{lc:two-return} and
Lemma~\ref{run:own-singleton}, respectively.

For a run of length two, the proper five-source $02010$ pays
$49/200$ on its three own arcs when $d\le1$. When $d=2$,
Corollary~\ref{run:two-positive} pays $1/2$ on the same own set.
Both choices have zero port loads and dominate the stated values.

For $a\ge3$, Lemma~\ref{comb:half-strips} gives
\begin{equation}\label{analytic:positive-half-table}
\begin{array}{c|c}
a&\text{strict half-port payment}\\\hline
3&51/100+7d/200\\
4&431/1000+7d/200\\
a\ge5&61/200+P\lfloor(a-3)/2\rfloor+7d/200.
\end{array}
\end{equation}
For every $d$ its minimum occurs at $a=4$ and dominates $p_d$.
Indeed the $a\ge5$ value increases by $P$ every two increments,
and its first two values already exceed the $a=4$ value.
All supporting sources have at most nine arcs, so they are proper
under $B\ge10$.
\end{proof}

\begin{lemma}\label{analytic:M-three}
Every mixed minimizing increment cycle has $M\le3$.
\end{lemma}
\begin{proof}
If $M\ge4$, then $B\ge3M\ge12$, so the preceding profiles are
proper. They satisfy
\[
 p_d\ge\frac{119}{250}-\frac{479}{2000}(2-d).
\]
There is equality at $d=1,2$; at $d=0$ the right side is
$-3/1000$. The negative half-port profile of
Lemma~\ref{comb:negative-strips} has value zero at $b=2$ and at
least $3/5$ at every allowed $b\ge4$. Each long negative run
has two incident slots, so $\sum_i(2-d_i)=2L$. Adding these
profiles gives
\begin{equation}\label{analytic:M-margin}
 C>\frac{119}{250}M
       +\left(\frac35-\frac{479}{1000}\right)L
   =\frac{119}{250}M+\frac{121}{1000}L.
\end{equation}
At a length-two negative run the positive half-port loads add to
at most one on the shared physical arc. At a longer negative run,
each positive half-load meets at most one negative half-load.
All other charges lie in disjoint own sets or negative interiors.
Thus the physical capacity is at most one, also at repeated
incidences. The positive profiles make the sum strict. For $M\ge4$
the right side is at least $238/125>C_*$, a contradiction.
\end{proof}

\subsection{The negative lengths and their multiplicity}
\begin{lemma}\label{analytic:b-nine}
Every negative run of a mixed minimizing cycle has length at most nine.
\end{lemma}
\begin{proof}
A negative length $b\ge10$ implies $B\ge11$, so all half-port
sources are proper. Let $V(b)$ be the negative payment of
Lemma~\ref{comb:negative-strips}. For $b\ge11$, $V(b)\ge9/5$:
this holds at $b=11,12,13$ and $V(b+3)-V(b)=3/5$ thereafter.
One positive profile therefore gives
$C>9/5+389/2000=3989/2000>C_*$.

At $b=10$, $V(b)=3/2$. If $M\ge2$, two positive profiles give
$C>3/2+2(389/2000)=1889/1000>C_*$. For $M=1$, the twist gives
$a\equiv2\pmod3$, so $a=2$ or $a\ge5$. Number the negative
run's arcs $0,\ldots,10$. Its nine interior arcs split into
\[
 [1,3],\quad[4,6],\quad[7,9],
\]
with complete five-sources $[0,4],[3,7],[6,10]$. These proper
lawful images of $10210$ give the weak total $9/5$, with unit
load even at the two ports. At $a=2$ add the zero-port $49/200$
own payment from $02010$. At $a\ge5$ a complete $012012$ source
inside the positive run gives a disjoint own pair above
$P>49/200$. Hence in either case
\[
 C>\frac95+\frac{49}{200}=\frac{409}{200}>C_*.
\]
Here the twist supplies exactly the zero-port positive payment
needed after both negative ports have been used.
\end{proof}

\begin{lemma}\label{analytic:L-one}
A mixed minimizing cycle has at most one negative run of length
greater than two.
\end{lemma}
\begin{proof}
Suppose $L\ge2$. Then $B\ge10$, so all half-port profiles are
proper, and $M\le3$. If $M=3,L=2$, the positive incidences have
$d=0,1,1$. Thus
\[
 C>2\frac35+\frac{389}{2000}+2\frac{473}{2000}
   =\frac{747}{400}>C_*.
\]
If $M=L=3$, the lower bound is larger:
$C>9/5+3(389/2000)>C_*$.

It remains $M=L=2$. Every allowed negative length other than five
has $V(b)\ge9/10$; this follows from the exceptional values at
four and six and the three increasing residue classes. If either
negative length differs from five, then
$C>9/10+3/5+2(389/2000)>C_*$. If both are five and a positive
length is at least three, \eqref{analytic:positive-half-table}
gives
\[
 C>\frac65+\frac{431}{1000}+\frac{389}{2000}
   =\frac{3651}{2000}>C_*.
\]
Both positive lengths must therefore belong to $\{1,2\}$. The
twist $a_1+a_2-10\equiv1\pmod3$ leaves only $(a_1,a_2)=(1,1)$.

In this configuration $-^5+-^5+$, each negative run has complete
eight-source $20210212$. Theorem~\ref{opt:pair650} pays more than
$13/20$ on its two middle interior arcs, avoiding both ports.
Each positive singleton uses its full $102021$ profile above
$389/1000$, charging its own pair and its two ports. These four
charged sets are disjoint: the two ports of a length-five negative
run are distinct. Every complete source is proper at $B=12$, and
\[
 C>2\frac{13}{20}+2\frac{389}{1000}
   =\frac{1039}{500}>C_*.
\]
\end{proof}

\begin{lemma}\label{analytic:no-four-six}
A mixed minimizing cycle with $B\ge8$ has no negative run of
length four or six.
\end{lemma}
\begin{proof}
There is at most one such run, and all other negative lengths are
two. Write its length as $b=4$ or $6$.

If $M=3$, the positive incidences have $d=1,1,2$ and pay at
least $2p_1+p_2=949/1000$. Adding the negative payment $299/250$
or $13/10$ exceeds $C_*$. Here $B\ge11$, so all sources are
proper.

Suppose $M=2$. At $b=4$ the twist gives
$a_1+a_2\equiv1\pmod3$, hence $B\ge10$; at $b=6$ it gives
$a_1+a_2\equiv0\pmod3$, hence $B\ge11$. If a positive length
is at least three, its half-port value with $d=1$ is at least
$233/500$, the minimum in \eqref{analytic:positive-half-table}.
The other positive run pays at least $473/2000$. Adding even the
smaller negative payment gives
$C>299/250+233/500+473/2000>C_*$. Both positive lengths are
therefore at most two. The twist leaves $(2,2)$ at $b=4$, and
$(1,2)$ or $(2,1)$ at $b=6$.

For $b=6$, the two separate negative singleton payments total
$13/10$ and avoid the ports. Use the singleton positive run's
full $1020212$ profile above $473/1000$, and the other positive
run's zero-port $49/200$ own profile. Only the singleton uses
the shared short-separator port. These charges give
\[
 C>\frac{13}{10}+\frac{473}{1000}+\frac{49}{200}
   =\frac{1009}{500}>C_*.
\]

For $(b,a_1,a_2)=(4,2,2)$, start the word at
$\mathtt{----++--++}$ and set $j_0=0$. Reconstruct all labels
from the integer walk in \eqref{fd:increment-twist}, including
negative indices. The following whole-source maps give three
proper placements:
\begin{equation}\label{analytic:four-short-table}
\begin{array}{c|c|c|c}
\text{ordered support}&\text{map}&\text{source}&
 \text{physical charge}\\\hline
(2,1,\ldots,-5)&2-j&10201020&\{9,0,1\}\\
(5,4,\ldots,-1)&2-j&2021020&\{2\}\\
(2,3,\ldots,9)&j&10201020&\{3,4,5\}.
\end{array}
\end{equation}
On each eight-source, Lemma~\ref{ua:seven-sign} uses its suffix
$0201020$ for the negative entrance of the reflected $12021$ and
gives more than $61/200$ on local positions $1:3$. These are the
two triples displayed above. The seven-source retains the two contact jets and
terminal $020$ required for its $299/250$ center payment.
Their physical charges are disjoint, so
\[
 C>2\frac{61}{200}+\frac{299}{250}
   =\frac{903}{500}>C_*.
\]

Finally let $M=1$. At $b=4$ the twist gives $a\equiv2\pmod3$.
The value $a=2$ has $B=6$. Every $a\ge5$ has positive half-port
payment at least $61/200+P$, which together with $299/250$
exceeds $C_*$. At $a=5,B=9$ this construction still uses only
proper six- and seven-sources. At $b=6$ the twist gives
$a\equiv1\pmod3$. The value $a=1$ has $B=7$. For $a\ge7$,
the positive half-port value is at least $61/200+2P$ and its sum
with $13/10$ exceeds $C_*$. At $a=4,B=10$, use the two full
ordinary terminal triples above $61/200$ each, the middle own
arc above $63/500$, and the two negative singleton payments
totaling $13/10$. The negative charges avoid both ports, so
\[
 C>2\frac{61}{200}+\frac{63}{500}+\frac{13}{10}
   =\frac{509}{250}>C_*.
\]
\end{proof}

\subsection{Positive lengths and pure cycles}
\begin{lemma}\label{analytic:a-eight}
Every positive run of a mixed minimizing cycle has length at most eight.
\end{lemma}
\begin{proof}
Suppose $a\ge9$. Then $B\ge11$, so the sources of at most nine
arcs used here are proper. If a negative run is longer than two,
its half-port profile pays at least $3/5$. The positive run's
half-port value is at least $61/200+3P$, so
$C>61/200+3P+3/5>C_*$. If all negative lengths are two and
$M\ge2$, this positive run pays at least $3/8+3P$ and another
positive run pays at least $119/250$. Their compatible half-port
sum again exceeds $C_*$.

It remains $M=1,b=2$. The two formal ports now coincide, so
we instead use the complete nine-source $012021201$ across the
short separator. Its incremental pattern $\mathtt{+++--+++}$
is supplied by the end and beginning of the same positive run.
Lemma~\ref{six:lem:H3} pays the weak value $2P$ on the last two
own arcs, the shared port, and the first two own arcs. Number
the own arcs $0,\ldots,a$. The full support is
$[a-3,a+5]$, under the whole phase $j\mapsto j-(a-3)$;
its physical charge is $\{a-1,a,a+1,0,1\}$ modulo $B=a+2$.
These five arcs are distinct. The three pairs
\[
 \{2,3\},\quad\{4,5\},\quad\{6,7\}
\]
have complete forward sources $[0,5],[2,7],[4,9]$ inside the
positive run and each cost more than $P$. They avoid the first
and last own pairs because $a\ge9$. Thus $C>2P+3P=C_*$,
where strictness comes from the three forward pairs.
\end{proof}

\begin{proof}[Proof of Theorem~\ref{analytic:structural-theorem}]
Lemmas~\ref{analytic:M-three}, \ref{analytic:b-nine} and
\ref{analytic:L-one} bound the number of positive runs, the
negative lengths and their multiplicity. Lemma~\ref{analytic:no-four-six}
removes the negative lengths four and six for $B\ge8$.
Lemma~\ref{analytic:a-eight} gives the positive-length bound.
\end{proof}

\begin{proposition}[Pure increment cycles]\label{analytic:pure}
Among pure increment cycles with $B\ge8$, only the all-minus
cycle of length eight can belong to an original global minimizer.
\end{proposition}
\begin{proof}
For a pure positive cycle, the twist gives $B\equiv1\pmod3$,
so $B\ge8$ implies $B\ge10$. Use the proper complete $012012$
source at every physical start with coefficient one half. Each
arc belongs to two charged pairs, so its physical load is one.
The strict forward-pair bound gives $C>BP/2\ge5P=C_*$.

For a pure negative cycle, $B\equiv2\pmod3$. Use every proper
complete five-source with coefficient one third. Each arc belongs
to three charged triples, giving load one and the weak bound
$C\ge B/5$. Apart from $B=8$, the first possible length is
eleven, which already exceeds the budget. In both averages a
source crossing the seam uses the original twisted labels and
its lawful phase. Only physical scalar costs are summed.
\end{proof}
\space
\subsection{Three positive runs}\label{run:three-structure}
We now determine the mixed cycles with three positive runs.
Theorem~\ref{analytic:structural-theorem} leaves either three negative
runs of length two or one longer negative run and two of length two.
All charges below are original physical costs. We retain the complete
proper sources and their whole lawful charts from
Section~\ref{run:catalogue}; each source keeps the full maximizing
data required by its proof. Independently established scalar payments
may be added without identifying their selected lifts.

\begin{proposition}\label{analytic:M3-short}
Suppose a mixed minimizing cycle has three positive runs and all its
negative runs have length two. Its positive lengths must be
$(1,1,2)$ in cyclic order.
\end{proposition}
\begin{proof}
Put $A=a_1+a_2+a_3$. The twist gives $A\equiv1\pmod3$ and hence
$B=A+6\ge10$. A useful charge across a short separator is
\[
 \mathtt{012021201},\qquad \text{charge }2:6,
 \qquad \text{payment }\ge2P.
\]
Its increment factor is $\mathtt{+++--+++}$. Between positive runs
of length at least three, it charges the last two own arcs on the
left, the shared port, and the first two own arcs on the right.
The complete nine-source retains both returned contexts and their
internal forward witnesses, as required by Lemma~\ref{six:lem:H3}.
The two terminal own pairs within any such positive run are disjoint.

First suppose none of the positive lengths is one. If all three are
at least three, use this nine-source at every separator. Their
physical charges are disjoint, so $C\ge6P>C_*$.
If exactly one length is two, pay $1/2$ on its three own arcs by
Corollary~\ref{run:two-positive}. Between the other two runs use
the nine-source, and at their two other ends use the returned
terminal triples. The charges are disjoint and give
\[
 C>\frac12+2P+2\frac{379}{1000}>C_*.
\]
Otherwise write the lengths as $(2,2,a)$. The twist forces
$a\equiv0\pmod3$, so $a=3$ or $a\ge6$.
For $a=3$, let $s$ be the start of that run and set $j_s=0$.
The complete source $[s-3,s+6]=1210120212$, under the whole phase
$j\mapsto j+1$, is $2021201020$. Its charge $2:7$ is
$[s-1,s+4]$, the four own arcs and both ports of the length-three
run. Both actual exterior returns are retained. The strict
$807/1000$ bound of Proposition~\ref{late:r24:W807}, together with
the two zero-port half-unit payments, gives
\[
 C>\frac12+\frac12+\frac{807}{1000}
   =\frac{1807}{1000}>C_*.
\]
For $a\ge6$, the two returned terminal triples and one disjoint
internal forward pair instead give
$C>1+2(379/1000)+P>C_*$. Here $B\ge13$, so even the ten-source
used at $a=3$ is proper.

Now suppose exactly one positive length is one. If its two neighbors
are at least three, Lemma~\ref{run:own-singleton} pays $3/7$ on its
own pair, with zero port loads. The two proper six-source witnesses
give the required physical anchor-density facts; they do not identify
covectors across their containing intervals. Use the nine-source
between the two neighbors and one returned triple at each of their
other ends. This yields
\[
 C>\frac37+2P+2\frac{379}{1000}>C_*.
\]
Otherwise the two possible cyclic orders are $(1,2,a)$ and $(1,a,2)$.
The twist gives $a\equiv1\pmod3$, hence $a\ge4$.
The own payments $3/7$ and $1/2$ leave every port available to the
longer run. Its two returned triples and either its middle-center
payment ($a=4$) or an internal forward pair ($a\ge5$) give, in
either order,
\[
 C>\frac37+\frac12+2\frac{379}{1000}+\frac{63}{500}
   =\frac{1586}{875}>C_*.
\]
These are the full source placements of
Proposition~\ref{run:three-pattern}; the two whole lawful placements
do not identify the two cyclic types.

Finally, with at least two singleton runs write the lengths as
$(1,1,a)$. The twist gives $a\equiv2\pmod3$. If $a\ge5$, each
singleton pays $9/25$ on its own pair, while the other run pays its
two returned triples and one internal pair. Thus
\[
 C>2\frac9{25}+2\frac{379}{1000}+P
   =\frac{919483}{500000}>C_*.
\]
The only value left is $a=2$. The cases above exhaust the numbers
of positive lengths equal to one and two. Every nine-source used
is proper, and the sole ten-source was used only at $B\ge13$.
\end{proof}

\begin{proposition}\label{analytic:M3-long}
A mixed minimizing cycle with three positive runs has no negative
run longer than two.
\end{proposition}
\begin{proof}
There is at most one long negative run. If its length is seven,
eight or nine, its half-port value is at least $9/10$. The three
positive incidences have $d=1,1,2$, so
Lemma~\ref{analytic:three-positive-payments} gives
\[
 C>\frac9{10}+2\frac{473}{2000}+\frac{119}{250}
   =\frac{1849}{1000}>C_*.
\]
It remains to consider negative lengths $(5,2,2)$. Write the
positive lengths, starting after the length-five run, as $(x,y,z)$.
The end runs have $d=1$ and the middle run has $d=2$.
The complete eight-source $20210212$ pays more than $13/20$ on
the two middle interior arcs of the length-five run, with zero
port loads, by Theorem~\ref{opt:pair650}.

If $x$ or $z$ is at least three, its half-port payment is at least
$233/500$, the minimum of Lemma~\ref{comb:half-strips} at $d=1$.
The other two positive profiles therefore give
\[
 C>\frac{13}{20}+\frac{233}{500}
       +\frac{473}{2000}+\frac{119}{250}
   =\frac{3657}{2000}>C_*.
\]
Thus $x,z\le2$. If $y\ge5$, its half-port payment is at least
$3/8+P$, giving
$C>13/20+2(473/2000)+3/8+P>C_*$. Consequently $y\le4$.
The twist is $x+y+z\equiv1\pmod3$.
For $x+z=2,3,4$, respectively, this leaves
$y=2$, $y\in\{1,4\}$, and $y=3$. The six possible triples are
therefore
\[
 (1,2,1),\quad(1,1,2),\quad(2,1,1),\quad
 (1,4,2),\quad(2,4,1),\quad(2,3,2).
\]
Their lengths are at least thirteen, so every source used below
is proper.

At $(1,2,1)$ give each end singleton its full one-return four-arc
profile, above $473/1000$ by Corollary~\ref{loc:E473}.
The middle length-two run pays $1/2$ on its own arcs and the negative
run pays its zero-port pair. The two long-separator ports are
distinct; each short-separator port is charged by only one singleton.
Hence
\[
 C>2\frac{473}{1000}+\frac12+\frac{13}{20}
   =\frac{262}{125}>C_*.
\]
At $(2,3,2)$ use the middle run's complete $W807$ profile, including
both actual exterior returns. The end runs use their zero-port
$49/200$ own payments from the complete five-source $02010$.
All these charges avoid the negative pair, so
\[
 C>\frac{807}{1000}+2\frac{49}{200}+\frac{13}{20}
   =\frac{1947}{1000}>C_*.
\]

At $(1,4,2)$ give coefficient one half to the singleton's full
$473/1000$ profile and to the adjacent returned terminal triple of
the middle run. Their common port then has load one.
Give coefficient one to the middle run's other returned triple
and its middle $63/500$ payment, to the length-two run's zero-port
$49/200$ own payment, and to the negative pair. The middle own
charges are disjoint, and every other physical load is at most one.
Thus
\[
 C>\frac12\frac{473}{1000}+\frac32\frac{379}{1000}
       +\frac{63}{500}+\frac{49}{200}+\frac{13}{20}
   =\frac{913}{500}>C_*.
\]
For $(2,4,1)$ interchange left and right by whole lawful reflection
and time reversal, retaining the complete sources. This gives a
separate compatible allocation for that cyclic type.

For the last two triples we record every source and charge.
Take $j_0=0$ and reconstruct the native labels from the printed
increment word, with $j_{i+13}=j_i+1$ modulo three. The maps in
the tables act on the entire ordered source and its full lift;
only physical costs are reduced modulo thirteen. The final column
is a strict lower payment.

For $(x,y,z)=(1,1,2)$ the word is
$\mathtt{-----+--+--++}$, and the placements are
\begin{equation}\label{run:three-B13-forward}
\begin{array}{c|c|c|c|c}
\text{support}&\text{map}&\text{source}&\text{charge}&\text{payment}\\\hline
(2,1,\ldots,-5)&2-j&10201020&\{12,0,1\}&61/200\\
(6,5,\ldots,-1)&1-j&20210212&\{2,3\}&13/20\\[0pt]
[3,9]&1+j&1020212&\{4,5,6,7\}&473/1000\\[0pt]
[5,12]&2+j&01020212&\{8,9\}&9/25\\[0pt]
[5,14]&2+j&0102021202&\{11\}&3/40.
\end{array}
\end{equation}
For $(x,y,z)=(2,1,1)$ the distinct cyclic word is
$\mathtt{-----++--+--+}$. Its placements are
\begin{equation}\label{run:three-B13-reverse}
\begin{array}{c|c|c|c|c}
\text{support}&\text{map}&\text{source}&\text{charge}&\text{payment}\\\hline
[3,10]&1+j&10201020&\{4,5,6\}&61/200\\
(6,5,\ldots,-1)&1-j&20210212&\{2,3\}&13/20\\
(2,1,\ldots,-4)&2-j&1020212&\{1,0,12,11\}&473/1000\\
(0,-1,\ldots,-7)&-j&01020212&\{10,9\}&9/25\\
(0,-1,\ldots,-9)&-j&0102021202&\{7\}&3/40.
\end{array}
\end{equation}
In both tables, Lemma~\ref{ua:seven-sign} uses the complete suffix
$0201020$ of the eight-source to establish the entrance sign of the
reflected $12021$. Its charge on local positions $1:3$ gives the
displayed $61/200$ payment. The ten-source keeps its
returned four-core, internal seven-arc sign and circle witness,
and uncharged final cap, as required by Lemma~\ref{late:r18:cap}.
The $9/25$ row charges only its own pair; its two anchors remain
uncharged by that row. Each table has five pairwise disjoint
physical charged sets. Therefore both give
\[
 C>\frac{61}{200}+\frac{13}{20}+\frac{473}{1000}
       +\frac9{25}+\frac3{40}
   =\frac{1863}{1000}>C_*.
\]
All six triples have been excluded.
\end{proof}

\begin{corollary}\label{analytic:M3-only}\label{run:three-only}
A mixed minimizing cycle with three positive runs has $B=10$,
negative lengths $(2,2,2)$ and positive lengths $(1,1,2)$ in cyclic
order. Its cyclic increment word is $\mathtt{--+--+--++}$.
\end{corollary}
\begin{proof}
Combine Propositions~\ref{analytic:M3-short} and~\ref{analytic:M3-long}.
\end{proof}
\space
\subsection{One and two positive runs}\label{run:census}
Theorem~\ref{analytic:structural-theorem} and
Lemma~\ref{analytic:a-eight} reduce the mixed case to at most three
positive runs, with
\[
 1\le a_i\le8,\qquad b_i\in\{2,5,7,8,9\},
 \qquad \#\{i:b_i>2\}\le1.
\]
The three-run case was settled in Corollary~\ref{analytic:M3-only}.
We now determine the one-run and two-run alternatives directly from
the twist and compatible local payments. No upper bound on the period
is assumed in this argument.

\subsubsection{The payments used at a separator}
We retain the physical own arcs and ports of Section~\ref{run:profiles}.
Put
\[
 P=\frac{180483}{500000},\qquad t_0=\frac{379}{1000},
 \qquad C_*=5P.
\]
The following placements explain the combinations used below. Their
complete-source qualifications are those of Section~\ref{run:catalogue}.
Number the own arcs of a positive run of length $a$ as $0,\ldots,a$.
For $a\ge3$, the ordinary terminal payments charge
\[
 \{-1,0,1\},\qquad \{a-1,a,a+1\},
\]
respectively, above $61/200$. Their complete six-sources are
$(3,2,\ldots,-2)$, with whole label map $j\mapsto-j$, and
$[a-3,a+2]$, with whole phase $3-a$. Set the control label at own
arc zero to zero; both sources then become $012021$.
Corollary~\ref{ek:short-sources} gives the $61/200$ bound on these
triples. If the facing negative run has length two, the actual return
extends the corresponding source to $0120212$, and the same corollary
improves its charge to more than $t_0$. These are alternative payments on the
same arcs, rather than additional charges.

When $B>9$, a negative run of length two whose neighboring positive
runs both have length at least three admits the complete nine-source
$012021201$, which
pays weakly $2P$. It charges the terminal own pair of the first run,
the shared separator port, and the initial own pair of the second.
Indeed, if the first run ends at own arc $t$, the source is
$[t-3,t+5]$, phased so that the label at $t$ is zero, and its charge
is $[t-1,t+3]$. Within a positive run of length $a\ge5$, the pair
$\{2,3\}$ has the complete source $[0,5]=012012$ and pays more than
$P$. More generally the pairs
\[
 \{2+2r,3+2r\},\qquad
 0\le r<\lfloor(a-3)/2\rfloor,
\]
are disjoint from both terminal own pairs and from one another.

A positive run of length two pays more than $49/200$ on its three
own arcs, with zero port loads: its complete five-source
$[-1,3]$, phased by $+2$, is $02010$.
A negative run of length five pays more than $13/20$ on its two
middle interior arcs, also with zero port loads. With its starting
arc numbered zero, its complete eight-source is
$[-1,6]=20210212$, after setting the initial control label to zero.
These two consequences are Theorems~\ref{f245:theorem}
and~\ref{opt:pair650}.

We also use the half-port profiles of
Lemmas~\ref{comb:half-strips} and~\ref{comb:negative-strips}.
If a positive run has exactly one adjacent negative run of length two,
its payment is at least
\[
 \frac{473}{2000};
 \qquad \frac{233}{500}\quad(a\ge3);
 \qquad \frac{17}{50}+P\lfloor(a-3)/2\rfloor\quad(a\ge5).
\]
At $a=3$ the value is at least $109/200$.
The first bound for a run of length two follows from its stronger
zero-port $49/200$ payment. The negative half-port values needed here
are
\[
 V(5)=\frac35,\qquad V(7)=\frac9{10},\qquad
 V(8)=V(9)=\frac65.
\]
At each separator these profiles have total physical load at most one.
When full terminal payments are used instead, their charged sets will
be kept disjoint explicitly. All exclusions below have $B\ge12$;
all their complete sources have at most nine arcs and are therefore
proper. In particular, these payments impose no new premise on the
short candidates of periods eight and nine.

\begin{proposition}[One positive run]\label{run:one-candidates}
If an original global minimizer has $B\ge8$ and exactly one positive
run, its run pair is one of
\[
 (6,2),\quad(3,5),\quad(6,5),\quad(2,7),\quad(3,8),\quad(1,9).
\]
\end{proposition}
\begin{proof}
Write the pair as $(a,b)$. The bound $a\le8$ and the twist
$a-b\equiv1\pmod3$ leave the following possibilities:
\[
\begin{array}{c|c}
 b&a\\\hline
 2&6\\
 5&3,6\\
 7&2,5,8\\
 8&3,6\\
 9&1,4,7.
\end{array}
\]
Here the condition $a+b\ge8$ has already removed the shorter pairs.
The half-port profiles exclude $(8,7),(6,8),(7,9)$, respectively, by
\[
 \frac{61}{200}+2P+\frac9{10}>C_*,\qquad
 \frac{61}{200}+P+\frac65>C_*,\qquad
 \frac{61}{200}+2P+\frac65>C_*.
\]

For $(a,b)=(5,7)$, the complete eight-source $10120121$ charges the
six positive own arcs above $609/1000$ by Theorem~\ref{loc:OwnSix609}.
Its support is $[-1,6]$, with own arc zero labeled zero.
The six negative interior arcs split into the disjoint triples
$[1,3]$ and $[4,6]$, each retaining its complete weak five-source.
They pay $6/5$ and are disjoint from the positive own arcs. Thus
\[
 C>\frac{609}{1000}+\frac65=\frac{1809}{1000}>C_*.
\]

For $(a,b)=(4,9)$, the positive run uses its two ordinary terminal
triples and its middle opposite-center payment $63/500$.
The latter charges own arc $2$, on the complete source
$[-2,6]$, sent to $102012021$ by whole phase $+2$, of
Theorem~\ref{pf:opposite126}.
Number the negative interior arcs $1,\ldots,8$.
Two complete weak five-sources charge $[2,4]$ and $[5,7]$, leaving
both ports $1,8$ free. All five charged sets are disjoint, so
\[
 C>2\frac{61}{200}+\frac{63}{500}+\frac65
   =\frac{242}{125}>C_*.
\]
These exclusions leave precisely the six stated necessary alternatives.
\end{proof}

\begin{proposition}[Two positive runs]\label{run:two-candidates}
If an original global minimizer has $B\ge8$ and exactly two positive
runs, there are thirteen necessary run configurations.
When both negative lengths are two, the unordered positive pair is
\[
 \{a_1,a_2\}=\{1,4\},\ \{2,3\},\ \{2,6\},\ \{4,4\}.
\]
Otherwise write the cyclic increments as $+^x-^b+^y-^2$, where
$b>2$ is the distinguished negative length. The ordered pair $(x,y)$
is one of
\[
\begin{array}{c|l}
 b&(x,y)\\\hline
 5&(1,1),(1,4),(4,1),(2,3),(3,2)\\
 7&(1,3),(2,2),(3,1)\\
 8&(1,1).
\end{array}
\]
In particular $b=9$ is impossible.
\end{proposition}
\begin{proof}
Suppose first that both negative lengths are two.
If both positive lengths are at least three and one is at least five,
use the fused nine-source at each separator and an internal forward
pair in the longer positive run. The two fusions use opposite
terminal own pairs in each run, and the internal pair avoids them.
Their disjoint charges give $C>2P+2P+P=C_*$.
Thus if both lengths are at least three, they lie in $\{3,4\}$;
the twist $a_1+a_2\equiv2\pmod3$ leaves only $(4,4)$.

If one positive length is one, the other is $1,4$ or $7$.
The first choice has $B=6$. For the choice seven,
Lemma~\ref{run:own-singleton} pays more than $3/7$ on the singleton's
own arcs. The longer run uses its two returned terminal triples
and the internal pairs $\{2,3\},\{4,5\}$, giving
\[
 C>\frac37+2t_0+2P>C_*.
\]
The two neighboring positive occurrences in the singleton lemma are
here the same physical run; only its two distinct density owners
are used. The longer run is charged once.
Finally, if one length is two and the other is at least two, the
twist and $a_i\le8$ leave the other length three or six.
This proves the four short-negative alternatives.

Now let $b>2$ be the unique long negative length.
Each positive run has exactly one short-negative incidence.
For $b=5$, if both positive lengths are at least three, use the
zero-port negative middle pair, the fused nine-source at the short
separator, and an ordinary terminal triple at each end facing the
long separator. These four charges are disjoint and give
\[
 C>\frac{13}{20}+2P+2\frac{61}{200}>C_*.
\]
One positive length is therefore at most two.
If it is one, the twist $x+y\equiv2\pmod3$ leaves the other length
$1,4$ or $7$. At seven the compatible half-port profiles give
\[
 C>\frac35+\frac{473}{2000}+\frac{17}{50}+2P>C_*.
\]
If the shorter length is two and neither length is one, the other
is three or six. At six, combine the shorter run's zero-port
$49/200$ payment, both terminal triples of the longer run, one of
its internal forward pairs, and the zero-port negative middle pair.
The ordinary triple faces the long separator, while the returned
triple faces the short one. Their ports and own pairs are disjoint
from the other charges, so
\[
 C>\frac{49}{200}+\frac{61}{200}+t_0+P+\frac{13}{20}>C_*.
\]
The five displayed ordered pairs for $b=5$ remain.

For $b=7$, any positive length at least five gives the compatible
half-port bound
\[
 C>\frac{17}{50}+P+\frac{473}{2000}+\frac9{10}>C_*.
\]
Thus $x,y\le4$. The twist $x+y\equiv1\pmod3$ leaves totals four
and seven. Total seven consists of lengths three and four, whose
half-port payments give
\[
 C>\frac{109}{200}+\frac{233}{500}+\frac9{10}>C_*.
\]
Total four gives exactly $(1,3),(2,2),(3,1)$.

For $b=8$ or $9$, a positive length at least three would give
\[
 C>\frac{233}{500}+\frac{473}{2000}+\frac65>C_*.
\]
Hence $x,y\le2$. At $b=8$ the twist leaves only $(1,1)$.
At $b=9$ it leaves $(1,2),(2,1)$.
For either of the latter configurations, the negative run admits
weak five-source charges $[2,4]$ and $[5,7]$, of total $6/5$,
avoiding both ports. The singleton positive run uses its complete
one-return four-arc payment above $473/1000$, while the length-two
run pays more than $49/200$ on its own arcs. Their charged sets
are disjoint: the shared short-separator port is used only by the
singleton. Therefore
\[
 C>\frac65+\frac{473}{1000}+\frac{49}{200}
   =\frac{959}{500}>C_*.
\]
The singleton payment is Corollary~\ref{loc:E473}; its full seven-source
retains the actual short return.

All cases are now exhausted. When both negative lengths are two,
cyclic rotation exchanges the two positive lengths. When $b>2$,
its distinguished position fixes the ordered pair $(x,y)$; the proof
covers both orders and does not identify them by reflection.
\end{proof}

\subsection{The twenty-one necessary types}
\begin{theorem}[The run classification]\label{run:21}\label{thm:main}\label{en:finite-reduction}
Every original global minimizer with $B\ge8$ has one of the twenty-one
essential increment types in Table~\ref{run:table}: six types with one
positive run, thirteen with two, one with three, and the all-minus
eight-arc type. In particular every original global minimizer has
$B\le13$, and any prescribed fundamental cut has at most $N\le14$
positive pieces.
\end{theorem}
\begin{proof}
Theorem~\ref{analytic:structural-theorem} puts every mixed type into
$M=1,2,3$. Propositions~\ref{run:one-candidates}
and~\ref{run:two-candidates} give the first nineteen types, while
Corollary~\ref{analytic:M3-only} leaves the single cyclic pattern
$(a_1,a_2,a_3)=(1,1,2)$ with all $b_i=2$.
Proposition~\ref{analytic:pure} leaves only the all-minus word with
$B=8$ in the pure case. These configurations translate directly into
Table~\ref{run:table}; the maximal listed length is thirteen.
The types with $B\le7$ already satisfy this bound. A prescribed cut
splits at most one physical positive arc, so $N\le B+1\le14$.
\end{proof}

Table~\ref{run:table} records necessary alternatives. The proof uses
the run inequalities for arbitrary lengths. The separate bounded
enumeration in \href{https://arxiv.org/abs/2611.21470v1}{version~1}
can be checked by the \path{development} task of the reader archive.

\begin{longtable}{@{}c r l l@{}}
\caption{The twenty-one necessary increment types grouped by their negative runs. A pair $(a,b)$ records a positive run and the following negative run.}\label{run:table}\\
\toprule Identifier&$B$&Cyclic increments&Run pairs\\\midrule\endfirsthead
\toprule Identifier&$B$&Cyclic increments&Run pairs\\\midrule\endhead
\bottomrule\endfoot
\multicolumn{4}{l}{\emph{All negative runs of length two}}\\*[2pt]
R03&8&$\mathtt{--++++++}$&$(6,2)$\\
R05&9&$\mathtt{--+--++++}$&$(1,2),(4,2)$\\
R06&9&$\mathtt{--++--+++}$&$(2,2),(3,2)$\\
R07&10&$\mathtt{--+--+--++}$&$(1,2),(1,2),(2,2)$\\
R15&12&$\mathtt{--++--++++++}$&$(2,2),(6,2)$\\
R17&12&$\mathtt{--++++--++++}$&$(4,2),(4,2)$\\
\multicolumn{4}{l}{\emph{One negative run of length five}}\\*[2pt]
R02&8&$\mathtt{-----+++}$&$(3,5)$\\
R04&9&$\mathtt{-----+--+}$&$(1,5),(1,2)$\\
R08&11&$\mathtt{-----++++++}$&$(6,5)$\\
R10&12&$\mathtt{-----+--++++}$&$(4,5),(1,2)$\\
R11&12&$\mathtt{-----++--+++}$&$(3,5),(2,2)$\\
R12&12&$\mathtt{-----+++--++}$&$(2,5),(3,2)$\\
R13&12&$\mathtt{-----++++--+}$&$(1,5),(4,2)$\\
\multicolumn{4}{l}{\emph{One negative run of length seven, eight, or nine}}\\*[2pt]
F41&9&$\mathtt{-------++}$&$(2,7)$\\
F35&11&$\mathtt{--------+++}$&$(3,8)$\\
F40&10&$\mathtt{---------+}$&$(1,9)$\\
F24b&13&$\mathtt{-------+++--+}$&$(1,7),(3,2)$\\
F24a&13&$\mathtt{-------+--+++}$&$(3,7),(1,2)$\\
F39&13&$\mathtt{-------++--++}$&$(2,7),(2,2)$\\
F36&12&$\mathtt{--------+--+}$&$(1,8),(1,2)$\\
\multicolumn{4}{l}{\emph{The pure all-minus type}}\\*[2pt]
R01&8&$\mathtt{--------}$&---\\
\end{longtable}
\space
\space

\section{From twenty-one types to six}\label{run:twenty-one-to-six}
For a mixed word in Table~\ref{run:table}, either every negative run
has length two or exactly one has length greater than two. We group
the words by that longer length. Qualified local costs and endpoint
comparisons exclude three types in the first group, five in the
length-five group and all seven in the last mixed group.
\subsection{Negative runs all of length two}\label{run:short}
The six words in this group are R03, R05, R06, R07, R15 and R17.
Their positive-run lengths are respectively $(6)$, $(1,4)$, $(2,3)$,
$(1,1,2)$, $(2,6)$ and $(4,4)$, up to cyclic rotation. We exclude the
last three.

\subsubsection{R07 and its shared arc}
The two length-one positive runs require a comparison on a common
maximizing lift. The two genuine cuts in
Section~\ref{new:glue} identify the lifts on their overlap and
give the actual covector twist. In the native word $0212101020$,
write $s$ for the own reciprocal on the shared arc, $p,v$ for the
two exterior owners and $x=s^2$. The shared law and the positive cap
durations yield the following exhaustive comparison.
\begin{enumerate}
\item For $s\le1/5$, the two actual binary means give an upper bound
on both mean factors. In the full-frame formulas of
Section~\ref{new:r07}, these bounds force $F_LF_R>(1+x)^2/4$,
whereas positivity of the last arc gives $F_LF_R<(1+x)^2/4$.
\item For $x\ge1/6$, the unsquared shared law gives
$\Gamma>83376585/113379904>H_*$ for the adjoint determinant
$\Gamma$, contrary to the harmonic ceiling $H_*$ on the same lift.
\item For $1/25<x<1/6$ and $p,v\le s$, the two cap slopes force
one short side and one long side. In the order with the short side
first, the remaining coordinates
$(p/s,v/s,r_L/s,r_R/s,x)$ lie in
$[3/11,1]^2\times[3/2,2]^2\times[1/8,1/6]$.
Here $r_L,r_R$ are the reciprocals on either side of the shared
arc. A closed cover contradicts the necessary collar, duration,
full-frame or cap-slope conditions throughout that box.
\item In the remaining mixed and upper ranges, positive deficits
and actual returned-source density bounds enclose the candidate
in \eqref{new:r07root}. The shared determinant, both exterior
cap slopes and the full signed frame give the necessary system
excluded by Lemma~\ref{new:r07finite}.
\end{enumerate}
Section~\ref{new:r07} proves the inclusions and contradictions
in these four ranges, including the endpoints of their comparisons.

\subsubsection{R15 and its central triple}\label{new:r15}
\begin{proposition}\label{new:r15thm}
R15 cannot be the essential type of an original global minimizer.
\end{proposition}
\begin{proof}
Its native word is $021202120120$.
The complete proper six-source $[0,5]+1=102010$ retains the
central $[1,5]+1=02010$ and the capped binary normality witness
$[0,4]$. Select its full nonzero lift before restriction.
The forward source $[1,7]+1$ and the whole reflected source
$(5,4,\ldots,-1)$ with map $j\mapsto2-j$ are $0201020$.
They give both actual unsorted H150 crossings.
The two switches adjacent to the center have positive third
reciprocals, so $yA,yB<1$. In the notation of
\eqref{new:centralguards} these inequalities give $Eb,Ea<2$.
The returned sources $[9,15]+2$ and $(9,8,\ldots,3)$ with map
$j\mapsto1-j$ are $0120212$. They give the terminal circles
for the central cap reciprocals $p=r_1$ and $q=r_5$.
At the terminal state $(u,v)$ in the first return chart the two
reciprocals are $p=(1-u)/v$ and $A=(u^2+v^2-1)/(2v)$.
The circle gives $v^2>(1-u)(3+u)$, hence $p<A$.
Also $u>0$ and $p^2+2pA=1-u-u(1-u)^2/v^2<1$.
The terminal state in the reflected return chart gives $q<B$
and $q^2+2qB<1$ by the same formulas.
The six-sources $[10,15]$ and $[2,7]+1$ are $201020$ with
physical owners $1$ and $5$, respectively; hence $p^2,q^2<1/6$.
All these supports are proper, and their state conclusions give
precisely \eqref{new:centralguards}. Theorem~\ref{new:centralfee}
therefore gives
\begin{equation}\label{new:r15budget}
 F=c_2+c_3+c_4>C_0=80483/125000.
\end{equation}

For each complete $0120212$ source, Lemma~\ref{cmp:forward-prefix-forcing}
applied to its first six arcs supplies the nonpositive entrance and the
entrance circle. The final return supplies the exit circle for the
$120212$ suffix on the same full lift. Proposition~\ref{six:prop:tail}
therefore gives cost greater than $2/5$ on each returned triple,
$\{11,0,1\}$ and $\{5,6,7\}$. The forward source
$[6,11]+2=012012$ gives $c_8+c_9>P$.
The two triples, the central charge $F$ and the forward pair are
pairwise disjoint and cover all physical arcs except arc $10$,
whose cost is positive. Thus
\[
 C>C_0+2\frac25+P=\frac{180483}{100000}=C_*,
\]
contradicting the strict original budget.
\end{proof}

\subsubsection{R17 and the entrance terms at both centers}
The native word is $021201210120$. Its two complete returned unions
charge $U_1=c_{11}+c_0+c_1+c_2+c_3$ and
$U_7=c_5+c_6+c_7+c_8+c_9$. Together with the distinct center arcs
$4,10$ these partition the physical costs. The weak union bounds
$U_1,U_7\ge2P$ and the strict center bounds $c_i>63/500$ put each
center in the window \eqref{new:r17window} by $C<C_*$.

On the common lift of Section~\ref{new:glue}, the exclusive integrals
in each union retain their entrance terms. Moving these terms into
the corresponding center expressions $G_4,G_{10}$ gives
$C>J_1+J_7+G_4+G_{10}$ by \eqref{new:r17regroup}.
Here $J_1,J_7$ are the remaining lower payments of the two unions.
Equation~\eqref{new:r17outer} gives $J_1,J_7>14/25$ and
Lemma~\ref{new:r17center} gives $G_4,G_{10}>343/1000$ on the
actual cost windows. Thus
\[
 C>2\frac{14}{25}+2\frac{343}{1000}
   =\frac{903}{500}>C_* .
\]
Each entrance term is part of an original physical cost and is counted
once in this sum. The regrouping and the closed comparisons are in
Section~\ref{new:r17}.

\begin{theorem}[Negative runs of length two]\label{run:short-survivors}
An original global minimizer with $B\ge8$ and a mixed increment
cycle whose negative runs all have length two has one of the three
necessary types R03, R05 and R06.
\end{theorem}
\begin{proof}
Theorem~\ref{run:21} gives the six possibilities above. The three
comparisons exclude R07, R15 and R17, respectively.
\end{proof}
\space
\subsection{One negative run of length five}\label{run:five}
The one- and two-positive-run propositions leave seven words with a negative run of length
five. With one positive run they are $-^5+^3$ and $-^5+^6$, namely
R02 and R08. With two positive runs, start immediately before the
length-five negative run and write
$$-^5+^a-^2+^b.$$
The remaining pairs are
$$(a,b)=(1,1),(1,4),(2,3),(3,2),(4,1),$$
corresponding to R04 and R10--R13. These are exactly the seven
members of this family in Theorem~\ref{run:21}; the bounds on all run lengths have already been proved in
Section~\ref{run:structure}.

The length-five negative run supplies a rooted prefix and its adjacent
one-sided leaf. The decisive information is the actual endpoint sign
and its transfer to the neighboring capped core. The following three
arguments retain that common geometry and give the full conditional
cost comparisons for R08, the pair R10/R13, and the pair R11/R12.
Their reusable one-sided construction is proved in
Section~\ref{late:r18:section}.

The auxiliary source comparisons are proved in
Sections~\ref{late:r10:section} and~\ref{late:r11:section}.

\subsubsection{R08 and the reverse-pair payment}\label{late:r08:section}
\begin{theorem}\label{late:r08:exclusion}
R08 does not occur as the positive-maximal-arc word of an original
global minimizer.
\end{theorem}
\begin{proof}
Use native word $02102120120$, $B=11$. Set
\[
 x=c_2+c_3,\quad Q=c_4+c_5+c_6,\quad R=c_{10}+c_0+c_1,
 \quad M=c_7+c_8+c_9,\quad P_*={180483}/{500000}.
\]
The descending six-source $(8,7,6,5,4,3)$ under full reflection
$j\mapsto1-j$ is $012021$. Corollary~\ref{ek:short-sources} gives
$Q>61/200$, and Theorem~\ref{loc:S405} gives $Q+c_7>81/200$.
Native $[-3,2]=012021$ gives $R>61/200$ and $R+c_9>81/200$
by the same two results. The proper forward-six sources
$[5,10]$ and $[6,11]$ give $c_7+c_8>P_*$ and $c_8+c_9>P_*$.

The prefix on $[-1,4]$ with $c_2\le13/20$ gives the positive
terminal root of its one-sided leaf exactly as in
Lemma~\ref{late:r18:prefix}. In Q's reversed chart, its endpoint
$w=b+iv$ maps back by
$F(w)=(-\overline w-3)/(1-\overline w)$. Positivity of its real
part is
\[
 \frac{b^2+2b-3+v^2}{(1-b)^2+v^2}>0,
\]
which supplies Q's missing terminal circle. Together with its
entrance circle this proves $Q>2/5$ by the circle-only tail
comparison. More explicitly, its determinant parameter satisfies
$\Gamma=3/4-\mathcal L/4$. Were $\mathcal L<59/1000$, one would have
$\Gamma>2941/4000>0$; together with the actual $K>0$, this contradicts
Lemma~\ref{cmp:loss-obstruction}. Thus the loss bound required by
Proposition~\ref{six:prop:tail} is obtained without assuming harmonic
calibration in advance. The combined support is $[-1,8]$, ten arcs, so one
genuine period suffices. The reflected prefix similarly gives
$c_3\le13/20\Longrightarrow R>2/5$, on the proper support
$[-3,6]$. In these applications determinant positivity follows
from the tail-loss comparison before the positive mean supplies
the future sheet and harmonic calibration.

If $c_2>13/20$, the unconditional tail bounds imply
$c_3<54483/100000<13/20$. Then the disjoint masks
$\{2\},\{4,5,6,7\},\{8,9\},\{10,0,1\}$ give
\[
 C>\frac{13}{20}+\frac{81}{200}+P_*+\frac25
   =\Cstar+\frac{1392}{125000}.
\]
The full reflected self-map $j'_i=1-j_{5-i}$ gives the same
contradiction for $c_3>13/20$. Hence any survivor has both
centres cheap, $Q,R>2/5$, and $M>P_*$. The disjoint partition
then gives
\[
 x<\Cstar-\frac45-P_* =\frac{80483}{125000}<\frac{13}{20}.
\]

The proper source $[-1,6]=20210212$ charges exactly the pair
$\{2,3\}$. Theorem~\ref{opt:pair650} gives $x>13/20$, contradicting
the preceding strict upper bound. Equivalently, the disjoint charges
$\{2,3\}$, $Q$, $R$ and either complete forward pair in $M$ give
\[
 C>\frac{13}{20}+\frac45+\frac{180483}{500000}
 =\frac{905483}{500000}
 =\Cstar+\frac{767}{125000}.
\]
No source matrix is identified with the lift used for a different scalar
charge; the two prefix forcings and final calibration are already
proved together inside Theorem~\ref{opt:pair650}.
\end{proof}
\space
\subsubsection{The complete physical partitions and reflection}

\begin{theorem}\label{late:r10:exclusion}
Neither $R10$ nor $R13$ is the positive-maximal-arc word of an original
global minimizer.
\end{theorem}
\begin{proof}
For $R10$ retain the following unconditional physical charges:
\[
 K=c_7+c_8+c_9>2/5,
\]
\[
 J=c_{10}>63/500,\qquad T=c_{11}+c_0+c_1>61/200.
\]
The sources and orientations are as follows. Descending $(11,10,9,8,7,6,5)$ under $j\mapsto-j$ is $0120212$.
Its actual $202$ core is $(9,8,7)$, its caps are $10,6$, and its
return is $5$. Its forward prefix supplies the negative entrance and
entrance circle, while the true return supplies the terminal circle.
Proposition~\ref{six:prop:tail} gives $K>2/5$.
The complete nine-source $[6,14]=210120102$ becomes $102012021$
under phase $+2$, retaining both actual outward forward-four witnesses.
Theorem~\ref{pf:opposite126} gives $J>63/500$ without a HIGH predicate.
The source of the ordinary $T$ bound is that of
Lemma~\ref{late:r10:circle-transfer}.

There are two exhaustive cases, according as $c_3\le13/20$ or
$c_3>13/20$.

Suppose first $c_3\le13/20$, including equality.
Lemma~\ref{late:r10:circle-transfer} gives $T>2/5$.
The complete proper source $[-1,6]=20210212$ charges its local
positions $3,4$, namely native arcs $2,3$. Theorem~\ref{opt:pair650}
gives $c_2+c_3>13/20$, on a genuine containing interval obtained
by a cut inside arc $-2$.
The original refined one-return Theorem~\ref{ek:E}, on the actual
source $[3,9]=0212101$, gives
\[
 c_4+c_5+c_6+c_7>E_*>{237\over500},\qquad
 E_*-{237\over500}={2033896299\over8759490110000}>0.
\]
Together with Lemma~\ref{late:r10:exclusive155}, the disjoint masks
\[
 \{2,3\},\quad\{4,5,6,7\},\quad\{8,9\},\quad
 \{10\},\quad\{11,0,1\}
\]
partition all twelve physical arcs. In particular shared arc $7$ is
paid only by the four-core. Consequently
\[
 C>{13\over20}+E_*+{31\over200}+{63\over500}+{2\over5}
 >{361\over200}=C_*+{17\over100000},
\]
a contradiction. The exact $E_*$ is retained; the last inequality
merely displays a smaller rational margin. No separate expensive
$c_2$ branch or weak-five charge is needed.

If $c_3>13/20$, the same unconditional rows give
\[
 c_2+c_4+c_5+c_6<{32383\over100000}<{13\over20}.
\]
Thus $c_2<13/20$ from unconditional estimates, and
Lemma~\ref{late:r10:one-sided} gives $c_5+c_6>283/1000$.
The resulting disjoint budgets are
\[
 c_2+c_4<{4083\over100000},\qquad
 {13\over20}<c_3<{69083\over100000}.
\]
Lemma~\ref{late:r10:height} gives the actual own height $h_E>1/12$.
The reverse prefix on $(6,5,4,3,2,1)$ under $j\mapsto1-j$ is
$202102$, with centre $3$. Its initial $202$ middle mean is
exactly the transformed own mean leaf $5$, not the one-sided mean
on arc $3$ of the incoming construction. Its positive outer root
$U$ satisfies $U>h_E$, and its following positive mode-$2$ tail
increases $u$. Thus this actual prefix has $u_L>U>h_E>1/12$.
Lemma~\ref{late:r10:enhanced-prefix} restores its strict exit sign.
Lemma~\ref{late:r10:circle-transfer}, on the genuine period
$[-3,8]$, then gives $T>2/5$ even though $c_3>13/20$.
The disjoint masks are now
\[
 \{3\},\quad\{5,6\},\quad\{7,8,9\},\quad
 \{10\},\quad\{11,0,1\},
\]
leaving only positive arcs $2,4$ uncharged. Consequently
\[
 C>{13\over20}+{283\over1000}+{2\over5}+{63\over500}+{2\over5}
 ={1859\over1000}=C_*+{5417\over100000}.
\]
This excludes both cases. Internal identities coupling means
or matrices used the full common lifts specified above. The final
additions concern only physical costs, so different scalar estimates
may be established on different genuine containing periods.

To transfer the entire conclusion to $R13$, use the lawful whole
reflected-time transformation and phase with label map
\[
 j'_i=1-j_{5-i},\qquad r\longmapsto i=5-r\pmod{12}.
\]
Its native word is $021021201210$. Its twist is still $+1$, since
\[
 j'_{i+12}=1-j_{5-i-12}=j'_i+1\pmod3.
\]
The transformation carries the entire path, full state and frame
endpoints, real winding, the state covector, matrix adjoint,
normal multiplier, and every positive maximal support arc. It
neither presumes reflection symmetry of a trajectory nor identifies
the two cyclic words by a reflection quotient. The predicates $c_2,c_3$
of $R10$ become respectively $c'_3,c'_2$ of $R13$.

For $c'_2\le13/20$, the reflected first-case masks are
\[
 \{2,3\},\quad\{0,1,10,11\},\quad\{8,9\},\quad
 \{7\},\quad\{4,5,6\}.
\]
They partition all twelve arcs and have total greater than $361/200$.
For $c'_2>13/20$, the reflected second case has masks
\[
 \{2\},\quad\{0,11\},\quad\{8,9,10\},\quad
 \{7\},\quad\{4,5,6\},
\]
omits only arcs $1,3$, and has total greater than $1859/1000$.
The common nine-source $[3,11]$ becomes the complete reversed
source $[-6,2]$, and the eleven-source $[-1,9]$ becomes $[-4,6]$.
Their genuine containing periods transform as whole periods;
source length, cap completeness and calibration are preserved at
the seam. The same two strict rational margins therefore exclude
$R13$ as well.
\end{proof}
\space
\subsubsection{The disjoint partition}
\begin{theorem}\label{late:r11:exclusion}
Neither R11 nor R12 occurs as the positive-maximal-arc word of an
original global minimizer.
\end{theorem}
\begin{proof}
First suppose $c_2>13/20$. Lemma~\ref{late:r11:seven} and
\eqref{late:r11:charges} give
\[
 c_3<C_*-\frac{13}{20}-\frac{91}{200}-2\frac{61}{200}
 =\frac{8983}{100000}<\frac{13}{20}.
\]
Thus Lemma~\ref{late:r11:cheap-circles} gives $R>2/5$.
The disjoint masks $\{2\}$, $\{7,8,9,10\}$,
$\{4,5,6\}$ and $\{11,0,1\}$ now imply
\[
 C>\frac{13}{20}+\frac{91}{200}+\frac{61}{200}+\frac25
 =\frac{181}{100}=C_*+\frac{517}{100000},
\]
a contradiction.

It remains $c_2\le13/20$, including equality. Then
Lemma~\ref{late:r11:cheap-circles} gives $Q>2/5$. The unconditional
bridge payment and $R>61/200$ force
\[
 c_3<C_*-\frac{91}{200}-\frac25-\frac{61}{200}
 =\frac{64483}{100000}<\frac{13}{20}.
\]
Consequently $P,Q,R>2/5$. The complete proper source
$[-1,6]=20210212$ has eight arcs; its local charged positions $3,4$
are exactly native arcs $2,3$. Theorem~\ref{opt:pair650} applies on a
genuine containing interval and gives $c_2+c_3>13/20$. These are
scalar costs on the same physical trajectory; their full lifts are
not identified. The resulting disjoint partition gives
\[
 C>c_2+c_3+P+Q+R>\frac{13}{20}+3\frac25
 =\frac{37}{20}=C_*+\frac{4517}{100000},
\]
again a contradiction. Arc 7 is left unpaid in this final partition.
No bridge or overlapping triple is charged twice.

Finally the full reflected reanchoring
\[
 j'_i=1-j_{5-i}
\]
takes R11 to R12, whose native word is $021021201020$, and preserves
the $+1$ twist. It carries the state, group, matrix, both state
covectors, endpoints and real winding. The physical index map is
$r\mapsto5-r$ modulo twelve. The centres 2 and 3 are interchanged,
and the three charged triples become
\[
 P'=\{7,8,9\},\qquad Q'=\{11,0,1\},\qquad
 R'=\{4,5,6\}.
\]
The bridge 7 maps to arc 10, and the pair $\{2,3\}$ is retained.
All strict costs and surpluses therefore transport to R12.
\end{proof}
\space

\begin{theorem}\label{run:five-survivors}
An original minimizing word in the length-five-negative-run family
must have type R02 or R04.
\end{theorem}
\begin{proof}
The family consists of the seven listed words. The preceding
source-qualified prefix and bilateral comparisons exclude R08 and
R10--R13. The remaining two words are R02 and R04.
\end{proof}
\space
\subsection{A single negative run of length seven, eight, or nine}
\label{run:long}

The remaining mixed group in Theorem~\ref{run:21} has one negative run of
length $7$, $8$ or $9$; every other negative run has length $2$.
Its three one-positive-run and four two-positive-run words are specified in Table~\ref{run:table}. Reflection is used only by transporting the full geometric construction, not by identifying cyclic types.
\begin{theorem}[Exclusion of the long-negative-run family]
\label{run:long-exclusion}
No original global minimizer belongs to the seven-word long-negative-run
family of Theorem~\ref{run:21}.
\end{theorem}

For F41 and F40 we determine which rooted response is positive.
For F35 and F36 one expensive centre suffices, since a common
complement avoids both possible centres. For F24a, F24b and F39,
we choose the five-arc window that pays for a disjoint triple
according to which centre supplies the required lower bound.

Throughout this section, the native controls are reconstructed from the
printed increment word by $j_0=0$ and $j_{i+B}=j_i+1$ modulo three.
A source $a:b$ lists all indices from $a$ to $b$ in their indicated
order. On a decreasing source the displayed label reflection accompanies
full time reversal. Physical costs alone are reduced modulo $B$.
All the sources below are complete and proper. The full geometric
endpoints, real winding and selected normal lifts are those of
Section~\ref{fd:section}.

\subsubsection{Rooted alternatives and their actual continuations}
Let $U,V$ be the horizontal endpoints of an actual whole mode-$2$ tail
rooted at a $0\to2$ switch in its lawful chart. Its rooted response is
\[
 \Psi(U,V)=2(1-U)\log\frac{1-U}{1-V}
                    -(V-U)(2+V-U).
\]
This is the full-duration response of \eqref{pc:Psi}, not an assigned
adjoint at a comparison endpoint. Theorem~\ref{pc:paired-response}
says that a complete $21021$ context forces a positive response in at
least one of the two lawful orientations at its middle connector.
The two alternatives must first belong to the stipulated common full
lift. A short rigid capped source supplies such a lift by
Corollary~\ref{pc:short-rigid-global}; its covector twist is proved
there, not assumed.

\begin{lemma}[A blocked continuation]\label{nx:02101}
A complete actual $02101$ source on a qualified full normal maximizing
lift has nonpositive rooted response at its first $0\to2$ switch.
\end{lemma}
\begin{proof}
The response identity \eqref{nx:root-response} converts a positive
response into the strict actual ingress of \eqref{hc:ingress}.
The retained terminal $101$ gives the opposite reflected-tail inequality
on this same lift, contradicting \eqref{sw:two-tail}.
The full binary derivation is in Section~\ref{nx:rooted-proofs}.
Equality of the response is retained.
\end{proof}

\begin{lemma}[Two successful continuations]\label{nx:rooted-continuations}
A complete $021020$ source with positive rooted response has central
cost greater than $299/250$. A complete $02102101$ source with positive
rooted response has central cost greater than $13/20$. Both statements
use a qualified common full lift, and the charged arc is local position $2$.
\end{lemma}
\begin{proof}
The positive response gives strict ingress at the actual central arc.
For $021020$, both central contact jets and the full terminal $020$
are retained, so Theorem~\ref{nx:central1196} applies.
For $02102101$, a contrary central cost at most $13/20$ gives the
prefix exit inequality. Lemma~\ref{hc:regeneration} propagates it along
the two following complete arcs to a fresh strict ingress; the actual
terminal $101$ then contradicts the preceding blocked-continuation
argument. Only the prefix estimate is iterated.
Section~\ref{nx:rooted-proofs} gives the equations and raw-switching
regularity at both joins.
\end{proof}

For the next two words the rigid sources, in native indexing, are
$[7,11]+2=12021$ at $B=9$ and $[7,12]+2=102021$ at $B=10$.
Corollary~\ref{pc:short-rigid-global}, using genuine cuts inside arcs
$2$ and $7$ and their translates, supplies one normalized lift on the
unfolding. Its restrictions retain every complete response source used
below. The scalar weak-five payments may use separate genuine lifts;
no covector identification is needed when their physical costs are added.

\begin{proposition}\label{nx:rooted9}
The word $\mathtt{-------++}$, of length $9$, is excluded.
\end{proposition}
\begin{proof}
The complete five-source $[1,5]=21021$ has connector $3$.
Its reverse source $(3,2,1,0,-1)$ with whole map $j\mapsto-j$
is $02101$, so Lemma~\ref{nx:02101} blocks that alternative.
The response in the forward source $[3,8]=021020$ is therefore
positive, giving $c_5>299/250$.

Resolve a second disjunction at connector $4$, using
$[2,6]+1=21021$. Its forward $[4,8]+1=02101$ is blocked.
The reverse source $(4,3,2,1,0,-1)$ with whole map $j\mapsto2-j$
is $021020$, giving $c_2>299/250$.
All response supports have five or six distinct physical arcs and are
retained complete in the indicated genuine periods. Each singleton
bound follows from its own resolved disjunction. Since the two charges
are distinct,
\[
 C>c_2+c_5>\frac{299}{125}>C_*.
\]
\end{proof}

\begin{proposition}\label{nx:rooted10}
The word $\mathtt{---------+}$, of length $10$, is excluded.
\end{proposition}
\begin{proof}
Again $[1,5]=21021$ has connector $3$. The reverse source
$(3,2,1,0,-1)$ with whole map $j\mapsto-j$ is $02101$ and is
blocked. The forward source $[3,10]=02102101$ has positive response,
so $c_5>13/20$.
The two proper five-sources $[0,4]+1=10210$ and
$[5,9]=10210$ charge $\{1,2,3\}$ and $\{6,7,8\}$ at least
$3/5$ each by Lemma~\ref{strip:five-cap-completion}.
Cuts in physical arcs $5$ and $0$ respectively retain these sources.
Their charges and $\{5\}$ are disjoint, so
\[
 C>\frac35+\frac{13}{20}+\frac35
   =\frac{37}{20}>C_*.
\]
Every response continuation has at most eight arcs, fewer than $10$;
its caps and terminal arcs have not been deleted.
\end{proof}

\subsubsection{The two length-eight negative runs}
\label{pf:M135}\label{pf:F35}\label{pf:M136}\label{pf:F36}
The remaining words with a negative run of length eight use the same
paired-centre geometry. Here a single maximum bound suffices because
one common scalar complement avoids both possible centres.

\begin{proposition}[A paired centre and its complement]
\label{pc:B11}\label{pc:B12}
The following two words are excluded:
\[
 \mathtt{--------+++},\qquad\mathtt{--------+--+}.
\]
Their lengths are $11$ and $12$, respectively.
\end{proposition}
\begin{proof}
Use the following rigid sources and genuine cuts:
\[
\begin{array}{c|c|c}
 B&\text{rigid source}&\text{cuts inside physical arcs}\\\hline
11&[9,13]+2=12021&2,9\\
12&[6,11]+1=102021&-1,6.
\end{array}
\]
Corollary~\ref{pc:short-rigid-global} gives a common normalized lift
and its proved twist. The rigidity overlap needs only positive cap
portions; each cost comparison below still retains its full proper
source. At $B=11$ the resulting periods retain
$[3,12]$, $[10,19]$, and, after translating the cut $9$ by $-11$,
$[2,6]$. At $B=12$ they retain $[0,10]$, $[7,17]$ and $[2,6]$.
These are genuine minimizing periods with the prescribed endpoints
and real winding.

In both words $[2,6]+1=21021$, with connector $4$. The ordered sources
\[
 (4,5,6,7,8,9),\quad j\mapsto j+1;
 \qquad (B+4,B+3,B+2,B+1,B,B-1),\quad j\mapsto-j
\]
are $021020$. Their ingress switches are $4|5$ and
$(B+4)|(B+3)$; their centres are physical arcs $6$ and $2$.
Their full terminal triples are $(7,8,9)$ and $(B+1,B,B-1)$.
The proved twist identifies connector $B+4$ with $4$ on the common
lift. Theorem~\ref{lc:paired} therefore gives
\[
 \max(c_2,c_6)>299/250.
\]
It does not give two singleton inequalities.

The complementary complete scalar sources are
\[
\begin{array}{c|c|c|c}
 B&\text{native source}&\text{cut}&\text{physical charge}\\\hline
11&[6,13]=02120102&5&\{7,8,9,10,0,1\}\\
12&[6,14]=021210102&5&\{7,8,9,10,11,0,1\}.
\end{array}
\]
Their costs are greater than $71/100$ by Theorem~\ref{sj:W630},
and at least $16/25$ by Theorem~\ref{loc:AA640}, respectively.
Both sources are proper and retain their two exterior caps.
Each complement avoids both possible centres $2,6$; its separately
chosen scalar lift need not be identified with the paired lift.
Paying only a successful centre gives
\[
 C>\frac{299}{250}+
 \begin{cases}71/100,&B=11,\\16/25,&B=12,
 \end{cases}
 =\begin{cases}953/500,&B=11,\\459/250,&B=12.
 \end{cases}
\]
Both values exceed $C_*$. The other centre and physical arcs $3,4,5$
remain unpaid.
\end{proof}

\subsubsection{A length-seven run with one short separator}
\label{pf:M098}\label{pf:F24}
All three remaining thirteen-arc words have one negative run of length
seven and one of length two. Their three positive-length pairs are
$(1,3)$, $(2,2)$ and $(3,1)$. They share the proper ten-source
\[
 [12,21]+2=2021021020,
\]
retained by a cut inside physical arc $11$.
Corollary~\ref{hc:ten-arc}, including its regeneration arcs and terminal
$020$, yields the alternative
\begin{equation}\label{pf:two-thirteen-max}
 \max(c_2,c_5)>13/20.
\end{equation}
Every branch below uses just one successful centre.

\begin{proposition}\label{run:long-seven-asymmetric}
The following two words are excluded:
\[
 \mathtt{-------+--+++},\qquad\mathtt{-------+++--+}.
\]
\end{proposition}
\begin{proof}
For $w_-:=\mathtt{-------+--+++}$ and
$w_+:=\mathtt{-------+++--+}$, the following two scalar rows are
available independently of that alternative:
\[
\begin{array}{c|c|c|c}
 &\text{ordered native source}&\text{canonical word}&\text{charge}\\\hline
w_-&[9,14]+2&101202&\{10,11,12,0\}\\
w_-&[5,11]&1020212&\{6,7,8,9\}\\
w_+&[6,11]+1&101202&\{7,8,9,10\}\\
w_+&(15,14,\ldots,9),\ j\mapsto-j&1020212&\{11,12,0,1\}.
\end{array}
\]
The two costs for either word exceed $103/250$ and $473/1000$
by Theorem~\ref{fee:four} and Corollary~\ref{loc:E473}. Their cuts are
respectively $8,4$ for $w_-$ and $5,8$ for $w_+$.
If $c_2>13/20$, use $[2,6]=10210$ in $w_-$, charging
$\{3,4,5\}$, or $[3,7]+1=10210$ in $w_+$, charging
$\{4,5,6\}$. Otherwise \eqref{pf:two-thirteen-max} forces
$c_5>13/20$; use $[1,5]+2=10210$ in either word, charging
$\{2,3,4\}$. The cuts for these five-sources are $1,2$ and $0$,
respectively. Each chosen triple costs at least $3/5$.
The successful centre, selected triple and two fixed scalar charges are
pairwise disjoint in each branch. Hence both words satisfy
\[
 C>\frac{103}{250}+\frac{473}{1000}+\frac35+\frac{13}{20}
   =\frac{427}{200}>C_*.
\]
The branch with $c_2=13/20$ belongs to the second alternative, so the
boundary is included.
\end{proof}

\begin{proposition}\label{nx:mixed13}
The word $\mathtt{-------++--++}$ is excluded.
\end{proposition}
\begin{proof}
Use the same ten-source and centre alternative
\eqref{pf:two-thirteen-max}. The complete ordered sources
$[5,14]$ and $(15,14,\ldots,6)$, with maps $j$ and $-j$ and cuts
$4,5$, are $1020102010$.
The prefix $10201020$ contains the seven-source at positions
$1,\ldots,7$ used by Lemma~\ref{ua:seven-sign}. Its actual unsorted
chart gives $Ea>1$. Reading prefix positions $(4,3,2,1,0)$ with
whole map $j\mapsto2-j$ gives the complete capped source $12021$
with entrance $1-Ea<0$. Theorem~\ref{loc:ST305} therefore charges
local positions $1,2,3$ strictly above $61/200$.
The two physical charges are $\{6,7,8\}$ and $\{12,0,1\}$.
Their supporting caps and sign witnesses remain in the two proper
complete ten-sources.

If $c_2>13/20$, the proper five-source $[2,6]=10210$, cut $1$,
pays at least $3/5$ on $\{3,4,5\}$.
Otherwise $c_5>13/20$, and $[1,5]+2=10210$, cut $0$,
pays at least $3/5$ on $\{2,3,4\}$.
The four selected physical charges are disjoint in either branch, so
\[
 C>2\frac{61}{200}+\frac35+\frac{13}{20}
   =\frac{93}{50}>C_*.
\]
The weak five-source remains weak; strictness is supplied by the selected
centre and both tails. Equality of a centre with $13/20$ is not lost.
\end{proof}

\begin{proof}[Proof of Theorem~\ref{run:long-exclusion}]
The single-run cases of lengths $9$ and $10$ are excluded by
Propositions~\ref{nx:rooted9} and~\ref{nx:rooted10}.
Proposition~\ref{pc:B11} excludes the two words of lengths $11$ and
$12$. Proposition~\ref{run:long-seven-asymmetric} excludes the two asymmetric
thirteen-arc words, and Proposition~\ref{nx:mixed13} excludes the
equal positive-length pair $(2,2)$. These are all seven words in
this family, as classified in Theorem~\ref{run:21}.
\end{proof}
\space
\subsection{The six remaining types}\label{new:uniform}
\begin{theorem}[The mixed alternatives and the all-minus word]\label{new:uniform-six}
Every original global minimizer with $B\ge8$ has one of the essential increment types R01--R06. These types have $B=8$ or $B=9$.
\end{theorem}
\begin{proof}
The preceding three comparisons exhaust the mixed words in
Table~\ref{run:table}, according to their negative runs. In order
of the four rows in the table below, the surviving counts are
$(3,2,0,1)$. For the pure word the all-length classification
leaves exactly R01. This gives the six necessary types in
Table~\ref{ten:table}.
\end{proof}
\begin{center}
\begin{tabular}{@{}lcl@{}}
\toprule
Negative runs & Candidates & Types left by the comparisons\\
\midrule
All of length two & $6$ & R03, R05, R06\\
One of length five & $7$ & R02, R04\\
One of length seven, eight, or nine & $7$ & None\\
Pure all-minus word & $1$ & R01\\
\bottomrule
\end{tabular}
\end{center}

\begin{center}
\begin{minipage}{\textwidth}
\centering
\makeatletter
\def\@captype{table}
\makeatother
\begin{tabular}{ccll}\toprule
Type&$B$&Cyclic increments&Exclusion\\\midrule
R01&8&\texttt{--------}&Theorem~\ref{eight:thm:main}\\
R02&8&\texttt{-----+++}&Theorem~\ref{mix:thm:r02}\\
R03&8&\texttt{--++++++}&Theorem~\ref{mix:thm:r03}\\
R04&9&\texttt{-----+--+}&Theorem~\ref{r04:exclusion}\\
R05&9&\texttt{--+--++++}&Theorem~\ref{nine:thm:r05}\\
R06&9&\texttt{--++--+++}&Theorem~\ref{r06:exclusion}\\\bottomrule
\end{tabular}
\caption{The six necessary long types and their final exclusions.}\label{ten:table}\label{five:table}
\end{minipage}
\end{center}
\space

\section{Exclusion of R01}
\subsection{The common geometric and covering principle}\label{eight:common}
For the three eight-arc types in Table~\ref{ten:table}, we use full signed matrix closure. In the native label convention they are
\[
\begin{array}{ccl}
\mathrm{R01}&02102102&\texttt{--------},\\
\mathrm{R02}&02102120&\texttt{-----+++},\\
\mathrm{R03}&02120120&\texttt{--++++++}.
\end{array}
\]
In each row the labels continue by $j_{i+8}=j_i+1$. Thus the final increment is read across the true twisted seam. The common feature is a full signed matrix closure, which reduces eight complete arc multipliers to five coordinates. The three resulting charts need not have the same domain or the same rejection predicates.

\subsubsection{The physical subset of a closed parameter domain}\label{eight:common-candidates}
In the eight-arc comparisons a candidate is an original global minimizer in Section~\ref{fd:section}, with the indicated eight positive maximal arcs. Its state lies in $K$, its physical cost satisfies $C\le C_{\rm oct}<C_*$, and its frame and state satisfy \eqref{fd:endpoints} with the prescribed real winding. The complete multiplier of each arc is $q_i>1$; its physical duration is $d_{j_i}\log q_i$. In particular $q_i$ is neither a duration nor a square-root parameter for a half-return.

Each chart below is obtained from the full state and frame equations. A lawful change of cut or whole reflection transports the physical trajectory and its complete source qualifications. It can select a representative in a smaller closed root domain. The root may also contain points with a zero limiting duration, extra winding, or no physical reconstruction. Inclusion of such points enlarges the search domain and is harmless: the required implication is from every strict physical candidate to at least one root box. No converse is used.

Every local cost or state inequality is applied with its complete supporting source, on a genuine containing interval. Proposition~\ref{fd:normality} permits joint normalization of the selected full maximizing lift. A scalar physical conclusion can subsequently be combined with another such conclusion. This operation does not identify independently selected covectors on different windows.

\begin{lemma}[Finite exclusion of the physical subset]\label{eight:common-cover}
Let $\mathcal S$ be a set of qualified physical candidates, and let each candidate have a coordinate representative in a finite union of closed root boxes. Suppose a finite covering construction has the following properties. Every split covers its parent, including the common boundary; every contraction retains the representatives of all candidates in the box; and every terminal box has a verified strict contradiction to a necessary condition for its represented candidates. Then $\mathcal S$ is empty.
\end{lemma}
\begin{proof}
If a candidate existed, choose one root box containing its representative. A contraction retains that point and a covering split places it in at least one child. Since the construction is finite and all branches terminate, the point reaches a terminal box. Its verified contradiction is incompatible with the candidate. Overlaps of closed boxes can repeat the point but cannot remove it.
\end{proof}

The geometric inclusion, the finite covering identities and the numerical signs are separate obligations in this lemma. File identities and counts bind the actual inputs to a checking record; they supply neither a geometric inclusion nor a sign proof. Unknown arithmetic outcomes and unfinished branches cannot serve as terminal contradictions.

For R01, some contradictions use a necessary local-minimum condition for the complete moving-state cost. For R02 and R03, the accepted terminal predicates use only necessary domains, physical costs, state inequalities and full closure values. Their Taylor or polynomial bounds may differentiate the explicit functions to enclose their values on a box. Such differentiation does not add a stationarity equation or a positive-semidefinite Hessian condition to those two candidate sets. 
\space
\subsection{The all-minus word and full closure}\label{eight:sec:foundations}
We exclude the strict essential eight-arc control word $02102102$ from
global minimizers in the full Reinhardt admissible class. Full signed
frame closure admits a smooth five-dimensional Gram chart. A
four-duration inequality makes its state return elliptic and gives a
stable expression for the complete moving-state cost. Physical changes
of cut and reflection, together with the compactness bound for
minimizers, reduce the problem to a finite rational cover. Rigorous
interval bounds then contradict the comparison bound or a necessary
local-minimum condition throughout this cover.

We use the original geometric minimizing class and normalization of
Section~\ref{fd:section}, with its cited control correspondence and
reconstruction inputs. The following exclusion concerns that class
before imposing a switching word or compact-set restriction.

\begin{theorem}\label{eight:thm:main}
An original globally minimizing Reinhardt control cannot have eight
positive maximal arcs with successive labels
\[
                 0,2,1,0,2,1,0,2,
\]
where the twisted continuation has label $j_{i+8}=j_i+1\pmod3$.
Equivalently, its strict essential type cannot be $B=8$ with all eight
cyclic increments negative. We call this type R01.
\end{theorem}

The assertion is about eight essential control arcs, not eight sides of
a polygon. We prove it directly, without using a census of other types or
an upper bound on their number of arcs. Its role in the present classification is to remove R01 from the six-type conclusion of Theorem~\ref{new:uniform-six}. None of the other whole-class exclusions is a premise of this argument.

\subsubsection{Coordinates and physical time}
Use \eqref{fd:domain}--\eqref{fd:controls}, with the local notation
\begin{equation}\label{eight:eq:matrices}
 w=u+iv,\qquad A=A_*,\qquad R=R_*,\qquad \cost=C.
\end{equation}
Thus the full measurable control triangle, physical time ($X^2=-I$),
state and signed frame endpoints, and prescribed real winding are
unchanged. The reconstructed area is $\sqrt3\cost$; neither trace closure
nor an angle only modulo $2\pi$ suffices. The inverse control map
\eqref{fd:controls} retains the original barycentric class without a time
change, with $(x,y)=(u,v)/\sqrt3$ and
$\mu=(\alpha_2,\alpha_0,\alpha_1)$ in \cite{HV2024}.

Theorem~\ref{fw:finiteness} gives the finite representation used to count
essential arcs. The present exclusion uses the stated eight-arc hypothesis.
Zero-duration pieces are deleted and adjacent equal controls
are merged in the \emph{twisted} extension, including its seam. The
remaining arcs are maximal and have positive time. A cut at a switch
gives exactly $B$ arcs; a cut inside an arc may give $B+1$ pieces.

\subsubsection{Two imported bounds and their consequences}
The imported compactness theorem and smoothed-octagon comparison are
Theorem~\ref{fd:compact} and \eqref{fd:budget}. In this notation
$\budget=\Cstar$ and $\cost\le C_{\rm oct}<\budget$.
The cuts \eqref{fd:compact-cuts} constrain original minimizers, not local
competitors. The rational check of the strict budget is given with
\eqref{fw:budget}.

\begin{lemma}\label{eight:lem:cap}
On an original minimizing trajectory,
$99/200\le v\le39/5$. Every complete vertex arc has multiplier
$1<q\le Q:=\Qcap$, where its physical duration is $d_j\log q$.
\end{lemma}
\begin{proof}
The height bounds are \eqref{fd:compact-bounds}, proved in
Lemma~\ref{fw:compact-bounds}. Rotate a complete arc cyclically into mode $0$. Invariance of $K$ makes
the same height bounds valid at both endpoints. Its exact flow is
$w(p)=p(w(1)+1)-1$, $1\le p\le q$, so $v_{\rm out}=qv_{\rm in}$ and
$q\le(39/5)/(99/200)=520/33$. Since $d_0=(1+u)/v$ is constant,
$d\log p/dt=1/d_0$; its elapsed time is $d_0\log q$, not $\log q$.
\end{proof}

\subsubsection{Full closure in duration variables}
Write $B_j=\mathsf B_j$ for the three control matrices in
Section~\ref{sec:framework}, so $d_j=\operatorname{tr}(B_jX)$ as in
\eqref{fw:reciprocals}.
On a vertex-$j$ arc, $d=\tr(B_jX)$ is constant and positive. For
$h=\log q$, let $E=\exp(hB_j)$ and $Q_j(q)=E^{-2}$.
Direct integration gives
\begin{equation}\label{eight:eq:arc-frame}
 X_{\rm out}=EX_{\rm in}E^{-1},\qquad
 M=I+d\sinh h\,X_{\rm in}-2d(\cosh h-1)X_{\rm in}B_j
   =X_{\rm in}Q_j(q)X_{\rm out}^{-1}.
\end{equation}
The formulas can also be checked by differentiation:
$p\,dM/dp=dMX(w(p))$ and $M(1)=I$.
Thus the physical-time and frame scales agree.

For the R01 labels $j_i=-i\pmod3$, put
\begin{equation}\label{eight:eq:N}
 \Nmat(q)=\begin{pmatrix}q&-q\\q&q^{-1}-q\end{pmatrix},\qquad
 P=\begin{pmatrix}1&-2\\1&0\end{pmatrix}.
\end{equation}
Since $B_{-i}=R^iB_0R^{-i}$ and
$Q_0(q)R=P\Nmat(q)P^{-1}$, multiplication gives
\[
 \prod_{i=0}^7Q_{j_i}(q_i)
 =P\left(\prod_{i=0}^7\Nmat(q_i)\right)P^{-1}R^{-8}.
\]
The frame telescope in \eqref{eight:eq:arc-frame}, together with the twisted
state endpoint, makes $G(T)=R$ equivalent to $\prod Q_{j_i}=R$.
As $R^3=-I$, this is exactly
\begin{equation}\label{eight:eq:signed-closure}
                  \Nmat(q_0)\cdots\Nmat(q_7)=-I.
\end{equation}
In particular the sign in this equation is not dispensable.
\space

\subsection{The Gram chart and its finite cover}
Full signed closure gives the five-dimensional Gram chart of
Theorem~\ref{eight:thm:chart}. The four-duration gap inequality
gives an elliptic state half-return, so the moving state is recovered
as its unique upper-half-plane fixed point
(Theorem~\ref{eight:thm:elliptic}).
Lemma~\ref{eight:lem:stationarity} varies the complete moving-state
cost in this chart with the actual endpoints and winding. Hence its
gradient vanishes and its Hessian is positive semidefinite at an
original interior minimizer.

A physical change of cut and whole reflection yield the canonical
representatives in Section~\ref{eight:sec:cover}. The multiplier cap
and angular inequalities place them in the closed middle and endpoint
roots there. Proposition~\ref{eight:prop:finite} verifies the complete
covers of these roots. The terminal predicates use necessary physical
and chart conditions, the octagonal cost bound and the local-minimum
conditions just derived. Appendix~\ref{app:eight-details} gives the
chart, the inclusions, the terminal rules and their arithmetic bounds.
The checking tasks are specified in Section~\ref{reader:eight}.

\begin{proof}[Proof of Theorem~\ref{eight:thm:main}]
Suppose an original minimizer has the strict R01 type. Cut at one of its
eight switches and normalize its first label. Equations
\eqref{eight:eq:arc-frame}--\eqref{eight:eq:signed-closure} give full signed closure,
and Theorem~\ref{eight:thm:chart} places it in the open Gram cube.
Its physical reconstruction is the unique elliptic moving-state branch
of Section~\ref{eight:sec:cost}. Lemma~\ref{eight:lem:stationarity}
supplies the zero-gradient and positive semidefinite Hessian conditions.
Choose a lawful representative satisfying \eqref{eight:eq:canonical},
with an arc of minimum multiplier first. It is still an original
minimizer, so it obeys the cap of Lemma~\ref{eight:lem:cap}, the strict budget
\eqref{fd:budget}, and the angular bounds of Lemma~\ref{eight:lem:angle}.
The closed middle and endpoint covers contain it. This contradicts
Proposition~\ref{eight:prop:finite}.
\end{proof}

A zero multiplier increment $q_j-1=0$ has zero physical time:
$Q_j(1)=I$ and $Q_j(a)Q_j(b)=Q_j(ab)$. Deleting such an arc and merging
equal neighbors preserves the original state, frame and winding, and
leaves at most seven essential arcs. One must not instead delete
$\Nmat(1)$, which is not the identity. These lower essential types are
outside the theorem. In particular, the two alternating finite boundary
points of the R01 chart can encode the four-arc smoothed-octagon
competitor; strict interior descent does not contradict its role in
\eqref{fd:budget}.
\space
\section{Exclusion of R02}\label{mix:r02}
\begin{theorem}\label{mix:thm:r02}
An original global minimizer cannot have eight positive maximal arcs
with native labels $\mathtt{02102120}$ and twisted continuation
$j_{i+8}=j_i+1$. This excludes the strict essential type R02,
with increments $\mathtt{-----+++}$.
\end{theorem}

The full signed frame gives \eqref{mix:signed}. We use the positive
chart of Lemma~\ref{mix:r02-chart} and the common physical cost
circuit \eqref{mix:numerators}--\eqref{mix:positive-cost}.
We next place the strict source in a closed root and specify the
necessary inequalities for the cover.

\subsection{Complete sources and a compact root}\label{mix:r02-domain}
Throughout the source tables, $[s,t]$ denotes the increasing list of unwrapped arc indices when $s\le t$ and the decreasing list when $s>t$. Labels at unwrapped indices use $j_{i+8}=j_i+1$; physical charges are reduced modulo eight. The column $\ell(j)$ gives the label map on the entire ordered source. A decreasing list is transported by whole reflection with time reversal, including its state, frame, covectors and charges.

\begin{table}[htbp]\centering\small
\begin{tabular}{@{}clcll@{}}\toprule
Row & Charged arcs & $b_j$ & Complete source & Ordered support; $\ell(j)$\\\midrule
0 & $1,2,3$ & $3/5$ & $10210$ & $[0,4]$; $j+1$\\
1 & $2,3,4$ & $3/5$ & $10210$ & $[1,5]$; $j+2$\\
2 & $0,1,6,7$ & $81/200$ & $012021$ & $[5,10]$; $j+2$\\
3 & $4,5,6,7$ & $81/200$ & $012021$ & $[0,-5]$; $-j$\\
4 & $0,5,6,7$ & $103/250$ & $101202$ & $[1,-4]$; $-j$\\
5 & $0,1,7$ & $61/200$ & $012021$ & $[5,10]$; $j+2$\\
6 & $4,5,6$ & $61/200$ & $012021$ & $[0,-5]$; $-j$\\
7 & $0,1,4,5,6,7$ & $71/100$ & $02120102$ & $[3,10]$; $j$\\
8 & $2,3$ & $13/20$ & $20210212$ & $[-1,6]$; $j$\\\bottomrule
\end{tabular}
\caption{The nine R02 inequalities $\sum_{i\text{ charged}}c_i\ge b_j$. Only rows 0 and 1 are intrinsically weak; the other strict bounds are safely weakened in \eqref{mix:necessary-values}.}\label{mix:r02-fees}
\end{table}

Rows 0--1 are Theorem~\ref{wk:five600}; rows 2--3 are Theorem~\ref{loc:S405}; row 4 is Theorem~\ref{fee:four}. For rows 5--6, the full $012021$ source gives a strictly negative
entrance for its $12021$ suffix by Lemma~\ref{cmp:forward-prefix-forcing}.
Corollary~\ref{ek:short-sources} applies Theorem~\ref{loc:ST305} on the
same lift to local positions $2,3,4$. Under the two maps in the table these
are $\{0,1,7\}$ and $\{4,5,6\}$, respectively. The supports of rows 0--6 are
proper. Rows 7 and 8 have eight-arc supports and use the whole-period permissions of Theorems~\ref{fee:W71} and~\ref{opt:pair650}, respectively. For row 7 the two tails and common switching derivative belong to one selected complete lift; for row 8 both caps and rooted contacts lie inside its genuine switch-cut interval. Neither application continues a state covector periodically or claims strict properness. All rows retain the normality and joint normalization stipulated by their individual theorems. They are scalar facts about the same physical candidate; their addition does not glue covectors selected on different intervals.

Whole reflection and reanchoring give the exact R02 involution
\begin{equation}\label{mix:r02-reflection}
 j'_r=1-j_{5-r}=j_r\pmod3,\qquad q'_r=q_{5-r},\qquad c'_r=c_{5-r}.
\end{equation}
The twisted seam is unchanged, as $j'_{r+8}=j'_r+1$. The lawful symmetry of Section~\ref{fd:section} transports the full endpoint and real winding. Thus we choose $c_2\le c_3$. Row 7 and the budget give
\begin{equation}\label{mix:r02-cheap}
 c_2+c_3<\frac{109483}{100000},\qquad
 c_2<\frac{109483}{200000}<\frac{13}{20},\qquad c_3>\frac{13}{40},
\end{equation}
the last inequality following from row 8 and the ordering.

The two proper sources $[-1,4]$ and $[6,1]$, with maps $j$ and $1-j$, are both $202102$. On each containing interval, the initial $202$ and its last complete arc supply the strict ingress of \eqref{hc:ingress} to the respective centre, physical arc 2 or 3. Write $d_j^{\rm nat}(i)=d_{j-j_i}(i)$, with the reciprocal index taken modulo three, for the native-label component at cut $i$. The original and reflected horizontal coordinates have the signs of $d_0^{\rm nat}-d_2^{\rm nat}$ and $d_1^{\rm nat}-d_2^{\rm nat}$, respectively. The crossing estimate of Lemma~\ref{hc:crossing}, applied first from cut 2 to 3 and then from cut 4 to 3, gives
\begin{equation}\label{mix:r02-crossing}
 d_2(2)>d_1(2),\quad d_2(3)>d_0(3),\quad
 d_2(3)>d_1(3),\quad d_2(4)>d_0(4).
\end{equation}
Here and below, unadorned reciprocal indices refer to the canonical cut chart of \eqref{mix:physical-reciprocals}. The contact floor \eqref{hc:cost-floor} gives
\begin{equation}\label{mix:r02-contact}
 c_i>W(q_i)\quad(i=2,3),\qquad
 W(x)={x\log x\over x-1}-1.
\end{equation}
For arc 2, \eqref{mix:r02-cheap} also permits Theorem~\ref{hc:prefix650}. At its exit the switching function $q_{\rm sw}=w\chi$ satisfies $q_{\rm sw}\ge0$ and $q_{\rm sw}'>kq_{\rm sw}$. Its derivative matches across the switch by \eqref{hc:C1-identity}. On the outgoing mode-$0$ side, \eqref{hc:variable-rate} gives $w'/w=k$ at the switch. Thus $\chi\ge0$ and $\chi'>0$ at the entrance to arc 3 on the same lift. Along arc 3 the derivative with respect to binary shear time is $\mu-u^2$, with $\mu=2K-3$ constant, and $u$ is strictly increasing. The first subsequent zero of $\chi$ can occur only at $u>\sqrt\mu>0$. At the actual tie at cut 4 this gives $d_0^{\rm nat}(4)>d_2^{\rm nat}(4)$, that is, $d_1(4)>d_0(4)$. The reflected argument, when $c_3\le13/20$, gives $d_0(2)>d_1(2)$. The proper source $[5,8]=1201$ gives $d_1(7)<d_0(7),d_2(7)$ by Lemma~\ref{fin:forward-force}.

Since $W'(x)=(x-1-\log x)/(x-1)^2>0$, the positive logarithm series at $13/5$ and $6$ gives
\[
 W(13/5)>{162628\over295245}>{109483\over200000},\qquad
 W(6)>{18405\over16807}>{109483\over100000}.
\]
Equations \eqref{mix:r02-cheap}--\eqref{mix:r02-contact} imply $q_2<13/5$ and $q_3<6$. Let $a_i$ be the active physical reciprocal on arc $i$. At cut 3, the third reciprocal is strictly larger than both $a_2,a_3$ by \eqref{mix:r02-crossing}. The Lorentz identity gives $a_i^2+2a_2a_3<1$ for both $i=2,3$. Compactness gives $a_2,a_3\ge1/15$. Since $(15/16)^2+2(15/16)/15=257/256>1$, both $a_2,a_3<15/16$. Along every positive physical arc,
\begin{equation}\label{mix:cost-log-upper}
 {dc\over d\log p}=a(a+b(p)+c(p))=1+a^2-b(p)c(p)<1+a^2.
\end{equation}
Consequently $c_3<(481/256)\log q_3$. If $q_3\le19/16$, then $\log x<(x-x^{-1})/2$ gives
\[
 c_3<{481\over512}\left({19\over16}-{16\over19}\right)
 ={50505\over155648}<{13\over40},
\]
contrary to \eqref{mix:r02-cheap}. Thus $q_3>19/16$ and $\alpha>u_3>105/256$.

With both active reciprocal bounds, $q_2q_3\le7/5$ would imply
\[
 c_2+c_3<{481\over256}\log(7/5)
 <{1443\over2240}<{13\over20},
\]
contrary to row 8. Thus $q_2q_3>7/5$.

We have proved the root inclusion
\begin{equation}\label{mix:r02-root}
 (r_0,r_1,r_2,r_3,r_4)\in
 [33/520,1]^2\times[5/13,1]\times[1/6,16/19]\times[33/520,1].
\end{equation}
The actual root uses outward binary endpoint enclosures of this exact box. Weakening its faces retains every strict candidate. The positivity guards in Lemma~\ref{mix:r02-chart}, the cap $q_i\le520/33$, \eqref{mix:r02-crossing}, the conditional cheap-centre comparison and the nine fees give only necessary intersections in this root. In particular the R03 density owner used next is not imported into R02.

\begin{proof}[Proof of Theorem~\ref{mix:thm:r02}]
The source analysis places a representative of every original
candidate in \eqref{mix:r02-root}, with the stated domain,
cost, state and full-closure conditions. The complete cover of this
root in Proposition~\ref{mix:finite-covers} uses closed splits and
conservative contractions. Each terminal gives a strict
contradiction to a necessary condition or an empty necessary
intersection. Lemma~\ref{eight:common-cover} therefore excludes
every such candidate.
\end{proof}
\space
\section{Exclusion of R03}\label{mix:r03}
\begin{theorem}\label{mix:thm:r03}
An original global minimizer cannot have eight positive maximal arcs
with native labels $\mathtt{02120120}$ and twisted continuation
$j_{i+8}=j_i+1$. This excludes the strict essential type R03,
with increments $\mathtt{--++++++}$.
\end{theorem}

The full signed frame gives \eqref{mix:signed}. We use the positive
chart of Lemma~\ref{mix:r03-chart} and the common physical cost
circuit \eqref{mix:numerators}--\eqref{mix:positive-cost}.
We next place the strict source in a closed root and specify the
necessary inequalities for the cover.

\subsection{Proper sources and the positive-tail root}\label{mix:r03-domain}
Table~\ref{mix:r03-fees} uses the same support and label conventions as Table~\ref{mix:r02-fees}. Every support has length five, six or seven, strictly less than eight.

\begin{table}[htbp]\centering\small
\begin{tabular}{@{}clcll@{}}\toprule
Row & Charged arcs & $b_j$ & Complete source & Ordered support; $\ell(j)$\\\midrule
0 & $1,2,3,4$ & $81/200$ & $012021$ & $[5,0]$; $1-j$\\
1 & $0,1,6,7$ & $81/200$ & $012021$ & $[5,10]$; $j+2$\\
2 & $2,3,4,5$ & $41/100$ & $1012012$ & $[1,7]$; $j+2$\\
3 & $0,5,6,7$ & $41/100$ & $1012012$ & $[9,3]$; $1-j$\\
4 & $4,5$ & $180483/500000$ & $012012$ & $[7,2]$; $-j$\\
5 & $5,6$ & $180483/500000$ & $012012$ & $[8,3]$; $1-j$\\
6 & $1,2,3$ & $2/5$ & $0120212$ & $[13,7]$; $2-j$\\
7 & $0,1,7$ & $2/5$ & $0120212$ & $[5,11]$; $j+2$\\
8 & $3,4,5$ & $47/400$ & $01201$ & $[6,2]$; $2-j$\\
9 & $4,5,6$ & $47/400$ & $01201$ & $[7,3]$; $-j$\\
10 & $5,6,7$ & $47/400$ & $01201$ & $[8,4]$; $1-j$\\
11 & $3,4,5,6$ & $1/3$ & $012012$ & $[7,2]$; $-j$\\
12 & $4,5,6,7$ & $1/3$ & $012012$ & $[8,3]$; $1-j$\\\bottomrule
\end{tabular}
\caption{The thirteen proper-source R03 inequalities. Their original strict lower bounds are weakened to closed necessary inequalities in \eqref{mix:necessary-values}.}\label{mix:r03-fees}
\end{table}

Rows 0--1 use Theorem~\ref{loc:S405}. For rows 2--3 the complete seven-source preserves the uncharged support arc as well as both caps; its $41/100$ estimate is Theorem~\ref{loc:Edge41}, also weaker than the own-four bound in Theorem~\ref{fee:four}. Rows 4--5 use Theorem~\ref{nx:forward-pair}. Rows 6--7 use Proposition~\ref{six:prop:tail}: the complete forward prefix supplies the negative entrance and entrance circle by Lemma~\ref{cmp:forward-prefix-forcing}, and the actual final return supplies its terminal circle. This $2/5$ payment is not a theorem about an isolated $0212$ string. Rows 8--10 use the complete capped estimate of Proposition~\ref{strip:forward1175}; rows 11--12 use Corollary~\ref{strip:pthird}. All retain their actual returns, contacts, boundary ties and one selected complete lift wherever means or switching rows are coupled. Independently established scalar consequences may be used together without identifying the covectors of different sources.

The complete proper source $[6,11]=201020$ has local density owner 3, physical arc 1. Theorem~\ref{new:sixdensity} therefore gives
\begin{equation}\label{mix:r03-density}
 d_0(1)^2<1/6,\qquad \delta/6-n_{10}^2>0.
\end{equation}
For cuts $i=4,5,6,7$, the proper source $[i-2,i+1]$ is a whole phase of $1201$. Lemma~\ref{fin:forward-force} says that its canonical second reciprocal is strictly smaller than both others.

These state comparisons improve the duration-product cuts to
\begin{equation}\label{mix:r03-products}
 q_4q_5>5/4,\qquad q_5q_6>5/4.
\end{equation}
We include the state estimate underlying this strengthening. In canonical coordinates a positive increment transports
\[
 \Phi_q(a,b,c)=\left((q-1)a+qb,{c-(q-1)a\over q},a\right),\qquad ab+ac+bc=1.
\]
Suppose $q,r\ge1$, $qr\le5/4$, all three successive states are positive, and their first and last second reciprocals are minimal. Then
\begin{equation}\label{mix:two-step-state}
 a^2<3/5,\qquad (a')^2<3/5,
\end{equation}
where $a'$ is active after the first step. To prove the first inequality, divide the initial coordinates by $b$, writing $t=a/b\ge1$, $s=c/b\ge1$. The last second-versus-third comparison gives
\[
 t\le T(q,p)={q(2p-q)\over q^2-2pq+2p},\qquad p=qr.
\]
The denominator is positive for $1\le q,r\le5/4$. For $q\le15/14$, differentiation in $p$ and then $q$ gives $t\le T(15/14,5/4)=30/19$, whence
\[
 a^2={t^2\over t+s+ts}\le{t^2\over1+2t}\le{900\over1501}<3/5.
\]
For $q\ge15/14$ one has $r\le7/6$. The other terminal comparison gives $s\ge Ht-J$, where
\[
 H=q^2/r^2-(q-1)^2,\qquad J=q^2(1-r^{-2}).
\]
Here $H-J\ge593/1372>0$, and $R(t)=t^2/[t+(1+t)(Ht-J)]$ is increasing for $t\ge1$: with $u=1-r^{-2}$, its derivative has positive numerator $2qt\{(1-qu)t-qu\}$, since $qu\le65/196<1/2$. At the endpoint,
\[
 a^2\le R(T)={(2r-1)^2\over3}\le{16\over27}<3/5.
\]
Finally $S(a,b,c)=(c,b,a)$ satisfies $\Phi_qS\Phi_q=S$. Applying the first inequality to the reversed two-step state chain, with multipliers $r,q$, gives the second inequality in \eqref{mix:two-step-state}. This is an identity for physical states and does not reverse or join local covectors.

If $qr\le5/4$, \eqref{mix:cost-log-upper} and \eqref{mix:two-step-state} now give
\[
 c_i+c_{i+1}<{8\over5}\log(qr)\le{8\over5}\log(5/4)
 <{9\over25}<{180483\over500000},
\]
contrary to the proper forward-pair fee. For pair $(4,5)$ the forcing cuts are 4 and 6 and its forward-pair source is $[2,7]=120120$; for pair $(5,6)$ they are 5 and 7 with source $[3,8]=201201$. These are whole phases of $012012$, with their charged local middle pair and both caps retained. This proves \eqref{mix:r03-products}.

The exact whole reflection for R03 is
\begin{equation}\label{mix:r03-reflection}
 j'_r=1-j_{2-r}=j_r\pmod3,\qquad q'_r=q_{2-r}.
\end{equation}
It preserves the full return and the $+1$ seam. For example the matrix
$C=\left(\begin{smallmatrix}-1&-3\\-1&1\end{smallmatrix}\right)$ satisfies $C^2=4I$, $CB_jC^{-1}=-B_{1-j}$ and $CRC^{-1}=R^{-1}$; hence the reversed inverse frame product again equals $R$. It fixes density owner 1 and permutes the thirteen fee rows in pairs $(0,1)$, $(2,3)$, $(4,5)$, $(6,7)$, $(8,10)$, $(11,12)$, with row 9 fixed. Thus we choose $q_3\le q_7$, or $A\le E$ in \eqref{mix:r03-positive}. Its exact symmetric identity is
\[
 F=(1+C)(1-AB-DE)-AEK>0.
\]
Since $K\ge1$, this yields
\[
 A^2K\le AEK<(1+C)(1-AB-DE)\le1+C=q_5^2\le Q^2,
\]
so $q_3^2<1+Q=553/33$. Lemma~\ref{eight:lem:cap} now proves the root inclusion
\begin{equation}\label{mix:r03-root}
 (r_3,r_4,r_5,r_6,r_7)\in
 [\sqrt{33/553},1]\times[33/520,1]^4.
\end{equation}
The first lower endpoint, like all rational endpoints, is rounded outward in the actual closed root. The exact positivity guards of Lemma~\ref{mix:r03-chart}, the cap, \eqref{mix:r03-density}, the four forward conditions, \eqref{mix:r03-products} and the thirteen rows are necessary cuts on this root. Their use does not require every formal point of the box to be physical.

\begin{proof}[Proof of Theorem~\ref{mix:thm:r03}]
The source analysis places a representative of every original
candidate in \eqref{mix:r03-root}, with the stated domain,
cost, state and full-closure conditions. The complete cover of this
root in Proposition~\ref{mix:finite-covers} uses closed splits and
conservative contractions. Each terminal gives a strict
contradiction to a necessary condition or an empty necessary
intersection. Lemma~\ref{eight:common-cover} therefore excludes
every such candidate.
\end{proof}
\space
\section{Exclusion of R04}\label{r04:section}
We read the R04 word in Table~\ref{ten:table} from its fourth increment.
Its two positive increments are separated. Minimizing the two positive
rectangles at fixed area preserves the complete signed frame and reduces
the cost to a function of four exterior reciprocals. This cost exceeds the
octagon cost throughout the open domain, and its continuous extension
takes the octagon cost at one explicit degenerate boundary point.

\begin{theorem}[Exclusion of R04]\label{r04:exclusion}
An original global minimizer cannot have nine essential positive arcs with native labels $021210102$ and increments $\mathtt{--+--+---}$, with the genuine extension $j_{i+9}=j_i+1$.
\end{theorem}

\subsection{Reciprocal coordinates and cost for nine arcs}\label{nine:source-interface}
For $W=(u,v)$ in the state domain of Section~\ref{fd:section}, put
\[
 d_0=\frac{1+u}{v},\qquad d_1=\frac{u^2+v^2-1}{2v},\qquad d_2=\frac{1-u}{v}.
\]
All three numbers are positive and satisfy $d_0d_1+d_1d_2+d_2d_0=1$. Conversely this relation, with positivity, gives $v=2/(d_0+d_2)$ and $u=(d_0-d_2)/(d_0+d_2)$. At a switch, two assigned active reciprocals $x,z$ determine the third as
\[
 \kappa(x,z)=\frac{1-xz}{x+z}>0.
\]
The assignment to the three control labels remains part of the state.

Let $r_i$ be the active reciprocal on complete arc $i$. Put $x=r_i$, $z=r_{i-1}$, $w=r_{i+1}$, $\alpha=\varepsilon_{i-1}$ and $\beta=\varepsilon_i$. The exact multiplier and physical duration are
\begin{equation}\label{nine:multiplier}
 q_i=\frac{(x+z)^\alpha(x+w)^\beta}{(1+x^2)^{(\alpha+\beta)/2}},\qquad \tau_i=x\log q_i>0.
\end{equation}
Indeed, in the mode-$0$ chart $d_0=x$ is constant, $d_1^{\rm out}+x=q_i(d_1^{\rm in}+x)$ and $d_2^{\rm out}+x=(d_2^{\rm in}+x)/q_i$. At either end use the assigned neighbor or $\kappa(x,\cdot)$; the identity $x+\kappa(x,z)=(1+x^2)/(x+z)$ gives all four cases of \eqref{nine:multiplier}. Positive endpoint states remain positive along the arc. These formulas give a necessary coordinate description of the actual trajectory, while full signed frame closure and real winding are retained separately.

Define the symmetric physical-cost potential
\begin{equation}\label{nine:potential}
 \phi(x,t)=3xt-1-(x^2+t^2)\log(x+t)
       +\frac{x^2\log(1+x^2)+t^2\log(1+t^2)}2,
\end{equation}
whose first derivative is
\[
 g_1(x,t)=2t+\frac{2xt}{x+t}-\frac{x}{1+x^2}
                   -x\log\frac{(x+t)^2}{1+x^2}.
\]
Integrating the physical density on one complete arc gives
\[
 x\{d_2^{\rm in}-d_2^{\rm out}+d_1^{\rm out}-d_1^{\rm in}\}-x^2\log q.
\]
At a switch from reciprocal $x$ to $t$, the paired endpoint
contributions equal $\varepsilon(3xt-1)$ by
$x[\kappa(x,t)-t]+t[\kappa(x,t)-x]=1-3xt$.
Summing the endpoint and logarithmic terms with
\eqref{nine:multiplier} around the twisted cycle gives
\[
 C=\sum_i\varepsilon_i\phi(r_i,r_{i+1}).
\]
We use the same reciprocal and cost conventions in the R05 and R06
comparisons.

\subsection{The strict physical source and the four exterior coordinates}
Let $\mathcal S_{04}$ be the class of original global-minimizing trajectories of this type, with their actual reciprocal states, original full maximizing lift (including its state covector), complete signed endpoints and real winding. The nine physical arcs are all strictly positive. In native order put
\begin{equation}\label{r04:native}
 (r_0,\ldots,r_8)=(A,B,Q_1,R_1,C,Q_2,R_2,D,E),
 \qquad r_{i+9}=r_i.
\end{equation}
Here $B$ denotes a reciprocal coordinate, not the essential count. The multiplier formula \eqref{nine:multiplier} gives
\begin{equation}\label{r04:durations}
\begin{gathered}
 q_0=\frac{1+A^2}{(A+B)(A+E)}=:p,\quad
 q_1=\frac{1+B^2}{(A+B)(B+Q_1)},\\
 q_2=\frac{Q_1+R_1}{B+Q_1},\quad
 q_3=\frac{Q_1+R_1}{C+R_1},\quad
 q_4=\frac{1+C^2}{(C+R_1)(C+Q_2)},\\
 q_5=\frac{Q_2+R_2}{C+Q_2},\quad
 q_6=\frac{Q_2+R_2}{D+R_2},\\
 q_7=\frac{1+D^2}{(D+R_2)(D+E)},\quad
 q_8=\frac{1+E^2}{(D+E)(A+E)}=:q.
\end{gathered}
\end{equation}
Since $\tau_i=r_i\log q_i>0$, every displayed multiplier is greater than one. In particular
\begin{equation}\label{r04:widths}
 X_1=Q_1-C>0,\quad Y_1=R_1-B>0,\quad
 X_2=Q_2-D>0,\quad Y_2=R_2-C>0.
\end{equation}
Write $I_i=X_iY_i$. Set
\begin{equation}\label{r04:four-coordinates}
\begin{gathered}
 H=\frac{(1+B^2)(1+D^2)}{(A+B)(D+E)},\qquad M=p^2+q^2-1,\\
 C^2=\frac{pqH}{M}-1,\\
 I_1=\frac{(B+C)^2}{2}
 -\frac{(1+B^2)^2(q^2-1)}{2(A+B)^2M},\\
 I_2=\frac{(C+D)^2}{2}
 -\frac{(1+D^2)^2(p^2-1)}{2(D+E)^2M}.
\end{gathered}
\end{equation}
These are identities at the same physical candidate. Their derivation from the signed frame and the fixed-area variation is given in Section~\ref{r04:frame-reduction}. All denominators are positive since $p,q>1$, so $M>1$. The positive root determines $C$ without a branch choice.

For positive $a,b,X,Y$ define
\begin{equation}\label{r04:rectangle}
\begin{gathered}
 k(u,v)=2-\frac{2uv}{(u+v)^2},\qquad
 \Gamma(a,b;X,Y)=\int_a^{a+Y}\int_b^{b+X}k(u,v)\,dv\,du,\\
 \mathcal E(a,b,I)=\min_{X>0}\Gamma(a,b;X,I/X).
\end{gathered}
\end{equation}
Lemma~\ref{r04:shape-unique} proves existence and uniqueness of this minimum for every $a,b,I>0$. It also proves that the unique stationary shape is its global minimum. The actual two shapes in \eqref{r04:widths} are stationary by the fixed-area physical variation. Thus they realize $\mathcal E(B,C,I_1)$ and $\mathcal E(C,D,I_2)$, respectively.

Define $\Omega_{04}$ as the following open domain of $z=(A,B,D,E)$: the four coordinates are positive, $p,q>1$, the positive root in \eqref{r04:four-coordinates} exists, $I_1,I_2>0$, and the two unique minimizing shapes satisfy
\begin{equation}\label{r04:omega}
\begin{gathered}
 X_1<\frac{1+B^2}{A+B}-B-C,\qquad
 Y_2<\frac{1+D^2}{D+E}-C-D,\\
 (B+C+Y_1)(C+D+X_2)<1+C^2.
\end{gathered}
\end{equation}
The three inequalities follow exactly from $q_1,q_7,q_4>1$. Consequently every member of $\mathcal S_{04}$ projects into $\Omega_{04}$.\label{r04:source-bridge} This is an inclusion into a larger domain: none of the original winding, state or maximizing qualifications is replaced by a system of weak algebraic equations.

\subsection{The exact physical cost}
The potential in \eqref{nine:potential} has mixed derivative $\phi_{uv}=k$. For $s=(A,B,C,D,E)$ with cyclic indices define
\begin{equation}\label{r04:cost-function}
 C_5=-\sum_{i=0}^4\phi(s_i,s_{i+1}),\qquad
 G=C_5+\mathcal E(B,C,I_1)+\mathcal E(C,D,I_2).
\end{equation}
The fundamental theorem of calculus gives
\begin{align*}
 \Gamma(B,C;X_1,Y_1)
 &=\phi(R_1,Q_1)-\phi(B,Q_1)-\phi(R_1,C)+\phi(B,C),\\
 \Gamma(C,D;X_2,Y_2)
 &=\phi(R_2,Q_2)-\phi(C,Q_2)-\phi(R_2,D)+\phi(C,D).
\end{align*}
Adding these to $C_5$ cancels precisely the compressed edges $(B,C)$ and $(C,D)$. Hence the once-counted physical nine-arc cost is
\begin{equation}\label{r04:cost-bridge}
\begin{split}
 C_{04}=G={}&-\phi(A,B)-\phi(B,Q_1)+\phi(Q_1,R_1)-\phi(R_1,C)\\
 &-\phi(C,Q_2)+\phi(Q_2,R_2)-\phi(R_2,D)\\
 &-\phi(D,E)-\phi(E,A)
\end{split}
\end{equation}
at every member of $\mathcal S_{04}$. In particular no additional $3I_i/2$ is charged after the rectangle integral has been inserted.

Put
\begin{equation}\label{r04:constants}
\begin{gathered}
 s=\sqrt2,\quad b_0=(5+4s)^{-1/2},\quad a_0=(1+s)b_0,\quad
 z_0=(a_0,b_0,b_0,a_0),\\
 m_0=\frac{16+8s-(6+4s)\log2}{5+4s}.
\end{gathered}
\end{equation}
The exact identity
\begin{equation}\label{r04:budget-identity}
 m_0=J_8=\frac{2(8-4\sqrt2-\log2)}{2\sqrt2-1}
\end{equation}
follows by cross multiplication, using $s^2=2$: both cross products are $16+24s-(10+8s)\log2$.

\begin{theorem}[Strict cost bound on the four-dimensional domain]\label{r04:global-bound}
For every $z\in\Omega_{04}$, the function in \eqref{r04:cost-function} satisfies $G(z)>m_0$.
\end{theorem}

The proof occupies Appendix~\ref{app:r04-source}. A rectangle comparison and shape estimates give lower bounds for the cost. The tail and side estimates place every possible counterexample in $[0,5]\times[0,1]^2\times[0,5]$, and symmetry allows us to take $A\le E$. This box lies in the single root of the $31$-run interval cover. Near $z_0$ the cover uses the analytic estimate on the cube of radius $11/10000$. Section~\ref{reader:r04} gives the verification procedure.

\begin{proof}[Proof of Theorem~\ref{r04:exclusion}]
A member of $\mathcal S_{04}$ would project into $\Omega_{04}$ and, by \eqref{r04:cost-bridge}, have $C_{04}=G$. The competitor budget \eqref{fd:budget}, Theorem~\ref{r04:global-bound} and \eqref{r04:budget-identity} give
\[
 J_8\ge C_{04}=G>m_0=J_8,
\]
a contradiction. Thus $\mathcal S_{04}=\varnothing$ and its qualified global closure $\overline{\mathcal S_{04}}$ is empty as well.
\end{proof}

At $z_0$ the continuous extension has $p=q=1$, $C=b_0$, $I_1=I_2=2b_0^2$ and $X_i=Y_i=\sqrt2b_0$. Its cost is $m_0$. These zero-duration arcs place it on the boundary of the relaxed domain, outside $\mathcal S_{04}$. The strict pointwise bound requires no uniform positive gap from $J_8$ on the open domain and does not classify arbitrary weak nine-tuples as original minimizing trajectories.
\space
\section{Exclusion of R05}\label{nine:section}
The two returns and the central own pair determine the source
coordinates and the necessary equations for this type. A closed
cover reduces the problem to a fixed collection of cells, where we
exclude those necessary values.

\begin{theorem}[Exclusion of R05]\label{nine:thm:r05}
An original global minimizer cannot have nine essential positive arcs with native labels $021210120$ and increments $\mathtt{--+--++++}$, with the genuine extension $j_{i+9}=j_i+1$.
\end{theorem}
The proof first identifies these original minimizers with the candidates to which the source reduction applies. It then excludes their necessary values on a complete finite cover. It does not assert that every positive nine-tuple, or every zero of an auxiliary residual system, is absent.

\subsection{The original R05 source}\label{nine:r05-source}
Write the actual active reciprocals of R05 in native order as
\begin{equation}\label{nine:native}
 (r_0,\ldots,r_8)=(E,a,Q,R,b,Z,c,y,f),\qquad r_{i+9}=r_i.
\end{equation}
Here $Q$ is a reciprocal coordinate, not the multiplier bound $\Qcap$. The outer mean owners are $E,Z$, the density owners are $a,b$, and the central own pair is $Q,R$.

\begin{lemma}[Qualification of the R05 sources]\label{nine:qualification}
Every original R05 global minimizer admits one jointly normalized full maximizing lift on the actual unfolding, with the matrix twist and both state-covector twists. Its proper returned source is native $[-1,6]$, with whole phase $+1$ equal to $01020212$. On this source the density owners are native $1,4$, the exterior mean owners are native $0,5$, and the own pair is native $2,3$. Moreover $a^2,b^2<1/6$.
\end{lemma}
\begin{proof}
Take two cuts strictly inside native arcs $7$ and $8$, in that order. Their overlap contains the next twisted translate $[9,14]=102021=1|0202|1$ of the complete capped source $[0,5]=021210=0|2121|0$. Its four-core has the two internal switching roots and both positive caps required by Lemma~\ref{pc:rigidity}. Between the two cuts is the actual switch from label $2$ to label $0$. Proposition~\ref{pc:two-cuts} therefore applies: first it proves normality, then identifies the whole jointly normalized lifts on the overlap and proves their matrix and state-covector twists. This supplies the full lift, rather than an identification of independently selected determinants.

The twisted continuation gives $j_{-1}\cdots j_6=20212101$, hence the stated phase-$+1$ eight-source. It is proper since $8<9$, so Lemma~\ref{as:placement} retains it in a genuine fundamental interval. In its local numbering the four-core is $[2,5]$, the density owners are $2,5$, and the own pair is $3,4$; subtracting one gives exactly the native dictionary in the statement. Both actual returns $[-1,1]$ and $[4,6]$, and both exterior caps, remain present. Section~\ref{pf:M034} supplies the chronological root jets and density profiles on this same lift.

For the sharper density bounds use the proper six-source $[7,12]=201020$, whose local arc $3$ is physical arc $1$, and the reversed source $[7,6,5,4,3,2]$ with the whole reflected label map $j\mapsto1-j$. The latter is also $201020$, and its local arc $3$ is physical arc $4$. Theorem~\ref{new:sixdensity} gives $a^2,b^2<1/6$. Each six-source is placed on its own genuine interval if necessary; only these physical state conclusions are exported. Their covectors are not spliced into the common returned-source lift.
\end{proof}

Define the actual uniform reciprocal-shear mean by
\begin{equation}\label{nine:mean}
 M(x;k,z)=\int_0^1\left(\frac{x-((1-t)k+tz)}{x+((1-t)k+tz)}\right)^2\,dt.
\end{equation}
On an actual positive window this is the uniform average over $[k,z]$. It is not a physical-time average. The source identifies three values at the same candidate:
\begin{equation}\label{nine:three-means}
 \mu=M(Q;a,R)=M(R;b,Q),\qquad \nu_a=M(E;a,f),\qquad \nu_b=M(Z;b,c).
\end{equation}
In particular, the two exterior functions retain their different owners even when their values happen to coincide.

\subsection{The source reduction and its closed domain}\label{nine:root-domain}
We use three previously established inputs from the complete-source proof collection described in Section~\ref{verify:r05-source-inputs}. The first consists of the continuous source, calibration, cap and domain comparisons; we call it D1. Its statements retain the actual owners, selected roots, complete lift and permitted limits. The second, D2, is the closed-root, continuation and half-$F$ covering theorem. The third, D3, derives the physical necessary-value equations and their smooth auxiliary extensions. Section~\ref{verify:r05-source-inputs} gives their mathematical interfaces; the electronic collection supplies the original physical-graph texts, the $29$ domain sources, both complete coefficient proofs, and the finite covers. D1 is a geometric implication and D2 is a covering theorem; neither follows from the final finite table being empty.

Set
\begin{equation}\label{nine:shape-coordinates}
 H=Q+R,\quad p=\frac{R-a}{H},\quad t=\frac{Q-b}{H},\quad s=p+t,\quad \eta_s=p-t,
 \quad e=\frac{\eta_s}{s},\quad \delta=\frac{R-Q}{H}.
\end{equation}
The positive inner windows make $s>0$. The actual calibration equality in \eqref{nine:three-means} selects its original positive branch; no alternate root of its cleared quadratic is admitted. In these coordinates
\begin{equation}\label{nine:shape-reconstruction}
 \begin{gathered}
 R=\tfrac H2(1+\delta),\quad Q=\tfrac H2(1-\delta),\quad
 a=\tfrac H2(1-s+\delta-\eta_s),\quad b=\tfrac H2(1-s-\delta+\eta_s),\\
 D=\frac{(a+b)^2}{2}-(R-a)(Q-b)
   =H^2\left(\frac12-s+\frac{s^2+\eta_s^2}{4}\right)>0.
 \end{gathered}
\end{equation}
The symbol $D$ in this section denotes this area parameter. With $h_k=(1+k^2)/2$, put $F=D/(h_ah_b)$ and $r=(h_a-h_b)/D$.

\begin{proposition}[Retained source-to-domain input]\label{nine:source-reduction}
Every original global minimizer of type R05, qualified as in Lemma~\ref{nine:qualification}, either contradicts one of the retained source comparisons or has its coordinates $(s,e,H,c,f)$ in the closed box
\begin{equation}\label{nine:B0}
 B_{05}=\left[\frac{49}{125},\frac{2341}{4000}\right]
 \times\left[-\frac{327}{1000},\frac{327}{1000}\right]
 \times\left[\frac{13}{10},2\right]
 \times\left[\frac3{20},\frac43\right]^2.
\end{equation}
The same candidate retains its complete caps, original normal full lift, all maximizing conditions, selected calibration and return roots, actual means, signed endpoints and real winding. The accepted closed-root and continuation covers reduce it further to at least one of the fixed $46\,383$ complete original cells.
\end{proposition}
\begin{proof}[Identification and composition of the retained inputs]
The source comparisons in D1 apply to the actual owners identified in Lemma~\ref{nine:qualification}. The complete state compactness theorem, transported by whole cyclic phases, bounds every endpoint reciprocal in $[1/15,400/99]$ and each multiplier by $\Qcap$. This includes the third components at the end switches. The two proper six-sources give the sharpened density bounds. The returned eight-source supplies its two genuine return jets, the common determinant and the unsquared signed sum of the two constants.

The complete source $[9,14]=102021$ in the overlap of Lemma~\ref{nine:qualification} is the twisted translate of native $[0,5]$. On this same lift, the capped four-core argument in Corollary~\ref{loc:four-calibration} first gives $\Gamma=\det\Lambda>0$. The binary invariant $K=(3+\mu)/2>0$ then gives $\operatorname{tr}(\Lambda X)>0$ at an actual point by Corollary~\ref{cal:one-point}, so the original harmonic bound gives $\Gamma<H_*$. This uses the original global budget \eqref{fd:budget}, with each physical arc charged once.

The two-cap and selected-calibration comparisons give
\[
 \frac{49}{125}<s<\frac{2341}{4000},\qquad |\eta_s|<\frac{327}{1000}s,
 \qquad \frac{13}{10}<H<2.
\]
Their shape bounds are whole-domain comparisons, including both signs of imbalance and balance. The lower scale bounds come from the actual own-pair payment and the two-cap floor, not a scale chosen independently of the candidate. In particular $a,b>3/20$ and $a+b>5/8$. Positive returns give $f>a$ and $c>b$; the selected-return and state-product comparisons give $c,f<4/3$. These are the five required coordinate inclusions. Replacing their strict endpoints by the closed endpoints of \eqref{nine:B0} only enlarges the domain.
The actual central coordinates lie in the domain \eqref{r05:central-coefficient-domain} of both coefficient comparisons. The two proper forward pairs and the same-source frame bound establish this inclusion in Section~\ref{verify:r05-source-inputs}.

The additional necessary conditions used during evaluation retain those same premises. They include the selected calibration cone, $|r|\le3/7$, $F>1/20$, the whole-sector exclusion $F\ge1/2$, the return collars $f-a,c-b>1/1000$, and the same-side owner and state-product bounds. The return collars are proved from the undivided frame equations on the full closed thickness $[0,1/1000]$; they are not inferred by extending a zero-width derivative test. The owner bounds use their own selected mean roots and true leaves. Thus a zero exterior leaf is not confused with a zero return. All omitted maximizing or winding conditions remain conjuncts even when an evaluator uses only a subset of them.

The matching closed covers in
Section~\ref{reader:r05:reduction} take precisely the retained records
from each stage to the next, preserve the same qualified candidate,
and leave \(46383\) complete original cells. This proves the last
inclusion of the proposition.
\end{proof}

\subsection{The complete caps and necessary value equations}\label{nine:system}
The finite argument works with an outer box of the following thirteen values at one candidate:
\begin{equation}\label{nine:lift}
 z=(a,b,Q,R,c,y,f,E,Z,D,\mu,\nu_a,\nu_b).
\end{equation}
The three parameters $\theta=(s,e,H)$ range over their entire original cell intervals. They are fixed when differentiating in $z$. On the auxiliary open domain
\[
 \Omega_{05}=(0,\infty)^{10}\times\mathbb R^3
\]
all formulas below are smooth. The last three coordinates are unrestricted real variables for this analytic extension. They still represent the actual means \eqref{nine:three-means} at a physical candidate. An auxiliary coordinate is not an additional physical degree of freedom.

Put $K=(3+\mu)/2$, $K_a=(3+\nu_a)/2$ and $K_b=(3+\nu_b)/2$. The original common signed lift gives
\begin{equation}\label{nine:lambda}
 \lambda_a=-\frac{KD}{2}-\frac{h_a^2K_a-h_b^2K_b}{2D},\qquad
 \lambda_b=-\frac{KD}{2}+\frac{h_a^2K_a-h_b^2K_b}{2D}.
\end{equation}
Indeed its unsquared equations are $\lambda_a+\lambda_b=-KD$ and $Kh_a^2K_a-\lambda_a^2=Kh_b^2K_b-\lambda_b^2$. Since $KD>0$, subtracting the latter equations determines the difference and gives \eqref{nine:lambda}. There is no independent square-root choice on either side. The determinant is
\[
 \Gamma=\frac K2(h_a^2K_a+h_b^2K_b)-\frac{K^2D^2}{4}
       -\frac{(h_a^2K_a-h_b^2K_b)^2}{4D^2}.
\]
Both negative squares remain in every use of this expression; $\Gamma$ is not a physical cost.

For one side $(k,P,O,\nu,\lambda)$ define
\[
 \upsilon_k=\frac{k(k+P)}2,\qquad \omega_k=\frac{kh_k}{k+O},\qquad \rho_k=1-\frac{k^2}{\omega_k},
\]
\[
 I_k=(1-\mu)(\omega_k-\upsilon_k)-2k^2\log(\omega_k/\upsilon_k)
                   +k^4(1/\upsilon_k-1/\omega_k).
\]
The complete cap residual is
\begin{equation}\label{nine:cap}
 \begin{split}
 \mathcal P(k,P,O,\nu,\lambda)={}&-h_k^2\frac{3+\nu}{2}-2\lambda \omega_k-K\omega_k^2+\frac12\\
 &+\frac{(\rho_k^2-\mu)\omega_k^2}{2k^2}
       -\frac{(1+k^2)\omega_kI_k}{2k^2}.
 \end{split}
\end{equation}
Use $(a,Q,E,\nu_a,\lambda_a)$ and $(b,R,Z,\nu_b,\lambda_b)$ to obtain $\mathcal P_a,\mathcal P_b$. The logarithmic integral term is part of each cap. For example, replacing $\omega_k$ by a variable $w$ gives $I_k'(w)=(1-k^2/w)^2-\mu$ and $I_k(\upsilon_k)=0$, which recover that term from the original primitive without truncation.

For the central functions set
\begin{equation}\label{nine:central-functions}
 \begin{gathered}
 U=cf-1+(c+f)y,\quad Y=1+y^2,\quad V=U+2Y=cf+(c+f)y+1+2y^2,\\
 \mathcal F=\frac{4U+2U^2/Y}{(1+c^2)(1+f^2)},\quad
 v_a=\frac{2y(c+y)(1+f^2)}V,\quad
 v_b=\frac{2y(f+y)(1+c^2)}V.
 \end{gathered}
\end{equation}
Here $\mathcal F_i$ denotes its partial derivative with respect to $i$. With $d=h_a-h_b$, the first ten necessary equations at the minimizing candidate are
\begin{equation}\label{nine:first-ten}
 \begin{gathered}
 D-(a+b)^2/2+(R-a)(Q-b)=0,\\
 M(Q;a,R)-\mu=M(R;b,Q)-\mu=0,\\
 M(E;a,f)-\nu_a=M(Z;b,c)-\nu_b=0,\\
 \mathcal F-D/(h_ah_b)=0,\\
 A_a:=2E(f-a)-(1-d/D)(1+a^2)+v_a=0,\\
 A_b:=2Z(c-b)-(1+d/D)(1+b^2)+v_b=0,\\
 \mathcal P_a=\mathcal P_b=0.
 \end{gathered}
\end{equation}
The first equation defines $D$ as in \eqref{nine:shape-reconstruction}. On the common lift of Lemma~\ref{nine:qualification}, integration between successive switching zeros on each of the two inner arcs gives the same $\mu$ by Lemma~\ref{loc:binary-normality}. The outer means $\nu_a,\nu_b$ are defined in \eqref{nine:three-means}. Full signed frame closure gives $\mathcal F=D/(h_ah_b)$ and $A_a=A_b=0$. The nonzero central off-diagonal entry and determinant one recover the remaining matrix entry on the strict source graph. The two actual returns on the same lift give $\mathcal P_a=\mathcal P_b=0$ by \eqref{nine:cap}. The actual frame path supplies its real winding.

We give the remaining three rows explicitly. Differentiation
takes place on the physical source graph. The shared potential
\eqref{nine:potential} gives the once-counted cost
\begin{equation}\label{nine:cost}
 \begin{split}
 C_{05}={}&-\phi(E,a)-\phi(a,Q)+\phi(Q,R)-\phi(R,b)-\phi(b,Z)\\
          &+\phi(Z,c)+\phi(c,y)+\phi(y,f)+\phi(f,E).
 \end{split}
\end{equation}
In particular the closing edge occurs once, and the negative terms cannot be discarded or optimized independently.

The retained source chain-rule theorem D3 differentiates this complete physical cost on its strictly positive whole-period frame graph, before restricting the resulting values to the actual means and complete caps. In its local calibrated family, the selected calibration and central roots vary smoothly at the actual positive point. Three signed frame entries, including a nonzero central off-diagonal entry, and determinant one give the complete frame equation. The nine positive durations and positive state coordinates persist in a neighbourhood. The resulting frame paths have fixed signed endpoints and the same real winding by continuity. Lemma~\ref{fw:balanced-reconstruction} therefore gives genuine Reinhardt competitors throughout this neighbourhood. Only geometric admissibility is needed for the nearby competitors. The original mean--frame proof and the physical chain rule are supplied in the electronic source collection, with the interface summarized in Section~\ref{verify:r05-source-inputs}. At the actual minimizing point its three central necessary values are $R_i=0$, $i=y,c,f$, where
\begin{equation}\label{nine:central-rows}
 \begin{gathered}
 Q_c=h_ah_b\left[K-\frac{3d^2}{2D^2}
           +\frac{d(h_b\nu_b-h_a\nu_a)}{2D^2}\right],\\
 B_y=g_1(y,c)+g_1(y,f)-\frac{3+\nu_a}{4}(v_a)_y-\frac{3+\nu_b}{4}(v_b)_y,\\
 B_c=g_1(c,Z)+g_1(c,y)-\frac{3+\nu_a}{4}(v_a)_c-\frac{3+\nu_b}{4}(v_b)_c-ZK_b,\\
 B_f=g_1(f,E)+g_1(f,y)-\frac{3+\nu_a}{4}(v_a)_f-\frac{3+\nu_b}{4}(v_b)_f-EK_a,\\
 R_i=B_i-\mathcal F_iQ_c.
 \end{gathered}
\end{equation}
All central derivatives here are ordinary partial derivatives of the explicit rational functions \eqref{nine:central-functions}. The $c$ equation uses $Z,\nu_b$ and the $f$ equation uses $E,\nu_a$. The signed contrast in $Q_c$ and both owner terms are essential. The source derivation of these physical necessary values is an input distinct from the analytic extension used below: freezing $\mu,\nu_a,\nu_b$ in \eqref{nine:cost} would not prove \eqref{nine:central-rows}.

Three further zero rows retain the original parameters:
\begin{equation}\label{nine:chart-rows}
 Q+R-H=0,\qquad a+b-(1-s)(Q+R)=0,\qquad
 (1-e)(R-a)-(1+e)(Q-b)=0.
\end{equation}
Thus \eqref{nine:first-ten}--\eqref{nine:chart-rows} give sixteen necessary equalities. They are used as a necessary subsystem of the original physical problem, with all original source conditions still retained.

\subsection{Necessary contractions and exact rational separation}
Section~\ref{nine:ambient} augments the sixteen physical equations
to eighteen or twenty-four and gives the whole-box derivative
identities used for their enclosures.\label{nine:separation}
Fix one complete original cell. The accepted value constructor starts from its full five-dimensional endpoints, the selected calibration branch, complete frame, actual means and original source bounds. It either proves a necessary contradiction or encloses the same candidate in a positive box $X\subset\Omega_{05}$. In the latter case it initially retains the entire parameter box for $\theta=(s,e,H)$. A later closed chart split restricts this box only by covering it with both closed children. A subsequent contraction always intersects the old box. It cannot replace one candidate's mean, root or owner by a separately favorable choice.

Here is the inclusion principle used for the first-order contractions. For a residual $T$ and fixed $\theta$, with center $m\in X$,
\[
 0=T(z;\theta)=T(m;\theta)
       +\sum_j\left(\int_0^1\partial_jT(m+t(z-m);\theta)\,dt\right)(z_j-m_j)
\]
at every candidate zero. Whole-box derivative intervals therefore give necessary linear interval rows. If a coefficient omits zero, isolating that coordinate and intersecting its existing interval retains the same zero. Sequential use remains valid because each previous step has preserved it. A sign-separated value residual or a strictly empty necessary intersection gives a contradiction. Failure to certify a denominator, a disjoint pair of purported enclosures of the same function, or an unresolved calculation is not an exclusion.

For the stronger rows, twice applying the fundamental theorem of calculus gives
\begin{equation}\label{nine:taylor}
 T(z;\theta)=T(m;\theta)+\nabla T(m;\theta)\cdot\xi
 +\int_0^1(1-t)\xi^t\nabla^2T(m+t\xi;\theta)\xi\,dt,
 \qquad \xi=z-m.
\end{equation}
The line segment remains in $X$, and every parameter value in the current retained parameter box is quantified; initially this is the full original parameter box. If $\mathbb H_{ij}$ encloses each Hessian entry throughout this product domain, the remainder is enclosed by
\begin{equation}\label{nine:remainder}
 \mathcal R_T=\frac12\sum_i\mathbb H_{ii}[\xi_i^2]
                  +\sum_{i<j}\mathbb H_{ij}[\xi_i][\xi_j].
\end{equation}
The square interval in the diagonal term has lower endpoint zero. The factor $1/2$ on the diagonal and its cancellation on the symmetric off-diagonal terms follow directly from the integral in \eqref{nine:taylor}; no common intermediate point for all entries is asserted.

Put $\rho_j=\max(m_j-\inf X_j,\sup X_j-m_j)$ and $\xi_j=\rho_ju_j$, $|u_j|\le1$. A zero-radius coordinate uses $u_j=0$, without division. Suppose $\mathbb T_0$ encloses $T(m;\theta)$, $\mathbb J_j$ encloses $\partial_jT(m;\theta)$, and write
\[
 A_j=\rho_j\operatorname{mid}\mathbb J_j,\qquad
 E_J=\sum_j\rho_j\operatorname{rad}\mathbb J_j.
\]
Then a common zero satisfies the two rational inequalities
\begin{equation}\label{nine:two-rows}
 A\cdot u\le-\inf\mathbb T_0-\inf\mathcal R_T+E_J,\qquad
 -A\cdot u\le\sup\mathbb T_0+\sup\mathcal R_T+E_J.
\end{equation}
The same construction with whole-box first derivatives and zero explicit remainder gives the first-order necessary rows. Eighteen equalities supply thirty-six such rows, and the original cost condition $C_{05}-C_*\le0$ supplies one further row, used once. Adding the two second-order rows for each of eighteen or twenty-four equalities gives respectively $73$ or $85$ rows. These are linear necessary inequalities at a candidate, rather than a positive-semidefiniteness requirement on the Hessian of its physical cost.

\begin{lemma}[Rational support and coordinate contraction]\label{nine:support}
Suppose the necessary rows are $a_j\cdot u\le b_j$, with $|u_i|\le1$. For nonnegative rational weights $\lambda_j$, put $a=\sum_j\lambda_ja_j$ and $b=\sum_j\lambda_jb_j$. If
\begin{equation}\label{nine:strict-margin}
 -\sum_i|a_i|-b>0,
\end{equation}
there is no candidate in the box. Otherwise, if $a_k>0$,
\[
 u_k\le\frac{b+\sum_{i\ne k}|a_i|}{a_k};
\]
if $a_k<0$ the same quotient is a lower bound.
\end{lemma}
\begin{proof}
Nonnegative combination gives $a\cdot u\le b$, whereas $a\cdot u\ge-\sum_i|a_i|$. This proves the strict contradiction. Isolating $a_ku_k$ and using $-\sum_{i\ne k}a_iu_i\le\sum_{i\ne k}|a_i|$ gives the contraction, with the inequality reversed when dividing by a negative coefficient.
\end{proof}
Every adopted support recomputes its nonlinear values, derivative coefficients, nonnegative rational weights and strict margin. Floating-point optimization proposes weights only. Rational arithmetic checks their actual use; outward floor and ceiling round the reconstructed coordinate endpoints. A zero margin is retained. The nonlinear residuals may exclude all auxiliary zeros in a box, but cost and known-domain filters exclude only the qualified physical subset to which their theorems apply. The proof uses the latter, narrower statement.

\subsection{Arithmetic and the complete original-cell cover}\label{nine:finite-cover}
The final original-cell arithmetic uses arbitrary-precision integers and intervals with endpoints on the scale $2^{-128}$. The earlier whole-sector half-$F$ theorem has its own $2^{-160}$ interval implementation and $96$-term logarithmic bound; it is an input to the overlay and is not re-evaluated by the final-cell kernel. Rational inputs, multiplication and division are rounded outward; square roots use integer square roots with an upward correction when needed. Positive logarithms are reduced by exact powers of two to $[1,2]$. For $z=(x-1)/(x+1)$ the eighty-term expansion
\[
 \log x=2\sum_{j=0}^{79}\frac{z^{2j+1}}{2j+1}+R,\qquad
 |R|\le\frac{2|z|^{161}}{161(1-|z|^2)}
\]
encloses the complete tail, with the endpoint range of $z$ itself rounded outward. The same method encloses $\log2$. The calibrated roots and actual means retain their separately proved convergent-series tails and positive-root interval bounds; all operations preserve the same actual values.

One can also enclose the complete mean and its first derivatives by writing
\[
 q=\frac{z-k}{2x+k+z},\qquad r=\frac{2x-k-z}{2x+k+z},\qquad
 M(x;k,z)=r^2+\sum_{n\ge1}q^{2n}\left[r^2+\frac{4nr}{2n+1}+\frac{2n-1}{2n+1}\right].
\]
This follows from the uniformly convergent expansion of $(1+tq)^{-2}$ on $[-1,1]$. If $|q|\le v<1$ and $|r|\le R_0$, the tails after $n=40$ of the value, $r$ derivative and $q$ derivative are bounded respectively by
\[
 \frac{(1+R_0)^2v^{82}}{1-v^2},\qquad
 \frac{2(1+R_0)v^{82}}{1-v^2},\qquad
 \frac{2(1+R_0)^2v^{81}(41-40v^2)}{(1-v^2)^2}.
\]
They can be intersected with the independent compact-integral enclosures because both contain the same functions. An empty intersection between two such enclosures signals an unsafe evaluation, not absence of a physical candidate.

The fixed original-cell list has exact identities consisting of its seed, continuation address, five-dimensional endpoints and original record. Their fixed evaluation profiles partition these identities as follows:
\begin{center}
\begin{tabular}{lr}\toprule
Evaluation profile&Complete original cells\\\midrule
Value and first-order contractions&$17\,766$\\
Eighteen-residual support&$12\,290$\\
Twenty-four-residual support&$16\,327$\\\midrule
Total&$46\,383$\\\bottomrule
\end{tabular}
\end{center}
Each original identity occurs once in this partition. The final verification reconstructs each complete original record and its mathematical contradiction using the appropriate accepted equations and arithmetic. Fifteen cells use separate complete-original support plans: thirteen use their full nonlinear calculations, and two use complete closed internal covers. Each such plan is checked on its entire original cell; an inherited acceptance label is not a mathematical terminal.

Each internal subdivision covers the current enclosing box. An interior grid point $l<m<h$ splits one coordinate interval $[l,h]$ into $[l,m]$ and $[m,h]$. For a native-coordinate split, both children retain the entire current parameter box. For a chart-coordinate split, both retain the entire current native enclosure, while their closed parameter intervals cover the current parameter box. If both descendants are excluded, so is the parent. Otherwise the original cell remains unresolved. Depth limits, failed guards, unknown labels, unsafe arithmetic and missing branches never count as success. Internal leaf counts are not added to the number of original cells.

\begin{lemma}[Exhaustion of the R05 original cells]\label{nine:exhaustion}
Under the source implications of Proposition~\ref{nine:source-reduction}, the union of the fixed $46\,383$ original cells contains no original R05 minimizing candidate.
\end{lemma}
\begin{proof}
The final ordered cell list has positions $0,\ldots,46\,382$. At each position the verification binds the recorded source identity and the complete original cell, including its inherited record number. Every original cell appears once, and none remains unresolved. Each has either a strict necessary contradiction on the whole cell or a complete internal closed cover whose every terminal has such a contradiction. The value construction encloses the lift of every qualified candidate in that cell. The contractions in Section~\ref{nine:separation} preserve the same candidate; each internal closed split preserves coverage. Finite induction along its recorded tree would therefore place an existing candidate in a terminal where a strict contradiction holds, which is impossible.

The final acceptance also binds the original source and expression versions, the full ordered inputs and their mathematical receipts. Its evaluation-profile counts agree with the disjoint original-identity partition above; the check is on complete identities and records, not just their numbers. The accepted source-to-domain implication remains a premise of this lemma, rather than a conclusion drawn from the acceptance index.
\end{proof}

\begin{proof}[Proof of Theorem~\ref{nine:thm:r05}]
Suppose an original R05 global minimizer exists. Its finite essential arcs are positive, with the genuine twist and original cost budget of Section~\ref{fd:section}. Lemma~\ref{nine:qualification} supplies the required complete sources and one normal jointly normalized full lift, including both state-covector twists. Its actual states, caps, selected root jets, three means and real winding therefore have precisely the qualifications of Proposition~\ref{nine:source-reduction}. That proposition places the same candidate in one complete original cell, unless an earlier source comparison has already contradicted it. Lemma~\ref{nine:exhaustion} excludes the former possibility as well. Thus no such original global minimizer exists.
\end{proof}
The argument concerns strict essential R05 arcs. Closed numerical boxes add harmless auxiliary boundary points, while any physical boundary use is restricted to the qualified limits in the source theorems. It grants no general stationary or Hessian condition at zero physical duration, $D=0$, a lost selected root, or a formal tuple without the original full lift and winding. The independent short octagon class is unaffected.
\space
\section{Exclusion of R06}\label{r06:section}
In this type the middle shared arc has the opposite cap orientation from the two exterior shared arcs. The proof retains that distinction and couples all three caps through one original full lift.

\begin{theorem}[Exclusion of R06]\label{r06:exclusion}
An original global minimizer cannot have nine essential positive arcs with native labels $021202120$ and increments $\mathtt{--++--+++}$, with the genuine extension $j_{i+9}=j_i+1$.
\end{theorem}

The argument has two parts. We first place every such minimizing trajectory in one fixed five-dimensional root, keeping its full signed endpoints, maximizing lift and original cost. We then exclude its necessary values on a complete closed cover. The mathematical reduction is given below and in Appendix~\ref{app:r06-source}; the independent finite verification is specified in Section~\ref{reader:r06}.

\subsection{The actual source and its positive durations}
Use the reciprocal coordinates of Section~\ref{nine:source-interface}. Write the native own reciprocals as
\begin{equation}\label{r06:native}
 (r_0,\ldots,r_8)=(d,e,b,m,D,E,B,f,c).
\end{equation}
Within this section and its mathematical appendix, $D$ and $B$ are ordinary reciprocal coordinates; in particular this $B$ is distinct from the essential-arc count. Put
\[
 h(t)=\frac{1+t^2}{2},\qquad h_1=h(e),\quad h_2=h(E),\quad H=h(m).
\]
Let $\mathcal S_{06}$ denote the original global-minimizing trajectories of this type, together with their actual labelled reciprocal data on a genuine period. Thus the nine maximal arcs are positive, the state and group endpoints satisfy the original signed twist, and the real winding and budget of Section~\ref{fd:section} are retained. The auxiliary equations introduced below do not define this source class.

For an arc with own reciprocal $x$, predecessor $p$, successor $n$ and adjacent increments $\alpha,\beta\in\{-1,1\}$, the physical multiplier is
\[
 Q=\frac{(x+p)^\alpha(x+n)^\beta}{(1+x^2)^{(\alpha+\beta)/2}},
 \qquad \tau=x\log Q>0.
\]
The same formula applies to negative increments: the physical time is not reversed. Substitution of the R06 word gives
\begin{equation}\label{r06:multipliers}
\begin{gathered}
 Q_d=\frac{d+c}{d+e},\quad
 Q_e=\frac{1+e^2}{(e+d)(e+b)},\quad
 Q_b=\frac{b+m}{b+e},\\
 Q_m=\frac{(m+b)(m+D)}{1+m^2},\quad
 Q_D=\frac{D+m}{D+E},\quad
 Q_E=\frac{1+E^2}{(E+D)(E+B)},\\
 Q_B=\frac{B+f}{B+E},\quad
 Q_f=\frac{(f+B)(f+c)}{1+f^2},\quad
 Q_c=\frac{(c+f)(c+d)}{1+c^2}.
\end{gathered}
\end{equation}
All nine are greater than one. In particular
\begin{equation}\label{r06:strict-widths}
 c>e,\qquad m>e,\qquad m>E,\qquad f>E.
\end{equation}
The actual third reciprocal at the $c,f$ switch is
\begin{equation}\label{r06:seam}
 s=\frac{1-cf}{c+f}>0,\qquad cf+s(c+f)=1.
\end{equation}
Since $Q_c=(c+d)/(c+s)$ and $Q_f=(f+B)/(f+s)$, also $d>s$ and $B>s$. Thus both actual rectangle areas
\begin{equation}\label{r06:areas}
 A=(c-e)(d-s),\qquad Z=(f-E)(B-s)
\end{equation}
are positive. These conclusions use positive physical arcs, without a boundary-face exclusion or an assumed numerical collar.

\begin{lemma}[Qualification on one original period]\label{r06:qualification}
Every member of $\mathcal S_{06}$ has one jointly normalized original full maximizing lift on native $[-1,7]$. Its matrix adjoint and state covector retain the genuine endpoint relations. In a common whole phase $+1$, its three five-arc supports are
\[
 [-1,3]:01020,\qquad [1,5]:02010,\qquad [3,7]:01020,
\]
with shared reciprocals $e,m,E$, respectively. All actual own and adjacent third reciprocals lie in $[1/15,5]$, and $e^2,E^2<1/6$.
\end{lemma}
\begin{proof}
The interval $[-1,7]$ has literal word $202120212$ and contains all nine complete physical arcs once. Native $[0,4]=02120$ becomes $10201$ under whole phase $+1$; its binary core $020$ has two distinct interior switches and both actual third-mode caps. The homogeneous argument of Lemma~\ref{pc:rigidity} first proves normality of a full maximizing lift selected on this very period. Dividing its objective multiplier, matrix adjoint and state covector by the same positive number gives the stated normalization. Appendix~\ref{app:r06-source} proves the free-state endpoint row for that same lift; no covector is copied across a seam.

The reversed support $[4,3,2,1,0,-1]$ with whole reflected label map $j\mapsto2-j$, and the forward support $[2,7]$ with whole phase $+1$, both become the proper six-source $201020$. Their physical shared centers are $e$ and $E$. Theorem~\ref{new:sixdensity} gives the two density bounds. Only physical state consequences are exported from these recuts. The compactness theorem~\ref{fd:compact}, transported by the lawful whole phases, bounds every actual reciprocal, including the third component $s$; its elementary bound $400/99<5$ gives the stated convenient interval. The middle support remains inner-sharing, so no exterior density bound is imposed on $m$.
\end{proof}

\subsection{The complete frame and a common cap system}
Define
\begin{equation}\label{r06:frame-quantities}
 R=(m-e)\left(b+\frac{m+e}{2}\right),\quad
 L=(m-E)\left(D+\frac{m+E}{2}\right),\quad
 F=H^2-RL,\quad G=2h_1Hh_2.
\end{equation}
The original full signed matrix equation gives
\begin{equation}\label{r06:frame-equations}
 F(1+s^2)=G,\qquad
 A=\frac{(s+e)^2}{2}-\frac{h_1^2L}{F},\qquad
 Z=\frac{(s+E)^2}{2}-\frac{h_2^2R}{F}.
\end{equation}
Appendix~\ref{app:r06-source} derives these identities from all four signed frame entries. In particular their use retains the original sign of the frame.

Use the fixed-interval mean $M$ of \eqref{nine:mean}, and set
\begin{equation}\label{r06:actual-means}
 \begin{gathered}
 \mu_0=M(d;e,c),\qquad \mu_1=M(b;e,m),\\
 \mu_2=M(D;E,m),\qquad \mu_3=M(B;E,f),\qquad K_i=\frac{3+\mu_i}{2}.
 \end{gathered}
\end{equation}
These are four different integrals at one candidate. For positive $t,x,y$ define the bridge function
\begin{equation}\label{r06:bridge}
\begin{split}
 \mathcal B(t;x,y;\mu,\nu)
 ={}&\frac12+h(t)\log\frac{(t+x)(t+y)}{1+t^2}\\
 &+\frac{h(t)}2
 \left(p+q+\frac{\mu-1}{p}+\frac{\nu-1}{q}\right),\qquad
 p=\frac{2t}{t+x},\quad q=\frac{2t}{t+y}.
\end{split}
\end{equation}
Write
\[
 B_e=\mathcal B(e;d,b;\mu_0,\mu_1),\quad
 B_m=\mathcal B(m;b,D;\mu_1,\mu_2),\quad
 B_E=\mathcal B(E;D,B;\mu_2,\mu_3).
\]
The three restrictions of the same full lift have common determinant and the two unsquared overlap equations
\begin{equation}\label{r06:overlaps}
 B_e-B_m+RK_1=0,\qquad B_E-B_m+LK_2=0.
\end{equation}
Eliminating the common determinants, using only $K_1,K_2>0$, gives
\begin{equation}\label{r06:targets}
 h_1^2K_0-R^2K_1-H^2K_2+2RB_m=0,\qquad
 h_2^2K_3-H^2K_1-L^2K_2+2LB_m=0.
\end{equation}
The unsquared equations are retained. Equality of three determinant values alone would not identify the matrices or either state covector.

\subsection{The enclosing root and its fourteen-variable system}
\begin{proposition}[Source-to-root reduction]\label{r06:source-root}
The projection $(e,E,b,D,m)$ of every member of $\mathcal S_{06}$ belongs to
\begin{equation}\label{r06:root-box}
 B_{06}=\left[\frac1{16},\frac{27}{64}\right]^2
       \times\left[\frac1{16},5\right]^2
       \times\left[\frac38,1\right].
\end{equation}
At the same candidate,
\begin{equation}\label{r06:sharper-domain}
 \frac{21}{50}<m<\frac{471}{500},\quad
 b>\frac{m+e}{2},\quad D>\frac{m+E}{2},\quad
 s<m,\quad F>h_1h_2.
\end{equation}
The complete original minimizing and full-lift qualifications remain in force.
\end{proposition}
\begin{proof}
The reserve identity in Lemma~\ref{r06:reserve} gives $s<m$ and $F>h_1h_2$ directly on the strict source. The actual return inequalities give the two owner quadratics in Lemma~\ref{r06:mean-root}. The analytic Lemma~\ref{r06:middle-less-than-one} then proves $m<1$. It uses only these quadratics, the compact and density bounds, and the strict frame reserve. The midpoint comparison and the lower and upper middle estimates of Appendix~\ref{app:r06-source} give the remaining assertions. This order is independent of any unit-sector coefficient certificate or boundary exclusion.

Finally $[1/15,1/\sqrt6]\subset[1/16,27/64]$, $[1/15,5]\subset[1/16,5]$ and $[21/50,471/500]\subset[3/8,1]$. Replacing strict source bounds by these closed outer intervals proves the projection claim without selecting an owner ordering.
\end{proof}

For the finite verification use fourteen coordinates
\begin{equation}\label{r06:fourteen}
 z=(e,E,b,D,m,s,c,d,B,f,u_0,u_1,u_2,u_3)
 \in\Omega_{06}:=(0,\infty)^{10}\times\mathbb R^4.
\end{equation}
At a physical candidate $u_i=\mu_i$; the auxiliary mean coordinates are unrestricted real variables in this smooth extension. Put $\widehat F=G/(1+s^2)$ and use $K_i=(3+u_i)/2$ in the bridges. Define twelve smooth residuals by
\begin{equation}\label{r06:twelve-rows}
\begin{aligned}
 f_0&=(H^2-RL)(1+s^2)-G,\\
 f_1&=(c-e)(d-s)-(s+e)^2/2+h_1^2L/\widehat F,\\
 f_2&=(f-E)(B-s)-(s+E)^2/2+h_2^2R/\widehat F,\\
 f_3&=cf+s(c+f)-1,\\
 f_4&=h_1^2K_0-R^2K_1-H^2K_2+2RB_m,\\
 f_5&=h_2^2K_3-H^2K_1-L^2K_2+2LB_m,\\
 f_6&=B_e-B_m+RK_1,\qquad f_7=B_E-B_m+LK_2,\\
 f_8&=M(d;e,c)-u_0,\qquad f_9=M(b;e,m)-u_1,\\
 f_{10}&=M(D;E,m)-u_2,\qquad f_{11}=M(B;E,f)-u_3.
\end{aligned}
\end{equation}
Every source point satisfies $f_i=0$. The use of $\widehat F$, rather than substitution of $H^2-RL$ off the frame, fixes the ambient functions whose derivatives and Hessians are enclosed. The mean in \eqref{nine:mean} is smooth even when its two endpoints coincide. Bounds valid only for actual returned means may intersect their values, but are not differentiated as ambient formulas.

With the potential $\phi$ in \eqref{nine:potential}, the original cost is
\begin{equation}\label{r06:cost}
\begin{split}
 C_{06}={}&-\phi(d,e)-\phi(e,b)+\phi(b,m)+\phi(m,D)-\phi(D,E)\\
          &-\phi(E,B)+\phi(B,f)+\phi(f,c)+\phi(c,d).
\end{split}
\end{equation}
The closing edge is included once. Set
\[
 J_8:=C_{\rm oct}=\frac{2(8-4\sqrt2-\log2)}{2\sqrt2-1},
\]
as in \eqref{fd:budget}. The thirteenth row is the one-sided inequality $C_{06}-J_8\le0$. There are thus twenty-four signed equality inequalities and one cost inequality, not two directions for the cost row.

\subsection{Exhaustion and conclusion}
\begin{theorem}[Complete finite exclusion]\label{r06:finite-cover}
The fixed closed binary partition of $B_{06}$ has $806\,875$ original leaves. Every original leaf is excluded for a point of $\mathcal S_{06}$ by the necessary systems and candidate-preserving comparisons above and in Appendix~\ref{app:r06-source}. The independently assembled cover has no retained or unvisited original leaf.
\end{theorem}
\begin{proof}
The exact midpoint parser reconstructs every original interval from the single root \eqref{r06:root-box}, with both closed children retained at each split. In a surviving leaf the full enclosing $c$ interval is divided into eight closed pieces, and the actual mean values lift each surviving piece to its complete fourteen-dimensional box. All contractions retain the same possible source point. A leaf is excluded only by a necessary strict sign, an empty necessary intersection, or a fully checked affine contradiction with a valid whole-box Taylor remainder. Further subdivisions likewise retain both complete closed children.

In the earlier verification, independent subtrees covered $34\,838$ original cells. Their inputs and proof chains were matched to the original tree, and the remaining $772\,037$ cells were independently checked. The \path{R06} task reconstructs the original tree from \eqref{r06:root-box} and checks the exclusion on each of its $806\,875$ leaves, including the leaves in those subtrees. Section~\ref{reader:r06} gives the arithmetic, row construction and exact acceptance rules for this verification.
\end{proof}

\begin{proof}[Proof of Theorem~\ref{r06:exclusion}]
Suppose an original global minimizer has the stated type. Its nine essential maximal arcs are positive. Lemma~\ref{r06:qualification} qualifies the same trajectory, its actual means and one complete normalized maximizing lift; Proposition~\ref{r06:source-root} places its projection in $B_{06}$. It lies in an original leaf and in at least one member of each subsequent complete closed cover. At every candidate-preserving step its same fourteen values survive unless a necessary condition is contradicted. Theorem~\ref{r06:finite-cover} leaves no such member. This contradiction proves the exclusion.
\end{proof}

In particular $\mathcal S_{06}=\varnothing$, and consequently its genuine closure in the bounded reciprocal coordinates is empty. This final implication concerns limits of that same strict original source class. It does not identify a larger weak algebraic envelope with that closure or impose new variational conditions at a zero-duration boundary.
\space
\section{Short chains and completion of the proof}\label{five:completion}
\begin{theorem}[Exclusion of all long minimizing types]\label{ten:assembly-proposition}
No original global minimizer has essential count $B\ge8$.
\end{theorem}
\begin{proof}
Theorem~\ref{new:uniform-six} places every such minimizer among
R01--R06. Theorems~\ref{eight:thm:main}, \ref{mix:thm:r02},
\ref{mix:thm:r03}, \ref{r04:exclusion}, \ref{nine:thm:r05} and
\ref{r06:exclusion} exclude these types in order.
\end{proof}

\subsection{Short-chain rigidity}
The domain and the exhaustive short word list are proved in
Section~\ref{short:integration}. We assemble the geometric comparisons
here, with proofs in Appendix~\ref{app:framework}. In the short-chain
formulas, $\mathscr J=C$ and $\mathscr J_8=J_8=C_{\rm oct}$. The comparison
budgets are $U=C_*$, $V=361/200$, and $K_{\rm sh}=29/16$.

\begin{proposition}[Short-chain rigidity]\label{short:rigidity}
An original Reinhardt global minimizer with $B\leq7$ is an affine
image of the smoothed regular octagon.
\end{proposition}
\begin{proof}
By Lemma~\ref{short:interface}, all the variations and complete-PMP
conditions below apply in the original geometric class.  The valid short words are listed in Section~\ref{short:words}.

\emph{At most four arcs.}
Lemma~\ref{sh:omitted} excludes the one-, two- and three-arc words
and the mixed four-cycle.  Its missing-label argument uses the
straight boundary components and the strict positive determinant
of their tangents; it also handles the repeated-arc issue at $B=3$.
For the coherent four-cycle, Proposition~\ref{sh:octagon-rigidity}
classifies the whole family rather than imposing symmetry.  The
balanced support octagon is an intersection of two congruent
squares.  Its hyperbolic rounding gives a one-parameter area
function, strictly increasing away from the regular octagon.
Thus $\mathscr J\geq\mathscr J_8$, with equality precisely for
its affine images.

\emph{Five arcs.}
For the word $01202$, Lemmas~\ref{sh:mixed-five-chart}
and~\ref{sh:h-crossing} give the exact two-parameter cost and
exclude all its interior critical points.  Proposition~\ref{sh:mixed-five}
then checks every boundary and escaping sequence, giving
$\mathscr J>\mathscr J_8$.  For $02102$,
Proposition~\ref{sh:reverse-five} gives the stronger bound
$\mathscr J>K_{\rm sh}$.  Both compare all actual strict chains of the
stated words.

\emph{Six arcs.}
Proposition~\ref{sh:six-010202} gives an exact three-parameter
closure family about every actual $010202$ chain and proves
that its cost has no critical point.  For $021020$,
Lemma~\ref{sh:six02-chart} gives a local exact closure family
with parameters $(a,b,t)$.  Lemma~\ref{sh:six02-finite} proves
\[
 \mathscr J\leq K_{\rm sh}\quad\Longrightarrow\quad
 \left.\frac{\partial\mathscr J}{\partial b}\right|_{a,t}>0.
\]
The octagonal competitor gives $\mathscr J\leq\mathscr J_8<K_{\rm sh}$,
so decreasing $b$ slightly improves the cost.  The full endpoint
conditions, strict star and real winding persist as verified in
Lemma~\ref{short:interface}.

\emph{Seven arcs.}
The representative $\texttt{+------}$ has a cut $0210202$.
Proposition~\ref{sh:reverse-seven} and
Lemma~\ref{sh:reverse-inequality} give a lower comparison
$\mathscr J\geq L>\mathscr J_8$ on its complete necessary domain,
including the noncompact tails and the octagonal boundary collar.
No uniform positive gap is needed.

For $\texttt{+++++++}$, Lemmas~\ref{sh:state-filter}
and~\ref{sh:state-monotone} use complete maximization and a
finite state graph at budget $V$ to prove $q_i<2$ and
$757/2000\leq a_i\leq67/50$.  The exact full-matrix closure
chart, a largest-exponent cut and Lemma~\ref{sh:pure-finite}
then give $\mathscr J>U>\mathscr J_8$, proving
Proposition~\ref{sh:forward-seven}.  The prerequisite budget
$V>U$ keeps the state-domain reduction independent of the final
cost exclusion.

For $\texttt{++++---}$, take the actual extended arc indices
$3,4,5,6,7,8$.  Modulo seven these are $3,4,5,6,0,1$, all
distinct, and their phased reduced word is $010212$.  They omit
positive arc $2$.  Cutting in the interior of that omitted arc
puts the six-arc support strictly inside one genuine fundamental
interval.  Lemma~\ref{sw:six} applies to its one complete
maximizing lift and excludes this last word.

Only the coherent four-arc family remains.  The octagonal
competitor and its strict equality classification force the
asserted affine smoothed regular octagon.
\end{proof}

\space

\begin{proof}[Proof of Theorem~\ref{five:main}]
An original global minimizer exists by the geometric inputs of
Section~\ref{fd:section}. Theorem~\ref{fw:finiteness} gives its finite
vertex-control representation. Theorem~\ref{ten:assembly-proposition}
excludes $B\ge8$. Proposition~\ref{short:rigidity}, with the domain identification of Lemma~\ref{short:interface}, therefore identifies every global minimizer with a nonsingular affine image of the smoothed regular octagon. The minimum density is $C_{\rm oct}/2$: the normalization fixes the minimum circumscribed hexagon area at $2\sqrt3$ and gives disk area $\sqrt3 C$. Substitution of \eqref{fd:budget} gives the stated constant. Conversely, nonsingular affine transformations preserve $\delta_L$, so every such image attains the same minimum. A disk has equality in the density inequality exactly when it is itself a global minimizer; this proves the asserted equality range. The octagon has $B=4$, and \eqref{fd:cut-count} gives $N\le B+1=5$ for every prescribed cut.
\end{proof}
\space

\clearpage\appendix
\section{Geometric and switching proofs}\label{app:framework}
\begin{proofappendix}
\section{The geometric problem and its full lift}\label{sec:framework}
\subsection{Normalization and admissible trajectories}
We use the geometric class, matrices, physical-time equations, cost and
full endpoints of Section~\ref{fd:section}, in particular
\eqref{fd:domain}--\eqref{fd:endpoints}. The prescribed real winding and
six-curve reconstruction remain part of that class. We verify the
control convention and the local assertions stated there.

To verify the control convention, put
\[
\mathsf B_0=\begin{pmatrix}1/2&1\\0&-1/2\end{pmatrix},\quad
\mathsf B_1=\begin{pmatrix}0&-1/2\\-1/2&0\end{pmatrix},\quad
\mathsf B_2=\begin{pmatrix}-1/2&1\\0&1/2\end{pmatrix}.
\]
Then $d_j=\tr(\mathsf B_jX)$ are the positive quantities
\begin{equation}\label{fw:reciprocals}
(d_0,d_1,d_2)=\left(\frac{1+u}{v},\frac{S}{2v},\frac{1-u}{v}\right),
\qquad DX(W)F_j=\frac{[\mathsf B_j,X]}{d_j}.
\end{equation}
The positive $d_j$ make the control changes \eqref{fd:controls} inverse,
measurable and vertex-preserving. In the coordinates of \cite{HV2024}, $(x,y)=(u,v)/\sqrt3$. For its control $\alpha$, take $\mu=(\alpha_2,\alpha_0,\alpha_1)$; conjugation by $D_*$ sends its control matrix to $-2\sum_j\mu_j\mathsf B_j/\sqrt3$. The common factor cancels in the commutator-to-trace quotient. Thus the complete measurable control class and the physical time are unchanged.

\subsection{The original compact set}
Let $\mathcal R(w)=(w-3)/(w+1)$. The compactness theorem of \cite[Definition 5.0.1 and Theorem 5.0.2]{HV2024} places every normalized global minimizer in
\[
K=\Omega\smallsetminus\bigcup_{j=0}^2\mathcal R^j
\bigl(\{v>9\sqrt3/2\}\cup\{v>15(1-u)\}\bigr).
\]
Its six cuts are \eqref{fd:compact-cuts}; they constrain the minimizing
trajectory, not a competing local replacement.

\begin{lemma}\label{fw:compact-bounds}
The rates on $K$ satisfy \eqref{fd:compact-bounds}.
\end{lemma}
\begin{proof}
The side cuts give the three individual bounds. Put $c_0=513/1000$. The cusp cuts imply $v<39/5$ and $(u\pm1)^2+v^2\ge c_0v$. If $v<1/2$, the side cut gives $u^2\ge q_0=1+2v/15-v^2>0$. Hence, with $\delta=1-|u|$,
\[
c_0v\le\delta^2+v^2\le(1-\sqrt{q_0})^2+v^2.
\]
Writing $d=c_0/2-1/15=1139/6000$, this implies $\sqrt{q_0}\le1-dv$, whose positive right side may be squared. Consequently
$v\ge c_0/(1+d^2)=18468000/37297321>99/200$.
For any two rates $x,y$, the three identities
\[
\frac{2ab}{a+b}=v,\quad
\frac{2ac}{a+c}=\frac{4v}{(u+1)^2+v^2},\quad
\frac{2bc}{b+c}=\frac{4v}{(u-1)^2+v^2}
\]
give $xy\le(39/10)(x+y)$. Also $(15-x)(15-y)\ge0$, so $xy\ge15(x+y)-225$. Combining them gives the asserted pairwise bound.
\end{proof}

\subsection{Five-endpoint replacements}
\begin{lemma}[Normalized reconstruction]\label{fw:balanced-reconstruction}
Let a trajectory of \eqref{fd:dynamics} with states in $\Omega$
and controls in $\Delta$ satisfy \eqref{fd:endpoints} and the prescribed
real winding. Then the six arcs $G(t)q_j$ form a regular convex boundary
traversed once. Its disk has minimum circumscribed centrally symmetric
hexagon area $2\sqrt3$ and satisfies $\delta_L(D)=|D|/(2\sqrt3)$.
\end{lemma}
\begin{proof}
By the twisted extension \eqref{fd:extension}, the arcs concatenate
with matching tangents to the regular closed curve
$\gamma(t)=G_e(t)q_0$, $0\leq t\leq6T$. Its nonzero tangent
$G_eX_eq_0$ has winding one. Indeed one can homotope the periodic
$X_e=X(W_e)$ to a constant within the upper half-plane parametrization,
and the real winding condition homotopes the group loop to the six
consecutive rotations of angle $\pi/3$ in the $D_*$ coordinates.

Write $\mathsf B=\sum_j\mu_j\mathsf B_j$ and $d=\tr(\mathsf BX)>0$.
The identity $X\mathsf BX=\mathsf B+dX$ gives
$X+X^{-1}X'=-2\mathsf B/d$. Consequently, for $\gamma_j=Gq_j$,
\[
 (\det(\gamma_0',\gamma_0''),\det(\gamma_2',\gamma_2''),
   \det(\gamma_4',\gamma_4''))
   =\frac{\sqrt3}{d}(\mu_1,\mu_2,\mu_0)\geq0
\]
almost everywhere, and the opposite arcs have the same curvatures.
The absolutely continuous tangent angle of $\gamma$ is nondecreasing
and changes by $2\pi$. Relative to any oriented tangent, the signed
normal coordinate first increases and then decreases back to zero.
Every tangent therefore supports the convex hull, which has nonempty
interior since the tangent winds once. At a point of the curve any
supporting line contains the nonzero tangent, by taking secants from
both sides. If two parameter values gave the same point, their
oriented tangents would agree, since the hull lies to their left.
The angle would then be constant along one of the
two intervening intervals. Positive speed prevents a return along
such an interval. Thus the closed curve is simple and, as it lies in
the boundary of its convex hull, traverses that boundary once.

The opposite arcs give central symmetry. For any boundary point $p$,
one of the six arcs gives $p=G(t)q_j$, and its simultaneous neighbors
give a boundary pair $y,y-p$ with $\det(p,y)>0$. Choose oriented
linear coordinates in which $p=(1,0)$. The two tangents at $y,y-p$
are linearly independent since the common tangent map is
$G(t)X(W(t))$ and the preimages are $q_{j+1},q_{j+2}$. Their common
positive height is therefore below the top supporting line. At this
interior height they are the endpoints of a horizontal section of
width one. The section width is concave in the height and is two at
height zero. It is therefore strictly greater than one at lower
positive heights and strictly less than one at all higher heights,
including the top. This is the only boundary pair of the form
$y',y'-p$ above the horizontal axis.

By the critical lattice characterization
\cite[Theorems~2.3.6 and~2.3.7]{HV2024}, the minimum covolume of a
$D$-admissible lattice is the minimum of $\det(p,y)$ over positively
ordered boundary triples $p,y,y-p$. Each triple is one of those
represented by the six arcs, so in the original coordinates its
determinant is $\sqrt3/2$. The reconstructed disk has a regular
$C^1$ boundary, so it is not a parallelogram. The same characterization
gives minimum circumscribed hexagon area $4(\sqrt3/2)=2\sqrt3$.
The minimum packing covolume is four times the minimum
$D$-admissible covolume, by central symmetry and rescaling.
\end{proof}

\begin{proof}[Proof of Lemma~\ref{fd:replacement}]
Translate the original tail in time and attach it to the new window.
Matching $G,W$ matches all five endpoint coordinates and the tangents.
The global endpoints and twist are unchanged. After extending the
comparison reparametrization to the whole interval, the matrix
quotient of the new and old paths lies in a real logarithm neighborhood
of $I$. Its logarithm is traceless, continuous, and zero at both
endpoints, so exponentiating a scalar multiple gives a fixed-endpoint
homotopy. Reparametrization is used only for this comparison. The
replacement still solves the physical-time equation.
Lemma~\ref{fw:balanced-reconstruction} gives the convex boundary and
its hexagon normalization. Green's formula for the six arcs gives
\[
|D|=\frac12\sum_{j=0}^5\int\det(\gamma_j,\gamma_j')\,dt
=\sqrt3\int\ell(W)\,dt.
\]
All terms outside the replacement cancel. This proves the assertion.
\end{proof}
All duration variations below are taken about positive arcs with positive caps. Small variations preserve positivity and remain in the open set $\Omega$. Signed analytic extensions of durations may be used for differentiation, but are not themselves competing geometric trajectories.

\subsection{The full free-time maximum principle}
Use the full lift of Definition~\ref{fd:full-lift}, with the matrix pairing
$\operatorname{tr}(\Lambda G^{-1}\delta G)$ and the free-time equations
\eqref{fd:PMP}. The switching rows and their positive factor $\beta$ are
\eqref{fd:switching-rows}--\eqref{fd:ties}; costates from independently
selected lifts are never combined.
Thus a maximizing $02,01,12$ tie has signs respectively
$(0,\le0,\ge0)$, $(\ge0,0,\le0)$, $(\le0,\ge0,0)$.

\begin{lemma}\label{fw:regularity}
Along every measurable maximizing lift, each $\sigma_i$ is $C^1$ with absolutely continuous first derivative. Its second derivative is the control-weighted sum of the supported vertex equations almost everywhere.
\end{lemma}
\begin{proof}
Let $L_j$ be differentiation along the full equations with mode $j$. From $\alpha_0=\beta(F_0-F_2)$ and its cyclic versions,
\[
L_j\sigma_i-L_k\sigma_i=f_i(H_j-H_k),\qquad
f_i=\frac{\alpha_i\nabla\beta}{\beta},\qquad L_jH_j=0.
\]
Hence $\Psi_i=L_j\sigma_i-f_iH_j$ is independent of $j$. On the maximizing trajectory $\sigma_i'=\Psi_i$ almost everywhere. The right side is continuous and absolutely continuous. Differentiating once more on a supported mode, where $H_j=0$, gives $L_j\Psi_i=L_j^2\sigma_i$. Averaging gives the assertion without differentiating the measurable control.
\end{proof}

\subsection{Symmetry and genuine reanchoring}
The whole-state, frame and lift transports are those of
Proposition~\ref{fd:lawful-symmetry}; the twisted state and group extension
is \eqref{fd:extension}, with no state-costate extension.
\begin{proof}[Proof of Lemma~\ref{fd:reanchor}]
The equations, endpoints, and periodicity of the cost follow from the covariance identities. The family
$G_e(\lambda s)^{-1}G_e(\lambda s+t)$, $0\le\lambda\le1$, has constant endpoints $I,R_*$ and preserves the real homotopy class. More geometrically, with $\Gamma(t)=G_e(t)q_0$,
\[
G_s(t)q_j=G_e(s)^{-1}\Gamma(s+t+jT).
\]
These are the same convex boundary, with a new starting point and an area-preserving linear transformation. Its area and circumscribed-hexagon normalization are unchanged. The original maximum principle applies afresh.
\end{proof}
A complete local support using fewer than $B$ physical arcs can be put inside such a genuine interval by placing the cut in a remaining positive arc. The numerical properness condition guarantees room for the support; it does not replace the geometric reanchoring argument.

\subsection{Finiteness and the competitor budget}
\begin{theorem}[Finite switching]\label{fw:finiteness}
The control of every normalized global minimizer is, up to a null set,
piecewise constant with finitely many pieces, each at a vertex of $\Delta$.
\end{theorem}
The proof is given in Section~\ref{fw:finiteness-proof}, after the
harmonic bound. It applies to an arbitrary minimizing disk.

\paragraph{Elementary scalar enclosures.}
For integers $N\ge1,n\ge0$, $x\ge1$ and $0\le w<n+2$, put
$v_x=(x-1)/(x+1)$ and
\begin{equation}\label{hc:log-bounds}
 L_N(x)=2\sum_{j=0}^{N-1}\frac{v_x^{2j+1}}{2j+1},\qquad
 U_N(x)=L_N(x)+\frac{2v_x^{2N+1}}{(2N+1)(1-v_x^2)}.
\end{equation}
Then $L_N(x)\le\log x\le U_N(x)$, strictly for $x>1$.
Evaluate $U_N-L_N$ as its displayed tail, without rounded subtraction.
As functions of $v_x$, the lower and upper error derivatives are
$2v_x^{2N}/(1-v_x^2)$ and
$4v_x^{2N+2}/((2N+1)(1-v_x^2)^2)$; both errors vanish at zero.
Also
\begin{equation}\label{strip:exp-upper}
 \mathcal E_n(w)=\sum_{j=0}^n\frac{w^j}{j!}
 +\frac{w^{n+1}}{(n+1)![1-w/(n+2)]}\ge e^w,
\end{equation}
strictly for $w>0$: bound $(n+1+j)!\ge(n+1)!(n+2)^j$.

Write $\mathcal A=|D|=\sqrt3\,C$. For a global minimizer the
comparison with the smoothed octagon gives
$\mathcal A\le\mathcal A_{\rm oct}=\sqrt3\,C_{\rm oct}$ and
$C\le C_{\rm oct}$ by \eqref{fd:budget}. In particular,
\begin{equation}\label{fw:budget}
C< C_*:=\Cstar,\qquad
\mathcal A<A_0:=\frac{781513607211}{250000000000},\qquad
\pi>p_0:=\frac{3141592653589}{10^{12}}.
\end{equation}
Here are exact bounds sufficient to verify the constants. Squaring the consecutive rationals
$1414213562373095/10^{15}$ and $1414213562373096/10^{15}$ brackets $\sqrt2$; the analogous numerators $1732050807568877,1732050807568878$ bracket $\sqrt3$. Moreover $L_{16}(2)<\log2<U_{16}(2)$ by \eqref{hc:log-bounds}.
The quotient $(8-4h-l)/(2h-1)$ decreases in both $h$ and $l$ on these intervals. Its rational upper bound proves the first two claims. Machin's identity $\pi=16\arctan(1/5)-4\arctan(1/239)$, with respectively $26$ and $10$ alternating-series terms and their next-term remainders, proves the third. The branch of the identity follows from $\tan(4\arctan(1/5))=120/119$: subtracting $\arctan(1/239)$ gives the angle $\pi/4$ in $(0,1)$.


\section{Switching restrictions and normality}\label{sec:switching}
Every letter in this section denotes an actual positive maximal control arc. The two exterior letters of a local word remain positive caps. Roots and crossing derivatives always belong to the chosen lift of the whole containing interval.

\subsection{An independent edge argument}
\begin{lemma}\label{sw:four}
On any closed window of an original global minimizer where the control
stays on one edge of $\Delta$ almost everywhere, it is a vertex control
with finitely many switches, up to a null set. The words $2020,0101,1212$
cannot occur on an original global minimizer.
\end{lemma}
\begin{proof}
By a cyclic phase, work on the edge with modes $0,2$. Here $u'=v>0$,
so use $u$ as independent variable. Put $r=v^{-2}$ and
$s(u)=\int_{u_-}^{u}r(\xi)\,d\xi$, where $u_-$ is the entrance
coordinate. With derivatives now taken in $u$, the equations are
\[
 r'=2mr,\qquad s'=r,\qquad
 -\frac1{1+u}\le m\le\frac1{1-u}.
\]
Here $m=d\log(v^{-1})/du$ is affine in the simplex control, with its
lower and upper endpoints given by modes $0$ and $2$, respectively.
For $Q_u=\left(\begin{smallmatrix}1&u\\0&1\end{smallmatrix}\right)$,
the frame equation gives
\[
 G(u)Q_u=G(u_-)Q_{u_-}
   \begin{pmatrix}1&0\\s(u)&1\end{pmatrix},\qquad
 \ell\,dt=\frac{1+3r+u^2r}{2}\,du.
\]
Fixing the two $u$ endpoints and the endpoint values of $r,s$ fixes
the full state and frame endpoints. A nearby edge trajectory with
these data reconstructs by $v=r^{-1/2}$ and $dt=du/v$. It stays in
$\Omega$ by openness, and the frame formula gives uniform closeness
in the $u$ parametrization. Lemma~\ref{fw:replacement} makes it an
admissible competitor. Since the integrals of $1$ and $r$ are fixed,
the original trajectory locally minimizes $\int u^2r\,du$ there.

Apply the maximum principle to this edge problem with its own
nonzero multipliers. Let $\lambda$ be the covector of $r$, $\zeta$
the constant covector of $s$, and $\lambda_0\ge0$ the cost multiplier.
The Hamiltonian is $2mr\lambda+r\zeta-\lambda_0 u^2r$, so that
\[
 \chi:=r\lambda,\qquad \chi'=rP(u),\qquad P(u)=\lambda_0u^2-\zeta.
\]
The right-hand side is continuous. If $P$ is nonzero as a polynomial,
its at most two zeros give at most three zeros of $\chi$ by Rolle's
theorem. If $P$ vanishes identically, then $\lambda_0=\zeta=0$ and
$\chi$ is a nonzero constant by nontriviality. Maximization gives mode
$2$ for $\chi>0$ and mode $0$ for $\chi<0$. This proves the first claim.

For the word $2020$, at the three successive switches $u_1<u_2<u_3$
the signs of $\chi$ give $P(u_1)\le0$, $P(u_2)\ge0$ and $P(u_3)\le0$.
If $\lambda_0>0$, strict convexity of $P$ contradicts these inequalities.
If $\lambda_0=0$, the constant $P$ gives either strict monotonicity
of $\chi$ or a nonzero constant, so three switching zeros are impossible.
Cyclic relabeling proves the other exclusions.
\end{proof}

\subsection{Binary identities and actual signs}
Write
\[
\Lambda=\begin{pmatrix}l&b\\m&-l\end{pmatrix},\quad
K_e=2ul+b-u^2m,\quad L=l-um,\quad\chi=vp_v,\quad ds=dt/v.
\]
On either mode $0$ or mode $2$, substitution in the full equations gives
\begin{equation}\label{sw:binary}
(K_e)_s=0,\quad L_s=K_e,\quad m_s=-2L,\quad
\chi_s=2K_e-\eta(3+u^2),\qquad u_s=v^2>0.
\end{equation}
At a $02$ tie, $\chi=0$; its maximizing signs are $\chi\ge0$ on mode $0$ and $\chi\le0$ on mode $2$. Set
\[
B_e=K_e-\eta(3+u^2)/2,\qquad w=S/v^2.
\]
Then $\sigma_0=w\chi$, $B_e'=-\eta uv$, and $w'=2/v$ on mode $0$, $w'=-2/v$ on mode $2$.

\begin{lemma}[Complete binary action and transport]\label{cmp:action}
On a binary block of one full active-zero lift, retain its objective
multiplier $\eta$, write $K=K_e$, and put
\[
 P_\eta=\eta u+2L-\chi,\qquad
 \mathsf B=\begin{pmatrix}L&K\\m&-L\end{pmatrix},\qquad
 N(S)=\begin{pmatrix}1&0\\S&1\end{pmatrix}.
\]
For physical cost $A$ and shear $S$ from a fixed initial point,
\begin{equation}\label{cmp:action-transport}
 2\eta A=\Delta P_\eta=\eta\Delta u+2KS-\Delta\chi,
 \qquad \mathsf B(S)=N(S)^{-1}\mathsf B(0)N(S).
\end{equation}
At an active mode, set $e=1+u$ for mode $0$ and $e=u-1$ for mode $2$.
The complete reconstruction is
\begin{equation}\label{cmp:full-reconstruction}
 \begin{gathered}
 l=L+um,\quad b=K-2uL-u^2m,\quad p_v=\chi/v,\\
 p_u=-\chi/e-\operatorname{tr}((\Lambda-\eta A_*)X)/v.
 \end{gathered}
\end{equation}
On a single mode, put $\varepsilon=1$ for mode $0$ and $-1$ for mode $2$;
$e=u+\varepsilon$ and the constant $\kappa=v/e\ne0$ are its ray data.
With arbitrary constants $K,B,M,C$, the same solution is
\begin{equation}\label{cmp:ray-primitive}
 \begin{gathered}
 L=B-\frac K{\kappa^2e},\qquad
 m=M+\frac{2B}{\kappa^2e}-\frac K{\kappa^4e^2},\\
 \chi=C+\frac{(4\eta-2K)/e+2\eta\varepsilon\log|e|-\eta e}{\kappa^2}.
 \end{gathered}
\end{equation}
On one constant-mode arc the four arbitrary initial data
$(K,L,m,\chi)$, equivalently $(K,B,M,C)$, give every active-zero
full adjoint, including singular matrix adjoints. Across a binary switch the actual tie $\chi=0$ is
additionally retained. These identities neither prove normality nor
identify lifts across a global seam.

If the lift is normal, after joint normalization set $P=P_1$.
It is continuous at actual
binary switches. Under a lawful chart change evaluate it on the whole
transported lift; under reversal use $P_\gamma=\epsilon_\gamma
(u_\gamma+2L_\gamma-\chi_\gamma)$, with $\epsilon_\gamma=-1$ for reversal
and $+1$ otherwise. For consecutive blocks with actual joins $q_i$,
\begin{equation}\label{cmp:action-cocycle}
 2\sum_i A_i=P_{\gamma_n}(q_n)-P_{\gamma_1}(q_0)
 +\sum_{i=1}^{n-1}[P_{\gamma_i}(q_i)-P_{\gamma_{i+1}}(q_i)].
\end{equation}
Thus for two overlapping blocks $2J=\Delta P_L+\Delta P_R-2c(I)$,
where $I$ is their shared physical arc. Different-chart endpoint
potentials need not agree.
\end{lemma}
\begin{proof}
Equation~\eqref{sw:binary} gives
$(P_\eta)_s=\eta(3+u^2+v^2)=2\eta v\ell$ and integrates the matrix
entries to $(L+KS,K,m-2LS-KS^2)$. For a single ray, substitute
$ds=de/(\kappa^2e^2)$ to obtain \eqref{cmp:ray-primitive}.
Reconstruction is invertible since
$v>0$. Conversely these formulas give the matrix adjoint and the
$p_v$ adjoint equation; differentiating the imposed active zero energy
then gives $v(p_u'+H_u)=0$, hence the remaining state-adjoint equation.
At a binary switch $\chi=0$ makes the two reconstructions of $p_u$
agree. Continuity, oriented telescoping and inclusion--exclusion prove
the remaining assertions.
\end{proof}

\begin{lemma}\label{sw:ABA}
On a normal lift normalized to $\eta=1$, the second switch of $202$ has $u>0$, and the first switch of $020$ has $u<0$. Their cyclic versions give the corresponding strict signs of
\[
u^2+v^2-2u-3\quad\text{and}\quad3-u^2-v^2-2u.
\]
\end{lemma}
\begin{proof}
At the two switches of $202$, $u_1<u_2$ and $B_e(u_1)\ge0\ge B_e(u_2)$. If $u_2\le0$, then
$B_e(u_2)-B_e(u_1)=-(u_2-u_1)(u_2+u_1)/2>0$, a contradiction. The other assertion follows with the signs reversed. Finally
\[
\Re\mathcal R(W)=\frac{u^2+v^2-2u-3}{(1+u)^2+v^2},\qquad
\Re\mathcal R^2(W)=\frac{3-u^2-v^2-2u}{(1-u)^2+v^2}
\]
has positive denominators and transports the result legally.
\end{proof}

\begin{lemma}[Signed capped-pair comparison]\label{sc:pair}
On one selected nonzero full free-time maximizing lift, retain an actual
positive source $1\mid02\mid1$ or $1\mid20\mid1$, with complete interior
arcs and both actual third-mode ties. Set $\delta=1$ or $-1$, respectively.
At its middle put
\[
 x=(1+u)/v,\quad y=(u^2+v^2-1)/(2v),\quad z=(1-u)/v,
 \qquad xy+xz+yz=1.
\]
For the positive physical mode-$0,2$ durations $d_x,d_z$, set
$P=e^{\delta d_x/x}$, $R=e^{\delta d_z/z}$. These are signed endpoint
parameters of the actual arcs; their durations remain positive.
With $U=\operatorname{tr}((\Lambda-\eta A_*)X)$, put
$\nabla_H=v\nabla$ and $V=\nabla_HU+J(vp)$, fixing $\Lambda,\eta$ in the gradient, with
$J(a,b)=(-b,a)$. At the middle $vp=(h,0)$, $h\ge0$, $U=-h$ and
$\delta V_y\ge0$. Define
\[
 \nu_x=\frac{(x^2-1,2x)}{1+x^2},\quad
 \nu_z=\frac{(z^2-1,-2z)}{1+z^2},\quad f(g)=\frac{g\log g}{g-1}-1,
\]
\[
 A_x=\frac{(P-1)(1+x^2)}2,\quad
 \beta_x=x\left[\frac{x(P-1)}{x+y}-f(P)\right],\quad
 b_x=2x\left[\frac{\log P}{P-1}-\frac{z+x/P}{x+z}\right],
\]
and define $A_z,\beta_z,b_z$ by $(x,z,P)\leftrightarrow(z,x,R)$.
The unnormalized actual cap rows are
\begin{equation}\label{sc:rows}
 h=A_x(V\cdot J\nu_x+2\eta\beta_x)
   =A_z(V\cdot J\nu_z+2\eta\beta_z).
\end{equation}
Put
\begin{gather*}
 D=x(P-1)+z(R-1)-(x+y)f(P)-(y+z)f(R),\\
 q_0=\frac{2xzD}{y(x+z)},\qquad
 r=\frac{2\{z/(P-1)+x/(R-1)\}}{y(x+z)^2}.
\end{gather*}
Then $V_y=\eta q_0-rh$, with $r>0$ if $\delta=1$ and
$r<-2/(1-xz)<-2$ if $\delta=-1$. The selected lift is normal.
After joint normalization, its actual cap lies in the algebraic segment
\begin{equation}\label{sc:interval}
 0\le h\le\bar h:=\frac{\delta q_0-M}{\delta r}\le h_*:=q_0/r,
 \qquad M=\max(0,\delta b_x,\delta b_z),\qquad \delta q_0\ge M.
\end{equation}
Writing $\ell=x+y+z$ and
$N=((x-z)/(x+z),y-(x^2+z^2)/(x+z))$, the exact polarization is
\begin{equation}\label{sc:polarization}
 \mathcal C=\frac34-\frac{N\cdot V}2-\frac{(\ell-N_y)h}2
 -\frac{|V|^2}4+\frac{hV_y}2,\qquad
 \operatorname{tr}(\Lambda X)=\ell-h.
\end{equation}
It is strictly concave on the cap line, so the actual determinant is
bounded below by the minimum of its values at $0,\bar h$, or more coarsely
by the minimum at $0,h_*$.
In the latter pair $hV_y=0$. These are algebraic comparison endpoints,
not independently normalized or calibrated trajectories.
The sum $I$ of the two complete physical costs is
\begin{equation}\label{sc:cost}
 \begin{split}
 I=\delta\bigl[&(x+z)D+f(P)+f(R)+x^2f(P)/P+xy(1-P^{-1})\\
 &+z^2f(R)/R+zy(1-R^{-1})\bigr].
 \end{split}
\end{equation}
After multiplication by $\delta$ every term is nonnegative and the six
non-residual terms are positive.
\end{lemma}
\begin{proof}
The binary tie gives $p_v=0$ and the active zero energy gives $U=-h$.
The third Hamiltonian inequality gives $2hv^2/(u^2+v^2-1)\ge0$.
Write $K=K_e$. Differentiating $U$ at fixed $\Lambda,\eta$ and using that energy yields
\begin{gather*}
 K=\frac{\eta(3+u^2)-vV_y}{2},\quad L=\frac{V_x+\eta u}{2},\quad
 m=-\frac\eta2+\frac{h-V_y/2}{v},\\
 l=L+um,\quad b=K-2uL-u^2m,\quad p=(h/v,0).
\end{gather*}
Thus $\chi_t=-V_y$; its actual maximizing crossing gives $\delta V_y\ge0$.
The mode-zero and mode-two cap endpoints and signed shears from the middle are
\[
 (u_A,v_A)=((1+u)/P-1,v/P),\quad s_A=-x(P-1)/v,
\]
\[
 (u_B,v_B)=(1-(1-u)/R,v/R),\quad s_B=z(R-1)/v.
\]
Integrating \eqref{sw:binary} with $ds=x^2du/(1+u)^2$ or
$z^2du/(1-u)^2$ gives
\[
 m_j=m-2Ls_j-Ks_j^2,\qquad
 \chi_A=x(P-1)(V_y-\eta b_x),\quad
 \chi_B=-z(R-1)(V_y-\eta b_z).
\]
Eliminating $p_u$ from the active energy and the third-mode tie gives
\[
 m_j+\eta/2=\frac{2K-\eta(3+u_j^2)}{2v_j^2}
       +\frac12(e_j/v_j^2+e_j^{-1})\chi_j,
 \quad e_A=1+u_A,\ e_B=u_B-1.
\]
Substitution reduces the two residuals to $(x+z)/2$ times the corresponding
row of \eqref{sc:rows}. No multiplier was normalized in this calculation.
The normal determinant is $2y(x+z)^2/[(1+x^2)(1+z^2)]>0$.
Using
\[
 y(x+z)^2(1,0)+z(1+x^2)\nu_x+x(1+z^2)\nu_z=0
\]
therefore gives the stated solved line and signs of $r$.
If $\eta=0$, the signs $h\ge0$, $\delta V_y\ge0$ force $h=V_y=0$;
either cap row gives $V_x=0$. The displayed reconstruction annihilates
all matrix and state covectors. Homogeneous adjoint uniqueness propagates
zero on the original interval, contradicting nontriviality. Normalize jointly.
The actual endpoint signs $\chi_A\ge0$, $\chi_B\le0$ now give
\eqref{sc:interval}, including its zero-height cases without division by $q_0$.

Direct matrix expansion gives the general identity
\[
 \det\Lambda=\frac{3\eta^2}4+\frac{\eta(\ell U-N\cdot\nabla_HU)}2
 +\frac{U^2-|\nabla_HU|^2}4.
\]
Use $U=-h$, $\nabla_HU=V-(0,h)$ after normalization to obtain
\eqref{sc:polarization}. If $V_x=p_0+\alpha h$,
its quadratic coefficient is $-(\alpha^2+r^2+2r)/4<0$, proving both endpoint
comparisons. Finally, integrating the physical density gives
\begin{align*}
 I/\delta={}&x[(x+y)(1-P^{-1})+(x+z)(P-1)-x\log P]\\
 &+z[(y+z)(1-R^{-1})+(x+z)(R-1)-z\log R].
\end{align*}
For $\delta=-1$ both integrals have reversed endpoint order; multiplication
by $\delta$ restores the actual positive physical costs. Expanding $D$ and
using $xy+xz+yz=1$ proves \eqref{sc:cost}.
Since $f(1)=0$ by continuity and
$f'(g)=(g-1-\log g)/(g-1)^2>0$ for $g\ne1$,
$f(P),f(R),1-P^{-1},1-R^{-1}$ all have sign $\delta$; the residual sign
follows from \eqref{sc:interval}.
\end{proof}

\subsection{Normality of the selected six-arc lift}
\begin{theorem}\label{sw:normality}
If an interval contains six consecutive positive maximal arcs and avoids the three words of Lemma~\ref{sw:four}, every selected nonzero full free-time maximizing lift on that interval is normal.
\end{theorem}
\begin{proof}
Suppose $\eta=0$. On a binary interval, \eqref{sw:binary} reduces to $\chi_s=2K_e$, with $K_e$ constant. A tie with a third-mode cap gives, after eliminating $p_u$ from the active Hamiltonian,
\begin{equation}\label{sw:abnormal-cap}
m=\frac{K_e}{v^2}+A_j\chi,\qquad
A_0=\frac12\left(\frac{1+u}{v^2}+\frac1{1+u}\right),\quad
A_2=-\frac12\left(\frac{1-u}{v^2}+\frac1{1-u}\right).
\end{equation}
The signed capped-pair comparison, Lemma~\ref{sc:pair}, excludes every
abnormal four-letter $ABCA$ with three distinct labels: a cyclic relabeling
makes it $1021$ or $1201$, both covered before normalization.

Two other short arguments suffice. In every $ABA$ block, the two distinct zeros of affine $\chi$ imply $K_e=\chi=0$ on the full block. In $ABACA$, the $ABA$ and $ACA$ blocks overlap on a positive $A$ arc, where two independent Hamiltonian differences vanish. Indeed
\[
\det(F_0-F_2,F_1-F_2)=\frac{-4v^5}{(u-1)(u+1)S}\ne0.
\]
Thus $p=0$ on that arc. Its adjoint and active Hamiltonian give $U=U_u=U_v=0$. The linear map $\Lambda\mapsto(U,U_u,U_v)$ has determinant $4/v^2$, so $\Lambda=0$ as well. Finally, for a binary core with at least two switches and third-mode caps at both ends, $K_e=\chi=0$ throughout the core; the two cap equations give $m=0$ at distinct ends. Since $L$ is constant and $m_s=-2L$, $m=L=0$. This excludes $ABCBA$ and $ABCBCA$.

It remains a short literal classification. After adjacent equal letters, the three hard four-letter words, and the patterns $ABCA,ABACA,ABCBA$ are excluded, the possible five-letter words are
\[
02021,\quad02121,\quad10102,\quad10202,\quad21010,\quad21210.
\]
The first, third, and sixth have no extension satisfying the same restrictions. The others extend only to
$021210$, $102021$, and $210102$, respectively. Each is an $ABCBCA$ word. This list can also be checked by fixing the first letter and distinguishing whether the third equals the first; each subsequent letter has two possibilities and the displayed shorter patterns remove all omitted branches. Thus six positive arcs are impossible for an abnormal lift. Homogeneous uniqueness propagates the vanishing found in each case, contradicting nontriviality.
\end{proof}
No compactness or quantitative cost bound enters this theorem. In particular, it supplies normality on the very interval whose lift is later used.

\subsection{A six-letter exclusion without a span bound}
\begin{lemma}\label{sw:six}
The words $202101,010212,121020$ are impossible on a full maximizing lift with positive arcs and caps. The conclusion holds for arbitrary positive durations.
\end{lemma}
\begin{proof}
Consider $202101$ and its switches $t_1<\cdots<t_5$. Normality follows from Theorem~\ref{sw:normality}; normalize jointly. On the mode $2$ tail after the left $202$, Lemma~\ref{sw:ABA} gives $u>0$. Let $y=-\chi$ and $z=y'=-2B_e/v$. The rate $b=v/(1-u)>0$ is constant there, and
\[
z'-bz=2u>0,\qquad y(t_2)=0,\qquad z(t_2)\ge0.
\]
For $d_-=t_3-t_2>0$, integration of the strict differential inequality yields
\[
0<y(t_3)<\frac{1-e^{-bd_-}}b z(t_3).
\]
Since $w'/w=-c$ on mode $2$, where $c=2v/S$, at the central tie
\begin{equation}\label{sw:tail}
\sigma_0<0,\qquad
\sigma_0'<\left(c-\frac{b}{1-e^{-bd_-}}\right)|\sigma_0|.
\end{equation}
For the right tail use time reversal together with
\[
\psi(W)=\frac{3-\overline W}{1+\overline W},\qquad
Q=\frac12\begin{pmatrix}-1&3\\1&1\end{pmatrix},\quad Q^2=I.
\]
The transported lift is $p_{\rm new}=-(D\psi)^{-t}p$, $\Lambda_{\rm new}=-Q\Lambda Q$, $\eta_{\rm new}=\eta$. Substitution gives the label permutation $(1\ 2)$, $b(\psi W)=c(W)$, $c(\psi W)=b(W)$, and
\[
(\sigma_0,\sigma_1,\sigma_2)_{\rm new}
=(-\sigma_1,-\sigma_0,-\sigma_2)_{\rm old}.
\]
The reversal cancels the derivative sign. Applying \eqref{sw:tail} to the reversed right $101$ with $d_+=t_4-t_3$ gives
\[
\sigma_1'<\left(b-\frac{c}{1-e^{-cd_+}}\right)|\sigma_1|.
\]
At the central $21$ tie, $q=\sigma_1=-\sigma_0>0$. Therefore
\begin{equation}\label{sw:two-tail}
\sigma_0'+\sigma_1'<-\left(\frac{be^{-bd_-}}{1-e^{-bd_-}}+
\frac{ce^{-cd_+}}{1-e^{-cd_+}}\right)q<0.
\end{equation}
But the actual crossing has $\sigma_2'\le0$, whereas the zero-sum identity gives $\sigma_0'+\sigma_1'=-\sigma_2'\ge0$. Cyclic transport proves the two remaining words.
\end{proof}
The final two-tail contradiction applies at any actual $2\to1$ tie
of one full maximizing lift with $q=\sigma_1=-\sigma_0>0$ and the
two displayed tail inequalities for positive rates and durations;
it does not require either tail to have arisen from an $ABA$ block.

\subsection{Twisted counts and necessary run data}
On the extended trajectory let $w_i$ be the consecutive maximal-arc labels. A physical period has $B$ essential switches and
\begin{equation}\label{sw:twist}
w_{i+B}=w_i+1\pmod3.
\end{equation}
This counts arcs intrinsically on the twisted circle. Cutting can split at most one such arc, proving \eqref{fd:cut-count}. The increments $\delta_i=w_{i+1}-w_i$ lie in $\{1,2\}$.
If both occur, write the alternating positive and reverse run lengths as $a_i,b_i$, $1\le i\le M$. The hard four-letter words forbid a reverse run of length one: its surrounding increments $121$ are exactly a forbidden four-letter pattern. The six-letter exclusion forbids a reverse run of length three: its surrounding increments $12221$ are exactly the corresponding six-letter pattern. Thus
\begin{equation}\label{sw:necessary}
a_i\ge1,\qquad b_i\ge2,\quad b_i\ne3,\qquad
\sum_i(a_i+b_i)=B.
\end{equation}
If $L=\sum_i a_i$, the twist also gives $L+2(B-L)\equiv1\pmod3$.
Sorting $a_i$ and $b_i$ separately gives an unordered necessary key $(aa,bb)$. Actual incidence information is forgotten, so this is an upper cover, not a realization theorem. Bounds used for a key must hold for all of its actual occurrence orders.


\section{A harmonic calibration of the matrix costate}\label{sec:calibration}
Fix one normal full lift on a genuine fundamental interval, normalized jointly to $\eta=1$. Its Casimir
\[
\mathcal C=\det\Lambda=-\tfrac12\tr(\Lambda^2)
\]
is constant by the Lax equation. In this section assume $\mathcal C>0$ and $\tr(\Lambda X)>0$ at one actual point. Both hypotheses must be established before the following bound is used.

\begin{lemma}[The future sheet]\label{cal:future}
Set $g=D_*^{-1}GD_*$, $\widehat X=D_*^{-1}XD_*$, and $\widehat\Lambda=D_*^{-1}\Lambda D_*$. With $H=\diag(1,-1)$ and $S_0=\left(\begin{smallmatrix}0&1\\1&0\end{smallmatrix}\right)$, write
\[
\widehat X=xJ+r_1H+r_2S_0,\qquad
Y=-\widehat\Lambda/\sqrt{\mathcal C}=zJ+y_1H+y_2S_0.
\]
Then
\[
x=\frac{3+u^2+v^2}{2\sqrt3v}\ge1,\quad
x^2-|r|^2=z^2-|y|^2=1,\quad z\ge1,
\]
and, with $\kappa=2\sqrt{\mathcal C}$,
\begin{equation}\label{cal:hamilton}
q=zx-y\cdot r=\frac{\tr(\Lambda X)}\kappa\ge1,
\qquad h_p:=\max_jp^tF_j=\sqrt3x-\kappa q\ge0.
\end{equation}
\end{lemma}
\begin{proof}
The determinant identities are direct. On the sheet $z\le-1$, Lorentz Cauchy gives $zx-y\cdot r<0$, contradicting the specified positive trace. Continuity preserves the sheet $z\ge1$, and on this sheet Lorentz Cauchy gives $q\ge1$. Explicitly write $x=\cosh s$, $|r|=\sinh s$, $z=\cosh t$, $|y|=\sinh t$; then $q\ge\cosh(t-s)$.

The fields have the strictly positive equilibrium weights
\[
\omega_0=\frac{(1+u)(v^2+(1-u)^2)}{4v^2},\quad
\omega_1=\frac S{2v^2},\quad
\omega_2=\frac{(1-u)(v^2+(1+u)^2)}{4v^2}.
\]
They sum to one and $\sum_j\omega_jF_j=0$. Therefore $h_p\ge0$, and the zero maximum Hamiltonian gives the equality in \eqref{cal:hamilton}.
\end{proof}

\begin{lemma}[The real rotation identity]\label{cal:rotation}
The legal three-copy extension of this same matrix lift satisfies
\begin{equation}\label{cal:rotation-eq}
\int_0^{3T}x\,dt=\mathcal A,\qquad
\pi=\int_0^{3T}\frac{x+q}{z+1}\,dt.
\end{equation}
No periodicity of the state costate is assumed.
\end{lemma}
\begin{proof}
The Lax equation gives $(g\widehat\Lambda g^{-1})'=0$, hence
$\widehat\Lambda(T)=R^{-1}\widehat\Lambda(0)R$ with $R=e^{\pi J/3}$.
Extend $\widehat\Lambda(t+kT)=R^{-k}\widehat\Lambda(t)R^k$. This is its continuous Lax transport on the extended state trajectory. Since $R$ commutes with $J$, the scalars $x,z,q$ are $T$-periodic. The scalar Hamiltonian inequality extends in this way without extending $p$. The area formula gives the first identity.

Choose $S\in\mathrm{SL}_2(\R)$ with $Y(0)=SJS^{-1}$, and set $f=g_e^{-1}S$. Then $Y=fJf^{-1}$, $f'=-\widehat X_ef$, and $f(3T)=-f(0)$. The original real homotopy class implies that the Iwasawa angle of $f$ changes by $-\pi$. To see the exact integer, homotope $g_e$ to $e^{\pi Jt/(3T)}$, invert and multiply by $S$, and then move $S$ continuously to $I$. The odd multiple of $\pi$ cannot change; at $S=I$ it is $-\pi$.

Write
\[
f=\begin{pmatrix}1&\xi\\0&1\end{pmatrix}
\begin{pmatrix}\sqrt h&0\\0&1/\sqrt h\end{pmatrix}e^{J\theta},\quad h>0.
\]
The closed hyperbolic projection has $\xi+ih=(y_1+i)/(z+y_2)$. Its connection form is
\[
-\tfrac12\tr(Jf^{-1}df)=d\theta-\tfrac12d\xi/h.
\]
On the positive hyperboloid, introduce the globally smooth one-form
\[
\alpha=-\frac{y_1dy_2-y_2dy_1}{z+1}.
\]
Direct differentiation gives
\[
\frac{d\xi}{h}-\alpha=d\left(2\arctan\frac{y_1}{1+z+y_2}\right).
\]
The denominator is positive, so the potential is a single-valued real function, including at the hyperboloid's pole. Its integral on the closed projection is zero. Since the left side of the connection formula equals $-q\,dt$ on $f'=-\widehat Xf$, integration yields
\[
\int_0^{3T}q\,dt=\pi+\tfrac12\int_0^{3T}\alpha(Y')\,dt.
\]
The Lax equations in these coordinates are
\[
z'=-2(y_1r_2-y_2r_1),\quad y_1'=2xy_2-2zr_2,
\quad y_2'=2zr_1-2xy_1.
\]
Consequently
\[
q-\tfrac12\alpha(Y')=\frac{x+q}{z+1}.
\]
This proves the second identity. The signed area contribution has been retained, so self-intersection or negative oriented area causes no difficulty.
\end{proof}

\begin{theorem}[Harmonic Casimir bound]\label{cal:harmonic}
For the chosen lift,
\begin{equation}\label{cal:bound}
\mathcal C\le\frac34\left(\frac{2\mathcal A}{\pi}-1\right)^2.
\end{equation}
For a global minimizer this implies
\[
\mathcal C<H_*:=\Hcap<H_0:=\Hzero.
\]
\end{theorem}
\begin{proof}
Let $w=y_1r_2-y_2r_1$. The Euclidean identity
$|y|^2|r|^2=(y\cdot r)^2+w^2$ and the two hyperboloid equations give the exact relation
\[
2xq(z+1)-(x+q)^2=z^2-1+w^2\ge0.
\]
All denominators below are positive. The function $q\mapsto2xq/(x+q)$ is increasing, and \eqref{cal:hamilton} gives $q\le\sqrt3x/\kappa$. Thus
\[
\frac{x+q}{z+1}\le\frac{2xq}{x+q}
\le\frac{2\sqrt3}{\sqrt3+\kappa}x.
\]
Integrating \eqref{cal:rotation-eq} gives $\pi\le2\sqrt3\mathcal A/(\sqrt3+\kappa)$, or
$\kappa\le\sqrt3(2\mathcal A/\pi-1)$. Its right side is positive because $\kappa>0$, so squaring is legitimate. Using \eqref{fw:budget} then yields
\[
\mathcal C<\frac{3(2A_0-p_0)^2}{4p_0^2}=H_*.
\]
The final comparison with $H_0$ is exact rational arithmetic. For the
weaker bound, $(x+q)/(z+1)\le(x+q)/2$ and
$q\le\sqrt3x/\kappa$ give
$\kappa(2\pi-\mathcal A)\le\sqrt3\mathcal A$ after integration.
Lemma~\ref{fw:balanced-reconstruction} gives a circumscribed hexagon
of area $2\sqrt3$, so $\mathcal A\le2\sqrt3<2\pi$. Dividing and
squaring gives $\mathcal C\le3\mathcal A^2/[4(2\pi-\mathcal A)^2]$
and the stated $H_0$.
\end{proof}

\begin{corollary}[One-point calibration criterion]\label{cal:one-point}
Fix one complete normal maximizing free-time lift on a genuine fundamental
interval of the original problem, with its original state and group
endpoints and prescribed real winding, and normalize the whole multiplier
jointly to $\eta=1$. At one actual point, in a lawful whole chart, write
$\Lambda=\left(\begin{smallmatrix}l&b\\m&-l\end{smallmatrix}\right)$
and put $K=2ul+b-u^2m$, $L=l-um$.
If $\mathcal C=\det\Lambda>0$, then
\[
 K>0\quad\Longleftrightarrow\quad m<0
 \quad\Longleftrightarrow\quad\operatorname{tr}(\Lambda X)>0.
\]
Any one of these signs licenses \eqref{cal:bound} on this same lift;
for an original global minimizer it gives $0<\mathcal C<H_*$.
For such a global minimizer, an independently proved floor
$\mathcal C\ge H_*$ contradicts $K>0$ on that lift. No assertion is made for an arbitrary local extremal.
\end{corollary}
\begin{proof}
Direct expansion at the specified actual state gives
$\mathcal C=-L^2-Km$ and
$\operatorname{tr}(\Lambda X)=(K-v^2m)/v$, independently of the active
control. Since $v>0$ and $Km=-\mathcal C-L^2<0$, the three signs are
equivalent. Lemma~\ref{cal:future} and Theorem~\ref{cal:harmonic} now
apply on the original interval. Their matrix-only three-copy transport
retains the real winding and assumes no periodic state covector.
\end{proof}

The preceding proof is valid for measurable maximizing controls. All derivatives are identities of absolutely continuous functions and almost-everywhere equations; none differentiates a control. Every application must establish normality, positivity of the determinant, and one positive trace for the same selected lift. A scalar local cost conclusion may be transferred to a global sum, but its costate and Casimir may not be transferred to another lift.

\subsection{Proof of finite switching}\label{fw:finiteness-proof}
\begin{proof}[Proof of Theorem~\ref{fw:finiteness}]
Fix an original global minimizer with a measurable control and select
one nonzero full maximizing lift on its genuine fundamental interval.
By continuity, $\max_j H_j=0$ holds at every time. Suppose that
$p(t_0)=0$ at some time $t_0$. Then $U(t_0)=0$. The map from the
matrix coordinates $(Z_{11},Z_{12},Z_{21})$ to $(U,U_u,U_v)$,
with the state fixed, has determinant $4/v^2$ by \eqref{fd:fields}.
If also $\nabla_W U(t_0)=0$, then $Z(t_0)=0$. In this case $\eta>0$,
since for $\eta=0$ the full homogeneous adjoint equations with
$p(t_0)=0$ and $\Lambda(t_0)=0$ would give the zero lift. After joint
normalization we would have $\Lambda(t_0)=A_*$, so
\[
 \det\Lambda(t_0)=\frac34,\qquad
 \operatorname{tr}(\Lambda(t_0)X(W(t_0)))=\ell(W(t_0))>0.
\]
Theorem~\ref{cal:harmonic} applies to this lift and gives
$\det\Lambda<H_*<3/4$. Thus $\nabla_W U(t_0)\ne0$ whenever $p(t_0)=0$.

We next prove that on each sufficiently short one-sided interval at
any time $t_0$, the control stays on a fixed edge of $\Delta$ almost
everywhere. The vectors $F_0,F_2$ have equal first coordinates and
distinct second coordinates, and $F_1$ has negative first coordinate.
Thus if $p(t_0)\ne0$, the three values $p(t_0)^t F_j(W(t_0))$
are unequal for at least one pair. A strict gap to the maximum
excludes a fixed vertex throughout a neighborhood of $t_0$.
If $p(t_0)=0$, integrating the adjoint equation gives
\[
 p(t)=-(t-t_0)\nabla_W U(t_0)+o(|t-t_0|).
\]
Indeed, the term linear in $p$ tends to zero, uniformly over the
controls, and $\nabla_W U$ is continuous. On each side the leading
nonzero covector gives a strict gap between a fixed pair of vertices.
The gap persists for sufficiently small $|t-t_0|$ on that side.
Maximization excludes the lower vertex of that pair there. The two
sides may use different edges.

Apply the first assertion of Lemma~\ref{sw:four} to the resulting closed
one-sided windows. At each $t_0$ their union gives a relative neighborhood
with finitely many vertex pieces, allowing a switch at $t_0$ itself.
A finite subcover of the compact fundamental interval gives a finite
vertex-control representation of the original trajectory.
\end{proof}


\space
\section{The short-chain domain and word list}
\label{short:integration}
We identify the short-chain variational domain with the normalized
minimizing class. The identification and the word list are used in
Proposition~\ref{short:rigidity}.

\subsection{Identification of the geometric and variational domains}
\label{short:interfaces}
The comparison constant is
\[
 J_8=C_{\rm oct}=\frac{2(8-4\sqrt2-\log2)}{2\sqrt2-1},\qquad
 U=C_*=\frac{180483}{100000},\quad V=\frac{361}{200},\quad
 K_{\rm sh}=\frac{29}{16}.
\]
The exact inequalities are $J_8<U<V<K_{\rm sh}$. In the explicit
closure calculations below we write $\mathscr J=C$ and
$\mathscr J_8=J_8$. The four-arc variation uses a scalar edge cost $J$
that differs from area by a positive factor and endpoint-fixed terms.

\begin{lemma}[Compatibility of the short-chain proof]\label{short:interface}
The normalized trajectory of an original global minimizer satisfies the geometric,
endpoint, winding and minimizing hypotheses of the short-chain
argument. Its essential count is the same $B$ in both arguments.
The complete-PMP statements used in the short proof may be applied
on their actual containing intervals without matching independently
chosen costates through the twisted seam.
\end{lemma}
\begin{proof}
The diagonal-conjugate coordinates used in the short-chain proof give
\[
\begin{gathered}
 X(u,v)=\frac1v\begin{pmatrix}u&-u^2-v^2\\1&-u\end{pmatrix},
 \qquad R_* =\frac12\begin{pmatrix}1&-3\\1&1\end{pmatrix},\\
 -1<u<1,\quad v>0,\quad u^2+v^2>1.
\end{gathered}
\]
These are the data in Section~\ref{sec:framework}. The group endpoint
is $G(T)=R_*$ and the state endpoint is
$W(T)=R_*^{-1}\!\cdot W(0)$, rather than merely a trace condition or
a state-return condition. Both use the fixed-endpoint real winding
class. The scalar normalization follows from \eqref{fd:cost}.
The physical circle-coordinate labels $0,2,4$ name the curved
opposite pairs. They are twice the reduced labels after one cyclic
phase, since mode $j$ curves the pair with index $2(j-1)$ in
Section~\ref{sec:framework}. Only the edge variation and the octagon
geometry retain physical labels. The closure charts use reduced labels.

The finite positive bang representation is supplied by
Theorem~\ref{fw:finiteness}. Zero-duration pieces are deleted and
equal adjacent controls are merged, including the twisted seam,
before either argument counts $B$. Both extend the reduced labels by
$j_{i+B}=j_i+1$ modulo three. Consequently their increment strings
have the same orientation and satisfy $\sum_i\epsilon_i\equiv1$
modulo three. Only cyclic rotation is used to identify the short
representatives below.

The admissible competitors are those of the original full control
triangle, with free terminal time and variable entrance state. No
upper bound on the number of arcs is imposed on that competitor
class. Lemmas~\ref{fw:replacement} and~\ref{fw:reanchor} apply to
fixed-endpoint local windows and genuine reanchoring, respectively.
For the six-arc arguments the variation is instead a whole-cycle
family satisfying exactly
\[
 G(0)=I,\quad G(T)=R_*,\quad
 W(T)=R_*^{-1}\!\cdot W(0),
\]
with $W(0)$ allowed to vary as its formulas prescribe. These are
the original endpoint conditions; they are not a claim that both
state endpoint values remain fixed. After rescaling the time
interval, the pointwise quotient of sufficiently close group paths
has a continuous real matrix logarithm, providing an endpoint-fixed
homotopy of the group paths. Thus the real winding class is
preserved. Positive durations and strict-star admissibility persist
by openness. Lemma~\ref{fw:balanced-reconstruction} gives the nearby
convex disk with the prescribed least circumscribed hexagon area.
Thus area comparisons preserve the density comparison even though
the entrance state varies.

The oriented four-arc variation proved in
Theorem~\ref{sh:edge-descent} concerns exactly the reduced words
$0101,1212,2020$.  The six-letter exclusion is
Lemma~\ref{sw:six}.  For the pure-forward state filter we use
Lemma~\ref{fin:forward-force} on an actual containing interval.
At the actual middle switch, its
reciprocals are
\[
 d_0=a_i,\quad d_2=a_{i-1},\quad
 d_1=\frac{1-a_i a_{i-1}}{a_i+a_{i-1}}.
\]
Thus $d_1<\min(d_0,d_2)$ is exactly
\[
 a_i^2+2a_i a_{i-1}>1,\qquad
 a_{i-1}^2+2a_i a_{i-1}>1.
\]
Strict star gives $0<a_i a_{i-1}<1$. Each capped $1201$ support is
placed in one genuine fundamental interval. Different cuts may
supply different full lifts; only the resulting scalar state
inequalities are combined. No periodic state-costate assumption is
made. This verifies the interfaces needed below.
\end{proof}

\subsection{Exhaustion of short words}\label{short:words}
For $B\ge4$, any cyclic factor $+-+$ uses four distinct physical
arcs and is prohibited by Lemma~\ref{sw:four}. For $B=3$ the formal
factor repeats an arc and must not be used as a local exclusion.
Keeping that case for a separate geometric argument gives the
following complete list, up to rotation:
\begin{center}
\begin{tabular}{cl}\toprule
$B$ & Representatives before the six-arc exclusion\\\midrule
1 & $\texttt{+}$\\
2 & $\texttt{--}$\\
3 & $\texttt{++-}$\\
4 & $\texttt{++++},\ \texttt{+---}$\\
5 & $\texttt{+++--},\ \texttt{-----}$\\
6 & $\texttt{++----},\ \texttt{+--+--}$\\
7 & $\texttt{+++++++},\ \texttt{++++---},\ \texttt{+------}$\\
\bottomrule
\end{tabular}
\end{center}
Indeed if $k$ signs are negative, $B-2k\equiv1\pmod3$.
For $B\ge4$ distribute the negative signs into cyclic runs, with
no singleton negative run between positive signs. Listing the
possibilities gives the table. The marked-string counts before
forbidden factors are $1,1,3,5,11,21,43$. The exact elementary
enumerator supplied as \path{short/check_short_words.py} checks
these counts, all representatives, and the six-arc placement below.
It does not purport to prove geometric feasibility of a word.

\space
\subsection{A four-arc variation with fixed complete endpoints}

The reduction of short chains begins with a variation that preserves both
components of the endpoint: the position matrix and the reduced state.
In the four-arc variation the physical labels of active opposite pairs are
$0,2,4$; dividing these labels by two gives the reduced labels $0,1,2$.
The sign of a transition is read in the fixed cyclic order of the controls.
The reversed oriented edge is a different variational problem.

\begin{theorem}[Oriented four-arc descent]\label{sh:edge-descent}
An actual positive four-arc window with the oriented physical word
$0,2,0,2$ has an arbitrarily small, strictly area-decreasing variation
with the same word and the same complete entrance and exit states.
Consequently, a cyclic sign factor $+-+$ cannot occur in a minimizing
chain when its four physical arcs are distinct.
\end{theorem}

Suppose that the physical word is $0$--$2$--$0$--$2$.  The two controls
lie on one edge of the control simplex.  After a cyclic relabeling of the
three control coordinates, the word is $A$--$B$--$A$--$B$, where the edge
parameter of \cite[Section~8.1]{HV2024} takes the values
\[
 u_{\rm edge}(A)=\frac12,
 \qquad
 u_{\rm edge}(B)=-\frac12.
\]
In the edge coordinates, put
\begin{equation}\label{sh:eq:edge-basic-coordinates}
 q=\frac1y,
 \qquad a=\sqrt3,
 \qquad y>0.
\end{equation}
We retain throughout the full strict-star condition
\begin{equation}\label{sh:eq:edge-full-star}
 |x|<\frac1{\sqrt3},
 \qquad
 x^2+y^2>\frac13.
\end{equation}
Only the first inequality is needed for the signs below; the second is part
of admissibility and persists by openness.  Fix the full entrance state
$z_-=(g_-,X_-)$, and write $x_0+iy_0$ for the upper-half-plane coordinate
of $X_-$.  Thus $g_-$, $x_0$, and
\[
 q(x_0)=q_0:=y_0^{-1}
\]
are fixed throughout the variation.
Since $x'=y>0$, we use $x$ as the independent variable.  The two bang
equations become
\begin{equation}\label{sh:eq:edge-bang-equations}
 \frac{d\log q}{dx}=p_A(x)=\frac{a}{1-ax},
 \qquad
 \frac{d\log q}{dx}=p_B(x)=-\frac{a}{1+ax}.
\end{equation}
Define
\begin{equation}\label{sh:eq:edge-rho-h}
 \rho=q^2,
 \qquad
 h(x)=p_A(x)-p_B(x)=\frac{2a}{1-a^2x^2}>0.
\end{equation}
Write $s$ for the shear state and
\[
 S:=\Delta s=s(x_4)-s(x_0)
\]
for its increment across the four-link window.  Since $x'=y$ and
$q=1/y$, we have
\begin{equation}\label{sh:eq:edge-shear-cost}
 ds=\rho\,dx,
 \qquad
 S=\int_{x_0}^{x_4}\rho\,dx,
 \qquad
 dJ=x^2\rho\,dx.
\end{equation}
Here $J$ is the edge cost after removal of a boundary term fixed by the
endpoint data; it differs from the Reinhardt area by a positive constant
factor and fixed endpoint terms
\cite[Eqs.~(8.1.3)--(8.1.5)]{HV2024}.

Fix entrance and exit coordinates $x_0<x_4$, and let
\[
 x_0<x_1<x_2<x_3<x_4
\]
be the switching positions.  For the unperturbed trajectory, put
\begin{equation}\label{sh:eq:edge-tail-masses}
 R_i=\int_{x_i}^{x_4}\rho\,dx,
 \qquad
 m_L=R_1-R_2,
 \qquad
 m_R=R_2-R_3.
\end{equation}
Thus $m_L,m_R>0$ are the $\rho$-masses of the two middle links.

Let $\xi_i=\dot x_i$, $h_i=h(x_i)$, and $\rho_i=\rho(x_i)$.  Moving a
switch changes $\log q$ on the entire downstream tail.  The three
downstream increments are
\begin{equation}\label{sh:eq:edge-switch-jumps}
 (b_1,b_2,b_3)
 =(h_1\xi_1,-h_2\xi_2,h_3\xi_3).
\end{equation}
The linearized conditions that the terminal values of $q$ and $S$ remain
fixed are
\begin{equation}\label{sh:eq:edge-linear-constraints}
 b_1+b_2+b_3=0,
 \qquad
 R_1b_1+R_2b_2+R_3b_3=0.
\end{equation}
Because $R_1>R_2>R_3$, these equations have rank two.  Their kernel is
\begin{equation}\label{sh:eq:edge-kernel-b}
 (b_1,b_2,b_3)
 =c_0\bigl(m_R,-(m_L+m_R),m_L\bigr),
\end{equation}
or, equivalently,
\begin{equation}\label{sh:eq:edge-kernel-xi}
 (\xi_1,\xi_2,\xi_3)
 =c_0\left(
 \frac{m_R}{h_1},
 \frac{m_L+m_R}{h_2},
 \frac{m_L}{h_3}
 \right).
\end{equation}

\begin{lemma}[Exact edge fiber]
\label{sh:lem:exact-edge-fiber}
Near the given four-link chain, the $A$--$B$--$A$--$B$ trajectories with
the fixed full entrance state and the same exit data $(S,x_4,q_4)$ form a
one-dimensional analytic manifold.  Equivalently, they have the same full
exit state $(g_4,X_4)$.  The tangent is
\eqref{sh:eq:edge-kernel-xi}, and both local directions keep all four links
positive.
\end{lemma}

\begin{proof}
Let $\mathbf x=(x_1,x_2,x_3)$, and let $q_{\mathbf x}$ solve the four bang
equations with $q_{\mathbf x}(x_0)=q_0$.  On the ordered open set
\[
 \Omega=\{x_0<x_1<x_2<x_3<x_4\},
\]
define the exact endpoint map
\[
 F(\mathbf x)=\left(
  \log q_{\mathbf x}(x_4),
  \int_{x_0}^{x_4}q_{\mathbf x}(x)^2\,dx
 \right).
\]
It is analytic.  If $d=(h_1,-h_2,h_3)$, direct differentiation gives
\begin{equation}\label{sh:eq:edge-exact-jacobian}
 DF(\mathbf x)=
 \begin{pmatrix}
  d_1&d_2&d_3\\
  2R_1d_1&2R_2d_2&2R_3d_3
 \end{pmatrix}.
\end{equation}
Since $R_1>R_2>R_3$, this matrix has rank two.  Hence the analytic
implicit function theorem makes $F^{-1}(F(\mathbf x^*))$ an exact
one-dimensional analytic manifold, not merely a linearized constraint
set.  Its tangent is \eqref{sh:eq:edge-kernel-xi}.  Because $\Omega$ is open,
both signs of a sufficiently small parameter retain all four positive
links.  The passage from the four flow durations to the four successive
$x$-endpoints is triangular with positive diagonal entries $x'=y$; it is
therefore a local analytic change of coordinates.

It remains to check that these edge data are precisely the full-state
endpoint data.  In the notation of \cite[Eq.~(8.1.2)]{HV2024},
\[
 g(t)=g(t_0)h(t_0)^{-1}h(t),
 \qquad
 h(s,x)=
 \begin{pmatrix}
  1&-x\\
  s&1-sx
 \end{pmatrix},
\]
and the additive constant in $s$ is a gauge.  More explicitly, if
$s_i=s(x_i)$ and $S=s_4-s_0$, then
\begin{equation}\label{sh:eq:edge-relative-position}
 h(s_0,x_0)^{-1}h(s_4,x_4)
 =
 \begin{pmatrix}
  1+x_0S&x_0-x_4-x_0x_4S\\
  S&1-x_4S
 \end{pmatrix}.
\end{equation}
Thus
\[
 g_4=g_-h(s_0,x_0)^{-1}h(s_4,x_4)
\]
is determined exactly by the fixed entrance data and $(S,x_4)$.  The
reduced state $X_4$ is determined by $x_4+iy_4=x_4+i/q_4$.  Consequently,
after the full entrance state is fixed,
\[
 (S,x_4,q_4)
 \quad\longleftrightarrow\quad
 (g_4,X_4)
\]
are equivalent endpoint data on this fixed control edge.  The endpoint image of edge-confined trajectories is three-dimensional;
its four-duration fixed-endpoint fiber therefore has dimension one.
\end{proof}

We now compute the area along this exact fiber.  Define the weighted
means
\begin{equation}\label{sh:eq:edge-weighted-means}
 \mu_L=\frac1{m_L}\int_{x_1}^{x_2}x^2\rho\,dx,
 \qquad
 \mu_R=\frac1{m_R}\int_{x_2}^{x_3}x^2\rho\,dx.
\end{equation}

\begin{lemma}[First and second edge variations]
\label{sh:lem:edge-variations}
Along the tangent \eqref{sh:eq:edge-kernel-xi} with $c_0\ne0$,
\begin{equation}\label{sh:eq:edge-first-variation}
 \dot J=2c_0m_Lm_R(\mu_L-\mu_R).
\end{equation}
If this derivative vanishes, write
$\mu_L=\mu_R=\lambda$.  The constrained second variation is
\begin{equation}\label{sh:eq:edge-second-variation}
\begin{split}
 Q=-2c_0^2\Bigg[&
 \frac{\rho_1}{h_1}(x_1^2-\lambda)m_R^2
 +\frac{\rho_2}{h_2}(\lambda-x_2^2)(m_L+m_R)^2\\
 &+\frac{\rho_3}{h_3}(x_3^2-\lambda)m_L^2
 \Bigg]<0.
\end{split}
\end{equation}
\end{lemma}

\begin{proof}
Let
\[
 T_i=\int_{x_i}^{x_4}x^2\rho\,dx.
\]
A switch increment $b_i$ changes $\rho=q^2$ on its downstream tail by
$2b_i\rho$ to first order.  Consequently
\[
 \dot J=2\sum_{i=1}^3b_iT_i.
\]
Substitution of \eqref{sh:eq:edge-kernel-b} and grouping adjacent tails gives
\[
 \sum_{i=1}^3b_iT_i
 =c_0m_R(T_1-T_2)-c_0m_L(T_2-T_3)
 =c_0m_Lm_R(\mu_L-\mu_R),
\]
which is \eqref{sh:eq:edge-first-variation}.

Suppose now that $\mu_L=\mu_R=\lambda$.  For completeness, we derive the
quadratic form rather than differentiating a chosen parametrization of
the fiber.  Set
\[
 \mathcal L=J-\lambda S-\nu\log q(x_4),
 \qquad
 U_i=\int_{x_i}^{x_4}(x^2-\lambda)\rho\,dx.
\]
Switch stationarity is $2U_i=\nu$ for $i=1,2,3$.  If
$\sigma=(1,-1,1)$ and $d_i=\sigma_i h_i$, direct differentiation gives
the Hessian at a stationary switch configuration as
\begin{equation}\label{sh:eq:edge-hessian-matrix}
 D^2\mathcal L
 =2\nu\,d\,d^{T}
 -2\diag\bigl(
 d_i(x_i^2-\lambda)\rho_i
 \bigr)_{i=1}^3.
\end{equation}
The equality $\mu_L=\mu_R=\lambda$ gives
\[
 U_1-U_2=m_L(\mu_L-\lambda)=0,
 \qquad
 U_2-U_3=m_R(\mu_R-\lambda)=0.
\]
Hence $U_1=U_2=U_3$, and choosing $\nu=2U_i$ makes
$D\mathcal L=0$ in all three switch variables.  Let $\mathbf x(r)$ be an
analytic parametrization of the exact fiber with
$\mathbf x'(0)=\xi$.  Both $S$ and $\log q(x_4)$ are constant on this
curve, so $\mathcal L(\mathbf x(r))$ differs from $J(\mathbf x(r))$ by a
constant.  Since $D\mathcal L=0$,
\begin{equation}\label{sh:eq:edge-exact-fiber-second-derivative}
 \frac{d^2}{dr^2}J(\mathbf x(r))\bigg|_{r=0}
 =D^2\mathcal L[\xi,\xi].
\end{equation}
Thus the quadratic form below is the actual second derivative along an
exact fixed-endpoint curve; no acceleration or exterior correction remains.
On the terminal-$q$ tangent space, $d\cdot\xi=\sum b_i=0$, so the rank-one
term vanishes.  Inserting
\eqref{sh:eq:edge-kernel-b} into the diagonal term yields precisely
\eqref{sh:eq:edge-second-variation}.

It remains to prove the sign.  Equality of the two means forces
\begin{equation}\label{sh:eq:edge-mean-order}
 x_1^2>\lambda>x_2^2,
 \qquad
 x_3^2>\lambda.
\end{equation}
If $x_2\leq0$, then $x^2$ is strictly decreasing on $(x_1,x_2)$, giving
the first two inequalities.  If $x_3^2\leq\lambda$, convexity of $x^2$
would make $x^2<\lambda$ throughout $(x_2,x_3)$, contradicting its mean.
If $x_2\geq0$, the right interval first gives
$x_2^2<\lambda<x_3^2$; if $x_1^2\leq\lambda$, the same endpoint-maximum
argument on the left interval gives a contradiction.  Thus every term in
brackets in \eqref{sh:eq:edge-second-variation} is strictly positive.
\end{proof}

\begin{proof}[Proof of Theorem~\ref{sh:edge-descent}]
Choose a nonzero tangent to the exact fiber in
Lemma~\ref{sh:lem:exact-edge-fiber}. The resulting two-sided variation fixes
the complete endpoint. If \eqref{sh:eq:edge-first-variation} is nonzero,
choose its decreasing sign.  If it vanishes,
Lemma~\ref{sh:lem:edge-variations} gives a strictly negative constrained
second variation.  The four-link replacement may therefore be inserted
into any longer chain without moving its complement.

To spell out the concatenation, let $z_{\rm in}$ and $z_{\rm out}$ be the
full states at the ends of the window.  The exact fiber keeps both fixed.
The exterior flow therefore restarts from precisely the same state
$z_{\rm out}$, and uniqueness of the state equation makes the entire
downstream trajectory, including its final full state, unchanged.  The
window duration may change, but terminal time is free, so the exterior is
only translated in time.  For a twisted-closed trajectory this preserves
both the position condition on $g$ and the twisted terminal condition on
$X$.  Strict-star admissibility and embeddedness persist for a sufficiently
small $C^1$ replacement.  No additional global closure equation or
second-order compensator is needed.

A global minimizer
must minimize on every edge-restricted subinterval with the same endpoint
data \cite[Prop.~8.1.3]{HV2024}; hence it cannot contain the
factor.  Cyclic label symmetry supplies the other two physical shifts.
\end{proof}

\begin{remark}[The weighted form of pair compression]
The literal duration velocity $(1,-1,-1,1)$ suggested by ``shorten the
middle pair and lengthen the outer pair'' generally leaves a shear defect.
Formula \eqref{sh:eq:edge-kernel-xi} is its exact state-dependent correction:
the three switches move in the same $x$-direction, but at weights forced
by the two tail-mass constraints.  No link is added and no exterior
compensator is used.
\end{remark}

\begin{remark}[Why the reverse phase is different]
\label{sh:rem:reverse-edge-phase}
For the reverse word $B$--$A$--$B$--$A$, the three switch jumps are
$(-h_1\xi_1,h_2\xi_2,-h_3\xi_3)$.  At stationarity the quadratic form
\eqref{sh:eq:edge-second-variation} changes sign and is strictly positive.
The first-order alternative still gives descent away from stationarity,
but the fixed-edge calculation does not exclude a stationary reverse
phase.  Reflection does not identify the two oriented factors: restoring
the positive boundary parameter also reverses time, so transition strings
are reversed, not sign-complemented.
\end{remark}

\subsection{The finite list of short sign cycles}

Let $B$ denote the number of essential positive arcs, after zero arcs
have been deleted and equal adjacent controls have been merged, including
at the twisted seam.  If $a_i\in\mathbb Z/3\mathbb Z$ is the active
label, extend it by $a_{i+B}=a_i+1$.  Then
\begin{equation}\label{sh:twist}
 \epsilon_i=a_{i+1}-a_i\in\{1,-1\},\qquad
 \sum_{i=0}^{B-1}\epsilon_i\equiv1\pmod3.
\end{equation}
For $B\geq4$, a cyclic occurrence of $+-+$ uses four different physical
arcs.  Cutting immediately before its first arc puts it inside one
fundamental interval.  Theorem~\ref{sh:edge-descent} then applies to
that interval, with its actual twisted labels.  Its complement remains
fixed.  No continuation of a state costate through the seam is involved.

The complete list for $B\leq7$ is given in
Section~\ref{short:words}.  The three-arc representative is retained
there until the following geometric argument: a formal $+-+$ at
$B=3$ repeats a physical arc and cannot be excluded by the four-arc
variation.

\begin{lemma}\label{sh:omitted}
A regular strict-star twisted-closed chain cannot omit one active label.
The mixed four-arc sign cycle $+---$ has no such realization.
\end{lemma}
\begin{proof}
Write the six boundary curves in the circle representation as
$\sigma_j=g e_j$, where $e_j=(\cos(j\pi/3),\sin(j\pi/3))$.
If physical label $4$ is omitted, $\sigma_4$ is straight throughout a
fundamental interval.  Regular $C^1$ matching puts its pieces on one
line.  The twist gives $\sigma_4'(T)=\sigma_5'(0)=-\sigma_2'(0)$,
whereas strict star gives
\[
 \det(\sigma_2',\sigma_4')
 =\det X\det(e_2,e_4)=\frac{\sqrt3}{2}\det X>0.
\]
This is impossible.  The classes $+$ and $--$ omit labels.  For the
three-arc class, choose the cut giving the reduced word $121$; it omits
label $0$.

For the mixed four-arc class use physical labels $(0,4,2,0)$ and
extend by $a_{r+4}=a_r+2$.  Around the seam the labels are
\[
\begin{array}{c|rrrr|rrrr}
 r&-2&-1&0&1&2&3&4&5\\\hline
 a_r&0&4&0&4&2&0&2&0.
\end{array}
\]
Thus $\sigma_2$ is straight on $(q_{-2},q_2)$, and $\sigma_4$ is
straight on $(q_2,q_6)=(q_{-2},q_2)+T$.  The identity
$\sigma_4(t+T)=-\sigma_2(t)$ makes their carrying lines parallel.
Taking the common tangent limits at $q_2$ contradicts the same positive
determinant as above.
\end{proof}

\subsection{All coherent four-arc chains}

We give the entire four-arc comparison, because comparison on a symmetric
subfamily would not suffice for equality rigidity.  A coherent four-arc
chain has physical labels $(0,2,4,0)$ and eight curved boundary arcs.
Between them lie eight positive straight segments: each active opposite
pair contributes two curved arcs under the sixfold twist, and maximality
of the control intervals separates consecutive curved arcs by a positive
inactive interval.  Consecutive inactive pieces are collinear by $C^1$
matching.  Distinct maximal faces lie on distinct support lines, since a
face of a convex disk is connected.  The four opposite support pairs
therefore bound a proper centrally symmetric octagon $P$.

Let $S_0,S_1,S_2,S_3$ be its four strips in cyclic order. Number the
extended arcs so that $G_e(t)q_0$ is curved on arc $n$ exactly when
$n\equiv0\pmod3$. At the start of arc $n$ its tangent belongs to
strip $\lfloor(n-1)/3\rfloor$ modulo four. The opposite faces are
identified, and $G_e(t)q_j=G_e(t+jT)q_0$ shifts the arc index by $4j$.
Thus for $n=0,1,2,3$ the three support pairs with indices $j=0,1,2$
omit respectively strips $0,3,2,1$.

At each switch, pull back its six tangent lines by the common map
$G(t)^{-1}$. The line at $q_i$ has equation $B(q_i,z)=B(q_i,q_i)$,
where $B(p,z)=\det(Xp,z)$ is symmetric because $\tr X=0$.
The reconstructed disk from Lemma~\ref{fw:balanced-reconstruction}
contains $0$ in its interior. The six lines are distinct because
$X$ is invertible and opposite supports are separated by $0$.
Each $q_i$ has a unique supporting line, so the other five support
inequalities are strict there. The sides therefore occur in the cyclic
order of the oriented tangents $Xq_i$, each with $q_i$ in its relative
interior. This order is counterclockwise since $\det X=1$.
If $z_+$ is the intersection of the lines at $q_i$ and $q_{i+1}$,
then symmetry of $B$ and $q_{i-1}=q_i-q_{i+1}$ give
\[
 B(q_i,2q_i-z_+)=B(q_i,q_i),\qquad
 B(q_{i-1},2q_i-z_+)=B(q_{i-1},q_{i-1}).
\]
Thus $q_i$ is the midpoint of its side. For a centrally symmetric
hexagon centered at the origin, the shoelace formula gives the area as
four times the determinant of two consecutive side midpoints in counterclockwise
order. Here the area is $4\det(q_i,q_{i+1})=2\sqrt3$, and $G(t)$
preserves it. Hence all four switch hexagons
$H_r=\bigcap_{j\ne r}S_j$ have this area. Omitting a strip adds
two opposite exterior triangles, so the four triangles have the
same positive area.

Normalize alternating strip pairs to $|x|,|y|\leq1$, and write four
successive vertices as
\[
 (1,a),\quad(b,1),\quad(-c,1),\quad(-1,d).
\]
Set $\alpha=1-a$, $\beta=1-b$, $\gamma=1-c$, $\eta=1-d$ and
$D_0=\alpha\gamma+\beta\eta$.  Twice the four exterior triangle areas
are
\[
 \alpha\beta,\quad\gamma\eta,\quad
 \frac{(a+d)^2\beta\gamma}{D_0},\quad
 \frac{(b+c)^2\alpha\eta}{D_0}.
\]
Their equality, with all factors positive, gives
$\alpha\beta=\gamma\eta$ and $(a+d)\beta=(b+c)\eta$.
Adding the latter identity to $(1-a)\beta=(1-c)\eta$ yields $b=d$,
and then $a=c$.  The remaining equality is
\[
 (1-a)^2+(1-b)^2=(a+b)^2,
 \qquad 1=a+b+ab.
\]
Consequently $P$ is linearly equivalent to the intersection of two
congruent squares,
\begin{equation}\label{sh:octagon-polygon}
 P_p=\{(x,y):|x|,|y|,|cx+dy|,|-dx+cy|\leq1\},\quad
 c^2+d^2=1,\quad c,d>0,\quad p=c-d.
\end{equation}
Changing the marking replaces $p$ by $-p$, so take $0\leq p<1$ and put
$s=c+d=\sqrt{2-p^2}$.

The rounding is also forced. On curved arc $0$ the adjacent
components $j=1,5$ lie at extended arc indices $4,20$, hence on the two
nonincident strips $1,2$. The same index calculation shifts by one
strip for each succeeding curved arc. The two inactive lines are the
asymptotes of the curved component \cite[Lemma~3.1.8]{HV2024}.
At the corner $x=1$, $cx+dy=1$, these are $y=1$ and $dx-cy=1$.
In coordinates
\[
 U=1-y,\quad V=1-dx+cy,\quad A_c=1+c-d,
\]
the incident sides become $V+cU=A_c$ and $U+cV=A_c$.
A hyperbola $UV=h$ is tangent to both precisely when
\begin{equation}\label{sh:hyperbola}
 h=\frac{A_c^2}{4c},\qquad
 (U,V)=\left(\frac{A_c}{2c},\frac{A_c}{2}\right),\quad
 \left(\frac{A_c}{2},\frac{A_c}{2c}\right).
\end{equation}
Both tangencies lie inside their sides.  The other type interchanges $c$
and $d$.  The two neighboring tangencies on $x=1$ have ordinates
$(s-1)/(2c)$ and $-(s-1)/(2d)$, so every intervening straight segment is
positive.  The graph $V=h/U$ is strictly convex and tangent to both
incident lines, hence clips exactly the prescribed corner.  These eight
arcs form a convex regular $C^1$ boundary.

Its unnormalized balanced determinant is
\[
 D_p=\frac{1-h_c}{d}=\frac{1-h_d}{c}
     =\frac{s+3}{2(s+1)}.
\]
Scaling by $\lambda_p^2=\sqrt3(s+1)/(s+3)$ restores determinant
$\sqrt3/2$.  Integration of $h/U$ between the endpoints in
\eqref{sh:hyperbola} shows that a $c$-type corner removes area
\[
 E_c=\frac{(1+c-d)^2}{4cd}
 \left(\log\frac1c-\frac{2(1-c)}{1+c}\right),
\]
and $E_d$ is obtained by exchange.  The polygon has area $8/(s+1)$;
there are four corners of each type.  Multiplication by
$\lambda_p^2/\sqrt3$ gives the normalized cost
\begin{equation}\label{sh:octagon-cost}
 \mathscr J(p)=2\delta(p),\qquad
 \delta(p)=
 \frac{4s(s-1)-(1+p^2)\log\dfrac2{1-p^2}
       +4p\operatorname{artanh}(p/s)}{(s-1)(s+3)}.
\end{equation}
The positive determinant reconstruction also proves that the regular
member $p=0$ is an admissible competitor.  Uniqueness of the prescribed
hyperbolas means that no further rounding parameter has been suppressed.

\begin{proposition}\label{sh:octagon-rigidity}
Every coherent four-arc chain has $\mathscr J\geq\mathscr J_8$.
Equality holds precisely for an affine image of the smoothed regular
octagon.
\end{proposition}
\begin{proof}
The function $\delta$ is even.  For $0<p<1$ set
$\lambda=\log(2/(1-p^2))$, $a=\operatorname{artanh}(p/s)$ and
$D=(s-1)(s+3)$.  Differentiation gives
\[
 \delta'(p)=\frac{2E}{sD^2},\qquad
 E=2(4+s-s^3)a-p(s^2+3)\lambda+\frac{ps(s-1)^3}{s+1}.
\]
If
\[
 B_0=(1-c)(1+d+cd)\log(1/d)
      -(1-d)(1+c+cd)\log(1/c),
\]
then $E=4B_0+ps(s-1)^3/(s+1)$.  Put $t=d/c\in(0,1)$ and
$q=\sqrt{1+t^2}$, so $c=1/q$, $d=t/q$.  The assertion $B_0>0$ is
\[
 \frac{\log(1+t^2)}{\log(1/t)}
 <\frac{2(q^3-q^2+t^3)}{(1-t)(2q^2+t)}.
\]
The positive series for $\operatorname{artanh}$ gives
\[
 \log(1/t)>\frac{2(1-t)}{1+t},\quad
 \log(1+t^2)<U(t):=
 \frac{t^2(t^4+12t^2+12)}{6(1+t^2)(2+t^2)}.
\]
Also $q^3-q^2+t^3\geq t^2(t+2(1+t^2)/(4+t^2))$.
Subtract $U(t)(1+t)/(2(1-t))$ from the resulting lower bound for the
right side.  The difference is
\[
 \frac{t^3 P(t)}{12(1+t^2)(2+t^2)(1-t)(2+t+2t^2)(4+t^2)},
\]
where
\[
 P(t)=48-24t+60t^2-20t^3-2t^5-11t^6-3t^7-2t^8
      >24+22t^2>0.
\]
Thus $B_0>0$, $E>0$ and $\delta'(p)>0$.  At $p=0$,
$\delta(0)=(8-4\sqrt2-\log2)/(2\sqrt2-1)$.
The equality case follows from the complete polygon and rounding
classification above.
\end{proof}

\subsection{The two five-arc classes}

For the algebra below it is convenient to use a fixed diagonal conjugate
of the circle coordinates.  The reduced state and twist are
\begin{equation}\label{sh:uv}
 X(u,v)=\frac1v\begin{pmatrix}u&-u^2-v^2\\1&-u\end{pmatrix},\quad
 R_* =\frac12\begin{pmatrix}1&-3\\1&1\end{pmatrix},\quad
 -1<u<1,\quad v>0,\quad u^2+v^2>1.
\end{equation}
The change is a conjugacy, so group closure, normalized area and the
fixed-endpoint homotopy class are unchanged.  Put
\[
 P(s)=\frac{1+s^2}{1-s^2},\qquad
 Q(s)=\frac{1+s^2}{2s},\qquad
 L(s,t)=\frac{4st}{(1-s^2)(1-t^2)}.
\]
All arc exponents below are the positive real exponents in the exact
vertex-control flow.  They exceed one on positive arcs.

\paragraph{Transfers at a switch.}
On a mode-$j_i$ arc of duration $\theta_i$, $d_{j_i}$ is constant by
\eqref{fw:reciprocals}. The frame equation and the commutator identity
give $(GX)'=-2GX\mathsf B_{j_i}/d_{j_i}$. Lift the labels to integers
with outgoing increments $\epsilon_i=j_{i+1}-j_i\in\{1,-1\}$, and
set $q_i=e^{\theta_i/d_{j_i}}$ and
$\mathsf P=\left(\begin{smallmatrix}1&-2\\1&0\end{smallmatrix}\right)$.
Since $\mathsf B_j=R_*^{-j}\mathsf B_0R_*^j$, the frames
$G_iX_iR_*^{-j_i}\mathsf P$ at successive switches differ by the right
multipliers
\begin{equation}\label{sh:group-transfers}
 \mathsf P^{-1}e^{-2\theta_i\mathsf B_0/d_{j_i}}
 R_*^{-\epsilon_i}\mathsf P
 =
 \begin{cases}
 M(q_i)=\begin{pmatrix}0&q_i\\-q_i^{-1}&q_i\end{pmatrix},
       &\epsilon_i=1,\\
 N(q_i)=\begin{pmatrix}q_i&-q_i\\q_i&q_i^{-1}-q_i\end{pmatrix},
       &\epsilon_i=-1.
 \end{cases}
\end{equation}
Their ordered product equals
$\mathsf P^{-1}R_*^{1-\sum_i\epsilon_i}\mathsf P$ by the full endpoint
conditions. In the state chart at the start of mode $0$ the flow is
$W\mapsto q(W+1)-1$. Rephasing by $R_*$ or $R_*^{-1}$ at the next
switch gives the respective state transfer matrices
\begin{equation}\label{sh:short-state-transfers}
 T(q)=\frac12\begin{pmatrix}q&q-4\\q&q\end{pmatrix},\qquad
 T_-(q)=\frac12\begin{pmatrix}q&q+2\\-q&2-q\end{pmatrix}.
\end{equation}

For the cost let $a=(1+u)/v$ on the phased mode-$0$ arc and let
$e_-,e_+$ be its endpoint values of $1+u$. Integrating
$\ell\,dt=\tfrac12(1+(3+u^2)/v^2)\,du$ gives
\begin{equation}\label{sh:short-arc-cost}
 c=\frac{1+a^2}{2}(e_+-e_-)
   +2a^2(e_-^{-1}-e_+^{-1})-a^2\log(e_+/e_-).
\end{equation}
Equivalently, with the entrance reciprocals $d_0,d_1,d_2$ in the
unphased chart, the rational term for mode $j$ is
\[
 (1-q^{-1})d_j(d_j+d_{j-1})+(q-1)d_j(d_j+d_{j+1}).
\]
The reciprocals of the two other modes change to
$qd_{j+1}+(q-1)d_j$ and $(d_{j-1}-(q-1)d_j)/q$, respectively.
These formulas follow by substituting the mode-$0$ state map and
using $\sum_{i<j}d_id_j=1$, then applying a cyclic phase.
In particular the rational term is
$d_j[(d'_{j+1}-d_{j+1})+(d_{j-1}-d'_{j-1})]$, where primes denote
exit values. When these terms are summed over the twisted cycle,
combine the two contributions at each switch. If $D_i$ denotes the
active reciprocal on arc $i$, the contribution at its outgoing switch
is $\epsilon_i(3D_iD_{i+1}-1)$ by the reciprocity identity. The
twisted seam has the same formula because the phase permutes the
three reciprocals. Thus the complete rational part is
\begin{equation}\label{sh:switch-rational-cost}
 \sum_i\epsilon_i(3D_iD_{i+1}-1).
\end{equation}

\begin{lemma}\label{sh:mixed-five-chart}
Every actual strict five-arc chain of reduced word $01202$ is represented
by $0<a,r<1$, $(1+a)(1+r)>2$, with exponents
\[
 (P(a),Q(a),Q(r),P(r),L(a,r)).
\]
Putting $x=2a/(1-a)$ and $y=2r/(1-r)$, its cost is
\begin{equation}\label{sh:mixed-five-cost}
 \mathscr J=2-\frac{G(x)+G(y)}{xy+2x+2y+3},\quad xy>2,
\end{equation}
where, for $t=x+1$,
\[
 G(x)=\frac{2(t^2+1)}{(t-1)^2}\log\frac{t^2+1}{2t}
       +(t^2-1)\log\frac{t^2+1}{t^2-1}-3.
\]
\end{lemma}
\begin{proof}
Let $Q_u=\left(\begin{smallmatrix}1&u\\0&1\end{smallmatrix}\right)$
and $N_b=\left(\begin{smallmatrix}1&0\\b&1\end{smallmatrix}\right)$.
On a $0$ or $2$ arc,
\begin{equation}\label{sh:shear}
 G_{\rm out}Q_{u_{\rm out}}=G_{\rm in}Q_{u_{\rm in}}N_b,
 \qquad b=\int\frac{du}{v^2}.
\end{equation}
In particular a suffix using only these two modes fixes the second
column of $GQ_u$. For a mode-$1$ arc from $(x,V)$ to $(z,W)$ with
exponent $q$, its frame transfer $C_1$ in the $Q_u$ coordinates has
second column
\begin{equation}\label{sh:mode-one-column}
\begin{aligned}
 (C_1)_{12}&=-VW\sinh(\log q),\\
 (C_1)_{22}&=\frac WV\{\cosh(\log q)-x\sinh(\log q)\}.
\end{aligned}
\end{equation}
Indeed $X(z,W)Q_z$ has second column $(-W,0)^t$. The relative frame is
\[
 -X(x,V)e^{-2\log(q)\mathsf B_1}X(z,W).
\]
With $\lambda=(q-1)/(q+1)$, the state map is
$w\mapsto(w-\lambda)/(1-\lambda w)$ on the complex upper half-plane.
Writing $D_m=(1-\lambda x)^2+\lambda^2V^2$, its exit height is
$W=V(1-\lambda^2)/D_m$.

Let $p,q>1$ be the first two exponents. The increment string is
$+++--$, so the product in \eqref{sh:group-transfers} is $I$.
The inverse of $N(f)N(g)$ has lower-right entry zero, since the
upper-left entry of $N(f)N(g)$ is zero. Hence
$(M(p)M(q)M(e))_{22}=0$ for the third exponent $e>1$,
which gives $p^2q^2=p^2+q^2$. Setting $a=q(1-1/p)$ gives
$0<a<1$, $p=P(a)$ and $q=Q(a)$.

For the entrance height write $u=u_0$, $w=v_0^2$ and
$\Delta_0=(1-u)^2+w$. The required transfer
$Q_{-u}R_*Q_{u_T}$ has second column
$(-2w/\Delta_0,2(1-u)/\Delta_0)^t$. The mode-$1$ entrance is
$(x,V)=(p(1+u)-1,pv_0)$. The upper entries in
\eqref{sh:mode-one-column} and in the required column coincide
exactly when $D_m=\lambda p^2\Delta_0$. Since
$(1+\lambda)^2=4\lambda p^2$, this gives
\[
 w=\frac{((1-\lambda x)/(p\sqrt\lambda))^2-(1-u)^2}{1-\lambda}
   =(1+u)(1+2a-u).
\]
The lower entries then coincide as well: the initial shear is
$(1+u)(1-1/p)/w$, and the lower numerator after multiplication by
that shear matrix is $(1+\lambda)^2-2\lambda p^2(1+u)$ over $D_m$.
Substitution into the state map and into $W_T=\Phi W_0$ gives
$v_2/(1+u_2)=v_T/(1-u_T)=v_0/[a(1+u_0)]=:B_0>0$.

Write the three suffix lines as
$v=B_1(1-u)$, $v=A_0(1+u)$, $v=B_0(1-u)$.
Their shear residual is exactly
\[
 b_{\rm actual}-b_{\rm required}
       =-\frac2{v_2}\left(\frac{u_2}{A_0}+\frac1{B_0}\right).
\]
It vanishes precisely when $A_0=-u_2B_0$, and forces $u_2<0$.
Set $r=-u_2\in(0,1)$.  Solving the linear equation for the entrance
state gives, with $D=(r-1)(a^3-a+1)-a^2(r+3)<0$,
\begin{align*}
 u_0&=-\frac{a^3r-a^3+a^2r-5a^2-ar+a-r+1}{D},\\
 v_0^2&=-\frac{4a^2(a-1)(a+1)^2(r-1)(ar-a-r-3)}{D^2}.
\end{align*}
The suffix intersections now give its three displayed exponents.
Their positivity reduces to $(1+a)(1+r)>2$.
Integrating the original arc cost, and putting $z=(1-a)(1-r)$, gives
\[
 \mathscr J=\frac{4(z+2)}{4-z}
 -\frac{z}{4-z}\left[f_0(a)+f_0(r)+\log L(a,r)\right],\quad
 f_0(s)=\frac{\log P(s)}{s^2}
       +\left(\frac{1+s}{1-s}\right)^2\log Q(s).
\]
Substitution of $x,y$ gives \eqref{sh:mixed-five-cost}.
The argument uses necessary identities of the full group closure; it
never eliminates an exceptional return by a nonsingularity assumption.
\end{proof}

\begin{lemma}\label{sh:h-crossing}
The function $h(x)=(x+2)G'(x)$ satisfies
\[
 h(x)>2\quad(0<x<\sqrt2),\qquad h(\sqrt2)=2,
 \qquad 0<h(x)<2\quad(x>\sqrt2).
\]
\end{lemma}
\begin{proof}
Put $t=x+1$, $L_0=\log((t^2+1)/(2t))$,
$M_0=\log((t^2+1)/(t^2-1))$ and $K_0=t^2-2t-1$.
Direct differentiation yields
\[
 h-2=A L_0+B M_0+C,
\quad A=-\frac{4(t+1)^2}{(t-1)^3},\quad B=2t(t+1),\quad
 C=-\frac{2K_0(2t^2+t+1)}{t(t-1)(t^2+1)}.
\]
The bounds $L_0<(t-1)^2/(2t)$ and $M_0>2/t^2$ give
$h>4(t+1)/(t(t^2+1))>0$.
For $1<t\leq2$, set $b=(t-1)/(t+1)$ and use
$L_0<2b^2+2b^6/(3(1-b^4))$.  Then
\[
 h-2>-\frac{t^3+3t^2-9t-7}{3t(t^2+1)}>0,
\]
since the numerator is increasing on $[1,2]$ and equals $-5$ at $2$.
For $t\geq2$ use $R=(t^2-1)/(2t)$ and
\[
 \frac{\log R}{K_0}>\frac2{t^2+2t-1},\qquad
 M_0<\frac2{t^2}+\frac2{3t^2(t^4-1)}.
\]
The first inequality is the integral inequality
$\log R/(R-1)>2/(R+1)$.  Substitution into
\[
 \frac{h-2}{K_0}=A\frac{\log R}{K_0}
 +\frac{2(t+1)(t^2-t+2)}{(t-1)^3}M_0
 -\frac{2(2t^2+t+1)}{t(t-1)(t^2+1)}
\]
gives a strictly negative upper bound whose numerator is $-2p(t)$ and
whose denominator is positive.  The positivity of $p$ follows from
\[
 p(v+2)=3v^7+36v^6+183v^5+508v^4+831v^3+782v^2+363v+46.
\]
Thus $(h-2)/K_0<0$.  At the sole zero $t=1+\sqrt2$ of $K_0$,
$L_0=M_0=\log\sqrt2$ and $A+B=C=0$.  This proves the assertion.
\end{proof}

\begin{proposition}\label{sh:mixed-five}
Every strict five-arc chain of word $01202$ satisfies
$\mathscr J>\mathscr J_8$.
\end{proposition}
\begin{proof}
An interior stationary point of \eqref{sh:mixed-five-cost} must have
$h(x)=h(y)$.  Since $xy>2$, Lemma~\ref{sh:h-crossing} then forces
$x,y>\sqrt2$.  But
\[
 G(\sqrt2)=(3+2\sqrt2)\log2-3>1,
\]
using $\log2>56/81$ and $\sqrt2>39/28$.
Since $G'>0$, the numerator of $\mathscr J_x$ is
\[
 (y+2)(G(x)+G(y)-h(x))+G'(x)>0,
\]
a contradiction.  On the boundary $xy=2$, put
$a=x/(x+2)$, $r=(1-a)/(1+a)$ and
\[
 c=\frac{1-a^2}{1+a^2},\qquad d=\frac{2a}{1+a^2},\qquad p=c-d.
\]
Substitution gives $\mathscr J(x,2/x)=2\delta(p)$, with the exact
function in \eqref{sh:octagon-cost}.  Hence the boundary value is at
least $\mathscr J_8$, with equality only at $x=y=\sqrt2$.
Finally $G(x)\to-1$ as $x\downarrow0$, and
$-1<G(x)<10\log(x+1)-1$ for $x\geq1$.  For an escaping sequence with
$xy\geq2$, $Q=\max(x,y)\to\infty$ and
\[
 \frac{|G(x)+G(y)|}{xy+2x+2y+3}
 \leq\frac{2+20\log(Q+1)}{2Q+3}\longrightarrow0.
\]
Thus the cost tends to $2$ at infinity.  Any interior point with cost
at most $\mathscr J_8$ would give an interior attained minimum, already
excluded by the stationary equations.
\end{proof}

\begin{proposition}\label{sh:reverse-five}
Every strict five-arc chain of word $02102$ satisfies
$\mathscr J>29/16>\mathscr J_8$.
\end{proposition}
\begin{proof}
All five increments are negative and the group product is $I$.
For three consecutive exponents $p,d,q$, the inverse of the remaining
two $N$ factors has lower-right entry zero. Direct multiplication gives
\[
 pdq\,(N(p)N(d)N(q))_{22}=(p^2-1)(q^2-1)-d^2.
\]
Thus $q_i^2=(q_{i-1}^2-1)(q_{i+1}^2-1)$ at every cyclic cut.
With $p_i=q_i^2-1$ we obtain
\[
 p_{i-1}p_{i+1}=1+p_i,\quad
 (p_0,\ldots,p_4)=
 \left(A,B,\frac{1+B}{A},\frac{1+A+B}{AB},\frac{1+A}{B}\right).
\]
Let $S=\sum q_i$, $T=\sum q_iq_{i+2}$ and $K_0=2+S+T$.
Use the state transfer $T_-$ of \eqref{sh:short-state-transfers}.
For its fivefold product, direct multiplication using the displayed
recurrence gives
\[
 \Delta=4\prod q_i-(\tr\prod T_-(q_i))^2
       =4K_0\prod(q_i-1)>0,\qquad
 d_i^2=\frac{q_i+1}{K_0(q_i-1)}.
\]
Thus the scalar-return case is excluded by an identity, rather than by
an assumed regularity of the endpoint map. In mode-$0$ coordinates
at its entrance, the two active reciprocals are $d_i,d_{i-1}$ for
modes $0,1$, while at its exit they are $d_i,d_{i+1}$ for modes $0,2$.
Using $d_0d_1+d_1d_2+d_2d_0=1$ gives
$e_-=2d_i(d_i+d_{i-1})/(1+d_i^2)$,
$e_+=2d_i/(d_i+d_{i+1})$. Substitution into
\eqref{sh:short-arc-cost} gives the rational term
$(q_i-1)d_i(2d_i+d_{i-1}+d_{i+1})$.
Its sum is $2+6/K_0$.  One verifies the cross terms by writing
$q=(p,r,c/a,pr/(ab),c/b)$, where
$a^2=p^2-1$, $b^2=r^2-1$, $c^2=p^2+r^2-1$, and taking all roots
positive.  The logarithmic coefficient on arc $i$ is $-d_i^2$, so
\[
 \mathscr J=2-\frac{\sum_i f(q_i)-6}{2+S+T},\qquad
 f(q)=\frac{q+1}{q-1}\log q.
\]
The positive logarithm series gives $f(q)<(q+q^{-1}+10)/6$.
Multiplying the recurrence gives $\prod(q_i-q_i^{-1})=1$, hence
$S-\sum q_i^{-1}\geq5$.  If $Q=(\prod q_i)^{1/5}$, convexity of
$z\mapsto\log(1+e^z)$ gives $Q^2\geq Q+1$, and therefore
$T\geq5Q^2>13$.  Also $T>2S-5$ by expansion of
$\sum(q_i-1)(q_{i+2}-1)>0$.  For $S\leq9$ use the first bound,
and for $S\geq9$ the second, to get $9T>7S+54$.  Consequently
\[
 \frac{S+\sum q_i^{-1}+14}{6(2+S+T)}<\frac3{16},
\]
which proves the cost bound.  Finally
$186\sqrt2+32\log2-285>209/2349>0$, using
$\sqrt2>41/29$ and $\log2>56/81$, is exactly $29/16>\mathscr J_8$.
\end{proof}

\subsection{The two six-arc classes}

The following use of a parameter family is local.  At an actual strict
chain, positive durations and strict star are open conditions.  An exact
nearby family with the same full group and twisted state endpoints
preserves the fixed-endpoint homotopy class of the group path: for close
paths their pointwise quotient has a real matrix logarithm, which gives
an endpoint-fixed homotopy.  The curvature identities in the original
model then give a regular convex boundary with the same winding.
Consequently, a nonzero derivative of the original cost in such a family
excludes local minimality.  This argument requires neither a globally
connected parameter domain nor a generic rank assumption.

\begin{proposition}\label{sh:six-010202}
An actual strict six-arc chain of word $010202$ is not a local area
minimum, even among chains with the same six controls.
\end{proposition}
\begin{proof}
The increment string is $+--+--$ and the group product is $-I$.
The first-column entries of $M(e)N(f)N(g)$ coincide. Thus for the
prefix $A=M(p)N(q)N(s_0)$, $A_{21}+A_{22}=0$.
Multiplication gives $pq s_0(A_{21}+A_{22})=p^2+q^2-p^2q^2$.
The first two exponents are therefore $P(a),Q(a)$, $0<a<1$.
The second-column calculation in Lemma~\ref{sh:mixed-five-chart}
applies because the suffix is again all in modes $0,2$.
Its lines can be written
\[
 v=B_0(1+u),\quad v=sB_0(1-u),\quad
 v=rB_0(1+u),\quad v=B_0(1-u),
\]
where $r,s>0$.  Its shear, multiplied by $B_0^2$, is
\[
 \frac1{1+u_2}+\frac1{1-u_T}
       +\frac1{rs}-\frac1r-\frac1s-1.
\]
The required shear consists of the first two terms.  Thus
$1-r-s-rs=0$, so $s=(1-r)/(1+r)$ and $0<r<1$.
The first suffix switch is $u_3=-r$.  Positive duration gives
$u_2=-z$ with $r<z<1$.  Set
\[
 k=\frac{(1-a^2)(1-z)}{2a},\qquad
 B_0^2=\frac{3+a+(1-a)z}{(1-a)(1-z)}.
\]
The six exponents are
\[
 \left(P(a),Q(a),\frac{1-r}{1-z},Q(r),P(r),
              \frac{2r}{(1+r)k}\right).
\]
Their strict domain is $0<a,z,r<1$, $r<z$, $k(1+r)<2r$,
with the additional suffix star condition $rsB_0^2>1$.
The entrance formulas are those of Lemma~\ref{sh:mixed-five-chart}
with $r$ there replaced by $z$.  Thus these expressions are an exact
three-parameter full-closure family around every actual point.
The mode-$1$ star check follows from
$u_{\theta}'=(u^2-v^2-1)/2$: at $u=\pm1$ this is $-v^2/2<0$,
so a trajectory with both endpoints in $(-1,1)$ cannot leave and return.
The $0$ and $2$ checks follow from their monotone horizontal coordinates;
the displayed suffix condition checks their smallest star value.

Define
\begin{align*}
 F(s)&=\frac{1+s^2}{s^2}\log P(s)
       +\frac{4s}{(1-s)^2}\log Q(s),\\
 T(s)&=F(s)+2\log(1-s),\qquad
 \omega=\frac{(1-a)(1-z)}{4-(1-a)(1-z)}.
\end{align*}
The full original arc integral simplifies to
\begin{equation}\label{sh:six01-cost}
 \mathscr J=2+\omega[6-F(a)-T(r)+2\log(1-z)].
\end{equation}
For a direct check, its rational part is $2+6\omega$ and the six
logarithmic weights are
$-\omega(a^{-2},((1+a)/(1-a))^2,1,((1+r)/(1-r))^2,r^{-2},1)$.
Their combination uses only
\[
 \frac{Q(a)Q(r)E_2E_5}{P(a)P(r)}
       =\frac{(1-r)^2}{(1-z)^2}.
\]
Since $G(x)=F(x/(x+2))-3$, Lemma~\ref{sh:h-crossing} gives
\[
 F'(s)=\frac{h(2s/(1-s))}{1-s},\qquad
 T'(s)=\frac{h(2s/(1-s))-2}{1-s}.
\]
Thus $T'$ vanishes only at $\rho=\sqrt2-1$.
At a stationary point, \eqref{sh:six01-cost} first forces $r=\rho$.
Put $t=1-a$, $q=1-z$, $H=4-tq$ and
$W=6-F(a)-T(r)+2\log q$.  Then
\[
 \mathscr J_z=-\frac{2t(2W+H)}{H^2},\qquad
 \mathscr J_a=-\frac{4qW}{H^2}-\frac{tq}{H}F'(a).
\]
The first equation forces $W=-H/2$, and the second then forces
$a=\rho$.  Since $z>r=\rho$, $0<q<1-\rho$.  The remaining equation
would be
\[
 \Psi(q)=16-4F(\rho)+4\log\frac q{1-\rho}-(1-\rho)q=0.
\]
But $\Psi'(q)>0$ and
$\Psi(q)<\Psi(1-\rho)=16-4F(\rho)-(1-\rho)^2<0$,
because $F(\rho)=G(\sqrt2)+3>4$.
There is no stationary point in the open actual domain.
\end{proof}

The other six-arc word requires a finite inequality.  We state its
polynomial and interval witnesses explicitly, so that the computational
input is a finite assertion about specified functions and domains.

\begin{lemma}\label{sh:six02-chart}
Every actual strict six-arc chain of word $021020$ belongs locally to the
following exact full-closure family, with $a,b,t\in(0,1)$:
\begin{gather}
 E=(P(a),L(a,b),P(b),L(b,t),P(t),Q(a)Q(t)/P(b)),\\
 R=1-ab-at-bt,\\
 S_1=b+t-a+abt,\quad S_3=a+b-t+abt,\quad S_5=a+b+t-abt,\\
 Z=2+b^2-b(a+t)+b^2at,\qquad
 \kappa=\frac{(1-b^2)R^2}{Z^2-S_5^2},\\
 (g_0,\ldots,g_5)=(a^{-1},S_1/R,b^{-1},S_3/R,t^{-1},S_5/R),\\
 \mathscr J=2+\kappa\left(6-\sum_{i=0}^5g_i^2\log E_i\right).
 \label{sh:six02-cost}
\end{gather}
At an actual point $R,\kappa,g_i>0$, $E_i>1$ and
$\kappa g_i g_{i-1}<1$ for all cyclic $i$.
\end{lemma}
\begin{proof}
Write $p,d,q,e,f,g$ for the six exponents. The increment string is
$----++$, so $A N(e)M(f)M(g)=-I$ for $A=N(p)N(d)N(q)$.
The first column of $N(e)M(f)M(g)$ is $(0,-f/(eg))^t$.
Consequently $A_{22}=0$, and direct multiplication gives
$pdq A_{22}=(p^2-1)(q^2-1)-d^2$. Hence $p=P(a)$, $q=P(b)$
and $d=L(a,b)$ for $0<a,b<1$. Write
\[
 A=\begin{pmatrix}A_0&B_s\\-1/B_s&0\end{pmatrix},
 \qquad B_s=\frac{p(q^2-1)}{dq}>0,\qquad
 A_0+B_s=-\frac p{dq}.
\]
The upper-left entry of the closed product gives $g=fB_s/e$.
Its upper-right entry then gives $e^2=(q^2-1)(f^2-1)$.
Since $f>1$, write $f=P(t)$ for the unique $0<t<1$. Positivity
gives $e=L(b,t)$ and $g=Q(a)Q(t)/P(b)$. Conversely, these relations
make the lower-left entry of the product zero and both diagonal
entries $-1$, so the signed group equation holds.

Here is the second-column reduction for a $021$ prefix and a suffix
in modes $0,2$. Let $u=u_0$, $w=v_0^2$, $k=p/d$,
$x=1-2/d+k(1+u)$ and $V=kv_0$, the mode-$1$ entrance data.
In \eqref{sh:mode-one-column}, $\lambda=b^2$ and
$D_m=(1-b^2x)^2+b^4k^2w$. The upper-right equation for the full
transfer is $-2b^2k^2w/D_m=-2w/((1-u)^2+w)$. It fixes the height
through the identity
\begin{equation}\label{sh:prefix-021-height}
 w=\frac{[2(2a/(1-a^2)+b)/P(a)-b(1+u)]^2-(1-u)^2}{1-b^2}.
\end{equation}
Here $[2(2a/(1-a^2)+b)/P(a)-b(1+u)]=(1-b^2x)/(bk)$.
The lower-right equation is now automatic. The shear in the first
two arcs is
$(1+u)(1-1/p)/w+(2-p(1+u))(d-1)/(p^2w)$.
Using it in \eqref{sh:mode-one-column}, the remaining condition is
$4k^2-4/d^2-(1-b^2)^2/b^2=0$, precisely the relation
$d^2=(p^2-1)(q^2-1)$. All divisors are positive on the actual chain.

Write the prefix exit as $(z,v_z)$, and denote the prescribed twisted
endpoint by $(u_T,v_T)=\Phi(u,v_0)$. The three state maps of the $020$
suffix are affine, with combined scale $eg/f=B_s$. They therefore
reach the prescribed endpoint exactly when
\[
 v_T=B_s v_z,\qquad
 1+u_T-B_s(1+z)=2g(1-f^{-1})=\frac{2Q(a)t}{P(b)}.
\]
The vertical equation follows from $D_m=b^2k^2\Delta_0$,
where $\Delta_0=(1-u)^2+w$, since
$v_T/v_z=4b^2k/(1-b^4)=B_s$.
For the horizontal equation put $c=(1-u)/(1+u)$, $H=a+b$,
$M_0=1-ab$, and $c_0=(1-b^2x)/(bk(1+u))$. By
\eqref{sh:prefix-021-height},
\[
 c_0-bc=\frac{2aM_0(1+c)}{1+a^2},\qquad
 c_0+bc=\frac{2(aM_0+Hc)}{1+a^2},\qquad
 \Delta_0=\frac{(1+u)^2(c_0-bc)(c_0+bc)}{1-b^2}.
\]
The mode-$1$ state map gives
$B_s(1+z)=4(1+u)c_0/(b\Delta_0)-4k/(1+b^2)$, whereas
$1+u_T=4(1-u)/\Delta_0$. Substitution in the horizontal equation yields
\[
 t=\frac{M_0c-aH}{Hc+aM_0},\qquad
 c=\frac{a(H+M_0t)}{M_0-Ht}.
\]
Here $R=M_0-Ht>0$ because $c>0$.
The entrance state is
\[
 u_0=\frac{(1-a^2)(1-bt)-2a(b+t)}{(1+a^2)(1-bt)},\quad
 v_0^2=\frac{4a^2R^2}{\kappa(1+a^2)^2(1-bt)^2}.
\]
Indeed $c=aS_5/R$ and $c_0=aZ/R$, so
$a^2(1+u)^2/w=(1-b^2)R^2/(Z^2-S_5^2)=\kappa$.
The reciprocal update after \eqref{sh:short-arc-cost} gives
$D_i=\sqrt\kappa\,g_i$ using these entrance data and the six exponents.
At every switch the third reciprocal is positive, so the two active
reciprocals have product less than one. The signed sum of neighboring
products in \eqref{sh:switch-rational-cost}, divided by $\kappa$, is
\[
 \frac{(S_5-S_1)/a-(S_1+S_3)/b+(S_5-S_3)/t}{R}=2.
\]
Thus the rational part is $2+6\kappa$, and the logarithmic weights
are $-\kappa g_i^2$.

Conversely, the displayed state equations and the signed group
product give both original endpoint conditions for the constructed
family. All divisors are nonzero at the actual point. The positive
durations and strict star condition therefore persist for sufficiently
small changes of all three parameters, and the local homotopy
preserves the real winding.
\end{proof}

\begin{lemma}[Six-arc derivative inequality]\label{sh:six02-finite}
On the actual domain of Lemma~\ref{sh:six02-chart},
\[
 \mathscr J\leq29/16\quad\Longrightarrow\quad
 \frac{\partial\mathscr J}{\partial b}>0,
\]
where $a,t$ are held fixed.
\end{lemma}
\begin{proof}
The formula is symmetric in $a,t$, so take $a\geq t$ and put
\[
 v=\frac{1-t}{1+t},\qquad
 r=\frac{1-a}{(1+a)v},\qquad b=vy.
\]
Then $0<v<1$, $0<r\leq1$.  The positive exponents $E_3,E_5>1$
imply respectively
\begin{equation}\label{sh:six-enlargement}
 1<y<Y:=\sqrt{\frac{1+r^2}{1+r^2v^4}}.
\end{equation}
We estimate on the larger closed $y$ interval.

First, the region $r\leq3/5$ has $\mathscr J>29/16$ by a polynomial
witness.  For $q\geq1$,
\[
 \log q\leq\mathcal U(q):=
 \frac{(q-1)(q^2+10q+1)}{6q(q+1)},
\]
while $\log q\leq q-1$ for all $q>0$.
Replace the first five logarithms in \eqref{sh:six02-cost} by
$\mathcal U(E_i)$ and the last by $E_5-1$.  Parameterize a rational
outer cover of \eqref{sh:six-enlargement} by
\[
 y=1+\frac{r^2(1-v^4)}{2(1+r^2v^4)}\zeta,\qquad0\leq\zeta\leq1.
\]
For $0<v<1$ and $0<r\leq3/5$ this enlargement has $vy<1$,
$y<1+r$ and positive weights. Clear denominators on that region and
remove positive factors. After $r=3x/5$ the resulting polynomial of
coordinate degrees $(66,32,12)$ extends to the closed cube $[0,1]^3$.
Its $67\cdot33\cdot13=28,743$ scaled integer Bernstein
coefficients are all strictly positive.  The finite witness consists of
this array, the cleared-denominator identity, and the inverse Bernstein
transform back to the rational bound.  This defines a check using only
integer arithmetic; in particular, the list is checked against the cost
formula rather than accepted as an unexplained table of signs.

It remains to take $1/2\leq r\leq1$; this slightly larger interval
simplifies the remaining witnesses.  Define
\begin{align*}
 T_0&=1+r-y+ryv^2,\\
 U_1&=y+r-1+ryv^2,\qquad U_3=y+1-r+ryv^2,\\
 V_0&=1+v^2[y(1+r)-r],\\
 A_0&=T_0+v[(y-1)(y-r)+r]+ry^2v^3,\\
 B_0&=2+v(1+r-y)+v^2[y^2+y(1+r)]+ryv^3+ry^2v^4,\\
 C_0&=\frac{1-v^2y^2}{A_0B_0}.
\end{align*}
All their required denominators are positive.  Set
\[
 f(s)=\frac{\log P(s)}{s^2},\quad
 M=\frac{\log(P(rv)P(v)/P(vy))}{v^2},\quad
 K_0=6-\frac{\log E_0}{a^2}-\frac{\log E_4}{t^2}-f(vy),
\]
where $E_0=(1+r^2v^2)/(2rv)$ and $E_4=(1+v^2)/(2v)$.
Then the exact stable cost is
\begin{equation}\label{sh:six-stable}
 \mathscr J=2+vC_0H,\quad
 H=T_0^2K_0-U_1^2\log E_1-U_3^2\log E_3-V_0^2M,
\end{equation}
with $E_1=y(1-r^2v^2)/(r(1-v^2y^2))$ and
$E_3=y(1-v^2)/(1-v^2y^2)$.
Holding $a,t$ fixed is the same as holding $v,r$ fixed.  Thus
\begin{equation}\label{sh:six-derivative}
 \mathscr J_b=(C_0H)_y,\qquad
 M_y=-\frac{4y}{1-v^4y^4}.
\end{equation}

The two noncompact ends are analytic.  For $v\leq1/100$, the five
bounds on the logarithmic contributions, grouping terms $0,4$ together,
are respectively
\[
 \frac{11}{200},\quad\frac{93}{2500},\quad\frac{27}{1000},
 \quad\frac{103}{10000},\quad\frac{11}{200}.
\]
They follow from $T_0>3/8$, $\kappa<51v/100$,
$a,t>49/50$, $\log P(x)\leq2x^2/(1-x^2)$ and the monotonicity of
$v\log(c/v)$ in this interval. The five strict upper bounds add to
$369/2000$, so the sum of the five contributions is strictly smaller
than $369/2000$,
so $\mathscr J>3631/2000>K_{\rm sh}$.
For $v\geq299/300$ put $\varepsilon=1-v$ and $\delta=1-b$.
Then $0<\delta\leq\varepsilon\leq1/300$,
$a<17/50$, $t<1/500$, $\kappa<2\delta/5$ and $g_1,g_3,g_5<21/10$.
For the even terms $0,4$ use
$\log P(s)/s^2=2\sum_{k\geq0}s^{4k}/(2k+1)\leq2/(1-s^4)<21/10$
for $s=a,t$. Their contributions are less than $(42/25)\varepsilon$.
Also $E_2<1/\delta$, $E_1,E_3<2/\delta$ and
$E_5<Q(t)^2\leq1/t^2\leq4/\varepsilon^2$.
Thus the other four contributions are bounded by
$[1/2+441/125+441/125]\varepsilon\log(2/\varepsilon)$,
using the monotonicity of $s\log(2/s)$ in this range. This gives
\[
 2-\mathscr J<\varepsilon
 \left[\frac{42}{25}+\frac{1889}{250}\log\frac2\varepsilon\right]
 <\frac{13643}{75000},
\]
using $\log600<7$.  Hence $\mathscr J>136357/75000>K_{\rm sh}$.
These estimates are uniform in the other parameters.

The remaining compact domain is parameterized by
\begin{equation}\label{sh:six-root}
 (v,r,z)\in[1/100,299/300]\times[1/2,1]\times[0,1],\qquad
 y=1+z(Y-1).
\end{equation}
On this root \path{SHORT6-ALL} verifies complete closed coverage
and one of the following alternatives on every box:
\[
 \mathscr J>K_{\rm sh},\qquad \mathscr J_b>0,\qquad
 \mathscr J_b+8(\mathscr J-K_{\rm sh})>0,\qquad
 (a,b,t)\in[9/25,47/100]^3 .
\]
On the last cube the task uses the separate bound
\(\mathscr J_b>1/1000\). Thus every terminal proves the claimed
implication. The rational endpoint grid and the signs are checked
by the task.

The outward arithmetic rule, stable logarithm evaluation and
whole-box derivative enclosures are specified in
Appendix~\ref{sh:arithmetic-six}.  Every recorded lower bound
is rational and strictly positive.
Together with the analytic tails and the polynomial region these
witnesses exhaust the necessary domain, proving the implication.
\end{proof}

At an actual point with cost at most $\mathscr J_8<K_{\rm sh}$,
Lemma~\ref{sh:six02-finite} gives a strict decrease by decreasing $b$
slightly in the exact family of Lemma~\ref{sh:six02-chart}.
Combined with Proposition~\ref{sh:six-010202} and the short-word table,
this already proves Proposition~\ref{short:rigidity} for $B\leq6$.
No compactness assertion for a class of arbitrary fixed-word competitors
is needed for this conclusion.

\subsection{A complete cost comparison for the reverse seven-arc word}

\begin{proposition}\label{sh:reverse-seven}
Every actual strict seven-arc chain of word $0210202$ satisfies
$\mathscr J>\mathscr J_8$.
\end{proposition}

We supply the reduction, the noncompact estimates and the finite
verification separately.  The proof concerns every feasible chain of
this word; complete PMP is not used.

The prefix $021$ has exponents $P(a),d,P(b)$, where
\begin{gather*}
 0<a,b<1,\qquad B_0=(1-a^2)(1-b^2),\qquad d=4ab/B_0>1,\\
 D=4ab-B_0>0,\qquad U_0=4ab+B_0.
\end{gather*}
The increment string is $----+--$, with sum $-5$, so
\eqref{sh:group-transfers} gives the complete group product
\[
 N(p)N(d)N(q)N(E_3)M(E_4)N(E_5)N(E_6)=I.
\]
The first-column entries of $M(e)N(f)N(g)$ agree.  The corresponding
condition on the inverse prefix is
\[
 \frac{(p^2-1)(q^2-1)-d^2}{pdqE_3}=0.
\]
Every denominator is positive, so it gives the asserted prefix relation.
The upper-right group equation gives \eqref{sh:prefix-021-height},
which also applies at a scalar state return. In terms of
$c=(1-u_0)/(1+u_0)$ set
\[
 \tau=\frac{(1-a^2)c-2a^2}{a(1+a^2)}.
\]
After the first arc, $(1-u_1)/(1+u_1)=a\tau$, so $\tau>0$.
Write the prefix exit and twisted endpoint as $(z,v_z)$ and
$(t,v_T)$. The four suffix rays have slopes $AB_1,sB_1,rB_1,B_1$,
where $AB_1=v_z/(1+z)$ and $B_1=v_T/(1-t)$. Let $s_{\rm req}$
be the required suffix shear. In the same $Q_u$ coordinates as
\eqref{sh:mode-one-column}, with $x,V$ the mode-$1$ entrance,
the upper-left entry gives
\[
 s_{\rm req}=\frac{(1-u_0)^2+v_0^2}{2v_0^2}
       \left(\frac{1-b^4(x^2+V^2)}{1-b^4}-\frac{1-u_0}{2}\right).
\]
Integrating $du/v^2$ on the four rays and subtracting the endpoint
terms gives
\[
 c_{\rm req}:=B_1^2s_{\rm req}-\frac1{A^2(1+z)}-\frac1{1-t}
 =\frac{A-r}{Ars}-\frac1r-\frac{1+A^{-2}}2.
\]
Put $\beta=A(2c_{\rm req}+1+A^{-2})/2$. Then
$s=(A-r)/(A+\beta r)$. Positivity of the interior arcs gives
$r<A$, $s<1$ and hence $\beta>-1$. With $Z=(1+z)/(1-z)$ and
$T_0=(1-t)/(1+t)$, positivity of the two exterior arcs gives
$s>AZ$, $r>T_0$. Thus $AZ<1$, $1+A\beta Z>0$, and
$T_0<r<R_{\max}=A(1-AZ)/(1+A\beta Z)$.

Put $X_s=2b(1+a^2)+D\tau$, $X_e=2a(1+b^2)+D\tau$ and
$V_0=1+ab+(a+b)\tau$. In the stated state maps, inserting
$c=[2a^2+a(1+a^2)\tau]/(1-a^2)$ and the height
\eqref{sh:prefix-021-height} gives
\begin{align*}
 \frac{R_{\max}}{T_0}-1
  &=\frac{(\tau^2-1)(1+a^2)(1-ab)}{(1-a^2)V_0},\\
 \beta+1&=\frac{2(1-ab)^2(U_0-D\tau^2)}{X_sX_e},&
 1-\beta&=\frac{2DV_0^2}{X_sX_e}.
\end{align*}
All denominators are positive from $0<a,b<1$, $D>0$ and $\tau>0$.
In particular $-1<\beta<1$ and
\begin{equation}\label{sh:reverse-domain}
 1<\tau<\sqrt{U_0/D}.
\end{equation}
Only necessity of this enlarged domain is used below.

Put
\begin{gather*}
 \kappa=\frac{B_0}{4+4(a+b)\tau+D\tau^2},\quad
 \sigma=\frac{2b(1+a^2)+D\tau}{(1-a^2)(1+b^2)},\quad
 \varepsilon=\frac{2a(1+b^2)+D\tau}{(1+a^2)(1-b^2)},\\
 h=\frac{2(1-ab)^2(U_0-D\tau^2)}{(1-a^4)(1-b^4)},\quad
 C=(\sigma+\varepsilon)/2,\\
 f_\sigma=\frac{(1-b^2)[2b(1+a^2)+D\tau]}{2b(1-ab)V_0},\quad
 f_\varepsilon=\frac{(1-a^2)[2a(1+b^2)+D\tau]}{2a(1-ab)V_0},\\
 \Phi(s,e)=s^2\log\frac{2s}{s+e}+e^2\log\frac{2e}{s+e}.
\end{gather*}
The original cost has the lower bound
\begin{equation}\label{sh:reverse-L}
\begin{split}
 L(a,b,\tau)=2+\kappa\bigg[&6-\frac{\log P(a)}{a^2}
 -\frac{\log P(b)}{b^2}-\tau^2\log d
 +\sigma^2\log f_\sigma+\varepsilon^2\log f_\varepsilon\\
 &-\Phi(\sigma,\varepsilon)-2\Phi(C,C+\sqrt h)\bigg]
 \leq\mathscr J.
\end{split}
\end{equation}
For the cost calculation, the first three active reciprocals,
divided by $\sqrt\kappa$, are $a^{-1},\tau,b^{-1}$.
Put $s_0=(AB_1)^{-1}$, $e_0=B_1^{-1}$ and
$p=(rB_1)^{-1}-s_0$, $q=(sB_1)^{-1}-e_0$. Then $p,q>0$ and
\[
 H:=pq=\frac{\beta+1}{AB_1^2},\quad
 (s_0,e_0,H)= (\sqrt\kappa\,\sigma,\sqrt\kappa\,\varepsilon,\kappa h),
 \qquad (1+z,1-t)=(f_\sigma,f_\varepsilon).
\]
These identities follow from the prefix and endpoint state maps
and the displayed formulas for the rays. For an explicit rational
cost calculation, the prefix uses \eqref{sh:short-arc-cost} and its
reciprocal update. On the edge suffix the rational integrand is
$\frac12+3/v^2-\frac12(d/du)[u(1-u^2)/v^2]$ on either mode.
Thus the rational suffix integral is
\[
 \frac{t-z}{2}+3s_{\rm req}
 +\frac12\left(\frac{z(1-z^2)}{v_z^2}
                -\frac{t(1-t^2)}{v_T^2}\right).
\]
Substituting the prefix data and $s_{\rm req}$ makes the total
rational part $2+6\kappa$.

At the intersection of two suffix rays with reciprocals $d_0,d_2$,
the values of $1+u,1-u$ are $2d_0/(d_0+d_2),2d_2/(d_0+d_2)$.
Consequently the sum of the four suffix logarithmic penalties is
\[
 \Phi(s_0,e_0+q)-\Phi(s_0+p,e_0+q)+\Phi(s_0+p,e_0)
       -s_0^2\log(1+z)-e_0^2\log(1-t).
\]
We bound this expression by comparing rectangles with the same area.
For $s_0,e_0,p,q>0$, $H=pq$ and $s_0+e_0\geq\sqrt H$, the rectangle
integral
\[
 \int_{s_0}^{s_0+p}\int_{e_0}^{e_0+q}\frac{se}{(s+e)^2}\,de\,ds
\]
is at most its value on the square with outer endpoints
$(s_0+e_0)/2$ and side $\sqrt H$.
Indeed pair $(u,v)$ with $(v,u)$ on the unit square, put
$T=u+v$, $\Delta=u-v$, $S=p+q$, $\delta=p-q$, and
\[
 D_0=s_0-e_0+\delta T/2,\quad D_1=S\Delta/2,
 \quad L_0=s_0+e_0+ST/2,\quad L_1=\delta\Delta/2.
\]
Completing the square gives
\[
 \frac12\left(\frac{(D_0+D_1)^2}{(L_0+L_1)^2}
 +\frac{(D_0-D_1)^2}{(L_0-L_1)^2}\right)
 \geq\frac{D_1^2}{L_0^2+L_1^2}.
\]
With $\eta=S/(2\sqrt H)\geq1$ and $c=(s_0+e_0)/\sqrt H\geq1$,
the latter is at least $\Delta^2/(c+T)^2$, because the difference of
cross-multiplied denominators is
$(\eta^2-1)(c^2-\Delta^2)+2cT\eta(\eta-1)\geq0$.
Use $se/(s+e)^2=1/4-(s-e)^2/(4(s+e)^2)$ and integrate.
Since $\Phi_{se}=-1-2se/(s+e)^2$, the three $\Phi$ terms in the
suffix penalty equal $\Phi(s_0,e_0)+H$ plus twice that integral.
Here $H/(s_0+e_0)^2=A(\beta+1)/(A+1)^2<1/2$. For the comparison
square, the corresponding three terms reduce to
$\Phi(s_0,e_0)+2\Phi((s_0+e_0)/2,(s_0+e_0)/2+\sqrt H)$,
since $\Phi$ vanishes on the diagonal. Dividing by $\kappa$ uses the
degree-two homogeneity of $\Phi$ and gives precisely
\eqref{sh:reverse-L}.

\begin{lemma}\label{sh:reverse-inequality}
The function in \eqref{sh:reverse-L} satisfies
$L>\mathscr J_8$ throughout \eqref{sh:reverse-domain} and $D>0$.
\end{lemma}
\begin{proof}
The symmetry in $a,b$ will always be retained.  The entropy function has
the convergent expansion
\begin{equation}\label{sh:entropy-series}
 \Phi(s,e)=\frac{(s-e)^2}{4}
 \left(3-\sum_{k\geq1}\frac{w^k}{k(k+1)(2k+1)}\right),\qquad
 w=\frac{(s-e)^2}{(s+e)^2}<1.
\end{equation}
This follows by twice integrating
$[(1+u)^2\log(1+u)+(1-u)^2\log(1-u)]''=6+2\log(1-u^2)$.
In particular $0\leq\Phi\leq3(s-e)^2/4$.
The stable combined logarithms are
\begin{equation}\label{sh:combined-log}
 \sigma^2\log f_\sigma+\varepsilon^2\log f_\varepsilon
       -\Phi(\sigma,\varepsilon)
 =\sigma^2\log(1-u_s)+\varepsilon^2\log(1-u_e),
\end{equation}
where
\[
 u_s=\frac{(1-a^2)(1+b^2)[(1-ab)\tau-a+b]}{2b(1+a^2)V_0},\qquad
 u_e=\frac{(1+a^2)(1-b^2)[(1-ab)\tau+a-b]}{2a(1+b^2)V_0}.
\]
Both belong to $(0,1)$ on the necessary domain.  This exact combination
is essential in the noncompact estimates: estimating its original
summands separately would lose their cancellation.

We first list and prove the uniform tail estimates.  For $\tau\geq2$
use the coordinates
\begin{equation}\label{sh:tau-coordinates}
 \xi^2=\frac{U_0-D\tau^2}{2B_0},\quad
 d=1+\frac{2(1-\xi^2)}{\tau^2-1},\quad
 b=\frac{d(1-a^2)}{\sqrt{4a^2+d^2(1-a^2)^2}+2a},
 \qquad0\leq\xi\leq1.
\end{equation}
They include both finite boundary faces.  If $1/3\leq a\leq1/2$ and
$\tau\geq128$, then $1/3\leq b<3/5$, $\kappa\leq8/(27\tau)$,
$\sigma,\varepsilon<2$, $f_\sigma,f_\varepsilon>1/(3(1+\tau))$,
$f(a)+f(b)<9/2$, and $\tau^2\log d<3$, where $f(s)=\log P(s)/s^2$.
The two entropy bounds total at most $6$, cancelling the positive $6$.
Consequently
\[
 L\geq2-\frac8{27\tau}\left(\frac{15}2+8\log(3(1+\tau))\right)
 >\frac{539}{288}>K_{\rm sh}.
\]
The logarithmic error decreases with $\tau$; $\log387<6$ completes the
rational comparison.
For $0<a\leq1/3$ and $\tau\geq2$, the same coordinates give
$\kappa\leq a/(1+\tau)$ and $1\leq\sigma,\varepsilon\leq1+2/(\tau-1)$.
If $a\leq1/64$, the sharper $\sigma<17/16$ yields
\[
 L\geq2-\frac a{1+\tau}
 \left[5+\frac{8431}{3840}\log\frac1a
          +\frac{2593}{256}\log(2(1+\tau))\right]
 >\frac{374983}{204800}>K_{\rm sh}.
\]
The error is increasing in $a$ and decreasing in $\tau$, so use
$a=1/64$, $\tau=2$, $\log64<21/5$, $\log6<9/5$.
For all $a\leq1/3$ and $\tau\geq128$, instead use
\[
 L\geq2-\frac a{1+\tau}
 \left[\frac{81}{40}+3+4\log\frac5{3a}
       +\frac{81}{64}\log\frac{6(1+\tau)^2}{a}\right]
 >\frac{238247}{123840}>K_{\rm sh}.
\]
The same monotonicities reduce this to $a=1/3$, $\tau=128$,
$\log5<2$ and $\log299538<13$.

For $1/3\leq a\leq1/2$ and $b\geq99/100$, put $s=1-b^2<1/50$.
Equation~\eqref{sh:combined-log} gives
$\sigma<4$, $\varepsilon<4/s$, $u_s<1/2$, $u_e<s^2$ and
$\kappa\leq s/8$.  Hence
\[
 L>2-\frac s8\left[39+\frac{11}{10}\log\frac2s
                         +2\log\frac8{3s}\right]
 >\frac{1491}{800}>K_{\rm sh}.
\]
Here $s\log(c/s)$ increases on the stated range, and both required
logarithms at $s=1/50$ are less than $5$.
For $a\geq1/2$, $b\geq39/40$, let $A=1-a^2\leq3/4$ and
$s=1-b^2<1/20$.  Cancelling the factors in
$\sigma^2u_s$ and $\varepsilon^2u_e$ before estimating gives respectively
bounds $7$ and $12$, while $u_s<1/4$, $u_e<s^2$ and
$\kappa<As/11$.  Thus
\[
 L>2-\frac{As}{11}
 \left[26+4\log\frac2A+\frac{11}{10}\log\frac2s
                   +2\log\frac4{As}\right]
 >\frac{4067}{2200}>K_{\rm sh}.
\]
The error increases separately in $A,s$; use $(A,s)=(3/4,1/20)$ and
$\log(8/3)<1$, $\log40<4$, $\log(320/3)<5$.
All logarithmic constants just used follow, for example, by comparing
the stated argument with a finite positive Taylor sum for the exponential.

It remains to control the low-angle corner with $a\leq1/3$, $\tau\leq2$.
Put
\begin{equation}\label{sh:low-angle-coordinates}
 r=\frac{1-b}{a(1+b)},\quad b=\frac{1-ar}{1+ar},\quad
 \tau=1+rz.
\end{equation}
Then $0<r<1$, and $D\tau^2\leq U_0$ gives
\[
 \frac{\tau^2-1}{\tau^2+1}
 \leq\frac{r(1-a^2)}{1-a^2r^2}\leq r,
 \qquad 0\leq z\leq\frac{\tau^2+1}{\tau+1}\leq\frac53.
\]
This covers all coupled approaches to $a=0$, $b=1$ and $\tau=1$.
After cancellation in \eqref{sh:combined-log}, direct positive-factor
estimates give, uniformly in the other two variables,
\begin{align*}
 a\leq1/256:\quad
 L&>2-\frac{ar}{2}
 \left[60+\frac{11}{10}\log\frac1{ar}
             +4\log\frac9{8r}+\frac{121}{100}\log\frac3a\right]
 >\frac{94493}{51200}>K_{\rm sh},\\
 r\leq1/100:\quad
 L&>2-\frac{ar}{2}
 \left[18+\frac{11}{10}\log\frac1{ar}
             +2\log\frac9{8r}+\frac{441}{100}\log\frac3a\right]
 >\frac{115217}{60000}>K_{\rm sh}.
\end{align*}
For the first inequality one uses
$\kappa\leq ar/2$, $\sigma<11/10$, $r\varepsilon\leq3$,
$u_e/r^2<17/20$, and $(1-u_s)/a>1/3$.
For the second, $\tau\leq61/60$, $\sigma,r\varepsilon<21/10$ and
$u_e/r^2<8/5$ suffice.  The bracket times $ar$ increases in both $a,r$.
Evaluation at $(1/256,1)$ and $(1/3,1/100)$, respectively, gives the
displayed rational bounds.

We now specify the finite part.  Put $b_0(a)=(1-a)/(1+a)$,
$\delta=b-b_0(a)$ and $\rho=\sqrt2-1$.
The three derivative witnesses on
\[
 1/3\leq a\leq1/2,\quad0\leq\delta\leq1/100,\quad1\leq\tau\leq2
\]
prove $L_\tau>1/1000$, $L_\delta|_{\tau=1}>1/1000$ and
$F''(a)>1/10$, where $F(a)=L(a,b_0(a),1)$.
These are respectively three-, two- and one-dimensional inequalities.
The identities $F(\rho)=\mathscr J_8$ and $F'(\rho)=0$ follow by
substitution and the symmetry $F(a)=F(b_0(a))$.  Integrating in the order
$a$, then $\delta$ at $\tau=1$, then $\tau$, gives
\begin{equation}\label{sh:reverse-collar}
 L\geq\mathscr J_8+\frac{(a-\rho)^2}{20}
             +\frac\delta{1000}+\frac{\tau-1}{1000}.
\end{equation}
This is a comparison of costs, not an assertion that $L$ has the same
feasible derivative as the original cost.

The other compact witnesses are listed below.  An entry ``cost'' means a
strict lower bound for $L$ (or a stated smaller expression); the three
collar rows retain their derivative meanings.
\begin{center}\small
\begin{tabular}{p{.58\textwidth}rr}
Domain and quantity&Initial leaves&Threshold\\\hline
$[1/2,3/4]^2\times[0,1]$ in $(a,b,\xi)$, cost&5,347&$K_{\rm sh}$\\
Collar $L_\tau$&388&$1/1000$\\
Collar $L_\delta$ at $\tau=1$&85&$1/1000$\\
$F''$ on $[1/3,1/2]$&256&$1/10$\\
$[1/3,1/2]\times[2,128]\times[0,1]$ in $(a,\tau,\xi)$&15,911&$K_{\rm sh}$\\
$[1/64,1/3]\times[2,128]\times[0,1]$ in $(a,\tau,\xi)$&7,206&$K_{\rm sh}$\\
$a\in[1/3,1/2]$, $b\in[b_0+1/100,99/100]$, $\xi\in[0,1]$&54,655&$\mathscr J_8$\\
$[1/2,39/40]\times[3/4,39/40]\times[0,1]$ in $(a,b,\xi)$&39&$K_{\rm sh}$\\
$[1/256,1/3]\times[1/100,1]\times[0,5/3]$ in $(a,r,z)$&4,011&$\mathscr J_8$.
\end{tabular}
\end{center}
In rows with $(a,b,\xi)$ use
$\tau=\sqrt{1+2B_0(1-\xi^2)/D}$; the variable $\sqrt h$ is then
$2(1-ab)\xi/\sqrt{(1+a^2)(1+b^2)}$, with no endpoint singularity.
The variable interval for $b$ in the seventh row is affinely
parameterized by $[0,1]$.
For the eighth row replace the final entropy term by $3h/2$.
For the last row use its sharper upper bound
\[
 2\Phi(C,C+\sqrt h)\leq\frac{3h}{2}
                           -\frac{h^2}{48(C+1)^2},
\]
which follows from the first subtracted term of
\eqref{sh:entropy-series} and $\sqrt h<2$.
Boxes with a strictly negative upper bound for the numerator of $h$
have no necessary-domain points.  Boxes wholly in $\tau\geq2$ may
invoke the already verified fifth and sixth rows and their tails.
No such reference is used in the fifth or sixth row, so the checking
order is acyclic.

The finite checking rule and the entropy remainder bounds are
specified in Appendix~\ref{sh:arithmetic-reverse}.  The witnesses
are exact binary covers and every leaf proves its stated strict
inequality; the completed cover has no unresolved leaf.

Finally the domains exhaust the necessary set.  If $\tau\geq2$ and
neither variable lies in $[1/3,1/2]$, they cannot both be at least $1/2$,
since then $U_0/D\leq25/7<4$.  The remaining asymmetric case is covered
by the sixth row and its tails.  For $\tau\leq2$, a variable in
$[1/3,1/2]$ is handled by the collar, the seventh row or its $b$ tail.
If both variables exceed $1/2$, use the first or eighth row and the high
tail.  Otherwise exchange them and apply the last row with its two
uniform tails.  All interfaces overlap as closed sets.  Strictness in
the collar follows from $\delta>0$ and $\tau>1$ for an actual chain.
This proves the lemma.
\end{proof}

\begin{proof}[Completion of the proof]
The full-closure reduction gives \eqref{sh:reverse-domain} and
$\mathscr J\geq L$.  Lemma~\ref{sh:reverse-inequality} supplies the
strict comparison.
\end{proof}

\subsection{Pure-forward seven-arc minimizers and the complete PMP filter}

The pure-forward case is different from Proposition~\ref{sh:reverse-seven}:
we use inequalities supplied by maximization over the full control
triangle.  Their role is to separate the candidate set from the
zero-duration octagonal boundary and to reduce a genuinely
four-dimensional domain to a finite verified inequality.

\begin{proposition}\label{sh:forward-seven}
A positive pure-forward seven-arc trajectory satisfying the original
complete free-time maximizing PMP conditions cannot be an original
global minimizer.
\end{proposition}

We first record the exact chart and state formula used in the finite
proof. With seven positive increments, \eqref{sh:group-transfers}
gives the full matrix identity
\[
 M(q_0)\cdots M(q_6)=I,\qquad q_i>1.
\]
The sign on the right is fixed by the original twist, rather than by a
projective trace equation.  Put $t_i=q_i^2-1$.  Multiplication and
elimination of the last three positive exponents give
\begin{equation}\label{sh:pure-chart}
\begin{gathered}
 F=1+t_0+t_2-t_0t_1t_2,\qquad
 E=t_0(1+t_1)+t_3F,\\
 t_4=\frac{(1+t_2)(1-t_0t_1)}E,\\
 t_5=\frac{(1-t_1t_2)E}{(1+t_1)(1+t_2)},\qquad
 t_6=\frac{(1+t_1)(1-t_2t_3)}E.
\end{gathered}
\end{equation}
The complete positive chart has $t_0,\ldots,t_3>0$ and
$t_0t_1,t_1t_2,t_2t_3<1$.  The denominators are positive, and the same
relations hold at every cyclic cut.  In particular,
\begin{equation}\label{sh:pure-cyclic}
 t_{i+3}t_{i+4}
 =\frac{(1-t_{i-1}t_i)(1-t_it_{i+1})}{1+t_i}.
\end{equation}
These identities are verified by multiplying the seven matrices;
conversely substitution of \eqref{sh:pure-chart} gives $I$ in all four
entries.  In logarithmic $t$ coordinates the three defining inequalities
are affine half-spaces, so the positive chart is connected.

For the state return use $T$ from \eqref{sh:short-state-transfers} and
its cyclic products
\begin{equation}\label{sh:state-return}
 P_i=T(q_{i-1})\cdots T(q_0)T(q_6)\cdots T(q_i),
\end{equation}
with cyclic interpretations at the endpoints.  Define
\[
 \tau_0=\tr P_0,\quad \Delta=4\prod q_i-\tau_0^2,\qquad
 n_i=2(P_i)_{21}+(P_i)_{11}-(P_i)_{22}.
\]
For a return matrix $P_i$ with lower-left entry $c_i$, the
upper-half-plane fixed point yields
$a_i=\operatorname{sgn}(c_i)n_i/\sqrt\Delta$.
We establish below that $\Delta>0$ and $c_i>0$.
On the original state, the two active reciprocals at the entrance
of the phased mode-$0$ arc are $a_i,a_{i-1}$ for modes $0,2$.
At its exit they are $a_i,a_{i+1}$ for modes $0,1$.
The reciprocity identity gives
$e_-=2a_i/(a_i+a_{i-1})$ and
$e_+=2a_i(a_i+a_{i+1})/(1+a_i^2)$.
Their ratio and the arc integral \eqref{sh:short-arc-cost} give
\begin{equation}\label{sh:state-formulas}
 q_i=\frac{(a_i+a_{i-1})(a_i+a_{i+1})}{1+a_i^2},\qquad
 \mathscr J=3\sum_i a_i a_{i+1}-7-\sum_i a_i^2\log q_i.
\end{equation}
In summing the rational terms, the neighboring contributions with
denominator $a_i+a_{i+1}$ combine by
$[a_i(1+a_i^2)+a_{i+1}(1+a_{i+1}^2)]/(a_i+a_{i+1})
=1+a_i^2-a_i a_{i+1}+a_{i+1}^2$.
Equivalently, $n_i=2q_i\partial_{q_i}\tau_0-\tau_0$ before
substitution of the closure chart.  This identity is obtained by
differentiating the product in \eqref{sh:state-return}; it verifies the
normalization of the cost numerator used below.

\begin{lemma}\label{sh:half-state}
On every finite closed chart obtained from \eqref{sh:pure-chart} with
$t_i\geq0$, adjacent products at most one, and $E>0$, the state return
has $\Delta>0$.
\end{lemma}
\begin{proof}
Put $r_i=1/q_i>0$. If
$(x_{i+1},y_{i+1})=(x_i,y_i)M(q_i)$, eliminating $x_i$
gives $r_{i-1}y_{i-1}+r_i y_{i+1}=y_i$. This is the eigenvalue-one
recurrence for the seven-cycle Jacobi matrix with edge weights $r_i$.
The identity monodromy gives a two-dimensional eigenspace.  Thus its
six-point principal and adjacent cofactors vanish.
The matrix $\left(\begin{smallmatrix}-1&2\\1&0\end{smallmatrix}\right)$
conjugates $T(q)/\sqrt q$ to $M(\sqrt q)$, which replaces
$r_i$ by $\sqrt{r_i}$. Reversing the order of the factors does not
change the trace, since $K M(q)=M(q)^t K$ for
$K=\left(\begin{smallmatrix}-1&1\\1&0\end{smallmatrix}\right)$.
The half-state monodromy therefore has the trace of
$P_0/\sqrt{\prod q_i}$. The adjacent constraints say
$r_i^2+r_{i+1}^2\geq1$.

For a six-point path, write its five squared full edges as
$a,b,c,d,e\in(0,1]$.  Its zero determinant is
\[
 D_6=1-a-b-c-d-e+ac+ad+ae+bd+be+ce-ace=0.
\]
If $X=a+b-1$, then
$Xd=(1-e)[X+c(1-a)]$.
Write $A=\sqrt a$, $B=\sqrt b$, $C=\sqrt c$, $D=\sqrt d$,
$E_0=\sqrt e$, and $X_h=A+B-1>0$.
If $X=0$ then $e=1$, and the half-path determinant is $X_hD>0$.
If $X>0$ and $a=1$, it is $X_h(D+E_0-1)>0$, since $d+e=1$.
Otherwise set $\alpha=(1-a)/X$ and
$\alpha_h=(1-A)/X_h$.  The identity
$X_h^2-X=2(1-A)(1-B)\geq0$ gives $\alpha_h<\sqrt\alpha$.
With $s=c\alpha$, $d=(1-e)(1+s)\leq1$, so
\[
 \frac D{1-E_0}\geq1+s+\sqrt{s(1+s)}>1+C\alpha_h.
\]
The half-path determinant factors as
$[X_h+C(1-A)](E_0-1+D/(1+C\alpha_h))$ and is positive.

For the remaining five-point path put $X=a+b-1$, $Y=c+d-1$ and
$R=bc-XY>0$, the last strict inequality following from the full adjacent
cofactor identity $R=\prod_{j\ne i}r_j/r_i$.
Its half determinant is
$D_{5,h}=X_hY_h-\sqrt{bc}<0$, and
$D_{5,h}+\sqrt R>0$.
For the first inequality, $X_h\leq\sqrt b$ and $Y_h\leq\sqrt c$.
Equality in their product requires $a=d=1$, which gives $XY=bc$
and contradicts $R>0$.
For the second, $X_h\geq\sqrt X$, $Y_h\geq\sqrt Y$, and, when $XY>0$,
$\sqrt{XY}+\sqrt{bc-XY}>\sqrt{bc}$; the case $XY=0$ is immediate.
Thus the periodic half adjacent cofactors are positive; the
antiperiodic ones have sign opposite the corresponding signed edge.

If the half return were scalar $\pm I$, all six-point cofactors would
vanish.  If its trace were $2$ and it were nonscalar, its adjugate would
be $\eta vv^T$, with $\eta>0$ and all adjacent coordinates of $v$
having the same sign.  Eigenvalue one would then be the Perron root,
contradicting the three-point principal-path root
$\sqrt{r_i+r_{i+1}}>1$.  At trace $-2$, the antiperiodic signs instead
make every signed edge term of $v^THv$ negative, contradicting $Hv=v$.
The positive chart is connected and its regular point
$q_i=2\cos(2\pi/7)$ has strictly elliptic half return.
All finite closed-chart points are limits of positive points (shrink the
$t_i$ and then add a smaller positive increment).  Continuity and the
excluded parabolic cases prove strict ellipticity throughout the stated
closed chart.
\end{proof}

On the connected positive exponent chart every $c_i=(P_i)_{21}$
is nonzero by ellipticity, so its sign is constant. At its regular
point $q_i=q_*=2\cos(2\pi/7)$,
$T(q_*)/\sqrt{q_*}$ has angle
$\theta=\arccos(\sqrt{q_*}/2)$ and positive lower-left entry.
Since $2\pi<7\theta<7\pi/3$, its seventh power has positive
lower-left entry by the factor $\sin(7\theta)/\sin\theta$.
Hence $c_i>0$ on the chart and $a_i=n_i/\sqrt\Delta$ on the
original state. The angle bound follows from
$1/2<\cos(2\pi/7)<1$, which gives $2\pi/7<\theta<\pi/3$.

\begin{lemma}\label{sh:state-filter}
For an actual pure-forward seven-arc global-minimizer candidate,
\begin{equation}\label{sh:PMP-pairs}
 a_i^2+2a_i a_{i-1}>1,\qquad
 a_{i-1}^2+2a_i a_{i-1}>1,\qquad 0<a_i a_{i-1}<1.
\end{equation}
Moreover $1/5<a_i<11/4$.
\end{lemma}
\begin{proof}
At each switch choose a genuine fundamental interval containing positive
pieces of the capped forward word $1201$.  Apply the complete
maximizing-PMP forward-force lemma, Lemma~\ref{fin:forward-force}, on
that interval.  In its normalized middle-switch coordinates,
\[
 d_0=a_i,\qquad d_2=a_{i-1},\qquad
 d_1=\frac{1-a_i a_{i-1}}{a_i+a_{i-1}}.
\]
The conclusion $d_1<\min(d_0,d_2)$ is precisely the first two
inequalities.  The last is strict star.  Different legitimate cuts may
have different lifts; only their scalar state conclusions are combined.

Put $\alpha_i=\arctan a_i$ and $x_i=2\alpha_i/\pi$.
The pair inequalities say $P_i:=x_i+x_{i+1}<1$ and
$F_{ij}:=2x_i+x_j>1$ for each oriented adjacent pair.
The state phase is
\begin{equation}\label{sh:phase}
 \phi=\frac{3\pi}{2}-\sum_i\alpha_i\in(0,\pi).
\end{equation}
Here the branch can be fixed on the entire connected exponent chart;
we do not assume that its geometrically feasible subset is connected.
Indeed, in the coordinates $b_j=\log t_j$, $0\leq j\leq3$, the full
positive group chart is exactly
\[
 b_0+b_1<0,\qquad b_1+b_2<0,\qquad b_2+b_3<0.
\]
It is convex and contains the regular point.  By
Lemma~\ref{sh:half-state}, the functions
\[
 \phi=\arccos\frac{\tau_0}{2\sqrt{\prod q_i}}\in(0,\pi),
 \qquad a_i=\frac{n_i}{\sqrt\Delta}
\]
are single-valued real analytic functions on this entire chart.
For this continuation the $a_i$ need not initially be declared positive.
The state identities extend throughout the chart in the denominator-free
form
\[
 q_i(\Delta+n_i^2)=(n_i+n_{i-1})(n_i+n_{i+1}).
\]
They follow directly from the cyclic state-return products; equivalently,
they hold near the regular point by the original state recursion and
then everywhere by the identity principle for real analytic functions
on a connected domain.  The regular trajectory is strictly in the star
domain, so that neighborhood is nonempty.  As $q_i>0$ and $\Delta>0$,
no adjacent sum $n_i+n_{i+1}$ vanishes.  Its sign is positive at the
regular point and hence is positive throughout the connected chart.
Consequently the logarithmic recursion in
\eqref{sh:state-formulas} has no hidden logarithm branch there.

Differentiating the trace product gives
$d\phi=-\tfrac12\sum_i a_i\,d\log q_i$.
If $K=\partial\log q/\partial a$, the recursion gives a symmetric
matrix with
\[
 K_{ii}=\frac1{a_i+a_{i-1}}+\frac1{a_i+a_{i+1}}
                   -\frac{2a_i}{1+a_i^2},\qquad
 K_{i,i+1}=\frac1{a_i+a_{i+1}},
\]
and $Ka=(2/(1+a_i^2))_i$.  Thus, with the ordinary real branch of
$\arctan$,
\[
 d\left(\phi+\sum_i\arctan a_i\right)=0
\]
on the whole connected exponent chart.  At the regular point
$q_*=2\cos(2\pi/7)$ one has
$a_i=\sqrt{q_*/(4-q_*)}$.  Write
$\theta=\arccos(\sqrt{q_*}/2)$, so that
$\arctan a_i=\pi/2-\theta$ and
$2\pi<7\theta<7\pi/3$.  Therefore its principal half-state phase is
$\phi=7\theta-2\pi$, and the constant is exactly $3\pi/2$.
This proves \eqref{sh:phase} for every actual candidate, regardless
of the connected components of the additional geometric constraints.
In particular $S:=\sum x_i<3$.
The two exact linear identities
\[
 -7x_0=P_1+P_5+S-2F_{01}-F_{23}-F_{54}-2F_{06}<-1,
\]
\[
 9x_0=3P_0+3P_6+3S-3F_{12}-F_{34}-F_{43}-3F_{65}<7
\]
give $\pi/14<\alpha_i<7\pi/18$ cyclically.
The lower bound implies $a_i>1/5$.  For the upper bound put
$y=\arctan(4/11)<\pi/6$; then
$\tan(3y)=1388/803<\sqrt3$, because
$3\cdot803^2-1388^2=7883>0$.
Hence $y<\pi/9$ and $a_i<\cot(\pi/9)<11/4$.
\end{proof}

\begin{lemma}\label{sh:state-monotone}
The extension of the cost in \eqref{sh:state-formulas} to all positive
state tuples is strictly increasing in each state variable.
\end{lemma}
\begin{proof}
Write $\mathscr J_{\rm ext}+7=\sum_i F(a_i,a_{i+1})$, where
\[
 F(a,b)=3ab-(a^2+b^2)\log(a+b)
       +\frac{a^2}{2}\log(1+a^2)+\frac{b^2}{2}\log(1+b^2).
\]
For $r=b/a>0$ and $v=1/a^2>0$,
\[
 \frac{F_a}a=3r-\frac{1+r^2}{1+r}-2\log(1+r)
                  +\log(1+v)+\frac1{1+v}
 >2\left(r+\frac r{1+r}-\log(1+r)\right)>0.
\]
The first inequality uses
$(\log(1+v)+(1+v)^{-1})'=v/(1+v)^2>0$.
The last bracket vanishes at zero and has derivative
$(1+r+r^2)/(1+r)^2>0$.  Symmetry gives $F_b>0$.
\end{proof}

The following finite state-graph calculation is a domain bound, and is
kept separate from the final cost exclusion.  Assume only the closed
versions of \eqref{sh:PMP-pairs}, $a_i\in[1/5,11/4]$, and
$\mathscr J_{\rm ext}\leq V$.  Divide the interval into closed cells
$C_j=[l_j/D,u_j/D]$, where $u_j=l_j+1$.
Retain a directed pair $(j,k)$ only if
\begin{equation}\label{sh:DP-pair}
 u_j^2+2u_ju_k\geq D^2,\qquad
 u_k^2+2u_ju_k\geq D^2,\qquad l_jl_k\leq D^2.
\end{equation}
Every discarded pair is impossible.  On a retained pair take a directed
lower bound for $F(l_j/D,l_k/D)$, multiply by $10^9$ and round down to an
integer $W_{jk}$; otherwise set $W_{jk}=+\infty$.
Lemma~\ref{sh:state-monotone} proves that these are lower bounds on
entire cells, although their lower corners need not satisfy group
closure.

Let $\otimes$ denote minimum-plus multiplication and write $W^m$ for
its powers.  For a cycle whose zeroth cell is $i$,
\[
 B_i=\min_j\{(W^2)_{ij}+(W^5)_{ji}\}-7\cdot10^9
\]
is a necessary cost lower bound.  At $D=100$ on the full state interval,
there are 255 cells and 13,896 retained pairs; the inequalities
$B_i\leq1805000000$ restrict every state to $[9/25,7/5]$.
On this already established interval repeat with $D=2000$.
There are 2,080 cells and 2,514,283 retained directed pairs, and the
necessary state range becomes
\begin{equation}\label{sh:DP-state}
 757/2000\leq a_i\leq67/50.
\end{equation}
For a proposed $q_0\geq2$, retain a central triple of cells $(i,j,k)$
only when
\[
 \frac{(u_i+u_j)(u_i+u_k)}{D^2+l_i^2}\geq2.
\]
This is a necessary upper bound, with no monotonicity assumption on
$q_0$ as a function of $a_0$.
The exact minimum of
$W_{ij}+W_{ik}+(W^5)_{jk}-7\cdot10^9$ over the remaining triples is
$1811660735$.  Dividing by $10^9$ gives a bound strictly larger than
$V=1.805$.  Thus every low-cost candidate has
\begin{equation}\label{sh:q-cutoff}
 q_i<2\quad\hbox{for every }i.
\end{equation}
All edge enclosures are completed before these integer minima are taken;
the remaining recurrence is finite integer arithmetic.  The retained
triple count is 731,743,421.  It describes a finite enlarged-state proof,
not a sampling of geometric trajectories.

Cut at a largest exponent and write $T_0=\max_i t_i=t_0$.
If $T_0\geq1$, then $T_0>11/20$ immediately.
Otherwise \eqref{sh:pure-cyclic} gives
\[
 T_0^2\geq\frac{(1-T_0^2)^2}{1+T_0},\qquad
 T_0^3-2T_0^2-T_0+1\leq0.
\]
The cubic is decreasing on $[0,11/20]$ and has value $91/8000$ at its
right endpoint.  Hence $q_0>\sqrt{31/20}$ in either case, and $E>11/20$.
This avoids a singular tail denominator at a zero-duration boundary.

For the interval divisions one may use a uniform rational discriminant
floor on the filtered candidate set.  Here is its quantitative proof.
Put $\rho_i=\sqrt{r_i}$ and let $C_i$ be the periodic half-state adjacent
cofactor in Lemma~\ref{sh:half-state}.  Since $q_i\leq2$, $r_i\geq1/2$.
More generally, if $r_i\geq\delta$, the notation in that proof gives
$X_h,Y_h\geq\delta^2/3$ and $R\geq\delta^6$.
Set $z=\sqrt{XY}$, $w=\sqrt R$ and $A=D_{5,h}+\sqrt R$.  Then
\[
 A\geq\delta^4/9-z,\qquad
 A\geq z+w-\sqrt{z^2+w^2}\geq2zw/3.
\]
The last bound uses $z^2+w^2=bc\leq1$, so that the rationalizing denominator is at most $1+\sqrt2<3$.
Splitting at $z=\delta^4/18$ gives $A\geq\delta^7/27$ and hence
\[
 \frac{\rho_i C_i}{\prod_j\rho_j}\geq\frac{\delta^8}{27}.
\]
Take $c=1/6912$ and write $b=\tr(P_0/\sqrt{\prod q_i})=2\cos\phi$,
$e=2-b$.  Summing the angular PMP inequalities gives
$\sum\alpha_i>7\pi/6$; by \eqref{sh:phase}, $0<\phi<\pi/3$,
so $0<e<1$.
Differentiating the seven-cycle characteristic identity gives
\[
 a_i\sin\phi=\frac{\rho_iC_i}{\prod_j\rho_j}-\frac e2
              \geq c-\frac e2.
\]
Since $a_i\leq7/5$ and $\sin\phi\leq\sqrt e$, either $e\geq c$,
or $(7/5)\sqrt e>c/2$.  In both cases the normalized discriminant obeys
\[
 \overline\Delta:=\frac\Delta{\prod q_i}=e(4-e)
 \geq\frac{25}{3121348608}>\frac1{125000000}.
\]
Moreover $\overline n_i:=n_i/\sqrt{\prod q_i}>1/32000$, using
$a_i\geq9/25$.  The same weaker floors hold for the raw $\Delta,n_i$.
These are necessary candidate bounds and may be intersected with value
enclosures after the filter is imposed.  They are not bounds on unrelated
points of a raw parameter rectangle and are not used to clip its
Taylor derivatives.

\begin{lemma}[Finite pure-forward inequality]\label{sh:pure-finite}
On the chart \eqref{sh:pure-chart}, every point satisfying the necessary
low-cost filters \eqref{sh:PMP-pairs}, \eqref{sh:DP-state},
\eqref{sh:q-cutoff}, the maximum cut and the positive state branch has
$\mathscr J>U$.
\end{lemma}
\begin{proof}
It suffices to prove positivity of the explicit numerator
\begin{equation}\label{sh:Phi}
 \Phi_0=3\sum_i n_i n_{i+1}-\sum_i n_i^2\log q_i-(7+U)\Delta.
\end{equation}
Use the closed root in $(q_0,q_1,q_2,q_3)$ with
$q_0\in[q_{\rm lo},2]$, $q_1,q_2,q_3\in[1,2]$, where
$q_{\rm lo}=$ \texttt{0x1.3eb7a9792f3c5p+0} is an exact binary rational below
$\sqrt{31/20}$.  Impose $q_j\leq q_0$, the tail nonnegativity and all
cyclic adjacent-product inequalities.  In numerator coordinates the
state filters are
\begin{gather*}
 n_i^2+2n_in_{i-1}-\Delta\geq0,\quad
 n_{i-1}^2+2n_in_{i-1}-\Delta\geq0,\quad
 \Delta-n_in_{i-1}\geq0,\quad3n_in_{i-1}-\Delta\geq0,\\
 (4489/2500)\Delta-n_i^2\geq0,\qquad
 n_i^2-(573049/4000000)\Delta\geq0,\qquad n_i>0.
\end{gather*}
Products of the matrices \eqref{sh:state-return} define all the
polynomials; there is no fitted polynomial kernel.

For reproducibility we describe the finite witness and its checking
rules.  Its 16 root slabs share exact binary endpoints and partition the
root in the $q_0$ direction.  Each node records its rational box, any
necessary-domain contraction and either a two-child split or one of the
strict terminal tests below.  There are 139,838 nodes, 69,911 splits and
69,927 terminal leaves; all 16 trees end without a pending box.
A terminal box must satisfy one of the following:
\begin{enumerate}
\item a necessary nonnegative constraint has a strictly negative upper
bound, or a positive denominator/branch requirement is impossible;
\item a directed lower bound for \eqref{sh:Phi} is strictly positive;
\item for specified nonnegative rational numbers $\lambda_k$ and
necessary constraints $g_k\geq0$, a directed lower bound for
$\Phi_0-\sum_k\lambda_kg_k$ is strictly positive.
\end{enumerate}
The last test is valid because it bounds a function no larger than
$\Phi_0$ at every candidate.  A heuristic may suggest the numbers
$\lambda_k$ but does not supply any exclusion.

Besides elementary interval contractions, the exact cyclic identity
allows the following one.  If upper bounds for
$t_{i-1},t_{i+1},t_{i+3},t_{i+4}$ are $A,B,C,D$ and $P=CD<1$, then
\[
 t_i\geq\frac{2(1-P)}{A+B+P+
 \sqrt{(A-B)^2+P(2(A+B)+P+4AB)}}.
\]
This is the rationalized first positive root of
$(1-Au)(1-Bu)-P(1+u)$. The polynomial is positive at zero and
nonpositive at each finite reciprocal $1/A,1/B$, so the root lies
at or below every finite reciprocal. For an actual $u=t_i$ below it,
both $1-Au$ and $1-Bu$ are positive. The numerator in
\eqref{sh:pure-cyclic} is then at least $(1-Au)(1-Bu)>P(1+u)$,
contradicting $t_{i+3}t_{i+4}\leq P$. A zero denominator in the root
formula means $A=B=P=0$, which also contradicts
\eqref{sh:pure-cyclic}. These contractions preserve candidates.
Derivative bounds used on the full uncontracted box must still hold
throughout that box.

The directed arithmetic and complete tree-replay rules are given
in Appendix~\ref{sh:arithmetic-forward}.  They verify every
terminal test on its whole box and the exact coverage of all
sixteen roots.  The state-graph inequalities and matrix identities
above are separate prerequisites.
Positivity at all terminal leaves, together with exact
tree coverage, proves \eqref{sh:Phi} on the necessary set.  Since
$\Delta>0$ by Lemma~\ref{sh:half-state}, the stated cost inequality
follows.
\end{proof}

\begin{proof}[Proof of Proposition~\ref{sh:forward-seven}]
If an actual candidate has cost greater than $V$, it cannot minimize.
Otherwise the state filter and finite state graph apply.  The maximum cut
puts it in the complete chart used by Lemma~\ref{sh:pure-finite},
which gives $\mathscr J>U>\mathscr J_8$, contradicting the octagonal
competitor.  The use of the larger prerequisite budget $V$ ensures that
the domain reduction does not assume the conclusion at budget $U$.
\end{proof}

\begin{remark}\label{sh:scope}
The proof combines several kinds of result.  The reverse five- and
seven-arc inequalities compare all feasible trajectories in those word
classes.  The six-arc arguments give actual decreases within their
closed word families.  The mixed and pure-forward seven-arc arguments
use complete maximization for an original global minimizer.  Thus
Proposition~\ref{short:rigidity} supplies the lower end of the global
candidate classification.  It does not assert the cost inequality
$\mathscr J\geq\mathscr J_8$ for every feasible seven-arc chain.
\end{remark}
\space
\label{revision:appendix-end-1}
\end{proofappendix}
\section{Fundamental local estimates}\label{app:local}
\begin{proofappendix}
\section{Cap transport and the cost of an own pair}\label{new:local}
The hypotheses in this section always include the original global
minimizer framework, genuine full endpoints and real winding. A
``proper source'' means fewer than $B$ complete positive maximal arcs.
A genuine cut in an omitted positive arc retains the source. The
selected nonzero maximizing lift is first proved normal and only then
jointly normalized. All means in an identity are actual means on that
same full lift. Algebraic comparison points need not be realizable.

We use the original harmonic ceiling
\[H_*:=\Hcap.\]

\subsection{The positive-cap identity}
Use an actual returned leaf on each side of a positive shared cap.
Write $u,v$ for their positive outer-root magnitudes, $h^2,\nu$ for
the corresponding positive actual weighted means, and $x$ for the
shared reciprocal squared. Put
\[
 A=1+x,\quad k=\frac{4x}{1+x},\quad
 t=\frac{k}{(1-u)(1-v)},\quad \Phi(t)=t-1-\log t.
\]
The actual domain is
\[
0<h<u<1,\quad0<\nu<v^2,\quad0<v<1,\quad0<t<1.
\]
Let $\lambda$ be the oriented integration constant of this actual cap,
and put $E=\lambda+3x/2$. The positive-cap identity is
\begin{equation}\label{new:capidentity}
 E=\frac A2\{\Phi(t)+(1-t)(u+v)+tuv\}
   +\frac A4\left\{\frac{u^2-h^2}{1-u}
                         +\frac{v^2-\nu}{1-v}\right\}.
\end{equation}
Here $E$ is a matrix-transport term. To derive the identity, set
$K=(3+h^2)/2$, $\widehat K=(3+\nu)/2$ and
$\beta=1-t(1-u)$. The actual binary primitive and cap equality
of Section~\ref{pf:M002} read
\begin{align*}
 \chi_\beta&=x\left\{\frac{h^2-1}{1-u}(t^{-1}-1)
                  -2\log t+(1-u)(t-1)\right\},\\
 f_\beta&=\frac{(h^2-\beta^2)x}{2(1-\beta)^2}
                      -\frac{A\chi_\beta}{2(1-\beta)},\\
 0&=\frac12-\frac{\widehat K A^2}{4}
       -\frac{2\lambda x}{1-\beta}
       -\frac{Kx^2}{(1-\beta)^2}-f_\beta.
\end{align*}
Solving for $\lambda$, using $4x/A=t(1-u)(1-v)$, gives
\eqref{new:capidentity}. In particular the identity preserves the
oriented constant; neither its square nor a separately normalized cap
has been substituted. No physical cap cost is added to the charge.

\begin{lemma}[Root-width contribution of one cap]\label{new:onecap}
In the preceding actual domain,
\[
 E+\frac{x(h^2+3u)}{2(1-u)}
 >\frac u2+\frac{(1-x)(u^2-h^2)}{4(1-u)}>\frac u2.
\]
\end{lemma}
\begin{proof}
Set $t_0=k/(1-u)$. The actual positive other root gives
$v=1-t_0/t>0$, so $t_0<t<1$. Substitution yields
\[
\Phi(t)+(1-t)(u+v)+tuv
=u+k-\log t-\frac{k}{t(1-u)}.
\]
For fixed actual $x,u$, the derivative of the last two terms is
$(t_0-t)/t^2<0$. Their value is therefore strictly larger than
its comparison value at $t=1$. Dropping only the positive other-leaf
deficit in \eqref{new:capidentity} gives
\[
E>\frac{A(2u-u^2-h^2)}{4(1-u)}-\frac{2xu}{1-u}.
\]
Adding the displayed density correction gives the first inequality.
The actual positive shared arc and positive outer roots imply $x<1/3$,
so the final remainder is strictly positive. No physical zero-length
cap or separate calibrated comparison lift has been introduced.
\end{proof}

\subsection{The own-pair theorem}
Consider the complete proper source $01020212$. Let
$-\xi<z<\rho$ be the three actual roots of its own pair, with
common weighted mean $h^2$, $K=(3+h^2)/2$, and positive shear $T$.
Set $w=\xi+\rho$ and let $x_L,x_R$ be the two actual anchor densities.
The physical cost on local arcs $\{3,4\}$ is exactly
\begin{equation}\label{new:owncost}
Q=\frac w2+KT.
\end{equation}
The original full-matrix transport gives
\begin{equation}\label{new:owntransport}
KT=E_L+E_R+
\frac{x_L(h^2+3\xi)}{2(1-\xi)}+
\frac{x_R(h^2+3\rho)}{2(1-\rho)}.
\end{equation}
The two cap defects in this identity are not two physical payments.

\begin{theorem}[Analytic own-pair bound]\label{new:own}
On this proper source,
\[
KT>\frac w2+
 \frac{(1-x_L)(\xi^2-h^2)}{4(1-\xi)}+
 \frac{(1-x_R)(\rho^2-h^2)}{4(1-\rho)},\qquad Q>w.
\]
Consequently the original source-qualified ceilings
$x_L,x_R<47/250$ imply $Q>9/25$.
If two additional qualified sources give $x_L,x_R<1/6$, then
$Q>3/7$. The latter density hypothesis is not a consequence of the
bare eight-letter source.
\end{theorem}
\begin{proof}
Apply Lemma~\ref{new:onecap} to the two actual caps in
\eqref{new:owntransport}. To obtain the rational constants, put
$a=\sqrt{x_L}$, $b=\sqrt{x_R}$ and
\[
U=\frac{1+\xi}{1-\xi},\qquad V=\frac{1+\rho}{1-\rho}.
\]
Direct integration of the two physical shear profiles gives
\[
T=abUV-\frac{a^2U+b^2V}{2}.
\]
If $a,b<\sqrt C$, scale both by the same factor until one reaches
$\sqrt C$. This strictly increases $T$, because the actual $T$ is
positive and the expression is homogeneous of degree two. On the
resulting boundary $a=\sqrt C$, the derivative in $b$ is
$V(\sqrt C U-b)>0$; the symmetric argument handles the other
boundary. Hence, with $P_0=\xi\rho$,
\begin{equation}\label{new:shearupper}
T<C\frac{w+2P_0}{1-w+P_0}.
\end{equation}
This does not assert coordinatewise monotonicity throughout the
unrestricted anchor square. It is a comparison using positivity and
homogeneity first.

Suppose $0<w\le w_0<2/3$. The derivative of the fraction in
\eqref{new:shearupper} with respect to $P_0$ is
$(2-3w)/(1-w+P_0)^2>0$. Since $P_0\le w^2/4$,
\[
T<C\frac{w(1+w/2)}{(1-w/2)^2}.
\]
The actual two leaf-quadratic inequalities give
$w>\sqrt{12h^2-3z^2}\ge3h$; thus
$K<(3+w^2/9)/2$. Combining this with $KT>w/2$ forces
\[
1<\frac{2KT}{w}<C\frac{(3+w^2/9)(1+w/2)}{(1-w/2)^2}.
\]
The last expression increases in $w$ on the relevant interval.
For $(C,w_0)=(47/250,9/25)$ its endpoint is
$5224332/5253125<1$. For $(C,w_0)=(1/6,3/7)$ it is
$2516/2541<1$. Both assumptions $w\le w_0$, including equality,
are impossible. Equation~\eqref{new:owncost} then proves the claims.
\end{proof}

\begin{corollary}[Two adjacent own pairs]\label{fee:adjacent720}\label{six:prop:adjacent}
A proper complete source $01020212101$ satisfies
\[
 c_3+c_4+c_5+c_6+c_7>\frac{18}{25}.
\]
\end{corollary}
\begin{proof}
Its source on $[0,7]$ is $01020212$ and charges $\{3,4\}$.
Its source on $[3,10]$, after whole phase $+1$, is the same source
and charges $\{6,7\}$. Both are proper and retain all their caps
and actual returns. Theorem~\ref{new:own} pays more than $9/25$
on each pair. Adding the two disjoint physical charges and the
positive cost of arc $5$ proves the assertion. This addition requires
no comparison of the two source multipliers.
\end{proof}

\subsection{A six-arc density source}
\begin{theorem}[Shortened density support]\label{new:sixdensity}
A complete proper $201020$ source, or a lawful whole image, gives
$x<1/6$ for the own reciprocal squared of its local arc $3$.
The conclusion is a physical state inequality; it adds no cost charge.
\end{theorem}
\begin{proof}
The first five arcs are a cyclic image of $12021$. Its binary core
has two distinct internal ties and both third-mode caps; full-lift
rigidity excludes abnormality. The last five arcs are $01020$ and
supply the two actual returned means, their roots and both midpoint
ties. Thus the normality witness and all three matrix rows below
belong to the stated six-arc support.

Put $R=1+x^{-1}$, $A=1+x$, and use the actual root variables
$p,q\in(0,1)$ and means
$0<\mu_L<(1-p)^2$, $0<\mu_R<(1-q)^2$.
The positive shared arc gives $Rpq>4$. The same matrix gives
\begin{align*}
V_0&=\log(R/4)+g(p)+g(q),&g(z)&=\log z+(z^{-1}-z)/2,\\
B_0&=(1-AV_0)/2+A(\mu_L/p+\mu_R/q)/4,\\
\Gamma&=A^2(3+\mu_L)(3+\mu_R)/16-B_0^2.
\end{align*}
At fixed $R,p,q$, both mean derivatives of $\Gamma$ are positive
on the entire closed mean rectangle. Indeed the first, multiplied
by $16p/A^2$, has its minimum at the upper mean corner. With
$k=4/R$, $t=pq$, that minimum is
\[
(1-q)(4-t)+k-t+4\log(t/k)>0;
\]
the last two terms vanish at $t=k$ and their derivative in $t$ is
$-1+4/t>0$. The other derivative is symmetric. Therefore
\[
\Gamma\ge\frac{9A^2}{16}-\frac{(1-AV_0)^2}{4}.
\]
On $4<R\le7$, the elementary bounds
\[
\log(R/4)\le V_0\le R/8-2/R,\qquad AV_0\le11/16<1
\]
permit squaring with the correct sign, giving $\Gamma\ge\phi(R)$,
where
\[
\phi(R)=\frac{9R^2}{16(R-1)^2}
 -\frac14\left(1-\frac{R\log(R/4)}{R-1}\right)^2.
\]
Writing $z=\log(R/4)$, one obtains
\[
8(R-1)^3\phi'(R)=-9R+4(R-1-Rz)(R-1-z).
\]
Both factors are positive: the preceding logarithm bound gives
$R-1-Rz\ge R+1-R^2/8\ge15/8$ and $R-1-z>2$.
Substitute the lower bound
$z\ge2(R-4)/(R+4)$ to bound the right side by
$-(R-4)^2(4R^2+R+4)/(R+4)^2\le0$.
Four positive atanh terms give
\[
\log(7/4)>\frac{381685152}{682050985},\qquad
\phi(7)>\frac{446891067930101121}{607599733733185600}>H_*.
\]
Thus this domain would have actual $\Gamma>H_*>0$.
An actual positive binary invariant supplies the future sheet on the
same full lift; the original harmonic theorem gives $\Gamma<H_*$.
This contradiction excludes $4<R\le7$, so $x<1/6$.
Only the actual lift is calibrated; the zero-mean comparison corner
is not. All closed comparison boundaries used above are retained.
\end{proof}

Theorem~\ref{new:own} with $9/25$ applies on every qualified proper
returned eight-source. The stronger $3/7$ version requires both
additional density guards. They can be proved on separate genuine
cuts and exported as physical state facts, but this does not identify
the matrices or covectors of those separate lifts. Joint matrix
arguments still need the full common-lift theorem.

\space
\section{Local physical-cost inequalities}\label{sec:local-costs}

The local estimates fall into four families: capped binary cores, forward
pairs and their physical arms, returned cores and shared-arc unions, and
central arcs coupled by actual switching responses. We state their usable
forms here. Sections~\ref{new:local} and~\ref{fee:analytic} prove the
cap-transport, four-arc, opposite-centre and nonpositive-tail estimates.
The remaining local derivations are in
Sections~\ref{lp:binary}--\ref{lp:finite}, including finite rational
or interval witnesses for compact residual domains. This distinction is
important: the global argument adds physical costs, whereas the proofs of
several local inequalities compare quantities belonging to one full lift.

\subsection{The source and the charge}\label{lc:source-convention}

A source $j_0\cdots j_{L-1}$ means $L$ consecutive complete positive
maximal arcs of the original minimizing trajectory, in chronological
order. Its positions are numbered from zero. A charge $I$ means
$c(I)=\sum_{i\in I}c_i$, where $c_i=\int_{\text{arc }i}\ell\,dt$ is the
entire physical cost. Thus a positive cap or a returned arc may be required
in the source without belonging to its charge. An interval $[a,b]$ in a
charge denotes the integer set $\{a,a+1,\ldots,b\}$.

Each result retains its stated scope: explicit local hypotheses take
precedence over the complete-source convention below. Statements and proof
steps on partial arcs or algebraic, analytic, comparison or relaxed domains
require only their stated data; comparison points carry no actual-source or
lift hypothesis unless specified. For complete sources we retain the
original state and frame endpoints, real
winding, compact geometric restrictions and measurable-control maximum
principle of Section~\ref{fd:section}. A source is first placed in a
genuine containing fundamental interval, and a nonzero full maximizing
lift in the sense of Definition~\ref{fd:full-lift} is selected there.
Normality is established before jointly dividing
the objective multiplier, the matrix adjoint and the state covector by
its positive value. These operations are never performed separately on
the individual leaves of a source. For a proper source, $L<B$, a cut inside
the preceding positive arc gives the needed containment. Where two
sources must share a lift across different cuts, that extra assertion is
proved explicitly below.

A lawful image means a whole cyclic symmetry, or a whole reflection
combined with reversal of physical time. It acts on the state, frame,
matrix adjoint and the state covector, and transports the charge with
its source. A bare permutation of control names is not an additional
symmetry. In particular, a displayed source is not a statement about an
isolated control word with freely assigned endpoint data.

\subsection{Capped cores and their common binary algebra}

The binary family uses actual roots, exterior maximizing jets and the
same matrix adjoint throughout each source. The derivation is given in
Section~\ref{lp:binary}; it includes the order in which normality,
positive determinant, positive trace and harmonic calibration are proved.
The bounds below charge complete physical arcs rather than a scalar
relaxation of a middle-root mean.

\begin{theorem}[Capped core costs]\label{lc:binary}
Under the complete-source hypotheses above, the following costs have the
stated bounds. Every listed source includes its uncharged exterior arcs.
\begin{enumerate}[label=(\roman*)]
\item The sources $02010$, $10121$, $21202$ have
$c([1,3])>49/200$.
\item The sources $102021$, $210102$, $021210$ have
$c([1,4])>389/1000$.
\item The source $10210$, and its lawful images, has
$c([1,3])\ge3/5$.
\item The source $102102$, and its lawful images, has
$c([1,4])>3/5>62/125$.
\item The source $012021$, and its lawful images, has
$c([1,4])>81/200$.
\end{enumerate}
The third assertion is weak. Its equality case is not discarded.
\end{theorem}
\smallskip\noindent Proofs: Theorems~\ref{f245:theorem}, \ref{loc:A389}, \ref{wk:five600}, Proposition~\ref{strip:r496}, and Theorem~\ref{loc:S405}.

\begin{theorem}[One returned four-core]\label{lc:one-return}
Put
\[
 E_*:=\frac{4154032208439}{8759490110000}.
\]
The right sources $1020212$, $2101020$, $0212101$ satisfy
$c([1,4])>E_*$. The left sources $0102021$, $1210102$, $2021210$
satisfy $c([2,5])>E_*$. The two exterior caps and the returned arc are
positive and uncharged.
\end{theorem}
\smallskip\noindent Proof: Theorem~\ref{ek:E}.

\begin{theorem}[A four-core with two returns]\label{lc:two-return}
For the complete source $01020212$ on an original minimizer, with both
actual returned mean leaves and both exterior caps retained in a genuine
fundamental interval, the charge on arcs $2,3,4,5$ satisfies
\[
 c([2,5])>74/125>11/20>261/500.
\]
The same statement holds for whole lawful cyclic images. Such a source
is automatically proper by Lemma~\ref{loc:automatic-properness}.
\end{theorem}
\smallskip\noindent Proof: Section~\ref{pf:M101} and
Theorem~\ref{loc:D522}. The correlated shear envelope is stated
independently in Lemma~\ref{pf:two-density-envelope}.

\subsection{Forward pairs and four physical arms}

A forward pair needs its complete surrounding caps even when only its
two central arcs are paid. The two-arm estimate has a different support,
although some of its numerical consequences coincide. Their proofs are
grouped in Section~\ref{lp:forward}.

\begin{theorem}[The forward central pair]\label{lc:forward-pair}
In the complete positive source $012012$, and its lawful images,
\[
 c_2+c_3>\frac{180483}{500000}>\frac{20901}{62500}.
\]
The four supporting arcs $0,1,4,5$ remain complete and uncharged.
\end{theorem}
\smallskip\noindent Proof: Proposition~\ref{fp:final}.

\begin{theorem}[Four own arcs]\label{lc:own-four}
The complete six-sources $101202$, $212010$, $020121$ satisfy
\[
 c([1,4])>103/250.
\]
The seven-source $1012012$ has this bound on $[1,4]$, and
$0120121$ has it on $[2,5]$, with all supporting arcs retained.
\end{theorem}
\smallskip\noindent Proof: Theorem~\ref{fee:four}, with the actual
return and capped-forward endpoint signs established there.

\begin{proposition}[Scalar costs used in the finite catalogue]\label{lc:scalar-substitutions}
On the complete positive source $01201$, and either whole cyclic image,
\[
 c([1,3])>47/400.
\]
On the complete positive source $10120121$, and either whole cyclic image,
\[
 c([1,6])>609/1000.
\]
All supporting arcs and the genuine full-lift conditions above are retained.
\end{proposition}
\smallskip\noindent Proofs: Proposition~\ref{strip:forward1175} and
Theorem~\ref{loc:OwnSix609}. These are the exact thresholds used by the
finite source catalogue.

\subsection{Returns, entrance geometry, and shared unions}

A returned tail retains more information than a generic three-core.
It is essential to say which extra geometric inequality is supplied by
its source.

\begin{theorem}[Nonpositive and returned tails]\label{lc:returned}
Consider the complete capped core $1\mid202\mid1$, whose actual
entrance is $u=-a\le0$.
\begin{enumerate}[label=(\roman*)]
\item On $12021$, the charge $[1,3]$ exceeds $61/200$.
\item On $120212$, the same charge exceeds $379/1000$.
\end{enumerate}
The return in (ii) is the complete arc at position $5$.
No entrance-circle or LOW-mean condition is required for (ii).
In particular it implies both weaker bounds $361/1000$ and $3/8$.
The full source $012021$ gives $c([2,4])>61/200$;
$0120212$ gives $c([2,4])>379/1000$.
The full source $202120212$ gives $c([4,6])>379/1000$;
its complete reflected-time image $010201020$ gives
$c([2,4])>379/1000$. In all cases the entire sign-producing
context, both caps and the stated actual return are retained.
\end{theorem}
\smallskip\noindent Proofs: Theorem~\ref{fee:ST305},
Theorem~\ref{cmp:returned-379}, and
Corollaries~\ref{ek:short-sources} and~\ref{ek:H7-K}.

\begin{theorem}[Shared-arc and shared-derivative unions]\label{lc:unions}
The following are once-counted physical charges of complete sources:
\begin{enumerate}[label=(\roman*)]
\item $02120102$ has $c([1,6])>71/100$ (hence $>63/100$);
\item $20212101020$ has $c([2,8])>153/200$;
\item $0102010201020$ has $c([2,8])>79483/100000$.
\end{enumerate}
All means, switching derivatives and matrix comparisons used within an
individual assertion belong to one selected full normal lift.
\end{theorem}
\smallskip\noindent Proofs: Theorems~\ref{fee:W71}, \ref{ua:AA765}, and \ref{ua:seven}.

\subsection{Actual central contacts and the terminal comparison}

The complete-terminal and prefix estimates have different sufficient
interfaces. Let a complete mode-$1$ arc be entered at $2\to1$ and exited
at $1\to0$ on one full normal lift. Its rate
$k=2v/(u^2+v^2-1)>0$ is constant. In the whole cyclic chart making it
mode $0$, put $\xi=1+u_{\rm new}$, $v_{\rm new}=k\xi$,
$\xi_L=A$, $\xi_R=AE$, $E=e^{kd}>1$ and
$\tau=\sqrt{1+k^2}$. The actual domain is $AE<2$, $A\tau^2>2$.
The raw switching response $q$ is given by the complete solution
\eqref{hc:normal-family}--\eqref{hc:response}. The scaled switching
functions $U,V$ in that solution satisfy the actual contact conditions
\begin{equation}\label{lc:contact-jets}
 V(A)=U(AE)=0,\qquad V_\xi(A)\le0,\qquad U_\xi(AE)\le0.
\end{equation}
They are one-sided derivatives from the central arc, in the increasing
coordinate $\xi$; $U,V$ are not physical state coordinates.

An actual preceding mode-$2$ tail of full duration $d_0>0$ supplies
\emph{strict ingress} if
\begin{equation}\label{lc:ingress}
 u_L>0,\qquad q_L<0,\qquad
 q_L'+kq_L<\frac{b_L}{1-e^{-b_Ld_0}}q_L,\qquad
 b_L=\frac{v_L}{1-u_L}>0.
\end{equation}
A complete initial $202$ supplies these data. Regeneration must prove
them again at the new actual entrance, with its own rate and duration.
The central proofs are collected in Section~\ref{lp:centers}.

\begin{lemma}[The actual terminal inequality]\label{lc:terminal}
Suppose this central exit is followed by a complete positive $020$ on
the same lift, with $u_R<0$. Then
\[
 q_R>0,\qquad q_R'<(k-a_R)q_R,\qquad a_R=\frac{v_R}{1+u_R}>0.
\]
\end{lemma}
\smallskip\noindent Proof: Lemma~\ref{hc:terminal}. Its weighting has the variable rate $w^{\prime}/w=2v/(u^2+v^2-1)$; only its actual junction value equals $k$. Identity~\eqref{hc:C1-identity} identifies the same raw derivative across the switch.

\begin{theorem}[A complete seven-arc center]\label{lc:direct-center}\label{nx:central1196}
A properly retained complete positive source $2021020$ satisfies
\[
 c_3>299/250.
\]
Only its complete central arc is charged. More generally, the same
conclusion holds for a complete central arc satisfying
\eqref{lc:ingress}, both actual contact jets \eqref{lc:contact-jets},
and the complete same-lift terminal $020$ in Lemma~\ref{lc:terminal}.
\end{theorem}
\smallskip\noindent Proof: Section~\ref{nx:central-proof}, using the common central equations and terminal identities of Section~\ref{lp:centers}.

\subsection{Prefix propagation and disjunctive central costs}

\begin{theorem}[Actual prefix propagation]\label{lc:prefix}
For a complete actual central arc entered at $2\to1$ and exited at
$1\to0$, retain one full normal maximizing lift and the actual strict
ingress \eqref{lc:ingress}. If its physical cost satisfies
$C_{\rm loc}\le13/20$, then
\[
 u_R<0,\qquad q_R\ge0,\qquad q_R'>kq_R.
\]
The actual contacts supply both jets \eqref{lc:contact-jets} by full
maximization. No terminal word is assumed. The same conclusion holds,
in particular, under the weaker budget $C_{\rm loc}\le74/125$.
If a complete actual terminal $020$ is available on the same lift,
its comparison contradicts either small budget.
\end{theorem}
\smallskip\noindent Proof: Theorem~\ref{hc:prefix650}; the automatic
contact conditions and the weaker-budget consequence are recorded in
Theorem~\ref{hc:H592}.

\begin{corollary}[A regenerated alternative]\label{lc:regeneration}
For the complete proper ten-source $2021021020$,
\[
 \max(c_3,c_6)>13/20.
\]
If its two reverse four-core sources retain the hypotheses of
Theorem~\ref{lc:binary}(iv), then
\[
 c([2,7])>13/20+62/125=573/500.
\]
\end{corollary}
\smallskip\noindent Proof: Lemma~\ref{hc:regeneration} and Corollary~\ref{hc:ten-arc}. If $c_3>13/20$, pay the disjoint reverse-four charge $[4,7]$; if $c_6>13/20$, pay $[2,5]$. The branches are not added.

A paired argument requires one common lift for the two rooted tails.
The following criterion supplies it from genuine fundamental intervals;
it is stronger information than the independent validity of two scalar
cost inequalities.

\begin{lemma}[A sufficient common-lift criterion]\label{lc:gluing}
Let $a<b<a+T$ be two actual cuts. Suppose $[b,a+T]$ contains a complete
positive capped-binary core with two internal switching roots, and
$(a,b)$ contains a switch of distinct modes. The full lifts on the two
genuine fundamental intervals can be jointly normalized, agree on their
overlap, and extend to one normalized full reduced maximizing lift on
the actual unfolding.
\end{lemma}
\smallskip\noindent Proof: Lemma~\ref{pc:rigidity} and Proposition~\ref{pc:two-cuts}.

\begin{theorem}[Paired rooted responses]\label{lc:paired}
Suppose one normal full maximizing lift contains the complete positive
source $21021$, with middle mode-$0$ connector. Take its right rooted
mode-$2$ tail, and take the corresponding opposite rooted tail after
the lawful full reversed-time transformation. For a tail whose entrance
and exit horizontal coordinates are $U,V$, define
\[
 \Psi(U,V)=2(1-U)\log\frac{1-U}{1-V}-(V-U)(2+V-U).
\]
At least one of the two full-duration responses is strictly positive.
If both corresponding ordered extensions are complete proper $021020$
sources, their two central physical costs satisfy
\[
 \max(c_-,c_+)>299/250.
\]
The extensions retain their actual endpoints, compact geometry and full
maximizing lifts. If reanchoring is needed for the selected extension,
its positive response is state-determined and the ingress is rederived
on that complete lift. The other center has no separate lower bound
asserted by this theorem.
\end{theorem}
\smallskip\noindent Proof: Theorem~\ref{pc:paired-response}, Lemma~\ref{pc:strict-ingress}, and Corollary~\ref{pc:max-center}.

\paragraph{Proof order and finite arithmetic.}
Sections~\ref{lp:binary}--\ref{lp:centers}, placed after the global
assembly, prove these statements from the full equations. The remaining
closed-domain sign comparisons are in Section~\ref{lp:finite}.
Polynomial coefficients or interval cells are used only after the
actual source has been shown to satisfy their full domain. They cannot
supply a missing cap, a missing return, or a common lift. The weak
proper-five inequality remains weak, and every maximum inequality is
used branch by branch.
\space
\section{Two estimates for complete returned supports}
\label{six:sec:returned}
These local statements use the independently established returned-tail
interface, loss obstruction, LOW bound and HIGH absorption. They do not
use the returned-word exclusions or the finite census
conclusion. Their individual source dependencies remain those stated below.

The endpoint information of a returned support is stronger than that of
an isolated binary core. We first use both endpoint circles to estimate
a returned triple, and then retain the density dependence of the two
cap defects of a returned own pair. The inequalities in this section
are statements about complete supports with the interfaces specified
below; arbitrary formal triples do not acquire these hypotheses merely
by having the same labels.

\subsection{The returned triple interface}
Consider the complete positive source $120212$ in a genuine containing
interval. In its capped $202$ triple write the entrance as $-a\le0$,
the internal roots as $r<s$, the exit as $b$, and the actual middle mean
as $\mu$. The actual root/mean geometry gives
$r^2\le\mu<b^2$ and $s>\sqrt\mu>0$
Lemma~\ref{cmp:loss-obstruction} and its root identities.
Its physical cost is $T$. Put
\[
 q(u)=v(u)^{-2},\quad \sigma=\int_{-a}^bq(u)\,du,
 \quad Z=\frac{3+\mu}{2}\sigma,
 \quad X=\frac{f(-a)}\sigma,\quad Y=-\frac{f(b)}\sigma,
\]
\[
 \mathcal L=(X+Y+Z)^2-4ZY-\mu,
 \qquad \Delta_0=\frac{59}{1000}.
\]
The function $f$ is the binary switching response in the original
normalized lift. We use the returned-tail interface and loss obstruction
proved in Appendices~\ref{app:local}--\ref{app:crossed}: the actual signed
strip has $Y,Z>0$, $-Z<X+Y<Z$, and $\mathcal L\ge\Delta_0$.
Its LOW estimate gives $T>407/1000$ when $\mu\le a^2$.
On the HIGH branch $a^2<\mu<b^2$, the returned-HIGH absorption result
states that $w=a+b\le1/2$ and $Z<1/4$ imply $\mathcal L<Z^2$.
This last implication includes the exact polynomial comparison described
in Appendix~\ref{six:app:computation}; it is not being replaced by a
numerical plot or a sample of mean values.

We additionally retain the two actual circles
\begin{equation}\label{six:eq:two-circles}
 v(-a)^2>(1-a)(3+a),\qquad v(b)^2>(1-b)(3+b).
\end{equation}
The first is an entrance hypothesis; the last follows from the complete
return. They are not interchangeable endpoint assertions.

\begin{proposition}\label{six:prop:tail}
A complete returned triple satisfying this interface and
\eqref{six:eq:two-circles} has $T>2/5$.
\end{proposition}
\begin{proof}
Suppose $T\le2/5$. The LOW estimate leaves only the HIGH branch. Write
$w=a+b$. The exact physical moment identity gives
\begin{equation}\label{six:eq:tail-moment}
 0\le a<b,\quad 0<w<4/5,\qquad
 Z<\frac{(3+\mu)(4/5-w)}6.
\end{equation}
Both logarithmic density slopes belong to
$[-2/(1+u),2/(1-u)]$. Forward comparison from the entrance and backward
comparison from the exit, using the same actual density, give
\[
 q(u)<\min\left\{\frac A{(1-u)^2},\frac B{(1+u)^2}\right\},\quad
 A=\frac{(1+a)^2}{(1-a)(3+a)},\quad
 B=\frac{(1+b)^2}{(1-b)(3+b)}.
\]
The crossing $u_*=(\sqrt B-\sqrt A)/(\sqrt B+\sqrt A)$ lies in
$(-a,b)$. At the left endpoint this follows from
$(1-a)(3+a)\ge3(1-a)^2$; the right endpoint is the reflected comparison.
The excluded case $a=b=0$ is the only simultaneous equality.
Integrating on the two sides of $u_*$ yields
\begin{equation}\label{six:eq:tail-envelope}
 \sigma<\sqrt{AB}-\frac{1+a}{2(3+a)}-\frac{1+b}{2(3+b)}=:S(a,b).
\end{equation}

For $0<w\le1/3$ put $p=ab\in[0,w^2/4]$. Then
\[
 S=\frac{1+w+p}{\sqrt F}-\frac{3+2w+p}{9+3w+p},\qquad
 F=9-6w-3w^2+(10+2w)p+p^2,
\]
\[
 S_p=\frac{4(1-3w-w^2+p)}{F^{3/2}}
       -\frac{6+w}{(9+3w+p)^2}\le\frac4{27}.
\]
Indeed, if the first numerator is positive, use
$1-3w-w^2+p\le1-3w$, $F\ge9(1-w)$ and
$1-3w\le(1-w)^{3/2}$. The last inequality follows by squaring:
$(1-w)^3-(1-3w)^2=w(3-6w-w^2)\ge0$.
If that numerator is nonpositive the claimed derivative bound is immediate.
At $p=0$ the integral representation in \eqref{six:eq:tail-envelope}
is increasing in $w$: one endpoint coefficient is fixed, the other
increases, and the positive integration interval extends. Hence
\[
 S(w,p)\le S(1/3,0)+1/243
 =\frac2{\sqrt{15}}-\frac{11}{30}+\frac1{243}
 <\frac{749}{4860}<\frac{31}{200}.
\]
Here $(31/60)^2-4/15=1/3600$ and the final gap is $43/48600$.
Since $\mu<w^2\le1/9$, this gives
\[
 Z<\frac{217}{900}<\frac14,\qquad
 \Delta_0-\left(\frac{217}{900}\right)^2=\frac{701}{810000}>0.
\]
For $w\ge1/3$, the right side of \eqref{six:eq:tail-moment} is bounded by
$h(w)=(3+w^2)(4/5-w)/6$, whose derivative is negative. Therefore
\[
 Z<\frac{98}{405}<\frac14,\qquad
 \Delta_0-\left(\frac{98}{405}\right)^2=\frac{2939}{6561000}>0.
\]
In particular the entire range $w\le1/2$ contradicts HIGH absorption.

It remains $1/2\le w<4/5$. Put $I_0=(r+a)/((1-r)(1+a))$.
The actual entrance-response inequality of
Lemma~\ref{cmp:tail-interface} is
\[
 X\le\frac{\mu-a^2}{2L_0}+\frac{2G_0\mu}{3+a},\quad
 G_0=\frac{I_0}{(1-a^2)D_0}\le\frac1{1-a^2},
\]\[
 L_0=(1+a)^2D_0,
\]
\[
 D_0=\frac{r+a}{(1-r)(1+a)}
      +\frac{(1+r)(b-r)}{(1-r)^2(1+b)}.
\]
At fixed $a,b$, $(D_0)_r=4(b-r)/[(1+b)(1-r)^3]\ge0$.
Since $r>-a$, we obtain $L_0\ge w(1-a)/(1+b)$.
Now $a<w/2<2/5$, $b<4/5$, and consequently
\begin{equation}\label{six:eq:tail-X}
 X/\mu<3+50/63=239/63<4.
\end{equation}
Also \eqref{six:eq:tail-moment} gives $Z<(3+\mu)/20$.
For $\mu\ge1/25$, the clipped signed-strip bound gives
\[
 \mathcal L<4Z^2-\mu<\frac{(3+\mu)^2}{100}-\mu
 \le\frac{3276}{62500}<\Delta_0,
\]
where the last gap is $823/125000$; the comparison function decreases
on $[1/25,16/25]$. For $0<\mu\le1/25$, \eqref{six:eq:tail-X} and $Y>0$
give the unclipped bound
\[
 \mathcal L<(Z+4\mu)^2-\mu
 <\left(\frac{3+81\mu}{20}\right)^2-\mu.
\]
The last expression is convex. Its endpoint values $9/400$ and
$896/15625$ are below $\Delta_0$, with gaps $73/2000$ and
$207/125000$, respectively. Every remaining HIGH point is therefore
impossible. The width seams were included throughout.
\end{proof}

\subsection{Density-dependent payment for an own pair}
Consider a proper complete source $01020212$ with both actual returned
mean leaves and both exterior caps. Charge only its local arcs $3,4$.
Write their cost as $Q$, their genuine height as $h>0$, their shear as
$T>0$, and put $K=(3+h^2)/2$. The original same-lift identities are
\begin{equation}\label{six:eq:own-interface}
 Q=w/2+KT,\quad w>3h,\quad
 N=e_L+e_R=K(T-E)-\frac{h^2}2(x_L+x_R),
\end{equation}
where $0<x_L,x_R<47/250$ and $e_L,e_R\ge0$.
Writing $\delta=3/4-\Gamma$ for the determinant gap of this same
harmonically qualified lift, the cap gaps obey
\[
 \frac{29}{2000}<\delta\le\frac34(1-3x_j)e_j+e_j^2.
\]
Here $E=x_La/(1-a)+x_Rb/(1-b)$ is the cap-density term,
where $-a$ and $b$ are the two exterior roots of the actual own pair.
For $0<h<1/3$ the actual leaf-radius and density correlation give
\[
 R_h=\frac{8h}{4-3h}<1,\quad
 k=\frac{1-R_h}{1+R_h},\quad E>kT,\quad x_L+x_R\ge LT,
\]
\[
 L=\frac{(1-R_h)^2}{R_h(1+R_h)}
 =\frac{(4-11h)^2}{8h(4+5h)}.
\]
These identities use no numerical own-pair estimate. The exact transport
is \eqref{pf:d8-transport}; the chronological density ceilings and cap
gaps are \eqref{pf:d8-density-ceilings} and Section~\ref{pf:M033}.
The two actual middle means give $Q=w/2+KT$, and
Lemma~\ref{loc:leaf-correlations} gives $w>3h$.
For $h<1/3$, apply Lemma~\ref{pf:own-shear-correlation} with
$(a,b,R,q)=(\xi,r,R_h,q_z)$, $q_z=v(z)^{-2}$.
It gives $E>kT$ and $T\le2q_zR_h/(1+R_h)$. The same actual
density identities give
\[
 x_L+x_R\ge2q_zk^2(1+z^2)\ge2q_zk^2
 \ge\frac{(1-R_h)^2}{R_h(1+R_h)}T=LT.
\]
All these quantities belong to the one selected full lift and its
continuous physical profile. In particular, no conclusion of
Proposition~\ref{pf:own-fifth} is used to establish this interface.

The following numerical consequence depends only on the displayed
interface. It also applies when that interface has been recovered by
a one-sided-leaf construction; it does not require an additional
unproved returned source.

\begin{proposition}\label{six:prop:own}
Suppose the identities and inequalities of the preceding interface
hold on one selected actual lift. Then $Q>283/1000$.
\end{proposition}
\begin{proof}
Put $e_0=47/2500$. At $e=e_0+x/20$, the cap expression minus
$29/2000$ is exactly
\[
 -\frac{291}{6250000}-\frac{73x}{25000}-\frac{11x^2}{100}<0.
\]
It increases strictly with $e\ge0$ for the stated density range.
Thus $e_j>e_0+x_j/20$, and $N>47/1250+(x_L+x_R)/20$.
For $h\ge1/3$, \eqref{six:eq:own-interface} already gives $Q>1/2$.
For $h<1/3$, set
\[
 A=K(1-k)-h^2L/2=\frac{h(368+88h+7h^2)}{16(4+5h)},
 \qquad H=A-L/20.
\]
The same identities give $N<AT$, so
\begin{equation}\label{six:eq:own-H}
 HT>47/1250.
\end{equation}
In particular $H>0$ is a conclusion, not an assumption used in division.
If $h\ge283/1500$, the strict width inequality proves the assertion.
Otherwise a contrary $Q\le283/1000$ implies
\[
 \frac{47}{1250}K<H\left(\frac{283}{1000}-\frac{3h}2\right).
\]
The opposite difference, multiplied by $16h(4+5h)>0$, is
\[
 F=\frac{21}2h^5+\frac{131523}{1000}h^4
 +\frac{1275373}{2500}h^3-\frac{830077}{10000}h^2
 -\frac{1601}{1250}h+\frac{283}{625}.
\]
Its exact decomposition is
\begin{align*}
 F={}&510h(h-11/100)^2+131(h^2-121/10000)^2\\
 &+\frac{3217}{100}h^2-\frac{373}{50}h+\frac{541}{1250}\\
 &+\frac{21}2h^5+\frac{523}{1000}h^4+\frac{373}{2500}h^3
  +\frac{77}{400}h^2+\frac{41}{5000}h+\frac{82029}{100000000}.
\end{align*}
The nonsquare quadratic has positive leading coefficient and discriminant
$-2569/62500<0$. All other terms are nonnegative, with a positive
constant. Hence $F>0$, a contradiction.
\end{proof}
\space
\section{A conditional fee for a central triple}\label{new:conditional}
The following estimate keeps state information from the neighboring
sources. It will be used both to exclude R15 and to obtain a lower
bound for R06. In this section put
\[
 C_0=\frac{80483}{125000}<\frac23.
\]
All costs retain the original physical normalization.

\begin{theorem}[Central three-arc fee]\label{new:centralfee}
Let a complete proper six-source be a lawful whole image of $102010$
on an original globally minimizing trajectory. Its last five arcs
are $02010$, with successive own reciprocals $p,A,y,B,q$.
Select one full maximizing lift on a genuine interval containing the
six-source. Suppose qualified sources on the same trajectory give
\begin{gather}\label{new:centralguards}
 Ea>\max(1,b),\qquad Eb>\max(1,a),\qquad Ea,Eb<2,\\
 p<A,\quad q<B,\quad p^2+2pA<1,\quad q^2+2qB<1,
 \quad p^2,q^2<1/6,\notag
\end{gather}
where
\[
 a=\frac{2y}{y+A},\qquad b=\frac{2y}{y+B},\qquad
 E=\frac{(y+A)(y+B)}{1+y^2}.
\]
Then the complete physical cost $F=\alpha+A_c+\beta$ of the three
interior arcs of $02010$ satisfies $F>C_0$.
\end{theorem}
\begin{proof}
The first five arcs of $102010$, reversed with the whole label map
$j\mapsto2-j$, are $12021$. Their capped binary core excludes
abnormality. Jointly normalize the objective multiplier, entire
matrix and the state covector on the selected full lift.
All central means and matrix identities below belong to this lift.
The additional sources in \eqref{new:centralguards} export only
physical state inequalities; their multipliers need not be identified
with the central lift. The original harmonic bound gives
$\Gamma<H_*$: if $\Gamma\le0$ this is immediate; if $\Gamma>0$,
the actual binary invariant $K>3/2$ gives the future sheet on this
same lift before Corollary~\ref{cal:one-point} is applied.

Assume, for contradiction, that $F\le C_0$.
We first exclude two exterior ranges without finite computation.
Put
\[
 f(s,x)=\frac{1-sx}{s+x},\qquad
 \mathcal H(t)=1-3t-t\log\frac{1+t}{4t},\qquad
 h(y)=\sqrt{1+y^2}-y.
\]
The two crossings imply
$A,B>g(y):=\max\{h(y),(1-y^2)/(2y)\}$.
The exact physical leaf and half-center integrals are
\begin{align*}
 L(s;l,u)&=s(u-l)+s(1+s^2)
       \left(\frac1{s+l}-\frac1{s+u}\right)
       -s^2\log\frac{s+u}{s+l},\\
 Z_y(s)&=y[s-f(y,s)]-y^2\log\frac{y+s}{\sqrt{1+y^2}}.
\end{align*}
Thus
$F=L(A;p,y)+Z_y(A)+Z_y(B)+L(B;q,y)$.
For $y\le1/\sqrt3$, the mandatory central interval gives
$F>\mathcal H(y^2)$.
For $y\ge1/\sqrt3$, define
\[
 D_y(s)=
 \begin{cases}Z_y(s)+L(s;s,y),&s\le y,\\Z_y(s),&s\ge y.\end{cases}
\]
The actual leaves contain these mandatory portions, so
$F>D_y(A)+D_y(B)$. For $s<y$,
\[
 D_y'(s)=s\left(2z-3+\frac{2z}{1+z}
                         -2\log\frac{1+z}{2}\right),\qquad z=y/s>1.
\]
The bracket vanishes at $z=1$ and has derivative
$2-2z/(1+z)^2>0$. For $s\ge y$,
$Z_y'(s)=y[1+(1-ys)/(y+s)^2]>0$.
Consequently $D_y$ is increasing and
$F>2D_y(h(y))=2\mathcal H(h(y)^2)$.
These are comparisons of the complete center and mandatory leaf
portions, with their common state preserved.

The function $\mathcal H$ decreases on $(0,1/3)$.
The positive exponential partial sums give
$\log(53/16)<5/4$ and $\log(97/64)<5/12$.
Since $h(65/72)=4/9$, the two closed ranges satisfy
\begin{align*}
 y\le2/7&\ \Longrightarrow\ F>32/49
       =C_0+\frac{56333}{6125000},\\
 y\ge65/72&\ \Longrightarrow\ F>158/243
       =C_0+\frac{192631}{30375000}.
\end{align*}
It remains to treat $2/7<y<65/72$.

\paragraph{The exact central matrix and leaf floors.}
Set
\begin{gather*}
 k=\frac{4y^2}{1+y^2},\quad D=4-k,\quad P_0=ab,\quad
 E=\frac{k}{P_0},\\
 z=\log E,\quad t=\frac{2z}{D},
 \quad r=1-a,\quad s=1-b.
\end{gather*}
The crossing guards give
\[
 Ea=k/b,\qquad Eb=k/a,\qquad
 k/2<a,b<\min\{k,\sqrt k\}.
\]
Indeed $Ea,Eb<2$ give the lower bounds; $Ea,Eb>1$ give
$a,b<k$, and $Ea>b$, $Eb>a$ give $a^2,b^2<k$.
These last two inequalities also imply $E>1$.

The actual mean excesses are
$X=(\mu_A-r^2)/(aD)>0$ and $Y=(\mu_B-s^2)/(bD)>0$.
The original binary matrix and physical cost identities give
\begin{gather}\label{new:centralmatrix}
 K_A=\frac{3+r^2+aDX}{2},\qquad
 K_B=\frac{3+s^2+bDY}{2},\\
 J_m=X+Y+\frac{2z}{D}+\frac{4(a+b)-4-k}{2D},\qquad
 \Gamma(X,Y)=\frac{4K_AK_B}{D^2}-J_m^2,\notag\\
 A_c=\frac{2(a+b)(E-1)-kz}{D},\qquad
 J_a=2a+2k/a-k,\qquad J_b=2b+2k/b-k.\notag
\end{gather}
For the left leaf put $Q=2-a$ and
$d=(A-p)/(A+p)-(A-y)/(A+y)>0$.
The terminal circle and the strict leaf-gap identity
\eqref{f245:leaf-gap}, with $Eab=k$, $1+r=Q$ and $J_A=J_a$, imply
\[
 \alpha>m_a^{\rm term}:=
       \frac{J_a(2-a)(k-a)}{2kD}.
\]
Indeed $p<\sqrt{1+A^2}-A$ implies
$d/(1+d/Q)>Qa(k-a)/(2k)$, while the leaf-gap cost exceeds
$J_ad/[aD(1+d/Q)]$.
When $k\ge4/7$, the additional density $p<1/\sqrt6$ gives
\[
 \alpha>m_a^{\rm dens}:=
       \frac{J_a(2-a)}{2D}\left(1-\sqrt{\frac D{6k}}\right).
\]
Use $m_a=m_a^{\rm term}$ for $k\le4/7$ and
$m_a=\max\{m_a^{\rm term},m_a^{\rm dens}\}$ thereafter.
The right leaf gives $\beta>m_b$ by the identical actual calculation.
The two floors for one leaf are combined by maximum, not by addition.

\paragraph{A signed cost-to-mean comparison.}
For comparison cost $c\in[0,C_0]$, define
\begin{gather*}
 \lambda_a=\frac{aD}{J_a(2-a)},\quad w=\lambda_a c,\\
 T_a(c)=\frac c{J_a}\left\{\frac r{1-w}
                   +\frac{(2-a)w}{3(1-w)^2}\right\},\\
 T_a^+(c)=\max\{0,T_a(c)\}.
\end{gather*}
The denominator is uniformly positive. For $k\le1$,
$\lambda_a\le1$. On the high domain $\lambda_a\le4/3$, since
\[
 a[4J_a(2-a)-3aD]\ge a^2(a-2)(7a-10)>0
\]
for $1\le k<2$, $a\le\sqrt k<\sqrt2<10/7$.
Hence $w\le8/9$ even on the larger cost interval $[0,2/3]$.

The actual leaf gap gives $d<(2-a)w/(1-w)$ at $c=\alpha$.
For $r\ge0$, decreasing-weight covariance gives
$\mu_A-r^2\le rd+d^2/3$.
For $r<0$, integration by parts gives the same bound with the
right side replaced by its positive part: for $0\le v\le d$,
\[
 \frac vd\left(rd+\frac{d^2}{3}\right)
       -\left(rv+\frac{v^2}{3}\right)=\frac{v(d-v)}3\ge0.
\]
The actual excess is positive, so $d>-3r$; the quadratic is increasing
on the branch used in the substitution. Thus $X<T_a^+(\alpha)$
in either sign. This argument does not assume $a\le1$.

The function $T_a^+$ is increasing and convex on the entire comparison
cost interval. For $r\ge0$ this follows from its two nonnegative
power series. For $-1<r<0$, its positive branch starts at
$w_0=-3r/(1-2r)$; for
$f_r(w)=w[r/(1-w)+(1+r)w/(3(1-w)^2)]$,
\[
 f_r'(w)=\frac{3r+(2-r)w}{3(1-w)^3},\qquad
 f_r''(w)=\frac{2[1+4r+(2-r)w]}{3(1-w)^4}.
\]
Both numerators are positive at $w_0$ and increase thereafter.
Joining this branch to zero retains convexity, including $a=1$.

\paragraph{Eliminating the allocation.}
The Hessian of \eqref{new:centralmatrix} in $X,Y$ has eigenvalues
$-P_0$ and $P_0-4$, both negative. Its coordinate derivatives are
nonpositive on the comparison quadrant. To check the latter, the
coefficient of $X$ in its expanded quadratic is
\[
 A_X=s-2t+\frac{(k-P_0)(2-b)}D.
\]
For $k\le1$, its maximum is at $a=b=k$, where
$DA_X\le(1-k)(4+k-k^2)+4\log k\le0$.
For $1\le k<2$, its maximum is at $a=\sqrt k,b=1$, where
$DA_X\le k-\sqrt k-2\log k\le0$.
The first one-variable expression increases to zero at $1$, and the
second decreases from zero. The assertion for $Y$ is symmetric.

Therefore $\Gamma(T_a^+(\alpha),T_b^+(\beta))$ is decreasing
and concave in the two comparison costs. Put $U=C_0-A_c$.
The actual inequalities give $\alpha>m_a$, $\beta>m_b$ and
$\alpha+\beta\le U$. If $A_c+m_a+m_b>C_0$ the budget is
impossible. Otherwise increasing total comparison cost to $U$ and
using the chord of the remaining allocation segment gives
\begin{equation}\label{new:allocation}
\Gamma_{\rm actual}\ge
\min\left\{
 \Gamma(T_a^+(m_a),T_b^+(U-m_a)),\quad
 \Gamma(T_a^+(U-m_b),T_b^+(m_b))\right\}.
\end{equation}
Both endpoints are required. If $\alpha+\beta=U$, no enlargement
is needed; the same chord argument applies. Thus this reduction
retains the boundary $F=C_0$. Artificial allocation endpoints need
not represent trajectories and are not assigned calibrated lifts.

All actual data with $2/7<y<65/72$ lie in the closed domains of
Lemma~\ref{new:r15finite} below by \eqref{new:centralguards}.
That lemma and \eqref{new:allocation} give
$\Gamma_{\rm actual}\ge H_*$, contradicting
$\Gamma_{\rm actual}<H_*$ on the selected full lift.
Together with the two exterior ranges this proves $F>C_0$.
\end{proof}

\begin{lemma}[Three closed comparison covers]\label{new:r15finite}
Every point in the following closed domains satisfying
$A_c+m_a+m_b\le C_0$ has both lower bounds in \eqref{new:allocation} at least $H_*$:
\begin{center}
\begin{tabular}{ccl}\toprule
$k$ range&$(a,b)$ range&High-domain restriction\\\midrule
$[16/53,4/7]$&$[k/2,k]^2$&none\\
$[4/7,1]$&$[k/2,k]^2$&none\\
$[1,16900/9409]$&$[k/2,\sqrt k]^2$&$a^2,b^2\le k$\\\bottomrule
\end{tabular}
\end{center}
\end{lemma}
\begin{proof}[Computer-assisted proof]
Here $A_c\ge0$: $E\ge1$, $a+b\ge k$ and
$\log E\le E-1$ give
$DA_c\ge[2(a+b)-k](E-1)\ge0$.
Thus the stated budget condition and nonnegative leaf floors give
the comparison cost interval used in \eqref{new:allocation}.
Write $a=kx$, $b=kv$ and start from $(x,v)\in[1/2,1]^2$.
The high-domain curved inequalities are retained.
The witnesses have respectively $390193$, $13434119$ and $9922659$
nodes. Their interval evaluators use $48$, $46$ and $46$ binary
fractional bits, with $72$, $64$ and $64$ positive logarithm terms.
For $E\ge1$, $\zeta=(E-1)/(E+1)$, the logarithm tail after $N$ terms
is bounded by $2\zeta^{2N+1}/[(2N+1)(1-\zeta^2)]$.
Integer square roots enclose radicals by exact squared comparisons.
All products and quotients round outward by signed integer floor
and ceiling; stored endpoints have bounds $2^{60}$ in the first
cover and $2^{61}$ in the other two. The corresponding intermediate
products fit in signed $128$-bit integers. A denominator or range
failure cannot accept a leaf.

Each terminal box either violates the curved domain, has
$A_c+m_a+m_b>C_0$ throughout, or proves both expressions in
\eqref{new:allocation} at least the outward upper enclosure of $H_*$.
The finite trees use closed children sharing the split face and
have no unresolved box. Comparison costs are clipped only to
$[0,2/3]$, which contains every feasible allocation; clipping does
not discard the equality $F=C_0$ or the faces $a=1,b=1$.
Finite induction over these trees proves the lemma.
\end{proof}
\space
\providecommand{\tightlist}{\setlength{\itemsep}{0pt}\setlength{\parskip}{0pt}}

\section{Binary cores, actual caps, and correlated costs}\label{lp:binary}
The assembly charges three- and four-arc cores in many overlapping sources. We first derive the common binary equations, root geometry and cap losses; the separate scalar bounds then differ only in which actual returns are available. This order makes the source of each extra term explicit.

We use the source convention of Section~\ref{lc:source-convention}.
Containment remains a geometric hypothesis: changes of cut use
Lemma~\ref{fd:reanchor}, and whole changes of chart use
Proposition~\ref{fd:lawful-symmetry}.

\subsection{Binary cores and their exterior maximizing jets}
\label{loc:binary-section}
For complete-source cost estimates, each source uses one full lift on its
containing interval; no state covector is continued periodically across the seam.

\begin{lemma}[Automatic properness of two complete sources]
\label{loc:automatic-properness}
Suppose one genuine fundamental interval of an original global minimizer
contains $L$ consecutive complete positive maximal arcs. Then $L\leq B$.
If their source is $012012$ or $01020212$, or a whole lawful cyclic or
reflected-time image, then $L<B$.
\end{lemma}
\begin{proof}
A cut at a switch leaves exactly $B$ complete maximal arcs in one period.
An interior cut splits one physical arc and leaves $B-1$ complete arcs.
Consequently $L\leq B$, and equality requires a switch cut with the source
occupying the entire switch-to-switch period. The actual twist then gives
$j_{i+B}=j_i+1\pmod3$ at the closing switch.

For $012012$ with $B=6$, the resulting cyclic increment word is
$\mathtt{+++++-}$. Its last positive increment, closing negative
increment and first positive increment form $+,-,+$.
For $01020212$ with $B=8$, the closed increment word is
$\mathtt{+--+--+-}$, with the same forbidden triple across the seam.
Each such triple is a cyclic phase of the complete four-arc word $0101$,
excluded by Lemma~\ref{sw:four} without a normality hypothesis. These
four-arc witnesses are proper since $4<6,8$.
A whole cyclic phase leaves the increments unchanged. A full reflection
with reversed physical time reverses their cyclic order and preserves the
$+1$ twist and the forbidden triple. Thus equality is impossible in every
stated image.
\end{proof}

The containing interval and complete maximal arcs are essential here.
This lemma makes no assertion for arbitrary endpoint-constrained
extremals, for fragments of maximal arcs, or for a multi-period interval.
It neither continues a costate across the seam nor assumes its periodicity.

\begin{lemma}[Binary reduction]\label{loc:binary-reduction}
On a connected sequence of modes $0$ and $2$, write
\[
 \begin{gathered}
 \Lambda=\begin{pmatrix}l&b\\m&-l\end{pmatrix},\quad
 k=2ul+b-u^2m,\quad L=l-um,\\
 \chi=vp_v,\quad
 q=v^{-2},\qquad ds=dt/v=q\,du.
 \end{gathered}
\]
For an active mode $j\in\{0,2\}$, the complete equations imply
\begin{align}
 k_s&=0,& L_s&=k,&m_s&=-2L,
 &\chi_s&=2k-\eta(3+u^2)-vH_j,\label{loc:raw-binary}\\
 \det\Lambda&=-L^2-km,&
 \operatorname{tr}(\Lambda X)&=(k-v^2m)/v.\label{loc:binary-casimir}
\end{align}
Let $e=1+u$ on mode $0$ and $e=u-1$ on mode $2$, and set
\[
 f=\left(k-\frac\eta2(3+u^2)\right)q
       +\frac12(eq+e^{-1})\chi.
\]
The inactive Hamiltonian identity is
\begin{equation}\label{loc:inactive-identity}
 H_1=\frac{2v^3}{u^2+v^2-1}(f-m-\eta/2)
       -\frac{1-u^2+v^2}{u^2+v^2-1}H_j.
\end{equation}
In particular, on the actual core $H_j=0$, the function
$g=m+\eta/2-f$ is nonnegative. It vanishes at an exterior mode-$1$
cap. Moreover
\begin{equation}\label{loc:binary-q}
 q_u=-\frac{2q}{1+u}\quad(j=0),\qquad
 q_u=\frac{2q}{1-u}\quad(j=2),\qquad
 H_0-H_2=\frac{2v\chi}{1-u^2}.
\end{equation}
\end{lemma}
\begin{proof}
The two active state fields both have $u_t=v$, so $u_s=v^2$ and
$ds=q\,du$. The Lax equation gives
\[
 l_t=(b+m(u^2+v^2))/v,\quad
 b_t=-2(l(u^2+v^2)+ub)/v,\quad
 m_t=2(um-l)/v.
\]
Substitution into the definitions of $k,L$ proves their three differential
identities and direct expansion proves \eqref{loc:binary-casimir}.
For the state adjoint, differentiate $vp_v$ using
$p_t=-\nabla_W H_j$, with $l,b,m$ held fixed in this gradient.
The terms involving $p_u$ are eliminated using the displayed expression for
$H_j$, leaving the fourth identity of \eqref{loc:raw-binary}.
The same elimination in $H_1$ gives \eqref{loc:inactive-identity}.
Finally, differentiating $v^{-2}$ using $dv/du=v/e$ proves the first two
identities of \eqref{loc:binary-q}; subtraction of the two original
Hamiltonians proves the third. These computations use the active Hamiltonian
only; the inactive Hamiltonian is never set to zero in the core interior.
Since $u^2+v^2-1>0$, full maximization gives $g\geq0$.
\end{proof}

\begin{lemma}[Equality-only capped-binary rigidity]\label{pc:rigidity}
Suppose that two full lifts with the same objective multiplier follow the
same actual trajectory on a connected interval containing a complete
positive binary core, two distinct internal binary switches, and positive
portions of its actual exterior arcs of the third mode. Both core-endpoint
third-mode ties belong to this interval; the far endpoints of the exterior
arcs need not. In a lawful whole chart the
binary modes are $0,2$.  Then the two lifts agree on the interval.
In particular, a full nonzero maximizing lift on a containing genuine
period is normal.  The contexts $12021$ and $102021$ both satisfy these
hypotheses.
\end{lemma}
\begin{proof}
Subtract the full adjoint equations and active zero energies. Their
homogeneous binary reduction, Lemma~\ref{loc:binary-reduction}, gives
\begin{equation}\label{pc:homogeneous-binary}
 K_s=0,\qquad\chi_s=2K,\qquad L_s=K,\qquad m_s=-2L.
\end{equation}
The two distinct actual internal ties force $K=\chi=0$.
The two actual third-mode ties in \eqref{loc:inactive-identity} then force
$m=0$ at both core endpoints; hence $L=m=0$ throughout. Reconstruction
and the active energy give $\Delta\Lambda=\Delta p=0$.
Explicitly, the homogeneous solution is
$\chi=c+2Ks$, $L=L_0+Ks$, $m=m_0-2L_0s-Ks^2$.
The two root rows and two cap rows have determinant
$4S_c(s_1-s_2)\ne0$ in $(K,c,L_0,m_0)$, where
$0<s_1<s_2<S_c$. No cap coefficient or matrix entry is divided by,
and the difference need not maximize.
Uniqueness for the full linear adjoint equation propagates zero through
all actual switches on the connected interval. Applying this argument to
an abnormal lift annihilates it, contrary to nontriviality. Thus normality
precedes joint normalization of every multiplier.
\end{proof}

\begin{corollary}[Reconstruction from a rigid capped core]
\label{pc:rigid-reconstruction}
For the actual state path in Lemma~\ref{pc:rigidity}, its normalized
full reduced lift is reconstructed as follows. Let $s_1<s_2$ be two
internal binary ties, $a,b$ the core endpoints, and $S=s_b-s_a>0$.
In the notation of Lemma~\ref{loc:binary-reduction},
\[
 k=\frac12\left(3+\frac{\int_{s_1}^{s_2}u(s)^2\,ds}{s_2-s_1}\right),
 \qquad \chi(s)=\int_{s_1}^s(2k-3-u(t)^2)\,dt.
\]
Compute $f_a,f_b$ from the displayed definition of $f$ with $\eta=1$,
and put
\[
 m_a=f_a-\tfrac12,\quad m_b=f_b-\tfrac12,\quad
 L_a=(m_a-m_b-kS^2)/(2S).
\]
For $\tau=s-s_a$, the full reconstruction is
\[
 L=L_a+k\tau,\quad m=m_a-2L_a\tau-k\tau^2,\quad
 l=L+um,\quad b=k-2uL-u^2m,
\]
\[
 p_v=\chi/v,\qquad
 p_u=\frac{-\operatorname{tr}((\Lambda-A_*)X)-F_j^vp_v}{v}.
\]
Here $s$ is shear time, and $j$ is the actual adjacent binary mode.
\end{corollary}
\begin{proof}
Integrate \eqref{loc:raw-binary} from the two internal zeros and use
\eqref{loc:inactive-identity} at the two actual caps. The final row is
the active zero energy, with $F_j^u=v>0$. All formulas therefore concern
the unique actual lift, including the state covector.
\end{proof}

\begin{lemma}[Normality and actual roots of a capped four-core]
\label{loc:binary-normality}
For the positive context $102021$, the selected complete lift is normal.
After jointly normalizing it to $\eta=1$, its charged $0202$ core has endpoints
$-\alpha,\beta$ and internal switches $-\xi,z,\zeta$ satisfying
\begin{equation}\label{loc:four-roots}
 \begin{gathered}
 -\alpha<-\xi<z<\zeta<\beta,\quad
 \alpha>\xi>h>0,\quad\beta>\zeta>h,\\
 |z|\leq h,\qquad h^2=2k-3,
 \quad \frac32<k<2.
 \end{gathered}
\end{equation}
If $S$ is the total core shear and $A$ its physical cost, then
\begin{equation}\label{loc:four-cost}
 2A=d+3S+M,\qquad d=\alpha+\beta,\quad
 M=\int u^2\,ds,
 \qquad A>2S.
\end{equation}
\end{lemma}
\begin{proof}
At each of the three internal switches $\chi=0$, while full maximization
requires $\chi\geq0$ on mode $0$ and $\chi\leq0$ on mode $2$.
Lemma~\ref{pc:rigidity} proves normality; normalize the complete lift jointly.

Put $\mu=2k-3$. On each of the two middle arcs the zero-to-zero equation is
\[
 \mu=\frac{\int u^2q\,du}{\int q\,du}>0.
\]
The weak one-sided maximizing jets give
$\xi^2\geq\mu$, $z^2\leq\mu$, $\zeta^2\geq\mu$.
Ordering implies $-\xi<0<\zeta$. Equality $\xi=\sqrt\mu$ would make
$\int_{-\xi}^{z}(\mu-u^2)q\,du$ strictly positive, since this interval
has positive length; equality on the right is excluded in the same way.
The open state domain gives $|u|<1$, proving \eqref{loc:four-roots}.
On the first and last arcs $|u|>h$, so their terminal switching values give
$\chi(-\alpha)>0>\chi(\beta)$.
Finally, $v^2ds=du$ gives \eqref{loc:four-cost}; its last inequality follows
from $\ell\,dt=(3+u^2+v^2)ds/2>2ds$.
\end{proof}

\begin{lemma}[One exterior-leaf calculation]\label{loc:outer-leaf-calculus}
In one jointly normalized full normal binary lift, consider an actual
first mode-$0$ leaf $[-t,-r]$ or last mode-$2$ leaf $[r,t]$, with
$0<h<r<t<1$, $h^2=2k-3$, and an internal tie at $\mp r$.
Write $q_t$ for its outer density, $f_t$ for its outer cap function,
and $D_{\rm out}$ for differentiation in shear directed towards the
outer endpoint. Then
\[
 \frac{(h^2-t^2)q_t}{2}<f_t<0,\qquad
 D_{\rm out}f_t<-t+\frac{t^2-h^2}{2(1+t)}<0.
\]
With $E_t=-f_t$ one also has
\[
 \frac{E_t}{q_t}>e(t,h):=\frac{t^2-h^2}{2}
 -\frac{(t-h)^2(t+2h)}{3(1-t^2)}.
\]
If the actual endpoint additionally satisfies
$q_t<[(1-t)(3+t)]^{-1}$, then
\[
 \frac{E_t}{q_t}>e_2(t,h):=\frac{t^2-h^2}{2}
 -\frac{2(t-h)^2(t+2h)}{3(1-t)(3+t)}.
\]
These lower bounds are positive respectively for
$h<t\leq259/500$ and $h<t\leq349/500$.
\end{lemma}
\begin{proof}
Use the outward coordinate $x=-u$ on the first leaf and $x=u$ on the
last. This is an integral substitution on the same lift. In both cases,
\[\begin{gathered}
 q(x)=\frac{\kappa^{-2}}{(1-x)^2},\quad
 \psi(x)=\int_r^x(y^2-h^2)q(y)\,dy,\\
 f(x)=\frac{(h^2-x^2)q(x)}2+B(x)\psi(x),\quad
 B(x)=\frac{q(x)(1-x)+(1-x)^{-1}}2.
\end{gathered}\]
If $\sigma$ is the leaf shear, then
$0<\psi(t)<(t^2-h^2)\sigma$ and
\[
 \frac{2\sigma B(t)}{q_t}
 =\frac{(t-r)(1+\kappa^{-2})}{1-r}<1.
\]
Indeed $\Omega$ gives $\kappa^{-2}<(1-t)/(1+t)$, and the remaining
positive difference is $(1-t)(1+r)$. Hence the stated bounds on $f_t$.
Differentiating gives
\[
 D_{\rm out}f=-x+\frac{f}{q(1-x)}
                  +\frac{(x^2-h^2)q(1-x)}2,
\]
so $f_t<0$ and $q_t(1-t)<(1+t)^{-1}$ prove the strict jet bound.
Finally, with $I=\int_r^t(x^2-h^2)(1-x)^{-2}\,dx$,
\[
 \frac{E_t}{q_t}=\frac{t^2-h^2}{2}
 -\frac{1-t}{2}\{1+(1-t)^2q_t\}I,\qquad
 0<I<\frac{(t-h)^2(t+2h)}{3(1-t)^2}.
\]
The coefficient of $I$ is smaller than $(1-t)/(1+t)$ by $\Omega$,
and smaller than $2(1-t)/(3+t)$ under the additional circle hypothesis.
These give $e,e_2$. Their positivity follows from
\[
 \frac{e(t,h)}{t^2-h^2}
 =\frac12-\frac{t}{3(1-t^2)}
   +\frac{2h^2}{3(1-t^2)(t+h)}>0
\]
on the first range, and from $(t-h)(t+2h)\leq t(t+h)$ and
$3(1-t)(3+t)-4t\geq139597/250000>0$ on the second.
\end{proof}

\begin{lemma}[The actual matrix apex]\label{loc:four-core-apex}
In the preceding notation let
$E_a=-f(-\alpha)$, $E_b=-f(\beta)$, $E=E_a+E_b$, and
$\mathcal C=\det\Lambda$. Then $E_a,E_b>0$ and
\begin{align}
 |E_a-E_b|&<kS^2,\label{loc:apex-asymmetry}\\
 \mathcal C&=\frac k2+\frac{kE}{2}
   -\frac{(E_a-E_b)^2}{4S^2}-\frac{k^2S^2}{4},\label{loc:cap-casimir}\\
 \mathcal C&>\frac k2-\frac{k^2S^2}{4},\label{loc:apex-floor}\\
 \mathcal C&>\frac34-\frac{9A^2}{64}\qquad(A\leq1).
 \label{loc:structural-casimir}
\end{align}
Whenever the harmonic cap $\mathcal C<C_0<3/4$ applies to this same lift,
\begin{equation}\label{loc:shear-completion}
 S>\frac{\sqrt{2k-4C_0}+E/S}{k}.
\end{equation}
For the value $C_0=\Hcap$, it also follows that $S>81/500$.
\end{lemma}
\begin{proof}
Apply Lemma~\ref{loc:outer-leaf-calculus} with
$(t,r)=(\alpha,\xi)$ and $(\beta,\zeta)$. It gives $E_a,E_b>0$ and
$f_s(0)>0>f_s(S)$ on the same lift. The maximizing gap
$g=m+1/2-f\geq0$ vanishes at the actual exterior ties, so its inward
jets imply $m_s(0)\geq f_s(0)>0>f_s(S)\geq m_s(S)$.
Thus $L(0)<0<L(S)$.
The unique zero of $L$ is an actual interior point, the apex of the concave
quadratic $m$. Its interpolation is
\[
 m(s)=(1-s/S)(-1/2-E_a)+(s/S)(-1/2-E_b)+ks(S-s).
\]
The endpoint derivatives give \eqref{loc:apex-asymmetry}; evaluating
$-L^2-km$ at an endpoint gives \eqref{loc:cap-casimir}.
At the actual apex, $m<-1/2+kS^2/4$, proving the first strict lower bound.
Using $S<A/2$ and the fact that $k/2-k^2A^2/16$ increases on
$3/2\leq k\leq2$, $A\leq1$, proves the second.

If $\mathcal C<C_0<k/2$, equations \eqref{loc:apex-asymmetry}
and \eqref{loc:cap-casimir} first imply $E<kS^2$. Thus $kS-E/S>0$.
The exact identity
\[
 \mathcal C-\left\{\frac k2-\frac14(kS-E/S)^2\right\}
   =\frac{E_aE_b}{S^2}>0
\]
proves \eqref{loc:shear-completion} on its positive square-root branch.
If $S\leq81/500$, the unconditional bound
\eqref{loc:apex-floor}, increasing in $k$ on this range, gives
\[
 \mathcal C>\frac34-\frac9{16}(81/500)^2>C_0.
\]
The last inequality is an exact rational comparison for $C_0=\Hcap$.
\end{proof}

\begin{lemma}[Leaf correlations and endpoint densities]\label{loc:leaf-correlations}
For an actual four-core,
\begin{gather}
 \xi^2-\xi z+z^2>3h^2,\qquad
 \zeta^2+\zeta z+z^2>3h^2,\label{loc:leaf-quadratics}\\
 d>3h,\qquad \alpha\beta>2h^2,\qquad
 (\alpha-\beta)^2>4h^2-d^2/3,\qquad M>h^2S,\label{loc:leaf-consequences}\\
 \frac{q(-\alpha)}S\geq\frac{1-\beta}{d(1+\alpha)},\qquad
 \frac{q(\beta)}S\geq\frac{1-\alpha}{d(1+\beta)}.\label{loc:densities}
\end{gather}
Moreover, for every $0<h<1$, one has
\begin{equation}\label{loc:radius}
 \xi,\zeta<\rho(h):=\frac{8h}{4-3h}.
\end{equation}
\end{lemma}
\begin{proof}
The two exact middle-leaf equations are
\[
 \int_{-\xi}^{z}\frac{u^2-h^2}{(1-u)^2}\,du=0,
 \qquad
 \int_z^\zeta\frac{u^2-h^2}{(1+u)^2}\,du=0.
\]
Subtract the comparison weight $(1+h)^{-2}$. The differences are
\[
 \frac{(u-h)(u+h)^2(2+h-u)}{(1-u)^2(1+h)^2},\qquad
 -\frac{(u-h)^2(u+h)(2+h+u)}{(1+u)^2(1+h)^2},
\]
respectively. Their integrals are strictly negative, since $|z|\leq h$
and both exterior portions have positive length. The unweighted integrals
are therefore positive, which is \eqref{loc:leaf-quadratics}.
Set $T=\sqrt{12h^2-3z^2}$,
$a_0=(T+z)/2$, $b_0=(T-z)/2$.
The positive roots give $\alpha>a_0$, $\beta>b_0$;
$T\geq3|z|$, $a_0b_0=3h^2-z^2\geq2h^2$, and
$a_0^2-a_0b_0+b_0^2=3h^2$.
Writing $\alpha=a_0+x$, $\beta=b_0+y$, with $x,y>0$, gives
\[
 \alpha^2-\alpha\beta+\beta^2
 =3h^2+(2a_0-b_0)x+(2b_0-a_0)y+x^2-xy+y^2>3h^2.
\]
These identities prove the first three consequences. Each middle arc has
mean $u^2=h^2$, and the two positive exterior arcs have $u^2>h^2$,
proving $M>h^2S$.

On both active modes
$-2/(1+u)\leq(\log q)_u\leq2/(1-u)$.
The upper comparison from $-\alpha$ gives
$q(u)\leq q(-\alpha)(1+\alpha)^2/(1-u)^2$.
Integrating over the entire core proves the first inequality in
\eqref{loc:densities}; integration backwards proves the second.

For the radius bound, the actual roots satisfy $h<\xi,\zeta<1$.
After the change of variable $u\mapsto-u$ on the mode-$2$ leaf,
both middle-mean identities have the form
$\int_z^r(h^2-u^2)(1+u)^{-2}\,du=0$, with $|z|\le h$ and $r>h$.
Thus $J_h(r):=\int_{-h}^r(h^2-u^2)(1+u)^{-2}\,du\ge0$,
and $J_h$ strictly decreases for $r>h$.
If $h\ge4/11$, then $\rho(h)\ge1>r$ already.
For $0<h<4/11$, use the exact primitive
$P_h(u)=-u+2\log(1+u)+(1-h^2)/(1+u)$.
Since $J_h(\rho(h))=P_h(\rho(h))-P_h(-h)$ vanishes at $h=0$ and
\[
 \frac{d}{dh}J_h(\rho(h))
 =\frac{54h^4(h-4)(5h-4)}{(h-1)(3h-4)^2(5h+4)^2}<0,
\]
we have $J_h(\rho(h))<0\le J_h(r)$, whence $r<\rho(h)$.
This proves both bounds without reflecting or renormalizing the lift.
\end{proof}

\begin{corollary}[Automatic calibration of a capped four-core]
\label{loc:four-calibration}\label{loc:future-order}
Suppose a genuine fundamental interval of an original global minimizer
contains the complete positive context $102021$, or a lawful whole cyclic
or reflected-time image. Every selected nonzero full maximizing free-time
lift is normal. After joint normalization, its determinant satisfies
$0<\mathcal C<H_*$, before any improved local cost bound is used.
\end{corollary}
\begin{proof}
Lemma~\ref{loc:binary-normality} proves normality and gives
$3/2<k<2$, $S<A/2$. The positive exterior arcs make the four-core
cost a strict part of the period cost, so $A<C_*<2$ by
\eqref{fw:budget}. Its actual apex
\eqref{loc:apex-floor} therefore gives
\[
 \mathcal C>\frac k2-\frac{k^2S^2}{4}
 >\frac k2\left(1-\frac{kA^2}{8}\right)>0.
\]
Apply Corollary~\ref{cal:one-point} with this actual $k>0$.
All caps and both original endpoints remain on the chosen interval;
this argument neither transfers a costate nor calibrates a larger union.
\end{proof}

\begin{lemma}[The signed internal root of a return]\label{loc:return-root}
On an actual $020$ triple in a normal lift, the first internal switch has
$u<0$. On an actual $202$ triple, the second internal switch has $u>0$.
The triples need not share their binary constant $k$ with any other core.
\end{lemma}
\begin{proof}
Normalize the complete lift once. The middle arc has two distinct roots of
$\chi$, so its own $\mu=2k-3$ is the strictly positive weighted mean of
$u^2$. For $020$, weak jets give $r^2\geq\mu$ and $s^2\leq\mu$
for its two ordered switches $r<s$, hence $r<0$. Equality
$r=-\sqrt\mu$ would give a strictly positive middle integral of
$\mu-u^2$, contradicting its zero value. For $202$ the two inequalities
are reversed and imply $s>\sqrt\mu>0$ by the same argument.
\end{proof}

\subsection{A signed cap strip for nonpositive-entry three-cores}
\label{loc:signed-strip-section}
\begin{lemma}[Signed strip]\label{loc:signed-strip}
Consider an actual positive capped context $12021$, charging its complete
$202$ core, whose canonical entrance is $u_0=-a\leq0$.
The selected complete lift is normal. Normalize it jointly and write its
core coordinates as $-a<r<s<b$, its shear as $\sigma$, and
$\mu=2k-3$. Then
\[
 0<\mu<b^2<1,\quad r^2\leq\mu,\quad s>\sqrt\mu>0,
 \qquad k=(3+\mu)/2>0.
\]
Let $f$ be the cap function of Lemma~\ref{loc:binary-reduction} and
$\tau=m+1/2$. The initial value $f_0$ is allowed to have either sign;
the terminal value satisfies $f_3<0$. The actual endpoint jets give
\[
 f_s(0)>0>f_s(\sigma),\qquad
 \tau_s(0)>0>\tau_s(\sigma).
\]
In the convention
\[
 X=f_0/\sigma\quad\hbox{(signed)},\qquad
 Y=-f_3/\sigma>0,\qquad Z=k\sigma>0,\qquad Q=X+Y,
\]
one has
\begin{equation}\label{loc:strip-casimir}
 -Z<Q<Z,\qquad
 \mathcal C=\frac34-\frac14\mathcal L,\qquad
 \mathcal L=(Q+Z)^2-4ZY-\mu=(Q-Z)^2+4ZX-\mu.
\end{equation}
No return-circle hypothesis is used. For the actual physical cost $T$,
\begin{equation}\label{loc:strip-cost}
 2T=w+3\sigma+\int u^2ds,\qquad w=a+b,\qquad T>2\sigma.
\end{equation}
If additionally $w<1$, then
\begin{equation}\label{loc:strip-X-upper}
 X\leq U(\mu):=
 \frac{(\mu-a^2)_+(1+b)}{2(1-a)w}+\frac\mu{1-a^2}.
\end{equation}
\end{lemma}
\begin{proof}
Lemma~\ref{pc:rigidity} proves normality from the two internal roots and actual exterior caps. Normalize the complete lift jointly.

The middle-arc mean and weak jets give $\mu>0$ and $r^2\leq\mu$.
The strict root calculation in Lemma~\ref{loc:return-root} gives
$s>\sqrt\mu>0$, hence $\mu<b^2<1$.

Full maximization gives $g=\tau-f\geq0$ on the core, with zero values
at the actual exterior caps. Therefore
$\tau_s(0)\geq f_s(0)$ and $\tau_s(\sigma)\leq f_s(\sigma)$.
On mode $2$, differentiation of the exact expression for $f$ gives
\[
 f_s=-u-\frac f{q(u-1)}+\frac{(\mu-u^2)q(u-1)}2.
\]
At $u=-a$, put $H=1+a$, $B=(\mu-a^2)/2$.
Substitution of the actual $f_0$ gives
\begin{equation}\label{loc:strip-initial-jet}
 f_s(0)=a+B(H^{-1}-q_0H)
                  -\tfrac12(1+(q_0H^2)^{-1})\chi_0.
\end{equation}
Here $\chi_0\leq0$. The remaining expression is affine in
$0<q_0<(1-a^2)^{-1}$; its two endpoint values are
\[
 a+B/(1+a)=\frac{a^2+2a+\mu}{2(1+a)}>0,\qquad
 \frac{a(1-\mu)}{1-a^2}\geq0.
\]
Thus $f_s(0)>0$, without a sign assumption on $f_0$ or a split at $a=0$.

Apply Lemma~\ref{loc:outer-leaf-calculus} to the final leaf with
$(t,r,h)=(b,s,\sqrt\mu)$. It gives
$(\mu-b^2)q(b)/2<f_3<0$ and
$f_s(\sigma)<-b+(b^2-\mu)/[2(1+b)]<0$.
The strict matrix jets follow from the actual endpoint gaps. Since
$\tau_{ss}=-2k$, its apex is an actual interior point.

The quadratic with endpoint values $\tau(0)=f_0$,
$\tau(\sigma)=f_3$ has derivatives $Z-Q$ and $-Z-Q$ at its ends.
The strict matrix jets give $-Z<Q<Z$.
Substituting these endpoint values in $-L^2-km$ gives both equal
expressions in \eqref{loc:strip-casimir}. Notice that $X$ retains its
signed initial-cap convention. The physical cost identity follows from
$v^2ds=du$, and $T>2\sigma$ from the open state domain.

Finally the lower logarithmic derivative bound
$(\log q)_u\geq-2/(1+u)$ gives
\[
 \sigma\geq q_0(1-a)w/(1+b),\qquad
 q_0/(2\sigma)\leq(1+b)/[2(1-a)w].
\]
The integral to the actual first root gives
$|\chi_0|\leq\mu\sigma$, and $\Omega$ gives
$|A_-(-a)|<(1-a^2)^{-1}$.
Insert these bounds in
$f_0=(\mu-a^2)q_0/2+A_-\chi_0$,
discarding its negative first term only when $\mu<a^2$.
This proves \eqref{loc:strip-X-upper}.
\end{proof}

\begin{lemma}[Permitted signed-loss relaxation]\label{loc:strip-relaxation}
For an actual three-core as above, if $X\leq0$, then
$\mathcal L<Z^2$. If $X\geq0$, at fixed actual $a,b,Z$ one has
\[
 \mathcal L<(U(\mu)+Z)^2-\mu.
\]
The expression on the right is convex on each of its two affine-$U$
pieces in $\mu$. Thus its maximum on $0\leq\mu\leq M$ is attained
among $0$, $M$, and the knot $a^2$ when that knot belongs to the interval.
\end{lemma}
\begin{proof}
If $X\leq0$, then $Y-X\geq|Q|$ and $|Q|<Z$, so
\[
 \mathcal L=Z^2+Q^2-2Z(Y-X)-\mu
 \leq Z^2+Q^2-2Z|Q|-\mu<Z^2.
\]
If $X\geq0$, the derivative in $Y$ of the first loss expression is
$2(X+Y-Z)<0$ throughout $0\leq Y\leq Y_{\rm actual}$, by the actual
strip. Lowering $Y$ to zero gives
$\mathcal L<(X+Z)^2-\mu\leq(U(\mu)+Z)^2-\mu$.
The last expression has nonnegative second derivative on each affine
piece. The comparisons vary a scalar $\mu$ after freezing the actual
$Z$; no equation coupling the varying scalar to the actual shear is
imposed. Any upper bound for $Z$ is substituted only after this reduction.
\end{proof}

\subsection{The generic four-core estimate}
\label{loc:A-section}
Throughout this section put
\[C_0=\Hcap.\]
\begin{lemma}[A common product envelope for binary cap losses]\label{loc:product-loss}
Let $a,b\ge0$, $0<d=a+b\le D<1$, $p=ab$, and $0\le h\le1/2$.
Using the exterior loss $e$ of Lemma~\ref{loc:outer-leaf-calculus}, put
\[
\begin{aligned}
 Q(d,h,p)&=\frac{(1-b)e(a,h)/(1+a)+(1-a)e(b,h)/(1+b)}d,\\
 \Psi(d,h)&=\frac{-64h^3+(12d^2+48d-48)h^2-3d^4-4d^3+12d^2}
 {12d(d+2)^2}.
\end{aligned}
\]
Define
\[
 c(D)=\frac{(2-D)^2}{12(1-D)^2},\quad
 t(D)=\frac{5D^3+13D^2+8D-4}{12(1+D)^3},\quad
 \Gamma(D)=1-c(D)-\max(t(D),0).
\]
On the whole scalar interval $0\le p\le d^2/4$,
\begin{equation}\label{loc:product-envelope}
 dQ_p\le-\Gamma(D),\qquad
 Q\ge\Psi(d,h)+\frac{\Gamma(D)}d(d^2/4-p).
\end{equation}
In particular
\[
 \begin{array}{c|c|c}
 D&\Gamma(D)&\text{strict lower bound}\\\hline
 3/10&299/588&1/2\\
 13/25&1231001/11852352&1/10\\
 9/25&2144981/5030912&9/25
 \end{array}
\]
No root, mean or physical-source condition is imposed on the comparison products.
\end{lemma}
\begin{proof}
Set $A_h=(1-h)^2(1+2h)$ and $B_h=(1+h)^2(1-2h)$; both lie in $[0,1]$.
Direct symmetric reduction, with $l=1+d+p$, gives
\[
 dQ_p=-1+A_h\frac{(2-d)^2}{12(1-d+p)^2}+B_hT_d(l),\quad
 T_d(l)=\frac{d(d+2)^2}{3l^3}-\frac{4+4d-d^2}{12l^2}.
\]
The first positive fraction is at most $c(d)$. Since
\[
 T_d'(l)=\frac{(4+4d-d^2)l-6d(d+2)^2}{6l^4},
\]
$T_d$ has no interior maximum on $l\ge1+d$ and tends to zero there;
hence $T_d(l)\le\max(T_d(1+d),0)=\max(t(d),0)$.
Both $c$ and $t$ increase on $(0,1)$, with derivatives
$(2-d)/(6(1-d)^3)$ and $(d^2+5d+10)/(6(1+d)^4)$.
Thus the derivative bound follows even when $T_d$ is negative.
Finally $Q(d,h,d^2/4)=\Psi(d,h)$; integrate the derivative bound to this
balanced product. The table is exact rational substitution.
\end{proof}

\begin{theorem}\label{loc:A389}
In each actual positive context $102021$, $210102$, or $021210$, the
physical cost of the four complete arcs in positions $1,2,3,4$ is strictly
larger than $389/1000$. The exterior arcs are uncharged.
\end{theorem}
\begin{proof}
It suffices to use the canonical context and the notation of
Section~\ref{loc:binary-section}. Suppose $A\leq B=389/1000$.
Corollary~\ref{loc:four-calibration} gives $0<\mathcal C<C_0$
on this same complete lift. Lemma~\ref{loc:four-core-apex} gives
\begin{equation}\label{loc:A-small-domain}
 S>81/500,\quad d<2B-3(81/500)=73/250<3/10,
 \quad 0<h^2<d^2/9<1/100.
\end{equation}

\paragraph{A correlated lower bound for the cap losses.}
The positive exterior losses and \eqref{loc:densities} give
\begin{equation}\label{loc:A-raw-loss}
 E/S>Q(d,h,\alpha\beta).
\end{equation}
Apply Lemma~\ref{loc:product-loss} with $D=3/10$, and retain the actual
leaf correlation $(\alpha-\beta)^2>4h^2-d^2/3$. Since
$d^2/4-\alpha\beta=(\alpha-\beta)^2/4$, this gives, with $\mu=h^2$
and $R(d)=(12-4d-3d^2)/(3(d+2)^2)$,
\begin{equation}\label{loc:A-loss-envelope}
 \frac ES>R(d)(d/4-\mu/d)
  +\frac{\max(0,4\mu-d^2/3)}{8d}
  +\frac{8\mu(1-2\sqrt\mu/d)}{3(d+2)^2}.
\end{equation}
The comparison keeps the actual common $d,\mu$ throughout.

\paragraph{Minimizing the shear envelope.}
For $0\leq\mu\leq d^2/12$ use
$\sqrt\mu/d\leq1/\sqrt{12}<3/10$ in the last term of
\eqref{loc:A-loss-envelope}; it is at least $16\mu/[15(d+2)^2]$.
For $d^2/12\leq\mu\leq d^2/9$, the function
$\mu(1-2\sqrt\mu/d)$ increases, and the last term is at least
$T(d)=4d^2/[45(d+2)^2]$.
Insert these two continuous lower envelopes in
\eqref{loc:shear-completion}. The derivative with respect to $\mu$,
multiplied by the positive $k^2$, is
\[
 \frac{4C_0-k}{2\sqrt{3+\mu-4C_0}}+H_{\rm lo}(d)
 \quad\hbox{or}\quad
 \frac{4C_0-k}{2\sqrt{3+\mu-4C_0}}+H_{\rm hi}(d),
\]
where
\begin{align*}
 H_{\rm lo}(d)&=\frac{15d^4+20d^3+120d^2+432d-720}{120d(d+2)^2},\\
 H_{\rm hi}(d)&=\frac{105d^4+148d^3+1320d^2+3600d-2160}{720d(d+2)^2}.
\end{align*}
Both functions increase on $0<d<3/10$. After multiplying by positive
factors their derivative numerators are
\[
15d^5+90d^4-40d^3-624d^2+2160d+1440,
\quad105d^5+630d^4-728d^3-4560d^2+6480d+4320;
\]
positivity follows by dropping the positive nonconstant terms and bounding
the negative terms at $3/10$.
The rational bounds
\[
 \begin{gathered}
 147/200<C_0<1471/2000,\quad 243/1000<\sqrt{3-4C_0},\\
 k<301/200,\qquad\sqrt{3+\mu-4C_0}<27/100
 \end{gathered}
\]
give the two derivative signs
\[
 \frac{721}{243}+H_{\rm lo}(3/10)
 =\frac{721}{243}-\frac{128653}{42320}<0,
\quad
 \frac{287}{108}+H_{\rm hi}(9/50)
 =\frac{287}{108}-\frac{67974937}{28514400}>0.
\]
The second sign is used only after $d>9/50$ has been proved below.
Thus, on that range, the minimum is at $\mu=d^2/12$ and
\begin{equation}\label{loc:A-J}
 S>J(d):=\frac{\sqrt{3+d^2/12-4C_0}+R(d)d/6+T(d)}{(3+d^2/12)/2}.
\end{equation}
The radical quotient increases there. The derivative of the remaining
rational term has the sign of
$5d^5-2d^4-240d^3-1120d^2-1296d+1440>0$ on $[0,3/10]$,
again by the endpoint bound on the negative terms. Hence $J$ increases.

\paragraph{The full-profile reserve.}
The original compact set gives $(u^2+v^2-1)/v\geq2/15$.
If $d<D<1$ and $V^2-2V/15+D^2-1<0$, this implies $v>V$
everywhere in the core. Put $y=(u^2+v^2-1)/v$.
On mode $0$,
$y_s=v(1+u)+v^3/(1+u)\geq2v^2$;
on mode $2$ the corresponding derivative is at most $-2v^2$.
On each of the four positive arcs integrate from the endpoint where $y$
is smaller. If its shear length is $s_i$, this gives
\[
 \int (u^2+v^2-1)\,ds>2VS/15+V^3\sum s_i^2
                    \geq2VS/15+V^3S^2/4.
\]
Consequently
\begin{equation}\label{loc:A-profile-cost}
 A>(2+V/15)S+V^3S^2/8.
\end{equation}
Since $q<V^{-2}$ throughout the actual core, direct integration gives
\[
 M=\int_{-\alpha}^{\beta}u^2q(u)\,du
 <\frac{\alpha^3+\beta^3}{3V^2}<\frac{d^3}{3V^2}.
\]
Combining this with
$d=S-M+\int(u^2+v^2-1)ds$ gives
\begin{equation}\label{loc:A-profile-width}
 d>(1+2V/15)S+V^3S^2/4-D^3/(3V^2).
\end{equation}

The following three rational rows finish the proof:
\[
\begin{array}{c|c|c|c|c}
s&D=2B-3s&V&d_0&s'\\\hline
81/500&73/250&5127/5000&4583/25000&4549/25000\\
4549/25000&5803/25000&651/625&10637/50000&9257/50000\\
9257/50000&11129/50000&5219/5000&21727/100000&18563/100000.
\end{array}
\]
For each row, direct rational arithmetic gives
\[
 1+2V/15-D^2-V^2>0,
\quad (1+2V/15)s+V^3s^2/4-D^3/(3V^2)>d_0.
\]
Writing $k_0=(3+d_0^2/12)/2$ and
$T_0=k_0s'-R(d_0)d_0/6-T(d_0)$, one also has
$T_0>0$ and $3+d_0^2/12-4C_0-T_0^2>0$.
These are explicit rational certificates for $J(d_0)>s'$.
The first row establishes $d>4583/25000>9/50$ before the
high-$\mu$ derivative argument is used. Each row therefore proves the next
lower bound for $S$ without circularity. Finally, at
$V=5219/5000$, $s=18563/100000$,
\[
 (2+V/15)s+V^3s^2/8-389/1000
   =\frac{2274532511176128313}{30000000000000000000000}>0.
\]
Equation~\eqref{loc:A-profile-cost} contradicts the closed contrary budget.
Legal cyclic transport gives the other two literal contexts.
\end{proof}

\begin{lemma}[A scalar endpoint principle]\label{loc:endpoint-minimum}
Let $\Psi(d_0,h)=c_0-ah^2-bh^3$, where $c_0,a,b>0$, and put
\[
 g(k)=\frac{\sqrt{2k-4C_0}+\Psi(d_0,\sqrt{2k-3})}{k}.
\]
On any closed interval where $k\geq3/2$ and $2k-4C_0>0$, the minimum
of $g$ is at an endpoint. The same is true after adding a nonnegative
constant to $\Psi$.
\end{lemma}
\begin{proof}
For $h=\sqrt{2k-3}>0$ direct differentiation gives
\[
 k^2g'(k)=\frac{4C_0-k}{\sqrt{2k-4C_0}}-(c_0+3a)-bh(k+3).
\]
The derivative of the right side is
$-k/(2k-4C_0)^{3/2}-b[(k+3)/h+h]<0$.
Thus $g'$ can change sign only from positive to negative, and $g$ has no
interior minimum. Continuity includes the endpoint $k=3/2$.
\end{proof}

\subsection{One-return four-core and circle-qualified returned-core estimates}
\label{ek:section}
Costs are measured in physical time. A context consists of complete positive
maximal arcs on a genuine containing interval. Select one full lift there;
its actual capped core proves normality before joint normalization.
No binary invariant is continued through a third-mode cap without its
lawful full transport.
We use
\[
 C_0=\frac{29025933167887355473205403}{39478417604337498218323684}.
\]
The binary identities, endpoint jets, actual four-core apex and leaf
correlations are those proved in
Lemmas~\ref{loc:binary-reduction}--\ref{loc:leaf-correlations}.

\subsubsection{The one-return four-core cost}
\begin{theorem}[Refined one-return estimate]\label{ek:E}
For the actual right contexts $1020212$, $2101020$, $0212101$, charge
positions $1,2,3,4$; for the actual left contexts $0102021$, $1210102$,
$2021210$, charge positions $2,3,4,5$.  Positions start at zero.
Under the complete-source hypotheses above, the charged cost satisfies
\[
 A>B_E:=\frac{4154032208439}{8759490110000}.
\]
The two exterior caps and the one returned arc are positive and uncharged.
\end{theorem}
\begin{proof}
We give the right-context proof and then explain the left transport.  Put
\[
 s_*={983\over5000},\quad r_*={8\over47},\quad x_*=1+r_*,
 \quad K_*=4-x_*^2,
 \quad d_{\max}={s_*K_*x_*\over x_*+s_*K_*}
 ={314171715\over875949011}.
\]
Thus $2B_E=d_{\max}+3s_*$ exactly.  Suppose $A\leq B_E<1$.
Use the actual variables of the capped four-core: endpoints
$-\alpha,\beta$, total shear $S$, $d=\alpha+\beta$, and
$h^2=2k-3$.  Corollary~\ref{loc:four-calibration} gives $0<\mathcal C<C_0$
on this selected lift.  The exact comparison
$3/4-(9/16)(81/500)^2>C_0$ gives $S>81/500$.
The physical identity $2A=d+3S+M$, with $M=\int u^2\,ds>0$, gives the
fresh initial domain
\[
 0<\alpha,\beta<d<I={2025476111709\over4379745055000}
 <{259\over500}<1,
 \qquad 0<h<d/3,\qquad p=\alpha\beta\leq d^2/4.
\]
This width bound follows directly from the stated contrary budget.

\paragraph{Cap loss and endpoint comparison.}
Let $E_a=-f(-\alpha)$, $E_b=-f(\beta)$ and $E=E_a+E_b$.
Lemma~\ref{loc:outer-leaf-calculus} gives the positive endpoint losses
$e(\alpha,h),e(\beta,h)$ on the fresh domain. Thus
\eqref{loc:densities} gives
\begin{equation}\label{loc:E-Q}
 {E\over S}>Q={ (1-\beta)e(\alpha,h)/(1+\alpha)
 +(1-\alpha)e(\beta,h)/(1+\beta)\over d}.
\end{equation}
Lemma~\ref{loc:product-loss}, with $D=13/25>259/500$, gives
$E/S>\Psi(d,h)$ throughout the fresh domain. Direct differentiation gives
\[
 12d^2(d+2)^3\Psi_d
 =d^2(24-28d-18d^2-3d^3)
 +12h^2(8+12d-6d^2-d^3)+64h^3(2+3d)>0
\]
on the indicated box.  For fixed $d_0$ write
$\Psi(d_0,h)=c_0-a_0h^2-b_0h^3$, where
\[
 c_0={d_0(12-4d_0-3d_0^2)\over12(d_0+2)^2},\quad
 a_0={4-4d_0-d_0^2\over d_0(d_0+2)^2},\quad
 b_0={16\over3d_0(d_0+2)^2}>0.
\]
All three constants are positive at the widths below.  Lemma~\ref{loc:endpoint-minimum} therefore puts the minimum of
$2\{\sqrt{3+h^2-4C_0}+\Psi(d_0,h)\}/(3+h^2)$ at a height endpoint.
This is a scalar comparison, not a second physical core.

\paragraph{One actual return and the fresh bootstrap.}
The final actual $212$ transforms to $020$.  Its own first root is negative
by Lemma~\ref{loc:return-root}, giving
$(1+\beta)^2+v(\beta)^2>4$ at the actual core exit.  This uses the
return's own binary constant.  Backward density integration yields
\[
 S<F(d,1+\beta),\qquad F(d,x)={dx\over(4-x^2)(x-d)}.
\]
For $1\leq x\leq1+d$, $F_x$ has sign $2x^3-dx^2-4d$, whose derivative
is positive.  Thus $F$ has no interior maximum, and
\[
 F(d,1+d)-F(d,1)={2d^2\over3(1-d)(3+d)}>0,
 \quad S<G(d)={d(1+d)\over3-2d-d^2},\quad G'>0.
\]
Also $F_d=x^2/[(4-x^2)(x-d)^2]>0$.  Successively apply
\[
\begin{array}{c|c|c|c}
s&D=2B_E-3s&d_0&s'\\\hline
81/500&2025476111709/4379745055000&73/250&93/500\\
93/500&1710134467749/4379745055000&317/1000&97/500\\
97/500&1605020586429/4379745055000&163/500&49/250\\
49/250&1578742116099/4379745055000&329/1000&983/5000.
\end{array}
\]
These are finite rational comparisons, specified completely as follows.
For each row check $s-G(d_0)>0$ and, at both $h=0,D/3$, put
$T_h=(3+h^2)s'/2-\Psi(d_0,h)$ and check
\[
 T_h>0,\qquad 3+h^2-4C_0-T_h^2>0.
\]
The first comparison gives $d>d_0$; the latter two choose the positive
square-root branch in the actual shear completion
$S>(\sqrt{3+h^2-4C_0}+E/S)/k$.
Monotonicity and the endpoint principle then give $S>s'$ throughout the
height interval.  Direct rational substitution verifies every displayed
comparison (the smallest squared residual exceeds $16/100000$).
We obtain
\[
 S>s_*,\qquad329/1000<d<d_{\max},\qquad0<h<h_{\max}:=d_{\max}/3.
\]
The width inequality is strict, even under the closed contrary budget.

\paragraph{Retain the asymmetry.}
Since $d<d_{\max}<9/25$, the third row of Lemma~\ref{loc:product-loss}
gives $\Gamma(9/25)/d>1$. Therefore
\[
 E/S>\Psi(d,h)+(\beta-d/2)^2.
\]
On this final domain $2-5d_{\max}>0$, so $F_x>0$ for all $x\geq1$ used
below.  Moreover
\[
 {1\over k}>{2\over3+h_{\max}^2}>{33\over50},\qquad
 {2\over3+h_{\max}^2}-{33\over50}
 ={50671561585382537\over14455169411838021175}>0.
\]

\paragraph{The closed low-radius case.}
If $h\leq2/25$, Lemma~\ref{loc:leaf-correlations} gives the actual last
root $r_0<8h/(4-3h)\leq r_*$.  For $\beta\leq r_*$,
\[
 S<F(d,1+\beta)<F(d_{\max},1+\beta)
 \leq F(d_{\max},1+r_*)=s_*,
\]
a contradiction.  The equality in the final scalar comparison does not
lose strictness.  For $\beta>r_*$, retain the actual last mode-$2$ leaf:
\[
 S<{J(\beta-d,\beta,r_0)\over4-(1+\beta)^2},\quad
 J(a,b,r)=(1-b)^2\left[
 {(1+r)^2\over(1-r)^2}\left({1\over1+a}-{1\over1+r}\right)
 +{1\over1-b}-{1\over1-r}\right].
\]
Here $J_r=4(r-a)(1-b)^2/[(1+a)(1-r)^3]>0$ and
$J_d=(1-b)^2(1+r)^2/[(1-r)^2(1+b-d)^2]>0$ after setting $a=b-d$.
Enlarge $r_0$ to $r_*$ and $d$ to $d_{\max}$, preserving all denominator
signs.  Multiplication of the resulting $s_*-S$ bound by the positive
$(b+3)(1-r_*)^2(1+b-d_{\max})$ gives
\[
 q(b)=-{20594857\over11045000}b^2
 +{18296134235204047\over9674856826495000}b
 -{323946560474377\over1209357103311875}.
\]
It satisfies $q(r_*)=0$, $q(d_{\max})=412481343731491/2418714206623750>0$,
and
\[
 q(b)-q(d_{\max}){b-r_*\over d_{\max}-r_*}
 =-{20594857\over11045000}(b-r_*)(b-d_{\max}).
\]
Thus $q(b)>0$ on $r_*<b\leq d_{\max}$, again a contradiction.

\paragraph{Two closed high-radius intervals.}
Use the two rows
\[
\begin{array}{c|c|c}
h_0&h_1&s\\\hline
2/25&11/100&79/400\\
11/100&h_{\max}&1973/10000.
\end{array}
\]
The same endpoint tests with $d_0=329/1000$, now at $h_0,h_1$, give
$S>s$ on each entire closed interval.  Put $D=2B_E-3s$.
Since $\alpha>h\geq h_0$, we have $\beta<D-h_0$.
For $\beta\leq D/2$, the exact positive comparison
$s-F(D,1+D/2)>0$ is contradictory.  Otherwise set
\[
 \sigma(b)=s+{33\over50}(b-D/2)^2,\quad
 \delta(b)=2B_E-3\sigma(b),\quad L=D/2,
 \quad U=D-h_0.
\]
The retained asymmetry and $\beta-d/2>\beta-D/2>0$ give
$S>\sigma(\beta)$ and $d<\delta(\beta)$.
On the full closed interval $[L,U]$, the three denominator reserves
\[
 \delta(b)\geq\delta(U)>0,\quad
 4-(1+b)^2\geq4-(1+U)^2>0,\quad
 1+b-\delta(b)\geq1+L-D>0
\]
are exact rational inequalities.  Therefore
$S<F(\delta(\beta),1+\beta)$ contradicts positivity of
\[
 P(b)=\sigma(b)[4-(1+b)^2][1+b-\delta(b)]-\delta(b)(1+b).
\]
Set $x=b-L$; in both rows $0\le x\le U-L<1/10$.
Writing $P(L+x)=\sum_{i=0}^6a_ix^i$, exact rational expansion gives
\[
 (a_1,\ldots,a_6)>
 (-1117/5000,207/50,13/10,9/10,-15/4,-4/3)
\]
componentwise, and $a_0>c_i$, where
$c_1=479/100000$, $c_2=41/12500$.
Dropping positive cubic and quartic terms, and using
$207/50-15/4000-1/7500>413/100$, gives
\[
 P(L+x)>c_i-\frac{1117}{5000}x+\frac{413}{100}x^2
 \ge c_i-\frac{1247689}{413000000}
 =\begin{cases}
 730581/413000000,&i=1,\\
 106951/413000000,&i=2.
 \end{cases}
\]
These exact positive bounds cover both whole closed intervals, including
all endpoints. All height cases are now exhausted.

For a left context, its actual initial return gives its entrance circle.
In the same actual integrals change variable $u\mapsto-u$ and reverse
limits.  This exchanges $\alpha,\beta$ and the two cap losses, placing
that one returned endpoint in the role just proved.  It neither assumes a
second circle nor identifies different binary invariants.  Lawful whole
cyclic transport gives the other four contexts and their stated charges.
\end{proof}

\begin{corollary}[A weaker comparison for overlapping cores]\label{loc:E473}
Every one-return charge in Theorem~\ref{ek:E} exceeds $473/1000$.
\end{corollary}
\begin{proof}
The exact constant in Theorem~\ref{ek:E} exceeds $473/1000$.
This corollary is used only in a few overlapping-core comparisons below.
\end{proof}

\begin{lemma}[Two actual return circles and cap losses]
\label{loc:double-return-loss}
Retain the complete source $01020212$ on one full normal lift, with
four-core endpoints $-\alpha,\beta$, mean $h^2$, density $q=v^{-2}$,
shear $S$ and losses $E_a,E_b$ of Lemma~\ref{loc:four-core-apex}.
Both actual returns give
\begin{equation}\label{loc:double-circle}
 q(-\alpha)<\frac1{(1-\alpha)(3+\alpha)},\qquad
 q(\beta)<\frac1{(1-\beta)(3+\beta)}.
\end{equation}
Moreover,
\begin{equation}\label{loc:double-return-loss-formula}
 \frac{E_a}{q(-\alpha)}>e_2(\alpha,h),\qquad
 \frac{E_b}{q(\beta)}>e_2(\beta,h),\qquad
 e_2(t,h)=\frac{t^2-h^2}{2}
       -\frac{2(t-h)^2(t+2h)}{3(1-t)(3+t)}.
\end{equation}
For $h<t\leq349/500$, these lower losses are positive. In that domain,
with $d=\alpha+\beta$,
\[
 \frac{E_a}{S}>\frac{(1-\beta)e_2(\alpha,h)}{d(1+\alpha)},\qquad
 \frac{E_b}{S}>\frac{(1-\alpha)e_2(\beta,h)}{d(1+\beta)}.
\]
No numerical four-core cost floor is used.
\end{lemma}
\begin{proof}
Apply Lemma~\ref{loc:return-root} to the actual initial $010$, in its
whole cyclic chart $202$, and to the final $212$, in its whole cyclic
chart $020$. Their strict roots give respectively
$v(-\alpha)^2>(1-\alpha)(3+\alpha)$ and
$v(\beta)^2>(1-\beta)(3+\beta)$, proving
\eqref{loc:double-circle}. The binary constants in the return charts
need not agree with the four-core constant; all restrictions belong to
the same jointly normalized full lift.

Apply Lemma~\ref{loc:outer-leaf-calculus} to both actual outer leaves,
using their separately proved circles. It gives
\eqref{loc:double-return-loss-formula}, including positivity for
$h<t\leq349/500$, without identifying their return-chart constants.
Only after establishing positivity multiply the two loss bounds by
the corresponding endpoint-density inequalities~\eqref{loc:densities}.
\end{proof}

\begin{theorem}[A coarse two-return comparison]\label{loc:D522}
In the three actual positive contexts $01020212$, $12101020$, $20212101$,
the physical cost of positions $2,3,4,5$ is strictly larger than
$261/500$. The two caps and both returns are uncharged.
\end{theorem}
\begin{proof}
The source retains all eight complete maximal arcs in one genuine
fundamental interval of the original minimizer, as required throughout
this section. Lemma~\ref{loc:automatic-properness} makes it proper.
The independent two-return theorem in Section~\ref{pf:M101} therefore
gives the stronger bound $A>74/125>261/500$ on these same four charges.
Its proof uses the actual circle and cap-loss identities of
Lemma~\ref{loc:double-return-loss}, not the present numerical conclusion.
Whole lawful cyclic transport gives the other two contexts. Thus the
original complete-source scope is retained without an additional
properness assumption.
\end{proof}

\subsection{Exact polynomial certificates for the four-core estimates}
\label{loc:AED-certificates}
The common product envelope replaces separate polynomial symmetrizations.
Its rational table and the bootstrap tables in the proofs are complete
finite certificates, not samples of physical trajectories.

The bootstrap tables in the proofs constitute equally explicit rational
certificates. A radical comparison $\sqrt{R}>T$ is checked only after
$T>0$, and is then exactly equivalent to $R-T^2>0$.
To make positivity readily
checkable without printing large denominators, the next table gives
conservative positive rational lower bounds for each group of residuals.
Every bound follows by substitution in these formulas and clearing positive
integer denominators.
\[
\begin{array}{l|c}
\text{group of residuals}&\text{common lower bound}\\\hline
\text{Four-core height residuals}&10^{-5}\\
\text{Four-core width residuals}&10^{-6}\\
\text{Four-core squared radicals}&10^{-7}\\
\end{array}\]
All unsquared branch thresholds in these groups exceed $1/10$.

\subsection{The three-arc estimate}\label{sec:f245}\label{f245:section}
This section proves the local estimate that will be used for a charged
three-arc core. The only physical contexts covered here are
\[
02010,\qquad10121,\qquad21202.
\]
The three middle arcs are charged; the first and last arcs are strictly
positive, actual, uncharged caps. All quantities below come from one
complete normal lift, with its multiplier jointly normalized once. In
particular, the two leaf invariants are not independently chosen costates.
The admissible reanchoring requirement of Lemma~\ref{fw:reanchor} is retained.
Write
\begin{equation}\label{f245:constants}
 U=\frac{49}{200},\qquad
 H_* =\frac{29025933167887355473205403}{39478417604337498218323684}.
\end{equation}
The symbol $\mathscr C$ below denotes the determinant of the matrix covector,
not the cost of a trajectory.

\begin{theorem}[Qualified three-arc bound]\label{f245:theorem}
For an actual context of one of the three types above, to which the full
normal PMP and the global harmonic calibration apply, the cost $F$ of
its three middle arcs satisfies $F>49/200$.
\end{theorem}

\subsubsection{Reconstruction from a single full lift}\label{f245:reconstruction}
It suffices first to treat $02010$. Denote its charged arcs by $A,c,B$, of
modes $2,0,1$, with costs $\alpha,A_c,\beta$. Put
\[
\Lambda=\begin{pmatrix}\lambda_1&\lambda_2\\m&-\lambda_1\end{pmatrix},
\quad K=2u\lambda_1+\lambda_2-u^2m,
\quad L=\lambda_1-um,\quad\chi=vp_v,\quad ds=dt/v.
\]
Differentiation of the full equations, using only the active Hamiltonian
$H_j=0$, gives on a $0/2$ edge
\begin{equation}\label{f245:edge}
 K_s=0,\quad L_s=K,\quad m_s=-2L,\quad
 \chi_s=2K-3-u^2,\qquad \mathscr C=-L^2-Km.
\end{equation}
Indeed $H_0-H_2=2v^2p_v/(1-u^2)$; maximization therefore gives
$\chi\ge0$ in mode $0$ and $\chi\le0$ in mode $2$.
At an actual $2\to0$ switch, $\chi'$ is strictly positive. If it vanished
at $u\ne0$, its common second derivative would be $-2u$, inconsistent
with the sign change. If it vanished at $u=0$, then $K=3/2$ and
$\chi'=-u^2/v<0$ on the neighboring nontrivial arcs, again impossible.
This argument allows tied switches and does not assume transversality.
The continuous maximum-zero relation at a switch follows from its
a.e. version by continuity of all three Hamiltonians.

On the middle mode-$0$ arc set $\xi=1+u$, $v=k\xi$, with endpoints
$\xi=a,aE$, where $E>1$. Its lawful cyclic image
$\mathcal R(w)=(w-3)/(w+1)$ is
\[
 U_s=1-\frac4{(1+k^2)\xi},\qquad
 V_s=\frac{4k}{(1+k^2)\xi}.
\]
(The subscript on $U_s,V_s$ here distinguishes state coordinates.)
At the exit put $b=4/((1+k^2)aE)$ and define
\begin{equation}\label{f245:variables}
 \begin{gathered}
 r=1-a,\quad s=1-b,\quad P=ab,\quad Q=a+b,\quad R=rs,
 \\
 D=4-EP,\quad z=\log E,\quad t=2z/D.
 \end{gathered}
\end{equation}
Then $k^2=D/(Eab)$, $a,b>0$, $aE<2$, and $bE<2$.
The actual positive leaf widths satisfy $0<d_A<a$, $0<d_B<b$.
Let
\begin{equation}\label{f245:mean-definition}
 e(w,d)=
 \frac{\int_0^d(2w\delta+\delta^2)(1+w+\delta)^{-2}\,d\delta}
 {\int_0^d(1+w+\delta)^{-2}\,d\delta},\qquad
 e_A=e(r,d_A),\quad e_B=e(s,d_B).
\end{equation}
The two actual roots on each leaf, together with \eqref{f245:edge}, imply
\[
 K_A=(3+r^2+e_A)/2,\qquad K_B=(3+s^2+e_B)/2,\qquad e_A,e_B>0.
\]
The strict inequalities are the central-facing switch argument just given.

Here is an explicit common-matrix reconstruction. Direct conjugation by
$R_* =\tfrac12\left(\begin{smallmatrix}1&-3\\1&1\end{smallmatrix}\right)$
gives, on the central arc,
\[
 K_B=\frac{4(K_A+2Lk^2\xi-k^4m\xi^2)}{(1+k^2)^2\xi^2}.
\]
Since $J_m=L+K_A/(k^2\xi)$ is constant, solving this equation gives
\begin{align*}
 L&=J_m-K_A/(k^2\xi),\\
 m&=2J_m/(k^2\xi)-K_A/(k^4\xi^2)-K_B(1+k^2)^2/(4k^4),\\
 \mathscr C&=K_AK_B(1+k^2)^2/(4k^4)-J_m^2.
\end{align*}
Integrating $\chi_s$ from the entry root gives
\[
 \chi(aE)=k^{-2}\bigl((2K_A-4)(E-1)/(aE)-a(E-1)+2z\bigr).
\]
At the exit $H_0=H_1=0$ gives
\[
 m+\tfrac12=(2K_A-3-u^2)/(2k^2\xi^2)
                +(1+k^2)\chi/(2k^2\xi).
\]
Substitution determines, with $x=e_A/(aD)$ and $y=e_B/(bD)$,
\begin{equation}\label{f245:matrix}
 J_m=x+y+t+\frac{4Q-4-EP}{2D},\qquad
 \mathscr C=\frac{4K_AK_B}{D^2}-J_m^2.
\end{equation}
Consequently the full determinant expansion is
\begin{align}\label{f245:determinant}
 \mathscr C={}&\frac34+sx+ry+\frac{R(4-P)t}{D}-(x+y+t)^2
 +Pxy\notag\\
 &+\frac{P(E-1)(2-b)x}{D}+\frac{P(E-1)(2-a)y}{D}
 +\frac{PG}{D^2},\\
 G={}&[6-2R-P(E-1)](E-1-z)+(2-P)z(E-1).\notag
\end{align}
This is an identity; in particular none of its remainder terms has yet
been discarded.

The actual leaf endpoint heights are $V_A^2=aD/(Eb)$ and
$V_B^2=bD/(Ea)$. Integration of the physical cost yields
\begin{equation}\label{f245:costs}
 \alpha=\frac{d_A}{2}+\frac{K_Ad_A(1+r)}{V_A^2(1+r+d_A)},\quad
 \beta=\frac{d_B}{2}+\frac{K_Bd_B(1+s)}{V_B^2(1+s+d_B)},\quad
 A_c=\frac{2Q(E-1)-EPz}{D}.
\end{equation}
Thus $d_A<2\alpha$, $d_B<2\beta$, and direct integration of
$\ell\,dt>du/2$ on the central arc and its cyclic image gives
$a(E-1)<2A_c$, $b(E-1)<2A_c$.
Put $J_A=2a+Eb(2-a)$ and $J_B=2b+Ea(2-b)$. The useful strict leaf bound
has the following short proof. For $q=1+r$,
\begin{equation}\label{f245:leaf-gap}
 \alpha-\frac{d_AJ_A}{aD(1+d_A/q)}
 =\frac{d_A^2}{2(q+d_A)}+
       \frac{e_Ad_Aq}{2V_A^2(q+d_A)}>0.
\end{equation}
The corresponding identity holds on the right leaf. This proof uses the
actual cost and positive mean and requires no extra endpoint geometry.
For $w\ge0$, the decreasing-weight covariance inequality gives
$e(w,d)\le wd+d^2/3$. Hence, defining
\[
 N(w,d)=(\max(w,0)+d/3)(1+d/(1+w)),
\]
we obtain $x\le\alpha N(r,d_A)/J_A$ and
$y\le\beta N(s,d_B)/J_B$; the negative-coordinate case follows also from
the signed mean proof below.

Set $\Phi(E)=E+1-2E\log E/(E-1)$, with $\Phi(1)=0$.
The central maximizing signs in the two actual charts give
\begin{equation}\label{f245:central-signs}
 e_A\ge a[\Phi(E)-r(E-1)],\qquad
 e_B\ge b[\Phi(E)-s(E-1)].
\end{equation}
The common integral identity
\[
 \Phi(E)=2(E-1)^2\int_0^1\frac{v(1-v)}{1+(E-1)v}\,dv
       \ge\frac{(E-1)^2}{3E}
\]
will also be used. Both signs are derived before exchanging scalar names;
no odd permutation of the physical modes is being invoked.

\subsubsection{The necessary third-control root condition}\label{f245:root}
At an actual $0/2$ tie $p_v=0$. The remaining Hamiltonian difference is
$H_1-H_0=-2v^3p_u/(u^2+v^2-1)$, so full maximization implies $p_u\ge0$.
The active Hamiltonian then gives
\[
 p_u=m+\tfrac12-(2K-3-u^2)/(2v^2).
\]
At the two central-facing roots the same matrix has
\begin{align*}
 L_A&=J_m-K_AEb/D,&m_A&=2J_mEb/D-K_AE^2b^2/D^2-4K_B/D^2,\\
 L_B&=-J_m+K_BEa/D,&m_B&=2J_mEa/D-K_BE^2a^2/D^2-4K_A/D^2.
\end{align*}
The second line is obtained by the actual cyclic conjugation above.
On the right leaf put $q_B=2-b$, $d=d_B$ and
\[
 \tau=\frac{dq_B}{V_B^2(q_B+d)},\qquad
 T_B=\frac{Ea}{D}+\tau=\frac{Ea(bq_B+2d)}{bD(q_B+d)}.
\]
Forward propagation gives $m_{\rm out}=m_B-2L_B\tau-K_B\tau^2$,
and hence the necessary inequality
\begin{equation}\label{f245:root-test}
 p_{\rm out}=\frac12-\frac{4K_A}{D^2}+2J_mT_B-K_BT_B^2
 -\frac{(e_B-2sd-d^2)q_B^2}{2V_B^2(q_B+d)^2}\ge0.
\end{equation}
Backward propagation on the left gives the analogous formula with
$m(t_0)=m_A+2L_A\tau_A-K_A\tau_A^2$ and
$\tau_A=d_A(2-a)/(V_A^2(2-a+d_A))$.
Thus both exterior tests are valid in their own charts before scalar sorting.
Only the right test is needed by the finite residual proof.

\subsubsection{An analytic treatment of negative coordinates}\label{f245:negative-section}
\begin{lemma}\label{f245:negative}
If at least one of $r,s$ is negative and
$0<F\le A\le U_-=129/500$, then $\mathscr C>H_*$.
\end{lemma}
\begin{proof}
The proof uses \eqref{f245:determinant}--\eqref{f245:central-signs}, not a
separate negative-coordinate theorem. All budgets in this proof refer to
the same three actual arcs.

\paragraph{Signed means and the retained determinant remainder.}
If $r=-\kappa<0$, put $q=1-\kappa>0$ and
$M(v)=v^2/3-\kappa v$. Since
\[
 (v/d)M(d)-M(v)=v(d-v)/3\ge0,
\]
we have $M(v)\le\max(0,M(d))$ on $[0,d]$. Integration by parts
with $\omega(v)=(q+v)^{-2}$ writes the mean in
\eqref{f245:mean-definition} as
\[
 e_A=\frac{d\omega(d)M(d)+\int_0^d[-v\omega'(v)]M(v)\,dv}
 {d\omega(d)+\int_0^d[-v\omega'(v)]\,dv}.
\]
Thus $e_A>0$ implies $d_A>3\kappa$ and
$e_A\le d_A^2/3-\kappa d_A$. No monotonicity of the signed quadratic
has been presumed. The central sign, the integral for $\Phi$, and
$(E-1)/\sqrt E\ge z$ give
\begin{equation}\label{f245:signed-mean}
 e_A\ge a(z^2/3+\kappa z),\qquad
 e_A\ge z^2/3+\kappa z,\qquad d_A\ge z+3\kappa.
\end{equation}
The last implication follows by factoring
$d_A^2/3-\kappa d_A-z^2/3-\kappa z
=(d_A+z)(d_A-z-3\kappa)/3$.
In particular $\kappa<2\alpha/3$ and
$a<a_m:=1+2U_-/3=293/250$. The analogous facts hold on the other side.

More generally the positive mean implies $-U_-<r,s<1$:
if a coordinate is negative this follows from the stronger depth bound
just proved. Therefore $P<(1+U_-)^2<2$ and
$P(E-1)<2U_-(1+U_-)$. Every coefficient used from $G$ in
\eqref{f245:determinant} is positive. Denoting its last four positive
terms by $\mathrm{Rem}$, the inequalities $E-1-z\ge z^2/2$ and
$E-1\ge z$ give
\begin{equation}\label{f245:retained-rem}
 \mathrm{Rem}\ge \frac P4
       [5-rs-P-P(E-1)/2]t^2.
\end{equation}
Also the central cost has the positive-tail estimate
\begin{equation}\label{f245:central-cubic}
 A_c\ge\frac{2Q-EP}{2}t+\frac{QD}{4}t^2+\frac{QD^2}{24}t^3.
\end{equation}
These follow by expanding $e^z-1$, with $z=Dt/2$.

\paragraph{Mixed signs: a uniform right-leaf coefficient.}
After the two chart derivations we may name the negative coordinate
$r=-\kappa<0$ and take $s=1-b\ge0$. The actual exit has
$\xi=aE\in(1,2)$ and $v_{\rm state}^2\le243/4$ by the original
compact set. Hence
\[
 b=\frac{4\xi}{\xi^2+v_{\rm state}^2}
 \ge\frac{16\xi}{4\xi^2+243}>b_0:=\frac{16}{247}.
\]
The rational function is increasing on $[1,2]$. Moreover
\[
 \frac{bD}{J_B}=\frac{b(4-\xi b)}{2b+\xi(2-b)}
 \le\frac{b(4-b)}{2+b},
\]
since its $\xi$ derivative is negative. The strict bound
\eqref{f245:leaf-gap}, with $\beta\le U_-$, therefore implies
\[
 d_B<h(b):=\frac{U_-b(4-b)(2-b)}{4-b^2-U_-b(4-b)}.
\]
The denominator decreases on $[0,1]$ to $3(1-U_-)>0$.
Since $N(1-b,d)$ increases in $d$ and $J_B>2+b$,
\[
 y<\beta\,\mathcal F(b),\qquad
 \mathcal F(b)=\frac{[1-b+h(b)/3][1+h(b)/(2-b)]}{2+b}.
\]
For $k_0=233/500$, exact reduction gives
\[
 k_0-\mathcal F(b)=
 \frac{P_4(b)}{500(371b^2+516b-2000)^2},
\]
where
\[
 P_4(b)=135570353b^4-146040824b^3-770234352b^2
                 +1105088000b-68000000.
\]
On $0\le b\le1$, $b^2\le b$ gives $P_4''(b)\le-789869412<0$.
Its endpoint values on the required interval are
\[
 P_4(16/247)=\frac{1173789821728512}{3722098081}>0,
 \qquad P_4(1)=256383177>0.
\]
Concavity therefore proves $y<k_0\beta$ on the entire closed interval.

\paragraph{Mixed signs: a triangular comparison region.}
Put $u_0=1-b/2\in[1/2,1)$. The actual chart bound gives $D>4u_0$.
Using \eqref{f245:signed-mean} and $d_A<2\alpha$ yields
\[
 \alpha>u_0t,\qquad
 0<\kappa<K_0:=2(\alpha-u_0t)/3,
\]
\[
 L_\kappa:=u_0t^2/3+\kappa t/2<x<
 X_\kappa:=\alpha^2/(3u_0)-\kappa\alpha/(2u_0).
\]
The bounds fit exactly:
\[
 X_\kappa-L_\kappa
  =\frac{(\alpha+u_0t)(K_0-\kappa)}{2u_0},\qquad
 X_{K_0}=L_{K_0}=\alpha t/3.
\]
All compared endpoints are nonnegative. Furthermore
$P(E-1)<2U_-b$, $P<a_mb$, and
$D>4-b(1+8U_-/3)\ge(9-8U_-)/3$ imply
\[
 \frac{4-P}{D}<R_0:=1+\frac{6U_-}{9-8U_-}=\frac{707}{578}.
\]
Equation~\eqref{f245:retained-rem}, with $P\ge b$ and $rs\le0$, gives
\[
 \mathrm{Rem}\ge\gamma b t^2,\qquad
 \gamma=(4-5U_-/3)/4=357/400.
\]
Keeping the beneficial term $sx$ in the determinant gives
\begin{equation}\label{f245:mixed-loss}
 \frac34-\mathscr C\le
 (t+x+k_0\beta)^2+\kappa(k_0\beta+R_0t)-sx-\gamma bt^2.
\end{equation}
Replacement of $y$ is valid because the preceding loss has derivative
$2(t+x+y)+\kappa>0$ in $y$. At fixed $\kappa$ the displayed expression
is convex in $x$; on either affine side $x=L_\kappa,X_\kappa$ it is
convex in $\kappa$. Applying the convex chord inequality twice reduces
it to the three points
\begin{equation}\label{f245:three-corners}
 (0,u_0t^2/3),\qquad(0,\alpha^2/(3u_0)),\qquad
 (K_0,\alpha t/3).
\end{equation}
This is a scalar comparison, not an assertion that these points are
physical trajectories.

\paragraph{Exact allocation of the two leaf costs.}
Equation~\eqref{f245:central-cubic}, $Q\ge1+b$, $D>4u_0$, and
$(2Q-EP)/2\ge1+b/2-bA_c$ give $A_c\ge g$, where
\[
 g=\frac{(1+b/2)t+(1+b)u_0t^2+(1+b)(2-b)^2t^3/6}{1+bt},
 \qquad H=g+u_0t.
\]
The following two identities have nonnegative right sides:
\begin{align*}
 (1+bt)(g-t)&=bt/2+(1-b)(1+b/2)t^2+(1+b)(2-b)^2t^3/6,\\
 H-(2t-t^2/2)&=\frac{t^2}{1+bt}
 \left[(1-b)(b+3/2)+\{b/2+(1+b)(2-b)^2/6\}t\right].
\end{align*}
Since $H<A$, first $t<2A/3$ and then
\begin{equation}\label{f245:mixed-t}
 t<\frac{A}{2-A/3}\le T_0:=43/319.
\end{equation}
Now fix $B=\alpha+\beta$ and decrease $\alpha$ to $u_0t$, keeping
$\beta=B-\alpha\ge0$. Each of the three losses in
\eqref{f245:three-corners} decreases strictly as $\alpha$ increases.
For the first two corners the derivatives are respectively
\[
 -2k_0(t+u_0t^2/3+k_0\beta),\quad
 2W(2\alpha/(3u_0)-k_0)-2s\alpha/(3u_0),
\]
where $W=t+\alpha^2/(3u_0)+k_0\beta>0$. Both are negative since
$2\alpha/(3u_0)\le4U_-/3<k_0$.
For the third put $\nu=t/3$, $\kappa=2(\alpha-u_0t)/3$ and
$\beta=B-\alpha$. Its derivative is
\[
 2t(\nu-k_0+R_0/3)+2\alpha\nu(\nu-k_0)
 +2k_0\beta(\nu-k_0+1/3)-k_0\kappa-s\nu<0.
\]
Every term is nonpositive and the first is negative: use
$\alpha,\beta,\kappa,s\ge0$, $\nu<T_0/3<k_0-1/3$ and
$k_0>(T_0+R_0)/3$.
At $\alpha=u_0t$ the three corners coincide. Increasing $B$ to
$A-g$ is then legitimate, since $A-H>0$ and the root of the square
is positive. We obtain
\[
 \frac34-\mathscr C< W_0^2-su_0t^2/3-\gamma bt^2,
 \qquad W_0=t+u_0t^2/3+k_0(A-H)>0.
\]
Set $q_0=1-2k_0=17/250$ and $c_0=(1+3k_0)/6=1199/3000$.
Direct subtraction shows $W_0\le k_0A+q_0t+c_0t^2$:
after division of this gap by $t^2$ and multiplication by $1+bt$,
the result is
\[
 \frac{(1-b)(6k_0b+9k_0-1)}6+
 t\left[\frac{b(3k_0-1+b)}6+
             \frac{k_0(1+b)(2-b)^2}6\right]\ge0.
\]
Also $su_0/3+\gamma b-1/3=b(\gamma-1/2+b/6)\ge0$.
The square roots being positive, we may square this bound. Define
\[
 M_0=\frac13-q_0^2-2k_0U_-c_0-2q_0c_0T_0-c_0^2T_0^2
     =\frac{105198654743}{473062500000}>0.
\]
Expansion and completion of the square imply
\[
 \frac34-\mathscr C<(k_0A)^2+2k_0Aq_0t-M_0t^2
 \le(k_0A)^2(1+q_0^2/M_0)
 \le L_m:=\frac{97015095513566129007}{6574915921437500000000}.
\]
The decisive rational margin is
\[
 (3/4-H_*)-L_m=
 \frac{3527766575724023280412656655948597000737}
 {383975261035380330465773933925280437500000000}>0.
\]
This proves the mixed case, including $s=0$.

\paragraph{Both coordinates negative.}
Write $r=-\kappa$, $s=-\lambda$, so $1<a,b<a_m$. In
\eqref{f245:retained-rem}, use $P\ge1$, $P\le a_m^2$,
$rs\le(a_m-1)^2$ and $P(E-1)/2\le U_-a_m$. It follows that
\[
 \mathrm{Rem}\ge\frac{411807}{500000}t^2>\frac23t^2.
\]
For negative $w$, $N(w,d)$ increases in $d$ and decreases in $w$.
At $w=-U_-$ and $d=2U_-$ it equals $2n_0$, where
\[
 n_0=\frac{U_-(1+U_-)}{3(1-U_-)}=\frac{27047}{185500}.
\]
Since $J_A,J_B>2$, we have $x\le n_0\alpha$, $y\le n_0\beta$.
The function $a+b-ab/2-U_-\min(a,b)$ increases separately on either
ordered half of $[1,a_m]^2$: its derivatives are at least
$1/2-4U_-/3$ and $1/2-U_-/3$, both positive. Thus
\[
 A_c\ge c_1t,\qquad c_1=3/2-U_-=621/500.
\]
Both compulsory debits $d_A,d_B\ge z+3\kappa,z+3\lambda$ give
$B:=\alpha+\beta>Dt/2>\ell_1t$, with
$\ell_1=1-2U_-/3=207/250$. Hence
\[
 t<\frac{A}{c_1+\ell_1},\qquad B\le A-c_1t.
\]
Discard the positive $rs$ term, retain $2t^2/3$ of the remainder, and
use $\kappa<2\alpha/3$, $\lambda<2\beta/3$. This gives
\[
 \frac34-\mathscr C\le(t+n_0B)^2+\frac{n_0}{3}B^2-\frac23t^2.
\]
It increases in $B\ge0$. After the residual substitution it is
\[
 Q_A(t)=q_2t^2+2n_0A h_1t+(n_0^2+n_0/3)A^2,
\]
where
\begin{align*}
 q_2&=(1-n_0c_1)^2+n_0c_1^2/3-2/3
 =\frac{2036665287275407}{25807687500000000}>0,\\
 h_1&=1-c_1(n_0+1/3)=\frac{37555313}{92750000}>0.
\end{align*}
Thus $Q_A$ increases for $t\ge0$. At $t=A/(c_1+\ell_1)$ the
residual equals $2A/5>0$. Homogeneity yields
\[
 \frac34-\mathscr C<Q_A\left(\frac{A}{c_1+\ell_1}\right)
 \le L_b:=\frac{7428370561341219727}{767940001171875000000},
\]
and
\[
 (3/4-H_*)-L_b=
 \frac{38589144975363733283697549133345658036971433}
 {7579264015334677252965623903889516796875000000}>0.
\]
The opposite mixed orientation follows by scalar exchange after the two
actual chart derivations. Together the cases exhaust every negative
coordinate, and all strict inequalities survive $F=A=U_-$.
\end{proof}

\subsubsection{The nonnegative region: analytic outer reduction}\label{f245:outer-section}
We now suppose $r,s\ge0$, and sort the already derived scalar names so
$0<b\le a\le1$. In this subsection suppose $F\le U$ and
$\mathscr C\le H_*$. All comparison inequalities precede calibration.
The width bound gives $P(E-1)<2Ub$. Since $P\le1$, $R+P\le1$,
$r,s\ge R$, and $(4-P)/D\ge1$, the remainder in
\eqref{f245:determinant} is nonnegative and
\begin{equation}\label{f245:rough-det}
 \mathscr C\ge3/4+Rh-h^2,\qquad h=x+y+t.
\end{equation}
The function $N(w,d)$ increases in $d$, and
$N(w,1-w)=2(1+2w)/(3(1+w))\le1$ for $0\le w\le1$.
Therefore $\alpha\ge J_Ax$, $\beta\ge J_By$, and
$\alpha+\beta\ge J_0(x+y)$, where $J_0=2Q-P$.
Put
\[
 J=J_0-2Ub,\quad D_*=4-P-2Ub,\quad
 \kappa_0=Q(1-b/2),\quad V_0=J_0/2+Ub.
\]
Then $D>D_*>4-2b>0$, $2Q-EP>J>0$, and expansion of the positive
exponential tail gives, for every integer $m\ge2$,
\begin{equation}\label{f245:central-tail}
 A_c\ge Jt/2+\sum_{n=2}^m
            \frac{QD_*^{n-1}t^n}{2^{n-1}n!}.
\end{equation}
In particular
\begin{equation}\label{f245:outer-budgets}
 F\ge J(x+y)+Jt/2+\kappa_0t^2,
 \qquad F\ge J_0h-V_0t+\kappa_0t^2.
\end{equation}

\paragraph{Compactness and the scalar envelopes.}
The original compact rate bound $k\le15$ gives
$EP=4/(1+k^2)\ge2/113$. Moreover, along the central ray,
\[
 \ell/k=\frac{(1+k^2)\xi/2-1+2/\xi}{k^2}
 \ge\frac{2\sqrt{1+k^2}-1}{k^2}>\frac4{31}.
\]
For the last comparison the expression decreases with $\sqrt{1+k^2}$;
at $k=15$ it is equivalent to $62\sqrt{226}>931$, whose squared gap
is $1983$. Integration of $d\log\xi=k\,dt$ gives $z<31A_c/4$.
By \eqref{strip:exp-upper}, $e^{z_0}<\mathcal E_6(z_0)<27/4$
for $z_0=1519/800$, with exact gap
\[
 66551571865629146457223/921257902080000000000000.
\]
It follows that $P>p_0:=8/3051$ and $Q>1/10$, since
$4p_0-1/100=149/305100>0$.

Put $h_0=243/2000$, $L_0=h_0^2$. The exact gap
\[
 3/4-H_*-L_0
 =\frac{89765086136957463911195871}
 {39478417604337498218323684000000}>0
\]
and \eqref{f245:rough-det} imply
\[
 h>H:=\frac{R+\sqrt{R^2+4L_0}}2>R.
\]
The increasing lower bound
\[
 W(R)=R+\frac{L_0R}{R^2+L_0}<H
\]
has residual $-L_0^3/(R^2+L_0)^2$ in $X^2-RX-L_0$;
its derivative has positive numerator $R^4+L_0R^2+2L_0^2$.
Two further global support lines are
\[
 H\ge\frac45(R+h_0),\qquad H\ge\frac{49}{50}R+\frac7{25}h_0.
\]
Their residuals are $-(2R-3h_0)^2/25$ and
$-(7R-48h_0)^2/2500$, respectively; both lines are positive.
The further global bound
$H\ge Z(R):=R/2+h_0(16h_0^2+3R^2)/(16h_0^2+R^2)$
is retained: its residual is $-R^6/[4(16h_0^2+R^2)^2]$.
Completing the square in the first bound of
\eqref{f245:outer-budgets}, over the full interval $0\le t\le h$, gives
\begin{equation}\label{f245:envelope}
 F\ge g_0(h),\qquad
 g_0(v)=
 \begin{cases}Jv/2+\kappa_0v^2,&v\le J/(4\kappa_0),\\
 Jv-J^2/(16\kappa_0),&v\ge J/(4\kappa_0).
 \end{cases}
\end{equation}
It is strictly increasing and its branches agree. Hence $F>g_0(H)$.
Useful ratios are
\begin{equation}\label{f245:ratios}
 2-U\le J/\kappa_0\le3-2U,\quad
 J_0\ge3b,\quad J_0\ge Q(2-Q/4),\quad
 \kappa_0\ge Q(1-Q/4).
\end{equation}
For the lower ratio,
$J-(2-U)\kappa_0=U(a-b)-Uab/2+(1-U/2)b^2$ increases in $a$
and equals $(1-U)b^2$ at $a=b$. For the upper ratio use
$J-2\kappa_0=b(b-2U)$ and $\kappa_0\ge b(2-b)$.
The other comparisons follow from $ab\le Q^2/4$ and $b\le Q/2$.

\paragraph{The remaining outer inequalities.}
For $1/5\le Q\le9/16$ set
$\rho=1-Q+p_0$, $H_d=\rho^2+L_0$, and
$H_n=\rho(\rho^2+2L_0)$. Since $R>\rho>0$,
$H>H_n/H_d$. Also put $d=1+p_0+2h_0/7>1$.
The second support line gives $H>49(d-Q)/50$.
Square completion in the second bound of
\eqref{f245:outer-budgets} therefore gives $4(4-Q)(F-U)>M(Q)$, where
\begin{equation}\label{f245:poly-medium}
 M(Q)=\frac{49}{50}Q(8-Q)(4-Q)(d-Q)
       -4Q(1+U/2)^2-4U(4-Q).
\end{equation}
Here $50M''/49=-64-24d+(72+6d)Q-12Q^2\le-353/8<0$.
Thus its minimum is at an endpoint, where
\[
 M(1/5)=\frac{3367552463}{25425000000}>0,\qquad
 M(9/16)=\frac{2712290783}{5554176000}>0.
\]
This excludes the medium-$Q$ region.

For $Q\ge9/16$ in the linear branch,
$H\ge H_b=(2-U)/4=351/800$. The identity $R=H-L_0/H$ gives
$R>2/5$, since $H_b-L_0/H_b-2/5=1327/260000>0$.
Since $R\le(1-b)^2$, we have $b<2/5$ and
\[
 \tfrac94\kappa_0-J_0
  =[(a-b)(1-b/2)+b(2-5b)]/4\ge0,\qquad
 V_0/J_0\le\delta_0:=1/2+U/3=349/600.
\]
Square completion yields $F>J_0(H-9\delta_0^2/16)$.
Here $J_0\ge1071/1024$, and
\[
 \frac{1071}{1024}\left(\frac{351}{800}-\frac{121801}{640000}\right)
   -U=\frac{9724729}{655360000}>0.
\]
This excludes the linear branch, including its seam.

For $1/10\le Q\le1/5$, the actual compactness correlation gives
\[
 J_B\ge J_A\ge2a+\frac4{113a}-\frac2{113}
 \ge4\sqrt{2/113}-2/113>\lambda_0:=257/500.
\]
The squared margin for the last comparison is
$1538319/3192250000>0$. Thus $\alpha+\beta\ge\lambda_0(x+y)$.
Define
\[
 \begin{gathered}
 y_K=\frac{8-113Q^2}{8+113Q^2},\quad
 \tau_0=y_K+y_K^3/3+y_K^5/5,\\
 D_Q=4-Q^2/4-UQ,\quad c_Q=Q(2-U-Q/4)/2,
 \end{gathered}
\]
\[
 C_6(Q,t)=c_Qt+\sum_{n=2}^6\frac{QD_Q^{n-1}t^n}{2^{n-1}n!}.
\]
Since $E\ge8/(113Q^2)$ and $D<4$, we have
$t>\log(E)/2\ge\operatorname{arctanh}(y_K)>\tau_0$.
Also $D>4-Q\ge19/5$ and $z<31U/4$ give
$t<T_{max}=1519/1520$. The function $\tau_0$ decreases in $Q$, and
\[
 T_{max}-\tau_0(1/10)
   =\frac{54593043297215119}{964267380585605360}>0.
\]
Write $t=\tau_0+(T_{max}-\tau_0)\theta$, $0\le\theta\le1$.
Since $J/2\ge c_Q$ and $D_*\ge D_Q$, equation
\eqref{f245:central-tail} gives
$F>\lambda_0H_n/H_d+C_6(Q,t)-\lambda_0t$.
For an explicit polynomial clearing, put
\[
 y_N=8-113Q^2,\quad y_D=8+113Q^2,\quad T_d=15y_D^5,
 \quad T_n=15y_Ny_D^4+5y_N^3y_D^2+3y_N^5,
\]
\[
 t_n=T_n+(T_{max}T_d-T_n)\theta.
\]
The polynomial
\begin{equation}\label{f245:poly-thin}
 H_dT_d^6[C_6(Q,t_n/T_d)-U]
       +\lambda_0H_nT_d^6-\lambda_0H_dt_nT_d^5
\end{equation}
has $518$ positive Bernstein coefficients, of bidegree $(73,6)$, on
$[1/10,1/5]\times[0,1]$. Its clearing denominator $H_dT_d^6$ is positive.
This excludes the thin-$Q$ region.

Finally suppose we are in the quadratic branch with $R\ge1/5$.
Then $R<H\le J/(4\kappa_0)\le251/400$.
At fixed $R$, write $a=1-R/(1-b)$,
$0\le b\le1-\sqrt R$. The second derivative of
$JH/2+\kappa_0H^2$ in $b$ is
\[
 \frac{H[H(R+(1-b)^3)+R]}{(b-1)^3}<0.
\]
Its minimum is therefore at $b=0$ or $a=b$.
At $b=0$, use $H^2=RH+L_0$ and the first support line to obtain
\begin{equation}\label{f245:poly-highzero}
 (1-R)(H+H^2)\ge f_0(R):=\frac45(1-R^2)(R+h_0)+(1-R)L_0.
\end{equation}
Since $f_0''=-(24R+8h_0)/5<0$, its minimum on
$[1/5,251/400]$ is at an endpoint. The endpoint excesses over $U$ are
$68609/5000000$ and $198016309/1600000000$, respectively.
At $a=b$, put $w=\sqrt R$, $v=1-w$, and
$A=(1-w)(2-(1-w)/2-U)+(1-w^2)w^2>0$.
The same identity and support line give
\begin{equation}\label{f245:poly-highdiag}
 (2v-v^2/2-Uv)H+(2v-v^2)H^2
 \ge f_d(w):=\frac45 A(w)(w^2+h_0)+(1-w^2)L_0.
\end{equation}
On $[447/1000,4/5]$,
$f_d'''=-3(20000w^3-1514w+755)/625<0$, since
$20000(447/1000)^2>1514$. Moreover
$f_d'(447/1000)=321815473402479/625000000000000>0$.
Concavity of $f_d'$ and its positive left endpoint value preclude
a negative-to-positive crossing; hence $f_d$ has no interior minimum.
Its endpoint excesses over $U$ are, respectively,
\[
 \frac{22593992898613871}{1250000000000000000},\qquad
 \frac{10231929}{100000000},
\]
both positive.
This interval contains every required $\sqrt R$; the squared endpoint
gaps are $191/10^6$ and $1/80$. All comparisons include their endpoints.

\begin{lemma}[Outer reduction]\label{f245:outer}
Under $F\le49/200$, any nonnegative-coordinate point with
$\mathscr C\le H_*$ lies in
\[
 0\le r<9/20,\qquad 1/30<b\le a=1-r,\qquad
 1<E<104/55,\qquad R<1/5.
\]
It lies in the quadratic branch of \eqref{f245:envelope}.
\end{lemma}
\begin{proof}
The intervals above include all seams $Q=1/5$, $Q=9/16$,
$R=1/5$, and the linear/quadratic join. Thus the only region not excluded
is quadratic with $R<1/5$. Since $(1-a)^2\le R$ and
$(1-11/20)^2-1/5=1/400$, we have $a>11/20$.
The actual cusp inequality is
\[
 b\ge\frac{16aE}{243+4(aE)^2}.
\]
Its right side increases for $11/20\le aE<2$, and at $11/20$
it exceeds $1/30$ by $1979/732630$. Thus $b>1/30$.
Finally $a(E-1)<2U$ gives $E<1+2U/(11/20)=104/55$.
\end{proof}

\subsubsection{The residual comparison and its certificate}\label{f245:residual-section}
The remaining step is computer-assisted. We give the analytic reduction,
the exact meaning of every acceptance test, and a reproducible integer
verification specification. The finite certificate is a cover of an
auxiliary comparison domain, not a collection of physical samples.
Put $\eta=d_B/b$ and consider the closed rectangle
\begin{equation}\label{f245:rectangle}
 \mathcal B=[0,9/20]\times[1/30,1]\times[1,104/55]\times[0,1]
\end{equation}
in the order $(r,b,E,\eta)$. It contains every actual residual point.
The zero-width faces are comparison faces only; the actual leaves remain
strictly positive.

\paragraph{Continuous means and monotonicity.}
For $w,d\ge0$, $q=1+w$, direct integration gives
\begin{equation}\label{f245:positive-mean}
 e(w,d)=q(2+d)-\frac{2q(q+d)}d\log(1+d/q)
 =wd+2d^2\int_0^1\frac{v(1-v)}{q+dv}\,dv,
\end{equation}
with $e(w,0)=0$. Consequently $e(w,d)>wd$ for $d>0$ and the mean
increases strictly with $d$. It is nondecreasing in $w$ on
$0\le w,d\le1$, because
\[
 \partial_we=d-2d^2\int_0^1\frac{v(1-v)}{(q+dv)^2}\,dv
 \ge d-d^2/(3q^2)\ge0.
\]
Let $S(\gamma)=\sum_{n\ge1}n\gamma^{2n}/(2n+1)$.
With $\gamma=d/(2q+d)$ the continuous series is
\begin{equation}\label{f245:mean-series}
 e(w,d)=wd+4q(1+\gamma)S(\gamma),\qquad
 0\le S-S_N\le\frac{\gamma^{2N+2}}{2(1-\gamma^2)}.
\end{equation}
For $\gamma_E=(E-1)/(E+1)$, \eqref{hc:log-bounds} gives
$z=L_M(E)+\varepsilon_z$, $0\le\varepsilon_z\le U_M(E)-L_M(E)$, and
\begin{equation}\label{f245:log-series}
 \Phi(E)=4(1+\gamma_E)S(\gamma_E),\qquad
 \mathrm{rem}:=E-1-z=\frac{2\gamma_E^2}{1-\gamma_E^2}
                   +\frac{\gamma_E\Phi(E)}{1+\gamma_E}.
\end{equation}
No formula divides by $d$ or $E-1$ in the verification.

The following corner rules enclose entire boxes.
The central cost increases separately in $a,b,E$ on
$0<a,b\le1$, $1\le E<2$. Its $a$ derivative has numerator
$2[(4+b^2E)(E-1)-2bE\log E]\ge0$; the $b$ case is symmetric.
Since $E(4-P)-(D+4\log E)=4(E-1-\log E)\ge0$, the numerator
$2Q(4-P)-P(D+4\log E)$ of its $E$ derivative is bounded below
by $(2Q-EP)(4-P)>0$, including at $E=1$.
Thus we evaluate its lower and upper corners using the noncancelling form
\[
 A_c=\frac{(2Q-EP)(E-1)+EP\,\mathrm{rem}}D.
\]
For $\nu=e_B/b$ the integral identity gives
\[
 \nu=(1-b)\eta+2b\eta^2
       \int_0^1\frac{v(1-v)}{2-b+b\eta v}\,dv.
\]
It increases in $\eta$ and decreases in $b$, since
$\partial_b\nu\le-\eta+2\eta^2/[3(2-b)^2]\le0$.
The remainder $e_B-b(1-b)\eta$ increases in both $b$ and $\eta$.
The exact range of $b(1-b)$ is determined by its endpoint values and
its maximum $1/4$ if the interval contains $1/2$.
Finally $V_A^2=4a/(Eb)-a^2$ and $aD=4a-Eba^2$ increase in $a$ and
decrease in $b,E$, because $Eab<2$.

Writing $\xi=aE$, the right cost is
\begin{equation}\label{f245:beta-integrand}
 \beta=\int_0^\eta I(b,\xi,v)\,dv,\quad
 I=\frac b2+\frac{\xi(2-b)^2[3+(1-b+bv)^2]}
                    {2(4-\xi b)(2-b+bv)^2}.
\end{equation}
It increases in $\xi$ and $\eta$. Its monotonicity in $b$ is the
explicit polynomial identity in Section~\ref{f245:beta-monotone}.
Consequently the lower and upper $(b,\xi,\eta)$ corners bound the
entire cost box. These rules are bounds on functions over a box, not
point acceptance tests.

\paragraph{A strict upper bound from one residual budget.}
For a box let $B$ be an upper bound for
$U-A_c-\beta$, decreased to $U$ if necessary. If $B\le0$, an actual
positive left leaf is impossible. Otherwise
$\alpha\le B$. For $w=r_{\min}$ and $V^2$ an upper bound for $V_A^2$,
\eqref{f245:positive-mean} and the exact cost give
\[
 \alpha>L(w,V^2,d_A),\qquad
 L(w,V^2,d)=\frac d2+
 \frac{(3+w^2+wd)d(1+w)}{2V^2(1+w+d)}.
\]
This expression increases in $w$ and decreases in $V^2$.
Put $q=1+w$, $h_c=1+wq/V^2$, and
$p_c=q+(3+w^2)q/V^2-2B>0$. The inequality $L<B$ is equivalent to
$h_cd^2+p_cd-2Bq<0$, and hence
\begin{equation}\label{f245:upper-root}
 d_A<\delta:=\frac{4Bq}{\sqrt{p_c^2+8h_cBq}+p_c}.
\end{equation}
Taking rigorous upper bounds, set
$\bar d=\min(a_{\max},\delta_{\max},2B)$ and
$u_e=e(r_{\max},\bar d)_{\max}$. Then
\begin{equation}\label{f245:strict-upper}
 e_A<u_e.
\end{equation}
Strictness comes from the positive actual leaf and
$e_A>rd_A$, including $r=0$; it is not a numerical tolerance.

\paragraph{Lower bounds and the continuous affine root.}
The central sign supplies
$l_0=\max(0,[a(\Phi-r(E-1))]_{\min})\le e_A$.
For a stronger bound define
\begin{align*}
 q_B&=2-b,\qquad\epsilon=b\eta/q_B,\qquad\nu=e_B/b,\\
 M_0&=\mathrm{rem}-r(E-1),\\
 T_0&=[4M_0-(D-2)\nu]/D,\\
 p_{i0}&=-2(E-1)a[\Phi-r(E-1)]/D^2,\\
 I_0&=\eta(1+\epsilon)T_0-
       \eta e_B(1+\epsilon/2)/q_B\\
 &\quad-b^2\eta^2/(2q_B)+K_B\eta^2Ea/D,\\
 p_0&=p_{i0}-EaI_0/[D(1+\epsilon)^2],\\
 A_r&=2[(E-1)+E\eta/(1+\epsilon)]/D^2.
\end{align*}
Substitution in \eqref{f245:root-test} gives the exact identity
\begin{equation}\label{f245:affine-root}
 p_{\rm out}=A_re_A+p_0.
\end{equation}
Only if the entire box has $(A_r)_{\min}>0$ may we divide and set
$l=\max(l_0,[-p_0/A_r]_{\min})$. Otherwise only $l_0$ is used.
In particular, division is forbidden on the face $E=1,\eta=0$,
where $A_r=0$.

For fixed actual $a,b,E,e_B$, equation~\eqref{f245:determinant} has
second derivative $-2/(a^2D^2)$ with respect to $e_A$.
For $l_e<u_e$ its exact chord identity is
\begin{equation}\label{f245:chord}
 \mathscr C(e)=\frac{(u_e-e)\mathscr C(l_e)+(e-l_e)\mathscr C(u_e)}{u_e-l_e}
       +\frac{(e-l_e)(u_e-e)}{a^2D^2}.
\end{equation}
Thus $l_e\ge u_e$ contradicts $l_e\le e_A<u_e$; otherwise two
bounds $\mathscr C(l_e),\mathscr C(u_e)\ge H_*$ imply the actual
$\mathscr C\ge H_*$. The determinant is always evaluated in full.
No mean, realization, budget, or maximizing condition is imposed on
an artificial endpoint $l_e$ or $u_e$.

\begin{lemma}[Finite residual lemma, computer-assisted]\label{f245:residual}
Every box in the complete closed certificate specified in
Section~\ref{f245:certificate} either contains no actual residual point
with $F\le U$, or implies $\mathscr C\ge H_*$ for every such point.
The same-lift calibration therefore excludes every actual residual point
in \eqref{f245:rectangle} with $F\le U$.
\end{lemma}
\begin{proof}
The cover has $44\,813$ leaves. Its eight acceptance predicates are:
\begin{enumerate}
\item $b_{\min}>a_{\max}$, contradicting scalar ordering.
\item $[r(1-b)]_{\min}\ge1/5$, invoking Lemma~\ref{f245:outer}, including equality.
\item $[b(243+4(aE)^2)-16aE]_{\max}<0$, contradicting the actual cusp.
\item $[U-A_c-\beta]_{\max}\le0$, contradicting $\alpha>0$.
\item $l_0\ge u_e$, contradicting \eqref{f245:strict-upper}.
\item $(A_r)_{\min}>0$ and $l\ge u_e$, giving the same contradiction.
\item $l_0<u_e$ and both full endpoint bounds at $l_0,u_e$ are at least $H_*$.
\item $(A_r)_{\min}>0$, $l<u_e$, and both full endpoint bounds at $l,u_e$ are at least $H_*$.
\end{enumerate}
The last two use \eqref{f245:chord}. Every recorded predicate is checked
on the entire reconstructed closed box, with directed integer bounds.
The prefix-free complete binary tree covers the outward enclosure of
\eqref{f245:rectangle}; both children retain every common split face.
Thus the predicates cover every actual residual point with $F\le U$.
\end{proof}

\begin{proof}[Completion of Theorem~\ref{f245:theorem}]
Suppose $F\le U$. If either $r$ or $s$ is negative,
Lemma~\ref{f245:negative} gives $\mathscr C>H_*$, since
$U<129/500$. Otherwise the outer reduction and the residual certificate
prove $\mathscr C\ge H_*>0$. At an actual state $K_A>0$, so Corollary~\ref{cal:one-point}
on this qualified complete lift gives $\mathscr C<H_*$.
This contradiction proves the claim. Lawful cyclic transport covers the
other two listed contexts without changing physical time or cost.
\end{proof}

\space
\section{Forward supports, homogeneous cores, and physical arm integrals}\label{lp:forward}
Forward runs require two kinds of payments: a central pair supported by surrounding caps, and larger disjoint physical arms. We derive the common cap support once, then use it for the pair and arm estimates. Reverse homogeneous cores use the same full-lift setup but their own signed series.

\subsection{Cap forcing for forward pairs}\label{sec:final-analysis}
\subsubsection{Strict forcing by a forward pair}
\begin{lemma}\label{fin:forward-force}
At the middle switch of an actual positive $1201$ word on one full maximizing lift, let
\[
x=(1+u)/v,\qquad y=S/(2v),\qquad z=(1-u)/v.
\]
Then $xy+xz+yz=1$ and $y<x,y<z$. This conclusion has no cost threshold.
\end{lemma}
\begin{proof}
Apply Lemma~\ref{sc:pair} on this same full lift and normalize jointly.
Use its forward specialization $\delta=-1$, $P=F^{-1}$, $R=E^{-1}$,
where $E=e^{t_2/z}>1$, $F=e^{t_0/x}>1$. Write $V=(a,b)$, so $b\le0$,
and identify its $\nu_x,\nu_z$ with $\nu_1,\nu_2$, putting $\nu_0=(1,0)$
and $\lambda_i=V\cdot J\nu_i$. With $w_0=1-xz=y(x+z)>0$,
the common cap rows become
\begin{equation}\label{fin:responses}
h+A_z\lambda_2=c_z,\qquad h+A_x\lambda_1=c_x,
\end{equation}
where
\begin{align*}
A_z&=\frac{(E-1)(1+z^2)}{2E},&
c_z&=\frac{z(E-1)}{E^2}\bigl((1+z^2)f(E)-w_0(E-1)\bigr),\\
A_x&=\frac{(F-1)(1+x^2)}{2F},&
c_x&=\frac{x(F-1)}{F^2}\bigl((1+x^2)f(F)-w_0(F-1)\bigr).
\end{align*}
Here the displayed $A_i$ are positive; they are the negatives of the signed
coefficients in Lemma~\ref{sc:pair}. Its solved line is
\begin{equation}\label{fin:b-response}
b=rh-Q,\qquad Q=rh-b\ge0,
\end{equation}
where
\begin{align*}
r&=\frac{2}{w_0(x+z)}\left(\frac{zF}{F-1}+\frac{xE}{E-1}\right)>2,\\
Q&=\frac{2xz}{w_0(x+z)}\left\{\frac{(1+x^2)f(F)}F+
\frac{(1+z^2)f(E)}E-w_0\left(\frac{F-1}F+\frac{E-1}E\right)\right\}.
\end{align*}

For $G>1$, define $g(G)=(G-1-\log G)/(G-1)$. The elementary inequality $G^2-1-2G\log G>0$ follows by differentiating twice. It gives
\[
f(G)<(G-1)/2,\quad g(G)>(G-1)/(2G),\quad f(G)/G=(G-1)/G-g(G).
\]
Substitution in $Q$ implies
\begin{equation}\label{fin:one-small}
0\le Q<\frac{xz}{w_0}\left((x-y)\frac{F-1}F+(z-y)\frac{E-1}E\right).
\end{equation}
In particular, $y<\max(x,z)$.

The exact endpoint responses in Lemma~\ref{sc:pair}, with the same
forward substitution, give the actual nonnegative cap-ray coefficients
\[
h_-=A_z(b_z-b),\quad h_+=A_x(b_x-b),
\]
\[
b_z=2z\left(f(E)-\frac{z(E-1)}{x+z}\right),\quad
b_x=2x\left(f(F)-\frac{x(F-1)}{x+z}\right).
\]
Thus $Q+b_x=rh+h_+/A_x\ge0$. The preceding strict scalar bounds give
\begin{align*}
\frac{w_0(Q+b_x)}{2x}
<&\frac{F-1}{2F}\{x(z-y)+(F-1)y(z-x)\}
+\frac{E-1}{2E}z(z-y).
\end{align*}
If $z\le\min(x,y)$ the right side is nonpositive, a contradiction. If $z\le y$ still held, it would force $x<z\le y$, contradicting \eqref{fin:one-small}. Hence $z>y$. Interchanging $(x,F)$ and $(z,E)$ proves $x>y$. Positive durations supply strictness even when a cap coefficient or first switching jet vanishes.
\end{proof}

\subsection{The fixed-state support geometry of a forward pair}
\label{loc:forward-support-section}
In this section $x,y,z$ denote the positive reciprocal coordinates at the
middle of an actual capped $1201$ pair; by Lemma~\ref{fin:forward-force},
$x>y$ and $z>y$. Set $\ell=x+y+z$.
Use the actual complete central durations $d_2,d_0>0$,
$E=e^{d_2/z}$, $F=e^{d_0/x}$, and the cap responses
$A_i,c_i$ from \eqref{fin:responses}.

\begin{lemma}[Fixed-state cap family and support]\label{loc:forward-support}
Use the common affine cap line of Lemma~\ref{sc:pair}, in the forward
notation $V(h)=(ch-p_*,rh-Q)$, $0\le h\le h_*:=Q/r$, $r>2$.
Its determinant, gradient $N=\nabla_H\ell$ and actual trace
are those of \eqref{sc:polarization}; direct state algebra gives
$|N|^2=\ell^2-3$. The determinant is strictly concave along this line.
For $i=x,z$, with corresponding duration $d_i=d_0,d_2$, let
\[
 \alpha_i=\frac{i-y}{i+y},\qquad \beta_i=\frac{\alpha_i d_i}2,
 \qquad C_x=\frac{(x+z)(x-y)}{2x},\quad
 C_z=\frac{(x+z)(z-y)}{2z}.
\]
Both algebraic endpoints $V(0),V(h_*)$ belong to the triangle with vertices
\begin{align}
 W_0&=(-(x+x^{-1})\beta_x,0)=(-C_xd_0,0),\nonumber\\
 W_2&=((z+z^{-1})\beta_z,0)=(C_zd_2,0),\nonumber\\
 W_{12}&=-(y+y^{-1})(-\beta_x\nu_2+\beta_z\nu_1),\label{loc:forward-triangle}
\end{align}
where $\nu_1,\nu_2$ are the normals in Lemma~\ref{fin:forward-force}.
In addition, $V(h_*)$ lies on the segment $[W_0,W_2]$, and
\begin{gather}
 |W_{12}|^2=(y+y^{-1})^2(\beta_x+\beta_z)^2
                      -4(1+y^2)\beta_x\beta_z,\label{loc:forward-norm}\\
 h_*\leq(d_0C_x+d_2C_z)\frac{d_0d_2}{d_0+d_2}
       \leq\frac{\ell(d_0+d_2)^2}{8}.\label{loc:forward-height}
\end{gather}
For $\mathscr L(W)=N\cdot W/2+|W|^2/4$ and
$H=(\ell-N_y)h_*/2$, the actual determinant therefore satisfies
\begin{equation}\label{loc:forward-support-bound}
 \mathcal C_{\rm actual}\geq\frac34-
 \max\{\mathscr L(W_{12}),\mathscr L(W_0)+H,\mathscr L(W_2)+H\}.
\end{equation}
\end{lemma}
\begin{proof}
The common cap-line and determinant assertions follow from
Lemma~\ref{sc:pair} with the forward specialization just stated.

Put $\gamma_i=c_i/(2A_i)$. Substituting the explicit cap responses,
and writing $q_i=d_i/i>0$, gives the exact positive defect
\[
 \beta_i-\gamma_i=i\{\delta_0(q_i)+\alpha_i\delta_1(q_i)\}>0,
\quad \delta_0(q)=\frac{\sinh q-q}{e^q-1},\quad
 \delta_1(q)=\frac{q-1+e^{-q}}2.
\]
Both functions are positive for $q>0$, by differentiating their
numerators. The cap rows therefore imply
$V\cdot J\nu_i\leq2\gamma_i<2\beta_i$ when $h\geq0$.
The two algebraic endpoints also have $V_y\leq0$.
Solving the intersections of these three half-planes gives precisely
\eqref{loc:forward-triangle}; their positive normal dependence proves
that the intersection is bounded. At $h_*$ its vertical coordinate is
zero. Direct scalar products, using $xy+xz+yz=1$, prove
\eqref{loc:forward-norm}.

The same positive normal relation, $\gamma_i<\beta_i$, and
$\log G>1-G^{-1}$ for $G>1$ give respectively
\[
 Q\leq\frac{d_0z(x-y)+d_2x(z-y)}{y(x+z)},\qquad
 r>\frac{2xz}{y(x+z)^2}(d_0^{-1}+d_2^{-1}).
\]
Division proves the first height bound when $Q>0$; if $Q=0$ it is
immediate. Since $0<C_x,C_z<\ell/2$ and
$d_0d_2\leq(d_0+d_2)^2/4$, the second follows.

The determinant is concave along the fixed-state cap line, so its minimum
is at $h=0$ or $h=h_*$. At both endpoints $hV_y=0$.
The convex function $\mathscr L$ is bounded above by its triangle vertex
values at $h=0$, and by its two horizontal vertex values at $h=h_*$.
Since $\ell>|N|$, $H\geq0$, so the two conclusions combine to give
\eqref{loc:forward-support-bound}. The algebraic endpoints and vertices
are never asserted to be physical lifts and are never calibrated.
\end{proof}

\begin{lemma}[A common three-vertex estimate]\label{loc:forward-three-vertices}
Write
\[
 R=x/y>1,\quad t=z/y>1,\quad
 \rho=\frac{R-1}{R+1},\quad\sigma=\frac{t-1}{t+1},\quad
 w=\rho+\sigma,\quad n=3-w-\rho\sigma.
\]
Let $k_i=\ell d_i$, $k=k_0+k_2$, and assume
$w<w_0<1$, $n>n_0>0$, $k<B<n_0$, and $R,t<3$.
Set
\[
 u_0=\frac{2w_0}{(1+w_0)n_0},\qquad
 l_0=\frac{w_0}{(1+w_0/2)n_0}.
\]
Then the three losses in \eqref{loc:forward-support-bound} obey
\begin{align}
 \mathscr L(W_{12})&<\frac{Bw_0^2}{4(n_0-B)},\label{loc:vertex-nonhorizontal}\\
 \mathscr L(W_i)+H&\leq u_0^2B^2/4+u_0B^2/6
 &&\text{if }N\cdot W_i\leq0,\label{loc:vertex-negative}\\
 \mathscr L(W_i)+H&\leq\frac{(D_0+H_0)^2}{4H_0}
 &&\text{if }N\cdot W_i\geq0,\label{loc:vertex-positive}
\end{align}
where $H_0=2u_0B^2/3$ and
$D_0=w_0^2B/(8n_0)+l_0^2B^2/4$.
\end{lemma}
\begin{proof}
The exact scalar identities, obtained by substituting $R,t$, are
\begin{gather*}
 (y+y^{-1})/\ell=4/n,\qquad
 N\cdot W_{12}/2=\rho\sigma k/n,\qquad
 |W_{12}|\leq2(\rho k_0+\sigma k_2)/n,\\
 N\cdot W_0/2=\frac{\rho(\sigma-\rho)k_0}{(1+\rho)n},\qquad
 N\cdot W_2/2=\frac{\sigma(\rho-\sigma)k_2}{(1+\sigma)n},\\
 c_0:=C_x/\ell=\frac{2\rho(1-\rho\sigma)}{(1+\rho)n},\qquad
 c_2:=C_z/\ell=\frac{2\sigma(1-\rho\sigma)}{(1+\sigma)n}.
\end{gather*}
Put $\delta=|\rho-\sigma|$. The first three formulas and
\eqref{loc:forward-norm} give
\[
 \mathscr L(W_{12})\leq
 \frac{B(w^2-\delta^2)}{4n_0}+\frac{B^2(w+\delta)^2}{4n_0^2}.
\]
For $b=B/n_0<1$, the bracket is bounded by the exact square
\[
 (w^2-\delta^2)+b(w+\delta)^2
 =\frac{w^2}{1-b}-(1-b)\left(\delta-\frac{bw}{1-b}\right)^2.
\]
This proves \eqref{loc:vertex-nonhorizontal}.
The scalar formulas give $c_0,c_2\leq u_0$;
when a vertex has nonnegative signed term, its own coefficient is at most
$l_0$, since its corresponding $\rho$ or $\sigma$ is at most $w/2$.
Also
\[
 \frac{\ell-N_y}{2\ell}
 =\frac{R^2+Rt+t^2}{(R+t)(R+t+1)}\leq2/3.
\]
This is equivalent to
$R^2-Rt+t^2-2R-2t\leq0$ on $[1,3]^2$;
separate convexity bounds it by its four vertex values $-3,-1,-1,-3$.
Thus \eqref{loc:forward-height} gives
$H\leq(2/3)(c_0k_0+c_2k_2)k_0k_2/k$.
Discarding a nonpositive signed term proves \eqref{loc:vertex-negative}.
For a nonnegative term put $\lambda=k_i/k$.
The inequality $\rho(\sigma-\rho)\leq w^2/8$ and its mate give
\[
 \mathscr L(W_i)+H\leq D_0\lambda+H_0\lambda(1-\lambda)
 \leq(D_0+H_0)^2/(4H_0),
\]
by completing the square. This proves the final estimate.
\end{proof}

\begin{corollary}[Ratio-free support bounds]\label{loc:forward-support-global}
Without assuming $x<1$, $z<1$, or an upper ratio bound, let
$T=d_0+d_2$. Every point in the support triangle satisfies
\[
 |V|\leq(y+y^{-1})T/2,\qquad
 N\cdot V/2\leq\frac{(x-y)(z-y)}{4y}T.
\]
The same cap family has $h_*\leq\ell T^2/8$.
\end{corollary}
\begin{proof}
The height assertion is already \eqref{loc:forward-height}.
At $W_{12}$, \eqref{loc:forward-norm} and
$0<\alpha_i<1$ give the norm bound. At $W_0,W_2$, use
$C_x,C_z<\ell/2<(y+y^{-1})/2$; convexity of the norm then proves it
throughout the triangle.
The signed bound at $W_{12}$ is equality, since
$N\cdot W_{12}=(x-y)(z-y)T/(2y)$.
At $W_0$ its value is $(z-x)(x-y)d_0/(4x)$.
If this is nonpositive the claimed bound is immediate; otherwise
$(z-x)/x<(z-y)/y$ since $x>y$, and the result follows.
The $W_2$ case is the exchanged scalar computation.
Linearity in $V$ completes the proof on the triangle.
\end{proof}

\subsection{Own-run and mixed-edge estimates}
\label{loc:own-section}
\begin{lemma}[The two-arm integral]\label{loc:two-arm-integral}
Let four complete consecutive modes $1,2,0,1$ occupy $[t_0,t_4]$,
with middle switch $t_2$. Suppose their actual endpoint states satisfy
\[
 z(t_0)<y(t_0),\qquad x(t_4)<y(t_4).
\]
At the middle use $\rho,\sigma,w$ from
Lemma~\ref{loc:forward-three-vertices}, so $0<\rho,\sigma<1$ and $0<w<2$.
Their physical cost $J$ satisfies
\begin{equation}\label{loc:arms-G}
 J>G(\rho,\sigma)=\Phi((1-\rho)(1-\sigma)),\qquad
 \Phi(q)=\frac{4-4q+q\log q}{4-q}.
\end{equation}
At fixed $w=\rho+\sigma$, the exact algebraic minimum is balanced:
\begin{equation}\label{loc:arms-g}
 G\ge\Phi((1-w/2)^2)>
 \underline G(w):=\frac{w(4-w)(6+w)}{2(6-w)(2+w)}.
\end{equation}
The balanced comparison is strictly increasing in $w\in(0,2)$;
its equality point is not asserted to be an actual trajectory.
\end{lemma}
\begin{proof}
The reciprocal dynamics, obtained directly from the original fields, are
\[
\begin{array}{c|ccc}
 &x'&y'&z'\\\hline
0&0&(x+y)/x&-(x+z)/x\\
1&-(x+y)/y&0&(y+z)/y\\
2&(x+z)/z&-(y+z)/z&0.
\end{array}
\]
On the actual right two arcs, $q=(x-y)/(x+y)$ has
$q'=-2/(x+y)<0$. It goes from $\rho>0$ to a negative value, hence has
one zero. On $0\le q\le\rho$, the sum $a(q)=x+y$ satisfies
$a(q)\ge2y(t_2)/(1-q)$: mode $0$ has constant $x$, and on mode $1$
its constant $y(t_3)>y(t_2)$ gives the same inequality.
Since $\ell a=1+a^2(3+q^2)/4$, that portion costs at least
$\rho/2+y(t_2)^2F(\rho)/2$, where
$F(b)=b+2\log(1-b)+4/(1-b)-4$.
The omitted cost is strictly positive by the strict endpoint sign.
For the disjoint left portion use $q=(z-y)/(z+y)$,
$q'=2/(z+y)>0$, and read its integral backwards from the middle;
this reverses an integration parameter, not the physical lift.
Thus $J>w/2+y^2[F(\rho)+F(\sigma)]/2$.
For $q=(1-\rho)(1-\sigma)$ one has $y^2=q/(4-q)$ and
\[
 F(\rho)+F(\sigma)=w+2\log q+4(2-w)/q-8,
\]
which proves the exact product formula.
Now $\Phi'(q)=(4\log q-8-q)/(4-q)^2<0$ on $(0,1)$,
and $q\le(1-w/2)^2$, with equality precisely when $\rho=\sigma$.
Finally $\log q>-(1-q)/\sqrt q$: with $a=\sqrt q<1$, the difference
$\log a-(a-a^{-1})/2$ vanishes at $1$ and has derivative
$-(a-1)^2/(2a^2)<0$. Substituting $a=1-w/2$ gives $\underline G(w)$.
\end{proof}

\begin{lemma}[Four charges with endpoint state signs]\label{loc:state-sign41}
Under the hypotheses of Lemma~\ref{loc:two-arm-integral}, one has
$J>41/100$.
\end{lemma}
\begin{proof}
Assume $J\leq a=41/100$. The balanced comparison in
\eqref{loc:arms-g} and the exact gap
\[
 \underline G(249/500)-41/100=2600387/1145333000>0
\]
force
\begin{equation}\label{loc:own41-gate}
 w<249/500<1/2,\qquad R,t<3,\qquad n>39/16.
\end{equation}

The actual central two-arc cost $I$ is strictly smaller than $J$.
Its durations sum to $T< a/\sqrt3$. The exact density profiles, measured
from the middle along each actual arc, are
\[
 \ell_0(v)=(x+y)e^{v/x}+(x+z)e^{-v/x}-x,
\quad
 \ell_2(v)=(y+z)e^{v/z}+(x+z)e^{-v/z}-z.
\]
They are strictly convex with initial slopes $(y-z)/x$, $(y-x)/z$.
Integrated tangent bounds give
\[
 I>\ell T-\tfrac12\{(z-y)d_0^2/x+(x-y)d_2^2/z\}>
 \ell T-T^2,
\]
using $R,t<3$. Consequently $k=\ell T<B=a+a^2/3=13981/30000$.

Apply the support estimates with $n_0=39/16$.
For the nonhorizontal vertex retain $w_0=249/500$; for the two horizontal
bounds use the weaker $w_0=1/2$, so $u_0=32/117$, $l_0=32/195$.
Their exact upper values are
\[
 L_{12}=\frac{866835981}{59144000000},\quad
 L_- =\frac{2150151971}{154001250000},\quad
 L_+ =\frac{1864266368476801}{133467750000000000}.
\]
The formulas for them are \eqref{loc:vertex-nonhorizontal}--
\eqref{loc:vertex-positive}, and direct rational comparison gives
$L_-,L_+<L_{12}$ and
\[
 59/4000-L_{12}=5538019/59144000000>0.
\]
Hence the actual determinant exceeds $2941/4000>0$.
Independently, \eqref{loc:forward-height} gives
$0\leq h/\ell\leq T^2/8<a^2/24<1$, so its actual trace $\ell-h$
is positive. The same complete lift is normal by
Lemma~\ref{fin:forward-force}. Only now apply the harmonic cap,
which gives $\mathcal C<C_0<2941/4000$, a contradiction.
\end{proof}

\begin{corollary}[Four own arcs]\label{loc:Own39}
In the actual positive contexts $101202$, $212010$, $020121$, the cost
of the four complete interior arcs is greater than $41/100$.
\end{corollary}
\begin{proof}
Use the full legal cycle taking $101202$ to $212010$.
Its four charges are $1201$. The actual left return $212$, transformed
to $020$, gives $z(t_0)<y(t_0)$ by Lemma~\ref{loc:return-root};
the actual right return $010$, transformed to $202$, gives
$x(t_4)<y(t_4)$. These are precisely the two hypotheses of
Lemma~\ref{loc:state-sign41}, which proves $J>41/100$.
Legal cycles preserve the two other contexts and their charged positions.
\end{proof}

\begin{theorem}[Mixed edge and long own-run resources]\label{loc:Edge41}
For each actual positive full-seven context $1012012$ and its legal phases,
charge the four complete arcs $[1,4]$; for $0120121$ and its legal phases,
charge $[2,5]$. Each charged cost exceeds $41/100$.
The respective uncharged supporting arcs $5$ and $1$ must be complete;
the two exterior caps must remain positive.
Consequently an ordinary forward run with at least seven increments has
own-arc cost greater than $41/50$ on its full eligible context.
\label{loc:OwnLong82}
\end{theorem}
\begin{proof}
For $1012012$, the charges have modes $0120$.
The legal $+1$ cycle makes them $1201$ with middle coordinates
$(x_{\mathrm{new}},y_{\mathrm{new}},z_{\mathrm{new}})=(z,x,y)$.
The left return $101$ gives $x(t_1)>y(t_1)$, so the first required
canonical endpoint sign is $z_{\mathrm{new}}<y_{\mathrm{new}}$.
The fully capped $2012$ on arcs $[3,6]$ gives
$z(t_5)<x(t_5)$ by Lemma~\ref{fin:forward-force} in its legal phase,
so the second sign is $x_{\mathrm{new}}<y_{\mathrm{new}}$.
Lemma~\ref{loc:state-sign41} applies. Its latter sign uses the complete
uncharged supporting arc $5$; a fragment of that arc would not suffice.

For $0120121$, use the $-1$ cycle on the four charges $2012$ to obtain
canonical $1201$ with coordinates $(y,z,x)$.
The fully capped $0120$ on $[0,3]$ gives the initial sign
$x<z$, and the terminal return $121$ gives the final sign $z>y$.
In the new chart these are exactly the two endpoint inequalities.
This uses the complete uncharged arc $1$ and proves the second estimate.

A forward run with $a$ increments has $a+1$ own arcs, with successive
labels increasing by one. Its first seven context arcs form the first
mixed pattern; its last seven form a legal phase of the second.
For $a\geq7$, the first four and last four own arcs are disjoint.
Apply both bounds on one eligible whole-run lift and add their physical
costs. Their supporting contexts may overlap but their charges do not.
The remaining own arcs, if any, have positive cost. The result is
$>2(41/100)=41/50$.
\end{proof}

\begin{theorem}[Six own arcs]\label{loc:OwnSix600}\label{loc:OwnSix609}
For each actual positive context $10120121$, $21201202$, $02012010$,
the physical cost of complete positions $1,\ldots,6$ exceeds $609/1000>3/5$.
All eight source arcs and the same full containing-interval hypotheses are retained.
\end{theorem}
\begin{proof}
Work with $10120121$. Its prefix $[0,6]=1012012$ and suffix
$[1,7]=0120121$ supply two capped middles at $t_3$ and $t_5$,
with the same full normal lift. Use canonical coordinates $(z,x,y)$ at
$t_3$ and $(y,z,x)$ at $t_5$.
The two-arm integrals of Lemma~\ref{loc:two-arm-integral} give
\begin{equation}\label{loc:own6-joint-integral}
 J_6>G(\rho_3,\sigma_3)+G(\rho_5,\sigma_5).
\end{equation}
Here the two inner arms really are disjoint: on the shared actual modes
$2,0$ between $t_3,t_5$, the quantity $q=(z-x)/(z+x)$ has
$q'=-2/(x+z)<0$. It starts positive and ends negative.
The right arm at $t_3$ uses precisely $q\geq0$; the left arm at $t_5$
uses $q\leq0$, with canonical variable $-q$. They meet only at the
unique zero. The outer arms lie in disjoint outer arc pairs, and their
strict endpoint signs leave positive omitted cost.

Suppose $J_6\leq A=609/1000$. Select one actual middle at which $G<A/2$,
which exists by \eqref{loc:own6-joint-integral}.
The same increasing balanced comparison and
$\underline G(7/20)-A/2=10451/10622000>0$ give
\[
 w<w_0=7/20,\quad R,t<27/13,\quad n>n_0=4191/1600.
\]
The central two arcs at the selected actual middle are a subset of the six
charges. The integrated tangent proof above, now using $R,t<27/13$,
gives
\[
 k=\ell T<B=A+\frac{14}{13}\frac{A^2}{6}
   =\frac{8782389}{13000000}.
\]

Use the support family of Lemma~\ref{loc:forward-support} and put
$u=2w_0/[(1+w_0)n_0]$, $l=w_0/[(1+w_0/2)n_0]$.
As before $c_0,c_2\leq u$, and the coefficient belonging to a vertex
with nonnegative signed term is at most $l$.
In addition $c_0+c_2\leq2l$: indeed
$\rho/(1+\rho)+\sigma/(1+\sigma)\leq w/(1+w/2)$ by concavity,
and $1-\rho\sigma\leq1$.
The nonhorizontal loss satisfies
\[
 \mathscr L(W_{12})<\frac{Bw_0^2}{4(n_0-B)}<59/4000;
\]
the exact margin is $25154971/6125936000>0$.
For a horizontal vertex with nonpositive signed term use
$\lambda^2(1-\lambda),\lambda(1-\lambda)^2\leq4/27$ in the
height bound to obtain
\[
 \mathscr L(W_i)+H\leq B^2(u^2/4+16l/81)<59/4000.
\]
Its exact margin is $87509319798137/3242851377187500000>0$.

For a vertex with nonnegative signed term, let $\lambda=k_i/k$.
The sharper own coefficient and the retained duration split give
\[
 \mathscr L(W_i)+H\leq P(\lambda):=
 d\lambda+e\lambda^2+\frac23B^2[l\lambda+u(1-\lambda)]\lambda(1-\lambda),
\quad d=\frac{w_0^2B}{8n_0},\quad e=\frac{l^2B^2}{4}.
\]
Put $a_3=2B^2(u-l)/3>0$, $p=d+B^2(5l-u)/6$ and
$q=B^2l(2/3-l/4)>0$. Completing the square gives, for $0\le\lambda\le1$,
\[
 P(\lambda)=\frac{a_3}4+p\lambda-q\lambda^2
             -a_3(1-\lambda)(\lambda-1/2)^2
 \le\frac{a_3}4+\frac{p^2}{4q}
 =\frac{2918301475265239513583}{205447452944760000000000}.
\]
Here $u=22400/113157$, $l=22400/196977$, and the exact margin is
\[
 \frac{59}{4000}-\left(\frac{a_3}4+\frac{p^2}{4q}\right)
 =\frac{112048455669970486417}{205447452944760000000000}>0.
\]

Thus the actual determinant is $>2941/4000>0$.
Independently $h/\ell\leq T^2/8<A^2/24$ gives positive
actual trace, with ratio
$(\ell-h)/\ell>1-A^2/24=7876373/8000000>0$. The same-lift harmonic cap is now applicable and gives
$\mathcal C<C_0<2941/4000$, a contradiction.
\end{proof}

\subsection{Forward and reverse homogeneous cores}\label{sec:linear-strips}

\subsubsection{Actual reverse caps and common scalar kernels}
\label{strip:reverse-kernel}

Throughout this analytic block, $C=\det\Lambda$ is the same-lift matrix
Casimir, while $A,I,J,W$ denote physical costs. An eligible actual
six-context consists of six positive vertex arcs in one genuine fundamental
minimizing interval, with all four charged interior arcs complete.
Reanchoring retains every complete charged, return and sign arc and positive
portions of the actual exterior caps; the charged span is strictly shorter
than the interval. Select one full lift there and prove normality from its
actual positive cap equations, as in the proof of Theorem~\ref{sw:normality}, before
joint normalization. The phrase ``selected-six premise'' adds no numerical
or selection hypothesis.

\begin{lemma}[The actual reverse-cap line]\label{strip:reverse-primitive}
On one complete normal maximizing free-terminal-time lift, jointly
normalized to $\eta=1$, retain an actual positive capped reverse pair
$1\mid02\mid1$ with both charged arcs complete. Its actual cap equations,
outer-ray residuals and determinant are respectively
\eqref{sc:rows}, \eqref{ctr:ratio-residual} and
\eqref{sc:polarization}. They imply a nonempty affine cap segment
$0\le h\le h_*$ and strict determinant concavity on it.
These identities require neither a numerical cost budget nor harmonic
calibration; the algebraic endpoints are not asserted to be trajectories.
\end{lemma}
\begin{proof}
Apply Lemma~\ref{sc:pair} with $\delta=1$ at the actual middle,
using $x=(1+u)/v$, $y=(u^2+v^2-1)/(2v)$ and $z=(1-u)/v$.
Put $P=e^{t_0/x}$, $R=e^{t_2/z}$ for the two complete physical durations.
Its $A_x,A_z,\beta_x,\beta_z,\nu_x,\nu_z$ give exactly the
reverse cap rows \eqref{sc:rows}. Write $D_{\rm cap}=D$ and
$r_0=r>0$ in its solved line $V_y=q_0-r_0h$; then
$D_{\rm cap}\ge0$, $q_0\ge0$ and $0\le h\le h_*=q_0/r_0$.
The common polarization and endpoint minimum retain the full matrix
and the state covector, without a cost budget or calibration.
Put $d(g)=g(g-1-\log g)/(g\log g-g+1)$. The actual cap residual implies
\begin{equation}\label{ctr:ratio-residual}
 y/z<\max(Pd(P),d(R)),\qquad
 y/x<\max(d(P),Rd(R)).
\end{equation}
Put $a(g)=\log g/(g-1)$. The actual endpoint signs in
\eqref{sc:interval} give
\[
 V_y\ge b_x=2x\left\{a(P)-\frac{z+x/P}{x+z}\right\}.
\]
Since $V_y=q_0-r_0h$ and actual $h\ge0$, one has $q_0\ge b_x$.
Set $k=y/z$ and $r=x/z$. Dividing this last inequality by positive factors gives the nonnegative outer-ray residual
\[
 r[P+k/P-(P+k)a(P)]
   +R+3k-k(P+1)a(P)-(k+1)Ra(R).
\]
Its two components are
\[
 f(P)[Pd(P)-k]/P,
 \qquad f(R)[d(R)-k]+k[2-(P+1)a(P)].
\]
The final bracket is strictly negative, since $a(P)>2/(P+1)$ by Jensen's inequality applied to
$a(P)=\int_0^1[1+(P-1)t]^{-1}dt$. This proves the first inequality in \eqref{ctr:ratio-residual}; the outgoing mode-$2$ cap similarly gives $V_y\ge2z\{a(R)-(x+z/R)/(x+z)\}$, proving the second inequality by the same algebra.

The remaining cap-segment, strict concavity and trace assertions are
\eqref{sc:interval}--\eqref{sc:polarization} on this same selected lift.
\end{proof}

\begin{lemma}[A common positive kernel]\label{strip:reverse-credits}
Put
\[
 f(g)=\frac{g\log g}{g-1}-1,\qquad
 d(g)=\frac{g(g-1-\log g)}{g\log g-g+1}.
\]
For every $s>0$,
\[
 \frac s2+\frac{5s^2}{60+s^2}
 < f(e^s)<\frac s2+\frac{s^2}{12}.
\]
For $0<s\le5/4$ one also has
\[
 \log d(e^s)<\frac{2s}{3}+\frac{s^3}{360}.
\]
Consequently $d(e^s)<e^{cs}$ for $0<s\le B$ in all three rows
\[
 (c,B)=(67/100,1),\quad(27/40,154/125),\quad(7/10,6/5).
\]
The quartic lower bound $f(e^s)>s/2+s^2/12-s^4/720$ also follows.
\end{lemma}
\begin{proof}
For $t=s/2>0$, write
$\phi(s)=f(e^s)-s/2=t\coth t-1$. The absolutely convergent identities
\begin{align*}
 (3+t^2)\sinh t-3t\cosh t
 &=\sum_{n\ge2}\frac{4n(n-1)t^{2n+1}}{(2n+1)!}>0,\\
 (15+t^2)t\cosh t-(15+6t^2)\sinh t
 &=\sum_{n\ge3}\frac{8n(n-1)(n-2)t^{2n+1}}{(2n+1)!}>0
\end{align*}
prove the upper and lower bounds for $\phi$ after division by
$3\sinh t$ and $(15+t^2)\sinh t$, respectively. The quartic consequence is
\[
 \frac{5s^2}{60+s^2}-\left(\frac{s^2}{12}-\frac{s^4}{720}\right)
 =\frac{s^6}{720(60+s^2)}>0.
\]
For the logarithmic bound set $A=e^s-1-s>0$ and
$B=e^{-s}-1+s>0$, so $d(e^s)=A/B$. Direct expansion gives
\[
 \left(\frac23+\frac{s^2}{120}\right)AB-A'B+AB'
 =\sum_{j\ge0}\frac{(8j^3+56j^2+46j+10)s^{2j+6}}
                         {60(2j+6)!}>0.
\]
Indeed the coefficient of $s^{2n}$ before the low-degree cancellations
is $(8n^3-16n^2-74n+160)/(60(2n)!)$; the terms through degree four
vanish, and $n=j+3$ gives the displayed positive series.
Dividing by $AB$ and integrating from zero, where $A/B\to1$, proves
$\log d(e^s)<2s/3+s^3/360$ for every $s>0$.
The reserves
$c-2/3-B^2/360$ for the three stated $(c,B)$ rows are respectively
$1/1800$, $23159/5625000$ and $11/375$, all positive.
\end{proof}

For $0\le s\le6/5$ put
\[
 L(s)=(1-s+s^2/2-s^3/6)/2,\quad
 P(s)=1-s/2+s^2/6-s^3/24,\quad N(s)=L(s)+P(s).
\]
\begin{lemma}[A two-vertex reciprocal credit]\label{strip:cone-credit}
Let $p,q,r>0$, $pq+pr+qr=1$, with all ordered ratios at most $\rho\ge1$.
For $L,P>0$,
\[
 Lp^2+Ppq\ge
 \min\left\{\frac{L+P}{1+2\rho},
            \frac{L+\rho P}{\rho(\rho+2)}\right\}.
\]
This is the exact minimum on the algebraic cone. If
$\rho P\ge(1+\rho)L$, the first expression is the minimum.
For the functions $L(s),P(s)$ just defined,
$0\le s\le6/5$, this condition holds for every $\rho\ge1$.
\end{lemma}
\begin{proof}
Write $a=q/p$, $b=r/p$. The quantity is
$(L+Pa)/(a+b+ab)$ and decreases with $b$.
For $a\le1$ the largest permitted $b$ is $\rho a$; the resulting
quotient has derivative with numerator
$-L(1+\rho+2\rho a)-\rho Pa^2<0$.
For $a\ge1$ the largest $b$ is $\rho$; the derivative has the
constant sign of $\rho P-(1+\rho)L$.
It is therefore enough to test $(a,b)=(1,\rho),(\rho,\rho)$.
Both are permitted, proving the minimum formula.

The indicated $L,P$ decrease to
$L(6/5)=29/250>0$, $P(6/5)=71/125>0$, and
\[
 P-2L=s(3s^2-8s+12)/24\ge0,
\]
since the quadratic has negative discriminant and positive leading term.
Thus $\rho P-(1+\rho)L\ge(\rho-1)L\ge0$.
\end{proof}

\subsubsection{A complete forward three-core}
\label{strip:forward-three}

\begin{proposition}[Forward three-core]
\label{strip:pf3}\label{strip:forward1175}
On an eligible full positive context $01201$, or either additive phase,
the sum $J=c_1+c_2+c_3$ of its three complete interior physical costs
satisfies $J>47/400>29/250$. The uncharged caps $0,4$ remain positive.
\end{proposition}
\begin{proof}
Use one complete normal lift, jointly normalized as in
Lemmas~\ref{fin:forward-force} and~\ref{loc:forward-support}.
Only the two genuinely capped middles of this five-arc context will be
used.  At either middle write
\[
 xy+xz+yz=1,\quad R=x/y>1,\quad t=z/y>1,\quad D=R+t+Rt,
 \quad \ell=(R+t+1)/\sqrt D.
\]
For its adjacent own durations $t_x,t_z$, put $q_i=t_i/i$ and
$T=t_x+t_z$.  Direct integration of the physical density gives
\begin{align*}
 A_x&=xy(e^{q_x}-1)+xz(1-e^{-q_x})+x^2(2\sinh q_x-q_x),\\
 A_z&=xz(1-e^{-q_z})+yz(e^{q_z}-1)+z^2(2\sinh q_z-q_z).
\end{align*}
With $N=\nabla_H\ell$, Lemma~\ref{loc:forward-support} supplies the support triangle and
$h_*\le\ell T^2/8$; we use its three-vertex determinant bound directly.

First we select one of the two existing middles.  Present the context
as $12012$.  Along its shared mode-$0$ arc the neighboring middle is
\[
 (x,Y,Z)=(x,F(x+y)-x,(x+z)/F-x),\qquad F>1.
\]
The unused minima at the two middles are $y$ and $Z$.
If both
$\min(x/y,z/y)$ and $\min(Y/Z,x/Z)$ are at least $m>1$,
put $s=\max(y,Z)$ and choose $w$ with
$(x+s)(x+w)=1+x^2$.  Then $x/s,w/s\geq m$.
Monotonic increase of the second coordinate and decrease of the third
show that the symmetric segment $(x,s,w)$ to $(x,w,s)$ lies in
the actual shared arc.  Its factor is $(x+w)/(x+s)>1$.
In coordinates $R=x/s,t=w/s$ its physical cost is
\[
 A_{\rm sym}=\frac R{R+t+Rt}
 \left[2(t-1)-R\log\frac{R+t}{R+1}\right].
\]
The elementary strict bound $\log g<g-1$ gives
\[
 A_{\rm sym}>L(R,t):=
 \frac{R(R+2)(t-1)}{(R+t+Rt)(R+1)}.
\]
This rational function increases in $t$, since
\begin{align*}
 L_t&=\frac{R(R+2)(2R+1)}{(R+1)(R+t+Rt)^2}>0,\\
 L(R,m)-\frac{m-1}{m+1}
 &=\frac{(m-1)(R-m)}{(R+1)(m+1)(R+m+Rm)}\ge0.
\end{align*}
Hence the actual shared cost exceeds $(m-1)/(m+1)$.

Suppose now $J\leq a\leq1/5$ and set $M=(1+a)/(1-a)$.
Since the shared arc has cost strictly less than $J$, the preceding
argument selects a genuine capped middle with
$\min(R,t)<M$.  Its two adjacent complete costs satisfy
$A_x+A_z<J\leq a$.  No third middle or extra cap has been introduced.

All terms in the two cost formulas are positive and $xz>1/3$.
Thus $e^{-q_i}>d:=1-aD/(Rt)>1-3a>0$.
The elementary inequalities
$1-e^{-q}>qe^{-q/2}$, $e^q-1>q$, and
$2\sinh q-q>q$ imply
\[
 \ell T<\frac a{\sqrt d},\qquad
 \ell T<B:=\frac a{\sqrt{1-3a}}.
\]
Define $k=(y+1/y)/\ell=(R+1)(t+1)/(R+t+1)$ and $g=G/\ell$, where $G=(x-y)(z-y)/y$.
The unbounded ratio sector is controlled without polynomial expansion:
$D/(Rt)\le2+1/R$ and $R+t+1\ge R+2$ give
\[
 (R+t+1)^2d\ge(R+2)^2[1-a(2+1/R)]>(R+1)^2(1-2a).
\]
Indeed the last difference is
\[
 2R+3-a(5R+10+4/R)\ge R+1-4/(5R)>0,
\]
and the factor $1-a(2+1/R)\ge1-3a$ is positive.
Applying it and its exchanged version in
Corollary~\ref{loc:forward-support-global} gives, with
$b=a/\sqrt{1-2a}$,
\[
 |V|<\frac{(M+1)b}{2},\qquad
 \frac{g\ell T}{4}<\frac{(M-1)b}{4}.
\]
Here the norm bound holds throughout the support triangle. The nonhorizontal
vertex loss is less than $(M-1)b/4+(M+1)^2b^2/16$.
For a horizontal vertex opposite in sign to $N_x$, discard its
nonpositive linear loss and use the same norm bound.
For a horizontal vertex of the same sign, retain both terms. If $R\ge t$,
that vertex is $W_2=(L_z,0)$, and
\[
 \frac{N_xL_z}{2}=\frac{t_z(t-1)(R-t)}{4t\sqrt D}
 \le\frac{g\ell T}{4},\qquad
 \frac{L_z}{\ell t_z}=\frac{(R+t)(t-1)}{2t(R+t+1)}
 <\frac{t-1}{2t}<\frac{M-1}{2}.
\]
For $t\ge R$ exchange the scalar names; no actual arc is reversed.
The height loss $H=(\ell-N_y)h_*/2$ is less than $B^2/8$, since
$|N|<\ell$. Thus the three-vertex bound gives
\[
 C>\min\{L_0,L_{h,o},L_{h,s}\},
\]
where
\begin{align*}
 L_0&=\tfrac34-(M-1)b/4-(M+1)^2b^2/16,\\
 L_{h,o}&=\tfrac34-(M+1)^2b^2/16-B^2/8,\\
 L_{h,s}&=\tfrac34-(M-1)b/4-(M-1)^2B^2/16-B^2/8.
\end{align*}
This also covers $h_*=0$ without dividing by its height.

Set $a=47/400$. Then
\[
 M=447/353,\qquad B^2=2209/103600,\qquad b^2=2209/122400.
\]
The rational bracket $b<134341/10^6$ is strict because
\[
 (134341/10^6)^2-b^2
 =\frac{18154993}{153000000000000}>0.
\]
Substituting this upper bracket in the negative linear terms in $b$,
and retaining the exact $b^2$, gives the respective determinant floors
\[
 \frac{28035851215757}{38130354000000},\qquad
 \frac{11717258087207}{15801218697600},\qquad
 \frac{11913796512899}{16136865500000}.
\]
The first is the smallest. Its difference from the harmonic ceiling
\[
 \mathcal H=
 \frac{29025933167887355473205403}{39478417604337498218323684}
\]
is the positive rational number
\[
 \frac{10483829210426125893586595821272197}
 {376331509653305185634762839376034000000}.
\]
Hence the actual determinant satisfies $C>\mathcal H>0$.
Independently, $h/\ell<T^2/8<B^2/24$, so the actual future trace
has ratio
$(\ell-h)/\ell>2484191/2486400>0$.
Only now apply the harmonic calibration to this same complete lift;
it gives $C<\mathcal H$, a contradiction.  The closed contrary
budget included equality, proving the strict assertion.
\end{proof}

\subsubsection{The central forward-pair estimate}\label{sec:forward-final}
\begin{proposition}\label{fp:final}
In a genuine complete positive forward context $012012$ on an original
minimizer, or a lawful symmetry image of the entire problem, the two
complete central physical arcs satisfy
\[
 c_2+c_3>\frac{180483}{500000}>\frac{20901}{62500}.
\]
All six actual arcs, including the positive exterior caps at positions
$0$ and $5$, are retained in one genuine fundamental interval.
\end{proposition}
\begin{proof}
Lemma~\ref{loc:automatic-properness} makes this six-arc source proper.
Choose the full maximizing lift on its genuine containing interval before
restriction. Lemmas~\ref{fin:forward-force} and~\ref{loc:forward-support}
establish normality from its actual complete caps; jointly normalize the
entire lift. Theorem~\ref{nx:forward-pair}, proved independently in
Section~\ref{nx:forward-proof}, applies to these same two central physical
charges. The exact difference between the two displayed constants is
$13275/500000>0$. No extra cap, multi-period interval or freely chosen
endpoint extremal has been introduced.
\end{proof}

\begin{corollary}[A four-charge consequence]\label{strip:pthird}
The four interior physical arcs of a complete $012012$ source under the
same hypotheses have cost greater than $1/3$.
\end{corollary}
\begin{proof}
Proposition~\ref{fp:final} gives $c_2+c_3>20901/62500>1/3$.
The other two interior arcs have positive cost.
\end{proof}

\subsubsection{Reverse windows in a genuine interval}
The next lemma allows the three-fifths estimate to be applied without
an additional omitted-arc assumption. It changes the genuine cut and
selects a fresh lift; it does not enlarge the analytic domain of
Theorem~\ref{wk:five600}.

\begin{lemma}[Completing the caps of a reverse five-window]
\label{strip:five-cap-completion}
Let a normalized original global minimizer of physical cost below
$\Cstar$ contain the actual consecutive pattern
$1\mid021\mid0$, with all three interior maximal arcs complete and
positive portions of both exterior maximal arcs in one genuine
fundamental interval. Then its three physical interior costs satisfy
$c_1+c_2+c_3\ge3/5$. Whole lawful images have the same conclusion.
No matrix adjoint is transferred from one cut to another.
\end{lemma}
\begin{proof}
The three complete charged arcs occupy a span strictly shorter than the
fundamental period, since both retained cap portions have positive
duration. Thus the number $B$ of maximal arcs in a period satisfies
$B\ge4$. If $B=4$, the two exterior occurrences would be the same
physical arc in consecutive periods. The actual twist would require
$j_4=j_0+1\pmod3$, whereas $10210$ has $j_0=1,j_4=0$.
This contradiction proves $B\ge5$.

Reanchor at the actual initial switch of the first exterior maximal
arc, using Lemma~\ref{fd:reanchor}, whose assertion is for every real cut,
including a switch. In its notation
$G_s(T)=G_e(s)^{-1}G_e(s)R_*=R_*$ and
$W_s(T)=\Phi W_s(0)$. The same fixed-endpoint homotopy retains the
real winding. At a switch the trajectory remains absolutely continuous
and the state equations hold almost everywhere; no differentiability
of the control at the cut is required. Thus it is a global minimizer
of the original free-terminal-time endpoint problem, not a problem
with a newly fixed value of $T$. The original maximum principle
therefore supplies its own full nonzero maximizing lift with the
zero free-time Hamiltonian. The ensuing five entire maximal
arcs are now contained in one genuine fundamental interval because
$5\le B$. When $B=5$ they fill the switch-to-switch period; no positive
uncharged interval is invented. Its four source ties, between successive
arcs of $10210$, are still strictly inside $(0,T)$, so all cap equations
use interior switch ties rather than any identification of endpoint
covectors. Theorem~\ref{wk:five600} requires five
complete positive arcs in one genuine fundamental interval, so includes
this equality case: neither its statement nor its proof assumes $5<B$.
Choose its own nonzero complete maximizing lift on this reanchored
original minimizing interval. Lemma~\ref{p32:five-normal} proves
normality before joint normalization; the proper-five proof establishes
its own positive determinant and future trace before harmonic calibration.
It gives $c_1+c_2+c_3\ge3/5$.

Reanchoring and whole lawful transport preserve complete physical costs,
the original state/group endpoint problem and real winding. Only these
scalar costs are identified with those of the initial interval. In
particular, no Casimir or state covector is continued across its cut.
\end{proof}

\subsubsection{A complete reverse four-core}
\label{strip:reverse-four}
\begin{proposition}[Reverse four-core]\label{strip:r496}
In an eligible full positive reverse context $102102$ or an additive
phase, on a genuine containing full fundamental interval with the one
complete lift and cap-derived normality specified in
Section~\ref{strip:reverse-kernel}, the four complete interior
arcs have physical cost
\[
 A=c_1+c_2+c_3+c_4>3/5>62/125.
\]
Both exterior caps remain positive.
\end{proposition}
\begin{proof}
The actual window $0,\ldots,4$ is $10210$ and retains the complete
charged arcs $1,2,3$. Lemma~\ref{strip:five-cap-completion} gives
$c_1+c_2+c_3\ge3/5$, while $c_4>0$.
\end{proof}

\subsubsection{Six charged arcs in a reverse eight-context}
\label{strip:reverse-six}
\begin{proposition}[Reverse six-core]\label{strip:r600}
On a full positive reverse eight-context $10210210$, or an additive
phase, charge exactly the complete arcs $1$ through $6$ and retain
positive caps $0,7$ in one genuine fundamental minimizing interval.
Assume also the separate eligibility of its contained reverse
six-context on $2,\ldots,7$. Then
\[
 A_6=c_1+\cdots+c_6\ge6/5>3/5.
\]
\end{proposition}
\begin{proof}
The windows $0,\ldots,4$ and $3,\ldots,7$ both read $10210$.
Their complete charged triples are $1,2,3$ and $4,5,6$.
Apply Lemma~\ref{strip:five-cap-completion} to both and add their
scalar inequalities. These two charge sets are disjoint; neither
choice of containing lift is imposed on the other.
\end{proof}

\subsubsection{A weighted reverse four-core}
\label{strip:weighted-reverse}
\begin{proposition}[Half-port reverse bound]\label{strip:weighted308}
In an eligible full reverse six-context $102102$, retaining positive
caps $0,5$ and its actual selected-six premise, the four complete
interior costs satisfy
\[
 W=\tfrac12c_1+c_2+c_3+\tfrac12c_4\ge3/5>77/250.
\]
\end{proposition}
\begin{proof}
Use the actual windows $0,\ldots,4=10210$ and
$1,\ldots,5=02102$. A whole phase $+1$ takes the latter to $10210$.
Lemma~\ref{strip:five-cap-completion} gives
$c_1+c_2+c_3\ge3/5$ and $c_2+c_3+c_4\ge3/5$.
Half the sum is exactly $W$, with weights
$(1/2,1,1,1/2)$ on its four complete physical charges.
\end{proof}

\subsubsection{A short--reverse--short ten-arc union}
\label{strip:j-union-section}
\begin{proposition}[Ten-arc union]\label{strip:j1010}
On the full positive context $12101021020212$ or an additive phase,
charge the ten complete arcs $2,\ldots,11$ once each and retain
returns $0,13$. Retain the actual selected-six premise on the
contained reverse context $4,\ldots,9$. Then
\[
 J=c_2+\cdots+c_{11}>1073/1000>101/100.
\]
\end{proposition}
\begin{proof}
The actual window $4,\ldots,8=10210$ has complete charged mask
$\{5,6,7\}$, with cost at least $3/5$ by
Lemma~\ref{strip:five-cap-completion}. The qualified one-return
source $7,\ldots,13=1020212$ keeps its actual caps $7,12$ and
return $13$; its complete charged mask $\{8,9,10,11\}$ has cost
strictly greater than $473/1000$ by Corollary~\ref{loc:E473}.
Both unit-weight masks are disjoint and contained in $\{2,\ldots,11\}$.
The other charged arcs have positive cost, so
$J>3/5+473/1000=1073/1000$. Only the physical cost inequalities are
added; no cap calculation uses an adjoint selected on a different cut.
\end{proof}

\subsection{Uniform state bounds from the critical hexagon}\label{p32:state}
The bounds in this section concern the actual geometric trajectory, before a
costate is selected. They therefore do not require a sign for its determinant.
Put
\[
 \begin{gathered}
 f_* =2\Cstar-3=\frac{30483}{50000},\qquad
 \delta=\frac{388657}{5000000},\\
 M=\frac{739}{200},\quad
 p_* =\frac{1767}{100000},\qquad [r]_+=\max(r,0).
 \end{gathered}
\]
\begin{lemma}[The area constraint]\label{p32:area}
At every actual state of a normalized global minimizer with physical cost
$C<\Cstar$, its three reciprocal coordinates satisfy
\begin{equation}\label{p32:area-eq}
 x,y,z>0,\qquad xy+yz+zx=1,\qquad
 F(x,y,z):=\sum_{\rm cyc}[xy-2z\sqrt{xy}]_+<f_*.
\end{equation}
\end{lemma}
\begin{proof}
The reciprocal identity follows from the state coordinates. The
critical-hexagon area inequality of Hales--Vajjha
\cite[Theorem 5.0.6 and Lemma 5.0.4]{HV2024}, in the normalization used here, is
\[
 |D|\ge\frac{3\sqrt3}{2}
       +2\sum_i[|T_i|-\sqrt{|T_i||T_{i,\mathrm{ext}}|}]_+.
\]
The triangle areas are $\sqrt3xy/4$ and $\sqrt3z^2$, cyclically.
The reciprocal order differs only by a permutation; the half-plane
coordinates in that source are $(u,v)/\sqrt3$. Since $|D|=\sqrt3C$,
division by $\sqrt3$ gives $C\ge3/2+F/2$, proving the assertion.
\end{proof}
\begin{proposition}\label{p32:state-bounds}
Every triple satisfying \eqref{p32:area-eq} obeys, cyclically,
\begin{equation}\label{p32:state-list}
 \delta<x<M,\quad xy>p_*,\quad \frac2{x+y}<\frac{301}{40},
 \quad \frac1x+\frac1y<\frac{18183}{1000}.
\end{equation}
Moreover
\begin{equation}\label{p32:state-vector}
 x+y+z<\frac{397}{100},\qquad
 \left|\left(\frac{x-z}{x+z},y-\frac{x^2+z^2}{x+z}\right)\right|
 <\frac{179}{50}.
\end{equation}
\end{proposition}
\begin{proof}
Order $x=m\le y\le z$. For the lower bound on $m$ only $m\le1/3$
needs consideration. Set $r=yz,s=xz,t=xy$. Then $r\ge s\ge t$,
$r+s+t=1$, $st=m^2r$, and
\[
 r=\frac{t(1-t)}{t+m^2},\qquad s=\frac{m^2(1-t)}{t+m^2},\qquad
 m^2\le t\le m(\sqrt{1+m^2}-m).
\]
This interval is precisely $m\le y\le z$. On it
\[
 r'=\frac{m^2(1-2t)-t^2}{(t+m^2)^2}\ge0,\quad
 s'=-\frac{m^2(1+m^2)}{(t+m^2)^2}<0,\quad t\le\frac13.
\]
The smallest positive-part argument is negative. The largest,
$A=r-2\sqrt{st}=r-2m\sqrt r$, is nonnegative and increasing, since
$r\ge(1-m^2)/2\ge4m^2$. The remaining argument
$B=s-2\sqrt{rt}$ decreases from a nonnegative to a negative value.
Where $B\le0$, $F=A$ increases. Where $B\ge0$,
\[
 F=1-t-2\sqrt{Q(t)},\qquad
 Q(t)=t(1-t)+\frac{2mt^{3/2}(1-t)}{t+m^2},
\]
and
\[
 Q'(t)=1-2t+
 \frac{m\sqrt t\{t(1-3t)+m^2(3-5t)\}}{(t+m^2)^2}>0.
\]
Consequently the minimum occurs at $B=0$. Writing $a=\sqrt r$ there,
\[
 t=(1-a)^2,\quad s=2a(1-a),\quad
 m^2=\frac{2(1-a)^3}{a},\quad F_{\min}(m)=a^2-2ma,
 \quad \frac23\le a<1.
\]
The function $2(1-a)^3/a$ decreases. For $a_0=17247967/20000000$,
\begin{align*}
 2(1-a_0)^3-\delta^2a_0&=
 \frac{1156868807273}{4000000000000000000000}>0,\\
 a_0^2-2\delta a_0-f_*&=
 \frac{20750537}{400000000000000}>0.
\end{align*}
Thus $m\le\delta$ would imply $F\ge F_{\min}(m)>f_*$.

For $p=xy$, discard its positive-part term and bound each remaining
positive part by its argument. Cauchy's inequality gives
\[
 F\ge1-p-2\sqrt{2p(1-p)}=:g(p).
\]
The function $g$ decreases on $(0,1/2]$, and
\[
 1-p_*-f_*>0,\qquad
 (1-p_*-f_*)^2-8p_*(1-p_*)=\frac{207601}{10000000000}>0.
\]
Hence $p>p_*$. The arithmetic--geometric mean inequality now gives
\[
 \frac2{x+y}\le\frac1{\sqrt p}<\frac{301}{40},\qquad
 z\le\frac{1-p}{2\sqrt p}<\frac{1-p_*}{2\sqrt{p_*}}<M;
\]
the last comparisons are verified by squaring positive quantities.
For $R=1/\delta$, $H=301/40$, $r=1/x$, $s=1/y$, we have
$r,s<R$ and $2rs/(r+s)<H$. Expanding $(R-r)(R-s)>0$ yields
\[
 r+s<\frac{R^2}{R-H/2}
 =\frac{2000000000000000}{109995466641651}<\frac{18183}{1000}.
\]
If the largest reciprocal is at most one, their sum is at most three;
otherwise $x+y+z<z+z^{-1}<M+M^{-1}<397/100$.
For $\ell=x+y+z$ and the vector $n$ in \eqref{p32:state-vector}, expansion
using $xy+yz+zx=1$ gives $|n|^2=\ell^2-3$.
Finally $(397/100)^2-3<(179/50)^2$.
\end{proof}
\begin{remark}\label{p32:quadratic-area}
For $a=\sqrt{yz},b=\sqrt{zx},c=\sqrt{xy}$, the closed area constraint is
$a^2+b^2+c^2=1$ together with
\[
 a^2-2bc\le f_*,\qquad 1-a^2-2a(b+c)\le f_*,\qquad\text{cyclically}.
\]
Indeed a sum of three positive parts is the largest of its eight subset
sums. The smallest of the three arguments is negative, so the three-term
sum imposes no additional condition. These six quadratic inequalities
may be retained in an interval box instead of selecting positive-part cases.
\end{remark}

\subsection{A three-fifths estimate on a five-arc source}\label{p32:weak-five}
Only one complete exterior arc on each side is retained in the following
statement. Its support, not merely its numerical constant, is needed in
short twisted periods.
\begin{theorem}\label{p32:weak-five-thm}\label{wk:five600}
Suppose a normalized original global minimizer with $C<\Cstar$ has
five consecutive complete positive maximal arcs $1\mid021\mid0$ in one
genuine fundamental interval. The sum $I$ of the three interior physical
costs satisfies $I\ge3/5$. The conclusion is invariant under the lawful
symmetries of the full problem.
\end{theorem}
We give the analytic reductions and specify the finite residual witness.
In this section $\mathcal C=\det\Lambda$ is the matrix invariant, whereas
$I,J_1,J_2,J_3$ denote physical costs. At the first interior reverse switch
let the reciprocals be $(a,b,c)$, with $ab+ac+bc=1$. The three own
exponents are $P,F,G>1$: on an arc whose active reciprocal is $d$ and
whose physical duration is $t$, its own exponent is $\exp(t/d)$.

\subsubsection{Normality and actual exponent bounds}
\begin{lemma}\label{p32:five-normal}
A nonzero full maximizing lift on a genuine interval containing the
complete positive subcontext $1021$ is normal.
\end{lemma}
\begin{proof}
This is the $\delta=1$ case of Lemma~\ref{sc:pair}, on the very
containing interval whose full lift has been selected.
\end{proof}
Divide all multipliers of this selected lift by the same positive $\eta$.
Define
\[
 r(X)=\frac{X\log X}{X-1}-1,\qquad
 \rho(s)=\frac{(1+s)\log(1+s)-s}{s^2},
\]
with $r(1)=0$, $\rho(0)=1/2$. Then $r(1+s)=s\rho(s)$ and
\begin{equation}\label{p32:rho-integral}
 \rho(s)=\int_0^1\frac{1-z}{1+sz}\,dz,\quad
 0<\rho(s)\le\tfrac12,\quad
 0\le\tfrac12-\rho(s)\le s/6.
\end{equation}
\begin{lemma}\label{p32:exponent-bounds}
If $I<3/5$, then
\[
 P,F,G<14/5,\qquad PF,FG<31/10,\qquad PFG<196/25,
 \qquad P+F-2<21/10.
\]
\end{lemma}
\begin{proof}
The physical pair identity \eqref{sc:cost}, with $\delta=1$ and
$(x,y,z;P,R)=(a,b,c;P,F)$, gives $J_1+J_2$ with nonnegative residual
and positive credits. Here the function called $f$ there is $r$.
Let $f(z)=r(e^z)$ and
$g(z)=f(z)+p_*(1-e^{-z})$. Direct differentiation gives
\[
 f''(z)=\frac{e^z\{(z-2)e^z+z+2\}}{(e^z-1)^3}>0,
\]
since the numerator bracket and its first derivative vanish at zero and
its second derivative is $ze^z>0$. Thus $g'>1/2$. Also
\[
 1-f'(z)-\tfrac12e^{-z}
 =\frac{e^{-z}(2z-3+4e^{-z}-e^{-2z})}{2(1-e^{-z})^2}\ge0,
\]
so $g'<1$. Consequently $g(A+B/2)<g(A)+g(B)$ for $A\ge0,B>0$.
For exact numerical comparisons use $L_8$ from \eqref{hc:log-bounds}.
The resulting positive reserves in $g(\log(14/5))-3/5$ and
$2r(44/25)-3/5$ are, respectively,
\[
 \frac{39476340070128859199786571}{3039261631380934619459800000},
 \qquad
 \frac{95546363226776403357876125063}{5222386528361254699941587114385}.
\]
Apply \eqref{sc:cost} to both overlapping actual pairs. The first
comparison bounds each exponent. Convexity gives
$r(P)+r(F)\ge2r(\sqrt{PF})$; since $(44/25)^2<31/10$, this bounds
$PF$, and likewise $FG$. The inequality for $g$ gives
$P\sqrt F<14/5$ and $G\sqrt F<14/5$, proving the triple-product
bound. Finally $PF=1+(P+F-2)+(P-1)(F-1)$.
\end{proof}

\subsubsection{A rational description of the actual cap}
Put
\[
 u=\frac{a-c}{a+c},\quad x=ac,\quad q=\frac{x}{1-u^2},\quad
 A=1+u,\quad C_0=1-u,\quad B=1-x.
\]
Thus $a=A\sqrt q$, $c=C_0\sqrt q$, $b=B/(2\sqrt q)$.
Let $s=P-1$, $t=F-1$, $\varepsilon=s+t$,
$\lambda=s/\varepsilon$, $\bar\lambda=1-\lambda$,
$e=A\lambda+C_0\bar\lambda$, $\rho_s=\rho(s)$ and $\rho_t=\rho(t)$.
The following expressions do not divide by either individual duration:
\begin{align}\label{p32:cap-coeff}
 W={}&qe-B/2+(1+A^2q)\lambda(1/2-\rho_s)
                 +(1+C_0^2q)\bar\lambda(1/2-\rho_t),\notag\\
 v_0={}&2Aq\lambda-(1+A^2q)\lambda\rho_s-C_0\bar\lambda W/e,\notag\\
 v_1={}&\frac{\lambda(A^2q-1)+\bar\lambda(1-C_0^2q)}{2e},\quad
 h_0=\frac{x\lambda\bar\lambda W}{e},\quad
 h_1=-\frac{qB\lambda\bar\lambda}{e}.
\end{align}
Write the nonnegative reverse crossing response as $Y=\varepsilon H$.
The original cap and the unused endpoint Hamiltonians imply
\begin{equation}\label{p32:H-interval}
 W\ge0,\quad
 \max\left(0,\frac{A\lambda(2\rho_s-C_0)}P,
                 \frac{C_0\bar\lambda(2\rho_t-A)}F\right)
 \le H\le\frac{AC_0W}{B}.
\end{equation}
The full cap vector and binary constant are
\[
 V_x=\varepsilon(v_0+v_1H),\quad V_y=\varepsilon\sqrt q H,\quad
 h=\varepsilon^2(h_0+h_1H),\quad K=(3+u^2-\varepsilon H)/2.
\]
The determinant is exactly
\begin{align}\label{p32:det-cap}
 \mathcal C(H)={}&\frac34-\frac\varepsilon2
       \left\{uv_0+H\left(\frac{B}{2e}-q\right)\right\}\notag\\
 &-\frac{\varepsilon^2}{4}\{(v_0+v_1H)^2+qH^2\}
 -\frac{\varepsilon^2(3+u^2)}2(h_0+h_1H)
 +\frac{\varepsilon^3H}2(h_0+h_1H).
\end{align}
Its quadratic coefficient is
$-\varepsilon^2(v_1^2+q)/4+\varepsilon^3h_1/2<0$.
Thus its minimum on any retained closed cap interval is at an endpoint,
without any assumption on the signs of $K$ or $\mathcal C$.

Here is a direct derivation connecting these expressions to the original
lift. Set $E=As+C_0t$ and
$W_{\max}=q\{2E-A^2r(P)-C_0^2r(F)\}-r(P)-r(F)$.
The two exterior ties give $0\le w\le W_{\max}$ and
\[
 Y=AC_0(W_{\max}-w)/B,\qquad
 V_x=\frac{(A^2q-1)Y}{2A}+2Aqs-(1+A^2q)r(P)-\frac{wC_0t}{E},
\]
\[
 K=(3+u^2-Y)/2,\quad L=(V_x+u)/2,\quad
 m=-\tfrac12+\frac{wAC_0qst}{E}-\frac{Yq}{2}.
\]
Substitute $Y=\varepsilon H$ and $\mathcal C=-L^2-Km$.
The unused endpoint Hamiltonians, divided only by positive state factors,
are the two lower bounds in \eqref{p32:H-interval}. Thus this interval
contains the cap of the actual selected lift; it is not a comparison
matrix to which calibration is independently applied.
The first two costs are
\begin{align}\label{p32:costs12}
 J_1&=\frac{(A^2q+AB/2+x)s+\{x+A^2q(1-\rho_s)\}s^2}{P},\notag\\
 J_2&=\frac{(C_0^2q+C_0B/2+x)t+\{x+C_0^2q(1-\rho_t)\}t^2}{F}.
\end{align}
At the next switch put
\[
 \begin{gathered}
 E_0=A+2t,\quad D_0=B-2C_0qt>0,\quad N_0=1+E_0^2q,\quad
 \\
 d=\frac{D_0}{2F\sqrt q},\quad e_0=E_0\sqrt q,\quad h_n=\frac{2D_0}{N_0}.
 \end{gathered}
\]
The next reciprocal triple is $(e_0,d,c)$, and the terminal triple is,
up to a cycle,
\begin{equation}\label{p32:terminal-state}
 X=(d+e_0)/G-d,\qquad d,\qquad Z=G(d+c)-d.
\end{equation}
The third cost and its derivative are
\begin{equation}\label{p32:cost3}
 J_3=d\{(d+e_0)(1-G^{-1})+(d+c)(G-1)-d\log G\},\quad
 J_3'=d\{(d+e_0)G^{-2}+d+c-dG^{-1}\}>0.
\end{equation}
The last sign follows from $G^2-G+1>0$.

\subsubsection{The double-collapse boundary}
\begin{lemma}\label{p32:collapse}
Under the original cap and endpoint maximizing inequalities,
$\varepsilon\le1/113$ implies $\mathcal C>\Hcap$.
\end{lemma}
\begin{proof}
Write $\bar e=\lambda a+\bar\lambda c$, $\ell=a+b+c$,
$\alpha=1/2-\rho_s$, $\beta=1/2-\rho_t$,
$\Delta=\lambda(a+b)\alpha+\bar\lambda(b+c)\beta$, and
$D_b=(\bar e-b)/2+\Delta$. The cap residual is $\varepsilon D_b$;
write its parameter as $\varepsilon z$, $0\le z\le D_b$. The cap
formula becomes $V=\varepsilon\bar V$, with
\[
 \bar V_y=\frac{2ac(D_b-z)}{b(a+c)},\quad
 \bar V_x=\frac{a^2-1}{2a}\bar V_y
 +(a+c)\lambda\{a-(a+b)\rho_s\}
 -\frac{(a+c)c\bar\lambda z}{\bar e},
\]
\[
 h=\frac{\varepsilon^2ac(a+c)\lambda\bar\lambda z}{\bar e}.
\]
Put $D=\delta$, $S=397/100$, $N=179/50$. For $E\ge\varepsilon$ set
\[
 B_0=(M-D)/2+ES/6,\quad R_0=1/2+E/6,\quad
 B_{\max}=\max(D+D^{-1},M+M^{-1}).
\]
Then $0\le D_b\le B_0$. For
$n=((a-c)/(a+c),b-(a^2+c^2)/(a+c))$ expansion gives
\[
 n\cdot\bar V=n_0-\frac{a(c-b)(a+b)\lambda\alpha}{b}
                   -\frac{c(a-b)(b+c)\bar\lambda\beta}{b},\quad
 n_0=-\frac{(a-b)(c-b)\bar e}{2b}
                    +\frac{zac(\bar e-b)}{b\bar e}.
\]
The outgoing endpoint inequality implies, for
$g_F=(a-b)\bar e/2-az$,
\[
 g_F\ge-\varepsilon B_1,\qquad B_1=aB_0+(a+b)S/6.
\]
Indeed its exact left side is
$g_F+ta(D_b-z)+a\Delta+b\bar\lambda(a+c)\beta$.
Interchanging $a,c$ and $\lambda,\bar\lambda$ gives the incoming inequality.
The endpoint defects give the uniform bound
\begin{equation}\label{p32:normal-vector-bound}
 n\cdot\bar V\le\varepsilon C_N,\qquad
 C_N=\frac{B_0+S/6}{D}+\frac{MS}{6}.
\end{equation}
Indeed, interchange the sides so $a\le c$. The identity
\[
 n_0=-\frac{c-b}{b}g_F-\frac{\lambda a(c-a)z}{\bar e}
 \le\frac{\varepsilon B_1(c-b)_+}{b}
\]
uses $g_F\ge-\varepsilon B_1$ and $g_F\le0$ when $c<b$.
For $c\ge b$, normalization gives
\[
 \frac{B_1(c-b)}b
 =B_0(b^{-1}-2a-c)+\frac S6(b^{-1}-2a-b)
 \le\frac{B_0+S/6}{D}.
\]
The two response corrections are at most $\varepsilon MS\lambda^2/6$
and $\varepsilon MS\bar\lambda^2/6$, respectively: use
$(b-c)_+,(b-a)_+\le b$, $a(a+b),c(b+c)\le MS$.
Their sum is at most $\varepsilon MS/6$.

At $z=0$ the original cap vector has the representation
\[
 V=(b+b^{-1})(-\beta_a\nu_c+\beta_c\nu_a),\quad
 \beta_a=as\{a/(a+b)-\rho_s\},\quad
 \beta_c=ct\{c/(b+c)-\rho_t\},
\]
where $\nu_a=((a^2-1),2a)/(1+a^2)$ and
$\nu_c=((c^2-1),-2c)/(1+c^2)$ are unit vectors.
Thus $|V|\le\varepsilon B_{\max}MR_0$ there.
At $z=D_b$, $V_y=0$ and $|V|\le\varepsilon S^2R_0$.
Affine dependence on $z$ and $\bar e^2\ge4ac\lambda\bar\lambda$ give
\[
 |V|\le\varepsilon W_0,\quad W_0=\max(B_{\max}M,S^2)R_0,
 \qquad h\le\varepsilon^2SMB_0/4.
\]
The exact invariant identity is
\[
 \mathcal C=\frac34-\frac{n\cdot V}{2}
 -\frac{(\ell-n_y)h}{2}-\frac{|V|^2}{4}+\frac{hV_y}{2}.
\]
Its last term is nonnegative. Hence
\[
 \mathcal C\ge\frac34-\varepsilon^2K(E),\qquad
 K(E)=C_N/2+W_0^2/4+(S+N)SMB_0/8.
\]
Rational substitution gives
$K(1/113)<186053/1000$ and
$3/4-186053/(1000\cdot113^2)>\Hcap$.
These are integer cross-multiplications, completing the proof.
\end{proof}

\subsubsection{Full reflection and terminal information}
The representative may be chosen with $P\le G$ by the full lawful
reflection of Proposition~\ref{fd:lawful-symmetry}, followed by two cyclic
phases. With $\psi(u,v)=(-u,v)$, put $T=\Phi^2\psi$ and $Q=R_*^{-2}\operatorname{diag}(1,-1)$:
\[
 Q=\tfrac12\begin{pmatrix}-1&-3\\-1&1\end{pmatrix},\qquad
 T(u,v)=\left(\frac{u^2+2u+v^2-3}{(u-1)^2+v^2},
               \frac{4v}{(u-1)^2+v^2}\right).
\]
Apply \eqref{fd:reflection} to the whole containing interval and then
these phases: at reversed time the matrix is $-Q\Lambda Q^{-1}$,
while both state-covector entries transform as $-DT^{-t}p$.
The frame is reanchored as in that formula, and $\eta$, maximizing
Hamiltonians, density, own reciprocals, original endpoints and real
winding are preserved. All caps and charged arcs move with the source;
no costate is periodically extended. The label map is $j\mapsto1-j$,
which sends the reversed word $10210$ to itself and exchanges $P,G$.

Let $K_n$ be the transported binary constant on the third charged arc,
$\mu_n=2K_n-3$. Its own mode-zero coordinates run from $h_n-1$ to
$h_nG-1$. Integrating the nonnegative switching function, which vanishes
at the actual final tie, gives
\begin{equation}\label{p32:terminal-mu}
 \mu_n\le M(G):=1+h_n^2G-2h_n\{1+r(G)\}.
\end{equation}
Since $r''<0$, $M$ is convex. For any justified $G\le g$, the chosen
$P\le G$ therefore gives
\begin{equation}\label{p32:terminal-convex}
 M(1)-\mu_n+\max\{M(P)-M(1),M(g)-M(1)\}\ge0.
\end{equation}
The actual incoming response supplies a stronger affine restriction.
Exact transport gives
\[
 E_F=\varepsilon\{FH+C_0\bar\lambda(A-2\rho_t)\}\ge0,
 \quad \chi_n=h_nqtE_F,\quad
 \chi_G=\chi_n-\frac{d^2(G-1)}{h_nG}\{M(G)-\mu_n\}.
\]
Thus the actual final tie is
\begin{equation}\label{p32:terminal-tie}
 M(G)-\mu_n=k(G)E_F,\qquad k(G)=\frac{16F^2q^2tG}{N_0^2(G-1)}.
\end{equation}
These formulas use both the shear
$L_2=L+KC_0qt$, $m_2=m-2LC_0qt-K(C_0qt)^2$
and the normally forced state-adjoint equation.
For stable evaluation define
\[
 T_s=q(A+4t+2t^2),\quad
 B_s=1-qA(4-A)+t\{2+2q(A^2-6)+tq(8A-12+2At)\},
\]
\[
 \delta_0=\frac{8qC_0\bar\lambda B_s-8v_0T_s+8\varepsilon h_0}{N_0^2},
 \qquad
 \delta_1=\frac{4T_s^2-8v_1T_s-4q+8\varepsilon h_1}{N_0^2}.
\]
Then $M(1)-\mu_n=\varepsilon(\delta_0+\delta_1H)$.
If $U$ encloses the maximum increment in \eqref{p32:terminal-convex},
the necessary affine cut is $\delta_0+\delta_1H+U/\varepsilon\ge0$.
For $g>1$, any nonnegative lower bound $k_l\le k(g)$ strengthens it to
\begin{equation}\label{p32:terminal-cut}
 \delta_0+\delta_1H+U/\varepsilon
       -k_l\{FH+C_0\bar\lambda(A-2\rho_t)\}\ge0.
\end{equation}
This uses the actual $E_F\ge0$. If $g\le1$ the divided cut is omitted,
retaining the collapsed boundary.
Interior maximization also gives, for
$l=\min(u,As/P)$ when $u\ge0$ and
$l=\min(-u,C_0t/F)$ when $u\le0$,
\begin{equation}\label{p32:interior-cut}
 \varepsilon H\ge|u|l-l^2/3.
\end{equation}
Indeed integrate $\chi_\sigma=\mu-u^2$ up to the actual binary zero.
On the nonnegative interval the square increases while its own weight
$(1+u)^{-2}$ decreases, so the weighted mean is bounded by the ordinary
mean $u^2-ul+l^2/3$. The negative outgoing interval gives the reflected
inequality. The case $l=0$ uses no division.
The weaker consequences $H\ge2u^2\lambda/(3P)$ for $u\ge0$ and
$H\ge2u^2\bar\lambda/(3F)$ for $u\le0$ are also valid.

\subsubsection{The compact residual and its finite witness}
\begin{lemma}\label{p32:weak-residual}
For a normal actual prefix with $I<3/5$ and $\varepsilon\ge1/113$,
the preceding necessary conditions imply $\mathcal C>\Hcap$.
\end{lemma}
\begin{proof}
We specify the finite calculation, including its domain and checking rule.
The closed initial box, in the order $(u,x,\varepsilon,\lambda)$, is
\begin{equation}\label{p32:root-box}
 [-u_*,u_*]\times[p_*,1-2p_*]\times[1/113,21/10]\times[0,1],
 \qquad u_*=(M-\delta)/(M+\delta).
\end{equation}
The state bounds and Lemma~\ref{p32:exponent-bounds} put every actual
contrary point in this box. The faces $\lambda=0,1$ are outer relaxations,
not positive-arc assumptions. The cap variable is handled by affine
restrictions and concavity; the final exponent is projected by necessary
inequalities, not chosen as an independent root.

\paragraph{State and pair-cost restrictions.}
Test \eqref{p32:area-eq}--\eqref{p32:state-list} at the common state
and at its two actual neighbors. Their triples are
\[
 \left(a,\frac{B-2Aqs}{2P\sqrt q},(2P-A)\sqrt q\right),
 \qquad(e_0,d,c).
\]
In particular,
\[
 s\le\min\left\{\frac{b-\delta}{a+\delta},
                    \frac{M-c}{a+c},\frac95\right\},\quad
 t\le\min\left\{\frac{b-\delta}{c+\delta},
                    \frac{M-a}{a+c},\frac95\right\}.
\]
Intersect these with $s=\varepsilon\lambda$,
$t=\varepsilon(1-\lambda)$. The exact costs give
$J_1+J_2\ge\ell\varepsilon\bar e/(1+\varepsilon)$,
$\bar e=\sqrt q e$, and
\[
 J_1+J_2=j_0+xj_x,\quad j_0=As/(2P)+C_0t/(2F),
\]
\[
 j_x=\frac{3+u^2}{2}\left(\frac{s}{C_0P}+\frac{t}{AF}\right)
 +\left(1+\frac{A(1-\rho_s)}{C_0}\right)\frac{s^2}{P}
 +\left(1+\frac{C_0(1-\rho_t)}A\right)\frac{t^2}{F}.
\]
Here $j_x>0$ for $s+t>0$, so $x<(3/5-j_0)/j_x$.
Endpoint lower bounds are
\[
 x\ge\frac{A-R_P}{3+u},\qquad x\ge\frac{C_0-R_F}{3-u},
\]
where
\[
 R=(1+A^2q)\lambda\alpha+(1+C_0^2q)\bar\lambda\beta,
\]
\[
 R_P=\frac{2AC_0(R+sW)+4BA\lambda\alpha}{e},\qquad
 R_F=\frac{2AC_0(R+tW)+4BC_0\bar\lambda\beta}{e}.
\]
They are the endpoint inequalities evaluated at the largest allowed $H$,
where their left sides are maximal. Every divisor is positive.

\paragraph{Projection of the terminal exponent.}
Let $R=3/5-J_1-J_2$. The product and terminal-state conditions give
\[
 G\le\min\left\{\frac{31}{10F},\frac{196}{25PF},\frac{14}{5},
 \frac{d+e_0}{d+\delta},\frac{M+d}{d+c},
 \frac{d+e_0}{d+p_*/d},\frac{301(d+e_0)}{80}\right\}.
\]
The remaining product $XZ\ge p_*$ gives
\[
 d(d+c)G^2-\{(d+e_0)(d+c)+d^2-p_*\}G+d(d+e_0)\le0.
\]
A negative discriminant proves emptiness; otherwise its larger root is
an upper bound, with no selection of a physical branch.
With $y=G-1$, \eqref{p32:cost3} also gives
\[
 J_3=\frac{d(d+e_0+c)y+\{dc+d^2(1-\rho(y))\}y^2}{G}
 \ge d(d/2+c)y+d(d/2+e_0)y/G.
\]
Thus an enclosing $g\ge G$ may be reduced to
\[
 \min\left(g,1+\frac{R}{d(d/2+c)+d(d/2+e_0)/g}\right).
\]
Four outward iterations are used. The inequality $r(F)+r(G)<3/5$
and the lower polynomial
\[
 r(G)\ge z+(1+z)(z^2/3+z^4/5),\quad z=(G-1)/(G+1),
\]
give a further enclosing upper bound by eight bisections, as interval
restrictions, not approximate roots. The chosen representative still
requires $P\le g$ and $J_3(P)\le R$.

\paragraph{Cap restrictions and determinant evaluation.}
Retain \eqref{p32:H-interval}, \eqref{p32:interior-cut},
\eqref{p32:terminal-convex}, and \eqref{p32:terminal-cut}.
For an affine restriction $c_0+c_1H\ge0$, divide only when the entire
interval for $c_1$ has one strict sign. If it straddles zero, either the
whole expression has negative upper bound, proving emptiness, or no
division is made. If $U<0$, evaluate $U/\varepsilon$ on the entire positive
$\varepsilon$ interval: replacing its denominator by the lower endpoint
would reverse the needed bound. On a nonempty cap interval evaluate
\eqref{p32:det-cap} at both endpoints. Its strict concavity then bounds
every actual cap. Intersecting algebraically identical factorizations
of $\delta_0,\delta_1$ can sharpen the enclosure but adds no premise.

\paragraph{Coverage and outward arithmetic.}
Apply the exact closed-tree test of Section~\ref{cv:cover} to the
proper-five witness on \eqref{p32:root-box}. Each path records a split
coordinate and a lower/upper bit at its rational midpoint; all terminals
are assigned. There are $240\,722$ leaves, $240\,721$ internal nodes
and maximum depth $35$.

Apply Lemma~\ref{cv:outward} on the $2^{-112}$ grid, including outward
integer square roots. Division requires a strictly positive denominator
interval after any sign reduction.
For $z=s/(2+s)\in[0,1)$ use the convergent identity
\begin{equation}\label{p32:rho-series}
 \rho(s)=\frac12\left(1-z+(1-z^2)
             \sum_{k=1}^{\infty}\frac{z^{2k-1}}{2k+1}\right).
\end{equation}
The first $36$ terms give a lower bound; the omitted contribution is
at most $z^{73}/150$, by replacing all subsequent denominators by $75$
and summing the geometric series. Since $\rho$ decreases, endpoint
evaluations enclose a whole interval.

For every reconstructed leaf, fresh evaluation of the preceding formulas
proves either impossibility of one necessary condition or a determinant
lower bound strictly greater than $\Hcap$ at both retained cap endpoints.
The witness supplies the complete tree and fresh leaf evaluations;
stored acceptance reasons, samples or totals do not replace these
finite computer-assisted inequalities. Every restriction retains actual
contrary points, so Theorem~\ref{cv:finite-rejection} proves the lemma.
\end{proof}

\begin{proof}[Proof of Theorem~\ref{p32:weak-five-thm}]
Select a nonzero full maximizing lift on the specified genuine interval.
Lemma~\ref{p32:five-normal} makes it normal; normalize jointly.
If $I<3/5$, the exponent, collapse and residual lemmas give
$\mathcal C>\Hcap>0$ on this same lift.
It remains to verify future trace before applying harmonic calibration.
For the unscaled cap parameter $\beta$ the original first-switch formulas
are
\[
 0\le\beta\le D\le as+ct,\quad
 h=\frac{ac(a+c)st\beta}{as+ct},\quad
 \operatorname{tr}(\Lambda X)=\ell-h.
\]
Moreover $J_1>acs$, $J_2>act$; for example
\[
 J_1-acs=\frac{a(a+b)s+a^2s^2(1-\rho_s)}{1+s}>0.
\]
Consequently
\[
 \frac h\ell\le acst<I\frac{st}{s+t}
 \le I\varepsilon/4<\frac{63}{200},
 \qquad \frac{\operatorname{tr}(\Lambda X)}\ell>\frac{137}{200}.
\]
No determinant sign entered this trace calculation. Harmonic calibration
now applies to the selected full lift, with its original group endpoint
and real winding, and gives $\mathcal C<\Hcap$, a contradiction.
Thus $I\ge3/5$. In particular no strict $3/5$ claim is extracted from
an argument that only excluded $I<3/5$.
\end{proof}

\space
\section{Returned tails, shared arcs, and common-derivative unions}\label{lp:returns}
A returned arc supplies a geometric circle or a switching sign at one actual endpoint. We retain these facts while combining tails, shared four-cores and common-derivative unions. In particular, two overlapping core costs are never added without a once-counted union argument.

\subsection{A signed root-mean comparison}
\begin{lemma}\label{fin:signed-mean}
If $-1<w<s<1$ and the $(1+u)^{-2}$-weighted mean $\mu$ of $u^2$ on $[w,s]$ satisfies $\mu>w^2$, then
\[
\mu<\frac{w^2+ws+s^2}{3}.
\]
\end{lemma}
\begin{proof}
The strict excess implies $-\sqrt\mu<w<\sqrt\mu$, and the defining mean forces $s>\sqrt\mu$. The function $u^2-\mu$ has just one sign change on the interval. Since $q(u)=(1+u)^{-2}$ strictly decreases,
$q(u)(u^2-\mu)<q(\sqrt\mu)(u^2-\mu)$ away from the zero. Integrating gives $\int_w^s(u^2-\mu)du>0$, which is the assertion.
\end{proof}
This argument does not assume that $u^2$ is monotone on the interval. For an actual $2\to0$ root the required excess is strict: $\chi_u=(\mu-u^2)/v^2$ and $\chi$ changes from nonpositive to nonnegative. If equality held at $w\ne0$, the common one-sided second derivative $-2w/v^2$ would produce the same sign on both sides. If $w=\mu=0$, $\chi_u<0$ off the root. Both are impossible.

\subsubsection{The adjacent-pair row}
Consider seven complete positive arcs $h_0,\ldots,h_6$ with word
$0201020$ on one full normal maximizing lift, normalized jointly on an eligible
genuine complete interval. Their first ordinary $F$ core
is $h_1,h_2,h_3$. The right analytical triple is the cyclic image of
$h_2,h_3,h_4=010$, namely $202$, with positive caps $h_1,h_5$.
The left analytical mean comes from the two roots of the actual leaf
$h_5$, supplied by $h_4,h_6$; no cap after $h_6$ is required.
All chart transformations carry the same complete lift.

The shared positive arc is $h_4$. Its entrance $h_3|h_4$ supplies the
$1\to0$ cap equation and its exit $h_4|h_5$ the binary tie.
The outside, right and left means below belong respectively to
$h_1,h_3,h_5$, whose two roots are supplied by
$(h_0,h_2)$, $(h_2,h_4)$ and $(h_4,h_6)$.
Thus every primitive, height join and conjugation below uses these
seven arcs. For the adjacent word $020102010$, apply this row to
positions $0,\ldots,6$ and, on the other side, to positions
$8,\ldots,2$ under full reflected time reversal with $j\mapsto-j$.
The latter has the unsorted parameters $(h,g,E)=(b_2,a_2,E_2)$.
On their shared actual arc $4$, write $\xi=1+u$, $v=k\xi$, $R=1+k^2$, and let $\xi_0<\xi_1$ be the endpoints. The right coordinates are
$u_R=1-4/(R\xi)$, $v_R=4k/(R\xi)$; the left coordinates are $u_L=1-\xi$, $v_L=k\xi$. Put
\[
p=\frac4{R\xi_0},\qquad q=\xi_1,\qquad t=Rpq/4>1,
\quad A=R/k^2.
\]
Actual root signs give $0<p,q<1$. The two middle means satisfy
$0<\mu_R<(1-p)^2$, $0<\mu_L<(1-q)^2$, with $K_R=(3+\mu_R)/2$, $K_L=(3+\mu_L)/2$.
The following transfer uses only the shared arc, its initial $1\to0$
cap, its terminal $0\to2$ root, and the actual adjoining middle
leaf. Thus it applies whenever those same endpoint data are retained.
The action and matrix transport in \eqref{cmp:action-transport} give
\begin{equation}\label{fin:shared-matrix}
\begin{split}
L&=B-\frac{K_L}{k^2\xi},\qquad
m=\bar m+\frac{2B}{k^2\xi}-\frac{K_L}{k^4\xi^2},
\\
K_R&=-\frac{4k^4}{R^2}\bar m,\\
B&=\frac12-\frac A2\log t
-\frac A4\left(\frac{1-\mu_R}{p}+\frac{1-\mu_L}{q}-p-q\right),\\
\mathcal C&=\frac{A^2(3+\mu_L)(3+\mu_R)}{16}-B^2.
\end{split}
\end{equation}
Apply \eqref{cmp:ray-primitive} on the actual shared mode-$0$ arc,
with $(\kappa,K,B,M)=(k,K_L,B,\bar m)$. The terminal tie $\chi(q)=0$
fixes $C$. At $\xi_0=4/(Rp)$, the actual entrance cap is
\[
 m(\xi_0)+\frac12=\frac{2K_L-3-(\xi_0-1)^2}{2k^2\xi_0^2}
 +\frac{R\chi(\xi_0)}{2k^2\xi_0}.
\]
Substituting the primitive and solving this row gives the displayed $B$.
Conjugating the reconstructed matrix by
$\frac12\left(\begin{smallmatrix}1&-3\\1&1\end{smallmatrix}\right)$ gives
$K_R=-4k^4\bar m/R^2$, and its determinant gives $\mathcal C$.
All three lines thus describe the same full matrix.

For the right $F$ side let its original unsorted parameters be $(h,g,E)$. The actual middle interval is $[1-g,1-p]$, so $0<p<g<2$. Matching its endpoint height gives
\[
j^2=\frac{k^2p^2(2-g)^2}{g^2(2-p)^2},\qquad
\nu=Eh=\frac4{g(1+j^2)},\quad D=4-\nu g,
\quad \Phi=\frac{2(g-p)(4-g-p)}{j^2g^2(2-p)^2}>0.
\]
The complete original matrix reconstruction at this adjacent $F$, or directly substitution of the first line of \eqref{fin:shared-matrix} in the conjugated invariant at the outer leaf, yields
\let\savedordinarylabel\label
\begin{equation}\label{fin:outside-row}
\begingroup\let\label\savedordinarylabel\label{ek:outside-row}\endgroup
\mu_{\rm out}=D^2\left(\frac{A^2(3+\mu_L)}{16}+B\Phi
+\frac{(3+\mu_R)\Phi^2}{4}\right)-3.
\end{equation}
For clarity, the entire propagation giving this row can be recorded in four short formulas. Set $K_R=(3+\mu_R)/2$, $K_L=(3+\mu_L)/2$. At the root $s=1-p$ in the right chart,
\[
L_s=\frac{K_R}{k^2p}-B,\qquad
m_s=\frac{2B}{k^2p}-\frac{K_R}{k^4p^2}-\frac{A^2K_L}{4}.
\]
Here the subscript indicates evaluation at the endpoint, not differentiation. The preceding middle mode-$0$ leaf has rate $jg/(2-g)$ and positive shear
\[
\tau=\frac{(g-p)(2-g)}{j^2g^2(2-p)}.
\]
Backward propagation to $w=1-g$, with $T_0=1/(k^2p)-\tau$, gives
\[
L_w=K_RT_0-B,\qquad m_w=2BT_0-K_RT_0^2-A^2K_L/4.
\]
At that point $v=jg$. The inverse cyclic conjugation yields the other invariant
\[
K_{\rm out}=\frac{4\{K_Rg^2-2L_wgv^2-m_wv^4\}}{(g^2+v^2)^2}.
\]
Substitute the preceding expressions and observe that
\[
\frac1{j^2g}-T_0=\Phi,\qquad D=\frac{4j^2}{1+j^2}.
\]
Then $2K_{\rm out}-3$ is exactly the right side of \eqref{fin:outside-row}. This gives the common-matrix row without an unproved compatibility assertion.
The outer root supplies $\mu_{\rm out}>(1-h)^2$ and $h<\nu$. The left side has the separately transported analogous formulas, with $(h,g,E)=(b_2,a_2,E_2)$ rather than scalar-sorting the occurrences.

\begin{lemma}[An outside-mean comparison]\label{fin:H150}
On an actual seven-arc source just described,
\[
 g<3/2,\qquad \nu=Eh>\frac{3g}{2g+p}>\max\{1,g\}.
\]
In particular these conclusions hold on either lawful analytical side
of the adjacent context.
\end{lemma}
\begin{proof}
Set $d=g-p>0$ and
\[
 u_R=(1-p)^2-(1-p)d+d^2/3.
\]
The actual root excess and Lemma~\ref{fin:signed-mean} give
$(1-g)^2<\mu_R<u_R$, including a negative entrance. If $g>1$,
\[
 0<u_R-(1-g)^2=d(3-2g-p)/3,
\]
so $2g+p<3$ and $g<3/2$; the latter is immediate for $g\le1$.
Define $r=q/t$, so $0<r<q<1$ and $Rpr=4$.
Keep all actual state and matrix parameters fixed, and put
\[
 u_L=(1-r)^2,\qquad
 B_f=\frac12-\frac A4\left(\frac{1-u_R}{p}+2-2r-p\right),
\]
\[
 M_f=D^2\left(\frac{A^2(3+u_L)}{16}+B_f\Phi
             +\frac{(3+u_R)\Phi^2}{4}\right)-3.
\]
The actual left root gives $\mu_L<(1-q)^2<(1-r)^2$, and
\[
 \frac{4(B_f-B)}A=\frac{u_R-\mu_R}p+
 \frac{(1-q)^2-\mu_L}q+2(\log t-q+r)>0,
\]
since $q-r<1-1/t<\log t$.
All coefficients in \eqref{fin:outside-row} are positive, so
$\mu_{\rm out}<M_f$. No scalar comparison is assigned a physical lift.

Define
\begin{align*}
P_0={}&g^3+2g^2p-7g^2-4gp^2+2gp+12g+4p^3-7p^2,\\
P_1={}&-g^4-3g^3p+11g^3+17g^2p-40g^2+3gp^3-21gp^2\\
&+4gp+48g+p^4-p^3+12p^2-24p.
\end{align*}

Put $\mathfrak a=P_0/[3pg(2-p)^2]$,
$\mathfrak b=P_1/[3pg^2(2-p)^3]$ and $v_0=1-\nu^{-1}$.
Substitute $j^2=4/(\nu g)-1$,
$k^2=j^2g^2(2-p)^2/[p^2(2-g)^2]$ and
$r=4/[(1+k^2)p]$ in the displayed $M_f$. Exact cancellation gives
\let\savedordinarylabel\label
\begin{equation}\label{fin:factor}
\begingroup\let\label\savedordinarylabel\label{loc:H-quadratic}\endgroup
 \frac{M_f}{\nu^2}=v_0^2+2\mathfrak a d v_0-\mathfrak b d^2,
 \qquad
 \frac{M_f-(1-\nu)^2}{\nu^2}
       =d(2\mathfrak a v_0-\mathfrak b d).
\end{equation}
This rational identity uses $0<p<\min(g,1)$, $g<2$, $\nu>0$ and
$\nu g<4$, as supplied by the actual source; no zero denominator is
introduced at a scalar face.

The coefficient signs follow from two short polynomial comparisons.
First, $P_0>0$ on $0<p<1$, $p\le g\le3/2$: it is concave in $g$,
since $(P_0)_{gg}=6g+4p-14\le-1$. Its endpoint values are
\[
 P_0(p,p)=3p(p-2)^2>0,\qquad
 P_0(3/2,p)=45/8+15p/2-13p^2+4p^3.
\]
The latter is concave in $p$, with second derivative $24p-26\le-2$
and positive endpoint values $45/8,33/8$. Thus $\mathfrak a>0$.

Define
\[
 W=3P_1-2(2-p)P_0-p(2-p)^3.
\]
Then
\[
 \mathfrak b-\frac{2\mathfrak a}{3g}-\frac1{9g^2}
       =\frac{W}{9pg^2(2-p)^3}.
\]
On $0\le p\le1$, $p\le g\le3/2$,
\[
 W_{ggg}=174-72g-42p\ge24,\qquad
 W_{gg}(3/2,p)=-1-3(1-p)-8p(1-p)\le-1.
\]
Thus $W$ is concave in $g$. Its endpoint values satisfy
\[
 W(p,p)=p(2-p)^2(7p+4)>0\quad(p>0),
\]
\[
 16W(3/2,p)-293
  =(1-p)(22+80p+408p^2-192p^3)\ge0.
\]
Hence $W>0$ on the actual domain, and $\mathfrak b>0$ follows
without a separate sign argument for $P_1$. The closed polynomial
faces do not extend the rational identity to $p=0$.

If $\nu\le1$, the actual $h<\nu$ gives
$\mu_{\rm out}>(1-h)^2>(1-\nu)^2$, whereas $v_0\le0$ and
\eqref{fin:factor} give $M_f<(1-\nu)^2$. Thus $\nu>1$.
Now $v_0>0$, and $0<\mu_{\rm out}<M_f$ in the same quadratic gives
\[
 v_0>d c_*,\qquad c_*:=\sqrt{\mathfrak a^2+\mathfrak b}-\mathfrak a>0.
\]
For $c\ge0$, one has $c_*>c$ exactly when
$\mathfrak b>2\mathfrak a c+c^2$.
The $W$ comparison therefore gives $v_0>d/(3g)$.
Consequently $\nu>3g/(2g+p)>1$. If $g>1$, the actual strict
root correlation $2g+p<3$ makes this lower bound larger than $g$;
if $g\le1$, the conclusion $\nu>g$ is already immediate.
No own-HIGH or cost hypothesis, harmonic calibration, or boundary
lift is used.
\end{proof}

\subsection{Shared-arc unions of four-cores}
\label{loc:shared-A-section}
Put $C_0=\Hcap$ and $z_0=243/2000$. In particular
\begin{equation}\label{loc:common-z0}
 3/4-C_0-z_0^2
 =\frac{89765086136957463911195871}
 {39478417604337498218323684000000}>0.
\end{equation}

\begin{lemma}[One shared physical arc]\label{loc:shared-A-matrix}
Suppose the last mode-$2$ arc of a canonical capped $0202$ core is also
the first mode-$0$ arc of a second canonical binary core under the legal
cycle $\widehat w=(w+3)/(1-w)$.
Let its first and last coordinates in the left chart be $r,\beta$,
where $0<r<\beta$, and suppose its end is a negative first internal
root in the right chart. On the shared arc write
\[
 v=\kappa(1-u),\quad x=\kappa^{-2},\quad
 L=\lambda+Kx/(1-u),\quad
 m=M-2\lambda x/(1-u)-Kx^2/(1-u)^2.
\]
Then, on the same normalized complete lift,
\begin{gather}
 \widehat K=-4M/(1+x)^2,\qquad
 \mathcal C=K\widehat K(1+x)^2/4-\lambda^2,\label{loc:shared-A-C}\\
 \widehat u=-1+\frac{4x}{(1+x)(1-u)},\qquad
 \widehat v^2=(1+\widehat u)^2/x,\qquad
 \widehat L=-\lambda-\frac{\widehat Kx}{1+\widehat u},\label{loc:shared-A-transform}\\
 0<x<1/3,\qquad
 B<\mathcal R(x):=\frac{1-9x^2}{1+5x},\qquad\mathcal R'(x)<0,
 \label{loc:shared-A-cost}
\end{gather}
where $B$ is the cost of this complete shared arc.
If $T_L$ is the left exclusive three-arc shear, then, whenever
$\mathcal C<K/2$,
\begin{equation}\label{loc:shared-left-branch}
 \lambda< -\sqrt{K/2-\mathcal C}+KT_L-\frac{Kx}{1-r}.
\end{equation}
If the right core is also a capped four-core, its exclusive shear $T_R$
and positive first-root magnitude $\widehat\xi$ satisfy the independent
right-end inequality
\begin{equation}\label{loc:shared-right-branch}
 \lambda< -\sqrt{\widehat K/2-\mathcal C}
       +\widehat K T_R-\frac{\widehat Kx}{1-\widehat\xi}.
\end{equation}
\end{lemma}
\begin{proof}
On the shared arc $ds=x\,d((1-u)^{-1})$.
Apply the matrix transport in \eqref{cmp:action-transport} to obtain the displayed expressions and
$\mathcal C=-\lambda^2-KM$.
Reconstruct $l=L+um$, $b=K-2ul+u^2m$ and conjugate
$\left(\begin{smallmatrix}l&b\\m&-l\end{smallmatrix}\right)$ by
$\frac12\left(\begin{smallmatrix}1&3\\-1&1\end{smallmatrix}\right)$.
The definitions of $\widehat K,\widehat L$ then give
\eqref{loc:shared-A-C} and \eqref{loc:shared-A-transform} directly.
This is a full-matrix conjugacy on one physical arc; the shear coordinates
in its two charts have not been identified.

At its actual end $\widehat u<0$, so
$\beta<(1-3x)/(1+x)=:d_x$.
Since $\beta>r>0$, this proves $0<x<1/3$.
The exact physical cost, with no clock change, is
\[
 B=\frac12\int_r^\beta\left\{1+\frac{x(3+u^2)}{(1-u)^2}\right\}du
 <\frac12\int_0^{d_x}\left\{1+\frac{x(3+u^2)}{(1-u)^2}\right\}du
 =1-3x+x\log\frac{4x}{1+x}.
\]
An antiderivative is
$u/2+x[4/(1-u)+2\log(1-u)-(1-u)]/2$.
For $t=4x/(1+x)\in(0,1)$,
$\log t<2(t-1)/(t+1)$: their difference vanishes at $1$ and has
derivative $-(t-1)^2/[t(t+1)^2]<0$.
This gives \eqref{loc:shared-A-cost}, whose derivative is
$-(45x^2+18x+5)/(1+5x)^2<0$.
The scalar integral endpoints do not extend the physical trajectory.

At the actual left entrance, its negative cap jet and loss give
$L(a)<-\sqrt{K/2-\mathcal C}$.
Integrate $L_s=K$ along only its exclusive triple, and subtract
$Kx/(1-r)$ at the shared entrance, proving
\eqref{loc:shared-left-branch}.
At the right shared end $\widehat u=-\widehat\xi$,
\eqref{loc:shared-A-transform} gives
$\widehat L=-\lambda-\widehat Kx/(1-\widehat\xi)$.
Forward integration along the right exclusive triple reaches its actual
positive exit jet, where
$\widehat L>\sqrt{\widehat K/2-\mathcal C}$ because its exit cap
loss is positive. Rearrangement gives \eqref{loc:shared-right-branch}.
\end{proof}

\begin{theorem}[The generic adjacent four-core union]\label{loc:AA640}
For the actual positive context $021210102$ and its legal phases, the cost
$J$ of the seven complete arcs $1,\ldots,7$, counted once, exceeds $16/25$.
\end{theorem}
\begin{proof}
The cores $A_L=[1,4]$ and $A_R=[4,7]$ share only arc $4$.
The actual windows $[0,6]=0212101$ and $[2,8]=1210102$ are respectively
a right and a left one-return context, so Corollary~\ref{loc:E473} gives
$A_L,A_R>E=473/1000$.
Let $e_L=[1,3]$, $e_R=[5,7]$ denote their exclusive costs and $B=c_4$.
Thus
\[
 J=A_R+e_L=A_L+e_R,\qquad B=A_L+A_R-J.
\]
Suppose $J\leq U=16/25$. Then $e_L,e_R<e=167/1000$ and $B>2E-U$.
The left and right canonical $0202$ charts are the two successive legal
cycles. They transport one complete lift, which is normal by the left core.
The leaf correlations give $h<2e_L/3$, $\widehat h<2e_R/3$,
since each exclusive triple has width greater than $3h$ or $3\widehat h$.
Consequently
\[
 3/2<K,\widehat K<k=(3+4e^2/9)/2=6777889/4500000,
 \qquad T_L<T=e/2=167/2000.
\]
The left apex gives
$\mathcal C>3/4-9U^2/64=:C_{\rm low}=1731/2500>0$.
Corollary~\ref{loc:four-calibration} gives $\mathcal C<C_0$ on this
same lift. Thus, for the actual $z=\sqrt{3/4-\mathcal C}$,
\[
 z_0<z<z_1=3U/8=6/25.
\]

The actual shared interval satisfies Lemma~\ref{loc:shared-A-matrix}.
Take $X=103/500$, $x_0=1/10$. The exact comparisons
\begin{align*}
 2E-U-\mathcal R(X)&=194/126875>0,\\
 C_{\rm low}-k^2(1+x_0)^2/4
 &=49726705145159/8100000000000000>0.
\end{align*}
give $x_0<x<X$.
By \eqref{loc:shared-left-branch}, $r>0$ and $K>3/2$,
\[
 \lambda<-(3/2)x-z+kT.
\]
Its right side is negative before squaring, because
\[
 (3/2)x_0+z_0-kT=\frac{1311592537}{9000000000}>0.
\]
Therefore
\[
 \mathcal C<Q(x,z):=k^2(1+x)^2/4-((3/2)x+z-kT)^2.
\]
On the whole closed rectangle $[x_0,X]\times[z_0,z_1]$, $Q_x$ decreases
in both variables and
$Q_x(X,z_1)=221145810103507/2250000000000000>0$.
Hence $\mathcal C<Q(X,z)$.
Put $w=(3/2)X-kT=1649092537/9000000000>0$.
Using the actual identity $\mathcal C=3/4-z^2$ gives
\[
 \mathcal C-Q(X,z)=3/4+w^2+2wz-k^2(1+X)^2/4
 >\frac{259602733384166413}{81000000000000000000}>0,
\]
where the last rational is its value at $z=z_0$.
This contradiction proves the result. The shared arc is counted once;
no exterior return beyond the full nine-context was used.
\end{proof}

\subsection{A six-charge four-core and returned-tail union}
\label{loc:DK-section}
\begin{theorem}\label{loc:DK610}
In the actual positive full-eleven context $01020212012$ and its legal
phases, the cost $J$ of complete arcs $2,\ldots,7$, counted once,
exceeds $61/100$. The completely transported opposite chronology
$01201020212$ has the same bound on complete arcs $3,\ldots,8$.
\end{theorem}
\begin{proof}
Select one full lift containing all eleven complete arcs; its capped
four-core proves normality before joint normalization. The complete
source $[0,7]=01020212$ is proper, so Section~\ref{pf:M101} gives
$A_D=c([2,5])>74/125$.
The full reversal of $[3,10]=20212012$, with $j\mapsto1-j$, is
$20120212$. Its first seven arcs supply the strict negative entrance
by Corollary~\ref{loc:ST305-transfer}; its final six are $120212$.
Thus Theorem~\ref{cmp:returned-379} gives $A_K=c([5,7])>379/1000$.
Both actual caps and the returned arc remain on the selected lift.

Write $P=c_2+c_3+c_4$, $B=c_5$, $Q=c_6+c_7$.
Then $A_D=P+B$, $A_K=B+Q$, $J=P+B+Q$.
Under $J\le U=61/100$, the stronger constituent bounds give
\begin{equation}\label{loc:DK-budgets}
 P<p=231/1000<249/1000,\quad Q<e=9/500<11/125,\quad B>361/1000.
\end{equation}
The four-core leaf correlations, apex and calibration give
\[
 h<2p/3,\quad K<k=(3+(2p/3)^2)/2,\qquad
 0<C_{\rm low}:=3/4-9U^2/64<\mathcal C<C_0.
\]
Hence $z=\sqrt{3/4-\mathcal C}$ lies between
$z_0=243/2000$ and $z_1=183/800$.

Use the shared-arc formulas of Lemma~\ref{loc:shared-A-matrix}.
The whole phase $+1$ makes the tail's charged $212$ into $020$;
its first root is $-\xi<-\widehat h<0$.
Its full-eight chart is $01020120$, whose reflected-time image
under $j\mapsto2-j$ is the above $20120212$.
The strict negative canonical $202$ entrance just proved therefore
makes this $020$ exit $\beta>0$; the weaker $\beta\ge0$ suffices.
Since $361/1000-\mathcal R(37/200)=59/38500>0$, the actual shared cost gives
\begin{equation}\label{loc:DK-xupper}
 x<X=37/200<219/1000.
\end{equation}
For the two right-exclusive arcs put $\mu=\widehat h^2=2\widehat K-3$.
Their actual roots and physical cost imply
\[
 w=\xi+\beta<W=2e=9/250,\quad T<\bar T=e/2=9/1000,\qquad
\]
\[
 0<\mu<\xi^2<w^2,\quad \widehat K<k_h=(3+W^2)/2.
\]
The actual start density is $x/(1-\xi)^2$. Logarithmic density
comparison along these two arcs gives
\[
 q_{\rm end}\le\frac{x(1+\xi)^2}{(1-\xi)^2(1-\beta)^2}
 \le\frac{x(1+w)^2}{(1-w)^2}<\bar q:=\frac{X(1+W)^2}{(1-W)^2}.
\]
The middle inequality follows by clearing positive square roots;
the numerator difference is $\beta(1-\beta+\beta\xi+\xi^2)\ge0$.
At the actual final cap let $f$ be its signed inactive value.
Integration after the second root gives $\chi_{\rm end}\le\mu T$, hence
\[
 f\le\mu A,\qquad A=\{\bar q+((1+W)\bar q+1)\bar T\}/2.
\]
Two exact comparisons are
\[
 1-\bar q(1+W)^2=559517659243/726012500000>0,
\]
\[
 1/4-k_hA=59121297316151933/726012500000000000>0.
\]
The actual exit branch remains to be fixed. Put $d=1+\beta$.
At the cap $m=f-1/2$, and the inward derivative of its nonnegative
gap gives $\widehat L_{\rm exit}\ge-f_s/2$.
The raw derivative and $f\ge(\mu-\beta^2)q/2$ give
\[
 f_s=-\beta-\frac f{qd}+\frac{(\mu-\beta^2)qd}{2}
 \le-\beta+\frac{\mu-\beta^2}{2}(qd-d^{-1})<0.
\]
Indeed $qd^2<1$ proves the last sign when $\mu\ge\beta^2$
(including $\beta=0$, since $\mu>0$); otherwise it is bounded by
$-\beta+\beta^2/(2d)<0$.
Thus $\widehat L_{\rm exit}>0$, and the determinant identity now gives
\[
 \widehat L_{\rm exit}^2=\widehat K/2-\widehat Kf-\mathcal C
 \ge z^2+\mu(1/4-\widehat KA)>z^2.
\]
Integrating $\widehat L_s=\widehat K$ from the shared end therefore yields
\begin{equation}\label{loc:DK-lambda}
 \lambda<-z+b-(3/2)x,\qquad b=k_h\bar T=1688229/125000000.
\end{equation}

Set $x_0=1/10$, $a=3/2$. The exact margin
\[
 C_{\rm low}-kk_h(1+x_0)^2/4=284287889871/25000000000000>0
\]
gives $x>x_0$. Also
$ax_0+z_0-b=32249271/125000000>0$, so the bound on $\lambda$
is negative before squaring. Consequently
$\mathcal C<Q_*(x,z):=kk_h(1+x)^2/4-(ax+z-b)^2$.
Put $F=3/4-z^2-Q_*$ and $A_2=a^2-kk_h/4>0$.
On the rectangle $[x_0,X]\times[z_0,z_1]$,
\[
 F_z=3x-2b\ge17061771/62500000>0,
\]
\[
 F_x(x,z_0)\le F_x(X,z_0)=-4694043906513/25000000000000<0.
\]
Thus
\[
 F(x,z)\ge F(X,z_0)=21844694839413531/250000000000000000>0,
\]
contradicting $\mathcal C<Q_*$.
Whole reflected/reversed cyclic transport gives the opposite chronology
and its stated charges. All physical charges are counted once.
\end{proof}

\subsection{Three- and four-arc tails}\label{sec:tails}

\subsubsection{One comparison kernel for capped tails}
All quantities in this subsection belong to one actual positive capped
$12021$ source in an eligible complete interval. Its $202$ core starts at
$-a\le0$ and has switches and exit $-a<r<s<\beta$.
Lemma~\ref{loc:signed-strip} proves normality before joint normalization.
Write
\[
 \begin{gathered}
 w=a+\beta,\quad q=v^{-2},\quad \sigma=\int q\,du,\quad
 \mu=2K-3,\\
 Z=K\sigma,\quad X=f(-a)/\sigma,\quad Y=-f(\beta)/\sigma.
 \end{gathered}
\]
Thus $r^2\le\mu<\beta^2$, $s>\sqrt\mu>0$, $K=(3+\mu)/2$,
$Y,Z>0$, $-Z<X+Y<Z$, and
\begin{equation}\label{K361:loss}
 \mathcal C=\tfrac34-\tfrac14\mathcal L,\qquad
 \mathcal L=(X+Y+Z)^2-4ZY-\mu=(X+Y-Z)^2+4ZX-\mu.
\end{equation}
The initial response $X$ remains signed.

\begin{lemma}[Response, shear and terminal reserve]\label{cmp:tail-interface}
Set
\[
 I=\frac{r+a}{(1-r)(1+a)},\quad
 D_0=I+\frac{(1+r)(\beta-r)}{(1-r)^2(1+\beta)},\quad
 L_0=(1+a)^2D_0,\quad G_0=\frac I{(1-a^2)D_0}.
\]
Then $0\le G_0\le(1-a^2)^{-1}$ and
\begin{equation}\label{K361:response}
 X\le U_\kappa(\mu),\qquad
 U_\kappa(\nu)=\frac{(\nu-a^2)_+}{2L_0}+\kappa G_0\nu,
\end{equation}
\begingroup
\makeatletter
\protected@edef\@currentlabel{\p@equation\theequation}
\label{ek:response}
\makeatother
\endgroup
with $\kappa=1$; if the actual entrance also satisfies
$v(-a)^2>(1-a)(3+a)$, take $\kappa=2/(3+a)$ and additionally
$\sigma<J_{\rm in}:=w(1+a)/[(1-\beta)(1-a)(3+a)]$.
An actual final return, giving the complete source $120212$, supplies
\begin{equation}\label{K361:circle}
 v(\beta)^2>(1-\beta)(3+\beta),
\end{equation}
and the sharper shear ceiling
\begin{equation}\label{K361:shear}
 \sigma<J_0:=\frac{(1+\beta)^2}{4-(1+\beta)^2}
 \left[\frac{(1-r)(r+a)}{(1+r)^2(1+a)}
       +\frac{\beta-r}{(1+r)(1+\beta)}\right].
\end{equation}
The terminal reserve, with or without that return, is
\begin{equation}\label{cmp:terminal-reserve}
 Y>\frac{(\beta^2-\mu)(1-a)(1-\beta)}{2w(1+\beta)(1-s)}
 \begin{cases}(1+s)/(1+\beta),&\text{without a return},\\
 (3+s)/(3+\beta),&\text{with the actual return}.
 \end{cases}
\end{equation}
\end{lemma}
\begin{proof}
On the first arc $q=c_q/(1-u)^2$; in either binary mode
$(\log q)_u\ge-2/(1+u)$. Thus $\sigma\ge c_qD_0$ and the
actual switching integral gives $0\le-\chi(-a)\le\mu c_qI$.
The cap identity is
$f(-a)=(\mu-a^2)q(-a)/2+A_-(-a)\chi(-a)$.
Its coefficient satisfies $|A_-(-a)|<(1-a^2)^{-1}$, improved to
$2/[(1-a^2)(3+a)]$ by the entrance circle. This proves the response,
including its positive part; discarding $u^2$ in a signed integrand
uses the nonnegativity of the actual integral, not pointwise positivity.
Forward density comparison from that circle proves $J_{\rm in}$.

The lawful whole cyclic map $(z+3)/(1-z)$ carries the actual final
$212$ to $020$. Its first root is negative by
Lemma~\ref{loc:return-root}, which is exactly \eqref{K361:circle}.
Backward density comparison, with the exact first-arc shape, gives
\eqref{K361:shear}. The returned leaf has its own binary invariant;
none is identified across the cap.

On the last arc write $q=c/(1-u)^2$. Then
$\sigma_3=c(\beta-s)/[(1-\beta)(1-s)]$ and
$0<-\chi(\beta)<(\beta^2-\mu)\sigma_3$.
The cap formula therefore gives
\[
 -f(\beta)>\tfrac12(\beta^2-\mu)q(\beta)
 \left[1-\frac{(\beta-s)(1+c)}{1-s}\right].
\]
The open state domain gives $c<(1-\beta)/(1+\beta)$; the return
improves it to $(1-\beta)/(3+\beta)$. Substitute either in the
bracket and use
$\sigma\le q(\beta)(1+\beta)w/(1-a)$ to obtain
\eqref{cmp:terminal-reserve}.
\end{proof}

\begin{lemma}[Signed-strip comparisons]\label{cmp:strip-envelope}
For $Z>0$, $-Z<X+Y<Z$, and arbitrary real bounds $X\le U$, $Y\ge V$,
\begin{equation}\label{cmp:unclipped-strip}
 (X+Y+Z)^2-4ZY-\mu\le(Z+U+V)^2-4ZV-\mu.
\end{equation}
The inequality is strict when $Y>V$. If $U,V\ge0$, the sharper bound is
\begin{equation}\label{cmp:sharp-strip}
 (X+Y+Z)^2-4ZY-\mu
 \le[Z+\min(U+V,Z)]^2-4ZV-\mu.
\end{equation}
This is the sharp supremum under the displayed constraints and is
strict when $Y>V$. In particular the actual tail satisfies
\[
 \mathcal L<\min\{(Z+U_\kappa(\mu))^2,4Z^2\}-\mu;
 \qquad X\le0\ \Longrightarrow\ \mathcal L<Z^2-\mu.
\]
At fixed actual $a,\beta,r,Z$, a nonnegative piecewise affine response
$U(\nu)$ reduces $(Z+U(\nu))^2-\nu$ on any mean interval to its
endpoints and affine knots. A shear bound may be substituted afterward
at each candidate only when it still bounds the actual $Z$.
Alternatively, if $Z<(3+\mu)c$ with $c$ fixed by actual state and
budget, $[(3+\nu)c+U(\nu)]^2-\nu$ has the same endpoint property.
No varying mean is assigned a physical trajectory or cost.
\end{lemma}
\begin{proof}
Put $Q=X+Y$, $A_1=U-X\ge0$, $A_2=Y-V\ge0$. The difference
between the right and left sides of \eqref{cmp:unclipped-strip} is
\[
 (A_1-A_2)^2+2A_1(Z+Q)+2A_2(Z-Q)\ge0,
\]
strictly positive if $A_2>0$. This identity permits signed bounds
without asserting coordinatewise monotonicity outside the actual strip.
For the sharp nonnegative comparison, at fixed $Q$ the loss decreases
with $Y$, whose least allowed value is $\max(V,Q-U)$. The resulting expression increases
for $Q\le U+V$ (derivative $2(Q+Z)>0$) and decreases for
$Q\ge U+V$ (derivative $2(Q-Z)<0$). Its supremum occurs at
$Q=\min(U+V,Z)$, proving \eqref{cmp:sharp-strip}.
For $X\le0$ use $U=V=0$. Each untruncated squared affine expression
is convex on its own piece, proving the endpoint assertions.
When retaining a positive reserve or clipping at $Z$, one must also
include the clipping crossings; global convexity is not asserted.
\end{proof}

\begin{lemma}[Coefficient ceilings and endpoint responses]
\label{cmp:coefficient-responses}\label{cmp:circle-responses}
In the common actual capped-tail setting, suppose $w<1$ and the first
coefficient in $q=c_q/(1-u)^2$ has a qualified bound $c_q\le C(a,\beta)$,
independent of the first root when $a,\beta$ are fixed. Put
\begin{gather*}
 A_C=\frac{1+C}{2(1+a)},\qquad
 U_- =\frac{(\mu-a^2)(1+\beta)}{2w(1-a)},\qquad
 U_C=\frac{(\mu-a^2)(1-\beta)}{2w(1+a)}+A_C\mu,\\
 R=\frac{(\beta^2-\mu)(1-a)(1-\beta)}{2w(1+\beta)^2}.
\end{gather*}
Then $Y>R>0$. In LOW, $0<\mu\le a^2$, the signed bound is
$X\le U_C$. In HIGH, $\mu\ge a^2$, one has
$X\le\max(U_-,U_C)$.

In particular, an actual entrance circle
$v(-a)^2>(1-a)(3+a)$ permits
\[
 C=C_{\rm in}:=\frac{(1+a)^2}{(1-a)(3+a)},\qquad
 A_C=\frac{2}{(1-a^2)(3+a)}.
\]
Write $U_+=U_C$ for this circle choice. Then $X\le U_+$ if
$r\ge\beta/3$. For $r\le\beta/3$ the actual mean satisfies
$\mu<13\beta^2/27$. These circle assertions require no returned terminal arc.
\end{lemma}
\begin{proof}
Retain the actual first coefficient $c_q$. With $I,D_0$ from
Lemma~\ref{cmp:tail-interface}, logarithmic density comparison gives
\[
 c_qD_0\le\sigma\le c_qD_+,\qquad
 D_+=\frac{w}{(1+a)(1-\beta)}=D_0\big|_{r=\beta}.
\]
The actual first switching integral is nonnegative, while
$|A_-(-a)|=(1+c_q)/[2(1+a)]\le A_C$. The cap identity gives
\[
 \frac{f(-a)}{c_q}\le n_\mu(r):=
 \frac{\mu-a^2}{2(1+a)^2}+A_C\mu I.
\]
Only after using that nonnegative actual integral is its $u^2$ moment
discarded. The coefficient ceiling stays fixed in every root comparison.

In LOW, if $n_\mu(r)<0$, the upper bound on $\sigma/c_q$ gives
$X\le n_\mu(r)/D_+\le n_\mu(\beta)/D_+=U_C$.
If $n_\mu(r)\ge0$, use the lower bound on $\sigma/c_q$.
The quotient $n_\mu(r)/D_0$ increases with $r$: its first numerator
is nonpositive and $(D_0)_r=4(\beta-r)/[(1+\beta)(1-r)^3]\ge0$.
The derivative of $I/D_0$ has sign opposite to
$(w-2)r^2+2(\beta-a)r+2a\beta-w$, whose maximum on the real line is
$2w(w-1-a\beta)/(2-w)<0$. This proves the same signed bound.

In HIGH, $X\le U(r):=n_\mu(r)/D_0$. Write
\[
 N=w+2a\beta+2(\beta-a)r-(2+w)r^2,
 \qquad D_0=\frac{N}{(1-r)^2(1+a)(1+\beta)}.
\]
The concave $N$ is positive on $[-a,\beta]$, since its endpoint values
are $w(1-a^2)$ and $w(1-\beta^2)$. The endpoints of $U$ are $U_-,U_C$.
For $H=\max(U_-,U_C)$, the quadratic $N(H-U)$ has leading coefficient
\[
 -H(2+w)-\frac{(\mu-a^2)(1+\beta)}{2(1+a)}+A_C\mu(1+\beta)
 \le-A_C\mu(1+a)\le0,
\]
because $H\ge U_C\ge A_C\mu$ and $a^2\le\mu$.
Its nonnegative endpoint values prove $U\le H$.

For the stated circle choice the coefficient formula simplifies as
claimed. Its additional root comparison is the exact identity
\begin{gather*}
 U_+-U=\frac{(\beta-r)
 \{a^2(1-a)(3+a)(\beta-r)+\mu H_r\}}{(1-a)(3+a)wN},\\
 H_r=a^2(\beta-r)+2a(1+\beta)+(3+2\beta)r-\beta\ge0
 \quad(r\ge\beta/3).
\end{gather*}
At the actual first root $r^2<\mu$, so the root-qualified
Lemma~\ref{fin:signed-mean} gives
$r>-\beta/2$ and $\mu<(r^2+r\beta+\beta^2)/3$.
For $r\le\beta/3$ the latter is at most $13\beta^2/27$.
Finally, the nonreturned line of \eqref{cmp:terminal-reserve}, with
$s>\sqrt\mu>0$, gives $Y>R$.
\end{proof}

\begin{lemma}[Signed LOW width contraction]\label{cmp:low-width-contraction}
Fix $h>0$, put $b=1-x$, and consider the full domain
\[
 0\le w\le h,\quad0\le x\le1,\quad
 0\le k\le\min(x^2,b^2).
\]
Suppose $z,u,\rho$ are continuous, differentiable in $w$, and satisfy
throughout that domain
\[
 z_w,u_w\ge0,\quad\rho_w\le0,\qquad
 z\ge\tfrac12,\quad(k-x^2)/2\le u\le k/2,
 \quad0\le\rho\le(b^2-k)/2.
\]
At $w=h$, suppose also that $z=z_0+z_1k$, $u=u_0+u_1k$ and
$\rho=\rho_0-\rho_1k$, with $z_1,\rho_1\ge0$.
Then
\[
 F=(z+u+\rho)^2-4z\rho-k
\]
is nondecreasing in $w$ and convex in $k$ at $w=h$.
The zero-shear expression $(u+\rho)^2-k$ is at most $1/4$.
Consequently, for any $c\ge1/4$, strict upper bounds by $c$ for
\[
 F(h,x,0),\qquad F(h,x,x^2)\ (x\le1/2),\qquad
 F(h,x,(1-x)^2)\ (x\ge1/2)
\]
bound the full positive-shear comparison by $c$.
\end{lemma}
\begin{proof}
The hypotheses imply $z+u-\rho\ge xb+k\ge0$ and
$z-u-\rho\ge(1-b^2)/2\ge0$. These are exactly the signs of
$F_z/2$ and $-F_\rho/2$; also $F_u=2(z+u+\rho)\ge0$.
Thus every term of $dF/dw$ is nonnegative. Moreover
$-1/2\le u+\rho\le1/2$, proving the zero-shear assertion.
The leading mean coefficient is
$(z_1+u_1-\rho_1)^2+4z_1\rho_1\ge0$, which gives the three endpoints.
In a physical application the signed-strip comparison and convexity in
its positive actual shear are applied first; only the resulting
normalized scalar functions enter this lemma.
\end{proof}

\begin{lemma}[A balanced LOW boundary]\label{cmp:low-balance}
Put $b=1-x$, $0\le x\le1$ and $0\le k\le\min(x^2,b^2)$.
For the two rows $i=C,R$ of the following table define
\begin{gather*}
 z_i=(3+k/9)j_i/2,\qquad
 u_i=\frac{(k-x^2)(2+x)+k(1+C_i)}{2(3+x)},\\
 \rho_i=(b^2-k)R_i,\qquad
 F_i=(z_i+u_i+\rho_i)^2-4z_i\rho_i-k.
\end{gather*}
\[
\begin{array}{c|c|c|c}
 i&j_i&C_i&R_i\\\hline
 C&\dfrac{9(3+x)}{(2+x)(3-x)(9+x)}
  &\dfrac{(3+x)^2}{(3-x)(9+x)}
  &\dfrac{(3-x)(2+x)}{2(4-x)^2}\\[5pt]
 R&\dfrac{9(4-x)}{(3-x)(2+x)(10-x)}
  &\dfrac{(4-x)^2(3+x)^2}{(2+x)(10-x)(3-x)^2}
  &\dfrac{3(3-x)(2+x)}{2(10-x)(4-x)}
\end{array}
\]
Both scalar envelopes attain their maximum at $(x,k)=(1/2,1/4)$:
\begin{equation}\label{cmp:low-balanced-values}
 F_C\le\frac{5361444}{11055625}<\frac12,\qquad
 F_R\le\frac{142503006}{276390625}<\frac{13}{25}.
\end{equation}
\end{lemma}
\begin{proof}
Write $z_i=z_0+z_1k$, $u_i=u_0+u_1k$ and
$\rho_i=\rho_0-\rho_1k$, where $z_1,\rho_1\ge0$.
The quadratic coefficient of $F_i$ in $k$ is
$(z_1+u_1-\rho_1)^2+4z_1\rho_1\ge0$.
We prove $(F_i)_k(x,0)>0$, and then show that
$F_i(x,x^2)$ increases for $x\le1/2$ while
$F_i(x,(1-x)^2)$ decreases for $x\ge1/2$.
These three signs reduce both the mean and entrance to the stated point.

Here are elementary polynomial sign checks for all six signs. Set
\[
 d_C=(3-x)(x+2)(x+3)(x+9),\qquad
 d_R=(10-x)(3-x)(x+2)(x+3).
\]
Let $P_{iK},P_{iL},P_{iR}$ be the respective products of
$(F_i)_k(x,0)$, $\frac{d}{dx}F_i(x,x^2)$ and
$-\frac{d}{dx}F_i(x,(1-x)^2)$ by the positive factors
\[
\begin{array}{c|ccc}
 i&D_{iK}&D_{iL}&D_{iR}\\\hline
 C&\dfrac{(4-x)^4d_C^2}{x+3}&(4-x)^5d_C^3&d_C^3\\[4pt]
 R&\dfrac{(4-x)^2(3-x)d_R^2}{x+3}
   &(4-x)^3(3-x)^2d_R^3&(3-x)^2d_R^3.
\end{array}
\]
Each product is a polynomial. For a polynomial
$P(t)=\sum c_jt^j$ on $[0,1]$, put
$H=c_1+\sum_{j\ge2}\min(c_j,0)$. Then
$P(t)\ge c_0+Ht\ge c_0+\min(H,0)$.
Exact expansion after the listed substitutions gives
\[
\begin{array}{c|c|r|r}
 P&x&c_0&H\\\hline
 P_{CK}&1-t&600696&-253818\\
 P_{CL}&t/2&1428513408&-260156608593/2048\\
 P_{CR}&(1+t)/2&456265575/128&39780873/64\\
 P_{RK}&1-t&56862&105366\\
 P_{RL}&(1-t)/2&4787820755585/16384&-2395333095267/16384\\
 P_{RR}&1-t/2&24879312&43228534629/4096
\end{array}
\]
All six grouped lower bounds are strictly positive. Thus the asserted
mean and boundary signs hold on their closed intervals. Evaluation at
the common endpoint proves \eqref{cmp:low-balanced-values}.
Only scalar functions have been varied.
\end{proof}

\begin{lemma}[The actual loss obstruction]\label{cmp:loss-obstruction}
Let one full normal lift of an original global minimizer be jointly
normalized on a genuine fundamental interval. Suppose that, at an actual
point in a lawful whole chart, $K>0$, and write
$\mathcal C=3/4-\mathcal L/4$. Then $\mathcal L<\Delta:=59/1000$
is impossible.
\end{lemma}
\begin{proof}
Such a loss would give $\mathcal C>2941/4000>0$.
Corollary~\ref{cal:one-point}, applied at the actual point on this same
lift, gives $\mathcal C<H_*$. The contradiction is the exact gap
\begin{equation}\label{cmp:harmonic-gap}
 \frac{2941}{4000}-H_*
 =\frac{573375701790091817085661}{39478417604337498218323684000}>0.
\end{equation}
\end{proof}

\begin{lemma}[A scalar low-mean criterion]\label{cmp:low-criterion}
Suppose $\mu\le a^2$ and a contrary budget is $2T\le V<1$.
Put
\[
 h_V(w)=\frac{(3+w^2)(V-w)}6,\qquad
 R(V)=\frac{3-V}{3(1+V)}-\frac{V}{4(1-V^2/4)^2}.
\]
If $R(V)>0$, then $\mathcal L<Z^2$. Consequently a uniform shear
bound $Z^2<\Delta:=59/1000$ excludes this contrary budget.
\end{lemma}
\begin{proof}
The physical identity \eqref{cmp:action-transport}, equivalently
\eqref{loc:strip-cost}, gives $0<w<V$ and $Z<h_V(w)$.
The low-mean response is $X\le G\mu$, $0\le G\le(1-a^2)^{-1}$.
For $\mu\le\min(a^2,\beta^2)$, maximizing separately over
$a\le w/2$ and $a\ge w/2$ gives
\[
 \frac{\mu}{(1-a^2)^2}\le\frac{w^2}{4(1-w^2/4)^2}.
\]
On the second half the derivative has sign
$-(1+a^2-2aw)<0$, since $1+a^2-2aw=1-w^2+(w-a)^2>0$.
The exact identity
\[
 V-\frac{2h_V(w)}{1-w^2}-\frac{w(3-V)}{3(1+V)}
 =\frac{4w[(V-w)(1-Vw)+(1-V)w(1+w)]}{3(1+V)(1-w^2)}>0
\]
and $w<V$ yield $2ZG+\mu G^2<V-wR(V)<1$.
Lemma~\ref{cmp:strip-envelope} now gives
$\mathcal L<Z^2+\mu(2ZG+\mu G^2-1)<Z^2$.
If $Z^2<\Delta$, then $\mathcal L<\Delta$, contradicting
Lemma~\ref{cmp:loss-obstruction} on the selected lift, where $K>0$.
\end{proof}

\begin{lemma}[A circle LOW reserve]\label{cmp:circle-low-width}
On one actual capped source $12021$ in the common eligible complete
interval, suppose $0<\mu\le a^2$ and
$v(-a)^2>(1-a)(3+a)$. If $w\le1/3$, then
\[
 \mathcal L<\frac{w^2}{2}.
\]
Without any width restriction, its three-arc physical cost satisfies
\begin{equation}\label{cmp:circle-low-cost}
 T>81/200.
\end{equation}
Neither assertion requires a returned terminal arc or a charged prefix.
\end{lemma}
\begin{proof}
For the narrow assertion use the signed $U_+$ and reserve $R$ from
Lemma~\ref{cmp:circle-responses}. The untruncated signed comparison
\eqref{cmp:unclipped-strip} is strict because $Y>R$.
It is convex in $Z$, so use the endpoints
$0$ and $(3+\mu)J_{\rm in}/2$.
Put $x=a/w$, $b=1-x$, and $k=\mu/w^2$; thus
$0\le k\le\min(x^2,b^2)$. At the positive endpoint the normalized
shear, response and reserve are
\begin{align*}
 j&=\frac{1+wx}{(1-wb)(1-wx)(3+wx)},& z&=(3+w^2k)j/2,\\
 u&=\frac{(k-x^2)(1-wb)}{2(1+wx)}
       +\frac{2wk}{(1-w^2x^2)(3+wx)},&
 \rho&=\frac{(b^2-k)(1-wx)(1-wb)}{2(1+wb)^2}.
\end{align*}
The normalized loss is
$F=(z+u+\rho)^2-4z\rho-k$.
For fixed $x,k$, both $z,u$ increase with $w$, and $\rho$ decreases.
For the only nonimmediate term,
\[
 \frac{d}{dw}\frac{w}{(1-w^2x^2)(3+wx)}
 =\frac{3+3w^2x^2+2w^3x^3}{[(1-w^2x^2)(3+wx)]^2}>0.
\]
Throughout $0\le w\le1/3$,
\[
 z\ge\tfrac12,\qquad (k-x^2)/2\le u\le k/4,
 \qquad 0\le\rho\le(b^2-k)/2.
\]
Indeed the positive coefficient in $u$ is at most $1/4$.
Thus the whole-width hypotheses of Lemma~\ref{cmp:low-width-contraction}
hold with $h=1/3$; it also bounds the zero-shear endpoint by $1/4$.

At $w=1/3$ the positive-shear expression is exactly $F_C$ in
Lemma~\ref{cmp:low-balance}. That lemma and the zero-shear bound give
\[
 \mathcal L<\frac{5361444}{11055625}w^2<\frac{w^2}{2}.
\]
The artificial face $w=0$ enters only the continuous normalized
comparison, never an actual cap quotient.

For the cost assertion suppose $2T\le V=81/100$.
Lemma~\ref{cmp:low-criterion} applies because
$R(V)=22982170633/202388176501>0$.
The LOW mean satisfies $\mu\le w^2/4$, so
$Z<(3+w^2/4)(V-w)/6$, a decreasing function on $0<w<V$.
At $w=33/100$, its squared gap below $\Delta=59/1000$ is
$87454079/250000000000>0$.
For smaller widths the narrow result gives the gap
$\Delta-(33/100)^2/2=91/20000>0$.
Thus $\mathcal L<\Delta$, contradicting Lemma~\ref{cmp:loss-obstruction}
for the actual $K>0$ and proving \eqref{cmp:circle-low-cost}.
\end{proof}

\begin{lemma}[A returned LOW reserve]\label{cmp:returned-low-width}
On one actual complete returned source $120212$ in the common eligible
tail interface, suppose $0<\mu\le a^2$ and $w\le1/3$. Then
\[
 \mathcal L<\frac{13}{25}w^2<\frac8{15}w^2.
\]
Both caps and the actual final return are retained on the same normal
lift. No entrance circle or contrary cost budget is required.
\end{lemma}
\begin{proof}
Write $b=\beta$. Backward density comparison throughout the actual core is
\[
 q(u)\le\frac{q(b)(1+b)^2}{(1+u)^2},\qquad
 q(b)<\frac1{(1-b)(3+b)}.
\]
The second inequality is supplied by the actual final return.
Evaluation at $u=-a$ and integration give respectively
\begin{gather*}
 c_q<\bar C:=\frac{(1+b)^2(1+a)^2}{(1-b)(3+b)(1-a)^2},\\
 \frac\sigma{q(b)}\le\frac{w(1+b)}{1-a},\qquad
 \sigma<\bar J:=\frac{w(1+b)}{(1-a)(1-b)(3+b)}.
\end{gather*}
These bounds include negative first roots and retain the same density
and full lift; no invariant is identified across the terminal cap.
The actual returned terminal bracket is
$(1-b)(3+s)/[(3+b)(1-s)]>3(1-b)/(3+b)$, since $s>\sqrt\mu>0$.
Together with Lemma~\ref{cmp:coefficient-responses}, this yields
\[
 X\le U:=\frac{(\mu-a^2)(1-b)}{2w(1+a)}
                  +\frac{(1+\bar C)\mu}{2(1+a)},\qquad
 Y>R:=\frac{3(b^2-\mu)(1-a)(1-b)}{2w(3+b)(1+b)}.
\]
Apply \eqref{cmp:unclipped-strip} with this signed $U$, then use
convexity in $0\le Z\le(3+\mu)\bar J/2$.

Put $x=a/w$, $B=1-x$ and $k=\mu/w^2$.
The normalized positive-shear comparison has
\begin{align*}
 C&=\frac{(1+wB)^2(1+wx)^2}{(1-wB)(3+wB)(1-wx)^2},&
 j&=\frac{1+wB}{(1-wx)(1-wB)(3+wB)},\\
 z&=(3+w^2k)j/2,&
 u&=\frac{(k-x^2)(1-wB)+wk(1+C)}{2(1+wx)},\\
 \rho&=\frac{3(B^2-k)(1-wx)(1-wB)}{2(3+wB)(1+wB)}.
\end{align*}
Throughout $0\le w\le1/3$ and $0\le k\le\min(x^2,B^2)$,
$z,u$ increase with $w$ and $\rho$ decreases.
For $u$, its nonpositive first term increases, as do $w/(1+wx)$ and $C$;
for the latter group $(1+wB)^2/(3+wB)$ with its remaining increasing factors.
At fixed width, the logarithmic derivative of $\bar C$ in $a$ is
\[
 \frac4{1-a^2}-\frac2{1+b}-\frac1{1-b}+\frac1{3+b}>\frac12.
\]
Hence $C\le(1+w)^2/[3(1-w)^2]\le4/3$, and on the whole domain
\[
 z\ge\tfrac12,\qquad(k-x^2)/2\le u\le7k/18\le k/2,
 \qquad0\le\rho\le(B^2-k)/2.
\]
At the width endpoint $z,u,\rho$ are affine in $k$, with the
nonnegative slopes required by Lemma~\ref{cmp:low-width-contraction}.
That lemma handles the zero-shear endpoint and contracts the positive
endpoint to $w=1/3$.
The artificial width zero is used only for this continuous normalized comparison.

At $w=1/3$ the positive-shear expression is exactly $F_R$ in
Lemma~\ref{cmp:low-balance}. Together with the zero-shear bound this gives
\[
 \mathcal L<\frac{142503006}{276390625}w^2<\frac{13}{25}w^2,
\]
which proves both stated inequalities without any additional source.
\end{proof}

\begin{corollary}[Low means with and without a return]\label{sj:low}\label{cmp:low-costs}
An actual capped $12021$ core with nonpositive entrance and its own
middle mean $\mu\le a^2$ costs $T>13/40$.
If its actual final return is also retained, so the full source is
$120212$, then
\begin{equation}\label{cmp:low-return}
 T>407/1000>101/250>379/1000>3/8.
\end{equation}
The returned assertion requires no entrance circle.
\end{corollary}
\begin{proof}
For a return suppose $2T\le V=407/500$.
Lemma~\ref{cmp:low-criterion} applies because
\[
 R(V)=\frac{69053086394231}{631400423219307}>0.
\]
The LOW mean gives $\mu\le w^2/4$, so the decreasing shear bound is
$Z<(3+w^2/4)(V-w)/6$ on $0<w<V$.
At $w=1/3$ its squared gap below $\Delta=59/1000$ is
$17353079/104976000000>0$.
For $w\le1/3$, Lemma~\ref{cmp:returned-low-width} gives the gap
\[
 \Delta-\frac{13}{25}\left(\frac13\right)^2=\frac{11}{9000}>0.
\]
Lemma~\ref{cmp:loss-obstruction}, with the actual $K>0$, excludes the
contrary budget.

Without a return use $V=13/20$, for which
$R(V)=6117863/22525371>0$, and split at $c=7/40$.
For $w\le7/20$, the compact cut first gives
$v>14/15$, then $v^2>3607/3600>1$, hence
$q<(17/15-u^2)^{-1}$. Therefore
$Z<j_\Omega(w):=(3+w^2)w/[2(17/15-w^2)]$, an increasing function
on that interval. The exact gaps
\[
 \Delta-j_\Omega(c)^2=\frac{1048196707}{896507168000},\qquad
 \Delta-h_V(c)^2=\frac{211782839}{147456000000}
\]
prove $Z^2<\Delta$ on both sides of $c$ and complete the proof.
\end{proof}

\begin{lemma}[Uniform returned-HIGH absorption]\label{cmp:high-absorption}
For one actual source in Lemma~\ref{cmp:tail-interface}, retain its
actual final return. If $\mu>a^2$, $w\le1/2$ and $Z<1/4$, then
$\mathcal L<Z^2$. No entrance circle is required.
\end{lemma}
\begin{proof}
The actual $2\to0$ root has $r^2<\mu$, as proved after
Lemma~\ref{fin:signed-mean}. That lemma on $[r,s]$, followed by
$s<\beta$ and $s>|r|$, gives
\[
 \mu<\frac{r^2+rs+s^2}{3}<M:=\frac{r^2+r\beta+\beta^2}{3}.
\]
Thus $M>\max(a^2,r^2)$, $r>-\beta/2$ and $a<\beta$.
With $I,D_0,L_0,J_0$ from the common interface, put
\[
 c_*:=\frac{(1+\beta)^2(1-r)^2}{(1-\beta)(3+\beta)(1+r)^2},
 \quad B_2:=\frac{r^3+a^3}{3(1+a)^2},\quad J_0=c_*D_0,
\]
\begin{equation}\label{cmp:high-ray}
 E:=\frac{M-a^2}{2L_0}
   +\frac{(1+c_*)(MI-B_2)}{2(1+a)D_0},\qquad
 V:=\frac{3(\beta^2-M)}{2(3+\beta)^2J_0}.
\end{equation}
Both are nonnegative on this actual HIGH domain. We claim, for
$t=\mu/M\in(0,1)$, that $X\le tE$ and $Y>tV$.

Write the first density as $q=c/(1-u)^2$. Matching at $r,s$ and
the actual terminal circle give
\[
 c\left[\frac{(1+r)(1-s)}{(1-r)(1+s)}\right]^2
 <\frac{1-\beta}{3+\beta},\qquad c<c_*.
\]
Hence $|A_-(-a)|=(1+c)/[2(1+a)]<(1+c_*)/[2(1+a)]$.
For $I_\nu=\int_{-a}^r(\nu-u^2)(1-u)^{-2}du$, the switching
integral and first-shape bound are $-\chi(-a)=cI_\mu\ge0$ and
$\sigma\ge cD_0$. Moreover
\[
 I_\mu\le tI_M\le t(MI-B_2),\qquad
 \mu-a^2\le t(M-a^2),
\]
since $\int_{-a}^r u^2(1-u)^{-2}du\ge B_2\ge0$.
The entrance cap identity therefore proves $X\le tE$, retaining
the first-arc second moment rather than discarding it.

The backward density calculation for $J_0$, before inserting the
terminal circle, gives
$\sigma/q(\beta)\le(1-\beta)(3+\beta)J_0$.
The actual terminal cap estimate in the common interface gives
\[
 -f(\beta)>\frac{(\beta^2-\mu)q(\beta)}2
 \frac{(1-\beta)(3+s)}{(3+\beta)(1-s)}.
\]
As $s>0$, this implies
$Y>3(\beta^2-\mu)/[2(3+\beta)^2J_0]\ge tV$.
No invariant of the returned leaf has been identified with the core.

Apply the signed-strip envelope with the nonnegative bounds $tE,tV$.
Unclipping its square only increases the bound, so
\begin{equation}\label{cmp:high-rational}
 \mathcal L\le Z^2+t\{2Z(E-V)+t(E+V)^2-M\}.
\end{equation}
This keeps $X$ signed and avoids any monotonicity assertion outside
$-Z<X+Y<Z$. It remains to prove
$2Z(E-V)+(E+V)^2<M$.

For the finite polynomial comparison use
\[
 0\le w\le\tfrac12,\quad 0\le x,y\le1,\qquad
 a=wx/2,\quad \beta=w(1-x/2),\quad r=w(y-x/2).
\]
Put \(J=(3+M)J_0/2\) and
\(\bar z=\{w/4+2(1-2w)J\}/(2-3w)\).
The two checked residuals are
\[
 M-(E+V)^2=\frac{w^2P_0}{3H^2},\qquad
 M-(E+V)^2-2\bar z(E-V)=\frac{w^2P_1}{12(2-3w)H^2}.
\]
Here
\(H=6(1+a)^2(1-\beta)(3+\beta)(1+r)^2(1+\beta)n\),
where \(n=1+2a\beta/w+2(\beta-a)r/w-(2+w)r^2/w\), with the
continuous polynomial interpretation at \(w=0\).
The quadratic \(wn\) is concave in \(r\) and takes the values
\(w(1-a^2)\), \(w(1-\beta^2)\) at \(-a,\beta\). Thus \(n\ge3/4\)
for \(w>0\), and the same bound holds at the limiting face.
The other factors and \(2-3w\) are positive on the comparison box.

The \path{local-algebra} task evaluates the exact constructors in
\path{check_returned_high.py}. It checks these clearing identities
and the strict signs of all Bernstein coefficients of \(P_0,P_1\)
in tensor degrees \((18,20,14)\), \((19,20,14)\), respectively.
The Bernstein partition of unity proves both displayed residuals positive
for $w>0$.

Now $\bar z$ is a convex combination of $1/4$ and $J$, with weights
$w/(2-3w)$ and $2(1-2w)/(2-3w)$. The actual return gives
$Z<KJ_0<J$, so $Z<\min(1/4,J)\le\bar z$.
If $E\le V$ the first residual proves absorption; if $E>V$ the
second does. Equation~\eqref{cmp:high-rational} proves the lemma.
The artificial face $w=0$ is checked only through the factored polynomial
identities. Actual HIGH has $M>\mu>a^2\ge0$, so no division by zero
occurs. The argument also allows $t=1$ in the relaxed comparisons.
\end{proof}

\subsubsection{A returned nonpositive tail}\label{loc:K361-section}
\begin{theorem}[Returned nonpositive core]\label{cmp:returned-379}
Consider an actual positive context $120212$ in one applicable complete
basic interval, with the full original hypotheses. If its complete
$202$ core enters at $u_0\le0$, its three charged arcs cost
$T>379/1000$. Its two caps and final returned arc are positive and
uncharged. No entrance circle is required.
\end{theorem}
\begin{proof}
Suppose \(2T\le V_*=379/500\) and set \(\Delta=59/1000\).
The LOW case contradicts \eqref{cmp:low-return}. In HIGH,
the complete return and the cost budget give
\[
 Z<\min\{j(w),h_{V_*}(w)\},\qquad
 j(w)=\frac{(3+w^2)w(1+w)}{2(1-w)(3+w)}.
\]
The first function increases and the second decreases. At
\(c=143/500\) each is below \(z_*=97/400\), and
\(\Delta-z_*^2=31/160000>0\).
Thus \(Z<z_*<1/4\). For \(w\le1/2\),
Lemma~\ref{cmp:high-absorption} gives \(\mathcal L<Z^2<\Delta\).

For \(w\ge1/2\), we have \(a<V_*/2\), \(\beta<V_*\)
and \(\mu<V_*^2<3/5\). The weak response
\eqref{loc:strip-X-upper} gives
\[
 X/\mu<\frac{1+V_*}{1-V_*/2}+\frac{1}{1-V_*^2/4}<4,
 \qquad Z<q_*(3+\mu),\qquad q_*=43/1000 .
\]
For \(X\le0\), the signed envelope gives \(\mathcal L<Z^2\).
For \(X>0\) apply the following bounds to that same envelope:
\[
\begin{array}{c|c|c}
 \text{range of }\mu &\text{upper bound for }\mathcal L
 &\text{least gap to }\Delta\\ \hline
 {[1/25,3/5]}&4q_*^2(3+\mu)^2-\mu&2394463/78125000\\
 {[0,1/25]}&[q_*(3+\mu)+4\mu]^2-\mu&42359/1000000
\end{array}
\]
In the first row the upper bound decreases. In the second it is
convex, so compare both endpoints. Both ranges give \(\mathcal L<\Delta\).
The actual \(K>0\) now gives a contradiction by
Lemma~\ref{cmp:loss-obstruction}.
\end{proof}

\subsubsection{Transfer of the nonpositive-tail estimate}
\begin{corollary}[Actual localized-tail transfer]\label{loc:ST305-transfer}
The right contexts $1201210$, $2012021$ and $0120102$ have cost greater than
$61/200$ on complete arcs $[3,4,5]$. Their left contexts
$0212012$, $1020120$ and $2101201$ have the same bound on $[1,2,3]$.
All seven arcs are positive and the proper-context guard is eight.
\end{corollary}
\begin{proof}
In the actual context $1201210$, its forward five-arc source $12012$
at indices $0,\ldots,4$ has the shared mode-$0$ arc at index $2$.
The actual forward comparison in Lemma~\ref{fin:forward-force} gives
$Y>Z>0$ at its terminal state in reciprocal coordinates.
Apply the lawful cyclic map $\mathcal R^2$ to the complete actual
$01210$ tail at indices $2,\ldots,6$. It becomes $12021$ and its
entrance coordinate is exactly
\[
 u_0=(Z-Y)/(Z+Y)<0.
\]
Theorem~\ref{loc:ST305} therefore charges only $[3,4,5]$ and retains
caps $[2,6]$. The forward source supplies a sign on the same lift;
it is not charged by this resource.
The full reflection and time reversal sends the word to $2101201$;
a legal cyclic addition gives $0212012$, with charged indices
$[1,2,3]$ and caps $[0,4]$. The cyclic copies give the stated six contexts.
Every seam retains one genuinely admissible reanchor of the complete
seven-arc context.
\end{proof}

\subsubsection{The actual forward-prefix entrance}\label{loc:positive-S}
\begin{lemma}[Forward-prefix forcing]\label{cmp:forward-prefix-forcing}
On one selected nonzero full free-time maximizing lift in a genuine
complete containing interval, retain the actual positive source $012021$.
Let $E$ be the physical cost of arc $1$ and $T$ that of arcs $2,3,4$.
At the entrance to its $202$ core, write
\[
 x=(1+u)/v,\qquad y=(u^2+v^2-1)/(2v),\qquad z=(1-u)/v.
\]
Then the selected lift is normal, and
\begin{equation}\label{Spos:prefix}
 x<\min(y,z),\qquad u=-a<0,\qquad
 v^2>(1-a)(3+a).
\end{equation}
All six arcs, including both exterior caps, are retained; no cost,
mean, or harmonic-height hypothesis is required.
\end{lemma}
\begin{proof}
The complete prefix $0120$ retains its two middle arcs $1,2$ and
its actual third-mode ties supplied by arcs $0,3$.
Apply the whole cyclic phase $+1$ to the state and full lift. It gives
$1201$ and reciprocal coordinates $(x',y',z')=(z,x,y)$ at the
same actual middle switch. Lemma~\ref{sc:pair} proves normality on
the selected lift before joint normalization; Lemma~\ref{fin:forward-force}
then gives $y'<\min(x',z')$. Thus $x<z,y$, and
\[
 u=(x-z)/(x+z)<0,\qquad
 2v(y-x)=v^2-(1-a)(3+a)>0.
\]
The containing interval was selected for the whole source, so this
restriction requires neither a new reanchor nor a second lift.
\end{proof}

\subsection{Complete-prefix interfaces at a nonpositive tail entrance}
\label{Sneg:interfaces}
Fix the actual complete context $012021$ on one normalized complete lift.
The cost $E$ is that of its mode-$1$ arc at position $1$; the cost $T$
is that of its three complete arcs $202$ at positions $2,3,4$.
The exterior positive caps $0,5$ are retained and uncharged.
At the tail entrance write $u=-a\leq0$, and at its exit write $u=\beta$.
Put
\[
 x=(1-a)/v,\qquad y=(v^2+a^2-1)/(2v),\qquad z=(1+a)/v,
 \qquad xy+xz+yz=1.
\]
Normality, the actual mean $0<\mu<\beta^2$, $\beta>0$, and the signed
strip follow from Lemma~\ref{loc:signed-strip}; set $K=(3+\mu)/2$.
Lemma~\ref{cmp:forward-prefix-forcing} gives $a>0$ and $y>x$.
For the prefix response calculation set
\[
 P=vp,\quad U=\operatorname{tr}((\Lambda-A_*)X),\quad
 V=\nabla_H U+JP,\quad J(\xi,\eta)=(-\eta,\xi),
\]
where $\nabla_H=v\nabla$ holds the matrix costate fixed, and put
\[
 \nu_0=(1,0),\quad\nu_1=\frac{(x^2-1,2x)}{1+x^2},\quad
 \nu_2=\frac{(z^2-1,-2z)}{1+z^2},\quad\lambda_i=V\cdot J\nu_i.
\]
At the actual $1\to2$ tie, $P=h\nu_2$, $U=-h$ and $h\ge0$.

\begin{lemma}[Actual complete-prefix credit]\label{Sneg:prefix}
For every such nonpositive entrance,
\[
 E>\frac{(\mu-a^2)_+}{1+a}.
\]
If $\mu>a^2$, then the entrance necessarily satisfies $y>x$.
\end{lemma}
\begin{proof}
The actual $1\to2$ cap has $P=h\nu_2$, $h\geq0$, and its outgoing
maximizing jet gives $\lambda_2\leq0$.
The cap formula for its vertical response is
\[
 \lambda_0=V_y=-\frac{\mu-a^2+p_c}{v},
 \qquad p_c=\frac{2hv}{1+z^2}\geq0.
\]
The positive zero-sum normal weights are
$w_0=y(x+z)^2$, $w_1=z(1+x^2)$, $w_2=x(1+z^2)$.
Thus, if $\mu>a^2$, their response balance gives
\[
 \lambda_1\geq\frac{w_0}{w_1}
                   \frac{\mu-a^2+p_c}{v}>0.
\]

Under the full lawful map $w\mapsto(w-3)/(w+1)$, the complete preceding
mode-$1$ arc becomes mode $0$. Its constant reciprocal coordinate is
$y$ in the old chart; its terminal horizontal coordinate is
$u_*=(y-x)/(y+x)$, and its terminal vertical response is $\lambda_1$.
For its actual duration $d>0$, the original response equations give
\[
 \dot\chi=-\chi/y-V_y,\qquad \dot V_y=2u,\qquad
 \chi(0)=0,\quad \chi(d)=2hy/(1+y^2)\geq0.
\]
Since $u(t)<u_*$ for $t<d$, the integrating factor yields
\[
 0\geq\int_0^d e^{t/y}V_y(t)\,dt,\qquad
 V_y(t)>\lambda_1-2u_*(d-t).
\]
If $u_*\leq0$, the latter is positive everywhere, a contradiction.
For $u_*>0$, the strict exponentially weighted mean of $d-t$ is less
than $d/2$. Equivalently, for $b=d/y>0$,
$(b-2)e^b+b+2>0$; this follows by twice differentiating, since its
second derivative is $be^b>0$ and its value and first derivative at
zero vanish. Hence $\lambda_1<du_*$.

On this same transformed physical arc, $\ell>y$ pointwise, since
\[
 \ell-y=\frac{(1-u')^2+v'^2}{2v'}>0.
\]
Therefore $E>yd>y\lambda_1/u_*$.
Using $1+x^2=(x+y)(x+z)$, its coefficient of
$\mu-a^2+p_c$ is
\[
 \frac{y w_0}{u_*w_1v}
 =\frac{y^2(x+z)}{z(y-x)v}
 >\frac1{1+a}.
\]
Indeed $zv=1+a$, and the reciprocal constraint gives the strict identity
\[
 y^2(x+z)-(y-x)=x(1-yz)=x^2(y+z)>0.
\]
This proves the assertion when $\mu>a^2$. Otherwise its right side is
zero and the actual complete prefix has positive cost.
Every inequality involved the original single normalized lift.
\end{proof}

\begin{lemma}[Circle shear and the wide-width bound]\label{Sneg:circle}
Suppose the actual entrance has $y>x$. Put $w=a+\beta$ and
$\sigma=\int_{\mathrm{tail}}dt/v$, $Z=K\sigma$.
Then
\[
 \sigma<J(a,\beta):=\frac{w(1+a)}
 {(1-\beta)(1-a)(3+a)}.
\]
If $E+T\leq A\leq1/2$, then $0<w<2A$ and
\[
 Z<h_A(w):=\frac{(3+w^2/4)(2A-w)}6,\qquad
 \mathcal L<4Z^2,
\]
where $\mathcal C=3/4-\mathcal L/4$ is the actual signed-strip
determinant.
\end{lemma}
\begin{proof}
The entrance condition is exactly $v^2>(1-a)(3+a)$, so
Lemma~\ref{cmp:tail-interface} supplies the shear ceiling $J=J_{\rm in}$.
The exact physical-cost identity and Lemma~\ref{Sneg:prefix} give
\[
 0<Z<
 \frac{3+\mu}{6}
 \left(2A-w-\frac{2(\mu-a^2)_+}{1+a}\right).
\]
In particular $w<2A$. Set $D=2A-w>0$.
If $\mu\leq w^2/4$, the required bound follows immediately.
Otherwise $\mu<\beta^2$ forces $\beta>w/2$, hence $a<w/2$ and
$\mu-a^2\geq\mu-w^2/4>0$.
The difference of the preceding numerator from $(3+w^2/4)D$ is
\[
 (\mu-w^2/4)D-\frac{2(3+\mu)(\mu-a^2)}{1+a}
 \leq(\mu-a^2)\left(D-\frac{2(3+\mu)}{1+a}\right)<0,
\]
because $D<1$ and $a<1$. This proves the $Z$ bound.
Finally $-Z<Q<Z$ and $Y,\mu>0$ in
$\mathcal L=(Q+Z)^2-4ZY-\mu$ imply $\mathcal L<4Z^2$.
\end{proof}

\subsubsection{A nonpositive-entrance four-arc estimate}\label{Sneg:section}
\begin{proposition}\label{loc:S405-nonpositive}
In an actual complete context $012021$, let $E$ be the physical cost
of arc $1$ and $T$ the cost of the three complete arcs $2,3,4$.
If the entrance coordinate of the $202$ tail is nonpositive, then
$E+T>81/200$. Both exterior caps are positive and uncharged.
\end{proposition}
\begin{proof}
Assume $E+T\le A=81/200$, and put $\Delta=59/1000$.
All comparisons use the actual complete normalized lift from
Section~\ref{Sneg:interfaces}, retaining $0<\mu<\beta^2$ and
$-Z<X+Y<Z$. A scalar comparison below does not replace its physical state.

\paragraph{The forced-circle LOW sector.}
The full prefix forces $y>x$ by Lemma~\ref{cmp:forward-prefix-forcing}.
If $\mu\le a^2$, Lemma~\ref{cmp:circle-low-width} gives
$T>81/200$, contradicting the budget since $E>0$.

\paragraph{The two HIGH root sectors.}
It remains that $\mu>a^2$. The actual full-prefix
Lemma~\ref{Sneg:prefix} forces $y>x$, hence the entrance circle, and gives
$E>(\mu-a^2)/(1+a)$. Set $w=a+\beta$. Then
\[
 \begin{gathered}
 0<Z<\min\{Z_p,B(\mu)\},\qquad
 Z_p=\frac{(3+\mu)w(1+a)}{2(1-\beta)(1-a)(3+a)},\\
 B(\nu)=\frac{3+\nu}{6}
       \left(2A-w-\frac{2(\nu-a^2)}{1+a}\right).
 \end{gathered}
\]
Lemma~\ref{Sneg:circle} gives $0<w<2A$,
$Z<h_A(w)=(3+w^2/4)(2A-w)/6$ and $\mathcal L<4Z^2$.
This decreasing width bound, with $W=23/40$, has the strict reserve
\[
 \Delta-4h_A(W)^2=10173728831/14745600000000>0.
\]
Hence $w<W$, and HIGH implies $0\le a<\beta$, so $a<W/2<1$.
Lemma~\ref{cmp:circle-responses} now reduces the actual response to
$U_+$ on $a^2\le\mu\le\beta^2$, and additionally $U_-$ on
$a^2\le\mu\le13\beta^2/27$. These two scalar domains cover both
root sectors; empty mean intervals contain no actual point.
All following substitutions use the signed identity
\eqref{cmp:unclipped-strip}.

\paragraph{Narrow widths, $w\le3/10$.}
Discard the terminal reserve and use $Z_p$.
For $k=a/w\in[0,1/2]$ and $q=\mu/w^2\ge k^2$, both
$Z_p/w$ and $U_-/w$ are nonnegative and increasing with $w$.
The derivative of the first term of $U_+/w$ is
$-(q-k^2)/[2(1+kw)^2]\ge-q/2$; that of its second term is at least
$6q/(3+kw)^2>q/2$, because $kw\le3/20$.
Thus $(Z_p+U_\pm)/w$ increases as well.
The normalized squared-total bound is convex in $q$, so compare at
$w=3/10$ and only at the three mean endpoints
$\mu=a^2$ (using $U_+$), $\mu=\beta^2$ (using $U_+$), and
$\mu=13\beta^2/27$ (using $U_-$).

Here is explicit scalar arithmetic for all three endpoints. Put
$a=3x/20$, $\beta=3(2-x)/20$, $0\le x\le1$, and denote their
means and totals by $(m_0,S_0),(m_+,S_+),(m_-,S_-)$, respectively,
where $S_j=Z_p(m_j)+U_j(m_j)$ and $U_0=U_+$.
Define positive denominators
\begin{gather*}
 D=(x+20)^2(20-3x)^2(14+3x)^2,\qquad D_-=6350400D,\\
 D_0=19600D(20+3x)^2,\qquad D_+=1225D(20+3x)^2.
\end{gather*}
The three polynomials are fully specified by
\begin{align*}
 P_0&=D_0(72/1225-S_0^2+m_0),\\
 P_+&=\frac{D_+}{x}(72/1225-S_+^2+m_+),\\
 P_-&=D_-(72/1225-S_-^2+m_-).
\end{align*}
The quotient continues polynomially at $x=0$. If
$P=\sum c_jx^j$, then $P\ge c_0+\sum_{j\ge1}\min(c_j,0)$ on $[0,1]$.
Exact expansion gives, respectively,
\[
 \begin{array}{c|r|r}
 &\deg P&c_0+\sum_{j\ge1}\min(c_j,0)\\\hline
 P_0&10&2970552730752\\
 P_+&7&499952198988\\
 P_-&10&545077762579
 \end{array}
\]
Thus each endpoint loss is at most $72/1225$, including $x=0$.
If the frozen normalized loss is negative, its smaller-width bound is
already negative. Otherwise multiply by $w^2\le(3/10)^2$.
The table and mean convexity therefore prove
$\mathcal L\le72/1225<\Delta$, with gap $11/49000$.

\paragraph{Critical widths, $3/10\le w\le d=13/40$.}
Put $\delta=\mu-a^2\ge0$ and $z=607/2500$. A fixed mixture of the
two actual shear bounds gives the useful affine ceiling
\begin{equation}\label{Sneg:mixed-shear}
 Z<z+a^2-11\delta/20.
\end{equation}
To prove it, write $J=2Z_p/(3+\mu)$ and
$G=J/6+(2A-w)/9$, $g=121/1500+a/100$.
On the enclosing rectangle $w\in[3/10,d]$, $a\in[0,13/80]$,
$J$ is separately convex in $w,a$: width convexity is immediate, and
\[
 J={w^2\over(2-w)(2+w)(1+a-w)}
   +{w\over2(2+w)(3+a)}+{w\over2(2-w)(1-a)}
\]
is a sum of positive convex reciprocal functions of $a$.
Thus $G-g$ has its maximum at a corner. The four corner reserves
$g-G$, in increasing $w$ then $a$, are
\[
 1/5250,\quad112793/623796800,\quad7/243000,\quad
 252253/16354324800.
\]
Consequently $G<g<1/12$. Combining the strict circle and budget bounds,
\begin{align*}
 Z&<(Z_p+2B(\mu))/3\\
  &=(3+a^2)G+\delta\left[G-{2(3+a^2)\over9(1+a)}\right]
                        -{2\delta^2\over9(1+a)}\\
  &<121/500+3a/100+a^2/12
                       +\delta[-293/500+203a/300].
\end{align*}
Here $1/(1+a)\ge1-a$ and only favorable terms were dropped.
The gap to \eqref{Sneg:mixed-shear} is therefore greater than
\[
 H=1/1250-3a/100+11a^2/12+\delta(9/250-203a/300).
\]
If the last coefficient is nonnegative, completing the square gives
$H\ge61/110000$. Otherwise use
$\delta\le \beta^2-a^2\le d^2-3a/5$, obtaining
\[
 H\ge496a^2/375-2363a/19200+1841/400000
       \ge2120451959/1218969600000>0.
\]
All denominators are positive on the enclosing rectangle; the affine
shear ceiling is positive at every actual datum because $Z>0$.

The existing response bounds now finish both root sectors without a
terminal reserve. Indeed
\[
 U_+\le(2/3)a^2+(23/12)\delta,\qquad U_-\le9\delta/4.
\]
For the first, its cap coefficient is at most
$5/3-40/93=115/93<5/4$, and its circle coefficient is at most $2/3$
since $1-3a-a^2>0$ for $a\le13/80$.
For the second, $\mu\le13\beta^2/27$ implies
$a<7w/17\le91/680<1/7$, giving $1/(2w)+1/[2(1-a)]\le9/4$.
Apply \eqref{cmp:unclipped-strip} with the applicable nonnegative
response and discard only the actual positive terminal reserve.
Writing $s=a^2$, $t=\delta$, the two loss bounds are
\[
 F_+=(z+5s/3+41t/30)^2-s-t,\qquad
 F_-=(z+s+17t/10)^2-s-t.
\]
They are convex in $(s,t)$. For $F_+$ use the polygon
$s,t\ge0$, $s\le d^2/4$, $s+t\le d^2$; its three nonzero vertices
are $(d^2/4,0),(0,d^2),(d^2/4,3d^2/4)$.
Their differences from $z^2$, after division by positive factors,
are bounded above respectively by
\[
 5/6+11/144-1,\quad41/60+(1681/900)(11/100)-1,\quad
 173/240+(29929/14400)(11/100)-1,
\]
all negative, using $z<1/4$ and $d^2<11/100$.
For $F_-$ use $s,t\ge0$, $s+t\le m_*=13d^2/27<51/1000$.
At its two nonzero vertices the corresponding upper bounds are
$1/2+51/1000-1<0$ and
$17/20+(289/100)(51/1000)-1<0$.
Thus both comparisons are at most $z^2$, and
$\mathcal L<z^2<\Delta$ with gap $301/6250000$.
The closed polygons include all mean and width boundaries and assert
no new physical state.

\paragraph{Wide widths, $13/40\le w\le23/40$.}
Set $f_\pm=(1\mp\beta)/[2w(1\pm a)]$ and
$c_*=2/[(1-a^2)(3+a)]\le2/3$; the latter uses
$1-3a-a^2\ge351/6400>0$ for $a\le23/80$.
Directly from the budget,
\begin{align*}
 B+U_-&\le A-w/2+[(2A-w)/6]\mu
                      +[f_--1/(1+a)]\delta,\\
 B+U_+&\le A-w/2+[(2A-w)/6+c_*]\mu
                      +[f_+-1/(1+a)]\delta.
\end{align*}
For $U_+$,
$f_+-1/(1+a)=1/(2w)-3/[2(1+a)]\le20/13-120/103$;
the total coefficient is at most $600361/535600<3/2$.
For $U_-$, $a\le w/2\le23/80<1/3$ gives
$f_--1/(1+a)\le g(w):=1/(2w)+3w/4-1/2$.
Here $g>0$ and $g'=-1/(2w^2)+3/4<0$ throughout the width interval.
The total coefficient is therefore at most $42527/31200<3/2$
(and at most $671/600$ for $w\ge2/5$).
Both chosen $\delta$-coefficient bounds, $20/13-120/103$ and $g(w)$, are positive.
Thus $\delta\le\mu$ gives
\[
 0<Z+U_\pm\le B+U_\pm\le A-w/2+(3/2)\mu.
\]
Discard the reserve. Convexity in $0\le\mu\le w^2$ leaves its two
endpoints. At $\mu=0$ the width function is convex; at $\mu=w^2$
its second derivative is $93/100-9w+27w^2>0$.
The four exact endpoint target gaps, with columns $\mu=0,w^2$, are
\[
 \begin{array}{c|cc}
 w=13/40&31/160000&39671/10240000\\
 w=23/40&7231/160000&136391/10240000.
 \end{array}
\]
All common width and mean boundaries are included.

Thus $\mathcal L<\Delta$ in every actual sector, contradicting
Lemma~\ref{cmp:loss-obstruction} for the actual $K>0$.
\end{proof}

\begin{theorem}[The four-arc scalar resource]\label{loc:S405}
Every actual positive six-arc context $012021$ has cost greater than
$81/200$ on its four complete interior arcs. The same holds for its
lawful cyclic images and full reflection/time-reversal images.
\end{theorem}
\begin{proof}
Lemma~\ref{cmp:forward-prefix-forcing} makes the actual entrance
strictly negative, so Proposition~\ref{loc:S405-nonpositive} applies.
All transformations transport
the complete state and full lift together and preserve physical cost.
\end{proof}

\subsubsection{The two outward entrances for a shared charge}\label{loc:W-unions}
In the complete context \(02120102\), the central capped pair
\(1201\) at positions \(2:5\) gives \(y<x,z\) at the shared state.
The two whole charts of Theorem~\ref{fee:W71} send its reciprocals
\((x,y,z)\) to \((y,z,x)\) and \((y,x,z)\).
In each chart the entrance is strictly negative and the entrance
circle holds, because the second reciprocal exceeds the first.
These are conditions on the same actual state. Theorem~\ref{fee:W71}
uses them and the common switching derivative to charge the two
disjoint triples \(1:3\) and \(4:6\).
\space
\subsection{Central geometry and exact transport}\label{ctr:section}

All sources in this section consist of complete positive maximal arcs of
an original minimizer, retained in one genuine containing interval by
Lemma~\ref{fd:reanchor}. Each argument selects one full lift on that interval;
its state covector is not extended periodically. Write
\[
 \Lambda=\begin{pmatrix}A&B\\m&-A\end{pmatrix},\qquad
 K=2uA+B-u^2m,\quad L=A-um,\quad \chi=vp_v,
 \quad ds=dt/v.
\]
Here and below $t$ denotes physical time; accumulated physical cost will be denoted by $\xi$.

Consider now the actual seven-arc configuration
\begin{equation}\label{ctr:gap}
                   020\mid1\mid202.
\end{equation}
Denote the central endpoints by $C,D$, its physical cost by $c>0$, and the two genuine binary means by $\mu_L,\mu_R$. At the endpoints, actual maximization gives
\begin{equation}\label{ctr:cap-signs}
 p_v(C)\ge0\ge p_v(D),\qquad H_0(C)=H_1(C),\qquad
 H_1(D)=H_2(D),
\end{equation}
and $H_1\ge H_0,H_2$ throughout the centre. With $j=-1$ at $C$ and $j=1$ at $D$, the cap equation is
\begin{equation}\label{ctr:cap-m}
 m_i+\frac12=\frac{\mu_i-u_i^2}{2v_i^2}+\chi_i\beta_i,
 \qquad \beta_i=\frac12\left(\frac{u_i-j}{v_i^2}+\frac1{u_i-j}\right).
\end{equation}
In particular $\chi_i\beta_i\ge0$.
Define
\[
 A_j=\frac{(u-j)^2-v^2}{4v^2(u-j)},\qquad
 B_j=\frac{(u-j)^2+v^2}{4(u-j)^2}>0.
\]
The actual one-sided maximizing derivatives are precisely
\begin{equation}\label{ctr:switch-rows}
 \begin{aligned}
 S_L&=L_C-u_C/2+A_{-1}(\mu_L-u_C^2)-B_{-1}\chi_C\ge0,\\
 S_R&=-L_D+u_D/2-A_1(\mu_R-u_D^2)+B_1\chi_D\ge0.
 \end{aligned}
\end{equation}
For clarity, differentiating $H_{\rm new}-H_{\rm old}$ along the central field gives
\[
 p\cdot(D\Delta F\,F_1-DF_1\,\Delta F)
       -\nabla\operatorname{tr}((\Lambda-A_*)X)\cdot\Delta F.
\]
Substitution of the appropriate cap tie and \eqref{ctr:cap-m} reduces this to $4v^2/(u^2+v^2-1)$ times the corresponding expression in \eqref{ctr:switch-rows}. Its multiplier is positive. Thus weak switches and triple ties are included, with no transversality assumption.

\subsubsection{A finite crossing-height bound}

\paragraph{Central coordinates and physical cost.}
The following identities apply more generally to any normal normalized
mode-$1$ arc with the actual endpoint cap equations, genuine means,
maximizing rows and full state-adjoint transport above. They require
neither the complete source \eqref{ctr:gap} nor a height or cost cut.
Its reciprocal $y=(u^2+v^2-1)/(2v)>0$ is constant, and its endpoint
$x=(1+u)/v$ values satisfy $a>b>0$, $ay<1$. Put
\[
 h=\sqrt{1+y^2},\quad k=h-y,\qquad
 u(t)=\frac{t^2+2ty-1}{1+t^2},\quad
 v(t)=\frac{2(t+y)}{1+t^2}.
\]
The physical cost, directly integrated from the original density, is
\begin{equation}\label{ctr:density}
 c=\int_b^a f(t)\,dt,\qquad
 f(t)=y-\frac{y^2}{y+t}+\frac{y(1+y^2)}{(y+t)^2}
     =y\left(1+\frac{1-yt}{(y+t)^2}\right)>y.
\end{equation}
Thus $c>\rho:=y(a-b)$. If $u_C\ge0\ge u_D$, then $b\le k\le a$;
the unique crossing has height $V=1/k$, since the mode-$1$ field has
$u'=-v(1-u^2+v^2)/(u^2+v^2-1)<0$.

\begin{lemma}[Finite height]\label{ctr:height}
On the actual complete positive source \eqref{ctr:gap}, retained on one
full maximizing lift, the central crossing is strictly interior.
If its physical cost satisfies $c\le1/10$, its height obeys $V^2>14/5$.
In particular the original hypotheses $c<73/1000$ and the actual
compact-set crossing cut imply $V^2>2$.
\end{lemma}
\begin{proof}
Write $z=(1-u)/v$, so $xy+xz+yz=1$. The internal complete sources
$2012$ and $0120$ retain their positive outer caps. Whole phases $+2$
and $+1$, respectively, carry them to $1201$ on the same full lift.
Lemma~\ref{fin:forward-force} therefore gives
\[
 y>z_C,\quad x_C>z_C,\qquad y>x_D,\quad z_D>x_D.
\]
Hence $u_C>0>u_D$. Put $Y=y^2$. Since $z_C=(1-ay)/(a+y)$,
these orders also give $a>(1-Y)/(2y)$ and $b<y$.
For $Y<1/3$, the interval $[y,(1-Y)/(2y)]$ lies strictly inside
$[b,a]$, so \eqref{ctr:density} gives
\begin{equation}\label{ctr:geometric-height}
 c>1-3Y-Y\log\frac{1+Y}{4Y}
 >\mathcal R(Y):=\frac{(1-3Y)(3Y+7)}{8(Y+1)}.
\end{equation}
Here $\log g<(g-g^{-1})/2$ for $g>1$: the difference vanishes at
$1$ and its derivative is $(g-1)^2/(2g^2)>0$. Also
$\mathcal R'(Y)=-(9Y^2+18Y+25)/[8(Y+1)^2]<0$.
At the crossing $q=V^2>1$ gives $Y=(q-1)^2/(4q)$, increasing in $q$.
If $Y\ge1/3$, then $q\ge3$. Otherwise $q\le14/5$ would imply
$c>\mathcal R(81/280)=1/10+647/808640$, a contradiction.
The rational bound \eqref{ctr:geometric-height} is independent of both
outward means and remains available at other explicitly checked cost cuts.
\end{proof}

\paragraph{Source-independent central transport and threshold.}
Retain the broader central cap interface above; no strict endpoint orders
from Lemma~\ref{fin:forward-force} are used in this algebra.
In the moving binary-root frames the actual central transport is
$P=H(a)^{-1}N(y(a-b)/2)H(b)$, where
\[
 H(t)=\begin{pmatrix}t/(t+y)&-2/(1+t^2)\\1/2&t(t+y)/(1+t^2)\end{pmatrix},
 \qquad N(\sigma)=\begin{pmatrix}1&0\\\sigma&1\end{pmatrix}.
\]
Write $P=\left(\begin{smallmatrix}p&n\\r&w\end{smallmatrix}\right)$. Its entries are
\begin{align*}
 p&=\frac{a^2b+2aby-b^2y+b+y}{(1+a^2)(b+y)},&
 n&=-\frac{2(a-b)(a+b+2y)}{(1+a^2)(1+b^2)},\\
 r&=\frac{y(a-b)(ab+1)}{2(a+y)(b+y)},&
 w&=\frac{-a^2y+ab^2+2aby+a+y}{(a+y)(1+b^2)}.
\end{align*}
Let $\mathsf B_C=\left(\begin{smallmatrix}L&K_L\\m&-L\end{smallmatrix}\right)$ and $\mathsf B_D=P^{-1}\mathsf B_CP$. Thus
\[
 K_R=w^2K_L+2nwL-n^2m,\quad
 L_D=(pw+nr)L+rwK_L-pnm,\quad m_D=p^2m-2prL-r^2K_L.
\]
Here $K_i=(3+\mu_i)/2$. Set
$K^*=(H(a)\mathsf B_CH(a)^{-1})_{12}$,
$\Delta U=(y-b)/(y+b)-(y-a)/(y+a)$, and
$F_D=2y/((y+b)(u_D-1))$. The integrated normalized central $\chi$
equation is the affine cost row
\begin{equation}\label{ctr:affine-cost}
 c=\frac{y(a-b)}2 K^*+\frac{\Delta U}2-\frac{F_D\chi_D}2.
\end{equation}
The lawful cyclic chart sending the centre to mode $0$ has initial
$\chi=0$, final $\chi=F_D\chi_D$, and shear $y(a-b)/2$.
Thus this is exactly \eqref{cmp:action-transport}, retaining the state
adjoint as well as matrix transport.

Put $m_0=(w^2K_L-K_R)/n^2$ and
\[
 \mathcal W=(a-b)y^2+y(2ab-b^2+1)+a^2b+a>0.
\]
Let $c_0^{\rm aff}$ be the right side of \eqref{ctr:affine-cost}
at $L=0,m=m_0$, using the cap equation for $\chi_D$. The exact solution is
\[
 \gamma=\frac{2y\mathcal W}{(a-b)(1+y^2)(a+b+2y)}>0,
 \quad L=(c-c_0^{\rm aff})/\gamma,\quad m=m_0+2wL/n.
\]
Substitution in \eqref{ctr:cap-m} and \eqref{ctr:switch-rows} gives affine
rows $S_L=l_L\mu_L+l_R\mu_R+l_cc+l_0$ and
$S_R=r_L\mu_L+r_R\mu_R+r_cc+r_0$, where
\[
 l_R=-\frac{(a+y)(1+a^2)(1+b^2)^2(a^2+2ay-by)}
 {8a(a-b)(b+y)(a+b+2y)\mathcal W}<0,
\]
\[
 r_R=\frac{(a+y)^3(1+b^2)^3(1-ay)}
 {8(a-b)(b+y)^2(1-by)(a+b+2y)\mathcal W}>0.
\]
The nonnegative combination $r_RS_L-l_RS_R$ eliminates $\mu_R$;
its $\mu_L$ coefficient is
\[
 -\frac{y(a+y)(1+a^2)^3(1+b^2)^3}
 {64a(a-b)(b+y)^2(1-by)(a+b+2y)\mathcal W}<0.
\]
Direct cancellation therefore gives the reusable threshold
\begin{equation}\label{ctr:threshold}
 \mu_L\le A(c-c_0),\qquad
 A=\frac{4(1+y^2)(a+y)^2}{y^2(1+a^2)^2}>0.
\end{equation}
There is no division by either crossing distance $a-k$ or $k-b$.
With $\varepsilon=y/h$, $r_1=(a+y)/h$, $s_1=(b+y)/h$ and $d=r_1-s_1$,
the threshold and the exact integral \eqref{ctr:density} are
\begin{align*}
 c_0=\frac{\varepsilon}{1-\varepsilon^2}\bigg[
 &d+s_1^{-1}-r_1^{-1}-\frac\varepsilon4(r_1-r_1^{-1})^2\\
 &-\varepsilon^2d(2+r_1^{-2})
 +\varepsilon^3(1-s_1^2/r_1^2)\bigg],\\
 c=\frac\varepsilon{1-\varepsilon^2}
 &(d+s_1^{-1}-r_1^{-1}-\varepsilon\log(r_1/s_1)).
\end{align*}
Equivalently,
\[
 \mu_L\le-\frac{4r_1^2\mathcal B}{(1+r_1^2-2\varepsilon r_1)^2},\qquad
 \mathcal B=\log(r_1/s_1)-\tfrac14(r_1-r_1^{-1})^2
 -\varepsilon d(2+r_1^{-2})+\varepsilon^2(1-s_1^2/r_1^2).
\]
These exact scalar identities require no source-specific height conclusion.
For the conditional low-height box $1<V^2\le2$, $y\ge1/15$,
$b\le k\le a$ and $0<\rho<1/20$, they imply the original bounds
\begin{equation}\label{ctr:dimensionless}
 1/16<1/\sqrt{226}\le\varepsilon\le1/3,
 \qquad 0<d<\frac{1-\varepsilon^2}{20\varepsilon}\le255/320<1.
\end{equation}
Indeed $d=\rho/(yh)=\rho(1-\varepsilon^2)/\varepsilon$.
\paragraph{The broader numerical consequence.}
The same normal cap/transport interface, without the complete seven-source,
still gives $V^2>2$ if $u_C\ge0\ge u_D$, $y\ge1/15$ and $c<73/1000$.
To prove this, suppose $1<V^2\le2$ and put $t=yk$,
$Y=y^2=t^2/(1-2t)$, so $3/50<t\le1/4$.
If $\rho\ge1/20$, choose a subinterval containing $k$ of length $1/(20y)$.
Its midpoint is at most $k+1/(40y)$, which lies below $1/y$.
The density is decreasing and strictly convex on $(0,1/y)$, so Jensen gives
$c>f(k+1/(40y))/(20y)>73/1000$. For the last inequality put
$w=t+1/40+Y$, $z=t-3/50\in[0,19/100]$. Its positive-denominator
numerator has the exact nonnegative decomposition
\begin{align*}
 &\bigl(-\tfrac{23}{50}w^2-Yw+Y(1+Y)\bigr)(1-2t)^2\\
 &\quad=\frac{77z^4}{50}+\frac{573159z}{50000000}+\frac{11141}{156250000}\\
 &\qquad+\left(\frac{19}{100}-z\right)
       \left(\frac{2133z^2}{1250}+\frac{40891z}{500000}\right)>0.
\end{align*}
Hence $\rho<1/20$ and \eqref{ctr:dimensionless} holds.
Write $\theta=(r_1-1)/d\in[0,1]$. Full reflection and time reversal
fixing mode $1$ preserve cost, height and the cap interface, exchange
the genuine means, and give
\[
 a'=\frac{1-by}{b+y},\quad b'=\frac{1-ay}{a+y},\quad
 r_1'=s_1^{-1},\quad s_1'=r_1^{-1},\quad
 \theta'=(1-\theta)(1+\theta d).
\]
If $\theta>2/3$, then $d<1$ gives $\theta'\le1-\theta^2<5/9$.
Reapplying the width bound in the reflected chart retains
\eqref{ctr:dimensionless}; thus arrange $0\le\theta\le2/3$.
The elementary bounds
\[
 \frac{\log(r_1/s_1)}d>\frac2{r_1+s_1}\ge1-d/6,\qquad
 \frac{(r_1-r_1^{-1})^2}{4d}\le\theta^2d\le4d/9,
\]
\[
 2+r_1^{-2}\le3,\qquad
 \frac{1-s_1^2/r_1^2}{d}\ge\frac{2-d}{(1+2d/3)^2}
\]
give $\mathcal B/d>F(\varepsilon,d)$, where
\[
 F(e,d)=1-11d/18-3e+e^2\frac{2-d}{(1+2d/3)^2}.
\]
This decreases in $d<1$, since the last fraction has derivative
$(2d-11)/[3(1+2d/3)^3]$. At $D(e)=(1-e^2)/(20e)$,
\[
 F(e,D(e))=\frac{N(e)}{360e(e^2-30e-1)^2},
\]
\[
 N(e)=(1-e)\{e^3(270142-727631e-15131e^2)
                  +e^2(10342-311/e-11/e^2)\}>0.
\]
Indeed on $[1/16,1/3]$ the two parentheses are at least
$233254/9$ and $2550$, respectively. Thus $\mathcal B>0$, and the
exact threshold forces $\mu_L<0$, contradicting its genuine positivity.
No strict crossing-distance hypothesis was used.

\subsubsection{Common central geometry and endpoint profiles}\label{ctr:profiles}
The identities here retain one complete jointly normalized full lift, both
actual cap contacts and the maximizing signs \eqref{ctr:cap-signs}. Throughout,
$v>0$, $|u|<1$, $S=u^2+v^2-1>0$, and the genuine endpoint means lie in $(0,1)$. They
apply to each source supplying that interface; no signed diagonal entrance,
HIGH condition, determinant sign or harmonic calibration is assumed.
Place physical time zero at the crossing and write the actual central interval as $[-\alpha,\beta]$, with $\alpha,\beta\ge0$. Put
\[
 y=\frac{q-1}{2\sqrt q},\quad x=t/y,\quad
 D(x)=1+(q+1)\sinh^2(x/2).
\]
Direct substitution in the original state equations gives
\begin{equation}\label{ctr:centered-state}
 u(t)=-\frac{(q+1)\sinh x}{2D(x)},\qquad
 v(t)=\frac{\sqrt q}{D(x)},\qquad
 \ell(t)=\frac{q+1}{\sqrt q}\cosh(t/y)-y.
\end{equation}
The odd cost primitive is
\[
 \mathcal J(t)=\frac{y(q+1)}{\sqrt q}\sinh(t/y)-yt.
\]
Thus the two partial costs $c_L=\mathcal J(\alpha)$ and $c_R=\mathcal J(\beta)$ add to $c$. Define
\[
 H(q)=\frac{4q}{(q-1)(q+3)},\quad
 r(q)=\frac{2q(q+1)}{(q-1)(q+3)}.
\]
Since $\mathcal J'\ge(q+3)/(2\sqrt q)$,
\[
 \alpha/y\le H(q)c_L,\quad \beta/y\le H(q)c_R.
\]
With $w=\cosh(t/y)$, the exact horizontal displacement per unit physical cost is
\begin{equation}\label{ctr:displacement-cost}
 \left|\frac{du}{d\xi}\right|
 =\frac{4q(q+1)((q+1)-(q-1)w)}
 {(q-1)((q+1)w-q+1)^2(2(q+1)w-q+1)}\le r(q).
\end{equation}
The numerator is positive on the actual state domain, since
$1-u^2+v^2=(q+1)[q+1-(q-1)w]/(2D^2)>0$.
It decreases and both positive denominator factors increase with $w\ge1$. 
\paragraph{The exact nonnegative identity.}
Let $S=u^2+v^2-1$ and define at every central point
$\mu(t)=2(2uA+B-u^2m)-3$, whose endpoint values are the genuine means. The original state adjoint and active Hamiltonian give
\[
 p_v'=\frac{\mu-u^2}{v^2}
 -\frac{4v^2(1-u^2)}{S^2}p_u
 +\frac{2uv(1-u^2+v^2)}{S^2}p_v.
\]
The two full maximizing inequalities are equivalent to
\[
 p_u\le\frac{u(1-u^2+v^2)}{2v(1-u^2)}p_v
       -\frac{S}{2v(1-u^2)}|p_v|.
\]
Let $\zeta\ge0$ be the right side minus $p_u$. The signed drift cancels exactly, and integration yields
\begin{equation}\label{ctr:nonnegative-budget}
 -I=p_v(C)-p_v(D)+\int\frac{2v}{S}|p_v|\,dt
        +\int\frac{4v^2(1-u^2)}{S^2}\zeta\,dt\ge0,
 \quad I=\int_C^D\frac{\mu-u^2}{v^2}\,dt.
\end{equation}
Every term on the right is nonnegative, including the endpoint terms.

\paragraph{An exact quadrature with its endpoint correction.}
Write the crossing matrix as $\left(\begin{smallmatrix}A_0&B_0\\m_0&-A_0\end{smallmatrix}\right)$. The original group increment from the crossing is
\[
 P(t)=\begin{pmatrix}
 1-y\sqrt q(\cosh x-1)&-y\sqrt q\sinh x\\
 y\sinh x/\sqrt q&1+y(\cosh x-1)/\sqrt q
 \end{pmatrix}.
\]
It satisfies $P'=PX$ and $\Lambda(t)=P(t)^{-1}\Lambda(0)P(t)$. Moreover
\[
 P(t)\frac{(u(t),1)^{\mathsf T}}{v(t)}
       =(-\sqrt q\sinh x,\cosh x/\sqrt q)^{\mathsf T},
\]
whence
\[
 \frac{\mu+3}{v^2}
 =\frac{2B_0}{q}\cosh^2x-4A_0\sinh x\cosh x-2qm_0\sinh^2x.
\]
Set $a=\alpha/y$, $b=\beta/y$, $X_0=a+b$, $d_0=B_0/q+qm_0$, and $F_0=q-1+1/q$. Exact integration of this three-dimensional hyperbolic span gives
\begin{align}\label{ctr:corrected-quadrature}
 \frac Iy={}&\frac{\tanh X_0}{2}
 \left(\frac{\mu_L-u_C^2}{v_C^2}
             +\frac{\mu_R-u_D^2}{v_D^2}\right)
 +(d_0-F_0)(X_0-\tanh X_0)\notag\\
 &+\frac{3(q^2-1)}{2q}
 \left(\tanh(X_0/2)-\frac{\tanh X_0}{2}\right)
 (\cosh a+\cosh b).
\end{align}
The endpoint $u_i^2$ terms are essential; dropping them would give incorrect numerical constants.

\paragraph{Crossing-matrix control without calibration.}
Put $p_i=P_{11}(t_i)>0$, $\rho_i=P_{21}(t_i)/p_i$,
$t_0=-\rho_L\rho_R\ge0$, and $z_0=\tanh a\tanh b$.
Put $\epsilon=1$; the ordinary profile retains both adverse endpoint-state terms.
The exact cap identities and conjugation give
\[
 m_0-2\rho_iA_0-\rho_i^2B_0>-a_i,\qquad
 a_i=\frac{1+\epsilon u_i^2/v_i^2}{2p_i^2}.
\]
Convex weights eliminating $A_0$ give $m_0>t_0B_0-a_{\max}$.
Positive endpoint means in the null-quadratic identity give
$B_0>q^2z_0m_0$. These weights allow a zero value at an endpoint crossing.
If $t_0q^2z_0<1$, combining the rows without a sign assumption on
$B_0,m_0$ gives
\[
 d_0>-\frac{q(1+z_0)a_{\max}}{1-t_0q^2z_0}.
\]

\begin{lemma}[Physical-cost endpoint profiles]\label{ctr:profile-budget}
Suppose $0<c\le c_0$, $q>1$ lies in a certified interval,
$z>0$, $\eta>0$, $\gamma\ge0$ and $b(q)\ge0$. Suppose the actual
geometry and cap rows give
\[
 zq\le v^2\le q,\quad |u_i|\le r(q)c_i,\quad
 d_0>-\kappa q,\quad v_i|\beta_i|\le b(q)\sqrt q,
\]
\[
 \tanh X_0\ge\eta X_0,\quad X_0-\tanh X_0\le X_0^3/3,\quad
 \tanh(X_0/2)-\tfrac12\tanh X_0\ge\gamma X_0^3/8.
\]
With the same actual $q$ throughout, define
\begin{align}\label{ctr:profile-constructors}
 K&=\frac{(1+\kappa)q-1+q^{-1}}3-\frac{3\gamma}{8}(q-q^{-1})>0,\notag\\
 M&=\frac{2qKH^2}{\eta}+\frac{r^2}{z},\qquad
 N=b(q)\left[\frac{q-1}{2}KH^3+\epsilon\frac{r^2}{z(q+3)}\right].
\end{align}
Then $\mu_L+\mu_R\le Mc^2$ and $e_L+e_R\le Nc^3$, where
$e_i=\chi_i\beta_i\ge0$. No signed diagonal entrance or calibration is needed.
\end{lemma}
\begin{proof}
The corrected quadrature, $\cosh a+\cosh b\ge2$, and the cap bound give
\[
 I/y\ge\tfrac12\tanh X_0\sum_i(\mu_i-u_i^2)/v_i^2-KX_0^3.
\]
As $I\le0$, $X_0\le Hc$ and
$\sum_i u_i^2/v_i^2\le r^2c^2/(zq)$, this gives the stated $M$.
Retain both adverse endpoint-state terms:
\begin{equation}\label{ctr:profile-I}
 0\le-I\le\mathfrak B(q)c^3,\qquad
 \mathfrak B=yKH^3+\epsilon\frac{r^2}{z\sqrt q(q+3)}.
\end{equation}
Here $\tanh X_0\le X_0$ and $yH/(2q)=1/(\sqrt q(q+3))$.
The nonnegative identity bounds the aggregate
$|p_v(C)|+|p_v(D)|$ by $-I$; multiplying once by $b(q)\sqrt q$
proves the aggregate $N$ bound, using $e_i=v_i|\beta_i||p_v(i)|$. No physical cost is paid twice.
\end{proof}

\paragraph{The ordinary opposite-centre profile.}
On the source-qualified interval $27/10\le q\le243/4$, use
\[
 c_0=63/500,\quad z=(97/100)^2,\quad \kappa=27/50,\quad
 \eta=397/400,\quad\gamma=197/200,\quad\epsilon=1,
\]
\[
 b(q)=105413/113296,\qquad
 K_{\mathrm O}=\frac{691q-1600+3373/q}{4800}>0.
\]
This is profile $\mathrm O$; its constructors $M,N$ remain exactly
\eqref{ctr:profile-constructors}.

Here are the complete elementary guards for this profile. The exact state
formulas give, for $q\ge2$,
$qH\le4$, $q(q+1)H^2/2\le8$ and $(q-1)H^2\le64/25$.
Use $(H_*,r_*,q_*,p_*)=(96/85,176/85,8/3,49/50)$. With $X_*=H_*c_0$,
\[
 X_0\le X_*,\quad |u|\le r_*c_0,\quad
 q-v^2\le\frac{8c^2}{1-16c^2/25},\quad
 P_{11}\ge1-\frac{(16/25)c^2}{1-X_0^2/2},
\]
\[
 z_0\le X_*^2/4,\quad
 t_0q^2z_0\le\frac{(64/625)c_0^4}{(1-X_*^2/2)^2p_*^2}<1/1000.
\]
The height bound follows from $q-v^2\le2q(q+1)\sinh^2(|t|/(2y))$,
$\sinh x\le x\cosh x$ and $\cosh^2x\le(1-x^2)^{-1}$; the group
bounds follow directly from the displayed $P(t)$ using
$\cosh X\le(1-X^2/2)^{-1}$ and $\sinh X\le X/(1-X^2/2)$.
At this budget the exact scalar inequalities are
\[
 \frac{8c_0^2}{1-16c_0^2/25}<(1-z)q_*,\quad
 1-\frac{(16/25)c_0^2}{1-X_*^2/2}>p_*,\quad
 \frac{(1+X_*^2/4)(1+\epsilon U_*)}{2p_*^2(999/1000)}<\kappa,
 \quad U_*=\frac{r_*^2c_0^2}{zq_*}.
\]
Thus the crossing-matrix bound is proved before calibration. For $X_0<3/20$,
$\tanh X\ge X-X^3/3$ and
\[
 \tanh(X/2)-\tfrac12\tanh X
 \ge\frac{X^3}{8}\frac{(1-X^2/12)^3}{1+X^2/4}
\]
prove the stated hyperbolic constants. Use $u_*=27/100$ in
\[
 v|\beta_j|\le\frac{\sqrt q}{2}
       \left(\frac{1+u_*}{\sqrt z\,q}+\frac1{1-u_*}\right).
\]
Evaluation at $q_*$ gives $105413/113296$.
The portable \texttt{check\_center\_profiles.py} reproduces these exact
scalar guards and the unchanged whole-interval profile certificates.

\subsubsection{Source-independent central transport and endpoint losses}
\label{ctr:generic-transport}
Retain one actual normal central arc and its genuine caps on the same complete
jointly normalized full lift, with the state domain and positive means of
Section~\ref{ctr:profiles}. The identities below need no opposite-source
word, HIGH condition, Casimir sign or calibration. Signed diagonal entrances
are required only where explicitly invoked, and must first be proved on this
actual lift.

\paragraph{Signed moving-root transport.}
In accumulated physical-cost coordinate $\xi$ the matrix
$\mathcal B=\bigl(\begin{smallmatrix}L&K\\m&-L\end{smallmatrix}\bigr)$
satisfies
\[
 \mathcal B_\xi=[\mathcal B,M_0],\qquad
 M_0=\begin{pmatrix}0&-A\\B&0\end{pmatrix},\qquad
 A=\frac{4v^4}{S(3+R^2)},\quad B=\frac2{3+R^2},\quad R^2=u^2+v^2.
\]
Let $P_\xi=PM_0$, $P(0)=I$, and
$P=\bigl(\begin{smallmatrix}p&n\\\rho&w\end{smallmatrix}\bigr)$.
On $0<\xi\le c$ with $A_0B_0c^2<2$, suppose
$0<A\le A_0$ and $0<B\le B_0$. The differential equations and a
first-zero argument give
\[
 0<p,w\le1,\quad n<0<\rho,\quad -n\le A_0\xi,\quad
 \rho\le B_0\xi,\quad p,w\ge1-A_0B_0\xi^2/2.
\]
The last lower bounds keep the diagonal entries positive on the stated
short interval, closing the first-zero argument. Since $\det P=1$,
exact endpoint conjugation gives
\[
 L_D=(pw+n\rho)L_C+\rho wK_C-pnm_C.
\]
Whenever the actual endpoint signs $L_C>0>L_D$ have been independently
established, put $D_L=L_C-L_D>0$. Since $0<L_C<D_L$ and $K_C>3/2$,
any lower cap bound $m_C>-M_C/2$ yields
\begin{equation}\label{ctr:generic-drop}
 (1+2n\rho)D_L
 <\frac{-pnM_C-3\rho w}{2}.
\end{equation}
It is used only after $1+2n\rho>0$ has been proved. The favorable $\rho w$ term is retained in the profile below.

\paragraph{Casimir identity and minimax.}
At each actual endpoint of the same lift,
\[
\begin{split}
 C_{\rm mat}={}&\frac34-L_i^2
 +\frac{\mu_i(v_i^2-3)}{4v_i^2}
 -\frac{\mu_i^2}{4v_i^2}
 +\frac{(3+\mu_i)u_i^2}{4v_i^2}
 -\frac{3+\mu_i}{2}e_i.
\end{split}
\]
Define
\begin{equation}\label{ctr:generic-cap-loss}
\begin{gathered}
 Z_i=\frac{\mu_i(3-v_i^2)_+}{4v_i^2}
     +\frac{\mu_i^2}{4v_i^2}+\frac{3+\mu_i}{2}e_i\ge0,\\
 C_{\rm mat}\ge3/4-(L_i^2+Z_i).
\end{gathered}
\end{equation}
Only the favorable $u_i^2$ term was dropped.
\begin{lemma}[Ordinary endpoint losses]\label{ctr:compensated-loss}
Let $a_1,a_2,z_1,z_2,Z\ge0$, $D>0$,
$a_1+a_2\le D$ and $z_1+z_2\le Z$. Then
\[
 \min_i(a_i^2+z_i)\le\mathcal F(D,Z):=
 \begin{cases}(D+Z/D)^2/4,&Z\le D^2,\\
 (D^2+Z)/2,&Z\ge D^2.
 \end{cases}
\]
This sharp bound is nondecreasing in both budgets, and is always at most
$(D+Z/D)^2/4$.
\end{lemma}
\begin{proof}
Saturate $a_1+a_2=D$ and $z_1+z_2=Z$. Write $a\ge b$ for the two
lengths, $e=a-b$, and $m=\min(a^2+z_1,b^2+z_2)$.
Weighted and ordinary averaging give
\[
 m\le ab+aZ/D=
 \frac{(D+Z/D)^2-(e-Z/D)^2}{4},\qquad
 m\le\frac{a^2+b^2+Z}{2}\le\frac{D^2+Z}{2}.
\]
For $Z\le D^2$ the first bound is attained by
$e=Z/D,z_1=0$; for $Z\ge D^2$ the second is attained by
$a=D,z_1=(Z-D^2)/2$. These feasible equality cases prove the stated
sharp piecewise bound. Nested feasible sets prove monotonicity.
For $Z\ge D^2$, the difference between the smooth upper bound and
$(D^2+Z)/2$ is $(Z^2-D^4)/(4D^2)\ge0$; for $Z\le D^2$
the sharp bound already equals the smooth expression.
\end{proof}

\begin{theorem}[Signed physical-cost profile]\label{ctr:signed-profile}
Retain the actual central arc, both cap contacts and the same complete
jointly normalized full maximizing lift of Section~\ref{ctr:profiles}.
Suppose a profile of Lemma~\ref{ctr:profile-budget} supplies
$M,N,z,\epsilon$ on $0<c\le c_0$, and that actual-source arguments have
independently established $L_C>0>L_D$.
Put
\begin{gather}
 A_q=\frac{4q^2}{(q-1)(q+3)},\quad B_-=\frac2{q+3},\quad
 B_+=\frac2{3+zq},\quad k=A_qB_+,\quad p_0=1-kc^2/2,\notag\\
 U=\frac{r^2c^2}{zq},\qquad
 D(q,c)=\frac{c\{A_q(1+\epsilon U)-3B_-p_0^2\}}
                    {2(1-2kc^2)}.\label{ctr:drop-profile}
\end{gather}
Require $1-2kc_0^2>0$. Then $0<L_C-L_D<D(q,c)$, and $D/c$ and $D$
increase with $c$. The ordinary endpoint-loss budgets are
\begin{equation}\label{ctr:ordinary-profile}
 Z(q,c)=\frac{Mc^2(3-zq)_+}{4zq}
       +\frac{M^2c^4}{4zq}+\frac{(3+Mc^2)Nc^3}{2},\qquad T=D^2.
\end{equation}
Here $T>0$ increases and $Z\ge0$ is nondecreasing with $c$.
For $0<\delta<3/4$, a whole-domain certificate at $c=c_0$ of
\begin{equation}\label{ctr:profile-reserve}
                    4\delta T-(T+Z)^2>0
\end{equation}
therefore implies $C_{\rm mat}>3/4-\delta>0$ on the actual lift.
\end{theorem}
\begin{proof}
The actual circle gives $R^2=1+2yv\le q$ and $v\le\sqrt q$.
Since $A=v^3/[y(2+yv)]$ increases with $v>0$, the exact generator obeys
$A\le A_q$ and $B_-\le B\le B_+$. The first-zero bounds of Section~\ref{ctr:generic-transport} give
$p,w\ge p_0>0$, $-n\le A_qc$, and
$\rho w\ge B_-p_0^2c$ by integrating $\rho'=Bw$.
Also $1+2n\rho\ge1-2kc^2>0$. The actual cap gives
$m_C>-(1+\epsilon U)/2$, so the exact signed conjugation, retaining
its favorable $\rho w$ term, proves \eqref{ctr:drop-profile}.
Since $A_q-3B_->0$, the positive numerator of $D/c$ increases and its
positive denominator decreases, proving both stated monotonicities.

The exact cap Casimir identity bounds the aggregate adverse loss by
\eqref{ctr:ordinary-profile}, using
$\sum\mu_i^2\le(\sum\mu_i)^2$ and
$\sum(3+\mu_i)e_i\le(3+\sum\mu_i)\sum e_i$.
Apply Lemma~\ref{ctr:compensated-loss} with
$a_1=L_C$, $a_2=-L_D$ and the two nonnegative cap losses.
The actual length sum is smaller than $D(q,c)\le D(q,c_0)$,
and their loss sum is at most $Z(q,c)\le Z(q,c_0)$.
Monotonicity of the sharp bound therefore permits enlargement of both
budgets to $c_0$. Its smooth upper bound is
$(T+Z)^2/(4T)$ at that budget, which the reserve makes strictly
smaller than $\delta$. No monotonicity of this smooth expression
in $T$ is assumed.
\end{proof}

\space

\space
\subsection{The common returned-terminal structure}\label{ret:section}
\begin{lemma}[The actual common returned-terminal identities]\label{ret:structure}
Retain the complete positive context $01201020120$, or a lawful whole
image, in a genuine containing fundamental interval of an original
minimizer. Its first constituent has support $[0,7]$, core $[3,5]$,
caps $2,6$, return $7$ and first-sign support $[0,4]$; its other
constituent has support $[3,10]$, core $[5,7]$, caps $4,8$, return $3$
and first-sign support $[6,10]$. All arcs are complete and maximal.
On one jointly normalized full lift, the bridge identities
\eqref{ret:domain}--\eqref{ret:matrix-bridge}, physical-cost comparisons
\eqref{ret:exclusive-cost}--\eqref{ret:J-cost}, mean comparison
\eqref{ret:Czero}, and conditional radius/shear bounds
\eqref{ret:radius-implication}--\eqref{ret:Q-bound} hold as proved below.
The radius/shear bounds require the stated actual height $0<h\le1/3$.
No numerical union floor or LOW-mean budget is assumed.
\end{lemma}
\begin{proof}

Choose one full lift on the stated containing interval. Either actual
capped binary triple supplies the two internal switches and two cap ties
of Lemma~\ref{pc:rigidity}, proving normality. Jointly normalize this lift
to $\eta=1$ before restricting to the two constituents.

Write $I$ for the shared cost of arc $5$, $E_R=c\{3,4\}$ and $E_L=c\{6,7\}$. Thus
\[
 U=E_R+I+E_L=T_R+E_L=T_L+E_R,
 \quad T_R=E_R+I,\quad T_L=E_L+I.
\]
Theorem~\ref{cmp:returned-379}, applied to the two full contexts, gives $T_R,T_L>379/1000$. They are scalar cost statements only; all subsequent differential reasoning uses the single lift just selected.

\subsubsection*{The actual two-chart bridge}
The right full context, cycled by $M=\frac12\left(\begin{smallmatrix}1&-3\\1&1\end{smallmatrix}\right)$, becomes $20120212$. Reversing the left full context and exchanging modes $0,2$ also gives $20120212$. These are the full lawful cyclic and reflected/time-reversed transports, including state, group, the state covector, all matrix entries, endpoints and winding. In the simple reflected chart
\[
 u_L=-u,\quad v_L=v,\quad
 p_L=-\operatorname{diag}(-1,1)p,\quad
 \Lambda_L=-S\Lambda S,\qquad S=\operatorname{diag}(1,-1),
\]
with reversed physical time. In both canonical charts the three-core is $202$. The complete first-sign supports and Corollary~\ref{loc:ST305-transfer} give nonpositive entrances. Lemma~\ref{loc:signed-strip} consequently supplies
\[
 -\alpha_j\le0,\quad -h_j\le r_j\le h_j<s_j<1,
 \quad \mu_j=h_j^2>0,\quad K_j=(3+\mu_j)/2,
 \quad L_j(-\alpha_j)<0.
\]
The inequalities allow weak internal jets and triple cap ties.

On the original shared mode-$0$ arc put $\xi=1+u$, $v=k\xi$ and $R=1+k^2$. Its endpoints satisfy $\xi_0<\xi_1$. In the two canonical charts the same states are
\[
 u_R=1-\frac4{R\xi},\quad v_R=\frac{4k}{R\xi},
 \qquad u_L=1-\xi,\quad v_L=k\xi.
\]
Set
\[
 p=1-s_R=\frac4{R\xi_0},\quad q=1-s_L=\xi_1.
\]
The actual positive duration and root signs imply
\begin{equation}\label{ret:domain}
 0<p,q<1,\quad Rpq>4,\quad R>4,
 \quad0<\mu_R<(1-p)^2,\quad0<\mu_L<(1-q)^2.
\end{equation}
These are exactly the two charts of \eqref{fin:shared-matrix}.
The original shared invariant is $K_L$, the terminal $0\to2$ root has
$\chi(\xi_1)=0$, and the initial $1\to0$ cap gives
$R\chi=2k(\ell-\operatorname{tr}(\Lambda X))$.
Thus the common transfer applies to this same matrix, using only retained
endpoint data; its state-adjoint row is the action identity
\eqref{cmp:action-transport}.
Define
\[
 A=\frac R{R-1},\quad z=\log(R/4),\quad
 S_\mu=\mu_L+\mu_R,\quad M_\mu=\mu_L/q+\mu_R/p,
\]
\[
 g(x)=\log x+(x^{-1}-x)/2,\quad V=z+g(p)+g(q),
 \quad B_0=(1-AV)/2,\quad B=B_0+AM_\mu/4.
\]
The diagonal constant in \eqref{fin:shared-matrix} is exactly $B$. Its two endpoint evaluations therefore give the same-lift identities
\begin{equation}\label{ret:matrix-bridge}
 \mathcal C=\frac{A^2(3+\mu_L)(3+\mu_R)}{16}-B^2,
 \quad L_R(s_R)=\frac{K_R}{(R-1)p}-B,
 \quad L_L(s_L)=\frac{K_L}{(R-1)q}-B.
\end{equation}

\subsubsection*{Three exact scalar comparisons}
On either exclusive two-arc piece let $Q=v^{-2}$. The mode equations give
$(\log Q)_u=2/(1-u)$ on mode $2$ and $-2/(1+u)$ on mode $0$.
Integrating backward from $Q(s)=1/((R-1)(1-s)^2)$ shows
\[
 Q(u)\ge\frac1{(R-1)(1-u)^2},\qquad 0\le u\le s.
\]
There is a strict comparison on a nonempty part of the actual mode-$0$ middle arc. Retaining the width contribution $\alpha_j/2$ below $u=0$, and integrating the positive-side density, gives
\begin{equation}\label{ret:exclusive-cost}
 E_j>\alpha_j/2+J_R(s_j),\qquad J_R(s):=
 \frac{Rs/2+\log(1-s)+2/(1-s)-2}{R-1}.
\end{equation}
In $p,q$ coordinates the shared cost is
\begin{equation}\label{ret:shared-cost}
 I=\frac{R(p+q)/2-2(p^{-1}+q^{-1})-\log(Rpq/4)}{R-1}.
\end{equation}
Adding it to $J_R(1-p)+J_R(1-q)$ cancels every endpoint term. Thus
\begin{equation}\label{ret:F-cost}
 \begin{gathered}
 U>F(R)+(\alpha_L+\alpha_R)/2,\\
 F(R)=\frac{R-4-\log(R/4)}{R-1},\qquad
 F'(R)=\frac{2+1/R+\log(R/4)}{(R-1)^2}>0.
 \end{gathered}
\end{equation}
No physical arc has been counted twice.

A second comparison retains the entrance jets. If $\Sigma_j$ is an exclusive shear, then $L_j(s_j)<K_j\Sigma_j$ and
$2E_j=\alpha_j+s_j+3\Sigma_j+\int u^2ds$. Hence
\[
 E_j>(\alpha_j+s_j)/2+3L_j(s_j)/(2K_j).
\]
Insert \eqref{ret:matrix-bridge} and add the shared cost once. With
\[
 \delta=\frac{\mu_L}{3+\mu_L}+\frac{\mu_R}{3+\mu_R}\in(0,1/2),
 \qquad W=V-AM_\mu/2,
\]
exact cancellation gives
\begin{equation}\label{ret:J-cost}
 U>J+(\alpha_L+\alpha_R)/2,\qquad J=W+B\delta
 =\frac1A-\frac{A-1}{2}M_\mu+B(\delta-2/A).
\end{equation}
This identity is valid before any sign for $B$ is established.

Finally, on the entire closed mean rectangle from \eqref{ret:domain}, the expression in \eqref{ret:matrix-bridge} increases in both means. For example
\[
 \frac{16}{A}\mathcal C_{\mu_L}
       =A(3+\mu_R)-8B/q.
\]
Its coefficients in $\mu_L,\mu_R$ are $-2A/q^2$ and $A(1-2/(pq))$, both negative. At their upper corner the expression is
\[
 \widehat D_L=A\{A_0+(4/q)\log(R/R_0)\}-4/q,
 \quad R_0=4/(pq),\quad A_0=(4-pq)(2-p)/q.
\]
At $R_0$ it equals $4(1-p)/q\ge0$. Its derivative is
\[
 \frac{4G}{q(R-1)^2},\qquad
 G=R-1-\log(R/R_0)-qA_0/4.
\]
Here $G'=1-1/R>0$ and
$G(R_0)=(4-pq)(1/(pq)-(2-p)/4)>0$. Since actual $R>R_0$, both mean derivatives are strictly positive. Therefore
\begin{equation}\label{ret:Czero}
 \mathcal C\ge C_{00}:=9A^2/16-B_0^2.
\end{equation}
This is an algebraic comparison at fixed actual $R,p,q$; it selects no second physical lift.
Since $g(1)=0$ and $g'(x)=-(1-x)^2/(2x^2)\le0$,
\[
 z\le V\le R/8-2/R.
\]
For the upper bound put $w=pq$ and use $p+q\le1+w$, obtaining
$V\le\log(Rw/4)+(1-w^2)/(2w)$; this decreases in $w>4/R$.

\subsubsection*{A scalar radius bound from the actual two-root equations}
For $0<h\le1/3$ define $d(h)=8h/(4-3h)$. We first prove the scalar implication
\begin{equation}\label{ret:radius-implication}
 -h\le b\le h<t<1,\quad
 \int_b^t\frac{u^2-h^2}{(1+u)^2}\,du\le0
                    \quad\Longrightarrow\quad t<d(h).
\end{equation}
The scalar radius argument in the proof of Lemma~\ref{loc:leaf-correlations}
shows that
$F_h(d(h))>0$, where $F_h(t)=\int_{-h}^t(u^2-h^2)(1+u)^{-2}\,du$,
and that $F_h$ increases for $t>h$.
Extending the lower endpoint in \eqref{ret:radius-implication} to $-h$
only adds a nonpositive integral. Thus $F_h(t)\le0$ excludes
$t\ge d(h)$, including weak input inequalities.
This uses that scalar integral proof, not a four-core hypothesis.

For any actual height $h\le1/3$, apply this to its middle mode-$0$ arc, where $Q$ is a positive multiple of $(1+u)^{-2}$ and the zero-to-zero mean gives equality in the integral: $s_j<d(h_j)$. On its actual first mode-$2$ arc,
$\chi(r)=0$ and $\chi(-\alpha)\le0$ give
\[
 \int_{-\alpha}^r\frac{h^2-u^2}{(1-u)^2}\,du\ge0.
\]
If $\alpha\le h$ the desired bound is immediate. Otherwise the scalar substitution $w=-u$ gives \eqref{ret:radius-implication} with lower endpoint $-r$, so $\alpha<d(h)$. This substitution does not reflect the physical trajectory or choose a new lift. The radius and shear comparisons in this paragraph require $h\le1/3$; no such bound is inferred without that hypothesis.

Lemma~\ref{pf:own-shear-correlation}, with
$(a,b,z,T,q)=(\alpha,s,r,\Sigma,Q(r))$, radius $d=d(h)$ and
$a=0$ allowed by continuity, recovers the exact exclusive shear.
The actual density is $x_R=Q(s)(1-s)^2=1/(R-1)$ and
$Q(r)=x_R[(1+s)/(1-s)]^2/(1+r)^2$. Its upper bound gives
$(R-1)\Sigma<2(1+d)(d-r^2)/[(1-d)^2(1+r)^2]$.
This decreases in $r\ge-h$, with derivative sign $-2(r+d)/(1+r)^3<0$.
At $r=-h$ it gives
\begin{equation}\label{ret:Q-bound}
 (R-1)\Sigma<Q(h):=
 \frac{2h(5h+4)(3h^2-4h+8)}{(1-h)^2(4-11h)^2}.
\end{equation}

\subsubsection*{A general logarithm enclosure}
For $x>1$ put $w=(x-1)/(x+1)$ and define
\begin{equation}\label{ret:logbounds}
 L_4(x)=2(w+w^3/3+w^5/5+w^7/7),\qquad
 U_4(x)=L_4(x)+\frac{2w^9}{9(1-w^2)}.
\end{equation}
Expanding $\log x=2\int_0^w(1-t^2)^{-1}\,dt$ and bounding its
positive geometric-series remainder gives $L_4(x)<\log x<U_4(x)$.
\end{proof}

\begin{corollary}[A returned-terminal union from its nine-source]
\label{opt:union722}\label{ret:U568}
Let a genuine containing fundamental interval of an original global
minimizer contain the eleven complete positive arcs $01201020120$,
or a lawful whole image. Its once-counted physical charge satisfies
\[
 U=c_3+c_4+c_5+c_6+c_7\ge\frac{180483}{250000}=2P,
 \qquad P=\frac{180483}{500000}.
\]
In particular $U>71/125+k/100$ for every $k\in\{0,1,2\}$,
without a LOW-mean hypothesis.
\end{corollary}
\begin{proof}
The full eleven-source has at least eleven physical arcs in its genuine
containing period. Its positions $1,\ldots,9$ have word $120102012$.
The whole cyclic phase $+2$ turns this into $012021201$.
This nine-source is proper, and its local charge $2,\ldots,6$ is
exactly the original charge $3,\ldots,7$.
Its two returned contexts, forward entrance witnesses, complete outgoing
leaves and ties are all retained: the canonical returned contexts are
positions $(1,\ldots,6)$ and $(7,\ldots,2)$ with reflected map $1-j$,
and the two forward witnesses are $[0,3]$ and $[5,8]$.
Choose its genuine containing interval and full maximizing lift before
restriction, proving normality before joint normalization.
The local conclusion of Lemma~\ref{six:lem:H3} and the contradiction
proved in Section~\ref{six:sec:first-contact}, namely
\eqref{six:eq:union-input}, gives $U\ge2P$.
Only this physical scalar inequality is exported; no adjoint is identified
with a lift chosen on another containing interval.
Finally
\[
 2P-\left(\frac{71}{125}+\frac2{100}\right)
 =\frac{33483}{250000}>0,
\]
which proves every older strict bound, including its former LOW cases.
\end{proof}

\space

\subsubsection{Short full sources and the actual seven-arc circle transfer}
\begin{corollary}[Short forward-four attachments]\label{ek:short-sources}
An actual proper positive $012021$, with all six arcs retained, supplies
$\cost(2)+\cost(3)+\cost(4)>61/200$.
An actual proper positive $0120212$, with all seven arcs retained, supplies
$\cost(2)+\cost(3)+\cost(4)>379/1000>3/8$.
Here proper means that the source length is less than the basic-period arc
count, so a cut inside the preceding positive arc places the entire source
inside a genuine complete basic interval.  All lawful whole cyclic and
reflected-time images have the transported physical charges.
\end{corollary}
\begin{proof}
Lemma~\ref{cmp:forward-prefix-forcing} on the full six-arc
subsource gives its strictly negative tail entrance and entrance circle.
The suffix $12021$ has its full capped $202$ at positions $2,3,4$
of the six-source. Theorem~\ref{loc:ST305} therefore gives the
$61/200$ bound on those positions.
In the seven-source the final positive arc~2 is additionally present,
giving the actual returned context $120212$. The actual negative
entrance therefore permits Theorem~\ref{cmp:returned-379}.  The caps, return and first-sign support
remain uncharged.  The complete lift is selected before these restrictions,
and its actual capped core establishes normality before normalization.
The reflected-time action reverses chronology and carries state, group,
all covectors and winding, rather than permuting labels in isolation.
\end{proof}

\begin{lemma}[The actual seven-arc entrance circle]\label{ek:H7-circle}
On an actual proper complete positive source $2021202$, retain the
original hard-word, full-lift and normality qualifications of the seven-arc
argument. Its terminal core~4,5,6 has an actual negative entrance $-a$
and satisfies $v(-a)^2>(1-a)(3+a)$. This implication has no numerical
cost budget or additional middle-mean hypothesis.
\end{lemma}
\begin{proof}
Use the full seven-arc chart $0201020$ and the unsorted parameters
$(h,g,E)$ of its first $F$ core $1,2,3$, with $\nu=Eh$.
The seven-arc row of Section~\ref{lp:returns} has exactly this support,
so Lemma~\ref{fin:H150} gives $\nu>\max\{1,g\}$.

The full reflected-time map is $u\mapsto-u$, $v\mapsto v$,
$j\mapsto2-j$, $\Lambda\mapsto-J\Lambda J$ and
$(p_u,p_v)\mapsto(p_u,-p_v)$, where $J=\diag(1,-1)$.
It reverses the seven indices and carries the original middle exit
to the actual entrance of the terminal core~4,5,6.

The actual first-$F$ central rate $\kappa_F$, distinct from the shared
rate $k$, gives
\[
 g=\frac4{(1+\kappa_F^2)\nu},\qquad
 a=\nu-1,\qquad v_K^2=\kappa_F^2\nu^2=\frac{4\nu}{g}-\nu^2.
\]
Thus $a>0$ and $v_K^2-(1-a)(3+a)=4(\nu-g)/g>0$. Every mean, endpoint and transported matrix lies
on this same seven-arc support; no harmonic calibration or numerical
cost hypothesis is used.
\end{proof}

\begin{corollary}[Proper-nine returned attachment]\label{ek:H7-K}
An actual proper positive source $202120212$, retaining all nine arcs
on one qualified complete lift, has
\[
 \cost(4)+\cost(5)+\cost(6)>379/1000>3/8.
\]
Its full reflected-time image $010201020$ has the same bound on local
positions $2,3,4$.  In the first orientation caps~3,7 and return~8, together
with sign-source arcs 0,1,2, are uncharged; in the second, caps~1,5,
return~0 and sign-source arcs 6,7,8 are uncharged.
\end{corollary}
\begin{proof}
The seven-arc source supplies the actual negative entrance; its full
returned attachment permits Theorem~\ref{cmp:returned-379}. The actual
final cap and returned arc are still needed: the seven-source alone
does not supply the returned-tail conclusion.  Normalization, means and all chart changes
refer to the same complete lift.  Whole legal transports preserve cost
and the stated physical masks.
\end{proof}

\subsection{Returned unions and a physical bridge}\label{sec:tail-bridge}
\subsubsection{Adjacent returned four-cores}
\begin{theorem}[Adjacent returned union]\label{ua:AA765}
On a full positive eleven-arc support $20212101020$, or either additive phase, the seven physical charges $2,\ldots,8$ satisfy
\[J>153/200.\]
\end{theorem}
\begin{proof}
The complete subsource $[1,9]=021210102$ supplies precisely the two
same-lift canonical four-cores and shared-arc geometry in the proof of
Theorem~\ref{loc:AA640}. Their charges are now $[2,5]$ and $[5,8]$.
The full eight-sources $[0,7]$ and $[3,10]$, under whole phases $+1,+2$,
are $01020212$. They are proper inside the eleven-source, so
Section~\ref{pf:M101} gives both costs $>E=74/125$.
Suppose $J\le V=153/200$. Each exclusive triple then costs
$<e=V-E=173/1000$, and the shared cost is $B>2E-V=419/1000$.

Use the selected complete lift, normalized by its actual four-core, and
retain the notation of that proof. Its leaf, shear, apex
and calibration arguments give
\[
 K,\widehat K<k=(3+4e^2/9)/2=6779929/4500000,\quad T_L<T=e/2,
\]
\[
 C_{\rm low}:=3/4-9V^2/64<\mathcal C<C_0.
\]
Thus $z=\sqrt{3/4-\mathcal C}$ satisfies
$z_0=243/2000<z<z_1=3V/8=459/1600$.
Set $x_0=2/25$, $X=9/50$. The exact comparisons
\begin{align*}
 2E-V-\mathcal R(X)&=877/19000>0,\\
 C_{\rm low}-k^2(1+x_0)^2/4
 &=1803636889503/312500000000000>0.
\end{align*}
give $x_0<x<X$ by \eqref{loc:shared-A-C}--\eqref{loc:shared-A-cost}.
The actual left branch~\eqref{loc:shared-left-branch} gives
$\lambda<-(3/2)x-z+kT$. Its right side is negative, since
$(3/2)x_0+z_0-kT=1000572283/9000000000>0$.
Consequently $\mathcal C<Q(x,z):=k^2(1+x)^2/4-((3/2)x+z-kT)^2$.
Here $Q_x$ decreases in both variables and
\[
 Q_x(X,z_1)=120789381432419/2025000000000000>0.
\]
Put $w=(3/2)X-kT=1257072283/9000000000>0$. On the whole rectangle,
\begin{align*}
 \mathcal C-Q(X,z)&=3/4+w^2+2wz-k^2(1+X)^2/4\\
 &>\frac{358129395870914563}{27000000000000000000}>0,
\end{align*}
where the rational bound is the value at $z=z_0$.
This contradicts $\mathcal C<Q(x,z)\le Q(X,z)$.
Both returned sources retain their actual caps and chronological returns;
only their scalar costs were imported. All shared-matrix comparisons use
one full lift, and the shared physical arc is charged once.
\end{proof}

\subsubsection{An exact bridge lemma}
\begin{lemma}\label{ua:bridge}
Suppose $0<\beta<1$, $d>0$, and $d^2\ge4(1-\beta)/(3+\beta)$. A decreasing-$u$ mode-1 arc containing $u\in[0,\beta]$ on
\[v=(d+\sqrt{d^2+4(1-u^2)})/2\]
has cost at least $47\beta/100$, for running cost $(3+u^2+v^2)/(2v)$ and dynamics
\[\frac{du}{dt}=-\frac{v(1-u^2+v^2)}{u^2+v^2-1}.\]
\end{lemma}
\begin{proof}
Put $D=d^2$, $R=\sqrt{D+4}$ and $q=\sqrt{D+4-4\beta^2}$.
The positive cost density, with $r=d/v\in(0,1)$, and its integral are
\[
 F(u,D)=\frac{r[4(1-r)/(1-u^2)+r]}{2(2-r)},\qquad
 I(D,\beta)=\frac{2d\beta}{d+q}-\frac D2 L,\quad
 L=\operatorname{arctanh}\frac{\beta(R-d)}{R+q}.
\]
The primitive is verified by differentiation, using $v_u=-2u/(2v-d)$.
Since $1-u^2=v(v-d)$, also $F=d(dv+4)/(2v^2(2v-d))$ and
\[
 \partial_v\log F=\frac d{dv+4}-\frac2v-\frac2{2v-d}<0.
\]
Thus $F$ increases with $u$, so $I(D,\beta)/\beta$ increases with $\beta$.
If $D\ge4/3$, then $r\ge2/3$ and
$F\ge f(r):=r(4-3r)/(2(2-r))$.
Since $f''=4/(r-2)^3<0$ and $f(2/3)=f(1)=1/2$, this gives
$I(D,\beta)\ge\beta/2>47\beta/100$.

If $0<D<4/3$, put $b=(4-3D)/(4+D)\in(0,1)$.
The hypothesis is equivalent to $\beta\ge b$, and $D=4(1-b)/(3+b)$.
Since $\sqrt{D+4-4b^2}=(b+2)d$ and $R/d=2/\sqrt{1-b}$, the primitive gives
\[
 I(D,b)=\frac{2b+(1-b)\log(1-b)/2}{3+b}.
\]
The positive series $(1-b)\log(1-b)=-b+\sum_{n\ge2}b^n/[n(n-1)]$ gives
\[
 \frac{I(D,b)}b\ge\frac{47}{100}+\frac{P(b)}{3+b},\qquad
 P(b)=\frac32+\sum_{n=2}^7\frac{b^{n-1}}{2n(n-1)}-\frac{47}{100}(3+b).
\]
Here $P''(b)\ge1/6$, since
$P''(b)-1/6=b(150b^3+140b^2+126b+105)/420\ge0$.
Taylor's inequality about $t=2/3$, followed by completing the square, yields
\[
 P(b)\ge P(t)-3P'(t)^2
 =\frac{1369}{1530900}-3\left(\frac{92}{42525}\right)^2
 =\frac{2122319}{2411167500}>0.
\]
Consequently $I(D,\beta)/\beta\ge I(D,b)/b>47/100$.
Any longer bridge has additional positive cost.
\end{proof}

\subsubsection{Three full relaxed residual domains}\label{ua:loss}
\begin{lemma}
For any row below, retain all $a,b\ge0$ with $a+hb\le V$, $b\le q$,
$-\min(a,b)\le r\le b$ and $r^2\le\nu\le b^2$. Put $w=a+b$,
$k=(3+q^2)/2$ and
$\bar Z=(3+b^2)\min\{(V-a-hb)/3,J_0\}/2$.
Then $(\bar Z+U_1(\nu))^2-\nu<\delta$ and
$\bar Z^2-r^2<\delta$; the $w=0$ extension is zero.
\[
\begin{array}{c|rrrrrr}
 &V&q&h&S&j&Z\\\hline
0&43383/50000&81/500&19/10&41/100&31/200&47/200\\
1&4501/5000&1513/10000&97/50&21/50&4/25&121/500\\
2&91/100&81/500&97/50&43/100&4/25&97/400
\end{array}
\]
Here $\delta=589/10000$ in row~0 and $\delta=59/1000$ in rows~1,2.
Rows~0 and~1 retain their original domains and conclusions. In particular,
row~0 includes $b\le7283/45000<81/500$, without an actual-mean restriction.
\end{lemma}
\begin{proof}
Direct differentiation gives
$(J_0)_r=-4(1+b)^2(a+r)/[(1+a)(1-b)(3+b)(1+r)^3]\le0$.
For $a\le b$, $J_0\le J_0(-a)=w(1+b)/[(1-a)(1-b)(3+b)]$.
For $a\ge b$, $J_0(-b)$ increases with $a$.
Set $W=2S/(1+h)$. If $a+hb\ge S$, the cost branch gives
$\bar Z\le k(V-S)/3<Z$.
Otherwise, when $a\le b$ one has $w\le W$, $a\le W/2$ and
$w/(3+b)\le W/(3+W/2)$; hence
\[
 \bar Z\le\frac{k(1+q)W}{(1-W/2)(1-q)(3+W/2)}<Z.
\]
When $a\ge b$ and $a+hb\le S$, use $J_0\le J_0(S-hb,b,-b)$.
The numerator of $j-J_0(S-hb,b,-b)$ is respectively
\[
\begin{aligned}
 P_0(b)&=298746\left(b-\frac{37819}{298746}\right)^2
 +\frac{37462337}{298746}+b^3(186860-22371b+5890b^2),\\
 P_1(b)&=19108\left(b-\frac{2443}{19108}\right)^2
 +\frac{280067}{19108}+b^3(11722-1459b+388b^2),\\
 P_2(b)&=38165\left(b-\frac{35591}{305320}\right)^2
 +\frac{89054179}{2442560}
 +22000(b-1/8)^2(b+1/4)\\
 &\qquad+b^3(1469-2922b+776b^2).
\end{aligned}
\]
Their positive denominators are $200(1-b)^3(b+3)(141-190b)$,
$25(1-b)^3(b+3)(71-97b)$ and
$25(1-b)^3(b+3)(143-194b)$.
All displayed remainders are positive on the respective $[0,q]$;
thus $J_0<j$ and $\bar Z<kj<Z$.
The rational comparisons for all three rows are
\[
 \frac{k(V-S)}3<Z,\qquad kj<Z,\qquad
 \frac{k(1+q)W}{(1-W/2)(1-q)(3+W/2)}<Z.
\]
The earlier bounds $\bar Z<M$ on $a\le b$ also remain in rows~0,1,
with $M=23/100,121/500$: splitting at $w=2/7$ gives respectively
$k(1+q)/[11(1-1/7)(1-q)]<M$ and $k[V-(1+h)/7]/3<M$.

Negative $U_1$ must be treated before squaring an upper response bound.
Write $I=(r+a)/[(1-r)(1+a)]$ and
$D_0=I+(1+r)(b-r)/[(1-r)^2(1+b)]$. Then
\[
 D_0-\frac a{1+a}
 =\frac{b+2br-(b+2)r^2}{(1+b)(1-r)^2}\ge0.
\]
The numerator is concave on $[-b,b]$, with endpoint values
$b(1-4b-b^2)$ and $b(1-b^2)$; $1-4q-q^2>0$ in every row.
For $a>0$ the response identity gives
$U_1\ge-a/[2(1+a)]\ge-V/[2(1+V)]$; for $a=0$ it is nonnegative.
Consequently the negative-response squared candidate is at most
$\max\{Z^2,[V/(2(1+V))]^2\}<\delta$ in every row.
We may henceforth assume $U_1\ge0$ when squaring upper bounds.

For $\nu\le a^2$, $U_1\le\nu/(1-a^2)$.
If $a\le1/2$, $(Z+4\nu/3)^2-\nu$ decreases on $[0,q^2]$.
If $a\ge1/2$, convexity in $\nu$ reduces
$[k(V-a-hb)/3+\nu/(1-a^2)]^2-\nu$ to $\nu=0,b^2$.
The ratio $(V-a)/(1-a^2)$ decreases with $a$, so
$2b/(1-a^2)\le2(V-1/2)/(h(3/4))<kh/3$.
Hence $k(V-a-hb)/3+b^2/(1-a^2)$ decreases with $b$ and is nonnegative;
both endpoint losses are at most $[k(V-1/2)/3]^2$.
The sufficient comparisons, valid in all three rows, are
\[
 \frac83(Z+4q^2/3)<1,\quad Z^2<\delta,\quad
 \frac{2(V-1/2)}{h(3/4)}<\frac{kh}3,\quad
 [k(V-1/2)/3]^2<\delta.
\]

For $a^2\le\nu\le b^2$, retain the full-root scalar envelope from
Lemma~\ref{cmp:coefficient-responses}:
\[
 U_1\le\max(U_-,U_+),\qquad
 U_- =\frac{(\nu-a^2)(1+b)}{2w(1-a)},\qquad
 U_+ =\frac{(\nu-a^2)(1-b)}{2w(1+a)}+\frac\nu{1-a^2}.
\]
Its root algebra uses only $a,b<1$ and no actual-mean restriction.
Convexity in $\nu$ reduces the squared envelopes to $a^2,b^2$;
the first endpoint is LOW. At the second keep the common entrance $a$:
\[
 C(b)=\frac{(3+b^2)(1+b)}{2(1-b)(3+b)},\qquad
 Z_J=\frac{C(b)(a+b)}{1-a}\ge\bar Z.
\]
The derivative of $Z_J+U_-$ in $a$ is
$(1+b)[C(b)-(1-b)/2]/(1-a)^2\ge0$, since
$C(b)-(1-b)/2=4b/[(1-b)(3+b)]$.
The exact positive-reciprocal decomposition
\[
 Z_J+U_+=\frac{C(b)(1+b)+b^2/2}{1-a}
 +\frac1{2(1+a)}-C(b)-\frac{1-b}{2}
\]
is convex in $a\in[0,b]$. Thus both sums are bounded by the endpoint
maximum $b\max(E_0(b),E_1(b))$, where
\[
 E_0=C+(1+b)/2,\qquad E_1=2C/(1-b)+b/(1-b^2).
\]
Each $E_i$ increases from $1$: $C$ is the product of the increasing
positive factors $(3+b^2)/2$, $(1+b)/(3+b)$ and $1/(1-b)$.
For $Q=81/500\ge q$ the exact rational comparisons give
\[
 q^2(E_i(q)^2-1)\le Q^2(E_i(Q)^2-1)<11/200<589/10000\le\delta
 \qquad(i=0,1).
\]
This proves HIGH on all three full domains, including its joins with LOW.
Finally $\bar Z^2-r^2<Z^2<\delta$ throughout; all root/mean endpoints
and the stated zero-width extension are retained.
\end{proof}

\subsubsection{A returned tail and its complete bridge}
\begin{theorem}[Returned tail with a complete bridge]\label{ua:four}
On one full normal maximizing lift of the original minimizing trajectory in a genuine containing interval, suppose a complete capped returned source $120212$ has its three-arc $202$ core beginning at $u=-a\le0$. Let $T$ be its physical cost. Let $B_1$ be the cost of the complete following mode-$1$ arc, beginning at the core's exit $u=\beta>0$. Retain the actual final mode-$2$ return. If that same mode-$1$ arc ends at horizontal coordinate at most zero, then
\[
T+B_1>\frac{91}{200}>\frac{4501}{10000}.
\]
\end{theorem}
\begin{proof}
Suppose $T+B_1\le91/200$ and put $V=91/100$.
The core cost identity gives $0<\beta\le a+\beta<2T<V<1$.
The actual return gives $v(\beta)^2>(1-\beta)(3+\beta)$.
On the following mode-$1$ arc the positive invariant
$d=(u^2+v^2-1)/v$ satisfies $d^2>4(1-\beta)/(3+\beta)$;
the complete arc decreases from $\beta$ to at most zero.
Thus Lemma~\ref{ua:bridge} applies before any bootstrap, while
Theorem~\ref{cmp:returned-379} independently gives $T>379/1000$.
Consequently
\[
 B_1\ge47\beta/100,\qquad
 \beta<38/235<81/500,\qquad a+97\beta/50<V.
\]
The first-switch response and returned shear identities require only
$a,\beta<1$, not the later numerical width restriction in the returned-tail
certificate. They bound the actual loss by the full relaxed expressions
in row~2 of the preceding residual implication.
That row gives $\mathcal L<59/1000$, contradicting
Lemma~\ref{cmp:loss-obstruction} on the actual full lift with $K>0$.
\end{proof}

\subsubsection{Seven consecutive charges}
\begin{theorem}[Seven-charge union]\label{ua:seven}
On a genuine containing interval of the original minimizing trajectory, with one full normal lift and all actual positive caps and returns retained, on the full thirteen positive arcs $0102010201020$, or an additive phase, one has
\[c_2+c_3+c_4+c_5+c_6+c_7+c_8>417/500>8291/10000>79483/100000.\]
\end{theorem}
\begin{proof}
The complete nine-sources $[0,8]$ and $[4,12]$ are both $010201020$.
They are proper since $9<13\le B$ by Lemma~\ref{loc:automatic-properness},
and stay on the given full normal lift. Corollary~\ref{ek:H7-K} gives
costs $T_2=c[2,4]>379/1000$ and $T=c[6,8]>379/1000$,
with nonpositive canonical entrances.
Apply whole reflected time reversal $j\mapsto2-j$ to this lift.
The charged order is $[8,7,6]\mid5\mid[4,3,2]=202\mid1\mid202$;
its two returned contexts are $[9,8,7,6,5,4]$ and $[5,4,3,2,1,0]$,
including caps $(9,5),(5,1)$ and returns $4,0$, respectively.
The complete bridge $5$ starts at the first tail's positive exit
and ends at the second tail's nonpositive entrance.
Thus Theorem~\ref{ua:four} gives $T+c_5>91/200$.
The disjoint physical pieces give
$T+c_5+T_2>417/500>8291/10000>81283/100000>79483/100000$.
\end{proof}

\subsubsection{A proper eleven-arc mixed source}
\begin{lemma}[A one-sided sign and its eight-arc payment]\label{ua:seven-sign}
On one full normal maximizing lift of the original trajectory, let
$h_0,\ldots,h_6$ be seven complete positive arcs with word $0201020$.
Write $(a,b,E)$ for the original, unsorted parameters of its first
three-arc core $h_1,h_2,h_3$ as in \eqref{f245:variables}.
Then $Ea>1$.  Equivalently, after full time reversal with
$j\mapsto2-j$, the word is $2021202$ and its actual state at the
boundary between positions 3 and 4 has horizontal coordinate
$1-Ea<0$.  No cap after the last arc of this reversed seven-word is
required for this state sign.

An actual proper complete source $10201020$ in an original global
minimizer has $c_1+c_2+c_3>61/200$, in local positions.
\end{lemma}
\begin{proof}
Normalize the given full normal lift jointly. The seven-arc row of
Section~\ref{lp:returns} is exactly this source, with
$(h,g,E)=(a,b,E)$: its outside, right and left means belong to
$h_1,h_3,h_5$, with roots supplied by $(h_0,h_2)$, $(h_2,h_4)$,
$(h_4,h_6)$, and shared arc $h_4$. No cap after $h_6$ is used.
Lemma~\ref{fin:H150} therefore gives $Ea>\max\{1,b\}$ directly.
The auxiliary bounds $a,b<3/2$ also remain valid without a cost budget:
if $a=1+\kappa>1$, the actual leaf has $0<d_A<a$ and
Lemma~\ref{fin:signed-mean} gives
$0<e_A<d_A^2/3-\kappa d_A$, hence $3\kappa<d_A<a$ and
$\kappa<1/2$; the right leaf gives the same bound for $b$.
Finally $h_2$ ends at $u=aE-1$. Full reflected time reversal with
$j\mapsto2-j$ carries this actual state to $1-aE<0$ at the
boundary between positions $3$ and $4$ of $h_6,\ldots,h_0$,
transporting the whole lift. No scalar sorting is used.

For $10201020$, choose one full normal lift on a genuine period
retaining all eight arcs by Lemma~\ref{as:placement}. Positions $1:7$
form $0201020$, so the exit of position $3$ has coordinate $Ea-1$.
Reverse the whole interval with $j\mapsto2-j$. Its ordered positions
$(4,3,2,1,0)$ form the complete $12021$, with actual entrance
$1-Ea<0$ at the switch from original position $4$ to $3$.
Theorem~\ref{loc:ST305} gives $61/200$ strictly on the reflected
$202$, which charges precisely original positions $1,2,3$.
\end{proof}

\space
\space
\section{Central contacts, prefix propagation, and paired responses}\label{lp:centers}
The central-arc payments are driven by actual ingress and outgoing switching jets. We integrate the complete central adjoint solution, distinguish a prefix from a complete-terminal estimate, and then derive regenerated and paired alternatives. The terminal weighting rate is variable, and its raw derivative is matched at the actual switch.

\subsection{Central costs from actual switch contacts}\label{hc:section}\label{sec:central-contact}

All arcs in this section are positive maximal arcs of an original minimizing
trajectory. Retain each stated source in a genuine containing fundamental
interval and select one full lift there. Apply Theorem~\ref{sw:normality}
before joint normalization; an isolated word is never assigned an independent
costate. For a proper source $L<B$, Lemma~\ref{as:placement} supplies the
interval by cutting inside the preceding positive arc. All chart changes
are whole cyclic symmetries or whole reflection with time reversal.

We prove two different statements. A contact-qualified prefix of small cost
has an outgoing differential sign, without a terminal-word assumption. A
complete seven-arc source has a stronger cost bound, using its actual terminal
three arcs. Only the first statement will be iterated.

\subsubsection{The complete central response}
Let the complete central arc have old mode $1$, with immediate actual
neighbors of modes $2$ and $0$. Denote its endpoints by $L,R$, its duration by
$d>0$, and its physical cost by $C_{\mathrm{loc}}$. On this arc
\[
 k=\frac{2v}{u^2+v^2-1}>0
\]
is constant. Apply the whole cyclic map which makes it mode $0$, and put
\begin{equation}\label{hc:parameters}
 \begin{gathered}
 \xi=1+u_{\rm new},\quad v_{\rm new}=k\xi,\quad
 \xi_L=A,\quad \xi_R=AE,\\
 E=e^{kd}>1,\quad \tau=\sqrt{1+k^2},\quad z=\log E=2h.
 \end{gathered}
\end{equation}
Thus $\xi'=k\xi$, and the actual open-domain restrictions are
\begin{equation}\label{hc:domain}
 AE<2,\qquad A\tau^2>2.
\end{equation}
The incoming and outgoing switches in this chart are respectively $1\to0$
and $0\to2$.

The complete family is Lemma~\ref{cmp:action}, specialized in
\eqref{cmp:ray-primitive} to
\begin{equation}\label{hc:normal-family}
 \varepsilon=\eta=1,\quad e=\xi,\quad\kappa=k,\qquad
 (K,B,M,C)=(K_0,-k^2C_1/2,C_2,C_0).
\end{equation}
The four constants remain arbitrary, and \eqref{cmp:full-reconstruction}
retains the state covector and all matrix entries, including singular
matrices. The whole cyclic map sends $p$ to $(D\Phi)^{-T}p$ and conjugates
$\Lambda$, preserving the full Hamiltonians and their differentiated
adjoint equations.

Let $q=\sigma_{0,\rm old}$ be the old switching response. After imposing the
actual incoming and outgoing contacts, two residual constants $x,d_f$ remain.
The scaled functions are
\begin{align}\label{hc:response}
 U(\xi)&=-2x\left(\frac1\xi-\frac1{AE}\right)
       +2\bigl(\log(\xi/A)-z\bigr)-(\xi-AE),\\
 V(\xi)&=d_f(\xi-A)+\tau^2\log(\xi/A)
                         +2\left(\frac1\xi-\frac1A\right),\nonumber\\
 \mathfrak R(\xi)&=k^4q
       =-\frac{\tau^2\xi-2}{\xi}U(\xi)-(2-\xi)V(\xi),
 \nonumber\\
 \xi\mathfrak R_\xi&=k^3q'.\nonumber
\end{align}
Here $U=k^2\chi_{\rm new}$ and $V=k^4Y_{\rm new}$; neither is the physical
horizontal state. The constant $d_f$ is not a duration. In this chart
\begin{equation}\label{hc:contact-signs}
 \sigma_{1,\rm new}=(2-\xi)V/k^4\le0,
 \qquad
 \sigma_{0,\rm new}=(\tau^2\xi-2)U/(k^4\xi)\ge0.
\end{equation}
The prefactors are positive by \eqref{hc:domain}. Continuity at the actual
switches gives $V(A)=U(AE)=0$. Full central maximization consequently gives
\begin{equation}\label{hc:contact-jets}
 V_\xi(A)\le0,\qquad U_\xi(AE)\le0.
\end{equation}
The second derivative is taken from the left and the first from the right,
in the increasing coordinate $\xi$. These are properties of actual
contacts, rather than constraints imposed on comparison endpoints.

The residual constants making both endpoint responses vanish are
\[
 x_*={A^2E\over2}-{AEz\over E-1},\qquad
 d_*={-\tau^2z/(E-1)+2/(AE)\over A}.
\]
Define $f_-,f_+$ as the values of $\xi\mathfrak R_\xi$ at $A,AE$ with
$(x,d_f)=(x_*,d_*)$, and put $Q=f_-+f_+$. The endpoint system is
nonsingular: its determinant with respect to $(x,d_f)$ is
$2(E-1)^2(AE-2)(A\tau^2-2)/(AE)$. Direct differentiation gives
\begin{align}\label{hc:affine-coefficients}
 \ell&={2\over A\tau^2-2}-{1\over E-1}-1,&
 m_0&={2-A\over(E-1)(2-AE)}>0,\\
 r&=1+{4-A(E^2+1)\over(E-1)(2-AE)},&
 s&=1+{2(AE\tau^2-E^2-1)\over E(E-1)(A\tau^2-2)},\nonumber\\
 n&=s+\ell={AE\tau^2-2\over E(E-1)(A\tau^2-2)}>0,&
 \rho&={r-1\over m_0}={4-A(E^2+1)\over2-A}.\nonumber
\end{align}
If $D_\pm=(\xi\mathfrak R_\xi)(A), (\xi\mathfrak R_\xi)(AE)$ and
$R_\pm=\mathfrak R(A),\mathfrak R(AE)$, the two complete affine rows are
\begin{equation}\label{hc:affine-rows}
 D_-=\ell R_-+m_0R_++f_-,\qquad
 D_-+D_+=rR_+-sR_-+Q.
\end{equation}

\subsubsection{Actual ingress, crossing, and the contact cost floor}
We use the following sufficient ingress data on the same full lift: an
actual preceding mode-$2$ tail of duration $d_0>0$, with
\begin{equation}\label{hc:ingress}
 u_L>0,\quad q_L<0,\quad
 q_L'+kq_L<\frac{b_L}{1-e^{-b_Ld_0}}q_L,
 \qquad b_L={v_L\over1-u_L}>0.
\end{equation}
The complete initial $202$ supplies these data. Its second root has
$U>0$ by Lemma~\ref{sw:ABA}. On the following complete actual mode-$2$
tail, $r=1-U<1<R(s)$ for every $s>0$ by Lemma~\ref{pc:endpoint}; hence
its full-duration response~\eqref{pc:Psi} is positive. Lemma~\ref{pc:strict-ingress}
gives \eqref{hc:ingress} on this same full lift, with its actual tail duration.
Later the same rooted-tail criterion will apply after actual regeneration.

The outgoing actual tie gives $q_R\ge0$. Set
\[
 X_-=-k^4q_L>0,\qquad Y_+=k^4q_R\ge0,\qquad
 a_0={b_L\over k(e^{b_Ld_0}-1)}+{1\over E-1}>0.
\]
Equations \eqref{hc:affine-rows} and \eqref{hc:ingress} imply
\begin{equation}\label{hc:strict-affine-ingress}
 f_-+a_0X_-+m_0Y_+<0.
\end{equation}
In particular $f_-<0$.

\begin{lemma}[Actual crossing]\label{hc:crossing}
Under \eqref{hc:ingress}, $u_R<0$ and
\begin{equation}\label{hc:crossing-interval}
 {2\over\tau E}<A<{2\over\tau},\qquad
 {65\over63}<\tau< T:={451\over30}.
\end{equation}
The crossing conclusion does not use a cost budget.
\end{lemma}
\begin{proof}
In old coordinates,
$u=(4-\tau^2\xi^2)/((2-\xi)^2+k^2\xi^2)$.
The numerator $N_-=A(E-1)f_-$ is strictly concave in $A$, since
\[
 (N_-)_{AA}=-2\tau^2(E^2-3E+2+z)<0.
\]
The bracket vanishes at $E=1$ and has derivative
$(E-1)(2E-1)/E>0$. If $u_R\ge0$, then $\tau>E$ and
$2/\tau^2<A\le2/(\tau E)$. At these two bounding values,
\[
 N_-={4(\tau^2-1)(Ez-E+1)\over E\tau^2},
 \qquad N_-={4F_0\over E^2\tau},
\]
where, with $J=(E+1)z-2(E-1)>0$,
\[
 F_0=EJ(\tau-E)^2+(E^2+E-1)J(\tau-E)
                  +(E^2-1)(Ez-E+1)>0.
\]
Here $z>2(E-1)/(E+1)$ and $Ez-E+1>0$ for $E>1$.
Concavity gives $f_->0$, a contradiction.
The actual signs $u_L>0>u_R$ yield the bounds on $A$.
At the resulting actual old-coordinate crossing, $1<v<8$, so
$k=2v/(v^2-1)>16/63$; the original compact cuts give $k\le15$.
Thus $\tau>65/63$ and $\tau\le\sqrt{226}<T$.
\end{proof}

The first actual contact jet in \eqref{hc:contact-jets} gives
$d_f\le d_c=(2-\tau^2A)/A^2$. As $R_+$ decreases with $d_f$,
\begin{align}\label{hc:contact-floor}
 Y_+&\ge Y_t=(2-AE)\left(\tau^2(E-1-z)-{2(E-1)^2\over AE}\right),\\
 f_-+m_0Y_t&={ (A\tau^2-2)(E-1)\over A}(A_t-A),
 &A_t&={2(Ez-E+1)\over(E-1)^2}.\nonumber
\end{align}
The strict ingress gives $f_-+m_0Y_t<0$, whether or not $Y_t$ is
positive. Hence
\begin{equation}\label{hc:A-floor}
 A>A_t>0.
\end{equation}

Integration of the physical density, with physical time
$dt=d\xi/(k\xi)$, gives
\begin{equation}\label{hc:physical-cost}
 C_{\mathrm{loc}}={\tau^2A(E-1)+4(1-E^{-1})/A-2z\over2(\tau^2-1)}.
\end{equation}
If $a_A=A(E-1)$ and $b_A=4(1-E^{-1})/A$, then
\[
 C_{\mathrm{loc}}={a_A\over2}+{a_A+b_A-2z\over2(\tau^2-1)}>{a_A\over2},
\]
since $a_A+b_A\ge4(\sqrt E-1/\sqrt E)>4z$.
Consequently the contact implies the strictly increasing floor
\begin{equation}\label{hc:cost-floor}
 C_{\mathrm{loc}}>W(E):={E\log E\over E-1}-1,
 \qquad W'(E)={E-1-\log E\over(E-1)^2}>0.
\end{equation}
For fixed $A,E,z$, \eqref{hc:physical-cost} strictly decreases with
$\tau^2$. None of these algebraic comparisons assigns a physical
trajectory or a budget to $A_t$.

\subsubsection{Positive series and finite polynomial certificates}
Use \eqref{hc:log-bounds} at $x=E$, abbreviating $L_N(E),U_N(E)$ as $L_N,U_N$.

The summed numerator $N_Q=AE(E-1)Q$ has
\[
 (N_Q)_{AA}=-2E\tau^2\bigl((E^2+1)z-(E^2-1)\bigr)<0.
\]
At the two crossing endpoints $2/\tau$ and $2/(\tau E)$ it is,
respectively, $32e^{3h}F/\tau$ and $32e^hF/\tau$, where
\begin{align}\label{hc:hyperbolic-force}
 F&=a(h)\tau^2+b(h)\tau+c(h),\\
 a(h)&=h\cosh h-\sinh h,\qquad c(h)=2\sinh^3h,\nonumber\\
 b(h)&=\cosh h\bigl(\sinh h\cosh h-h\cosh(2h)\bigr).\nonumber
\end{align}
Thus $F>0$ implies the actual $Q>0$.
The complete negative series is
\begin{equation}\label{hc:negative-series}
 -b(h)=\sum_{n\ge1}d_nh^{2n+1},\qquad
 d_n={(4n-1)9^n+4n+1\over4(2n+1)!},\quad
 d_1={4\over3},\quad d_2={6\over5}.
\end{equation}
For $n\ge2$, $d_{n+1}/d_n<33/98$. To check this for every index, after
clearing positive denominators and writing $n=j+2$, the two coefficients
of $9^n$ and $1$ are
\[
 528j^3+4356j^2+8022j,\qquad
 528j^3+4620j^2+12874j+11200.
\]
Both are nonnegative and the second is strictly positive for $j\ge0$.
The positive terms give
$a\ge h^3/3+h^5/30$ and $c\ge2h^3+h^5$. If the tail after $d_1$ is
bounded by $\kappa h^5$, then
\begin{equation}\label{hc:force-bracket}
 F>{h^3\over3}\bigl(D(\tau)-B_\kappa(\tau)h^2\bigr),\quad
 D(t)=(t-2)^2+2,\quad B_\kappa(t)=3\kappa t-t^2/10-3.
\end{equation}
More generally, for $0<h<H<\sqrt{98/33}$ the whole tail factor is
bounded by $(6/5)/(1-(33/98)H^2)$.

\begin{lemma}[Uniform short-factor force]\label{hc:short-factor-force}
Every actual crossing above with $1<E\le31/10$ has $Q>0$, without a
cost budget or an upper bound on $\tau$.
\end{lemma}
\begin{proof}
The logarithmic upper bound gives
\[
 \frac87-U_2(31/10)=\frac{14612893}{1495585700}>0,
 \qquad h=\tfrac12\log E<\frac47.
\]
The all-index tail bound just proved permits $\kappa=27/20$, since
\[
 \frac{27}{20}-\frac{6/5}{1-(33/98)(4/7)^2}
 =\frac{15}{8548}>0.
\]
The force bracket in \eqref{hc:force-bracket} is affine in $h^2$.
At $h^2=0$ it is $D(\tau)>0$; at $h^2=16/49$ it is
\[
 P(\tau)=\frac{253\tau^2-1304\tau+1710}{245}>0,
 \qquad
 1012\cdot245P(\tau)=(506\tau-1304)^2+30104.
\]
Thus it is positive between these endpoints, including either sign or
zero value of $B_{27/20}(\tau)$. Hence $F>0$, and concavity of $N_Q$
on the actual crossing interval proves $Q>0$.
\end{proof}

Here is the exact convention for every polynomial certificate in this
section. If a polynomial $p$ of degree $n$ is considered on $[a,b]$, define
\begin{equation}\label{hc:bernstein-generator}
 a_i=[x^i]p(a+(b-a)x),\qquad
 \beta_j(p;[a,b])=\sum_{i=0}^j{\binom ji\over\binom ni}a_i
 \quad(0\le j\le n).
\end{equation}
Then
$p(a+(b-a)x)=\sum_{j=0}^n\beta_j\binom njx^j(1-x)^{n-j}$.
The basis is nonnegative and sums to one, so positivity of all listed
coefficients proves positivity on the entire closed interval. This finite
rational formula specifies every coefficient, without a sampling rule or
an unrecorded subdivision.

\subsubsection{A contact-qualified prefix of cost at most $13/20$}
\begin{theorem}[Prefix exit estimate]\label{hc:prefix650}\label{hc:H650}
Suppose the complete actual central arc and its adjacent controls satisfy
\eqref{hc:domain}, the full maximizing contacts above, and the strict actual
ingress \eqref{hc:ingress}. If $C_{\mathrm{loc}}\le13/20$, then
\[
 u_R<0,\qquad q_R\ge0,\qquad q_R'>kq_R.
\]
No terminal word is required for this conclusion.
\end{theorem}
\begin{proof}
Put $\alpha_0=13/20$. The floor \eqref{hc:cost-floor} and the exact reserve
\[
 {E L_3(E)\over E-1}-1-\alpha_0
 \bigg|_{E=31/10}={7066595\over463424804}>0
\]
force $E<31/10$. Lemma~\ref{hc:short-factor-force} gives $Q>0$.

It remains to control $\mathfrak G=Q-\rho f_-$ when $\rho<0$.
It is affine in $\tau^2$, with coefficient $AH/(2-A)$, where
\[
 H=-A(E-1)(E^2-E+1)+2E(E+1)z-4E^2+6E-2.
\]
This decreases strictly in $A$ and increases in $z$. If $E\ge2$, then
\[
 H(A_t)={2E(E-2)(z-E+1)\over E-1}\le0,
\]
and $A>A_t$ gives $H<0$, also at $E=2$. If $1<E\le2$ and $\rho<0$,
then $A>A_r=4/(E^2+1)$; moreover
\[
 H(A_r,U_1)={(E-1)^3(E^2-5)\over3(E^2+1)}<0.
\]
Thus $H<0$ in every case in question. Since actual $\tau<2/A$,
\begin{equation}\label{hc:G-cross}
 \begin{gathered}
 \mathfrak G>{2\over E}\left({2J(z)\over A}-W_0\right),\\
 J(z)=(E+1)(E^2+1)z-(E-1)(2E^2-E+3),\\
 W_0=(E-1)^3(E+1).
 \end{gathered}
\end{equation}
Here $J>0$, because
\[
 J(L_2)={(E-1)^2(2E^3+E^2+8E+1)\over3(E+1)^2}>0.
\]
The evaluation at $\tau=2/A$ is only an affine comparison, with no
assumed state or lift at that value.

The same fixed-$A$ comparison for physical cost gives
\begin{equation}\label{hc:cross-cost}
 C_{\mathrm{loc}}>C_\times(A,z)
  ={A[2(E-E^{-1})-Az]\over4-A^2}.
\end{equation}
For $0<A<2$, this increases strictly in $A$ and decreases in $z$.
Indeed, with $R_E=2(E-E^{-1})>4z$, the numerator of its $A$ derivative is
$R_E(4+A^2)-8Az=(R_E-2z)(4+A^2)+2z(A-2)^2>0$.

Set $t=(E-1)/(E+1)$ and use the bounds $L_3<z<U_3$.
If $\mathfrak G\le0$, \eqref{hc:G-cross} gives
$A>a_g=2J(L_3)/W_0$; independently the contact gives
$A>a_t=A_t(E,L_3)$. In terms of $t$ these are
\[
 a_g={(1-t)(6t^5+16t^3-5t+15)\over30t},\qquad
 a_t={(1-t)(3t^4+3t^3+5t^2+5t+15)\over15}.
\]
Both are positive for $0<t\le21/41$. For $a=a_g$ or $a_t$, put
\[
 \Pi(a)=a[2(E-E^{-1})-aU_3]-\alpha_0(4-a^2),\qquad
 E={1+t\over1-t}.
\]
The two certificates below prove $\Pi(a_g)>0$ for
$0<t\le257/657$ and $\Pi(a_t)>0$ for
$257/657\le t\le21/41$. In either case an actual point has $0<a<A<2$;
monotonicity in \eqref{hc:cross-cost} then implies
$C_{\mathrm{loc}}>C_\times(a,U_3)>\alpha_0$, a contradiction.
No admissibility of the scalar comparator is asserted.
We have proved $\mathfrak G>0$ whenever $\rho<0$.

Finally, the right affine row is
\begin{equation}\label{hc:exit-affine}
 k^3(q_R'-kq_R)=f_++(\rho-1)m_0Y_++nX_-.
\end{equation}
For $\rho\ge1$, $f_+=Q-f_->0$ makes this positive. For $\rho<1$, write
$\delta=-f_--a_0X_--m_0Y_+>0$. Exact substitution gives
\[
 k^3(q_R'-kq_R)
 =\mathfrak G+\bigl(n+(1-\rho)a_0\bigr)X_-+(1-\rho)\delta>0.
\]
For $0\le\rho<1$ one has $\mathfrak G=Q-\rho f_->0$ directly, and the
remaining range was just proved. This establishes the exit estimate.
\end{proof}

\begin{table}[ht]
\centering\small
\begin{tabular}{llll}
Polynomial & $t$ interval & Degree & Every Bernstein coefficient exceeds\\
$378000t^2(1+t)\Pi(a_g)$ & $[0,257/657]$ & $18$ & $53$\\
$94500(1+t)\Pi(a_t)$ & $[257/657,21/41]$ & $16$ & $146$
\end{tabular}
\caption{Two complete rational prefix certificates. Substitute the displayed
$a_g,a_t,L_3,U_3,E=(1+t)/(1-t)$ and use
\eqref{hc:bernstein-generator} in $t$. The positive clearing factors produce
polynomials of the stated degrees; exact coefficient evaluation proves the
strict rational lower bounds. The face $t=0$ is a polynomial comparison
face only.}\label{hc:650-certificates}
\end{table}

\subsubsection{The actual terminal three arcs}
\begin{lemma}[Exponential terminal comparison]\label{hc:terminal}
Suppose the outgoing arc is the complete first mode-$0$ arc of an actual
positive $020$ on the same lift, and $u_R<0$. Then $q_R>0$ and
\begin{equation}\label{hc:terminal-inequality}
 q_R'<(k-a_R)q_R,\qquad a_R={v_R\over1+u_R}>0.
\end{equation}
The raw $q,q'$ are those at the actual central exit.
\end{lemma}
\begin{proof}
Let $e>0$ be the first tail duration. Its first root has $u(e)<0$ by
Lemma~\ref{sw:ABA}; since $u'=v>0$, the whole tail has $-1<u<0$.
The full normal equations give
\[
 a={v\over1+u}>0\text{ constant},\quad
 \chi''+a\chi'=-2u>0,\quad \chi(e)=0,\quad\chi'(e)\le0.
\]
Backward integration gives $\chi>0$ and $\chi'<0$. With $y=-\chi'>0$,
$y'+ay=2u<0$, and hence
\[
 0<\chi_R=\int_0^e y(t)\,dt
       <{-\chi_R'\over a}(1-e^{-ae}),\qquad
 {\chi_R'\over\chi_R}<-{a\over1-e^{-ae}}<-a.
\]
The weighting $q=w\chi$, $w=S/v^2$, has the variable logarithmic rate
\begin{equation}\label{hc:variable-rate}
 {w'\over w}=c(t)={2v\over S},\qquad
 c'=-{2v^2((1+u)^2+v^2)\over(1+u)S^2}<0.
\end{equation}
Only its endpoint value $c(R)=k$ is used. Thus
$q_R'=kq_R+w_R\chi_R'<(k-a_R)q_R$.
There is no assumption that $c(t)$ remains equal to $k$ on the tail.
The same raw response is $C^1$ at both actual central switches.
Differentiation of the full Hamiltonians gives
\begin{equation}\label{hc:C1-identity}
 D_j(q)-D_1(q)={3u^2+v^2-3\over v^2}(H_1-H_j),\qquad j=0,2,
\end{equation}
which vanishes at the respective actual $1\to0$ or $2\to1$ tie. The scalings in
\eqref{hc:response}, $\mathfrak R=k^4q$ and
$\xi\mathfrak R_\xi=k^3q'$, therefore apply to these same terminal data.
\end{proof}

\subsubsection{Eliminating the actual terminal response}
The following identities apply to any complete central arc with strict
ingress, both actual contacts and the same-lift terminal $020$. They
isolate the shared algebra used in the complete-terminal estimate.

Write $\Gamma=a_R/k=AE/(2-AE)>0$ and
$\mathcal E=k^3(q_R'-kq_R)$. Lemma~\ref{hc:terminal} gives
\begin{equation}\label{hc:terminal-scaled}
 \mathcal E+\Gamma Y_+<0.
\end{equation}
The endpoint row \eqref{hc:exit-affine} gives the exact cancellation
\begin{equation}\label{hc:cancellation}
 (\rho-1)m_0+\Gamma={1\over E-1},\qquad
 \mathcal E+\Gamma Y_+=f_++{Y_+\over E-1}+nX_-.
\end{equation}
Because strict ingress gives $f_-<0$, positivity of $Q=f_-+f_+$ is
already incompatible with the actual terminal inequality.

Independently, differentiation of \eqref{hc:response} at $\xi=AE$ gives
\begin{align}\label{hc:exit-floor}
 \mathcal E+\Gamma Y_+
 &=-(\tau^2AE-2)U_\xi(AE)
 +(2-AE)\{V(AE)-AEV_\xi(AE)\}\\
 &\ge(2-AE)\left(\tau^2z-{4(E-1)\over AE}\right).\nonumber
\end{align}
The first term is nonnegative by the actual outgoing jet. In the second
term the expression in braces is
$-d_fA+\tau^2(z-1)+4/(AE)-2/A$; the incoming jet gives
$d_f\le(2-\tau^2A)/A^2$, proving the inequality.

If the final bracket is nonpositive, put
\[
 T_b={4(E-1)\over AEz}\ge\tau^2>1,\qquad
 D_b=4(E-1)-AEz>0.
\]
The physical cost \eqref{hc:physical-cost} decreases with $\tau^2$,
so
\begin{align}\label{hc:eliminated-cost}
 C_{\mathrm{loc}}&\ge C_b(A)
 ={2A(E-1)^2-AEz^2+2(E-1)z\over D_b},\\
 C_b'(A)&={2(E-1)(4(E-1)^2-Ez^2)\over D_b^2}>0.\nonumber
\end{align}
Here $z<(E-1)/\sqrt E$, which follows from
$2\sinh(z/2)>z$. The actual incoming contact supplied $A>A_t$.
Since $D_b$ is positive on $[A_t,A]$, one has
$C_{\mathrm{loc}}>C_b(A_t)$. This is an algebraic elimination of an
actual inequality; no physical lift or cost budget is assigned to $A_t$.
The complete low- and high-factor proof using these identities is
Section~\ref{nx:central-proof}.

\begin{corollary}[A convenient central payment]\label{hc:full922}\label{hc:H922}
A proper complete positive source $2021020$ has $c_3>461/500$.
The same weaker bound holds under the sufficient interface of
Theorem~\ref{lc:direct-center}.
\end{corollary}
\begin{proof}
Theorem~\ref{lc:direct-center} gives the stronger complete-terminal
bound $299/250$, and $299/250>461/500$.
\end{proof}

\subsubsection{Regeneration on an actual ten-arc source}
The prefix theorem is useful because its exit sign regenerates strict
ingress without assigning a terminal word to the first center.

\begin{lemma}[Actual return and regenerated ingress]\label{hc:regeneration}
Suppose an actual mode-$1$ central arc satisfies the conclusion of
Theorem~\ref{hc:prefix650}, and its actual continuation begins with
$021$. Then the next mode-$1$ arc has the strict ingress
\eqref{hc:ingress}, with its own rate and incoming-tail duration.
If its actual outgoing control is $0$, the full contact assumptions of
Theorem~\ref{hc:prefix650} hold anew at this second center.
\end{lemma}
\begin{proof}
On the first following mode-$0$ arc use the original binary quantities
\[
 K_e=2ul+b-u^2m,\qquad B_e=K_e-(3+u^2)/2,\qquad
 \chi=vp_v,\qquad w=S/v^2.
\]
The full equations give
\[
 K_e'=0,\qquad B_e'=-uv,\qquad \chi'=2B_e/v,\qquad w'=2/v.
\]
At the old central exit $w'/w=k$ on the outgoing side, and $q'$ is
continuous by \eqref{hc:C1-identity}. Since $q=w\chi\ge0$ and
$q'>kq$, we have $\chi\ge0$, $\chi'>0$ and
$B_e>0$, including the possible case $q=0$.
At the next actual $0\to2$ tie, $\chi$ returns to zero. It must first
decrease somewhere, so $B_e<0$ somewhere on this actual mode-$0$ arc.
Here $u'=v>0$, and $B_e$ increases while $u<0$. The negative value
therefore occurs after $u>0$, and thereafter $B_e$ decreases. The actual
next tie has strictly $u>0$ and $B_e<0$.

On the following complete mode-$2$ arc, its actual root has $U>0$.
Again $1-U<1<R(s)$ for every $s>0$, so the full-duration response
\eqref{pc:Psi} is positive. Lemma~\ref{pc:strict-ingress} supplies \eqref{hc:ingress} with the
second center's own rate $k_2=2v_L/S_L$ and actual incoming-tail duration.
The raw response and derivative match at the actual junction by
\eqref{hc:C1-identity} and its lawful cyclic versions.

This proves ingress only. The additional contact obligations are supplied
separately by the second center's actual immediate controls $2,1,0$.
In its own whole cyclic chart these are $1,0,2$, so full maximization
and continuity give the actual ties and both one-sided signs in
\eqref{hc:contact-signs}. The duration and cost belong to its complete
actual central arc. Its crossing and compact-rate bounds are derived
anew by Lemma~\ref{hc:crossing}; no rate, mean, invariant or costate is
copied from the first center.
\end{proof}

\begin{corollary}[Ten-arc maximum and protected union]\label{hc:ten-arc}
Retain the complete actual positive source $2021021020$, indexed $0$ through
$9$, properly in one genuine containing interval. Then
\begin{equation}\label{hc:ten-max}
 \max(c_3,c_6)>13/20,\qquad c_3+c_6>13/20.
\end{equation}
If its two full reverse four-core sources have the original hypotheses of
Proposition~\ref{strip:r496}, then additionally
\begin{equation}\label{hc:ten-union}
 \sum_{j=2}^7c_j>{13\over20}+{62\over125}={573\over500}.
\end{equation}
Here $c_j$ is the physical cost of the complete actual arc $j$.
\end{corollary}
\begin{proof}
Assume both center costs at most $13/20$. The initial $202$ supplies the
first center's strict ingress. Its actual contacts qualify
Theorem~\ref{hc:prefix650}, whose exit criterion applies without
assuming an immediate terminal $020$ there. The actual return through
arcs $4,5$ regenerates the second center's ingress by
Lemma~\ref{hc:regeneration}; the actual triples $2,3,4$ and $5,6,7$
supply both centers' contact hypotheses separately. The second application
of the prefix theorem contradicts the complete actual terminal $020$
on arcs $7,8,9$ by Lemma~\ref{hc:terminal}. Hence at least one center
is expensive, proving \eqref{hc:ten-max}.

If $c_3>13/20$, the full source $[3,8]=102102$ gives the reverse
four-core payment $\sum_{j=4}^7c_j>62/125$, disjoint from arc $3$;
arc $2$ stays positive and unpaid. If $c_6>13/20$, the full source
$[1,6]=021021$ becomes $102102$ under one whole cyclic phase and gives
$\sum_{j=2}^5c_j>62/125$, disjoint from arc $6$; arc $7$ stays positive.
Each source retains its complete caps and satisfies the required
normal-lift, determinant/trace and calibration hypotheses.
Each alternative gives \eqref{hc:ten-union}. They are exhaustive branches
of one disjunction and are not added to one another. In particular,
\eqref{hc:ten-max} does not assert two unconditional singleton bounds.
\end{proof}

\begin{remark}\label{hc:scope}
Theorem~\ref{lc:direct-center} also gives the convenient complete-seven-arc costs
$C_{\mathrm{loc}}>37/50$ and $C_{\mathrm{loc}}>13/20$ whenever those
complete sources are retained. Its stronger constant cannot replace the
prefix threshold in Theorem~\ref{hc:prefix650} or in the iteration above:
its proof uses the complete actual terminal $020$ at the center to which
it is applied. Every named placement must still retain its full proper
support, actual contacts and sign witnesses, and must charge only its
stated physical arcs. Overlap of supports does not identify their
multipliers or permit repeated payment of the same physical cost.
\end{remark}

\subsection{A lower-budget consequence of actual strict ingress}
\label{pc:h592-section}

\begin{theorem}[All lower budgets and retained regeneration]
\label{hc:H592}\label{pc:lower-budgets}
Under the same complete actual central-arc and strict-ingress hypotheses
as Theorem~\ref{hc:prefix650}, every budget $0<\beta\le13/20$ gives
\begin{equation}\label{pc:h592-exit}
 C_{\rm loc}\le\beta\quad\Longrightarrow\quad
 u_R<0,\qquad q_R\ge0,\qquad q_R'>kq_R.
\end{equation}
No terminal word is needed. If a complete positive same-lift terminal
$020$ follows, then $C_{\rm loc}>\beta$. In particular these conclusions
hold at the useful budget $74/125$ and every smaller budget.
\end{theorem}
\begin{proof}
Budget inclusion gives the exit criterion by Theorem~\ref{hc:prefix650};
Lemma~\ref{hc:terminal} contradicts it when the terminal is present.
The raw derivative agrees at the actual junction by
\eqref{hc:C1-identity}. Every regenerated application still needs its own
actual ingress, complete duration, proper support and contacts; a cost
contradiction additionally needs its own same-lift terminal. Budget
inclusion supplies none of these hypotheses and turns no maximum bound
into two singleton payments.
\end{proof}

\subsection{Paired rooted responses on a common full lift}
\label{pc:section}

The estimates in this section couple two rooted tails through one actual
matrix costate.  A scalar resource can be added to another resource without
identifying their lifts.  The argument below needs more: the two matrix
reconstructions and the state covector must belong to one jointly
normalized full maximizing lift.  We first give a sufficient, explicitly
checkable way of obtaining such a lift.  The ensuing local theorem does
not use a sign of the matrix determinant, a harmonic estimate, or a cost
budget.  Its cost consequence is a bound on the \emph{maximum} of two
specified complete-arc costs.

\subsubsection{Two genuine cuts and equality-only rigidity}
\label{pc:gluing}

Write $W=(u,v)$, $S=u^2+v^2-1$, and
\[
 \begin{gathered}
 \Omega=\{(u,v):v>0,\ -1<u<1,\ S>0\},\\
 R_*={1\over2}\begin{pmatrix}1&-3\\1&1\end{pmatrix},\qquad
 \Phi(w)={w+3\over1-w},\quad w=u+iv.
 \end{gathered}
\]
The actual unfolded geometric trajectory satisfies
\begin{equation}\label{pc:geometric-twist}
 G(t+T)=G(t)R_*,\qquad W(t+T)=\Phi(W(t)),\qquad
 j_{i+B}=j_i+1\pmod3.
\end{equation}
These are geometric relations; no periodicity of a state covector is
included in them.  We write the full reduced lift as
$(\eta,\Lambda,\symbf p)$, with $\symbf p=(p_u,p_v)^t$ and
$\Lambda$ traceless.  Its equations are those of \eqref{fd:PMP}; in the
present notation,
\begin{equation}\label{pc:full-pmp}
 \begin{split}
 H_j&=\symbf p^tF_j+\tr((\Lambda-\eta A_*)X),\qquad
 A_*={1\over2}\begin{pmatrix}0&3\\-1&0\end{pmatrix},\\
 \Lambda'&=[\Lambda,X],\qquad
 \symbf p'=-DF_\rho^t\symbf p
   -\nabla_W\tr((\Lambda-\eta A_*)X),\qquad
 \max_jH_j=H_\rho=0.
 \end{split}
\end{equation}
The free-time equality, both state-adjoint equations, and the matrix
adjoint are all retained.

\begin{proposition}[A sufficient two-cut gluing criterion]\label{pc:two-cuts}
Let $a<b<a+T$ be two actual cuts of an original minimizer.  Suppose that
$[b,a+T]$ contains a complete positive capped-binary core as in
Lemma~\ref{pc:rigidity}, and that $(a,b)$ contains an actual switch of two
distinct modes.  Then the genuine full lifts on $[a,a+T]$ and $[b,b+T]$
can be jointly normalized to $\eta=1$, agree on their overlap, and extend
to one normalized full reduced maximizing lift on the actual unfolding.
Its matrix and state-covector twists are
\begin{equation}\label{pc:full-twist}
 \Lambda(t+T)=R_*^{-1}\Lambda(t)R_*,\qquad
 \symbf p(t+T)=D\Phi(W(t))^{-t}\symbf p(t).
\end{equation}
This conclusion does not assert normality of arbitrary short windows or
an enlarged-interval transversality condition.
\end{proposition}
\begin{proof}
Genuine reanchoring, preserving all original group/state endpoint
coordinates and real winding, supplies a full lift on each period.
Lemma~\ref{pc:rigidity} makes each lift normal.  Normalize the whole lift,
not just the matrix, and restore both to the same unfolded chart before
subtracting.  Constant left group reanchoring leaves the reduced pairing
unchanged.  A cyclic change of chart transports all components by
\begin{equation}\label{pc:cyclic-transport}
 W\mapsto\Phi^k(W),\quad
 \Lambda\mapsto R_*^{-k}\Lambda R_*^k,\quad
 \symbf p\mapsto D\Phi^k(W)^{-t}\symbf p.
\end{equation}
The equality-only lemma identifies the two normalized restrictions on
$[b,a+T]$.  Hence they glue on $J=[a,b+T]$.

On $J$ the spatial matrix $P=G\Lambda G^{-1}$ is constant, by
$G'=GX$ and $\Lambda'=[\Lambda,X]$.  Equation
\eqref{pc:geometric-twist} therefore proves the matrix part of
\eqref{pc:full-twist} on $J\cap(J-T)=[a,b]$.  This calculation does not
require $P$ to commute with $R_*$, and it permits a zero or repeated matrix
eigenvalue.  Tracelessness excludes a scalar-identity gauge.

Pull back the already existing shifted lift on this self-overlap:
\[
 \widehat\Lambda(t)=R_*\Lambda(t+T)R_*^{-1}=\Lambda(t),\qquad
 \widehat{\symbf p}(t)=D\Phi(W(t))^t\symbf p(t+T).
\]
The full Hamiltonian identity follows from
$D\Phi F_j=F_{j+1}\circ\Phi$,
$X(\Phi W)=R_*^{-1}X(W)R_*$ and $A_*R_*=R_*A_*$.
Differentiating this identity gives the full state-adjoint covariance,
including derivatives of $D\Phi$; thus the pullback solves
\eqref{pc:full-pmp}.  At the actual switch in $(a,b)$ both lifts satisfy
the two adjacent active zero energies by continuity.  Their matrix and
objective parts already agree, so their state-covector difference
annihilates both adjacent fields.  Directly from \eqref{fd:fields},
\begin{align*}
 \det(F_0,F_1)&={v^3((1+u)^2+v^2)\over(1+u)S},\\
 \det(F_0,F_2)&=-{2v^3\over1-u^2},\\
 \det(F_1,F_2)&={v^3((1-u)^2+v^2)\over(1-u)S}.
\end{align*}
All are nonzero in $\Omega$.  The state-covector difference is zero there,
and its homogeneous equation $\delta\symbf p'=-DF_\rho^t
\delta\symbf p$ makes it zero on $[a,b]$.  This proves the second
part of \eqref{pc:full-twist}, rather than assuming periodic costates.

Only now extend one existing full-period restriction by
\eqref{pc:full-twist}.  The proved seam matching gives continuity;
Hamiltonian covariance gives the equations and maximization on every
translate.  Initial-value uniqueness identifies the extension with the
original glued lift on $J$.  A second global normalized lift would agree
on the same rigid core, hence everywhere.  No minimization problem on
$J$ and no extra endpoint transversality have been introduced.
\end{proof}

\begin{corollary}[A short rigid source determines the global lift]
\label{pc:short-rigid-global}
Let $\tau_i$ be the actual unfolded switching times, with physical arc
$i=(\tau_i,\tau_{i+1})$. Suppose $[s,t]$ is an actual complete positive
capped-binary source as in Lemma~\ref{pc:rigidity}, with
$t-s+1\le B$. There is a unique normalized full reduced maximizing
lift on the unfolding, and it satisfies \eqref{pc:full-twist}.
Every genuine fundamental-interval lift retaining such a complete
binary core, both cap ties and positive portions of its third-mode caps
is normal and, after joint normalization and lawful alignment, is a
restriction of this lift. If two such closed physical binary cores are
disjoint, this identifies every genuine fundamental-interval lift.
\end{corollary}
\begin{proof}
Geometric period translation preserves positive maximal arcs and their
order, so $\tau_{i+B}=\tau_i+T$. This uses the physical count $B$;
a prescribed cut may instead give $N=B+1$ pieces by splitting one arc.
Choose $a$ inside arc $t-B$ and $b$ inside arc $s$. Since
$t-B\le s-1$, the actual distinct-mode switch $\tau_s$ lies in
$(a,b)$, while $a+T$ lies inside arc $t$. Thus $a<b<a+T$, and
$[b,a+T]$ retains the whole binary core $[s+1,t-1]$, both endpoint
ties and positive cap portions. This includes $t-s+1=B$, when the
cuts are in adjacent arcs $s-1,s$. Genuine reanchoring and
Proposition~\ref{pc:two-cuts} now supply the full lift and its proved
twist. Any other normalized lift agrees on a rigid core and hence on
its connected domain by full-adjoint uniqueness. Whole reflected time
reversal transports this construction as well.

For the last assertion, a cut can meet at most one of the two disjoint
closed binary cores. The other core and positive portions of both caps
remain in that genuine period, including when the cut is a switch.
Apply the preceding uniqueness assertion.
\end{proof}
Partial caps suffice here only for rigidity. Each scalar source theorem
still requires its stated complete arcs and other hypotheses, and each
harmonic calibration remains on a specified genuine complete period
with its original endpoints, real winding and qualifying signs.

\subsubsection{Full rooted-tail response and endpoint maximization}
\label{pc:tail}

Use the normalized binary variables and equations of
Lemma~\ref{loc:binary-reduction}, writing its $k$ as $K$ and its
$u^2+v^2-1$ as $S$.
\label{pc:binary-variables}
\begin{equation}\label{pc:binary-equations}
 K_s=0,\quad L_s=K,\quad m_s=-2L,\quad\chi_s=2K-3-u^2.
\end{equation}
Its cap function $f_j$ satisfies $m+1/2-f_j\ge0$, with equality
at an actual third-mode cap, by \eqref{loc:inactive-identity}.

A root of a binary switching equality has $\chi=0$.  On an active mode2
arc binary maximization gives $\chi\leq0$; the third-mode inequality must
also be retained.

Consider one complete positive mode2 tail rooted at an actual $0\to2$
switch.  Write $u(0)=U$, $r=1-U$, $b=v(0)/r>0$, its duration $d>0$, and
$E=e^{bd}>1$.  Thus
\begin{equation}\label{pc:tail-state}
 u(t)=1-r e^{-bt},\qquad v(t)=br e^{-bt}.
\end{equation}
Define its full-duration rooted response by
\begin{equation}\label{pc:Psi}
 \Psi(U,V)=2(1-U)\log{1-U\over1-V}
                 -(V-U)(2+V-U),\qquad V=u(d),
\end{equation}
and set
\begin{equation}\label{pc:R-J}
 R(s)={2(s-1+e^{-s})\over(1-e^{-s})^2},\qquad
 J_r(e)=2-{r\log e\over e-1}+{r^2\over2e},\quad
 J_r(1)=2-r+{r^2\over2}.
\end{equation}
A direct substitution of $V=1-r/E$ shows
$\Psi\leq0$ if and only if $R(\log E)\leq r$.

\begin{lemma}[Endpoint characterization on a failed tail]\label{pc:endpoint}
If $\Psi\leq0$ on the actual whole tail, then $U<0$, and its entire binary
maximizing inequality $\chi\leq0$ is equivalent to
\begin{equation}\label{pc:K-endpoint}
 K\leq J_r(E)<J_r(1)={3+U^2\over2}.
\end{equation}
Under this condition $\chi<0$ at every interior point.  The endpoint
inequality remains weak and allows an endpoint triple Hamiltonian tie.
\end{lemma}
\begin{proof}
For $s>0$,
\[
 \begin{gathered}
 R'(s)={4e^{-s}(\sinh s-s)\over(1-e^{-s})^3}>0,\qquad
 \lim_{s\downarrow0}R(s)=1,\\
 e^s-R(s)={2(\sinh s-s)\over(1-e^{-s})^2}>0,
 \end{gathered}
\]
since $(\sinh s-s)'=\cosh s-1>0$. Thus $1<R(s)<e^s$ for $s>0$,
and failure gives $r>1$, hence $U<0$.  Integrating \eqref{pc:binary-equations} from the actual root,
with $e=e^{bt}$, gives
\begin{equation}\label{pc:chi-exact}
 \chi(e)={ (2K-4)(e-1)/r+2\log e-r(1-e^{-1})\over b^2}
       ={2(e-1)\over b^2r}\bigl(K-J_r(e)\bigr).
\end{equation}
No initial derivative has been set to zero.  Differentiation gives
\[
 J_r'(e)={r\over2e^2}\bigl(R(\log e)-r\bigr).
\]
It is strictly negative for $1<e<E$, even if it vanishes at $E$.
Hence $J_r(E)<J_r(e)<J_r(1)$ in the interior.  The endpoint of
\eqref{pc:chi-exact} proves necessity of \eqref{pc:K-endpoint}; the same
formula proves its sufficiency and strict interior negativity.
\end{proof}

\subsubsection{A local paired-response theorem}
\label{pc:local}

The lawful reversed-time transformation used below is
\begin{equation}\label{pc:reverse}
 \psi(w)={3-\overline w\over1+\overline w},\quad j\mapsto-j,\quad
 Q_*={1\over2}\begin{pmatrix}-1&3\\1&1\end{pmatrix},\quad
 \Lambda\mapsto-Q_*\Lambda Q_*,\quad
 \symbf p\mapsto-D\psi(W)^{-t}\symbf p.
\end{equation}
The objective multiplier is unchanged.  Here $Q_*^2=I$,
$X(\psi W)=-Q_*X(W)Q_*$, $Q_*A_*Q_*=-A_*$, and
$D\psi F_j=-F_{-j}\circ\psi$.  These identities prove equality of the
full Hamiltonians under reversal; their derivatives give both transformed
state-adjoint equations.  The map is reflection $(u,v)\mapsto(-u,v)$
composed with $\Phi^2$, and preserves the physical cost.  Both $\Phi$ and $\psi$ preserve the six original compact cuts: their full actions cyclically permute the three rate cuts and three cusp-height cuts, while reflection interchanges the two side labels.  Thus these lawful whole transports do not impose the compact domain on an arbitrary fractional-linear chart.

\begin{theorem}[Paired rooted response]\label{pc:paired-response}
Suppose that one normal full maximizing lift of an actual minimizing
trajectory contains five complete positive consecutive arcs with whole
canonical word $21021$.  Jointly normalize the full lift to $\eta=1$.
Let the middle mode0 arc be the connector.  Its right neighbor is a rooted
mode2 tail; apply \eqref{pc:reverse} to obtain the opposite rooted mode2
tail from its left neighbor.  The exterior two arcs are their actual
following mode1 caps.  With the response \eqref{pc:Psi} evaluated over
each full tail, one has
\begin{equation}\label{pc:response-disjunction}
 \Psi_+>0\quad\hbox{or}\quad\Psi_->0.
\end{equation}
There is no determinant-sign hypothesis, no assumption that a local cap
matrix is invertible, and no extra periodic transversality in this theorem.
\end{theorem}
\begin{proof}
Assume both responses fail. Write the connector endpoints as
$(p-1,ap)$ and $(q-1,aq)$, where $a=v/(1+u)>0$, $A=1+a^2$ and $p<q$.
Put
\begin{equation}\label{pc:symmetric-variables}
 h={4\over Ap},\qquad z={q\over p}>1.
\end{equation}
The reversed connector has endpoints $h/z,h$ and the same rate $a$.
The two rooted tails start at $U_+=q-1$ and $U_-=h-1$.
Lemma~\ref{pc:endpoint} therefore gives
\begin{equation}\label{pc:symmetric-domain}
 0<q,h<1,\quad z>1,\quad p=q/z,\quad
 A={4z\over qh},\quad a^2={4z\over qh}-1>3.
\end{equation}
Their full multipliers $E,F>1$ and rates $aq/(2-q),ah/(2-h)$ satisfy
\begin{equation}\label{pc:two-failures}
 R(\log E)\le2-q,\qquad R(\log F)\le2-h.
\end{equation}

\paragraph{One cap row and its whole reverse.}
Let $K_+,K_-$ be the binary constants of this same transported matrix;
they need not agree. At the left connector endpoint set
\begin{align}
 S_c&={1/p-1/q\over a^2},&
 \chi_A&={q-p-2\log z+(4-2K_+)(1/p-1/q)\over a^2},
                                                        \label{pc:connector-chi}\\
 m_A&={K_+-(3+(p-1)^2)/2\over a^2p^2}
                +{A\chi_A\over2a^2p}-{1\over2},&
 L_A&={(A^2p^2/4)K_--K_++a^4p^2m_A\over2a^2p}.
                                                        \label{pc:connector-matrix}
\end{align}
The right actual root gives $\chi_B=0$, and integration gives the first
row. The actual left $1\to0$ tie gives $m_A$; the full reversed-matrix
identity $K_-=4(K_++2a^2pL_A-a^4p^2m_A)/(A^2p^2)$ gives $L_A$.
These reconstruct the full matrix without division by its determinant.
At the forward root,
\begin{equation}\label{pc:root-matrices}
 L_B=L_A+K_+S_c,\qquad m_B=m_A-2L_AS_c-K_+S_c^2.
\end{equation}
For a rooted tail with data $(K,L,m;U,v,e)$, put $r_*=1-U$, $b_*=v/r_*$ and
\begin{equation}\label{pc:tail-operation}
 S_*={e-1\over b_*^2r_*},\quad V_*=1-r_*/e,\quad v_*=v/e,
 \qquad \mathscr A=m-2LS_*-KS_*^2+\tfrac12-f_*.
\end{equation}
Here $f_*$ is the mode-$2$ cap function of Lemma~\ref{loc:binary-reduction}
at $(V_*,v_*,\chi(e))$, with the complete response~\eqref{pc:chi-exact}.
The actual forward terminal tie and its whole reversed counterpart are
\begin{equation}\label{pc:actual-cap-pair}
 \begin{gathered}
 \mathscr A(K_+,L_B,m_B;q-1,aq,E)=0,\\
 \text{and its whole swap }(q,E,K_+)\leftrightarrow(h,F,K_-).
 \end{gathered}
\end{equation}
Both occur at actual endpoints of complete tails; all remaining binary
and third-mode inequalities are still hypotheses of the full lift.

Let $\mathscr F_+$ be the first cap row multiplied by
$D_0^2q/(4h)>0$, where $D_0=4z-qh>0$. Direct substitution gives
\begin{align}\label{pc:negative-matrix}
 \partial_{K_+}\mathscr F_+&=-{(E-1)(2z-h)\over q}<0,\\
 \partial_{K_-}\mathscr F_+&=-{z[(2-q)(E-1)+q(z-1)]\over h}<0.\nonumber
\end{align}
No cap-system inverse, rank assumption or determinant sign is used.
Write $J(x,e)=J_{2-x}(e)$ and $J_0(x)=J(x,1)=2-x+x^2/2$.
The endpoint lemma gives $K_+\le J(q,E)<J_0(q)$ and
$K_-\le J(h,F)<J_0(h)$. Evaluate the first row at
$(J(q,E),J_0(h))$ and the reversed row at $(J(h,F),J_0(q))$,
calling the results $\mathcal F_+,\mathcal F_-$. Negativity of both
coefficients gives
\begin{equation}\label{pc:scaled-signs}
 \mathcal F_+(z)\le0,\qquad\mathcal F_-(z)\le0.
\end{equation}
Thus each comparison needs the opposite root bound, rather than its
full-duration response formula inside the cap row.

\paragraph{A common endpoint comparison.}
In the forward row abbreviate $R=R(\log E)$, $r=2-q$, and set
\begin{align}
 D&={r(E-1)^2\over2E^2q}\bigl[(qE+2)R-2Eq\bigr],\nonumber\\
 V&={q-h\over q}D+
 {hr^2(E-1)^2\over2q^2E}
          \left[{R\over E}-q(E+1-R)\right],\nonumber\\
 S_1&=D+r(E-1)(1-h),\qquad
 C_2=q(1-h/2),\qquad C_1=rE-2.\label{pc:residual-coefficients}
\end{align}
The complete cap substitution gives the exact identity
\begin{equation}\label{pc:residual-identity}
 \mathcal F_+(z)=V+S_1(z-1)+C_2(z-1)^2+C_1(z\log z-z+1).
\end{equation}
It can be checked by collecting $z^2,z\log z,z,1$: the coefficients are
$C_2,C_1,S_1-2C_2-C_1,V-S_1+C_2+C_1$, respectively.
The logarithm bound $L_1<\log E$ gives $R>2E/(E+1)$ and hence
$D>r^2(E-1)^2/(Eq(E+1))>0$. Since $C_1>-q$,
\[
 \mathcal F_+''(z)=2C_2+C_1/z>q(1-h)>0\qquad(z\ge1).
\]
Indeed $C_1/z>-q$ whether $C_1$ is negative or nonnegative.
Thus $\mathcal F_+'(1)=S_1>0$, and the same holds after the whole swap.
Since the actual $z>1$, \eqref{pc:scaled-signs} forces
\begin{equation}\label{pc:boundary-negative}
 \mathcal F_+(1)<0,\qquad\mathcal F_-(1)<0.
\end{equation}
Only scalar coefficients are held fixed here: $z=1$ is no physical
zero-duration trajectory, maximizing comparator or endpoint variation.

Retain the sharper scalar bound $R>R_0=(11E+1)/(3(E+3))$.
Indeed, for
\[
 g(E)={(E-1)(17E^2+8E-1)\over6E^2(E+3)},\qquad
 (\log E-g(E))'={(E-1)^4\over E^3(E+3)^2}>0,
\]
the difference vanishes at $E=1$, and substitution in the formula for
$R$ gives $R>R_0$. Also $R<E$ by Lemma~\ref{pc:endpoint}.
The quadratic $\Delta_E(s)=s/E-(2-s)(E+1-s)$ is concave. Its endpoint
values on $[R_0,E]$ are
\[
 \Delta_E(E)=E-1>0,\qquad
 \Delta_E(R_0)={(E-1)^2[15(E-1)^2+14(E-1)+8]\over9E(E+3)^2}>0.
\]
Consequently the common scalar comparison is
\begin{equation}\label{pc:Delta}
 \Delta_E(R)={R\over E}-(2-R)(E+1-R)>0.
\end{equation}
If $h\le q$, the forward failure gives $q\le2-R$ and $E+1-R>0$, so
\begin{equation}\label{pc:V0}
 V={q-h\over q}D+
 {hr^2(E-1)^2\over2q^2E}\left[{R\over E}-q(E+1-R)\right]>0.
\end{equation}
This contradicts $\mathcal F_+(1)=V<0$. If $q\le h$, the whole swapped
argument contradicts $\mathcal F_-(1)<0$ instead. Equality is included;
the two weak failures cannot coexist.
\end{proof}

\subsubsection{From a successful response to a complete central cost}
\label{pc:ingress-application}

The five-arc hypothesis alone does not supply a numerical central charge.
For that purpose the two candidate extensions must be complete actual
six-arc sources.  The forward ordered arcs from the connector and the
backward ordered arcs in \eqref{pc:reverse} must each be $021020$.
Their respective third arcs are the centers and their last three arcs are
the full terminal $020$ configurations.  Each source must be proper in
the physical period, with every arc and cap positive and complete.

\begin{lemma}[A successful rooted response supplies strict ingress]
\label{pc:strict-ingress}
On an actual complete rooted mode2 tail, $\Psi>0$ gives the strict incoming
conditions for the actual-cost $74/125$ central estimate: the incoming
central state has $u_L>0$, its raw switching response $q_{\rm raw,L}<0$,
and
\begin{equation}\label{pc:strict-row}
 q_{\rm raw,L}'+\kappa q_{\rm raw,L}
      <{b\over1-e^{-bd}}q_{\rm raw,L},\qquad \kappa={2v_L\over S_L}.
\end{equation}
All quantities are from the same full lift at the actual tail endpoint.
\end{lemma}
\begin{proof}
The original raw switching response is $q_{\rm raw}=\sigma_0=(S/v)p_v$.  Set $y=-\chi\geq0$, $w=S/v^2>0$, and
$\mathcal D(t)=(1-e^{-bt})y'(t)-by(t)$.  The full binary equations give
\[
 y''-by'=2u,\qquad {w'\over w}=-{2v\over S},\qquad y(0)=0.
\]
Integrating with \eqref{pc:tail-state} gives
\[
 b\mathcal D(d)=2(bd-1+e^{-bd})-(1-U)(1-e^{-bd})^2.
\]
The homogeneous initial derivative cancels; it has not been suppressed.
Thus $\Psi>0$ is equivalent to $\mathcal D(d)>0$.
By Lemma~\ref{pc:endpoint}, $1-U<R(\log E)<E$, so
$u_L=1-(1-U)/E>0$.
If $y(d)=0$, the left-hand nonnegativity of $y$ would give $y'(d)\leq0$,
contradicting $\mathcal D(d)>0$.  Therefore $y(d)>0$.
Finally $q_{\rm raw}=-wy$, so $q_{\rm raw,L}<0$ and
$q_{\rm raw}'+(2v/S)q_{\rm raw}=-wy'$.  Multiplying
$\mathcal D(d)>0$ by the positive factor $w/(1-e^{-bd})$ proves
\eqref{pc:strict-row}.  By Lemma~\ref{fw:regularity}, $q_{\rm raw}=\sigma_0$ is $C^1$ across the actual $2\to1$ switch.  The tail-side value and derivative are therefore exactly the center-side entrance data.
\end{proof}

\begin{corollary}[One expensive center]\label{pc:max-center}
Under Theorem~\ref{pc:paired-response}, suppose both ordered extensions
are proper complete $021020$ sources with their original compact geometry.
If their two physical center costs are $c_-,c_+$, then
\begin{equation}\label{pc:H592-max}
 \max(c_-,c_+)>{299\over250}.
\end{equation}
\end{corollary}
\begin{proof}
The budget-free response disjunction~\eqref{pc:response-disjunction}
selects a side with $\Psi>0$. If the chosen common lift contains that
side's entire extension, use it. Otherwise, reanchor outside the proper
six-arc extension and choose a complete normal maximizing lift there.
The value of $\Psi$ depends only on the actual states and full tail duration,
so its positivity is unchanged. Lemma~\ref{pc:strict-ingress} then proves
the strict actual ingress anew on this selected lift; no covector from a
different interval is transferred.

The complete central arc has actual $2\to1$ and $1\to0$ contacts.
In its mode-$0$ chart the full maximizing signs give
$V(A)=0$, $V_\xi(A)\le0$, $U(AE)=0$, $U_\xi(AE)\le0$ and $q_R\ge0$.
Thus both actual contact jets are retained. Its entire positive terminal
$020$ is present on the same selected lift. The sufficient-interface
form of Theorem~\ref{lc:direct-center} applies and gives
$c_\pm>299/250$ on the selected side. The other side is not selected
independently, so the conclusion is a maximum, not two singleton bounds.
In particular the weaker payment $74/125$ is available whenever it is
convenient in a disjoint cost sum.
\end{proof}

\subsection{Short central consequences used by the coarse resource argument}
\label{sec:central-singletons}
The following statement records only the source interfaces needed by the
coarse packing argument. Its numerical estimate follows from the stronger
strict-ingress theorem; no separate lower-threshold computation is needed.
The strict-ingress theorem in particular gives the consequence
\begin{equation}\label{h:criterion}
 C_{\rm loc}\le39/100\quad\Longrightarrow\quad
 u_R<0,\quad q_R\ge0,\quad q_R'>kq_R.
\end{equation}
Indeed $39/100<74/125$, and all hypotheses of Theorem~\ref{hc:H592}
are retained. This implication is used only at that strict-ingress interface.
\begin{theorem}[Short central sources]\label{loc:H390}
For complete actual positive sources retained on one genuine full lift:
\begin{enumerate}
\item $2021020$ has $c_3>39/100$;
\item $202102101$ has $c_3>39/100$ and $c_5>39/100$;
\item $2021021020$ has $\max(c_3,c_6)>39/100$.
\end{enumerate}
All exterior arcs remain, and indices start at zero.
\end{theorem}
\begin{proof}
The first conclusion follows from Corollary~\ref{hc:full922}; the third
follows from Corollary~\ref{hc:ten-arc}. For the second, the actual
initial $202$ gives $U>0$ at its $0\to2$ switch by
Lemma~\ref{loc:return-root}. The full mode-$2$ tail therefore has
$\Psi>0$ by the contrapositive of Lemma~\ref{pc:endpoint}.
The complete suffix $[1,8]=02102101$ remains on this same normal full
lift, with its actual root, two central contacts, positive continuation
and terminal $101$. Lemma~\ref{nx:rooted-continuations} consequently gives
$c_3>13/20>39/100$; its local charged position $2$ is original position $3$.
Whole reflected-time transport~\eqref{pc:reverse} takes $202102101$
to itself and this charge to $8-3=5$, giving the other inequality.
\end{proof}

\space
\section{Closed-domain finite witnesses for the local comparisons}\label{lp:finite}
The preceding reductions leave specified compact sign problems. This section supplies their exact polynomial constructors, full closed-domain coverage rules and outward remainder bounds. These witnesses prove the residual signs only after the geometric and same-lift domain hypotheses have been established.


\subsection{Exact arithmetic for the three-arc certificate}\label{f245:certificate}
This section specifies the finite mathematical witness for
Lemma~\ref{f245:residual}. Its proof content consists of rational
inequalities on whole boxes, together with the analytic soundness
argument in Section~\ref{sec:f245}. Neither a hash nor a recorded success
flag is an inequality. The supplied witness is verified afresh, including
all coverage data and every assigned predicate. No physical integration,
optimization, inversion of a trajectory, or cyclic-word enumeration is
part of this verification.

\subsubsection{Bernstein reconstruction}\label{f245:bernstein}
For $p(a+(b-a)X)=\sum_{j=0}^np_jX^j$, use the exact tensor
conversion~\eqref{cv:bernstein-coefficients}, with coefficients $\beta_i$,
and Lemma~\ref{cv:bernstein}. The coordinatewise inverse check is
\[
 p_j=\binom nj\sum_{i=0}^j(-1)^{j-i}\binom ji\beta_i.
\]

The polynomial \eqref{f245:poly-thin} is completely specified by its
formula, rational constants, bidegree $(73,6)$ and affine box
$[1/10,1/5]\times[0,1]$. It needs no input search or subdivision.
Substitution into \eqref{cv:bernstein-coefficients} gives
$\beta_{ij}>1$ for all $74\times7=518$ coefficients.
The electronic witness includes these exact coefficients; its
reconstruction checks their formula and the inverse polynomial identity.
This is the only sizable polynomial expansion in the outer argument.
The other outer comparisons \eqref{f245:poly-medium},
\eqref{f245:poly-highzero}, and \eqref{f245:poly-highdiag} use the
printed support-line and derivative proofs, with no coefficient table.

\subsubsection{A printed proof of right-cost monotonicity}\label{f245:beta-monotone}
Put $q=2-b+bv$, $M=4-\xi b$, and $u=1-b$.
Differentiating \eqref{f245:beta-integrand} gives
\[
 \partial_bI=\frac{2N}{M^2q^3},\qquad N=A\xi^2+B\xi+4q^3,
\]
where the following factorizations give the exact coefficients:
\begin{align*}
 A&=(2-b)^3+b(2-b)(6-b)v+b^2(2-b)v^2+b^3v^3,\\
 B&=2(2-b)(1-v)(3b^2v+b^2-2bv-4b+4)-4A.
\end{align*}
The discriminant identity
\[
 16Aq^3-B^2=4(2-b)^2(1-v)P_m
\]
holds identically, where
\begin{align*}
 P_m={}&(1-v)(4v^3+3v^2+18v+3)\\
 &+4u\{4v^4+2v^2(1-v)+3v^2+6v+3\}\\
 &+2u^2(1-v)(12v^3+v^2+8v+9)\\
 &+4u^3(1-v)^2(4v^2+2v+3)
   +u^4(1-v)^3(4v+3).
\end{align*}
These polynomial identities are verified by expansion. On $0\le b,v\le1$
every displayed summand is nonnegative, so $A\ge1$ and
$P_m\ge3(1-v)+12(1-b)\ge3(1-v)^4$. Therefore
\[
 4AN=(2A\xi+B)^2+4(2-b)^2(1-v)P_m\ge0.
\]
For $0\le\xi\le2$, $M\ge2$ and $q\ge1$, proving the required
nondecreasing property. It is strict for $v<1$; at $v=1$,
$N=8(\xi-2)^2$. Integration over any positive interval in $v$
proves strict monotonicity of $\beta$ in $b$ when $\eta>0$.
At $\eta=0$ the cost is zero, continuously.

\subsubsection{The binary witness and closed coverage}
The data file \path{independent_SMALL_COMPLETE_RESULT.json} contains
an array of $44\,813$ pairs $(p,\ell)$. Here $p$ is a nonempty finite
string of zeros and ones, and $\ell$ is exactly one of the eight
predicate names below. Only the following data fields have mathematical
meaning: \texttt{bits}, \texttt{root\_bounds},
\texttt{certified\_leaves}, \texttt{open\_leaves},
\texttt{pending}, and \texttt{axis\_rule}.
The two unresolved arrays must be empty and the axis rule is the fixed
cyclic order $r,b,E,\eta$. The root coordinates, divided by $2^{72}$,
are exactly
\[
\begin{array}{c|r|r}
 &\text{lower integer}&\text{upper integer}\\\hline
 r&0&2125064917291340346164\\
 b&157412216095654840456&4722366482869645213696\\
 E&4722366482869645213696&8929565713062601858626\\
 \eta&0&4722366482869645213696
\end{array}
\]
Each endpoint is an outward rounding of \eqref{f245:rectangle}, by less
than $2^{-72}$. At depth $k$ split coordinate $k\bmod4$ at the integer
$m=\lfloor(L+H)/2\rfloor$. A zero selects $[L,m]$ and a one selects
$[m,H]$. Both children are closed. Repeating this instruction reconstructs
each box without a floating-point operation.

Apply the exact path reconstruction and full-tree test of
Section~\ref{cv:cover}. The maximum actual depth is $43$, not the stored
depth allowance; the complete tree has $89\,625$ nodes.

The assigned predicate counts are
\[
\begin{array}{l|r}
 \texttt{central\_endpoints}&24578\\
 \texttt{root\_endpoints}&1627\\
 \texttt{central\_mean\_budget}&5055\\
 \texttt{exterior\_mean\_budget}&5073\\
 \texttt{no\_left\_budget}&6998\\
 \texttt{cusp\_B}&1320\\
 \texttt{high\_R}&83\\
 \texttt{b\_above\_a}&79
\end{array}
\]
The counts alone prove nothing; every assigned test must pass again.

\subsubsection{Directed arithmetic and deterministic verification}
All arithmetic intervals use the grid $2^{-160}$; multiply serialized
$2^{-72}$ endpoints by $2^{88}$ to embed them exactly. Apply the
integer outward operations of Lemma~\ref{cv:outward}, including the
zero-crossing square and enclosing integer square roots. A denominator
interval containing zero or an undefined operation makes verification fail.

Use $N=26$ in \eqref{f245:mean-series} and $M=32$ in
\eqref{f245:log-series}, with the displayed nonnegative geometric tails.
All other operations follow the expressions in Section~\ref{sec:f245}
in their printed order. The positive denominators on the enlarged root
satisfy $a,b>0$, $E<2$, $aE<2$, $D>2$, $q_B\ge1$ and
$1+\epsilon\ge1$. The affine coefficient $A_r$ is checked separately.
Use the proven monotone corners for $A_c,\beta,\nu,V_A^2,aD$ and the
right-mean remainder; never impose actual realization conditions on a
comparison endpoint. Set
\[
 B=\min([U-A_c-\beta]_{\max},U_+),\qquad
 U_+=2^{-160}\lceil2^{160}U\rceil.
\]
In the determinant predicates compare the lower integer enclosure with
$\lceil2^{160}H_*\rceil$, which is stronger than comparison with $H_*$.

For clarity, a complete checker consists of the following finite steps.
\begin{enumerate}
\item Read the rational constants and the witness; verify its domain,
legal labels, empty unresolved lists and complete closed coverage.
\item Reconstruct one leaf box and compute the label's expressions by
these directed rules and monotone corners.
\item For a root label reject if the affine coefficient lower bound is
nonpositive. For an endpoint label reject if its lower mean is not less
than its strict upper mean. Reject any undefined interval operation.
\item Check exactly the assigned inequality from the list in the proof
of Lemma~\ref{f245:residual}; never relabel a failed leaf.
\item Accept the finite witness only after all $44\,813$ leaves pass.
\end{enumerate}
Lemma~\ref{cv:outward}, the predicate implications in
Section~\ref{f245:residual-section}, and Theorem~\ref{cv:finite-rejection}
establish the lemma independently of the implementation's success message.
The electronic witness includes a direct implementation of these
formulas, the exact witness, and separate exact symbolic and
Bernstein-reconstruction checks. Its inputs are read-only. The residual
checker requires only Python's standard integer library; the optional
symbolic identity and coefficient reconstruction checks use exact
rational polynomial arithmetic.

\subsubsection{What the computer-assisted step does and does not remove}
The negative-coordinate proof is analytic; its right-leaf quartic is
positive by strict concavity and two endpoint values. The outer reduction uses four
explicit finite polynomial identities. The remaining nonnegative
comparison uses the certified closed cover just specified. In particular,
this proof does not claim to have replaced that cover by a short
algorithm-free inequality. Simple monotonicity of the whole determinant
cannot do so: at $a=b=1$ its mean derivative is negative, while in other
parts of the rectangle it can be positive. Concavity justifies two
endpoints; it does not justify dropping one of them everywhere. A
stronger local monotonicity statement can compress a cover but still
requires proof of its remaining graph inequalities.

\subsection{An affine absorption principle}
The following scalar comparison is useful whenever a response has an
explicit affine upper bound; it requires no subdivision.

\begin{lemma}[Affine absorption]\label{Sneg:affine-absorption}
Let $0\le x\le1$, $z\ge b\ge0$, $n\ge0$, and
$A_0=z(z-b)-n\ge0$. If $0\le f\le(1-x)z+bx$, then
\[
 f^2-(1-x)^2n\le z^2-n-x(2-x)A_0\le z^2-n.
\]
\end{lemma}
\begin{proof}
Put $\lambda=z-b\ge0$. The upper comparison is nonnegative,
so it may be squared. Since
$\lambda^2-n=A_0-b\lambda\le A_0$, expansion gives
$(z-x\lambda)^2-(1-x)^2n
\le z^2-n-2xA_0+x^2A_0$.
\end{proof}
\space
\space
\section{Analytic physical-cost estimates}\label{fee:analytic}
The following estimates retain the complete source and charge conventions
of Section~\ref{lc:source-convention}. Every matrix comparison is made
on one full normalized lift; algebraic endpoints in a proof are used
only for comparison. The four-arc and central estimates below use the
analytic cap identities. The nonpositive-tail proof also uses two
explicit polynomial minorants whose coefficient bounds are displayed.

\subsection{Four consecutive own arcs}\label{fee:four-section}
In state coordinates $(u,v)$ put
\begin{equation}\label{fee:four:eq:coordinates}
 x=\frac{1+u}{v},\qquad y=\frac{u^2+v^2-1}{2v},\qquad z=\frac{1-u}{v}.
\end{equation}
On the open state domain, $x,y,z>0$, $xy+xz+yz=1$, and the physical density is $\ell=x+y+z\ge\sqrt3$. 
\begin{theorem}\label{fee:four}
Let four complete consecutive control arcs have labels $1201$. Suppose their complete supporting context gives, in the canonical reciprocal coordinates, the endpoint signs
\begin{equation}\label{fee:four:eq:end-signs}
 z(t_0)<y(t_0),\qquad x(t_4)<y(t_4).
\end{equation}
Then their total physical cost satisfies
\[
J>\frac{103}{250}=0.412.
\]
Consequently the complete sources $101202$ and $1012012$ charge positions $1,\ldots,4$ at this value, and $0120121$ charges positions $2,\ldots,5$ at this value, together with their whole lawful images. Every supporting arc required by those sources remains complete and positive.
\end{theorem}

We use the cap equations of the full maximum principle on the selected complete lift. At the central forward switch, the two actual positive caps give normality and $x,z>y$ by the forward cap equations. Set
\begin{equation}\label{fee:four:eq:rho}
 \rho=\frac{x-y}{x+y},\quad \sigma=\frac{z-y}{z+y},\quad w=\rho+\sigma,\quad n=3-w-\rho\sigma.
\end{equation}
Thus $0<\rho,\sigma<1$. Let $d_0,d_2>0$ be the two central physical durations and $T=d_0+d_2$.

The matrix invariant is $\mathcal D=\det\Lambda$, distinct from physical cost. The analytic harmonic comparison for a complete normal lift states that positive determinant and positive future trace imply
\begin{equation}\label{fee:four:eq:harmonic}
 \mathcal D<\mathcal H:=\frac{29025933167887355473205403}{39478417604337498218323684}<\frac{2941}{4000}.
\end{equation}
It follows from the original frame endpoints and real winding, and the smoothed-octagon competitor. We verify both signs on this same selected lift before applying Corollary~\ref{cal:one-point}.

\paragraph{The physical arms and cap support.}
\begin{lemma}[The physical-arm integral]\label{fee:four:lem:arms}
Under \eqref{fee:four:eq:end-signs},
\begin{equation}\label{fee:four:eq:arms}
 J>\Phi((1-\rho)(1-\sigma))\ge\Phi((1-w/2)^2)>g(w),
\qquad g(w)=\frac{w(4-w)(6+w)}{2(6-w)(2+w)}.
\end{equation}
The function $g$ is strictly increasing on $(0,2)$.
\end{lemma}
\begin{proof}
On the right arm use $q=(x-y)/(x+y)$; along its two actual modes, $q'=-2/(x+y)$. It travels from $\rho>0$ to a negative endpoint by \eqref{fee:four:eq:end-signs}. Along $0\le q\le\rho$, its sum $x+y$ is at least $2y(t_2)/(1-q)$. Since $\ell(x+y)=1+(x+y)^2(3+q^2)/4$, this portion costs at least $\rho/2+y(t_2)^2F(\rho)/2$, where
\[
F(s)=s+2\log(1-s)+\frac4{1-s}-4.
\]
The left arm gives the same expression with $\sigma$, on a disjoint physical portion; the omitted portions have positive cost. Since $y(t_2)^2=q_0/(4-q_0)$ for $q_0=(1-\rho)(1-\sigma)$, addition gives
\[
J>\Phi(q_0),\qquad \Phi(q)=\frac{4-4q+q\log q}{4-q}.
\]
Now $\Phi'(q)=(4\log q-8-q)/(4-q)^2<0$ on $(0,1)$ and $q_0\le(1-w/2)^2$. The strict elementary inequality $\log q>-(1-q)/\sqrt q$ gives the last inequality in \eqref{fee:four:eq:arms}. It follows by differentiating $\log s-(s-s^{-1})/2$ for $s=\sqrt q<1$. Finally
\[
g'(w)=\frac{p(w)}{2(6-w)^2(w+2)^2}>0,\quad p(w)=w^4-8w^3-20w^2-48w+288;
\]
on $0<w<2$, the polynomial $p$ decreases, since $p'(w)<-16w^2-40w-48<0$, and $p(w)>p(2)=64$.
\end{proof}

\begin{lemma}[Three cap losses]\label{fee:four:lem:losses}
Write $k=\ell T$, $R=x/y$ and $t=z/y$. Suppose
\[
w<w_0<1,\qquad n>n_0>0,\qquad k<B<n_0,\qquad x/y,z/y<3.
\]
Define
\begin{equation}\label{fee:four:eq:loss-defs}
 \begin{gathered}
 u_0=\frac{2w_0}{(1+w_0)n_0},\qquad l_0=\frac{w_0}{(1+w_0/2)n_0},\\
 H_0=\frac{2u_0B^2}{3},\qquad D_0=\frac{w_0^2B}{8n_0}+\frac{l_0^2B^2}{4},\\
 L_{12}=\frac{Bw_0^2}{4(n_0-B)},\quad L_- =\frac{u_0^2B^2}{4}+\frac{u_0B^2}{6},\quad L_+=\frac{(D_0+H_0)^2}{4H_0}.
 \end{gathered}
\end{equation}
The actual determinant satisfies
\begin{equation}\label{fee:four:eq:det-loss}
 \mathcal D>\frac34-\max\{L_{12},L_-,L_+\}.
\end{equation}
Different valid pairs $(w_0,n_0)$ may be used for the nonhorizontal and horizontal losses. Moreover the actual trace is $\ell-h$ with $0\le h/\ell\le T^2/8$.
\end{lemma}
\begin{proof}
We give the cap-support calculation explicitly. Solving the two linear cap equations leaves an affine segment for the two-dimensional response vector $V$. Its endpoints lie in the triangle with vertices
\[
 W_0=(-C_xd_0,0),\quad W_2=(C_zd_2,0),\quad
 W_{12}=-(y+y^{-1})(-\beta_x\nu_2+\beta_z\nu_1),
\]
where
\[
C_x=\frac{(x+z)(x-y)}{2x},\quad C_z=\frac{(x+z)(z-y)}{2z},\quad \beta_x=\rho d_0/2,\quad\beta_z=\sigma d_2/2.
\]
The unit cap normals satisfy $\nu_1\cdot\nu_2=(y^2-1)/(1+y^2)$. Their explicit coordinates are
\[
\nu_1=\frac{(x^2-1,2x)}{1+x^2},\qquad\nu_2=\frac{(z^2-1,-2z)}{1+z^2}.
\]
Their inner product is the value just stated. Write $\mathsf J(a,b)=(-b,a)$. The analytically integrated full-PMP cap rows are
\[
h+A_i(V\cdot\mathsf J\nu_i)=2A_i\gamma_i,\qquad A_i=\frac{(1-e^{-q_i})(1+i^2)}2>0,
\]
with $i=x,z$, $d_x=d_0$, $d_z=d_2$, $\nu_x=\nu_1$ and $\nu_z=\nu_2$. Their exact positive defects define $\gamma_i$ explicitly:
\[
\beta_i-\gamma_i=i\left\{\frac{\sinh q_i-q_i}{e^{q_i}-1}+\alpha_i\frac{q_i-1+e^{-q_i}}2\right\}>0,
\]
with $q_i=d_i/i$, $\alpha_x=\rho$, $\alpha_z=\sigma$. Their three inward half-planes give precisely these vertices. The cap rows determine an affine line $V_y=rh-Q$ with $r>2$, $Q\ge0$, and $0\le h\le h_*=Q/r$; at the final algebraic endpoint $V_y=0$. Thus $V(h_*)$ lies on $[W_0,W_2]$. These are the exact cap equations and inequalities used from the normal maximum principle, rather than a finite numerical assertion.

Put $N=((x-z)/(x+z),\,y-(x^2+z^2)/(x+z))$, so $|N|^2=\ell^2-3$. The exact polarization on this actual normalized family is
\[
\mathcal D=\frac34-\frac{N\cdot V}{2}-\frac{(\ell-N_y)h}{2}-\frac{|V|^2}{4}+\frac{hV_y}{2}.
\]
If $V_x=\alpha h+\beta$, its quadratic coefficient is $-(\alpha^2+r^2-2r)/4<0$. Concavity of the determinant on the affine cap segment, followed by convexity of $\mathscr L(V)=N\cdot V/2+|V|^2/4$, gives
\[
 \mathcal D\ge\frac34-\max\{\mathscr L(W_{12}),\mathscr L(W_0)+H,\mathscr L(W_2)+H\},\quad H=(\ell-N_y)h_*/2.
\]
For the height bound the same two solved cap rows, using $\log G>1-G^{-1}$, give
\[
 Q\le\frac{d_0z(x-y)+d_2x(z-y)}{y(x+z)},\qquad
 r>\frac{2xz}{y(x+z)^2}(d_0^{-1}+d_2^{-1}).
\]
Division (with $h_*=0$ immediate if $Q=0$) gives
\[
h_*\le(d_0C_x+d_2C_z)\frac{d_0d_2}{T}\le\ell T^2/8,
\]
proving the trace assertion. All vertices here are algebraic comparisons; only the actual selected lift will be harmonically compared.

For completeness set $k_i=\ell d_i$, $\delta=|\rho-\sigma|$, and $b=B/n_0<1$. Substitution gives
\[
 \mathscr L(W_{12})\le\frac{B(w^2-\delta^2)}{4n_0}+\frac{B^2(w+\delta)^2}{4n_0^2}.
\]
The identity
\[
(w^2-\delta^2)+b(w+\delta)^2=\frac{w^2}{1-b}-(1-b)\left(\delta-\frac{bw}{1-b}\right)^2
\]
proves the $L_{12}$ bound. At the horizontal vertices put $c_i=C_i/\ell$. The exact scalar identities are
\[
 c_0=\frac{2\rho(1-\rho\sigma)}{(1+\rho)n},\qquad
 c_2=\frac{2\sigma(1-\rho\sigma)}{(1+\sigma)n},
\]
\[
 N\cdot W_0/2=\frac{\rho(\sigma-\rho)k_0}{(1+\rho)n},\qquad
 N\cdot W_2/2=\frac{\sigma(\rho-\sigma)k_2}{(1+\sigma)n}.
\]
They give $c_i\le u_0$. A nonnegative linear term forces its corresponding $\rho$ or $\sigma$ to be at most $w/2$, giving $c_i\le l_0$. Also
\[
\frac{\ell-N_y}{2\ell}\le\frac23,
\]
because $R^2-Rt+t^2-2R-2t\le0$ on $[1,3]^2$, as follows from separate convexity and its four negative corner values. Thus
\[
H\le\frac23(c_0k_0+c_2k_2)\frac{k_0k_2}{k}.
\]
Dropping a nonpositive linear vertex term yields $L_-$. In the other case write $\lambda=k_i/k$ and use $\rho(\sigma-\rho)\le w^2/8$ and its mate; then
\[
\mathscr L(W_i)+H\le D_0\lambda+H_0\lambda(1-\lambda)\le\frac{(D_0+H_0)^2}{4H_0}.
\]
This proves $L_+$. The strict bounds on $w,n,k$ give \eqref{fee:four:eq:det-loss}; zero cap height is included without division by that height.
\end{proof}

\paragraph{The resulting numerical bound.}
\begin{proposition}\label{fee:four:prop:criterion}
Let $0<a<1$, $0<w_0\le1/2$, $g(w_0)>a$, and put $B=a+a^2/3$. Choose $B<n_0\le3-w_0-w_0^2/4$. If the three cap losses at valid upper parameters satisfy
\[
\max\{L_{12},L_-,L_+\}<\frac{59}{4000},
\]
then the cost in Theorem~\ref{fee:four} exceeds $a$.
\end{proposition}
\begin{proof}
Assume $J\le a$. Lemma~\ref{fee:four:lem:arms} gives $w<w_0\le1/2$. Hence $\rho,\sigma<1/2$, so $x/y,z/y<3$, and
\[
n=3-w-\rho\sigma\ge3-w-w^2/4>n_0.
\]
The central two-arc cost $I$ is strictly less than $J$. Its two density profiles, measured from their common middle, are
\[
 \ell_0(s)=(x+y)e^{s/x}+(x+z)e^{-s/x}-x,\qquad
 \ell_2(s)=(y+z)e^{s/z}+(x+z)e^{-s/z}-z.
\]
They are strictly convex and have initial slopes $(y-z)/x$ and $(y-x)/z$. Integrated tangent bounds therefore give
\[
I>\ell T-\tfrac12\{(z-y)d_0^2/x+(x-y)d_2^2/z\}>\ell T-T^2.
\]
Since $I\ge\sqrt3T$ and $I<J\le a$, we obtain $T<a/\sqrt3$ and
\[
k=\ell T<a+a^2/3=B.
\]
The cap losses force $\mathcal D>2941/4000>0$. Independently the actual trace is positive, since $h/\ell\le T^2/8<a^2/24<1$. The harmonic comparison \eqref{fee:four:eq:harmonic} on this same full lift now gives a contradiction.
\end{proof}

\begin{proof}[Proof of Theorem~\ref{fee:four}]
Take
\[
a=\frac{103}{250},\quad w_0=\frac{249}{500},\quad n_0=\frac{2439}{1000},\quad B=\frac{87859}{187500}.
\]
The two domain reserves are
\[
 g(w_0)-a=\frac{309721}{1145333000}>0,\qquad
 3-w_0-w_0^2/4-n_0=\frac{999}{1000000}>0.
\]
For the nonhorizontal vertex,
\[
L_{12}=\frac{5447345859}{369453500000},\qquad \frac{59}{4000}-L_{12}=\frac{1046633}{184726750000}>0.
\]
For the horizontal vertices use the weaker valid bounds $w<1/2$, $n>39/16$. Then $u_0=32/117$, $l_0=32/195$, and
\[
 \begin{split}
 \frac{59}{4000}-L_-&=\frac{122238280763}{192501562500000}>0,\\
 \frac{59}{4000}-L_+&=\frac{54090241339849439}{83417343750000000000}>0.
 \end{split}
\]
Proposition~\ref{fee:four:prop:criterion} proves $J>103/250$. For $101202$, a whole phase makes the four charged modes $1201$;
the two actual returns give the two endpoint signs. For $1012012$,
whole phase $+1$ makes its charges $1201$. Its initial return $101$
gives the first sign, and the complete capped forward context $2012$
on $[3,6]$ gives the second by Lemma~\ref{fin:forward-force}.
For $0120121$, phase $-1$ makes the charges $1201$; the capped
forward prefix $0120$ gives the first sign and the actual return
$121$ gives the second. In particular the uncharged arcs $5$ and
$1$ in the two seven-sources, respectively, remain complete.
\end{proof}

\subsection{The opposite centre}
\label{fee:opposite-section}

\begin{theorem}\label{fee:opposite126}
Retain the complete positive source $1\mid020\mid1\mid202\mid1$,
its two actual outward forward-force contexts and their negative
canonical entrances, and one complete normalized maximizing lift on
a genuine fundamental interval of an original minimizer. The central physical charge is $c=c_4>63/500$.
\end{theorem}

The full source supplies the signed cap data of
Section~\ref{pf:M084}. We retain both adverse endpoint-state losses
in the ordinary signed profile of Section~\ref{ctr:profiles}.

\begin{proof}
Assume $0<c\le c_0=63/500$ and set
\[
 s=c_0^2,\qquad z=9409/10000,\qquad \eta=397/400,
 \qquad b=105413/113296<15/16.
\]
Let $q$ be the squared central height at its $u=0$ crossing, as in
\eqref{ctr:centered-state}. The source geometry and compact cuts of
Section~\ref{pf:M084} give $27/10<q\le243/4$.
The ordinary profile, at this same actual $q$, has constructors
\[
 E=(q-1)(q+3),\quad H=4q/E,\quad r=2q(q+1)/E,
 \quad K=(691q-1600+3373/q)/4800,
\]
\[
 M=2qKH^2/\eta+r^2/z,\qquad
 N=b\{(q-1)KH^3/2+r^2/[z(q+3)]\}.
\]
It gives $\sum\mu_i\le Mc^2$ and $\sum e_i\le Nc^3$ without
calibration. Put
\[
 \begin{gathered}
 A=4q^2/E,\qquad B_-=2/(q+3),\qquad B_+=2/(3+zq),\\
 k=AB_+,\qquad p_0=1-ks/2,\qquad U=r^2s/(zq).
 \end{gathered}
\]
The signed-transport and endpoint-loss identities
\eqref{ctr:drop-profile}--\eqref{ctr:ordinary-profile} permit
enlargement to the budget $c_0$:
\begin{align*}
 D(q)&=\frac{c_0\{A(1+U)-3B_-p_0^2\}}{2(1-2ks)},\\
 Z(q)&=\frac{Ms(3-zq)_+}{4zq}
       +\frac{M^2s^2}{4zq}+\frac{(3+Ms)Nc_0^3}{2}.
\end{align*}
Here $D$ bounds the positive signed length sum and $Z$ bounds the
aggregate adverse endpoint loss. Both are algebraic bounds on the
same actual lift; they select no comparison trajectory.
If $D\le\overline D$, $Z\le\overline Z$ and
$\overline Z\le\overline D^2$, the compensated endpoint
comparison of Lemma~\ref{ctr:compensated-loss} gives
\begin{equation}\label{fee:opposite:minimax}
 \Gamma>\frac34-\frac14
       \left(\overline D+\frac{\overline Z}{\overline D}\right)^2.
\end{equation}
We will prove that the parenthesis is less than $243/1000$ everywhere.

\paragraph{Uniform profile bounds.}
The original elementary guards give $qH<4$, $r<21/10$ and
$K<691q/4800$. Consequently
\[
 M<\frac{691}{150\eta}+\frac{441}{100z}<\frac{19}{2}.
\]
For $q\ge16/5$ one has $r<2$ and $zq>3$, so
\[
 M<\frac{691}{150\eta}+\frac4z<9,
 \qquad qN<\frac{15}{16}\left(\frac{691}{150}+\frac4z\right)<\frac{25}{3}.
\]
For the second estimate use
$q(q-1)KH^3/2<(691/9600)(qH)^3<691/150$.
It follows that on this entire range
\begin{equation}\label{fee:opposite:Zfar}
 Z(q)<\frac1q\left\{\frac{81c_0^4}{4z}
                +\frac{25}{6}(3+9s)c_0^3\right\}<\frac{4}{125q}.
\end{equation}
The last positive reserve is
$8934453850739/23522500000000000$.

On $27/10\le q\le16/5$, the elementary functions instead satisfy
\[
 H<9/8,\quad K<7/20,\quad A<31/10,
 \quad q/(q+3)\le16/31.
\]
Indeed $H$ decreases, $K$ increases here, and $A<31/10$ is the
concave quadratic comparison $-9q^2+62q-93>0$.
Hence, even discarding $b<1$,
\[
 qN<\frac{11}{5}\frac7{20}\frac{31}{10}\frac{81}{64}\frac12
        +\frac{16}{31}\frac{441}{100z}<4.
\]
Also $k<3/2$ and $r^2/z<5$. With $q_*=27/10$, the positive
transport numerator and its positive denominator give
\[
 D(q)<\frac{c_0\{(31/10)(1+5s/q_*)
                  -(30/31)(1-3s/4)^2\}}{2(1-3s)}<\frac3{20}.
\]
Using $M<19/2$ and $qN<4$ in the endpoint loss gives
\[
 Z(q)<\frac{(19/2)s(3-zq_*)+(19/2)^2s^2}{4zq_*}
       +\frac{(3+(19/2)s)4c_0^3}{2q_*}<\frac{69}{5000}.
\]
The two final positive rational reserves are respectively
\[
 \frac{4149287786661}{2952353200000000},\qquad
 \frac{2108048371161}{29403125000000000}.
\]
Thus $\overline D=3/20$, $\overline Z=69/5000$ work on the near
range, with $\overline Z<\overline D^2$ and
\begin{equation}\label{fee:opposite:near}
 \overline D+\overline Z/\overline D=121/500<243/1000.
\end{equation}

\paragraph{A monotone compensation for the transport budget.}
For $q\ge16/5$, retain the favorable term in the numerator of $D$.
Since $r<2$ and $p_0^2\ge1-ks$, one has
\[
 D(q)\le D_1(q):=
 \frac{c_0\{A(1+4s/(zq))-3B_-(1-ks)\}}{2(1-2ks)}.
\]
Define $F(q)=q(2c_0-D_1(q))$. Direct cancellation gives
\[
 F(q)=\frac{c_0q\{a q^3+d q^2+e q-81z\}}{z(q+3)B(q)},
\]
where
\[
 \begin{gathered}
 a=7z^2-40sz,\quad d=12z^2+21z-144sz-24s,
 \quad e=36z-27z^2-72s,\\
 B(q)=zq^3+(2z+3-16s)q^2+(6-3z)q-9.
 \end{gathered}
\]
The denominator is positive for $q\ge2$: its two middle coefficients
exceed $4$ and $3$, so $B(q)>4q^2-9>0$.
We prove that $F$ increases throughout $q\ge2$.
Differentiation gives
\[
 F'(q)=\frac{c_0\sum_{j=0}^6 a_jq^j}{z(q+3)^2B(q)^2},
\]
with the following short coefficient formulas:
\begin{align*}
 a_6&=z(640s^2-168sz-96s+23z^2),\\
 a_5&=6z(640s^2-152sz-176s+16z^2+23z),\\
 a_4&=-9(-768s^2z+40sz^2+376sz+16s+2z^3-64z^2-23z),\\
 a_3&=-108(-24sz^2+4sz+4s+2z^3+z^2-8z),\\
 a_2&=81(104sz+16s+3z^3-16z^2-2z),\\
 a_1&=486(8s+3z^2-4z),\qquad a_0=2187z.
\end{align*}
Use only $47/50<z<941/1000$ and $3/200<s<2/125$.
The first, second and last coefficients are positive. The other four
have the elementary bounds
\[
 a_4>600,\qquad a_3>500,\qquad a_2>-1000,\qquad a_1>-500.
\]
For example, lower bounds proving these four assertions are
\begin{align*}
 a_4&>576(47/50)^2+207(47/50)-3888(2/125)-18>600,\\
 a_3&>864(47/50)-216(941/1000)^3-108(941/1000)^2\\
 &\qquad-432(2/125)(1+941/1000)>500,\\
 a_2&>81\{104(3/200)(47/50)+16(3/200)+3(47/50)^3\\
 &\qquad\qquad-16(941/1000)^2-2(941/1000)\}>-1000,\\
 a_1&>486\{8(3/200)+3(47/50)^2-4(941/1000)\}>-500.
\end{align*}
Consequently for $q\ge2$ the derivative numerator is greater than
\[
 600q^4+500q^3-1000q^2-500q
   =200q^2(3q^2-5)+500q(q^2-1)>0.
\]
This proves monotonicity without a parameter subdivision or a stored
coefficient sign list.

Three rational endpoint comparisons now give
\[
 F(16/5)>3/10,\qquad F(16)>3/5,\qquad F(25)>7/10.
\]
Their exact positive reserves are
\[
 \frac{4683633759921}{70665290089250},\qquad
 \frac{1452850562103}{21117972419875},\qquad
 \frac{34884475}{3165789776}.
\]
Therefore the following transport and endpoint-loss budgets hold:
\[
 \overline D=2c_0-d/q,\qquad \overline Z=4/(125q),\qquad
 \begin{array}{c|c}
 q\text{ interval}&d\\\hline
 {}[16/5,16]&3/10\\
 {}[16,25]&3/5\\
 {}[25,243/4]&7/10.
 \end{array}
\]
In all three ranges $\overline D>3/20$ and
$\overline Z\le1/100<\overline D^2$.

\paragraph{The endpoint maximum and the actual invariant.}
Put $x=1/q$ and $u=2c_0$. On each of the three ranges,
\[
 \overline D+\frac{\overline Z}{\overline D}
     =u-dx+\frac{(4/125)x}{u-dx}.
\]
This function is convex in $x$, with second derivative
$2ud(4/125)/(u-dx)^3>0$. Its maximum on each full closed interval
is therefore at an endpoint. The six positive reserves below
$243/1000$ are, in chronological interval order,
\[
 \frac{54587}{2532000},\quad\frac{4387}{3732000},\quad
 \frac{16453}{858000},\quad\frac{107}{11400},\quad
 \frac{93}{7000},\quad\frac{1179317}{3549987000}.
\]
Together with \eqref{fee:opposite:near}, these exhaust the entire source
domain. Applying \eqref{fee:opposite:minimax} gives
\[
 \Gamma>\frac34-\frac14(243/1000)^2
       =\frac{2940951}{4000000}
       =\frac{1470471}{2000000}+\frac9{4000000}>0.
\]
Each genuine adjacent mean gives the actual canonical
$K_i=(3+\mu_i)/2>0$. On this same full lift the one-point criterion
therefore proves the positive future sheet before harmonic comparison.
The harmonic theorem gives
$\Gamma<H_*<1470471/2000000$, a contradiction.
The contrary budget included equality, proving the strict claim.
\end{proof}

\subsection{A nonpositive three-arc tail}\label{fee:ST305-section}
\begin{theorem}\label{fee:ST305}\label{loc:ST305}
On one full lift of an original global minimizer, let $12021$ be an
actual positive capped context contained in a genuinely eligible complete
interval. If its complete $202$ core enters at $u_0\le0$, then its
physical core cost is $T>61/200$.
Consequently every localized seven-arc tail listed in
Corollary~\ref{loc:ST305-transfer} has this strict cost bound.
\end{theorem}
\begin{proof}
Use the common actual capped-tail interface, without an entrance circle
or final return. LOW has the stronger bound in Corollary~\ref{sj:low}.
Suppose $2T\le V:=61/100$, and write $b=\beta$, $w=a+b$,
$\Delta=59/1000$, $h=47/200$. Thus $a^2<\mu<b^2$, $a<w/2$,
$0<w<V$ and $Z<(3+\mu)(V-w)/6$. The actual middle root excess and
Lemma~\ref{fin:signed-mean} give $r>-b/2$ and
$\mu<M:=(r^2+rb+b^2)/3$. We prove $\mathcal L<\Delta$.
For $w\le7/20$, the compact-cut argument of Corollary~\ref{sj:low} gives
\begin{equation}\label{fee:ST305:density}
 q(u)<(c-u^2)^{-1},\qquad c=17/15.
\end{equation}
All cap, reserve and shear comparisons below are first applied to this
actual source; scalar endpoints assert no new trajectories.

\paragraph{All widths $w\ge h$.}
The qualified coefficient ceiling $C=(1+a)^2/(c-a^2)$ gives
$A_C=(1+C)/[2(1+a)]\le1$ for $w\le7/20$: its cleared difference is
$(c-1)(1+2a)-2a^2(1+a)>0$ for $a\le7/40$.
For $7/20\le w\le2/5$, the open state domain gives $A_C\le25/24$.
Write $A=1$ or $25/24$ in these respective ranges. Lemmas
\ref{cmp:coefficient-responses} and \ref{cmp:strip-envelope}, dropping
only the positive terminal reserve, reduce the loss by convexity in
$\mu\in[a^2,b^2]$ to the responses $U_-$ and
$(\mu-a^2)(1-b)/[2w(1+a)]+A\mu$ at those two endpoints.
At $\mu=a^2$, the loss is at most $z_0^2<\Delta$, where
$z_0=(V-w)/2\le3/16$: indeed, with $d=A+(V-w)/6\le343/300$,
\[
 1-Vd-d^2/25\ge281563/1125000>0.
\]
At $\mu=b^2$, let
\[
 P_-={ (w-2a)(1+b)\over2(1-a)},\quad
 P_A={ (w-2a)(1-b)\over2(1+a)}+Ab^2,\quad
 S={ (3+b^2)(V-w)\over6}+P.
\]
Both choices satisfy $S_a\le-1$. For $P_A$ this follows from
$2(w-a)(1+a)^2-w(1+w/2)\ge w^2(1+w/2)/2>0$;
its left side decreases in $a\le w/2$ when $w\le2/5$.
Thus $(S^2-b^2)_a<0$ wherever $S^2-b^2\ge\Delta$, since then
$S>b$. A first-crossing argument reduces to $a=0$, where both losses
are bounded by
\[
 T_A(w)=\left[V/2+(V/6+A-1/2)w^2-w^3/6\right]^2-w^2.
\]
This decreases on $[0,2/5]$: the bracket is at most $6119/15000$,
and $2(V/6+A-1/2)\le193/150$, whose product is less than one.
The starting gaps are
\[
 \Delta-T_1(h)={526723439\over409600000000}>0,\qquad
 \Delta-T_{25/24}(7/20)={253602311\over6400000000}>0.
\]
Finally $H(w)=(3+w^2)(V-w)/6$ decreases on $[0,V]$;
for $w\ge2/5$, $\mathcal L<4H(w)^2$ and
$\Delta-4H(2/5)^2=62941/6250000>0$.

\paragraph{The lower-root sector $w\le h$, $r\le b/3$.}
The actual signed mean gives $\mu<\lambda b^2$, $\lambda=13/27$,
so $a<7w/17<3/20$. The endpoint responses imply
$X\le q_1\mu$, $q_1=1/(2w)+10/17$, using $A_C\le1$ and
$1/[2(1-a)]\le10/17$. The ordinary terminal reserve gives
$Y>b^2-\mu>(14/13)\mu$: the needed denominator guard is
$1-h-2h(1+h)^2=192577/4000000>0$.
Also $Z<z:=9/40$. Indeed the density bound
$(3+\lambda w^2)w/[2(c-w^2)]$ increases and the budget bound
$(3+\lambda w^2)(V-w)/6$ decreases; at $w=33/200$ their gaps
below $z$ are $14061/53093200$ and $220003/144000000$.
Apply \eqref{cmp:unclipped-strip} with $U=q_1\mu$ and
$R=(14/13)\mu$. Convexity in $Z\in[0,z]$ and then in
$\mu\in[0,\lambda w^2]$ leaves $z^2<\Delta$ and
\[
 F(w)={135424w^4\over210681}+{4784w^3\over12393}
       -{131219w^2\over247860}+{13w\over120}+{81\over1600}.
\]
The $Z=0$ endpoint is nonpositive because
\[
 1-\lambda[1/2+(10/17+14/13)h]^2
 =\frac{17399411}{28177500}>0.
\]
Bounding the positive higher powers in $F$ at $h$ gives
$F\le81/1600+13w/120-d_0w^2$, $d_0=23595911/58522500$.
Completing the square leaves the positive gap
$\Delta-81/1600-(13/120)^2/(4d_0)=12954557/11797955500$.

\paragraph{A polynomial switching moment.}
It remains that $w\le h$ and $b/3\le r<b$.
Put $z=r+a$, $A=1+a$, $d=c-w^2$, and retain $I,D=D_0$.
For $-a\le u\le r$, define
\[
 p(u)=A^{-2}+2(u+a)A^{-3}
       +\frac{(3+a-2r)(u+a)^2}{A^3(1-r)^2},\qquad
 H=\int_{-a}^r(M-u^2)p(u)\,du.
\]
This is an explicit rational moment: the integrand is quartic. The exact error is
\[
 p(u)-(1-u)^{-2}
 =\frac{(u+a)^2(r-u)\{2(1-u)(1-r)+A(2-u-r)\}}
 {A^3(1-r)^2(1-u)^2}\ge0.
\]
Since $M\ge\max(a^2,r^2)$, its weight is nonnegative. Set
\[
 \begin{gathered}
 C_0=\frac{d+(1-r)^2}{2Ad},\quad
 P_0=\frac{M-a^2}{2A^2D},\quad E=P_0+C_0H/D,\\
 J=\frac{w-z^2/A}{d},\quad k=\min\{J/2,(V-w)/6\}.
 \end{gathered}
\]
The actual positive switching integral and the same coefficient ceiling
as in Lemma~\ref{cmp:tail-interface} give $X\le(\mu/M)E$:
pointwise $\mu-u^2\le(\mu/M)(M-u^2)$, and the entrance term obeys the
same inequality. Also $\sigma<J$, so $Z<(3+\mu)k$.
The signed-strip identity, dropping only $Y>0$, and mean-ray convexity
leave $9k^2$ and $[(3+M)k+E]^2-M$.
The split at $w=41/250$ bounds the first by $z_0^2$, where
$z_0=9/40$; the respective gaps for $(3+w^2)k<z_0$ are
$139299/207457000$ and $137/187500000$.
Use
$\mathcal B(\nu)=24289/100000+2\nu-26\nu^2/5$.
It is positive and $\mathcal B(\nu)^2<\Delta+\nu$ for $0\le\nu\le1/16$.
For $t=25\nu\in[0,1]$, the exact identity is
\begin{align*}
 3\cdot10^{10}[\Delta+\nu-\mathcal B(\nu)^2]
 ={}&133437(1-t)^4+4(8665437)t(1-t)^3\\
 &+6(5405885)t^2(1-t)^2+4(338781)t^3(1-t)+1371453t^4.
\end{align*}
On $[1/25,1/16]$, $B\ge3/10$ and $0<\mathcal B'\le8/5$;
the gap $\Delta+\nu-\mathcal B^2$ has second derivative at least $28/25$
and derivative $10757/3125000$ at $1/25$.
It therefore suffices to prove $G_i=\mathcal B(M)-Z_i-E>0$, with
$Z_d=(3+M)J/2$ and $Z_b=(3+M)(V-w)/6$, using the applicable width bound.

\paragraph{Entrance elimination.}
All entrance derivatives hold $w,z$ fixed; put $y=z/w$. Write
$n=1+2z(w-z)/w-r^2$, so
\[
 D=\frac{wn}{(1-r)^2A(1+b)},\qquad
 \frac{D_a}{D}=\frac{2r}{n}-\frac2{1-r}-\frac1A+\frac1{1+b}.
\]
For $u=zs-a$, $0\le s\le1$, put $g=M-u^2$ and $K=w+z-2zs$.
Then $g_a=-K$ and
\[
 wK-g=\tfrac23(w-z)^2+2z(w-z)(1-s)+z^2(1-s)^2+aK\ge0.
\]
Also $p_a\le-2p/A$: for its quadratic coefficient
$c_H=(3A-2z)/[A^3(A-z)^2]$,
$c_{H,a}+2c_H/A=-2(3A^2-3Az+z^2)/[A^4(A-z)^3]<0$.
Hence $H_a\le-(1/w+2/A)H$.
Since $n\ge1-w^2$, $r\ge-a\ge-w/2$, and $d\ge1$,
\[
 \frac{(H/D)_a}{H/D}\le\frac{6w^2-1}{w(1-w^2)},\qquad
 \frac{(C_0)_a}{C_0}=\frac{2(1-r)}{d+(1-r)^2}-\frac1A\le w/2.
\]
Thus $(C_0H/D)_a<0$ for $w\le h$; the sufficient numerator
$2-13w^2+w^4$ is positive even at $w=7/20$.
For $z=0$ the moment vanishes and no logarithmic division is used.

Only the elementary entrance response remains. Put $F_i=\mathcal B(M)-Z_i-P_0$,
$N=M-a^2$, $v=1-r$ and
\[
 R=\frac{v^2(1+b)}{An},\qquad P_0=\frac{NR}{2w},\qquad
 \ell=\frac{R_a}{R}=\frac2v-\frac1{1+b}-\frac1A-\frac{2r}{n}.
\]
Since $0\le r\le b\le w$, $n\ge1-r^2\ge v^2$ gives
$R\le1+w$ and $0\le\ell\le4w/(1-w^2)$. Direct differentiation yields
\[
 \begin{aligned}
 \frac{R_{aa}}R&=\ell^2-\frac2{v^2}-\frac1{(1+b)^2}
 +\frac1{A^2}+\frac2n+\frac{4r^2}{n^2}\\
 &\le C(w):=-1-\frac1{(1+w)^2}+\frac2{1-w^2}
             +\frac{20w^2}{(1-w^2)^2}\\
 &\le C(h)<\frac{17}{10}.
 \end{aligned}
\]
Here $N\ge0$, $N_a=-(w+z)$ and $N\le w^2$, so
$(P_0)_{aa}\le17h(1+h)/20$.
Also $\mathcal B(M)_{aa}\ge4-312h^2/5$, while
$(Z_d)_{aa}\le J\le h$ and $(Z_b)_{aa}=(V-w)/3\le h$.
Consequently
\[
 (F_i)_{aa}\ge4-312h^2/5-h-17h(1+h)/20
       =11563/160000>0.
\]
Reducing $a$ toward zero preserves $M\ge\max(a^2,r^2)$ and $r\ge b/3$.
It therefore suffices to bound both $(F_i)_a$ at $a=0$.

There put $y=z/w\in[1/3,1]$ and $m=(1+y+y^2)/3$. Exactly,
\[
 \begin{gathered}
 n=1+2wy(1-y)-w^2y^2\le1+w/2,\quad
 R=1-w+\frac{2w(1-y)^2}{n},\\
 R_a=w-\frac{4w^2y(1-y)^2}{n^2}\le w.
 \end{gathered}
\]
It follows that $(F_d)_a-w/40\ge L(w,y)$, where
\begin{align*}
 L(w,y)={}&\frac{(1+y)(1-5w)}2+\frac{52w^3(1+y)m}5
       -\frac{w^2m}2+\frac{w(1+y)(1-y)^2}{1+w/2}-\frac w{40}\\
 &+\frac{w^2}{2(c-w^2)}
       \bigl[(1+y)(1-wy^2)-(3+w^2m)y^2\bigr].
\end{align*}
The bracket is at most $1+y$ and decreases with $w$. Thus
\[
 L_w\le(1+y)\left[-\frac52+\frac{156h^2}5+\frac49
                  +\frac{ch}{(c-h^2)^2}\right]-\frac1{40}<0.
\]
At $w=h$, collecting powers of $y$ gives the coefficientwise lower bound
\[
 L(h,y)\ge P(y):=\frac{-y^4+496y^3-426y^2-384y+356}{2000}.
\]
The five coefficient differences are positive rational numbers.
On $[1/3,1]$, $P''\ge26/375>1/16$, while
$P(7/8)=104607/8192000>1/80$ and
$0<P'(7/8)=181/51200<1/250$.
Taylor's inequality and completing the square give
$P>1/80-8/250^2>1/100$. Hence $(F_d)_a>w/40+1/100$ at $a=0$.

For the other shear, still at $a=0$, put
\[
 Q_0=\frac{(Z_d)_a-(Z_b)_a}{w}
 =\frac{(1+y)(V-w)}6+
 \frac{w}{2(c-w^2)}\bigl[-(1+y)+w(1+y)y^2+(3+w^2m)y^2\bigr].
\]
Its derivative satisfies
$(Q_0)_y\ge(V-w)/6+w(3y-1/2)/(c-w^2)>0$.
At $y=1/3$ it is
\[
 \frac{2(V-w)}9-
 \frac{w(1-4w/27-13w^2/243)}{2(c-w^2)},
\]
which decreases with $w$: the numerator of the subtracted quotient's
 derivative is
$c(1-8w/27-13w^2/81)+w^2+13w^4/243>0$.
Its value at $h$ exceeds $-1/40$ by $14513729/4191685200$.
Thus $(F_b)_a=(F_d)_a+wQ_0>1/100$ at $a=0$.
The positive second derivatives extend both entrance inequalities
throughout the original mean domain. Since $(C_0H/D)_a\le0$, both
barrier gaps $G_i=F_i-C_0H/D$ increase and reduce to $a=0$.

\paragraph{The zero-entrance functions.}
The preceding entrance reduction leaves
$0\le w\le h=47/200$, $1/3\le y\le1$, and $a=0$. Put
\[
c=17/15,\quad r=wy,\quad m=(1+y+y^2)/3,\quad M=w^2m,
\]
\[
n=1+2wy-(2w+w^2)y^2,\quad d=c-w^2,\quad
 D=\frac{wn}{(1-r)^2(1+w)},\quad C=\frac{d+(1-r)^2}{2d}.
\]
Define the polynomial comparison density and its moment by
\[
p(u)=1+2u+\frac{(3-2r)u^2}{(1-r)^2},\qquad
 H=\int_0^r(M-u^2)p(u)\,du,
\]
and let
\[
E=\frac{M/2+CH}{D}.
\]
The formula at $w=0$ is interpreted by cancellation and continuity.
We prove $E_w<2/3$ and positivity of both boundary gaps.
All variables in this calculation are algebraic comparators; no new
physical lift is selected at a boundary value.

\paragraph{A derivative bound from moments.}
Put
\[
e_0=\frac{y(1+y)}3,\quad
 e_1=\frac{2y^2+2y^3-y^4}{12},\quad
 e_2=\frac{5y^3+5y^4-4y^5}{45}.
\]
These are the moments $e_j=\int_0^y(m-v^2)v^j\,dv$. In particular
\[
e_1\le\frac{ye_0}{2},\qquad e_2\le\frac{y^2e_0}{3};
\]
the respective differences are $y^4/12$ and $4y^5/45$.
Set
\[
B=(1-r)^2(e_0+2we_1)+(3-2r)w^2e_2,\quad
 A=\frac{1+w}{n},\quad F=CAB,\quad
 L=\frac{(1-r)^2(1+w)}n.
\]
Direct integration gives
\[
E=\frac{wmL}{2}+w^2F.
\]
We have $d>1$, $0<C<1$, and
\[
C_w=\frac{(1-r)(w-cy)}{(c-w^2)^2}<0,
\]
because $w\le h<c/3\le cy$. Also $r\le h<1/3$ and
\begin{align*}
 B_w&=2(1-r)\{-ye_0+(1-3r)e_1+3we_2\}\\
 &\le-(1-r)(1+r)ye_0<0.
\end{align*}
Thus $B\le e_0$.

We first prove $F<27/40$. For $y\le1/2$, $e_0\le1/4$ and
$n\ge1-w^2$, so
\[
F\le\frac1{4(1-w)}\le\frac1{4(1-h)}<\frac{27}{40}.
\]
For $y\ge1/2$, the two moment bounds imply
\[
B\le e_0(1-r+r^3/3).
\]
Moreover
\[
n-(1+w)(1-r)=w(1-y)(2y-1)+w^2y(1-y)\ge0.
\]
Hence
\[
F\le\frac23\left(1+\frac{h^3}{3(1-h)}\right)<\frac{27}{40}.
\]
The last two reserves are $2131/6120$ and $125677/27540000$.
Next
\[
A_w=\frac{(1-y)^2+y^2(1+w)^2}{n^2}>0.
\]
Since $n-(1-w)(1+wy)=w(1-y)(1+2y+wy)\ge0$ and $(1-y)^2+y^2(1+w)^2\le(1+wy)^2$, we have $A_w\le(1-w)^{-2}$. Hence $C_w,B_w<0$, $C<1$ and $B\le e_0$ give
\[
F_w\le e_0A_w\le\frac{2}{3(1-h)^2}<\frac65,
\]
with final reserve $21362/351135$.

The exact identity
\[
L=1-w+\frac{2w(1-y)^2}{n}
\]
and $n_w\ge-2w$, $n\ge1-w^2$ give
\[
L+wL_w\le1-2w+\frac{4w(1-y)^2}{(1-w^2)^2}
 \le1-2w+\frac92w(1-y)^2.
\]
Here $(1-h^2)^2-8/9=53437129/14400000000>0$.
Therefore
\[
E_w\le\frac m2+w\left\{\frac{27}{20}+\frac65w-m
                  +\frac94m(1-y)^2\right\}.
\]
The bracket is positive. With $C_{\mathrm{bd}}=27/20+6h/5=204/125$, the last
expression is at most
\[
U(y)=m\{1/2-h+(9h/4)(1-y)^2\}+hC_{\mathrm{bd}}.
\]
It is convex, since the exact identity is
\[
U''(y)=\frac{71}{1200}+\frac{47}{400}(3y-1)(6y-1)>0.
\]
Its endpoint reserves below $2/3$ are
\[
\frac23-U(1/3)=\frac{3578}{84375},\qquad
 \frac23-U(1)=\frac{1361}{75000}.
\]
Consequently $E_w<2/3$ throughout the required rectangle. 

\paragraph{The boundary signs.}
Write $V=61/100$ and
\[
J=\frac{w-r^2}{c-w^2},\quad
 Z_d=(3+M)J/2,\quad Z_b=(3+M)(V-w)/6,
\]
\[
\mathcal B(M)=\frac{24289}{100000}+2M-\frac{26M^2}{5},\quad
 G_i=\mathcal B(M)-Z_i-E.
\]
At the width barrier
\[
W(y)=\frac{337}{2000}-\frac{7y}{400}+\frac{153y^2}{4000},
\]
let
\[
N_i(y)=90\{17-15W(y)^2\}n(W(y),y)G_i(W(y),y),\qquad i=d,b.
\]
These definitions involve only rational arithmetic and the displayed
quartic moment integral. Each $N_i$ is a polynomial of degree 24.
For $0\le y\le1$, direct expansion gives the lower bounds
\[
1000N_d(y)>L_d(y),\qquad1000N_b(y)>L_b(y),
\]
where
\begin{align*}
L_d={}&-1720+31464y-43147y^2-10534y^3+19366y^4\\
 &+17232y^5-15284y^6+8456y^7-5392y^8,\\
L_b={}&10529-14880y-12474y^2+40542y^3-31138y^4\\
 &+16967y^5-10179y^6+2319y^7-1589y^8.
\end{align*}
For clarity, the entire arithmetic needed for this truncation is given
here. If $N_i=\sum c_{i,k}y^k$, the coefficient lower bounds
$1000c_{i,k}>l_{i,k}$ for $k=0,\ldots,8$ are the two integer rows
\[
\begin{array}{c|rrrrrrrrr}
i&0&1&2&3&4&5&6&7&8\\\hline
d&-1387&31464&-43147&-10534&19366&17232&-15284&8456&-5392\\
b&10576&-14880&-12474&40542&-31138&16967&-10179&2319&-1589
\end{array}
\]
Each unscaled coefficient reserve $c_{i,k}-l_{i,k}/1000$ exceeds
$1/100000$. Their minimum values, by row, are
\[
\frac{68220585961351}{3200000000000000000},\qquad
\frac{338561509273248403461}{6400000000000000000000000}.
\]
The two total negative tails are exactly
\begin{align*}
\sum_{k=9}^{24}(-c_{d,k})_+
 &=\frac{1090765939242063781497125327}{3276800000000000000000000000}
 <\frac{333}{1000},\\
\sum_{k=9}^{24}(-c_{b,k})_+
 &=\frac{150985480245549429849736467}{3276800000000000000000000000}
 <\frac{47}{1000}.
\end{align*}
For $0\le y\le1$, discard positive tail terms and bound $y^k\le1$
in negative terms. Subtracting these two tail bounds from the constant
terms of the integer rows gives exactly $L_d,L_b$. This supplies the
claimed minorants. 

We prove their signs by elementary calculus. Put
\[
g(y)=-15284y^2+8456y^3-5392y^4.
\]
It is strictly concave on $[0,1]$: $g''=-30568+50736y-64704y^2<0$
by its negative discriminant. For $1/3\le y\le3/4$, use the chord
through $0,3/4$ to get
\[
L_d(y)\ge P_-(y):=
 \frac{33003y^5+77464y^4-42136y^3-172588y^2+125856y-6880}{4}.
\]
Now
\[
P_-'''=3(165015y^2+154928y-21068)>0\quad(y\ge1/3),
\]
and $P_-'(3/4)=-5426745/1024<0$. A function whose derivative is convex
and negative at the right endpoint has no interior strict minimum: the
nonpositive sublevel set of its derivative is an interval containing
that endpoint. Thus its minimum is at an endpoint. Here
\[
P_-(1/3)=1249597/324>0,\qquad P_-(3/4)=5116913/4096>0.
\]
For $3/4\le y\le1$, the chord of $g$ through $3/4,1$ gives
\[
L_d(y)\ge P_+(y):=
 \frac{-18817y^5+116329y^4-42136y^3-172588y^2+125856y-6880}{4}.
\]
Its third derivative is the concave quadratic
\[
P_+'''=-282255y^2+697974y-63204,
\]
positive at the two endpoints with values $4824129/16$ and $352515$.
Hence $P_+''\ge P_+''(3/4)=1466629/64>20000$. Taylor's theorem and
completing the square, centered at $9/10$, give
\[
P_+(y)\ge P_+(9/10)-\frac{P_+'(9/10)^2}{40000}
 =\frac{320197948869839}{2560000000000}>0.
\]
Therefore $L_d>0$ on the whole required interval.

For $L_b$, put $D_b=-L_b'$. Explicitly,
\begin{align*}
D_b={}&14880+24948y-121626y^2+124552y^3-84835y^4\\
 &+61074y^5-16233y^6+12712y^7.
\end{align*}
Its third derivative has the positive decomposition
\[
D_b'''=24\{y^2(111230y^2-81165y+20000)
              +(132685y^2-84835y+31138)\}>0.
\]
Both quadratics have positive leading coefficient; the values of
$4ac-b^2$ are $2310642775$ and $9329204895$.
On $[1/3,2/3]$, convexity of $D_b'$ and its negative endpoint values
$-17274068/729$, $-11948228/729$ give $D_b'<0$, and thus
$D_b\ge D_b(2/3)=10856960/2187>0$.
On $[2/3,1]$, $D_b''\ge D_b''(2/3)>100000$. Taylor's theorem centered
at $3/4$ gives
\[
D_b(y)\ge D_b(3/4)-\frac{D_b'(3/4)^2}{200000}
 =\frac{54709251469071}{13107200000}>0.
\]
Consequently $L_b$ decreases and $L_b(y)\ge L_b(1)=97>0$.
Both boundary gaps $G_d(W(y),y)$ and $G_b(W(y),y)$ are therefore
strictly positive.

\paragraph{The width comparison and the actual loss.}
At $a=0$, the boundary satisfies
\[
 \frac{2003}{12000}\le W(y)\le\frac{757}{4000}<\frac15.
\]
For $w\le1/5$, write $D=w\delta$ and
$E=wm/(2\delta)+w^2Ce/\delta$, where
\[
\begin{split}
 \delta={}&\int_0^y\frac{dv}{(1-wv)^2}
 +\frac{(1+wy)^2}{(1-wy)^2}\int_y^1\frac{dv}{(1+wv)^2},\\
 e={}&\int_0^y(m-v^2)p(wv)\,dv.
\end{split}
\]
Integration gives $\delta=\delta_1+\delta_2$, where
\[
 \delta_1=\frac{y}{1-wy},\qquad
 \delta_2=\frac{(1-y)(1+wy)}{(1-wy)^2(1+w)}.
\]
Their sum is $n/[(1-wy)^2(1+w)]$, as required by the definition of $D$.
At fixed $y$ we have $w(\delta_1)_w/\delta_1=wy/(1-wy)$. For $y<1$,
\[
 \frac{w(\delta_2)_w}{\delta_2}
 =\frac{wy}{1+wy}+\frac{2wy}{1-wy}-\frac{w}{1+w}
 \le \frac{2wy}{1-wy}\le\frac12 .
\]
Thus $w\delta_w\le\delta/2$, also at $y=1$ where $\delta_2=0$.
Since $C=\tfrac12+(1-wy)^2/(2d)$ and $d=c-w^2$, we likewise have
$wC_w\ge-C/2$. On $0\le v\le y$, $m-v^2\ge0$ and $p(wv)$ is
nondecreasing in $w$, so $e_w\ge0$. Both terms of $E$ therefore
increase for $0<w\le1/5$ and extend continuously to $w=0$.
Direct differentiation gives
\[
 (Z_d)_w-4wm\ge
 \frac{3-(6+8c)w+3w^2-4w^3}{2c}\ge\frac{14}{425}>0.
\]
Since $\mathcal B'(M)M_w\le4wm$, this proves $(G_d)_w<0$.
On the other side of the boundary, the derivative estimate already
proved gives
\[
 (G_b)_w\ge2\frac{13}{27}\frac{2003}{12000}
       \left(2-\frac{52h^2}{5}\right)
       +\frac12-\frac{V^2}{18}-\frac23
       =\frac{112897679}{2700000000}>0.
\]
Hence $G_d(w,y)>0$ for $w\le W(y)$ and $G_b(w,y)>0$ for
$w\ge W(y)$. Together with entrance monotonicity these inequalities
exclude the remaining root sector. The wider widths and lower roots
were handled above. Thus the whole actual HIGH domain has
$\mathcal L<\Delta$, contradicting Lemma~\ref{cmp:loss-obstruction}.
The LOW branch was already excluded by its larger cost bound.
Therefore $T>61/200$.
\end{proof}
\subsection{A shared derivative and six disjoint charges}\label{fee:W71-section}
\begin{theorem}\label{fee:W71}\label{sj:W630}
In an actual positive full context $02120102$, indexed $0,\ldots,7$, the physical cost $W$ of exactly arcs $1,\ldots,6$ satisfies
\[
                         W>\frac{71}{100}>\frac{63}{100}.
\]
The entire source is retained in a genuine containing interval of an
original global minimizer, with one full lift shared by both tails.
\end{theorem}
\begin{proof}
The common switching derivative forces one actual tail to have a LOW
mean. The central forward pair supplies both negative entrances and
both circles, so the retained circle-LOW cost and ordinary nonpositive
cost add on disjoint triples. Use the selected full lift and its $C^1$
switching regularity (Lemma~\ref{fw:regularity}); its state covector is
not continued periodically.

\paragraph{The two tails on one full lift.}
The right capped context is actual $[3,4,5,6,7]=20102$, transformed by $j\mapsto j-1$ to $12021$. The reversed left context $[4,3,2,1,0]=02120$ becomes $12021$ by $j\mapsto1-j$ with the full reflected time reversal. Their charged cores are respectively $[4,5,6]$ and $[3,2,1]$. Thus their costs $T_R,T_L$ add to $W$, while exterior arcs $0,7$ remain uncharged.

Either actual capped triple proves normality of this selected complete
lift by the two-root, two-cap argument of Lemma~\ref{loc:signed-strip},
whose normality proof does not require the entrance sign. Normalize this
single lift jointly before the following transports.

At the actual central seam set $w=u+iv$, $\Lambda=\left(\begin{smallmatrix}l&b\\m&-l\end{smallmatrix}\right)$ and
\[
 x=\frac{1+u}{v},\quad y=\frac{u^2+v^2-1}{2v},\quad z=\frac{1-u}{v}.
\]
Then $x,y,z>0$, $xy+xz+yz=1$. Put $d=x+z$, $e=1-xz=dy>0$ and
\[
 R=\tfrac12\begin{pmatrix}1&-3\\1&1\end{pmatrix},\qquad S=\operatorname{diag}(1,-1).
\]
The actual right and reversed-left transports are
\begin{align*}
 w_R&=(w-3)/(w+1),&\Lambda_R&=R\Lambda R^{-1},\\
 w_L&=(-\overline w-3)/(1-\overline w),&\Lambda_L&=-RS\Lambda SR^{-1}.
\end{align*}
State covectors transform cotangentially, with an additional minus sign under time reversal. Both transformed matrices have determinant $\mathcal C=\det\Lambda$. The entrance coordinates are
\[
 u_R=\frac{y-x}{y+x},\quad v_R=\frac2{x+y},\qquad
 u_L=\frac{y-z}{y+z},\quad v_L=\frac2{y+z}.
\]
The complete central pair $[2,3,4,5]=1201$ gives $y<x,z$
by Lemma~\ref{fin:forward-force} on this same selected lift.
The full canonical reciprocal triples are $(y,z,x)$ and $(y,x,z)$.
Their entrances are therefore negative, and each has its second
reciprocal larger than its first. The identity
$2v_j(y_j-x_j)=v_j^2-(1-a_j)(3+a_j)$ gives both actual entrance
circles $v_j^2>(1-a_j)(3+a_j)$, with $a_j=-u_j>0$.
For each tail use the increasing coordinate $t=u_j$ along its $202$ core and write its two internal roots as $r_j<s_j$. Let $q_j=v_j^{-2}$ along that core, $ds=q_j(t)dt$, and define the \emph{actual} middle mean
\[
 \mu_j=\frac{\int_{r_j}^{s_j}t^2q_j(t)dt}{\int_{r_j}^{s_j}q_j(t)dt},\qquad
 K_j=(3+\mu_j)/2.
\]
It is positive. If the entrance is nonnegative, strict increase of $t$ gives $\mu_j>u_j^2$. These means, states, and matrices all belong to the same selected lift.

\paragraph{The common switching derivative.}
At the central $2\to0$ tie, $p_v=0$; write $p=(P,0)$. The inactive inequality $H_1\le0$ gives $P\ge0$. The common active Hamiltonian gives
\[
 \operatorname{tr}(\Lambda X)=\ell-Pv,\qquad
 \ell=\frac{3+u^2+v^2}{2v},\qquad h=Pv^2\ge0.
\]
Here $X=v^{-1}\left(\begin{smallmatrix}u&-u^2-v^2\\1&-u\end{smallmatrix}\right)$; in a binary chart $L=l-um$. The switching function is
\[
 \sigma_0=\frac{u^2+v^2-1}{v}p_v.
\]
It changes from nonpositive to nonnegative. By its $C^1$ regularity, $\sigma_0'\ge0$ at the switch. Direct differentiation of the full adjoint equation at the tie gives $p_v'=2m+1-2P$, so
\begin{equation}\label{fee:W71:sj:jet}
                       m+\tfrac12-P\ge0.
\end{equation}
This uses a weak derivative inequality and includes triple ties; no strict crossing or inactive gap is assumed.

For completeness, matrix reconstruction makes the consequence explicit. The two equations
\[
 K_j=2u_j(\Lambda_j)_{11}+(\Lambda_j)_{12}-u_j^2(\Lambda_j)_{21}
\]
and the active trace equation uniquely determine $(l,b,m)$.
The absolute value of their coefficient determinant is
\[
 \frac{8d^3(1-xz)}{(1+x^2)^2(1+z^2)^2}>0.
\]
They give
\begin{align*}
 m&=m_0-\frac{xzd^2h}{4e},\\
 m_0&=\frac{2xzd(x^2+xz+z^2+1)-K_Lx(1+z^2)^2-K_Rz(1+x^2)^2}{4de},\\
 P&=hd^2/4.
\end{align*}
Substituting these expressions into \eqref{fee:W71:sj:jet}, and using $y=e/d$, yields exactly
\begin{equation}\label{fee:W71:sj:meanjet}
 \begin{split}
 x(1+z^2)^2(\mu_L-u_L^2)&+z(1+x^2)^2(\mu_R-u_R^2)\\
                         &+2d^3h\le0.
 \end{split}
\end{equation}
The coefficients are strictly positive. Therefore at least one tail satisfies $\mu_j\le u_j^2$. Such a tail cannot have nonnegative entrance, since its actual middle mean would be strictly larger than the entrance square. Its entrance is consequently $-a<0$ and $\mu_j\le a^2$.

The selected LOW tail has its actual entrance circle. Hence
Lemma~\ref{cmp:circle-low-width} gives $T_j>81/200$.
The other actual capped tail has negative entrance, so
Theorem~\ref{fee:ST305} gives $T_{\rm other}>61/200$.
The two charged triples are disjoint and exhaust arcs $1,\ldots,6$.
Therefore
\[
 W=T_L+T_R>81/200+61/200=71/100>63/100.
\]
Both circles and the common mean-jet belong to the selected full
lift. The exterior arcs $0,7$ remain uncharged.

\end{proof}
\space
\label{revision:appendix-end-2}
\end{proofappendix}
\section{Further complete-source estimates}\label{app:sources}
\begin{proofappendix}
\section{Returned and coupled union geometry}\label{pf:F01}
The common-lift seven-charge inequality is Theorem~\ref{ua:seven}, on
its full thirteen-source $0102010201020$. 
The shared-core identities in Sections~\ref{pf:M002}--\ref{pf:M006}
also supply reusable inputs to returned-core and own-pair estimates.
\space
\subsection{Four coupled cores and common normalization}\label{pf:F02}\label{pf:M008}
Consider a complete chain with source $20212101020212101$, indexed
$-1,\ldots,15$, on one proper containing interval. Its four returned
eight-sources start at $-1,2,5,8$ and are whole phases of $01020212$.
Their four-core charges are $[1,4]$, $[4,7]$, $[7,10]$ and $[10,13]$.
Choose one genuine containing interval and
one full maximizing lift. Capped-core normality precedes normalization
of that lift. In particular, the disjoint first and third four-cores
have total cost below $C_*$, so at least one has cost $A<C_*/2<1$.
The structural apex inequality gives, on this same lift,
\[
 \mathcal C>\frac34-\frac9{64}A^2>\frac{63}{100}.
\]
The positive future trace follows before applying harmonic calibration.
Consequently $\mathcal C<H_*$, where
\[
 H_* = \frac{29025933167887355473205403}
              {39478417604337498218323684}<\frac{1471}{2000}.
\]
The inner once-counted physical charge is
\[
 C_{\mathrm{in}}=\sum_{i=2}^{12}c_i.
\]
We derive a local lower bound for this charge without using a
complementary numerical budget. The identities below apply to the stated complete common-lift source.

\subsection{The exact cap and transport identities}\label{pf:M002}
On an overlap, use the left core's canonical $0202$ chart. Write
$K=(3+h^2)/2$, $\widehat K=(3+\widehat h^2)/2$, let $r>h$ be the
left core's right root, and let $-\xi<-\widehat h$ be the right
core's left root. The shared density parameter is $x=1/\kappa^2$;
its integration constant is $\lambda$. If $\beta$ is the shared
arc's final coordinate, the exact state and determinant equations are
\begin{equation}\label{pf:shared-geometry}
 0<x<\frac13,\qquad
 1-\beta=\frac{4x}{(1+x)(1-\xi)},\qquad
 \mathcal C=\frac{K\widehat K(1+x)^2}{4}-\lambda^2.
\end{equation}
With $q(u)=x/(1-u)^2$, the binary equation and the actual third-mode
cap equality give
\begin{align*}
 \chi(\beta)&=x\int_r^\beta\frac{h^2-u^2}{(1-u)^2}\,du,\\
 m(\beta)&=-\frac{\widehat K(1+x)^2}{4}
       -\frac{2\lambda x}{1-\beta}-\frac{Kx^2}{(1-\beta)^2},\\
 f(\beta)&=\frac{h^2-\beta^2}{2}q(\beta)
       -\frac{1+x}{2(1-\beta)}\chi(\beta),\\
 m(\beta)+\frac12-f(\beta)&=0.
\end{align*}
These follow by integrating $\chi_s=h^2-u^2$, with $\chi(r)=0$,
and imposing the actual cap, rather than a determinant relaxation.
An antiderivative of the integrand is
\[
 \frac{2K-4}{1-u}-2\log(1-u)+(1-u).
\]
For an interior core with roots $-\xi_i,z_i,r_i$, incoming and
outgoing densities $x_-,x_+$, and actual own-pair shear $T_i$, exact
state propagation and matrix transport give
\begin{align}
 \frac{x_+}{x_-}
 &=\left(\frac{1+\xi_i}{1-\xi_i}
          \frac{1+z_i}{1-z_i}\frac{1-r_i}{1+r_i}\right)^2,
                                                        \label{pf:state}\\
 \lambda_-+\lambda_+
 &=K_i\left(T_i-\frac{x_-}{1-\xi_i}-\frac{x_+}{1-r_i}\right).
                                                        \label{pf:transport}
\end{align}
The incoming sign is fixed by whole-chart conjugacy. Both adjacent
root means equal $h_i^2$, namely the weighted means of $u^2$ on
$[-\xi_i,z_i]$ with density $(1-u)^{-2}$ and on $[z_i,r_i]$ with
density $(1+u)^{-2}$. These equations therefore share actual roots
and one matrix Casimir.

\subsection{A nonnegative cap defect}\label{pf:M003}
Put
\[
 t=\frac{4x}{(1+x)(1-r)(1-\xi)},\qquad
 \Phi(t)=t-1-\log t,\qquad e=\lambda+\frac32x.
\]
Positive shared length gives $0<t<1$. Substitution of the preceding
primitive in the cap equality, followed by
$t(1-r)(1-\xi)=4x/(1+x)$, gives the exact identity
\begin{align}\label{pf:positive-cap}
 e={}&\frac{1+x}{2}\bigl[\Phi(t)+(1-t)(r+\xi)+tr\xi\bigr]\notag\\
 &+\frac{1+x}{4}
 \left(\frac{r^2-h^2}{1-r}
       +\frac{\xi^2-\widehat h^2}{1-\xi}\right)\ge0.
\end{align}
Every summand is nonnegative by the actual root ordering. This argument
does not prescribe the sign of $\lambda$.

\subsection{Midpoint inequalities and the high-density sign}\label{pf:M004}
At a middle switch put
\[
 F_i=\frac{h_i^2-z_i^2}{2}q_i(z_i)\ge0,
 \qquad \Delta_i^2=\frac{K_i}{2}-\mathcal C.
\]
The maximizing cap inequality is $L_i^2+K_iF_i\le\Delta_i^2$.
For an overlap between cores $i$ and $k=i+1$, set
\[
 d_L=\frac{x}{1-r}-T_i^+,
 \qquad d_R=\frac{x}{1-\xi}-T_k^-.
\]
The exact midpoint matrices have
$L_i=\lambda+K_i d_L$ and $L_k=-\lambda-K_kd_R$.
Substitution of~\eqref{pf:shared-geometry} cancels $\lambda^2$ and gives
\begin{align}
 2\lambda d_L+K_i d_L^2+\frac{K_k(1+x)^2}{4}+F_i&\le\frac12,
                                                        \label{pf:mid-left}\\
 2\lambda d_R+K_k d_R^2+\frac{K_i(1+x)^2}{4}+F_k&\le\frac12.
                                                        \label{pf:mid-right}
\end{align}
Both inequalities concern the same $x$ and $\lambda$.

\paragraph{High-density sign.}\label{pf:M005}
The leaf-quadratic inequalities for the actual middle roots are
\[
 r^2+rz+z^2>3h^2,\qquad
 \xi^2-\xi\widehat z+\widehat z^2>3\widehat h^2.
\]
Together with $|z|\le h<r$ and $|\widehat z|\le\widehat h<\xi$,
they imply $z>-r/2$ and $\widehat z<\xi/2$. Exact density transport
now yields
\[
 d_L=\frac{x(1-2r-r^2+2z)}{(1-r)^2(1+z)},\qquad
 d_R=\frac{x(1-2\xi-\xi^2-2\widehat z)}{(1-\xi)^2(1-\widehat z)}.
\]
If both are nonpositive, these root bounds imply $r,\xi>3/10$ and
hence $(1-r)(1-\xi)<49/100$. But positive shared length gives
$(1-r)(1-\xi)>4x/(1+x)$. For
$x\ge x_*=2/\sqrt3-1>3/20$, the latter bound exceeds
$12/23>49/100$, a contradiction. At least one of $d_L,d_R$ is
positive. On that side, $K,\widehat K\ge3/2$ and
$3(1+x)^2/8\ge1/2$ in~\eqref{pf:mid-left} or
\eqref{pf:mid-right} force $\lambda<0$. No division by an
uncontrolled-sign expression is used.

\subsection{The density ceiling}\label{pf:M006}
Write $s=1-t$, $R=r+\xi$, $P=r\xi$, $d=R-P$, and
$\delta=(1-3x)/(1+x)=d+s(1-d)$. From~\eqref{pf:positive-cap},
\begin{align*}
 \frac{4e}{1+x}-(\delta^2-h^2-\widehat h^2)
 &\ge2\Phi(t)+2sR+2tP+r^2+\xi^2-\delta^2\\
 &=2\Phi(t)+2sd+R^2-\delta^2\\
 &\ge2\Phi(t)+2sd^2-s^2(1-d)^2\ge0.
\end{align*}
The final inequality uses $0\le d\le1$ and
$\Phi(1-s)\ge s^2/2$. Substitution in the determinant identity gives
\begin{equation}\label{pf:gap}
 g:=\frac34-\mathcal C\le a(x)e+e^2,
 \qquad a(x)=\frac{3(1-3x)}4.
\end{equation}
Dropping the nonnegative height terms directly in that determinant also
gives
\begin{equation}\label{pf:second-gap}
 g\le b(x)-3xe+e^2,
 \qquad b(x)=\frac{3(1-3x)^2}{16}.
\end{equation}
Suppose $x\ge47/250$. The high-density sign just proved gives
$0\le e\le3x/2$. On this interval the first right side increases
and the second decreases. Their intersection is
\[
 e_*(x)=\frac{(1-3x)^2}{4(1+x)}\le\frac32x.
\]
Thus $g\le Q(x):=a(x)e_*(x)+e_*(x)^2$. Both nonnegative factors
decrease on $[0,1/3]$, and exact rational comparison gives
\[
 Q(47/250)<\frac34-H_*.
\]
This contradicts $\mathcal C<H_*$. Therefore every overlap has
$x<47/250$. Notice that this proof uses neither a positive lower
cutoff on $x$ nor a small-height collar.

\subsection{Once-counted transport and the four-core contradiction}
\label{pf:M001}
Write the own-pair charges as
$P_i=(\xi_i+r_i)/2+K_iT_i$, and let $B_j$ be the cost of shared arc
$j$. For an own pair, exact density propagation dominates the comparison
density $x_{\rm out}/(1-u)^2$ on both actual arcs. Indeed, on the
second arc the quotient is
\[
 \left(\frac{1+r}{1-r}\frac{1-u}{1+u}\right)^2\ge1,
\]
and continuity at the middle switch propagates this bound to the first
arc. Integration gives
\[
 T\ge\frac{x_{\rm out}(\xi+r)}{(1-r)(1+\xi)}
       \ge\frac{x_{\rm out}r}{1-r}.
\]
The analogous incoming comparison gives
$T\ge x_{\rm in}\xi/(1-\xi)$. Apply these inequalities to the two
end cores and~\eqref{pf:transport} to the two middle cores. With
\[
 D_j=B_j+\frac{r_j+\xi_{j+1}}2+
 \frac{3x_j}{2}\left(\frac{r_j}{1-r_j}
                   +\frac{\xi_{j+1}}{1-\xi_{j+1}}\right),
\]
physical addition, with each shared arc included once, gives
\begin{equation}\label{pf:four-assembly}
 C_{\mathrm{in}}\ge D_1+D_2+D_3+e_1+2e_2+e_3+\frac{\xi_1+r_4}{2}.
\end{equation}
The discarded middle-core terms are
$\frac12h_i^2[x_{i-1}/(1-\xi_i)+x_i/(1-r_i)]\ge0$.

To bound a single $D$, put $A=1-r$, $B=1-\xi$, $S=A+B$,
$p=AB$, and $k=4x/(1+x)$. Then $t=k/p$, $k\le t\le1$,
$S\le1+p$, and direct integration of the shared arc gives
\[
 B_{\rm cost}=\frac{1+x}{2}S(1-t)+x\log t.
\]
Consequently
\begin{align*}
 D&=1-3x+S\left(\frac{x}{2}-\frac{(1+x)t}{8}\right)+x\log t\\
 &\ge 1-3x-\frac{(1+x)t}{8}
             +\frac{2x^2}{(1+x)t}+x\log t=:Q_1(x,t).
\end{align*}
The coefficient of $S$ is nonpositive, and
\[
 \partial_tQ_1=-\frac{((1+x)t-4x)^2}{8(1+x)t^2}\le0.
\]
It follows that
\[
 D\ge F(x):=\frac{7-18x-9x^2}{8(1+x)},\qquad
 F'(x)=-\frac{25+18x+9x^2}{8(1+x)^2}<0.
\]
Since $x<47/250$ and all $e_j\ge0$, equation~\eqref{pf:four-assembly}
yields
\[
 C_{\mathrm{in}}>3F(47/250)=\frac{206119}{198000}.
\]
This is a reusable once-counted bound on the complete common-lift chain.
\space
\subsection{Both actual returns and their density profiles}\label{pf:F09}\label{pf:M034}
Use the complete chronological source $01020212$, with charged
four-core $[2,5]$ and own pair $[3,4]$, on one genuine containing
normal lift. The left $010$ return becomes $202$ under whole phase
$+2$; its actual mean leaf $1$ runs from an inner root to the positive
outer root $a_L$. The right $212$ return becomes $020$ under whole
phase $+1$; its actual mean leaf $6$ runs from a negative outer root
to an inner root. These are direct forward chronological charts.

At either return, the weak maximizing root jets place the inner root
inside the square-mean strip and the outer root outside it.
Equality at a nonzero sign-changing root would give a nonzero second
derivative of the switching function and prevent the required sign
change; equality at zero would force the positive mean to vanish.
The weighted zero integral also excludes equality of the outer-root
magnitude with the mean height. Thus these two actual leaves supply
the strict root inequalities used in Section~\ref{pf:M005} on each
side of the four-core.

The derivations in Sections~\ref{pf:M002}--\ref{pf:M006} are local
to the shared arc, its two actual adjacent mean leaves, and the two
maximizing midpoint gaps. They do not use a second leaf beyond either
of those midpoints. Apply them on the left with the phase-$+2$ leaf
just specified, and on the right with the phase-$+1$ leaf. In both
cases the half-shear formulas, root inequalities, and positive shared
length imply that high density forces a negative shared $\lambda$.
The positive cap identity and the two gap inequalities then give,
on this single harmonically qualified lift,
\begin{equation}\label{pf:d8-density-ceilings}
 0<x_L,x_R<47/250,
 \qquad e_j=\lambda_j+3x_j/2\ge0.
\end{equation}
This proves both chronological density ceilings; it does not infer
the incoming assertion by reversing a one-sided scalar theorem.

Let $-\xi<z<r$ be the own-pair roots, $K=(3+h^2)/2$, and $T$ its
actual shear. Matrix transport gives
\begin{equation}\label{pf:d8-transport}
 KT=e_L+e_R+x_L\left(\frac K{1-\xi}-\frac32\right)
              +x_R\left(\frac K{1-r}-\frac32\right).
\end{equation}
With $M(u)=(1+u)/(1-u)$, state continuity gives
\[
 \frac{x_R}{x_L}=\left(\frac{M(\xi)M(z)}{M(r)}\right)^2.
\]
Writing $x=x_R$ and $q=x\bar q$, the four actual profile pieces are
\[
\bar q(u)=
\begin{cases}
 [M(r)/(M(z)M(\xi))]^2/(1+u)^2,&-\alpha\le u\le-\xi,\\
 [M(r)/M(z)]^2/(1-u)^2,&-\xi\le u\le z,\\
 M(r)^2/(1+u)^2,&z\le u\le r,\\
 1/(1-u)^2,&r\le u\le\beta.
\end{cases}
\]
Hence $S=xS_0$, with $S_0=\int\bar q\,du$, and
$A_D=(\alpha+\beta)/2+xW_0$, where
$W_0=\frac12\int(3+u^2)\bar q\,du>0$.
The same actual chart transformations give
\[
 1-a_L=\frac{4x_L}{(1+x_L)(1-\alpha)},\qquad
 1-a_R=\frac{4x_R}{(1+x_R)(1-\beta)}.
\]
The shared-length ratios are
$t_L=(1-\alpha)/(1-\xi)$ and $t_R=(1-\beta)/(1-r)$.

For later use retain the complete scale dependence of the exterior
cap losses. Put $b_a=\bar q(-\alpha)$, $b_b=\bar q(\beta)$, and
\[
 I_a=\int_\xi^\alpha\frac{t^2-h^2}{(1-t)^2}\,dt>0,
 \qquad I_b=\int_r^\beta\frac{t^2-h^2}{(1-t)^2}\,dt>0.
\]
The exact exterior integrals yield $E_a=xA_a-x^2B_a$, where
\[
 A_a=\frac{b_a}{2}\bigl[\alpha^2-h^2-(1-\alpha)I_a\bigr],
 \qquad B_a=\frac{b_a^2(1-\alpha)^3I_a}{2}>0,
\]
and analogously for $b$. With
$\Delta_A=A_a-A_b$, $\Delta_B=B_a-B_b$,
$\Sigma_A=A_a+A_b$, and $\Sigma_B=B_a+B_b$,
the cap-Casimir identity gives
\[
 \mathcal C=c_0+c_1x-c_2x^2,
\]
\[
 c_0=\frac K2-\frac{\Delta_A^2}{4S_0^2},\quad
 c_1=\frac{K\Sigma_A}{2}+\frac{\Delta_A\Delta_B}{2S_0^2},\quad
 c_2=\frac{K\Sigma_B}{2}+\frac{\Delta_B^2}{4S_0^2}
                         +\frac{K^2S_0^2}{4}>0.
\]
For fixed actual geometry, this is concave in the actual outgoing
scale. A valid closed scale interval can therefore be tested at its
endpoints. All actual return, cap, state, and signed apex conditions
remain constraints; independent arbitrary choices of geometry and
scale are not asserted to be realizable.

\subsection{An actual shear correlation and height floor}\label{pf:M033}\label{pf:M032}
The complete subsource at positions $1,\ldots,6$ is $102021$.
Corollary~\ref{loc:four-calibration}, using only the global budget
$C_*<2$, therefore supplies $0<\mathcal C<H_*<1471/2000$ on
the same selected full lift of every proper returned eight-source. The first cap gap implies
\[
 \frac34-\mathcal C\le\frac34(1-3x_j)e_j+e_j^2
                       \le\frac34e_j+e_j^2.
\]
Since
\[
 \frac{29}{2000}-\left(\frac34\frac{47}{2500}
                         +\left(\frac{47}{2500}\right)^2\right)
 =\frac{291}{6250000}>0,
\]
we have $e_L,e_R>47/2500$.

\begin{lemma}[Two-arc shear and density correlation]\label{pf:own-shear-correlation}
Let $0<a,b\le R<1$, $-a\le z\le b$, $q>0$, and define
\begin{align*}
 \frac Tq&=\frac{(1-z)(a+z)}{1+a}+\frac{(1+z)(b-z)}{1+b},\\
 x_L&=q(1-z)^2\left(\frac{1-a}{1+a}\right)^2,\\
 x_R&=q(1+z)^2\left(\frac{1-b}{1+b}\right)^2.
\end{align*}
Put $E=x_La/(1-a)+x_Rb/(1-b)$, $k=(1-R)/(1+R)$ and
$X=\max(x_L,x_R)$. Then
\[
 E\ge k(T+2qz^2),\qquad
 T\le\frac{2q(R-z^2)}{1+R},\qquad
 q\le X\left(\frac{1+R}{1-R}\right)^2.
\]
The first inequality is strict if $a,b<R$.
\end{lemma}
\begin{proof}
The ratio $(1-s)/(1+s)$ decreases, while the exact identity is
\[
 \frac{a(1-z)^2}{1+a}+\frac{b(1+z)^2}{1+b}-\frac Tq=2z^2.
\]
This proves the first bound. The derivatives of $T/q$ in $a,b$ are
$(1-z)^2/(1+a)^2$ and $(1+z)^2/(1+b)^2$; increasing both roots
to $R$ proves the second. For the third, solve the right density
identity if $z\ge0$ and the left if $z\le0$.
\end{proof}

Apply the lemma to the actual own pair with $a=\xi$, $b=r$,
$q=q_z:=q(z)=v(z)^{-2}>0$.
Its common-matrix identity \eqref{pf:d8-transport} becomes
$N:=e_L+e_R=K(T-E)-h^2(x_L+x_R)/2$. If $h\le1/14$,
Lemma~\ref{loc:leaf-correlations} gives $\xi,r<R:=\rho(h)=8h/(4-3h)<1$;
use $k=(1-R)/(1+R)$.
Consequently, using both actual density ceilings,
\[
 N<K(1-k)T<\frac{4K\rho(h)^2(47/250)}{(1-\rho(h))^2}
 \le\frac{442928}{12403125}<\frac{94}{2500}.
\]
The factors increase with $h$; the last reserve is
$46859/24806250>0$. This contradicts the two strict defects, proving
$h>1/14$ on the original source and lift.

We also retain the one-sided envelope on its full domain:
$0<h\le1/3$, $\xi,r<\rho(h)$, $|z|\le h$, $N>0$ and
$P=M(\rho(h))<2$, where $M(s)=(1+s)/(1-s)$. Write $R=\rho(h)$.
The same profiles give $x_L+x_R\ge2q_zk^2(1+z^2)$, so the lemma yields
\[
 N\le\frac{q_z}{(1+R)^2}
 \{6R^2-3(1+2R-R^2)z^2+h^2(R^2+2R-1-2z^2)\}.
\]
Here $R<1/3$, so the $h^2$ bracket is negative. Put
$A=6R^2$, $B=3(1+2R-R^2)$; then $0<A<B$.
Since $N>0$, $A-Bz^2>0$, allowing the one-sided density bound
$q_z\le x_R M(R)^2/(1+z)^2$. The completed square
\[
 \frac{AB}{B-A}-\frac{A-Bz^2}{(1+z)^2}
 =\frac{(A+Bz)^2}{(B-A)(1+z)^2}\ge0
\]
therefore proves the stronger strict version of the old envelope,
\[
 N<x_RF(P),\qquad F(P)=\frac{3(P-1)^2(P^2+2P-1)}{4(2P-1)}.
\]
The function $F$ increases for $P>1$: its positive factors increase, since
$\partial_P[(P^2+2P-1)/(2P-1)]=2P(P-1)/(2P-1)^2>0$.
No second density ceiling is used in this one-sided statement.

The leaf-quadratic inequalities give
\[
 \xi^2-\xi z+z^2>3h^2,\qquad r^2+rz+z^2>3h^2.
\]
With $W=\sqrt{12h^2-3z^2}\ge3h$, their positive roots imply
$\xi>(W+z)/2$, $r>(W-z)/2$, hence $\xi+r>3h$.
The exact own-pair cost is consequently
\begin{equation}\label{pf:d8-own-cost}
 P_{\mathrm{own}}=\frac{\xi+r}{2}+e_L+e_R
       +\frac{x_L(h^2+3\xi)}{2(1-\xi)}
       +\frac{x_R(h^2+3r)}{2(1-r)}.
\end{equation}
It follows already that $P_{\mathrm{own}}>2533/17500$.

A second useful correlation follows directly from the actual profiles:
\[
 T\le\frac{2\rho xM(\rho)^2M(h)^2}{(1-\rho)^2},
 \qquad KT>1/30.
\]
For $h\le1/20$, $x\le1/50$, the first right side multiplied by
$K$ is at most
\[
 \frac{1309566797}{90920340900}<\frac1{30}.
\]
Thus this closed rectangle is excluded. The assertion concerns the
actual own mean and density, not an arbitrary dummy mean.

\subsection{A refined defect and its cost consequence}\label{pf:M031}
Retain the following refined defect on the same full eight-source.
Put $a=47/2500$ and
$F_x(e)=\frac34(1-3x)e+e^2$. Direct algebra gives
\[
 \frac{29}{2000}-F_x(a+x/20)
 =\frac{291}{6250000}+\frac{73x}{25000}+\frac{11x^2}{100}>0.
\]
Convexity of $F_x$ on $[0,a+x/20]$ bounds it by its endpoint values,
both below $29/2000$. The actual gap therefore forces
$e_j>a+x_j/20$.

At the two actual endpoints, compactness $v\le9\sqrt3/2$ gives
\[
 x_L\ge\frac{4(1-\xi)^2}{243},\qquad
 x_R\ge\frac{4(1-r)^2}{243}.
\]
The factors depending on the actual endpoints are essential.
These endpoint factors and the refined defect are independent local
identities used in subsequent returned-core comparisons.

The sufficient bound $P_{\mathrm{own}}>1501/10000$ is an
immediate consequence of Proposition~\ref{pf:own-fifth}, since
$1/5>1501/10000$. That proposition uses the density, transport and
height statements above, not this weaker numerical consequence.

\paragraph{Applications.}\label{pf:M035}\label{pf:M036}
The defect and endpoint comparisons above are local source statements.
Their physical-cost consequences may be assembled by
Theorem~\ref{as:capacity}.
\space
\section{Shared-cap alternatives}\label{pf:shared-caps}\label{pf:F10}
\subsection{Rigidity on a capped three-core}\label{pf:M042}
The complete source $12021$ has two distinct internal binary switches
and third-mode caps. Lemma~\ref{pc:rigidity} therefore proves normality
and equality of any two jointly normalized full lifts on their connected
overlap. Corollary~\ref{pc:rigid-reconstruction}, with $k=K$, gives the
explicit reconstruction from its actual state path. Whole chart
alignment transports the matrix and the state covector; gluing adds
no endpoint or calibration conditions on an enlarged union.

\subsection{A shared cap costs more than one eighteenth}\label{pf:M040}
Retain a common normalized lift on which the preceding tail ends at
$\beta>0$, the following tail begins at $-\alpha<0$, and
\[
 L_p>\sqrt{H_1/2-\mathcal C}>3/25,
 \qquad L_q<0,\qquad
 0<\mathcal C<\overline C:=1470471/2000000.
\]
In the pre-cyclic cap chart its true mean is $0<\eta<1$,
$K=(3+\eta)/2$, and $-1<p<q<1$, with
$p<0$ and $p^2>\eta>q^2$. The actual state signs are
$U_p=\beta>0>U_q=-\alpha$. We prove that its complete physical
cost $E$ exceeds $1/18$.

Put $w=(1-u)^{-1}$, $k=\kappa>0$,
$D=1+k(2w-1)^2$, and $A=(2w-1)D/(4w)$. The exact same-cap
matrix reconstruction is
\[
 Z=B+2Kkw^2,\qquad H=\frac{4(Z^2+\mathcal C)}{KD^2},\qquad
 L=HA-Z.
\]
It implies
\begin{align*}
 H_w&=-32kwL/D^2,\\
 L_w&=-8k(2w-1)L/D+HA_w-4Kkw,\\
 A_w&=[1+k(2w-1)^2(4w+1)]/(4w^2)>0.
\end{align*}
The full cost is
$E=(q-p)/2+Kk(w_q-w_p)$. Use the increasing comparison parameter
\[
 e(w)=\frac{u(w)-p}{2}+Kk(w-w_p),\qquad
 e_w=\frac1{2w^2}+Kk>0,\qquad e(w_q)=E.
\]
It need not equal partial physical cost. Root ordering gives
$-p>|q|$, hence $q<E$ and
$1/2<w<1/(1-E)$ when $E<1$.
Let $w_*$ be the first zero of $L$; it exists because
$L_p>0>L_q$. All positive-$L$ comparisons stop there.

Define $a=8k(2w-1)/(De_w)$ and $b=4Kkw/e_w$.
For $1/2<w\le10/9$, we have $a<5/2$.
Indeed set $r=2w-1$, $b_0=1/(2w^2)$, and $t=r/w\le11/10$;
then
\[
 a\le\frac{8r}{K+b_0r^2+2r\sqrt{Kb_0}}
 <\frac{8t}{(2-t)(3/2+5t/3+t^2/2)}<\frac52.
\]
For the last inequality, the difference after multiplication is
$(90-41t-20t^2-15t^3)/12\ge49/800>0$.

Suppose $E\le1/18$, so $w<18/17$.
Dropping the positive $HA_w$ term gives
$L_e>-(5/2)L-4/(1-E)$. Gronwall to the first zero and
$e^x\le(1-x)^{-1}$ imply
\[
 L_p<\frac{4E}{(1-E)(1-5E/2)}\le\frac{144}{527}<\frac7{25}.
\]
Thus the same actual determinant satisfies
$\mathcal C>3/4-(7/25)^2=1679/2500>12E$.

The exact state on the cap is
$U(w)=[1-k(4w^2-1)]/D$.
Its outgoing negative sign gives $k>289/1007>1/4$;
the open-state condition gives $k<1/(2w-1)$.
On this strip $A_w/D^2\ge1/3$. To verify it, multiply by
$12w^2D^2$: the resulting quadratic in $k$ is concave, so it
suffices to use $k=1/4$ and $k=1/(2w-1)$.
At the latter endpoint the difference is
$(-8w^3+12w-3)/(24w^3)>0$ on $[1/2,18/17]$.
At $k=1/4$ it is
\[
 \frac{P_6(w)}{3w^2(4w^2-4w+5)^2},\quad
 P_6=-16w^6+32w^5-56w^4+88w^3-61w^2+15.
\]
Substitute $w=(1/2+(18/17)t)/(1+t)$, $t\ge0$.
The coefficients of $4\cdot17^6(1+t)^6P_6(w)$ are
\[
 772402208, 3771140192, 7631480812, 8532151624,
 5593973879, 1811941056, 70330092,
\]
all positive. This proves the whole closed-interval sign, including
its endpoint at infinity. Consequently
$HA_w\ge4\mathcal C/(3K)>2\mathcal C/3>8E$.

Put $N=2Kk(2w-1)-w^{-2}$.
The incoming positive state, $p<0$, and $\eta<p^2$ imply
\[
 N(w_p)<\frac{4(w_p-1)}{w_p(2w_p+1)}<0,\qquad
 N_e=\frac{4Kk+2/w^3}{Kk+1/(2w^2)}<8.
\]
Thus $N<8e\le8E$, and the exact equation becomes
\[
 L_e=-aL-2+\frac{HA_w-N}{e_w}>-\frac52L-2.
\]
At the first zero this gives
\[
 L_p<\frac45(e^{5/36}-1)\le\frac{221}{1854}<\frac3{25}.
\]
The rational bound uses $\mathcal E_1(5/36)$ in \eqref{strip:exp-upper}.
The final gap is $37/46350>0$, a contradiction. Hence $E>1/18$.

\space
\subsection{Same-cap responses and HIGH-following comparison}\label{pf:M045}\label{pf:F11}
Retain the single normalized lift, positive determinant $\mathcal C$,
cap coordinates, strict root and state signs of Section~\ref{pf:M040}.
The preceding terminal and following entrance satisfy
\[
 L_p>3/25,\quad L_q<0,\quad H_{\mathrm{pre}}>3/2,\quad
 L_p^2>H_{\mathrm{pre}}/2-\mathcal C,\quad m_p<-1/2.
\]
If the following actual mean is HIGH, its initial-cap identity gives
\[
 m_q+\tfrac12=\tfrac12(\mu-\alpha^2)\rho
 +\tfrac12[(U-1)\rho+(U-1)^{-1}]\chi>0,
 \qquad U=-\alpha<0,\quad\rho>0,\quad\chi\le0.
\]
Thus $m_p<-1/2<m_q$ in one canonical chart. Under exactly this interface,
Section~\ref{pf:M061} proves $E>21/200$, hence also $E>1/10$.
Its proof uses only the following identities, not the weaker conclusion.

Use $w,k,D,A,H,L,e$ from Section~\ref{pf:M040} and set
\[
 a=\frac{8k(2w-1)}{De_w},\quad
 N=2Kk(2w-1)-w^{-2},\quad R=1+2Kkw^2,\quad\omega=D/R>0.
\]
The exact same-cap responses, with no intervening change of frame, are
\begin{equation}\label{pf:highcap-response}
 L_e=-aL-2+\frac{HA_w-N}{e_w},\qquad m_e=-\omega L.
\end{equation}
The second follows from $m_w=-DL/(2w^2)$.
For the forcing comparison $A_w/D^2\ge1/3$, multiplication by the
positive $12w^2D^2$ gives the numerator
\[
 3[1+k(2w-1)^2(4w+1)]-4w^2[1+k(2w-1)^2]^2.
\]
It is concave in $k$. The endpoint $k=(2w-1)^{-1}$ reduces to the
cubic $-8w^3+12w-3$; the endpoint $k=1/4$ reduces to the polynomial
$P_6$ of Section~\ref{pf:M040}. Section~\ref{pf:M061} checks both
on its entire comparison strip.
Finally, whenever $0<a<5/2$ and $w<5/4$,
\[
 \frac{D_w}{De_w}=a/2<5/4,\qquad
 \frac{R_w}{Re_w}=\frac{4wx}{(1+x)^2}\le w<5/4,
 \qquad x=2Kkw^2>0.
\]
Both terms are positive, so $|(\log\omega)_e|<5/4$.
These generic response and weight identities retain the full actual
matrix and state adjoint of the cap.

\space
\section{Caps adjoining a nonreturned three-core}
\subsection{The actual cap interface}\label{pf:M058}\label{pf:F14}
We retain the local cap and tail comparisons for an actual returned core
before a cap and a capped, possibly nonreturned, core after it. Every
statement below keeps its own complete-source hypotheses. 

\paragraph{The cap interface.}
The cap comparison uses an original minimizing trajectory on one genuine
containing interval, with its original endpoints and real winding. Choose
one complete normal maximizing lift, jointly normalized, and one fixed
canonical chart. Its determinant and actual preceding binary
constant are positive. Corollary~\ref{cal:one-point} therefore gives
$\mathcal C<H_*<1470471/2000000<1471/2000$ before the cap
comparison. In the pre-cyclic mode-$2$
chart the cap density is $k/(1-u)^2$, its mean parameter is
$K=(3+\eta)/2\in(3/2,2)$, and its endpoints satisfy
\[
 -1<p<q<1,\qquad p<0,\qquad p^2>\eta>q^2.
\]
In the canonical chart of the adjoining cores its endpoint states and
matrix entries satisfy
\begin{gather*}
 U_p=\beta>0,\quad U_q=-\alpha<0,\quad
 L_p>\frac3{25},\quad L_q<0,\\
 m_p< -\frac12,\quad
 L_p^2>\frac{H_{\rm pre}}2-\mathcal C,\quad
 H_{\rm pre}>\frac32.
\end{gather*}
The same-cap comparison of Section~\ref{pf:M040} gives the unconditional cap bound $E>1/18$.
Only these actual endpoint data and the full seam equations enter that
comparison; a subsequent return is unnecessary. If the actual following
mean satisfies $\mu_{\rm post}>\alpha^2$, its entrance identity gives
\[
 m_q+\frac12
 =\frac{(\mu_{\rm post}-\alpha^2)\rho_q}{2}
  +\frac{(U_q-1)\rho_q+(U_q-1)^{-1}}2\chi_q>0.
\]
Indeed, $\rho_q>0$, $U_q<0$ and $\chi_q\le0$. This recovers the
same matrix entry from $m_p<-1/2$ to $m_q>-1/2$.
The stronger cap comparison proved in Section~\ref{pf:M061} then applies.
The signed-strip lemma~\ref{loc:signed-strip}, applied to a complete
negative-entry $12021$, supplies $3/2<H_{\rm post}<2$ and
$L_{\rm post,initial}<0$ with both caps retained. It has no return
hypothesis.

\subsection{A strict $21/200$ HIGH-following cap bound}\label{pf:M061}
Under exactly the actual cap hypotheses of Section~\ref{pf:M045},
including $m_p<-1/2<m_q$ in one canonical chart, the complete physical
cap cost satisfies $E>21/200$. We give the enlarged-domain comparison;
its numerical constants are established on that entire domain.

Use $w,k,D,A,H,L,e,a,N,R,\omega$ and the exact responses
\eqref{pf:highcap-response}. In particular $e_w>0$ and $e(w_q)=E$.
This increasing comparison parameter is not partial physical cost:
the latter has derivative $1/(2w^2)+(3+u^2)k/2$, whereas
$e_w=1/(2w^2)+Kk$.
Suppose $E\le21/200$. Root ordering gives $q<E$.
For $q>0$, the negative exit and $p<-q$ give
\[
 E>q\left(1+\frac{3(1-q)}{(1+q)^2(3-q)}\right)>5q/3,
\]
since $3/(3-q)>1$ and
$35800/48841-2/3=9718/146523>0$.
The conclusion $q<3E/5$ also holds for $q\le0$. Thus
\[
 \frac12<w<W:=(1-3E/5)^{-1}\le\frac{1000}{937}<\frac{10}{9},
 \quad k>\frac{877969}{3122031}>\frac14,
 \quad k<\frac1{2w-1}.
\]
It follows that $e_w>1627969/2000000>4/5$.
The AM--GM damping calculation of Section~\ref{pf:M040} applies
on $w<10/9$ and gives $0<a<5/2$.

The forcing inequality $A_w/D^2\ge1/3$ holds on this larger strip.
As in Section~\ref{pf:M045}, its numerator is concave in $k$.
At $k=(2w-1)^{-1}$, the cubic $-8w^3+12w-3$ has positive endpoint
values $2$ and $67657141/822656953$ on $[1/2,1000/937]$.
At $k=1/4$, use the polynomial $P_6$ of Section~\ref{pf:M040}.
Under $w=(1/2+(1000/937)t)/(1+t)$, the coefficients of
$4\cdot937^6(1+t)^6P_6(w)$ are
\[
\begin{gathered}
21656462794215814688,\quad105370132207929454816,\\
212455881278683048268,\quad237156580490629921000,\\
155499680007928507351,\quad49897734978496407240,\\
1369645344670652540.
\end{gathered}
\]
They are all positive. Hence $HA_w>2\mathcal C/3$.

The incoming state and mean give explicitly
\[
 k<(4w_p^2-1)^{-1},\quad
 K<2-w_p^{-1}+(2w_p^2)^{-1},\quad
 N(w_p)<\frac{4(w_p-1)}{w_p(2w_p+1)}<0.
\]
The remaining comparison interval satisfies $k(w_q-w)\le E/K$,
so the outgoing state yields
\[
 k(4w^2-1)>1-(16/3)WE\ge377/937>2/5,
 \qquad Kk>\frac3{5(4w^2-1)}.
\]
For $w\ge4/7$, $N_e<7$ follows from
$7e_w-N_w=3Kk+(7w-4)/(2w^3)>0$.
For $1/2<w<4/7$, this expression exceeds
\[
 \frac{158w^3-80w^2-35w+20}{10w^3(4w^2-1)}>0.
\]
After $w=(1/2+(4/7)t)/(1+t)$, its numerator times $(1+t)^3$
has positive coefficients $9/4,7,789/98,1152/343$.
Thus $N<7e$ and the master inequality is
\begin{equation}\label{pf:highcap-master}
 L_e>-aL-b_{\mathcal C}(e),\qquad
 b_{\mathcal C}(e)=2+(35/4)(e-2\mathcal C/21)_+.
\end{equation}
On positive components replace $a$ by $5/2$; on negative components
drop $-aL>0$. The positive-weight calculation in
Section~\ref{pf:M045} gives $|(\log\omega)_e|<5/4$, hence
\begin{equation}\label{pf:highcap-weight}
 \omega_{\max}/\omega_{\min}<e^{21/160}<160/139.
\end{equation}

Let $e_0=11/200$, $P=\int L_+\,de$ and $Q=\int(-L)_+\,de$.
If $\mathcal C\ge7/10$, then $b_{\mathcal C}=2$ through $1/15$,
and until the first zero
\[
 L(e)>\frac{23}{25}e^{-5e/2}-\frac45
 \ge h(e):=\frac{144-2760e+3450e^2-2875e^3}{1200}.
\]
The cubic bound follows from the nonnegative fourth-order remainder
in the exponential. Its value at $e_0$ is $138107/76800000>0$;
the decreasing exponential comparison therefore stays positive
through $e_0$. Hence the first zero is later, and
\[
 P>\int_0^{e_0}h(e)de=\frac{67076339}{20480000000}.
\]
Each subsequent negative component starts at a zero later than $e_0$.
The master inequality bounds all of them by the common nonnegative
envelope $2(e-e_0)+(35/8)(e-1/15)_+^2$. Thus
\[
 Q<(E-e_0)^2+\frac{35}{24}(E-1/15)_+^3
 \le\frac{2677169}{1036800000},
\]
\[
 P-\frac{160}{139}Q>
 \frac{69857236801}{230584320000000}>0.
\]
If $0<\mathcal C<7/10$, the actual preceding terminal determinant
gives $L_p>1/5$, and $b_{\mathcal C}\le467/160<3$.
Until the first zero,
$L(e)>(7/5)e^{-5e/2}-6/5\ge1/5-7e/2$,
which at $e_0$ is $3/400>0$. Consequently
\[
 P>913/160000,\qquad Q<3/800,\qquad
 P-\frac{160}{139}Q>30907/22240000>0.
\]
Again the negative-area estimate allows every later sign component.
In both cases $\int\omega L\,de>0$, so the exact response gives
$m_q-m_p<0$, contradicting the actual recovery.
This proves $E>21/200$ on precisely the stated cap interface.

\subsection{The geometry of an opposite centre}\label{pf:M063}
Consider the complete positive source
$1\mid020\mid1\mid202\mid1$ on one normalized complete lift.
Let $C,D$ denote the entrance and exit of its central mode1 arc,
whose physical cost is $c$. Its complete seven-arc subsource $020\mid1\mid202$ meets
Lemma~\ref{ctr:height}: the actual $2012$ and $0120$ force
$y>z_C$, $x_C>z_C$, $y>x_D$ and $z_D>x_D$ on this same lift.
Use $x=(1+u)/v$, $y=(u^2+v^2-1)/(2v)$ and $z=(1-u)/v$.
With $Y=y^2$ and $a=x_C>b=x_D$, the source-qualified rational
bound \eqref{ctr:geometric-height} gives $c>\mathcal R(Y)$ for $Y<1/3$;
in particular $c\le1/10$ implies the actual crossing height $V^2>14/5$.
This uses neither outward mean nor the two additional exterior caps.

For use with the full central transport, the exact entrance data are
\[
 u(t)=\frac{t^2+2ty-1}{1+t^2},\quad
 v(t)=\frac{2(t+y)}{1+t^2},\quad
 \alpha_L=u(a),\quad\alpha_R=-u(b),
\]
\[
 Y_L=a^2,\quad Y_R=\left(\frac{1-by}{b+y}\right)^2,
 \qquad Q_i(U)=Y_i/(1-U)^2.
\]
The $Y_i$ are density coefficients, not initial densities.
Full left reflection gives
$(L_L,m_L,\chi_L)=(-L_C,m_C,-\chi_C)$; the right data are
$(L_R,m_R,\chi_R)=(L_D,m_D,\chi_D)$.
The true means and both $K_i=(3+\mu_i)/2$ are actual, and
$\mathcal C=-L_i^2-K_im_i$ is the same determinant at both ends.
These data remain joined by the central matrix transport and the
state-adjoint cost identity \ref{ctr:affine-cost}; they are not
independent entrance parameters.

\space
\section{Returned-terminal source qualification and packing}\label{pf:F15}
\subsection{The independent nine-source geometry}\label{pf:M075}
Retain a proper complete positive source $012021201$, with charged
positions $2,\ldots,6$. Its ascending returned context
is local $(1,2,3,4,5,6)$; the descending context is $(7,6,5,4,3,2)$
with whole reflected map $j\mapsto1-j$. Both are $120212$ and share
local terminal arc $4$. The full forward sources on local $[0,3]$
and $[5,8]$ give the two strict negative entrances on the same selected
full lift. All nine source arcs are positive and proper. Reanchoring
in an omitted positive arc and capped-three-core normality supply a
genuine complete jointly normalized lift for each union separately.

After whole phase $+1$, local $[1,7]$ is $2010201$, the common
returned-terminal configuration in Lemma~\ref{ret:structure}.
The eleven-source's exterior arcs establish entrance signs; here the two
complete forward sources establish those same signs. Both actual returns,
all caps, signed-strip roots, full matrix and state adjoint, compact cuts
and complete endpoint problem are retained. The two-chart bridge,
actual mean identities and physical-cost comparisons therefore apply
on this selected lift without a numerical union estimate. These are the
independent geometric identities used by the local transfer proof of
Lemma~\ref{six:lem:H3}; they do not use its numerical conclusion.

\subsection{Three returned unions}\label{pf:M078}
Consider $\mathtt{--+++--+++--++++}$, with $B=16$, native controls
$0212010201210120$ and twist $j_{i+16}=j_i+1$.
The complete sources $[-3,5]$, $[2,10]$ and $[7,15]$, with whole
phases $0,+2,+1$, respectively, are $012021201$.
Their five-arc charges are
\[
 \{15,0,1,2,3\},\qquad\{4,5,6,7,8\},\qquad\{9,10,11,12,13\}.
\]
Each source is proper and retains both actual returned contexts, both
forward entrance witnesses and all outgoing leaves and ties of the
local H3 theorem. Cuts inside physical arcs $12,1,6$, respectively,
provide genuine containing periods. Apply the local nine-source
bound~\eqref{six:eq:union-input} on each selected full lift separately.
The charge sets are disjoint and omit only positive arc $14$, so
\[
 C>3\frac{180483}{250000}
  =\frac{541449}{250000}=C_*+\frac{180483}{500000}>C_*.
\]
No opposite-centre estimate, LOW-mean alternative or identification of
these three full lifts is needed.
\space
\subsection{Opposite-centre transport and the endpoint minimax}\label{pf:M084}
This subsection gives the mechanism for the unconditional
opposite-centre estimate. Select one genuine complete minimizing
interval containing the full positive source
$1\mid020\mid1\mid202\mid1$, with the original state and group endpoints,
real winding, all six compact cuts and the full measurable-control maximum
principle. The complete forward sources $2012$ and $0120$ give
$u_C>0>u_D$ in the central mode-$1$ frame. The outward left chart is
full reflection and time reversal, not an isolated label change.
The signed-strip proof~\ref{loc:signed-strip} treats both signs of
$(\mu-\alpha^2)/2$ and gives negative outward entrance $L$ without a
HIGH assumption. Thus $L_C>0>L_D$. Its two internal binary zeros and
both positive cap ties prove normality; normalize the entire lift once.

With $j=-1$ at $C$ and $j=1$ at $D$, the actual cap identities are
\[
 \chi_C\ge0\ge\chi_D,
 \qquad e_i:=\beta_j\chi_i\ge0,
 \qquad
 m_i+\frac12=\frac{\mu_i-u_i^2}{2v_i^2}+e_i,
 \qquad
 \beta_j=\frac12\left(\frac{u-j}{v^2}+\frac1{u-j}\right).
\]
The genuine means satisfy $0<\mu_i<1$ and $K_i=(3+\mu_i)/2$.
The cap identity gives $m_i>-(1+u_i^2/v_i^2)/2$.

\paragraph{Geometry and exact central coordinates.}
Let $y=(u^2+v^2-1)/(2v)>0$ be the central constant reciprocal and
$Y=y^2$. The mandatory interval and density calculation of
Section~\ref{pf:M063} gives, for $Y<1/3$,
\[
 c>\mathcal R(Y):=\frac{(1-3Y)(3Y+7)}{8(1+Y)},
 \qquad
 \mathcal R'(Y)=-\frac{9Y^2+18Y+25}{8(1+Y)^2}<0.
\]
At the actual $u=0$ crossing put $q=V^2$, so
$y=(q-1)/(2\sqrt q)$ and $Y=(q-1)^2/(4q)$. The compact cuts give
$q\le243/4$. The lawful cyclic and reflected charts preserve those cuts.
The exact centered state, full nonnegative adjoint identity, cap conjugation
and corrected quadrature are retained in Section~\ref{ctr:profiles}.
Its profiles use the same actual $q$ throughout; the independent mean
budget there requires neither the signed entrances nor calibration.

The source-independent transport and endpoint-loss results of
Section~\ref{ctr:generic-transport} now apply with these actual signed entrances.

The source-independent signed profile of Theorem~\ref{ctr:signed-profile}
applies with these proved signs; its circle, transport and endpoint-loss
constructors retain the actual $q$ and the separately qualified profile.

\paragraph{The two required rational certificates.}\label{ctr:essential-certificates}
All scalar guards in Section~\ref{ctr:profiles} are evaluated separately
at the printed budgets. Their height envelopes have elementary proofs:
write $E=(q-1)(q+3)$. For $1<a\le3$ and $q\ge a$,
\begin{align*}
 H'(q)&=-4(q^2+3)/E^2<0,\\
 r(a)-r(q)&=\frac{2(q-a)[(3-a)q+3(a+1)]}{(a-1)(a+3)E}\ge0.
\end{align*}
The numerators of $4-qH$, $8-q(q+1)H^2/2$, and
$64/25-(q-1)H^2$, after clearing their positive denominators, are
respectively positive multiples of
$2q-3$, $3(q-2)^3+16(q-2)^2+16(q-2)+1$, and
$(q-2)(4q^2+3q+18)$. Thus the stated $q\ge2$ bounds follow.
Each quadratic numerator of $K$ has positive leading coefficient and
negative discriminant, so $K,M,N>0$. On the opposite-height domain,
$A_q=qH<4$, $B_+<3/8$ because $3z(27/10)-7=62129/100000>0$, and
\[
 1-2kc^2>1-3(63/500)^2=238093/250000>0.
\]
Also $A_q-3B_-=[4(q-3/4)^2+15/4]/[(q-1)(q+3)]>0$.
Hence $p_0,D,Z$ and all their defining denominators are positive directly.
The optional inequalities $T-Z>0$ are unnecessary for the endpoint-loss
lemma. 

Only the following numerator signs remain. Write each indicated rational
function in coprime integer polynomials $P/Q$ with positive-leading $Q$.
Its denominator is positive by the preceding factorizations. For $P(q)=\sum_{j=0}^np_jq^j$ on $[l,h]$, the coefficient of $t^k$ in
$(1+t)^nP((l+ht)/(1+t))$ is exactly
\[
 \sum_{j=0}^n p_j
 \sum_{v=\max(0,k-n+j)}^{\min(j,k)}
 \binom jv l^{j-v}h^v\binom{n-j}{k-v}.
\]
Every such coefficient in these two rows is at least $1$.
Here $R=4\delta T-(T+Z)^2$, $\delta=29529/2000000$, and each row
uses its own profile's unchanged budget and constructors:
\begin{center}\small
\begin{tabular}{cclc}
Profile&Function&Interval&$\deg P$\\\hline
O&$R$, ordinary&$[27/10,3/z]$&34\\
O&$R$, ordinary&$[3/z,243/4]$&34\\
\end{tabular}
\end{center}
The 70 exact positive coefficients certify the closed endpoints as well
as the interiors. The portable \texttt{check\_center\_essential.py}
reconstructs these coefficients and factorizations. These are rational
identities and sign checks, independent of either opposite-centre conclusion.

A resulting floor $C_{\rm mat}>\overline C>0$ and the actual
$K_i>0$ license Corollary~\ref{cal:one-point} on the selected
source-qualified complete lift. It gives
$C_{\rm mat}<H_*<1470471/2000000$.

\space
\section{Own pairs and bilateral response}\label{pf:own-bilateral}
\subsection{A one-fifth bound for an own pair}\label{pf:M086}\label{pf:F18}

\begin{proposition}\label{pf:own-fifth}
Let $01020212$ be a complete positive chronological source on an original
minimizing trajectory, or a whole cyclic image of this source. Retain both
exterior caps and the returned mean arcs at positions $1$ and $6$.
Then the cost $P=c_3+c_4$ of its own pair satisfies $P>283/1000>249/1000>1/5$.
For a cyclic application this source is required to be proper, so $B>8$.
\end{proposition}
\begin{proof}
The complete $01020212$ source contains $102021$ at positions
$1,\ldots,6$. Corollary~\ref{loc:four-calibration} proves normality,
positive determinant and future trace before harmonic calibration on
one genuine containing lift. Its two chronological return means give
the density and transport identities of Sections~\ref{pf:M034}
and~\ref{pf:M033}. These are exactly the interface
\eqref{six:eq:own-interface}, as derived there without a numerical
own-pair bound. Proposition~\ref{six:prop:own} therefore gives
$P>283/1000$. It charges the same complete arcs $3,4$ and retains
both caps and returns; the displayed weaker bounds follow immediately.
\end{proof}

\space
\section{Two-return four-core estimates}\label{pf:two-return}
\subsection{A correlated density envelope and a four-core lower bound}
\label{pf:M096}\label{pf:F23}

\begin{lemma}[Correlated two-density shear envelope]\label{pf:two-density-envelope}
Let an actual positive binary core have endpoints $-\alpha,\beta$,
where $0<\alpha,\beta<1$, and let its first and last density profiles be
$x_L/(1+u)^2$ and $x_R/(1-u)^2$.  Suppose that
\begin{equation}\label{pf:two-density-ceilings}
 x_L<X,\qquad x_R<X,\qquad X>0.
\end{equation}
Put $d=\alpha+\beta$, $p=\alpha\beta$, and
$S=\int_{-\alpha}^{\beta}q(u)\,du$, where $q=v^{-2}$ is the actual
continuous density.  Then
\begin{equation}\label{pf:density-envelope}
 S<X\frac{d+2p}{1-d+p}.
\end{equation}
If $d<2/3$, one also has
\begin{equation}\label{pf:Phi-envelope}
 S<\Phi(d):=X\frac{2d(d+2)}{(2-d)^2},\qquad
 \Phi'(d)=\frac{4X(3d+2)}{(2-d)^3}>0.
\end{equation}
\end{lemma}
\begin{proof}
On each actual binary piece the logarithmic derivative of $q$ is
$-2/(1+u)$ or $2/(1-u)$. Continuity at the actual switches makes
$q$ positive and absolutely continuous on the full core. Its endpoint
scales are exactly
$q(-\alpha)(1-\alpha)^2=x_L$ and
$q(\beta)(1-\beta)^2=x_R$.
Apply Lemma~\ref{cv:correlated-profile} to this same actual profile.
Its endpoint-density conclusion is~\eqref{pf:density-envelope}.
The additional reduction to~\eqref{pf:Phi-envelope} uses $d<2/3$;
this restriction is retained, and direct differentiation gives the
stated positive derivative.
\end{proof}

For the complete two-return source $01020212$, the chronological
calculation in Section~\ref{pf:M034} supplies these ceilings with
$X=47/250$ on its harmonically qualified common lift.  The left return
$010$ at positions $0,1,2$ is a whole phase of $202$, with mean leaf~$1$;
the right return $212$ at $5,6,7$ is a whole phase of $020$, with mean
leaf~$6$.  Both are actual chronological mean leaves, so the incoming
ceiling is not inferred from an unqualified reversed flow.  The envelope
lemma and those local return calculations use no four-core numerical
cost floor.

\begin{proposition}\label{pf:four-eleven-twentieths}
On a proper complete positive source $01020212$, or a whole cyclic
image, the cost $A$ of arcs $2,3,4,5$ satisfies $A>11/20$.
Both returned mean leaves and both exterior caps are retained.
\end{proposition}
\begin{proof}
This is the identical proper eight-source and physical charge of
Section~\ref{pf:M101}.  A cut inside an omitted positive arc retains all
eight complete arcs on a genuine containing interval, with both actual
returns and exterior caps.  The stronger four-core estimate there gives
\[
 A>\frac{74}{125}>\frac{11}{20}.
\]
Its proof uses the local chronological return ceilings, the correlated
envelope above, and its own cap-loss and closed-domain comparisons.  It
does not use this weaker bound or the fixed-word outer bootstrap.
\end{proof}

\paragraph{Applications.}\label{pf:M097}
The same complete-source payment is available at every qualified
occurrence; costs can be combined subject to Theorem~\ref{as:capacity}.
\space
\subsection{A complete two-return four-core}\label{pf:M101}
We first prove a reusable source bound, independently of the fixed-word
estimates below. Its inputs are binary normality, the actual apex,
harmonic calibration, both chronological return-density ceilings, and
the profile and cap-loss lemmas proved in the following subsections.
Neither a prior four-core numerical floor nor any outer-pair cost bound
is used. In particular, this theorem supplies the weaker $11/20$
corollary rather than depending on it.

Consider a proper positive chronological source $01020212$, or a
whole lawful additive phase of it, on the minimizing trajectory.
All eight source arcs are complete and the number of physical arcs
is greater than eight. Charge only local arcs $2,3,4,5$, the four-core
$0202$. Its physical cost $A$ satisfies
\[
 A>\frac{74}{125}.
\]
Reanchor inside an omitted positive arc to choose a genuine complete
minimizing interval containing the source.

\paragraph{Normality and calibration.}
Suppose $A\le74/125<1$. The actual capped source
$1\mid0202\mid1$ proves normality by
Lemma~\ref{loc:binary-normality}; normalize the entire multiplier
tuple once. Its actual matrix apex~\ref{loc:four-core-apex} gives
\[
 \mathcal C>\frac34-\frac{9A^2}{64}
 \ge\frac{175179}{250000}>0.
\]
Corollary~\ref{loc:four-calibration} gives, on this same lift,
$\mathcal C<C_0:=H_*<\overline C:=1470471/2000000<3/4$.
Both chronological return ceilings are those proved in
\eqref{pf:d8-density-ceilings}:
$x_L,x_R<X=47/250$. The left return is local $0,1,2$, sent
from $010$ to $202$ by phase $+2$, and its genuine positive mean
leaf is local~$1$ in the inner-to-outer direction.
The right return is local $5,6,7$, sent from $212$ to $020$ by
phase $+1$, with mean leaf~$6$ in the outer-to-inner direction.
Their shared charged arcs are respectively~$2$ and~$5$.
The two signs and ceilings are proved in these actual chronological
charts, preserving their own mean equations and common determinant.
They require no prior four-core cost floor.

\paragraph{The actual contrary domain.}
Write the outer endpoints as $-\alpha,\beta$, inner switches as
$-\xi<z<r$, and the genuine mean as $h^2$. Put
\[
 d=\alpha+\beta,\quad p=\alpha\beta,\quad K=(3+h^2)/2,
 \quad S=\int q\,du,\quad M_2=\int u^2q\,du,\quad q=v^{-2}.
\]
The exact four-root and physical-cost identities give
\[
\begin{gathered}
 0<h<\xi<\alpha<1,\quad h<r<\beta<1,\quad |z|\le h,\\
 2A=d+3S+M_2,\quad M_2>h^2S,\quad d>3h,\quad p>2h^2,\\
 p>h(d-h),\quad p\le d^2/4,\quad p<d^2/3-h^2.
\end{gathered}
\]
The last upper bound is the actual endpoint-asymmetry consequence.
The structural Casimir estimate gives $S>81/500$: otherwise
$\mathcal C>3/4-(9/16)(81/500)^2=\overline C+9/4000000$.
Therefore $d<349/500$.
Positive shear completion gives
\[
 A>d/2+KS>3h/2+\sqrt{3+h^2-4\overline C}.
\]
At $h=1/5$ the unsquared remaining budget is $73/250>0$,
and the squared excess is $6897/500000>0$; the right side
increases with $h$. Hence $h<1/5$.
The independent correlated-density Lemma~\ref{pf:two-density-envelope}
gives the raw two-density envelope
\[
 S<E_1:=X\frac{d+2p}{1-d+p}.
\]
For $d\le43/100$ its derivative in $p$ is positive, so
$p\le d^2/4$ gives the symmetric upper envelope
\eqref{pf:Phi-envelope}. It increases
with $d$ and at $43/100$ is below $81/500$ by
$32157/12324500$. Thus $d>43/100$. The symmetric reduction
is used only on this small-width branch.

\paragraph{The two actual loss functions.}
Define
\[
\begin{gathered}
 e_2(t,h)=\frac{t^2-h^2}{2}
 -\frac{2(t-h)^2(t+2h)}{3(1-t)(3+t)},\\
 \ell_a=\frac{(1-\beta)e_2(\alpha,h)}{d(1+\alpha)},\quad
 \ell_b=\frac{(1-\alpha)e_2(\beta,h)}{d(1+\beta)},\\
 Q=\ell_a+\ell_b,\quad P=\ell_a\ell_b.
\end{gathered}
\]
Both losses are positive on the contrary domain. Indeed
$(t-h)(t+2h)\le t(t+h)$ reduces positivity to
$3(1-t)(3+t)-4t>0$, whose decreasing endpoint value at
$t=349/500$ is $139597/250000>0$.
The structural cap-loss Lemma~\ref{loc:double-return-loss}, using the
actual return-circle bounds~\eqref{loc:double-circle}, gives this same $e_2$.
The exact cap integrals therefore give actual normalized losses
$u=E_a/S>\ell_a>0$, $v=E_b/S>\ell_b>0$.
Their apex identities are
\[
 u+v<KS,\qquad
 \mathcal C=K/2-(KS-u-v)^2/4+uv.
\]
The functions $Q$ and $P$ are symmetric in $\alpha,\beta$, hence
rational in $d,p,h$. Their displayed endpoint formulas are also
unambiguous exact constructors for every polynomial used below:
expand and eliminate $\alpha+\beta=d$, $\alpha\beta=p$.
No replacement by a weaker single-loss function is made.

Section~\ref{pf:M112} proves the same-profile upper bound
\[
 S<E_2:=X\left(\frac{2-d}{1-d+p}-2B_R\right),\quad
 R=h(26/15+3h/4+h^2),\quad B_R=\frac{1-2R-R^2}{(1-R)^2}.
\]
Section~\ref{pf:M110} strengthens the actual product to
$uv>\lambda^2P$, where $\lambda=2(1+p)/(2+2p-d)>1$.
The positive branch $KS>u+v>Q$ is established before taking roots;
thus every contrary point satisfies
\[
 KS>Q+\sqrt{3+h^2-4\overline C+4\lambda^2P}.
\]
Put $T=74/125-d/2$. Since $KS<T$ and $KS>2\lambda\sqrt P$,
the decreasing quadratic comparison in Section~\ref{pf:M110} gives
\[
 B:=\overline C-K/2+T^2/4>0,\qquad B^2>T^2\lambda^2P.
\]
These are necessary inequalities for the same actual lift, without
any converse assertion about a solution of the scalar inequalities.

\paragraph{Exhaustion.}
Section~\ref{pf:M109} proves $d<3/5$ under the weaker contrary
budget $A\le297/500$. Its proof uses no four-core cost floor.
Section~\ref{pf:M033} supplies the independent genuine own-height floor
$h>1/14$ on this same full eight-source. Its proof uses the actual
return defects, the radius theorem and the two density ceilings, with
no four-core numerical floor. Sections~\ref{pf:M102}--\ref{pf:M104}
then exclude the remaining domain by a monotone complete product comparison
and six exact polynomial comparisons. No mean from a returned chart
is substituted for this own height.

\subsection{One mean-root profile controls both halves}\label{pf:M112}
This lemma retains a stronger same-profile density envelope for the
two-return source. It holds under its actual source hypotheses whenever
$A\le303/500$. The same normality, apex and future-sheet argument
then calibrates one full lift. The structural estimates give
$d<363/500$, $h<1/5$, and $x_L,x_R<X=47/250$.
For the height bound, use
$A>3h/2+\sqrt{3+h^2-4\overline C}$; at $h=1/5$ this already
exceeds $303/500$, with a positive squaring branch.

\paragraph{A scalar zero-facing root barrier.}
Let $r_0(h)>h$ solve
\[
 \int_0^{r_0}\frac{u^2-h^2}{(1+u)^2}\,du=0.
\]
The integral is negative at $h$ and increases thereafter. Put
$R=h(26/15+3h/4+h^2)$; for $0<h\le1/5$ one has
$h<R<1/2$. Its integral value is
\[
 I=R+(1-h^2)\left(1-\frac1{1+R}\right)-2\log(1+R).
\]
With $q=R/(2+R)$, the positive atanh remainder gives
$\log(1+R)<2[q+q^3/3+q^5/(5(1-q^2))]$.
The resulting rational lower bound for $I$ factors as
\[
 \frac{2h^3(60h^2+45h+104)P_0(h)}
 {75(60h^3+45h^2+104h+60)(60h^3+45h^2+104h+120)^3},
\]
where
\[
\begin{split}
 P_0(h)={}&6480000h^{10}+11340000h^9+48573000h^8+77834250h^7\\
 &+128081925h^6+173342700h^5+149309880h^4+125763360h^3\\
 &+72622928h^2+166400h+96000.
\end{split}
\]
All coefficients are positive. Hence $I>0$, and the unique root
satisfies $h<r_0<R$.

The two actual own mean equations, with their same central switch $z$,
are
\[
 \int_z^r\frac{u^2-h^2}{(1+u)^2}\,du=0,
 \qquad \int_{-z}^{\xi}\frac{u^2-h^2}{(1+u)^2}\,du=0.
\]
The second is a substitution in the actual left integral, not a
reversal of the physical PMP lift. If $z\ge0$, then $r\le r_0<R$;
if $z\le0$, then $\xi\le r_0<R$. Only one of these root bounds
is asserted on each branch.
With $M(t)=(1+t)/(1-t)$ and $q_z=q(z)$, continuity gives
\[
 x_L=\frac{q_z(1-z)^2}{M(\xi)^2},\qquad
 x_R=\frac{q_z(1+z)^2}{M(r)^2}.
\]
When $z\ge0$, the right bound and $(1-z)^2\le(1+z)^2$ imply
both estimates below; when $z\le0$, use the left bound instead:
\[
 q_z(1-z)^2<XM(R)^2,\qquad q_z(1+z)^2<XM(R)^2.
\]
This retains the same actual $z,q_z$ and both linked endpoint densities.

\paragraph{Integration of the two half-profile envelopes.}
For $u\le0$, the first mode-$0$ coefficient is $x_L<X$ and the
following mode-$2$ coefficient is $q_z(1-z)^2<XM(R)^2$.
If the third mode-$0$ arc begins before zero, monotonicity of $M$
bounds it by the continuation of the same mode-$2$ upper curve.
Thus
\[
 q(u)<X\max\{(1+u)^{-2},\ M(R)^2(1-u)^{-2}\}\qquad(u\le0).
\]
These curves meet at $u=-R$. If $\alpha\ge R$, integration gives
\[
 \int_{-\alpha}^0q\,du<X\left(\frac1{1-\alpha}-B_R\right),
 \qquad B_R=\frac{1-2R-R^2}{(1-R)^2}.
\]
If $\alpha<R$, the actual endpoint comparison gives the stronger
bound $X\alpha(1+\alpha)/(1-\alpha)^2$. It is at most the displayed
extension because $B_t'=-4t/(1-t)^3<0$ and the two formulas agree
at $t=R$. The positive half is identical using the actual last
mode-$2$ coefficient and the second same-profile bound. Adding gives
\[
 S<X\left(\frac{2-d}{1-d+p}-2B_R\right)=E_2.
\]
All denominators are positive; no symmetry of $\alpha,\beta,z$ is used.

\paragraph{The two retained cap losses.}
On $d\le363/500$, positivity of $e_2$ follows from
$3(1-363/500)(3+363/500)-4(363/500)>0$.
The actual cap integrals give $u>\ell_a$, $v>\ell_b$.
Apply Lemma~\ref{cv:guarded-loss} with
$(\kappa,t,a,b,W)=(K/2,KS,\ell_a,\ell_b,P)$.
Its guard is the actual apex inequality $KS>u+v$; its monotone
loss lowering and positive-root conclusion give
\[
 \mathcal C>K/2-(KS-Q)^2/4+P,
 \qquad KS>Q+\sqrt{3+h^2-4\overline C+4P}.
\]
This establishes all retained-loss and profile premises without a
four-core numerical cost improvement.

\subsection{A joint endpoint-density bound for the loss product}\label{pf:M110}
The two endpoint density coefficients belong to one actual four-core
profile. This correlation strengthens the product of its cap losses.
Assume the complete source, full-lift calibration and positive-loss
domain of Section~\ref{pf:M112}; the absolute upper bounds on
$x_L,x_R$ are not needed in this lemma.
At the actual endpoints,
\[
 q(-\alpha)=\frac{x_L}{(1-\alpha)^2},\qquad
 q(\beta)=\frac{x_R}{(1-\beta)^2},\qquad x_L,x_R>0.
\]
The cap integrals give
$u>q(-\alpha)e_2(\alpha,h)/S$ and
$v>q(\beta)e_2(\beta,h)/S$, strictly.

On each binary arc the logarithmic slope is one of the two endpoints
$-2/(1+u)$ and $2/(1-u)$, so Lemma~\ref{cv:correlated-profile} applies
to this one continuous actual profile. With
$e_\alpha=e_2(\alpha,h)$, $e_\beta=e_2(\beta,h)$, its exact product is
\[
 P=\frac{(1-d+p)e_2(\alpha,h)e_2(\beta,h)}{d^2(1+d+p)}
   =\ell_a\ell_b.
\]
The lemma therefore gives
\[
 uv>\frac{(1+\sqrt Z)^2}{4}P>\lambda^2P,\qquad
 Z=\frac{1+d+p}{1-d+p},\quad
 \lambda=\frac{2(1+p)}{2+2p-d}>1.
\]
The common profile, rather than independently chosen endpoint boxes,
is the hypothesis retaining this strict loss-product correlation.

\paragraph{Two correctly branched Casimir comparisons.}
Apply Lemma~\ref{cv:guarded-loss} with
\[
 (\kappa,t,a,b,W)=(K/2,KS,\ell_a,\ell_b,\lambda^2P).
\]
It gives
\[
 \mathcal C>K/2+KS\lambda\sqrt P-K^2S^2/4,
 \qquad KS>Q+\sqrt{3+h^2-4\overline C+4\lambda^2P}.
\]
Under $A\le A_*$ the physical moment identity gives
$KS<T_*=A_*-d/2$. The same lemma therefore gives
\[
 B:=\overline C-K/2+T_*^2/4>T_*\lambda\sqrt P>0,
 \qquad B^2>T_*^2\lambda^2P.
\]
Any larger positive ceiling for $KS$ is permitted by that lemma.
These are necessary comparisons on the actual lift; in particular,
the squared test is used only after its positive branch is proved.

\subsection{Two product-coordinate monotonicities}\label{pf:M111}
Let $Q,P$ be the exact rational functions constructed from the cap
losses in Section~\ref{pf:M101}, and set $W=\lambda^2P$.
The correlated product satisfies $W_p>0$ on the actual real-endpoint domain
$0<h<\alpha,\beta$, $\alpha+\beta=d\le3/5$,
including the symmetric limit $p=d^2/4$.

\paragraph{Strict log-concavity of one loss.}
For $0\le h<t\le3/5$, write $e(t):=e_2(t,h)=(t-h)f(t)$, where
\[
 f(t)=\frac{8h^2+(9-10t-3t^2)(h+t)}{6(1-t)(3+t)}>0,
\]
since $9-10t-3t^2\ge48/25$. Direct differentiation gives
\[
 f''(t)=-\frac{(1-h)(1+2h)}{3(1-t)^3}
         -\frac{(h+3)(3-2h)}{3(3+t)^3}<0.
\]
Both numerators are positive on the stated domain. Consequently
$(\log e)''=-1/(t-h)^2+f''/f-(f'/f)^2<0$.

For fixed $d,h$,
\[
 W=V e(\alpha)e(\beta),\qquad
 V=\frac{\lambda^2(1-d+p)}{d^2(1+d+p)}.
\]
Direct differentiation gives
\[
 \frac{V_p}{V}=
 \frac{d[(1+p)^2-d(1+p)+d^2]}
 {[(1+p)^2-d^2](1+p)(1+p-d/2)}>0.
\]
Its numerator bracket is $(1+p-d/2)^2+3d^2/4$.
When $\alpha>\beta$, strict log-concavity gives
\[
 \partial_p\log(e(\alpha)e(\beta))
 =\frac{(\log e)'(\beta)-(\log e)'(\alpha)}{\alpha-\beta}>0.
\]
At equality its limit is $-(\log e)''(d/2)>0$.
This proves $W_p>0$ over the entire stated real-endpoint interval.

\paragraph{A shared derivative identity.}
Put $k=1-d+p$, $l=1+d+p$, $m=9+3d+p$ and
$a_Q=(1-3h^2+2h^3)/12$, $b_Q=(2-3h^2-h^3)/6$,
$c_Q=(27-9h^2-2h^3)/48$.
Direct differentiation of the symmetric rational loss sum gives
\begin{equation}\label{pf:Qp-identity}
 dQ_p=\frac{a_Q(2-d)^2}{k^2}+\frac{2db_Q(d+2)}{l^2}
       -\frac{4c_Q(d+2)(d+6)}{m^2}-1.
\end{equation}
For $0\le h\le1/5$, these decreasing polynomials give
$28/375\le a_Q\le1/12$, $0<b_Q\le1/3$ and
$208/375\le c_Q\le9/16$.
On the full rational box $43/100\le d\le11/20$,
$0\le h\le1/5$, $0\le p\le121/1600$, use the increasing ratios
$(2-d)/(1-d)$ and $d(d+2)/(1+d)^2$ for the positive terms,
and endpoint bounds on the subtracted term:
\[
 dQ_p\le\frac{841}{972}+\frac{374}{961}
 -\frac{4(208/375)(243/100)(643/100)}
 {(9+33/20+121/1600)^2}-1<-\frac1{40}.
\]
No real-endpoint hypothesis is needed.

\subsection{Eliminating large total width}\label{pf:M109}
Under the same actual source hypotheses and $A\le297/500$,
one has $h<1/5$, $d<351/500$ and
$h(d-h)<p\le d^2/4$. Suppose $d\ge3/5$.
The whole closed box
\[
 3/5\le d\le351/500,\quad0\le h\le1/5,\quad
 h(d-h)\le p\le d^2/4
\]
satisfies the exact inequality
\[
 Q>351/1000-d/2-h^2.
\]
Here is its finite polynomial certificate. Parametrize the entire box by
\[
 d=3/5+(51/500)D,\quad h=H/5,\quad
 p=h(d-h)+(d/2-h)^2Z,\qquad 0\le D,H,Z\le1.
\]
This covers every point since $d/2-h\ge1/10$.
Multiply $Q-351/1000+d/2+h^2$ by its positive denominator
$3000d(1-d+p)(1+d+p)(9+3d+p)$; all factors are positive
because $0<d<1$ and $p\ge0$.
The resulting numerator has tensor degree $(8,8,4)$ and $405$
Bernstein coefficients, with exact minimum
$1280998719/5000000>0$. The displayed substitution, rational
constructor for $Q$ and coefficient formula determine the whole test.

Put $\delta=3-4\overline C=29529/500000$.
Then $\delta-(243/1000)^2=9/1000000>0$ and
$\delta+(1/5)^2<1/4$. Consequently
$\sqrt{\delta+h^2}-h^2$ increases on $0\le h\le1/5$ and
is greater than $243/1000$.
The actual shear completion therefore yields
\[
 KS>Q+\sqrt{\delta+h^2}>594/1000-d/2.
\]
Since $A=d/2+KS+(M_2-h^2S)/2>d/2+KS$, this contradicts
$A\le297/500$. Hence $d<3/5$, also under the stronger
contrary budget $A\le74/125$.

\subsection{A complete product comparison}\label{pf:M102}
Put $b_0=3-4\overline C=29529/500000$, $b=b_0+h^2$,
$K=(3+h^2)/2$, $X=47/250$ and $A_*=74/125$.
The own-height floor of Section~\ref{pf:M033} applies directly on this
same eight-source. Indeed its complete positions $1,\ldots,6$ are
$102021$, so Corollary~\ref{loc:four-calibration} supplies
$\mathcal C<H_*<1471/2000$ without a larger coupled source.
The actual return-cap gaps give $e_L,e_R>47/2500$ as in that section.
Its own-pair transport and radius contradiction apply to the same
roots $-\xi<z<r$, own shear $T_{\rm own}$ and both chronological
ceilings $x_L,x_R<X$, proving $h>1/14$ without a four-core cost floor.
Hence every actual contrary point satisfies
\begin{equation}\label{pf:d592-regional-domain}
 \frac1{14}<h<\frac15,\quad \frac{43}{100}<d<\frac35,\quad
 \max\{2h^2,h(d-h)\}<p\le\frac{d^2}{4},\quad
 p<\frac{d^2}{3}-h^2.
\end{equation}
This is the genuine four-core own mean: its two actual zero-to-zero
leaves are those used in that transport argument. The returned mean
heights and return-chart constants remain distinct.
The actual shear has the positive ceilings
\[
 KS<T_c=A_*-d/2,\qquad KS<T_s=KX\frac{2d(d+2)}{(2-d)^2},
\]
since $\partial_pE_1=X(2-3d)/(1-d+p)^2>0$ for $d<3/5$.
Write $Q_s=Q(d,h,d^2/4)$.

\paragraph{The complete product-coordinate comparison.}\label{pf:complete-product-monotonicity}
On $43/100\le d\le3/5$, $1/14\le h\le1/5$ and
$h(d-h)<p\le d^2/4$, put
\[
 F=Q+\sqrt{\delta+h^2+4W},\quad L=21/100+h/2,\quad T=-LQ_p-2W_p.
\]
The actual cap losses give $W>0$ and $W_p>0$ by
Section~\ref{pf:M111}. Retain $k,l,m$ from~\eqref{pf:Qp-identity}
and put $v=2+2p-d$. Here $k<1$, $m\ge9$ and $b_Q>0$.
The coefficient bounds depend only on $h$, so twice differentiating
\eqref{pf:Qp-identity} gives
\[
 dQ_{ppp}>2744/3125-143/4050>0.
\]

The only poles of $W$ as a rational function of $p$ are the distinct
roots of $l,m,v$, the last of multiplicity at most two; its polynomial
part has degree two. Define, after cancellation,
\[
 a_l=(lW)|_{p=-1-d},\quad a_m=(mW)|_{p=-9-3d},\quad
 b_v=(v^2W)|_{p=d/2-1},\quad
 a_v=\tfrac12(v^2W)_p|_{p=d/2-1}.
\]
Partial fractions therefore give
\[
 W=\Pi_2(p)+a_l/l+a_m/m+a_v/v+b_v/v^2,\qquad \deg\Pi_2=2.
\]
All four residues depend only on $d,h$. Put
\begin{align*}
 D_0&=600d^2(1+h)(1+h-d)^2(1+d-h)^2(3+d-h)^2
                      [2-d+2h(d-h)]^2,\\
 D_1&=900d^2(2-d)^2(2+d)^4(6+d)^4(d^2-2d+4)^3.
\end{align*}
The following polynomial constructors have positive Bernstein
coefficients after $d=43/100+17x/100$, $h=1/14+9y/70$,
with $x,y\in[0,1]$:
\[
\begin{array}{c|c|r}
 \text{constructor}&\text{degree}&\text{count}\\ \hline
 405d^2(d+4)(-a_l-8/5)&(4,6)&35\\
 9d^2(d+4)(7d+16)^2(1700-a_m)&(6,6)&49\\
 -324d^2(7d+16)^2a_v&(5,6)&42\\
 -432d(7d+16)b_v&(4,6)&35\\
 D_0T(d,h,h(d-h))&(9,12)&130\\
 D_1T(d,h,d^2/4)&(18,6)&133
\end{array}
\]
All clearing factors are positive. These $424$ exact comparisons use
\eqref{cv:bernstein-coefficients}, without subdivision.
The first four rows give $a_l<-8/5$, $a_m<1700$, $a_v,b_v<0$.
Since $l\le169/100$ and $m\ge1029/100$,
\begin{align*}
 W_{ppp}&=-6a_l/l^4-6a_m/m^4-48a_v/v^4-192b_v/v^5\\
 &>6\left[\frac{8/5}{(169/100)^4}
            -\frac{1700}{(1029/100)^4}\right]>0.
\end{align*}
Thus $T_{pp}=-LQ_{ppp}-2W_{ppp}<0$. The last two rows give its
positive endpoint values, so concavity proves $T>0$ throughout the
real-product interval. Finally
\[
 \delta+h^2-L^2=\tfrac34(h-7/50)^2+129/500000>0.
\]
Hence $-LQ_p>2W_p>0$ and
$F_p<Q_p+2W_p/L=-T/L<0$. Consequently
$F(d,h,p)\ge F(d,h,d^2/4)$, strictly off balance.

\subsection{Six balanced boundary comparisons}\label{pf:M103}
Set $\omega=(2-d)/(4(1-d))\in(0,1)$ and $B=B_R$.
The two already proved shear bounds combine exactly to give
\[
 KS<T_m:=K\{\omega E_1+(1-\omega)E_2\}
       =2KX\{\omega(1+B)-B\}.
\]
The cancellation retains the same $d,h,p$; $E_2$ itself decreases in $p$
and is not maximized at balance. Since $E_1$ increases in $p$,
the earlier ceilings $T_s,T_c$ also remain valid.
At $p=d^2/4$ put
\[
 q=Q_s,\quad \lambda_s=\frac{2(1+d^2/4)}{2+d^2/2-d},\qquad
 G_T=\delta+h^2+\lambda_s^2q^2-(T-q)^2.
\]
The balanced product is $P=q^2/4$. Write $D=d^2-2d+4$ and
\[
 Z=4500000d(d+2)^2(d+6)^2D^2,\quad
 U=60h^3+45h^2+104h-60=60(R-1).
\]
For $T=T_c,T_s,T_m$, respectively, the positive clearing factors are
$Z$, $Z(2-d)^4$, and $Z(1-d)^2U^4$.
Multiply $G_T$ by its factor and substitute each interval below
linearly from $[0,1]$. These exact constructors have positive Bernstein
coefficients at the displayed actual bidegrees:
\[
\begin{array}{c|c|c|c|r}
 T&d&h&\text{degree}&\text{count}\\\hline
 T_m &[43/100,53/100]&[1/14,1/10]&(12,18)&247\\
 T_c &[53/100,3/5]&[1/14,1/10]&(10,6)&77\\
 T_m &[43/100,261/500]&[1/10,13/100]&(12,18)&247\\
 T_c &[261/500,3/5]&[1/10,13/100]&(10,6)&77\\
 T_s &[43/100,129/250]&[13/100,1/5]&(14,6)&105\\
 T_c &[129/250,3/5]&[13/100,1/5]&(10,6)&77
\end{array}
\]
Formula~\eqref{cv:bernstein-coefficients} specifies all $830$ positive
coefficients. The factors are nonzero: $d>0$, $D=(d-1)^2+3>0$,
$d<1$, and $R<1/2$. The six closed rectangles cover the original
$(d,h)$ box, including their shared boundaries.

\subsection{Exhaustion by the complete comparison}\label{pf:M104}
At every actual contrary point the monotonicity just proved gives
\[
 KS>F(d,h,p)\ge q+\sqrt{\delta+h^2+\lambda_s^2q^2}.
\]
Select its rectangle and ceiling $KS<T$. The positive residual $G_T$
gives $T-q<\sqrt{\delta+h^2+\lambda_s^2q^2}$, irrespective of the
sign of $T-q$, a contradiction. Thus the original complete-source
bound $A>74/125$ follows with all cap, chronological-return,
common-lift and calibration hypotheses unchanged.

\space
\section{Regeneration and further central geometry}\label{pf:central-applications}
\subsection{The direct-centre packing}\label{pf:M129}\label{pf:F27}
The complete-source central estimates are used with their retained
contacts, terminal arcs and selected full lift.
\space
\subsection{Regeneration from a complete nine-source}\label{pf:M130}\label{pf:F29}
\begin{lemma}[Two singleton bounds from a nine-arc source]\label{pf:nine-singletons}

For any \(0<\alpha\le13/20\), an actual proper complete \texttt{202102101} has both \(c_3>\alpha\) and \(c_5>\alpha\), where the indices refer to that whole source. These are two individual statements, not alternatives of a maximum disjunction.
\end{lemma}
\begin{proof}
Choose a genuine full minimizing interval containing all nine complete
arcs. Proposition~\ref{fd:normality} gives normality of the selected
full lift before joint normalization.
The initial $202$ has second root $U>0$ by Lemma~\ref{sw:ABA}; hence
$1-U<1<R(s)$ for $s=bd>0$, where $b$ and $d$ are its own
positive rate and full tail duration. Its rooted response
\eqref{pc:Psi} is strictly positive by Lemma~\ref{pc:endpoint}.
The complete suffix at indices $1,\ldots,8$ is $02102101$ on this same
lift. Lemma~\ref{nx:rooted-continuations} therefore gives
$c_3>13/20\ge\alpha$: its charged local position $2$ is original arc $3$.
Full time reversal with $j\mapsto-j$ preserves the entire nine-word
and sends index $3$ to $5$. Lawful transport gives $c_5>13/20\ge\alpha$.
These are separate singleton conclusions, with no rate or costate
identified across the two orientations.
\end{proof}

Both singleton bounds can be used on their distinct physical charges;
they are not a maximum alternative.
\space
\subsection{An opposite centre with both sign contexts}
\label{pf:M132}\label{pf:F31}\label{pf:M133}
\begin{theorem}\label{pf:opposite126}
On one genuine containing interval of an original minimizer, retain the
complete positive source $1|020|1|202|1$, its two complete outward
forward-force contexts, and one full maximizing lift. Suppose these
contexts give the negative entrances of both outward capped triples in
their canonical charts. Then the central physical arc has cost
$c>63/500$. No HIGH condition on either actual mean is required.
\end{theorem}
\begin{proof}
Apply Theorem~\ref{fee:opposite126}. The two outward contexts, their
negative entrances and their matrix means belong to the selected
full lift specified here, as required there.
\end{proof}
\space
\space
\section{Detailed proofs for the terminal and finite-domain exclusions}
\label{nx:supplement}

\subsection{Rooted terminal obstructions and continuations}
\label{nx:rooted-proofs}

At an actual complete $0\to2$ root and its positive mode-$2$ tail, let
$U,V$ be the horizontal endpoints, $b=v/(1-u)>0$, and $d>0$ its full duration.
The full binary equations give
\begin{equation}\label{nx:root-response}
 \begin{split}
 \Psi(U,V)&=2(1-U)\log\frac{1-U}{1-V}-(V-U)(2+V-U),\\
 D_0&=(1-e^{-bd})y'(d)-by(d),\qquad y=-\chi,\\
 bD_0&=2(bd-1+e^{-bd})-(1-U)(1-e^{-bd})^2,\\
 \Psi&=b(1-U)D_0.
 \end{split}
\end{equation}
These identities retain the actual initial derivative, whose homogeneous
contribution cancels in $D_0$.

\paragraph{Proof of Lemma~\ref{nx:02101}.}
\begin{proof}
If its response were positive, the endpoint reconstruction and
$1<R(bd)<e^{bd}$ in Section~\ref{pc:tail} would give the strict ingress
\eqref{hc:ingress} at its actual $2\to1$ tie. The complete terminal $101$
supplies the opposite reflected-tail inequality on this same lift.
Here $q=\sigma_1=-\sigma_0>0$, and both rates and full durations are
positive, so the two-tail contradiction \eqref{sw:two-tail} applies.
The conclusion is weak: equality of the rooted response is retained.
\end{proof}
\paragraph{Proof of Lemma~\ref{nx:rooted-continuations}.}
\begin{proof}
Equation~\eqref{nx:root-response} supplies the strict actual ingress.
The two actual central contacts and full maximizing signs, in the cyclic
chart making the centre mode $0$, give
\[
 V(A)=0,\quad V_\xi(A)\le0,\qquad
 U(AE)=0,\quad U_\xi(AE)\le0,\qquad q_R\ge0.
\]
The first source has the complete terminal $020$ and therefore satisfies
all hypotheses of Theorem~\ref{nx:central1196}.

For the second, suppose its centre costs at most $13/20$. The prefix exit
estimate gives $u_R<0$, $q_R\ge0$, and $q_R'>kq_R$.
Lemma~\ref{hc:regeneration} propagates these facts along the following
complete mode-$0$ and mode-$2$ arcs, on the same lift, to a fresh strict
ingress. This propagation uses no further cost hypothesis. Its actual
terminal $101$ now gives exactly the two-tail contradiction in the proof
of Lemma~\ref{nx:02101}. Only the $13/20$ prefix theorem was iterated.
\end{proof}

\subsection{The complete central cost $299/250$}
\label{nx:central-proof}

We retain all actual-source hypotheses of Theorem~\ref{nx:central1196}.
The symbols $A,E,\tau,z=\log E=2h$, $X_->0,Y_+\ge0$, $f_-,f_+$ and $Q$
have the meanings in Section~\ref{hc:section}. In particular the complete
four-constant full normal family, including the state covector, yields
\eqref{hc:response}; the actual contact derivatives are
$V_\xi(A)\le0$ and $U_\xi(AE)\le0$. The open-domain conditions are
$AE<2$, $A\tau^2>2$. The strict ingress yields $f_-<0$ and the actual
crossing interval
\[
 \frac2{\tau E}<A<\frac2\tau.
\]
No upper bound on $\tau$ is needed in this argument.

\paragraph{The whole low-factor range.}
If $1<E\le31/10$, Lemma~\ref{hc:short-factor-force} gives $Q>0$
on this actual crossing, with no restriction on $\tau$ or the physical
cost. Hence $f_+=Q-f_->0$. The exact cancellation
\eqref{hc:cancellation} contradicts the negative actual terminal
expression~\eqref{hc:terminal-scaled}. Thus the complete source
necessarily has $E>31/10$, independently of its cost. The upper-factor
comparison below covers the full closed interval $[31/10,13/2]$.

\paragraph{A fresh factor cap from the physical cost.}
The actual contact floor and physical-cost comparison
\eqref{hc:A-floor}--\eqref{hc:cost-floor} give
$A>A_t>0$ and $C_{\mathrm{loc}}>W(E)$, with $W$ strictly increasing.
Use the same positive-series logarithmic bounds $L_N,U_N$ of
\eqref{hc:log-bounds}.
At $E=13/2$ direct rational evaluation gives
\[
 \frac{E L_6(E)}{E-1}-1-\frac{299}{250}
 =\frac{134750995111}{13455175781250}>0.
\]
Therefore the contrary physical budget $C_{\mathrm{loc}}\le299/250$
implies $31/10<E<13/2$.

\paragraph{Both contacts on the high-factor range.}
Both actual contact jets give the exit floor \eqref{hc:exit-floor}.
If its bracket were nonpositive, the complete elimination
\eqref{hc:eliminated-cost}, including $D_b>0$ on $[A_t,A]$, would give
$C_{\mathrm{loc}}>C_b(A_t)$.
Define the polynomial
\begin{align}\label{nx:P1196}
 P(E,z)&=(E-1)^2D_b(A_t)\left(C_b(A_t)-\frac{299}{250}\right)\\
 &=\frac1{125}\bigl(750E^3z-1098E^3-250E^2z^3
       +549E^2z^2-2049E^2z+3294E^2\nonumber\\
 &\hspace{28mm}-250Ez^2+1549Ez-3294E-250z+1098\bigr).
 \nonumber
\end{align}
The appearance of $A_t$ is an algebraic elimination, not the assignment
of an actual cost or lift to a comparison trajectory.

\paragraph{Two elementary polynomial comparisons.}
The polynomial $P_z$ is concave in $z$, since $P_{zzz}=-12E^2$.
At $z=0$ and $E=3+s$ its descending coefficients are
$6,4701/125,1901/25,6206/125$, all positive. The upper logarithmic
bound needs only its first term:
\[
 P_z(E,U_1(E))=\frac{E-1}{750(E+1)^2}
 \bigl[E^4(3223-125E)+2534E^3+E(7240E-6821)+1125\bigr]>0.
\]
Indeed $E\le13/2$ gives $3223-125E\ge4821/2$, and $E>1$
gives $7240E-6821>419$. Concavity therefore proves $P_z>0$
throughout $0\le z\le U_1(E)$.
For the lower logarithmic bound put $s=E-3\ge1/10$. Exact expansion gives
\begin{equation}\label{nx:bernstein-rule}
\begin{split}
 P(E,L_2(E))&=\frac{2(E-1)^3}{3375(E+1)^9}\,p(s),\\
 p(s)={}&12177s^9+353312s^8+4545624s^7+33802536s^6\\
 &+158618396s^5+480189480s^4+914569008s^3\\
 &+1000579904s^2+480280320s-1285632>0.
\end{split}
\end{equation}
The last two terms already have sum at least $46742400>0$;
all other coefficients are positive. These displayed identities certify
the entire closed factor interval, without subdivision or numerical logs.
Thus $P(E,\log E)>P(E,L_2(E))>0$, so $C_b(A_t)>299/250$,
contrary to the actual cost budget. The bracket in \eqref{hc:exit-floor}
is positive; its positive factor $2-AE$ contradicts the strict negative
terminal expression. This proves Theorem~\ref{nx:central1196}.
The low and high arguments overlap at $E=31/10$, and the latter includes
$E=13/2$.
Lemma~\ref{hc:terminal} already proves the raw $C^1$ matching and the
variable terminal rate \eqref{hc:variable-rate}; only its actual endpoint
value $c(R)=k$ enters the retained cancellation.

\subsection{The weighted bilateral four-core estimate}
\label{nx:weighted-proof}

Here all binary quantities are obtained from one selected normalized
original lift. Let $A=A_1+A_2+A_3+A_4$, $M=A_2+A_3$, and
$W=A_1/2+A_2+A_3+A_4/2$, so $2W=A+M$.
Write the outer cap and three root coordinates as
\[
 -\alpha<-\xi<z<\zeta<\beta,\qquad
 \alpha>\xi>h>0,\quad \beta>\zeta>h,\quad |z|<h,
 \qquad h^2=2K-3.
\]
In shear time $ds=v^{-2}\,du$ the complete normal equations are
\[
 L_s=K,\qquad m_s=-2L,\qquad \chi_s=h^2-u^2,
 \qquad C_{\rm mat}=-L^2-Km.
\]
With $Q=v^{-2}$ and $e=1+u$ on mode $0$, $e=u-1$ on mode $2$, put
\[
 f=\frac{h^2-u^2}{2}Q+\frac{eQ+e^{-1}}2\chi,
 \qquad g=m+\frac12-f\ge0.
\]
The actual exterior mode-$1$ contacts satisfy $g=0$.

On the left complete leaf set $p=1-\alpha$, $q=1-\xi$,
$E=q/p>1$, $\lambda=(v/(1+u))^{-2}$, and
\[
 S_1=\frac{\lambda(E-1)}q,\qquad
 I=q-p-2\log E+(1-h^2)(p^{-1}-q^{-1}).
\]
The actual root equation gives $\chi(-\alpha)=\lambda I$. Therefore
\[
 f_a=\frac{\lambda\{h^2-(p-1)^2\}}{2p^2}
       +\frac{\lambda(1+\lambda)I}{2p},\qquad
 f_1=\frac{\lambda(h^2-\xi^2)}{2q^2}<0.
\]
Define the strict response margin
\[
 \Delta=(1+\lambda)q(E\log E-E+1)-2\lambda(E-1)^2>0.
\]
The actual response predicate is equivalently
$(1+1/\lambda)pR(\log E)>4$, where
$R(s)=2(s-1+e^{-s})/(1-e^{-s})^2$.
The inequality $R(\log E)<E$ follows from
$E-E^{-1}-2\log E>0$, whose derivative is $(E-1)^2/E^2$.
Consequently
\[
 \lambda<\frac q{4-q}<1.
\]
Expanding the cap expression and propagating the actual matrix give the
two equal expressions
\begin{align*}
 f_a+KS_1^2-f_1
 &=-\xi S_1-\frac{S_1(1-\lambda)(\xi^2-h^2)}{2q}
                  -\frac{\lambda\Delta}{q^2}\\
 &=g_1+2L_1S_1.
\end{align*}
Thus
\begin{equation}\label{nx:weighted-left}
 L_1=-\frac\xi2-
 \frac{(1-\lambda)(\xi^2-h^2)}{4q}
 -\frac{g_1}{2S_1}-\frac{\Delta}{2q(E-1)}.
\end{equation}
At the right leaf use $q_R=1-\zeta$, $p_R=1-\beta$,
$E_R=q_R/p_R$, $\lambda_R=(v/(1-u))^{-2}$ and its strict
$\Delta_R>0$. The sign of the endpoint switching value is reversed;
substitution in $f$ and forward matrix propagation yield
\begin{equation}\label{nx:weighted-right}
 L_3=\frac\zeta2+
 \frac{(1-\lambda_R)(\zeta^2-h^2)}{4q_R}
 +\frac{g_3}{2S_4}+\frac{\Delta_R}{2q_R(E_R-1)}.
\end{equation}
Both equations concern the same original matrix; symmetry of the
physical trajectory has not been assumed.

Set $a=-L_1=\xi/2+\epsilon_L$, $b=L_3=\zeta/2+\epsilon_R$,
$r=\xi+\zeta$ and $P=a+b=KS_{\rm mid}$.
Equations~\eqref{nx:weighted-left}--\eqref{nx:weighted-right} imply
\[
 \epsilon_L,\epsilon_R>0,\qquad
 g_1<2S_1\epsilon_L<2S_1a,\quad
 g_3<2S_4\epsilon_R<2S_4b.
\]
Each complete middle arc has exact mean $u^2=h^2$; integration of its
physical cost therefore gives
\begin{equation}\label{nx:middle-excess}
 M=\frac r2+P,\qquad M-r=\epsilon_L+\epsilon_R>0.
\end{equation}

Suppose $W\le w:=61/200$. The original strict unweighted four-core bound
$A>389/1000$ gives
\[
 M=2W-A<\frac{221}{1000}<\frac29,\qquad
 r<M,\qquad A<2W\le\frac{61}{100}<1.
\]
Corollary~\ref{loc:four-calibration} on this actual capped-four
lift gives $0<C_{\rm mat}<H_*<92/125$.
For $\delta=K/2-C_{\rm mat}$ it follows that
\begin{equation}\label{nx:delta-lower}
 \delta>0,\qquad 4\delta=2K-4C_{\rm mat}>\frac7{125}.
\end{equation}

Let $D=v(z)^{-2}$ be the single actual middle height. Exact leaf matching
gives
\[
 S_{\rm mid}=DT_{\rm mid},\qquad
 T_{\rm mid}=
 \frac{(1-z)(\xi+z)}{1+\xi}
 +\frac{(1+z)(\zeta-z)}{1+\zeta}<r,
\]
since
$r-T_{\rm mid}=(\xi+z)^2/(1+\xi)+(\zeta-z)^2/(1+\zeta)>0$.
The exact mean equations give
$\zeta(h,z)=\xi(h,-z)$ and
\[
 \xi_z=\frac{(h^2-z^2)(1+\xi)^2}
                {(1-z)^2(\xi^2-h^2)}>0.
\]
If $z\ge0$, then $\zeta\le\xi$ and $\zeta<1/9$.
The right strict response and its own leaf height give
\[
 D<\frac{(1+\zeta)^2}
           {(1+z)^2(1-\zeta)(3+\zeta)}
 \le\phi(\zeta)<\phi(1/9)=\frac{25}{56},
 \quad \phi(x)=\frac{(1+x)^2}{(1-x)(3+x)}.
\]
Here $\phi$ is strictly increasing on $(0,1)$.
For $z\le0$ use the left side and $\xi\le\zeta$ instead.
Thus $h<r/2<1/9$ and $K=(3+h^2)/2<8/5$.
It follows that
\[
 P=KDT_{\rm mid}<\frac57r<\frac34r,
 \qquad P<\frac35M,
\]
where the second inequality uses \eqref{nx:middle-excess}.

At the two outer binary roots the same determinant identity gives
\[
 \delta=a^2+Kf_1+Kg_1<a^2+2KS_1a,
 \qquad
 \delta=b^2+Kf_3+Kg_3<b^2+2KS_4b.
\]
Divide by positive $a,b$ and add. Since $\delta>0$ and
$1/a+1/b\ge4/(a+b)=4/P$, we get
\[
 4\delta<P^2+2KS_{\rm out}P.
\]
The original open-domain cost density is strictly larger than twice
shear density; hence
$S_{\rm out}<(A-M)/2=W-M\le w-M$.
Therefore
\begin{align*}
 4\delta
 &<\frac9{25}M^2+\frac{48}{25}M(w-M)\\
 &\le\frac{192w^2}{325}=\frac{11163}{203125}.
\end{align*}
The last inequality is exactly
\[
 \frac{192w^2}{325}
 -\left\{\frac9{25}M^2+\frac{48}{25}M(w-M)\right\}
 =\frac{39}{25}\left(M-\frac{8w}{13}\right)^2\ge0.
\]
It contradicts \eqref{nx:delta-lower}, proving
Theorem~\ref{nx:weighted305}. All strict terms came from positive complete
leaves, the strict responses or the original open domain; no equality of
left and right states was used.

\subsection{The forward-pair domain and determinant certificate}
\label{nx:forward-proof}

Put $A_*=180483/500000$ in this subsection; this scalar is not the cost
matrix in the full PMP equations. We prove
Theorem~\ref{nx:forward-pair} by assuming that the two central complete
costs have sum $I\le A_*$.
Retain actual complete arcs $0,1,2,0,1,2$ and both exterior caps.
The central $1201$ context proves normality of the selected full lift.
At its actual middle switch let $x,z>y>0$ be the reciprocal variables,
$xy+xz+yz=1$. Whole reflection and time reversal may arrange
\[
 R=x/y\ge t=z/y>1,\qquad S=R+t,\qquad D=S+Rt,
 \qquad (x,y,z)=\frac{(R,1,t)}{\sqrt D}.
\]
Let $E,F>1$ be the exact factors of the two charged central arcs, whose
durations are $z\log E$ and $x\log F$. At the switch write
$vp=(h_0,0)$ and let the full response vector be $(a_0,b_0)$.
The actual ray and crossing conditions give $h_0\ge0$, $b_0\le0$.
The full matrix, not only the response vector, is reconstructed by
\[
 u_0=\frac{x-z}{x+z},\quad v_0=\frac2{x+z},\quad
 K=\frac{3+u_0^2-b_0v_0}{2},\quad L=\frac{a_0+u_0}{2},\quad
 m=-\frac12+\frac{h_0-b_0/2}{v_0},
\]
with $l=L+u_0m$ and $\beta=K-2u_0l+u_0^2m$.
Both state covectors are retained: $p_u=h_0/v_0$, $p_v=0$ at the
middle, and thereafter $p_v=\chi/v$ with $p_u$ supplied by the active
zero Hamiltonian.

\paragraph{The genuine exterior-cap inequalities.}
At the right end of the central pair the reciprocal triple is $(x,X,Y)$,
where
\[
 X=F(x+y)-x,\qquad Y=(x+z)/F-x.
\]
The actual neighboring forward context, in its whole cyclic chart,
gives $0<Y<\min(X,x)$. On the next complete positive arc let $d>0$,
$G=e^{d/X}>1$ and let $\lambda_R$ denote its transverse response in
that chart. Its genuine exterior cap equation and the positive-defect
identity imply
\[
 \lambda_R\le \frac{c_X(G)}{A_X(G)}
     =2\gamma_X<\frac{X-Y}{X+Y}\,d.
\]
The actual next endpoint remains in the open state domain, so
$G<(X+x)/X$ and
$d<X\log(1+x/X)<x$. Hence
\begin{equation}\label{nx:exterior-forward}
 \lambda_R<x\frac{X-Y}{X+Y}.
\end{equation}
No supporting-arc cost is assigned to the pair budget. The left analogue
is obtained by full reflection and time reversal of the same lift.

The following explicit formulae make the domain and the arithmetic
certificate unambiguous. Put
\[
 f(g)=\frac{g\log g-g+1}{g-1},\qquad
 k(g)=\frac{f(g)}{g-1}
     =\int_0^1\frac{1-s}{1+(g-1)s}\,ds.
\]
Scale $h_0=u/\sqrt D$, $b_0=v/\sqrt D$, and define
\begin{align*}
 r_0&=\frac{2D}{S^2}\left(\frac{tF}{F-1}+\frac{RE}{E-1}\right),\\
 q_0&=\frac{2Rt}{S^2}\left\{
 \frac{(D+R^2)f(F)}F+\frac{(D+t^2)f(E)}E
   -S\left(\frac{F-1}F+\frac{E-1}E\right)\right\},\\
 \beta_x&=2R\left(f(F)-\frac{R(F-1)}S\right),
 \qquad
 \beta_z=2t\left(f(E)-\frac{t(E-1)}S\right),\\
 A_x&=\frac{(F-1)(D+R^2)}{2FD},\qquad
 c_x=\frac{R(F-1)\{(D+R^2)f(F)-S(F-1)\}}{F^2D},\\
 A_z&=\frac{(E-1)(D+t^2)}{2ED},\qquad
 c_z=\frac{t(E-1)\{(D+t^2)f(E)-S(E-1)\}}{E^2D},\\
 v(u)&=r_0u-q_0,\qquad
 a(u)=\frac{(R^2-D)v(u)}{2RD}
               -\frac{F(c_x-u)}{R(F-1)}.
\end{align*}
These are the two actual central cap equations solved for the full
response. Direct integration of the full response drifts on the right
arc gives
\begin{align*}
 a_R(u)&=a(u)+\frac{2R(R^2-D)(F-1)}{DS}
           -\frac{2R^2\log F}{D}+\frac{2RS(1-F^{-1})}{D},\\
 v_R(u)&=v(u)+\frac{4R^2(F-1)}S-2R\log F.
\end{align*}
With
\[
 X_R=F(R+1)-R,\qquad Y_R=S/F-R,
\]
whole cyclic covariance converts \eqref{nx:exterior-forward} to
\[
 L_R(u):=\frac{2DY_Ra_R(u)+(Y_R^2-D)v_R(u)}{D+Y_R^2}
 <B_R:=\frac{R(X_R-Y_R)}{X_R+Y_R}.
\]
The left functions are obtained by the simultaneous exchange
$(R,t,E,F)\mapsto(t,R,F,E)$. This exchange includes the sign change
of the horizontal response under the full reflection.
Their slopes are strictly negative:
\[
 L_R'(u)=-\kappa_R,\qquad
 \kappa_R=
 \frac{2D\{(E-1)F(F+t)+ER(F-1)+EF\}}
 {(E-1)S\{R(F-1)^2+F^2+t\}}>0.
\]
Every denominator and summand is positive on the actual domain.

\paragraph{Exact height elimination.}
The central rays and signs give
\[
 0\le u\le U:=\frac{\min(q_0,q_0+\beta_x,q_0+\beta_z)}{r_0}.
\]
The actual adjacent first jets give the weak lower bounds
\begin{align*}
 u&\ge j_x:=c_x+A_x\{2R\log F-4R^2(1-F^{-1})/(R+1)\},\\
 u&\ge j_z:=c_z+A_z\{2t\log E-4t^2(1-E^{-1})/(t+1)\}.
\end{align*}
Set
\[
 \ell_R=\frac{L_R(0)-B_R}{\kappa_R},\quad
 \ell_L=\frac{L_L(0)-B_L}{\kappa_L},\quad
 L=\max(0,j_x,j_z,\ell_R,\ell_L).
\]
Every actual height lies in the closed comparator interval $[L,U]$.
Feasibility of the retained affine conditions requires
\[
 U\ge\max(0,j_x,j_z),\qquad L_R(U)<B_R,\qquad L_L(U)<B_L.
\]
Because both exterior responses strictly decrease, these conditions are
also sufficient for this \emph{affine relaxation}; they do not assert
physical existence. A singleton interval is retained unless its active
strict exterior inequality excludes it. No division by $q_0$ is used;
$q_0=0$ remains included.

The exact central cost, in a positive form, is
\begin{align}\label{nx:forward-cost}
 I={}&\frac{t(E-1)}D\left\{1+\frac SE+
                  \frac{t(E-1)}E(1-k(E))\right\}\\
 &+\frac{R(F-1)}D\left\{1+\frac SF+
                  \frac{R(F-1)}F(1-k(F))\right\}.
 \nonumber
\end{align}
The remaining geometric conditions are
\begin{equation}\label{nx:forward-geometry}
 \begin{gathered}
 R\ge t>1,\quad E,F>1,\\
 0<Y_R<\min(X_R,R),\qquad
 0<Y_L<\min(X_L,t),\\
 X_L=E(t+1)-t,\qquad Y_L=S/E-t.
 \end{gathered}
\end{equation}
The full determinant, reconstructed from the above matrix entries, is
\begin{equation}\label{nx:forward-det}
 \mathcal C(u)=\frac34
 -\frac{n_xa(u)+n_yv(u)/D}{2}
 -\frac{(\ell_0-n_y)u}{2D}
 -\frac{a(u)^2+v(u)^2/D}{4}+\frac{uv(u)}{2D},
\end{equation}
where
$n_x=(R-t)/S$, $n_y=1-(R^2+t^2)/S$, $\ell_0=S+1$.
Its quadratic coefficient is
$\{2r_0-r_0^2-D(a')^2\}/(4D)<0$, since $r_0>2$.
The adverse final term in \eqref{nx:forward-det} is retained; it need
not vanish at either clipped endpoint.
Lemma~\ref{cv:affine-fibre} therefore gives
$\mathcal C_{\rm actual}\ge\min\{\mathcal C(L),\mathcal C(U)\}$.

\paragraph{Analytic removal of the only singular corner.}
From the exact cost and \eqref{nx:forward-geometry}, using
$\log g<g-1$ and $1-\sqrt a>(1-a)/2$ for $0<a<1$, one obtains
\[
 I>\frac{2Rt-R-t}{2D},\qquad
 I>\frac{(2t+1)(R-t)}{2D}.
\]
For completeness the underlying strict factor bounds are
$E>S/(2t)$, $E^2>S/(t+1)$, $F^2>S/(R+1)$,
$E<S/t$, $F<S/R$; substitution in the positive cost formula gives the
two displayed estimates. They imply
\[
 t<T:=\frac{1+2A_*}{1-A_*}=\frac{860966}{319517}<3,
\]
\[
 \frac Rt<\phi(t):=
 \frac{2t+1+2A_*}{2(1-A_*)t+1-2A_*}
 \le\phi(1)=\frac{103387}{43226}=:P,
\]
since $\phi'<0$. The exact product
$TP=44506345921/6905720921<7$ gives $R<7$, and
$E<1+P<17/5$, $F<2$.

Put $p_g=(g-1)/g$ and $d_g=(g-1)/2-f(g)$.
The positive integral identities give
\[
 \frac{(g-1)^2}{6g}\le d_g\le\frac{(g-1)^2}{6},\qquad
 q_0=\frac{2Rt}{S}\left\{
 \frac{R-1}{2}p_F-\frac{R+1}{F}d_F
 +\frac{t-1}{2}p_E-\frac{t+1}{E}d_E\right\}.
\]
If $F\to1$, the geometry forces $t\to1$; on that face $q_0<0$ unless
$E=1$. If $E\to1$, the geometry forces $R=t=1$, and again $q_0<0$
unless $F=1$. Since actual feasibility requires $q_0\ge0$, only the
all-one corner can be a feasible factor-singular limit.

Suppose $1<t\le R\le1+\eta$ with $\eta=1/100$.
The same identity and $d_g/g\ge p_g^2/6$ imply
\[
 \frac{p_E^2+p_F^2}{3}\le\frac\eta2(p_E+p_F),\quad
 p_E+p_F\le3\eta,\quad
 0\le q_0\le\frac32(1+\eta)\eta^2<\frac1{5000}.
\]
Consequently $0\le u<1/10000$ and $-1/5000<v(u)\le0$.
In the exact expression for $a(u)$, its first term has magnitude
at most $q_0/2$, its last $Fu/[R(F-1)]$ is at most
$q_0S^2/(2DRt)\le q_0$, and
\begin{align*}
 \left|\frac{Fc_x}{R(F-1)}\right|
 &=\frac SD\left|\frac{(R+1)f(F)-(F-1)}F\right|\\
 &\le\frac SD\left\{\frac\eta2p_F+
                    \frac{(R+1)(F-1)p_F}{6}\right\}\\
 &<\frac3{20000}+\frac1{2000}.
\end{align*}
Here $F-1\le3/97<1/30$. Thus $|a(u)|<1/1000$.
Using $|n_x|,|n_y|\le1$, $\ell_0-n_y<4$, $D\ge3$ in the full
\eqref{nx:forward-det} gives
\begin{align*}
 \mathcal C(u)>{}&\frac34
 -\frac{1/1000+(1/5000)/3}{2}
 -\frac{4(1/10000)}6\\
 &-\frac{(1/1000)^2+(1/5000)^2/3}{4}
 -\frac{(1/10000)(1/5000)}6
 >\frac{749}{1000}>H_*.
\end{align*}
This bound applies to every central-feasible height, including zero.

On the complement $R\ge101/100$, the factor inequalities give
$E-1>1/1760>1/2000$. If $F-1\le10^{-8}$, then
$t-1<24(F-1)$ and the $q_0$ identity gives $p_E^2\le45(F-1)$.
But $E>S/(2t)$ gives
$p_E>(R-t)/S>(1/100-24\cdot10^{-8})/10>9/10000$, a contradiction.
Thus the whole remaining physical domain is contained in the closed box
\begin{equation}\label{nx:forward-box}
 [101/100,7]\times[1,3]\times[2001/2000,17/5]
 \times[100000001/100000000,2]
\end{equation}
in $(R,t,E,F)$, with all the necessary conditions above retained.

\paragraph{The complete finite inequality.}
The certificate proves on this closed box, subject to
\eqref{nx:forward-geometry}, \eqref{nx:forward-cost} with $I\le A_*$,
and the exact height feasibility conditions, that
\begin{equation}\label{nx:forward-finite}
 \mathcal C(L)>H_*,\qquad \mathcal C(U)>H_*.
\end{equation}
A closed binary tree consists of $129453$ nodes, $64726$ splits and
$64727$ terminal leaves. Each split is at the integer midpoint of the
longest root-relative dyadic coordinate, with lowest-index tie breaking.
The saved leaves specify their path, exact dyadic box, terminal predicate
and, for positive determinant leaves, a rational reserve.
The disjoint list of dispositions is: $36041$ determinant-positive,
$18430$ cost-excluded, $9951$ geometry-excluded, $187$ jet-excluded,
$117$ exterior-excluded, and $1$ height-excluded. There are no unresolved
leaves. The smallest recorded determinant reserve is
\[
 \frac{1178565330864301390370460635519}
 {170141183460469231731687303715884105728}>0.
\]
This is an interval reserve for determinants, not a uniform reserve for
physical costs.

Every inequality is evaluated with exact integer outward rounding.
The producer uses $128$-bit dyadics and $48$ Cayley terms; the independent
checker uses $192$-bit dyadics and exact rational sums of $96$ logarithm
terms. The exact logarithm enclosure is $L_N(g)<\log g<U_N(g)$
from \eqref{hc:log-bounds} for $g>1$.
The independent implementation obtains $k(g)$ from
$(g\log g-g+1)/(g-1)^2$ with rational arithmetic, independently of the
producer's cancellation-free positive $k$ series.
All denominator intervals must avoid zero. Pointwise minima and maxima
are enclosed without choosing a branch. Weak cost, height, ray and jet
faces are rejected only by strict separation. Strict geometry and
exterior conditions can be rejected at equality. The complete root tree
is reconstructed, both children of every split are required, and every
saved leaf box must agree exactly with its reconstructed box. Thus the
arithmetic and the coverage together establish
\eqref{nx:forward-finite}; a list of successful evaluations alone would
not suffice.

Finally concavity proves $C_{\rm actual}>H_*>0$ on the original full
lift. The original ratio-free height estimate and the physical duration
bound give
\[
 h_0\le\ell_{\rm mid}T^2/8,
 \qquad T\le A_*/\sqrt3,
 \qquad \operatorname{tr}(\Lambda X)
 =\ell_{\rm mid}-h_0
 \ge\ell_{\rm mid}(1-A_*^2/24)>0.
\]
The original harmonic theorem, with its original endpoints and real
winding, now gives $C_{\rm actual}<H_*$, a contradiction. The proof
uses harmonic calibration only at this final actual-lift step.

\space
\section{One-sided leaves and cap payments}
\label{late:r18:section}
Write $c_i$ for the complete physical cost of arc $i$. Indices are unfolded
before any change of chart: for a word of length $B$ the convention is
$j_{i+B}=j_i+1$ modulo three. All containing intervals below are genuine
fundamental intervals of the original minimizing trajectory. In particular,
a matrix calculation uses one full maximizing lift, normalized jointly
only after normality has been established.

\subsection{Replacing a missing return by a one-sided leaf}
\begin{lemma}[The one-sided transfer construction]\label{late:r18:prefix}
Retain the complete actual source $20210212101$, indexed from $-1$
to $9$, on a genuine containing interval of an original minimizer.
Use one full lift, with the capped four-core on $3,\ldots,8$ supplying
normality and calibration. If $c_2\le13/20$, the actual own pair on
arcs $5,6$ satisfies $c_5+c_6>283/1000$.

The two-circle returned-tail payment $>2/5$ also holds when the terminal
circle is supplied by the exit of an actual one-sided prefix rather
than by a final return. This assertion retains both complete caps,
the actual entrance circle, and the loss and harmonic qualifications
on one genuine containing interval.
\end{lemma}
\begin{proof}
The source on $[-1,9]$ is $20210212101$. Choose its full lift on one
containing genuine period. The actual capped four-core of $E$ gives
normality by Corollary~\ref{loc:four-calibration}. Normalize this same lift
jointly. That lemma then supplies positive determinant, the future sheet,
and the harmonic comparison, before any lower bound for $\delta$ is used.
Its prefix $[-1,4]$ is
$202102$, with central arc 2. Suppose $c_2\le13/20$.
The prefix estimate \ref{hc:prefix650}, together with the raw derivative
identity at the actual tie, gives at the beginning of the following
mode-zero arc
\[
 z<0,\qquad \chi(z)\ge0,\qquad \chi_t(z)>0.
\]
Indeed the raw switching function is $w\chi$, where $w>0$, and its
logarithmic derivative at the junction is the coefficient appearing in
the strict prefix inequality. The binary invariant $K$ on arcs 3 and 4
satisfies $\chi_s=h^2-u^2$, with $h^2=2K-3$. Hence $-h<z<0$.
If $r$ is the actual tie at the end of arc 3, then $\chi(r)=0$.
The switching function initially increases and cannot return to zero
before $u>h$. Thus $r>h>0$, and the original state domain gives $h<1$.
The exact one-sided integral is
\begin{equation}\label{late:r18:onesided}
 \int_z^r\frac{u^2-h^2}{(1+u)^2}\,du\ge0.
\end{equation}
The initial value $\chi(z)$ has not been set to zero.

At the actual preceding mode-one cap tie, the inactive identity gives
\[
\begin{gathered}
 f(z)=F(z)+A_0(z)\chi(z),\\
 F(z)=\frac{h^2-z^2}{2}q(z)>0,\qquad A_0(z)=\frac{(1+z)q(z)+(1+z)^{-1}}2>0.
\end{gathered}
\]
The determinant identity is $L(z)^2+Kf(z)=K/2-\Gamma$. Consequently
\begin{equation}\label{late:r18:capgap}
 L(z)^2+KF(z)\le K/2-\Gamma.
\end{equation}
This is an inequality at the actual cap, not a midpoint identity at a
fictitious switch. Since $(1+u)^{-2}$ strictly decreases, comparison
on the two sides of the single sign change at $h$ in
\eqref{late:r18:onesided} gives
\[
 (1+h)^{-2}\int_z^r(u^2-h^2)\,du
 >\int_z^r\frac{u^2-h^2}{(1+u)^2}\,du\ge0.
\]
Thus $r^2+rz+z^2>3h^2$. Together with $h^2>z^2$ this implies
$r^2+rz-2z^2>0$, hence $z>-r/2$.

Arc 4 is the complete shared mode-two arc from the actual binary root
$r$ to the actual third-mode tie. After cyclic change of chart it is
the incoming exterior arc of the E four-core. Let $y>0$ be its density,
$\lambda_{\rm in}$ its integration constant, and
$e_{\rm in}=\lambda_{\rm in}+3y/2$. The shared transport identities
\ref{pf:M002}--\ref{pf:M003} use the actual zero at $r$, the complete
shared arc, its terminal cap equation, and the matrix conjugacy to E;
none uses $\chi(z)=0$. They therefore give the same nonnegative cap
identity and $e_{\rm in}\ge0$. The left gap used in the density
comparison is supplied by \eqref{late:r18:capgap}, and the left shear
is the actual shear from $z$ to $r$:
\[
 L(z)=\lambda_{\rm in}+Kd_L,\qquad
 d_L=\frac{y}{1-r}-T_+
 =\frac{y(1-2r-r^2+2z)}{(1-r)^2(1+z)}.
\]
The right gap and the right root quadratic concern E's genuine own
middle switch and are unchanged. The high-density argument of
\ref{pf:M004}--\ref{pf:M005} uses precisely $z>-r/2$ on the left.
It implies $\lambda_{\rm in}<0$ in its high-density range and yields,
with $\delta=3/4-\Gamma$,
\[
 \delta\le\frac{3(1-3y)}4e_{\rm in}+e_{\rm in}^2,\qquad
 \delta\le\frac{3(1-3y)^2}{16}-3ye_{\rm in}+e_{\rm in}^2.
\]
Their monotone intersection gives $0<y<47/250$. The same-lift bound
$\delta>29529/2000000$ then gives the original defect consequences
\[
 e_{\rm in}>\frac{47}{2500},\qquad
 e_{\rm in}>\frac{47}{2500}+\frac y{20}.
\]
The genuine right return provides the corresponding right inequalities.

The incoming endpoint circle is also an actual consequence. Under the
cyclic map to E's $0202$ chart, the inequality $r>0$ becomes
$q_E>p_E$ in reciprocal coordinates, so at its entrance $u_E=-\alpha$
\[
 v_E(-\alpha)^2>(1-\alpha)(3+\alpha).
\]
E's right circle follows from its actual return on arcs $7,8,9$.
The density, defect, and circle inputs for the own-pair proof are
therefore all present. The own mean leaves remain E's actual
zero-to-zero leaves, with their original leaf-radius and density
correlations. In particular, the two interfaces in
\eqref{six:eq:own-interface} hold for these same leaves. The original height argument \ref{pf:M033}
gives $h_E>1/14$, and Proposition~\ref{six:prop:own} applies to arcs
$5,6$, proving the first assertion.

The second assertion assumes an actual terminal circle. To supply
that circle from a one-sided prefix, consider descending arcs
$(6,5,4,3,2,1)$ with the full reflected map $j\mapsto1-j$;
its central arc is $3$. Require the strict exit sign of
Theorem~\ref{hc:prefix650}, which in particular holds when
$c_3\le13/20$. The preceding one-sided-leaf argument then proves
a positive horizontal coordinate at the physical junction $1|2$
in that chart.
For the tail on $[-1,1]$, this is exactly
\[
 v(b)^2>(1-b)(3+b)
\]
at its actual terminal point. For the circle-only tail assertion, its nonpositive entrance,
entrance circle and complete capped context are retained as in the
statement. The applications specify their own complete sign witnesses.

In the proof of Proposition~\ref{six:prop:tail}, the terminal return
is used only through this terminal circle and its cap consequences.
In particular the coefficient of the actual last arc satisfies
$c=q(b)(1-b)^2<(1-b)/(3+b)$; the backward shear bound and the LOW
and HIGH profile bounds follow from the same inequality. No extra
returned mean is used in those comparisons. The two-circle proof
therefore applies to this capped tail. The H7 source, capped tail and opposite prefix must be retained on
one genuine containing interval, as required in the statement.
Normality, the actual loss obstruction and harmonic calibration are
applied on that selected full lift; the specific placements in the
applications verify this containment.
This proves $c_{-1}+c_0+c_1>2/5$ under the stated containing-interval qualifications.
\end{proof}

\subsection{A cap with only following entrance data}
\begin{lemma}\label{late:r18:cap}
A complete proper ten-source $0102021202$ in an original minimizing
trajectory has $c_6>3/40$. The assertion does not require a cap or
terminal return beyond its last arc.
\end{lemma}
\begin{proof}
The returned D source on $[0,7]$ supplies its true height $h>1/14$,
the left endpoint circle of cap 6, and, on one normalized full lift,
\[
 \alpha>h,\quad m_p<-\tfrac12,\quad L_p>0,\quad
 L_p^2>K_{\rm pre}/2-\Gamma,\quad
 0<\Gamma<\frac{1470471}{2000000},\quad L_p>\frac3{25}.
\]
The H7 suffix $[3,9]=2021202$ supplies the negative entrance $-a$
and the other endpoint circle of this cap. To justify the needed
following jet without a missing far cap, put $q=v^{-2}$, $ds=q\,du$,
$\mu=2K_{\rm post}-3$. The two actual internal binary ties give
\[
 \mu=\frac{\int u^2\,ds}{\int ds},\qquad0<\mu<1.
\]
At the initial mode-two point maximization gives $\chi_0\le0$.
Write $\tau=m+1/2$ and let $f$ be the inactive function. At the cap
tie $\tau_0=f_0$, while $g=\tau-f\ge0$ on the following arc.
With $H=1+a$, $B_0=(\mu-a^2)/2$, differentiation gives
\[
 f_s(0)=a+B_0(H^{-1}-q_0H)
       -\frac{1+(q_0H^2)^{-1}}2\chi_0.
\]
The last term is nonnegative. The remaining affine expression is
positive at both ends of $0\le q_0\le(1-a^2)^{-1}$, where its
values are $(a^2+2a+\mu)/(2(1+a))$ and
$a(1-\mu)/(1-a^2)$. Thus $\tau_s(0)>0$ and $L_q<0$.
The same initial inactive identity gives
\[
 m_q+\tfrac12=\frac{\mu-a^2}{2}\rho
 +\frac{(-a-1)\rho+(-a-1)^{-1}}2\chi_0
 >-\frac{a^2\rho}{2},\qquad\rho=v(-a)^{-2}.
\]
These are entrance statements, independent of a later cap.
The whole phase $+1$ sends the actual arcs $5,6,7$ from $212$ to $020$.
In this cap chart write $p_{\mathrm c}<q_{\mathrm c}$ for the two
switching abscissae and $s_{\mathrm c}$ for shear time. By
\eqref{loc:raw-binary}, the switching variable $\chi_{\mathrm c}$
satisfies $(\chi_{\mathrm c})_{s_{\mathrm c}}=\eta_{\mathrm c}-u_{\mathrm c}^2$,
where $\eta_{\mathrm c}=2K_{\mathrm c}-3$ and $K_{\mathrm c}$ is the
binary invariant. Both ties have $\chi_{\mathrm c}=0$, so
$\eta_{\mathrm c}$ is the mean of $u_{\mathrm c}^2$ on the intervening
arc and lies in $(0,1)$. The maximizing signs of $\chi_{\mathrm c}$ on
$020$ give $p_{\mathrm c}^2\ge\eta_{\mathrm c}\ge q_{\mathrm c}^2$.
If either comparison is an equality, the derivative at that tie
vanishes and its abscissa is nonzero. Since $(u_{\mathrm c})_{s_{\mathrm c}}=v_{\mathrm c}^2$
on both binary modes, the second derivative $-2u_{\mathrm c}v_{\mathrm c}^2$
is continuous across the tie and nonzero there. This makes $\chi_{\mathrm c}$ have the same strict sign
on both sides near the tie, a contradiction. Consequently
$p_{\mathrm c}^2>\eta_{\mathrm c}>q_{\mathrm c}^2$, and
$p_{\mathrm c}<0$ follows from $p_{\mathrm c}<q_{\mathrm c}$.
These are the root conditions of Proposition~\ref{six:prop:cap},
with the determinant and terminal sign supplied above.
For clarity, suppose its cost $E\le3/40$. The two endpoint circles
give $E\ge47(\alpha+a)/100$, hence
$a<29/329<9/100$ and $m_q+1/2>-(81/20000)\rho$.
The cap-master estimates used in that proposition give, for
$P=\int L_+\,de$, $Q=\int(-L)_+\,de$,
\[
\begin{array}{c|cc}
 &P&Q\\\hline
\Gamma\ge7/10 & >13/4000 & <401/1000000\\
0<\Gamma<7/10 & >1/200 & <3/5000.
\end{array}
\]
They precede any following HIGH-recovery argument and use only the
cap's actual domain and the endpoint matrix signs just proved.
The same weight estimates give
\[
 \int\omega L\,de
 >\rho\left(\frac{29}{20}P-\frac85Q\right)
 >\frac{81}{20000}\rho.
\]
The two lower brackets are $40709/10000000$ and $629/100000$.
Integration of $m_e=-\omega L$, starting with $m_p<-1/2$,
contradicts the entrance lower bound. Hence $E>3/40$.
\end{proof}

\subsection{A weighted response and a forward pair}
\begin{theorem}[A weighted bilateral core]\label{nx:weighted305}
Let a complete positive $1|0202|1$ source belong to one selected jointly
normalized full maximizing lift of an original minimizer. Suppose both
endpoint rooted responses, in their lawful full-source orientations, are
strictly positive on that same lift. If $A_1,\ldots,A_4$ are the four
complete interior costs, then
\[
 \frac12A_1+A_2+A_3+\frac12A_4>\frac{61}{200}.
\]
The halves are coefficients on complete arcs; neither the trajectory nor
the physical objective is truncated or reweighted.
\end{theorem}
\begin{proof}
Section~\ref{nx:weighted-proof} derives quantitative identities at the two
actual caps. They force a positive excess in the middle cost and bound
both inactive gaps. Under the contrary weighted budget these identities,
the original four-core estimate $A>389/1000$, and the same-lift harmonic
bound give simultaneously
\[
 4\delta>\frac7{125},\qquad
 4\delta<\frac{11163}{203125},
 \qquad
 \frac7{125}-\frac{11163}{203125}=\frac{212}{203125}>0.
\]
The full algebraic proof there uses no parameter cover.
\end{proof}

\begin{theorem}[A complete forward pair]\label{nx:forward-pair}
For a proper complete positive $012012$ source on an original minimizer,
with both exterior caps and a selected full normal maximizing lift
retained, its central costs satisfy
\[
 c_2+c_3>\frac{180483}{500000}.
\]
The indices in this display are local source positions.
\end{theorem}
\begin{proof}
Section~\ref{nx:forward-proof} derives the necessary four-dimensional
domain, including two strict inequalities from the actual exterior caps,
and eliminates the remaining height by a closed interval on which the
full matrix determinant is concave. An analytic corner estimate and the
complete rational interval certificate prove that both endpoint
determinants exceed $H_*$. Concavity gives the same inequality for the
actual lift. Its trace is positive by the ratio-free cap-height estimate;
the original harmonic theorem on that lift gives the contradiction
$\det\Lambda<H_*$. The cover includes zero height, zero ray, weak first
jets, and all endpoint min/max branches. No comparator is calibrated as a
physical trajectory.
\end{proof}

\space
\section{A returned six-arc union with a shared switching jet}\label{opt:r24-application}
\subsection{A six-arc returned-tail payment}

We use the following returned-tail estimate with both actual circles.
\begin{lemma}[Two-circle returned tail]\label{late:r24:two-circle-tail}
Let the complete positive source $120212$ lie in a genuine containing
fundamental interval of an original global minimizer, with one full
jointly normalized maximizing lift. Its capped $202$ core has entrance
$u=-a\le0$, internal roots $r<s$, exit $b$, and its actual middle
mean $\mu$. Retain its complete positive caps and final return.
Suppose also that its entrance satisfies
\[
 v(-a)^2>(1-a)(3+a).
\]
The actual return supplies the terminal circle
$v(b)^2>(1-b)(3+b)$. Then the cost $T$ of the three complete $202$
arcs satisfies $T>2/5$. This holds for either $\mu\le a^2$ or
$\mu>a^2$, including the weak allowed switching jets.
\end{lemma}
\begin{proof}
This is Proposition~\ref{six:prop:tail}, with its complete returned-tail
interface and both circles \eqref{six:eq:two-circles}. The returned-HIGH
absorption and original harmonic loss obstruction are retained there.
In the LOW case the stronger original bound
\eqref{cmp:low-return} is $T>407/1000$.
\end{proof}

\begin{proposition}\label{late:r24:W807}
Suppose the complete positive source $2021201020$, indexed $0,\ldots,9$,
is proper and belongs to an original global minimizer. Then its six
physical charges $2,\ldots,7$ satisfy
\[
 c_2+c_3+c_4+c_5+c_6+c_7>\frac{807}{1000}.
\]
\end{proposition}
\begin{proof}
Choose a genuine containing fundamental interval. Either capped binary
triple proves normality by Lemma~\ref{pc:rigidity}; normalize the whole
lift jointly. The interior source $[1,8]$ is $02120102$.
The left returned source, read in decreasing order $(5,4,3,2,1,0)$,
has word $021202$; full reflected time reversal with label map $1-j$
gives $120212$. The right returned source $[4,9]=201020$, with phase
$-1$, also gives $120212$. Their charged triples are respectively
$\{2,3,4\}$ and $\{5,6,7\}$.

At the central switching point write the positive reciprocal coordinates
as $x,y,z$, so $xy+xz+yz=1$, and put $d=x+z$. The two actual middle
means satisfy the common weak switching-jet identity
\eqref{fee:W71:sj:meanjet}:
\[
 x(1+z^2)^2(\mu_L-u_L^2)
 +z(1+x^2)^2(\mu_R-u_R^2)+2d^3h_0\le0,
 \qquad h_0\ge0.
\]
The coefficients are strictly positive; consequently at least one of
these two actual means satisfies $\mu_j\le u_j^2$. Equality is
included, and no side is fixed in advance.

The complete central $1201$ source gives $y<x,z$ by
Lemma~\ref{fin:forward-force}. The actual outward entrance coordinates
are
\[
 u_R=\frac{y-x}{y+x},\quad v_R=\frac2{x+y},\qquad
 u_L=\frac{y-z}{y+z},\quad v_L=\frac2{y+z}.
\]
Thus both entrances are negative. If $a_j=-u_j$, the reciprocal identity
gives
\[
 v_R^2-(1-a_R)(3+a_R)=\frac{4(z-y)}{x+y}>0,
 \qquad
 v_L^2-(1-a_L)(3+a_L)=\frac{4(x-y)}{y+z}>0.
\]
Both entrance circles are therefore actual state inequalities, and
Lemma~\ref{cmp:tail-interface} supplies both actual terminal circles
from the complete returns. The side selected by the weak jet is LOW,
so \eqref{cmp:low-return} gives $T_j>407/1000$. The other side,
whether LOW or HIGH, satisfies Lemma~\ref{late:r24:two-circle-tail}
and costs more than $2/5$. The two physical charges are disjoint, and
hence
\[
 T_L+T_R>\frac{407}{1000}+\frac25=\frac{807}{1000}.
\]
The means and the selecting jet are from the same lift. The theorem does
not charge the exterior caps or final returns a second time.
\end{proof}

\space
\section{Complete-source half-port profiles}\label{comb:large}
We prove the local profiles used in the all-length run-structure argument.  Throughout, the
trajectory is an original normalized global minimizer, so its physical
cost satisfies $\cost<\Cstar=180483/100000$.

In a mixed increment cycle write the maximal positive and negative runs
as the alternating occurrence cycle
\[
 a_1,b_1,\ldots,a_M,b_M,\qquad
 a_i\ge1,\quad b_i\ge2,\quad b_i\ne3.
\]
The negative restrictions follow from the switching restrictions proved
before this section.  A positive run with $a$ increments owns the arcs
$0,\ldots,a$ and has facing ports $-1,a+1$.  A negative separator with
$b$ increments contributes the interior arcs $1,\ldots,b-1$, with ports
$1,b-1$.  These own sets and interiors partition the physical period.
The two ports of a separator coincide exactly when $b=2$.  The two
incident slots are retained also when $M=1$ or $2$.

\subsection{Complete sources and simultaneous charges}
Every source used below has at most nine complete positive arcs.  For
$B\ge10$ it is proper: cut inside the preceding arc and apply
Lemma~\ref{as:placement} to obtain a genuine containing minimizing
interval with its original endpoints and real winding.  Normality is
established on that complete interval before joint normalization of the
objective multiplier, matrix adjoint and the state covector.  A whole
cyclic phase or reflected time reversal transports all of these data,
not merely the labels.  Different scalar rows may select different
lifts; only their physical costs are added by
Theorem~\ref{as:capacity}.

We use the constants
\[
\begin{gathered}
 A=\frac{389}{1000},\qquad E=\frac{473}{1000},\qquad D=\frac{74}{125},\\
 f=\frac{49}{200},\qquad o=\frac{41}{100},\qquad t=\frac{61}{200},\\
 k=\frac38,\qquad j=\frac{63}{500},\qquad p=\frac{180483}{500000}.
\end{gathered}
\]
The complete source and local charge of each are
\[
\begin{array}{c|c|c@{\qquad}c|c|c}
 &\text{source}&\text{charge}&&\text{source}&\text{charge}\\\hline
 A&102021&1{:}4&E&1020212&1{:}4\\
 D&01020212&2{:}5&f&02010&1{:}3\\
 o&101202&1{:}4&t&012021&2{:}4\\
 k&0120212&2{:}4&k&202120212&4{:}6\\
 j&102012021&4&p&012012&2,3
\end{array}
\]
All these bounds are strict.  The binary and own-four rows are
Theorems~\ref{lc:binary}, \ref{lc:one-return}, \ref{lc:two-return}
and~\ref{lc:own-four}; the conservative $E$ is
Corollary~\ref{loc:E473}.  The six- and seven-source attachments are
Corollary~\ref{ek:short-sources}.  For the nine-source returned attachment,
Corollary~\ref{ek:H7-K} supplies the stronger bound $379/1000>k$
without an entrance-circle or mean predicate.  The full opposite source
qualifies $j$ by Lemma~\ref{as:source-qualification} and
Theorem~\ref{pf:opposite126}, without a HIGH predicate.  The strong pair
is Theorem~\ref{nx:forward-pair}.  These are local consequences of the
original full problem and harmonic comparison; none assumes a bound on
$B$ or a period classification.

\begin{lemma}[Positive half-port profiles]\label{comb:half-strips}
Suppose $B\ge8$ and the supporting windows used in the proof are
proper, as happens for every run when $B\ge10$.
Let $d\in\{0,1,2\}$ be the number of incident separators of length two.
A positive run admits a strict charged profile of value at least
$u(a)+v(a)d$, with load at most one on every own arc and at most one half
on each facing port, where
\[
\begin{array}{c|cc}
a&u(a)&v(a)\\\hline
1&389/2000&21/500\\
2&23/200&13/100\\
3&51/100&7/200\\
4&431/1000&7/200\\
a\ge5&61/200+p\lfloor(a-3)/2\rfloor&7/200
\end{array}
\]
Both incident bonuses can be counted simultaneously.
\end{lemma}
\begin{proof}
Set the label at own arc $0$ to zero.  For $a=1$, take half the
$A,E,D$ core according as $d=0,1,2$.  All three charge $[-1,2]$.
The native full supports are $[-2,3]=210102$,
$[-2,4]=2101020$ for a right return, and
$[-3,4]=12101020$ for two returns; phase $+2$ gives the displayed
sources.  The left-return case is the full reflection of the right.
Their half-values are $389/2000,473/2000,37/125$, respectively;
they dominate $u(1)+v(1)d$, the last by $7/400$.

For $a=2$ and $d\le1$, use the own $f$ core on $[0,2]$;
its source $[-1,3]$, phased by $+2$, is $02010$.
This dominates $23/200+13d/100$.  If $d=2$, the full native source is
$[-3,5]=121012101$.  Phase $+1$ gives $202120212$ with charge
$[1,3]$; reading the entire source backward with map $\ell\mapsto-\ell$
gives the same canonical source with charge $[-1,1]$.
Half of each $k$ row gives $3/8$, with own loads $1/2,1,1/2$ and
half-load on each port.

For every $a\ge3$ a left $t$ charges $[-1,1]$ and a right $t$
charges $[a-1,a+1]$.  The right full source is $[a-3,a+2]$,
which phase $-(a-3)$ sends to $012021$.  The left full source
$[-2,3]$, read backward with $\ell\mapsto-\ell$, is the same source.
If the facing separator is two, extend the respective support to
$a+3$ or $-3$; the actual returning increment gives $0120212$ and
upgrades that identical charge to $k$.  Thus half of each terminal
row contributes $t+(k-t)d/2$ and half-loads on the ports.

At $a=3$, add half the own-four row $[-1,4]=101202$, charging
$[0,3]$.  Its half-load fills exactly the unused own capacity, giving
$t+o/2=51/100$.  At $a=4$, add the single middle charge at arc $2$:
the complete native source $[-2,6]=210120102$, phased by $+2$, is
$102012021$.  This gives $t+j=431/1000$ on disjoint charges.
Finally, for $a\ge5$ the interval $[2,a-2]$ contains
$\lfloor(a-3)/2\rfloor$ disjoint pairs
\[
 \{2+2r,3+2r\},\qquad 0\le r<\lfloor(a-3)/2\rfloor.
\]
Their complete sources $[2r,5+2r]$ lie in $[0,a]$ and phase to
$012012$.  They avoid both terminal charged triples.  These explicit
charges prove all load bounds.  Each construction contains a
positive-weight strict row.
\end{proof}

\begin{lemma}[Negative half-port profiles]\label{comb:negative-strips}
Suppose $B\ge8$ and the supporting windows used in the proof are
proper, as happens for every separator when $B\ge10$.
A negative separator admits a weak profile of value $V(b)$ with load at
most one internally and at most one half on each port, where
\[
 V(2)=0,\qquad V(4)=\frac{299}{250},\qquad V(6)=\frac{13}{10},
\]
\[
 V(b)=\begin{cases}
 b/5-1/2,&b\equiv1\pmod3,\\
 b/5-2/5,&b\equiv2\pmod3,\\
 b/5-3/5,&b\equiv0\pmod3,
 \end{cases}\qquad (b\ge5,\ b\ne6).
\]
The $b=4,6$ profiles are in fact strict and avoid both ports.
\end{lemma}
\begin{proof}
With separator arc $0$ phased to zero, maximality gives the full
source $[-1,5]=2021020$ when $b=4$.  Its actual rooted ingress,
central contact jets and complete terminal $020$ qualify
Theorem~\ref{lc:direct-center}, which gives $299/250$ on arc $2$.
For $b=6$ the full source is $[-1,7]=202102101$.
Lemma~\ref{pf:nine-singletons}, at $\alpha=13/20$, gives two
individual bounds on separator arcs $2,4$, not alternatives of a
maximum; their sum is $13/10$ on disjoint charges.

For every other allowed $b$, put $L=b-1=3m+r$ with $r\in\{0,1,2\}$.
A triple starting at $s$ charges $[s,s+2]$ and retains its complete
source $[s-1,s+3]$.  These are proper five-sources phase-equivalent to
$1\mid021\mid0$, so Theorem~\ref{wk:five600} gives each the weak
value $3/5$.  Use the following families:
\[
\begin{array}{c|c|c}
r&\text{triple starts}&\text{coefficient}\\\hline
0&1+3q\ (0\le q<m),\quad2+3q\ (0\le q<m-1)&1/2\\
1&1+3q,\ 2+3q\ (0\le q<m)&1/2\\
2&2+3q\ (0\le q<m)&1
\end{array}
\]
All immediate caps lie in $[0,b]$.  Within each half-weight family the
cores are disjoint, so overlap between the two families has load at
most one.  Each port lies in at most one half-weight family; the
whole-weight family avoids both.  The total triple weights are
$m-1/2,m,m$, respectively, giving exactly the displayed values.
\end{proof}

\space
\label{revision:appendix-end-3}
\end{proofappendix}
\section{Returned unions and first contact}\label{app:returned}
\begin{proofappendix}
\renewcommand{\Cstar}{C_*}
\section{Returned unions and the qualified first-contact argument}\label{six:sec:problem}
This appendix proves the proper-nine returned-union payment used in
the common local-cost catalogue. It works in the original geometric
class of Section~\ref{fd:section}, with its full state and signed frame
endpoints, real winding and strict physical budget. Put
$P=C_*/5=180483/500000$.

The labels, twisted continuation, complete physical sources and charged
sets are those of Sections~\ref{fd:section} and~\ref{as:section}.
All matrix identities in a coupled-source argument refer to one selected
full maximizing lift. Separate proper sources can export physical state
or cost inequalities without identifying their covectors.

\subsection{Lift coordinates and local-source notation}
\label{six:sec:lift-notation}
Use the full lift of Definition~\ref{fd:full-lift} and equations
\eqref{fd:PMP}, writing $\eta_0=\eta$, $\pi=p$, and $U_*=U$ to distinguish
the state covector from local endpoint coordinates. The multiplier is
nonzero as a whole. Normalization requires $\eta_0>0$ and divides every
component by that same value. The matrix pairing remains
$\operatorname{tr}(\Lambda G^{-1}\delta G)$, and no periodicity of either
state-covector component is imposed on a local source.

In any fixed binary chart write
\[
 \Lambda=\begin{pmatrix}l&b_\Lambda\\m&-l\end{pmatrix},\quad
 K_{\rm bin}=2ul+b_\Lambda-u^2m,\quad L=l-um,
\]\[
 \chi=v\pi_v,\quad q_d=v^{-2},\quad ds=dt/v=q_d\,du.
\]
Thus $m$ is the lower-left entry of the actual matrix adjoint and
$L$ is its displayed shear-adjusted diagonal entry, not a free
response function. On modes zero and two,
\[
 (K_{\rm bin})_s=0,\quad L_s=K_{\rm bin},\quad m_s=-2L,
 \qquad \Gamma:=\det\Lambda=-L^2-K_{\rm bin}m.
\]
The normalized binary switching response is
\[
 f=\left(K_{\rm bin}-\frac{3+u^2}{2}\right)q_d
       +\frac12(\epsilon q_d+\epsilon^{-1})\chi,
 \qquad
 \epsilon=1+u\ (j=0),\quad \epsilon=u-1\ (j=2).
\]
This defines the $f$ used in the returned-tail interface. Its source
is Appendices~\ref{app:local}--\ref{app:crossed}. A lawful orientation-preserving cyclic chart transformation
$W\mapsto M\cdot W$ transports the state covector by the inverse
transpose Jacobian and the matrix by
$\Lambda\mapsto M\Lambda M^{-1}$; both are transported together. Time-reversing reflection additionally
uses the original sign changes and reversed physical time, not this
orientation-preserving rule alone.
Symbols such as $B,R,K$ used later for local scalar parameters are
local to their displayed formulas and are not the arc count or the
fixed endpoint matrix unless expressly indicated.

The source names $D$ and $W$ refer respectively to the complete
words $01020212$ and $02120102$, with positions numbered from zero.
The $D$ four-core charges positions $2,3,4,5$; its exterior arcs
remain part of the eight-arc supporting source. The $W$ source
charges its six middle positions $1,\ldots,6$. The returned
extension and circle witnesses needed for any particular $W$
application remain separate supporting hypotheses; the word alone
does not assert them. These conventions are those of the original
source catalogue and Appendices~\ref{app:local}--\ref{app:crossed}. Here the letter naming a $W$
source is unrelated to the state pair $W=(u,v)$.
\space
\section{Statement and original-domain interfaces}\label{six:sec:hypotheses}
H1 collects the geometric framework of Section~\ref{fd:section} and
Appendix~\ref{app:framework}, together with the local estimates cited
in the proof of Lemma~\ref{six:lem:H1}. The first-contact reduction H2
and the nine-arc transfer H3 use the analytic and finite comparisons in
Sections~\ref{six:app:physical-cap}--\ref{six:app:reduction-proofs}.

\begin{description}
\item[H1: geometric and local interfaces.]
An original global minimizer obeys the normality, compactness and
replacement results stated in Section~\ref{fd:section}. On each genuine
containing interval one may choose a full free-time maximizing lift;
normality is established before joint normalization of its multipliers.
Whenever used below, positive determinant and future sheet on that
same lift justify harmonic calibration. The complete returned-tail,
proper-own-pair, cap-response and shared-arc interfaces stated in
Sections~\ref{six:sec:returned}, \ref{six:sec:cap} and \ref{six:sec:adjacent}
are valid on their indicated actual supports. This includes the LOW
and HIGH absorption implications, the proper forward bound
\eqref{six:eq:forward-input}, and the complete four- and seven-charge
bounds stated in Theorems~\ref{ua:four} and~\ref{ua:seven}.
The named interfaces include the tail moment and entrance-response
identities, the cap master differential comparison on its stated
small-cost domain, and the shared-arc minimization and actual-mean
orientation used inside their proofs. Their exact lift conventions
are fixed in Section~\ref{six:sec:lift-notation}; these imported
implications are not inferred merely from the switching words.
No such implication is assumed for an arbitrary local PMP trajectory.

\end{description}

\begin{lemma}[The original-domain interfaces]\label{six:lem:H1}
For an original global minimizer, the interfaces designated H1 hold
on the complete supports and in the actual lifts specified above.
\end{lemma}
\begin{proof}
The global representation, compactness, full free-time lift and
lawful reanchoring are set out in Section~\ref{fd:section} and
Appendix~\ref{app:framework}. Normality and harmonic comparison
are applied in the order of Theorem~\ref{sw:normality} and
Corollary~\ref{cal:one-point}; all charges are the physical charges
of the original trajectory. The returned-tail response and reserve
are Lemma~\ref{cmp:tail-interface}. Their actual loss obstruction,
LOW estimate and HIGH absorption are, respectively,
Lemma~\ref{cmp:loss-obstruction}, equation~\eqref{cmp:low-return}
and Lemma~\ref{cmp:high-absorption}.

The proper-own transport, density and defect inequalities are
\eqref{pf:d8-transport}, \eqref{pf:d8-density-ceilings}, and
Sections~\ref{pf:M003}--\ref{pf:M006}, with the actual shear
correlation of Lemma~\ref{pf:own-shear-correlation}.
The cap equations used below are those in
Sections~\ref{pf:M040}, \ref{pf:M045}, \ref{pf:M058}
and~\ref{pf:M061}; the strengthened cap estimate is proved below,
not assumed in H1. The shared-arc identities are
Lemma~\ref{loc:shared-A-matrix} and the common-normalization
construction in Section~\ref{pf:F02}; the subsequent proper-own
and shared-arc comparisons retain the actual means of
Section~\ref{pf:F09} and
Sections~\ref{pf:M101}--\ref{pf:M104}. The shared-minimization
formula is proved in Section~\ref{pf:M001}. The density-dependent
own-pair interface is derived in Section~\ref{six:sec:returned}
from the preceding structural identities, without an older numerical
own-pair bound. The older four- and seven-charge
applications use \ref{ua:four} and \ref{ua:seven}; their actual
entrance-circle and union qualifications are
\ref{ek:H7-circle} and \ref{loc:W-unions}.
The forward-source estimate used in \eqref{six:eq:forward-input}
is Proposition~\ref{fp:final}, applied on its complete forward support.
Thus each imported interface is a previously established local
statement on its indicated original support. No identity or
inequality is inferred for an arbitrary formal switching tuple.
\end{proof}

\begin{lemma}[H2: qualified first-contact reduction]\label{six:lem:H2}
Under H1 and the analytic and finite lemmas specified in
Sections~\ref{six:app:physical-cap}--\ref{six:app:reduction-proofs},
for the actual R09 returned union of cost $U<2P$, the original
endpoint bootstrap covers every candidate and gives
$5/9<p,q\le11/15$. The original common-bridge comparison gives both
strict terminal rows and $K<0$ at the common logarithm
$\mathsf X=\log(\mathsf C/4)$, with
$\mathsf C=196517/46875$. Its finite downward continuation retains
the listed true-profile, jet, kappa, endpoint and row-independent
guards. The O and K boundary exclusions, reversed-endpoint-order
reduction and first-root argument yield a qualified first T boundary.
The root-aware two-face reduction excludes its lower geometric face
and gives, on the other face, a common-width comparison satisfying
\eqref{six:eq:S-strip}, $K+4F_p<0$ and $F_q<0$ with the exact
functions of Section~\ref{six:sec:first-contact}.
This implication includes all endpoint equalities and branch ties
specified in its source domains. It does not require periodic state
covectors or noncentral state monodromy.

\end{lemma}

\begin{lemma}[H3: local dependency restriction]\label{six:lem:H3}
Under the same framework and identified lemmas,
the successful implication in H2 also applies under the contrary
$U<2P$ to the complete proper nine-source $012021201$, with charges
$\{2,3,4,5,6\}$, its two actual returned contexts, two forward entrance
witnesses, all retained outgoing leaves and ties, and one common full
lift. Only the local budgets and the global-minimizer qualifications
in H1 may be imported. In particular no Casimir lower bound derived
from R09's three exterior forward pairs is imported. The support
identification in Section~\ref{six:sec:first-contact} specifies this
restriction; equality $U=2P$ is not part of the contrary domain.
\end{lemma}

Appendix~\ref{six:app:reduction-proofs} combines the physical
implications with the scalar comparisons and locates the detailed
analytic proofs in the separate \path{h2-h3-analytic-sources} supplement.
Appendix~\ref{six:app:computation} specifies the finite comparisons.
These interfaces concern the qualified sources of original global minimizers.

\space
\section{The common returned union and first contact}
\label{six:sec:first-contact}
The preceding cost estimates alone do not resolve the common returned
union. Its two returned triples overlap, and the shared arc must not
be paid twice. The original approach retains the joint bridge and the
two exact profile means, rather than summing independent endpoint
relaxations. Here we specify the scalar comparison reached by that
approach and explain its final reduction. The physical implication to
this scalar system is Lemma~\ref{six:lem:H2}, proved in Appendix~\ref{six:app:reduction-proofs}.

\subsection{One common logarithm}
To avoid confusion with the full physical cost $C$, write
\[
 \mathsf C=\frac{196517}{46875},\quad
 \mathsf X=\log(\mathsf C/4),\quad
 p_*=\frac{11}{15},\quad u_*=\frac{180483}{250000}.
\]
At a common width $w>1$ let
\[
 L(w)=\frac{\log w}{w-1},\qquad
 m(z,w)=1+\frac{(2-z)^2}{w}-2(2-z)L(w).
\]
Both means use this same logarithm. Put $y=w^2-1$ and
\[
 D_{\rm g}=\tfrac12(1/p+1/q-p-q),\quad
 e_z=1-\frac{2-z}{w+1},\quad
 g_z=\max\{0,1-(2-z)/w\},\quad \mathsf R=\frac{\mathsf C}{pq}.
\]
The two normalized terminal complements are
\begin{align*}
 F_q={}&t_q-d_{11}m(q,w)-d_{12}m(p,w),\\
 F_p={}&t_p-d_{21}m(q,w)-d_{22}m(p,w),
\end{align*}
where
\begin{gather*}
 d_{11}=\frac1{2q}+\frac{py}{\mathsf Cq},\quad
 d_{12}=\frac1{2p}+\frac{\mathsf Cq}{16py},\quad
 d_{21}=\frac1{2q}+\frac{\mathsf Cp}{16qy},\quad
 d_{22}=\frac1{2p}+\frac{qy}{\mathsf Cp},\\
 t_q=\mathsf X-1+D_{\rm g}+\frac32e_q+\frac{pq}{\mathsf C}
            +\frac{3py(e_q^2-1)}{\mathsf Cq},\\
 t_p=\mathsf X-1+D_{\rm g}+\frac32e_p+\frac{pq}{\mathsf C}
            +\frac{3qy(e_p^2-1)}{\mathsf Cp}.
\end{gather*}
The corresponding action counts one contribution from each of the two
profiles, including when their endpoint parameters coincide. Put
$z_1=p$ and $z_2=q$. Then
\begin{align*}
 K={}&2(\mathsf R-1)\left(1-u_*+\frac{g_q+g_p}{2}\right)
       -2\mathsf X-2D_{\rm g}\\
 &+\sum_{i=1}^2
 \frac{6y-3(2w-2+z_i)^2/(1+g_{z_i})+(2-z_i)(w-1)m(z_i,w)}{z_i^2}.
\end{align*}
The first-contact comparison retains $B_4=K+4F_p<0$ and
$F_q<0$ at its equal-width base. Consequently the function
\begin{equation}\label{six:eq:S24}
 S=K+4F_p+24F_q
\end{equation}
is strictly negative there. In particular no independent sign of
$F_p$ is assumed.

The comparison domain after the other base face is excluded is
\begin{equation}\label{six:eq:S-strip}
 \frac{13}{20}<p\le p_*,\qquad \frac59<q\le p,\qquad
 \frac{19}{15}\le w\le\frac{27}{20},\qquad w\ge2-p.
\end{equation}
The equal-width replacement and later changes of $p,w$ are scalar
comparisons. They do not assert that the varied tuple remains an
original state, satisfies group closure, or is a new physical extremal.

\subsection{Decrease in the endpoint parameter}
Appendix~\ref{six:app:derivatives} proves the common-profile derivative bounds
\[
 K_p<5/2,\qquad (F_q)_p<-1/2.
\]
Here the $K_p$ bound is the ordinary fixed-$q,w$ derivative when
$w>2-p$, and the derivative from the side of increasing $p$ when
$w=2-p$. No two-sided derivative across that geometry seam is
asserted. For completeness, the additional own-row derivative is obtained by
collecting the exact displayed $F_p$ at fixed $q,w$:
\[
 F_p=A_{\rm own}p+B_{\rm own}+D_{\rm own}/p,
\]
\begin{align*}
 A_{\rm own}={}&-\tfrac12+\frac3{2(w+1)}
 +\frac q{\mathsf C}\left(1+\frac{3(w-1)}{w+1}\right)
 -\frac{1/2+qy/\mathsf C}{w}
 -\frac{\mathsf C\,m(q,w)}{16qy},\\
 D_{\rm own}={}&\tfrac12-\frac{12qw(w-1)}{\mathsf C(w+1)}
 -(1/2+qy/\mathsf C)m(0,w).
\end{align*}
The strip bounds give $A_{\rm own}<1/2$:
use $(w-1)/(w+1)<3/20$, $q/\mathsf C<11/60$ and
$w+1\ge34/15$. Also
\[
 \frac{12qw(w-1)}{\mathsf C(w+1)}
 <\frac{2079}{4700}<\frac9{20},\qquad
 \frac12+\frac{qy}{\mathsf C}<\frac23,\qquad 0<m(0,w)<\frac23.
\]
Thus $D_{\rm own}>-71/180>-2/5$ and
\[
 (F_p)_p<\frac12+\frac{160}{169}<\frac32.
\]
Together these estimates give
\begin{equation}\label{six:eq:S-decrease}
 S_p<5/2+4(3/2)-24(1/2)=-7/2<0.
\end{equation}
At $w=2-p$, equation~\eqref{six:eq:S-decrease} has the same
right-derivative interpretation. The scalar expressions match
continuously at the seam. Increasing $p$ to $p_*$ preserves the
strip and immediately enters $w>2-p$ if the starting point is on
the seam; integrating the derivative bound on the smooth portion
and taking the continuous endpoint limit proves that $S$ decreases
strictly whenever $p<p_*$.

\subsection{The three boundaries}
Each true profile has $m_w<0$ and $m_{ww}>0$ throughout the strip.
For the latter assertion, with $a=2-z$,
\[
 m_{ww}=2a[a/w^3-L''(w)],\quad
 L''(w)=2\int_0^1\frac{t^2}{[1+(w-1)t]^3}\,dt.
\]
The elementary bound $(1+x)^{-3}\le1-3x+6x^2$ gives
$L''(w)\le164/375<1/2$, while $a/w^3>1/2$.
The original normalized-row width concavity with the opposite mean
fixed remains strict when that true mean is restored: the additional
terms are $-2d'm_w-dm_{ww}<0$, since $d'<0$, $m_w<0$ and $d>0$.
The exact action is also strictly concave within each geometry sector.

At $p=p_*$, however, the action has an upward derivative jump at
$w=2-q$. Thus only separate-sector concavity is available.
The remaining boundaries are precisely
\[
 w=19/15,\quad 5/9\le q\le p_*;\qquad
 w=27/20,\quad 5/9\le q\le p_*;
\]
\[
 w=2-q,\quad 13/20\le q\le p_*.
\]
Appendix~\ref{six:app:three-edges} proves the following three estimates:
\begin{equation}\label{six:eq:S-edge-inputs}
 \left.S\right|_{w=19/15}>\frac3{40},\qquad
 \left.S\right|_{w=27/20}>0,\qquad
 \frac{d^2}{dw^2}S(p_*,2-w,w)<-\frac{3851}{28500}.
\end{equation}
The endpoints of the last curve belong to the first two boundaries.
Consequently its concavity and endpoint positivity make the whole
seam positive. Concavity on each width sector then gives $S>0$
throughout \eqref{six:eq:S-strip}, using \eqref{six:eq:S-decrease}.
This contradicts the necessary strict sign in \eqref{six:eq:S24}.

\begin{theorem}[Common-log scalar comparison]\label{six:thm:scalar}
With the exact constants and functions defined above, every tuple in
\eqref{six:eq:S-strip} satisfies $K+4F_p+24F_q>0$.
In particular no such tuple satisfies both $K+4F_p<0$ and $F_q<0$.
\end{theorem}
\begin{proof}
The endpoint derivative \eqref{six:eq:S-decrease} reduces to $p=p_*$.
The three boundary signs proved in Appendix~\ref{six:app:three-edges}
and separate-sector concavity proved in Appendix~\ref{six:app:derivatives}
then give positivity on every closed width sector. At a seam the
same positive boundary value is retained. These sectors cover the
entire strip, including its equalities.
\end{proof}

\subsection{The nine-source support}
The standard eleven-source is $01201020120$. The nine-source
$012021201$, after phase $+1$, is exactly its positions $1,\ldots,9$.
Its ascending returned context occupies local $(1,2,3,4,5,6)$.
The descending context occupies $(7,6,5,4,3,2)$ and becomes $120212$
under full reflection with $j\mapsto1-j$.
The two triples charge $\{2,3,4\}$ and $\{4,5,6\}$, so their union
charges $\{2,3,4,5,6\}$ only once.
Complete forward sources at $[0,3]$ and $[5,8]$ supply the two actual
negative entrance witnesses.

The outgoing leaves used in the later first-contact argument are
standard eleven-source positions $7$ and $3$, with ending ties $7/8$
and $2/3$. They become nine-source positions $6$ and $2$, with ties
$6/7$ and $1/2$; all remain inside the support. The two deleted
exterior arcs were entrance witnesses, not these delayed leaves.
This support check is necessary for transferring the argument and
must be combined with the dependency restriction proved in
Appendix~\ref{six:app:reduction-proofs}. In particular a lower Casimir bound obtained
from R09's three exterior forward pairs cannot be imported into a
general nine-source. The intended transferred conclusion is only the
weak bound \eqref{six:eq:union-input}, not an exclusion of its equality case.

The successful local implication is recorded as H3. It retains the
original global-minimizer, compactness and harmonic qualifications;
it is not a general local-control theorem. With those hypotheses,
the support restriction above and the scalar contradiction imply
exactly the weak union estimate used in the applications.
\space
\section{A cap between a four core and a returned tail}
\label{six:sec:cap}
The cap considered here follows a complete $D$ four-core and precedes
a returned $120212$ tail. Its role is different from that of the bridge
following the tail. Both supporting sources are retained on one proper
containing interval, with one full jointly normalized maximizing lift.

In the canonical decreasing-$U$ chart the cap runs from $\alpha>h>1/14$
to $-a<0$. Both endpoints have their actual circle inequalities.
At the preceding endpoint the four-core supplies
\[
 0<\Gamma<\frac{1470471}{2000000},\quad K_{\rm pre}>3/2,
 \quad m_p<-1/2,\quad
 L_p>\sqrt{K_{\rm pre}/2-\Gamma}>3/25.
\]
At the following endpoint $L_q<0$. These are actual endpoint identities
and signs, not a calibration imposed on an independently chosen matrix.
In the cap's own mode-two chart the root and mean qualifications are
\[
 -1<p<q<1,\qquad p<0,\qquad p^2>\eta>q^2,\qquad
 K=(3+\eta)/2\in(3/2,2).
\]
The symbols $m_p,L_p,m_q,L_q$ are the values of $m=\Lambda_{21}$
and $L=\Lambda_{11}-U\Lambda_{21}$ at the two physical endpoints,
in the fixed canonical chart of the adjoining cores. They are not
entries from independently chosen lifts. To distinguish charts, let
$(u,v)$ denote the cap's pre-cyclic mode-two chart and $(U,V)$ the
canonical decreasing-$U$ chart. They satisfy
$(U,V)=R\cdot(u,v)$, with simultaneous transport
$\Lambda_{\rm can}=R\Lambda_{\rm pre}R^{-1}$ and the cotangent
Jacobian rule of Section~\ref{six:sec:lift-notation}.
In the pre-cyclic chart the actual density is
\[
 v^{-2}=\frac{k}{(1-u)^2},\qquad
 k=\frac{(1-u)^2}{v^2}>0.
\]
The mode-two equation makes $k$ constant along this cap. It is a
state-density parameter, distinct from $K_{\rm bin}$ or the
preceding core's $K_{\rm pre}$. With $w=(1-u)^{-1}$ the cyclic
state formulas are
\[
 U=\frac{1-k(4w^2-1)}{1+k(2w-1)^2},\qquad
 V=\frac{4w\sqrt{k}}{1+k(2w-1)^2}.
\]
These conventions are specified in the original
Sections~\ref{pf:M058} and~\ref{pf:M040}.
Thus the derivative $m_e$ below differentiates this same canonical
matrix entry along the physical cap, reparametrized by $e$.

The generic cap response equations of
Sections~\ref{pf:M040}, \ref{pf:M045} and~\ref{pf:M061}
apply before their final HIGH recovery step. We use only those equations
and the estimates derived from the endpoint conditions just listed.

\begin{proposition}\label{six:prop:cap}
The complete cap under these qualifications has physical cost $E>3/40$,
irrespective of the following tail's actual mean.
\end{proposition}
\begin{proof}
Suppose $E\le3/40$. On the cap,
$d=(U^2+V^2-1)/V>0$ is invariant. The two circles give
\[
 d^2>\frac{4(1-\alpha)}{3+\alpha},\qquad
 d^2>\frac{4(1-a)}{3+a}.
\]
Lemma~\ref{ua:bridge} applied to the two halves of this cap gives
$E\ge(47/100)(\alpha+a)$. The density is even in $U$, so the
negative half is integrated using $U\mapsto-U$ on the same cap. Thus
\begin{equation}\label{six:eq:cap-small-entrance}
 a<\frac{15}{94}-\frac1{14}=\frac{29}{329}<\frac9{100}.
\end{equation}
Let $\rho=V(-a)^{-2}$, and let $\mu>0$ and $\chi_q\le0$ be the genuine
following middle mean and first switching response. The first-cap identity is
\[
 m_q+\frac12=\frac{(\mu-a^2)\rho}{2}
 +\frac{[(-a-1)\rho+(-a-1)^{-1}]\chi_q}{2}.
\]
Its second term is nonnegative, hence
\begin{equation}\label{six:eq:cap-lower-response}
 m_q+\frac12>-\frac{81}{20000}\rho.
\end{equation}
No comparison between $\mu$ and $a^2$ was needed.

For the generic cap chart put
\[
 w=(1-u)^{-1},\quad k>0,\quad D=1+k(2w-1)^2,\quad
 R=1+2Kkw^2,\quad \omega=D/R,
\]
\[
 e_w=\frac1{2w^2}+Kk>0,\qquad e(w_p)=0,\quad e(w_q)=E.
\]
The increasing comparison parameter $e$ is not identified with the
physical cost of a partial cap. Its final value is the complete cost.
The exact response equation is
\begin{equation}\label{six:eq:cap-response}
 m_e=-\omega L.
\end{equation}
On a cap of the assumed size the original master comparison gives
\[
 \frac12<w<\frac{1000}{937}<\frac{10}{9},\qquad
 L_e>-\lambda L-b_\Gamma(e),\quad 0<\lambda<\frac52,
\]
\[
 b_\Gamma(e)=2+\frac{35}{4}(e-2\Gamma/21)_+,\qquad
 |(\log\omega)_e|<\frac54.
\]
In particular
\begin{equation}\label{six:eq:cap-weight}
 \frac{\omega_{\max}}{\omega_{\min}}<e^{3/32}<\frac{32}{29}.
\end{equation}
The last inequality follows from $e^x<(1-x)^{-1}$ for $0<x<1$.

Set $P_+=\int L_+\,de$, $P_-=\int(-L)_+\,de$, and $e_0=11/200$.
If $\Gamma\ge7/10$, the positive master comparison through $e_0$ gives
\[
 P_+>\frac{67076339}{20480000000}>\frac{13}{4000}.
\]
The first zero lies after $e_0$, and every subsequent negative component
is bounded by the negative envelope. Consequently
\[
 P_-<(E-e_0)^2+\frac{35}{24}(E-1/15)_+^3
 \le\frac1{2500}+\frac7{8294400}<\frac{401}{1000000}.
\]
For $\Gamma<7/10$ the preceding determinant gives $L_p>1/5$ and
$b_\Gamma<3$. Until the first zero,
$L(e)>1/5-7e/2$, so
\[
 P_+>\frac{913}{160000}>\frac1{200},\qquad
 P_-<\frac32(E-e_0)^2\le\frac3{5000}.
\]
If $E<e_0$, the positive comparison already contradicts $L_q<0$.
Thus these integrals are taken on the stated complete interval whenever
needed.

It remains to compare the response weight with the actual following
entrance density. At the terminal state the same chart gives
\[
 U=\frac{1-k(4w^2-1)}D=-a,\qquad
 V=\frac{4w\sqrt k}D,\qquad
 \rho=\frac{D^2}{16w^2k}.
\]
Solving for the actual $k,D$ gives
\[
 k=\frac{1+a}{(2w-1)[2w+1-a(2w-1)]},\qquad
 D=\frac{4w}{2w+1-a(2w-1)}.
\]
Since $K<2$, positive denominators imply
\[
 \frac{\omega_q}{\rho}>
 \frac{4w(1+a)[2w+1-a(2w-1)]}
 {(4w^2-1)-a(2w-1)^2+4w^2(1+a)}>\frac85.
\]
For the last inequality the numerator minus $8/5$ times the denominator is
\[
 -\frac{24}{5}w^2+4w+\frac85
 +a\left[\frac85(w+1)-4wa(2w-1)\right].
\]
The first quadratic is concave on $[1/2,10/9]$ and has positive endpoint
values $12/5,16/135$. The bracket is positive because its first term is
at least $12/5$, while the subtracted term is less than $44/81$, using
\eqref{six:eq:cap-small-entrance}.

By \eqref{six:eq:cap-weight},
\[
 \int\omega L\,de
 \ge\omega_{\min}P_+-\omega_{\max}P_-
 >\rho\left(\frac{29}{20}P_+-\frac85P_-\right).
\]
Indeed first factor $\omega_{\min}$, use
$P_+-(32/29)P_->0$, and then
$\omega_{\min}>(29/32)\omega_q>(29/20)\rho$.
This avoids the unjustified assertion $\omega_{\max}\le\omega_q$.
For the two determinant ranges the bracket is respectively greater than
\[
 \frac{29}{20}\frac{13}{4000}-\frac85\frac{401}{1000000}
 =\frac{40709}{10000000}>\frac{81}{20000},
 \qquad
 \frac{29}{20}\frac1{200}-\frac85\frac3{5000}
 =\frac{629}{100000}>\frac{81}{20000}.
\]
Integrating \eqref{six:eq:cap-response} and using $m_p<-1/2$ now gives
$m_q<-1/2-(81/20000)\rho$, contradicting
\eqref{six:eq:cap-lower-response}.
\end{proof}
\space
\section{Adjacent returned pairs}\label{six:sec:adjacent}
Retain the complete positive proper source $01020212101$.
Its returned sources on $[0,7]$ and $[3,10]$, with whole phase $+1$
on the latter, are both $01020212$. Their own-pair charges are
$\{3,4\}$ and $\{6,7\}$. Both contain their actual returns and
positive caps; because $11<B$, each can be retained in a genuine
containing interval.

Corollary~\ref{six:prop:adjacent} therefore gives
\[
 c_3+c_4+c_5+c_6+c_7>\frac{18}{25}>\frac{161}{250}.
\]
It adds the two disjoint $9/25$ own-pair payments and the positive
shared cost $c_5$. This scalar addition needs no identification of
separately chosen matrices, covectors or cap defects. The complete
eleven-source and its physical charge are the same in every use of
the weaker $161/250$ consequence.
\space
\section{The eleven-source first-contact comparison}
\label{six:sec:applications}
The following two applications use the qualified first-contact union.
The complete-forward
estimate used below is
\begin{equation}\label{six:eq:forward-input}
 c_2+c_3>P=\frac{180483}{500000}
\end{equation}
for an actual complete proper $012012$ source, including both exterior
caps; see Proposition~\ref{fp:final}.
The second input concerns a complete nine-source $012021201$. Its
charged local positions are $\{2,3,4,5,6\}$ and its bound is
\begin{equation}\label{six:eq:union-input}
 c_2+c_3+c_4+c_5+c_6\ge2P.
\end{equation}
The latter input retains both returned contexts, both complete forward
entrance witnesses, all actual outgoing leaf ties used by the
first-contact argument, and the genuine containing interval.
We state explicitly where \eqref{six:eq:union-input} is needed; the cost
accounting does not prove that input.

\subsection{The complete eleven-source}
Its native labels and twisted extension are
\[
 (j_0,\ldots,j_{10})=(0,2,1,2,0,1,2,0,1,2,0),
 \qquad j_{i+11}=j_i+1\pmod3.
\]
The genuine fundamental interval from switch $7$ to switch $18$
carries $01201020120$. Its eleven-arc length is allowed by the
containing-interval hypothesis; its internal caps and sign ties do
not require an exterior switch or a periodic state covector.
The union charges physical positions $\{10,0,1,2,3\}$, with cost $U$.
The three complementary forward pairs are $\{4,5\}$, $\{6,7\}$
and $\{8,9\}$. Their complete supports are respectively
$[2,7]=120120$, $[4,9]=012012$, and $[6,11]=201201$;
whole phase shifts $-1$, $0$, and $+1$ take these to $012012$.
They provide the strict forward payment, so
\[
 U=C-(c_4+c_5)-(c_6+c_7)-(c_8+c_9)
    <\Cstar-3P=2P=180483/250000.
\]
This uses no double counting of the two overlapping returned triples.

Under H2 the resulting actual candidate reaches the comparison strip
\eqref{six:eq:S-strip} with $S=K+4F_p+24F_q<0$.
The scalar theorem of Section~\ref{six:sec:first-contact} gives $S>0$,
a contradiction. This argument is about the actual full word;
it excludes the qualified actual candidate, rather than every
auxiliary covariance owner.

For clarity, the logical join behind H2 uses a new common-log root,
not a substitution into a root constructed with another logarithm.
After whole reflection orders the widths, the inherited reversed-order
argument excludes the alternative endpoint order. At fixed endpoints
and fixed smaller profile one lowers only the larger profile width.
On the specified comparison interval, its own normalized row $F$ is
strictly concave, the excluded O/K lower boundary has $F(l)>0$, and
the starting point has $F(u_0)<0$. There is a unique crossing $u_T$:
at a zero, concavity gives
\[
 F'(u_T)\le\frac{F(u_T)-F(l)}{u_T-l}<0,
\]
and afterwards the derivative cannot increase. In the original row
normalization this gives $m_p=0$, $m_p'>0$, while the opposite row
retains $m_q>0$. The inherited root-aware estimate
$(K-4m_p)'>1/6$ on the interval from the equal-width base to $u_T$
is combined with the retained $K(u_T)<0$ and
$B_4(u_T)=K(u_T)$, where $B_4=K-4m_p$. Integration gives
\[
 B_4(v_0)<K(u_T)-(u_T-v_0)/6<0
\]
at the equal-width base $v_0$. The other row remains positive,
so $S=(K-4m_p)-24m_q<0$. It is not necessary, and not justified,
to impose $m_p=0$ at the equal-width base itself.
The guard preservation, O/K exclusions and root-aware estimate in
this paragraph are precisely the implications of Lemma~\ref{six:lem:H2};
their finite prerequisites are not reproved by the concavity argument.

Under H3 the same contradiction applies to a proper nine-source
whenever its five-arc union has $U<2P$. Therefore every such actual
source has $U\ge2P$, establishing \eqref{six:eq:union-input} under the
stated hypotheses. No strict lower bound at equality is inferred.

\space

\section{Scalar derivative identities}\label{six:app:derivatives}
Use the constants and strip of Section~\ref{six:sec:first-contact}; as in
Appendix~\ref{six:app:three-edges}, write $C=\mathsf C$.
These computations justify the derivative inputs to the scalar
comparison. They do not supply the physical-to-scalar implication.

\subsection{The true branch}
For $z\in[5/9,11/15]$, put $a=2-z$, $s=1-z$,
$r=1-a/w$ and $\mu=m(z,w)$. The density is proportional to
$(1-x)^{-2}$ on $[-s,r]$, and direct differentiation gives
\[
 \mu'=(r^2-\mu)/(w-1).
\]
At $w=27/20$, if $r\ge0$, then $r\le5/81$, $s\ge4/15$,
and $r/s\le25/108<1/4$. The density ratio is $400/729$, so
\[
 \mu\ge\frac{400}{729}\frac{s^2-sr+r^2}{3}
       \ge\frac{325}{2187}s^2>r^2.
\]
If $r\le0$ the strict inequality follows from $x^2>r^2$
in the interior of the integration interval. For $r>0$, at any
zero of $J=\mu-r^2$ one has $J'=-2ra/w^2<0$.
Thus a zero can only cross from positive to negative as $w$
increases; the proved positive upper value rules out every such
zero below it. Hence $\mu>r^2$ and $\mu'<0$ throughout the strip.
The computation of $m_{ww}>0$ in Section~\ref{six:sec:first-contact}
then applies to both profiles.

\subsection{Endpoint derivatives}
At fixed $q,w$ the exact other-row complement is
$F_q=Ap+B+D/p$, where, with $y=w^2-1$ and
$e=1-(2-q)/(w+1)$,
\begin{align*}
 A&=-(1+1/w)/2+q/C+y[3e^2-3-m(q,w)]/(Cq)-Cq/(16yw),\\
 D&=2(L-1/w)-[Cq/(16y)]m(0,w).
\end{align*}
Here $A<-47/54+11/60$. Also
$2(L-1/w)>(w-1)/w^2\ge60/361$,
$m(0,w)<203/323<2/3$ and $Cq/(16y)<693/2176$.
Thus $D>-1/20$, and
\[
 (F_q)_p<-47/54+11/60+20/169
        =-1/2-6269/91260<-1/2.
\]
Remove the two nonnegative mean terms from $K$ and call the result
$K_0$. On the open positive-$p$ geometry $w>2-p$, exact
fixed-width differentiation, with $b=2-p$, gives
\[
 (K_0)_p=-\frac{\mathsf R}{p}[3-2u_*+g_q-2/w]
       -1/w+1/p^2+1-[3w(2+b)-12]/p^3.
\]
The bracket is greater than $-109/4750$ and
$\mathsf R/p<2448/125$, so its contribution is less than $1/2$.
The other three contributions are bounded above by
$-2/3$, $89/25$, and $-1984/2025$.
Their sum with $1/2$ is $9775/4050<5/2$.
The $q$-mean term is independent of $p$. The $p$-mean term is
$b(w-1)m(p,w)/p^2$; its positive coefficient decreases in $p$,
and its mean derivative is $2(L-b/w)<0$ by the trapezoidal bound
$L<(1+1/w)/2$ and $w+1-2b\le-11/60$.
It follows that $K_p<(K_0)_p<5/2$ for $w>2-p$.
At the closed face $w=2-p$, the positive-geometry formula
agrees continuously with the exact action and gives its derivative
from increasing $p$. The displayed bound holds for this right
derivative, not in general for a two-sided derivative through the
geometry seam. For fixed $q,w$, the allowed segment from $p$ to
$p_*$ lies entirely on this side of the seam, so continuity and
integration give the endpoint comparison used in the main text.

\subsection{Normalized-row width concavity}
Fix the opposite mean $\mu>0$ and vary the own width $u$ at
fixed endpoints $p,q\in[5/9,11/15]$. Put
\[
 b=2-p,\ t=u-1,\ Y=u^2-1,\ r=1-b/u,\ \nu=m(p,u),
 \quad\gamma=q/(Cp),\quad\beta=Cp/(16q).
\]
The own row has the exact form, with a width-independent $B_0$,
\[
 F=B_0-\mu/(2q)+(3/2)e_p-\beta\mu/Y
          +\gamma Y(3e_p^2-3-\nu)-\nu/(2p),
 \qquad e_p=1-b/(u+1).
\]
Twice differentiating, using $\nu'=(r^2-\nu)/t$ once at each
stage, gives
\begin{align*}
 F''={}&-\frac{2\beta\mu(3u^2+1)}{Y^3}
       -\frac{b(3+12\gamma b)}{(u+1)^3}
       -2\gamma r(1+b/u^2)\\
      &-\frac{(\nu-r^2)+rbt/u^2}{pt^2}.
\end{align*}
For $r\ge0$ all terms are nonpositive and some strictly negative.
For $r<0$, the density random variable $Z$ on $[-(1-p),r]$
satisfies
\[
 \nu-r^2=\mathbb E[(r-Z)(-Z-r)]>-2r\mathbb E[r-Z],\qquad
 \mathbb E[r-Z]=b(L-1/u)>bt/[u(u+1)].
\]
Thus the last numerator exceeds $(-r)bt^2/[u^2(u+1)]>0$.
For the possibly positive third term group its $\gamma$-part
with that of the second term. On $9/8\le u<b\le13/9$,
\[
 (-r)(1+b/u^2)\le35903/59049<5/8,
 \qquad 6b^2/(u+1)^3>13/20.
\]
The grouped bracket is therefore less than $-1/40$.
This proves $F''<0$ in both geometries.

If the opposite mean is restored to its true width-dependent value,
the extra second-derivative terms are $-2d'\mu'-d\mu''<0$,
where $d>0$, $d'<0$, $\mu'<0$, $\mu''>0$.
Thus both normalized rows are strictly concave at common width.

\subsection{Action sector concavity}
For an endpoint $q$ put $a=2-q$. Its contribution to the action,
apart from a width-independent constant, is
\[
 G_q(w)=(\mathsf R-1)g_q+
 [6(w^2-1)-3(2w-a)^2/(1+g_q)+a(w-1)m(q,w)]/q^2.
\]
The identity $[(w-1)m(q,w)]'=r^2$ yields
\[
 G_q''=\begin{cases}
 [-12+2a^2(w-a)/w^3]/q^2,&w<a,\\
 2a[a(w-a)-q^2(\mathsf R-1)]/(q^2w^3),&w>a.
 \end{cases}
\]
Both are negative. For the second, $w<2$ gives
$a(w-a)<aq<1$, while $q^2(\mathsf R-1)>50/27$.
The reflected side is identical. At $w=a$ the first derivative
jumps upward from $0$ to $(\mathsf R-1)/a+3a/q^2$.
Thus $K$, and hence $S$, is strictly concave on each sector,
but no concavity across that seam is asserted.

\subsection{Monotonicity of the original bridge determinant}
\label{six:app:bridge-mean}
The early harmonic comparison uses one common bridge, not independent
endpoint determinants. Here write $C_0=196517/46875$ and let
$0<p,q\le4/5$, $7<R<8$ and $Rpq>C_0$. Put
\[
 X=\log(Rpq/4)>0,\qquad A=R/(R-1),
\]
\[
 B=\frac12-\frac{AX}{2}
   -\frac A4\left(\frac{1-\nu}{p}+\frac{1-\mu}{q}-p-q\right),
 \qquad c=\frac{A^2(3+\mu)(3+\nu)}{16}-B^2.
\]
On the complete support rectangle
$0\le\mu\le(1-q)^2$, $0\le\nu\le(1-p)^2$, differentiation gives
\[
 \frac{16}{A^2}c_\mu=3+\nu-\frac{8B}{Aq}.
\]
Substituting $B$ and then using the upper bound for $\mu$ yields
\[
 \frac{16}{A^2}c_\mu\ge
 \nu-1-\frac{2p}{q}+\frac{2(1-\nu)}{pq}
                 +\frac4{Rq}+\frac{4X}{q}.
\]
The coefficient of $\nu$ on the right is $1-2/(pq)<0$.
Its upper support bound therefore gives
\[
 \frac{16}{A^2}c_\mu\ge
 \frac{4(1-p)}q-2p+p^2+\frac4{Rq}+\frac{4X}{q}
 \ge\frac1{25}+\frac4{Rq}+\frac{4X}{q}>0.
\]
Here $4(1-p)/q\ge1$ and $-2p+p^2\ge-24/25$.
Reflection gives $c_\nu>0$ as well. In particular, when the actual
means are at least $m=3/400$, lowering both to $m$ at fixed actual
$R,p,q$ does not increase the same bridge determinant $c$; it decreases
$c$ if either mean exceeds $m$. Since $m<1/25$, these
constant floors remain inside both support intervals.

The abbreviated assertion $B<1/2$ alone does not imply either mean
derivative sign; the estimates above provide the needed strict
argument. The intermediate mean values are scalar comparison values,
not new physical profiles. Subsequent radius transport and endpoint
replacement retain their separate source hypotheses and order: lower
the means first at fixed actual endpoints and radius, transport radius
within its proved interval, and only then use the endpoint comparison
on the target face. No global-minimizer compactness bound is imposed
on those formal intermediate tuples.
\space
\section{The three scalar boundaries}\label{six:app:three-edges}
We prove the three estimates in \eqref{six:eq:S-edge-inputs} directly
from the functions in Section~\ref{six:sec:first-contact}. Throughout this
appendix write $P=11/15$, $C=\mathsf C$, $X=\mathsf X$ and
$h=69517/250000$. This temporary $C$ is a scalar comparison constant,
not the total physical cost.

\subsection{The lower boundary}
Set $w=19/15$, $L=\log(w)/(w-1)$ and
\[
 M=79/19-4L,\qquad \nu=34/15-38L/15,\qquad
 k=1/2+1496/(3375C).
\]
The positive series for $\operatorname{artanh}$ and the trapezoidal
upper bound for $1/t$ give
\[
 15/17+20/4913<L<17/19,\quad
 0<M<49/80,\quad 3/5<k<49/80,\quad \nu>0.
\]
Collecting the displayed original action and row, with both geometry
terms nonpositive, gives
\[
 K=A_Kq+B_K+D_K/q+E/q^2,\qquad F_q=A_Fq+B_F+D_F/q,
\]
where the constants $B_K,B_F$ will not be needed, and
\begin{align*}
 A_K&=15/19,\\ D_K&=30hC/11-1361/285-(8/15)M,\\
 E&=208/75+(8/15)M,\\
 A_F&=11/68+253/(255C)-(15/19)k-C\nu/(16P(136/225)),\\
 D_F&=1/2-1672/(1275C)-kM.
\end{align*}
In particular $A_F<-619/9690$. Put $J_j=K+jF_q$.
Writing $J_{29}=Aq+B+D/q+E/q^2$ gives $A<0$ and
\[
 PD+2E=2hC+108397/8550-533368/(19125C)
       +(152/225-(319/15)k)M>119/180.
\]
Indeed the four terms exceed, respectively,
$419/180$, $63/5$, $-20/3$, and $-38/5$; for the last use
$k,M<49/80$. Now
$q^3J_{29}'=Aq^3-Dq-2E<0$: if $D\ge0$ this is immediate;
if $D<0$, use $q\le P$ and $PD+2E>0$.

On the diagonal, put $D_\Sigma=15/11+136/(225C)+225C/2176$.
Direct substitution, with no change of logarithms, gives
\begin{align*}
 J_{29}(P)={}&2h(225C/121-1)+27X-208/165-534/121\\
 &+29(3277/11220)-29(3529/(3825C))
       +(152/121-29D_\Sigma)\nu.
\end{align*}
The first six terms exceed $7/5$: termwise lower bounds are
\[377/100,\quad63/50,\quad-127/100,\quad-443/100,\quad846/100,\quad-639/100.\]
Here $419/100<C<21/5$, $h>5/18$ and $X>469/10000$ suffice.
Also $\nu<1552/73695<211/10000$ and $D_\Sigma<397/200<2$.
The last coefficient is negative with magnitude less than $227/4$.
Consequently
\[
 J_{29}(P)>7/5-(227/4)(211/10000)=8103/40000.
\]
For $J_{28}$, $A_{28}<0$ as before. Moreover $M>11/19$,
$kM>1/3$ and $1672/(1275C)>3/10$, so $D_F<-2/15$.
Thus $D_{28}>D_{29}$, and the same argument proves $J_{28}'<0$.
At $q=P$, positivity of both means implies
\[
 F_q<t_q=X+3277/11220-3529/(3825C)<1/8,
\]
since $X<1/21$, $3277/11220<293/1000$ and
$3529/(3825C)>219/1000$. Hence
$J_{28}(q)\ge J_{28}(P)>3103/40000>3/40$.

For completeness, the logarithm bound is elementary. For $x>1$,
\[
 \frac{d}{dx}\left(\frac{x^2-1}{2x}-\log x\right)
 =\frac{(x-1)^2}{2x^2}>0,
\]
so $\log x<(x^2-1)/(2x)$. At $x=C/4$, it suffices that
$21C^2-8C-336<0$. This polynomial is increasing for $C>4$;
using $C<1677/400$ gives the strict upper value
$-67491/160000<0$. Thus $X<1/21$ as used above.

It remains to restore the other row. For any common width set
\[
 H=w-1,\quad y=w^2-1,\quad D_g=C/(16y)-y/C,\quad
 m_1=2L-4/w,\quad m_0=1+4/w-4L.
\]
Direct subtraction of the original rows yields
\[
 F_p-F_q=(p-q)\{A+B(1/p+1/q)\},
\]
\[
 A=\frac3{2(w+1)}-\frac{6H^2}{C(w+1)}-D_gm_1,
 \qquad B=\frac{12wH}{C(w+1)}-D_gm_0.
\]
For example, the mean contribution follows by expanding
$m(z,w)=z^2/w+m_1z+m_0$; the common $X$ cancels exactly.
On the lower boundary $0<D_g<9/20$, $m_1<0$, $m_0<2/3$,
and $12wH/(C(w+1))=152/(85C)>2/5$.
Thus $A>0$ and $B>1/10$. Since $p=P\ge q$,
\[
 S=J_{28}+4(F_p-F_q)\ge J_{28}>3/40.
\]

\subsection{The upper boundary}
Now $w=27/20$ and $5/9\le q\le P$. The positive Cayley series gives
\[
 X>18034/384017> X_0:=469/10000,\qquad
 L>1869280/2180283>L_0:=17147/20000.
\]
The coefficient of $X$ in $S$ is $26$. The mean coefficients are
$\ell_q-4d_{21}-24d_{11}$ and $\ell_p-4d_{22}-24d_{12}$,
where $\ell_z=(2-z)(w-1)/z^2$. They are negative:
$\ell_q\le819/500$, $\ell_p=399/484<1$, whereas each row sum
is at least $210/11$. Both means decrease with their common $L$.
Therefore $S>S_0:=S(X_0,L_0)$.
Write
\[
 \begin{gathered}
 \mu_0=A+Bq+Dq^2,\quad n_0=200767/12150000,\\
 A=72031/135000,\quad B=-337031/270000,\quad D=20/27.
 \end{gathered}
\]
Introduce the following exact rational constants:
\begin{gather*}
 H_0=69517/125000+5/81,\quad d=1/2+3619/(6000C),\\
 E=1/2+55C/987,\quad V=987/(880C),\quad f=375C/3619,\\
 t_1=13/94+748/(705C),\quad W=-1/2-29723/(15510C).
\end{gather*}
For $q\ge q_0:=13/20$ the exact collection of $q^2S_0$ is
$P_+(q)=\alpha+\beta q+\gamma q^2+\delta q^3$, with
\begin{align*}
 \alpha&=21/10+(7/10)A,\\
 \beta&=CQ_b-R_b/C+T_b,\\
 Q_b&=(15/11)(H_0-13/27)-220A/987,\\
 R_b&=49896/1175+3619A/250,\\
 T_b&=179/20+(7/20)(2B-A)-14A,\\
 \gamma&=CQ_g+R_g/C+T_g,\\
 Q_g&=100/99-220B/987,\qquad R_g=6468/1175-3619B/250,\\
 T_g&=26X_0-H_0-104/165-783/484+13/27-85678/7755\\
 &\hspace{1cm}+14/27-(287/20)B-(8841/484)n_0,\\
 \delta&=4W+24t_1-(4E+24d)D-(4V+24f)n_0.
\end{align*}
These are explicit rational expressions, so the following enclosing
intervals are checked by multiplication by positive denominators:
\begin{gather*}
 13/200<Q_b<17/250,\quad2509/50<R_b<251/5,\\
 258/625<T_b<8557/20000,\\
 161/125<Q_g<129/100,\quad2357/100<R_g<11793/500,\\
 59159/10000<T_g<59163/10000,\\
 6435/10000<d<644/1000,\quad7334/10000<E<7341/10000,\\
 267/1000<V<268/1000,\quad434/1000<f<436/1000,\\
 3909/10000<t_1<3916/10000,\quad-958/1000<W<-955/1000,\\
 165/10000<n_0<166/10000.
\end{gather*}
Together with $419/100<C<21/5$ they imply
\[
 -113/10<\beta<-56/5,\quad169/10<\gamma<17,\quad
 -83/10<\delta<-41/5.
\]
Thus $P_+'$ is concave, and its endpoint lower bounds are
$P_+'(13/20)>599/4000$ and $P_+'(11/15)>12/125$.
It follows that $P_+$ increases.

For $q\le q_0$ the negative-geometry collection is a cubic $P_-$
whose coefficients differ from those of $P_+$ by
\[
 (273/200,\ 65C/99-3/20,\ -94/27-100C/99,\ 20/27).
\]
In particular $\beta_-<-8$, $\gamma_-<10$, $\delta_-<-7$.
Hence $P_-''<20-42q\le-10/3$ and
$P_-'(5/9)<-91/27$. Thus $P_-$ decreases. The two cubics
agree at $q_0$, so both sectors reduce to that common value.

At $q_0$, $\mu_0=7031/200000<9/250$ and $n_0<1/60$.
Their negative coefficients allow these two upper bounds to give a
further strict lower comparison $\bar S=A_CC-B_C/C+D_0$, where
\begin{align*}
 A_C&=300H_0/143-(4400/12831)(9/250)-(5850/3619)(1/60),\\
 B_C&=4(386399/310200)+24(165253/91650)\\
 &\hspace{1cm}+(7238/325)(9/250)+(12831/4400)(1/60),\\
 D_0&=26X_0+26(13031/17160)-466/47-H_0-213/169-783/484\\
 &\hspace{1cm}-(3451/169)(9/250)-(8841/484)(1/60).
\end{align*}
Positive-denominator multiplication gives
\begin{gather*}
 A_C>2501/2000>5/4,\qquad
 0<B_C<491094/10000<4911/100,\\
 D_0>32547/5000>13/2.
\end{gather*}
Consequently
\[
 S_0(q_0)>\bar S>419/80-4911/419+13/2=561/33520>0.
\]
Dividing the positive sector cubics by $q^2$ proves $S>S_0>0$
on the whole upper boundary.

\subsection{The moving seam}
Here $q=2-w$, $19/15\le w\le27/20$ and $b=19/15$.
All derivatives in this subsection move both $q$ and $w$.
Put $H=w-1$, $y=w^2-1$,
$\mu=1+w-2wL$ and $\nu=1+b^2/w-2bL$.
Exact collection gives $S=R-A_\mu\mu-A_\nu\nu$, where
\begin{align*}
 A_\mu&=[14+P(C/(4y)+24y/C)]/q-wH/q^2,\\
 A_\nu&=14/P+(q/P)(4y/C+3C/(2y))-bH/P^2,\\
 R&=(C/(Pq)-1)(2h+1-b/w)+26X-22\\
 &\quad+13(1/P+1/q-P-q)+(36-6b)/(w+1)+28Pq/C\\
 &\quad+12qH[-2b+b^2/(w+1)]/(CP)
       -72PwH(w+2)/(Cq(w+1))\\
 &\quad+(3w^2-6)/q^2+(3bw-6)/P^2.
\end{align*}
Partial fractions give the exact derivative
\[
 R''=36/q^4+Z/q^3+2b[1-C/(2P)]/w^3+Y/(w+1)^3+48b/(CP),
\]
where $Z=4Ch/P+C+2-384P/C<-109/2$ and
$Y=72-12b-96P/C-144b^2/(CP)<-104/3$.
The first two terms together are less than $13/4$:
$(36-(109/2)q)/q^4$ decreases on this interval and at $q=13/20$
is $92000/28561<13/4$. The remaining three terms are less
than $-7/4$, $-8/3$ and $99/5$, respectively. Thus $R''<559/30$.

The density representation of the seam mean and the positive
artanh series give
\[
 \frac{H^2}{3w^2}<\mu<\frac{H^2}3,\qquad
 \mu'=\frac{H^2-\mu}{wH}>\frac8{57},\qquad \mu>\frac7{500},
\]
\[
 L<\frac2{w+1}+\frac{H^2}{6w(w+1)},\qquad
 \mu''=\frac{2(1+1/w-2L)}{H^2}>\frac4{3w(w+1)}>\frac{21}{50}.
\]
To bound the coefficient derivatives use
\[
 \frac1{qy}=\frac1{3q}+\frac1{2H}-\frac1{6(w+1)},
 \qquad \frac yq=-w-2+\frac3q.
\]
These give $A_\mu>45/2$, $A_\mu'>32$, $A_\mu''>24$.
For explicit sufficient comparisons, its four defining terms are
bounded below by $1909/100$, $6/5$, $69/20$, and $-9/8$.
For the derivative, the terms from $14/q$ and $24Py/(Cq)$
exceed $26$ and $19$, the negative derivative from $PC/(8H)$
has magnitude less than $11/2$, and $(wH/q^2)'<15/2$;
the omitted derivative terms are positive.
For the second derivative the two $y$ terms are positive,
$(14/q)''>70$ and $(wH/q^2)''<5+25+16=46$.
Therefore
\[
 (A_\mu\mu)''>24(7/500)+64(8/57)+(45/2)(21/50)
               =42/125+512/57+189/20.
\]

The other profile satisfies $\nu>0$, $\nu'<0$, $\nu''>0$.
The last sign follows directly from
$\nu''=2b[b/w^3-L'']$, $L''\le164/375<1/2<b/w^3$;
the first derivative is the true-profile width sign retained in
Section~\ref{six:sec:first-contact}.
We also have $A_\nu>0$, $A_\nu'<0$, $A_\nu''>0$.
Indeed use $qy=-w^3+2w^2+w-2$ and
$q/y=1/(2H)-3/(2(w+1))$. The first derivative of $4qy/(CP)$
is less than $2$, that of $q/y$ is less than $-3$,
and its coefficient $3C/(2P)>8$; the remaining linear derivative
is negative. For the second derivative,
$4(-6w+4)/(CP)>-6$ and $1/H^3-3/(w+1)^3>22$.
Consequently $(A_\nu\nu)''>0$ term by term. Combining yields
\[
 S''<559/30-(42/125+512/57+189/20)=-3851/28500<0.
\]
Its two endpoints were proved positive above. This proves the third
boundary assertion, and hence all three scalar signs.
\space
\section{From a physical returned leaf to the closed endpoint cover}
\label{six:app:physical-cap}
This appendix supplies the physical inequalities used to enter the
closed endpoint cover in H2. Sections~1--2 of the reader map in
\path{h2-h3-analytic-sources/README.md} give the proofs for the individual
pieces and their domain joins. The hypotheses here are
those of H1: one genuine complete normal maximizing lift, its actual
positive middles and subsequent leaves, and the original primitive,
transfer, mean and cost identities. It does not construct a physical
leaf from a relaxed parameter tuple. In this appendix put
\[
 C_0=\frac{196517}{46875},\quad T=\frac{11}{15},\quad
 b_0=\frac{3953}{62500}.
\]
The parameters here come from the actual shared mode-zero arc,
not from a freely chosen tuple. Write $\xi=1+u$, $v=k_s\xi$,
$R=1+k_s^2$, and let $\xi_0<\xi_1$ be its two endpoints.
The complete cyclic and reflected/time-reversed contexts have states
\[
 (u_R,v_R)=\left(1-\frac4{R\xi},\frac{4k_s}{R\xi}\right),
 \qquad (u_L,v_L)=(1-\xi,k_s\xi).
\]
Set $p=1-s_R=4/(R\xi_0)$, $q=1-s_L=\xi_1$,
$\mu=\mu_L$ and $\nu=\mu_R$, where $s_j$ is the terminal
root of the actual $202$ core in that chart and $\mu_j$ its
true weighted middle mean. The common diagonal bridge constant is
\[
 A_b=\frac R{R-1},\quad
 g(x)=\log x+\frac{x^{-1}-x}{2},\quad
 V_b=\log(R/4)+g(p)+g(q),
\]\[
 B=\frac{1-A_bV_b}{2}+\frac{A_b}{4}
           \left(\frac\mu q+\frac\nu p\right).
\]
In particular $B$ here is not the number of arcs. The same actual
matrix satisfies
$\Gamma=A_b^2(3+\mu)(3+\nu)/16-B^2$.
These are the definitions and matrix identity of Section~\ref{pf:F02};
reflection transports both covectors and all matrix entries with
reversed time before this identity is used.

The non-circular initial scalar reduction uses the once-counted budget
$U<2P$ and gives
\begin{equation}\label{six:eq:initial-cap-box}
\begin{gathered}
 7<R<54/5<13,\quad 14/25<p,q<9/10,\quad Rpq>C_0,\quad B>b_0,\\
 0<\mu<(1-q)^2,\quad0<\nu<(1-p)^2.
\end{gathered}
\end{equation}
For clarity, the order of that reduction is: the signed strip and
same-lift harmonic comparison exclude $4<R\le7$; compactification
and cost give $R<17$ and $B>0$; the returned $379/1000$ bounds give
the exclusive-cost ceiling and shared-cost floor. These imply $R<13$
and $Rpq>C_0$, respectively. Alternating the exclusive-cost and
bridge comparisons gives $p,q>8/15$, $R<11$, $p,q>5/9$,
$R<54/5$ and $p,q>14/25$. The actual cap condition proved below,
with signed means, excludes $p\ge9/10$ or $q\ge9/10$ without using
these last lower endpoint bounds. The source locators and exact
endpoint reserves are in the accompanying initial-domain supplement.
These imported scalar comparisons remain part of H2; the arguments
below do not replace their geometric hypotheses.

\subsection{A later actual zero is sufficient}
\begin{lemma}\label{six:lem:later-zero}
Suppose $d>0$, $R_0=1+d$, $0<e_0<e_1$, and
\[
 F(e)=ae+b-2R_0\log e-4/e.
\]
If $F(e_0)\ge0$ and $F(e_1)=0$, then
\[
 \Delta=R_0^2-4a>0,\qquad
 ae_0-R_0<0\quad\hbox{or}\quad
 ae_0^2-2R_0e_0+4<0.
\]
\end{lemma}
\begin{proof}
Write $Q(e)=e^2F'(e)=ae^2-2R_0e+4$.
If $F'\ge0$ on $(e_0,e_1)$, its continuous numerator, a
nonzero polynomial, is positive on a nonempty subinterval. Integration
would give $F(e_1)>F(e_0)$, a contradiction. Thus $Q(t)<0$ for
some $e_0<t<e_1$. If $a\le0$, both asserted conclusions follow
immediately. If $a>0$, $Q$ has two distinct positive roots
$t_-<t_+$, since their product and sum are $4/a$ and $2R_0/a$.
In particular $\Delta>0$ and $e_0<t_+$. If $e_0<R_0/a$ the
first alternative holds. Otherwise $e_0$ lies between the vertex
and $t_+$ and $Q(e_0)<0$. This also treats the vertex itself.
\end{proof}
On the physical mode-zero ray, $e$ increases strictly with time and
$F(e)=2d^2e\,g(e)$. Actual maximization gives $F(e_0)\ge0$;
the specified later mode-one tie gives an actual $e_1>e_0$ with
$F(e_1)=0$. Thus Lemma~\ref{six:lem:later-zero} applies without
assuming that this is the first zero, that any zero is simple, or
that the derivative at the actual exit is negative. An earlier
contact or a triple tie does not invalidate the argument.

\subsection{Positive clearing and the precise cap predicate}
Set
\[
 t=R-1,\quad K=(3+\mu)/2,\quad
 w=\frac{2-q}{1-r}>1,\quad x=2w-2+q,\quad
 y=w^2-1,\quad S=x^2+tq^2.
\]
Here $r$ is the actual middle's outgoing coordinate. The original
same-bridge transfer reads
\[
 d=\frac{tq^2}{x^2},\quad e_0=\frac{x}{w},\quad
 \zeta=-\frac{2y}{tq^2},\qquad
 M_{\rm new}=-\frac{R^2(3+\nu)}{8t^2}+2B\zeta-K\zeta^2.
\]
Consequently, with $a=(1+d)^2+2d^2M_{\rm new}$, define
\[
 A=4S^2-R^2(3+\nu)q^4-32Ky^2-32Btq^2y,
 \qquad a=\frac{A}{4x^4}.
\]
Direct substitution gives
\[
 \Delta=\frac{S^2-A}{x^4},\quad
 ae_0-R_0=\frac{A-4wxS}{4wx^3},\quad
 Q(e_0)=\frac{A-8wxS+16w^2x^2}{4w^2x^2}.
\]
All denominators are positive. Put
\[
 N_D=S^2-A,\quad N_E=A-4wxS,\quad
 N_Q=A-8wxS+16w^2x^2,\quad Z=S-4wx.
\]
The exact necessary condition is $N_D>0$ and
$(N_E<0\ \hbox{or}\ N_Q<0)$. Since
$N_E=SZ-N_D$, $N_Q=Z^2-N_D$, and $S>Z$,
this is equivalent to
\begin{equation}\label{six:eq:strict-physical-cap}
 N_D-\max(0,Z)^2>0.
\end{equation}
Indeed, for $Z\le0$ the first alternative is automatic; for $Z>0$
the union of the two alternatives is $N_D>Z^2$, because $Z^2<SZ$.
No division by $Z$ occurs, and its zero seam is retained.

\subsection{The signed mean bound has its own physical proof}
Write $s=1-q$. On the actual interval $[-s,r]$ the density is
$\omega(u)=(1-u)^{-2}$ and $\mu$ is the weighted mean of $u^2$.
The strict outgoing jet gives $\mu>r^2$, while support gives
$\mu<s^2$; hence $|r|<s$. If $r\ge s/2$, then
$u^2-r^2$ changes from positive to negative at $-r$ on this
interval. The weight is strictly increasing, so
\[
 \int_{-s}^r(u^2-r^2)\omega(u)\,du
 <\omega(-r)\int_{-s}^r(u^2-r^2)\,du
 =\omega(-r)\frac{(s-2r)(s+r)^2}{3}\le0,
\]
contradicting $\mu>r^2$. Therefore $r<s/2$. Put
$m=(s^2-sr+r^2)/3$. The identities
\[
 m-r^2=(s-2r)(s+r)/3>0,\qquad
 s^2-m=(2s-r)(s+r)/3>0
\]
show that $u^2-m$ has precisely one positive-to-negative sign
change on $[-s,r]$, including when $r<0$. Its unweighted integral
is zero. The same weight comparison proves
\[
 \mu<\mu_S:=\frac{s^2-sr+r^2}{3},\qquad
 0<\eta:=\frac{(w-1)(1+q)}{3(1-q)}<1.
\]
Here $\eta$ is a profile coordinate, not a multiplier normalization.
Closed comparison boxes may include $\eta=0,1$ without declaring
those artificial boundary points physical positive middles.

\subsection{A product gap supplies the saved quantitative collar}
\begin{lemma}\label{six:lem:quantitative-collar}
Under \eqref{six:eq:initial-cap-box}, the physical transfer and strict
order condition above imply
\[
 \frac{1}{50000000}<w-1\le\frac{11}{13}.
\]
\end{lemma}
\begin{proof}
The actual bridge identity and strict mean supports give
\[
 B\le\frac12-\frac{R}{2(R-1)}
       \{\log(Rpq/4)+2-p-q\}.
\]
For $v=pq>4/R$, $p+q\le1+v$ and the function
$\log(Rv/4)+1-v$ is strictly increasing for $0<v<1$.
Thus $0<B<3/[2(R-1)]<1/4$.
Also $q>4/13$, and $|r|<s$ gives $w<(2-q)/q<11/2$.

Along the formal segment $1\le z\le w$, hold the actual
$R,q,\mu,\nu,B$ fixed and replace $w$ by $z$ in the polynomials
above. Then
\[
 t<12,\quad K<2,\quad x<10,\quad
 0\le y<117/4,\quad S<112,\quad z<11/2.
\]
Differentiation, with primes denoting $z$ derivatives, gives
\begin{align*}
 A'&=32Sx-128Kzy-64Btq^2z,\\
 N_E'&=A'-4(xS+2zS+4zx^2),\\
 N_Q'&=A'-8(xS+2zS+4zx^2)+32zx^2+64z^2x.
\end{align*}
Using $q<9/10$, the termwise absolute bounds are
$|N_E'|<484608/5<100000$ and
$|N_Q'|<760448/5<160000$.
These are polynomial bounds, so the formal endpoint $z=1$
requires neither a mean quotient nor physical realizability.

Let $\gamma=Rpq-4$ and $d_R=(1-p)^2-\nu>0$.
Evaluation at $z=1$, using the actual bridge identity, gives
\begin{align*}
 N_E(1)&=Rq^3\{4(1-p)+(2-p)\gamma+Rq d_R\}
             >\frac{448\gamma}{2197},\\
 N_Q(1)&=q^2\gamma\{\gamma+2Rq(1-p)\}+R^2q^4d_R
             >\frac{16\gamma^2}{169}.
\end{align*}
At least one of $N_E(w),N_Q(w)$ is strictly negative. The
fundamental theorem of calculus therefore gives respectively
\[
 w-1>\frac{448\gamma}{2197\cdot100000}
 \quad\hbox{or}\quad
 w-1>\frac{\gamma^2}{1690000}.
\]
Now $\gamma>C_0-4=36068/187500>19/100$.
The first bound exceeds $19/50000000$ since $448/2197>1/5$;
the second exceeds
\[
 \frac{361}{16900000000}
 =\frac{1}{50000000}+\frac{23}{16900000000}.
\]
Both prove the strict lower collar. Finally,
$w-1=3\eta(1-q)/(1+q)\le11/13$ follows from
$q\ge14/25$ and $\eta\le1$.
\end{proof}
This supplies the actual qualification for the saved candidate
intersection. It uses neither a qualitative compactness collar nor
an exterior-pair Casimir floor. In particular it is proved after the
initial scalar reduction and is not fed back into that reduction.

\subsection{How the closed cover uses these predicates}
After inserting the actual bridge, the coefficients of $\mu,\nu$
in $N_D$ are respectively $16y^2+8Rqy>0$ and
$R^2q^4+8Rq^2y/p>0$. Raising the opposite mean to $(1-p)^2$
therefore increases $N_D$. The derivative in $p$ of this upper
expression is
\[
 -2R^2q^4(1-p)-16Rq^2y(1-p)/p<0.
\]
Lowering $p$ to $C_0/(Rq)$ is thus a necessary upper comparison,
not a construction of a new physical state. Denote the resulting
cap by $H_3=N_{D,3}-\max(0,Z)^2$. The original radius derivative
identity, with $L_0=\log(C_0/4)>0$, is
\begin{align*}
 \partial_R(N_{D,3}-Z^2)
 &=\frac{8q(q-2)(w-1)}{3w^2}\mathcal P
   -2(C_0-4)q^3-16q^2yL_0+8qy(\mu-\mu_S),\\
 \mathcal P&=4qw^2(w-1)+q(2w+1)+(w+2)(w-1)(w+1)>0.
\end{align*}
It is strictly negative also at $w=1$. On $Z\le0$, the derivative
of $H_3$ adds $2q^2Z\le0$; the two formulas agree at $Z=0$.
Thus transport to a valid lower radius face is safe, provided the
own upper mean is at most $\mu_S$. All cost-face denominators
must be positive on the whole comparison rectangle.

There are exactly three initial rejection types: an upper radius
strictly below its lower face; an upper bound $N_{D,3}\le0$;
or an upper bound $N_{D,3}-Z^2\le0$ together with a
\emph{same-box lower bound} $Z\ge0$. The last guard cannot be
replaced by a point value. Each contradicts the strict actual cap.
The twelve closed root rectangles have $q$ cuts
$14/25,2/3,4/5,9/10$ and profile cuts $0,1/4,1/2,3/4,1$.
Each split replaces a rectangle by its two closed midpoint children;
induction preserves the exact union and every seam. Sixteen terminal
rectangles are rejected and twenty-two retained. These counts alone
prove no sign. The closed-cover proof map in the analytic supplement
gives each rejection's polynomial, radius face and guards. The tests
of Appendix~\ref{six:app:computation} verify their finite inequalities.

The non-circular next steps are
\[
\begin{gathered}
 p,q\le17/20\ \Longrightarrow\ \mu,\nu\ge333/73600
 \ \Longrightarrow\ p,q\le13/16,\\
 \mu,\nu\ge261/38912\ \Longrightarrow\
 \text{diagonal and A/B/C comparisons}\ \Longrightarrow\ p,q\le4/5.
\end{gathered}
\]
The seven diagonal pieces leave two pieces; they do not exclude all
nine. A/B/C respectively use the true-rising tail, a low-radius
harmonic comparison, and the weak-profile cap. Only after their
closed union is established is the mean floor $3/400$ available.

For $q\in[T,4/5]$, the complete profile ownership is
\[
\begin{array}{c|l}
 {}[0,1/8],[1/8,1/4],[1/4,3/8]&
 \text{early closed cap owners }r1010,r1011,r1110\\
 {}[3/8,1/2]&r1111\text{ product faces and tail}\\
 {}[1/2,5/8]&r1210\text{ tail, eight quarters and N1}\\
 {}[5/8,3/4]&r1211\text{ complete strong/weak strip}\\
 {}[3/4,7/8]&r1310\text{ strong/weak delay and own budget}\\
 {}[7/8,1]&r1311\text{ rising tail, three delays and own budget}.
\end{array}
\]
For the first three subintervals, grouped in the first row, the named rectangles cover the low-profile
part meeting this endpoint range; other early rectangles cover the
higher-$q$ portions. The supplementary closed-cover map specifies
all exact coordinate intersections. In particular all rows own their
$q=T$ faces, without requiring $q>T$.

Here the strong-companion branch is $p>T$, not $p\ge T$.
Writing $R_b=C_0/(Tq)$, the weak branch $p\le T$ has $R>R_b$
by strict product. It must use its separate weak tail/delay and
cannot borrow a reflected high-profile floor. For $r1111$, the
low-radius faces join at $Q=5C_0/28$, and the all-companion weaker $e_0$ comparison
face C owns $R=R_b$. For $r1210$, four closed radius quarters on each side
of $Q$ cover the entire lower-radius interval; seven are closed
by the joint comparison and the eighth by N1. For $r1211$,
the old $[T,Q]$ and new $[Q,3/4]$ small-endpoint comparisons
are distinct inputs. On $[3/4,4/5]$, the weak, high-radius,
strong-tail and harmonic owners precede the remaining
H1260/Fsmall conjunction. The later $r1310$ and $r1311$
strong branches use only already established reflected strips;
their weak branches use product bounds without those floors.
Their delay and budget comparisons are separate obligations.

Thus the stated source comparisons imply that
the closed owners exclude every actual $q\ge T$ in the initial
box. Whole reflection on the same lift excludes $p\ge T$ as well.
In particular, the closed conclusion $p,q\le T$ used in the main
text holds. The formal closed S24 domain is not trimmed
on the strength of this observation.

\subsection{Locality of the endpoint comparison}
The initial-domain, cap-predicate and closed-cover comparisons have
separate source, sign and covering obligations. The sixteen initial
rejection signs use fifteen arrays with $948$ exact coefficients;
the closed-owner plan supplies their union. H1260, high radius and
the highest profile retain their own domains and analytic joins, as
specified in Appendix~\ref{six:app:computation}. These comparisons
must be composed with those qualifications.

Throughout, the local nine-source argument assumes $U<2P$ for
contradiction, but retains the original global minimizer and its
same-lift harmonic and compactness qualifications. It does not use
an exterior-pair Casimir floor or impose periodic state-covector
data on an isolated source. The resulting nine-source estimate
remains $U\ge2P$. The endpoint cover does not itself close the
remaining low-endpoint region: the guarded common-log and
first-contact argument is still needed. Appendix~\ref{six:app:reduction-proofs} completes the H2 reduction and
H3 transfer within H1; this endpoint cover alone does not prove them.
\space
\section{Finite comparisons for returned sources}
\label{six:app:computation}
The returned-source argument uses exact polynomial signs and closed
endpoint covers after its physical source reductions. We list their
mathematical domains and joins here. The accompanying dependency index
binds each comparison to its constructor and checking route.

\subsection{Closed parts of the endpoint comparison}
\label{six:app:owner-implementation}
Let \(T_0=11/15\), \(Q_0=196517/262500\) and
\(C_0=196517/46875\). The following branches precede the
analytic continuation in Appendix~\ref{six:app:reduction-proofs}.
The task entries \path{H2-01} through \path{H2-21} and
\path{H3-01} through \path{H3-05} supply the closed covers, exact
polynomials, leaf predicates and coefficients for these comparisons.
Section~\ref{reader:local} identifies the separate finite inputs
used by the returned-HIGH and forward-pair bounds.

\begin{center}\small
\begin{tabular}{@{}p{.25\textwidth}p{.67\textwidth}@{}}
\toprule Domain or boundary & Condition retained at the join\\\midrule
\(q\in[T_0,Q_0]\) & The distinct mean, logarithm, event, prefix,
 physical-primitive and two-owner comparisons keep their original
 source-qualified branches.\\
\(q\in[Q_0,3/4]\) & The fixed \(34\)-owner cover uses its original
 closed interval. The strong large-\(p\) inequalities retain the
 Sturm sign and true-mean transport. The reflected profile gives
 the highest-\(\eta\) continuation.\\
\(q\in[3/4,4/5]\) & The weak radius starts at
 \(R_b=C_0/(T_0q)\), with \(R=R_b+\rho(8-R_b)\),
 \(\eta\in[5/8,3/4]\), \(\rho\in[0,1]\) and
 \(C_0/(Rq)<p\le T_0\). The strong integration-by-parts test,
 weak \(R\in[8,9]\) test and monotone tail \(R\ge9\) keep their
 separate hypotheses. The interface satisfies \(15C_0/(11T_0)<8\).\\
\bottomrule
\end{tabular}
\end{center}

The same-bridge mean transfer in Appendix~\ref{six:app:bridge-mean}
uses the positive outgoing leaf, its terminal root and the shared
arc under \(U<2P\). These are the data used in the branch joins.

The both-weak \(\kappa\) test and the O/K boundary comparisons retain
their complete closed strips. The symmetric floor128 cover gives
\(G<0\) on all 64 leaves of chart A and rules out 44 leaves of
chart B by the necessary guards. On 19 further B leaves it gives
positive cost at the covariance upper mean. With strict \(G>0\)
and the original profile guards, J27 supplies the same sign on
the last leaf, B111011. Appendix~\ref{six:app:reduction-proofs}
uses these signs in the vertical comparison with the true mean.
The sign on the line \(\eta=973/1000,\ q\in[11/15,4/5]\)
extends upward using \(G_r=-2rI_0<0\), \(r_\eta>0\) and
\(r\ge11/155\).

The H1260 polynomial uses positive homogeneous clearing on
\(q\in[3/4,4/5]\), \(\eta\in[5/8,3/4]\), \(\rho\in[0,1]\).
The exact inverse identities and the positive Bernstein signs are
checked before using \(H=J(r_F)>0\). The next comparison explains
how this test and the actual maximizing leaf give the positive
action used by Fsmall.

For the moving-cost input put $z=5/53$ and
\[
 L_{\rm hi}=2\sum_{j=0}^5\frac{z^{2j+1}}{2j+1}
       +\frac{2z^{13}}{13(1-z^2)}.
\]
Then $\log(29/24)<L_{\rm hi}$. With $a=6/5$, $v=29/24$,
\[
 N_{\rm lo}=1+\frac{a^2}{v}-\frac{2aL_{\rm hi}}{v-1}
 =\frac{78409276064951567888403}{6726771099843633746732825}>0.
\]
This input is bound to the six raw/H1260/Fsmall comparisons.

For the delay comparison use the natural closed box
\[
 q\in[3/4,4/5],\quad\eta\in[5/8,3/4],\quad
 k_b\in[k_{\rm right},C_0/T_0],\qquad
 k_{\rm right}=21q/4+C_0/(4T_0).
\]
Put $w=1+3\eta(1-q)/(1+q)$, $k_b=Rq$,
$x=2w-2+q$, $S=x^2+(R-1)q^2$, $r_0=S/x^2$, $e_0=x/w$, and
\[
 a=2-q,\quad y=w^2-1,\quad z_{\rm cap}=2-aw,\quad
 r_F=\frac{z_{\rm cap}}{2y},\qquad r_{\rm end}=a/x.
\]
Here $q\le4/5$, $k_b<6$ and $w\ge29/24$ give
\[
 S-4wx=q(k_b-4)-4(w^2-1)
 <\frac85-\frac{265}{144}<0.
\]
For the source-bound lower comparison $A$, let
$J_0=A-2r_0+4/e_0$ and $J(r)=J_0+r(A-r_0)$. Then
\[
 H=J(r_F)=(1+r_F)J_0+r_F(r_0-4/e_0)>0,
 \qquad r_F>0.
\]
Since $r_0-4/e_0<0$, this implies $J_0>0$, and affinity gives
$J>0$ throughout $[0,r_F]$. The complete physical leaf,
$F_{\rm in}\ge0$ and source-bound entry inequality still supply
its actual delay. The moving-action proof also gives
$0<r_F<r_{\rm end}$ and uses the positive primitive derivative and
$r_F-\log(1+r_F)>0$. No positivity test is used to invent a leaf.
The resulting qualified H2/H3 implications are assembled in
Appendix~\ref{six:app:reduction-proofs}.
\space
\section{Proof of the qualified reduction and its local transfer}
\label{six:app:reduction-proofs}
We prove Lemmas~\ref{six:lem:H2} and \ref{six:lem:H3} using the
interfaces in H1. The separate \path{h2-h3-analytic-sources}
supplement contains the analytic derivations, and its \path{README.md}
locates the hypotheses, exact formulas and boundary joins used at each
step. The finite signs and closed covers are those specified in
Appendix~\ref{six:app:computation}, checked by the fixed 1118R01
reader archive.

\subsection{The actual tuple and the common logarithm}
Choose the complete jointly normalized lift before restricting it.
The returned contexts give $U=E_L+I+E_R$ and
$T_L=E_L+I,T_R=E_R+I>379/1000$. Under $U<2P$, therefore,
\[
 E_L,E_R<85733/250000,\qquad I>9017/250000.
\]
Appendix~\ref{six:app:physical-cap} proves the physical cap/order and
quantitative collar links. Its initial-domain and closed-cover
lemmas, with the source and sign arguments of
Appendix~\ref{six:app:owner-implementation}, give $5/9<p,q\le11/15$.
The same local budget gives $p+q>5/4$ and $pq>341/900$, as in
Section~3 of the reader map. All strong-companion and weak-companion
branches are retained; equality does not acquire a strict opposite
profile floor.

The actual-leaf argument in Section~3 of the reader map starts from
the complete leaf's maximizing entry and later prescribed zero. The
strict binary jet, the qualified LOW174 physical-floor lemma, the
once-counted short-exit comparison and the kappa64 lemma give the
two strict kappa bounds, $w,u>9/8$, and the endpoint guards. These
are used in that order; neither the physical floor nor the kappa
lemma is justified by a later radius projection.

For the $q$-side write $y=w^2-1$, $x=2w-2+q$,
$X_R=\log(Rpq/4)$ and $M=\mu/q+\nu/p$, where $\mu=m(q,w)$ and
$\nu=m(p,u)$ are the two actual means. Substituting the common
bridge $B=\tfrac12-R(X_R+D_{\rm g}-M/2)/(2(R-1))$ into the
strict terminal discriminant and dividing by $16Rq^2y>0$ gives
\[
 T_{\rm term}=X_R+D_{\rm g}-M/2-1+b_R+a_RR+(c_R+1)/R,
\]
with
\[
 a_R=-\frac{\nu q^2}{16y},\quad b_R=\frac{3(x^2-q^2)}{8y},
 \quad c_R=\frac{3(x^2-q^2)^2-16(3+\mu)y^2}{16q^2y}.
\]
The coefficients $a_R,b_R,c_R$ are constant when the profiles are fixed.
Here $R_0=\mathsf C/(pq)<R_{\rm act}<25/2$. The strict profile
bounds give $c_R<-1$ and $-3/100<a_R<0$ on the fixed-profile
segment, so $T_{\rm term}'=1/R+a_R-(c_R+1)/R^2>1/20$.
The strict negative terminal row at the actual radius therefore stays
negative at $R_0$. The reflected comparison uses the other exact
width and mean. The identity $x^2-q^2=4ye_q$ in the $q$-row
gives the coefficients of Section~\ref{six:sec:first-contact} after
putting $u=w$, including the term $pq/\mathsf C$.

For the action, the actual once-counted cost and the same bridge
give a strict inequality after lowering the two physical exit
widths to their nonnegative geometric lower bounds. At the fixed
profiles, let $g_q=\max(0,1-(2-q)/w)$ and
$g_p=\max(0,1-(2-p)/u)$. The resulting action $K_{\rm exact}$
satisfies
\[
 \partial_R K_{\rm exact}
 =2\left(1-u_*+\frac{g_q+g_p}{2}\right)-\frac2R>
 \frac{89}{350},
\]
since $R\ge R_0>7$ and $1-u_*>27/100$. This gives the action
$K<0$ at the same $\mathsf X=\log(\mathsf C/4)$ as the two strict
terminal rows. Section~3 of the reader map gives the actual width
comparison and both radius projections.

For side $j\in\{p,q\}$ let $M_j$ denote the negative of the normalized
terminal complement multiplied by its positive row scale $s_j$, where
$s_q=16\mathsf Cpq(w^2-1)$ and $s_p=16\mathsf Cpq(u^2-1)$.
Write $m_j(w,u)=M_j(\mathsf X)/s_j$ for the normalized terminal margin
at the common logarithm, with the $q$-width listed first.
The true profile means in these rows are $m(q,w)$ and $m(p,u)$.
The rational logarithm bounds satisfy $X_-<\mathsf X<X_+$.
At the same tuple,
\[
 M_j(X_-)=M_j(\mathsf X)+s_j(\mathsf X-X_-),\qquad
 K(X_+)=K(\mathsf X)-2(X_+-\mathsf X).
\]
Thus a weak terminal inequality $M_j(\mathsf X)\ge0$ supplies the
strict terminal inequality needed by the earlier comparison at $X_-$,
and the action stays negative at $X_+$. This implication is used
pointwise, with the applicable mean guards.

\subsection{Guarded continuation and a new comparison root}
The finite comparison domain is larger than the physical domain.
Only its listed guards are continued; no new physical exit or
frame closure is asserted along the path. Section~4 of the reader
map locates the proofs of preservation
of the density covariance bound, true mean branch, M9 and mean
floors, weak kappa width floor, paired lower-profile disjunction
and sum, and the action and row-independent filters. An own-row
upper filter is required only where that row is nonnegative, and
is restored at its zero. At endpoint $T=11/15$, the closed
all-companion strips give the required profile \emph{upper} bound;
whereas the strict-companion comparison at $p>T$ supplies a lower
bound.
With the larger width $u$ decreasing, the fixed-side row improves
as its opposite true mean increases, while the action decreases.
The initial connected admissible interval can end at width equality
($O$), an active kappa equality ($K$), or the zero of the moving
terminal comparison ($T$). The two geometry formulas agree at
their common boundary. For a tied endpoint each active equality is
retained. In particular the sign required in the opposite row is
strict at the lower end of the interval.

The direct O boundary lemma derives its scalar discriminant sign
from the action, balances the two sides, and relays the guards to
the closed midpoint cover. The strict terminal rows at $X_-$ give
$G>0$ at the covariance upper mean $\mu_{\rm cov}$. For B111011,
the adjacent leaf B111010 gives $G<0$ on its lower face.
The J27 bound uses the positive distance from that face to the first
zero of $G$. Together floor128 and J27 give
$k(w,\mu_{\rm cov}(w))>0$ on the guarded part of chart B.
Here $k(w,\lambda)$ is $q^2/2$ times the symmetric action at $X_+$
evaluated with mean $\lambda$.

Fix $q$. At $X_-$ let $\tau(w)$ be the common terminal right-hand
side divided by the sum of the mean coefficients,
\[
 d_s=\frac1q+\frac{w^2-1}{\mathsf C}
               +\frac{\mathsf C}{16(w^2-1)}>0.
\]
The sign of $G$ is the sign of $\mu_{\rm cov}-\tau$.
The lower line of the residual B strip, at chart height $5/8$,
has $G<0$. The remaining cover guards persist as $w$ decreases to
this line. The component of $G>0$ containing the comparison width
has a left endpoint $w_0$
with $G(w_0)=0$. By continuity along this component,
$k(w_0,\tau(w_0))=k(w_0,\mu_{\rm cov}(w_0))\ge0$.
On the component, the terminal expression gives
$-3/4<\tau<3/64$ and $\tau'>-153/64$. With $a=2-q$ the action
derivative is
\[
 \frac{d}{dw}k(w,\tau(w))=
 a\left(\frac{\mathsf C-q^2}{w^2}+3+\tau\right)+a(w-1)\tau'
 >\frac{3305}{1280}.
\]
Thus
$k(w,\tau(w))>0$ at the comparison point. Adding the two mirrored
strict terminal rows gives an average of the two means greater
than $\tau(w)$. The coefficient of this mean in $k$ is
$(2-q)(w-1)>0$, which makes the action comparison at $X_+$ positive,
contradicting its required negative sign.

The K and KK boundary lemmas use the strict $G<0$ sign on the
entire A chart, with separate exact-kappa and fixed-width paths.
Section~4 of the reader map locates the original zero-delay and
vertical completion proofs as well as the K and KK comparisons.

Order the widths as $u\ge w$. The reversed-endpoint lemma
excludes $p\le q$ by lowering $u$ to $w$, temporarily omitting
the own row. At equality of widths, exact row dominance restores
that row and its filter, giving the already excluded O boundary.
For $p=q$ the two rows coincide; a simultaneous kappa tie is
covered by K/KK. Thus $p>q$ and $u>w$ remain.

Put $F(u)=-m_p(w,u)$ with the other profile fixed.
The normalized-row lemma gives $F''<0$
and the common-X change adds only a constant to the old normalized
row. Put $l=\max(w,u_\kappa(p))$. If $F(l)\le0$, the weak own
row and strict other row/action restore the old-domain subset
and contradict O or K. Hence $F(l)>0$, whereas $F(u_0)<0$
at the initial point. Continuity and strict concavity give a new
root $u_T\in(l,u_0)$ with
\[
 F'(u_T)\le\frac{F(u_T)-F(l)}{u_T-l}<0.
\]
At this root $m_p(w,u_T)=0$ and the true mean $m(p,u_T)$ is positive.
The opposite row and action retain strict margins, and both kappa
bounds remain strict because $u_T>l$. The physical leaf's prescribed
endpoint was used earlier to give the actual terminal sign.

\subsection{The two faces and the exact scalar conclusion}
The complete continuous-kappa-cap lemma excludes $u_T\le4/3$:
the coordinate and face concavity arguments reduce its overcover
to the top, diagonal and sum faces, each with a strict rational
reserve. On the sum face the second derivative is
$F_{pp}+F_{qq}-2F_{pq}<0$, using the separately proved
$F_{pq}>0$; coordinate concavity alone would not suffice. The equality
$u_T=4/3$ uses the original active row and strict mean caps, not a
fictitious positive extension length. The nonpositive-$p$-profile
lemma then excludes $4/3<u_T\le2-p$, including equality.
The actual jet gives $u_T<22/13$. The smaller-width action lemma
gives $w<27/20$. Section~5 of the reader map locates the proofs of
these comparisons and the closed cap faces.

Put $v_0=\max(w,2-p)$ and vary only the own true width $v$ on
$[v_0,u_T]$. Write $K(v)$ for the action at this fixed $q$-profile.
The exact identities give
\[
 K'(v)>55/6,\qquad
 \frac{d}{dv}m_p(w,v)<9/4,\qquad
 \frac{d}{dv}\bigl(K(v)-4m_p(w,v)\bigr)>1/6.
\]
At $v_0=2-p$ the derivative is from the positive-geometry side.
Since $m_p(w,u_T)=0$ and $K(u_T)<0$, integration yields
\[
 K(v_0)-4m_p(w,v_0)<K(u_T)-(u_T-v_0)/6<0.
\]
The other margin satisfies $m_q(w,v_0)>m_q(w,u_T)>0$.
To exclude the
lower geometric face let $b=2-p$ and suppose $w\le b$. Since
$p>q$ and $p+q>5/4$, we have $p>5/8$. At $v_0=b$
the transported $q$-row is $m_q(w,b)>0$. Keep the positive
opposite mean $m(p,b)$ fixed and put $F_b(t)=-m_q(t,b)$ for
$9/8\le t\le b$. Since $q<p$ the moving $q$-profile has
negative support throughout this interval. Its normalized row
is strictly concave. The two endpoint comparisons give
$F_b(9/8)>9/400$ and $F_b(b)>71/40000$, as proved in the
lower-face sources in Section~5 of the reader map. The chord
bound gives $F_b(w)>0$, a contradiction which includes $w=b$.

Thus $v_0=w>2-p$. At this equal-width tuple the row conventions are
$F_p=-m_p(w,w)$ and $F_q=-m_q(w,w)$.
The surviving domain lies in the closed enlargement
\eqref{six:eq:S-strip}, and the inequalities are
$K+4F_p<0$, $F_q<0$. Both concern this one tuple. This
proves Lemma~\ref{six:lem:H2}. The scalar proof in
Section~\ref{six:sec:first-contact} then contradicts this conclusion.

\subsection{Why the whole chain restricts to a proper nine-source}
For $012021201$, the two actual contexts are positions
$(1,2,3,4,5,6)$ and $(7,6,5,4,3,2)$, the latter under full
reflection $j\mapsto1-j$. The complete internal forward witnesses
$[0,3]$ and $[5,8]$ provide the two entrances. The shared arc and
its endpoint ties define the same bridge and means above on the
selected lift. The own middle is position 5, its outgoing leaf
position 6, and its prescribed tie is 6/7. On the reflected side
they are positions 3, 2 and the reversed 1/2 tie. Thus neither
complete leaf nor its endpoint control lies outside the source.
Section~\ref{pf:M075} and Lemma~\ref{ret:structure} give these
physical identities independently of any numerical union conclusion.
The union bound uses the full first-contact argument here, with
these local structural identities as its inputs.

For contradiction suppose its once-counted union costs $U<2P$.
The two returned costs give the same exclusive ceilings and shared
floor used above. The signed strip gives $\Gamma>3/4-U^2>0$.
The positive binary invariant puts this determinant on the future
sheet, so the same-lift harmonic comparison and original compactness
bounds remain available on the genuine containing interval. In every
endpoint comparison the budget is this local $U$. The upper
determinant bound still comes from the original minimizing period
on the same lift. The radius and product bounds then use the
local exclusive ceilings, the shared floor and the two profile
means. The endpoint cover uses those bounds, the signed cap and
positive leaf lengths. The strict leaf guards and both radial
comparisons use the same two complete leaves and the once-counted
budget. Once the common logarithm and guards are obtained, the
continuation, endpoint comparisons and the two faces depend only
on the resulting scalar inequalities. Section~6 of the reader map
traces these dependencies, including the strict, weak and
coincident cases in the endpoint cover. In particular none uses
the exterior Casimir floor $147/250$, exterior forward payments,
period-eleven closure or periodic covectors. This proves
Lemma~\ref{six:lem:H3} under H1 and hence the weak bound $U\ge2P$.
\space

\space
\label{revision:appendix-end-4}
\end{proofappendix}
\section{Reverse pairs and crossed cores}\label{app:crossed}
\begin{proofappendix}
\section{A source-qualified reverse-pair payment}\label{opt:pair650-section}
The additional information needed here is a pair of strict prefix
forcings on the same full lift. Positivity of two scalar cap responses
alone is insufficient.

\begin{proposition}\label{late:r08:cubic}
Consider an actual capped reverse pair $1|02|1$, with a normal full
maximizing lift normalized to $\eta=1$. At its middle let
\[
 p,q,r>0,\quad pq+pr+qr=1,\quad q\ge p,r,\qquad
 F=e^{t_0/p},\quad G=e^{t_2/r},\quad t_0,t_2>0.
\]
Let $x$ be the sum of the physical costs of the two complete interior
arcs $0,2$, and let $T=t_0+t_2$. In the
bilateral-cap notation put
\[
 a=(p^{-1}-p)/2,\quad b=(r^{-1}-r)/2,\quad
 f(Y)=\frac{Y\log Y}{Y-1}-1,
\]
\begin{align*}
 \beta_p&=p\left[\frac{p(F-1)}{p+q}-f(F)\right],\\
 \beta_r&=r\left[\frac{r(G-1)}{r+q}-f(G)\right],\\
 A&=(p+p^{-1})\beta_p,\qquad D=(r+r^{-1})\beta_r.
\end{align*}
Assume the actual cap rows and strict prefix forcings are
\begin{align}
 V_x+aV_y+\frac{h}{p(F-1)}&=A,&
 -V_x+bV_y+\frac{h}{r(G-1)}&=D,\label{late:r08:rows}\\
 A&>\frac{hF}{p(F-1)}\ge0,&
 D&>\frac{hG}{r(G-1)}\ge0.\label{late:r08:forcing}
\end{align}
If $x\le13/20$, the determinant of this same actual lift satisfies
\[
 \frac34-\Gamma<\frac{T^3}{5}+\frac4{49}T^4
 \le\frac{x^3}{15\sqrt3}+\frac{4x^4}{441}<\frac1{80}.
\]
In particular $\Gamma>59/80$. No harmonic calibration is assumed
in deriving this lower bound.
\end{proposition}
\begin{proof}
Put $s=\log F$, $t=\log G$, $\Sigma=s+t$. The positive responses give
$q/p<d(F)$ and $q/r<d(G)$, where
\[
 d(Y)=\frac{Y(Y-1-\log Y)}{Y\log Y-Y+1}.
\]
Since $q$ is the largest reciprocal, these inequalities control every
ordered ratio. The exact physical-cost decomposition is
\begin{align*}
 x={}&(p+r)D_{\rm cap}+f(F)+f(G)
 +\frac{p^2f(F)}F+pq(1-F^{-1})\\
 &+\frac{r^2f(G)}G+rq(1-G^{-1}),\qquad D_{\rm cap}\ge0.
\end{align*}
The kernel bound \ref{strip:reverse-credits} gives
$f(e^u)>u/2+5u^2/(60+u^2)$. At $u=10/9$ its right side is
$1465/2232>13/20$, so $s,t<10/9$, including the case $x=13/20$.
The logarithmic ratio bound with $(c,B)=(27/40,154/125)$ and the
exponential bound \ref{strip:exp-upper} give all ratios $<17/8$.

For $0\le u\le6/5$ set
\[
 L(u)=\frac{1-u+u^2/2-u^3/6}{2},\qquad
 P(u)=1-u/2+u^2/6-u^3/24.
\]
The kernel bound $f(e^u)>u/2$ and the positive Taylor remainder
$e^{-u}>1-u+u^2/2-u^3/6$ prove the first credit below by
multiplication; the cubic polynomial is positive on $[0,6/5]$.
Integrating that same exponential bound from $0$ to $u$ proves
the second credit. These decreasing polynomials therefore supply
$f(e^u)e^{-u}>uL(u)$, $1-e^{-u}>uP(u)$ and
\[
 p^2L(u)+pqP(u)\ge\frac{L(u)+P(u)}{1+2\rho}
\]
on the ratio cone of size $\rho$, by \ref{strip:cone-credit}.
Using $L(10/9)=607/4374>1/8$, $P(10/9)=1297/2187>7/12$
and $\rho=17/8$ in the full cost gives
\[
 x>\frac{40}{63}\Sigma+\frac{\Sigma^2}{25}.
\]
Its value at $\Sigma=1$ exceeds $13/20$, so $\Sigma<1$.
The logarithmic ratio bound with $c=67/100$ now gives $\rho=2$,
because $e^{67/100}<4564450/2330000<2$.
Using $L(1)=1/6$, $P(1)=5/8$ gives
\[
 x>\frac{79}{120}\Sigma+\frac5{122}\Sigma^2.
\]
At $\Sigma=15/16$ the right side is $20401/31232>13/20$.
Hence $\Sigma<15/16$. Finally
\[
 e^{201/320}<\frac{1242605721}{662937600}<\frac{15}{8},
\]
so every reciprocal ratio is $<15/8$.

Write $\ell=p+q+r$. The ratio cone gives $\ell<9/5$.
Here is the elementary geometric verification: with
$\alpha=p/q$, $\beta=r/q$, $\rho=15/8$, and
$K_0=(\rho+2)^2/(2\rho+1)$, the separately concave polynomial
\[
 K_0(\alpha+\beta+\alpha\beta)-(1+\alpha+\beta)^2
\]
is nonnegative at the four corners of $[1/\rho,1]^2$.
Its corner values are $0$, $3(\rho-1)^2/(2\rho+1)$ and, at each of the two mixed corners,
$(\rho-1)^3(\rho+1)/(\rho^2(2\rho+1))$.
Thus $\ell^2\le961/304<81/25$.
Also $q\ge1/\sqrt3$ and the physical density bound gives
$T\le x/\sqrt3<\sqrt3/4$.

Set $\mu=4/7$, $Q_0=8/11$, $\gamma=256/3117$.
For $0<u<15/16$ the kernel and exponential estimates yield
\[
 e^u-1-u<Q_0u^2,\qquad f(e^u)>u/2+\gamma u^2.
\]
Let $d=q-p\ge0$, $e=q-r\ge0$. Expanding $A$ and retaining its
negative linear term gives, with $v=q/p$ and $w=r/p\le v$,
\begin{equation}\label{late:r08:response}
 A<\frac{p+r}{p}\left(\mu t_0^2-\frac{dt_0}{2}\right),
 \qquad
 D<\frac{p+r}{r}\left(\mu t_2^2-\frac{et_2}{2}\right).
\end{equation}
Indeed the quadratic coefficient before this relaxation is
$Q_0-\gamma(1+v)\le19304/34287<4/7$.
Since both responses are positive, $d<2\mu t_0$ and $e<2\mu t_2$.
Moreover $k_0=1/(2\sqrt3t_0)>1/(2x)\ge10/13>\mu$, and
\[
 \frac{A}{t_0^2}<(1+v)\left[\mu-k_0(1-v^{-1})\right]\le2\mu.
\]
The last inequality follows by subtracting $2\mu$:
the difference is $(v-1)[\mu-k_0(v+1)/v]\le0$.
The upper bracket is positive because the preceding response bound
is positive. Exchanging $p,r$ gives
\begin{equation}\label{late:r08:quadratic}
 0<A<2\mu t_0^2,\qquad0<D<2\mu t_2^2.
\end{equation}

Put $S=p+r$, $N_x=(p-r)/S$, $N_y=q-(p^2+r^2)/S\ge0$,
$k_p=1/[p(F-1)]$, $k_r=1/[r(G-1)]$.
Solving \eqref{late:r08:rows} gives
\[
 V(h)=\frac{A-hk_p}{a+b}(b,1)
       +\frac{D-hk_r}{a+b}(-a,1).
\]
Using $1-pr=q(p+r)$, the signed coefficients are exactly
\[
 \frac{bN_x+N_y}{a+b}=\frac{p(q-r)}{qS}\ge0,\qquad
 \frac{-aN_x+N_y}{a+b}=\frac{r(q-p)}{qS}\ge0.
\]
Thus $N\cdot V(h)\le N\cdot V(0)$ for the actual $h\ge0$.
Multiplying \eqref{late:r08:response} by these nonnegative
coefficients yields
\[
 N\cdot V(0)\le\frac\mu q(et_0^2+dt_2^2)-\frac{deT}{2q}.
\]
Set $U=d/(2\mu t_0)$, $W=e/(2\mu t_2)$ in $[0,1)$.
The right side is
\[
 \frac{2\mu^2}{q}t_0t_2\{t_0W(1-U)+t_2U(1-W)\}
 \le\frac{2\mu^2}{q}t_0t_2\max(t_0,t_2).
\]
The last step holds throughout the closed unit square since the
two nonnegative weights sum to $U+W-2UW\le1$.
Consequently
\[
 \tfrac12N\cdot V(h)\le\sqrt3\mu^2t_0t_2\max(t_0,t_2).
\]

The strict forcings \eqref{late:r08:forcing} also imply
\[
 h<Ap(1-F^{-1})<At_0,\quad h<Dr(1-G^{-1})<Dt_2,
 \quad h<2\mu\min(t_0^3,t_2^3).
\]
They make both coefficients in the solved vector positive.
The bilateral geometry gives each coefficient vector norm at most
one: for example
\[
 (a+b)^2-(b^2+1)=(a-r)(a+r^{-1})\ge0,\quad
 a-r=\frac{(p+r)(q-p)}{2p}\ge0.
\]
Hence $|V(h)|\le A+D<2\mu(t_0^2+t_2^2)$.
All these bounds concern the actual height, not an algebraic
endpoint of the cap segment.

The polarization identity \ref{sc:polarization}, with $h,V_y\ge0$,
now gives
\begin{align*}
 \tfrac34-\Gamma
 &=\tfrac12N\cdot V+\tfrac12(\ell-N_y)h
      +\tfrac14|V|^2-\tfrac12hV_y\\
 &<\sqrt3\mu^2t_0t_2\max(t_0,t_2)
       +\tfrac95\mu\min(t_0^3,t_2^3)
       +\mu^2(t_0^2+t_2^2)^2.
\end{align*}
For $z=\min(t_0,t_2)/T\in(0,1/2]$, divide the right side by $T^3$.
It becomes
\[
 F(z)=cz(1-z)^2+b_0z^3+d_0(1-2z+2z^2)^2,
 \quad c=\sqrt3\mu^2,\ b_0=\tfrac95\mu,\ d_0=\mu^2T.
\]
With $u=2z$, its derivative has the exact decomposition
\[
 \frac{d}{du}F(u/2)=A_0(1-u)^2+(B_0-A_0)u^2+d_0u^3,
\]
where $A_0=c/2-2d_0>0$ and $B_0=3(c+b_0)/8-3d_0$.
Since $T<\sqrt3/4$ and $b_0\ge\sqrt3\mu$,
$B_0-A_0\ge3\sqrt3\mu(1-\mu)/8>0$.
Thus the maximum on the whole closed interval occurs at $z=1/2$,
including the algebraic zero-duration boundary. It follows that
\[
 \tfrac34-\Gamma<\frac{\sqrt3\mu^2+(9/5)\mu}{8}T^3
                     +\frac{\mu^2}{4}T^4
 <\frac{T^3}{5}+\frac4{49}T^4.
\]
The last coefficient comparison uses $\sqrt3<7/4$.
Finally $x\le13/20$ gives $T^4<1/50$ and $T^3<27/500$,
so the loss is $<1523/122500<1/80$; the final reserve is
$33/490000$. No division by a possibly zero $h$, $d$, $e$ or
$V_y$ has occurred.
\end{proof}

\begin{theorem}[A complete reverse-pair source]\label{opt:pair650}
Let the complete positive source $20210212$, with local positions
$0,\ldots,7$, lie in one genuine containing fundamental interval of
an original global minimizer. This includes a proper occurrence and
a whole eight-arc occurrence in a genuine switch-cut period.
Then its two central complete physical costs satisfy
\[
 c_3+c_4>\frac{13}{20}.
\]
The statement also holds for a whole lawful cyclic image or a whole
reflection combined with time reversal. All eight source arcs, including
the actual caps and both preceding $202$ triples in their respective
orientations, are retained in one genuine fundamental interval with the
original state and group endpoints and real winding.
\end{theorem}
\begin{proof}
Suppose $x=c_3+c_4\le13/20$. Since both costs are positive,
each is strictly less than $13/20$. For a proper occurrence choose a
genuine cut inside the positive arc preceding the source. For an
eight-arc period cut instead at the switch beginning position $0$,
using Lemma~\ref{fd:reanchor}; all eight arcs then remain complete.
Every root, cap tie and contact used below lies strictly inside this
interval. Select one nonzero full maximizing lift there. The capped pair on positions
$2,3,4,5$ is $2\mid10\mid2$, which whole phase $+2$ takes to
$1\mid02\mid1$. Lemma~\ref{sc:pair} proves normality of this selected
lift. Only then jointly normalize its objective multiplier, matrix
adjoint and the state covector. In this actual pair chart, the
maximum conditions and Lemma~\ref{sc:pair} with $\delta=1$ give
that the scaled state covector is $(h,0)$, $h\ge0$, $V_y\ge0$,
and $0\le h\le h_*$, where $h_*$ is the affine-cap endpoint
defined in \eqref{sc:interval}. These signs hold before any determinant or harmonic
calibration is asserted.

The forward prefix $(0,1,2,3,4,5)$ is $202102$ with centre $3$.
The opposite prefix $(7,6,5,4,3,2)$, under the full reflected-time map
$j\mapsto1-j$, is also $202102$, with centre $4$.
In each chart its initial complete $202$ and positive following tail
supply precisely the strict actual ingress \eqref{hc:ingress}.
Both actual maximizing contacts are retained. Apply
Theorem~\ref{hc:prefix650} twice on the same selected full lift,
transported as a whole for the second application. The two exits are
the same physical middle junction $3|4$.

At each exit the conclusions are $u<0$, $q\ge0$ and $q'>kq$.
The raw switching function $q=w\chi$, with $w>0$, is $C^1$ at the
junction by \eqref{hc:C1-identity}. On the outgoing mode-$0$ side,
the one-sided logarithmic derivative of $w$ at the junction is
$2v/S=k$. Thus $q\ge0$ gives $\chi\ge0$, and $q'-kq>0$ gives
$2K-3-u^2>0$ in each oriented chart.
For the fixed-chart reciprocal coordinates $(x_f,y_f,z_f)$ at this
junction, the two horizontal inequalities give $z_f>x_f,y_f$.
The whole phase making the pair $1\mid02\mid1$ has the dictionary
\[
 (p,q,r)=(y_f,z_f,x_f),\qquad q>p,r,\qquad pq+pr+qr=1.
\]
In the notation of Proposition~\ref{late:r08:cubic}, the full matrix
conjugacies give the two strict forcing values
\[
 -\frac{2(p+r)\{2rV_x+(r^2-1)V_y+2h\}}{(1+r^2)^2}>0,
 \qquad
 -\frac{2(p+r)\{-2pV_x+(p^2-1)V_y+2h\}}{(1+p^2)^2}>0.
\]
For clarity the cyclic state chart is represented by
$R=\tfrac12\left(\begin{smallmatrix}1&-3\\1&1\end{smallmatrix}\right)$;
the opposite chart is reflection and time reversal followed by this
cyclic map. Conjugate the entire matrix and transport each state covector
by the inverse Jacobian with the reversal sign. Substitution of the
actual cap reconstruction in Lemma~\ref{strip:reverse-primitive}
gives the displayed formulas. The two binary invariants are not
identified with one another.

Writing $a=(p^{-1}-p)/2$ and $b=(r^{-1}-r)/2$, these inequalities are
\[
 V_x>h/p-aV_y,\qquad V_x<bV_y-h/r.
\]
Insert them in the actual cap rows \eqref{late:r08:rows} to obtain
\[
 A>\frac{hF}{p(F-1)}\ge0,\qquad
 D>\frac{hG}{r(G-1)}\ge0.
\]
Thus every hypothesis of Proposition~\ref{late:r08:cubic} holds on
this same actual normalized lift, including the strict forcings when
$h=0$. It yields $\Gamma>59/80>0$, without prior harmonic calibration.

On the initial actual $202$ in positions $0,1,2$, the two genuine
binary switching zeros give
\[
 2K-3=\frac{\int u^2\,ds}{\int ds}>0.
\]
The interval has positive duration, so its denominator is positive.
In particular $K>0$ for this same matrix in a lawful whole chart.
Corollary~\ref{cal:one-point} supplies the future sheet only after
positive determinant has been proved. The harmonic theorem on this
genuine full interval now yields
\[
 \Gamma<H_*<\frac{1470471}{2000000}<\frac{59}{80},
 \qquad
 \frac{59}{80}-\frac{1470471}{2000000}
 =\frac{4529}{2000000}>0,
\]
a contradiction. No costate is extended periodically, no fictitious
switching zero is inserted, and no algebraic cap endpoint is calibrated
as a trajectory. The conclusion is strict also when the hypothesized
pair cost equals $13/20$.
\end{proof}
\space
\section{The half-unit estimate for a doubly crossed core}
\label{late:r24:section}

The estimate proved in this section concerns a complete five-arc source,
not a free three-arc endpoint problem. Its two exterior arcs have positive
physical duration. All switching equations are equations of a full
maximizing free-time lift: the cost multiplier, the whole matrix adjoint,
and the state covector are retained. Normalization is performed jointly,
after normality has been established. Throughout this section set
\[
 C_* =\frac{180483}{100000},\qquad
 H_* =\frac{29025933167887355473205403}
 {39478417604337498218323684},\qquad
 \overline H=\frac{1470471}{2000000}.
\]
The geometric, endpoint, winding, and harmonic-calibration hypotheses are
those of Section~\ref{fd:section} and Corollary~\ref{cal:one-point}.
In particular, a genuine fundamental interval is always an interval of
the original minimizing trajectory, with its original state and frame
endpoints and its prescribed real winding.

\subsection{Exact reconstruction of the crossed core}

\begin{definition}\label{late:r24:crossed-source}
A \emph{doubly crossed three-arc core} is the middle triple of an actual
complete positive source $02010$, $10121$, or $21202$, on one complete
normal maximizing lift of an original global minimizer. In its lawful
$02010$ chart, its original, unsorted parameters $(a,b,E)$ satisfy
\begin{equation}\label{late:r24:crossings}
 E>1,\quad a,b>0,\quad Ea,Eb<2,\qquad
 Ea>\max(1,b),\quad Eb>\max(1,a).
\end{equation}
The complete source is retained on a genuine containing fundamental
interval. Its two positive caps, all active and inactive Hamiltonian
conditions, and both actual leaf means belong to that selected lift.
\end{definition}

Write $\alpha,c,\beta$ for the three charged physical costs, in order.
The central mode is $0$, and the two charged leaves have modes $2,1$.
Put
\begin{equation}\label{late:r24:parameters}
\begin{gathered}
 h=E-1,\quad z=\log E,\quad P=ab,\quad Q=a+b,\quad R=Eab,\quad D=4-R,\\
 f(E)=\frac{2(E\log E-E+1)}{(E-1)^2}.
\end{gathered}
\end{equation}
Here $R$ denotes $Eab$, and not the product $(1-a)(1-b)$ occurring in
\eqref{f245:variables}. The inequalities \eqref{late:r24:crossings}
imply $D>0$ and
\begin{equation}\label{late:r24:product-crossings}
 R>\max(a,b,a^2,b^2).
\end{equation}

For $v\in\{a,b\}$ define $q_v=2-v$, $r_v=1-v$ and
$V_v^2=v^2D/R$. If $d$ is its actual leaf width, then $0<d<v$.
The exact mean excess and physical leaf cost are
\begin{align}
 e_v(d)
 &=\frac{\displaystyle\int_0^d
 (2r_v s+s^2)(q_v+s)^{-2}\,ds}
 {\displaystyle\int_0^d(q_v+s)^{-2}\,ds}
 =d\left[q_v-f\left(1+\frac d{q_v}\right)\right],
 \label{late:r24:mean}\\
 A_v(d)
 &=\frac d2+
 \frac{[3+r_v^2+e_v(d)]dq_v}{2V_v^2(q_v+d)}
 \label{late:r24:leaf-cost}\\
 &=\frac d2+\frac{q_v^2}{2V_v^2}
 \left[d-2\log\left(1+\frac d{q_v}\right)
                  +\frac{4d}{q_v(q_v+d)}\right].\nonumber
\end{align}
Thus $\alpha=A_a(d_A)$, $\beta=A_b(d_B)$, and
\begin{equation}\label{late:r24:central-cost}
 c=\frac{2Qh-Rz}{D},\qquad
 e_A=e_a(d_A)>0,\quad e_B=e_b(d_B)>0.
\end{equation}
These are the integrals in \eqref{f245:mean-definition} and
\eqref{f245:costs}; they do not replace the true means by an envelope.
The central maximizing signs are
\begin{equation}\label{late:r24:central-signs}
 e_A\ge ah(a-f(E)),\qquad e_B\ge bh(b-f(E)).
\end{equation}
Indeed the expression $\Phi(E)-(1-a)h$ in
\eqref{f245:central-signs} equals $h(a-f(E))$.

For completeness, the matrix reconstruction can be written without any
choice of independent leaf costates. In a binary chart set
\[
 \Lambda=\begin{pmatrix}l&b_0\\m&-l\end{pmatrix},\qquad
 K=2ul+b_0-u^2m,\qquad L=l-um,\qquad \chi=vp_v,
 \qquad ds=dt/v.
\]
The full adjoint equation and the active equation $H_j=0$ give
\begin{equation}\label{late:r24:edge}
 K_s=0,\quad L_s=K,\quad m_s=-2L,\quad
 \chi_s=2K-3-u^2,\qquad \det\Lambda=-L^2-Km.
\end{equation}
The two actual roots of each leaf determine
\[
 K_A=\frac{3+(1-a)^2+e_A}{2}>\frac32,\qquad
 K_B=\frac{3+(1-b)^2+e_B}{2}>\frac32.
\]
Integrating on the central ray and using the actual exit tie gives
\begin{equation}\label{late:r24:matrix}
 J_m=\frac{e_A/a+e_B/b+2z+2Q-2-R/2}{D},\qquad
 \Gamma:=\det\Lambda=\frac{4K_AK_B}{D^2}-J_m^2.
\end{equation}
No positivity assumption on $\Gamma$ is used in these identities.
With $x=e_A/(aD)$ and $y=e_B/(bD)$, write
\begin{align}
 C_0&=\frac{[3+(1-a)^2][3+(1-b)^2]
                -(2z+2Q-2-R/2)^2}{D^2},
 \label{late:r24:central-determinant}\\
 DL_x&=(1-b)(4-P)+Ph-4z,\qquad
 DL_y=(1-a)(4-P)+Ph-4z.\nonumber
\end{align}
Then the exact determinant is
\begin{equation}\label{late:r24:determinant}
 \Gamma=C_0+L_xx+L_yy-x^2-(2-P)xy-y^2.
\end{equation}
The symbol $C_0$ is an algebraic part of the same matrix determinant;
it is not the determinant of a newly asserted zero-leaf trajectory.

\subsection{Consequences of a contrary half-unit budget}

\begin{lemma}\label{late:r24:cheap-domain}
Suppose \eqref{late:r24:crossings}--\eqref{late:r24:central-signs}
hold for the true means and costs, and
\begin{equation}\label{late:r24:contrary-budget}
 \alpha+c+\beta\le\frac12.
\end{equation}
Then
\[
 1<E<\frac94,\qquad R<\frac32,\qquad D>\frac52,
 \qquad 0<P<\frac32,
 \qquad L_x,L_y<0.
\]
These conclusions precede harmonic calibration.
\end{lemma}
\begin{proof}
The signed weighted-mean argument in the proof of Lemma~\ref{f245:negative}
gives, if $a>1$,
\[
 0<e_A\le\frac{d_A^2}{3}-(a-1)d_A,
 \qquad d_A>3(a-1),\qquad
 \alpha>\frac32(a-1).
\]
The same statement holds on the other leaf. Consequently
\[
 \alpha+\beta>\frac32(Q-2)_+.
\]
If $P\ge3/2$, then $Q\ge2\sqrt P$ and
$3(\sqrt{3/2}-1)>1/2$, a contradiction. Hence $P<3/2$.

At fixed $P,Q$,
\[
 D^2c_E=2Q(4-P)-4P(z+1)+EP^2
       >(4-P)(2Q-EP).
\]
The difference is $4P(E-1-\log E)>0$. At fixed $E,P$, both $c$ and
$\tfrac32(Q-2)_+$ are nondecreasing in $Q$. If $R\ge3/2$, first lower
$Q$ to $2\sqrt P$ and then lower $E$ to $3/(2P)$. These comparisons
have positive denominators and $2Q-EP>0$: writing $t=\sqrt P$, the
original upper state inequalities give $Et<2$, and the same upper bound
persists while $E$ decreases. The resulting lower bound for the total
cost is
\[
 J(t)=\frac{6/t-4t-\tfrac32\log((3/2)/t^2)}{5/2}
       +3(t-1)_+,\qquad 0<t<\sqrt{3/2}.
\]
Its derivative is $(-12/t^2+6/t-8)/5<0$ for $t\le1$, and
$(7+6/t-12/t^2)/5\ge1/5$ for $t\ge1$. Therefore
\[
 J(t)\ge J(1)=\frac45-\frac35\log\frac32>\frac12,
\]
again a contradiction. This proves $R<3/2$.

At fixed $E,b$, the numerator of $D^2c_a$ is
$8h-4Ebz+2Eb^2h>8(h-z)>0$; the other partial derivative is positive
as well. Lowering both coordinates to $1/E$ is therefore permissible as
a central algebraic comparison and yields
\[
 c\ge B_0(E):=\frac{4(E-1)-\log E}{4E-1}.
\]
The numerator of $B_0'$ is $8+1/E+4\log E>0$, and
$B_0(9/4)>1/2$ since $\log(9/4)=2\log(3/2)<1$.
Thus $E<9/4$.

It remains to prove $L_x<0$; the opposite coefficient has the same
proof with the separately derived scalar names exchanged. If $b\ge1$,
then $(1-b)(4-P)\le0$ and $Ph<2h/E<4\log E$. If $b<1$, write
\[
 N(a,b)=DL_x=4(1-b)+ab(E-2+b)-4\log E.
\]
Its derivative in $b$ is $-4+a(E-2+2b)<-4+aE<-2$.
The crossed lower bound permits decreasing $b$ to $\max(1,a)/E$.
For $a\le1$, the resulting expression increases with $a$ and is at most
\[
 N(1,1/E)=\frac{h(5E-1)}{E^2}-4\log E<0.
\]
The last sign follows by differentiating
$4\log E-5+6/E-1/E^2$, whose derivative is
$2(2E-1)(E-1)/E^3>0$ and whose value at $E=1$ is zero.
For $a>1$, the relevant boundary is $b=a/E$. Here $E<2$ and $a^2<R<2$.
The numerator of its derivative in $a$ is
\[
 3a^2+2aE(E-2)-4E
 <6+2E^2-8E=2((E-2)^2-1)<0.
\]
Decreasing $a$ to $1$ therefore bounds this case by the same negative
quantity. This proves both coefficient signs.
\end{proof}

\begin{definition}\label{late:r24:scalar-domain}
An \emph{exact scalar comparison} retains the central parameters,
\eqref{late:r24:mean}--\eqref{late:r24:determinant}, the nonnegative
true-mean branches, and the central signs. Its widths may be zero only
where the corresponding minimum mean excess is zero. Unless a closed
face is stated explicitly, it retains the strict state inequalities
\eqref{late:r24:crossings}. It need not satisfy an exterior-root test or
represent a maximizing trajectory.

For fixed $(E,a,b)$, define the true central minimum as follows. If
$v\le f(E)$, let $d_{v,\min}=m_v=0$. If $v>f(E)$, let
$d_{v,\min}\in(0,v)$ be the positive-mean root
\[
 e_v(d_{v,\min})=vh(v-f(E)),\qquad
 m_v=A_v(d_{v,\min}).
\]
Only states for which the stated roots exist and the indicated budget
is feasible are included. An actual state satisfying
\eqref{late:r24:contrary-budget} has these properties.
\end{definition}

\begin{lemma}[Exact allocation reduction]\label{late:r24:allocation}
At a state in Lemma~\ref{late:r24:cheap-domain}, put $\ell=1/2-c$.
Then the two exact scalar comparisons
\[
 (\alpha,\beta)=(m_a,\ell-m_a),\qquad
 (\alpha,\beta)=(\ell-m_b,m_b)
\]
exist whenever an actual contrary-budget point exists, and
\begin{equation}\label{late:r24:allocation-bound}
 \Gamma_{\rm actual}\ge
 \min\{\Gamma(m_a,\ell-m_a),\Gamma(\ell-m_b,m_b)\}.
\end{equation}
The determinant on the right is evaluated using the exact inverse-cost
means at the unchanged central state.
\end{lemma}
\begin{proof}
For one leaf, let $w=d/q_v$. Direct differentiation gives
\[
 e_v'(d)=q_v+\frac{f(1+w)-2}{1+w},\qquad
 e_v''(d)=\frac{2[1-f(1+w)]}{q_vw(1+w)}>0,
\]
with continuous values $e_v(0)=0$ and $e_v'(0)=1-v$. Thus $e_v$ is
increasing on its positive-mean branch. Also
\[
 A_v'(d)=\frac12+\frac{q_v^2}{2V_v^2}
 \left(1-\frac2{q_v+d}+\frac4{(q_v+d)^2}\right)>0,
 \qquad
 A_v''(d)=\frac{q_v^2(q_v+d-4)}{V_v^2(q_v+d)^3}<0
\]
for $d<v$. The first sign follows from $(q_v+d-1)^2+3>0$.

By \eqref{late:r24:product-crossings}, $V_v^2<v(4-v)$. At the coordinate
barrier $d=v$ the exact cost therefore satisfies
\begin{equation}\label{late:r24:inverse-barrier}
 A_v(v)>1-\frac{q_v^2\log(2/q_v)}{4-q_v^2}>\frac12.
\end{equation}
For the last inequality use $1+2\log X<X^2$ for $X=2/q_v>1$.
If a closed lower-crossing face weakens the height comparison to
$V_v^2\le v(4-v)$, the first comparison in
\eqref{late:r24:inverse-barrier} becomes weak but the final comparison
remains strict. Every budget below $1/2$ thus has a unique inverse width
strictly below $v$, and on a nonnegative-mean branch the cost inequality
$A_v(d)>d/2$ also gives $d<1$.

The central cost is positive: $z<h$ and $2Q>R$ imply
$2Qh-Rz>h(2Q-R)>0$. Hence $\ell<1/2$. The actual leaves show that
$m_a+m_b\le\ell$. On the retained branches the functions
\[
 X(A)=\frac{e_a(A_a^{-1}(A))}{aD},\qquad
 Y(B)=\frac{e_b(A_b^{-1}(B))}{bD}
\]
are increasing and convex. By Lemma~\ref{late:r24:cheap-domain}, the
polynomial \eqref{late:r24:determinant} is decreasing in each
nonnegative mean. Its Hessian has eigenvalues $-P$ and $-(4-P)$, so it
is jointly concave. Its composition with the two convex maps $X,Y$ is
therefore concave. Increasing unused comparison budget to $\ell$ can
only decrease the determinant. On the closed allocation interval
$[m_a,\ell-m_b]$, concavity places the minimum at an endpoint, proving
\eqref{late:r24:allocation-bound}. A singleton interval needs no
interior argument. No state or costate is assigned to either artificial
comparison.
\end{proof}

\subsection{Scalar estimates for the direct sectors}

The following scalar estimates, and the scalar inputs explicitly identified
below, are used with their stated domains. Their independent analytic
proofs are given in Section~\ref{late:r24:scalar-proofs}. We retain separate statements so
that the two allocation endpoints and the different closed state faces
are not identified with one another.

\begin{lemma}[Direct sector estimates]\label{late:r24:direct-sectors}
Use the true means, exact costs, central signs, and contrary budget of
Definitions~\ref{late:r24:crossed-source} and
\ref{late:r24:scalar-domain}, with the strict original state inequalities.
The following conclusions hold.
\begin{enumerate}[label=(\roman*)]
 \item If $E\ge26/15$, then $\Gamma>1471/2000$.
 This estimate also holds for the exact scalar comparisons with
 $a,b\ge1/E$, $a,b<f(E)$, $1<E<9/4$, $R<3/2$, and total cost at most
 $1/2$ that occur on the closed high-$E$ face of the profiles below.
 \item If $a,b\ge1$, then $\Gamma>737/1000$.
 \item If $0<a\le1<b$, then
 $\Gamma>1471/2000+1/97500$. The same conclusion holds for
 $0<a<1=b$, provided the crossed and upper inequalities remain strict.
 It holds for arbitrary retained exact allocations, including slack and
 a zero unsigned mean. The other orientation has the identical statement
 with the original scalar leaf names exchanged.
\end{enumerate}
\end{lemma}

The three assertions are proved in
Theorems~\ref{late:r24:a:high-e}, \ref{late:r24:a:fully-signed},
\ref{late:r24:a:mixed}, and Corollary~\ref{late:r24:a:unit}, respectively.
Their proofs use the zero-leaf estimate of
Theorem~\ref{late:r24:a:zero-leaf} and the central comparison of
Lemma~\ref{late:r24:a:r11-central}; neither uses the profile conclusions
being assembled here.
The first proof compensates a leaf loss of at most
$109\ell/1000$ by a central gain greater than $\ell/9$ before the
zero-leaf cost face. The second retains the signed mean and bounds the
whole cost simplex by
\[
 \Gamma>\frac34-\left(\frac{c^2}{12}
       +\frac8{25}c\ell^2+\frac{14}{375}\ell^3
       +\frac{16}{225}\ell^4\right)>\frac{737}{1000}.
\]
The third retains the common signed depth and proves, jointly rather
than with independent leaf budgets,
\[
 (-L_x)x+(-L_y)y<\frac{37}{4875},\qquad
 \left(\frac34-C_0\right)+x^2+(2-P)xy+y^2<\frac{69}{10000}.
\]
Their sum gives the exact identity
\[
 \frac34-\frac{37}{4875}-\frac{69}{10000}
 =\frac{1471}{2000}+\frac1{97500}.
\]
The coordinate-one statement retains strict crossings; it is not a
statement about the closed circle $Ea=1$ or $Eb=1$.

\subsection{The unsigned zero-minimum endpoint}

In this subsection $1<E<26/15$ and $0<a,b<1$. We consider an allocation
endpoint with $a\le f(E)$, $\alpha=e_A=0$, and
$\beta=1/2-c\ge m_b$. Every occurrence of a mean below uses its exact
inverse-cost value.

The interval and exact determinant concavity in
Lemma~\ref{late:r24:a:zero-profile-inputs} reduce the comparison to
the faces in \eqref{late:r24:zero-faces}.

\begin{lemma}[Zero-minimum boundary estimates]\label{late:r24:zero-boundaries}
On the closed scalar domain of Lemma~\ref{late:r24:a:zero-profile-inputs}, each
face in \eqref{late:r24:zero-faces} satisfies $\Gamma>\overline H$.
More explicitly, the following statements are used.
\begin{enumerate}[label=(\roman*)]
 \item On the circle $a=1/E$, $1/E\le b\le f(E)$, with the $B$ mean
 zero and the $A$ leaf assigned the exact budget $1/2-c$, one has
 $\Gamma>\overline H$. This includes the cap $b=f(E)$, and includes
 the opposite orientation obtained by exchanging scalar names.
 \item On $a=f(E)$, $e_A=0$, the exact full-budget determinant is
 increasing in $b$ on the feasible continuation down to $b=1/E$.
 Thus (i) bounds this cap face.
 \item If $e_A=0$ and $e_B=bh(b-f(E))>0$, lowering $a$ to $1/E$
 with the true minimum $B$ width fixed does not increase its total
 cost and does not increase the determinant. Any new budget slack
 may then be filled using the exact $B$ inverse. On the resulting
 circle $a=1/E$, $b>f(E)$, $\alpha=0$,
 $\beta=1/2-c\ge m_b$, the exact determinant exceeds $\overline H$.
 \item If both leaf budgets vanish, the retained signs give
 $a,b\le f(E)$, and $C_0>1471/2000$.
 \item The $E=26/15$ face is covered by the exact-comparison statement
 in Lemma~\ref{late:r24:direct-sectors}(i).
\end{enumerate}
The only possible closed $E=1$ corner has $k=a=b=1$; its exact
comparison determinant exceeds $299/400$.
\end{lemma}

For (i), use Lemmas~\ref{late:r24:b:opposite-circle}
and~\ref{late:r24:b:positive-a-cap}. Part (ii) is
Lemma~\ref{late:r24:b:cap-transfer}. Part (iii) follows from
Lemmas~\ref{late:r24:b:saturation-transfer}
and~\ref{late:r24:b:first-circle}. Part (iv) is
Theorem~\ref{late:r24:a:zero-leaf}, and part (v) uses
Theorem~\ref{late:r24:a:high-e}. The closed $E=1$ corner is included
in the proof of Lemma~\ref{late:r24:a:zero-profile-inputs}.
The two circle orientations are kept distinct in these proofs.
Part (iii) is a statement about the exact determinant and exact costs.
It is not a comparison of an exact endpoint with a smaller determinant
envelope.

\begin{proposition}\label{late:r24:zero-minimum}
Every unsigned low-$E$ full-budget allocation endpoint whose selected
central minimum is zero satisfies $\Gamma>\overline H$.
\end{proposition}
\begin{proof}
Take the interval through the given point in
Lemma~\ref{late:r24:a:zero-profile-inputs}. Its endpoints are covered by
Lemma~\ref{late:r24:zero-boundaries}, including their intersections.
Strict concavity of the exact determinant in
Lemma~\ref{late:r24:a:zero-profile-inputs} gives the same strict lower bound in
the interior. A singleton is already a boundary case. If the interval
starts from an original strictly crossed point, $k=Ea>1$ precludes an
$E=1$ endpoint; the stated corner estimate also handles the enlarged
closed domain. The other allocation orientation follows by exchanging
the two already derived scalar laws.
\end{proof}

\subsection{The unsigned positive-minimum endpoint}

The remaining unsigned endpoint has $a>f(E)$ and
\begin{equation}\label{late:r24:positive-endpoint}
 e_A=ah(a-f(E)),\qquad \alpha=m_a,\qquad
 \beta=\frac12-c-m_a\ge m_b,
 \quad 1<E<26/15,\quad 0<a,b<1.
\end{equation}
The following three separate analytic estimates are needed for this case.

The proof uses the curvature of the true minimum cost, the comparison
when both costs are at their exact minima, and the circle comparison
with both leaf budgets. These are Lemmas~\ref{late:r24:c:G-curvature},
\ref{late:r24:c:both-minima} and~\ref{late:r24:c:two-positive-circle},
respectively. The last retains the common budget of the two leaves.

\begin{lemma}[Exact positive-minimum concavity]\label{late:r24:positive-concavity}
Fix $R,b$ at a point of \eqref{late:r24:positive-endpoint}, and let
$k=R/b$. Along each feasible open profile component, the exact
full-budget determinant is strictly concave in
\[
 \xi=\frac{\mu_A-1}{k}=a-\frac{2\log E}{E-1}.
\]
\end{lemma}
\begin{proof}
Put
\[
 t=e_B(\beta)/b,\quad K_0=3+(1-b)^2,\quad
 C_1=k+2b-R/2,\qquad \beta=1/2-G(\mu_A).
\]
At fixed $R,b$, the $B$ height is fixed. The exact $B$ mean is an
increasing convex function of its budget. By
Lemma~\ref{late:r24:c:G-curvature}, $\beta$ is convex in $\xi$, and
hence $t_{\xi\xi}\ge0$. No sign of $t_\xi$ is assumed.

The saturation identity gives the exact numerator
\begin{equation}\label{late:r24:positive-matrix}
 Z=D^2\Gamma=(4+k\xi)(K_0+bt)-(\xi+C_1+t)^2.
\end{equation}
Let $W=\xi+C_1+t-b(4+k\xi)/2$. In original scalar variables,
\[
 2W=2t-DL_y+(2-ab)e_A/a>0.
\]
The closed-rectangle estimate \eqref{late:r24:b:closed-coefficients} gives \(L_y<0\), including at \(Eb=1\). Differentiating \eqref{late:r24:positive-matrix} twice yields
\[
 Z_{\xi\xi}=-2\bigl[t_\xi^2+(2-R)t_\xi+1\bigr]
             -2Wt_{\xi\xi}<0,
\]
because, for every real $v$,
\[
 v^2+(2-R)v+1=(v+1-R/2)^2+R(1-R/4)>0.
\]
This differentiates the full determinant, including its two-mean
product.
\end{proof}

\begin{proposition}\label{late:r24:positive-minimum}
Every unsigned low-$E$ full-budget allocation endpoint whose selected
central minimum is positive satisfies $\Gamma>\overline H$.
\end{proposition}
\begin{proof}
Use the orientation \eqref{late:r24:positive-endpoint}. At fixed $R,b,k$,
\[
 \frac{d\mu_A}{da}=E(a-f(E))>0,
\]
so $\xi$ increases strictly with $a$ on this branch. The state bounds
are
\[
 a\le\min(1,R),\qquad a\ge15k/26,\qquad a>f(k/a).
\]
Indeed $Eb\ge1$ is exactly $a\le R$; the bound $Ea=k<2$ is fixed,
and $Eb<2$ follows from $E\le26/15$ and $b<1$.
The equation $a=f(k/a)$ has at most one solution because
$(Ef(E))'=2(1-f(E))/(E-1)>0$. Since $k>1$, the profile cannot reach
$E=1$. The fixed value $R<3/2$ supplies no further endpoint.

A point satisfying $\beta=m_b$ is already a minimum-cost boundary.
Otherwise consider the maximal connected open component containing it
on which $\beta>m_b$ and the varying state bounds are strict. It is an
interval in $a$, and therefore in $\xi$; global connectedness is not
required. Exact inverse continuity and
\eqref{late:r24:inverse-barrier} supply finite continuous endpoint
values. Were all inequalities strict at an endpoint, the component
would extend. Every endpoint is therefore one of the following.
At simultaneous boundary equalities, use (iv) first whenever $Eb=1$,
with the closed range $f(E)\le b_c\le1$ of
Lemma~\ref{late:r24:c:two-positive-circle}; this includes $a=f(E)$,
zero remaining circle-side budget, and $a=R=1$. Next use (iii)
whenever $a=1$ and $R>1$. All remaining uses of the zero-minimum
proposition in (i) and (ii) then have $a,b<1$ and strict crossed
inequalities.
\begin{enumerate}[label=(\roman*)]
 \item $a=f(E)$. Here the selected minimum vanishes, and
 Proposition~\ref{late:r24:zero-minimum} applies.
 \item $\beta=m_b$. If $b>f(E)$, use
 Lemma~\ref{late:r24:c:both-minima}. If $b\le f(E)$, then
 $\beta=e_B=0$, and use Proposition~\ref{late:r24:zero-minimum} in
 the opposite orientation. This includes $b=f(E)$.
 \item $a=1$ and $R>1$. Since $b<1$, the inequalities
 $Ea=E>1>b$ and $Eb=R>1=a$ are still strict. The exact costs and
 central signs are retained, so the coordinate-one assertion of
 Lemma~\ref{late:r24:direct-sectors}(iii) applies.
 \item $a=R\le1$, equivalently $Eb=1$. Exchange the two derived
 scalar leaf laws. The resulting parameters have
 $a_c=1/E$, $f(E)\le b_c\le1$, the $B_c$ leaf at its true
 (possibly zero) minimum,
 and the $A_c$ leaf assigned the exact remaining budget. This is
 Lemma~\ref{late:r24:c:two-positive-circle}. In particular the
 intersection $a=1,R=1$ belongs here, not to (iii).
 \item $E=26/15$. No cheap positive-minimum point reaches this
 face: the central coordinate derivatives are positive and the
 corner estimate used in Lemma~\ref{late:r24:direct-sectors}(i)
 gives
 \[
 c\ge c(26/15,f(26/15),15/26)>1/2.
 \]
\end{enumerate}
A width cannot reach its coordinate barrier at an endpoint, by
\eqref{late:r24:inverse-barrier}. A vanishing $A$ width is case (i),
and a vanishing $B$ width is case (ii). Thus the list is exhaustive.
Every possible endpoint has exact determinant greater than
$\overline H$. Lemma~\ref{late:r24:positive-concavity} proves the
same bound at every interior point of its own component. Equality sets
and singleton components are already boundary cases. The opposite
minimum allocation follows by exchanging the two scalar laws.
\end{proof}

\subsection{The half-unit theorem and its calibration}

\begin{theorem}[Doubly crossed three-arc estimate]\label{late:r24:F500}
Every source of Definition~\ref{late:r24:crossed-source} satisfies
\[
 \alpha+c+\beta>\frac12.
\]
The theorem uses only the actual lift for the final calibration.
\end{theorem}
\begin{proof}
Assume \eqref{late:r24:contrary-budget}. If $E\ge26/15$, apply
Lemma~\ref{late:r24:direct-sectors}(i), including the equality. For
$E<26/15$ and $a,b\ge1$, apply part (ii). For precisely one coordinate
above $1$, or one equal to $1$ with the other below $1$, apply part
(iii) in the corresponding scalar orientation. These cases cover
$\max(a,b)\ge1$.

The remaining case is $0<a,b<1$, $1<E<26/15$. Apply the exact
allocation reduction. Each of its two endpoints has its own selected
minimum. If that minimum is zero, Proposition~\ref{late:r24:zero-minimum}
applies; if positive, Proposition~\ref{late:r24:positive-minimum}
applies. Thus both endpoints exceed $\overline H$, even when just
one minimum is zero. Equation~\eqref{late:r24:allocation-bound} yields
\[
 \Gamma_{\rm actual}>\overline H>H_*>0.
\]
The exact comparison of the constants has positive numerator
\[
\begin{aligned}
 &1470471\cdot39478417604337498218323684\\
 &\quad-2000000\cdot29025933167887355473205403\\
 &=1877293054396185839935164>0.
\end{aligned}
\]

Only now apply the global calibration. On the actual leaf,
$K_A>3/2$ was already proved from the actual roots, independently of
the determinant lower bound. The same-matrix identity
$\Gamma=-L^2-K_Am$ gives $m<0$, and hence
\[
 \operatorname{tr}(\Lambda X)=\frac{K_A-v^2m}{v}>0.
\]
Corollary~\ref{cal:one-point} and Theorem~\ref{cal:harmonic}, on this
same jointly normalized full lift and genuine fundamental interval,
give $\Gamma<H_*$. This is a contradiction. None of the scalar
comparison states is calibrated, and no periodicity of the state
covectors is assumed.
\end{proof}

\space
\section{Scalar boundary estimates for the crossed-core theorem}
\label{late:r24:scalar-proofs}
The scalar comparisons retain the exact costs and true means of
the crossed core. The direct sectors and both zero-minimum orientations
are treated first. The positive-minimum argument then uses implicit
curvature, the exact comparison with both costs at their minima,
and the closed circle with both budgets retained. Polynomial signs
in these analytic arguments use the displayed factorizations or
the exact Bernstein conversion \eqref{cv:bernstein-coefficients}.
\begingroup
\renewcommand{\paragraph}[1]{\par\smallskip\noindent\textit{#1}\quad}
\subsection{Central comparisons and the direct scalar sectors}
\label{late:r24:a:section}

We retain the exact scalar notation
\eqref{late:r24:parameters}--\eqref{late:r24:determinant}. In particular
$R=Eab$, $D=4-R$, $P=ab$, $Q=a+b$, $h=E-1$, and $z=\log E$.
Write $E_0=26/15$ and $l=1/E$. All assertions in this subsection are
inequalities for these scalar formulas. They do not assign a maximizing
lift to a comparison point and do not use harmonic calibration.
Subscripts $A,B$ refer to the leaves with state coordinates $a,b$;
in particular $e_A=e_a$, $e_B=e_b$, $m_A=m_a$, and $m_B=m_b$.

\begin{lemma}[Elementary logarithm and mean bounds]
\label{late:r24:a:elementary}
For $E>1$,
\begin{gather}
 \frac hE<z<\frac h{\sqrt E},\qquad
 \frac{2h}{E+1}<z<\frac{h(E+1)}{2E},\qquad
 z<\frac{E-E^{-1}}2,
 \label{late:r24:a:log-bounds}\\
 \frac{2}{E+1}<f(E)<1,\qquad
 \frac3{E+2}<f(E)<\frac{E+5}{3(E+1)},\qquad
 f(E)<\frac2{\sqrt E+1},
 \label{late:r24:a:f-bounds}\\
 f'(E)=\frac{2-(E+1)f(E)}{E(E-1)},\qquad
 -\frac13<f'(E)<0,\qquad
 (Ef(E))'=\frac{2(1-f(E))}{E-1}>0.
 \label{late:r24:a:f-derivatives}
\end{gather}
Also, with $u=1-E^{-1}\in(0,1)$,
\begin{equation}\label{late:r24:a:series}
 z=\sum_{n\ge1}\frac{u^n}{n},\qquad
 f(E)=1-\sum_{n\ge1}\frac{2u^n}{(n+1)(n+2)}.
\end{equation}
\end{lemma}
\begin{proof}
The lower logarithm bounds follow by differentiation from $E=1$.
The derivative of $\sqrt E-E^{-1/2}-\log E$ is
$(\sqrt E-1)^2/(2E^{3/2})>0$. The trapezoidal inequality for the
strictly convex function $1/t$ on $[1,E]$ gives the second upper bound;
the last upper bound follows either from it or by differentiation.
The defining formula for $f$ gives
\[
 f(E)=2\int_0^1\frac{1-t}{1+ht}\,dt,
 \qquad
 1-f(E)=2h\int_0^1\frac{t(1-t)}{1+ht}\,dt.
\]
These identities prove $0<f<1$, $-1/3<f'<0$, and, by Jensen's
inequality for the probability density $2(1-t)$ of mean $1/3$,
$f>3/(h+3)$. The bound $f>2/(E+1)$ follows upon substituting
$z>2h/(E+1)$ in its defining formula. Similarly
$z<h/\sqrt E$ gives $f<2/(\sqrt E+1)$.

For the upper bound and a further estimate used below, apply Jensen
to the probability density \(6t(1-t)\), of mean \(1/2\), in the
second integral. Bounding its denominator above by \(E\) also gives
\(E(1-f)\ge h/3\). Direct differentiation gives \eqref{late:r24:a:f-derivatives}.
Finally, integrate the geometric series on compact subintervals of
$0\le u<1$ and substitute into the defining formula for $f$ to obtain
\eqref{late:r24:a:series}.
\end{proof}

The exact central determinant has the useful expanded numerator
\begin{equation}\label{late:r24:a:N0}
 C_0=\frac34+\frac{N_0}{D^2},\qquad
 N_0=Q(2Ph-8z)-P^2h(E+1)+P(4h+2Ez)+8z-4z^2.
\end{equation}
For example, its first product before subtracting the square in
\eqref{late:r24:central-determinant} is
$P^2-2PQ+4Q^2-4P-8Q+16$; subtraction and collection give
\eqref{late:r24:a:N0}.

\begin{lemma}[The corner and the rational splitting point]
\label{late:r24:a:corner}
Put
\[
 c_f(E)=c(E,l,f(E)),\qquad
 g(E)=\left(\frac hE\right)^2
 \frac{E(E+1)f(E)^2+2f(E)-4}{(4-f(E))^2}.
\]
Then $C_0(E,l,f(E))=3/4-g(E)$, $c_f$ is strictly increasing for
$E>1$, $c_f(E)\to0$ as $E\downarrow1$, and $g$ is positive and
strictly increasing on $1<E<9/4$. At $E_0=26/15$,
\begin{equation}\label{late:r24:a:split-constants}
 \frac{409}{500}<f(E_0)<\frac{819}{1000},\qquad
 \log E_0<\frac{551}{1000},\qquad
 g(E_0)<\frac{29}{2000},\qquad c_f(E_0)>\frac12.
\end{equation}
\end{lemma}
\begin{proof}
Substitution of $a=l$, $b=f$, and $2Ez=h(2+hf)$ into
\eqref{late:r24:a:N0} gives the determinant identity and
\[
 c_f(E)=\frac hE\,
 \frac{2+(2E-1)f-hf^2/2}{4-f}.
\]
Differentiation using \eqref{late:r24:a:f-derivatives} gives
\begin{equation}\label{late:r24:b:corner-derivatives}
 c_f'(E)=\frac{T(E,f)}{E^2(4-f)^2},\qquad
 g'(E)=\frac{hS(E,f)}{E^3(4-f)^3},
\end{equation}
where \(T(E,f)=(E-2)f^2+(4-14E)f+16E+4\) and \(S\)
is given below.
For $1<E\le2$ its numerator decreases in $f\in(0,1)$ and exceeds
its value $3E+6$ at $f=1$. For $E\ge2$, omit the nonnegative
quadratic term and use $f<1$ to obtain the lower bound $2E+8$.
This proves monotonicity; the limit at $1$ follows from the formula.

The bracket in $g$ vanishes at $f=2/(E+1)$ and increases with $f$,
so $g>0$. In the formula for \(g'\) of \eqref{late:r24:b:corner-derivatives},
the numerator uses
\[
 S(E,f)=-(2E^2+E+1)(Ef^3+6f^2)
       +4(4E^2+4E+7)f-32.
\]
This polynomial is strictly concave in $f>0$. At the endpoints of
$[2/(E+1),1]$ its values are
\[
 S(E,1)=(E-1)(5-2E)(E+2)>0,
\quad
 S\left(E,\frac2{E+1}\right)
 =\frac{8E(E-1)(2E+1)^2}{(E+1)^3}>0
\]
for $1<E<9/4$. Thus $g'>0$.

At \(E_0\), use \(x=11/41\) and the three-term lower sum
\(L=2(x+x^3/3+x^5/5)\). The series tail gives the rational chain
\[
 \frac{214489}{390000}<L<z_0
 <L+\frac{2x^7}{7(1-x^2)}
 <\frac{143033}{260000}<\frac{551}{1000}.
\]
Since \(f(E_0)=(780z_0-330)/121\), this gives the two claimed
bounds for \(f\). The fraction defining \(g(E_0)\) increases in
\(f\), with derivative numerator \(2f[4E_0(E_0+1)+1]\).
Using \(f<819/1000\) gives
\(g(E_0)<1234073797/85503530450<29/2000\).
The sign \(c_f(E_0)>1/2\) is the sign of
\(-2/E_0+f(2E_0-3/2-z_0)\), which is greater than
\(-15/13+(409/500)(4247/3000)=81299/19500000\).
\end{proof}

\begin{theorem}[Complete capped zero-leaf comparison]
\label{late:r24:a:zero-leaf}
For every $E>1$, $l\le a,b\le f(E)$, and $c\le1/2$, the exact
central expression satisfies
\[
 C_0>\frac{1471}{2000}.
\]
In the high-$E$ part $E\ge E_0$, the circle point of cost $1/2$ is
\begin{equation}\label{late:r24:a:bc}
 a=l,\qquad b=b_c(E)=\frac4{E(4E-3-2z)}.
\end{equation}
Whenever the cheap capped square is nonempty, this point exists with
$l\le b_c<f$. Among its cost-$1/2$ points it has smallest $C_0$,
and $C_0(E,l,b_c(E))$ strictly increases with $E$ on its feasible
high-$E$ interval.
\end{theorem}
\begin{proof}
First, $Ef^2<4E/(\sqrt E+1)^2<4$ by
Lemma~\ref{late:r24:a:elementary}; hence $D>0$ on the entire capped
square before any restriction on $E$ is made. The globally valid
coordinate derivative is
\begin{equation}\label{late:r24:a:coordinate-cost}
 D^2c_a=2h[4+Eb^2-2Eb(z/h)]
 >2h[3+(\sqrt E b-1)^2]>0,
\end{equation}
and the same holds for $c_b$. Thus
\[
 c\ge B_0(E):=\frac{4h-z}{4E-1},\qquad
 (4E-1)^2B_0'(E)=8+1/E+4z>0.
\]
The trapezoidal bound gives $\log(3/2)<5/12$, so
$\log(9/4)<5/6$ and $B_0(9/4)>25/48>1/2$. Therefore $E<9/4$.
In particular $P<1$ and $D>4-E>7/4$.

\paragraph{The fixed-product boundary reduction.}
At fixed $E,P$, the coefficient $(N_0)_Q=2Ph-8z$ is negative, since
$z>h/E$ and $E<9/4$. The maximal sum in the square is
\[
 Q\le\min(l+P/l,f+P/f).
\]
The first value is attained on the circle $a=l$ when $P\le lf$;
the second is attained on the cap $a=f$ when $P\ge lf$, with scalar
names exchanged if necessary. Increasing $Q$ until the first boundary
of the square or of $c\le1/2$ is reached decreases $C_0$.
The only possible resulting faces are therefore the circle, the cap,
and the cost level $c=1/2$. A collapsed interval is retained.

On the circle put $u=h/E$ and $q=z/u$. Then
\[
 C_0-\frac34=\frac{u^2[b^2+A b+B]}{(4-b)^2},\quad
 A=-2+\frac{6(1-q)}u<-5,\quad B=8q-4q^2\le4.
\]
The numerator of its derivative in $b$ is
$b(8+A)+4A+2B$. It is negative if $A\le-8$; if $-8<A<-5$, it is
less than $3-20+8<0$. Thus $C_0$ strictly decreases along the circle.

For fixed $a,E$, write $N_0=A_2b^2+A_1b+A_0$, where
\[
 A_2=ah[2-a(E+1)],\quad
 A_1=2a^2h+4ah+(2Ea-8)z,\quad
 A_0=8(1-a)z-4z^2.
\]
The derivative numerator is
$(8A_2+EaA_1)b+4A_1+2EaA_0$. On $a=f$ its slope is
\[
 ah(Ea-4)[(E+1)a-2]<0.
\]
There is no interior minimum on a cap segment. Moreover its value at
$b=f$ is
\[
 \frac hE[(E+1)f-2]
 [E^2f^3-2(E^2+2E-1)f^2-4f+16].
\]
The bracket exceeds $14-2E^2-4E$, which is at least $7/8$ for
$1<E\le7/4$. Hence on that range the entire cap increases in $b$.
For $E\le E_0<7/4$, omit the cost restriction: the preceding
fixed-product comparison and the two boundary monotonicities put the
minimum over the square at $(l,f)$. Lemma~\ref{late:r24:a:corner}
then gives
\[
 C_0\ge\frac34-g(E_0)>\frac{1471}{2000}.
\]

\paragraph{High-$E$ exclusion of the cap and product bound.}
Assume $E\ge E_0$. Coordinate-cost monotonicity and increasing
$c_f$ exclude $a=f$ or $b=f$ from the cheap domain. Put
$v=\sqrt{EP}$. By $Q\ge2\sqrt P$,
\[
 c\ge\frac{4hv/\sqrt E-zv^2}{4-v^2}.
\]
Its derivative numerator in $v\in(0,2)$ is
\[
 \frac{4h}{\sqrt E}(4+v^2)-8zv
 >4z[(v-1)^2+3]>0.
\]
If $P\ge1/E$, the cost is at least $(4h/\sqrt E-z)/3$.
The derivative of the latter function is
$[2(E+1)-\sqrt E]/(3E^{3/2})>0$. At $E_0$, using
$\sqrt{15/26}>3/4$ and $z_0<551/1000$, it is greater than
$1649/3000>1/2$. Consequently
\begin{equation}\label{late:r24:a:high-product}
 P<1/E,\qquad R<1,\qquad D>3.
\end{equation}

\paragraph{Monotonicity on the high-$E$ cost boundary.}
On $c=1/2$, substitute
$Q=[2+E(z-1/2)P]/(2h)$ into \eqref{late:r24:a:N0}. The resulting
quadratic $N_0=A_cP^2+B_cP+C_c$ has coefficients
\[
 A_c=E(z-1/2)-E^2+1<1-3E/2<0,
\]
\[
 B_c=4E-2+\frac{2Ez}{h}(E-2z)>0,
 \qquad C_c=\frac{8z(E-2)}h-4z^2.
\]
Indeed $z<5/6$ and $E\ge26/15$ imply $E-2z>1/15$.
The derivative numerator of $C_0$ in $P$ is
\[
 T=(8A_c+EB_c)P+4B_c+2EC_c
 >\frac{8A_c}E+4B_c+2EC_c
\]
\[
 =8E-12+\frac8E+8z+\frac{8Ez}{h}J(E),
 \qquad J(E)=3E-4-(E+1)z.
\]
The inequality uses $A_c<0$, $B_c>0$ and $P<1/E$ term by term.
Now $J'=2-z-1/E>0$ and
\[
 J(E_0)>\frac65-\frac{41}{15}\frac{551}{1000}
 =-\frac{4591}{15000}>-\frac13.
\]
Since $z/h<1$,
\[
 T>\frac{16}{3}E-12+\frac8E\ge\frac{1088}{585}>0;
\]
the middle expression increases for $E\ge E_0$.
Thus $C_0$ increases strictly in $P$ on the cost boundary.

For precision, in the ordered square $l\le a\le b\le f$, the cost
level is a connected decreasing graph. Along it,
\[
 \frac{db}{da}=-\frac{2h+E(c-z)b}{2h+E(c-z)a},\qquad
 \frac{dP}{da}=\frac{2h(b-a)}{2h+E(c-z)a}\ge0.
\]
The denominator is $Dc_b>0$. The cap is absent, so the graph runs
from its circle endpoint to its diagonal endpoint, possibly as a
singleton. Its least product, and therefore its least $C_0$, occurs
on the circle. Solving the circle cost equation gives
\eqref{late:r24:a:bc}. Its denominator is positive because
$4E-3-2\log E$ equals $1$ at $E=1$ and has derivative $4-2/E>0$.
Since $B_0(E)\le1/2<c_f(E)$, one has $l\le b_c<f$.
The feasible high-$E$ values form an interval starting at $E_0$:
$B_0$ increases, and $B_0(E_0)<1/2$ follows from $4E_0-7<0<2z_0$.

\paragraph{Monotonicity of the circle--cost curve.}
Put $u=h/E$, $\rho=z-u$, $v=\rho/u$, and
$B_*=1+E(E-1/2)b$. On the circle,
\[
 C_0-\frac34=\frac{N}{(4-b)^2},\qquad
 N=u^2(b^2-2b+4)-6b\rho-4\rho^2.
\]
On $b=b_c(E)$,
\[
 b'=-b/E-(E-1/2)b^2=-bB_*/E,\qquad
 u'=E^{-2},\qquad \rho'=u/E.
\]
Differentiation gives
\[
 \frac{dC_0}{dE}=\frac{uJ_c}{E^2(4-b)^3},
\]
\begin{align*}
 J_c={}&(4-b)[2b^2-(4+6E)b+8-8hv]\\
 &+EbB_*[6v(4+b)+8uv^2-6ub].
\end{align*}
The series gives $v>u/2$. On $E_0\le E<9/4$, $1/E\le b\le1$,
one has $Eb\ge1$, $B_*\ge E+1/2$, and
\[
 (J_c)_v>-8h(4-b)+6(E+1/2)(4+b)
 =44-8E+(14E-5)b>0.
\]
Hence $J_c>J_0$, where
\[
 J_0=(4-b)[2b^2-(4+6E)b+8-4h^2/E]
       +hbB_*[12-3b+2u^2].
\]
We verify $J_0>0$ on the entire wedge. Let
$L_*=3E^2-3E+1/2$ and $V_*=12+2u^2$. At fixed $b$,
\[
 (J_0)_E=\mathscr D+4hbB_*u/E^2>\mathscr D,
\]
\begin{align*}
 \mathscr D
 & =-(4-b)(6b+4-4/E^2)+b(1+L_*b)(V_*-3b)\\
 &=-16+16/E^2-(6+4/E+2/E^2)b
                  +(3+L_*V_*)b^2-3L_*b^3.
\end{align*}
Its derivative in $b$ is
\[
 \mathscr D_b=-(6+4/E+2/E^2)+6b+L_*b(2V_*-9b)>0.
\]
Indeed the negative bracket exceeds $-10$, $L_*b\ge L_*/E>3/2$,
and $2V_*-9b\ge15$. Therefore
\[
 \mathscr D(E,b)\ge26-69/E+49/E^2-23/(2E^3)+E^{-4}.
\]
The derivative numerator of this last expression is
$69E^3-98E^2+(69/2)E-4>0$ for $E\ge E_0$, since both
$E^2(69E-98)$ and $(69/2)E-4$ are positive. At $E_0$ its value is
$23077/57122>0$. Explicitly the numerator before cancellation is
\[
 52\cdot456976-2070\cdot17576+22050\cdot676
 -598\cdot3375+101250=369232,
\]
with denominator $913952$. Hence $(J_0)_E>0$ throughout the wedge.

At $E=E_0$ the cubic coefficients of $J_0$ in $b$ are
\begin{align*}
 A_0&=32-16h^2/E,\\
 A_1&=-36-12E+4h^2/E+2h^3/E^2,\\
 A_2&=15+3E+hE(E-1/2)(12+2u^2),\\
 A_3&=-2-3hE(E-1/2).
\end{align*}
They satisfy $A_0>27$, $A_1>-56$, $A_2+A_3>32$, and $A_3<0$.
For the third inequality, $13+3E_0>18$ and
$9hE_0(E_0-1/2)=47619/3375>14$ suffice. Thus
\[
 J_0(E_0,b)>27-56b+32b^2
 =32(b-7/8)^2+5/2>0\qquad(0\le b\le1).
\]
On the other entrance boundary,
\[
 J_0(E,1/E)=30-2E-\frac{87}{2E}
             +\frac{37}{2E^2}-\frac3{E^3}>\frac32,
\]
using $2E<9/2$, $87/(2E)<26$, $37/(2E^2)>3$, and $3/E^3<1$.
For any wedge point, decrease $E$ at fixed $b$ until $E=E_0$ or
$E=1/b$. The path remains in the wedge. Its positive entrance value
and $(J_0)_E>0$ prove $J_0>0$, hence $J_c>0$.

It follows that the circle--cost determinant increases with $E$. At
$E_0$, $b_c(E_0)<f(E_0)$ and strict circle decrease give
\[
 C_0(E,l,b_c(E))\ge C_0(E_0,l_0,b_c(E_0))
 > C_0(E_0,l_0,f(E_0))>\frac{1471}{2000}.
\]
The fixed-product reduction has no other high-$E$ face: the cap is
absent, the circle decreases to $b_c$, and the cost boundary increases
from that same endpoint. This proves the theorem, including the
circle--diagonal singleton and the $E=E_0$ seam.
\end{proof}

\begin{lemma}[Quantitative central descent]
\label{late:r24:a:central-descent}
Let $E\ge E_0$, $l\le a,b\le f(E)$ and $c\le1/2$. Put
$\ell=1/2-c$ and $C_c(E)=C_0(E,l,b_c(E))$. If $\ell>0$, then
\[
 C_0>C_c(E)+\ell/9.
\]
At $E=E_0$ the stronger conclusion is
$C_0>C_c(E_0)+\ell/6$. These bounds include a path starting on either
circle. If $\ell=0$, then $C_0\ge C_c(E)>1471/2000$.
\end{lemma}
\begin{proof}
The preceding theorem gives $P<1/E$, $R<1$ and $D>3$. At fixed
$E,P$, increasing $Q$ until the cost boundary or circle is reached has
rate
\[
 -\frac{dC_0}{dc}=\frac{4z/h-P}{D}>\frac{3}{4E}>\frac13.
\]
No cap is reached at cost at most $1/2$. If the cost boundary is reached
first, its determinant is at least $C_c(E)$ by
Theorem~\ref{late:r24:a:zero-leaf}. Otherwise continue on the circle
by increasing its free coordinate $b$ to $b_c(E)$.
On that circle, with $u=h/E$ and $\rho=z-u$,
\[
 -\frac{dC_0}{dc}
 =\frac{6\rho(b+4)+8\rho^2-6bu^2}{(4-b)k_c},\qquad
 k_c=8h-2u-4\rho>0.
\]
Since $\rho>u^2/2$,
\[
 -\frac{dC_0}{dc}>\frac{3u^2}{k_c}
 >\frac{3u(1-u)}{6+2u^2}>\frac19.
\]
For the middle step use
$k_c<u(6+2u^2)/(1-u)$. The last step is
$27u-29u^2-6>0$ on $11/26\le u<5/9$; this concave quadratic has
values $157/676$ and $4/81$ at the two endpoints. Integration over
the one or two positive-cost path portions gives the first assertion.

At $E_0$, $u=11/26$, the stronger series bound
$\rho>u^2/2+u^3/3+u^4/4$ gives
\[
 \frac{6\rho(b+4)+8\rho^2-6bu^2}{4-b}
 >3u^2+\frac83u^3+\frac{18}{7}u^4.
\]
Here $b\ge15/26>4/7$, so $(b+4)/(4-b)>4/3$ and $4-b<7/2$.
Using $u^2>7/40$, $u^3>7/100$ and $u^4>3/100$, the right side
exceeds $3313/4200>7/9$. Meanwhile
\[
 k_c<\frac{88}{15}-\frac{11}{13}-\frac{121}{338}
 =\frac{14}{3}-\frac7{1690}<\frac{14}{3}.
\]
The circle rate is therefore greater than $1/6$. The fixed-product
rate was already greater than $1/3$, proving the stronger assertion.
If $\ell=0$, the cost-boundary comparison applies directly without
asserting a positive-length descent.
\end{proof}

\begin{lemma}[A central determinant bound]
\label{late:r24:a:r11-central}
Suppose
\begin{gather*}
 1<E<7/4,\quad 0<a,b<6/5,\quad Ea,Eb<2,\\
 Ea>\max(1,b),\quad Eb>\max(1,a),\quad 0<c\le2/5.
\end{gather*}
Then $C_0>3/4-c^2/12$.
\end{lemma}
\begin{proof}
Set
\[
 M=N_0+\frac{(2Qh-EPz)^2}{12}
   =D^2\left(C_0-\frac34+\frac{c^2}{12}\right).
\]
We prove $M>0$. At fixed $E,P$,
\[
 M_Q=2Ph-8z+\frac h3(2Qh-EPz).
\]
The state crossings imply $a,b>1/E$ and $1/E<a/b<E$.
The two geometric upper bounds for $Q$ are
\[
 Q_{\rm circ}=1/E+EP,\qquad
 Q_{\rm ratio}=(E+1)\sqrt{P/E}.
\]
If $v=\sqrt{EP}>1/\sqrt E$, their difference is
$(v-1)(Ev-1)/E$. The circle bound is smaller for $P\le1/E$,
and the ratio bound is smaller for $P\ge1/E$.

For $E\le3/2$, increase $Q$ to this smaller bound. Since
$P<36/25$, its circle endpoint has $Q<2$, and its ratio endpoint has
$v<\sqrt{54/25}<3/2$ and $Q=v+v/E<3$. The entire path therefore has
$2Qh-EPz<6h\le3$, while $z>h/E\ge2h/3$. Thus
\[
 M_Q<h\left(\frac{72}{25}-\frac{16}{3}+1\right)
 =-\frac{109}{75}h<0.
\]
The fixed denominator satisfies $D>4-54/25>0$. No cost restriction
at the algebraic endpoint is needed.

For $E\ge3/2$, the argument below first proves $P<1/E$.
Increase $Q$ only to the lesser of $Q_{\rm circ}$ and
\[
 Q_{\rm cost}=\frac{(2/5)D+EPz}{2h}.
\]
This path preserves $c\le2/5$, so $2Qh-EPz<8/5$. Consequently
\[
 M_Q<h\left(\frac{72}{25}-\frac{32}{7}+\frac8{15}\right)
 =-\frac{608}{525}h<0.
\]
We now establish all endpoint inequalities for these comparisons.

\paragraph{The low-$E$ circle.}
On $a=1/E$, let $u=h/E$, $q=z/u$, and $A=2E-q$. Substitution gives
\begin{equation}\label{late:r24:a:r11-V}
 \frac{M}{u^2}=V(b,q)
 =b^2-2b+\frac{6b(1-q)}u+8q-4q^2
       +\frac{(2+Ab)^2}{12}.
\end{equation}
The series in \eqref{late:r24:a:series} implies
\[
 1+u/2<q\le q_U(u):=1+u/2+u^2/3+\frac{u^3}{4(1-u)}.
\]
Both $q_U$ and $(q_U-1)/u$ increase. On every logarithm-comparison
segment used below, $u\le3/7$ and
$q\le q_U(3/7)=1027/784<21/16$. Hence $0<A\le5/2$.
For $0<b\le1$,
\[
 V_b\le-3+A/3+A^2/6\le-27/24<0,
 \qquad
 V_q=-6b/u+8-8q-b(2+Ab)/6<0.
\]
If $E\le3/2$, then $u\le1/3$, so it suffices to bound
$V(1,q_U(u))$. For $v=q-1$ the exact expression is
\[
 V(1,q)=3+\frac{(3-u)^2}{12(1-u)^2}
 -\left[\frac6u+\frac{3-u}{6(1-u)}\right]v-\frac{47}{12}v^2.
\]
Use $v\le u(6-2u-u^2)/[12(1-u)]$. Multiplication by
$1728(1-u)^2$ gives the lower polynomial
\[
 \mathscr R(u)=1296-4752u+3060u^2+288u^3
                       +352u^4-188u^5-47u^6.
\]
The least degree-six Bernstein coefficient on \([0,1/3]\), by
\eqref{cv:bernstein-coefficients}, is \(48241/729>0\).
Thus \(M>0\) on the low-\(E\) circle.

\paragraph{The low-$E$ ratio boundary.}
Let $v=\sqrt{EP}\in[1,3/2)$ and $k=(E+1)/E$. The boundary is
$(a,b)=(v,v/E)$ and
\[
 M(v)=A_4v^4+A_3v^3+A_2v^2-8kzv+8z-4z^2,
\]
\[
 A_4=-hk/E+z^2/12<0,\quad
 A_3=2hk/E-khz/3,\quad
 A_2=4h/E+2z+k^2h^2/3.
\]
In $M'(v)$, after separating the positive $z^2$ term, the coefficient
of $z$ is $-khv^2+4v-8k<0$. Replace $z$ by its strict upper bound
$hk/2$ and omit only that positive $z^2$ term. At $v=1$ this yields
\[
 M'(1)>\frac{h^2[(E+1)^2-12h]}{6E^2}>0,
\]
since $E^2-10E+13\ge1/4$ on $[1,3/2]$. At $v=3/2$, writing
$p=1/E\in[2/3,1]$, it yields
\[
 \frac{M'(3/2)}h>-1+7p-4p^2-\frac{k^2h}{8}
 \ge\frac{17}{9}-\frac14>0.
\]
The exact second derivative at $1$ is
\[
 M''(1)=z^2+(4-2kh)z+8h/E+2k^2h^2/3>0
\]
because $kh\le5/6$. Since $M''$ is a downward quadratic, positive
at $1$, it has at most one zero to the right, crossed from positive
to negative. Thus $M'$ first increases and then possibly decreases;
its positive endpoint values imply $M'>0$ on $[1,3/2]$.
Consequently $M(v)\ge M(1)>0$, the last point being the preceding
circle boundary with scalar names exchanged.

\paragraph{The high-$E$ product restriction.}
For $E\ge3/2$, the same diagonal cost calculation used in
Theorem~\ref{late:r24:a:zero-leaf} shows that $P\ge1/E$ would imply
\[
 c\ge\frac{4h/\sqrt E-z}{3}.
\]
Its derivative in $E$ is positive. At $E=3/2$, the rational bounds
$\sqrt{2/3}>13/16$ and $\log(3/2)<5/12$ make it greater than
$29/72>2/5$. Hence $P<1/E$, as required above.

\paragraph{The cost-qualified high-$E$ circle.}
On the circle,
\[
 c=\frac{2h/E+b(2h-z)}{4-b}.
\]
At fixed \(b\) its numerator increases with \(E\), since the derivative
is \(2/E^2+b(2-1/E)>0\). The bound
\(z<2h/(E+1)+h^3/[6E(E+1)]\) is the first term of the
artanh series with its geometric remainder. At \(E=3/2,8/5\) it
implies, respectively, \(c(3/2,2/3,47/50)>2/5\) and
\(c(8/5,5/8,19/25)>2/5\), by substitution in the displayed
circle formula. Since that formula is increasing in \(b\), the
following \(b_*\) are upper bounds for cheap circle points.
\[
\begin{array}{c|c|c|c|c}
 E & b_* & q_* \text{ for } q
 & v_* \text{ for } (q-1)/u & \text{lower bound for }V\\ \hline
 {[3/2,8/5]} &47/50&1607/1280&109/160&15107/4000000\\
 {[8/5,7/4]} &19/25&21/16&3/4&25237/120000
\end{array}
\]
Here the logarithm bounds follow from \(q_U\) at the right
endpoints. For each row with left endpoint \(E_-\), evaluate
\eqref{late:r24:a:r11-V} at \(b_*\), using \(V_b<0\), and bound it
below by
\[
 b_*^2-2b_*-6b_*v_*+4-4(q_*-1)^2
 +\frac{[2+(2E_--q_*)b_*]^2}{12}.
\]
Both square bases are positive. The rational expression is at least
the stated lower bound in each row, including the common seam.

\paragraph{The high-$E$ cost face.}
If the circle endpoint would have cost greater than $2/5$, stop the
$Q$ comparison at $c=2/5$. On that plane substitute
\[
 Q=\frac{c(4-EP)+EPz}{2h},\qquad c=2/5.
\]
The coefficient of $P^2$ in $M$ is
\[
 E(z-c)-h(E+1)+\frac{c^2E^2}{12}
 <-h-Ec+\frac{c^2E^2}{12}<0.
\]
Thus $M$ is concave in $P$. Circle cost increases strictly in $b$;
at $b=0$ it is $h/(2E)<2/5$, and at $b=1$ it exceeds $2/5$.
There is a unique $b_{\rm cap}\in(0,1)$ of cost $2/5$.
The present case has $P\ge P_{\rm cap}=b_{\rm cap}/E$ and $P<1/E$.
Enlarge this scalar interval to $[P_{\rm cap},1/E]$.
At the first endpoint the preceding circle estimate applies.
At the second, artificial endpoint $D=3$; putting $w=z/h$ gives
\begin{align*}
 M={}&\frac{6}{5E}+\frac{h(3E-1)}{E^2}
       +h(10+1/E)w-\frac{24}{5}w-4Eh w^2+\frac3{25}.
\end{align*}
For $w>1/E$ its derivative in $w$ is less than
$h(2+1/E)-24/5<0$. Since $z/h$ decreases with $E$ and
$z/h<5/6$ for $E\ge3/2$, replace $w$ by $5/6$ to obtain
\[
 M>H(E):=-\frac{25E^2}{9}+\frac{100E}{9}
              -\frac{419}{50}-\frac{109}{30E}+\frac1{E^2}.
\]
Here
\[
 H'(E)=\frac{50}{9}(2-E)+\frac{109E-60}{30E^3}>0,
 \qquad H(3/2)=53/900>0.
\]
Both concavity endpoints are positive, so $M>0$ on the whole cost
face. This covers all comparison branches and proves the lemma.
\end{proof}

\begin{theorem}[High-$E$ positive-leaf compensation]
\label{late:r24:a:high-e}
Suppose $E\ge E_0$, $a,b\ge1/E$, $Ea,Eb<2$, $D>0$, and the
exact nonnegative-mean leaf branches satisfy their central signs,
$0\le d_A<a$, $0\le d_B<b$, and
$\alpha+c+\beta\le1/2$. Then
\[
 \Gamma>\frac{1471}{2000}.
\]
Zero leaf budgets, total slack, and either lower-crossing equality are
allowed. Thus the statement applies to the exact closed comparison
faces as well as to actual positive leaves.
\end{theorem}
\begin{proof}
Coordinate-cost monotonicity \eqref{late:r24:a:coordinate-cost} and
Lemma~\ref{late:r24:a:corner} exclude $a\ge f$ or $b\ge f$ even if
an initial coordinate was outside the zero-leaf square: decrease that
coordinate to $f$ and the other to $l$, preserving $D>0$.
Therefore
\[
 l\le a,b<f(E)\le f(E_0)<819/1000<5/6.
\]
The product comparison \eqref{late:r24:a:high-product} now gives
$R<1$, $D>3$. Also
\[
 c\ge B_0(E)\ge B_0(E_0)>7147/17800>2/5,
 \qquad 0\le\ell:=1/2-c<1/10.
\]
The difference between $7147/17800$ and $2/5$ is $27/17800$.
Since $B_0(E)\le1/2$ and $z<5/6$,
$4E-7\le2z<5/3$, so $E<13/6$. We may also use
\[
 11/26\le u:=1-1/E<5/9,\qquad h<5/4,\qquad a,b>4/9.
\]
Coordinate-cost monotonicity and the circle cost equation imply
$b\le b_c(E)$ and the analogous bound for $a$. Since
$z<(E-E^{-1})/2$,
\[
 b\le b_c(E)<B(E):=\frac4{T(E)},\qquad
 T(E)=3E^2-3E+1,\qquad B(E_0)=300/361<5/6.
\]

\paragraph{Uniform small-budget leaf estimates.}
For the $A$ leaf write $r=1-a$, $q=2-a$, $V^2=aD/(Eb)$ and $d=d_A$.
Here $r>1/6$, $q>7/6$ and
$V^2\le a(4-a)\le95/36$.
The exact cost is bounded below by
$d/2+(3/2)dq/[V^2(q+d)]$, an increasing function of $d$.
At $d=1/10$ its cost/width ratio is at least
\[
 \frac12+\frac32\frac{35}{38}\frac{36}{95}
 =\frac12+\frac{189}{361}>1.
\]
Thus $\alpha<1/10$ forces $d<1/10$. Set
$\tau=35/38$ and $\kappa_0=31/30$. The true mean obeys
$e_A\le rd+d^2/3$; hence, for a positive budget,
\begin{equation}\label{late:r24:a:high-mean-ratio}
 \frac{x}{\alpha}
 <\frac{2(\kappa_0-a)}{aD+3Eb\tau}<\frac3{10}.
\end{equation}
The last inequality uses $aD/2>2/3$, $Eb\ge1$,
$2/3+105/76>2$, and $\kappa_0-a<53/90<3/5$.
For zero budget $x=0$; the corresponding right-leaf assertion is
$y<3\beta/10$ at positive budget and $y=0$ otherwise.

\paragraph{The linear mean losses.}
Let
\[
 M_x=-DL_x=4(z-1+b)+ab(2-E-b).
\]
It increases in $b$ and decreases in $a$ on the enlarged rectangle
needed to compare with $(a,b)=(5/6,l)$, since
$(M_x)_b=4+a(2-E-2b)>0$ and
$(M_x)_a=b(2-E-b)<0$. The series gives
\[
 M_x\ge4(z-u)-\frac56u^2
 >\frac76u^2+\frac43u^3+u^4
 \ge\frac{937266}{2741856}>\frac13.
\]
The last reserve over one third is $23314/2741856$, so $M_xD>1$.

At fixed $E,b$, put $H=(\kappa_0-a)M_x/D$. Exact differentiation gives
\[
 \left(\frac{M_x}{D}\right)_a
 =\frac{4bh[b-1+fh/2]}{D^2},
\]
\[
 D^2H_a=-M_xD+(\kappa_0-a)4bh[b-1+fh/2]<0.
\]
If the bracket is positive, its last term is less than
$(53/90)(25/6)(17/48)=22525/25920<1$; otherwise the sign is
immediate. Moreover $(aD+3Eb\tau)_a=4-2Eab>2$.
Thus the sufficient upper bound for $(-L_x)x/\alpha$ obtained from
\eqref{late:r24:a:high-mean-ratio} decreases in $a$. Lower $a$ to $l$;
this path remains cheap by coordinate-cost monotonicity. At its endpoint,
\[
 M_x=4z-4+3b+(2b-b^2)/E<A(E)+3b,
 \qquad A(E)=2E-4-1/E<0.
\]
Therefore
\[
 \frac{(-L_x)x}{\alpha}
 <\frac{31E-30}{15(4-b)}
       \frac{A(E)+3b}{4+b(3E^2\tau-1)}.
\]
Both factors increase in $b$. For the second the derivative numerator is
$12-A(E)(3E^2\tau-1)>0$; its numerator is positive because it exceeds
$M_x$. Replace $b$ by $B(E)$.
The first factor becomes
\[
 O(E)=\frac{(31E-30)T(E)}{180E(E-1)}.
\]
Its derivative exceeds $[93-31/(E-1)^2]/180>0$, and hence
\[
 O(E)<O(13/6)=22969/32760<5/7.
\]
For the second factor to be below $7/50$, clearing positive
denominators gives the polynomial condition
\[
 p(E)=-300E^4+(984+84\tau)E^3-634E^2-550E+50>0.
\]
As $\tau>9/10$, it suffices to bound
\[
 p_0(E)=-300E^4+(5298/5)E^3-634E^2-550E+50.
\]
Its second derivative is less than
$-3600E^2+6408E-1268$, which is $-4884/5$ at $E_0$ and decreases
thereafter. Thus $p_0$ is strictly concave. Its two endpoint values are
\[
 p_0(26/15)=97494/50625>0,\qquad
 p_0(13/6)=312456/6480>0.
\]
Their numerator sums are, respectively,
\[
 -137092800+279352944-96431400-48262500+2531250=97494,
\]
\[
 -42841500+69838236-19286280-7722000+324000=312456.
\]
The positive endpoint chord proves $p_0>0$ throughout.
Consequently $(-L_x)x<\alpha/10$ for positive $\alpha$.
Starting with the separately derived $L_y$ formula gives
$(-L_y)y<\beta/10$ for positive $\beta$. Zero budgets contribute zero.

\paragraph{The shared budget and the central compensation.}
Let $s=\alpha+\beta\le\ell$. Since $P>0$,
$x^2+(2-P)xy+y^2\le(x+y)^2$. The uniform statement, including
$s=0$, is
\[
 \Gamma\ge C_0-s/10-(9/100)s^2
          \ge C_0-\ell/10-(9/100)\ell^2.
\]
If $\ell>0$, then $\ell<1/10$ makes the last expression strictly
larger than $C_0-(109/1000)\ell$. Lemma
\ref{late:r24:a:central-descent} gives
\[
 \Gamma>C_c(E)+\left(\frac19-\frac{109}{1000}\right)\ell
       =C_c(E)+\frac{19}{9000}\ell>\frac{1471}{2000}.
\]
If $\ell=0$, both budgets vanish, and
$\Gamma=C_0\ge C_c(E)>1471/2000$ by
Theorem~\ref{late:r24:a:zero-leaf}. This also treats a zero-total-leaf
comparison with slack without asserting the false inequality $C_0>C_0$.
\end{proof}

\begin{lemma}[Signed-leaf width bounds]\label{late:r24:a:signed-width}
Consider a state with the strict crossings \eqref{late:r24:crossings},
the budget \eqref{late:r24:contrary-budget}, and true positive leaf
means satisfying the central signs \eqref{late:r24:central-signs}.
A leaf of coordinate $v\ge1$ and cost $A$ satisfies
\[
 v<5/4,\qquad d<7A/5,
 \qquad d\ge z+3(v-1).
\]
\end{lemma}
\begin{proof}
Write $\kappa=v-1\ge0$. The signed true-mean bound is
\[
 0<e\le d^2/3-\kappa d.
\]
For clarity it follows directly from the weighted mean even when its
quadratic integrand is not monotone. With
$M(s)=s^2/3-\kappa s$ and $\omega(s)=(q+s)^{-2}$, integration by
parts expresses $e$ as the weighted average of $M(d)$ and $M(s)$ with
positive weights $d\omega(d)$ and $-s\omega'(s)ds$. Since
$(s/d)M(d)-M(s)=s(d-s)/3\ge0$, positivity of that average implies
$M(d)>0$ and $e\le M(d)$. Thus $\kappa<d/3$.
The cost gives $d<2A<1$, whence $v<4/3$, $q=2-v>2/3$,
$V^2<v(4-v)<32/9$, and $q/(q+d)>2/5$. Therefore
\[
 A/d>\frac12+\frac32\frac25\frac9{32}
 =\frac{107}{160}>\frac23.
\]
It follows that $d<3A/2\le3/4$ and $v<5/4$. Reinserting
$q>3/4$, $V^2<55/16$, and $q/(q+d)>1/2$ gives
\[
 A/d>\frac12+\frac{12}{55}=\frac{79}{110}>\frac57.
\]

The central sign is $e\ge v[\Phi(E)+\kappa h]$, where
\[
 \Phi(E)=2h^2\int_0^1\frac{t(1-t)}{1+ht}\,dt
 \ge\frac{h^2}{3E}\ge\frac{z^2}{3}.
\]
Since $v\ge1$ and $h\ge z$, it gives
$e\ge z^2/3+\kappa z$. Subtracting the lower bound from the upper
bound yields
\[
 (d+z)(d-z-3\kappa)\ge0.
\]
As $d+z>0$, the final asserted width inequality follows. The same
formulas include $\kappa=0$.
\end{proof}

\begin{theorem}[The fully signed sector]
\label{late:r24:a:fully-signed}
Suppose an actual contrary-budget state has $a,b\ge1$, the original
strict crossed and upper state inequalities, true positive means,
central signs, and $\alpha+c+\beta\le1/2$. Then
\[
 \Gamma>\frac{737}{1000}.
\]
\end{theorem}
\begin{proof}
Lemma~\ref{late:r24:cheap-domain} gives $R<3/2$ and $D>5/2$.
Put $\ell=\alpha+\beta$, $\kappa=a-1$, and $\lambda=b-1$.
Since $a,b\ge1$, $E\le R<3/2$. Coordinate descent to $(1,1)$ gives
\[
 c\ge\frac{4h-Ez}{4-E}>h,
\]
whose strict difference from $h$ is $E(h-z)/(4-E)>0$.
Lemma~\ref{late:r24:a:signed-width} gives
\begin{equation}\label{late:r24:a:signed-sum}
 \ell>\frac{10}{7}z,\qquad
 \kappa+\lambda<\frac7{15}\ell-\frac23z.
\end{equation}
If $h\ge1/4$, then $z\ge\log(5/4)>1/5$, and the total cost
would exceed $1/4+2/7=15/28>1/2$. Thus $E<5/4$.
The original crossing $a<Eb$ and $b<5/4$ imply
$\kappa-\lambda<hb<(5/4)h$. Combining this with
\eqref{late:r24:a:signed-sum} and $\ell<1/2-h$ gives
\[
 \kappa<\frac7{60}+\frac{47}{120}h-\frac z3
 <\frac7{60}+\frac h8<\frac{71}{480}<\frac15,
\]
where $z>h/E>4h/5$ was used. The other original crossing yields
the same bound for $\lambda$, so $a,b<6/5$.

We also need the central-cost hypothesis of
Lemma~\ref{late:r24:a:r11-central}. Since $P\ge Q-1$ and $P\ge1$,
\[
 c<\frac{h(2Q-P)}D\le\frac{h(Q+1)}D<\frac{97}{75}h;
\]
here $Q<2+(7/15)\ell\le67/30$. On the other hand,
$\ell>(10/7)h/(1+h)$ implies
\[
 h<\frac{7\ell}{10-7\ell}.
\]
If $c\ge4/15$, then $\ell\le7/30$, and these two bounds give
\[
 c<\frac{97}{75}\frac{49}{251}
 =\frac{4753}{18825}<\frac{5020}{18825}=\frac4{15},
\]
a contradiction. Thus $0<c<4/15<2/5$. All hypotheses of
Lemma~\ref{late:r24:a:r11-central} have now been proved, so
$C_0>3/4-c^2/12$.

The true signed means and the width bounds imply
\[
 x<\frac{98}{375}\alpha^2<\frac4{15}\alpha^2,
 \qquad y<\frac4{15}\beta^2.
\]
The exact coefficients are
\[
 -L_x=\lambda\frac{4-P}{D}+\frac{4z-Ph}{D},\qquad
 -L_y=\kappa\frac{4-P}{D}+\frac{4z-Ph}{D}.
\]
Since $(4-P)/D<6/5$, $4z-Ph<3h$, $h<c$,
$\lambda<7\beta/15$ and $\kappa<7\alpha/15$, one obtains
\[
 -L_x<\frac65c+\frac{14}{25}\beta,\qquad
 -L_y<\frac65c+\frac{14}{25}\alpha.
\]
The coefficients $L_x$ and $L_y$ are negative by
Lemma~\ref{late:r24:cheap-domain}.
Using $\alpha^2+\beta^2\le\ell^2$ and
$\alpha\beta\le\ell^2/4$ in the exact determinant gives
\[
 \Gamma>\frac34-B(c,\ell),\qquad
 B=\frac{c^2}{12}+\frac8{25}c\ell^2
                  +\frac{14}{375}\ell^3+\frac{16}{225}\ell^4.
\]
All coefficients are positive. Since $\ell\le1/2-c$ and
$\ell\le1/2$,
\begin{align*}
 B&\le\frac{c^2}{12}
 +\left(\frac8{25}c+\frac{41}{1125}\right)(1/2-c)^2\\
 &=\frac{41}{4500}+\frac{49}{1125}c
            -\frac{901}{4500}c^2+\frac8{25}c^3=:F(c).
\end{align*}
The derivative has the sign of $2160c^2-901c+98$, whose
discriminant is $-34919<0$. Therefore $F$ increases, and
\[
 F(c)<F(4/15)=\frac{12713}{1012500}<\frac{13}{1000}.
\]
The last cross-product reserve is $449500>0$.
Hence $\Gamma>3/4-13/1000=737/1000$.
All signed-depth equalities, including $a=b=1$, have been covered by
the same true-mean formulas. No saturation of the leaf budget was used.
\end{proof}

\begin{theorem}[The mixed sector with a common budget]
\label{late:r24:a:mixed}
Suppose
\[
 0<a\le1<b,\quad E>1,\quad Ea>b,\quad Eb>1,\quad
 Ea,Eb<2,\quad R<3/2,
\]
with the remaining original crossed inequalities retained. Use the exact
true leaf means on their nonnegative or positive branches, with
\[
 e_A\ge\max(0,ah(a-f)),\qquad e_B\ge bh(b-f)>0,
 \qquad \alpha+c+\beta\le1/2.
\]
The exact widths lie below their coordinate barriers. Zero $A$ width
is permitted only when its central minimum is zero. Then
\[
 \Gamma>\frac{1471}{2000}+\frac1{97500}.
\]
The theorem permits arbitrary exact allocations and slack at the retained
central state.
\end{theorem}
\begin{proof}
Write $r=1-a\ge0$, $\kappa=b-1>0$, $p=d_A$, $s=d_B$,
$q_A=2-a$, $q_B=2-b$. The exact costs imply $p,s<1$ at positive
widths. The crossing gives $R>b^2$, and hence
\[
 1<b<\sqrt{3/2}<5/4,\quad q_B>3/4,\quad
 V_B^2=\frac{b^2(4-R)}R<4-b^2<3.
\]

\paragraph{The signed minimum and the shared width budget.}
Coordinate-cost monotonicity gives
\[
 c>c_1(E):=c(E,1/E,1)
 =\frac{2h(1+1/E)-z}{3}>\frac{h(1+2/E)}3.
\]
The last expression increases with $E$ and equals $308/585>1/2$
at $E_0$. Thus $E<E_0$ before the next comparison is used.
Let $H_b=bh(b-f(E))$. At width $h$,
\[
 \frac{H_b-e_b(h)}h
 =T(b):=b^2+b-2-bf(E)+f(1+h/(2-b)).
\]
Here $T(1)=0$, and \eqref{late:r24:a:f-derivatives} gives
\[
 T'(b)=2b+1-f(E)+f'(1+h/q_B)h/q_B^2
 >2-176/405>0.
\]
Therefore $e_b(h)<H_b$. On the positive-excess branch $e_b$ increases
strictly, so every retained signed width satisfies $s>h$. This remains
valid if the tested width $h$ initially lies before that branch: a
retained width at most $h$ would put $h$ on the increasing positive
branch and give the opposite inequality.

Its invariant is at least
\[
 K_{\min}(b)=\frac{4+Eb^2-2Ebz/h}{2}>3/2.
\]
Since $V_B^2<q_B(2+b)$,
\[
 \beta>\frac h2+
 \frac{K_{\min}(b)h}{(2+b)(q_B+h)}.
\]
The positive numerator of the ratio increases with $b$, with derivative
$E(b-z/h)>0$, while its denominator has derivative $h-2b<0$.
Evaluation at $b=1$ gives the exact identity
\[
 \beta>\frac h2+\frac{K_{\min}(1)h}{3E}=c_1(E).
\]
Thus $1/2\ge c+\beta>2c_1(E)$. Since
$h(1+2/E)/3=33/130>1/4$ at $E=13/10$, it follows that
\begin{equation}\label{late:r24:a:mixed-h}
 E<13/10,\qquad c>11h/13,\qquad \beta>11h/13,
 \qquad 0<h<13/44.
\end{equation}

For the unsigned leaf,
$V_A^2=4a/(Eb)-a^2<4a-a^2\le3$ and $q_A\ge1$.
At positive width, its exact cost first gives $\alpha/p>3/4$,
so $p<4\alpha/3\le2/3$. Reinserting this gives
$\alpha/p>4/5$, hence $p<5\alpha/4$. At zero budget use the
non-strict endpoint $p=0\le5\alpha/4$.
For the signed leaf, positivity and Jensen's bound in
\eqref{late:r24:a:f-bounds} give
\[
 q_B>f(1+s/q_B)>\frac3{3+s/q_B},\qquad
 s>3\kappa,\qquad q_B>1-s/3,
\]
\[
 0<e_B<s^2/3-\kappa s.
\]
Indeed $1-f(1+w)<w/(3+w)<q_Bw/3=s/3$ on this positive branch.
Its exact cost therefore satisfies
\[
 \frac\beta s>\frac12+\frac12\frac{1-s/3}{1+2s/3}
 =\frac{1+s/6}{1+2s/3}.
\]
Starting from $s<1$ gives $s<10\beta/7\le5/7<3/4$;
reinsertion gives $s<4\beta/3$. Combining these estimates with
\eqref{late:r24:a:mixed-h}, define
\begin{equation}\label{late:r24:a:mixed-triangle}
 P_h=\frac58-\frac{55h}{26},\qquad
 W_h=\frac23-\frac{44h}{39},\qquad
 0\le p<P_h,\quad p+s<W_h.
\end{equation}
For example, $\alpha<1/2-22h/13$ gives the first bound;
$p+s<(4/3)(\alpha+\beta)<(4/3)(1/2-11h/13)$ gives the second.
They remain strict when $p=0$, since the signed leaf and central cost
bounds are strict.

\paragraph{A central comparison valid on this entire sector.}
At fixed $E,P$, $(N_0)_Q=2Ph-8z<0$ because $P<3/(2E)$.
Increase $b$ and decrease $a=P/b$ until
\[
 b_*=v=\sqrt{EP},\qquad a_*=v/E.
\]
Since $v>b>1$ and $a_*<a\le1$, the sum $Q$ increases. Thus $C_0$
decreases along this central algebraic comparison. Write $v=1+t$;
then $0\le t\le h$ because $a_*\le1$. Put $u=h/E$,
$\rho=z-u$, and $A_*=2u-u^2$. The exact numerator at that endpoint
is
\begin{align}
 N={}&3u^2-6\rho-4\rho^2
 +(6u^2-12\rho+8u\rho)t\notag\\
 &+(6u+2\rho)t^2-2A_*t^3-A_*t^4.
 \label{late:r24:a:mixed-central-poly}
\end{align}
For an expansion check, before replacing $v$ by $1+t$ it is
\[
 -A_*v^4+2A_*v^3+(6u+2\rho)v^2
 -8(2-u)(u+\rho)v+8(u+\rho)-4(u+\rho)^2.
\]
The positive logarithm series gives
\[
 \rho\le u^2/2+\frac{u^3}{3(1-u)},
\]
\[
 3u^2-6\rho\ge-\frac{2u^3}{1-u},\qquad
 6u^2-12\rho+8u\rho\ge-\frac{4u^4}{3(1-u)}.
\]
The quadratic-through-quartic part in
\eqref{late:r24:a:mixed-central-poly} is nonnegative, since
$t\le h<3/10$ and its bracket exceeds
$6u-2u(69/100)>0$. Also $\rho<h^2/(2E)$ by
$z<(E-E^{-1})/2$. It follows that
\[
 N>-\left(\frac{2h^3}{E^2}+\frac{h^4}{E^2}
                    +\frac{4h^5}{3E^3}\right).
\]
As $2+h+4h^2/(3E)<2E$ for $h<3$, and $D^2>25/4$, we obtain
\begin{equation}\label{late:r24:a:mixed-central}
 C_0>\frac34-\frac8{25}\frac{h^3}{1+h}.
\end{equation}
No cost condition on the crossed-boundary comparison point was used.

\paragraph{The coupled linear losses.}
The original crossing gives $0<\kappa<h$ and $Er<h-\kappa$.
The exact means imply
\[
 \frac{e_A}a\le
 \frac{(h-\kappa)p+Ep^2/3}{1+\kappa},\qquad
 \frac{e_B}b<\frac{s^2/3-\kappa s}{1+\kappa}.
\]
The first comparison is strict at positive $p$; at $p=0$ both sides
are zero. The coefficients satisfy
\[
 -DL_x<(4-P)(h+\kappa),\qquad
 -DL_y<(4-P)(h-r).
\]
Both $L_x,L_y<0$ at the retained strict central state. This follows
from the coefficient proof in Lemma~\ref{late:r24:cheap-domain},
whose state part uses only $R<3/2$ and the crossed and upper bounds.
Moreover
\[
 \frac{4-P}{D}=1+\frac{Ph}{D}<1+\frac35u<\frac{74}{65},
 \qquad K_c:=\frac{4-P}{D^2}<\frac{148}{325}.
\]
Substitute the true-mean bounds and replace $h-r$ by $h$ in a
nonnegative term. Dropping the positive denominator $1+\kappa$ is
legitimate: before expansion its numerator is
\[
 (h+\kappa)[(h-\kappa)p+Ep^2/3]
                 +h[s^2/3-\kappa s]>0.
\]
The resulting bracket is
\[
 h^2p+(E/3)hp^2+(h/3)s^2
             +\kappa(Ep^2/3-hs)-\kappa^2p.
\]
For $p>0$ complete the square on $\kappa\ge0$:
\[
 \kappa(Ep^2/3-hs)-\kappa^2p
 \le\frac{(Ep^2/3-hs)_+^2}{4p}\le\frac{E^2p^3}{36}.
\]
For $p=0$ the left side is $-\kappa hs<0$, so the same upper bound
holds directly. Since $E<13/10$,
\begin{equation}\label{late:r24:a:mixed-linear}
 (-L_x)x+(-L_y)y<\frac{148}{325}\mathscr T,
 \quad
 \mathscr T=h^2p+\frac{13}{30}hp^2+\frac h3s^2
                                      +\frac{169}{3600}p^3.
\end{equation}

\paragraph{A uniform polynomial bound for the common width triangle.}
The polynomial $\mathscr T$ increases with $s$. Replace $s$ by
$W_h-p$; its second derivative in $p$ is
$(23/15)h+(169/600)p>0$. Hence its maximum over $0\le p\le P_h$
is at an endpoint. Put $\tau=44h/13\in[0,1]$. Then
\[
 h=13\tau/44,\quad P_h=5(1-\tau)/8,\quad
 W_h=(2-\tau)/3,\quad W_h-P_h=(1+7\tau)/24.
\]
At $p=0$,
\[
 \mathscr T=\frac{13\tau(2-\tau)^2}{1188}
 \le\frac{104}{8019}<\frac1{60}.
\]
The maximum of $\tau(2-\tau)^2$ is at $\tau=2/3$.
At $p=P_h$, the exact expansion in a nonnegative cubic basis is
\[
 \mathscr T=a_0(1-\tau)^3+a_1\tau(1-\tau)^2
              +a_2\tau^2(1-\tau)+a_3\tau^3,
\]
\[
 a_0=845/73728<23/2000,
 \quad a_1=845/16896+13/76032<51/1000,
\]
\[
 a_2=845/15488+208/76032<58/1000,
 \quad a_3=13/1188<11/1000.
\]
For the second parts of the coefficients use
$\tau(1+7\tau)^2=\tau(1-\tau)^2+
16\tau^2(1-\tau)+64\tau^3$. The two middle sums are
$7631/152064$ and $23959/418176$, so their displayed rational
comparisons can also be checked directly. Thus
\[
 \mathscr T<U(\tau):=\frac{23+33\tau-19\tau^2-15\tau^3}{2000}.
\]
With $X=\tau-1/2$, $|X|\le1/2$, the exact expansion gives
\begin{align*}
 U&=\frac{263}{16000}+\frac{11}{8000}X
       -\frac{83}{4000}X^2-\frac3{400}X^3\\
 &\le\frac{263}{16000}+\frac{11}{8000}X-\frac{17}{1000}X^2
 \le\frac{263}{16000}+\frac{121}{4352000}<\frac1{60}.
\end{align*}
The final comparison uses $1/60-263/16000=11/48000$.
Equation~\eqref{late:r24:a:mixed-linear} therefore proves
\begin{equation}\label{late:r24:a:mixed-linear-final}
 (-L_x)x+(-L_y)y<37/4875.
\end{equation}

\paragraph{The joint central and quadratic loss.}
Summing the two true-mean bounds and dropping only favorable negative
terms gives
\[
 x+y<\frac25\left[hp+\frac{13}{30}p^2+\frac{s^2}{3}\right].
\]
Again substitute $s=W_h-p$ and maximize the convex polynomial in $p$.
The endpoint $P_h$ dominates $0$, because their difference is
\[
 P_h[h+(23/30)P_h-(2/3)W_h]
 =P_h[5/144+(61/468)h]\ge0.
\]
The larger value is
\[
 Z(h)=hP_h+(13/30)P_h^2+(W_h-P_h)^2/3
 =587/3456-(2771/5616)h+(5423/36504)h^2.
\]
On $0\le h\le3/10$,
\[
 Z(h)-\left(17/100-4h/9\right)
 =-13/86400-(275/5616)h+(5423/36504)h^2
\]
\[
 \le-13/86400-(1606/365040)h<0.
\]
The comparison affine function is positive on this interval.
Since $P>0$,
\[
 x^2+(2-P)xy+y^2\le(x+y)^2
 <\frac4{25}(17/100-4h/9)^2.
\]
Combine this with \eqref{late:r24:a:mixed-central}. The sufficient
joint deficit is
\[
 B(h)=\frac8{25}\frac{h^3}{1+h}
                   +\frac4{25}(17/100-4h/9)^2.
\]
It is convex on $[0,3/10]$: the second derivative of
$h^3/(1+h)=h^2-h+1-(1+h)^{-1}$ is
$2-2/(1+h)^3\ge0$. Its endpoint values are
\[
 B(0)=289/62500<69/10000,
\]
\[
 B(3/10)=216/32500+121/562500
 =50173/7312500<69/10000.
\]
The last cross-product reserve is $2832500>0$. Consequently
\[
 \Gamma>\frac34-\frac{69}{10000}-\frac{37}{4875}
 =\frac{1471}{2000}+\frac1{97500}.
\]
Indeed the residual difference is
$19/2500-37/4875=(741-740)/97500$.
All bounds retained the same two leaf budgets. They apply unchanged
under allocation slack, and do not assume $a\le f$ or a saturated
central sign.
\end{proof}

\begin{corollary}[The strict-crossing coordinate-one face]
\label{late:r24:a:unit}
The conclusion of Theorem~\ref{late:r24:a:mixed} also holds when
$0<a\le1=b$, provided the original strict crossed and upper state
inequalities, $R<3/2$, exact central signs and shared budget are retained.
It therefore also holds for $a=1$, $0<b<1$, in the opposite scalar
orientation. It does not assert the same conclusion on a new crossed
circle equality.
\end{corollary}
\begin{proof}
Set $b=1$ and $\kappa=0$. Strict $Ea>1$ gives $R=Ea>1$,
$V_B^2=D/R<3$, and $a>1/E$. Therefore $c>c_1(E)$, and the same
preliminary comparison forces $E<E_0$. The exact identity
\[
 e_1(h)=h(1-f(E))=H_1>0
\]
shows that the retained signed-side width has $s\ge h$, with equality
allowed at its true minimum. Even at equality its exact cost satisfies
\[
 \beta=\frac s2+\frac{K_Bs}{V_B^2(1+s)}
 >\frac s2+\frac{K_{\min}(1)s}{3(1+s)}
 \ge\frac h2+\frac{K_{\min}(1)h}{3(1+h)}=c_1(E).
\]
The first comparison remains strict because $s\ge h>0$ and $V_B^2<3$.
Thus all four inequalities \eqref{late:r24:a:mixed-h} persist.

The unsigned width bootstrap uses only $V_A^2<3$, $q_A\ge1$.
For the other leaf the true mean obeys $0<e_B<s^2/3$ and
$q_B=1>1-s/3$; both exact cost-ratio steps give $s<4\beta/3$.
Therefore the same strict common-width triangle
\eqref{late:r24:a:mixed-triangle} holds, including $p=0$.
The central comparison remains valid: its crossed endpoint has
$b_* =\sqrt R>1$, $a_*=b_*/E<a\le1$, so
$0<t=b_*-1\le h$. The exact polynomial
\eqref{late:r24:a:mixed-central-poly} and its bound therefore apply.

The coefficient inequalities reduce to
\[
 -DL_x<(4-P)h,\qquad -DL_y<(4-P)(h-r),
 \qquad Er<h.
\]
The negative-coefficient proof still applies, or follows by lowering
$a$ to $1/E$ and using
$h(5E-1)/E^2-4\log E<0$. The coupled-depth term is now exactly zero,
which is bounded by the same nonnegative majorant $E^2p^3/36$.
All subsequent polynomial estimates in the mixed proof depend only on
\eqref{late:r24:a:mixed-h}, \eqref{late:r24:a:mixed-triangle},
$E<13/10$, and these coefficient bounds. They prove
\eqref{late:r24:a:mixed-linear-final} and the same joint deficit
below $69/10000$. The identical strict rational conclusion follows.
The opposite orientation uses the two already derived symmetric scalar
laws, with no physical relabeling assertion.
\end{proof}

\begin{lemma}[The zero-minimum exact profile]
\label{late:r24:a:zero-profile-inputs}
Fix $R,b$ and put $k=R/b$, $D=4-R$. Suppose that
\[
 0<R<3/2,\qquad 1\le k<2,\qquad 0<b\le1.
\]
Consider the scalar comparison states
\[
 1<E\le E_0,\quad a=k/E,\quad 1/E\le a\le f(E),\quad 1/E\le b\le1,
\]
\[
 Ea\ge\max(1,b),\quad Eb\ge\max(1,a),\quad Ea,Eb<2.
\]
Set $e_A=\alpha=0$ and $\beta=1/2-c\ge m_B$, using the exact
nonnegative-mean inverse for the $B$ leaf and its exact central-sign
minimum $m_B$. With $u=1-1/E$, define
\[
 A=1-a,\quad K_0=3+(1-b)^2,\quad t=e_B/b,\quad
 S=2\log E+2a+2b-2-R/2,
\]
\[
 W_0=S-3b/2,\qquad A_0=a^2K_0-3b.
\]
The feasible set is one interval, possibly empty or a singleton.
On every nontrivial interior of this interval,
\begin{equation}\label{late:r24:a:profile-inputs}
 c_u>0,\qquad t_u\le0,\qquad t_{uu}>-E^2,
 \qquad A_0<W_0,\qquad W_0>0.
\end{equation}
The exact determinant
\[
 \Gamma=\frac{(3+A^2)(K_0+bt)-(S+t)^2}{D^2}
\]
is strictly concave in \(u\). Its interval endpoints lie on the faces
\begin{equation}\label{late:r24:zero-faces}
 Eb=1,\qquad a=f(E),\qquad e_B=bh(b-f(E))>0,
 \qquad\beta=0,\qquad E=E_0.
\end{equation}
Intersections and singleton intervals are retained. The only possible
$E=1$ endpoint is $k=b=1$, where the continuous comparison has
\begin{equation}\label{late:r24:a:profile-unit-corner}
 \Gamma>\frac{299}{400}>\frac{1471}{2000}.
\end{equation}
No assertion about a physical outer root at a comparison state is needed.
\end{lemma}
\begin{proof}
\paragraph{The fixed exact inverse.}
Along the prescribed path
\[
 V_B^2=\frac{b^2D}{R}=\frac{bD}{k}
\]
is constant. Consequently the exact functions $e_B(d)$ and
$A_b(d)$ are fixed functions of the single width $d$.
Write $q=2-b\ge1$. Their derivative formulas in
Lemma~\ref{late:r24:allocation} give
\[
 e_B''(d)=
 \frac{2[1-f(1+d/q)]}{q(d/q)(1+d/q)}>0\quad(d>0),
\]
\[
 e_B'(0)=1-b\ge0,\qquad A_b'>0,
 \qquad A_b''<0.
\]
The expressions extend continuously to $d=0$. Thus the inverse
composition $t(\beta)=e_B(A_b^{-1}(\beta))/b$ is increasing
and convex, with nonnegative first derivative also when $b=1,d=0$.
For completeness, the derivative of this composition is
$e_B'/(bA_b')$, and its second derivative is
\[
 \frac{e_B''A_b'-e_B'A_b''}
 {b(A_b')^3}\ge0.
\]
Since \(V_B^2\le b(4-b)\), \eqref{late:r24:closed-inverse-barrier}
gives \(d<b\) for every \(\beta\le1/2\), including the lower-crossing
faces. In particular \(q+d<2\).

\paragraph{Curvature of the exact inverse.}
The central cost is
\[
 c(u)=\frac{2ku+2b(E-1)-R\log E}{D},
\]
so
\begin{equation}\label{late:r24:a:profile-cost-derivatives}
 c_u=\frac{2k+bE(2E-k)}D>0,\qquad
 c_{uu}=\frac{bE^2(4E-k)}D>0.
\end{equation}
As $\beta=1/2-c$, it follows that
\[
 t_u=-t_\beta c_u\le0,\qquad
 t_{uu}=t_{\beta\beta}c_u^2-t_\beta c_{uu}
       \ge-t_\beta c_{uu}.
\]
The integral formula for $f$ gives $1-f(1+w)\le w/3$ for $w\ge0$.
Therefore
\[
 e_B''(d)\le\frac{2}{3q},\qquad
 e_B'(d)\le1-b+\frac{2d}{3q}=:M.
\]
In the exact cost derivative, the inequality $q+d<2$ gives
\[
 1-\frac2{q+d}+\frac4{(q+d)^2}\ge1.
\]
It follows that
\[
 A_b'(d)\ge\frac12+\frac{q^2}{2V_B^2}
 =\frac{2w_0}{bD},\qquad
 w_0=k+(1-k)b=1+(k-1)(1-b)\ge1,
\]
and hence
\begin{equation}\label{late:r24:a:profile-tbeta}
 0\le t_\beta\le\frac{DM}{2w_0}.
\end{equation}

For a positive width, the exact cost, $K_B\ge3/2$,
$V_B^2\le b(4-b)\le3$, and $q/(q+d)>1/2$ give
$\beta/d>3/4$. Coordinate-cost monotonicity gives
\[
 c\ge B_0(E)=\frac{4(E-1)-\log E}{4E-1}
             >\frac{3(E-1)}{4E-1}.
\]
Consequently
\[
 d<\frac{4\beta}{3}<\frac{7\beta}{5}
 <\delta(E):=\frac{7(5-2E)}{10(4E-1)}.
\]
If $\beta=0$, the conclusion $d=0<\delta(E)$ holds directly;
no strict comparison of two zero multiples is used. Put
\[
 C=\frac{2\delta(E)}3=\frac{7(5-2E)}{15(4E-1)}.
\]
Since $q\ge1$, we have $bM\le b(1-b+C)$. The vertex
$(1+C)/2$ lies below $l=1/E$: after multiplication by positive
factors this is
\[
 100E-46E^2-30>0.
\]
This concave polynomial is $24$ at $1$ and $1154/225$ at $E_0$.
As $b\ge l$, it follows that
\[
 bM(4E-1)
 \le\frac{(E-1)(4E-1)}{E^2}
        +\frac{7(5-2E)}{15E}<2.
\]
The last inequality is equivalent to
$-16E^2+40E-15>0$, whose concave polynomial has endpoint values
$9$ and $1409/225$ on the same interval. Combining this with
\eqref{late:r24:a:profile-cost-derivatives} and
\eqref{late:r24:a:profile-tbeta}, using $k,w_0\ge1$, gives
\[
 t_{uu}\ge-\frac{bM E^2(4E-k)}{2w_0}>-E^2.
\]
This is a bound for the exact inverse, including its nonlinear dependence
on the remaining budget.

\paragraph{The state inequality in both coordinate orders.}
Let $\rho=\log E-u>0$. Since $k\ge1$ and $Eb<2$,
\[
 W_0=\frac{(k-1)(4/E-b)}2+2\rho>0.
\]
If $a\le b$, then
\[
 A_0\le b^2K_0-3b=b(b-1)(b^2-b+3)\le0<W_0.
\]
For $b\le a\le f$, let $G(a,b)=W_0-A_0$, at fixed $E$.
The two second partial derivatives are
\[
 G_{aa}=-2K_0<0,\qquad G_{bb}=-2a^2<0.
\]
For fixed $b$, concavity reduces the minimum over $b\le a\le f$
to $a=b$ or $a=f$. On the latter edge concavity in $b$ reduces it
to $b=l$ or $b=f$. The diagonal endpoints were just proved.
At the remaining corner $(a,b)=(f,l)$, the identity defining $f$ gives
\[
 G(f,l)=\frac32(l+f)+\frac{fu^2}{1-u}-f^2(3+u^2).
\]
Here $0<u\le11/26<1/2$ and $f>2/3$. The derivative of the right
side with respect to $f$ is less than $2-4<0$ throughout the
comparison interval. The positive series in
Lemma~\ref{late:r24:a:elementary} gives
\[
 f< F:=1-u/3-u^2/6.
\]
Substitution of $F$ in the preceding decreasing expression, followed
by multiplication by $36(1-u)$, yields
\[
 u^2(15+21u-13u^2-9u^3+3u^4+u^5)>0.
\]
Indeed its parenthesis exceeds
$15-13/4-9/8=85/8$. Thus $A_0<W_0$ for either coordinate order.
The corners in this argument are used only to prove a scalar inequality,
not to construct a feasible lift.

The estimate \eqref{late:r24:b:closed-coefficients} gives
\(W_*=S+t-b(3+A^2)/2=t-DL_y/2>0\).
Writing \(Z=D^2\Gamma\), use \(A_u=k\), \(S_u=2E-2k\) and
\(S_{uu}=2E^2\) to obtain
\begin{align*}
 Z_{uu}
 &=2k^2(K_0+bt)+4bkA t_u-2(S_u+t_u)^2
      -4E^2(S+t)-2W_*t_{uu}\\
 &<2E^2\{A_0-W_0-(1-ba^2)t-bA^2/2\}
       +4bkA t_u-2(S_u+t_u)^2<0.
\end{align*}
The first strict sign uses \(t_{uu}>-E^2\). The final one follows
from \(A_0<W_0\), \(t_u\le0\) and \(0<a,b\le1\).
The exact inverse extends continuously to zero width.

\paragraph{The entire feasible interval and its endpoints.}
The bounds $R<3/2$, $Ea=k<2$, and $Ea\ge\max(1,b)$ are fixed.
Since $a\le f<1$, the other lower crossing is just $Eb\ge1$.
It is preserved as $E$ increases, and $Eb<2$ is automatic from
$E\le E_0$ and $b\le1$. The condition $a\le f$ is equivalent to
$k\le Ef(E)$; the right side is strictly increasing by
\eqref{late:r24:a:f-derivatives}. These are lower-end interval
conditions. Meanwhile $\beta$ strictly decreases. The exact central
minimum excess is
\[
 H_B(E)=\max\{0,b(E-1)(b-f(E))\}.
\]
It is nondecreasing: on its positive part,
\[
 H_B'=b\left[b-\frac{2-f}{E}\right]>0,
\]
because $b>f$ and $(E+1)f>2$. The increasing exact mean as a
function of budget decreases along the profile. Consequently
$\beta\ge0$ and $e_B(\beta)\ge H_B$ are upper-end interval
conditions. Their intersection with the lower-end conditions and
$1<E\le E_0$ is therefore one interval.

The lower endpoints are $Eb=1$ or $a=f(E)$; the upper endpoints
are $E=E_0$, $\beta=0$, or positive equality $e_B=H_B>0$.
If $e_B=H_B=0$, strict increase of the true mean away from zero
forces $d=\beta=0$, so it adds no sixth face. At $\beta=0$ the
central signs also give $b\le f$, the genuine two-zero-leaf domain.
The inverse-width barrier excludes an additional width endpoint.
If $k>1$, then $k\le Ef(E)$ excludes $E=1$. If $k=1$ but $b<1$,
then $E\ge1/b>1$. Thus only $k=b=1$ can reach $E=1$, at an
intersection of the stated lower faces.

At that corner, $a=b=1$, $D=3$, $c=0$, $\beta=1/2$,
$q=1$, and $V_B^2=3$. The exact cost satisfies
\[
 A_b(d)\ge d/2+\frac{d}{2(1+d)}.
\]
Its right side is $8/15>1/2$ at $d=2/3$, so $d<2/3$.
The integrated mean bound gives $e_B\le d^2/3$ and hence
\[
 y=e_B/3<4/81<1/20,\qquad
 \Gamma=3/4-y^2>3/4-1/400=299/400.
\]
This proves the strict corner inequality directly.
\end{proof}
\space
\subsection{True minimum widths and the unsigned circle boundaries}

All comparisons in this subsection use the exact functions
\eqref{late:r24:mean}--\eqref{late:r24:determinant}. A closed circle
condition is a scalar condition, not a claim that the comparison is a
maximizing trajectory. Write $E_0=26/15$, $h=E-1$, and $f=f(E)$.
The integral bounds in Lemma~\ref{late:r24:a:elementary}
give \(1/E<f<1\), \((E+1)f>2\), and \(Ef\le E-h/3\).

On the closed rectangle $1<E\le E_0$, $1/E\le a,b\le1$, both
central cost derivatives are positive. Indeed
\[
 D^2c_a=8h-4Eb\log E+2Eb^2h>8(h-\log E)>0,
\]
and the other formula follows by exchanging the already defined scalar
names. Moreover the coefficient signs extend to this closed rectangle.
For example
\[
 N:=DL_y=4(1-a)+ab(E-2+a)-4\log E,
 \qquad N_a=-4+b(E-2+2a)<-2.
\]
At $a=1/E$ its coefficient of $b$ is $h^2/E^2>0$. Hence
\begin{equation}\label{late:r24:b:closed-coefficients}
 DL_y\le\frac{h(5E-1)}{E^2}-4\log E<0,
 \qquad L_x<0.
\end{equation}
The last strict inequality follows because the opposite difference
vanishes at $E=1$ and has derivative $2(2E-1)(E-1)/E^3>0$.
The inverse barrier in Lemma~\ref{late:r24:allocation} extends to
the equality height, because
\begin{equation}\label{late:r24:closed-inverse-barrier}
 V_v^2\le v(4-v)\quad\Longrightarrow\quad
 A_v(v)\ge1-\frac{(2-v)^2\log(2/(2-v))}{v(4-v)}>\frac12.
\end{equation}
The strict inequality is \(1+2\log X<X^2\) at \(X=2/(2-v)>1\).
Thus every nonnegative budget below $1/2$ used below has its unique
exact width below the coordinate $v$. At zero budget the width and
excess are zero; statements involving their ratios are used only at
positive budget and extended by continuity.

\begin{lemma}[Strong true-minimum estimate]
\label{late:r24:b:strong-minimum}
Let $E>1$ and $f(E)<v\le1$. Put $s=1-v$, $q=2-v$, and let $d>0$
solve the true central-minimum equation
\[
 e_v(d)=vh(v-f(E)).
\]
Then
\begin{equation}\label{late:r24:b:strong-width}
 vh\ge qd+3s,
\end{equation}
strictly if $v<1$. If $v=1$, then $d=h$. For the circle height
$V_v^2=v(4-v)$ every positive width satisfies
\begin{equation}\label{late:r24:b:circle-leaf-cost}
 3e_v(d)<(qd+3s)A_v(d).
\end{equation}
Consequently the exact central-minimum cost is strictly greater than
$3(v-f(E))$. This cost conclusion persists whenever
$V_v^2\le v(4-v)$ and for every larger retained true mean.
\end{lemma}
\begin{proof}
Set $S(w)=1-f(1+w)$. Its exact integral and differential identities are
\[
 S(w)=2w\int_0^1\frac{t(1-t)}{1+wt}\,dt,\quad
 S'(w)=2\int_0^1\frac{t(1-t)}{(1+wt)^2}\,dt
       =\frac{w-(w+2)S(w)}{w(1+w)}.
\]
Thus $S'>0$, $S'(0)=1/3$, and
\begin{equation}\label{late:r24:b:S-bounds}
 0<S(w)<\frac w{3+w},\qquad
 S(w)\le\frac w3-\frac{w^2}{6(1+w)}\quad(w>0).
\end{equation}
The second bound follows by inserting
$(1+wt)^{-1}\le1-wt/(1+w)$ in the integral. The exact excess is
$e_v(d)=d[s+S(d/q)]$.

Suppose first $0<s<1$. For $p\ge0$ define
\[
 L=q^2p+3s,\qquad H=L/v,\qquad
 B=s+\frac{qp[s+S(p)]}{L},\qquad W(p)=B-S(H).
\]
All denominators are positive. The strict first inequality in
\eqref{late:r24:b:S-bounds} gives $W(0)=s-S(3s/v)>0$.
At a hypothetical zero of $W$, differentiating $B$ and substituting
$S(H)=B$ gives
\[
 B'=\frac{q[3s(s+S(p))+pLS'(p)]}{L^2},\qquad
 S'(H)H'=\frac{q^2[L-(L+2v)B]}{L(v+L)}.
\]
Using the differential identity for $S(p)$ therefore yields
\begin{equation}\label{late:r24:b:first-contact}
 W'=\frac{q[C+ZS(p)]}{L^2(1+p)(v+L)},
\end{equation}
where the unexpanded expressions are
\begin{align*}
 C={}&(v+L)[3s^2(1+p)+Lp]
       -q(1+p)vL(q^2p+s)+q^2p(1+p)(L+2v)s,\\
 Z={}&-3s(v+L)+q^2p(vp-L).
\end{align*}
Here \(Z=-s[3(1+2s)+6q^2p+q^2(3+s)p^2]<0\).
By \eqref{late:r24:b:S-bounds}, replace \(S(p)\) by \(p/(p+3)\)
to bound \(C+ZS(p)\) below. The coefficients of
\(((p+3)C+pZ)/s\), in order of increasing degree in \(p\), are
\[
\begin{array}{c|l}
0 & 9s^2(s+2)\\
1 & 6s(2s^3+8s^2+9s+2)\\
2 & 4s^4+15s^3+16s^2+6s+1+3q^2(s^3+7s^2+6s)\\
3 & q^2(4s^3+16s^2+14s+2)\\
4 & q^5
\end{array}
\]
Every coefficient is positive. Thus every zero in
\eqref{late:r24:b:first-contact} would have $W'>0$, impossible at a
first zero reached from $W(0)>0$. Hence $W(p)>0$ for all $p\ge0$.

Take $p=d/q$. The resulting inequality is
\[
 vH[S(H)-s]<qp[s+S(p)]=e_v(d)=vh[S(h)-s]>0,
 \qquad H=\frac{qd+3s}{v}.
\]
The function $w\mapsto vw[S(w)-s]$ strictly increases on its positive
branch, which is an upper interval since $S'>0$. If $H\ge h$, the
last display is impossible. Therefore $H<h$, proving the strict case
of \eqref{late:r24:b:strong-width}. For $s=0$ the defining equation is
$dS(d)=hS(h)$. Strict increase of $wS(w)$ gives $d=h$, without using
the contact formula at its singular value $L=0$.

For the circle cost, $V_v^2=4-q^2$ and
\[
 \frac{A_v(d)}d=\frac12+
 \frac{[3+s^2+e_v(d)]q}{2V_v^2(q+d)}
 >\frac{2q+d}{2(q+d)},
\]
because $3+s^2+e_v(d)-V_v^2=2s+2s^2+e_v(d)>0$.
The second estimate in \eqref{late:r24:b:S-bounds} gives
\[
 \frac{3e_v(d)}d\le3s+\frac{d(2q+d)}{2q(q+d)}.
\]
The difference between $(qd+3s)(2q+d)/(2(q+d))$ and this upper bound is
\[
 \frac{sd[2q^2-q+(q+1)d]}{2q(q+d)}\ge0.
\]
The strict cost comparison proves \eqref{late:r24:b:circle-leaf-cost},
including $s=0$. At the central minimum,
\[
 3vh(v-f)=3e_v(d)<(qd+3s)A_v(d)\le vhA_v(d).
\]
Division by $vh>0$ proves the cost assertion. At fixed width the
positive nonconstant part of $A_v$ is proportional to $1/V_v^2$, so
weakening the height only increases cost. Also $e_v$ and $A_v$ are
strictly increasing on the retained positive branch. On an original
crossed state, $R>v$ gives
$V_v^2=v^2(4-R)/R<v(4-v)$, so this transfer applies to each original
leaf separately. No assertion for $v>1$ is made.
\end{proof}

\begin{lemma}[Positive circle-side cap]
\label{late:r24:b:positive-a-cap}
Suppose $1<E\le E_0$, $a=1/E$, $b=f(E)$, $e_B=0$, and let the
$A$ leaf have its exact true mean at budget
$\alpha_0=1/2-c_f(E)\ge0$, where $c_f(E)=c(E,1/E,f(E))$.
Then $\Gamma>1471/2000$.

More precisely, put $D_0=4-f$, $u=h/E$, $x_0=Ee_A/D_0$,
$C_{\rm corner}=3/4-g(E)$, and
\[
 N_p=27-6h-32h^2-12h^3,\qquad
 Z_p=24h^3+56h^2+62h+27.
\]
One has
\begin{equation}\label{late:r24:b:cap-rational-cost}
 \alpha_0\le\alpha_L:={N_p\over2E(h+3)(4h+9)}.
\end{equation}
For $0<h\le1/2$,
\begin{equation}\label{late:r24:b:cap-low-bounds}
 -L_{x,0}<\frac{2h}3,\qquad
 x_0<\frac{17-20h}{300},\qquad
 g(E)<\frac{h^3}{9E^2},
\end{equation}
and the resulting total deficit is at most
$29/2000-1/202500$. For $1/2\le h\le11/15$ and $\alpha_0>0$,
\begin{equation}\label{late:r24:b:cap-high-bounds}
 -L_{x,0}<H_c:=\frac{h(5-h)}8,\qquad
 \frac{x_0}{\alpha_0}<R_c:=\frac{4h}{16h+5},\qquad
 \alpha_0<\frac9{16}\left(\frac34-h\right).
\end{equation}
Writing $r_c(h)=h^2(5-4h)/16$, these upper bounds satisfy
\begin{equation}\label{late:r24:b:cap-high-gain}
 H_cR_c\alpha_0+R_c^2\alpha_0^2<r_c(h)\alpha_0,
 \qquad g'(E)>r_c(h)c_f'(E).
\end{equation}
There is a unique $E_*\in(1,E_0)$ with $c_f(E_*)=1/2$; positive
cap budget implies $E<E_*$.
\end{lemma}
\begin{proof}
Here $q=1+u$, $V_A^2=D_0/(E^2f)$, and the exact inverse gives
$0\le d<a$, $q+d<2$. The determinant and cost are exactly
\[
 \Gamma=\frac34-g(E)+L_{x,0}x_0-x_0^2,\quad
 \alpha_0=\frac d2+\frac{K_Adq}{V_A^2(q+d)},\quad
 K_A=\frac{3+u^2+e_A}2.
\]
For positive budget, $K_A\ge(3+u^2)/2$ and
$e_A\le ud+d^2/(3q)$ imply
\begin{equation}\label{late:r24:b:cap-Z}
 \frac{\alpha_0}{d}\ge\frac Z{4D_0},\qquad
 Z=2D_0+(3+u^2)E^2fq,\qquad
 \frac{x_0}{\alpha_0}\le\frac{4h}{Z}
       +\frac{16ED_0\alpha_0}{3qZ^2}.
\end{equation}
No upper mean has been substituted into the physical cost in this
argument.

At fixed \(E\), the derivative of the corner cost \(c_f\)
in \(f\) has positive numerator \(8E-2-4hf+hf^2/2\).
The coefficient of \(f\) in \(Z\) is
\((1+10h+16h^2+8h^3)/E>0\). The Jensen lower endpoint gives
\[
 D_0\le\frac{4h+9}{h+3},\quad Z\ge\frac{Z_p}{E(h+3)},
 \quad\alpha_0\le\alpha_L.
\]
Thus \(N_p\ge0\) on the feasible cap. Formula
\eqref{late:r24:b:cap-Z} gives
\[
 \frac{x_0}{\alpha_0}\le
 \frac{4hE(h+3)}{Z_p}+\frac{8E^3N_p}{3(1+2h)Z_p^2}.
\]
The coefficient
\(-L_{x,0}=[-f^2+(5+3h+2h^2)f-4]/[E(4-f)]\) increases with
\(f\), since its derivative numerator is \((f-4)^2+12h+8h^2\).
The upper Jensen endpoint gives
\begin{equation}\label{late:r24:b:cap-Pm}
 -L_{x,0}<\frac{hP_m}{3E(h+2)(11h+18)},\qquad
 P_m=6h^3+57h^2+122h+72.
\end{equation}
To prove \eqref{late:r24:b:cap-high-bounds}, use these rational
bounds and clear the positive denominators. The first two rows
below give the coefficient and ratio estimates, with
\(B=8h^3-13h^2-6h+12\); the last gives the cap-budget bound even on
\([1/3,3/4]\).
\[
\begin{array}{l|c|r|r}
\text{polynomial} & h\text{ interval} &\text{degree}&
 \min \text{ Bernstein coefficient}\\ \hline
-36+56h+81h^2-36h^3-33h^4 &[1/2,3/4]&4&91/16\\
3h(1+2h)Z_pB-2E^3(16h+5)N_p &[1/2,11/15]&8&1971/8\\
-135+516h-29h^2-408h^3-144h^4 &[1/3,3/4]&4&152/9
\end{array}
\]
For each row substitute \(h=h_-+(h_+-h_-)t\) and apply
\eqref{cv:bernstein-coefficients}. In particular
\begin{equation}\label{late:r24:c:extended-cap-budget}
 \alpha_0<\frac9{16}(3/4-h)\qquad(1/3\le h\le3/4)
\end{equation}
on the feasible cap.

These three estimates also prove the first inequality in
\eqref{late:r24:b:cap-high-gain}. After division by
$h^2\alpha_0$, the proposed loss bound becomes
\[
 \frac{5-h}{2(16h+5)}+
 \frac{16\alpha_0}{(16h+5)^2}<\frac{5-4h}{16}.
\]
Subtracting the first term on the left from the right leaves
$4(h-5/16)(3/4-h)/(16h+5)$. The second term is less than
$9(3/4-h)/(16h+5)^2$, and
\[
 4(h-5/16)(16h+5)=64h^2-25/4\ge39/4>9.
\]
In particular the actual loss satisfies
$(-L_{x,0})x_0+x_0^2<r_c(h)\alpha_0$.

For the same-curve gain, use
\eqref{late:r24:b:corner-derivatives}. The second inequality of
\eqref{late:r24:b:cap-high-gain} is equivalent to
\[
 P=16S-(5-4h)hE(4-f)T>0.
\]
Its cubic and quadratic coefficients in $f$ are
\begin{align*}
 &-64-149h-108h^2-27h^3-4h^4,\\
 &-384-510h-288h^2+6h^3+72h^4,
\end{align*}
both negative on $h\le3/4$, so $P$ is strictly concave in $f>0$.
At the lower and upper Jensen endpoints, clearing respectively
$(h+3)^3$ and $27(h+2)^3$ and removing the common factor $h$ gives
\begin{align*}
 P_L={}&243-1134h-323h^2+2942h^3
             +3776h^4+1696h^5+256h^6,\\
 P_U={}&1944-8100h-7462h^2+25681h^3
             +43326h^4+23895h^5+4532h^6.
\end{align*}
The least degree-six Bernstein coefficients on \([1/2,3/4]\)
are respectively \(256\) and \(88449/32\).
The two positive endpoints and concavity in \(f\) prove the gain.

By Lemma~\ref{late:r24:a:corner}, $c_f(E_0)>1/2$,
$g(E_0)<29/2000$, and $g$ increases up to $E_0$. Since $c_f$ tends
to zero as $E\downarrow1$, the asserted unique $E_*$ exists.
For positive cap budget $E<E_*$. Since
$r_c'(h)=h(5-6h)/8>0$ on $h<5/6$, integration along this same
corner curve gives
\[
 g(E_*)-g(E)>r_c(h)[c_f(E_*)-c_f(E)]
             =r_c(h)\alpha_0.
\]
The high-$h$ determinant is therefore greater than
$3/4-g(E_*)>1471/2000$.

It remains to prove \eqref{late:r24:b:cap-low-bounds} and the low-$h$
conclusion. For $0<h\le1/2$, \eqref{late:r24:b:cap-Pm} gives
$-L_{x,0}<2h/3$, since the cross-multiplied difference is
\[
 2E(h+2)(11h+18)-P_m=16h^3+45h^2+30h>0.
\]
The width estimate in \eqref{late:r24:b:cap-Z} gives
\[
 d\le\frac{2N_p}{(h+3)Z_p}<\frac{2-3h}{3E}.
\]
After removing \(h>0\), the cleared polynomial is
\(57+49h-170h^2-264h^3-72h^4\), with least degree-four
Bernstein coefficient \(3/2\) on \([0,1/2]\).
Since \(D_0\ge(11h+18)/(3(h+2))\), the true mean yields
\[
 x_0<\frac{(h+2)(2-3h)(2+8h+15h^2)}
             {9E(1+2h)(11h+18)}<\frac{17-20h}{300}.
\]
For the last bound the cleared numerator is
\(354-495h-7743h^2+16146h^3+9540h^4\).
Divide \([0,1/2]\) into four equal closed intervals. Formula
\eqref{cv:bernstein-coefficients} gives a least coefficient
\(1295/128\) across these four degree-four tests.

To sharpen $g$, the probability density $6t(1-t)$ has mean $1/2$
and variance $1/20$. The second $t$ derivative of $(1+ht)^{-1}$ is
at least $2h^2/E^3$. Strong convexity therefore gives
\[
 S(h)\ge\frac{2h}{3(h+2)}+\frac{h^3}{60E^3},\qquad
 f\le\frac{h+6}{3(h+2)}-\frac{h^3}{60E^3}.
\]
The desired inequality $g<h^3/(9E^2)$ is equivalent to $R_g(f)<0$,
where
\[
 R_g(f)=(18+26h+9h^2)f^2+(18+8h)f-36-16h.
\]
It increases with $f$, its value at the upper Jensen endpoint is
$h^3(14+9h)/(9(h+2)^2)$, and its derivative for
$f\ge3/(h+3)$ is at least $(162+198h+62h^2)/(h+3)$.
The preceding improvement of the upper endpoint suffices because
\begin{align*}
 &9(h+2)^2(162+198h+62h^2)-60E^3(h+3)(14+9h)\\
 &\hspace{8mm}=3312+2940h-4662h^2-7506h^3-3522h^4-540h^5>0.
\end{align*}
On $0\le h\le1/2$ the sum of the negative terms is at most
$9363h^2\le9363/4<3312$.

Combining these bounds yields
\[
 \Gamma>\frac34-F_c(h),\qquad
 F_c(h)=\frac{h^3}{9E^2}
            +\frac{289+2720h-3600h^2}{90000}.
\]
This function increases on $[0,1/2]$. Up to $17/45$ both derivative
terms are nonnegative. Above $17/45$, use $h>3/8$, $E\le3/2$ to
bound the first derivative below by $1/72$ and the second below by
$-11/1125$; their sum is positive. Hence
\begin{equation}\label{late:r24:c:cap-loss-endpoint}
 F_c(h)\le F_c(1/2)=\frac{11741}{810000}
                    =\frac{29}{2000}-\frac1{202500}.
\end{equation}
This proves the stated uniform margin. If $\alpha_0=0$, use
$d=x_0=0$ and the zero-leaf conclusion directly, without dividing
by the budget. All seams $h=1/2$, $b=f$, and $E=E_0$ are retained
whenever their budgets are feasible.
\end{proof}

\begin{lemma}[Opposite circle allocation]
\label{late:r24:b:opposite-circle}
Let $1<E\le E_0$, $a=1/E$, $a\le b\le f(E)$, and $e_B=0$.
Give the $A$ leaf its exact true mean at budget $\alpha=1/2-c\ge0$.
Then $\Gamma>1471/2000$. At fixed $E$ the exact full-budget
determinant strictly decreases as $b$ increases along the feasible
interval. Its endpoint is $b=\min(f(E),b_c(E))$, where
$c(E,1/E,b_c(E))=1/2$.
\end{lemma}
\begin{proof}
Put $u=1-a$, $q=1+u$, $\rho=\log E-u$, and let $d$ be the true
$A$ width. Here $e_B=0$, and the budget $\alpha$ determines the
excess $e_A(d)$ at $a=1/E$, including $e_A(0)=0$. Define
\[
 v=\frac b{4-b},\quad w=Ee_A(d)+2\rho,\quad
 k_0=3+u^2,\quad K_0=k_0-2a\rho,\quad
 \kappa=2h-\log E+\frac u2,
\]
\[
 H(d)=K_A(d)E^2\frac{dq}{q+d},\qquad T=\kappa+H(d)>0.
\]
The exact full budget is equivalently
\begin{equation}\label{late:r24:b:opposite-budget}
 c=\frac u2+\kappa v,\qquad
 \alpha=\frac d2+vH(d),\qquad
 v=\frac{a-d}{2T}.
\end{equation}
Here $H'>0$, directly from the positive derivative of the exact
nonconstant leaf-cost term. Differentiation gives
\begin{equation}\label{late:r24:b:opposite-slope}
 v_d=-\frac{T+(a-d)H'}{2T^2}<0,\qquad
 -w_v=\frac{2Ee_A'(d)T^2}{T+(a-d)H'}
       \le\frac{Ee_A'(d)(a-d)}v.
\end{equation}
The exact determinant in these variables is
\begin{equation}\label{late:r24:b:opposite-determinant}
 \Gamma=\frac{(K_0+aw)(1+3v^2)}4
              -\frac{[w+(w+6)v]^2}{16}.
\end{equation}
Its $w$ derivative is strictly negative, including the closed circle.
One direct check is
\[
 -(4-b)L_x=4\log E-4+3b+\frac{2b-b^2}{E}.
\]
This increases in $b\le1$. At $b=1/E$ its value vanishes at $E=1$
and its $E$ derivative is $(E^2-1)(4E-3)/E^4>0$.
Thus $L_x<0$, and adding any nonnegative $A$ mean preserves
$\Gamma_w<0$.

Define $Q=8(v\Gamma_v-u\Gamma_w)$. Differentiating
\eqref{late:r24:b:opposite-determinant} yields
\begin{align}
 Q_w&=u(1+2v-11v^2)-6v-2wv(1+v)<0,
 \label{late:r24:b:opposite-Q}\\
 Q_0:=Q(v,2\rho)&=2u(\rho-a)+L_1v+J_1v^2,\nonumber\\
 L_1&=6u+4u\rho-12\rho-4\rho^2,\nonumber\\
 J_1&=18u^2+2u\rho-24\rho-4\rho^2.\nonumber
\end{align}
For the first sign, $v\ge(1-u)/(3+u)$ and $u<1/2$ imply
\[
 u-(6-2u)v\le\frac{11u-u^2-6}{3+u}<0.
\]
Both $L_1$ and $J_1$ are positive. Indeed, with $\zeta=Ef$,
$\rho=u^2\zeta/2$, and Lemma~\ref{late:r24:a:elementary} implies
$\zeta<3/2$ for $h\le3/4$. Therefore
\[
 \frac{J_1}{u^2}=18-12\zeta+u\zeta-u^2\zeta^2>0,
 \qquad L_1\ge6u-9u^2-\frac94u^4>0.
\]
Consequently $Q_0$ increases in $v$ and is largest at $b=f$.
At that cap its exact value is
\[
 Q_0\bigl(f/(4-f)\bigr)
   =-\frac{2h\mathcal R(h,f)}{E^2(4-f)^2},
\]
where
\begin{equation}\label{late:r24:b:opposite-R}
 \mathcal R(h,b)=(3h+3h^2+2h^3)b^3
 +(4+6h+12h^2)b^2-(20+12h+8h^2)b+16.
\end{equation}
We prove this quantity positive in both required ranges.
For $1/2\le h\le3/4$, put $s=1-f<h/(h+3)$,
$A_3=3h+3h^2+2h^3$, and $A_2=4+6h+12h^2$. The expansion is
\[
 \mathcal R(h,1-s)=R_1-R_{b,1}s+(A_2+3A_3)s^2-A_3s^3,
\]
\[
 R_1=h(2h^2+7h-3)>0,\qquad
 R_{b,1}=-12+9h+25h^2+6h^3.
\]
The last two terms have nonnegative sum. If $R_{b,1}\le0$ the
claim follows. Otherwise
\[
 \mathcal R(h,f)>R_1-R_{b,1}\frac h{h+3}
 =\frac{h(3+9h-12h^2-4h^3)}{h+3}>0,
\]
since its bracket is at least $3-(9/4)h\ge21/16$.
For $0<h\le1/2$, $\mathcal R_{bb}>0$ and
$\mathcal R_b(h,1)\le-1/2$, so $\mathcal R$ decreases in
$0<b\le1$. Substitution of the upper Jensen endpoint gives
\[
 27(h+2)^3\mathcal R\left(h,\frac{h+6}{3(h+2)}\right)
 =h(216+852h+390h^2-33h^3+3h^4+2h^5)>0.
\]
Here $390h^2-33h^3>0$. Thus $Q_0<0$ for the entire range.

Suppose first $E\ge3/2$. The true mean has
\[
 e_A'(0)=u,\qquad e_A''(0)=\frac2{3q},\qquad
 e_A''(d)=\frac4q\int_0^1
       \frac{t(1-t)}{(1+(d/q)t)^3}\,dt,
\]
so $e_A''$ strictly decreases. For
$\phi(d)=e_A'(d)(a-d)$ one has
$\phi''=e_A'''(a-d)-2e_A''<0$ and
$\phi'(0)=2a/(3q)-u\le0$ since $u\ge1/3$.
Hence $\phi(d)\le au$, and \eqref{late:r24:b:opposite-slope}
gives $w_v\ge-u/v$. Using $\Gamma_w<0$ and
\eqref{late:r24:b:opposite-Q},
\[
 \frac{d\Gamma}{dv}\le\Gamma_v-\frac uv\Gamma_w
       =\frac{Q(v,w)}{8v}\le\frac{Q_0}{8v}<0.
\]

For $1<E\le3/2$, retain the actual mean. Put
\[
 t=Ee_A(d),\quad \eta=Ee_A'(d)(a-d),\quad m=\eta-u,
 \quad A_1=Q_w(v,2\rho),\quad Z_0=-8\Gamma_w(v,2\rho)>0.
\]
The exact expansion, including all $m$ terms, is
\begin{equation}\label{late:r24:b:actual-compensation}
 P_\eta:=8(v\Gamma_v-\eta\Gamma_w)
 =Q_0+A_1t-v(1+v)t^2+m[Z_0+(1+v)^2t].
\end{equation}
The bound $e_A''\le2/(3q)$ implies
\[
 m\le A_md-B_md^2,\qquad
 A_m=\frac2{3q}-Eu=\frac{2-5u-3u^2}{3aq}\ge0,
 \qquad B_m=\frac2{3aq}.
\]
If $m>0$, then
\begin{equation}\label{late:r24:b:eta-bounds}
 m\le A_md,\qquad
 m\le\frac{A_m^2}{4B_m}\le\frac a{6q}\le\frac16,
 \qquad t\ge\frac{ud}a.
\end{equation}
At $u=1/3$, $A_m=0$, so no positive-$m$ case occurs at the seam.
If $m\le0$, the positive bracket in
\eqref{late:r24:b:actual-compensation} gives
$P_\eta\le Q(v,w)\le Q_0<0$.

For $m>0$, the coefficient estimates are
\[
 A_1\le-K(u),\qquad K(u)=\frac{6-11u+u^2}{3+u},\qquad
 Z_0\le\frac{16u}{9a}.
\]
The last bound follows because $-L_x$ increases with $b$ on this
circle, $4-b\ge3$, and Lemma~\ref{late:r24:b:positive-a-cap}
gives $-L_{x,0}<2h/3$; in fact
$Z_0=8(-L_x)/(4-b)$. Since $v\le1/3$, the $mt$ term in
\eqref{late:r24:b:actual-compensation} is at most $8t/27$.
Thus
\[
 P_\eta\le Q_0-[K(u)-8/27]t+\frac{16u}{9a}A_md.
\]
The strict inequality $K-8/27-(16/9)A_m>0$ is proved by clearing
the positive denominator $27(3+u)(1-u^2)$: its numerator is
\[
 42-97u+113u^2+353u^3-27u^4>0\qquad(0\le u\le1/3).
\]
Here $42-97u\ge29/3$ and the remaining terms have nonnegative
sum. By \eqref{late:r24:b:eta-bounds},
\[
 P_\eta\le Q_0-\left(K-\frac8{27}-\frac{16}9A_m\right)
                  \frac{ud}a<0.
\]
Equation \eqref{late:r24:b:opposite-slope} now gives
$d\Gamma/dv\le P_\eta/(8v)<0$. At $d=0$ use directly
$P_\eta=Q_0<0$.

Finally the continuation to the claimed endpoint is feasible.
Central cost increases with $b$, so the available $A$ budget decreases
continuously; $a=1/E<f$ imposes no positive $A$ minimum, while
$e_B=0$ is allowed exactly up to $b=f$. The inverse barrier keeps
$d<a$, all lower crossings remain closed and all upper states remain
strict. The interval is therefore
$[a,\min(f,b_c)]$ whenever nonempty. At $b_c\le f$,
$\alpha=d=e_A=0$ and Theorem~\ref{late:r24:a:zero-leaf} applies.
At $f<b_c$ the endpoint is exactly
Lemma~\ref{late:r24:b:positive-a-cap}. Equality $b_c=f$, a singleton
interval, the seam $E=3/2$, and zero budget are included. The two
orientations have not been exchanged in this proof.
\end{proof}

\begin{lemma}[Exact cap transfer]
\label{late:r24:b:cap-transfer}
On the exact full-budget cap
\[
 1<E\le E_0,\qquad a=f(E),\quad e_A=0,\qquad
 1/E\le b\le1,\qquad \beta=1/2-c\ge m_b,
\]
the determinant increases strictly with $b$ at positive width and
extends continuously to zero width. Decreasing $b$ to $1/E$ preserves
feasibility and yields
\[
 \Gamma(E,f,b;\text{exact full }B\text{ budget})
 \ge\Gamma(E,f,1/E;\text{exact full }B\text{ budget})
 >\frac{1471}{2000}.
\]
\end{lemma}
\begin{proof}
Fix $E$, put $k=Ea=Ef$, $r=1-b$, $q=1+r$, $d=d_B$,
$p=q+d$, $\xi=r+d$, $e=e_B(d)$, and $K=(3+r^2+e)/2$.
By Lemma~\ref{late:r24:a:elementary},
\[
 1<k\le\frac{2E+1}3\le\frac{67}{45}<\frac32,
 \quad D=4-kb>\frac52,\quad
 V_B^2=\frac{bD}k\le b(4-b)\le3.
\]
The inverse barrier gives $d<b\le1$. At positive budget,
$q/(q+d)>1/2$ and $K>3/2$, so
\[
 \frac\beta d>\frac34,\qquad d<\frac43\beta<\frac75\beta.
\]
Since central coordinate monotonicity gives
$c\ge[4h-\log E]/(4E-1)>3h/(4E-1)$,
\[
 d<\delta(E):=\frac{7(5-2E)}{10(4E-1)},\qquad
 0\le r\le u:=1-1/E\le11/26,
\]
\begin{equation}\label{late:r24:b:cap-width-guard}
 \xi<u+\delta(E)
 =\frac7{10}-\frac{(E-1)(E+5)}{5E(4E-1)}<\frac7{10}.
\end{equation}
These bounds need only a nonnegative full budget, not yet the central
sign along a continuation.

Temporarily fix a positive true excess $e$ and differentiate its true
width with respect to $b$, retaining $k$. Set
\[
 m=\frac{(\partial_b d)_e}{d},\qquad
 v=\frac ed=q-f(1+d/q),\qquad w=d/q.
\]
Differentiating this exact mean law gives
\begin{equation}\label{late:r24:b:cap-width-derivative}
 m=\frac{b+(2+w)v}{q(2r+d-v)}
  >\frac{1+rd/q^2}{r+d}\ge\frac1\xi.
\end{equation}
The denominator is positive on the increasing mean branch. The strict
inequality follows because $v>r$ and the displayed quotient increases
with $v$.

With $\lambda=k/(bD)$, the exact cost is
$A=d/2+K\lambda dq/p$. At fixed $e$,
$K_b=-r$ and $\lambda_b=\lambda(-1/b+b\lambda)$.
Keeping the width derivative gives exactly
\begin{equation}\label{late:r24:b:cap-cost-derivative}
 \frac{A_b}{d}=\frac m2+C\lambda+\frac{Kbq}{p}\lambda^2,
 \qquad
 C=K\left(\frac{mq^2}{p^2}-\frac d{p^2}
                       -\frac q{bp}\right)-\frac{rq}p.
\end{equation}
The true mean is below the outer squared endpoint, so
\[
 \frac32<K<\frac{3+\xi^2}2<\frac{349}{200},\qquad
 \frac{Kbq}p>\frac32\frac{1-(11/26)^2}{17/10}
                =\frac{8325}{11492}.
\]
To control $C$, use $m\ge1/\xi$ and put
\[
 T=\frac q{bp}-\frac{q^2/\xi-d}{p^2}
  =\frac{2\xi^2+(1+r^2)\xi-(1-r)(1+r)^2}
         {(1-r)\xi(1+\xi)^2}.
\]
The assertion $T<11/20$ is equivalent to
\begin{align*}
 P(r,\xi)={}&11(1-r)\xi^3-(18+22r)\xi^2
 -(9+11r+20r^2)\xi+20(1-r)(1+r)^2>0.
\end{align*}
For fixed $r$, $P_\xi$ is a convex quadratic, negative at $0$ and
at $7/10$; its latter value is
$-1803/100-(5797/100)r-20r^2$. Thus $P$ decreases in $\xi$.
At the right endpoint,
\[
 1000P(r,7/10)=8653-2253r-34000r^2-20000r^3,
\]
which decreases in $r\ge0$. At $r=11/26$, its numerator over
the denominator $17576$ is
\[
 152085128-16753308-106964000-26620000=1747820>0.
\]
Consequently
\[
 C>-\frac{349}{200}\frac{11}{20}-\frac{11}{26}
   =-\frac75+\frac{893}{52000}>-\frac75.
\]
Equations \eqref{late:r24:b:cap-width-derivative} and
\eqref{late:r24:b:cap-cost-derivative} imply
\[
 \frac{A_b}d>\frac57-\frac75\lambda
                    +\frac{8325}{11492}\lambda^2>0,
\]
because the quadratic minimum is
\[
 \frac57-\frac{563108}{832500}
       =\frac{220744}{5827500}>0.
\]
This is positivity of the exact cost derivative at fixed positive
excess, not a derivative inferred from a minimum-cost lower bound.

On $a=f$, the central cost derivative simplifies to
\[
 c_b=\frac{2h[3+(1-a)^2]}{D^2}>0.
\]
For $F(b,e)=c(b)+A(b,d(b,e))$, therefore,
$F_e=A_d/e_d>0$ and $F_b>0$. The equation $F(b,e)=1/2$ gives
$e_b=-F_b/F_e<0$. Since $(bD)_b=4-2kb>1$, it also gives
\begin{equation}\label{late:r24:b:cap-y-decreasing}
 y=\frac e{bD},\qquad y_b<0.
\end{equation}
Put $\Delta=(E+1)a-2>0$ and $D_1=4-k$. Direct differentiation
of the exact central determinant and coefficient gives
\begin{equation}\label{late:r24:b:cap-determinant-derivatives}
 (C_0)_b=\frac{u\Delta(D_1D-2ku\Delta)}{D^3},\qquad
 (L_y)_b=-\frac{2ku\Delta}{D^2}.
\end{equation}
For an explicit algebraic check, set $A_0=4-2a+a^2$,
$M=D_1/2$, and $S_0=2\log E+2a-2$. The numerator of $C_0$ is
$A_0(b^2-2b+4)-(Mb+S_0)^2$; the cap identity gives
$S_0=2(k-1)-u\Delta$ and
$A_0=D_1+k(k-1-u\Delta)$. Its quotient derivative numerator is
$-kD_1u\Delta b+4D_1u\Delta-2ku^2\Delta^2$, proving the first
formula in \eqref{late:r24:b:cap-determinant-derivatives}.

Now
\[
 \Gamma_b=(C_0)_b+(L_y)_by+(L_y-2y)y_b.
\]
The last term is positive by \eqref{late:r24:b:closed-coefficients}
and \eqref{late:r24:b:cap-y-decreasing}. Since the true mean satisfies
$e<1-r^2=bq$, the first two terms equal
\[
 \frac{2ku\Delta}{D^2}
 \left(\frac{D_1}{2k}-\frac{u\Delta}{D}-y\right)>0.
\]
Indeed the bracket exceeds
\[
 \frac{(2-k)(8-kb)-2ku\Delta}{2kD}>0,
\]
because $(2-k)(8-kb)>13/4$ while $2ku\Delta<9/8$, using
$k<3/2$, $b\le1$, $u<1/2$, and $\Delta<h<3/4$.

To continue to $b=1/E$, first construct the exact inverse throughout
the central-cost sublevel interval, using its barrier. Central cost
decreases as $b$ decreases, so the budget stays nonnegative and below
$1/2$. The true full-budget excess increases by the derivative just
proved. The required excess $\max(0,bh(b-f))$ decreases when $b$
decreases, since its positive-branch derivative is $h(2b-f)>0$.
Therefore a central-sign-feasible starting point remains feasible.
Throughout,
\[
 Ea=k>1\ge b,\qquad Eb\ge1>a,\qquad Ea,Eb<2,
 \qquad R=kb<3/2.
\]
Zero-width starting points are included continuously and become
positive-width points as soon as the budget increases, unless the
interval is a singleton. The endpoint is thus reached without any
additional width or sign boundary.

At this endpoint exchange only the scalar names already defined:
\[
 (a_{\rm new},b_{\rm new},e_{A,\rm new},e_{B,\rm new})
       =(1/E,f,e_B,0).
\]
The central cost and common determinant formulas are unchanged. The
old $B$ height is $(4-f)/(E^2f)$, exactly the new $A$ height, and its
budget is $1/2-c_f$. Its true mean law has coordinate $1/E$; the
new zero $B$ leaf has coordinate $f$. Thus
Lemma~\ref{late:r24:b:positive-a-cap} applies with all costs and
means matched, including zero budget. This scalar equality does not
permute the modes of an actual trajectory.
\end{proof}

\begin{lemma}[Saturated circle slope and endpoint]
\label{late:r24:b:saturated-slope}
Let $E>1$, $h=E-1$, $u=h/E$, $\rho=\log E-u$, and on the circle
$a=1/E$ set $e_A=0$, $e_B=bh(b-f)$, with $f\le b\le1$.
Define the smaller determinant
\begin{equation}\label{late:r24:b:saturated-surrogate}
 \widetilde\Gamma_{\min}(E,b)
 =\frac34+\frac{u^2(b^2-2b+4)-6b\rho
                     -(2\rho+h(b-f))^2}{(4-b)^2}.
\end{equation}
The exact determinant equals this expression plus
$bu^2h(b-f)/(4-b)^2\ge0$. At fixed $E$,
\begin{equation}\label{late:r24:b:saturated-derivative}
 0<-\partial_b\widetilde\Gamma_{\min}<\frac{h^2}3,
\end{equation}
and consequently
\begin{equation}\label{late:r24:b:saturated-displacement}
 \widetilde\Gamma_{\min}(E,b)
 \ge C_{\rm corner}(E)-\frac{h^2}3(b-f).
\end{equation}
Equality holds at $b=f$; the last inequality is strict for $b>f$.
These derivative and displacement assertions do not require a
full-budget hypothesis.

If in addition $1<E\le E_0$, $f<b\le1$, and the true minimum
$B$ cost satisfies $c(E,1/E,b)+m_b=1/2$, then
$\widetilde\Gamma_{\min}>1471/2000$.
\end{lemma}
\begin{proof}
Since $hf=2E\rho/h$, one has
$2\rho+h(b-f)=hb-2\rho/h$. Thus
\[
 \widetilde\Gamma_{\min}
 =\frac34+\frac{-Ab^2-Bb+C}{(4-b)^2},
\]
\[
 A=\frac{h^3(E+1)}{E^2},\qquad
 B=2(u^2+\rho),\qquad C=4u^2-\frac{4\rho^2}{h^2}.
\]
Differentiation gives the exact positive expression
\begin{equation}\label{late:r24:b:saturated-exact-slope}
 -\partial_b\widetilde\Gamma_{\min}
 =\frac{8Ab+2b(u^2+\rho)+8\rho+8\rho^2/h^2}{(4-b)^3}.
\end{equation}
It increases in $0<b\le1$, since its numerator is positive and
affine increasing while its denominator decreases.

Put $v=\rho/h^2$. The identities and rational bound needed at $b=1$
are
\[
 \rho'=\frac h{E^2},\qquad
 v'=\frac{1/E^2-2v}{h},\qquad
 v<\frac{2E+1}{6E^2}.
\]
For the last inequality, differentiate $h^2/(2E)-\rho$ to obtain
$h^2/(2E^2)$. Integration from $0$ to $h$ and bounding the positive
integrand below using its final denominator give
$h^2/(2E)-\rho>h^3/(6E^2)$, which is the stated bound.
We claim
\begin{equation}\label{late:r24:b:rho-comparison}
 8v^2+10v<1+\frac6{E^2}.
\end{equation}
For $Q=1+6/E^2-10v-8v^2$, the convergent logarithm expansion near
$h=0$ gives
\[
 v=\frac12-\frac23h+\frac34h^2+O(h^3),\qquad
 Q=\frac{17}{18}h^2+O(h^3)>0.
\]
At any hypothetical first zero, substitution of
$8v^2+10v=1+6/E^2$ into its derivative gives
\[
 hE^3Q'=4E^3+2E+12-4E(4+5E^2)v.
\]
By the strict rational bound for $v$, multiplying this expression's
lower bound by $3E$ gives
\[
 12E^4-20E^3-4E^2+20E-8
   =4(E-1)^2(E+1)(3E-2)>0.
\]
Thus every zero has positive derivative, contradicting first contact
from positive values. This proves \eqref{late:r24:b:rho-comparison}
for every $E>1$. At $b=1$, the quantity in
\eqref{late:r24:b:saturated-exact-slope}, divided by $h^2$, is
\[
 \frac{8-6/E^2+10v+8v^2}{27}<\frac13.
\]
Its monotonicity in $b$ proves
\eqref{late:r24:b:saturated-derivative}. At $b=f$, the mean excess
$e_B=bh(b-f)$ vanishes and the surrogate equals $C_{\rm corner}$. Integration proves
\eqref{late:r24:b:saturated-displacement}, with its precise equality
case.

For the full-budget positive-saturation endpoint we need a second
same-corner gain:
\begin{equation}\label{late:r24:b:saturated-corner-gain}
 g'(E)>\frac{h^2}9c_f'(E)\qquad(1<E\le7/4).
\end{equation}
Using \eqref{late:r24:b:corner-derivatives}, this is equivalent to
$P=9S(E,f)-hE(4-f)T(E,f)>0$. Its full expansion in $f$ is
\[
 P=-A_3f^3-A_2f^2+A_1f-A_0,
\]
\begin{align*}
 A_3&=36+82h+63h^2+17h^3,\\
 A_2&=216+276h+132h^2+18h^3,\\
 A_1&=540+492h+276h^2+72h^3,\\
 A_0&=288+80h+144h^2+64h^3.
\end{align*}
It is strictly concave in positive $f$. At the lower Jensen endpoint,
clearing $(h+3)^3$ gives
\[
 h(1458+1539h-369h^2-1142h^3-504h^4-64h^5)>0.
\]
For $0<h\le3/4$, the negative part of its bracket is at most
\[
 h^2\left(369+1142\frac34+504\frac9{16}
                       +64\frac{27}{64}\right)
 =1536h^2\le864<1458.
\]
At the upper Jensen endpoint, clearing $27(h+2)^3$ gives
\[
 h(11664+18144h+1152h^2-11476h^3-6813h^4-1151h^5)>0.
\]
The negative part of this bracket is at most
$(11476+6813+1151)(3/4)^3=32805/4<11664$.
Concavity in $f$ proves \eqref{late:r24:b:saturated-corner-gain}.

Now let $\beta=m_b>0$ and $c+\beta=1/2$ as in the statement.
Lemma~\ref{late:r24:b:strong-minimum} gives $\beta>3(b-f)$,
including $b=1$. Since $c$ increases with $b$,
\[
 b-f<\frac\beta3=\frac{1/2-c(E,1/E,b)}3
                   \le\frac{1/2-c_f(E)}3.
\]
It follows that
\[
 \widetilde\Gamma_{\min}
 >\frac34-\left[g(E)+\frac{h^2}9(1/2-c_f(E))\right].
\]
In particular $c_f(E)<1/2$, so $E<E_*$ from
Lemma~\ref{late:r24:b:positive-a-cap}. On $(1,E_*]$ put
$K(E)=g(E)+(h^2/9)(1/2-c_f(E))$. Equation
\eqref{late:r24:b:saturated-corner-gain} gives
\[
 K'=g'-\frac{h^2}9c_f'+\frac{2h}9(1/2-c_f)>0.
\]
By the corner estimates in Lemma~\ref{late:r24:a:corner},
\[
 K(E)<K(E_*)=g(E_*)<g(E_0)<\frac{29}{2000}.
\]
Hence $\widetilde\Gamma_{\min}>1471/2000$. The seam $b=f$ has
zero minimum mean; when it is a full-budget saturation endpoint,
$\beta=0$, and the zero-leaf theorem applies directly rather than
using the strict positive-minimum inequality.
\end{proof}

\begin{lemma}[Closed circle with zero $A$ budget]
\label{late:r24:b:first-circle}
Suppose
\[
 1<E\le E_0,\qquad a=1/E,\qquad \alpha=e_A=0,\qquad
 1/E\le b\le1,\qquad \beta=1/2-c\ge m_b.
\]
Give $B$ its exact true mean at this full budget, including the
zero-width comparison when $\beta=0$. Then
$\Gamma>1471/2000$. In fact the explicitly smaller quantity
\begin{equation}\label{late:r24:b:first-envelope}
 \widetilde\Gamma
 =\frac34+\frac{u^2(b^2-2b+4)-6b\rho-(2\rho+t)^2}{(4-b)^2},
 \qquad u=1-1/E,\quad\rho=\log E-u,\quad t=e_B/b,
\end{equation}
satisfies the same bound throughout the circle domain.
\end{lemma}
\begin{proof}
Write $q=2-b$, $D=4-b$, and $k_0=b^2-2b+4$. Direct substitution
in the exact determinant gives
\begin{equation}\label{late:r24:b:first-exact}
 \Gamma=\widetilde\Gamma+\frac{bu^2t}{D^2}
       \ge\widetilde\Gamma.
\end{equation}
In particular the exact coefficient of $y=t/D$ is
$(bu^2-4\rho)/D$; replacing it by $-4\rho/D$ constructs the smaller
envelope in \eqref{late:r24:b:first-envelope}. All following profile
derivatives and endpoint bounds concern that same envelope.

Fix $b$. The exact $B$ mean and cost depend only on $b,d$, since
$V_B^2=bD$ is independent of $E$. The central cost and its two
$u$ derivatives are
\begin{equation}\label{late:r24:b:first-central}
 c=\frac{2u+b(2u/(1-u)+\log(1-u))}{D},\quad
 c_u=\frac{2+b(2E^2-E)}D,\quad
 c_{uu}=\frac{bE^2(4E-1)}D.
\end{equation}
At fixed $b$ the lower crossing condition is $E\ge1/b$. The budget
$\beta$ decreases strictly as $E$ increases. For $0\le\beta\le1/2$,
the inverse-width barrier \eqref{late:r24:closed-inverse-barrier}
gives a continuous true mean that strictly increases with $\beta$. The
required excess
\[
 H_B(E)=\max(0,b(E-1)(b-f(E)))
\]
is nondecreasing, since its derivative on the positive part is
\[
 H_B'=b\left[b-\frac{2-f}{E}\right]>0
\]
by $(E+1)f>2$. Thus the budget and mean conditions cut out a single
interval, possibly empty or a singleton. If nonempty, its closure
starts at $E=1/b$ and ends at $E=E_0$, zero budget, or positive
central saturation $e_B=H_B>0$. When $b=1$, we include $E=1$ by
continuity of the scalar formula and its inverse budget. If
$e_B=H_B=0$, the width and budget are zero. The inverse barrier
excludes a separate width endpoint.

Every interior point of a nontrivial feasible interval has
$\beta>0$, because the budget decreases strictly and is nonnegative
throughout the interval. On this interior, exact inverse convexity
in Lemma~\ref{late:r24:allocation} gives $t_\beta>0$ and
$t_{\beta\beta}\ge0$. Thus
\[
 t_{uu}=t_{\beta\beta}c_u^2-t_\beta c_{uu}
       \ge-t_\beta c_{uu}.
\]
The true mean derivative satisfies
$e_B'(d)\le1-b+2d/(3q)$. From the exact cost derivative and
$q+d<2$,
\[
 A_B'(d)\ge\frac12+\frac{q^2}{2bD}
           =\frac2{bD},\qquad
 t_\beta\le\frac D2\left(1-b+\frac{2d}{3q}\right).
\]
At this positive width, the exact circle cost gives $\beta/d>3/4$
as in \eqref{late:r24:b:cap-width-guard}, and hence
\[
 d<\frac75\beta<\frac{7(5-2E)}{10(4E-1)}.
\]
Using $b(1-b)\le1/4$ and $b/q\le1$ gives
\begin{align*}
 b\left(1-b+\frac{2d}{3q}\right)E(4E-1)
 &\le\frac{E(4E-1)}4+\frac{7E(5-2E)}{15}\\
 &=\frac{E^2}{15}+\frac{25E}{12}<4.
\end{align*}
The last polynomial increases, and its difference from $4$ at $E_0$
is exactly $1273/6750$. It follows that $t_{uu}>-2E$, and in
particular
\begin{equation}\label{late:r24:b:first-curvature-guard}
 t_{uu}+2E^2>0.
\end{equation}
Twice differentiating the whole envelope, retaining the nonlinear
inverse chain terms, gives
\begin{equation}\label{late:r24:b:first-curvature}
 D^2\widetilde\Gamma_{uu}
 =2k_0-6bE^2-2(t_u+2h)^2
                -2(2\rho+t)(t_{uu}+2E^2)<0.
\end{equation}
Indeed $E\ge1/b$, $b\le1$, and
$bk_0-3=(b-1)(b^2-b+3)\le0$, so
$k_0\le3/b\le3bE^2$; the final term is strictly negative by
\eqref{late:r24:b:first-curvature-guard} and $\rho>0$.
This proves strict concavity of the smaller fixed-budget profile,
without asserting concavity of the full $\Gamma$.

We next bound each endpoint of this same profile. At its lower
crossing endpoint $E=1/b$, one has $b=1-u$, $D=3+u$,
$k_0=3+u^2$, and $0<u\le11/26<1/2$. The logarithm series gives
\[
 \rho=\sum_{n\ge2}\frac{u^n}n
 \le\frac{u^2}2+\frac{u^3}{3(1-u)}
 =\frac{u^2(3-u)}{6(1-u)}.
\]
Consequently the numerator of $C_0-3/4$ is bounded below by
\[
 u^2(3+u^2)-6(1-u)\rho-4\rho^2
 \ge\frac{u^3(9-18u-3u^2+8u^3)}{9(1-u)^2}>0.
\]
The cubic bracket decreases on $[0,1/2]$ and is $1/4$ at its right
endpoint. Thus $C_0>3/4$ for $u>0$.

The true mean has $e_B\le ud+d^2/3$. The exact cost, $K_B\ge3/2$,
and $q+d<2$ give, with $N=9-u-2u^2$,
\[
 \beta\ge\frac{dN}{4(1-u)(3+u)},\qquad
 \beta<\frac{3-5u}{2(3+u)},\qquad
 d<\frac{2(1-u)(3-5u)}N.
\]
Substitution into the increasing upper mean proves
$y=e_B/[(1-u)(3+u)]<1/20$: the positive denominator clears to
\[
 9-39u+991u^2+2051u^3-1152u^4+12u^5>0.
\]
The cubic and quartic terms together are nonnegative on $[0,1/2]$,
and the quadratic has discriminant $1521-35676=-34155<0$.
Also $\rho/(3+u)$ strictly increases; its derivative numerator is
$h(3+u)-\rho>0$, since $\rho\le uh$. At $u_0=11/26$,
Lemma~\ref{late:r24:a:corner} gives $\log(26/15)<551/1000$.
The rational difference
\[
 3(3+u_0)-80(551/1000-u_0)=\frac{23}{650}>0
\]
therefore proves $4\rho/(3+u)<3/20$. Hence at this diagonal endpoint
\begin{equation}\label{late:r24:b:first-diagonal}
 \widetilde\Gamma=C_0-\frac{4\rho}{D}y-y^2
 >\frac34-\frac3{20}\frac1{20}-\frac1{400}
 =\frac{37}{50}.
\end{equation}
At the limiting point $E=b=1$ one has $C_0=3/4$, $\rho=0$,
and $\beta=1/2$. Here $A_b(d)\ge d/2+d/[2(1+d)]$, which is
$8/15$ at $d=2/3$. Hence $d<2/3$ and $y\le d^2/9<1/20$.
This gives $\widetilde\Gamma>299/400>37/50$ directly.

At any zero-budget endpoint, $t=0$ and the central sign gives
$b\le f$. Thus $\widetilde\Gamma=C_0>1471/2000$ by
Theorem~\ref{late:r24:a:zero-leaf}.

At an upper endpoint with $E=E_0$ and $\beta>0$, the logarithm bound gives
\[
 c\ge c(E_0,1/E_0,1/E_0)
 =\frac{44-15\log(E_0)}{89}>\frac25.
\]
Thus $\beta<1/10$ and $d<4\beta/3<2/15$.
The corner-cost inequality $c_f(E_0)>1/2$ implies $b<f(E_0)<1$.
Also $b\ge15/26>4/7$. Dividing the true mean upper bound by the
exact cost lower bound gives
\[
 \frac y\beta\le
 \frac{1-b+d/3}{bD/2+(3/2)q/(q+d)}<\frac14.
\]
For explicit rational estimates the numerator is below
$3/7+2/45=149/315$. Since $b>1/2$, $bD>7/4$, and
$(3/2)q/(q+d)>45/34>5/4$, the denominator exceeds $17/8$.
The desired comparison follows from $4\cdot1192<5355$.
The same logarithm bound gives $\rho(E_0)<16/125$, so
$4\rho/D<1/5$. Thus
\begin{equation}\label{late:r24:b:first-e0-loss}
 C_0-\widetilde\Gamma
 <\frac\beta{20}+\frac{\beta^2}{16}
 <\frac9{160}\beta<\frac\beta9.
\end{equation}
For clarity, the central gain that pays this loss can be seen on the
same circle. At fixed $E$, its exact derivatives are
\[
 -C_{0,b}=\frac{6\rho(b+4)+8\rho^2-6bu^2}{(4-b)^3},\qquad
 c_b=\frac{k_c}{(4-b)^2},\qquad
 k_c=8h-4\log E+2u>0.
\]
Since $\rho>u^2/2$,
\[
 -\frac{dC_0}{dc}>\frac{3u^2}{k_c}
 >\frac{3u(1-u)}{6+2u^2}>\frac19
 \quad\text{at }u=11/26.
\]
For the middle bound, $k_c=8h-2u-4\rho<8h-2u-2u^2$,
whose product with $1-u$ is $u(6+2u^2)$.
The last comparison is $27u-29u^2-6>0$, whose value at
$u=11/26$ is $157/676$.
Because the circle cost is increasing and $c_f(E_0)>1/2$, one can
increase $b$ to the unique circle point $b_c\le f$ with $c=1/2$.
The central drop exceeds $\beta/9$, so
\eqref{late:r24:b:first-e0-loss} and
Theorem~\ref{late:r24:a:zero-leaf} give
\[
 \widetilde\Gamma>C_0(E_0,1/E_0,b_c)>\frac{1471}{2000}.
\]
At a positive saturation endpoint, $b>f$ and $c+m_b=1/2$.
Lemma~\ref{late:r24:b:saturated-slope} then gives
$\widetilde\Gamma=\widetilde\Gamma_{\min}>1471/2000$.
Concavity on the interior and continuity at the endpoints give,
on every nontrivial feasible interval,
\[
 \Gamma(u)\ge\widetilde\Gamma(u)
 \ge\min\{\widetilde\Gamma(u_{\rm left}),
                       \widetilde\Gamma(u_{\rm right})\}
 >\frac{1471}{2000}.
\]
For a singleton the endpoint estimate gives the same bound.
The endpoint estimates apply also when several endpoint conditions hold.
\end{proof}

\begin{lemma}[Positive-$B$ saturation transfer]
\label{late:r24:b:saturation-transfer}
Suppose
\[
 1<E\le E_0,\qquad 1/E\le a\le f(E)\le b\le1,
 \qquad e_A=0,\quad e_B=bh(b-f),
\]
\[
 c+m_b\le1/2,\qquad R=Eab<3/2,
\]
with the true central-minimum $B$ width. At fixed $E,b$ and these
fixed means, the exact determinant strictly increases with $a$.
Lowering $a$ to $1/E$ decreases both central and $B$ costs.
Filling any remaining circle budget with the exact $B$ inverse yields
\begin{align*}
 \Gamma(E,a,b;0,bh(b-f))
 &\ge\Gamma(E,1/E,b;0,bh(b-f))\\
 &\ge\Gamma(E,1/E,b;\text{exact full }B\text{ budget})
 >\frac{1471}{2000}.
\end{align*}
The seam $b=f$, zero width, and equality in the budget are included.
\end{lemma}
\begin{proof}
Write $\mu_A=(1-a)^2+e_A$ and $\mu_B=(1-b)^2+e_B$.
Set $v=\log E/h$, $k=Eb$, and
\[
 B=3+\mu_B=4+Eb^2-2Ebv,\quad
 M=\frac{4-k}2,\quad N=b(E+1)-2v,\quad D=4-ka.
\]
The identity $2Ev=2+hf$ and the central-minimum mean give
\[
 DJ_m=Ma+N,\qquad 3+\mu_A=4-2a+a^2,
\]
so the exact common-matrix determinant is
\begin{equation}\label{late:r24:b:saturation-determinant}
 \Gamma=\frac{B(4-2a+a^2)-(Ma+N)^2}{(4-ka)^2}.
\end{equation}
Let $T=4M+kN$. Differentiation gives an affine numerator
\[
 L(a):=D^3\Gamma_a
 =2M(2B-T)a+8(k-1)B-2NT.
\]
The quadratic term in $a$ cancels: for numerator coefficients
$A_2,A_1,A_0$, the quotient derivative numerator is
$(8A_2+kA_1)a+4A_1+2kA_0$. Moreover
\[
 2B-T=kh\left(\frac{2-f}{E}-b\right)<0,
\]
because $Eb\ge Ef>2-f$. Since $M>0$, $L$ decreases with $a$,
and it suffices to bound $L(f)$.

Put $W=Eb-(2-f)>0$. The identities
$T=2B+bhW$ and $2(k-1)-N=hW/E$ give the exact factorization
\[
 L(f)=\frac{hW}{E}P(b),\qquad
 P(b)=16-[4+2f(3E-1)]b+E(Ef-2h)b^2.
\]
This polynomial decreases on $0<b\le1$. If its quadratic coefficient
is nonpositive this is immediate; otherwise its derivative is at most
\[
 P_b(1)=-4+2f(E^2-3E+1)-4Eh<0,
\]
since $1<E\le E_0<2$. Hence
\[
 P(b)\ge P(1)=12-2Eh+f(E^2-6E+2)
 >14-E^2-4E=9-6h-h^2>0.
\]
Here $f<1$, $E^2-6E+2<0$, and $h<3/4$ give the strict signs.
Thus $L(a)\ge L(f)>0$, proving $\Gamma_a>0$, including the
one-sided endpoint interpretations. At $b=f$, $W=(E+1)f-2>0$
although the minimum width is zero.

At fixed $E,b$, the true mean law, its minimum width $d$, and $K_B$
are independent of $a$. The exact minimum cost is
\[
 m_b=\frac d2+\frac{K_BEad(2-b)}{bD(2-b+d)},\qquad
 \partial_am_b=\frac{4EK_Bd(2-b)}{b(2-b+d)D^2}\ge0.
\]
The central cost also increases with $a$. Thus lowering $a$ preserves
$c+m_b\le1/2$, the fixed width $d<b$, and $R<3/2$.
All scalar state inequalities persist:
\[
 Ea\ge1\ge b,\qquad Eb\ge Ef>1\ge a,\qquad Ea,Eb<2.
\]
At $a=1/E$, the available budget $\beta_{\rm full}=1/2-c$ is at
least the new minimum cost and is strictly below $1/2$. The inverse
barrier remains strict even though $V_B^2=b(4-b)$. Increasing the
budget therefore gives a unique width below $b$, increases the true
excess, and preserves its central lower bound. By
\eqref{late:r24:b:closed-coefficients}, $L_y<0$ at this closed
circle, so $\Gamma=C_0+L_yy-y^2$ decreases as the nonnegative
excess is increased. These are the two displayed comparisons; equality
in the second is allowed if no slack remains. Their endpoint has exactly
the orientation and domain of Lemma~\ref{late:r24:b:first-circle}.
\end{proof}
\space
\subsection{Exact positive-minimum curvature and its boundary estimates}
\label{late:r24:c:section}

All means and costs in this subsection are those of
\eqref{late:r24:mean}--\eqref{late:r24:leaf-cost}, and every determinant
is the exact expression \eqref{late:r24:determinant}, unless a smaller
expression is explicitly displayed. We use the strong true-minimum
inequality of Lemma~\ref{late:r24:b:strong-minimum}, the positive-$A$
cap estimates of Lemma~\ref{late:r24:b:positive-a-cap}, and the
fixed-$E$ saturated-circle estimate of
Lemma~\ref{late:r24:b:saturated-slope}. These estimates concern scalar
comparisons and do not assign a maximizing lift to a comparison point.

\begin{lemma}[Curvature of the exact minimum cost]
\label{late:r24:c:G-curvature}
Fix $R,b$ and $k=R/b=Ea$. Retain
\[
 1<E<26/15,\qquad 0<a,b\le1,\qquad a>f(E),\qquad R<3/2,
\]
the strict crossed and upper state inequalities
\eqref{late:r24:crossings}, and the true central-minimum width
$0<d=d_{a,\min}<a$. Suppose $G=c+m_a\le1/2$. With
\[
 \mu=(1-a)^2+ah(a-f(E))
     =1+Ea^2-\frac{2Ea\log E}{h},
\]
regarded as the local parameter at fixed $R,b$, one has
\[
 G''(\mu)<0.
\]
More precisely, the outer endpoint $p=2-a+d$ satisfies
$G''<p''/2<0$. The assertion includes the $a=1$ face, with its
one-sided interpretation when that face is an endpoint of the unsigned
profile. It makes no connectedness assertion about the full feasible
set.
\end{lemma}
\begin{proof}
Put
\[
 \lambda=\frac RD,\quad r=1-a,\quad q=1+r=2-a,\quad
 s=r+d,\quad p=1+s=q+d,\quad w=k-a=ah,
\]
\[
 e=\mu-r^2>0,\qquad j=s^2-\mu,\qquad K=\frac{3+\mu}{2}.
\]
Both $k$ and $\lambda$ are fixed. The same number $\mu$ is the
weighted-square mean, with weight $(1+t)^{-2}$, on the two intervals
$[-r,k-1]$ and $[r,s]$. Indeed, direct integration on an interval whose
translated endpoints are $q,p$ gives denominator $(p-q)/(qp)$ and
numerator
\[
 (p-q)-2\log(p/q)+\frac{p-q}{qp}.
\]
The two mean identities are therefore
\begin{equation}\label{late:r24:c:paired-means}
 \mu-1=ak\left(1-\frac{2\log(k/a)}{k-a}\right),\qquad
 \mu-1=qp\left(1-\frac{2\log(p/q)}{p-q}\right).
\end{equation}
The first is precisely the central-sign equality and the second is the
true leaf mean. Since $0\le r<s$, one has $r^2<\mu<s^2$, so $j>0$.

Eliminating the logarithm in the central cost by the first identity in
\eqref{late:r24:c:paired-means} gives
\[
 c=\frac w2+\frac{Kwa}{V^2k},\qquad
 m_a=\frac d2+\frac{Kdq}{V^2p},\qquad
 V^2=\frac{a^2D}{R}.
\]
For example, multiplication of the first proposed formula by $D/h$
gives $aD/2+Kb=2(a+b)-Eab\log(E)/h$, which verifies it directly.
Consequently
\begin{equation}\label{late:r24:c:G-exact}
 G=\frac{k+p-2}{2}+K\lambda T,\qquad
 T=\frac{2-q^2/p}{a^2}-\frac1k.
\end{equation}
Here the two nonconstant cost factors add exactly as
\[
 \frac{w}{ak}+\frac{dq}{a^2p}
 =\frac1a-\frac1k+\frac q{a^2}-\frac{q^2}{a^2p}=T.
\]

We first establish a uniform outer-coordinate bound using the exact
cost. The integral formula
\[
 f(1+h)=2\int_0^1\frac{1-t}{1+ht}\,dt
\]
and strict Jensen inequality give
$f(E)>3/(E+2)>2/(E+1)$. Hence $k-1>r$. If $r>0$, the integral of
$(t^2-\mu)/(1+t)^2$ over $[-r,r]$ is strictly negative, whereas its
integral over $[-r,k-1]$ vanishes. Its integral over $[r,k-1]$ is
therefore positive. Starting at $r$, its primitive first decreases up
to $\sqrt\mu$ and then strictly increases. The unique positive-width
zero of this primitive is $s$, by the second exact mean equation.
Thus $k-1>s$. If $r=0$, uniqueness gives $k-1=s$ instead.

The crossing $Eb>1$ gives
\[
 V^2=\frac{4a}{Eb}-a^2<a(4-a)\le3.
\]
Moreover $a/k=1/E>1/2$, $q/p>1/2$ since $d<a$ and $q+a=2$, and
$K>3/2$. Thus the exact costs give
\[
 c>\frac34w,\qquad m_a>\frac34d.
\]
Since $w+d=k+s-1\ge2s$, the weak budget $G\le1/2$ implies the
strict bounds
\begin{equation}\label{late:r24:c:outer-guard}
 0\le r<s<\frac13,\qquad a>\frac23,\qquad p<\frac43,
 \qquad 3d-a=3s-2r-1<0.
\end{equation}

Let $M(r,s)$ denote the true weighted-square mean on $[r,s]$, and put
$H=\mu-rs=e-rd$. The exact mean identity gives
\[
 H=dS(d/q),\qquad
 S(x)=1-f(1+x)=2x\int_0^1\frac{t(1-t)}{1+xt}\,dt.
\]
The probability density $6t(1-t)$ has mean $1/2$. Strict Jensen
inequality and the elementary upper bound on the same integral yield
\begin{equation}\label{late:r24:c:H-bounds}
 S(x)>\frac{2x}{3(2+x)},\quad S(x)<\frac x3,
 \qquad
 H>\frac{2d^2}{3(p+q)}>\frac{d^2}{4}.
\end{equation}
The last inequality uses $p+q<8/3$ from
\eqref{late:r24:c:outer-guard}. Direct endpoint differentiation gives
\begin{equation}\label{late:r24:c:outer-Hessian}
 \begin{gathered}
 M_r=\frac{ep}{dq}>0,\qquad M_s=\frac{jq}{dp}>0,\\
 M_{rr}=\frac{2pH}{d^2q},\qquad
 M_{ss}=\frac{2qH}{d^2p},\qquad
 M_{rs}=1-\frac{2H}{d^2},\\
 \det\operatorname{Hess}M=\frac{4H}{d^2}-1>0.
 \end{gathered}
\end{equation}
The diagonal entries are positive; hence the true outer mean is
strictly jointly convex on the whole guarded domain.

All subsequent primes denote differentiation in $\mu$ at fixed $k$.
Differentiating the central mean in $a$ gives
\[
 \mu_a=\frac{ke}{aw}>0,\qquad
 \mu_{aa}=2\mu_a\left(\frac1w+\frac re\right)>0.
\]
Thus, with $L=1/w+r/e$,
\begin{equation}\label{late:r24:c:a-derivatives}
 a'=\frac{aw}{ke}>0,\qquad a''=-2(a')^2L<0.
\end{equation}
Twice differentiating $M(1-a,p-1)=\mu$ gives
\[
 0=\operatorname{Hess}M[(-a',p'),(-a',p')]
       +M_r(-a'')+M_sp''.
\]
The first two terms are positive, and $M_s>0$, so $p''<0$. The first
derivative identity gives a further useful strict bound:
\[
 p'=\frac{1+M_ra'}{M_s}>\frac1{M_s}>\frac32.
\]
Indeed $j<d(r+s)$ and hence
$M_s<q(r+s)/p<r+s<2/3$.

It remains to control the total derivatives of $T$ in
\eqref{late:r24:c:G-exact}; the sign $p''<0$ alone is insufficient.
Set $t=p'/a'>0$. The first derivative identity and
\eqref{late:r24:c:a-derivatives} give
\[
 t=\frac ej\left(\frac{kdp}{awq}+\frac{p^2}{q^2}\right).
\]
We claim
\begin{equation}\label{late:r24:c:cancellation}
 qt>dpL.
\end{equation}
Multiplication of $qt/(dp)-L$ by $j>0$ shows that this is equivalent to
\begin{equation}\label{late:r24:c:Xi}
 \Xi:=\frac{e-j}{w}+\frac ea+\frac{ep}{dq}-\frac{jr}{e}>0.
\end{equation}
In obtaining this exact expression we used $k=a+w$.
By \eqref{late:r24:c:H-bounds},
\[
 e<rd+\frac{d^2}{3q}\le rd+\frac{d^2}{3},\qquad
 e-j=2e-d(r+s)<0.
\]
At fixed $r,d,a,q,p$, use $j=d(r+s)-e$ to obtain
\[
 \Xi_e=\frac2w+\frac1a+\frac p{dq}
                  +\frac{rd(r+s)}{e^2}>0,\qquad
 \Xi_w=-\frac{e-j}{w^2}>0.
\]
Lemma~\ref{late:r24:b:strong-minimum}, in its strong-width form, gives
$w\ge W_0:=qd+3r$. Also
\eqref{late:r24:c:H-bounds} gives
$e>e_0:=rd+2d^2/[3(p+q)]$. Lowering $e$ to $e_0$ preserves
$e-j<0$. It is therefore enough to prove
$\Xi(e_0,W_0)>0$.

Here the required cancellation has an explicit polynomial form. Put
\begin{align*}
 U&=p+q=2q+d,\\
 C&=3rU+2d=6rq+(3r+2)d,\\
 B&=3U-4=2+6r+3d,\\
 \tau&=aq=1-r^2.
\end{align*}
Then $e_0=dC/(3U)$ and $j_0=d(C+dB)/(3U)$, and multiplication by
the positive denominator $3UaqCW_0$ gives exactly
\begin{equation}\label{late:r24:c:Xi-polynomial}
 (3UaqCW_0)\Xi(e_0,W_0)=d\mathcal P(r,d),
\end{equation}
where
\[
 \mathcal P=W_0\{2C^2-3\tau(2r+d)C+2\tau dB\}-dB\tau C.
\]
Expanding \(\mathcal P\) in \(d\), its four coefficients are
\[
\begin{array}{c|l}
0&108r^3q^3\\
1&78r^2+360r^3+534r^4+288r^5+36r^6\\
2&22r+123r^2+260r^3+189r^4+36r^5\\
3&2+14r+39r^2+36r^3+9r^4
\end{array}
\]
They are nonnegative, with a positive cubic coefficient. At \(r=0\)
the polynomial is \(2d^3>0\) and the clearing factors remain positive
since \(C=2d\) and \(W_0=d\). Thus
\eqref{late:r24:c:cancellation} holds on the whole stipulated domain.

Treating $a,p$ as independent variables, the exact partial derivatives
of $T$ are
\[
 \begin{gathered}
 T_a=-\frac{4d}{a^3p}<0,\qquad
 T_p=\frac{q^2}{a^2p^2}>0,\qquad
 T_{aa}=\frac{4(3d-a)}{a^4p}<0,\\
 T_{ap}=-\frac{4q}{a^3p^2},\qquad
 T_{pp}=-\frac{2q^2}{a^2p^3}=-\frac{2T_p}{p}.
 \end{gathered}
\]
The possibly positive term $T_a\,a''$ is controlled by the mixed term:
\[
 T_a\,a''+2T_{ap}a'p'
 =\frac{8(a')^2}{a^3p^2}(dpL-qt)<0.
\]
Therefore
\[
 T''<T_{aa}(a')^2+T_{pp}(p')^2+T_p\,p''<0.
\]
Moreover $T'=T_a\,a'+T_p\,p'<T_p\,p'$, so
\[
 T'+KT''<T_p\left[p'-\frac{2K}{p}(p')^2+Kp''\right]<0.
\]
Indeed $p'>3/2$, $p<4/3$, and $2K>3$ imply
$p'>p/(2K)$, while $p''<0$. Since $K'=1/2$ and $\lambda>0$,
\eqref{late:r24:c:G-exact} now gives
\[
 G''=\frac{p''}{2}+\lambda(T'+KT'')<\frac{p''}{2}<0.
\]
The implicit derivatives exist because $\mu_a,M_s>0$. Every strict
estimate survives the weak budget $G=1/2$. At $a=1$ the strong-width
lemma gives $d=h$, and the preceding formulas apply directly. The
proof retains $e>0$ and $E>1$; no expression containing $1/e$ or $1/h$
is evaluated on a zero-minimum or $E=1$ boundary.
\end{proof}

\begin{lemma}[Two true positive minima]
\label{late:r24:c:both-minima}
Suppose
\[
 1<E\le26/15,\qquad f(E)<a,b\le1,\qquad
 e_A=ah(a-f(E)),\qquad e_B=bh(b-f(E)),
\]
and both widths are their true central-sign minima, with
$c+m_a+m_b\le1/2$. Then
\[
 \Gamma>\frac{21467}{29160}>\frac{1471}{2000}>\overline H.
\]
The same strict bound holds when $a=f(E)$ or $b=f(E)$, with the
corresponding true minimum width, excess, and cost set to zero.
\end{lemma}
\begin{proof}
Write $v=\log(E)/h$, $Q=a+b$, and
\[
 \mu_A=1+Ea^2-2Eav,\quad \mu_B=1+Eb^2-2Ebv,
 \quad \kappa_A=3+\mu_A,\quad \kappa_B=3+\mu_B,
\]
\[
 N=(E+1)(a+b-2v)+2-\frac{Eab}{2}.
\]
The identity $2+hf=2Ev$ gives
$(1-a)^2+ah(a-f)=\mu_A$ and its $B$ counterpart. It also gives
\[
 \frac{e_A}{a}+\frac{e_B}{b}+2\log E+2Q-2-\frac R2=N.
\]
Thus the original determinant is exactly
\begin{equation}\label{late:r24:c:both-exact}
 \Gamma=\frac{\kappa_A\kappa_B-N^2}{D^2}.
\end{equation}
The scalar state inequalities needed for the minimum-cost comparison
hold automatically here: $f>3/(E+2)>1/E$, so $Ea>1\ge b$ and
$Eb>1\ge a$, while $Ea,Eb<2$ and $D>0$. In particular
$V_A^2<a(4-a)$ and $V_B^2<b(4-b)$.

We first reduce the range of $E$ without using a determinant bound.
The central cost is increasing in each coordinate, because
\[
 c_a=\frac{2h[4+Eb^2-2Ebv]}{D^2}>0,
 \qquad 4+Eb^2-2Ebv>4-E+E(b-1)^2>0.
\]
The same holds for $c_b$. Thus $c\ge c_{ff}$, the central cost at
$a=b=f$. The strict integral bound
$\log E<h-h^2/(2E)$ gives
\[
 c_{ff}>\frac{h[4f-(E+1)f^2/2]}{4-Ef^2}.
\]
The fraction in brackets increases with $f$: its derivative numerator
is $4[4-(E+1)f+Ef^2]>0$, since
$4-f-Ef(1-f)>3-E/4>0$. Substituting $f>3/(E+2)$ yields
\begin{equation}\label{late:r24:c:minimum-central-cost}
 c>\mathcal L(E):=
 \frac{3(E-1)(5E+13)}{2(4E^2+7E+16)}.
\end{equation}
The derivative numerator of $\mathcal L$, after removal of its
positive constant factor, is $3E^2+264E+219>0$, and
$\mathcal L(8/5)=105/208>1/2$. Hence every feasible state satisfies
$E<8/5$, or $0<h<3/5$.

On this smaller domain, $\log E<h$ gives
\[
 c_{ff}>h\frac{4f-Ef^2}{4-Ef^2}>\frac{3h}{4}.
\]
The second inequality is equivalent to $16f-12-Ef^2>0$.
This expression increases in $f\le1$; at $f=3/(E+2)$ its numerator
is $3(16-3E-4E^2)>0$, since the bracket decreases to $24/25$ at
$E=8/5$. We have proved
\begin{equation}\label{late:r24:c:three-quarter-central}
 c>3h/4.
\end{equation}

We next bound the zero-excess diagonal directly. Put
\[
 S=1-f,\qquad m=(E+1)f-2,\qquad D_f=4-Ef^2,
 \qquad \kappa=3+S^2.
\]
The integral estimates for $f$ and $1-f$ give
\[
 \frac3{h+3}<f<\frac{h+6}{3(h+2)},\qquad
 0<m<\frac h3,\qquad 0<S<\frac{h}{h+3},\qquad D_f\ge3-h.
\]
At $a=b=f$, direct substitution gives
\[
 \kappa=D_f+fm,\qquad J_mD_f=\frac{D_f}{2}+\frac{E+1}{E}m.
\]
Expanding the difference of squares therefore gives
\begin{equation}\label{late:r24:c:diagonal-identity}
 C_0(E,f,f)-\frac34=\frac{mhJ_f}{E^2D_f^2},
\end{equation}
where
\[
 \begin{aligned}
 J_f&=-E^2f^3+E^2f^2+(-E^2+4E+1)f-2(E+1)\\
    &=-h^2+(2h^2-3)S-2E^2S^2+E^2S^3.
 \end{aligned}
\]
For example the expansion follows by using
\[
 D_f\left(2f-\frac{E+1}{E}\right)
 +m\left(f^2-\frac{(E+1)^2}{E^2}\right)=\frac h{E^2}J_f.
\]
Dropping only favorable terms in $-J_f$ gives
\[
 \frac34-C_0(E,f,f)
 <\frac{h^2}{3E^2(3-h)^2}
   \left[h^2+\frac{3h}{h+3}
                    +\frac{2E^2h^2}{(h+3)^2}\right].
\]
This remains a valid upper bound if $J_f$ is positive. The right
side increases in $h$: $h/[E(3-h)]$ has derivative numerator
$3+h^2>0$, and each bracketed term increases. At $h=3/5$ it equals
\[
 \frac{25}{3072}\frac{451}{450}
 =\frac{451}{55296}<\frac1{120},
 \qquad 54120<55296.
\]
Hence
\begin{equation}\label{late:r24:c:diagonal-bound}
 C_0(E,f,f)>\frac34-\frac1{120}.
\end{equation}

Differentiate \eqref{late:r24:c:both-exact} at fixed $E$, retaining
both central-sign equalities. The two cancellation identities
\[
 (a-v)D+b\kappa_A=4(a+b)-v(4+R),
\]
\[
 (E+1-Eb/2)D+EbN=(E+1)\kappa_B
\]
give
\begin{equation}\label{late:r24:c:both-derivatives}
 \Gamma_a=\frac{2\kappa_B\mathcal T}{D^3},\qquad
 \Gamma_b=\frac{2\kappa_A\mathcal T}{D^3},
\end{equation}
with the same scalar
\[
 \mathcal T=-h^2(a+b)+[2(E^2+1)-ER]v-\frac{(E+1)D}{2}.
\]
On the entire algebraic square $f\le a,b\le1$, the identity
$2Ev=E+1-hS$ yields
\[
 \mathcal T_a=-h^2+\frac{EbhS}{2}<0,
 \qquad \mathcal T_b=-h^2+\frac{EahS}{2}<0.
\]
Thus $\mathcal T\ge\mathcal T(1,1)$, and
\[
 -\mathcal T(1,1)=\frac{h^3}{E}
       +hS\left(\frac E2+\frac1E\right)
 <h^2\frac{3E^2+2E-2}{2E(E+2)}
 \le\frac{37}{48}h^2.
\]
The last rational function increases: its derivative numerator over
the squared original denominator is $8(E^2+E+1)>0$, and its value
at $E=8/5$ is $37/48$. On the same square,
\[
 \kappa_A,\kappa_B\le3+S^2+hS\le\frac{563}{180},
 \qquad D\ge4-E>\frac{12}{5}.
\]
If $\mathcal T\ge0$, the desired derivative bound is immediate.
If $\mathcal T<0$, the preceding positive bounds applied to
\eqref{late:r24:c:both-derivatives} give
\begin{equation}\label{late:r24:c:both-slope-bound}
 \Gamma_a,\Gamma_b\ge
 -h^2\frac{2(563/180)(37/48)}{(12/5)^3}
 =-h^2\frac{2603875}{7464960}>-\frac{h^2}{2}.
\end{equation}
Only the explicit determinant is being compared on this square; the
intermediate points need not obey the original budget or represent
maximizing trajectories.

Integrating \eqref{late:r24:c:both-slope-bound} along the coordinate
segments from $(f,f)$ to $(a,b)$ and using
\eqref{late:r24:c:diagonal-bound} gives
\[
 \Gamma>\frac34-\frac1{120}
               -\frac{h^2}{2}\bigl[(a-f)+(b-f)\bigr].
\]
Every identity and bound above also holds for $f\le a,b\le1$.
On a zero-minimum face the unsigned mean law has its unique
nonnegative zero at width zero, and its exact cost is zero. The
minimum-cost inequalities of
Lemma~\ref{late:r24:b:strong-minimum} apply at the original state,
including either coordinate equal to one. Using equality $0=0$
for a coordinate equal to $f$, they give
\[
 (a-f)+(b-f)\le\frac{m_a+m_b}{3}
       \le\frac{1/2-c}{3}<\frac{1/2-3h/4}{3}.
\]
Therefore
\[
 \Gamma>\frac34-\frac1{120}
                   -\frac{h^2(1/2-3h/4)}6.
\]
The subtracted cubic $h^2/12-h^3/8$ has derivative
$h(4-9h)/24$ and maximum $4/729$ at $h=4/9$. Consequently
\[
 \Gamma>\frac34-\frac1{120}-\frac4{729}
 =\frac{21467}{29160}>\frac{1471}{2000},
\]
where the last comparison has numerator
$21467\cdot2000-1471\cdot29160=39640>0$.
The weak budget $c+m_a+m_b=1/2$ causes no loss of strictness.
The zero-minimum edges have been included by their exact zero-width
formulas and the weak first cost comparison, followed by the strict
central-cost bound. No strict inequality is inferred from continuity
alone. When both coordinates equal $f$,
\eqref{late:r24:c:diagonal-bound} is stronger.
The $a=1$ and $b=1$ faces were included in the derivative square and
in the strong-minimum lemma. The face $E=8/5$ and all larger retained
$E$ are infeasible by \eqref{late:r24:c:minimum-central-cost}.
\end{proof}

\begin{lemma}[The circle with two retained leaf budgets]
\label{late:r24:c:two-positive-circle}
Retain the exact scalar mean and cost laws on
\[
 1<E\le26/15,\qquad a=1/E,\qquad f(E)\le b\le1,
 \qquad e_B=bh(b-f(E)),\qquad \beta=m_b,
\]
with the true $B$ minimum and the exact remaining $A$ budget
\[
 \alpha=\frac12-c(E,1/E,b)-m_b\ge0.
\]
Then $\Gamma>1471/2000>\overline H$. The statement includes $b=1$,
$b=f(E)$, $\alpha=0$, and the zero-budget cap. The circle $Ea=1$ is
a closed scalar comparison face, not a strict-crossed physical source.
\end{lemma}
\begin{proof}
Write
\[
 f=f(E),\quad \delta=b-f,\quad D_0=4-f,\quad D=4-b=D_0-\delta,
 \quad u=\frac hE,\quad \rho=\log E-u,
\]
\[
 c_0=c(E,1/E,f),\qquad \alpha_0=\frac12-c_0,\qquad
 C_{\rm corner}=C_0(E,1/E,f)=\frac34-g(E).
\]
The integral inequality $f>1/E$ gives $Eb>1>a$, $Ea=1\ge b$,
$Eb<2$, $R=Eab=b\le1$, and $D,D_0\ge3$. The $A$ height is
\[
 V_A^2=\frac{D}{E^2b}<a(4-a),
\]
whereas the $B$ height is the equality value $V_B^2=b(4-b)$.
By \eqref{late:r24:closed-inverse-barrier}, both inverse widths
stay below their respective coordinate bounds.
The $A$ minimum is zero
because $a<f$. At $\delta=0$ the $B$ minimum and width are zero.
No exterior-root test is used for these comparisons.

The central cost strictly increases with $b$. Indeed, with
$v=\log(E)/h<1$ and $a=1/E$,
\[
 c_b=\frac{2h[4+Ea^2-2Eav]}{D^2}
     =\frac{2h[4+1/E-2v]}{D^2}>0.
\]
Thus $\alpha_0\ge\alpha\ge0$, so its cap inverse is defined.
Let $e_0$ denote the true $A$ excess at cap budget $\alpha_0$, and put
\[
 x_0=\frac{Ee_0}{D_0},\qquad
 x=\frac{Ee_A}{D},\qquad y=\frac{e_B}{bD}=\frac{h\delta}{D}.
\]
If $\delta>0$, then
Lemma~\ref{late:r24:b:strong-minimum} gives $m_b>3\delta$, including
$b=1$, and hence
\begin{equation}\label{late:r24:c:circle-budget}
 \alpha\le\alpha_0-3\delta,\qquad
 0\le\delta\le\min\left(1-f,\frac{\alpha_0}{3}\right).
\end{equation}
The first comparison is strict for $\delta>0$. In particular a
positive displacement is feasible only when $\alpha_0>0$. If
$\alpha_0=0$, necessarily $\delta=0$ and $\alpha=0$, and the
zero-budget cap is covered by Lemma~\ref{late:r24:b:positive-a-cap}.
We henceforth take $\alpha_0>0$.

At fixed $E$, the $A$ mean law is unchanged as $b$ varies, while
$V_A^2=(4/b-1)/E^2$ strictly decreases. Thus at any fixed true width
and mean its exact cost increases with $b$. Inverting the increasing
cost function shows that the $A$ mean excess at a fixed budget is no
larger at $b\ge f$ than at the cap. By Lemma~\ref{late:r24:allocation},
the cap mean excess is increasing and convex as a function of exact
cost. At zero cost the excess vanishes and the full mean equals
$(1-a)^2$. The chord from $(0,0)$ to $(\alpha_0,e_0)$ therefore
bounds the excess above. Using \eqref{late:r24:c:circle-budget} yields
\begin{equation}\label{late:r24:c:exact-chord}
 e_A\le\frac{\alpha}{\alpha_0}e_0,
 \qquad
 x\le\frac{D_0}{D}
        \left(1-\frac{3\delta}{\alpha_0}\right)x_0.
\end{equation}
This uses the true inverse and the common budget of the two leaves.
In particular, no upper mean is inserted into a lower cost estimate.
All factors in \eqref{late:r24:c:exact-chord} are nonnegative. Two
consequences used below are
\begin{equation}\label{late:r24:c:chord-factors}
 \frac{D_0}{D}\le1+\frac\delta3,
 \qquad
 \frac{1-3\delta/\alpha_0}{1-\delta/D_0}
 \le1-\left(\frac3{\alpha_0}-\frac1{D_0}\right)\delta
 \le1-\frac{17}{3}\delta.
\end{equation}
For the second, write the middle quotient as
$1-(3/\alpha_0-1/D_0)\delta/(1-\delta/D_0)$ and use
$\alpha_0\le1/2$, $D_0\ge3$. Its final upper bound stays positive
because $\delta\le\alpha_0/3\le1/6$.

Set $e_A=0$ temporarily, keeping the exact $B$ minimum.
Denote the resulting exact determinant by $\Gamma_{B,\min}$. The
explicit smaller expression from
Lemma~\ref{late:r24:b:saturated-slope} is
\[
 \widetilde\Gamma_{B,\min}
 =\frac34+
 \frac{u^2(b^2-2b+4)-6b\rho-(2\rho+h\delta)^2}{D^2}.
\]
It is smaller in the stated direction, since
\[
 \Gamma_{B,\min}-\widetilde\Gamma_{B,\min}
     =\frac{bu^2h\delta}{D^2}\ge0.
\]
For fixed $E$, that lemma proves
$0<-\partial_b\widetilde\Gamma_{B,\min}<h^2/3$ on
$f\le b\le1$, independently of any full-budget condition.
Consequently
\begin{equation}\label{late:r24:c:saturated-baseline}
 \Gamma_{B,\min}\ge\widetilde\Gamma_{B,\min}
 \ge C_{\rm corner}-\frac{h^2\delta}{3},
\end{equation}
with strictness in the last comparison for $\delta>0$ and equality
at $\delta=0$.

Restoring the actual value of $e_A$ gives
\begin{equation}\label{late:r24:c:effective-loss}
 \Gamma=\Gamma_{B,\min}-C_{\rm eff}x-x^2,\qquad
 C_{\rm eff}=(-L_x)+(2-ab)y.
\end{equation}
The coefficient can be expanded exactly. Put
\[
 M_0=-D_0L_{x,\rm cap},\qquad
 L=2E+1+\frac2E-\left(1+\frac1E\right)f.
\]
The original coefficient on the circle is
\[
 -DL_x=4\log E-4+3b+\frac{2b-b^2}{E}.
\]
It increases with $b\le1$. At $b=1/E$ it is zero at $E=1$ and has
$E$ derivative $(E^2-1)(4E-3)/E^4>0$. Thus $L_x<0$ throughout
this closed circle, without an appeal to a strict-crossed coefficient
statement. Also $ab<1$, so $C_{\rm eff}>0$. Expanding at $b=f+\delta$
gives
\begin{equation}\label{late:r24:c:effective-coefficient}
 C_{\rm eff}=\frac{M_0+L\delta-\delta^2}{D}.
\end{equation}
Indeed the coefficient of $\delta$ in
$-DL_x+(2-b/E)h\delta$ is
$3+(2-2f)/E+2h-hf/E=L$, and its quadratic coefficient is
$-(1+h)/E=-1$.

If $H$ is an upper bound on $-L_{x,\rm cap}$, then
both $H>0$ and $L\ge2E+1/E>0$. Equation
\eqref{late:r24:c:effective-coefficient} and $D\ge3$ give
\[
 0<C_{\rm eff}\le H+\frac{H+L}{3}\delta.
\]
Hence every nonnegative bound $X_\delta\ge x$ yields
\begin{equation}\label{late:r24:c:circle-envelope}
 \Gamma\ge\frac34-
 \left[g(E)+\frac{h^2\delta}{3}
          +(H+C_1\delta)X_\delta+X_\delta^2\right],
 \qquad C_1\ge\frac{H+L}{3}.
\end{equation}
The strict baseline comparison is available when $\delta>0$; at
$\delta=0$ the strict cap estimates below give strictness.

We record two elementary bounds on $L$. For $0<h\le1/2$, substitution
of $f>3/(E+2)$ gives $L<4$: after multiplication by $E(E+2)$, the
positive difference is
\[
 -2E^3-E^2+7E-1.
\]
It decreases on $[1,3/2]$ and has value $1/2$ at $3/2$. For
$1/2\le h\le3/4$, use $f>4/5$ and $1/E\le2/3$ to obtain
\begin{equation}\label{late:r24:c:L-bounds}
 L<2E+\frac15+\frac6{5E}\le\frac92.
\end{equation}

For \(h\le1/2\), set \(X=(17-20h)/300\) and \(H=2h/3\),
using \eqref{late:r24:b:cap-low-bounds}.
With \(F_0=HX+X^2\) the earlier bound is
\(g(E)+F_0<F_c(h)\), where \(F_c\) increases through \(h=1/2\)
and satisfies \eqref{late:r24:c:cap-loss-endpoint}.
For \(1/2\le h\le11/15\), take \(H=h(5-h)/8\) and
\(R_c=4h/(16h+5)\) from \eqref{late:r24:b:cap-high-bounds},
so \(x_0<R_c\alpha_0\).
The common cap-budget bound is
\eqref{late:r24:c:extended-cap-budget}.

\emph{The range $0<h\le1/3$.}
Take $H=2h/3$, $C_1=38/27$, and use \(X=(17-20h)/300\) as above. Since $H\le2/9$ and $L<4$,
$C_1>(H+L)/3$. By \eqref{late:r24:c:chord-factors}, we may use
\[
 X_\delta=X(1-q_0\delta),\qquad q_0=17/3.
\]
This is nonnegative throughout the enlarged feasible interval in
\eqref{late:r24:c:circle-budget}. The loss after $g(E)$ in
\eqref{late:r24:c:circle-envelope} is the quadratic
\[
 F_0+A\delta-B\delta^2,
\]
where
\[
 F_0=HX+X^2,\qquad
 A=\frac{h^2}{3}+C_1X-q_0HX-2q_0X^2,
 \qquad B=q_0X(C_1-q_0X)>0.
\]
Direct expansion gives
\[
 A=\frac{35122-180040h+433200h^2}{810000}.
\]
The numerator is convex, and its endpoint values on $[0,1/3]$ are
$35122$ and $23242$, both less than $40500$. Hence $A<1/20$.
Writing $v=q_0X$, one has
\[
 \frac{527}{2700}\le v\le\frac{289}{900}<\frac{C_1}{2}.
\]
The function $v(C_1-v)$ is increasing on this interval, so
\[
 B\ge\frac{527\cdot3273}{2700^2}
       =\frac{1724871}{7290000}>\frac15.
\]
Completion of the square, also valid when $A\le0$, gives
\[
 A\delta-B\delta^2
 \le\frac{\max(A,0)^2}{4B}<\frac1{320}.
\]
On $[0,1/3]$, the cap bound gives \(F_c(h)\le F_c(1/3)=251/22500\). Thus
\[
 g(E)+F_0+A\delta-B\delta^2
 <\frac{251}{22500}+\frac1{320}
 =\frac{29}{2000}-\frac{79}{360000}.
\]
Equation~\eqref{late:r24:c:circle-envelope} proves
\begin{equation}\label{late:r24:c:small-circle-margin}
 \Gamma>\frac{1471}{2000}+\frac{79}{360000}
       \qquad (0<h\le1/3).
\end{equation}

\emph{The range $1/3\le h\le1/2$.}
Take $H=2h/3$, $C_1=13/9$, with \(X\) as above. Again
$C_1>(H+L)/3$. By the exact chord and its first factor bound, a
nonnegative upper bound for $x$ is
\[
 X_\delta=X\left(1-\frac{3\delta}{\alpha_0}\right)
                      \left(1+\frac\delta3\right).
\]
It is decreasing and at most $X$, and
\[
 X_\delta'=-X\left(Q+\frac{2\delta}{\alpha_0}\right)
          \le-QX,\qquad
 Q=\frac3{\alpha_0}-\frac13
  >\frac{16}{3(3/4-h)}-\frac13,
\]
by \eqref{late:r24:c:extended-cap-budget}. Here and in the next
range, primes on $X_\delta$ mean differentiation in $\delta$ of the
displayed upper bound, not differentiation of a true minimum cost.
Set
\[
 F_\delta=\frac{h^2\delta}{3}
              +(H+C_1\delta)X_\delta+X_\delta^2.
\]
Its derivative obeys
\[
 F_\delta'=\frac{h^2}{3}+C_1X_\delta
                    +(H+C_1\delta+2X_\delta)X_\delta'
 \le\frac{h^2}{3}+C_1X-QHX<0.
\]
For the final sign, substitute the lower bound for $Q$ and
$X=(17-20h)/300$. Multiplication of the sufficient positive
difference by $2700(3-4h)>0$ gives exactly
\[
 P_{\rm mid}(h)=3440h^3-6044h^2+3738h-663.
\]
Its least degree-three Bernstein coefficient on
\([1/3,1/2]\) is \(1049/27>0\).
Hence $F_\delta\le F_0=HX+X^2$, and the cap bound \eqref{late:r24:c:cap-loss-endpoint} give
\begin{equation}\label{late:r24:c:middle-circle-margin}
 \Gamma>\frac{1471}{2000}+\frac1{202500}
       \qquad (1/3\le h\le1/2).
\end{equation}

\emph{The range $1/2\le h\le11/15$.}
Use \(H,R_c\) as above. Here $H<2/5$ and $L<9/2$, so
$C_1=49/30$ satisfies $C_1>(H+L)/3$. A nonnegative upper bound is
\[
 X_\delta=R_c(\alpha_0-3\delta)
                              \left(1+\frac\delta3\right).
\]
It is decreasing and at most $R_c\alpha_0$. Moreover
\eqref{late:r24:c:extended-cap-budget} gives $\alpha_0<9/64$, so
\[
 X_\delta'=R_c\left(\frac{\alpha_0}{3}-3-2\delta\right)
       \le-\frac{189}{64}R_c.
\]
For the same definition of $F_\delta$ as in the middle range,
\[
 F_\delta'\le\frac{h^2}{3}+C_1R_c\alpha_0
                      -\frac{189}{64}R_cH<0.
\]
Indeed division by $R_c>0$ gives
\[
 \frac{189}{64}H-\frac{h^2}{3R_c}
     =\frac{h(2195-2615h)}{1536},\qquad
 C_1\alpha_0<\frac{147}{160}\left(\frac34-h\right).
\]
Their difference
\[
 P_{\rm high}(h)=\frac{h(2195-2615h)}{1536}
                 -\frac{147}{160}\left(\frac34-h\right)
\]
is a concave quadratic. Its values at $1/2$ and $3/4$ are
$1819/30720$ and $935/8192$, respectively, both positive.
It is therefore positive on the larger interval $[1/2,3/4]$.
Thus \(F_\delta\le F_0=HR_c\alpha_0+R_c^2\alpha_0^2\).
By \eqref{late:r24:b:cap-high-gain}, \(F_0<r_c(h)\alpha_0\).
The feasible cap has \(E<E_*\). Integrating
\(g'>r_c(h)c_f'\) to the first zero of the cap budget
gives \(g(E_*)-g(E)>r_c(h)\alpha_0\), since \(r_c\) increases.
Consequently $g(E)+F_\delta<g(E_*)<29/2000$, and
\eqref{late:r24:c:circle-envelope} proves
$\Gamma>1471/2000$ also on this last range.

The three ranges include their overlap points $h=1/3$ and $h=1/2$
and exhaust $0<h\le11/15$. Every denominator used above is positive
on its stated branch; in particular $3/4-h>0$ throughout the required
range. The endpoint $b=1$ is included directly: $D=3$, the
strong-minimum lemma includes its zero signed-depth case, and the
$B$ equality-height inverse barrier remains strict. At $\delta=0$
the statement is the positive-$A$ cap and every bound agrees with
its cap value. At $\alpha=0$ the true $A$ inverse has $x=0$; the
same nonnegative upper bounds still apply. If $\alpha_0=0$, the
separate zero-budget case at the start applies and no division by
$\alpha_0$ is made. The enlarged intervals for $\delta$ are used
only to bound the displayed scalar upper expressions. Neither a
physical variation nor a derivative of $m_b$ has been inferred from
$m_b>3\delta$.
\end{proof}
\space
\endgroup
\space
\label{revision:appendix-end-5}
\end{proofappendix}
\section{Auxiliary comparisons for the reduction to six types}\label{app:run-details}
\begin{proofappendix}
\section{Comparisons with all negative runs of length two}
\subsection{Full-lift matching across the true seam}\label{new:glue}
We apply Proposition~\ref{pc:two-cuts} and Corollary~\ref{pc:short-rigid-global} on the following actual supports.
A complete capped binary core with two distinct internal binary ties
and its two third-mode cap ties excludes abnormality and identifies
the entire normalized lift by equality-only rigidity. This includes
the matrix and the state covector. On the overlap of two genuine
periods, $G\Lambda G^{-1}$ is constant; it first gives the matrix
twist. At an actual switch between distinct modes, the difference
of the state covectors annihilates two independent control fields.
It therefore vanishes in both coordinates. Adjoint uniqueness then
extends the full lift. Matrix twist alone would not prove this step.

For R07, the actual support $[0,5]=0|2121|0$, phased by $+1$,
is $1|0202|1$. Cuts inside arcs $-5$ and $0$ retain its closed core
$[\tau_1,\tau_5]$ and positive caps on their overlap, and the distinct
switch at $\tau_0$ lies between the cuts. A second source
$[6,10]=1|020|1$ has closed core $[\tau_7,\tau_{10}]$ disjoint
from the first on the period circle. Hence every genuine cut retains
at least one complete closed core, which supplies the stronger
compatibility of all required source lifts.

For R17, use $[4,8]+1=1|202|1$, with core
$[\tau_5,\tau_8]$, cuts inside $-4$ and $4$, and switch $\tau_4$.
The second source $[10,14]+2=1|202|1$ has disjoint closed core
$[\tau_{11},\tau_{14}]$. The same argument identifies all required
source lifts and their state covectors before any cross-seam
mean, signed entrance or determinant is compared.

\subsection{The four ranges for R07}\label{new:r07}
R07 has native word $0212101020$ and own reciprocals
\[
 (w,p,q_L,r_L,s,r_R,q_R,v,b,z).
\]
Use the common full lift of Section~\ref{new:glue}.
The sources $[0,5]+1$ and $[3,8]+2$ are both $102021$.
Their actual third-mode ties are $(\tau_1,\tau_5)$ and
$(\tau_4,\tau_8)$, respectively. Each source has its own cap at one
end of the shared arc. Corollary~\ref{loc:four-calibration}
gives $0<\Gamma:=\det\Lambda<H_*$ on this lift.
Put $x=s^2$, $H=(1+x)/2$ and $t_L=s+r_L$, $t_R=s+r_R$.
The positive shared arc gives $t_L,t_R>2s$, $t_Lt_R<2H$ and $x<1/3$.
Write $c_L=p,c_R=v$. With the divided mean
\[
 M(d;b,h)=\int_0^1\left(\frac{d-b-(h-b)t}{d+b+(h-b)t}\right)^2dt,
\]
the two mean equations are
$\mu_i=M(q_i;c_i,r_i)=M(r_i;s,q_i)$ and $K_i=(3+\mu_i)/2$.
They follow from the zero-to-zero integrals of
Lemma~\ref{loc:binary-normality} by changing to reciprocal shear.
The positive middle shears give $r_i>c_i,\ q_i>s$.
The signs of the two exterior internal roots give $q_i>c_i$ and $r_i>s$.

\paragraph{The full frame and the signed cap values.}
Put $A=(1+p^2)/2$, $B=(1+v^2)/2$ and
\[
 L=\frac{(p+s)^2}{2}-(r_L-p)(q_L-s),\qquad
 U=\frac{(v+s)^2}{2}-(r_R-v)(q_R-s),\qquad
 \Delta=H^2-LU.
\]
Write $\mathcal H_c=\operatorname{diag}(h_c,h_c^{-1})$ with
$h_c=(1+c^2)/2$, $J_a=\left(\begin{smallmatrix}a&1\\-1&0\end{smallmatrix}\right)$
and $N_o=\left(\begin{smallmatrix}1&0\\o&1\end{smallmatrix}\right)$.
The switch factors of Section~\ref{r06:source-details} are
$S(c,d)=\mathcal H_c^{-1}J_{(c+d)^2/2}$ for a forward switch and
$S(d,c)^{-1}$ for a reverse switch. The increment word
$\mathtt{--+--+--++}$ has sum $-2$, so the signed product is
$R_*^{1-(-2)}=-I$ in these bases. Its factors in order are
\[
\begin{aligned}
 &S(p,w)^{-1}S(q_L,p)^{-1}S(q_L,r_L)S(s,r_L)^{-1}\\[-2pt]
 &\quad{}\cdot S(r_R,s)^{-1}S(r_R,q_R)S(v,q_R)^{-1}
          S(b,v)^{-1}S(b,z)S(z,w).
\end{aligned}
\]
Define $D=((z+w)^2-(p+w)^2)/2$ and $O=((z+b)^2-(v+b)^2)/2$.
The identities $J_a^{-1}J_b=N_{b-a}$ and $N_aJ_b^{-1}=J_{b-a}^{-1}$
combine the two inner pairs with their following inverse factors to
produce $J_L^{-1},J_U^{-1}$, and combine $S(b,v)^{-1}S(b,z)$ to produce
$N_O$. Cancelling the end $J$ factors gives
\[
 \mathcal H_pJ_L^{-1}\mathcal H_sJ_U^{-1}\mathcal H_vN_O\mathcal H_z^{-1}
 =-\begin{pmatrix}1&D\\0&1\end{pmatrix}.
\]
The lower-right, lower-left and upper-right entries give
\begin{equation}\label{r07:full-frame}
 \Delta>0,\qquad h_z=ABH/\Delta,\qquad O=B^2L/\Delta,\qquad D=A^2U/\Delta.
\end{equation}
Positive arcs $0$ and $8$ give $z>p$ and $z>v$, respectively. Hence
$D,O>0$ and $L,U>0$.

In the left and right binary charts let $V_i$ be the state height and
$\mathcal L_i=l_i-u_i m_i$ the slope variable.
On the shared arc $V_LV_R=2/H$. The whole phase change in the moving
adjoint bases of Lemma~\ref{loc:binary-reduction} is
$\left(\begin{smallmatrix}sV_L/2&-1/H\\1/2&s/(HV_L)\end{smallmatrix}\right)$.
Conjugating the matrices with rows $(\mathcal L_i,K_i)$
and $(m_i,-\mathcal L_i)$ gives on this normalized lift
\[
 \mathcal L_L=\lambda+sK_L/V_L,\qquad
 \mathcal L_R=-\lambda-sK_R/V_R,\qquad
 \Gamma=H^2K_LK_R-\lambda^2.
\]
The inner cap heights are $V_R(\tau_4)=t_L/H$, $V_L(\tau_5)=t_R/H$.
Put
\begin{gather*}
 Y_L=Ap/(p+w),\quad Y_R=Bv/(v+b),\qquad Z_L=Y_L-L,\quad Z_R=Y_R-U,\\
 B_L=sH/t_R,\qquad B_R=sH/t_L.
\end{gather*}
The inner outward slopes are $\lambda+K_iB_i$. On the left,
subtracting the middle shears $q_L(r_L-p)/2$ and $r_L(q_L-s)/2$
from $\mathcal L_L(\tau_4)$ with $(\mathcal L_L)_\sigma=K_L$
gives $-\mathcal L_L(\tau_2)=K_L[p(p+q_L)/2-L]-\lambda$.
The outward leaf shear replaces $p(p+q_L)/2$ by $Y_L$. On the right,
adding the two middle shears from $\tau_5$ to $\tau_7$ gives
$\mathcal L_R(\tau_7)=K_R[v(v+q_R)/2-U]-\lambda$, and the outward leaf
replaces $v(v+q_R)/2$ by $Y_R$. Thus the outer outward slopes are
$K_iZ_i-\lambda$.

\paragraph{The low singleton range $s\le1/5$.}
For either core, write its internal roots as $-\xi<c<\zeta$.
By \eqref{loc:four-roots}, $c^2\le\mu<\min(\xi^2,\zeta^2)$.
The integrals of $u^2-\mu$ over $[-\xi,c]$ with weight $(1-u)^{-2}$
and over $[c,\zeta]$ with weight $(1+u)^{-2}$ both vanish.
If $c\ge0$, extend the first integral to $[-1,c]$. The added
$[-1,-\xi]$ part is strictly positive. If $\mu\ge3-4\log2$,
the $[-1,0]$ part is $(3-4\log2-\mu)/2\le0$ and the $[0,c]$ part
is nonpositive by $c^2\le\mu$, a contradiction. If $c<0$, extend the
second integral to $[c,1]$ and split at zero. The added $[\zeta,1]$
part is strictly positive, and the same two nonpositive bounds apply
on $[0,1]$ and $[c,0]$. Thus $\mu_i<3-4\log2$. Since
$\log2>2(1/3+1/81)=56/81$, $K_i<\overline K:=3-2\log2<13/8$.

Lemma~\ref{loc:outer-leaf-calculus} gives $E=-f>0$ at each cap.
The inward maximizing jets give strictly positive outward slopes,
as in the proof of Lemma~\ref{loc:four-core-apex}. At each cap
$m=-1/2-E$, so $\Gamma+\mathcal L^2=K_i(1/2+E)$.
Inserting the outward slopes and the common determinant gives, with
$d_i=1/2-H^2K_{\rm other}$,
\[
 K_iZ_i^2-2\lambda Z_i=d_i+E_{i,\rm out},\qquad
 K_iB_i^2+2\lambda B_i=d_i+E_{i,\rm in},
 \qquad E_{i,\rm out},E_{i,\rm in}>0.
\]
Here $H\le13/25$ and $H^2K_i<1/2$.
If $Z_i\le0$, positivity of the outward cap value
$P_{i,\rm out}=K_iZ_i-\lambda$ gives
$0<P_{i,\rm out}\le-\lambda$; the same determinant identity then
forces $H^2K_{\rm other}>1/2$, a contradiction. Hence $Z_i>0$.
Multiply the two displayed rows by $B_i,Z_i$ and add. The shared
constant cancels, giving $K_iZ_iB_i>d_i$.
Eliminating $w$ by $(z-p)w=D-h_z+A$ in~\eqref{r07:full-frame} gives
\[
 (p+w)Z_L=-\frac{(A-L)(z-p)}2
     -\frac{A^2H(H-B)}{\Delta(z-p)},
\]
and its right counterpart. Positive shared time gives
$1>r_Lr_R+s(r_L+r_R)>\max(2sp+s^2,2sv+s^2)$, so $A-L,B-U>0$.
Since $Z_L,Z_R>0$ and $z>p,v$, these rows force $v>s$, $p>s$.
Positivity of the two exterior numerators further gives
$U(A-L)>H(B-H)$ and $L(B-U)>H(A-H)$.
If either $L-(A-H)$ or $U-(B-H)$ were nonpositive, equality
would contradict the two strict rows immediately; in the negative
case both would be negative and their product would give $H^2<LU$, contrary to
$\Delta>0$. Hence $L>(p^2-x)/2$, $U>(v^2-x)/2$.
The cap-product floors consequently give
\[
 Z_L>c_0(v+s),\quad Z_R>c_0(p+s),\qquad
 c_0=\frac{1/(2\overline K)-H^2}{sH}.
\]
The last positive native arc gives
$(z+w)(z+b)>1+z^2$. The same frame has
$(z+w)^2=(p+w)^2+2A^2U/\Delta$ and
$(z+b)^2=(v+b)^2+2B^2L/\Delta$.
Square the last-arc inequality and substitute these identities.
With $X=Y_L^2/p^2$ and $Y=Y_R^2/v^2$, cancellation of positive factors
gives $(\Delta/X+2U)(\Delta/Y+2L)>4H^2$. Expand and use
$\Delta=H^2-LU$ to obtain
\[
 F_LF_R<H^2,\qquad
 F_L=\frac{2Y_L^2}{p^2}-L,\quad
 F_R=\frac{2Y_R^2}{v^2}-U.
\]
For $c_0\ge\kappa\in\{1/2,1/3\}$ the preceding deficit floors
and an increasing quadratic give
\[
 F_L>2\kappa(v+s)-\frac x2+
       \frac{[2\kappa(v+s)-x]^2}{2p^2},
\]
and the symmetric right inequality. The endpoint rational comparisons
$4/13-1591/5184>0$ and $4/13-572/1875>0$ give
$c_0>1/2$ for $s\le1/6$ and $c_0>1/3$ for $s\le1/5$.
In the first range, the two factors are positive and their product
exceeds
\[
 (v+v^2/(2p^2))(p+p^2/(2v^2))
 \ge(p+1/2)(v+1/2)>H^2.
\]
In the second range put $d=s(2/3-s)\ge1/12>1/16$.
Cauchy's inequality and $t+d/t\ge2\sqrt d$, $t=\sqrt{pv}$,
give $\sqrt{F_LF_R}>2/9+(4/3)\sqrt d>5/9>13/25\ge H$.
Both conclusions contradict the last-arc product. This excludes the
entire low range, without an anchor-order assumption.

\paragraph{The high singleton range $x\ge1/6$.}
On the left chart of the shared arc the values of $1-u$ at the
internal root and at its cap are $a_1=2s/t_L$ and $a_2=st_R/H$.
The integrated outward switching term at this cap is
\[
 \psi=x\bigl[(1-\mu_L)(a_2^{-1}-a_1^{-1})
                   +2\log(a_2/a_1)-a_2+a_1\bigr] .
\]
The cap identity of Section~\ref{pf:M002} equates
$-H^2K_R-2\lambda x/a_2-K_Lx^2/a_2^2+1/2$ with
$x[\mu_L-(1-a_2)^2]/(2a_2^2)+H\psi/a_2$.
Solving for $\lambda$ gives the unsquared law
\begin{equation}\label{r07:lambda}
 \lambda=\frac H{2s}\left\{(2-K_L)t_L+(2-K_R)t_R
               -2x\left(\frac1{t_L}+\frac1{t_R}\right)\right\}
        +H\log\frac{2H}{t_Lt_R}-\frac12.
\end{equation}
The shared-leaf cap identity \eqref{pf:positive-cap}, with strictly
positive shared duration, gives $\lambda>-3x/2$.
The determinant $H^2K_LK_R-\lambda^2$ increases in each mean
factor throughout the rectangle $[3/2,K_L]\times[3/2,K_R]$:
for example its first derivative is greater than
$(3H/2)(H-st_L)>0$.
At the zero-mean comparison denote the value of \eqref{r07:lambda} by
$\lambda_0$. Its derivative in each $t_i$ is
$H(t_i-2s)^2/(4st_i^2)>0$. The lower value has both $t_i=2s$.
For an upper bound increase the product to $2H$. On that product
boundary the coefficient of $t_L+t_R$ is $H/(4s)-s/2>0$,
and the sum is maximal at $(2s,H/s)$ or its exchange. Thus
\[
 \lambda_{00}:=H\log\frac{1+x}{4x}-\frac12
       <\lambda_0<\frac{1-6x-15x^2}{16x}<0.
\]
Only after these signs are proved may one square, obtaining
$\Gamma>9H^2/4-\lambda_{00}^2$.
Set $z_0=(1-3x)/(1+5x)\in(0,3/11]$ and
$\ell=2z_0+2z_0^3/3<\log((1+x)/(4x))$.
Then $\Gamma>3/4-D_2(z_0)$, where
\[
 D_2(z)=\frac{8z^3+16z^4+(44/3)z^5+(112/9)z^6
                       +(32/9)z^7+(16/9)z^8}{(3+5z)^2}.
\]
Every $z^n/(3+5z)^2$, $n\ge3$, increases. The endpoint gives
\[
 \Gamma>\frac{83376585}{113379904}>
       \frac{1470471}{2000000}>H_*,
\]
contrary to harmonic comparison on this actual lift.

\paragraph{The double-lower range $p,v\le s$.}
The cap derivative used here also applies in the mixed range. Write
$\eta_i=(q_i-c_i)/(q_i+c_i)$ and
$C_i=c_i(c_i+q_i)/2-L_i$, where $L_L=L$, $L_R=U$.
At the internal root, the outward slope of the maximizing gap is
\begin{equation}\label{r07:cap-slope}
 \mathscr D_i=\eta_i+\frac{(\eta_i^2-\mu_i)(1-c_i^2)}{2(1-\eta_i)}
                 -2K_i C_i+2\lambda .
\end{equation}
Indeed the outward derivative of $m$ is $-2(K_iC_i-\lambda)$.
At that root $f=(\mu_i-\eta_i^2)q/2$ with $q=c_i^2/(1-\eta_i)^2$.
The outward derivative in Lemma~\ref{loc:outer-leaf-calculus}
then gives \eqref{r07:cap-slope} for $m+1/2-f$.

Here $\mathscr D_i<0$ on the actual outer leaf. The exterior
three-core sources are $[6,10]=1|020|1$ and
$[8,12]+2=1|202|1$, for $b$ and $w$, respectively.
With $d_L=w,d_R=b$, the strict internal root on the $c_i$ side
and the weak one on the $z$ side give
$(d_i-c_i)/(d_i+c_i)>|(d_i-z)/(d_i+z)|$ by
Lemma~\ref{loc:return-root} and its mean comparison.
In particular $d_i>c_i$. If $d_i<z$ the displayed comparison gives
$d_i^2>c_i z$ by cross multiplication, and if $d_i\ge z$ it follows
from $z>c_i$. With $a_0=c_i(c_i+q_i)/2$ and outward shear
$\sigma$, the outward state coordinate $u$ satisfies
$q(1-u)=a_0+\sigma$ and $u'=c_i^2/(a_0+\sigma)^2$.
Differentiating the outward $f$ formula in
Lemma~\ref{loc:outer-leaf-calculus} gives
$f''=1-\mu_i-(1+c_i^2)/(a_0+\sigma)$. Since $m''=-2K_i$, the gap
satisfies $g''=(1+c_i^2)/(a_0+\sigma)-4$. At the actual cap
$a_0+\sigma=c_i(1+c_i^2)/(2(c_i+d_i))$, so $g''=2(d_i/c_i-1)>0$.
This is the minimum of $g''$ on the leaf. The initial gap is
nonnegative, the terminal gap is zero, and the leaf has positive
length. Strict convexity therefore gives $\mathscr D_i<0$.

For a lower anchor $c_i\le s$, suppress the side index and set
\[
 R=r_i/s,\quad Q=q_i/s,\quad T=1+R,\quad e=(R-1)/(R+1),\quad
 \alpha=c_i/s,\quad \eta=(Q-\alpha)/(Q+\alpha).
\]
The root comparison below gives $1<Q\le R$. When $Q<R$,
both endpoint squares on $[Q,R]$ are below the common mean
for their respective integrals. Removing this interval from
$M(Q;\alpha,R)$ and adding it to $M(R;1,Q)$ gives
\[
 m(e)<\mu<m(\eta),\qquad
 m(v)=v^2\int_0^1\frac{(1-t)^2}{(1+vt)^2}\,dt .
\]
Indeed each inserted or removed integrand is monotone on that
interval and is bounded by the middle-root square.
When $Q=R$, calibration gives $\alpha=1$ and $e=\eta$.
In both cases $m'(v)>0$ for $v>0$ and $m(v)<v^2/3$ give
$e\le\eta$ and $\mu<\eta^2/3$. Moreover
$0<L_i<(c_i+s)^2/2\le2x<H$.

For each side write $P_i=(1-\mu_i)T_i/2-2/T_i$ and
$F_i=P_i-2\log T_i$. In \eqref{r07:cap-slope} replace only $\lambda$
by $\lambda_i^*=H[P_i+\log((1+x)/(xT_i^2))]-1/2$.
Calling the resulting expression $\mathscr D_i^*$, the actual
shared law gives $\mathscr D_i=\mathscr D_i^*+H(F_{\rm other}-F_i)$.
The coefficient of the formal mean in $\mathscr D_i^*$ is
\[
 -\left(\frac{1-c_i^2}{2(1-\eta)}+C_i+HT\right)<0,
\]
since $C_i=a_0-L_i$ and $T>2$. Replacing $\mu$ by $\eta^2/3$
strictly lowers $\mathscr D_i^*$. Put $\beta=(1-\eta)/(1+\eta)$
so that $\alpha=\beta Q$ and $\beta R\le1$. Increase $Q$ to $R$
at fixed $R,\eta,x$, carrying $\alpha=\beta Q$ in this scalar expression.
Its $Q$ derivative divided by $x$ is
\[
 2\beta Q\left[K_0-\frac{\beta\eta^2}{3(1-\eta)}\right]-2K_0R<0,
 \qquad K_0=(3+\eta^2/3)/2.
\]
The same derivative applies when the other side's contribution is
held fixed. Put $y=x/(1+x)$ and write the resulting coefficients as
\[
 \begin{gathered}
 A_0=\eta+\frac{\eta^2}{3(1-\eta)},\qquad
 P_0=(1-\eta^2/3)T/2-2/T,\qquad E_0=A_0+P_0-1,\\
 d=1-\eta+K_0[1+2R+(\beta-2)R^2]
                  -\frac{(1+\beta^2R^2)\eta^2}{3(1-\eta)} .
 \end{gathered}
\]
After division by $2H$, the substituted $\mathscr D_i^*$ is
$E_0+dy-\log(T^2y)$. The bracket $C=1+2R+(\beta-2)R^2$ is at most $2R$,
and $\partial_\eta(K_0 C)\le2\eta R/3-3R^2/(1+\eta)^2<0$.
Its value at $\eta=0$ thus gives
$d<1+3(1+2R-R^2)/2<T^2$.

The function $\mathscr B=E_0+d/T^2$ is increasing in $\eta$ for
$e\le\eta<1$. In fact $\partial_\eta\mathscr B=R W/T^2$, where
\[
\begin{aligned}
 W={}&(R+2)\left[1+\frac{\eta(2-\eta)}{3(1-\eta)^2}
                         -\frac{\eta(R+1)}3\right]\\
 &-\frac{R(9+12\eta+3\eta^2+7\eta^3+3\eta^4)}{3(1+\eta)^3}
      +\frac{\eta}3 .
\end{aligned}
\]
At fixed $\eta$, this is concave in $R\in[1,(1+\eta)/(1-\eta)]$.
After multiplication by $3(1-\eta)^2(1+\eta)^3$ at the first end
and $3(1-\eta)^3(1+\eta)^2$ at the second, the values are
\[
 4\eta(4+\eta-3\eta^2+6\eta^3-2\eta^5),\qquad
 4\eta[(1+\eta)(1-\eta)^2+\eta^3(4-\eta-\eta^2)].
\]
Both are positive. Substitution of $\eta=e$ gives
$\mathscr B(R,e)=B_i^*:=e^2(1+3e+e^2+e^3)/[6(1-e)]>0$.
When $s+r_i\le\sqrt{1+x}$, we have $y\le1/T^2$.
Since $d<T^2$, the scalar expression is decreasing up to this
closed endpoint. Hence $\mathscr D_i^*>2HB_i^*$ and the
actual negative derivative gives $F_i-F_{\rm other}>2B_i^*$.
Adding the two inequalities excludes two short sides, including
threshold equality. The shared product excludes two long sides.
Exactly one is strictly short and one strictly long.

Name the short and long normalized widths $T<V$. On the actual short
side, the other mean has coefficient $-HV/2$ in $\mathscr D_i$.
The root bound $\mu_{\rm other}<((V-2)/V)^2$ permits replacement
of its $P_{\rm other}$ by $2-4/V$ for a lower bound.
Carry out the same short-mean and $Q$ substitutions as above.
For the short-side derivative they give
\[
 \frac{\mathscr D_i}{2H}>
 A_0+P_0/2+dy-2/V-\log(TVy).
\]
On $0\le y\le1/T^2$, $\partial_\eta(E_0+dy)>0$.
At $y=1/T^2$ this was proved by $W>0$. At $y=0$ it follows
from $T\le2/(1-\eta)$, giving
$E_0'\ge1+\eta^2/[3(1-\eta)^2]$. The derivative is affine in $y$.
As $P_0'=-\eta T/3<0$, the displayed lower bound also increases in
$\eta$. Now $y<1/(TV)$ by positive shared duration and $d<T^2<TV$.
Lower $\eta$ to $e$ and raise $y$ to $1/(TV)$. Direct substitution
gives the strictly lower value
\[
 \frac{2e(3-2e)}{3(1-e)}
 -\frac{e(e^3+2e^2+5e+12)}{3(1-e)V}.
\]
Because $\mathscr D_i<0$, we obtain
\[
 V<V_{\rm crit}(e):=\frac{e^3+2e^2+5e+12}{2(3-2e)},\qquad
 e^3+e^2+3e-1<0.
\]
The latter follows by subtracting $T=2/(1-e)$ from $V_{\rm crit}(e)$.
The cubic is increasing and positive at $37/125$.
The function $V_{\rm crit}$ is increasing on $(0,1)$ and
$V_{\rm crit}(6/25)^2<7<V^2$ by the long-side condition and $x<1/6$.
At the root of the cubic, $V_{\rm crit}=2/(1-e)$. These comparisons give
\[
 50/19<T<V<125/44,\qquad
 1936/13689<x<1/6 ,
\]
using $x>1/(V^2-1)$ in the last lower bound.

For $a=p/s,\ c=v/s\le1$, \eqref{r07:full-frame} gives $h_z>A,B$.
Thus $LU>H(H-B)$. Since
$L<x(1+a)^2/2$, $U<x(1+c)^2/2$, this implies
$x(1+a)^2(1+c)/(1+x)>1-c$, including $c=1$.
The left side is less than $4(1+c)/7$. Hence $c>3/11$,
and the other full-frame row gives $a>3/11$. The remaining candidate
is therefore contained in
$[3/11,1]^2\times[3/2,2]^2\times[1/8,1/6]$
in the coordinates $(c_i/s,c_{\rm other}/s,r_i/s,r_{\rm other}/s,x)$,
with the short side named first.
On this root \path{local-r07} checks the exact partition on the
\(2^{-96}\) lattice with directed arithmetic at scale \(2^{-128}\).
The tests contradict the strict collars, positive duration, full-frame
or actual cap-slope conditions. The certified mean roots and logarithm
lower bounds retain the arithmetic contract of this root.

\paragraph{The remaining mixed and upper ranges.}
Put
\[
 a=p/s,\quad c=v/s,\quad R=r_L/s,\quad S=r_R/s,
 \quad Q=q_L/s,\quad V=q_R/s.
\]
The proper returned eight-sources $[-1,6]+1$ and $[2,9]+2$ are both
$01020212$. Their endpoint density pairs are $(p^2,s^2)$ and $(s^2,v^2)$.
The retained returns and caps give $p^2,s^2,v^2<47/250$ by
\eqref{pf:d8-density-ceilings} on the common harmonic lift.
The low and high exclusions already give $1/25<x<1/6$.
The full frame gives $a,c>3/13$: using its positive deficits,
\[
 \frac{x(1+a)^2(1+c)}{1+x}>1-c,\qquad
 \frac{(s+p)^2}{1+x}<5/8,
\]
where the latter has the rational numerator reserve $29/2000$.
The actual calibrated root satisfies $Q\ge\min(a,1)R$, as shown below.
For $R\ge6$ this contradicts the left positive deficit, since
\[
 (6-a)(6a-1)>(1+a)^2/2\quad(3/13\le a<1),
\]
and $5(6-a)>(1+a)^2/2$ for $1\le a\le11/5$.
The upper density bounds also give $a,c<11/5$. Hence every remaining
candidate lies in the closed outward enlargement
\begin{equation}\label{new:r07root}
 (a,c)\in[3/13,11/5]^2,\quad(R,S)\in[1,6]^2,
       \quad x\in[1/25,1/6].
\end{equation}

For the root enclosure put $F(q;a,R)=M(q;a,R)-M(R;1,q)$.
At a zero on the full physical ray $q>\max(1,a)$, let $\mu$ be
the common mean, $f_a=((q-a)/(q+a))^2$ and
$f_R=((q-R)/(q+R))^2$.
If $q<R$, the second integrand is decreasing, so $\mu>f_R$;
the first has no interior maximum, so $f_a>\mu$.
If $q>R$, the first is decreasing and gives $f_a>\mu>f_R$
directly. The endpoint $q=R$ is a root only when $a=1$.
Homogeneity gives
\[
 qM_d(q;a,R)=\frac{R(\mu-f_R)+a(f_a-\mu)}{R-a}>0,
 \qquad \partial_qM(R;1,q)=\frac{f_R-\mu}{q-1}<0.
\]
Thus every zero is an upward crossing, so there is at most one.
The sign of $F(R)$ is the sign of $1-a$, because the
integrand on the two intervals ending at $R$ is strictly decreasing.
At $q=aR\ge1$, scaling identifies the first mean with $M(R;1,R/a)$.
On the added part between $R/a$ and $aR$ the integrand is at most
its equal values at those two endpoints, which are at most the mean
over the shorter interval. The integrand is strictly smaller at
interior points of the added part. Thus the longer mean is strictly
smaller, even when the shorter interval is a point. Hence
$F(aR)>0$ if $a>1$ and $F(aR)<0$
if $a<1$ and $aR\ge1$. Thus for $a>1$ the root lies in $(R,aR)$.
For $a<1$ it lies in $(1,R)$ exactly when $F(1)<0$.
At $F(1)=0$ the first mean has a strict left endpoint excess
over its actual common value $f_R$. Its $q$ derivative is positive
by homogeneity. The derivative of the second mean is half the
derivative of its integrand at $1$ and is negative. Thus also
$F_q(1)>0$, and the continuation clips to $1$ when $F(1)\ge0$.
When $a=1$ the root is $R$. The sign at $aR$ on the lower branch,
or $q>1$ when $aR<1$, proves $q\ge\min(a,1)R$.
For completeness, at a nonendpoint root
\[
 F_a=\frac{\mu-f_a}{R-a}<0,\qquad
 F_R=\frac{f_R-\mu}{R-a}-M_d(R;1,q)<0.
\]
The second homogeneity identity gives $M_d(R;1,q)>0$ by the
same endpoint comparison. Thus implicit differentiation makes
the root increasing in both $a,R$. The positive boundary derivative
and the bounded brackets give continuous clipping at $q=1$.
For $R=a>1$ the first mean is a point mean, with $F(a)<0<F(a^2)$.
At its root $f_R=\mu$, the first mean has positive own-reciprocal
derivative and negative integration-endpoint derivatives, each
endpoint derivative being half the integrand derivative.
The second mean has positive own-reciprocal derivative by its
strict other endpoint excess. These signs give the same implicit
monotonicity on approach to the diagonal. On $R=1,\ a\le1$, for
$q>1$ the first mean is at least $f_R$ and the second is strictly
smaller than $f_R$. The clipped value on that edge is $1$.
Consequently any actual root in a parameter box obeys
\[
 \widehat q(a_-,\max(a_-,R_-))\le q_{\rm actual}
       \le\widehat q(\min(a_+,R_+),R_+).
\]
Bisection changes a root endpoint only at a certified strict sign;
an uncertain sign retains the whole current bracket.

\begin{lemma}[R07 closed necessary-system cover]\label{new:r07finite}
No actual original minimizing R07 candidate lies in
\eqref{new:r07root} after the double-lower exclusion.
\end{lemma}
\begin{proof}[Computer-assisted proof]
We specify the correlated rows which carry the finite contradiction.
Put
\[
 l=(1+a)^2/2-(R-a)(Q-1),\quad
 u=(1+c)^2/2-(S-c)(V-1),\quad L=xl,\quad U=xu,
\]
$A=(1+a^2x)/2$, $B=(1+c^2x)/2$, $K_L=(3+\mu_L)/2$,
$K_R=(3+\mu_R)/2$, with the two calibrated means.
The shared law above computes the same unsquared $\lambda$.
Set $D_0=H^2K_LK_R$, $k=K_LL$, $j=K_RU$,
$b(y)=3(1+y)/4$ and $c(y)=9y(1-y)/8-3/4$.
The determinant is $\Gamma=D_0-\lambda^2<1471/2000$ since
$\Gamma<H_*$. Transport across an own pair with shear
$[q_i(r_i-c_i)+r_i(q_i-s)]/2$ gives the sum of its two
slot constants as $-K_i L_i$. Indeed one subtracts
$c_i(c_i+q_i)/2+s(s+r_i)/2$ from this shear, giving $-L_i$.
Thus the three constants are $\lambda,-k-\lambda,-j-\lambda$.
The three actual defects are
$\lambda+3x/2$, $-k-\lambda+3a^2x/2$ and
$-j-\lambda+3c^2x/2$. At their respective densities $y$ they are
strictly greater than the positive solution of
$e^2+3(1-3y)e/4=29/2000$. Here \eqref{pf:gap} applies at each
of the three actual densities, $e\ge0$ by \eqref{pf:d8-density-ceilings},
and $3/4-\Gamma>29/2000$. This solution is
\[
 \phi(y)=\frac{29/1000}{
 \sqrt{[3(1-3y)/4]^2+29/500}+3(1-3y)/4}.
\]
For a slot constant $\ell$, substitute $e=\ell+3y/2$ in
\eqref{pf:gap} to obtain $\Gamma+\ell^2+b(y)\ell+c(y)\ge0$.
Retain these strict floors and the resulting closed rows
\begin{gather*}
 b(x)\lambda+D_0+c(x)\ge0,\\
 [2k-b(a^2x)]\lambda+D_0+k^2-b(a^2x)k+c(a^2x)\ge0,\\
 [2j-b(c^2x)]\lambda+D_0+j^2-b(c^2x)j+c(c^2x)\ge0.
\end{gather*}
These affine rows are never divided by a coefficient of unknown sign.
The actual cap slopes are strictly negative by \eqref{r07:cap-slope}
and the complete outer-leaf argument. In the normalized variables
the left slope is
\[
 \eta+\frac{(\eta^2-\mu_L)(1-a^2x)}{2(1-\eta)}
 -(3+\mu_L)x\{a(a+Q)/2-l\}+2\lambda,
 \qquad \eta=(Q-a)/(Q+a).
\]
The right formula exchanges $(a,R,Q,l,\mu_L)$ with
$(c,S,V,u,\mu_R)$.

Equation~\eqref{r07:full-frame} determines the exterior reciprocals:
\[
 \Delta=H^2-LU>0,\quad J=ABH/\Delta,\quad
 Z=z/s=\sqrt{(2J-1)/x},
\]
\[
 W=w/s=\frac{A^2u/\Delta-(J-A)/x}{Z-a},\qquad
 B_o=b/s=\frac{B^2l/\Delta-(J-B)/x}{Z-c}.
\]
The actual three-core means on $[8,12]+2$ and $[6,10]$ are
$\mu_w=M(W;a,Z)$ and $\mu_b=M(B_o;c,Z)$. The weak root
jets at their $z$ ends and the strict roots at the other ends give
\[
 \mu_w\ge\left(\frac{Z-W}{Z+W}\right)^2,\qquad
 \mu_b\ge\left(\frac{Z-B_o}{Z+B_o}\right)^2,
 \quad W>a,\quad B_o>c.
\]
For the $w$ overlap, the same-lift determinant row
\eqref{pf:shared-geometry} is
$\Gamma=A^2K_L(3+\mu_w)/2-(k+\lambda)^2$.
Equating it with $D_0-\lambda^2$ and dividing by $K_L>0$ gives
\[
 A^2(3+\mu_w)/2-\Delta K_R-L(k+j+2\lambda)=0,
\]
The $b$ overlap gives the right counterpart. The actual outer
cap durations give $q_L<(1-pw)/(p+w)$ and $q_R<(1-vb)/(v+b)$.
The third reciprocals at the exterior switches give $wz,bz<1$.
The final native duration gives $(z+w)(z+b)>1+z^2$.

The remaining \path{local-r07} cover uses directed integer
enclosures at scale \(2^{-96}\) on the fixed dyadic partition.
Its divided-mean formula bounds the absolute even-moment coefficients
by \(4\), hence its tail after \(24\) terms by
\(4\rho^{50}/(1-\rho^2)\), where
\(\rho=(h-b)/(2d+b+h)\). The logarithm uses \(32\) positive
atanh terms and the geometric tail. The task checks complete closed
splits and one necessary contradiction at each terminal.
A strict row can be rejected at zero upper endpoint; a closed
affine row needs a strictly negative upper endpoint.
\end{proof}
The low, high, double-lower and mixed/upper ranges exhaust all actual
candidates. Their boundaries have been retained throughout, so R07
is excluded.

\subsection{The entrance terms for R17}\label{new:r17}
The native word is $021201210120$, with shared owners $1,7$ and
opposite centers $4,10$. Section~\ref{new:glue} identifies the entire
full lift across its true twist. Define the physical unions
\[
 U_1=c_{11}+c_0+c_1+c_2+c_3,\qquad
 U_7=c_5+c_6+c_7+c_8+c_9.
\]
There is the exact partition
\begin{equation}\label{new:r17partition}
C=U_1+U_7+c_4+c_{10}.
\end{equation}
The proper H3 nine-sources inside the actual eleven-sources give
$U_1,U_7\ge U_0=180483/250000$. The ordinary complete opposite
centers give $c_4,c_{10}>63/500$. Therefore each center separately
satisfies
\begin{equation}\label{new:r17window}
63/500<c_i<117483/500000\qquad(i=4,10).
\end{equation}
For each upper bound, the other center contributes its independent
$63/500$ floor; the two centers need not have equal costs.

\paragraph{Exclusive costs and the four actual means.}
The actual pair $\{11,0\}$ carries mean $\mu_0$ and entrance
$\lambda_{10R}$; $\{2,3\}$ carries $\mu_2,\lambda_{4L}$;
$\{5,6\}$ carries $\mu_6,\lambda_{4R}$; and
$\{8,9\}$ carries $\mu_8,\lambda_{10L}$.
All four entrance slopes are strictly positive by the complete
signed-strip source theorem. The exact original exclusive integral
contains $K\Sigma=L_{\rm exit}+\lambda$, where
$K=(3+\mu)/2$. Retaining this equality, instead of dropping its
entrance term, gives for the same physical exclusive pair
\[
E>\frac\alpha2+\frac{1-p}{2}+\frac{3x}{2p}
       -\frac{3\mathcal B}{3+\mu}+\frac{3\lambda}{3+\mu}.
\]
The strictness comes from the actual positive shear-square integral;
$\alpha\ge0$ is the original source entrance width. Adding the two
exclusive pairs and the shared physical arc gives the original
$J$ cancellation with both entrance terms retained. Hence set
\begin{align*}
G_4&=c_4+\frac{3\lambda_{4L}}{3+\mu_2}
           +\frac{3\lambda_{4R}}{3+\mu_6},\\
G_{10}&=c_{10}+\frac{3\lambda_{10L}}{3+\mu_8}
           +\frac{3\lambda_{10R}}{3+\mu_0}.
\end{align*}
Equation~\eqref{new:r17partition} and the exact cost algebra imply
\begin{equation}\label{new:r17regroup}
 C>J_1+J_7+G_4+G_{10}.
\end{equation}
The $G_i$ are regroupings of terms already present in disjoint
physical costs. They are not new fees, nor bounds for the center
arc $c_i$ alone. The discarded four $\alpha/2$ terms are nonnegative.

\paragraph{The outer bound.}
The actual normal capped sources at owners $1$ and $7$ imply
$x<1/6$. Their normality is established on the required full lift;
the density conclusion is not an extra payment. For either outer
profile the complete retained-cost expression is $J$.
The scalar argument below excludes $x\le1/10$, where
$J>511/880>14/25$. The remaining closed root is
\[
x\in[1/10,1/6],\qquad p,q\in[4/11,1].
\]
For clarity, the exact function to be bounded is the returned-union
quantity in \eqref{ret:J-cost}. Put $A_0=1+x$ and, for its two
actual means $m,n$,
\begin{gather*}
 0<m<(1-p)^2,\qquad0<n<(1-q)^2,\qquad pq>4x/A_0,\\
 \mathcal V=\log\frac{A_0}{4x}+g(p)+g(q),\qquad
 g(t)=\log t+(t^{-1}-t)/2,\qquad M=m/p+n/q,\\
 \mathcal B=(1-A_0\mathcal V)/2+A_0M/4,\qquad
 \delta=\frac m{3+m}+\frac n{3+n},\\
 \Gamma=\frac{A_0^2(3+m)(3+n)}{16}-\mathcal B^2,\qquad
 J=\mathcal V-A_0M/2+\mathcal B\delta
   =\frac1{A_0}-\frac{xM}2-\mathcal B(2/A_0-\delta).
\end{gather*}
Here $\mathcal B$ is the same actual transport constant retained
in each exclusive cost. The mean rectangle gives
$M\le(1-pq)^2/(pq)$ and therefore
$1/A_0-xM/2\ge93/112>14/25$. Thus $J\le14/25$ forces
$\mathcal B>0$ and $J>\mathcal V-A_0M/2$.
For
$f_x(t)=\log t+A_0-x/(2t)-(2+x)t/2$,
\begin{align*}
 f_x'(t)&=\frac{(1-t)((2+x)t+x)}{2t^2}\ge0,\\
 f_x(p)+f_x(q)-f_x(pq)
 &=(1-p)(1-q)\left(\frac{x}{2pq}+\frac{2+x}{2}\right)\ge0.
\end{align*}
It follows that
$J>(7-18x-9x^2)/[8(1+x)]$, proving the stated exclusion
of $x\le1/10$. Also $p,q\ge4/11$ on the remaining root.

For the actual two-mean rectangle consider
\[
\mathcal L=J+60x(\Gamma-H_*).
\]
Its Hessian is strictly negative definite throughout that rectangle.
Put
\[
 K=-\operatorname{diag}(p,q)\nabla^2\mathcal L
 \operatorname{diag}(p,q).
\]
Direct differentiation gives
$K_{12}\ge1867/7200>0$ and
$K_{11}-K_{12},K_{22}-K_{12}\ge119/3600>0$.
The bound $A_0\mathcal V\le21/16$ gives
$\mathcal B\ge-5/32$ on the entire comparison rectangle; its
positivity at the actual counterexample is not imposed at its corners.
The minimum of $\mathcal L$ is therefore bounded below
by the minimum at the four corners. Only the actual determinant is
calibrated. A four-corner lower bound
$\mathcal L>14/25$ contradicts $J\le14/25$ and
$\Gamma<H_*$ at the actual point.

The directed-integer cover uses $50$ fractional bits
and $56$ logarithm-series terms. Its finite tree has $165373$ nodes and $82687$ leaves, with all
four corners checked at each accepted leaf and no unresolved box. It proves
\begin{equation}\label{new:r17outer}
J_1,J_7>14/25.
\end{equation}
Indeed the four-corner inequality gives
$\mathcal L_{\rm actual}>14/25$, whereas
$\Gamma_{\rm actual}<H_*$ gives
$\mathcal L_{\rm actual}<J_{\rm actual}\le14/25$.

\paragraph{The center bound and its compact root.}
The center supports are $[0,8]+1$ and $[6,14]+2$, both
$102012021$. In the mode-$1$ coordinates of
Section~\ref{ctr:section}, let $a>b>0$ be the actual endpoint
reciprocals, $y$ the central own reciprocal and $q=V^2$ its crossing
height squared. Put
\[
 y=\frac{q-1}{2\sqrt q},\quad h=\frac{q+1}{2\sqrt q},\quad
 r=\frac{a+y}{h},\quad s=\frac{b+y}{h},\quad
 \rho=y(a-b),\quad\theta=\frac{r-1}{r-s}.
\]
The complete forward contexts $2012$ and $0120$ give $r>1>s>0$ and
$y^2+2ay>1$, $b<y$ by the endpoint orders in
Lemma~\ref{ctr:height}. The positive third reciprocal at the first
end gives $ay<1$. Both signed entrance slopes are strictly positive.
The full matrix is reconstructed from
$H(a)^{-1}N(y(a-b)/2)H(b)$ and the actual state-adjoint cost row
\eqref{ctr:affine-cost}. The endpoint cap equalities
\eqref{ctr:cap-m} and maximizing rows \eqref{ctr:switch-rows}
remain affine in the two actual returned means.

Original compactness gives $q\le243/4$ at this actual crossing.
If $q\le5/2$, then $y^2\le9/40$, so the decreasing rational
bound in \eqref{ctr:geometric-height} gives
$c>3991/15680>39/160$, contradicting \eqref{new:r17window}.
Whole reflection preserves the physical
cost, determinant and $G$, exchanges the two means and slopes, and
sends $(r,s)$ to $(1/s,1/r)$ and $\rho$ to $\rho/(rs)$.
Choose the orientation with $rs\ge1$, recomputing $\rho$ there.
The exact physical cost in these coordinates is
\[
 c=\rho\left\{1+\frac1{rs}
       -\frac{q-1}{q+1}\frac{\log(r/s)}{r-s}\right\}.
\]
The density in \eqref{ctr:density} is greater than $y$ on $[b,a]$
by $ay<1$, so $c>\rho$. The displayed cost and $rs\ge1$ give $c<2\rho$.
The strict positive center floor supplies the closed comparison root
\begin{equation}\label{new:r17root}
q\in[5/2,243/4],\quad
\rho\in[63/1000,117483/500000],\quad
\theta\in[1/2,1],\quad(\mu_L,\mu_R)\in[0,1]^2.
\end{equation}
Write $Y_0=q-1$, $Q_0=q+1$, $d=r-s=4q\rho/(Y_0Q_0)$,
$A_0=Q_0r-Y_0$ and $B_0=Q_0s-Y_0$. Thus the first two actual
reciprocals are $a=A_0/(2\sqrt q)$ and $b=B_0/(2\sqrt q)$.
The endpoint orders above give the following necessary rows:
\[
 Y_0^2+2A_0Y_0>4q,\quad B_0<Y_0,\quad A_0Y_0<4q,\quad B_0>0.
\]
The identity $rs=1+(2\theta-1)d-\theta(1-\theta)d^2$ gives the
orientation row $2\theta-1-\theta(1-\theta)d\ge0$.
For $q<3$, the first two endpoint orders give, respectively,
$\theta\rho>(q+1)(3-q)/(8q)$ and
$(1-\theta)\rho>(q-1)(3-q)/(4q)$. The program uses these as
necessary contractions, with the trivial lower bound zero where
the right sides are nonpositive.

On the larger comparison rectangle $\rho\in[1/16,6/25]$ with the
same $q,\theta$ ranges, one has $0<d\le16/35$,
$r\ge1$, $s\ge27/35$, $A_0\ge2$, and $6/5\le B_0\le2$.
These bound the state and logarithm denominators away from zero.
In the cap rows $u_C+1=2A_0(A_0+Y_0)/(4q+A_0^2)>0$ and
$u_D-1=2(B_0Y_0-4q)/(4q+B_0^2)<0$, hence $\beta_C>0>\beta_D$.
In the affine matrix solve, put $\delta=A_0-B_0=Q_0d>0$.
The off-diagonal entry is
$-8q\delta(A_0+B_0+2Y_0)/[(4q+A_0^2)(4q+B_0^2)]<0$.
The other denominator factor has positive denominator
$\delta Q_0^2(A_0+B_0+2Y_0)$ and numerator
\[
 2Y_0\bigl[\delta Y_0^2+Y_0\{B_0(B_0+2\delta)+4q\}
                  +A_0^2B_0+4qA_0\bigr]>0 .
\]
Thus these formulas extend through all the geometric
mean-value comparison segments.

For a fixed geometric box, enclose coefficients from above in the
six affine necessary cap, switching and slope rows. Because the
means are nonnegative, these rows define a single polygon containing
every actual mean pair in that box. Intersections of every pair of
nonparallel supporting lines are outward enclosed. A candidate
intersection is discarded only by a strict violation; exact parallelism
is recognized by exact mantissa products, while ambiguous determinants
force subdivision. This includes zero-dimensional and line-segment
polygons through their supporting-line intersections.

Let $h_L,h_R$ be the maximal coordinates of this fixed containing
polygon and put $\alpha=3/(3+h_L)$, $\beta=3/(3+h_R)$.
At the actual point, positivity of the actual slopes gives
\[
G\ge W:=c+\alpha\lambda_L+\beta\lambda_R.
\]
The right side is affine in the two means. The determinant concavity
follows from the full cap and transport rows. In the whole mode-$0$
chart let $(u_-,v_-)$ and $(u_+,v_+)$ be the central entrance and
exit states, and let $K_c,\ell,m$ be the central invariant, outgoing
slope and outgoing lower-left matrix entry.
The outgoing cap equation has
\[
 m=D_cK_c-Z_c,\qquad D_c=v_+^{-2}+2\beta_+\Sigma>0,
\]
where $\Sigma>0$ is the central shear,
$\beta_+=[(1+u_+)/v_+^2+(1+u_+)^{-1}]/2>0$, and
\[
 Z_c=\frac{3+u_+^2}{2v_+^2}
       +\beta_+[2c-(u_+-u_-)]+\frac12.
\]
Thus
\[
 \Gamma=\frac{Z_c^2}{4D_c}
       -D_c\left(K_c-\frac{Z_c}{2D_c}\right)^2-\ell^2.
\]
Let $d,e$ be the preceding and following own reciprocals in this
whole chart, $Q=(y+d)(y+e)/(1+y^2)>1$ and
$b_0=(1-ye)/(y+e)>0$. Put
\[
 t_L=d(y+d)/2+\Sigma,\quad t_R=(1-u_+)/v_+^2,\quad
 \mathfrak h_L=(1+d^2)/2,\quad \mathfrak h_R=(1+b_0^2)/2.
\]
The returned-invariant rows are
\[
 \mathfrak h_i^2K_i=(t_i^2-D_c)K_c-2t_i\ell+Z_c
 \qquad(i=L,R).
\]
They have a nonsingular coefficient matrix, since
\[
 t_L-t_R=\frac{(1+y^2)^2(Q^2-1)}{2(y+e)^2}>0,
 \qquad (t_L-t_R)(t_Lt_R+D_c)>0.
\]
Thus the two means map affinely and nonsingularly to $(K_c,\ell)$.
With $J_c$ denoting the Jacobian, the determinant Hessian is
$-2J_c^{\mathsf T}\operatorname{diag}(D_c,1)J_c<0$ at fixed
geometry and resulting actual cost. The central invariant $K_c$
is distinct from the returned factors $(3+\mu_i)/2$. Therefore both $\Gamma$ and
$W+8(\Gamma-\overline H)$ can be bounded below on the whole
polygon by their vertex minima, where
$\overline H=1470471/2000000>H_*$. The polygon, its weights and
vertex enclosures are kept fixed in the geometric mean-value
estimate. No concavity claim for the original rational $G$ is used.
Comparison-vertex slopes need not be positive; positivity is used
only at the actual point to obtain $G\ge W$.

\begin{lemma}[Complete weighted-center certificate]\label{new:r17center}
Under the actual complete-source hypotheses and the cost window
\eqref{new:r17window},
\[
G>343/1000.
\]
\end{lemma}
\begin{proof}[Computer-assisted proof]
The actual determinant is less than $H_*$ if it is nonpositive,
since $H_*>0$. If it is positive, either actual returned factor
$K_i=(3+\mu_i)/2>0$ gives the positive trace required by
Corollary~\ref{cal:one-point} on the common lift and hence $\Gamma<H_*$.
On a box, failure of a necessary geometric or cost row or an
empty containing polygon excludes the actual point. A vertex lower
bound $\min\Gamma>\overline H$ also excludes it by concavity.
Otherwise either polygon-wide vertex bound
$\min W>343/1000$ or $\min\{W+8(\Gamma-\overline H)\}>343/1000$
suffices.
The first expression is affine in the means and the second is
concave. At the actual point the latter expression is less than
$W\le G$, because $\Gamma<H_*<\overline H$.

The binary80 interval program covers the entire root
\eqref{new:r17root} by the six closed $q$ slabs
$[5/2,7/2]$, $[7/2,5]$, $[5,8]$, $[8,16]$, $[16,32]$ and
$[32,243/4]$. It verifies the preceding
necessary rows and the fixed-polygon vertex tests. The whole-root
analytic guards justify all interval denominators and analytic
extensions. Arithmetic uses $64$ significand bits with directed
\texttt{nextafterl} endpoint enlargement; exact products of the
$64$-bit mantissas fit in unsigned $128$-bit integers for the
parallel-line test. The centered logarithm kernel
$\operatorname{atanh}(z)/z$ uses terms through index $24$ with
value remainder $|z|^{50}/[51(1-z^2)]$ and derivative remainder
$|z|^{49}/(1-z^2)$. Unsafe intervals become unresolved or wide;
they cannot certify exclusion.

The \path{local-r17} task checks the exact root union, all
closed child endpoints, the vertex tests and the strict positive
reserves. The finite covers of all six slabs have no unresolved box.
\end{proof}
The actual opposite centers $4$ and $10$ both have the required
complete sources and cost windows. Equations~\eqref{new:r17regroup},
\eqref{new:r17outer} and Lemma~\ref{new:r17center} yield
\[
C>2\frac{14}{25}+2\frac{343}{1000}
 =\frac{903}{500}=C_*+\frac{117}{100000}.
\]
This contradicts the original budget and excludes R17. The full
maximizing-PMP and source hypotheses remain essential.

\space
\section{Exclusion of the two twelve-arc classes}
\label{late:r10:section}

Throughout this section the geometric hypotheses, full endpoint conditions,
and switching-type convention are those of the original problem. Put
\[
 C_*={180483\over100000},\qquad
 \delta_0={29529\over2000000},\qquad
 E_*={4154032208439\over8759490110000}.
\]
Physical arc indices are taken modulo $12$, whereas source intervals are
read on the actual twisted lift. For $R10$ that lift has
\[
 (j_0,\ldots,j_{11})=021021210120,
 \qquad j_{i+12}=j_i+1\pmod3.
\]
Write $c_i>0$ for the cost of its complete maximal arc $i$ and
$C=\sum_{i=0}^{11}c_i$. A surviving minimizing owner would satisfy $C<C_*$.
All phase changes and reflections below are the full transformations of
Proposition~\ref{fd:lawful-symmetry}; they transport the state, frame,
matrix adjoint, the state covector, and objective multiplier.

\subsection{Source placement and the conditional interfaces}

Every coupled calculation below is made on a genuine containing period,
with its original state and frame endpoints and real winding.
Lemma~\ref{fd:reanchor} supplies the reanchoring. The complete capped
four-core $[3,8]=021210$, transformed by phase $+1$ to $102021$, supplies
normality and the automatic calibration of
Corollary~\ref{loc:four-calibration}. Joint normalization of that full lift,
followed by the actual apex and future-sheet conclusions of that corollary,
gives
\begin{equation}\label{late:r10:calibration}
 0<\Gamma<H_*<{1470471\over2000000},\qquad
 \delta={3\over4}-\Gamma>\delta_0.
\end{equation}
In particular, no state covector is extended periodically.

\begin{lemma}[Transferred own-pair interface]\label{late:r10:one-sided}
On the complete source $[-1,9]=20210212101$, if $c_2\le13/20$, the
actual $E$ own pair $\{5,6\}$ has both chronological density ceilings
$0<x_L,x_R<47/250$ and both cap gaps
\begin{equation}\label{late:r10:own-gaps}
 e_j\ge0,\qquad
 \delta\le {3(1-3x_j)\over4}e_j+e_j^2,
 \qquad e_j>{47\over2500}+{x_j\over20}.
\end{equation}
Its two actual zero-to-zero mean leaves and their common invariant are
unchanged. In particular $c_5+c_6>283/1000$.
\end{lemma}
\begin{proof}
The eleven-source is identical to the complete source used in the
one-sided construction in the proof of Lemma~\ref{late:r18:prefix}.
Here retain it in the genuine twelve-arc period $[-1,10]$, together
with the capped calibration source $[3,8]$. The prefix $[-1,4]$ has
centre $2$, its actual outgoing arc is $3$, the shared incoming arc
of $E$ is $4$, the true own leaves are $5,6$, and the right return is
$7,8,9$. Thus that construction applies without any change of its
local hypotheses. In particular its cap equation retains the
nonnegative initial $\chi$, its one-sided integral yields the missing
midpoint gap, and its actual full matrix conjugacy gives
\eqref{late:r10:own-gaps}. The one-sided invariant on arc $3$ is not
identified with the own invariant of $E$.
Proposition~\ref{six:prop:own}, with these recovered interfaces, gives
the stated own-pair cost. Only this local construction is used, not
any exclusion of another switching type.
\end{proof}

\begin{lemma}[Circle transfer]\label{late:r10:circle-transfer}
For $R10$, put $T=c_{11}+c_0+c_1$. Then $T>61/200$ always.
If $c_3\le13/20$, then $T>2/5$. More generally the latter conclusion
holds whenever the actual prefix read on $(6,5,4,3,2,1)$ under
$j\mapsto1-j$ has the strict exit sign of Theorem~\ref{hc:prefix650}.
\end{lemma}
\begin{proof}
The complete source $[-3,2]=012021$ gives the negative entrance and
entrance circle for the capped $202$ tail on $[-1,1]$, whose two caps
are retained in $[-2,2]=12021$; see
Lemma~\ref{cmp:forward-prefix-forcing}. Theorem~\ref{loc:ST305} gives
the unconditional payment.
The stated reverse prefix is literally $202102$, with centre $3$.
Its exit sign, together with \eqref{hc:C1-identity}, gives the one-sided
leaf argument in the proof of Lemma~\ref{late:r18:prefix}:
at the actual junction $1|2$ its reflected horizontal coordinate is
strictly positive. In the original tail chart, at $w=b+iv$, the whole
transformation is
\[
 F(w)={-\overline w-3\over1-\overline w}.
\]
The inequality $\operatorname{Re}F(w)>0$ is exactly
$v^2>(1-b)(3+b)$, the terminal circle. This is intrinsic state data.
For the tail comparison choose the genuine period $[-3,8]$, of exactly
$12$ complete arcs. It retains the whole coupled source $[-3,6]$ and
the capped calibration source $[3,8]$. Thus both circles, both actual
caps and the calibration are available on one full lift. If the exit
sign was established on $[-1,10]$, only its state-circle consequence
is carried to this period; no external normalized costate is imported.
Proposition~\ref{six:prop:tail}, in the circle-only form proved in
Lemma~\ref{late:r18:prefix}, now gives $T>2/5$.
For $c_3\le13/20$ the required exit sign follows from
Theorem~\ref{hc:prefix650}, including equality at the cost threshold.
\end{proof}

\subsection{An exclusive payment at the right tail}

\begin{lemma}[Two exclusive arcs]\label{late:r10:exclusive155}
Every surviving $R10$ owner satisfies $c_8+c_9>31/200$.
\end{lemma}
\begin{proof}
Retain the proper nine-source $[3,11]=021210120$ on the genuine
period $[3,14]$. The one-return source $[3,9]$, in phase $+1$, is
$1020212$; its $0202$ core is $\{4,5,6,7\}$, with own means $5,6$
and true right returned mean $8$. The capped source $[6,10]$, in
phase $+2$, is $10201$; its $020$ core is $\{7,8,9\}$ and its
middle mean is arc $8$. Their shared physical arc is $7$, the last
mode-$2$ arc of the first chart and the first mode-$0$ arc of the
second. Their adjacent actual mean leaves are $6$ and $8$.
The capped four-core supplies \eqref{late:r10:calibration} on this
one full lift.

In the first chart write $K_L=(3+h_L^2)/2$, and let $r>h_L$ be
its positive outer root. In the second chart put
\[
 K=(3+\mu)/2,\qquad -s<-\sqrt\mu,
 \qquad |z|\le\sqrt\mu,
\]
where $-s,z$ are the two internal roots. Let $\beta>z$ be the
terminal cap coordinate. Lemma~\ref{loc:shared-A-matrix} and
\eqref{pf:positive-cap}, \eqref{pf:gap}, applied to the actual means
$6,8$, give, for the shared density $x$ and constant $\lambda$,
\begin{equation}\label{late:r10:shared-defect}
 \begin{gathered}
 \Gamma={K_LK(1+x)^2\over4}-\lambda^2,\quad
 0<x<1/3,\quad e=\lambda+3x/2\ge0,\\
 \delta\le {3(1-3x)\over4}e+e^2\le3e/4+e^2.
 \end{gathered}
\end{equation}
The local scope of these identities is precisely one true adjacent
mean on each side and their two midpoint gaps, as in
Section~\ref{pf:F09}; a right four-core is not required. Since
\[
 {29\over2000}-{3\over4}{47\over2500}
                     -\left({47\over2500}\right)^2
 ={291\over6250000}>0,
\]
\eqref{late:r10:calibration} and monotonicity for $e\ge0$ imply
\begin{equation}\label{late:r10:defect-floor}
 e>47/2500.
\end{equation}
No high-density ceiling or second determinant gap is needed here.

Descending arcs $(11,10,9,8,7,6,5)$ under $j\mapsto-j$ form
$0120212$. The entrance conclusion of
Lemma~\ref{cmp:forward-prefix-forcing} for this actual source becomes,
in the forward $020$ chart,
\begin{equation}\label{late:r10:terminal-circle}
 \beta>0,\qquad q(\beta)<{1\over(1-\beta)(3+\beta)},
 \qquad q=v^{-2}.
\end{equation}
Indeed the reverse chart is exactly the whole binary reflection of
the phase-$+2$ chart, since $2-(j+2)=-j$ modulo $3$.
This uses an actual entrance circle in reverse chronology; it does
not insert a terminal returned mean into the forward $020$ core.

Write $P=c_8+c_9$, $w=s+\beta$, and let $T>0$ be the shear over
precisely these two arcs, with $ds=dt/v=q\,du$.
The mean and physical-cost identities give
\[
 0<\mu<s^2<w^2,\qquad
 P={w\over2}+{1\over2}\int(3+u^2)\,ds
   >w/2+3T/2,\qquad P>2T.
\]
Suppose $P\le31/200$, and put $V=31/100$. Then
\begin{equation}\label{late:r10:exclusive-budget}
 0<w<V,\quad \mu<w^2<V^2<1/10,\quad
 T<\min\{31/400,(V-w)/3\},\quad K<30961/20000.
\end{equation}
In particular $KT<959791/8000000<1/2$.

We first determine the sign of the actual terminal matrix entry $L$.
At the last mode-$0$ cap put $d=1+\beta$ and abbreviate terminal
values by $q,\chi,f,L$. The actual maximizing equations give
\[
 \chi\ge0,\quad g=m+1/2-f=0,\quad
 f={(\mu-\beta^2)q\over2}+{dq+1/d\over2}\chi,
\]
\[
 g_s=-2L-f_s\le0,\qquad
 f_s=-\beta-{f\over qd}+{(\mu-\beta^2)qd\over2}.
\]
These are the mode-$0$ cap identities used in
Theorem~\ref{loc:DK610}; the inequality follows from $g\ge0$ on the
inward side. By \eqref{late:r10:terminal-circle} and
$\beta<V<\sqrt2-1$, one has $qd^2<1$. Hence
\[
 f_s\le-\beta+{\mu-\beta^2\over2}(qd-1/d)<0.
\]
For $\mu\ge\beta^2$ the second term is nonpositive; for
$\mu<\beta^2$ it is strictly less than $\beta^2/(2d)<\beta$.
It follows that $L\ge-f_s/2>0$, before any square root is taken.

The real root $z$ has $\chi(z)=0$, so integration of
$\chi_s=\mu-u^2$ on the last arc gives $0\le\chi\le\mu T$.
Retain this possibly nonzero terminal value. The cap determinant is
\begin{equation}\label{late:r10:cap-determinant}
 L^2=\delta+\mu/4-Kf.
\end{equation}
Put $D=(1-\beta)(3+\beta)$ and $B=\mu-\beta^2+d\mu T$.
If $B\le0$, the formula for $f$ gives $f\le\mu T/(2d)$,
and therefore
\[
 L^2\ge\delta+\mu\{1/4-KT/(2d)\}>\delta>\delta_0.
\]
If $B>0$, the same terminal circle gives
$f\le qB/2+\mu T/(2d)<B/(2D)+\mu T/(2d)$.
Substitution into \eqref{late:r10:cap-determinant} gives exactly
\[
 L^2>\delta+{3\beta^2-2\mu\beta-\mu^2\over4D}
             -{(3+\mu)\mu T\over dD}.
\]
Here $d^2+D=4$, so the last term includes the full contribution of
$\chi$. Complete the square as
$3\beta^2-2\mu\beta-\mu^2=3(\beta-\mu/3)^2-4\mu^2/3$.
Moreover
\[
 D>{22839\over10000}>{9\over4},\qquad
 dD=3+\beta-3\beta^2-\beta^3
 \ge\min\{3,3+V-3V^2-V^3\}
 ={2991909\over1000000}>{299\over100}.
\]
The minimum assertion follows from concavity on $[0,V]$.
Thus $(3+\mu)/(dD)<310/299<21/20$. Both signs of $B$ yield
\begin{equation}\label{late:r10:L-bound}
 L^2>\delta_0-{4w^4\over27}-{21\over20}w^2T.
\end{equation}
By \eqref{late:r10:exclusive-budget} and
$w^2(V-w)\le4V^3/27$,
\[
 {4w^4\over27}+{21\over20}w^2T
 <{4V^4\over27}+{7V^3\over135}
 ={7864824\over2700000000}.
\]
The exact reserve
\[
 \delta_0-(21/200)^2-{7864824\over2700000000}
 ={2231826\over2700000000}>0
\]
therefore proves $L>21/200$.

Finally, the full shared transform \eqref{loc:shared-A-transform}
and the actual equation $L_s=K$ give
\[
 L(-s)=-\lambda-{Kx\over1-s},\qquad
 KT=L-L(-s)=L+e+{x(\mu+3s)\over2(1-s)}>L+e>{619\over5000}.
\]
This contradicts \eqref{late:r10:exclusive-budget}, since
\[
 {619\over5000}-{959791\over8000000}
 ={30609\over8000000}>0.
\]
Thus $P>31/200$. Its charge excludes shared arc $7$ and both caps;
every coupled identity used only the common nine-source.
\end{proof}

\subsection{A height floor for the actual own pair}

\begin{lemma}[The refined height]\label{late:r10:height}
If $c_2\le13/20$, the actual height $h_E$ of the own mean on
arcs $5,6$ satisfies $h_E>1/12$.
\end{lemma}
\begin{proof}
Use the same full lift and densities supplied by
Lemma~\ref{late:r10:one-sided}. For either density $x$ and defect $e$,
put $F_x(e)=3(1-3x)e/4+e^2$. It increases for $e\ge0$ and
$0<x<1/3$, and
\[
 \delta_0-F_x(19/1000+x/20)
 ={307\over2000000}+{67x\over20000}+{11x^2\over100}>0.
\]
Thus $e>19/1000+x/20$.

Suppose $0<h=h_E\le1/12$. Write the actual roots as $-a,z,b$,
and put
\[
 R={8h\over4-3h},\quad k={1-R\over1+R},\quad
 X=\max(x_L,x_R),\quad S=x_L+x_R,\quad K={3+h^2\over2}.
\]
Lemma~\ref{loc:leaf-correlations} gives $h<a,b<R$, $|z|\le h$.
For $M(t)=(1+t)/(1-t)$ the exact density ratio is
\[
 {x_R\over x_L}=\left({M(a)M(z)\over M(b)}\right)^2.
\]
Since $z\le h<b$ this ratio is less than $M(R)^2$; using
$-z\le h<a$ proves the reciprocal bound. Consequently
$S>(1+k^2)X$.

Let $T$ be the actual own shear and
$D=x_La/(1-a)+x_Rb/(1-b)$. The exact transport and
Lemma~\ref{pf:own-shear-correlation} give
\[
 N:=e_L+e_R=K(T-D)-h^2S/2
 \le K(1-k)T-h^2S/2
 \le A(h)X-h^2S/2,
\]
where the last inequality uses
$T\le2q(z)R/(1+R)$, $q(z)\le X/k^2$, and
\[
 A(h)={4KR^2\over(1-R)^2}
 ={128h^2(3+h^2)\over(4-11h)^2}.
\]
The refined defects also give $N>19/500+S/20$, so
\begin{equation}\label{late:r10:height-contradiction}
 {19\over500}<A(h)X-(1/20+h^2/2)S<F(h)X,
 \quad F(h)=A(h)-(1/20+h^2/2)(1+k^2).
\end{equation}
This function increases on $[0,1/12]$. Indeed
$A(h)=h^2a(h)$, where $a(h)$ is positive, increasing, and at least
$24$, so $A'(h)\ge48h$. Also $\sigma(h)=1+k^2$ decreases and
$\sigma\le2$. Thus
$F'=A'-h\sigma-(1/20+h^2/2)\sigma'\ge46h\ge0$.
At $h_0=1/12$ one has
\[
 R_0={8\over45},\quad A(h_0)={3464\over12321},\quad
 \sigma(h_0)={4178\over2809},\quad
 F(h_0)={558222323\over2768775120},
\]
\[
 {19\over500}-{47\over250}F(h_0)
 ={133828918\over1384387560000}>0.
\]
Since $X<47/250$, this contradicts
\eqref{late:r10:height-contradiction}, including its endpoint $h=h_0$.
\end{proof}

\subsection{The strengthened central-prefix estimate}

\begin{lemma}[Exit under the refined ingress]\label{late:r10:enhanced-prefix}
Suppose an actual contact-qualified $202102$ prefix has central cost
$c<69083/100000$ and ingress $u_L>1/12$. Then its actual exit satisfies
\[
 u_R<0,\qquad q_R\ge0,\qquad q_R'>kq_R.
\]
Here $q$ is the raw switching response of
\eqref{hc:response}, not the reciprocal state density used above.
\end{lemma}
\begin{proof}
Use the actual full central response in \eqref{hc:parameters}--
\eqref{hc:affine-rows}, with $A>0$, $E>1$, $z=\log E$ and
$\tau^2=1+k^2$. The strict ingress, both contact jets and physical
cost imply
\begin{equation}\label{late:r10:prefix-domain}
 AE<2,\quad A\tau^2>2,\quad
 A>A_t={2(Ez-E+1)\over(E-1)^2},\quad
 a:=A(E-1)<2c<7/5.
\end{equation}
Lemma~\ref{hc:crossing} already gives $u_R<0$, without a cost
restriction, and the actual outgoing tie gives $q_R\ge0$.
Since
\[
 u_L={4-\tau^2A^2\over4-4A+\tau^2A^2}>1/12,
 \qquad \tau^2<{4(11+A)\over13A^2},
\]
the monotonicity of \eqref{hc:physical-cost} in $\tau^2$ gives
\begin{equation}\label{late:r10:C12}
 c>C_{12}(A,E):=
 {2A((11+A)(E-1)+13(1-E^{-1}))-13A^2\log E
       \over44+4A-13A^2}.
\end{equation}
The denominator is $(2-A)(13A+22)>0$.

First $C_{12}$ increases strictly in $A$ on $0<A<2/E$.
Put $e=E-1$, $a_1=22e+26e/E$, and $a_2=2e-13z$.
Its $A$-derivative has positive denominator and numerator
\[
 44a_1+88a_2A+(13a_1+4a_2)A^2.
\]
If $a_2\ge0$ this is positive. Otherwise $z<e$ gives
$a_2>-11e$, whence $13a_1+4a_2>0$, and, using $A<2/E$,
the first two terms exceed $968e-792e/E>0$.

\smallskip
\noindent\emph{The whole prefix has $E<47/20$.}
Write $D(A)=44+4A-13A^2$ and $W(E)=Ez/(E-1)-1$.
At the contact comparator $A=A_t$ direct substitution gives
\[
 C_{12}(A_t,E)=W(E)\{1+2R(A_t)/E\},\qquad
 R(A)={13(A^2-2A+4)\over2D(A)}.
\]
For $0<A<2$ one has $R(A)>491/1000$, since
\[
 1000D(A)(R(A)-491/1000)
 =12883A^2-14964A+4396,
\]
whose leading coefficient is positive and whose negative-discriminant
reserve is $4\cdot12883\cdot4396-14964^2=2613376>0$.
With $b=491/500$, the lower envelope $W(E)(1+b/E)$ is strictly
increasing for $E>1$. Its derivative, multiplied by
$E^2(E-1)^2$, is
\[
 (E-1)\{E^2+b(2E-1)\}-(1+b)E^2\log E.
\]
The trapezoid bound $\log E<(E-1)(E+1)/(2E)$ makes this greater
than $(E-1)^2\{(1-b)E+2b\}/2>0$.

For clarity the needed logarithm comparisons are entirely rational.
For $t=(E-1)/(E+1)$ let
\[
 L_N(E)=2\sum_{m=0}^{N-1}{t^{2m+1}\over2m+1},\qquad
 U_N(E)=L_N(E)+{2t^{2N+1}\over(2N+1)(1-t^2)}.
\]
The positive series for $\log E$ gives $L_N<\log E<U_N$.
At $E_0=47/20$, $t=27/67$ and
\[
 t>{402985\over1000000},\quad t^2>{1623969\over10000000},\quad
 t^3>{654435\over10000000},
\]
\[
 t^5>{106278\over10000000},\quad
 t^7>{17259\over10000000},\quad
 t^9>{2802\over10000000}.
\]
Multiplication of the first two bounds successively proves the
remaining power bounds. Substitution gives
$L_5(E_0)>854405/1000000>534/625$. Consequently
\[
 W(E_0)>{2741\over5625},\qquad
 W(E_0)(1+b/E_0)>{4566506\over6609375}
 >{69083\over100000},
\]
with last difference $55146875/660937500000>0$.
Since $A>A_t$, the two monotonicities prove
\begin{equation}\label{late:r10:E-bound}
 E<47/20.
\end{equation}
In particular $E<31/10$, so the original
Lemma~\ref{hc:short-factor-force} gives $Q=f_-+f_+>0$.
No extension of its domain is required. This $Q$ is an affine
response coefficient, unrelated to the cost $c_8+c_9$.

\smallskip
\noindent\emph{The remaining affine coefficient is positive.}
Use the notation of \eqref{hc:affine-coefficients}, and put
$G=Q-\rho f_-$. Only $\rho<0$ needs further work, since $f_-<0$.
In this sector the argument preceding \eqref{hc:G-cross} gives
\[
 H=-A(E-1)(E^2-E+1)+2E(E+1)z-4E^2+6E-2<0.
\]
The exact coefficient of $\tau^2$ in $G$ is $AH/(2-A)$.
The strengthened upper bound on $\tau^2$, compared with $4/A^2$,
therefore gives
\begin{equation}\label{late:r10:G12}
 G>{2\over E}(2J/A-W_0)-{4H\over13A},
\end{equation}
where
\[
 J=(E+1)(E^2+1)z-(E-1)(2E^2-E+3),\qquad
 W_0=(E-1)^3(E+1).
\]
Multiplying the right side of \eqref{late:r10:G12} by $13EA/4$
gives exactly $J_{12}-aP_{11}(E)/2$, where
\begin{align*}
 J_{12}&=(E+1)(11E^2+13)z-(E-1)(22E^2-11E+39),\\
 P_{11}(E)&=11E^3-11E^2-15E+13.
\end{align*}

For $1<E\le3/2$ the weaker original comparison
\eqref{hc:G-cross} already gives $G>0$. In fact $A<2/E$, and
\[
 J>{(E-1)^2(2E^3+E^2+8E+1)\over3(E+1)^2}
\]
by the lower logarithm $L_2$. Hence $EJ>W_0$ follows from
\[
 P_0(E)=-E^4-5E^3+8E^2+7E+3>0.
\]
This polynomial is strictly concave on $[1,3/2]$, with endpoint
values $12$ and $153/16$.

Assume now $G\le0$. Then $3/2<E<47/20$.
On this interval $P_{11}>0$: its value at $3/2$ is $23/8$ and
its derivative is positive and increasing. Put
$t=(E-1)/(E+1)<27/67<81/200$. Substitution of $L_3$ gives the
exact identity
\begin{align*}
 120(1-t)^3\{J_{12}(L_3)-(39/80)P_{11}(E)\}=P(t),\\
 P(t)=117-363t+507t^2-1569t^3-640t^4
                 +6144t^5-384t^6+2304t^7.
\end{align*}
We prove $P>0$ on the whole interval $[0,81/200]$.
Its derivative is
\[
 P'=-363+1014t-4707t^2-2560t^3
                      +30720t^4-2304t^5+16128t^6.
\]
For $0\le t\le1/3$, the positive higher terms are at most
$(10240/3+1792/9)t^2$; the remaining quadratic coefficient is
negative, so $P'\le-363+1014/3=-25$.
For $1/3\le t\le3/8$ use
\[
 30720t^4\le4320t^2,\qquad
 16128t^6\le(5103/16)t^2,\qquad
 -2560t^3\le-(2560/3)t^2.
\]
Dropping the other negative term gives
\[
 P'\le-363+1014t-(44227/48)t^2
 \le69/4-44227/432<0.
\]
For $3/8\le t\le81/200$ use
\[
 \begin{gathered}
 30720t^4\le(629856/125)t^2,\qquad 16128t^6<434t^2,\\
 -2560t^3\le-960t^2,\qquad -2304t^5\le-(243/2)t^2.
 \end{gathered}
\]
Their total quadratic coefficient is less than $-315$.
The resulting quadratic $-363+1014t-315t^2$ increases on this
interval and at $81/200$ equals $-159915/40000<0$.
Thus $P'<0$ everywhere. At $t=81/200$ the exact rational bounds
\[
 507t^2>83.16,\quad6144t^5>66.94,\quad2304t^7>4.11,
\]
\[
 363t=147.015,\quad1569t^3<104.232,\quad
 640t^4<17.219,\quad384t^6<1.695
\]
prove $P(81/200)>1.049>0$; all terminating decimals here are
exact rationals. This proves $P>0$ throughout.
Since $z>L_3$ and its coefficient in $J_{12}$ is positive,
$J_{12}>(39/80)P_{11}$. From $G\le0$ and the strict lower
comparison \eqref{late:r10:G12} it follows that
\begin{equation}\label{late:r10:a-bound}
 a>{2J_{12}\over P_{11}}>39/40.
\end{equation}
The denominator was shown positive before this division.

\smallskip
\noindent\emph{The forced value of $a$ contradicts the cost.}
At fixed $a=A(E-1)$ put $e=E-1$ and
\[
 D_a=44e^2+4ae-13a^2,\qquad
 B_a=a^2-2a\log E+4e^2/E.
\]
Then
\[
 C_{12}(a/e,E)={a\over2}+{13a\over2}{B_a\over D_a}.
\]
For $39/40\le a<7/5$ and $E\ge3/2$, this decreases strictly
in $E$ wherever $D_a>0$. Here is the derivative comparison.
For $F=D_a'B_a-D_aB_a'$, the logarithm trapezoid bound gives
\begin{align*}
 {E^2F\over4}\ge{}&-{11a^3\over2}
 +(48a^2-9a^3/2)e+(a^3+56a^2-26a)e^2\\
 &+(21a^2-44a)e^3+(44-22a)e^4.
\end{align*}
The first two terms are at least $(263/20)a^2>0$, using
$e\ge1/2$ and $a\le7/5$. On the same $a$ interval,
\[
 a^3+56a^2-26a>28,\quad
 a(44-21a)\le484/21<24,\quad 44-22a\ge66/5.
\]
The remaining bracket therefore exceeds
$28-24e+(66/5)e^2>0$, whose negative-discriminant reserve is
$4512/5$. Hence $F>0$, as asserted.

First decrease $a$ to $39/40$ at fixed $E$, using the already
proved $A$ monotonicity. This preserves $0<A<2/E$. Then increase
$E$ to $47/20$ with $a=39/40$ fixed; both $A=a/(E-1)$ and
$AE$ decrease, so the denominator remains positive. Thus
\[
 c>C_{12}(A,E)
  >C_{12}((39/40)/(E-1),E)
  >C_{12}(13/18,47/20).
\]
For the upper logarithm at $E_0$, use $t=27/67<403/1000$,
$t^3<0.065451$, $t^5<0.010631$, $t^7<0.001727$ and
$1-t^2>0.837$. Then
\[
 U_3(E_0)<0.806+0.043634+0.0042524+0.0006<171/200.
\]
Since the cost decreases in its logarithm argument, direct rational
substitution gives
\[
 C_{12}(13/18,47/20)
 >{39\over80}+{454753143\over2198754000}
 >{69083\over100000},
\]
with last difference
$768049218000/219875400000000>0$. This contradicts the assumed
cost. Therefore $G>0$ whenever $\rho<0$.

Finally retain the original affine exit row \eqref{hc:exit-affine}.
For completeness set $X_-=-k^4q_L>0$, $Y_+=k^4q_R\ge0$, with
$m_0,n,a_0>0$ as in \eqref{hc:affine-coefficients} and
\eqref{hc:strict-affine-ingress}. Then
\[
 k^3(q_R'-kq_R)=f_++(\rho-1)m_0Y_++nX_-.
\]
For $\rho\ge1$, $f_+=Q-f_->0$ proves positivity. For
$\rho<1$, put $\varepsilon=-f_--a_0X_--m_0Y_+>0$.
The same expression is exactly
\[
 G+\{n+(1-\rho)a_0\}X_-+(1-\rho)\varepsilon>0.
\]
For $0\le\rho<1$ use $G=Q-\rho f_->0$; for $\rho<0$ use the
preceding argument. This proves the strict exit assertion in every
case, including $\rho=0$.
\end{proof}
\space
\section{A bilateral-pair exclusion of R11 and R12}
\label{late:r11:section}
Throughout this section the trajectory is an original global minimizer,
with positive maximal arcs and the original group and state endpoints
and real winding. Its cost satisfies $C<\Cstar$. The native word of R11
is $021021202120$, with $j_{i+12}=j_i+1$ modulo three. An interval
$[i,k]$ denotes the complete arcs $i,\ldots,k$ on this unfolded word.
Write
\begin{equation}\label{late:r11:charges}
 P=c_8+c_9+c_{10},\qquad Q=c_4+c_5+c_6,\qquad
 R=c_{11}+c_0+c_1.
\end{equation}
These triples are disjoint, and
$C=c_2+c_3+c_7+P+Q+R$.

\subsection{Physical sources and calibration}
\begin{lemma}\label{late:r11:qualification}
The genuine switch-cut period $[3,14]$ admits a jointly normalized full lift
with
\[
 0<\Gamma<H_*<\frac{1470471}{2000000}.
\]
No improved numerical bound on the remaining three-arc charge is needed
for this qualification.
\end{lemma}
\begin{proof}
The source $[1,5]=21021$ becomes $10210$ under whole phase $+2$.
The independently qualified weak-five theorem~\ref{wk:five600} gives
$V:=c_2+c_3+c_4\ge3/5$. Put $W=P+R$. Positivity of
$c_5+c_6+c_7$ alone therefore gives
\[
 W<C-\frac35<\Cstar-\frac35<\frac54.
\]
The full source $[7,14]=02120102$ is retained in $[3,14]$. This interval
has exactly twelve complete arcs and is a fundamental period cut at a
switch. Reanchoring preserves the original endpoints and real winding.

The actual capped source $[7,11]=02120$, read backwards under
$j\mapsto1-j$, is $12021$. Its two roots and two caps establish
normality before the matrix and the state covector are jointly
normalized. The central forward-four $[9,12]=1201$ in $W$ supplies
both negative entrances, without harmonic calibration. Choose the
cheaper actual $W$ tail on this same selected lift. Its cost satisfies
$T\le W/2<5/8$. At its terminal point the signed strip
of Lemma~\ref{loc:signed-strip} gives
\[
 \frac32<K<2,\qquad m<-\frac12,\qquad
 0<L<K\sigma<2\sigma<T.
\]
Consequently
\[
 \Gamma=-L^2-Km>\frac34-T^2>\frac{23}{64}>0.
\]
The positive actual invariant $K$ supplies the future sign in
Corollary~\ref{cal:one-point}. The future-sheet and harmonic
results~\ref{cal:future} and~\ref{cal:harmonic} now apply to this
same genuine period. No costate is continued periodically across a seam.
\end{proof}

\begin{lemma}\label{late:r11:seven}
The physical charges in \eqref{late:r11:charges} satisfy
\[
 P>\frac25,\quad P+c_7>\frac{91}{200},\quad
 Q,R>\frac{61}{200}.
\]
\end{lemma}
\begin{proof}
On the full reflected-time chart $j\mapsto1-j$, the actual chain is
\[
 (10,9,8)\mid7\mid(6,5,4)=202\mid1\mid202.
\]
The first triple, of cost $P$, has caps 11 and 7 and return 6.
Its entrance circle is furnished by $[9,12]=1201$, and its terminal
circle by its return. The second triple, of cost $Q$, has caps 7 and 3
and its negative entrance and entrance circle from the seven-source
$(10,\ldots,4)$, by Lemma~\ref{ek:H7-circle}. The complete support
$[3,12]$ and the full $W$ source both lie in the period $[3,14]$.
Apply Lemma~\ref{late:r11:qualification} on that period.

Proposition~\ref{six:prop:tail} gives $P>2/5$.
The complete reversed source $(11,10,9,8,7,6)$ under $1-j$ is
$120212$. Its charged triple is $P$, its following complete mode-one
bridge is arc $7$, and its actual return is arc $6$. The bridge ends
at the negative entrance of $Q$, as proved by the retained H7 source.
Thus Theorem~\ref{ua:four}, on this same qualified interval, gives
$P+c_7>91/200$.
Theorem~\ref{loc:ST305} gives $Q>61/200$ and, on the other actual
capped $W$ triple, $R>61/200$.

\end{proof}

\subsection{The two cheap centre implications}
\begin{lemma}\label{late:r11:cheap-circles}
For the actual R11 trajectory,
\[
 c_2\le\frac{13}{20}\ \Longrightarrow\ Q>\frac25,
 \qquad
 c_3\le\frac{13}{20}\ \Longrightarrow\ R>\frac25.
\]
\end{lemma}
\begin{proof}
Suppose first $c_2\le13/20$. The actual prefix $[-1,4]=202102$
has centre 2. By \ref{hc:prefix650} and \ref{hc:C1-identity}, its
outgoing native mode-zero arc 3 starts at $z<0$ with
$\chi(z)\ge0$ and $\chi_t(z)>0$. Write $\mu_c=2K_c-3$ for its
binary invariant. Since $\chi_s=\mu_c-u^2$, we have
$\mu_c>z^2$. At its actual binary tie $3|4$, let the horizontal
coordinate be $r_0$. The switching function first increases and can
reach zero only after $u>\sqrt{\mu_c}$, so
\begin{equation}\label{late:r11:prefix-root}
 0<\sqrt{\mu_c}<r_0<1.
\end{equation}
The exact one-sided relation is
$\int_z^{r_0}(u^2-\mu_c)(1+u)^{-2}\,du\ge0$.
In particular $\mu_c$ has not been identified with a two-zero mean.

The reflected $Q$ triple $(6,5,4)$ already has its entrance circle
and caps 7 and 3. Its terminal junction is the same physical $3|4$.
The full state involution between its chart and the native chart is
\[
 \mathcal T(w)=\frac{-\overline w-3}{1-\overline w},\qquad
 \Re\mathcal T(b+iv)=
 \frac{b^2+2b-3+v^2}{(1-b)^2+v^2}=r_0>0.
\]
Thus $Q$ has its terminal circle as well. The union of this H7 source
and the prefix is $[-1,10]$, exactly one genuine twelve-arc period.
On this selected genuine period, the capped-three-core argument proves
normality before joint normalization. Its actual loss and determinant
satisfy $\Gamma=3/4-\mathcal L/4$. If $\mathcal L<59/1000$,
then $\Gamma>2941/4000>0$. The actual positive binary invariant
$K=(3+\mu)/2$ supplies the future sheet, and the harmonic ceiling
contradicts this determinant bound. Thus $\mathcal L\ge59/1000$
is established on this lift independently of the calibration on
$[3,14]$. The circle-only transfer in the proof of
Lemma~\ref{late:r18:prefix} now gives $Q>2/5$.

If $c_3\le13/20$, instead use descending native arcs $(6,\ldots,1)$
under full reflected time reversal $j\mapsto1-j$. The centre is 3;
the outgoing reversed native arc 2 gives a positive coordinate at
$1|2$. The same involution supplies the terminal circle for
$[-1,1]=202$, which represents $R$. Its capped context is
$[-2,2]=12021$, and the forward-four $[-3,0]=0120$ gives its entrance
circle. This full union $[-3,6]$ has ten arcs. On its own selected
normal lift the same loss-obstruction, future-sheet and circle-only
argument gives $R>2/5$.

\end{proof}
\space
\label{revision:appendix-end-6}
\end{proofappendix}
\section{Closure charts and finite covers for eight arcs}\label{app:eight-details}
\begin{proofappendix}
\section{A complete Gram chart}
All duration indices in this section are modulo eight and all pair
indices modulo four. Full closure has a stronger consequence than a
single scalar equation.

\begin{lemma}\label{eight:lem:gram}
Suppose $q_j\ge1$ and \eqref{eight:eq:signed-closure} holds. Set
$\alpha_i=q_{2i}^2+q_{2i+1}^2-1$ and
\[
 \Gram=\begin{pmatrix}
 q_0&0&0&-q_7\\q_1&q_2&0&0\\0&q_3&q_4&0\\0&0&q_5&q_6
 \end{pmatrix}.
\]
Then $\Gram^t\Gram-I$ is positive semidefinite of rank two. After
congruence by $\diag(\alpha_i^{-1/2})$, it has diagonal entries one
and consecutive off-diagonal entries $c,s,c,-s$, where
\begin{equation}\label{eight:eq:cs}
 c=\frac{q_1q_2}{\sqrt{\alpha_0\alpha_1}}
  =\frac{q_5q_6}{\sqrt{\alpha_2\alpha_3}}>0,\quad
 s=\frac{q_3q_4}{\sqrt{\alpha_1\alpha_2}}
  =\frac{q_7q_0}{\sqrt{\alpha_3\alpha_0}}>0,\quad c^2+s^2=1.
\end{equation}
\end{lemma}
\begin{proof}
Consider the antiperiodic recurrence
$q_{j-1}z_{j-1}+q_jz_{j+1}=z_j$, $z_{j+8}=-z_j$.
Its transfer on $(z_j,q_{j-1}z_{j-1})^t$ is
\[
 T(q)=\begin{pmatrix}1/q&-1/q\\q&0\end{pmatrix}
      =(S_0^{-1}\Nmat(q)S_0)^t,\qquad
 S_0=\begin{pmatrix}1&1\\1&0\end{pmatrix}.
\]
By full signed closure its antiperiodic solution space has dimension two.
Splitting odd and even coordinates gives the equations
$\Gram v=u$, $\Gram^tu=v$. Consequently
$\ker(\Gram^t\Gram-I)$ has dimension two. Its opposite-coordinate
principal block $\diag(\alpha_0,\alpha_2)$ is positive definite.
The Schur complement must vanish by rank two, proving positive
semidefiniteness. In the order $(0,2,1,3)$ the normalized matrix is
$\left(\begin{smallmatrix}I&U\\U^t&I\end{smallmatrix}\right)$;
its vanishing Schur complement says $U^tU=I$.
The signs of its four nonzero entries make $U$ exactly
$\left(\begin{smallmatrix}c&-s\\s&c\end{smallmatrix}\right)$.
This proves every assertion.
\end{proof}

\begin{theorem}\label{eight:thm:chart}
The full positive locus $q_j>1$ satisfying \eqref{eight:eq:signed-closure}
has a smooth bijective chart of differential rank five,
\[
                  (t,x_0,x_1,x_2,x_3)\in(0,1)^5.
\]
Its forward formulas are
\begin{gather}
 a_i=\begin{cases}t&i\text{ even},\\1-t&i\text{ odd},\end{cases}
 \quad u_i=1-a_ix_i,\quad v_i=\frac{a_i}{u_{i+1}},\notag\\
 D_i=1-x_i+(1-a_i)x_ix_{i+1},\quad
 z_i=\frac{a_iD_i}{u_{i+1}},\label{eight:eq:cube}\\
 q_{2i}^2=\frac{u_i}{z_i}
   =1+\frac{(1-a_i)(1-x_{i+1})}{a_iD_i},\qquad
 q_{2i+1}^2=\frac{v_i}{z_i}=\frac1{D_i}.\label{eight:eq:qchart}
\end{gather}
Every square root is positive. The inverse is
\begin{equation}\label{eight:eq:inverse}
 t=\frac{q_1^2q_2^2}{\alpha_0\alpha_1},\qquad
 x_i=\frac{1-q_{2i}^2/\alpha_i}{a_i}.
\end{equation}
\end{theorem}
\begin{proof}
On the cube, $u_i>1-a_i>0$ and
$D_i=1-x_iu_{i+1}\in(0,1)$. Also $z_i=u_i+v_i-1$ and
\[
 \alpha_i=z_i^{-1},\quad
 q_{2i}^2/\alpha_i=u_i,\quad q_{2i+1}^2/\alpha_i=v_i,\quad
 v_iu_{i+1}=a_i.
\]
These identities make the normalized Gram matrix the rank-two matrix
of Lemma~\ref{eight:lem:gram}, with $c=\sqrt t$ and $s=\sqrt{1-t}$.
The same parity splitting gives a two-dimensional space of antiperiodic
solutions of the recurrence. Since its unrestricted solution space also
has dimension two, its monodromy is \emph{exactly} $-I$.
The transpose conjugacy above now proves \eqref{eight:eq:signed-closure}.
This argument proves full closure, rather than merely its trace.

Conversely, Lemma~\ref{eight:lem:gram} gives $0<t=c^2<1$. Set
$u_i=q_{2i}^2/\alpha_i$, $v_i=q_{2i+1}^2/\alpha_i$.
Strict durations give $u_i,v_i<1$, and
$v_{i-1}u_i=1-a_i$ gives $u_i>1-a_i$. Thus the inverse in
\eqref{eight:eq:inverse} lies in the open cube. Substitution recovers
\eqref{eight:eq:cube}--\eqref{eight:eq:qchart} uniquely.
Both the squared forward map and this inverse are rational and smooth
near the image, with nonzero denominators. Differentiating their
composition gives a five-dimensional left inverse of the forward
Jacobian. Positive square roots preserve its rank.
\end{proof}

The same argument covers finite $q_j\ge1$ by allowing $0\le x_i\le1$
while retaining $0<t<1$ and $D_i>0$. On this domain,
\begin{equation}\label{eight:eq:boundary}
 q_{2i+1}=1\iff x_i=0,\qquad q_{2i}=1\iff x_{i+1}=1.
\end{equation}
A face $D_i=0$, namely $x_i=1$, $x_{i+1}=0$, is a duration pole,
not a finite zero-time arc. Our use of closed boxes will not impose
interior stationarity on the finite zero-time faces.

\section{The gap inequality and ellipticity}
The state return involves half durations. Its discriminant cannot be
inferred from the rank of the full Gram matrix alone.

\begin{lemma}[Four-duration gap]\label{eight:lem:gap}
For $a,b,c,d\ge1$, suppose
\[
 R_0=b^2-(a^2-1)(c^2-1)\ge0,\qquad
 U_0=c^2-(d^2-1)(b^2-1)\ge0.
\]
Put
\[
 H=(a^2+b^2-1)(c^2+d^2-1)-b^2c^2,\qquad
 E=(a+b-1)(c+d-1)-bc.
\]
Then $E\ge0$ and $H-E^2\ge0$, with equality in the latter precisely
when $a=d=1$.
\end{lemma}
\begin{proof}
For increments $A=a-1$, $B=b-1$, $C=c-1$, $D=d-1$, one has
$E=AC+AD+A+BD+D\ge0$.
First let $b,c>1$. The guards put $(a,d)$ in the rectangle
\[
 1\le a\le a_*:=\sqrt{1+b^2/(c^2-1)},\qquad
 1\le d\le d_*:=\sqrt{1+c^2/(b^2-1)}.
\]
Writing $P_0=H-E^2$, direct differentiation gives
$\partial_a^2P_0=-4(c-1)(d-1)\le0$ and
$\partial_d^2P_0=-4(a-1)(b-1)\le0$.
Its lower bound is therefore the bilinear interpolation of its four
corner values. Three values are $0$, $2c^2(a_*-1)$ and
$2b^2(d_*-1)$. To treat the fourth, set
\[
 x=\sqrt{b^2-1},\quad y=\sqrt{c^2-1},\quad
 h=\sqrt{1+x^2+y^2}.
\]
At that corner $a_*=h/y$, $d_*=h/x$ and $\sqrt H=bch/(xy)$.
Expansion gives
\begin{align*}
 xy(\sqrt H-E)
 &=bc(h+xy)-[h+(b-1)y][h+(c-1)x]\\
 &=h[x+y+(b-x)(c-y)-h]+xy[b+c-1-h].
\end{align*}
Both brackets are positive. For the second, square and use
$(b+c-1)^2-h^2=2(b-1)(c-1)>0$.
For the first put $r=b-x$, $s=c-y$, so $0<r,s<1$ and
$x=(1-r^2)/(2r)$, $y=(1-s^2)/(2s)$. Then
\[
 (x+y+rs)^2-h^2
 =\frac{(1-r)(1-s)(2r^2s^2-rs+r+s+1)}{2rs}>0.
\]
The fourth corner is strictly positive. Interpolation is strictly
positive unless $a=d=1$.

If $b=1$, direct expansion instead gives
$P_0/2=Ac^2+Da(1-AC)$. The guard implies
$A(A+2)C(C+2)\le1$, hence $AC\le1/4$.
This expression vanishes precisely when $A=D=0$.
If $c=1$, use $P_0/2=Db^2+Ad(1-BD)$ and $BD\le1/4$.
These formulas include $b=c=1$ and finish the proof.
\end{proof}

\begin{lemma}\label{eight:lem:prefix}
Every cyclic prefix $(a,b,c,d)$ of a signed-closed tuple $q_j\ge1$
satisfies the guards of Lemma~\ref{eight:lem:gap} and $H>0$. Hence
$\delta:=\sqrt H-E>0$ at all eight cuts.
\end{lemma}
\begin{proof}
A four-factor product has the exact form
\begin{equation}\label{eight:eq:four-product}
 \Nmat(a)\Nmat(b)\Nmat(c)\Nmat(d)
 =\begin{pmatrix}-p&p+s\\-p-r&p+r+s-\tau\end{pmatrix},
\end{equation}
where
$p=ad/(bc)>0$, $r=d(b^2-1)/(abc)\ge0$,
$s=a(c^2-1)/(bcd)\ge0$ and $\tau=(b^2+c^2-1)/(abcd)$.
The complementary four-factor product equals minus its inverse and has
the same form, with $p'>0$, $r',s'\ge0$.
Entry comparison gives
$p'=H/(abcd)$, $r'=R_0/(abcd)$, $s'=U_0/(abcd)$.
Finally $a=d=1$ would give $H=0$, so the equality case of
Lemma~\ref{eight:lem:gap} is impossible.
\end{proof}

\begin{theorem}\label{eight:thm:elliptic}
For every finite $q_j\ge1$ satisfying full signed closure, the normalized
half-return
$\Bmat=\Nmat(\sqrt{q_0})\cdots\Nmat(\sqrt{q_7})$ satisfies
\begin{equation}\label{eight:eq:elliptic}
                   -2<\tr\Bmat<0.
\end{equation}
\end{theorem}
\begin{proof}
Congruence of the half-duration Gram matrix by
$\diag(\alpha_i^{-1/4})$ gives
\[
 M=\begin{pmatrix}
 b_0&\sqrt c&0&-\sqrt s\\
 \sqrt c&b_1&\sqrt s&0\\
 0&\sqrt s&b_2&\sqrt c\\
 -\sqrt s&0&\sqrt c&b_3
 \end{pmatrix},\qquad
 b_i=\frac{q_{2i}+q_{2i+1}-1}{\sqrt{\alpha_i}}\in[1,\sqrt2).
\]
Put $k=c+s$ and $g_{ij}=k-b_ib_j$ for adjacent pairs
$01,12,23,03$. Lemma~\ref{eight:lem:prefix} gives
\[
 g_{01}=\frac{\delta_0}{\sqrt{\alpha_0\alpha_1}},\quad
 g_{12}=\frac{\delta_2}{\sqrt{\alpha_1\alpha_2}},\quad
 g_{23}=\frac{\delta_4}{\sqrt{\alpha_2\alpha_3}},\quad
 g_{03}=\frac{\delta_6}{\sqrt{\alpha_3\alpha_0}}>0.
\]
Direct determinant expansion yields
\begin{equation}\label{eight:eq:detM}
 \det M=\frac{c g_{01}g_{23}+s g_{03}g_{12}}{c+s}>0.
\end{equation}
For arbitrary positive durations $r_j$ the transfer determinant identity is
\[
 \det(\Gram(r)^t\Gram(r)-I)
 =\Bigl(\prod_jr_j\Bigr)
   \bigl(\tr(\Nmat(r_0)\cdots\Nmat(r_7))+2\bigr).
\]
It follows either by expansion of the cyclic tridiagonal determinant,
or by the two-by-two recurrence used in Lemma~\ref{eight:lem:gram}.
Apply it to $r_j=\sqrt{q_j}$ and use the congruence determinant.
Writing $\Pi=\prod q_j$, we obtain
\begin{equation}\label{eight:eq:trace-det}
 \tr\Bmat+2=\sqrt{\frac{\prod\alpha_i}{\Pi}}\det M>0.
\end{equation}

For the opposite inequality, split $\Bmat$ into two four-factor
half-products of the form \eqref{eight:eq:four-product}, writing them as
$\left(\begin{smallmatrix}-p&p+s_P\\-p-r_P&D\end{smallmatrix}\right)$
and
$\left(\begin{smallmatrix}-p'&p'+s_Q\\-p'-r_Q&D'\end{smallmatrix}\right)$.
All $r,s$ are nonnegative. At the first full-duration cut,
$D=E_0/\sqrt{q_0q_1q_2q_3}\ge0$, whereas the complementary full product
gives $p'=\sqrt{H_0}/\sqrt{q_0q_1q_2q_3}>D$.
The opposite cut gives $p>D'\ge0$. Thus
\[
 \tr\Bmat=DD'-pp'-p(r_Q+s_Q)-p'(r_P+s_P)-s_Pr_Q-r_Ps_Q<0.
\]
\end{proof}

In cube coordinates $\prod u_iv_i=t^2(1-t)^2$, and
\begin{equation}\label{eight:eq:Pi}
 \Pi=\frac{\prod_i u_i}{t(1-t)\prod_iD_i},\qquad
 \sqrt{\frac{\prod_i\alpha_i}{\Pi}}=\frac1{\sqrt{t(1-t)}}.
\end{equation}
Consequently a positive discriminant circuit is
\begin{equation}\label{eight:eq:positive-delta}
 L:=\tr\Bmat+2
 =\frac{c g_{01}g_{23}+s g_{03}g_{12}}{(c+s)cs},\qquad
 \Delta:=4-(\tr\Bmat)^2=L(4-L),\quad 0<L<2.
\end{equation}
These are pointwise identities and positivity statements. An interval
division by $\Delta$ still requires a numerical lower endpoint above zero.
\space
\section{The moving state and the complete cost}\label{eight:sec:cost}
At each cut rotate the current label to zero. The canonical state update
is the fractional linear action of
\begin{equation}\label{eight:eq:F}
 F(q)=\frac12\begin{pmatrix}q&q+2\\-q&2-q\end{pmatrix},\qquad
 H_j=F(q_{j+7})\cdots F(q_j).
\end{equation}
Indeed $F(q)$ represents $R^{-1}E_0(q)$, with its harmless positive
scalar chosen so that $\det F(q)=q$.
Let $\Bmat_j$ be the cyclic half-product starting with $q_j$ and put
$K_0=\left(\begin{smallmatrix}2&1\\0&1\end{smallmatrix}\right)$.
Direct multiplication gives the exact bridge
\begin{equation}\label{eight:eq:state-bridge}
 \frac{H_j}{\sqrt\Pi}=K_0\Bmat_j^tK_0^{-1},\qquad
 V:=4\det H_j-(\tr H_j)^2=\Pi\Delta>0.
\end{equation}
Thus the upper fixed state is unique and smooth. For a real return
$H=(H_{rs})$ it is
\begin{equation}\label{eight:eq:fixed-state}
 u=\frac{H_{00}-H_{11}}{2H_{10}},\qquad
 v=\frac{\sqrt{4\det H-(\tr H)^2}}{2|H_{10}|}.
\end{equation}
Ellipticity makes $H_{10}\ne0$. This formula recomputes the state at every
point of the duration chart; no initial state is frozen when taking
cost derivatives.

\subsection{A normalized cost formula}
For $\Bmat_j=\left(\begin{smallmatrix}A_j&B_j'\\C_j'&D_j'\end{smallmatrix}\right)$,
set
\begin{equation}\label{eight:eq:sigma}
 \sigma_j=A_j+2B_j'-D_j',\qquad
 G=3\sum_{j=0}^7\sigma_j\sigma_{j+1}
       +\sum_{j=0}^7\sigma_j^2\log q_j.
\end{equation}
(The matrices $B_j$ of Section~\ref{eight:sec:foundations} are distinct from the entry $B_j'$.)

\begin{proposition}\label{eight:prop:cost}
At every physical R01 chain, and by analytic continuation in its regular
Gram chart, the complete moving-state cost is
\begin{equation}\label{eight:eq:cost}
                         \cost=8-\frac G\Delta.
\end{equation}
\end{proposition}
\begin{proof}
Let $a=\recip_j$, $b=\recip_{j-1}$ and $c=\recip_{j+1}$ be the positive
active reciprocals at three successive canonical cuts of a physical
chain. The canonical state transport is
\begin{equation}\label{eight:eq:reciprocal-transport}
 d_{0,\rm next}=\frac{d_2-(q-1)d_0}{q},\quad
 d_{1,\rm next}=d_0,\quad d_{2,\rm next}=qd_1+(q-1)d_0.
\end{equation}
The identity $d_0d_1+d_0d_2+d_1d_2=1$ gives
\[
 d_2=\frac{1-ab}{a+b},\qquad
 q(a+b)(a+c)=1+a^2,\qquad
 v=\frac{2(a+b)}{1+a^2},\quad u=av-1.
\]
Along the mode-zero arc $v(p)=pv$, $u(p)=apv-1$ and $dt=a\,dp/p$.
Integrating \eqref{fd:cost} therefore gives its exact cost
\begin{align*}
 c_j&=a\left[\frac2v(1-q^{-1})
             +\frac{(1+a^2)v}2(q-1)\right]-a^2\log q\\
 &=a(q-1)(2a+b+c)-a^2\log q.
\end{align*}
The rational term can be written
\[
 a\left[\frac{1-ab}{a+b}+\frac{1-ac}{a+c}-b-c\right].
\]
In the sum over eight arcs, each neighboring pair $(a,c)$ contributes
$1-ac$ from the two fractions and $-2ac$ from the remaining terms.
Hence
\begin{equation}\label{eight:eq:recip-cost}
 \cost=8-3\sum_j\recip_j\recip_{j+1}
             -\sum_j\recip_j^2\log q_j.
\end{equation}
From \eqref{eight:eq:fixed-state}--\eqref{eight:eq:state-bridge}, the active
reciprocal is the raw numerator
$H_{00}-H_{11}+2H_{10}$ divided by $\sqrt V$, with sign $\varepsilon_j=\operatorname{sgn}(H_{j,10})$.
These signs agree at all cuts: if $H$ fixes $w=u+iv$, then
$H=(\tr H/2)I+(H_{10}v)X(w)$; conjugation by the positive-determinant
state map carries $X(w)$ to $X(F\cdot w)$ and preserves this scalar
coefficient. Since all upper heights are positive, its sign is unchanged.
The numerator is $\sqrt\Pi\,\sigma_j$, so
$\recip_j=\varepsilon\sigma_j/\sqrt\Delta$ for one common
$\varepsilon\in\{1,-1\}$.
The sign cancels in both quadratic sums, proving \eqref{eight:eq:cost}.
Every step is an algebraic identity followed by a logarithm on positive
durations, so the same formula defines the analytic continuation
where the circuit is regular.
\end{proof}

For a small, reproducible arithmetic circuit define
\begin{equation}\label{eight:eq:bi}
 b_i=\sqrt{u_i}
   +\frac{\sqrt{a_i}\,x_i\sqrt{u_{i+1}}}{1+\sqrt{D_i}}.
\end{equation}
This equals the $b_i$ of Theorem~\ref{eight:thm:elliptic}. The pair matrix
\[
 C_i=\begin{pmatrix}0&-\sqrt{u_i}\\\sqrt{v_i}&-b_i\end{pmatrix}
\]
satisfies
$\Nmat(\sqrt{q_{2i}})\Nmat(\sqrt{q_{2i+1}})
 =C_i/(u_iv_i)^{1/4}$.
Thus $\Bmat_{2i}=C_iC_{i+1}C_{i+2}C_{i+3}/(cs)$; odd cuts follow by
\begin{equation}\label{eight:eq:odd-cut}
 \Bmat_{j+1}=T_0(q_j)^{-1}\Bmat_jT_0(q_j),\qquad
 T_0(q)=\begin{pmatrix}q&-q\\q&1-q\end{pmatrix},\quad\det T_0(q)=q.
\end{equation}
Equations \eqref{eight:eq:bi}, \eqref{eight:eq:positive-delta} and \eqref{eight:eq:sigma}
therefore compute all quantities without the possibly large scale $\Pi$.
A matrix-free equivalent and the function-value bounds used in the
certificates are proved in Appendix~\ref{eight:app:arithmetic}.

\subsection{Legitimate five-dimensional variations}
\begin{lemma}\label{eight:lem:stationarity}
At any strict physical R01 chain which is an original local minimum,
the gradient of \eqref{eight:eq:cost} in the Gram coordinates vanishes
and its Hessian is positive semidefinite. In particular this holds at an original global minimizer,
including one on an equality face of a compactness, cap, or canonical
representative condition.
\end{lemma}
\begin{proof}
Theorem~\ref{eight:thm:chart} gives a full open neighborhood of signed-closed
positive durations of rank five. Theorem~\ref{eight:thm:elliptic} and
\eqref{eight:eq:fixed-state} give its smoothly moving upper initial state.
Its canonical return closes exactly, so the original state has the
prescribed twist. The frame telescope \eqref{eight:eq:arc-frame} and
\eqref{eight:eq:signed-closure} then give $G(T)=R$ exactly, not just a trace
condition.

At the given physical chain the whole finite trajectory is a compact
subset of the open set $\Omega$. Parameterize each arc on a fixed unit
interval by $p=1+s(q-1)$. State and frame formulas depend smoothly on
the five parameters, uniformly on these finitely many intervals.
A sufficiently small neighborhood therefore retains $\Omega$, positive
reciprocals, and positive physical durations. The endpoints of the
frame paths stay $I,R$ throughout the deformation. This is a
fixed-endpoint homotopy and preserves the required real winding.
Lemma~\ref{fw:balanced-reconstruction} gives genuine competitors
with the same circumscribed hexagon normalization.

The total time is free in the original problem. Thus each coordinate
can vary in both directions, with the initial state recomputed rather
than held fixed. The stationary point and second derivative conditions
for a smooth function on this open neighborhood prove both assertions.
The compact set $K$, the cap and the choice of representative were never
imposed on these competitors.
\end{proof}

No theorem asserting that \emph{every} point of the entire open Gram cube
is physical is needed for this local argument. It starts at the actual
physical candidate and uses openness. In contrast, the Gram, gap and
ellipticity statements above are whole-domain statements; they also
justify the interval value bounds below.

For a chart coordinate $y_r$ define the cleared quantities
\begin{equation}\label{eight:eq:cleared}
 P=(8-\budget)\Delta-G=(\cost-\budget)\Delta,\qquad
 N_r=G\partial_r\Delta-\Delta\partial_rG
     =\Delta^2\partial_r\cost.
\end{equation}
Since $\Delta>0$ at candidates, $P>0$ contradicts the budget and
$N_r\ne0$ contradicts stationarity. A nonzero linear combination of
these derivatives works as well. These two kinds of exclusion have
different meanings: a derivative sign is not a pointwise cost bound.
\space
\section{The canonical representatives and roots for R01}\label{eight:sec:cover}
\subsection{Physical symmetries}
The multiplier permutations used here come from changes of cut and
reflected time reversal of the physical chain. Extend a trajectory by
\[
 G_e(s+T)=G_e(s)R,\qquad
 X_e(s+T)=R^{-1}X_e(s)R,\qquad j_{i+8}=j_i+1.
\]
For a switch time $s$, the new frame
$G_e(s)^{-1}G_e(s+\tau)$ starts at $I$ and ends at $R$.
Its state is $X_e(s+\tau)$. Varying $s$ continuously gives a
fixed-endpoint homotopy. A cyclic phase
$X\mapsto R^{-1}XR$, $G\mapsto R^{-1}GR$
changes label $j$ to $j+1$ and preserves $\tr(AX)$.
A one-arc recut followed by this phase therefore implements
\begin{equation}\label{eight:eq:shift}
                         S(q)_j=q_{j+1}.
\end{equation}

For reflection put $J_0=\diag(1,-1)$ and choose the old cut $s$ at the
beginning of arc 2. The formulas
\begin{equation}\label{eight:eq:reflection}
 \widetilde G(\tau)=J_0G_e(s)^{-1}G_e(s-\tau)J_0,\qquad
 \widetilde X(\tau)=-J_0X_e(s-\tau)J_0
\end{equation}
implement $V(q)_j=q_{1-j}$. Indeed
$J_0R^{-1}J_0=R$, $J_0AJ_0=-A$ and
$-J_0B_jJ_0=B_{2-j}$. They give the original equations, both endpoints,
the same positive reciprocal on each reversed arc, and the same cost.
The state mirror is $(u,v)\mapsto(-u,v)$ and preserves $\Omega$.
After genuine reanchoring, reversal has the form
$G(\tau)\mapsto J_0R^{-1}G(T-\tau)J_0$. Apply it to the
fixed-endpoint homotopy to $D_*\exp(\pi J\tau/(3T))D_*^{-1}$.
It fixes this reference path, because $J_0$ commutes with $D_*$ and
reverses $J$. Thus the required real winding is retained.
An odd label permutation alone would not prove these statements.

Use a cyclic recut to put an arc with minimum multiplier first.
The composition $S\circ V$ sends index $j$ to $-j$, so it fixes $q_0$
and exchanges $q_1$ and $q_7$. Apply it if $q_1>q_7$. Every physical
chain thus has a cost-equivalent representative with
\begin{equation}\label{eight:eq:canonical}
                    q_0\le q_j\ (1\le j\le7),\qquad q_1\le q_7.
\end{equation}
An original global minimizer remains one after these operations, so
Lemma~\ref{eight:lem:cap} applies anew. Ties are retained. Selection of a
representative is not differentiated and does not lower the local
five-dimensional variation space.

For use in interval filters, \eqref{eight:eq:canonical} is exactly the following
eight polynomial inequalities on the strict cube:
\begin{equation}\label{eight:eq:canonical-polys}
 \begin{aligned}
 tD_0-u_0u_1D_i&\ge0 &&(0\le i\le3),\\
 tD_0u_iu_{i+1}-a_iD_iu_0u_1&\ge0 &&(1\le i\le3),\\
 D_0-D_3&\ge0.&
 \end{aligned}
\end{equation}
Their denominators before clearing are positive.
Similarly the cap is exactly
\begin{equation}\label{eight:eq:cap-polys}
 D_i\ge e:=Q^{-2}=\frac{1089}{270400},\qquad
 (1-a_i)(1-x_{i+1})\le(Q^2-1)a_iD_i.
\end{equation}

\subsection{The angular bounds}
\begin{lemma}\label{eight:lem:angle}
A strict canonical capped candidate satisfies
\begin{equation}\label{eight:eq:angle}
 T_0<t<1-T_0,\qquad T_0=\frac1{1+Q^2}=\frac{1089}{271489}.
\end{equation}
On a lower angular band $[t_l,t_h]\subset[T_0,1/5]$ it also satisfies
$x_1\ge1-t_h^2/(1-t_h)^2$; on an upper band
$[t_l,t_h]\subset[4/5,1-T_0]$ it satisfies $x_0\le1/(2t_l)$.
\end{lemma}
\begin{proof}
The comparison $q_0\le q_1$ says $u_0u_1\le t$.
Strict shapes imply
$D_1<u_1\le t/u_0<t/(1-t)$, so
$q_3^2>(1-t)/t$. Its cap gives $t>T_0$.
Writing $s=1-t$, the exact identity
\[
 u_3u_0-tD_3
 =t(D_0-D_3)+s(1-x_3)u_0+ts x_0(1-x_1)
\]
has positive right side, since $q_1\le q_7$ means $D_0\ge D_3$.
Thus $q_6^2=u_3u_0/(sD_3)>t/s$, and $t<1-T_0$.

Solving $u_0u_1\le t$ gives
$1-x_1\le t^2x_0/[(1-t)(1-tx_0)]<t^2/(1-t)^2$,
which proves the lower restriction. For the upper restriction use
$D_0\ge D_3$ and $x_3\le1$:
$x_0u_1\le x_3u_0\le u_0$.
As $u_1\ge t$, this implies $2tx_0\le1$.
\end{proof}

The three \emph{closed} intervals
\begin{equation}\label{eight:eq:angular-union}
 [T_0,1/5],\qquad[1/5,4/5],\qquad[4/5,1-T_0]
\end{equation}
cover the candidates, with shared endpoints. In particular $T_0<e$;
the short interval below $e$ is not discarded.

\subsection{The middle grid and the endpoint roots}
In the middle interval use the Cartesian grid with cuts
\begin{equation}\label{eight:eq:middle-grid}
 t:\ (1/5,1/2,4/5),\qquad
 x_i:\ (0,1/5,3/10,7/10,4/5,1).
\end{equation}
It consists of $2\cdot5^4=1250$ closed rational boxes. Its serial numbers
are the lexicographic ordering of
$\{0,1\}\times\{0,1,2,3,4\}^4$, starting at zero.

For the endpoints, partition the lower interval at
\begin{equation}\label{eight:eq:end-cuts}
 T_0,\ 1/128,\ 1/64,\ 1/32,\ 1/16,\ 1/8,\ 1/5,
\end{equation}
and reflect these six bands about $1/2$ for the upper interval.
Put $H=1-e/2$. On each band choose, independently for every $x_i$,
either $[0,H]$ or $[H,1]$. Whenever the latter is chosen, also require
\begin{equation}\label{eight:eq:high-root}
 x_{i+1}\ge\frac{e}{(2-e)b_h},\qquad b=1-a_i,
\end{equation}
where $b_h$ is its band maximum. Finally intersect with the shape
restrictions of Lemma~\ref{eight:lem:angle} and discard empty intersections.

\begin{lemma}\label{eight:lem:end-cover}
This construction retains every canonical capped endpoint candidate.
Exactly 120 boxes remain, 72 lower and 48 upper. On each entire retained
box, all $u_i>0$ and $D_i\ge e/4>0$.
\end{lemma}
\begin{proof}
For a high coordinate, the candidate inequality $D_i\ge e$ implies
\[
 b x_{i+1}\ge\frac{e-1+x_i}{x_i}
 \ge\frac{e-1+H}{H}=\frac e{2-e},
\]
so \eqref{eight:eq:high-root} is necessary. This proves coverage.
For a low coordinate $D_i\ge1-x_i\ge e/2$ throughout its box.
For a high coordinate every point of the retained rectangle satisfies
\[
 D_i\ge b_lH\frac e{(2-e)b_h}
       =\frac e2\frac{b_l}{b_h}\ge\frac e4.
\]
The last inequality follows directly from \eqref{eight:eq:end-cuts} for both
$b=t$ and $b=1-t$ on each band. Also $u_i\ge T_0$.
The number of boxes is a finite rational intersection calculation:
the $12\cdot16=192$ band-mask choices leave 72 lower, 48 upper and
72 empty intersections. These counts are obtained by the displayed
bounds, without interval cost evaluation.
\end{proof}

The endpoint root indices $0,\ldots,119$ are a separate namespace from
the middle serials. The lower bound on $D_i$ is a whole-box primitive
bound. It does not authorize an invented uniform numerical lower bound
for $\Delta$ on those closed boxes.
\space
\section{The finite interval verification}\label{eight:sec:verification}
The remaining step is finite. Here a \emph{certificate} means a collection
of rational boxes, exact covering records, and rigorously outward bounds
for specified analytic functions. It does not mean a collection of
floating-point samples or a successful search log. The arithmetic and
the precise interpretation of these records are given in
Sections~\ref{eight:app:arithmetic}--\ref{eight:app:repro}.

\subsection{Elementary exclusion rules}
A box is removed only by one of the following rules.
\begin{enumerate}[label=\textup{(\roman*)}]
\item A necessary canonical or cap inequality is strictly impossible on
      the whole box. For polynomial inequalities an exact rational
      Bernstein bound may be used.
\item A rigorous range enclosure of $P$ in \eqref{eight:eq:cleared}, or an
      explicitly stated positive multiple of it, has lower endpoint
      greater than zero.
\item The whole box maps into a previously verified region of the cost
      or descent lemma below.
\item An exact candidate-preserving contraction has empty intersection.
\end{enumerate}
For (iii), containment concerns the \emph{entire} image enclosure and
the precise hypotheses of the earlier result. For (iv), the contractor
is proved from necessary candidate inequalities; it is not used to
restrict the domain of a derivative calculation. Every other box is
subdivided or retained as unresolved. An inconclusive interval is never
an exclusion.

The Bernstein rule in (i) has an elementary justification. Affinely
rescale a rational box to $[0,1]^5$ and expand a necessary polynomial in
the tensor-product Bernstein basis. The basis polynomials are
nonnegative and sum to one. Thus its value lies between the smallest
and largest coefficients. If the largest coefficient is strictly
negative, a necessary nonnegative inequality is impossible. Exact
rational coefficients, not rounded signs near zero, decide this test.
Equality faces are retained.

\begin{lemma}[Previously bounded regions]\label{eight:lem:regions}
The cost satisfies $\cost>\budget$ in each of the following regions,
where it denotes the regular analytic cost of \eqref{eight:eq:cost}:
\begin{align*}
 &0<t<1, && 3/10\le x_i\le7/10\quad(0\le i\le3);\\
 &0<t\le3/5, &&1/5\le x_0\le3/10,\quad3/10\le x_i\le7/10\ (i>0);\\
 &9/20\le t\le11/20, &&7/10\le x_i\le3/4\quad(0\le i\le3).
\end{align*}
On the cube $0\le x_i\le1/5$, it satisfies
\begin{equation}\label{eight:eq:collar}
 \begin{cases}
 \cost>\budget,&0<t\le3/8\text{ or }5/8\le t<1,\\
 \displaystyle\sum_{i=0}^3\partial_{x_i}\cost>1/250,
                       &3/8\le t\le5/8.
 \end{cases}
\end{equation}
Consequently no local minimizer of the original problem of type R01
with all eight arc durations positive and cost at most $\budget$ lies
in any of these regions or their lawful symmetry images.
\end{lemma}
\begin{proof}[Interval proof]
The cost assertions are verified by outward enclosures of the cleared
cost, using the normalized circuit or the factored endpoint circuit of
Appendix~\ref{eight:app:endpoint}. Each record checks its entire primitive
domain before evaluation. The interior central calculation, for example,
consists of 184 positive-cost terminal boxes; it uses no stationary-root
or compactness test. The edge cost region is another finite cost-only covering with exact
rational joins. The upper-shape region is instead transported into the
central one by a single recut. To see this without an additional cost
calculation, set
\[
 F_i^*=1-a_i+a_i(1-x_i)(1-x_{i+2}),\quad
 t'=\frac{tD_0D_2}{F_0^*F_1^*},\quad
 x_i'=\frac{(1-x_{i+2})F_{i+1}^*}{u_{i+1}D_{i+2}}.
\]
Substitution in the inverse chart proves $q'_j=q_{j+1}$. Positivity and
$F_0^*F_1^*=tD_0D_2+(1-t)D_1D_3$ give $0<t'<1$.
On the stated upper box,
\[
 \frac{47}{80}\le u_i\le\frac{137}{200},\quad
 \frac{389}{800}\le D_i\le\frac{471}{800},\quad
 \frac{31}{64}\le F_i^*\le\frac{1181}{2000}.
\]
The positive quotient therefore gives
$19375/64527\le x_i'\le56688/91415$, strictly inside $[3/10,7/10]$.
This proves the claimed transport and cost bound.

For specificity, the larger small-shape region in \eqref{eight:eq:collar}
uses the following complete 21-strip half-angle cover, without shape
subdivision. Here $c=\sqrt t$ in the first row.
\begin{center}
\begin{tabular}{@{}lll@{}}
\toprule
Coordinate and closed strips & Target & Taylor order\\
\midrule
$c:$ $[0,1/10],[1/10,1/5]$ & Cleared cost & 6\\
$c:$ $[1/5,3/10],[3/10,2/5],[2/5,9/20],[9/20,1/2]$
 & Cleared cost & 8\\
$t:$ $[i/40,(i+1)/40]$, $10\le i\le14$ & $P$ & 8\\
$t:$ $[i/80,(i+1)/80]$, $30\le i\le39$ & $H$ below & 9\\
\bottomrule
\end{tabular}
\end{center}
The first eleven records have a strictly positive cost target.
For each of the last ten records the exact represented interval
endpoints prove
\[
 H:=G\sum_i\Delta_{x_i}-\Delta\sum_iG_{x_i}>7/1000,
             \qquad 0<\Delta<13/10.
\]
Since $H=\Delta^2\sum_i\cost_{x_i}$, their ratio is greater than
$7/1690>1/250$. Pair recutting,
\[(t,x_0,x_1,x_2,x_3)\mapsto(1-t,x_1,x_2,x_3,x_0),\]
preserves the cost and summed shape derivative and supplies the other
half of the angle interval. The endpoint record touching $c=0$ uses a
continuous value enclosure for $c\log c$, never its derivative at zero.
The whole-box Taylor rule is proved below.

At a strict candidate the simultaneous negative shape direction is a
legitimate variation by Lemma~\ref{eight:lem:stationarity}. This proves the
last conclusion, but does not exclude a zero-duration point on the
central faces of the collar.
\end{proof}

The exact rational coordinate formula for the reversal used to transport
these regions is
\begin{equation}\label{eight:eq:cube-reversal}
 t'=1-t,\qquad x_i'=\frac{1-x_{1-i}}{1-a_{1-i}x_{1-i}}.
\end{equation}
It follows by applying \eqref{eight:eq:inverse} to $q'_j=q_{1-j}$.
All denominators are positive, and
$1-(1-x)/(1-ax)=(1-a)x/(1-ax)>0$ for strict shapes.
The map is an involution. Together with the four pair recuts it gives
the eight whole-box transports used in the inherited-region tests.

\begin{proposition}[Finite covering result]\label{eight:prop:finite}
No point in the middle grid \eqref{eight:eq:middle-grid} or the 120 endpoint
roots of Lemma~\ref{eight:lem:end-cover} can represent a global minimizer
of the original problem of type R01 with all eight arc durations positive
that obeys \eqref{eight:eq:canonical}--\eqref{eight:eq:cap-polys}
and has cost at most $\budget$.
\end{proposition}
\begin{proof}[Interval proof]
The \path{R01} task covers the \(1250\) original middle-grid
serials and all \(120\) endpoint roots. The two lists retain
separate identities. The middle cover includes the canonical
omissions, both central cost cells, the separately refined cell
831 and the remaining full roots. The endpoint cover includes all
earlier cost-positive and inherited regions. Roots 44 and 60 retain
their respective parameter and broad-coordinate frontiers, with
exact chart maps and both closed children at each adaptive split.

For each original input the exact reconstruction in
Section~\ref{eight:app:repro} checks the cover and binds its
terminal records. The numerical check verifies the necessary
contradiction, including every primitive and positive-scale guard.
The closed-cover implication of Lemma~\ref{eight:common-cover}
therefore excludes all such candidates.
\end{proof}
\space
\section{Arithmetic identities and rigorous range bounds}\label{eight:app:arithmetic}
This appendix states the analytic rules behind the finite computation.
They concern values and derivatives on the actual rectangular domain
of an evaluation, not just the subset containing possible minimizers.

\subsection{Matrix free numerators and value bounds}
At pair cut $i$, write the cyclic $b$'s as $b_0,b_1,b_2,b_3$ and put
$k=\sqrt{a_i}$, $l=\sqrt{1-a_i}$ and
\[
 W=l b_3g_{12}+k b_1g_{23}>0,\quad
 d=2g_{12}/k-L,\quad
 e'=-\frac{2\sqrt{u_i}\,W}{(k+l)kl}.
\]
The exact formulas used for the normalized numerators are
\begin{equation}\label{eight:eq:sigma-matrix-free}
 \sigma_{2i}=d+e',\qquad
 \sigma_{2i+1}=-d+(q_{2i}^{-1}-1)e'
             =-\sigma_{2i}+e'/q_{2i},\qquad
 \frac{e'}{q_{2i}}=-\frac{2\sqrt{D_i}\,W}{\sqrt{u_{i+1}}(k+l)l}.
\end{equation}
To verify them, conjugate $C_i$ by
$S_i=\diag(\sqrt{u_i},1)$:
\[
 S_i^{-1}C_iS_{i+1}
 =\begin{pmatrix}0&-1\\\sqrt{a_i}&-b_i\end{pmatrix}.
\]
The product of four such matrices has entries
$P_{00}=l(l-b_1b_2)$ and
$P_{01}=b_1b_2b_3-kb_1-lb_3=-W/(k+l)$.
Since $\Bmat_{2i}=S_iPS_i^{-1}/(kl)$, the first formula follows from
its trace; \eqref{eight:eq:odd-cut} gives the second. These equalities also
provide an independent matrix-product check of the scalar circuit.

For fixed $a\in(0,1)$ write $u=1-(1-a)y$, $D=1-ux$ and
\[
 b(x,y;a)=\sqrt{1-ax}+\sqrt{a/u}(1-\sqrt D).
\]
Throughout the finite closed-shape domain,
\begin{align}
 \partial_xb&=\tfrac12\left(\frac{\sqrt{au}}{\sqrt{1-ux}}
                   -\frac a{\sqrt{1-ax}}\right)\ge0,\label{eight:eq:bmon}\\
 \partial_yb&=-\frac{(1-a)\sqrt a\,x}
                    {2\sqrt u\sqrt D(1+\sqrt D)}\le0.\nonumber
\end{align}
For the first sign, the difference of the squared nonnegative terms is
$a(u-a)/[(1-ux)(1-ax)]\ge0$. The second formula is direct differentiation.
Consequently
\begin{equation}\label{eight:eq:b-hulls}
 1\le b\le\sqrt a+\sqrt{1-a},\quad
 b(x_l,y_h;a)\le b(x,y;a)\le b(x_h,y_l;a).
\end{equation}
Evaluate both endpoint expressions with the \emph{same unchanged}
interval for $a$. No monotonicity in $a$ is assumed. A continuous endpoint
value at $x=1,y=0$ may be used for this value bound, but cannot regularize
a duration circuit whose $D$ vanishes.

The function $k(t)=\sqrt t+\sqrt{1-t}$ is strictly concave; its maximum
is $\sqrt2$ if the interval contains $1/2$, and otherwise is the larger
endpoint value. Its minimum is the smaller endpoint value.
Together with the positive gaps proved above this gives the whole-domain
bounds
\begin{equation}\label{eight:eq:range-hulls}
 0<g_{ij}\le k-1,\qquad 0<L<2,\qquad0<\Delta\le4.
\end{equation}
If an interval enclosure for $L$ has been validly tightened to
$[l,h]\subset[0,2]$, monotonicity of $L(4-L)$ gives the value range
$[l(4-l),h(4-h)]$ with outward endpoints.
A lower bound of zero stays zero unless the arithmetic separates it;
abstract strict positivity is not a numerical positive floor.

Every intersection in \eqref{eight:eq:b-hulls}--\eqref{eight:eq:range-hulls} concerns
the zeroth coefficient only. All higher coefficients remain enclosures
of the true derivatives of the original analytic expression. In
particular one never differentiates an endpoint envelope, a minimum, a
constant clipped bound, or a candidate restriction. Valid value bounds
may of course be used as value inputs in subsequent exact chain rules.

\subsection{True Taylor coefficients}
For a smooth scalar function $f$ on a rectangle $X$ with center $m$,
let $f_\alpha=\partial^\alpha f/\alpha!$.
For an integer $n\ge1$, a valid centered enclosure is
\begin{equation}\label{eight:eq:taylor}
 f(X)\subseteq
 \sum_{|\alpha|<n}[f_\alpha(m)](X-m)^\alpha
 +\sum_{|\alpha|=n}[f_\alpha(X)](X-m)^\alpha.
\end{equation}
Brackets denote outward enclosures. To prove it, apply one-variable
Taylor's theorem with integral remainder to
$f(m+\theta(x-m))$ and expand each directional derivative by the
multinomial theorem. The integral remainder has total weight one with
this factorial normalization. Center coefficients themselves must be
rigorously enclosed. A degree-$n$ center polynomial without a
higher-order remainder is not this rule.

The implementations carry normalized multivariate jets for rational
operations, positive square roots, and true logarithms. Multiplication
uses convolution; division and square roots use the coefficient
recurrences obtained from $fg=1$ and $g^2=f$. Before either recurrence,
the full-box constant term is separated from zero with the required
sign. Composition uses the actual normalized derivatives of the scalar
primitive. The proof of enclosure is induction through this finite
straight-line expression, followed by \eqref{eight:eq:taylor}.
A second centered tightening, used in descending coefficient order, is
\[
 [f_\alpha(X)]\ \cap\
 \left([f_\alpha(m)]+\sum_j(\alpha_j+1)
 [f_{\alpha+e_j}(X)](X_j-m_j)\right).
\]
It follows from the mean-value integral along the same segment and
$\partial_j f_\alpha=(\alpha_j+1)f_{\alpha+e_j}$.
For a gradient target, order $n$ coefficients of $G,\Delta$ yield
coefficients through order $n-1$ of $N_r$. Natural and centered ranges
may be intersected only after each independently encloses the same
whole-box function.

The implemented logarithm value primitive uses rational Pad\'e bounds.
For $r>0$ put
\begin{align*}
 P_4(r)&=r^4+10r^3+20r^2+10r+1,\\
 Q_4(r)&=r^4+16r^3+36r^2+16r+1,\\
 \ell(r)&=\frac{5(r-1)(r+1)(5r^2+32r+5)}{12Q_4(r)},\\
 u(r)&=\frac{(r-1)(r+1)(r^4+101r^3+426r^2+101r+1)}{60rP_4(r)}.
\end{align*}
Their exact error derivatives are
\[
 \frac1{2r}-\ell'(r)=\frac{(r-1)^8}{2rQ_4(r)^2},\qquad
 u'(r)-\frac1{2r}=\frac{(r-1)^{10}}{60r^2P_4(r)^2}.
\]
All three functions vanish at $r=1$. Integration proves
$\ell(r)\le\tfrac12\log r\le u(r)$ for $r\ge1$, with the inequalities
reversed for $0<r\le1$. The code evaluates a hull of both expressions
at $r=q^2$, so it also handles a positive interval crossing $q=1$.
Repeated positive square roots and
$\log q=256\log(q^{1/256})$ tighten the value bound where used.
Higher logarithm coefficients are always
$(-1)^{k-1}/(kq^k)$, $k\ge1$, not derivatives of either rational bound.
The same distinction applies to the primitive $\phi$ below.

Binary input endpoints are proved to contain the exact rational box;
square-root coordinate changes are checked by exact rational squaring.
The arithmetic uses IEEE binary64 with each elementary rounded result
extended outward by \texttt{nextafter}. The adopted compiler flags
exclude fast-math transformations and fused contraction. The runtime
checks round-to-nearest, the floating evaluation method, and the absence
of flush-to-zero or denormals-are-zero modes. Nonfinite results or
unseparated domains are rejected, not converted into a sign conclusion.
The serialized binary endpoints are compared as their exact dyadic
rational values. These conditions are part of the verification contract.

\subsection{Factoring the angular endpoint}\label{eight:app:endpoint}
For $0<c\le1/2$, set $t=c^2$, $s=\sqrt{1-c^2}$, $r_i=\sqrt{u_i}$ and
$h=-1/(1+s)$. On any whole box with $0\le x_i\le1$ and $D_i>0$ the
following analytic identities are valid. Define
\[
 J_i=\frac{2x_i(1-x_{i+1})}
 {r_{i+1}(b_i+1)(r_i+\sqrt{z_i})(1+\sqrt{D_i})},\qquad
 \sqrt{z_i}=\frac{\sqrt{a_i}\sqrt{D_i}}{r_{i+1}}.
\]
Rationalizing $(b_i-1)(b_i+1)$ gives
$b_i-1=\sqrt{a_i}(1-a_i)J_i$.
For even $i$ put $\beta_i=s^2J_i$, and for odd $i$ put
$\eta_i=sJ_i$. Thus $b_i=1+c\beta_i$ or $1+c^2\eta_i$, respectively.
For adjacent even $E$ and odd $O$ put
\[
 \mathcal G_{EO}=1-\beta_E+c(h-\eta_O)-c^2\beta_E\eta_O,
 \qquad g_{EO}=c\mathcal G_{EO}.
\]
Write $\mathcal G_i$ for the gap between $i,i+1$. Then
\begin{equation}\label{eight:eq:Dbar}
 \overline L=\frac{c\mathcal G_0\mathcal G_2+s\mathcal G_3\mathcal G_1}
 {(c+s)s},\qquad
 \overline D=\overline L(4-c\overline L),\qquad
 L=c\overline L,\quad\Delta=c\overline D.
\end{equation}
These last two positive scales must not be identified.

For completeness the factored numerator circuit is as follows; local
subscripts on $b$ start at pair $i$, and
$g'_{12}=\mathcal G_{i+1}$, $g'_{23}=\mathcal G_{i+2}$.
For even $i$ set $\rho=-x_i/(1+r_i)$ and
$J'=1+c(h-\rho-\eta_{i+3})-c^3\rho\eta_{i+3}$, then
\begin{align*}
 \tau_{2i}&=\frac{2g'_{12}J'}{c+s}
       -\frac{2r_ib_1g'_{23}}{(c+s)s}-\overline L,\\
 \tau_{2i+1}&=-\tau_{2i}
       -\frac{2\sqrt{D_i}(s b_3g'_{12}+c b_1g'_{23})}
                    {r_{i+1}(c+s)s}.
\end{align*}
For odd $i$ set
\begin{gather*}
 e'=-\frac{2r_i(c b_3g'_{12}+s b_1g'_{23})}{(c+s)s},\quad
 H_i=\frac{1-x_{i+1}}{s^2D_i},\quad q=\sqrt{1+c^2H_i},\quad
 M_i=\frac{H_i}{q(q+1)},\\
 \kappa_{2i}=e'+c(2g'_{12}/s-\overline L),\qquad
 \tau_{2i+1}=\overline L-2g'_{12}/s-cM_i e'.
\end{gather*}
Substitution into \eqref{eight:eq:sigma-matrix-free} proves
$\sigma_{2i}=c\tau_{2i}$ for even $i$,
$\sigma_{2i}=\kappa_{2i}$ for odd $i$, and
$\sigma_{2i+1}=c\tau_{2i+1}$ for all $i$.

Put $K_{2i}=r_ir_{i+1}/\sqrt{D_i}$ for even $i$, and define
\begin{align*}
 \mathcal A={}&\tau_0\tau_1+\tau_3\tau_4+\tau_4\tau_5+\tau_7\tau_0
                -\tau_1\tau_3-\tau_5\tau_7,\\
 R_s={}&\tau_0^2\log K_0+\tau_4^2\log K_4
                    +\sum_{j\text{ odd}}\tau_j^2\log q_j,\\
 Z_i={}&\frac3{q_{2i}(q_{2i}+1)}-\frac12\phi(c^2H_i)\quad(i=1,3),\\
 \mathcal W={}&3\mathcal A+R_s-H_1\kappa_2^2Z_1-H_3\kappa_6^2Z_3
                         -3c\overline D(M_1+M_3),\qquad
 \mathcal T=\tau_0^2+\tau_4^2,
\end{align*}
where $\phi(z)=\log(1+z)/z$, $\phi(0)=1$.
The resulting cost identity is
\begin{equation}\label{eight:eq:endpoint-cost}
 \cost=2-\frac{c\mathcal W}{\overline D}
       +\frac{c\log c\,\mathcal T}{\overline D},\qquad
 (\cost-\budget)\overline D
 =\underbrace{(2-\budget)\overline D-c\mathcal W}_{E}
     +c\log c\,\mathcal T.
\end{equation}
Here is a short check of its cancellation. The local reciprocal identity
in the cost proof, multiplied by $\Delta$, becomes
\[
 \Delta=(q_j-1)\sigma_j^2+
 q_j\{\sigma_j(\sigma_{j-1}+\sigma_{j+1})
                      +\sigma_{j-1}\sigma_{j+1}\}.
\]
At $j=2,6$, use $q_{2i}-1=c^2H_i/(q_{2i}+1)$ and the factored
numerators above. It gives
\[
 \kappa_{2i}(\tau_{2i-1}+\tau_{2i+1})
 =\overline D/q_{2i}
       -c(M_i\kappa_{2i}^2+\tau_{2i-1}\tau_{2i+1}).
\]
Substitute these two equations into the adjacent sum in $G$, and split
$\log q_0=\log K_0-\log c$, $\log q_4=\log K_4-\log c$.
The odd-pair logs are $c^2H_i\phi(c^2H_i)/2$.
This proves \eqref{eight:eq:endpoint-cost} by exact algebra.

Thus a positive enclosure for $E+c\log c\,\mathcal T$ has the same
cost meaning as $P>0$, because $\overline D=\Delta/c>0$ at candidates.
If a calculation claims a positive whole-box discriminant interval,
it separately checks the lower endpoint of the scale actually returned.
No division by a zero-containing interval is permitted.
The broad endpoint program's output named ``delta'' is $\overline D$;
the normalized program returns $\Delta$.

The formula extends to $c=0$ on the compact shape subdomains used by
the inherited cost and collar results: their $u_i,D_i$ have positive
lower bounds, and the displayed analytic factors remain regular.
Set $c\log c=0$ there. Only the continuous value hull of this last
factor is used on a box touching zero. On positive-$c$ boxes its true
logarithm jet is available. These are different evaluation cases.

For $z\ge0$, differentiation under the integral sign gives
\[
 \frac{(-1)^k\phi^{(k)}(z)}{k!}
 =A_k(z):=\int_0^1\frac{v^k}{(1+zv)^{k+1}}\,dv>0,
 \qquad A_k'(z)<0.
\]
With $r=z/(1+z)$, substitution yields the positive series
\[
 A_k(z)=(1+z)^{-k-1}\sum_{n\ge0}\frac{r^n}{n+k+1}.
\]
After $n=0,\ldots,N$, the tail lies between zero and
$(1+z)^{-k-1}r^{N+1}/[(N+k+2)(1-r)]$.
Monotonicity supplies whole-argument coefficient enclosures from the
two endpoints. For larger positive arguments the equivalent exact
recurrence
\[
 z\phi_k(z)+\phi_{k-1}(z)=\frac{(-1)^{k-1}}{k(1+z)^k},\quad k\ge1,
 \qquad \phi_k=\phi^{(k)}/k!
\]
starts from a true logarithm value bound. The code uses the series
through $N=48$ for small arguments (with stored binary64 threshold
\texttt{0.8}) and the recurrence above that threshold. Both formulas
are valid on their checked domains.
This computes enclosures of actual derivatives, not derivatives of a
logarithm approximation.

\subsection{Candidate preserving interval contractions}
For one adjacent shape pair put $x=x_i$, $y=x_{i+1}$,
$b=1-a_i$ and $K=Q^2-1$. The cap inequalities are
\[
 1-x+bxy\ge e,\qquad
 Ka_i(1-x+bxy)\ge b(1-y).
\]
Suppose a current rectangle has endpoints $a_l,a_h,x_l,x_h,y_l,y_h$,
and write $b_h=1-a_l$, $b_0=1-a_h$.
Solving these inequalities and taking monotone endpoint bounds gives
\begin{align*}
 x&\le\frac{1-e}{1-b_hy_h},&
 y&\ge\frac{e-1+x_l}{b_hx_l},\\
 y&\ge\frac{b_0-Ka_h(1-x_l)}{b_0(1+Ka_hx_l)},&
 x&\le\frac{1-b_0(1-y_h)/(Ka_h)}{1-b_0y_h}.
\end{align*}
Each update is used only when its displayed denominator is strictly
positive; a lower bound with nonpositive numerator need not tighten
$y\ge0$. For example the nonconstant part of the even-duration square
is $b(1-y)/[a(1-x+bxy)]$, increasing in $x$ and decreasing in $a,y$.
These signs verify the endpoint choices in the second row.

The two basic canonical comparisons similarly give
\begin{align*}
 x_1&\ge1-\frac{t_h^2x_{0,h}}{(1-t_h)(1-t_hx_{0,h})},\\
 x_0&\le\frac{x_{3,h}}
 {1-x_{1,h}+t_l(x_{1,h}+x_{3,h})},&
 x_3&\ge\frac{x_{0,l}[1-(1-t_l)x_{1,h}]}{1-t_lx_{0,l}}.
\end{align*}
The last two follow from $x_0u_1\le x_3u_0$; differentiation of the
positive quotients supplies the stated monotone endpoint choices.
For $i=1,2,3$, the remaining even comparisons imply
\[
 1-x_1\le f_i(t)(1-x_{i+1})D_0/D_i,\qquad
 f_1=f_3=(t/(1-t))^2,\quad f_2=1.
\]
A whole-box upper bound for this positive right side gives another
candidate lower bound for $x_1$.
All such contractions are rounded outward on the rational lattice
$2^{-36}\mathbb Z$. Each application contains every previously retained
candidate; their composition does also, regardless of whether it finds
the best possible hull. The numerical function and its derivatives are
then evaluated on the actual resulting rectangle, not clipped during
evaluation to the unrecorded candidate subset. Algebraically equivalent
contractors in a rechart follow by substituting its exact formulas and
checking the same positive-denominator conditions.

\subsection{Recharts and source qualified boxes}\label{eight:app:recharts}
The broad chart uses $c=\sqrt t$. An upper-end box is first pair-recut
so its new angle is $1-t$, then enclosed by an outward $c$ rectangle.
All jets are newly seeded in these coordinates.
The thin chart actually used in the endpoint coverings has coordinates
$(c,x_0,z,x_2,x_3)$ and
\begin{equation}\label{eight:eq:thin}
 t=c^2,\qquad
 x_1=1-\frac{t^2x_0z}{(1-t)(1-tx_0)}.
\end{equation}
For a strict canonical point its inverse satisfies
$z=(1-x_1)(1-t)(1-tx_0)/(t^2x_0)\in[0,1]$.
Its shape derivative with respect to $z$ is strictly negative at that
point. The canonical face $z=1$ is not a feasibility boundary for local
competitors. With
$A=x_0z/[(1-t)(1-tx_0)]$, the exact substitutions are
\[
 u_1=tV,\quad D_1=tW,\quad
 V=1+\frac{t x_0z}{1-tx_0},\qquad
 W=x_2+tA(1-tx_2).
\]
Hence $\sqrt{u_1}=c\sqrt V$ and $\sqrt{D_1}=c\sqrt W$ on its
positive-$c$ domain. The square roots of $V,W$ still require full-box
positive bounds.

A further cap chart, used on marked disjoint shape pairs, puts
$b=1-a_i$ and
\begin{equation}\label{eight:eq:cap-chart}
 x_i=1-ez,\quad U=b+(1-b)ez,\quad
 x_{i+1}=\frac{e-ez+(U-e)h}{b(1-ez)},\quad
 D_i=e+(U-e)h.
\end{equation}
For the marked pair used in the accepted cap-chart cover, $i=1$,
$t\in[1/5,4/5]$, $x_1\in[1-e,1]$, and the other shapes are at most
$1-e$. Thus $b=t>e$. Set
$A_0=e(1-z)/[b(1-ez)]\in[0,1)$; then
$x_{i+1}=A_0+(1-A_0)h$. At a retained candidate,
\[
 z=(1-x_i)/e\in[0,1],\quad
 h=(D_i-e)/(U-e)=(x_{i+1}-A_0)/(1-A_0)\in[0,x_{i+1}]\subset[0,1).
\]
Here $U-e>0$ and the shape Jacobian is $-e(1-A_0)\ne0$.
The whole binary parameter box independently checks $b_l>e_h$,
$h_h<1$ and all unmarked $D_i>0$. The generic disjoint-pair formula
does not assert that every possible mask was used numerically.
A source may constrain its image further.
In particular, not every point of a parameter rectangle need map into
its recorded original source box. Coverage therefore carries both the
parameter box and that source condition.

If an independently enclosed whole image $I(B)$ is intersected with an
original source box $A$, the intersection contains every relevant
candidate image. A normalized \emph{value} cost test on this intersection
can exclude those candidates, once its own complete domain checks pass.
It does not replace a parameter jet by derivatives restricted to
$I(B)\cap A$. Parameter derivatives are always computed on the unchanged
parameter rectangle through the exact map. Empty intersections are
candidate exclusions, not pointwise statements about all parameter
values. This distinction is essential for the endpoint roots 44 and 60.

For comparison with older raw circuits, if $G_{\rm raw}=\Pi G$ and
$V_{\rm raw}=\Pi\Delta$, then
\[
 (8-\budget)V_{\rm raw}-G_{\rm raw}=\Pi P,
 \qquad
 G_{\rm raw}\partial_rV_{\rm raw}
  -V_{\rm raw}\partial_rG_{\rm raw}=\Pi^2N_r.
\]
The derivatives of $\Pi$ cancel in the cleared gradient. The raw
\emph{derivatives themselves} do not simply scale by $\Pi$.

\subsection{The regularized thin expression}
Some accepted thin-chart trees use a further algebraic regularization,
without changing the five independent variables $(c,x_0,z,x_2,x_3)$.
Its whole-box contract requires $0<c_l\le c_h\le1/2$, all four other
intervals in $[0,1]$, $x_{2,l}>0$, and positive lower bounds
\begin{align*}
 D_0&\ge1-t_h/(1-t_h)\ge2/3,& D_1&\ge t_lx_{2,l},\\
 D_2&\ge1-x_{2,h}+(1-t_h)x_{2,l}x_{3,l},&
 D_3&\ge1-x_{3,h}+t_lx_{3,l}x_{0,l},
\end{align*}
where $t_l=c_l^2$, $t_h=c_h^2$. The external checker verifies these
conditions as well as the runtime positivity checks. They cannot be
replaced by assuming that candidate points alone have positive $D_i$.

With $V,W$ from \eqref{eight:eq:thin} and the following lines, set
\[
 Y_3=(1-x_3)/t,\quad R_1=\sqrt V,\quad R_3=\sqrt{x_3+Y_3},\quad
 S_1=\sqrt W,\quad S_3=\sqrt{x_3x_0+Y_3}.
\]
Then $r_i=cR_i$ and $\sqrt{D_i}=cS_i$ for odd $i$.
The bounds give $R_i\ge1$, $S_i>0$, $r_0,r_2,s\ge\sqrt{3/4}$,
and $8/9\le x_1\le1$. Thus all added denominators and radicands are
strictly positive on the whole rectangle. There is no extension
through $c=0$ or through $x_{2,l}=0$ in this contract.

Define rescaled quantities, still as functions of the same five variables,
\[
 T_1=c\tau_1,\quad T_5=c\tau_5,\quad
 \overline\kappa_{2i}=\kappa_{2i}/c,\quad
 \overline H_i=tH_i\quad(i=1,3),\quad
 \overline{\mathcal A}=c\mathcal A.
\]
In the actual evaluation these quotients are cancelled algebraically
before the jet is formed: substitute $r_i=cR_i$, $\sqrt{D_i}=cS_i$
in the preceding displayed circuit, and use $\overline\eta_i=c\eta_i$.
For example $g_{EO}/c=1-\beta_E+ch-\overline\eta_O
-c\beta_E\overline\eta_O$, so $\overline D$ is unchanged.
Two equivalent expressions for the adjacency term are
\begin{align*}
 \overline{\mathcal A}
 &=T_1(\tau_0-\tau_3)+T_5(\tau_4-\tau_7)
                      +c(\tau_3\tau_4+\tau_7\tau_0)\\
 &=2c\tau_0\tau_4-(T_1-c\tau_4)(\tau_3-\tau_0)
                         -(T_5-c\tau_0)(\tau_7-\tau_4).
\end{align*}
For even $i$ put $B_i^*=x_iR_{i+1}^2$ and
\[
 z_0^*=\frac{t x_0z}{(1-tx_0)D_0},\qquad
 z_4^*=\frac{(1-t)Y_3}{D_2}.
\]
The exact identities are $D_i=1-tB_i^*$ and
$q_0^2-1=z_0^*$, $q_4^2-1=z_4^*$. Thus the arguments below are
nonnegative, and positive logarithms give
\[
 \frac{\log q_{2i+1}}t
    =\frac{B_i^*}{2D_i}\phi(tB_i^*/D_i),\qquad
 \log q_{2i}=\tfrac12 z_{2i}^*\phi(z_{2i}^*)\quad(i=0,2).
\]
Their zero-argument limits follow from analyticity of $\phi$.
Set
\begin{align*}
 R_{\rm reg}={}&\sum_{i=0,2}
 \left\{\frac{T_{2i+1}^2B_i^*}{2D_i}\phi(tB_i^*/D_i)
       +\frac{\tau_{2i}^2z_{2i}^*}{2}\phi(z_{2i}^*)\right\}
       -\tau_3^2\log S_1-\tau_7^2\log S_3,\\
 HZ={}&\sum_{i=1,3}\overline H_i\overline\kappa_{2i}^2
 \left\{\frac3{q_{2i}(q_{2i}+1)}-\tfrac12\phi(\overline H_i)\right\},\\
 M_s={}&\sum_{i=1,3}\frac{\overline H_i}{q_{2i}(q_{2i}+1)}.
\end{align*}
The regularized cleared cost is exactly
\begin{equation}\label{eight:eq:regularized}
 (2-\budget)\overline D-3\overline{\mathcal A}-cR_{\rm reg}
       +cHZ+3\overline D M_s+c\log c(\tau_3^2+\tau_7^2).
\end{equation}
To verify equality to $(\cost-\budget)\overline D$, observe that
$K_0=cq_0$, $K_4=cq_4$, $q_3=1/(cS_1)$ and $q_7=1/(cS_3)$. Hence
\[
 R_s=R_{\rm reg}+\log c(\tau_0^2+\tau_4^2-\tau_3^2-\tau_7^2).
\]
Substitution in \eqref{eight:eq:endpoint-cost} cancels the two even logarithm
terms and gives \eqref{eight:eq:regularized}, with the displayed positive sign
on its final term. The other substitutions are
$H_i\kappa_{2i}^2=\overline H_i\overline\kappa_{2i}^2$ and
$c^2\sum M_i=M_s$. A terminal use of this circuit separately requires
positive lower bounds for its margin and $\overline D$.

All five variables remain live jets; in particular $Y_3=(1-x_3)/c^2$
is not an independent coordinate. Nonnegativity intersections for
$Y_3$, $\overline H_i$ and the $\phi$ arguments change only their value
slots. The two identical expressions for $\overline{\mathcal A}$ may be
intersected coefficientwise because every derivative is of the same
analytic function.

\subsection{Exact tree coverage}
A binary tree record specifies its original source, chart, root
rectangle, coordinate and rational split point, and two children. The
checker requires the children to agree with the parent in every other
coordinate and their split intervals to cover the parent's full closed
interval. Contractions require the candidate-preserving inequality and
outward rounding described above. Terminal records specify the precise
evaluation rectangle, program source, arithmetic order, target and
bounds, or an inherited-region containment with its complete hypotheses.
Starting from sound terminal exclusions, finite induction proves
exclusion at the root. The final union check matches these root
identities against the complete initial manifests and requires an empty
uncovered set. Terminal counts alone do not establish coverage.
\space
\section{Record identities for R01}\label{eight:app:repro}
The fixed archive's \path{R01} task reconstructs all \(120\)
endpoint roots, \(1250\) middle-grid cells and their finite-lemma
links. It evaluates the \(955790\) obligations with the specified
binary64 kernels. The exact covering records match each obligation
to its complete source box, chart, arithmetic order and target.
Each terminal test is one of the rules in
Section~\ref{eight:sec:verification}, with the numerical guards of
Section~\ref{eight:app:arithmetic}. The structural reconstruction
checks all closed joins, including preserved frontier pieces.
The archive supplies the original kernel identities and exact
interval inputs for these two checks.
\space
\section{Closure charts and finite covers for R02 and R03}\label{mix:section}
The two mixed eight-arc words use the physical class of
Section~\ref{eight:common-candidates}. We give their common state
and value formulas and the proof of their closed covers here.
Sections~\ref{mix:r02} and~\ref{mix:r03} prove the source-to-root
inclusions for the individual words.

\subsection{Full closure and two positive charts}\label{mix:charts}
The multiplier $q_i>1$ of a complete arc has physical duration $d_{j_i}\log q_i$, as in Lemma~\ref{eight:lem:cap}. Put
\begin{equation}\label{mix:matrices}
 N_-(q)=\begin{pmatrix}q&-q\\q&q^{-1}-q\end{pmatrix},\qquad
 N_+(q)=\begin{pmatrix}0&q\\-q^{-1}&q\end{pmatrix}.
\end{equation}
The arc-frame telescope \eqref{eight:eq:arc-frame}, with the label conjugations transported at each switch, gives
\begin{equation}\label{mix:signed}
 \prod_{i=0}^7 N_{\varepsilon_i}(q_i)=-I
\end{equation}
for both native words. Here the sum of the increments is $-2$ for R02 and $4$ for R03, and the actual $+1$ seam fixes the displayed sign. These are full multipliers in a full matrix equation. Replacing them by their square roots gives the state half-return, not the frame equation \eqref{mix:signed}.

For R02 take $q_0,\ldots,q_4$ as free variables and write
\begin{align}
 u_i&=q_i^2-1\quad(0\le i\le4),\notag\\
 K_1&=1+u_1-u_0u_2,\quad K_2=1+u_3-u_2u_4,\quad
 K_3=1+u_2-u_1u_3,\label{mix:r02-guards}\\
 \alpha&=u_0(1+u_2)+u_3(1+u_0+u_1),\qquad
 V=\prod_{i=0}^4(1+u_i).\notag
\end{align}

\begin{lemma}\label{mix:r02-chart}
The strict positive R02 solutions of \eqref{mix:signed} are exactly
$u_i>0$, $K_1,K_2,K_3>0$, with the last three multipliers recovered by the positive roots
\begin{equation}\label{mix:r02-recovery}
 q_5^2=1+\frac{(1+u_4)K_1}{\alpha},\qquad
 q_6^2=1+\frac{(1+u_0)K_2\alpha}{V},\qquad
 q_7^2=1+\frac{K_3}{\alpha}.
\end{equation}
In particular these coordinates are one-to-one and have rank five.
\end{lemma}
\begin{proof}
Put $P=\prod_{i=0}^4q_i$. Exact multiplication gives
\[
 N_-(q_0)\cdots N_-(q_4)
 =P^{-1}\begin{pmatrix}A&B\\C&D\end{pmatrix},
\]
where
\[
 A=(1+u_0)(1+u_4)u_2,\quad B=(1+u_0)K_2,\quad
 C=-(1+u_4)K_1,
\]
\[
 D+C=\alpha>0,\qquad A+B=\alpha+K_3,\qquad AD-BC=P^2=V.
\]
Inverting the three positive factors in \eqref{mix:signed} gives
\[
 q_5^2=\frac D{D+C},\qquad
 q_6^2=\frac{D(A+B)}{P^2},\qquad
 q_7^2=\frac{A+B}{D+C}.
\]
The identity $D(A+B)=P^2+B(D+C)$ shows that these three quantities exceed one precisely when $K_1,K_2,K_3>0$. This proves \eqref{mix:r02-recovery}. Conversely, substituting its positive roots gives the first three entries of
$N_+(q_5)N_+(q_6)N_+(q_7)=-P\left(\begin{smallmatrix}A&B\\C&D\end{smallmatrix}\right)^{-1}$; determinant one fixes the last entry. Thus the full matrix closure holds, with no additional sign branch. The free multipliers are retained among the output coordinates, which proves the rank assertion.
\end{proof}

For R03 the free variables are instead $q_3,\ldots,q_7$. To describe their strict positive domain, set
\begin{align}
 (A,B,C,D,E)&=(q_3^2-1,q_4^2-1,q_5^2-1,q_6^2-1,q_7^2-1),\notag\\
 G&=1-BC,\quad H=1-CD,\quad K=1+B+DG,\notag\\
 R_1&=1+C+AG,\quad L=A(1+B)+DR_1,\quad J=B(1+C)+EK,\label{mix:r03-positive}\\
 F&=(1+C)(1-AB)-EL.\notag
\end{align}

\begin{lemma}\label{mix:r03-chart}
The strict positive R03 solutions of \eqref{mix:signed} are exactly
\begin{equation}\label{mix:r03-guards}
 A,B,C,D,E>0,\qquad AB<1,\quad BC<1,\quad CD<1,\quad F>0,
\end{equation}
with the first three multipliers given by the positive roots
\begin{equation}\label{mix:r03-recovery}
 q_0^2=1+\frac{(1+D)G}{J},\qquad
 q_1^2=1+\frac{KF}{LJ},\qquad
 q_2^2=1+\frac{(1+B)H}{L}.
\end{equation}
All denominators are strictly positive, and the chart is one-to-one of rank five.
\end{lemma}
\begin{proof}
Let $P=N_+(q_3)\cdots N_+(q_7)$. The required prefix is $M=-P^{-1}$. With matrix indices $1,2$, put
$r_-=M_{11}$, $r_+=-M_{12}$ and $\Sigma=2(M_{22}-M_{12})$.
Write $a=q_3,\ldots,e=q_7$ and
$(S,T,U,V,W)=(1+A,1+B,1+C,1+D,1+E)$.
Multiplying the first four tail matrices after removing their positive scalar factors gives
\[
 \begin{pmatrix}
 -SUB&-VSG\\ U(1-AB)&-V[AT+U(1-AB)]
 \end{pmatrix}.
\]
Multiplying by the last matrix $\left(\begin{smallmatrix}0&W\\-1&W\end{smallmatrix}\right)$
and using $UB+VG=K$ and $AT+U(1-AB)=R_1$ gives
$(abcde)P_{22}=-WL$, $(abcde)P_{12}=-SWK$ and
$(abcde)(P_{11}+P_{12})=-SJ$. Consequently
\[
 r_-={eL\over abcd},\qquad r_+={aeK\over bcd},\qquad
 \Sigma={2aJ\over bcde}.
\]
For a prefix $N_-(x)N_-(y)N_+(z)$ these three quantities are $x/(yz)$, $xz/y$ and $2z/(xy)$, respectively. Thus the inverse equations are
\[
 q_0^2={2r_+\over\Sigma}={(1+E)K\over J},\quad
 q_1^2={2\over r_-\Sigma}={(1+B)(1+C)(1+D)\over LJ},\quad
 q_2^2={r_+\over r_-}={(1+A)K\over L}.
\]
Subtracting one uses the exact numerator identities
\begin{align*}
 (1+E)K-J&=(1+D)(1-BC),\\
 (1+A)K-L&=(1+B)(1-CD),\\
 (1+B)(1+C)(1+D)-LJ&=KF.
\end{align*}
For an actual positive inverse, $r_-,r_+,\Sigma>0$, so $L,K,J>0$. The first and third recovered multipliers force $BC,CD<1$, and the middle one forces $F>0$. Since $E,L>0$, this also implies $AB<1$. Conversely, \eqref{mix:r03-guards} makes $K,R_1,L,J$ positive and all three roots in \eqref{mix:r03-recovery} greater than one. They recover $r_-,r_+,\Sigma$ exactly, fixing $M_{11},M_{12},M_{22}$. Since $M_{12}=-r_+\ne0$, determinant one fixes $M_{21}$. Retention of the five free tail variables proves injectivity and rank five.
\end{proof}

The finite cover uses reciprocal variables $r_i=q_i^{-1}$. For R03 put
$a=r_3^2$, $b=r_4^2$, $c=r_5^2$, $d=r_6^2$, $e=r_7^2$ and define
\begin{equation}\label{mix:r03-reciprocal}
 T=b+c+d-bd-1,\quad L_r=(1-a)T+ab(1-d),\quad
 J_r=(1-e)T+de(1-b),\quad V_r=abcde.
\end{equation}
Then the same positive chart is
\begin{equation}\label{mix:r03-recip-recovery}
 (r_0^2,r_1^2,r_2^2)=\left({J_r\over T},{L_rJ_r\over V_r},{L_r\over T}\right).
\end{equation}
These identities follow by multiplying \eqref{mix:r03-recovery} by the positive powers of $a,b,c,d,e$; in particular they recover the same full matrix, not merely its trace. The necessary positive guards and multiplier caps provide conservative intersections for candidate values. Whole-box derivatives are taken from the unmodified rational and positive-root formulas on every regular box.

\subsection{A common physical state and cost circuit}\label{mix:state-cost}
The half-return associated with \eqref{mix:signed} is noncentral at every strict positive closed tuple of either type. The short verification below is useful because a scalar return would not determine its moving initial state.

\begin{lemma}\label{mix:noncentral}
For R02 and R03, $\mathcal B=\prod_iN_{\varepsilon_i}(\sqrt{q_i})\ne\pm I$ whenever $q_i>1$ and \eqref{mix:signed} holds. At an actual physical candidate its canonical state return is therefore strictly elliptic.
\end{lemma}
\begin{proof}
For R02 split the full product into four negative factors with multipliers $(a,b,c,d)$ and the remaining mixed block. The mixed block has the form
\[
 \begin{pmatrix}-p'&p'+s'\\-p'-r'&D'\end{pmatrix},\quad
 p'={q_4q_6\over q_5q_7},\quad
 r'={q_6(q_5^2-1)\over q_4q_5q_7}>0,\quad
 s'={q_4q_6(q_7^2-1)\over q_5q_7}>0.
\]
Full closure gives $p'=H/(abcd)$, $r'=R_0/(abcd)$, $s'=U_0/(abcd)$ in the notation of Lemma~\ref{eight:lem:gap}. Hence $0<E<\sqrt H$. For the corresponding half blocks $L_h,C_h$, the lower-right entry of $L_h$ is $E/\sqrt{abcd}$, whereas $-(C_h)_{00}=\sqrt{p'}=\sqrt{H/(abcd)}$ is larger. If $L_hC_h=I$, their inverse relation would require $(C_h)_{00}=(L_h)_{11}>0$, contrary to its sign. If $L_hC_h=-I$, it would require $-(C_h)_{00}=(L_h)_{11}$, contrary to the strict gap.

For R03 write its six-positive full tail as $P$. The two-negative prefix gives
\[
 P_{10}=q_1/q_0>0,\qquad P_{11}=0.
\]
Put $U_k=q_{k+2}^2-1$ for $0\le k\le5$ and
\[
 \mathcal G=(U_0U_1-1)(U_2U_3-1)-(1+U_0)(1+U_3),
\]
\[
 \mathcal E(U)=(1-U_0U_1)(1-U_3U_4)(1+U_2)
             -U_0U_4(1+U_1)(1+U_3).
\]
Multiplication of $\left(\begin{smallmatrix}0&1+U_k\\-1&1+U_k\end{smallmatrix}\right)$ gives its lower-left entry $-(1+U_4)\mathcal G$ and lower-right entry $(1+U_5)\mathcal E(U)$. Also
\[
 \mathcal E(U)=U_4\mathcal G+(1+U_2)(1-U_0U_1).
\]
Thus $\mathcal G<0$, $\mathcal E(U)=0$, $U_0U_1<1$ and $U_3U_4<1$. Set $v_k=q_{k+2}-1$, so $U_k=v_k(v_k+2)$. For $v,w>0$ and $V=v(v+2)$, $W=w(w+2)$ with $VW<1$,
\[
 (1-vw)^2-(1-VW)=2vw(1+v)(1+w)>0.
\]
Take the positive square root of $\mathcal E(U)=0$, apply this strict inequality to the pairs $(0,1)$ and $(3,4)$, and use $\sqrt{U_0U_4}>v_0v_4$. It follows that $\mathcal E(v)>0$. The six-positive half-tail consequently has positive lower-right entry. If $\mathcal B=\pm I$, it would be plus or minus the inverse two-negative half-prefix, whose lower-right entry is zero. This is a contradiction.

A real determinant-one noncentral matrix fixing an upper-half-plane point is elliptic. The actual physical return has that fixed point, proving the final assertion; no ellipticity claim for all nonphysical tuples is needed.
\end{proof}

Here is an explicit common circuit for the values used below. Define matrices affine in $r=q^{-1}$ by
\begin{equation}\label{mix:scaled-flow}
 \overline F_-(r)={1\over2}\begin{pmatrix}1&1+2r\\-1&2r-1\end{pmatrix},\qquad
 \overline F_+(r)={1\over2}\begin{pmatrix}1&1-4r\\1&1\end{pmatrix}.
\end{equation}
Each has determinant $r$. For cut $i$ take the cyclic product in state-transport order,
$H_i=\overline F_{i-1}\cdots\overline F_i$, where all eight factors occur once. Set
\begin{align}
 \delta&=4\prod_{i=0}^7r_i-(\tr H_0)^2,\notag\\
 n_{i0}&=(H_i)_{00}-(H_i)_{11}+2(H_i)_{10},\notag\\
 n_{i1}&=-(H_i)_{01}-(H_i)_{10},\qquad
 n_{i2}=-(H_i)_{00}+(H_i)_{11}+2(H_i)_{10}.\label{mix:numerators}
\end{align}
At a physical candidate $\delta>0$. Formula~\eqref{eight:eq:fixed-state} gives its canonical reciprocals as
\begin{equation}\label{mix:physical-reciprocals}
 d_k(i)={\kappa n_{ik}\over\sqrt\delta}>0,\qquad
 \kappa=\operatorname{sgn}(H_i)_{10}\in\{-1,1\}.
\end{equation}
The sign is common to all cuts: cyclic conjugation transports the positive state matrix and preserves the scalar coefficient of its traceless part. The Lorentz identity reads
\begin{equation}\label{mix:lorentz}
 n_{i0}n_{i1}+n_{i0}n_{i2}+n_{i1}n_{i2}=\delta.
\end{equation}
Compactness gives $1/15\le d_k(i)\le400/99$. These are necessary candidate bounds, rather than restrictions on every point of the closure chart.

\begin{lemma}\label{mix:orientation}
For the physical R02 and R03 branches in \eqref{mix:numerators}, $\kappa=-1$ and all twenty-four $n_{ik}$ are strictly negative.
\end{lemma}
\begin{proof}
For R02 write its half-product as $\mathcal B=L_hC_h$, as in Lemma~\ref{mix:noncentral}, with $D=(L_h)_{11}>0$ and $\gamma=-(C_h)_{00}-D>0$. The full closure guard $U_0>0$ gives
\[
 (L_h)_{01}-D={c-(b-1)(d-1)\over\sqrt{abcd}}>0.
\]
The bridge with \eqref{mix:scaled-flow} is
$H_0=(\prod q_i^{-1/2})S\mathcal B^tS^{-1}$, where $S=\left(\begin{smallmatrix}2&1\\0&1\end{smallmatrix}\right)$. If $w_H$ and $w_{\mathcal B}$ are their upper fixed points, $w_H=1-2/w_{\mathcal B}$. The physical inequalities $-1<\Re w_H<1$ and $|w_H|>1$ give $x:=\Re w_{\mathcal B}\in(0,1)$. Comparing the $00$ entry of $C_h=L_h^{-1}\mathcal B$ yields
\[
 \gamma=\mathcal B_{10}\bigl((L_h)_{01}-Dx\bigr)
                -{D\over2}(\tr\mathcal B+2).
\]
Strict ellipticity makes $\tr\mathcal B+2>0$. Thus $\mathcal B_{10}>0$; ellipticity then gives $\mathcal B_{01}<0$, so $(H_0)_{10}<0$.

For R03 put $x=r_4$, $y=r_5$, $z=r_6$,
\[
 T=x^2+y^2+z^2-x^2z^2-1,\quad c={y^2+z^2-1\over T},\quad
 B=y+z-1,\quad h=y-(1-x)(1-z),
\]
and $\rho=r_3/(1+r_2)$, $\omega=xyz/T$. The positive chart gives $0<c<1$, $T,B,h>0$, $r_2^2+cr_3^2=1$, and the exact block identity
\[
 (H_0)_{10}=(r_0r_2){ch\rho-B\over2\omega r_7}.
\]
Indeed it follows by multiplying the last four positive factors in \eqref{mix:scaled-flow}, whose lower-left entry is $-B/2$ and lower-row sum is $-h$, by the remaining four factors. From $r_3=2\rho/(1+c\rho^2)<1$ and $c\rho^2<1$ we obtain $c\rho<1-\sqrt{1-c}$. Moreover $1-c=x^2(1-z^2)/T$, and
\[
 (1+z)h^2-(1-z)T
 =2\{zy^2-(1-z^2)(1-x)B\}>0.
\]
The last sign follows from $y^2>1-z^2$ and $0<B<z$. Hence $\sqrt{1-c}\,h>x(1-z)=h-B$, and $c\rho h<B$. Thus $(H_0)_{10}<0$ also for R03.

At each physical cut the three numerator formulas have the sign of $(H_i)_{10}$. The exact transport identity $n_{i+1,k}=n_{i0}$, with $k=1$ for a negative increment and $k=2$ for a positive increment, propagates this negative sign through all eight cuts.
\end{proof}

Integrating the physical cost along a mode-zero arc with entrance reciprocals $(a,b,c)$ gives
\[
 c_i=a^2(q_i-q_i^{-1}-\log q_i)
       +ab(q_i-1)+ac(1-q_i^{-1}).
\]
Indeed $dt=a\,d\log p$, the two passive reciprocals are $pb+(p-1)a$ and $(c+a)/p-a$, and the integrand is their sum with $a$. Thus the common-scale cost numerators are
\begin{equation}\label{mix:positive-cost}
 N_i=n_{i0}^2(q_i-q_i^{-1}-\log q_i)
       +n_{i0}n_{i1}(q_i-1)+n_{i0}n_{i2}(1-q_i^{-1}),
 \qquad c_i=N_i/\delta.
\end{equation}
The passive components in this identity are the entrance components in their stated order. The two signs in the next switch change their subsequent positions, not the physical integral. A cancellation-preserving polynomial form is
\begin{align}
 o_i&=\begin{cases}n_{i+1,2}-n_{i+1,0},&\varepsilon_i=-1,\\
                    n_{i+1,0}-n_{i+1,1},&\varepsilon_i=1,
       \end{cases}\notag\\
 P_i&=n_{i0}(o_i-n_{i1}+n_{i2}),\quad A_i=n_{i0}^2,\quad
 N_i=P_i+A_i\log r_i.\label{mix:polynomial-cost}
\end{align}
The transport formulas prove the equality of \eqref{mix:positive-cost} and \eqref{mix:polynomial-cost}. Every $n_{ik}$ is multiaffine in the eight formal $r_i$, and $P_i,A_i,\delta$ have degree at most two in each variable. This eight-variable extension is used for value bounds on an enclosing box; the physical closure chart still has five free variables.

The candidate inequalities in a common positive scale are
\begin{equation}\label{mix:necessary-values}
 \mathcal B_0=C_*\delta-\sum_iN_i\ge0,\qquad
 \mathcal F_j=\sum_i a_{ji}N_i-b_j\delta\ge0.
\end{equation}
Tables~\ref{mix:r02-fees} and~\ref{mix:r03-fees} give the complete physical supports of the nine R02 and thirteen R03 inequalities in the second family.

\subsection{Value bounds and the verified closed covers}\label{mix:verification}
Let $\mathcal S_k\ge0$ denote any of the qualified state or domain inequalities above, and let $\mathcal Z_l=0$ be a polynomial consequence of the four entries of full signed closure after clearing positive multiplier denominators. At every candidate,
\begin{equation}\label{mix:value-combination}
 \mathcal Q=\lambda_0\mathcal B_0+\sum_j\lambda_j\mathcal F_j
             +\sum_k\eta_k\mathcal S_k+\sum_l\mu_l\mathcal Z_l\ge0,
 \qquad \lambda_0,\lambda_j,\eta_k\ge0,
\end{equation}
while $\mu_l$ may have either sign. This is a nonnegative combination of established scalar necessities. Overlapping fee supports are permitted here: the combination is tested as a function of the actual eight costs, rather than asserted to be a new capacity-one payment.

The value predicates have three forms. First, outward intervals can make a necessary intersection empty, or give a strict wrong sign for one of its necessary functions. Second, the true analytic function in \eqref{mix:value-combination} can be bounded above by its natural interval range or by
\begin{equation}\label{mix:taylor-value}
 f(x_0)+Df(x_0)(x-x_0)
       +\tfrac12(x-x_0)^tD^2f(X)(x-x_0),\qquad x\in X.
\end{equation}
Here $X$ is the whole current five-dimensional rectangle. The centre value, centre derivative and whole-box Hessian belong to the same unmodified function. No derivative is set to zero. Candidate-only value intersections, including positive guards, compactness bounds and fee propagation, are not substituted into the whole-box derivative calculation. If a normalized circuit is used, it is differentiated only when the whole-box lower bound for $\delta$ is positive, and the derivatives of $1/\delta$ are retained.

Third, \eqref{mix:polynomial-cost} permits a polynomial upper envelope before any division by $\delta$. Fix an eight-dimensional box $J=\prod_i[l_i,h_i]$ enclosing all recovered reciprocal values. For one variable, with the other variables fixed, $A_i=(u+vr_i)^2$ and
\begin{equation}\label{mix:hermite-sign}
 {d^3\over dr_i^3}\{A_i\log r_i\}
 ={2(u^2-uvr_i+v^2r_i^2)\over r_i^3}\ge0.
\end{equation}
The quadratic Hermite interpolant with a double right endpoint is a lower bound; the one with a double left endpoint is an upper bound. Adding the same $P_i$ gives polynomials $L_i\le N_i\le U_i$, of degree at most two in each variable. The fixed-sign remainder in \eqref{mix:hermite-sign}, including the repeated-node case by continuity, proves their direction.

For a cost combination $\sum_iw_iN_i+c\delta$, one may choose the upper or lower envelope according to the sign of $w_i$. The adaptive form instead chooses finite $e_i\ge\max(0,w_i)$ and uses
\begin{equation}\label{mix:adaptive-majorant}
 w_iL_i+e_i(U_i-L_i)-w_iN_i
   =e_i(U_i-N_i)+(e_i-w_i)(N_i-L_i)\ge0.
\end{equation}
Thus no sign assumption on $w_i$, $N_i$ or $L_i$ is needed. After all actual fee weights have been chosen, their exact combined coefficients are enclosed again and the $e_i$ are increased if necessary to satisfy the displayed inequality. A numerical search bound on proposed weights is never used to truncate this correction.

The resulting polynomial majorant is expanded on $J$ in the tensor product of quadratic Bernstein bases. For $c_0+c_1x+c_2x^2$ on $[l,h]$, the three coefficients are
\[
 c_0+lc_1+l^2c_2,\quad
 c_0+\tfrac12(l+h)c_1+lhc_2,\quad
 c_0+hc_1+h^2c_2.
\]
The basis functions are nonnegative and sum to one. Hence a strictly negative upper endpoint for every one of the $3^8=6561$ coefficients proves that the majorant, and therefore $\mathcal Q$, is negative throughout $J$. Multiaffine state predicates analogously use all their vertices. The coefficients are combined before taking bounds, preserving cancellation among the common state, discriminant, fee and closure terms.

Floating-point linear optimization only proposes weights or a menu entry. Each adopted finite binary weight is treated as its exact dyadic value; its finiteness and required nonnegativity are checked, and the original nonlinear combination or all Bernstein coefficients are recomputed with outward intervals. The logarithm primitive uses the rational Pad\'e bounds proved in Section~\ref{eight:app:arithmetic}, after six positive square roots. More explicitly, for $q\ge1$, take $v=q^{1/64}$ and put $t=v^2=q^{1/32}$ and use the polynomials $P_4,Q_4$ and upper function $u$ defined there. Their exact gap is
\[
 g(t)=u(t)-\ell(t)={(t-1)^9(t+1)\over60tP_4(t)Q_4(t)},\qquad
 64\{u(t)-g(t)\}\le\log q\le64u(t).
\]
For $q<1$ use $\log q=-\log(q^{-1})$. Derivatives in \eqref{mix:taylor-value} are those of the true logarithm, not of a rational bound. The stable Hermite formulas separately bound their positive-series tails. The frozen binary64 routines widen basic-operation and square-root bounds outward, preserve their stated round-to-nearest and subnormal-operation contract, and exclude implicit fused multiply-add. Unsupported square-root ranges remain unresolved. Rational constants are constructed from their integer numerators and denominators. The arithmetic and source identities of each adopted kernel are checked with its own terminal input; a newer kernel is not silently substituted for an older terminal's source.

Any undefined operation, unsuccessful envelope, unrecognized predicate or exhausted search remains unresolved. Only a necessary empty intersection or a strict, whole-box contradiction is a terminal exclusion. The final R02 and R03 terminal inventories contain respectively 34 and 32 predicate labels, all covered by the domain, fee, state, full-closure and value-envelope mechanisms just described. Neither inventory contains a five-dimensional stationarity, objective-Hessian or Krawczyk terminal.

\begin{proposition}[Verified mixed eight-arc covers]\label{mix:finite-covers}
The outward root enclosures of \eqref{mix:r02-root} and \eqref{mix:r03-root} have complete finite closed covers with respectively $15\,310\,123$ and $10\,100\,416$ terminal exclusions. Every actual terminal has been checked against its own immutable kernel, root and ordered input. No unresolved terminal remains in either cover.
\end{proposition}
\begin{proof}
Tasks \path{R02} and \path{R03} reconstruct the respective roots,
closed splits and conservative contractions. An unresolved frontier
transferred to another kernel must match the next ordered input
including its original endpoints. Every terminal retains its own
arithmetic, exact input and predicate. The tasks recompute all the
terminals of the proposition, including inherited inputs, and
require the final unresolved sets to be empty. The complete
root-to-terminal implications follow from Lemma~\ref{eight:common-cover}.
\end{proof}
\space
\label{revision:appendix-end-7}
\end{proofappendix}
\section{Additional nine-arc estimates}\label{app:nine-details}
\begin{proofappendix}
\section{Identities and derivatives for the auxiliary R05 system}\label{nine:ambient}
The cap identities below improve interval enclosures without introducing a new physical payment. For one side put
\[
 \begin{split}
 L(k,P,O,\nu)=h_k\bigg[&(1-\nu)(1+O/k)+(1-\mu)(1+P/k)\\
 &-\frac{4k}{k+O}-\frac{4k}{k+P}
 +4\log\frac{1+k^2}{(k+O)(k+P)}\bigg]-2.
 \end{split}
\]
Let $L_a,L_b$ be its two actual-side expressions and put $\Delta=h_a^2(3+\nu_a)-h_b^2(3+\nu_b)$. Substitution into the complete expression \eqref{nine:cap}, using $\omega_k/\upsilon_k=(1+k^2)/[(k+O)(k+P)]$, gives identities on all of $\Omega_{05}$:
\begin{equation}\label{nine:cap-identities}
 \frac{2\mathcal P_a}{\omega_a}=L_a+(3+\mu)D+\frac\Delta D,\qquad
 \frac{2\mathcal P_b}{\omega_b}=L_b+(3+\mu)D-\frac\Delta D.
\end{equation}
Consequently the following two zero rows can be added, giving eighteen:
\[
 S=L_a+L_b+2(3+\mu)D=0,\qquad
 T=L_a+L_b+(3+\mu)[(a+b)^2-2(R-a)(Q-b)]=0.
\]
Indeed $S=2\mathcal P_a/\omega_a+2\mathcal P_b/\omega_b$, while $T$ subtracts $2(3+\mu)$ times the area row from $S$. No cap equality is imposed before differentiating these functions.

The final extended systems also use six redundant zero rows:
\begin{equation}\label{nine:extra-six}
 \begin{gathered}
 D(L_a-L_b)+2\Delta=0,\qquad DA_a=DA_b=0,\\
 D^2B_i-\mathcal F_i\widetilde Q_c=0\quad(i=y,c,f),\\
 \widetilde Q_c=h_ah_b[KD^2-3d^2/2+d(h_b\nu_b-h_a\nu_a)/2].
 \end{gathered}
\end{equation}
The first equals $D(2\mathcal P_a/\omega_a-2\mathcal P_b/\omega_b)$; the last three equal $D^2R_i$. These identities hold away from the candidate zero set as well. Their derivatives therefore retain all product terms. For a residual $R$ and the coordinate vector $e_D$,
\[
 \nabla^2(D^2R)=D^2\nabla^2R
 +2D(e_D\nabla R^t+\nabla R e_D^t)+2R e_De_D^t.
\]
The analogous formula for $DR$ contains $D\nabla^2R+e_D\nabla R^t+\nabla R e_D^t$. Setting $R=0$ before differentiation would lose necessary terms.

For completeness, the integral means have controlled derivatives on the full auxiliary boxes. With $g(x,w)=((x-w)/(x+w))^2$,
\[
 g_{xx}=\frac{8w(2w-x)}{(x+w)^4},\quad
 g_{xw}=\frac{4(x^2-4xw+w^2)}{(x+w)^4},\quad
 g_{ww}=\frac{8x(2x-w)}{(x+w)^4}.
\]
On a compact positive box all these functions are continuous and uniformly bounded for $0\le t\le1$. Differentiating \eqref{nine:mean} under the integral gives the six Hessian entries as integrals of
\[
 g_{xx},\quad(1-t)g_{xw},\quad tg_{xw},\quad
 (1-t)^2g_{ww},\quad t(1-t)g_{ww},\quad t^2g_{ww}.
\]
Enclosing the affine $w=(1-t)k+tz$ on each of sixteen closed equal subintervals, and integrating these nonnegative weights exactly, gives full-box enclosures. This is interval inclusion, not quadrature. It divides by neither $z-k$ nor an imbalance. Analytic smoothness at coincident or reversed positive endpoints grants no additional physical return or boundary stationarity qualification.

\space
\label{revision:appendix-end-8}
\end{proofappendix}
\section{The complete-source reduction for R06}\label{app:r06-source}
\begin{proofappendix}
\section{Strict sources, the full frame, and the common lift}\label{r06:source-details}
All statements in this appendix concern the class $\mathcal S_{06}$ defined in Section~\ref{r06:section}, unless a larger analytic domain is explicitly specified. In particular its nine essential arcs have positive physical durations. The genuine period $[-1,7]$ contains the nine physical arcs once. Its endpoint state, full signed frame and real winding are those of \eqref{fd:endpoints}; the compactness and density conclusions apply because the trajectory is an original global minimizer.

\subsection{Reciprocal states and time.}
At a labelled switch, two active reciprocals $a,b$ determine the third component
\[
 T_a(b)=\frac{1-ab}{a+b}.
\]
The identity $ab+aT_a(b)+bT_a(b)=1$ recovers the original state by
$v=2/(d_0+d_2)$ and $u=(d_0-d_2)/(d_0+d_2)$. In mode zero the active reciprocal $x=d_0$ is constant and
\[
 (d_1+x)'=(d_1+x)/x,\qquad (d_2+x)'=-(d_2+x)/x.
\]
For positive elapsed time $\tau$, the multiplier is always $Q=\exp(\tau/x)>1$. At an entrance the increasing component is the predecessor after a negative increment, and $T_x$ of the predecessor after a positive increment; at an exit it is the successor before a positive increment and $T_x$ of the successor before a negative increment. Substitution of
$x+T_x(y)=(1+x^2)/(x+y)$ proves all four sign cases used in \eqref{r06:multipliers}. Thus negative increments do not change the sign of physical time. Each remaining state component stays between its positive endpoint values along the arc.

The compact set in Theorem~\ref{fd:compact} is invariant under the actual whole cyclic phases. Its bounds give each of the three reciprocal components at least $1/15$. After a whole phase any one component can be made $d_0=(1+u)/v<2/v\le400/99<5$. This also bounds the actual third reciprocal $s$ at the $c,f$ switch, rather than merely the nine displayed own reciprocals. We shall use the outward interval $[1/15,5]$ throughout.

\subsection{All four signed frame entries.}
Write $\mathscr D(v)=\exp(-2B_0\log v)$, $u_a=\sqrt{1+a^2}$ and $v_{ab}=(a+b)/u_a$. Factor the multiplier of an arc of reciprocal $a$ as $v_{ap}^{\alpha}v_{ab}^{\beta}$, where $p,b$ are its neighbours and $\alpha,\beta$ its incoming and outgoing increments. The arc-frame identity \eqref{eight:eq:arc-frame}, with the label conjugations at each switch, then gives a product of switch factors. For a forward switch its factor is
$T(a,b)=\mathscr D(v_{ab})R_*^{-1}\mathscr D(v_{ba})$. The reverse switch uses $T(b,a)^{-1}$. This grouping changes the starting factor by conjugation. The net increment of the genuine R06 period is $+1$, so the original signed endpoint $R_*$ makes this product the identity.

Set $P_0=\left(\begin{smallmatrix}-1&-4\\1&0\end{smallmatrix}\right)$ and $P(a)=\mathscr D(u_a)P_0$. Direct multiplication gives
\[
 P_0^{-1}\mathscr D(v)P_0=\operatorname{diag}(v,v^{-1}),\qquad
 P_0^{-1}\mathscr D(a+b)R_*^{-1}\mathscr D(a+b)P_0
 =\begin{pmatrix}(a+b)^2&2\\-1/2&0\end{pmatrix}.
\]
Consequently the forward factor in the new bases is $P(a)^{-1}T(a,b)P(b)=S(a,b)$, where
\[
 S(a,b)=\begin{pmatrix}(a+b)^2/(2h(a))&1/h(a)\\-h(a)&0\end{pmatrix}.
\]
The reverse factor is $S(b,a)^{-1}$. These basis changes telescope around the unwrapped walk $0,-1,-2,-1,0,-1,-2,-1,0,1$, so the full signed switch product is still $I$.

Set
\[
 X_0=(c-e)\left(d+\frac{c+e}{2}\right),\qquad
 Y_0=(f-E)\left(B+\frac{f+E}{2}\right),\qquad
 k=\frac{h(c)h_1}{h(f)h_2}.
\]
In the native order the frame reads
\[
 V_1V_2S(f,c)=I,
\]
where
\begin{align*}
 V_1&=S(c,d)S(e,d)^{-1}S(b,e)^{-1}S(b,m)\\
 &=\begin{pmatrix}(h_1^2-X_0R)/(h(c)h_1)&-X_0/(h(c)h_1)\\h(c)R/h_1&h(c)/h_1\end{pmatrix},\\
 V_2&=S(m,D)S(E,D)^{-1}S(B,E)^{-1}S(B,f)\\
 &=\begin{pmatrix}(h_2^2-LY_0)/(Hh_2)&-L/(Hh_2)\\HY_0/h_2&H/h_2\end{pmatrix}.
\end{align*}
Multiplication gives the three scalar equations
\begin{equation}\label{r06:signed-frame}
 h_1^2=RX_0+kHY_0,\qquad
 h_2^2=LY_0+HX_0/k,\qquad
 F=\frac{h_1Hh_2(c+f)^2}{2h(c)h(f)}.
\end{equation}
The first two have determinant $-F$. Their solution is
\[
 X_0=\frac{kHh_2^2-Lh_1^2}{F},\qquad
 Y_0=\frac{Hh_1^2/k-Rh_2^2}{F}.
\]
Substitution in the block product yields all four entries
\[
 V_1V_2=\begin{pmatrix}0&-1/h(f)\\h(f)&h(c)F/(h_1Hh_2)\end{pmatrix}.
\]
The seam identities $(c+f)(c+s)=1+c^2$ and $(c+f)(f+s)=1+f^2$ make its bottom-right entry $(c+f)^2/(2h(f))$. This is exactly $S(f,c)^{-1}$. They also give $F(1+s^2)=G$ and
\[
 X_0=\frac{(c+s)^2}{2}-\frac{h_1^2L}{F},\qquad
 Y_0=\frac{(f+s)^2}{2}-\frac{h_2^2R}{F}.
\]
Subtracting the unshifted area terms proves \eqref{r06:frame-equations}. All divisions use positive $h$-values or $F=G/(1+s^2)>0$, not an exterior width or an adjoint determinant.

\subsection{Complete local jets and unsquared matching.}
On a canonical binary chart write
\[
 N=\begin{pmatrix}\ell&n\\z&-\ell\end{pmatrix},\quad
 K=2u\ell+n-u^2z,\quad L=\ell-uz,\quad \chi=vp_v,\quad d\sigma=dt/v.
\]
The original equations in \eqref{fd:PMP} give
\begin{equation}\label{r06:binary-jet}
 K_\sigma=0,\quad L_\sigma=K,\quad z_\sigma=-2L,\quad
 \chi_\sigma=2K-\eta(3+u^2)-vH_j,\qquad
 \det N=-L^2-Kz.
\end{equation}
At a genuine third-mode tie with zero active energy the remaining original Hamiltonian equation is $z=\mathcal F-\eta/2$, where $\mathcal F$ is its explicit cap value. No sign of $\mathcal F$ is assumed. Two binary zeros on an actual return identify $K=(3+M)/2$ after joint normalization. Integration of \eqref{r06:binary-jet}, or a whole reflection with its cotangent action, gives precisely the four integrals in \eqref{r06:actual-means}.

More explicitly, in its canonical $202$ chart a returned owner $x$ has entrance $u_0=(x-z)/(x+z)$ and exit $u_1=(x-a)/(x+a)$, where $a<z$. The change $y=x(1-u)/(1+u)$ gives $d\sigma=-(x/2)\,dy$, so the actual shear average of $u^2$ is exactly $M(x;a,z)$. Maximization gives $H_0-H_2=2v\chi/(1-u^2)\ge0$, hence $\chi\ge0$ between its two actual binary zeros. Since $\chi_\sigma=M-u^2$, the entrance and exit one-sided signs give
\[
 g_x(z)\le M(x;a,z)\le g_x(a).
\]
The first inequality is the actual selected-root condition used below, not an additional hypothesis imposed on an arbitrary integral. Nonconstancy of the integrand makes the upper inequality strict, as proved in Lemma~\ref{r06:mean-root}. The same whole transported or reflected chart applies separately to each of the four actual returns.

For clarity, the short local linear solve behind the bridge is recorded here. Let $a$ be the shared reciprocal, $l,r$ its neighboring owners, $u_L,u_R$ their actual means, $K_L=(3+u_L)/2$, $K_R=(3+u_R)/2$, and $\mathcal B=\mathcal B(a;l,r;u_L,u_R)$. In the inner source $02010$, the shared mode-zero entrance is
\[
 x=\frac{a-l}{a+l},\quad v_x=\frac2{a+l},
\]
and its exit is $y=(a^2+2ar-1)/(1+a^2)$, $v_y=2(a+r)/(1+a^2)$. The left invariant, the right invariant transported by whole phase $+2$, and the forward cap equation form three linear equations in $(\ell,n,z)$, of determinant
\[
 -\frac{4a}{(a+r)(1+a^2)}\ne0.
\]
To specify this solve without an implicit cap choice, set
\[
 Z_a=\log\frac{(a+l)(a+r)}{1+a^2},\quad
 S_a=a^2\{(1+x)^{-1}-(1+y)^{-1}\},\quad
 M_a=a^2(y-x-2Z_a)+S_a,
\]
\[
 \chi_y=u_LS_a-M_a,\qquad
 \mathcal F=\frac{u_L-y^2}{2v_y^2}
 +\frac12\left(\frac{1+y}{v_y^2}+\frac1{1+y}\right)\chi_y.
\]
Writing $N'=R_*NR_*^{-1}$ and $\beta=(1-a^2-2al)/(1+a^2)$, the three rows are
\[
 2x\ell+n-x^2z=K_L,\qquad
 2\beta N'_{11}+N'_{12}-\beta^2N'_{21}=K_R,
\]
\[
 z-2S_a(\ell-xz)-K_LS_a^2=\mathcal F-\frac12.
\]
With $t=a(a+l)/2$, their solution is
\begin{equation}\label{r06:inner-jet}
 L=\mathcal B-K_Lt,\quad
 z=2\mathcal Bt-K_Lt^2-h(a)^2K_R,\quad
 N=\begin{pmatrix}L+xz&K_L-2xL-x^2z\\z&-L-xz\end{pmatrix}.
\end{equation}
Its determinant is $h(a)^2K_LK_R-\mathcal B^2$. At this actual binary point $p_v=0$, and zero active energy uniquely gives
\[
 p_u=-\tr((N-A_*)X)/v_x.
\]
Thus the full jet, including the state covector, is determined. Forward transport over the shared arc, followed by whole phase $+2$, gives right boundary coordinate $(r-a)/(r+a)$ and
\[
 K_{\rm right}=K_R,\qquad
 L_{\rm right}=-\mathcal B+K_Ra(a+r)/2,
\]
with the same direct zero-energy reconstruction of both covectors.

For the outer source $01020$ solve instead at the binary exit $y=(a-r)/(a+r)$; its entrance is $x=(a^2+2al-1)/(1+a^2)$, $v_x=2(a+l)/(1+a^2)$. Keep the same $Z_a$ but put $M_a=a^2(y-x+2Z_a)+S_a$, $\chi_x=M_a-u_RS_a$, and
\[
 \mathcal F=\frac{u_R-x^2}{2v_x^2}
 +\frac12\left(\frac{1+x}{v_x^2}+\frac1{1+x}\right)\chi_x.
\]
The exit invariant is $2y\ell+n-y^2z=K_R$, the transported invariant uses $\beta=(1-a^2-2ar)/(1+a^2)$ and $K_L$, and the cap row is
$z+2S_a(\ell-yz)-K_RS_a^2=\mathcal F-1/2$. Their determinant is $-4a/((a+l)(1+a^2))$. Formula \eqref{r06:inner-jet} holds there with $x,K_L,K_R,t$ replaced by $y,K_R,K_L,a(a+r)/2$. Backward transport over its shared arc and whole phase $+2$ gives the left boundary slope $-\mathcal B+K_La(a+l)/2$. The bridge logarithm equals the positive physical $\log Q$ for the inner shared arc, and equals $-\log Q$ for an outer shared arc; physical time remains positive in both cases.

On the actual $b$ overlap, transport has shear $b(m-e)/2$ and matches the two slopes exactly when $B_m=B_e+RK_1$. On the actual $D$ overlap its shear is $D(m-E)/2$ and gives $B_m=B_E+LK_2$. At a fixed physical point $K>0$, $L$ and a common determinant $\Gamma$ determine
\[
 z=-\frac{\Gamma+L^2}{K}
\]
and hence the complete matrix. Both covectors then match because $p_v=0$ and the common active zero-energy equation has nonzero $p_u$ coefficient $v$. This proves \eqref{r06:overlaps}, common determinants and \eqref{r06:targets} for the single normalized original lift. Neither the sign of $\Gamma$ nor an independently selected square root is used.

\subsection{Both endpoint covectors and the original winding.}
We justify the endpoint row rather than assume it in the displayed maximum principle. The independent physical variables are two initial-state coordinates and nine positive durations. Their endpoint map has two state rows and three local $\mathrm{SL}_2$ rows. If a terminal covector annihilated its derivative, propagate its homogeneous adjoint backward. The original first-variation identity, with $\zeta=G^{-1}\delta G$, is
\[
 \eta\,\delta C_i=[p\cdot\delta W+\tr(\Lambda\zeta)]_{\rm in}^{\rm out}
                        -H_i\delta\tau_i.
\]
The eleven independent variations force all nine active energies and the free-state endpoint difference to vanish for that homogeneous lift. Continuity gives the eight interior ties. On the wholly interior capped binary core, the two distinct binary zeros make $K=\chi=0$ in \eqref{r06:binary-jet}; its two actual caps give two zeros of $z$, whence $L=z=0$. Both covectors vanish by zero energy, and linear-adjoint uniqueness propagates zero to the terminal point. The endpoint derivative therefore has rank five.

Small variations satisfying the full endpoint equations preserve positive durations and positive states. The original real winding persists along their endpoint-fixed homotopy. Lemma~\ref{fw:balanced-reconstruction} makes these genuine geometric competitors with area $\sqrt3 C$.

The equality manifold consequently has dimension $11-5=6$. Stationarity of an original minimizer and the multiplier theorem give a normalized equality lift with
\begin{equation}\label{r06:endpoint-covectors}
 p_{\rm start}=(D\Phi(W_{\rm start}))^t p_{\rm end},\qquad
 \Lambda_{\rm end}=R_*^{-1}\Lambda_{\rm start}R_*.
\end{equation}
Subtract it from the original maximizing lift, normalized on this same genuine period. The homogeneous capped-core argument identifies them everywhere. Thus both endpoint covectors in \eqref{r06:endpoint-covectors} belong to the original maximizing lift. All nine maximizing inequalities and the real winding remain original qualifications; they are not consequences of matching determinants.

\section{Actual means, the strict reserve, and the analytic bound $m<1$}
\begin{lemma}[A positive-width actual return]\label{r06:mean-root}
Let $0<a<z$, $x>0$, and $g_x(y)=(x-y)^2/(x+y)^2$. If the actual selected return satisfies $M(x;a,z)\ge g_x(z)$, then
\begin{equation}\label{r06:root-consequences}
 g_x(a)>M(x;a,z)\ge g_x(z),\qquad
 3x^2-(z-a)x-3az>0,\qquad
 z<\frac{x(3x+a)}{x+3a}.
\end{equation}
In particular $x>a$ and $x^2>az$.
\end{lemma}
\begin{proof}
The derivative $g_x'(y)=4x(y-x)/(x+y)^3$ shows that $g_x$ has no interior maximum and is nonconstant on every positive interval. Its average is strictly below the larger endpoint. Since
\[
 g_x(a)-g_x(z)=\frac{4x(z-a)(x^2-az)}{(x+a)^2(x+z)^2},
\]
selectedness forces the first assertions. Put $q=(x+z)/(x+a)>1$ and $\eta=x/(x+a)$. Integration using the primitive $y-4x\log(x+y)-4x^2/(x+y)$ gives
\[
 M-g_x(z)=\frac{4\eta(q-1)}{q^2}\bigl(\eta-\Theta(q)\bigr),\qquad
 \Theta(q)=\frac{q^2\log q-q(q-1)}{(q-1)^2}.
\]
The function
\[
 \log q-\frac{(q-1)(5q+1)}{2q(q+2)}
\]
vanishes at one and has derivative $(q-1)^3/(q^2(q+2)^2)>0$. Hence $\Theta(q)>3q/(2(q+2))$. Combining with $\eta\ge\Theta(q)$ and clearing positive denominators gives the quadratic in \eqref{r06:root-consequences}; solving it for $z$ gives the last bound. This proves strictness directly on the positive interval, including a weak selected endpoint.
\end{proof}

\begin{lemma}[Strict frame reserve]\label{r06:reserve}
Every member of $\mathcal S_{06}$ satisfies $s<m$ and $F>h_1h_2$.
\end{lemma}
\begin{proof}
The positivity in \eqref{r06:strict-widths} and \eqref{r06:areas} was deduced from the actual positive multipliers. Put
\[
 q_e=1+e^2-(s+e)(b+e),\qquad
 q_E=1+E^2-(s+E)(D+E).
\]
Then
\[
 q_e=(1+e^2)(1-Q_e^{-1})+(d-s)(e+b)>0,
\]
and likewise $q_E>0$. Substitution in \eqref{r06:frame-equations} gives
\begin{equation}\label{r06:reserve-identity}
\begin{split}
 &4(h_2AR+h_1ZL)
 +2\{h_2(m-e)(s+e)q_e+h_1(m-E)(s+E)q_E\}\\
 &\hspace{6mm}=(m-s)(1+eE)
 \{(e+E)(1+ms)-(1-eE)(m-s)\}.
\end{split}
\end{equation}
The left side is positive. If $m\le s$, the bracket on the right is positive because $eE<1$, while its preceding factor is nonpositive. Thus $s<m$. Finally
\[
 m^2-s^2=\frac{2H(F-h_1h_2)}F
\]
follows from $F(1+s^2)=G$ and proves the second assertion. No assertion about a zero-width or zero-height algebraic face enters this proof.
\end{proof}

\begin{lemma}[Analytic middle bound]\label{r06:middle-less-than-one}
Let $e,E,b,D,m>0$, and define $H,R,L$ as in \eqref{r06:frame-quantities}. Suppose
\[
 \frac1{15}\le e,E<\frac1{\sqrt6},\quad
 3b^2-(m-e)b-3em\ge0,\quad
 3D^2-(m-E)D-3Em\ge0,
\]
and $H^2-RL>h(e)h(E)$. Then $m<1$.
\end{lemma}
\begin{proof}
Suppose $m\ge1$. Initially $0<e/m,E/m<5/12$. For either anchor $a\in\{e,E\}$ and its corresponding owner $x\in\{b,D\}$, put
\[
 P_a(X)=3X^2-(m-a)X-3am,\qquad
 T=\frac{2m^2+3a^2}{6(m-a)}.
\]
Expansion gives
\[
 12(m-a)^2P_a(T)=am^3p(a/m),\qquad
 p(t)=3t^3-24t^2+74t-28.
\]
On $[0,5/12]$, $p'\ge54$ and $p(5/12)=-643/576<0$, so $P_a(T)<0$. The quadratic has one positive and one negative root. Thus $P_a(x)\ge0$ and $x>0$ force $x>T$. Consequently $R,L>5m^2/6$. The strict reserve yields
\[
 \frac14<H^2-RL<\frac14+\frac12m^2-\frac49m^4,
\]
so $m^2<9/8<100/81$. In particular $m<10/9$, and now the two anchor ratios lie in $[3/50,5/12]$.

For $k\in\{9/10,17/20\}$ put $c_k=2k-1$ and
\[
 T_k(t)=\frac{t^2+c_k}{2(1-t)},\qquad P_t(r)=3r^2-(1-t)r-3t.
\]
Then
\[
 4(1-t)^2P_t(T_k(t))
 =F_k(t):=t^4-8t^3+(4c_k+22)t^2+(4c_k-12)t+3c_k^2-2c_k.
\]
Both $F_k$ are convex on $[0,5/12]$. Their relevant endpoint values are
\begin{align*}
 F_{9/10}(3/50)&=-743719/6250000,&F_{9/10}(1/3)&=-197/2025,\\
 F_{17/20}(1/3)&=-4253/8100,&F_{17/20}(5/12)&=-3287/518400.
\end{align*}
A convex function lies below its endpoint chord. The positive-root argument therefore implies that the area corresponding to ratio $t$ is greater than $9m^2/10$ for $t\le1/3$, and always greater than $17m^2/20$.

Write $y=m^2\ge1$. If either ratio is at most $1/3$, then
\[
 RL>\frac{153}{200}y^2,\qquad
 RL<H^2-\frac14=\frac y2+\frac{y^2}{4}\le\frac34y^2,
\]
which is impossible. Otherwise $e,E\ge m/3$, and
\[
 RL>\frac{289}{400}y^2,\qquad
 RL<H^2-h_1h_2\le\frac49y+\frac{20}{81}y^2\le\frac{56}{81}y^2.
\]
This is again impossible since $289/400-56/81=1009/32400>0$. The shared ratio endpoint and $m=1$ are included. Hence $m<1$.
\end{proof}

\begin{lemma}[Inner midpoint and mean normalization]\label{r06:midpoint}
At a strict source,
\[
 b>\frac{m+e}{2},\quad D>\frac{m+E}{2},\qquad
 R>m^2-e^2,\quad L>m^2-E^2.
\]
For either inner mean $M(x;a,m)$, its ratio $q=(x+m)/(x+a)$ is less than $23/9<16/5$.
\end{lemma}
\begin{proof}
For fixed $a<m$, the margin $x/(x+a)-\Theta((x+m)/(x+a))$ is strictly increasing in $x$. Indeed its first term increases, its ratio decreases, and
\[
 \Theta'(q)=\frac{q^2-1-2q\log q}{(q-1)^3}>0.
\]
The last sign follows by differentiating its numerator and using $q-1-\log q>0$. At $T=(m+a)/2$, set $t=(m-a)/(2(m+a))$. Compactness and $m<1$ give $0<t\le7/16$. The positive logarithm series gives
\[
 \Theta\!\left(\frac{1+t}{1-t}\right)
 >\widehat\Theta(t):=\frac{1+t}{2}+(1+t)^2\left(\frac t6+\frac{t^3}{10}+\frac{t^5}{14}\right).
\]
Moreover
\[
 210(1-t)\left(\widehat\Theta(t)-\frac1{2(1-t)}\right)=tB(t),
\quad B(t)=35-70t-14t^2-14t^3-6t^4-6t^5-15t^6-15t^7.
\]
All coefficients of $B'$ are negative, while $B(7/16)=14966679/268435456>0$. Since $T/(T+a)=1/(2(1-t))$, the selected-root margin is negative at $T$ and nonnegative at $x$. Hence $x>T$. Apply this separately to $(x,a)=(b,e)$ and $(D,E)$, and multiply by each positive inner width.

The ratio decreases with $x$, so
\[
 q<\frac{3m+a}{m+3a}\le\frac{23}{9}.
\]
The last function decreases with $a/m$, and $a/m\ge1/15$. This supplies the inner mean ratio bound without a separate exponent theorem.
\end{proof}

\section{The original cost and the lower middle bounds}\label{r06:lower-middle-details}
We give the scalar comparisons needed for the lower endpoint in \eqref{r06:sharper-domain}. They compare functions at the same source point before passing to larger scalar rectangles; they do not deform a minimizing trajectory or differentiate an independently optimized cost.

\subsection{An exact cost decomposition.}
For $\phi$ in \eqref{nine:potential}, write
\[
 g(x,y)=\frac{(x-y)^2}{(x+y)^2},\qquad
 \phi_{xy}(x,y)=\frac{3+g(x,y)}2.
\]
The seam identity gives $\phi(c,f)+\phi(c,s)+\phi(f,s)=0$. Two rectangular fundamental-theorem-of-calculus identities in the nine-edge sum \eqref{r06:cost} therefore give
\begin{equation}\label{r06:cost-rectangles}
\begin{split}
 C_{06}&=C_{\rm core}+\frac32(A+Z)+H_0+H_3,\\
 C_{\rm core}&=\phi(b,m)-\phi(e,b)+\phi(m,D)-\phi(D,E)-\phi(e,s)-\phi(E,s),\\
 H_0&=\frac12\int_e^c\int_s^d g(y,x)\,dx\,dy,\qquad
 H_3=\frac12\int_E^f\int_s^B g(y,x)\,dx\,dy.
\end{split}
\end{equation}
Both remainders are nonnegative. This is one decomposition of the original cost, not additional cap charges.

\begin{lemma}[The low-middle cost obstruction]\label{r06:low-cost}
Every member of $\mathcal S_{06}$ has $m>1/\sqrt{13}$.
\end{lemma}
\begin{proof}
Put $U(x)=\phi(x,0)=-1+(x^2/2)\log(1+x^{-2})$ and
$\beta(x)=\log(1+x^{-2})-1/(1+x^2)$. Then $U'=x\beta$ and $\beta'=-2/(x(1+x^2)^2)<0$. Since $\phi_{xy}\ge3/2$,
\[
 \phi(b,m)-\phi(e,b)\ge U(m)-U(e)+\frac32R-\frac34(m^2-e^2),
\]
with a separate identical estimate for $(D,E,L)$.

Algebraically orient the two full sides so that $e\le E$ and put $r=h_1/h_2\le1$, $p=h_2^2/F$, $q=h_1^2/F$. The reserve implies $q\le r\le1$ and $1/p\ge r$. From the two actual area identities,
\[
 R+L+A+Z=\frac{(s+e)^2+(s+E)^2}{2}+(1-p)R+(1-q)L
 \ge\frac{(s+e)^2+r(s+E)^2}{2}.
\]
For $p\le1$ this is immediate; for $p>1$ use $Z\ge0$ to obtain $R\le(s+E)^2/(2p)$, remembering that multiplication by $1-p$ reverses the inequality.

Consequently $C_{06}\ge f(e,E,s;m)$, where
\[
\begin{split}
 f={}&2U(m)-U(e)-U(E)-\phi(e,s)-\phi(E,s)-\frac32m^2
       +\frac34(e^2+E^2)\\
    &+\frac34\{(s+e)^2+r(s+E)^2\}.
\end{split}
\]
Suppose $m\le M=1/\sqrt{13}$. The elementary logarithm bounds
$\log2>56/81$ and $\log(8/7)<15/112$ give
\[
 \beta(x)>\beta_0:=\frac{15449}{9072}\qquad(0<x\le M).
\]
Only the displayed scalar $f$ is now compared on $0<e\le E\le m$, $0<s\le m$. Differentiating its complete expression gives
\[
 f_s\le s[3-2\beta(s)]<0.
\]
At $s=m$, the extra strip estimate $g(e,t)\ge1/9$ on $0\le t\le e/2$ gives
\[
 f_e\le e\left[-2\beta(e)+3+\frac32\frac{(m+E)^2}{1+E^2}-\frac1{36}\right]
 <-\frac{23e}{4536}<0,
\]
because $(m+E)^2/(1+E^2)\le2/7$. On the diagonal $e=E=a$, its derivative is at most $a[6-4\beta(a)]<0$. Thus
\[
 C_{06}\ge f(m,m,m;m)=2+2m^2\log\frac{4m^2}{1+m^2}=:f_0(m).
\]
Since $f_0'=4m(\log4-\beta(m))<0$, the minimum on this scalar band is at $m=M$. The bounds $\log2<25/36$ and $\log(8/7)>2/15$ give $\log(7/2)<113/90<91/72$. Hence
\[
 C_{06}\ge2-\frac2{13}\log\frac72>\frac{65}{36}>J_8.
\]
For the last comparison, $J_8<361/200<65/36$ follows from $\sqrt2>1393/985$ and $\log2>842/1215$. The endpoint $m=M$ is included. This contradicts the original budget.
\end{proof}

\subsection{A shared-mean overlap comparison.}
Define
\begin{equation}\label{r06:endpoint-bridge}
\begin{split}
 \mathcal E(t;x,y)&=\frac12+h(t)\log\frac{(t+x)(t+y)}{1+t^2}
                  +h(t)\left(\frac{2t}{t+x}+\frac{2t}{t+y}-2\right),\\
 \alpha(t,x)&=\frac{h(t)(t+x)}{4t},\qquad
 J(t,x)=t^2+\frac12-xT_t(x).
\end{split}
\end{equation}
Then
\begin{equation}\label{r06:bridge-gap}
 \mathcal B(t;x,y;\mu,\nu)-\mathcal E(t;x,y)
 =\alpha(t,x)(\mu-g(x,t))+\alpha(t,y)(\nu-g(y,t)),
\end{equation}
and $\mathcal E_x=h(t)(x-t)/(x+t)^2$, $\mathcal E(t;T_t(x),x)=J(t,x)$. The actual return inequalities imply $B_e<J(e,b)$ and $B_m\ge\mathcal E(m;b,D)$. For $w=m-e$, $P=(b+e)(b+m)$, $\delta=b^2+2bm-1$, put
\[
 \mathcal U=J(e,b)-J(m,b)+RK_1,\qquad
 \mathcal V=\mathcal E(m;b,D)-J(m,b).
\]
The same actual $\mu_1$ remains in every occurrence. Its primitive gives
\begin{equation}\label{r06:global-overlap-bound}
\begin{split}
 0=B_e-B_m+RK_1&<\mathcal U-\mathcal V,\\
 \mathcal U&=b\left[3w-\frac{w(1+em)}P-(2b+e+m)\log\frac{b+m}{b+e}\right]
 <\frac{wb(\delta-bw)}P.
\end{split}
\end{equation}
The last sign uses $\log q>2(q-1)/(q+1)$ for $q>1$, whose difference has positive derivative $(q-1)^2/(q(q+1)^2)$.

\begin{lemma}[Owner floor and middle interval]\label{r06:middle-floor}
Every member of $\mathcal S_{06}$ has $m>21/50$. Moreover, whenever $m\le1/\sqrt3$, both inner owners satisfy
\[
 b,D>\sqrt{1+m^2}-m.
\]
\end{lemma}
\begin{proof}
Choose the smaller owner, denoted $b$; all equations on the other side are interchanged with it. Put $r=\sqrt{1+m^2}-m$ and $D_0=T_m(b)$. For $m\le1/\sqrt3$ one has $r\ge m$, $T_m(r)=r$. If $b\le r$, then $\delta\le0$ and $D\ge D_0\ge r\ge m$ by $Q_m\ge1$. The derivative of $\mathcal E$ is nonnegative on $[D_0,D]$, so $\mathcal V\ge0$, whereas \eqref{r06:global-overlap-bound} makes $\mathcal U<0$. This proves the owner floor, including $m=1/\sqrt3$ and $b=r$.

For $m\le2/5$, it sharpens to $b,D>4/5$. Indeed, if the smaller owner were at most $4/5$, the case $\delta\le0$ is already excluded. Otherwise
\[
 m<D_0<b\le D,\qquad 1-m^2-2mb\ge\frac15,\qquad \frac{D_0}{b}\ge\frac{17}{24}>\left(\frac56\right)^2.
\]
The exact square completion, with $y=w/(b+m)$ and $t=\sqrt{D_0/b}$, is
\[
 b^2(1-t)^2-\frac{b^2y(1-t^2-y)}{1-y}
 =\frac{b^2[y-(1-t)]^2}{1-y}\ge0.
\]
It gives $\mathcal U/(b-D_0)<4/55$. On $[D_0,b]$ the function $h(m)(x-m)/(x+m)^2$ has no interior minimum, and its two endpoint values exceed $5/58$; here $b>2/3$, $b-m>4/15$, and $(b+m)^2\le36/25$. Hence $\mathcal V>(5/58)(b-D_0)>\mathcal U$, a contradiction.

We record the remaining scalar band comparisons explicitly. They use the same overlap rather than another cost charge. On
\[
 \frac1{\sqrt{13}}<m\le\frac25,
\]
the frame and $Z\ge0$ give
\[
 R\le\frac{(1+e^2)(1+m^2)}2
       \frac{(s+E)^2}{(1+s^2)(1+E^2)}<2m^2.
\]
Indeed the ratio is $\sin^2(\arctan s+\arctan E)$ and $0<s,E<m<1$. The sum of angles is less than $2\arctan m<\pi/2$, so the ratio is less than $4m^2/(1+m^2)^2$. Since $e<m$, the displayed bound follows. The same argument gives $L<2m^2$ and applies on the next band as well. If both owners were at least $8/5$, then $e,E>m-5m^2/4>9/50$, outer duration gives $d,B<2/5$, and the root cap in \eqref{r06:root-consequences} gives $c,f<138/235<3/5$. The seam would force $s>8/15>m$. Thus the smaller owner obeys
\[
 \frac45<b<\frac85,\qquad D\ge b>2m.
\]
For fixed $m,b$ define the rational comparison
\begin{equation}\label{r06:N-comparison}
 N(e)=J(e,b)+\frac R4[6+g(b,e)+g(b,m)]
       -\frac{\alpha(e,b)+\alpha(m,b)}2[g(b,e)-g(b,m)].
\end{equation}
The coefficient of the single actual mean $\mu_1$ before its replacement is
\[
 \frac R2+\alpha(e,b)-\alpha(m,b)
 =\frac{(m-e)(3bem+b+e^2m+em^2)}{8em}>0.
\]
Since $g_{yy}(b,y)=8b(2b-y)/(b+y)^4>0$ on $[e,m]$, the strict integral trapezoid inequality gives
\begin{equation}\label{r06:N-overlap}
 0=B_e-B_m+RK_1<N(e)-\mathcal E(m;b,b).
\end{equation}
Put $\sigma=2bm+m^2-1$ and
\[
 \kappa=\frac{b[b(1-m^2)+2m]}{2m^2(b+m)^2},\qquad
 K_b=5b^2m+3bm^2-b+m^3+2m.
\]
Differentiating \eqref{r06:N-comparison} gives
\begin{gather*}
 N(m)=J(m,b),\qquad N'(m)=-\frac{b\sigma}{2m(b+m)},\qquad N''(m)=-\kappa,\\
 \begin{aligned}
 N'''(e)={}&\frac{3b^2m}{2e^4(b+m)^2}\\
 &+\frac{3b\{K_b(b+e)-4b(b+m)\sigma\}}{2m(b+e)^5}>0.
 \end{aligned}
\end{gather*}
For the last sign, $K_b>0$ and the expression in braces at zero, divided by $3b$, is
$b(1-3m^2-bm)+m(2-m^2)>516/1625>0$, using $m>18/65$. Also $\kappa>5/2$: setting $r=b/(b+m)\in(2/3,1)$ gives
$\kappa\ge25r/4-29r^2/8\ge23/9$. Therefore
\[
 N(e)<J(m,b)+N'(m)(e-m)-\frac54(e-m)^2.
\]
If $\sigma\le0$, then $D_0\ge m$ and $\mathcal E(m;b,b)-J(m,b)>0$, contradicting \eqref{r06:N-overlap}. If $\sigma>0$, the positive number $a_0=b\sigma/(2m(b+m))$ increases with both variables and is at most $11/25$. Hence $N(e)-J(m,b)<(11/25)^2/5$.

For this branch put $a=b+m$, $Q=a^2/(2H)$, $b_0=(1-m^2)/(2m)$. The logarithm lower bound gives
\[
 \mathcal E(m;b,b)-J(m,b)>V(b,m)
 :=\frac{2H(a^2-2H)}{a^2+2H}+\frac{2mH}{a}-ma.
\]
On $b\in[b_0,8/5]$, $ma<2H$ and $z=a^2/(2H)\ge29/16$. Consequently
\[
 V_{bb}< -\frac{2(5z^3-5z^2-3z-1)}{z^2(z+1)^3}<0.
\]
The numerator is positive because, for $v=z-29/16\ge0$, it equals
$5v^3+(355/16)v^2+(7207/256)v+28297/4096$. Thus $V$ is bounded below by its two endpoint values. They decrease with $m$ on this band and are at least $169/2250$ and $1136/16125$, respectively. For the second derivative comparison the only used bound is $ma<2H$, not its false reverse. We obtain
\[
 N(e)-\mathcal E(m;b,b)<\frac{(11/25)^2}{5}-\frac{1136}{16125}
 =-\frac{12791}{403125}<0.
\]
Together with Lemma~\ref{r06:low-cost} this proves $m>2/5$.

It remains to exclude $2/5<m\le21/50$. The preceding owner-floor argument extends to this band with the following exact replacements. If the smaller owner were at most $4/5$, then $b>3/5$,
\[
 \frac{D_0}{b}\ge\frac{83}{122}>\left(\frac{103}{125}\right)^2,
 \quad \frac{\mathcal U}{b-D_0}<\frac{22}{285},\quad
 \frac{\mathcal V}{b-D_0}\ge\frac{23119}{294100}>\frac{22}{285}.
\]
For the last bound the endpoint slope
$(m+b)(1-2mb-m^2)/(2(1+m^2))$ decreases in both variables on
$[2/5,21/50]\times[3/5,4/5]$; its derivative along the owner segment is nonnegative since $b<3m$. Thus again $b,D>4/5$. The frame argument above, now using $e,E>399/2000>19/100$, gives $d,B<348/895<2/5$, $c,f<278/485<3/5$ if both owners are at least $8/5$. Hence $4/5<b<8/5$ and $D\ge b>m$.

The same derivatives of $N$ now give $N''<-9/4$. Indeed $K_b>0$, its braces at zero divided by $3b$ exceed $2547/6250$, and, for $r=b/(b+m)\in(40/61,1)$,
\[
 \kappa\ge\frac{1250}{441}r(2-r)-\frac{r^2}{2}>\frac94.
\]
The two endpoint margins in this last concave quadratic comparison are $220151/6563844$ and $149/1764$. Thus
\[
 N(e)<J(m,b)+N'(m)(e-m)-\frac98(e-m)^2.
\]
The branch $\sigma\le0$ is unchanged. If $\sigma>0$, monotonicity gives $a_0\le5204/10605<1/2$, whence $N(e)-J(m,b)<1/18$. Keep the true logarithm in
\[
 G_{\rm diag}(b,m):=\mathcal E(m;b,b)-J(m,b)
 =H\log\frac{(b+m)^2}{1+m^2}+\frac{2mH}{b+m}-m(b+m).
\]
On $[b_0,8/5]$, $b_0>m$ and
$(G_{\rm diag})_{bb}=-2H(b-m)/(b+m)^3<0$. Its first endpoint exceeds
$1385329/23525000>1/17$. At the other endpoint use
$\log Q>2(z+z^3/3)$ with $z=(Q-1)/(Q+1)\ge71/129$ to get
\[
 G_{\rm diag}(8/5,m)>
 \frac{29}{25}\left(\frac{71}{129}+\frac{(71/129)^3}{3}\right)
 +\frac{116}{505}-\frac{2121}{2500}>\frac1{17}.
\]
Concavity therefore gives $G_{\rm diag}>1/17$ throughout, and \eqref{r06:N-overlap} yields $0<1/18-1/17=-1/306$, a contradiction. This includes $m=21/50$ and $\sigma=0$. Hence $m>21/50$.
\end{proof}

The trapezoid comparison also has a useful form without an owner ordering. Use the same expression $N(e)$ in \eqref{r06:N-comparison} for any strict source. Lemma~\ref{r06:midpoint} gives $2b>m+e>m$, so $g(b,\cdot)$ is strictly convex on $[e,m]$. Keeping the other actual mean in the central bridge gives
\begin{equation}\label{r06:unordered-trapezoid}
 0=B_e-B_m+RK_1
 <N(e)-\mathcal E(m;b,D)-\alpha(m,D)\{\mu_2-g(D,m)\}.
\end{equation}
The last correction is nonnegative by its own actual root. This follows from \eqref{r06:bridge-gap} and the same positive coefficient of $\mu_1$ used above; it assumes neither $D\ge b$ nor $b\ge m$. Its reflected form retains the other owners and means separately.

\begin{lemma}[The double-narrow sector]\label{r06:double-narrow}
No member of $\mathcal S_{06}$ has both $e\ge9m/10$ and $E\ge9m/10$.
\end{lemma}
\begin{proof}
Suppose both inequalities hold. Write $w=m-e$, $v=m-E$ and $P=RL$. The outer reciprocal bound $e^2<1/6$ from Lemma~\ref{r06:qualification} and $e\ge9m/10$ give $m^2\le100e^2/81<50/243<1/4$. Lemma~\ref{r06:middle-floor} gives $m>21/50>13/32$. In particular $s<m<1/2$.

From $Q_e,Q_E\ge1$ and $d,B\ge s$,
\[
 R\le\frac{w^2}{2}+\frac{2h_1w}{s+e},\qquad
 L\le\frac{v^2}{2}+\frac{2h_2v}{s+E}.
\]
Expanding $Q_m$ and using these bounds only in its positive linear terms gives
\[
 Q_m\le\frac{P/(wv)+3wv/4+h_1v/(s+e)+h_2w/(s+E)}{2H}.
\]
First $F<H^2$ and the frame imply
$1+s^2>(1+81m^2/100)^2/(1+m^2)$, hence $s>3m/4$. Consequently
\[
 Q_m<\frac{P}{2Hwv}+\frac2{33}+\frac3{800},\qquad
 \frac{P}{Hwv}>\frac{24701}{13200}>\frac{37}{20}.
\]
This improves the seam bound to $s>81m/100$. To check it on the whole sector, put $x=w/m$, $y=v/m$, $t=m^2$, $\kappa=37/20$ and $\rho=81/100$. Subtract $1+\rho^2t$ from
\[
 \frac{(1+(1-x)^2t)(1+(1-y)^2t)}{1+t-2\kappa xyt},
\]
which is a strict lower bound for $1+s^2$. After multiplication by its positive denominator, the difference is $t$ times
\[
 (1-x)^2+(1-y)^2-1+2\kappa xy-\rho^2
 +t\{(1-x)^2(1-y)^2-\rho^2+2\kappa\rho^2xy\}.
\]
On $[0,1/10]^2$ the first line is at least $9/10000$ and the bracket in the second line is nonnegative. This proves the improved bound without assuming equal widths.

Set
\[
 a_0=\frac{171m}{100},\qquad
 \gamma=\frac{9564797}{6840000},\qquad
 A_*=\gamma m^2-\frac{10}{171}.
\]
The two ratios $h_1/h_2,h_2/h_1$ are less than $25/24$. Combining the previous duration ceilings with $s+e,s+E>a_0$ and $F>h_1h_2$ gives
\[
 \frac{h_1^2L}{F},\frac{h_2^2R}{F}
 <\frac{m^2}{192}+\frac{20H}{171},\qquad A,Z>A_*.
\]
The coefficient identity is $(171/100)^2/2-1/192-10/171=\gamma$. On $m\ge13/32$, $A_*>0$.

Let $X=c-e$, $Y=f-E$, $a=s+e$, $z=s+E$, $j=\sqrt{1+s^2}$ and $u=\sqrt{az}$. The same seam gives
\[
 XY\le(j-u)^2,\qquad
 (e+d)(E+B)\ge\left(u+\sqrt{AZ/(XY)}\right)^2.
\]
Both follow by subtracting the nonnegative square in the arithmetic--geometric mean identity. The function $u+v/(j-u)$ increases with $u,v$ and decreases with $j$ on $0<u<j$, $v>0$. Hence, with $t_0=\sqrt{1+m^2}$,
\begin{equation}\label{r06:narrow-lower-product}
 (e+d)(E+B)>\left(a_0+\frac{A_*}{t_0-a_0}\right)^2.
\end{equation}
All replacement denominators are positive since $a_0<171/200<1<t_0$.

The inner roots give $b,D>9m/10$, so $(e+b)/(m+b)>18/19$ and $(E+D)/(m+D)>18/19$. Using all three durations gives the opposite bound
\[
 (e+d)(E+B)<\left(\frac{19}{18}\right)^2(1+m^2).
\]
To contradict it, put $k_0=19/18$ and $\tau=1079/1000<t_0$. The difference
$A_*-(k_0t_0-a_0)(t_0-a_0)$ has the strict lower bound
\[
 D(m)=\{\gamma-(171/100)^2-k_0\}m^2
 +(k_0+1)(171/100)\tau m-10/171-k_0.
\]
This quadratic is concave on $[13/32,1/2]$, and
\[
 D(13/32)=\frac{8509439}{11673600000}>0,\qquad
 D(1/2)=\frac{6246473}{45600000}>0.
\]
Thus $a_0+A_* /(t_0-a_0)>k_0t_0$ throughout this interval. Squaring positive quantities in \eqref{r06:narrow-lower-product} contradicts the upper product bound. The supposed double-narrow source cannot exist.
\end{proof}

\section{The upper middle bound and the complete five-dimensional root}
\begin{lemma}[High-middle exclusion]\label{r06:middle-ceiling}
Every member of $\mathcal S_{06}$ satisfies $m<471/500$.
\end{lemma}
\begin{proof}
Suppose $M\le m<1$, where $M=471/500>15/16$. Set $t=m^2$, $t_0=225/256$, $x=e^2$, $y=E^2$. The midpoint and reserve give
\[
 0<-3t^2+2t+(4t-1)(x+y)-5xy.
\]
This polynomial increases in each anchor square on $[0,1/6]^2$. It follows that
\[
 x,y>\frac{18t^2-16t+1}{24t-11}>\frac9{100}t.
\]
For the last inequality, its cleared numerator is
$N(t)=(396/25)t^2-(1501/100)t+1$, with
$N(t_0)=715/16384>0$, $N'(t_0)=10267/800>0$, and $N''=792/25>0$. Thus $e,E>3m/10$.

Let
\[
 r(m,a)=\frac{m-a+\sqrt{m^2+34ma+a^2}}6,\qquad
 \mathcal A_m(a)=(m-a)\left(r(m,a)+\frac{m+a}{2}\right).
\]
Lemma~\ref{r06:mean-root} gives $R\ge\mathcal A_m(e)$ and $L\ge\mathcal A_m(E)$. Hence
\[
 0<F-h_1h_2\le\Psi_m(e,E):=H^2-\mathcal A_m(e)\mathcal A_m(E)-h(e)h(E).
\]
For $u=a/m\in[3/10,1]$ write $r=mf(u)$ and $\mathcal A_m=m^2A(u)$. Then
\[
 f''(u)=-\frac{48}{(u^2+34u+1)^{3/2}},\quad f'(1)=\frac13,
\]
and $(A'(u)+3u/2)'\le-1/6$. Its value at $3/10$ is
$-1/12+41\sqrt{1129}/33870<0$. Thus $\mathcal A_m'(a)<-3a/2$.
Also $\mathcal A_m(a)\ge m^2-a^2$: substitute $(1+u)/2$ in the positive-root quadratic to obtain $(1-u)(1-5u)/4\le0$. Therefore
\[
 \partial_e\Psi_m\ge e\left(\frac32t-\frac56\right)>0,
\]
and likewise for $E$. The maximum over the two independent anchors is at $a_*=1/\sqrt6$.

At this corner $\partial_m\mathcal A_m(a_*)>3m/2$, giving
\[
 \frac d{dm}\Psi_m(a_*,a_*)<m(3/2-2m^2)<0.
\]
Finally $(408/1000)^2<1/6<(4083/10000)^2$, and the owner quadratic at $m=M$, $a=408/1000$, $b=143/200$ equals $-1143/1000000$. Hence $r(M,a_*)>143/200$ and
\[
 \mathcal A_M(a_*)>\frac{890333}{1200000}.
\]
Consequently
\[
 \Psi_M(a_*,a_*)<
 \left(\frac{1+M^2}{2}\right)^2
 -\left(\frac{890333}{1200000}\right)^2-\left(\frac7{12}\right)^2
 =-\frac{8035455761}{36000000000000}<0.
\]
Monotonicity excludes the whole band $M\le m<1$, including its lower endpoint. No equality of the actual anchors or means was assumed.
\end{proof}

Lemmas~\ref{r06:middle-less-than-one}, \ref{r06:midpoint}, \ref{r06:middle-floor} and \ref{r06:middle-ceiling}, together with the compact and proper-six bounds, prove Proposition~\ref{r06:source-root}. The proof order begins with strict physical multipliers and the reserve, then proves $m<1$ analytically before using the midpoint. It requires neither a coefficient proof for a unit sector nor a theorem assigning physical status to arbitrary zero-duration algebraic points.

\section{Fourteen-dimensional enclosing boxes and further necessary tests}\label{r06:enclosure-details}
The twelve smooth functions in \eqref{r06:twelve-rows} and the original cost inequality are evaluated at the same candidate. Their ambient domain is the positive orthant in the first ten variables, with four real auxiliary mean variables. The fixed-interval integral defining $M$ is smooth there: on a compact positive box its denominator is uniformly positive for every integration parameter, so every derivative is obtained by differentiating under the integral. This remains true at equal or reversed endpoints. Actual selected-root restrictions are additional source-value bounds and are never substituted as identities for ambient derivatives.

\subsection{Closed reconstruction and the eight outer fibres.}
From a five-dimensional cell put $J=G-F$ and reconstruct the positive source seam by $s=\sqrt{J/F}$. Since $G\ge1/4$ and $1/15\le s\le5$,
\[
 F\ge\frac1{104},\qquad J=\frac{Gs^2}{1+s^2}\ge\frac1{904}.
\]
These are necessary value bounds. Define
\[
 U_0=T_e(b),\quad U_3=T_E(D),\qquad
 C(a,x)=\frac{x(3x+a)}{x+3a}.
\]
The actual outer points satisfy
\begin{align*}
 d&\in[s,U_0]\cap[e,5],&B&\in[s,U_3]\cap[E,5],\\
 c&\in[e,C(e,U_0)]\cap[T_s(C(E,U_3)),T_s(E)]\cap[1/15,5],
 &f=T_s(c).
\end{align*}
Here $C$ increases with its owner and $T_s$ is a decreasing involution. The two area identities may further intersect these intervals, dividing by an exterior width or height ceiling only when the entire denominator interval is positive. Otherwise the corresponding division is omitted and the undivided equation retained.

For the full remaining $c$ interval $[l,u]$, the eight closed intervals
$[l+j(u-l)/8,l+(j+1)(u-l)/8]$, $j=0,\ldots,7$, cover it exactly, including a singleton. Each includes its complete $d,B,f,s$ and four actual-mean enclosures. Outward endpoint rounding can add overlap but cannot remove a face. A successful test on one piece does not exclude the original five-dimensional cell; every piece must be excluded.

For an inner mean put $\eta=2x/(2x+a+m)$ and $\tau=(m-a)/(2x+a+m)$. The quadratic in Lemma~\ref{r06:mean-root}, the ratio bound in Lemma~\ref{r06:midpoint} and compactness give
\[
 \frac{40}{91}\le\eta\le\frac{75}{76},\qquad0\le\tau\le\frac{11}{21}.
\]
These intervals can refine values of an actual inner mean. They are not restrictions on a rectangular path used to differentiate the ambient $M$. The directed series and whole-box derivative constructions are given in Section~\ref{reader:r06}.

\subsection{Target deficits and undivided projections.}
At a positive-width selected point, for $a\le y\le z$,
\[
 g_x(a)-g_x(y)=\frac{4x(y-a)(x^2-ay)}{(x+a)^2(x+y)^2}.
\]
By Lemma~\ref{r06:mean-root}, $x^2-ay\ge x^2-az>x(z-a)/3$. Integrating the endpoint difference with $(x+y)^2\le(x+z)^2$ gives
\[
 g_x(a)-M(x;a,z)\ge\frac{2x^2(z-a)^2}{3(x+a)^2(x+z)^2}.
\]
For an outer side let $a=e$, $x=d$, $z=c$, $U=U_0$, $p=U-s$, $w=c-e$, with the other side treated separately. A valid sum upper bound is
\[
 W=\min\{26/5,\ U+1/s,\ U+C(a,U)\}>0.
\]
Indeed $(5-x)(5-z)\ge0$ and $xz\le1$ give $x+z\le26/5$, while the seam and root cap give the other two bounds. Define the actual target upper residual
\[
 T_0^{\rm up}=\frac{h_1^2}{2}\{3+g_{U_0}(e)\}-H^2K_2-R^2K_1+2RB_m.
\]
The target equation shows $T_0^{\rm up}=h_1^2(g_{U_0}(e)-\mu_0)/2\ge0$, and the endpoint deficit gives
\[
 T_0^{\rm up}\ge\frac{h_1^2w^2}{12W^2},\qquad
 0\le A\le pw.
\]
It follows without dividing by $p$ or $A$ that
\[
 12W^2p^2T_0^{\rm up}-h_1^2A^2\ge0.
\]
If $X=\sqrt{12}W\sqrt{T_0^{\rm up}}/h_1$ and $Y$ is the corresponding right expression, the same seam gives
$(s+e+X)(s+E+Y)\ge1+s^2$. These are consequences at one tuple; their two sides cannot use independently chosen means or seams.

A useful linear form also follows. The selected-root margin increases with the owner, as in Lemma~\ref{r06:midpoint}. On the selected interval of owners from $x$ to $U$, the inequalities $g_x(a)>M\ge g_x(z)$ give $M_a<0$, $M_z\le0$. Scale invariance gives $xM_x+aM_a+zM_z=0$, so $M(U;a,z)\ge M(x;a,z)$.

At fixed $U$, divide the averaged endpoint difference by $z-a$. The result is
\[
 \frac{4U}{(U+a)^2}\int_0^1 t\,
 \frac{U^2-a(a+t(z-a))}{(U+a+t(z-a))^2}\,dt .
\]
On $a\le y\le z_*=C(a,U)$ the numerator $U^2-ay$ decreases to $U(U^2-a^2)/(U+3a)>0$, while the denominator $(U+y)^2$ increases. The displayed ratio is therefore bounded below by its value at $z=z_*$. At that point set $q=(U+z_*)/(U+a)=4U/(U+3a)>1$ and $\eta=U/(U+a)=3q/[2(q+2)]$. The mean primitive in Lemma~\ref{r06:mean-root} gives the deficit
$4\eta[\log q/(q-1)-1+\eta(q-1)/q]$. After division by $z_*-a$, comparison with the coefficient in the next display reduces to
\[
 \log q>\frac{(q-1)(7-q)}{2(q+2)}+\frac{(q-1)^3}{6q^2}.
\]
The difference vanishes at one and its derivative is
$(q-1)^3(q^2+10q+4)/[3q^3(q+2)^2]>0$. Combining the two comparisons yields
\[
 g_U(a)-M(x;a,z)\ge
 \frac{(U-a)(U+3a)}{8U(U+a)^2}(z-a).
\]
Thus, with $k=h_1^2(U_0-e)(U_0+3e)/(16U_0(U_0+e)^2)$,
\[
 T_0^{\rm up}\ge kw,\qquad pT_0^{\rm up}-kA\ge0.
\]
The joint-seam version is obtained by multiplying the two inequalities
$k_0(s+c)\le k_0(s+e)+T_0^{\rm up}$ and its right analogue, then using $(s+c)(s+f)=1+s^2$. The original target inequalities are retained when $p$ or $k$ is zero. All these cuts follow by multiplying nonnegative quantities and do not create positive collars.

The lower-seam projection is also undivided. With $p_0=U_0-s$, $p_3=U_3-s$, $v_0=A+ep_0$, $v_3=Z+Ep_3$, it is
\[
 p_0p_3-v_0v_3-s(v_0p_3+v_3p_0)\ge0.
\]
Its difference from $p_0p_3[1-cf-s(c+f)]$ equals
\[
 (p_0c-v_0)p_3(f+s)+(p_3f-v_3)(v_0+sp_0)\ge0.
\]
This proves necessity also on $p_0=0$ or $p_3=0$, where division would lose a face.

\subsection{The original physical stationary rows.}
They are pointwise necessary rows on the strict source, separate from the twelve ambient functions. Put $q=(e,E,b,D,m,c)$ and differentiate the original $F=H^2-RL$ and $G$ before any interval intersection. On the actual frame graph,
\[
 s_j=\frac{G_j-(1+s^2)F_j}{2sF},\qquad
 f_j=-\frac{(1+c^2)s_j+(1+s^2)\delta_{jc}}{(c+s)^2}.
\]
Let $g_i$ denote the native partial derivatives of the original nine-edge cost, with its nine reciprocal variables independent, and set $w=c-e$, $v=f-E$. The exact mixed-derivative identity $\phi_{xy}=(3+g)/2$ gives
\[
 g_d=wK_0,\quad g_b=(m-e)K_1,\quad
 g_D=(m-E)K_2,\quad g_B=vK_3.
\]
For the ordered coordinates $q=(e,E,b,D,m,c)$ define
$g_{\rm retained}=(g_e,g_E,g_b,g_D,g_m,g_c)$. All subscripts $j$ on $s,f,A,Z,w,v$ below are ordinary partial derivatives of the actual six-coordinate physical graph; in particular $w_j=\delta_{jc}-\delta_{je}$ and $v_j=f_j-\delta_{jE}$. Differentiating its actual chart $d=s+A/w$, $B=s+Z/v$ cancels these widths and yields the six regular necessary values
\begin{equation}\label{r06:physical-six}
 E_j=g_{\rm retained}(j)+g_f f_j
 +K_0\{w s_j+A_j-(d-s)w_j\}
 +K_3\{v s_j+Z_j-(B-s)v_j\}=0.
\end{equation}
The full signed-frame formulas reconstruct the same endpoints on an open neighbourhood of the strict point. Positive durations and positive state coordinates persist, and the endpoint-fixed homotopy preserves the real winding. Lemma~\ref{fw:balanced-reconstruction} gives the genuine geometric competitors needed for these stationarity equations. The derivatives in \eqref{r06:physical-six} are of the original cost composed with the physical chart. The whole-box contractions use the ambient Jacobians and Hessians of \eqref{r06:twelve-rows}.

Finally, after the complete closed-cover exclusion has established $\mathcal S_{06}=\varnothing$, its actual closure is empty. None of the preceding source-to-root arguments assumes that conclusion or asserts that every weak algebraic tuple is physically realizable. In particular the closed numerical enlargements do not supply a maximizing lift, real winding or boundary variational freedom to a tuple that did not have them.
\space
\label{revision:appendix-end-9}
\end{proofappendix}
\section{The four-dimensional cost bound for R04}\label{app:r04-source}
\begin{proofappendix}
\section{The physical frame and the fixed-area variations}\label{r04:frame-reduction}
We first justify the reduction used in Section~\ref{r04:section} on the original physical class. This also explains why a mean equation is necessary at a minimizer, rather than an additional definition of that class.
Put $J_1=B+C$, $J_2=C+D$ and write
\[
 Q_1=C+J_1\alpha,\quad R_1=B+J_1\beta,\quad
 Q_2=D+J_2\gamma,\quad R_2=C+J_2\theta.
\]
All four shape parameters are positive. With $p,q,M$ as in \eqref{r04:four-coordinates}, put
\[
 y=\frac{1+B^2}{(A+B)J_1},\qquad
 u=\frac{1+D^2}{J_2(D+E)},\qquad
 k_1=\alpha\beta,\quad k_2=\gamma\theta.
\]
In the normalization of Section~\ref{fd:section}, the three constant matrices are
\[
 B_0=\begin{pmatrix}1/2&1\\0&-1/2\end{pmatrix},\quad
 B_1=\begin{pmatrix}0&-1/2\\-1/2&0\end{pmatrix},\quad
 B_2=\begin{pmatrix}-1/2&1\\0&1/2\end{pmatrix}.
\]
With the signed twist $R_*$ from that section, define
\[
 Q_j(v)=\frac{v+v^{-1}}2I-(v-v^{-1})B_j,
\]
\[
 \mathcal K_1(y;k)=\frac1y
 \begin{pmatrix}y^2+k&1-k-y^2\\k&1-k\end{pmatrix},
 \qquad
 \mathcal K_0(k)=\begin{pmatrix}1-k&-k\\k&1+k\end{pmatrix}.
\]
Substitution of the nine multipliers \eqref{r04:durations} into the product in native label order gives the finite rational identity
\begin{equation}\label{r04:compressed-product}
 P:=\prod_{i=0}^8 Q_{j_i}(q_i)
 =Q_0(p)\mathcal K_1(y;k_1)Q_1\left(\frac{1+C^2}{J_1J_2}\right)
       \mathcal K_0(k_2)Q_0(u)Q_2(q).
\end{equation}
Multiplying the displayed two-by-two matrices gives $P=R_*$ precisely at the three values
\begin{equation}\label{r04:frame-three}
 \frac{1+C^2}{J_1J_2}=\frac{pqyu}{M},\qquad
 2Mk_1=M-y^2(q^2-1),\qquad
 2Mk_2=M-u^2(p^2-1).
\end{equation}
For necessity put $T=(1+C^2)/(J_1J_2)$, $\xi=yT/(pqu)$, and
$W=2yT\mathcal K_1(y;k_1)Q_1(T)\mathcal K_0(k_2)$. Moving the exterior matrices to the other side of $P=R_*$ gives
\[
 W=\xi\begin{pmatrix}
 u^2(2p^2-1)&2p^2u^2-4p^2-2q^2-u^2+2\\
 u^2&2q^2+u^2-2
 \end{pmatrix}.
\]
With indices $1,2$, the contrast $W_{11}-W_{12}-W_{21}+W_{22}$ is $4y^2$ in the block product and $4M\xi$ on the right. Thus $\xi=y^2/M$ and $T=pqyu/M$. The differences $W_{11}-W_{21}$ and $W_{21}-W_{22}$ give
$y^2(1-2k_2)=\xi u^2(p^2-1)$ and $1-2k_1=\xi(q^2-1)$, respectively. Substitution of \eqref{r04:frame-three} into \eqref{r04:compressed-product} gives the four entries $(1/2,-3/2,1/2,1/2)$, proving sufficiency. Multiplication by $J_1J_2$ and $J_i^2$, respectively, gives \eqref{r04:four-coordinates} and $I_i=J_i^2k_i$. Thus the signed product is independent of the individual shape parameters when the two products $k_1,k_2$ are fixed.

For completeness this canonical product represents the full physical frame, not merely its trace. Write $X$ for the traceless state matrix; then $X^2=-I$, $B_j^2=I/4$ and $\tr(XB_j)=r_i$ on arc $i$. The state flow gives $X'= [B_j,X]/r_i$. The identity
\[
 XB_j+B_jX=\tr(XB_j)I
\]
shows that the physical transition on that arc is $X_{\rm in}Q_j(q_i)X_{\rm out}^{-1}$. Consecutive factors telescope. The genuine state twist consequently gives
\[
 G_9=X_0PR_*^{-1}X_0^{-1}R_*,
 \qquad G_9=R_*\quad\Longleftrightarrow\quad P=R_*.
\]
Both implications retain the sign of the endpoint matrix.

At an actual strict source point, vary $(A,B,D,E,\alpha,\theta)$ in a small connected open ball. Compute $C,k_1,k_2$ from \eqref{r04:four-coordinates} and \eqref{r04:frame-three}, then set $\beta=k_1/\alpha$ and $\gamma=k_2/\theta$. The positive root and nine positive durations persist by openness, as does the positivity of the third reciprocal $(1-r_ir_{i+1})/(r_i+r_{i+1})$ at every switch. The exact mode flow joins these states and keeps all three reciprocals positive along each arc. The label extension preserves the original state twist, and \eqref{r04:compressed-product} preserves the full frame endpoint. The resulting fixed-endpoint homotopy preserves the real winding. Lemma~\ref{fw:balanced-reconstruction} gives genuine Reinhardt competitors with variable entrance state and free durations. The nearby competitors need only the geometric admissibility conditions.

Fixing the four exterior variables now fixes $C,I_1,I_2$. Varying $\alpha$ varies $X_1>0$ with $Y_1=I_1/X_1$; varying $\theta$ does the same for the second block. The potential identity \eqref{r04:cost-bridge}, or its uncropped rectangle form before minimization, gives the derivative of the original cost in these directions. Fermat's theorem forces each derivative to vanish. This argument requires neither a rank assumption on a frame Jacobian nor a new maximizing lift at nearby competitors. It proves the required stationary-shape conditions for every original strict minimizer.

\section{Existence and uniqueness of the positive-area shape}\label{r04:shape-section}
For positive $x,a,z$ define the fixed-interval mean
\[
 m(x;a,z)=\int_0^1 g(x,(1-t)a+tz)\,dt,
 \qquad g(x,y)=\left(\frac{x-y}{x+y}\right)^2.
\]
This definition is smooth also when $a=z$. Since $k=(3+g)/2$, equality of the two far-edge $k$ means is equivalent to equality of the corresponding $g$ means. For the rectangle in \eqref{r04:rectangle} put $Q=b+X$, $R=a+Y$ and
\[
 \mu_X=\frac1Y\int_a^{a+Y} k(u,b+X)\,du,\qquad
 \mu_Y=\frac1X\int_b^{b+X} k(a+Y,v)\,dv.
\]
Since $k=(3+g)/2$, the equation $\mu_X=\mu_Y$ is exactly
$m(Q;a,R)=m(R;b,Q)$. Differentiating an integral on fixed unit intervals gives
\begin{equation}\label{r04:shape-derivative}
 \frac{d}{dt}\Gamma(a,b;e^t,Ie^{-t})=I(\mu_X-\mu_Y).
\end{equation}

\begin{lemma}[The positive-area shape]\label{r04:shape-unique}
For every $a,b,I>0$ the far-edge mean equation has exactly one solution $X>0$, $Y=I/X$. Its derivative difference changes from negative to positive at that solution, which is the unique global minimum of $\Gamma$. The minimizing widths are smooth functions of $a,b,I$ on their positive domain.
\end{lemma}
\begin{proof}
Write $Q=b+x$, $R=a+I/x$ and
$f(x)=m(Q;a,R)-m(R;b,Q)$. At a zero let the common mean be $\mu$. If $Q\ge R$, the first mean exceeds $g(Q,R)$ since $a<R\le Q$ and $g(Q,y)$ decreases for $y<Q$. If $b\ge R$, the second mean would be smaller than $g(R,Q)$, a contradiction. The case $R\ge Q$ is symmetric. Thus at every zero
\[
 a,b<\min(Q,R),\qquad d:=\mu-g(Q,R)>0.
\]
For $z>a$ integration by parts gives
\[
 (x-a)m_x(x;a,z)-[m(x;a,z)-g(x,z)]
 =\frac{4a}{z-a}\int_a^z\frac{(x-y)^2}{(x+y)^3}\,dy>0.
\]
Consequently $u=m_x(Q;a,R)>0$ and $v=m_x(R;b,Q)>0$ at a zero. Endpoint differentiation, with $r=R-a=I/x$, yields
\[
 f_Q=u+\frac d{x}>0,\qquad f_R=-v-\frac d r<0,
 \qquad f'(x)=u+\frac r x v+\frac{2d}x>0.
\]
These are signs at zeros, not a monotonicity assertion away from zeros. Each zero is a strict negative-to-positive crossing; two zeros would force an intervening zero with nonpositive left derivative. Hence there is at most one.

Existence follows from explicit positive brackets. The identity
\[
 1-m(x;a,z)=\frac{4x}{z-a}\int_a^z\frac{y}{(x+y)^2}\,dy
\]
implies, for $R\ge256\max(a,Q)$,
\[
 1-m(Q;a,R)\ge\frac Q{R-a}\log\frac R{\max(a,Q)}>\frac{4Q}R,
 \qquad 1-m(R;b,Q)\le\frac{4Q}R.
\]
We used $\log256=8\log2>4$. Therefore $f<0$ there. Interchanging the two coordinates gives $f>0$ when $Q\ge256\max(b,R)$. In particular
\[
 x_{\rm lo}=\min\left\{1,\frac I{256(a+b+1)}\right\},\qquad
 x_{\rm hi}=256(a+b+I+1)
\]
satisfy $f(x_{\rm lo})<0<f(x_{\rm hi})$. Continuity proves existence. Uniqueness and the crossing sign make $f$ negative before its root and positive afterwards. Equation~\eqref{r04:shape-derivative} proves the global minimum claim. The positive derivative at the root permits the implicit function theorem, proving smooth dependence. Nothing here assigns a unique shape when $I=0$.
\end{proof}

The same comparison proves that both far endpoints exceed both base endpoints at an optimal shape. We shall use this strict crossing property below. Also $\Gamma\to2I$ at either end of the logarithmic shape parameter, whereas $3I/2<\Gamma<2I$ at every finite positive rectangle.

\section{A sharp rectangle comparison}\label{r04:rectangle-section}
\begin{lemma}\label{r04:rectangle-comparison}
Let $S=a+b$, $I=XY$ and $L=\sqrt I$. Then
\[
 \Gamma(a,b;X,Y)\ge\Gamma(S/2,S/2;L,L).
\]
Equality holds exactly when $a=b$ and $X=Y$. Consequently, with
$\mathsf S(S,I)=\Gamma(S/2,S/2;\sqrt I,\sqrt I)$,
$\mathcal E(a,b,I)\ge\mathsf S(a+b,I)$.
\end{lemma}
\begin{proof}
For $0<\rho<1$ let $K_\rho=\{(u,v)>0:\rho u\le v\le u/\rho\}$. For $R=[a,a+Y]\times[b,b+X]$ we prove
\begin{equation}\label{r04:cone-bound}
 |R\smallsetminus K_\rho|\ge\rho\left(L-\frac{S(1-\rho)}{2\rho}\right)_+^2.
\end{equation}
Here $t_+=\max(t,0)$. Assume $X\ge Y$, and put
\[
 c=\frac{S(1-\rho)}\rho,\quad T=X+Y-c,\quad
 h_1=X+b-a/\rho,\quad h_2=Y+a-b/\rho.
\]
Then $h_1+h_2=T$. Define $F_1(h)$ as the area cut from $[0,Y]\times[0,X]$ by $\eta+\xi/\rho<h$, and let $F_2$ exchange $X,Y$. Translation and reflection of the two disjoint regions outside the cone give $|R\smallsetminus K_\rho|=F_1(h_1)+F_2(h_2)$. Their derivatives are
\[
 F_1'(h)=\bigl[\min(Y,\rho h)-\max(0,\rho(h-X))\bigr]_+,
\]
with the exchanged formula for $F_2$. If $T\le0$ the desired bound is zero. Otherwise monotonicity restricts the minimum of $Q(h)=F_1(h)+F_2(T-h)$ to $0\le h\le T$. For $Y\le h\le X$, $F_1'(h)\ge\rho Y\ge F_2'(T-h)$. For $h\ge X$, set $v=T-h\le Y$; either $F_1'(h)=\rho X\ge\rho v=F_2'(v)$, or $Q'(h)=(1-\rho)Y+\rho c>0$. Hence a minimum occurs with $h\le Y$. There $F_1'(h)=\rho h$, while
\[
 Q'(h)=\begin{cases}
 \rho(2h-T),&T-h\le Y,\\
 \rho(h-Y),&Y\le T-h\le X/\rho,\\
 -(1-\rho)X-\rho c,&T-h\ge X/\rho.
 \end{cases}
\]
Integrating these piecewise affine derivatives gives
\begin{equation}\label{r04:cone-minimum}
 \min Q=\begin{cases}
 0,&c\ge X+Y,\\
 \rho(X+Y-c)^2/4,&X-Y\le c\le X+Y,\\
 \rho Y(X-c),&0\le c\le X-Y.
 \end{cases}
\end{equation}
The middle line dominates \eqref{r04:cone-bound} by $X+Y\ge2L$. In the last line only $c<2L$ needs checking, and the difference divided by $\rho$ is $c(L-Y-c/4)\ge0$. Indeed if $Y\le L/2$, then $c<2L\le4(L-Y)$; otherwise $c\le X-Y=(L-Y)(L+Y)/Y\le3(L-Y)$.

Writing $r=\min(u/v,v/u)$, the kernel is $3/2+(1-r)^2/[2(1+r)^2]$. Integration of its nonnegative level sets gives
\begin{equation}\label{r04:cone-integral}
 \Gamma(R)=\frac32 I+
 \int_0^1\frac{2(1-\rho)}{(1+\rho)^3}|R\smallsetminus K_\rho|\,d\rho.
\end{equation}
The right side of \eqref{r04:cone-bound} is precisely the cone-exterior area of the diagonal square. This proves the inequality. If $X>Y$, the last comparison is strict for a positive interval of small $c>0$. If $X=Y$ but $a\ne b$, then $h_1\ne h_2$ and the unique minimum in the middle line occurs at $h_1=h_2=T/2$, again giving strictness on an interval. These observations prove the equality statement.
\end{proof}

Set $\mathsf S(S,0)=0$ for $S>0$, the continuous extension by the bound $\mathsf S(S,I)\le 2I$. Direct integration gives, for $c=S/2$,
\begin{equation}\label{r04:square-formula}
 \mathsf S(S,I)=3I-2(c+\sqrt I)^2\log\frac{2(c+\sqrt I)}{2c+\sqrt I}
 +2c^2\log\frac{2c+\sqrt I}{2c}.
\end{equation}
The identities $a_0^2+2a_0b_0=1$ and $(a_0+b_0)^2=1+b_0^2=2a_0^2$ now verify the boundary value used in the main text:
\[
 G(z_0)=5+3a_0^2-18a_0b_0-2a_0^2\log2
 =2+\frac{6-(6+4\sqrt2)\log2}{5+4\sqrt2}=m_0.
\]

\section{The square series and its derivative remainders}\label{r04:square-series}
For $S,I>0$ put $\sigma=\sqrt I/(S+\sqrt I)$ and $w=\sigma^2$. Then
\begin{equation}\label{r04:positive-series}
 \mathsf S(S,I)=I\left\{\frac32+
 \sum_{m=1}^{\infty}\frac{w^m}{m(2m+1)(2m+2)}\right\}.
\end{equation}
To prove it, substitute $u=S/2+\sqrt I\xi$, $v=S/2+\sqrt I\eta$ in the integral. The excess over $3I/2$, divided by $I$, is
\[
 \frac{\sigma^2}{2}\int_0^1\int_0^1
 \frac{(\xi-\eta)^2}{[1+\sigma(\xi+\eta-1)]^2}\,d\xi\,d\eta.
\]
Expand $(1+\sigma t)^{-2}$ uniformly for $|t|\le1$. Reflection about $(1/2,1/2)$ cancels the odd moments. With $t=\xi+\eta-1$, $z=\xi-\eta$, the remaining moments are
\[
 \int_0^1\int_0^1(\xi-\eta)^2(\xi+\eta-1)^{2n}\,d\xi\,d\eta
 =\frac23\int_0^1t^{2n}(1-t)^3\,dt
 =\frac4{(2n+1)(2n+2)(2n+3)(2n+4)}.
\]
This gives \eqref{r04:positive-series}. If $I\le S^2/2$, then $w\le(\sqrt2-1)^2$. Write $H_N(w)=\sum_{m=1}^N w^m/[m(2m+1)(2m+2)]$, $H=\lim H_N$ and $a_N=[(N+1)(2N+3)(2N+4)]^{-1}$. For $N\ge3$,
\begin{equation}\label{r04:series-remainders}
\begin{aligned}
 0\le H-H_N&\le a_Nw^{N+1}/(1-w),\\
 0\le H'-H_N'&\le(N+1)a_Nw^N/(1-w),\\
 0\le H''-H_N''&\le N(N+1)a_Nw^{N-1}/(1-w),\\
 0\le H'''-H_N'''&\le w^{N-2}/[4(1-w)].
\end{aligned}
\end{equation}
The first three bounds follow from decreasing coefficients; for the third use that $(m-1)/[(2m+1)(2m+2)]$ decreases for $m\ge3$. The last uses $(m-1)(m-2)/[(2m+1)(2m+2)]<1/4$. Thus a finite computation encloses the full function and its first three derivatives.

\section{Global and conditional shape bounds}\label{r04:shape-bounds}
\begin{lemma}\label{r04:shape-three-halves}
For the unique minimizing shape, if $a>b>0$ then
$X-Y>3(a-b)/2$. If in addition $I\le(a+b)^2/2$, then $X-Y<3(a-b)$. Exchanging the coordinates gives the corresponding bounds for $a<b$; if $a=b$ then $X=Y$.
\end{lemma}
\begin{proof}
Set $\delta=a-b>0$. First test the shape $XY=I$, $X-Y=3\delta/2$. Put $U=a+Y$, so $b+X=U+\delta/2$. Substitution in the two far-edge means gives
\[
 2(\mu_Y-3/2)=\int_{-Y-\delta}^{\delta/2}
 \frac{z^2}{X(2U+z)^2}\,dz,\qquad
 2(\mu_X-3/2)=\int_{\delta/2-Y}^{\delta/2}
 \frac{(\delta-z)^2}{Y(2U+z)^2}\,dz.
\]
The numerator measures have equal mass $Y^2/3+Y\delta/2+\delta^2/4$. Before the common interval only the first has mass; on the common interval its density is strictly smaller since $z^2\le(\delta-z)^2$ and $X>Y$. Their cumulative mass difference is therefore nonnegative and positive on an interval. The weight $(2U+z)^{-2}$ is strictly decreasing. Integration by parts proves $\mu_Y>\mu_X$. The unique crossing in Lemma~\ref{r04:shape-unique} puts the optimum at larger $X$, proving the lower bound.

For the upper bound test $X-Y=3\delta$, and normalize by
\[
 \begin{gathered}
 T=a+b+X+Y,\qquad d=\delta/T,\qquad y=Y/T,\qquad x=X/T=y+3d,\\
 u=(a+Y)/T=\tfrac12-d,\qquad v=\tfrac12+d.
 \end{gathered}
\]
Positivity and the area upper bound give $0<d<(1-2y)/4$, $0<y\le1/(2+\sqrt2)<3/10$ and $x<3/4$. Put $g_0(z)=-\log(1-z)/z$. Direct integration yields
\begin{equation}\label{r04:shape-trial-three}
 \frac{\mu_X-\mu_Y}{2d}
 =3u\frac{g_0(x)-g_0(y)}{x-y}
 +2\left(\frac1{1-y}-g_0(y)\right)
 -\frac{3u^2}{(1-x)(1-y)}.
\end{equation}
Jensen's inequality, with probability densities $2t$ and $2(1-t)$ respectively, gives
\[
 \frac{g_0(x)-g_0(y)}{x-y}\ge
 \frac1{2(1-2x/3)(1-2y/3)},\qquad
 2\left(\frac1{1-y}-g_0(y)\right)\ge\frac y{(1-y)(1-y/3)}.
\]
The resulting lower bound in \eqref{r04:shape-trial-three} has positive denominator
$4(3-y)(1-y)(3-2y)(1-3d-y)(3-6d-2y)$. On putting $d=(1-2y)t/4$, its numerator is
\[
 q_0(1-t)^3+3q_1t(1-t)^2+3q_2t^2(1-t)+q_3t^3,
\]
where
\begin{align*}
 q_0&=3y(2-y)(16y^2-27y+9),\\
 q_1&=(27-18y-141y^2+104y^3-4y^4)/4,\\
 q_2&=(27+18y-285y^2+268y^3-132y^4+32y^5)/8,\\
 q_3&=3y(2y+1)(9-33y+22y^2-8y^3)/8.
\end{align*}
All are positive for $0<y\le3/10$: the nontrivial factors in $q_0,q_3$ are at least $117/50,108/125$, and dropping positive terms gives $4q_1\ge11097/1250$, $8q_2\ge351/1250$. Thus the trial mean difference is positive. The unique minimum lies at smaller $X$, proving the upper bound.
\end{proof}

\begin{proposition}[An asymmetry term]\label{r04:asymmetry}
For all $a,b,I>0$,
\[
 \mathcal E(a,b,I)\ge
 \mathsf S\left(a+b,I+\frac9{16}(a-b)^2\right)
 -\mathsf S\left(a+b,\frac9{16}(a-b)^2\right).
\]
Equality holds exactly when $a=b$.
\end{proposition}
\begin{proof}
Put $S=a+b$, $d=|X-Y|$, $h=X+Y=\sqrt{d^2+4I}$. The piecewise minimum \eqref{r04:cone-minimum} gives the stronger cone bound
\[
 |R\smallsetminus K_\rho|\ge\frac\rho4[(h-c)_+^2-(d-c)_+^2],
 \qquad c=S(1-\rho)/\rho.
\]
It is the difference of the exterior areas of diagonal squares of side $h/2$ and $d/2$. Equation~\eqref{r04:cone-integral} therefore gives
$\Gamma\ge\mathsf S(S,I+d^2/4)-\mathsf S(S,d^2/4)$. For fixed $S$, the function $J\mapsto\mathsf S(S,J)$ is strictly convex: each cone contribution $\rho(\sqrt J-c/2)_+^2$ is convex, and where $J>c^2/4$ its second derivative is $\rho c/(4J^{3/2})>0$. Hence $T\mapsto\mathsf S(S,I+T)-\mathsf S(S,T)$ strictly increases. Lemma~\ref{r04:shape-three-halves} gives $d^2/4>9(a-b)^2/16$ when $a\ne b$, proving strictness. The symmetric case is the square.
\end{proof}

\begin{lemma}[Conditional shape bounds]\label{r04:conditional-shape}
Suppose $a>b>0$ and $0<I\le(a+b)^2/2$. If $(a-b)/(a+b)\le1/2$, then $X-Y<2(a-b)$. If also
\[
 \frac{a-b}{a+b}\le\frac3{10},\qquad
 \frac15\le\frac I{(a+b)^2},
\]
then $X-Y>17(a-b)/10$.
\end{lemma}
\begin{proof}
We give the finite rational coefficient checks, including the polynomials they certify. By homogeneity set $a+b=1$, $k=a-b$. At the trial shape $X-Y=ck$, put $Y=v$, $X=v+ck$ and
\[
 T=1+2v+ck,\quad d=k/T,\quad y=v/T,\quad x=(v+ck)/T,
 \quad u=\frac{1-(c-1)d}{2}.
\]
The integrated mean difference is
\begin{equation}\label{r04:general-shape-trial}
 \frac{\mu_X-\mu_Y}{2d}
 =cu\frac{g_0(x)-g_0(y)}{x-y}
 +(c-1)\left(\frac1{1-y}-g_0(y)\right)
 -\frac{cu^2}{(1-x)(1-y)},
\end{equation}
where $g_0(z)=\sum_{n\ge0}z^n/(n+1)$.
For $c=2$, let
\[
 F_6=2u\sum_{n=1}^6\frac{\sum_{j=0}^{n-1}x^jy^{n-1-j}}{n+1}
 +\sum_{n=1}^6\frac n{n+1}y^n
 -\frac{2u^2}{(1-x)(1-y)},
\]
\[
 P_+(k,v)=420(v+1)(2k+v+1)(2k+2v+1)^6F_6.
\]
This is an explicit polynomial of total degree eight. All omitted series terms are positive, so it suffices to show $P_+>0$ on $0\le k\le1/2$, $v>0$, $H=1/2-2kv-v^2\ge0$. Here $v\le1/\sqrt2<71/100$.

For a polynomial on $[a,b]\times[c,d]$, substitute $k=a+(b-a)t$, $v=c+(d-c)r$ and write the result as $\sum q_{ih}t^ir^h$. Its tensor Bernstein coefficients of bidegree $(m,n)$ are
\begin{equation}\label{r04:bernstein}
 B_{j\ell}=\sum_{i\le j,\ h\le\ell}
 q_{ih}\frac{\binom ji}{\binom mi}\frac{\binom\ell h}{\binom nh}.
\end{equation}
The \path{R04-guards} task starts from
\([0,1/2]\times[0,71/100]\), subdivides into closed rectangles
and checks \(H<0\) or \(P_+>0\).
On the initial corner's terminal rectangle the coefficient of
the corner \((0,0)\) may be zero. The other Bernstein coefficients
are positive, so \(v>0\) still gives a positive value.
Thus the trial mean difference is positive on the required region,
and the unique crossing gives \(X-Y<2(a-b)\).

For the lower bound set $c=17/10$ and define
\begin{align*}
 U(k,v)={}&cu\left(\sum_{n=1}^8
 \frac{\sum_{j=0}^{n-1}x^jy^{n-1-j}}{n+1}+\frac{x^8}{1-x}\right)\\
 &+(c-1)\left(\sum_{n=1}^8\frac n{n+1}y^n+\frac{y^9}{1-y}\right)
 -\frac{cu^2}{(1-x)(1-y)},\\
 P_-(k,v)={}&-5040(v+1)(17k+10v+10)(17k+20v+10)^8U(k,v).
\end{align*}
The geometric tails bound the omitted positive sums in \eqref{r04:general-shape-trial} from above, so the trial difference divided by $2d$ is at most $U$. The displayed expression for $P_-$ is a polynomial of total degree ten. The required region has \(0\le k\le3/10\),
\(0<v<71/100\) and \(1/5\le v^2+17kv/10\le1/2\).
The second shape check of \path{R04-guards} starts from
\([0,3/10]\times[0,71/100]\). Each closed terminal either
contradicts one of the area bounds or satisfies \(P_->0\), by
\eqref{r04:bernstein}. Thus \(U<0\) on the required region,
giving \(X-Y>17(a-b)/10\).
The exact polynomials, subdivision and sign checks for both trials
are supplied under \path{mathematical_inputs/r04/original/global_bound/}.
\end{proof}

\section{An explicit bound for the exterior tail}\label{r04:tail-section}
\begin{proposition}\label{r04:tail-bound}
If $z\in\Omega_{04}$ and $\max(A,E)\ge5$, then
$C_5>54433/30000>m_0$ and hence $G>m_0$.
\end{proposition}
\begin{proof}
The reflection $(A,B,D,E)\mapsto(E,D,B,A)$ preserves the relaxed domain and its cost; its precise action is recorded below. Assume $T=A\ge5$ and put
\[
 b=TB,\ c=TC,\ d=TD,\ e=TE,\ t=b+e,\ s=1-t,\ \ell=\log T.
\]
The inequalities $p,q>1$ give $b+e+be/T^2<1$ and $d+e+de/T^2<1$. Thus $0\le s\le1$ and $d\le b+s$. The first side condition says $C+X_1<x=(1-AB)/(A+B)$. Lemma~\ref{r04:shape-three-halves} gives $C+X_1>B$, whether $B\le C$ or $B>C$. Therefore $B^2+2AB<1$, $b<1/2$ and
\[
 c<Tx=\frac{1-b}{1+b/T^2}\le1-b.
\]
In particular $B,C,D,E\le1/T$.

The terms containing $A^2$, together with the constant and $-3A(B+E)$, are at least
\[
 F_T(t)=5-3t+T^2\log\frac{1+t/T^2}{1+1/T^2}\ge2+2s,
\]
because $F_T(1)=2$ and $F_T'\le-2$ for $0\le t\le1$. The remaining exterior logarithms contribute at least $(b^2+e^2)\ell/T^2$. The three internal logarithms have combined $\ell$ coefficient $-(b^2+2c^2+2d^2+e^2)/T^2$, so the net coefficient is exactly $-2(c^2+d^2)/T^2$. For $h,k>0$,
\[
 (h^2+k^2)\log(h+k)\ge-1/(2\exp(1))>-3/16.
\]
For $h+k<1$ this follows by minimizing $r^2\log r$; for $h+k\ge1$ it is immediate. Also $\sum_{z\in\{B,C,D,E\}}z^2\log(1+z^2)\le4/T^4$. Consequently
\begin{equation}\label{r04:tail-base}
 C_5\ge2+2s-
 \frac{3(bc+cd+de)+2(c^2+d^2)\ell+9/16}{T^2}-\frac4{T^4}.
\end{equation}
Using $c\le1-b$, $d\le b+s$ and $e=1-b-s$ gives
\begin{align*}
 bc+cd+de&\le3b(1-b)+s(2-3b)-s^2,\\
 c^2+d^2&\le1-2b+2b^2+2bs+s^2.
\end{align*}
In the first inequality, after replacing $c$ the coefficient of $d$ is $2(1-b)-s\ge0$, justifying the second replacement. Hence, since $\ell>3/2$, $s\le1$ and $b\le1/2$,
\begin{align*}
 3(bc+cd+de)+2(c^2+d^2)\ell
 &\le2\ell+b(1-b)(9-4\ell)+s(6-9b+4b\ell)+s^2(2\ell-3)\\
 &\le\max\{2\ell,9/4+\ell\}+s(3+4\ell).
\end{align*}
The function $(3+4\log T)/T^2$ decreases for $T\ge5$ and is less than $29/75<2$, using $\log5<5/3$. Thus the last slack is absorbed by $2s$ in \eqref{r04:tail-base}. Both branches in
\[
 \Delta(T)=\frac{9/16+\max\{2\log T,9/4+\log T\}}{T^2}+\frac4{T^4}
\]
decrease, and
\[
 C_5\ge2-\Delta(T),\qquad
 \Delta(T)\le\Delta(5)<\frac{215}{1200}+\frac4{625}=\frac{5567}{30000}.
\]
The elementary constants follow from
\[
 \log5>2\left(\frac23+\frac{(2/3)^3}{3}\right)=\frac{124}{81}>\frac32,
 \quad \exp(5/3)>\sum_{k=0}^4\frac{(5/3)^k}{k!}=\frac{10009}{1944}>5.
\]
Finally $\log2>56/81$, from its positive series, and $\sqrt2<283/200$ give
\[
 m_0<\frac{-1408+1720\sqrt2}{567}<\frac{5129}{2835}<\frac{54433}{30000}.
\]
Both rectangle costs are positive, completing the proof.
\end{proof}

\section{Exact polynomial tests for the two areas}\label{r04:area-section}
Use the following polynomials in the four exterior coordinates:
\begin{align*}
 U&=A+B,&V&=D+E,&T&=A+E,\\
 P_0&=1+A^2,&Q_0&=1+E^2,&R_B&=1+B^2,&R_D&=1+D^2,\\
 N&=P_0^2V^2+Q_0^2U^2-U^2V^2T^2,
 &\mathsf F&=P_0Q_0R_BR_D-N.
\end{align*}
Then
\begin{equation}\label{r04:area-rational}
\begin{gathered}
 C^2=\mathsf F/N,\\
 2I_1=(B+C)^2-\frac{R_B^2(Q_0^2-V^2T^2)}N,\qquad
 2I_2=(C+D)^2-\frac{R_D^2(P_0^2-U^2T^2)}N.
\end{gathered}
\end{equation}
On the strict domain $N>0$, since $P_0>UT$ and $N=(P_0^2-U^2T^2)V^2+Q_0^2U^2$. Define
\[
 K_1=[A(B+E)+BE]^2+B^2E^2+B^2+E^2,
\]
\[
 \mathcal Q_1=TK_1V^2+P_0Q_0(B^2-E^2)V-2BEUQ_0^2,
 \quad \mathcal Q_2(A,B,D,E)=\mathcal Q_1(E,D,B,A).
\]
\begin{proposition}\label{r04:area-sign}
For positive exterior variables with $p,q>1$ and $\mathsf F>0$,
\[
 \operatorname{sgn}I_i=\operatorname{sgn}\mathcal Q_i\quad(i=1,2),
 \qquad I_1<2BC,\quad I_2<2CD.
\]
\end{proposition}
\begin{proof}
For the first area set $\mathsf J_1=R_B^2(Q_0^2-V^2T^2)-B^2N-\mathsf F$. Exact expansion gives
\begin{equation}\label{r04:area-square}
 2NI_1=2B\sqrt{\mathsf FN}-\mathsf J_1,\qquad
 4B^2\mathsf FN-\mathsf J_1^2
 =4P_0R_B^2V(1-AD-AE-DE)\mathcal Q_1.
\end{equation}
The factor outside $\mathcal Q_1$ is positive. To fix the square-root branch, hold $A,B,E$ fixed and temporarily treat $V$ as independent, with $D=V-E$. The quadratic $\mathcal Q_1(V)$ has positive leading coefficient and negative constant. Its unique positive root $V_*$ satisfies
\[
 \mathcal Q_1(Q_0/T)=2B^2P_0Q_0^2/T>0,\qquad 0<V_*<Q_0/T.
\]
Writing $\mathsf J_1=j_2V^2+j_1V+j_0$ gives
\begin{align*}
 j_0&=2BU(1-AB)Q_0^2>0,\\
 j_1&=2EP_0R_BQ_0>0,\\
 j_2&=-2B^2P_0^2-2EP_0R_BT-2BU(1-AB)T^2<0.
\end{align*}
Let $h_2=TK_1$, $h_1=P_0Q_0(B^2-E^2)$, $h_0=-2BEUQ_0^2$. Their Sylvester determinant is
\[
 \det\begin{pmatrix}
 h_2&h_1&h_0&0\\0&h_2&h_1&h_0\\
 j_2&j_1&j_0&0\\0&j_2&j_1&j_0
 \end{pmatrix}
 =4B^3UP_0^2Q_0^4K_1\mathscr L,
\]
where $\mathscr L$ is specified without a long expansion by the following rational identity. Put $\chi=B^2$, $\rho=BE$ and $\tau=A(B+E)/(1-BE)$. Since $p>1$, $\chi>0$, $0<\rho,\tau<1$, and
\[
 -\mathscr L=\frac{\chi^3c_3+\chi^2c_2+\chi c_1+\rho^3}{\chi(\chi+\rho)},
\]
\begin{align*}
 c_3&=1-(1-\rho)\tau,\\
 c_1&=\rho[(1-\rho)\tau\{1+(1-\rho)(1-\tau)\}+\rho(2\rho+1)],\\
 c_2&=\rho(2\rho+1)+(1-\rho)^2\tau[(4\rho+\tau)(1-\tau)+\rho\tau^2].
\end{align*}
Every term is nonnegative and $\rho^3>0$. Thus the determinant is nonzero and $\mathsf J_1(V_*)\ne0$. The parameter domain $AB+AE+BE<1$ is connected: the displayed transformation has inverse
\[
 B=\sqrt\chi,\qquad E=\rho/\sqrt\chi,\qquad
 A=\frac{\tau(1-\rho)\sqrt\chi}{\chi+\rho}.
\]
The root and its $\mathsf J_1$ value vary continuously. At $A=B=E=1/4$,
\[
 V_*=\frac{17}{4\sqrt{21}},\qquad
 \mathsf J_1(V_*)=\frac{4913}{32768}\left(\frac{17}{\sqrt{21}}-1\right)>0.
\]
Hence $\mathsf J_1(V_*)>0$ everywhere. Concavity and $\mathsf J_1(0)>0$ give $\mathsf J_1>0$ on $[0,V_*]$.

If $\mathcal Q_1>0$, \eqref{r04:area-square} gives $|\mathsf J_1|<2B\sqrt{\mathsf FN}$, so $0<I_1<2BC$. If $\mathcal Q_1=0$, the positive branch gives $I_1=0$. If $\mathcal Q_1<0$, then $V<V_*$ and $\mathsf J_1>2B\sqrt{\mathsf FN}$, so $I_1<0<2BC$. Reflection proves the second assertion. The sign conclusion uses the branch argument, not merely the squared identity.
\end{proof}

\section{Central necessary inequalities and their branches}\label{r04:central-section}
Put
\begin{equation}\label{r04:central-coordinates}
 \xi=\frac{1-AB}{A+B},\qquad \eta=\frac{1-DE}{D+E},\qquad
 n=M(A+E)^2,\qquad \nu=\xi\eta,
\end{equation}
\[
 \zeta=C^2=\frac{(1+\xi^2)(1+\eta^2)}n-1,\qquad
 \delta=\xi^2+\eta^2+2-n=2(1+AE-A\xi-E\eta).
\]
The side bounds give $X_1<\xi-C$ and $Y_2<\eta-C$, hence $\xi,\eta>C>0$. The minimizing shapes cross the diagonal by Lemma~\ref{r04:shape-unique}: $Q_1,R_1>\max(B,C)$ and $Q_2,R_2>\max(C,D)$. The central inequality therefore gives $4C^2<1+C^2$, or $\zeta<1/3$.
The area formulas imply
\begin{equation}\label{r04:central-area}
 I_1=\frac{(1+\xi^2)\delta}{2n}-B(\xi-C),\qquad
 I_2=\frac{(1+\eta^2)\delta}{2n}-D(\eta-C).
\end{equation}
Thus $\delta>0$. Also exact expansion gives
\[
 1+\zeta-\delta=\frac{(\xi^2-\zeta)(\eta^2-\zeta)}{1+\zeta}>0,
 \qquad 0<\delta<1+\zeta<4/3.
\]
\begin{proposition}\label{r04:central-tests}
Every point of $\Omega_{04}$ satisfies
\begin{gather*}
 \beta_0=\nu+\zeta(1+\delta/2)+(\zeta-1)\delta^2/(4n)\ge0,\\
 \mathcal R_0=1-\delta-\nu\delta^2/(2n)+\delta^4/(16n^2)\ge0,\\
 0<\delta<1,\qquad \nu-\frac{\delta^2}{4n}\ge\sqrt{n\zeta}>0,
 \qquad \nu>\frac12.
\end{gather*}
\end{proposition}
\begin{proof}
Let $L_1=\xi-C$, $L_2=\eta-C$. The side bounds and $X_iY_i=I_i$ give
\[
 [(B+C)L_1+I_1][(C+D)L_2+I_2]
 \le L_1L_2(B+C+Y_1)(C+D+X_2)
 <(1+C^2)L_1L_2.
\]
Call the left side minus the last expression $\mathscr P$. Substituting \eqref{r04:central-area} gives the identities
\[
 \mathscr P=\alpha_0C-\beta_0,\qquad
 \alpha_0=(\xi+\eta)\left(1+\frac{(1+\nu)\delta}{2n}\right)>0,
\]
\[
 \alpha_0^2\zeta-\beta_0^2
 =-\frac{(1+\xi^2)^2(1+\eta^2)^2}{n^2}\mathcal R_0.
\]
It follows that $\beta_0\ge0$ and $\mathcal R_0\ge0$. Conversely those two inequalities imply $\mathscr P\le0$; the sign of $\beta_0$ is essential to this implication.

The frame also gives
\begin{equation}\label{r04:central-alternate}
 \nu^2=1-\delta+n\zeta,\qquad
 \mathcal R_0=\left(\nu-\frac{\delta^2}{4n}\right)^2-n\zeta.
\end{equation}
Write $h=\delta^2/(4n)>0$. Suppose $\delta\ge1$. Since $\nu\ge\zeta\ge\delta-1$ and $\xi^2+\eta^2\ge2\nu$, we have $n\ge2+2\nu-\delta\ge\delta$. Consequently $h\le\delta/4\le1/3$, whereas $n\zeta\ge\nu^2$. Equation~\eqref{r04:central-alternate} forces $h\ge\nu+\sqrt{n\zeta}\ge2\nu$; its other branch would give $h\le0$. Hence $\zeta\le\nu\le1/6$ and
\[
 \beta_0\le\zeta(\nu+1+\delta/2)-(1-\zeta)\sqrt{n\zeta}
 \le\frac{11}6\zeta-\frac56\sqrt\zeta<0,
\]
using $121<150$. This contradiction proves $\delta<1$.

If $h\ge\nu$, the same identity forces $h\ge\nu+\sqrt{n\zeta}$. But $n>1+2\nu$ and hence $\nu\le h<1/[4(1+2\nu)]$, which implies $\zeta\le\nu<1/5$. It follows that
\[
 \beta_0<\frac{17}{10}\zeta-\frac45\sqrt\zeta<0
\]
by $289<320$, another contradiction. Thus $h<\nu$ and the positive square-root branch in the statement is the only one possible.

Finally $p,q>1$ give $\xi>E$, $\eta>A$. The exact identity
\begin{equation}\label{r04:nu-half}
 \nu=1-\frac\delta2+(\xi-E)(\eta-A)
\end{equation}
then proves $\nu>1/2$. This is the additional necessary cut used in the numerical kernel.
\end{proof}

The completely polynomial version used by some contractions is also explicit. With $L=UV$, $X_0=1-AB$, $Y_0=1-DE$ and
\[
 \mathsf J_c=(1+AE)L-AX_0V-EY_0U,
\]
we have $n=N/L^2$, $\zeta=\mathsf F/N$, $\delta=2\mathsf J_c/L$ and $\nu=X_0Y_0/L$. Define
\begin{align*}
 \mathsf R_c&=L\mathsf J_c^4-2X_0Y_0N\mathsf J_c^2+N^2(L-2\mathsf J_c),\\
 \mathsf B_c&=X_0Y_0N^2+\mathsf FN(L+\mathsf J_c)+(\mathsf F-N)\mathsf J_c^2L.
\end{align*}
Then $\mathcal R_0=\mathsf R_c/(LN^2)$ and $\beta_0=\mathsf B_c/(LN^2)$. Thus all of
\[
 \mathsf J_c>0,\quad 2\mathsf J_c<L,\quad
 \mathsf R_c\ge0,\quad\mathsf B_c\ge0,\quad
 X_0Y_0N>L\mathsf J_c^2
\]
are necessary. They define a relaxation, not an equivalence with the original physical domain. The supplied exact symbolic scripts verify the displayed polynomial identities; the sign and branch arguments above justify their use.

\section{The quantitative local inequality}\label{r04:local-section}
Use $z_0$ from \eqref{r04:constants}, $\rho=11/10000$, and the square lower function
\[
 L(z)=C_5+\mathsf S(B+C,I_1)+\mathsf S(C+D,I_2).
\]
Then $G\ge L$ and $L(z_0)=m_0$.
\begin{proposition}\label{r04:local-bound}
If $z\in\Omega_{04}$ and $\|z-z_0\|_\infty\le11/10000$, then $G(z)>m_0$.
\end{proposition}
We specify all conditions checked by the finite local calculation. Its input is the rational decimal data in \texttt{r04\_conditional\_local\_certificate.json}; its checker is \texttt{verify\_r04\_conditional\_local.py}. These data are witnesses for the inequalities below, rather than numerical optimizer outputs admitted without verification.

Put
\[
 R(d,I)=\frac{d+\sqrt{d^2+4I}}2,\qquad
 L_1=\frac{1+B^2}{A+B}-B-C,\qquad
 L_2=\frac{1+D^2}{D+E}-C-D.
\]
For $s_1,s_2\in\{-1,1\}$ use the branch $s_1(B-C)\ge0$, $s_2(D-C)\ge0$, with
\[
 c_i=\begin{cases}17/10&s_i=1,\\2&s_i=-1,\end{cases}
 \qquad
 c_i'=\begin{cases}2&s_i=1,\\17/10&s_i=-1.\end{cases}
\]
First the checker verifies on the full cube that each base imbalance ratio is at most $3/10$ and each area divided by the square of its base sum is at least $1/5$. At feasible endpoints its upper bound $1/2$ follows from $p,q>1$. Lemma~\ref{r04:conditional-shape} therefore applies there. The following six smooth functions are nonnegative at feasible points of the branch:
\begin{align*}
 F_1&=1+A^2-(A+B)(A+E),&F_2&=1+E^2-(D+E)(A+E),\\
 F_3&=s_1(B-C),&F_4&=s_2(D-C),\\
 F_5&=L_1-R(c_1(B-C),I_1),&F_6&=L_2-R(c_2(D-C),I_2).
\end{align*}
Four more are $1+C^2-(B+C+Y)(C+D+X)$, with independent choices
\[
 Y\in\{I_1/L_1,R(c_1'(C-B),I_1)\},\qquad
 X\in\{I_2/L_2,R(c_2'(C-D),I_2)\}.
\]
Indeed $Y_1>I_1/L_1$, $X_2>I_2/L_2$ follow from the side bounds; the other lower bounds follow from the conditional shape differences, with the signs above. All ten functions vanish at $z_0$. The same formulas extend smoothly across the entire local cube. A Taylor segment need not be feasible or remain within its selected sign branch.

For a nonzero $h=z-z_0$ select $j\in\{1,2,3,4\}$ and $\sigma\in\{-1,1\}$ such that $\sigma h_j=\|h\|_\infty$. Then the six linear forms
\[
 H_{k,\tau}h=\sigma h_j+\tau h_k\ge0
 \quad(k\ne j,\ \tau\in\{-1,1\})
\]
are nonnegative. The four sign branches and eight maximum-coordinate cones give $32$ cases. For each, the data specify nonnegative rational $\lambda_i$, $\mu_{k,\tau}$, positive rational $\nu$ and a rational vector $v$. Set $P=L-\sum_{i=1}^{10}\lambda_iF_i$. Outward interval arithmetic verifies
\begin{align*}
 \left\|\nabla P(z_0)-\sum_{k,\tau}\mu_{k,\tau}H_{k,\tau}
 -\nu\sigma e_j\right\|_1&\le\varepsilon,\\
 \nabla^2P(z_0)&\succeq-\kappa I_4,\\
 \sum_{i,j,k=1}^4\sup_{\|z-z_0\|_\infty\le\rho}
 |\partial_i\partial_j\partial_kP(z)|&\le M_3.
\end{align*}
For the Hessian, Gershgorin row bounds are applied to $\nabla^2P(z_0)-vv^{\mathsf T}$. Take $\kappa\ge0$ as the maximum of zero and the outward upper bounds for the off-diagonal absolute row sum minus the diagonal entry. Since $vv^{\mathsf T}$ is positive semidefinite, this gives the stated lower bound for $\nabla^2P(z_0)$. The final strict check is
\begin{equation}\label{r04:local-budget}
 \varepsilon+2\rho\kappa+\rho^2M_3/6<\nu.
\end{equation}
All $32$ cases satisfy it; the largest certified upper bound for the left-to-right ratio is less than $0.896$.

\begin{proof}[Proof of Proposition~\ref{r04:local-bound}]
Put $r=\|h\|_\infty>0$. Nonnegative multipliers and $F_i(z)\ge0=F_i(z_0)$ imply
$L(z)-L(z_0)\ge P(z)-P(z_0)$. Taylor's theorem and $\|h\|_2^2\le4r^2$ give
\[
 P(z)-P(z_0)\ge\nu r-\varepsilon r-2\kappa r^2-M_3r^3/6>0
\]
by $r\le\rho$ and \eqref{r04:local-budget}. The center itself is outside $\Omega_{04}$ because $p=q=1$. Thus every point asserted in the proposition has $G\ge L>m_0$.
\end{proof}

The local checker uses $50$-digit interval arithmetic and the square series through $m=28$. The remainder bounds in Section~\ref{r04:square-series} enclose the full tail and its first three derivatives. In explicit notation the first three tails are bounded by
\[
 \frac{\alpha_{29}w^{29}}{1-w},\qquad
 \frac{29\alpha_{29}w^{28}}{1-w},\qquad
 \frac{28\cdot29\alpha_{29}w^{27}}{1-w},
 \quad \alpha_{29}=[29\cdot59\cdot60]^{-1},
\]
and the third derivative tail by $w^{26}/[4(1-w)]$. The center zeros used above follow from its exact algebraic identities. An interval merely containing zero would not establish them.

\section{The complete cover and the strict global bound}\label{r04:cover-section}
We explain how the analytic statements qualify the finite cover, including the lower functions used on boxes crossing a zero-area face. Section~\ref{reader:r04} gives the executable interfaces and data identities.

If $G(z)\le m_0$ at $z\in\Omega_{04}$, Proposition~\ref{r04:tail-bound} gives $0<A,E<5$. The diagonal-crossing property and the side bounds give $B^2+2AB<1$, $D^2+2DE<1$, hence $0<B,D<1$. Thus every possible counterexample lies in the closed box
\begin{equation}\label{r04:root-box}
 Q_{04}=[0,5]\times[0,1]\times[0,1]\times[0,5].
\end{equation}

The domain $\Omega_{04}$ is invariant under $(A,B,D,E)\mapsto(E,D,B,A)$. This fixes $C$, exchanges $I_1,I_2$ and sends the unique shapes to
\[
 (X_1',Y_1',X_2',Y_2')=(Y_2,X_2,Y_1,X_1).
\]
The side inequalities exchange, the two central factors exchange, and both $G$ and $L$ remain unchanged. Thus the cover may restrict to $A\le E$. This reflection is used only on the relaxed domain and its costs; it makes no assertion about reflection of the original prescribed real winding.

Here are the whole-box mathematical tests. If a smooth lower function $f$ satisfies $G\ge f$ at feasible points and the fixed smooth formulas $g_k$ are necessary nonnegative functions there, then for fixed $\lambda_k\ge0$ a box with center $c$ and half-widths $r_i$ is excluded by
\begin{equation}\label{r04:dual-box-test}
 f(c)-\sum_k\lambda_kg_k(c)
 -\sum_i r_i\sup_Q\left|\partial_if-\sum_k\lambda_k\partial_ig_k\right|>m_0.
\end{equation}
The multipliers are proposed numerically and then interpreted as exact binary rationals for interval verification. No optimization status is a proof of this inequality.

A second test uses an explicit positive semidefinite $P=\sum_j\alpha_jv_jv_j^{\mathsf T}$, $\alpha_j\ge0$. If $\nabla^2 f-P\succeq\ell I$ on the box with $\ell\le0$, then for any fixed vector $q$
\begin{equation}\label{r04:psd-box-test}
 f(c+h)\ge f(c)-q^{\mathsf T}Pq/2
 +[\nabla f(c)+Pq]^{\mathsf T}h+\ell\|h\|_2^2/2.
\end{equation}
This follows by completing the square in $P$ and using the Taylor integral remainder. It also applies to $C^{1,1}$ functions when the Hessian bounds hold almost everywhere.

For square costs write $\mathsf S(S,I)=S^2\psi(I/S^2)$. Its convexity gives the fixed-tangent lower bound
\[
 \mathsf S(S,I)\ge\psi'(r)I+[\psi(r)-r\psi'(r)]S^2
 \qquad(r\ge0).
\]
The same fixed $r$ must be used in the box value, center value, gradient and Hessian. Likewise conditional shape constants chosen from a whole-box test are kept fixed when evaluating at its center.

To enclose boxes whose area interval crosses zero, the kernel uses for $S>0$ the extension
\[
 \widetilde{\mathsf S}(S,I)=\frac32I+J[\widehat H(w)-3/2],
 \qquad J=\max(I,0),\quad w=\frac{J}{(S+\sqrt J)^2},
\]
where $\widehat H=3/2+H$ is the full positive series. It equals the square cost for $I\ge0$ and equals $3I/2$ for $I<0$. It is locally $C^{1,1}$ on $S>0$. For $I>0$, put $d=S+\sqrt I$, $l=\sqrt I$. Its derivatives are
\begin{align*}
 \widetilde{\mathsf S}_I&=\widehat H+IS\widehat H'/d^3,&
 \widetilde{\mathsf S}_S&=-2I^2\widehat H'/d^3,\\
 \widetilde{\mathsf S}_{II}&=S(2S+l/2)\widehat H'/d^4+IS^2\widehat H''/d^6,\\
 \widetilde{\mathsf S}_{IS}&=-I(4S+l)\widehat H'/d^4-2I^2S\widehat H''/d^6,\\
 \widetilde{\mathsf S}_{SS}&=6I^2\widehat H'/d^4+4I^3\widehat H''/d^6.
\end{align*}
The second derivatives have finite one-sided limits at zero. They vanish for $I<0$; taking both envelopes gives the almost-everywhere Hessian bound required for \eqref{r04:psd-box-test}. The explicit series remainders bound complete functions, not just truncated polynomials.

When $I_i,L_i>0$ on the whole box, the side condition also implies the exact trial-mean test
\[
 \mu_X(B,C;L_1,I_1/L_1)-\mu_Y(B,C;L_1,I_1/L_1)\ge0
\]
and its reflection. This uses the positive sign after the unique shape root, without approximating that root. The relevant mean is explicit:
\[
 \frac1Y\int_a^{a+Y}k(u,b)\,du
 =2-2q\varphi(r)+\frac{2q^2}{1+r},\quad
 q=\frac b{a+b},\quad r=\frac Y{a+b},
\]
where $\varphi(r)=\log(1+r)/r$, $\varphi(0)=1$. Its integral $\varphi(r)=\int_0^1(1+rt)^{-1}\,dt$ gives the required monotonicity and alternating envelopes at zero. All other necessary tests come from the area and central polynomials proved above, the positive durations, and the shape bounds; they are not sufficient definitions of the original trajectory class.

The single binary rational threshold for cost deletion is
\begin{equation}\label{r04:threshold}
 q_*=\frac{8128224356559707}{4503599627370496}.
\end{equation}
An exact rational check encloses $\sqrt2$ between two fifty-digit rationals and uses
\[
 \log2=2\sum_{j=0}^{40}\frac1{(2j+1)3^{2j+1}}+R,
 \qquad0<R<\frac2{83\,3^{83}(1-1/9)}.
\]
It gives $m_0<1.804828365994315<1.804828365994316<q_*$ and $q_*-m_0>2.4\cdot10^{-15}$. Thus a certified lower bound at least $q_*$ proves strict inequality over $m_0$.

\begin{proposition}[Complete directed interval exclusion]\label{r04:finite-cover}
For $z\in\Omega_{04}\cap Q_{04}$ with $A\le E$, the closed cover gives $G(z)>m_0$.
\end{proposition}
\begin{proof}
The task \path{R04-full} starts from a single closed box
$[0,20]\times[0,1]^2\times[0,20]$, which contains $Q_{04}$, and uses
the symmetry restriction $A\le E$. In each of the $31$ runs the
kernel recomputes the interval tests. It uses outward binary64
arithmetic and directed MPFR bounds for square roots and logarithms.
The build disables fast-math transformations and fused multiply-add
contraction. A test can discard a box if a necessary condition fails,
if the cost lower bound reaches the threshold~\eqref{r04:threshold},
or if the box lies inside the local cube of Proposition~\ref{r04:local-bound}.
Both local deletion radii used by the kernel are strictly less than
$11/10000$. An inconclusive test leaves the box available for further
tests or closed subdivision.

At the joins the structural check accounts for each pending box by
transfer, closed subdivision or the analytic tail estimate. For example,
the $19$ boxes left by the first run are subdivided into $1\,216$ inputs
to the next stage. The $203$ whole boxes set aside for the tail each have
lower endpoint at least $5$ in $A$ or in $E$. Proposition~\ref{r04:tail-bound}
therefore applies throughout their intersection with $\Omega_{04}$,
including the faces at $5$. The exact transfer check matches the
remaining boxes to the next inputs. At the end there are no pending or
gray boxes. Since every contraction preserves possible counterexamples
and the closed children cover their parent, induction over the runs
and their subdivisions excludes every counterexample in the stated domain.
\end{proof}

\begin{proof}[Proof of Theorem~\ref{r04:global-bound}]
If $z\in\Omega_{04}$ satisfied $G(z)\le m_0$, the tail and side
bounds would place it in $Q_{04}$. Reflection $(A,B,D,E)\mapsto(E,D,B,A)$
preserves $\Omega_{04}$, $Q_{04}$ and $G$, so we could take $A\le E$.
Proposition~\ref{r04:finite-cover} then gives $G(z)>m_0$.
\end{proof}
\space
\label{revision:appendix-end-10}
\end{proofappendix}
\section{Verification algorithms and their mathematical inputs}\label{app:verification}
\begin{proofappendix}
\section{Mathematical inputs and the reader archive}\label{verify:interfaces}
The finite checks use the fixed archive \path{reinhardt-reader-algorithms--1118R01.tar.xz}.
A check verifies a specified algebraic identity, bound, recurrence or
closed cover. The corresponding proof identifies the complete
physical source, shows that it lies in the checked domain, and applies
the checked inequality. The proofs of the first-contact reduction and
its nine-arc transfer use the separate \path{h2-h3-analytic-sources}
supplement. Its guide uses the numbering of version~1122R06 to locate the
original derivations. Sections~\ref{six:app:physical-cap}--\ref{six:app:reduction-proofs}
of the present paper give the hypotheses, arguments and supplementary
comparisons used here. The fixed 1118R01 archive provides the finite checks.
The common rules for outward evaluation and closed covers are in
Section~\ref{cv:section}.

In the archive, \path{CLAIM-MAP.json} gives the claim identifiers and
their statement locations in the bundled version~1 manuscript.
For the finite-switching input listed under \texttt{G01}, the present
manuscript supplies an analytic proof in Theorem~\ref{fw:finiteness}.
\path{TASKS.json} defines the execution inputs and acceptance conditions,
and \path{FILE-REFERENCE-MAP.json} locates the fixed data and sources.
Sections~\ref{reader:algorithms} and~\ref{verify:r05-source-inputs} give
the task entries, with the R05 coefficient and covering tasks in the latter.
An expected count identifies the comparison to be reproduced. Covering
and exclusion require the reconstructed closed splits and the verified
predicates at every terminal.
\space
\section{Finite certificates and their mathematical interpretation}
\label{cv:section}

The finite calculations in this paper have three logically separate parts.
First, an analytic argument associates to each hypothetical minimizing source
parameters satisfying explicitly stated necessary conditions. Second, a finite
certificate covers the resulting parameter domain and rejects every retained
cell, or supplies an exact linear combination of valid payment inequalities.
Third, the analytic classification and complete source placements associate
those conclusions to the cyclic words being removed. A saved success flag, a file digest, or a count of rejected cells
cannot replace any of these three parts. We give the verification principles
used here so that the electronic witnesses have a precise mathematical meaning.

\subsection{Affine fibres and guarded loss comparisons}
\label{cv:comparison}

\begin{lemma}[Affine-fibre Casimir comparison]\label{cv:affine-fibre}
Fix the geometric parameters. Suppose
$L(z)=L_0+\ell\cdot z$, $K(z)=K_0+k\cdot z$ and
$m(z)=m_0+r\cdot z$ are affine, and put
\[
 C(z)=-L(z)^2-K(z)m(z),\qquad
 Q=\ell\ell^{\mathsf T}+\tfrac12(kr^{\mathsf T}+rk^{\mathsf T}).
\]
If $Q$ is positive semidefinite and a qualified actual parameter belongs
to $\operatorname{conv}(V)$ for a finite set $V$, then
$C(z)\ge\min_{v\in V}C(v)$. The affine formulas must hold on that convex
hull; its vertices need not be physical. Endpoint bounds $C(v)\ge H$
contradict an actual strict calibration $C(z)<H$, whereas a weak calibration
$C(z)\le H$ requires strict endpoint bounds $C(v)>H$.
\end{lemma}
\begin{proof}
Expansion gives $C(z)=c+b\cdot z-z^{\mathsf T}Qz$. If
$z=\sum_i\lambda_iv_i$, $\lambda_i\ge0$, $\sum_i\lambda_i=1$, then
\[
 C(z)-\sum_i\lambda_iC(v_i)
 =\sum_i\lambda_i(v_i-z)^{\mathsf T}Q(v_i-z)\ge0.
\]
This proves both assertions. In one dimension the finite set consists
of the two ends of a closed enclosing interval, including a singleton.
\end{proof}

An affine matrix family is not automatically Casimir-concave: the check
$Q\ge0$ is indispensable. Affine cuts retain the entire closed feasible
fibre. Division by a cut's slope is allowed only after its sign is proved;
a zero slope is handled without division. Strict and weak contacts,
including singleton fibres, are distinguished as in Lemma~\ref{cv:contraction}.

\begin{lemma}[Guarded cap-loss comparison]\label{cv:guarded-loss}
Suppose $t>u+v>0$, $u>a\ge0$, $v>b\ge0$, $uv>W>0$, and
\[
 C=\kappa-(t-u-v)^2/4+uv\le\overline C,\qquad Q_0=a+b.
\]
On $[0,u]\times[0,v]$ the loss expression increases in each variable.
If $\Delta=4\kappa-4\overline C+4W\ge0$, then
$t>Q_0+\sqrt\Delta$. For every $T\ge t$,
\[
 B:=\overline C-\kappa+T^2/4>T\sqrt W>0,
 \qquad B^2>T^2W.
\]
\end{lemma}
\begin{proof}
The partial derivatives are $(t-u'+v')/2$ and $(t+u'-v')/2$,
both positive on the stated rectangle. Since $t>u+v>Q_0$,
\[
 C>\kappa-(t-Q_0)^2/4+W,
\]
which selects the positive square-root branch. Put $w=\sqrt{uv}$.
Then $t>u+v\ge2w>2\sqrt W$ and
$C\ge\kappa+tw-t^2/4>\kappa+t\sqrt W-t^2/4$.
The last quadratic decreases for all arguments at least $t$; replacing
$t$ by $T$ proves the positive inequality before it is squared.
\end{proof}

\begin{lemma}[Correlated endpoint profile]\label{cv:correlated-profile}
Let $0<\alpha,\beta<1$, and let $q>0$ be absolutely continuous on
$[-\alpha,\beta]$, with
$-2/(1+u)\le(\log q)'\le2/(1-u)$ almost everywhere. Put
\[
 x_L=q(-\alpha)(1-\alpha)^2,\quad x_R=q(\beta)(1-\beta)^2,
 \quad d=\alpha+\beta,\quad p=\alpha\beta,
\]
\[
 M(t)=\frac{1+t}{1-t},\quad Z=M(\alpha)M(\beta),\quad S=\int_{-\alpha}^{\beta}q(u)\,du.
\]
Then
\begin{align*}
 q(u)&\le\min\left\{\frac{x_LM(\alpha)^2}{(1-u)^2},
                   \frac{x_RM(\beta)^2}{(1+u)^2}\right\},\\
 S&\le\sqrt{x_Lx_R}Z-\tfrac12\{x_LM(\alpha)+x_RM(\beta)\}
 \le\sqrt{x_Lx_R}(Z-\sqrt Z).
\end{align*}
If $x_L,x_R<X$, then
\[
 S<X\frac{d+2p}{1-d+p};\qquad
 d<2/3\ \Longrightarrow\ S<\frac{2Xd(d+2)}{(2-d)^2}.
\]
If $e_\alpha,e_\beta>0$ and
$u>q(-\alpha)e_\alpha/S$, $v>q(\beta)e_\beta/S$, then
\[
 uv>\frac{(1+\sqrt Z)^2}{4}P>\lambda^2P,\quad
 P=\frac{(1-d+p)e_\alpha e_\beta}{d^2(1+d+p)},\quad
 \lambda=\frac{2(1+p)}{2+2p-d}.
\]
\end{lemma}
\begin{proof}
Integrate the upper logarithmic slope from the left and the lower slope
backward from the right. At the endpoints the resulting bounds give
$x_L\le x_RM(\beta)^2$ and $x_R\le x_LM(\alpha)^2$.
Thus the two envelope curves cross inside the interval at
$(b-a)/(b+a)$, where $a=\sqrt{x_L}M(\alpha)$,
$b=\sqrt{x_R}M(\beta)$. Integration of their minimum gives the first
bound on $S$; AM--GM on the subtracted bracket gives the second.
Increasing both endpoint coefficients to $X$ gives the strict raw
ceiling. Its derivative in $p$ is
$X(2-3d)/(1-d+p)^2>0$ when $d<2/3$, so $p\le d^2/4$ gives the
symmetric ceiling only under that guard.
Multiplication of the strict loss inequalities and the second $S$ bound
gives the stated factor $(1+\sqrt Z)^2/4$. Finally, for
$t=d/(1+p)\in(0,1)$,
\[
 Z-\left(\frac{2+t}{2-t}\right)^2
 =\frac{2t^3}{(1-t)(2-t)^2}>0.
\]
Taking positive square roots proves $(1+\sqrt Z)/2>\lambda$.
\end{proof}

\subsection{Outward arithmetic and analytic functions}
\label{cv:arithmetic}

An interval is a closed set $[l,u]$ with rational endpoints and $l\leq u$.
In fixed-point implementations stored integers $L,U$ represent the endpoints
$L2^{-p},U2^{-p}$ for a specified precision $p$. Integer floor and ceiling
operations must be used with their mathematical meanings also on negative
integers. All input constants are rationals, integers, or separately proved
rational enclosures.

\begin{lemma}[Outward evaluation]\label{cv:outward}
Suppose that each elementary operation maps its input intervals to an interval
containing every defined value of that operation. Suppose further that division
is performed only after the denominator interval has been proved to avoid
zero, and that every nonalgebraic primitive has a proved enclosure on its stated
domain. Then evaluation of any expression built from these operations encloses
its value at every point of the input box where the expression is defined.
\end{lemma}
\begin{proof}
Induct on the expression tree. At a leaf the assertion is the input enclosure.
At an internal node, the induction hypothesis encloses every operand and the
assumed elementary-operation rule encloses the result. Repeated occurrences of
a variable may enlarge the interval but cannot invalidate containment. For
fixed-point addition and subtraction the endpoint formulas are exact. For
multiplication, take the minimum and maximum of the four endpoint products,
then round the lower endpoint downward and the upper endpoint upward. Division
is handled by the four endpoint quotients when the denominator interval avoids
zero. Squaring an interval that contains zero has lower endpoint zero. A square
root requires a nonnegative argument interval and enclosing integer square
roots. These rules have the required containment property.
\end{proof}

The response function needed below can be evaluated without a floating-point
logarithm. Set
\[
 \kappa(a)=\int_0^1\frac{1-t}{1+at}\,dt
 =\begin{cases}
 ((1+a)\log(1+a)-a)/a^2,&a>0,\\
 1/2,&a=0.
 \end{cases}
\]
Differentiation under the integral sign shows that $\kappa$ is decreasing for
$a\geq0$. Writing $z=a/(a+2)$ gives, for $0\leq z<1$,
\begin{equation}\label{cv:kseries}
 \kappa(a)=\frac12-
 \sum_{j=0}^{\infty}\frac{z^{2j+1}}{(2j+1)(2j+3)}.
\end{equation}
One obtains this identity either by integrating a geometric series after the
Cayley substitution or by substituting
$\log(1+a)=2\sum_{j\geq0}z^{2j+1}/(2j+1)$ in the displayed formula for
$\kappa$. Absolute convergence justifies both manipulations for $z<1$.
After $m$ terms, the omitted positive tail $R_m$ satisfies
\begin{equation}\label{cv:ktail}
 \frac{z^{2m+1}}{(2m+1)(2m+3)}
 \leq R_m\leq
 \frac{z^{2m+1}}{(2m+1)(2m+3)(1-z^2)}.
\end{equation}
Indeed, the left side is the first omitted term, and replacing all later
positive denominators by the first omitted denominator bounds the sum by a
geometric series. The identity and bounds include $a=0$ by continuity. Thus
rational endpoint arithmetic supplies certified values of $\kappa$; monotonicity
gives its interval extension. These general arithmetic bounds remain
available independently of the archived inward-pair evaluator.

\begin{lemma}[Interval derivatives and centered forms]\label{cv:mean-value}
Let $f$ be continuously differentiable on a box $C$, let $c\in C$, and suppose
intervals $G_j$ enclose $\partial_j f$ throughout $C$. Then
\[
 f(C)\subseteq f(c)+\sum_{j=1}^dG_j(C_j-c_j).
\]
The conclusion remains valid when $f(c)$ is replaced by an outward enclosure.
For an exact expression formed from sums, products, reciprocals and
logarithms, propagate value--gradient pairs using
\[
 \begin{array}{c|c}
 \text{expression}&\text{gradient}\\ \hline
 u+v&\nabla u+\nabla v\\
 uv&v\nabla u+u\nabla v\\
 u^{-1}&-u^{-2}\nabla u\\
 \log u&u^{-1}\nabla u
 \end{array}
\]
Coordinates have unit gradients and constants have zero gradients.
The resulting gradient intervals are valid if primitive values and gradients
are enclosed on the whole box, reciprocal denominators avoid zero, and
logarithm arguments are positive. Additional smooth primitives may be used
with separately proved value and gradient enclosures.
\end{lemma}
\begin{proof}
For $x\in C$, the segment from $c$ to $x$ stays in $C$. Integrating the
identity $\frac{d}{dt}f(c+t(x-c))=\sum_j\partial_j f(c+t(x-c))(x_j-c_j)$
gives the inclusion. Derivative enclosure by automatic differentiation follows
by induction on the expression tree, applying the exact sum, product and chain
rules and Lemma~\ref{cv:outward} at every node. A removable singularity must be
replaced by its proved differentiable extension before this argument is used.
An endpoint or remainder bound for a primitive is not that primitive:
its derivative cannot replace the exact derivative. Intersections of two
valid whole-box enclosures of the same value or gradient remain valid;
an empty intersection on a nonempty analytic box signals a failed
calculation, not a physical exclusion.
\end{proof}

\subsection{Finite covers and conditional contractions}
\label{cv:cover}

Let $\mathcal X$ denote the set of qualified physical sources under discussion,
and let $\Phi:\mathcal X\to\mathbb R^d$ be the parameter map supplied by the
analytic reduction. The necessary domain is a set $D$ containing
$\Phi(\mathcal X)$; it may be a proper subset of a rational root box $Q$.
We never require every point of $D$ or $Q$ to be physically realizable.

\begin{theorem}[Sound finite rejection]\label{cv:finite-rejection}
Assume the following statements.
\begin{enumerate}
\item The analytic source argument proves $\Phi(\mathcal X)\subseteq D\subseteq Q$.
\item A finite rooted tree has root box $Q$. At each internal node with box
$C$, its children are closed boxes whose union contains $C\cap D$.
Alternatively, a finite forest is permitted if its roots cover $Q\cap D$.
\item Each terminal cell $C$ has an independently checked rejection. This is
an analytic exclusion of $C\cap D$, or an outward evaluation of a necessary
condition showing that no point of $C\cap D$ can satisfy that condition.
\item Every node needed for coverage has been accounted for, and there are no
unresolved terminal obligations.
\end{enumerate}
Then $\mathcal X$ is empty.
\end{theorem}
\begin{proof}
Suppose $x\in\mathcal X$. Its image $\Phi(x)$ belongs to a root box and to
$D$. At every internal node containing that image, the child-cover hypothesis
places it in at least one child. Since the tree is finite, this process reaches
a terminal cell containing $\Phi(x)$. The terminal rejection contradicts the
necessary conditions satisfied by $\Phi(x)$. The same argument applies to a
forest. Shared subdivision boundaries cause no difficulty, since children are
closed and at least one child contains every boundary point.
\end{proof}

For binary midpoint subdivision, coverage is established by checking the exact
integer endpoints, the selected coordinate, and the midpoint at every split.
A compact path representation is sufficient: the root is numbered $1$ and
children of $v$ are numbered $2v,2v+1$. A checker reconstructs each box from its
path, checks that the terminal paths are distinct and prefix-free, and checks
that every nonterminal prefix has both children. These conditions establish a
finite full tree with precisely the listed leaves. The identity
$\sum_v2^{-\operatorname{depth}(v)}=1$ is an additional combinatorial check;
it is not, by itself, proof that the geometric boxes cover the root. In
particular, combinatorial tree mass is not interchangeable with geometric
volume when integer rounding makes child widths unequal.

\begin{lemma}[Safe contraction]\label{cv:contraction}
Suppose a box $C$ is replaced by $C'$ after proving
$C\cap D\subseteq C'\subseteq C$. Rejection of $C'\cap D$ rejects $C\cap D$.
An empty contracted set rejects the original cell. A singleton intersection
must be retained unless it is separately rejected.
\end{lemma}
\begin{proof}
Any point of $C\cap D$ would also be a point of $C'\cap D$. The conclusion is
therefore the contrapositive of the containment assertion.
\end{proof}

A contraction is conditional on the necessary inequalities, not a claim about
all algebraic points of $C$. For example, from $\alpha(x)t\leq\beta(x)$ one
may divide only on a subcell where the sign of $\alpha$ has been established.
If its interval contains zero, retain the unsimplified inequality or use a
separately proved zero-slope case. Likewise, if an analytic argument excludes
$a+b<\varepsilon$, subsequent calculations may intersect with the closed
condition $a+b\geq\varepsilon$; the analytic argument must cover the omitted
actual sources. A strict necessary condition $f>0$ is contradicted by an upper
enclosure at most zero. A weak necessary condition $f\geq0$ requires an upper
enclosure strictly below zero. These distinctions apply to all boundary tests.

\begin{lemma}[Exact grouping of a remainder]\label{cv:grouping}
Let a finite family of unresolved boxes have the form $S_i\times I_i$, with
$S_i\subset\mathbb R^{d-1}$ and closed intervals $I_i$. Partition the indices
only by exact equality of $S_i$. For each group replace its intervals by the
connected components of their exact union. The resulting grouped boxes have
exactly the same union as the original unresolved boxes.
\end{lemma}
\begin{proof}
Within a fixed group, distributivity gives
$\bigcup_i(S\times I_i)=S\times\bigcup_i I_i$. Sorting rational endpoints
and merging only overlapping or touching intervals computes that union
exactly. Taking the union over the groups proves the assertion. Strict gaps
between components remain gaps. No analogous identity permits pooling
intervals attached to unequal projected boxes.
\end{proof}

For a grouped two-level cover, every final outer leaf must carry every component
of its original interval union. Each inner tree must cover its entire component.
When an outer cell is split, both children inherit the full original component
union. Numerical terminals attached to discarded ancestors cannot substitute
for the terminals of a final outer leaf. This is a direct application of
Theorem~\ref{cv:finite-rejection}, first to each inner tree and then to the
outer tree.

\subsection{Polynomial and capacity certificates}
\label{cv:bernstein-capacity}

\begin{lemma}[Bernstein enclosure]\label{cv:bernstein}
Let a polynomial $p$ on a rational box be transformed affinely to
$[0,1]^d$, and write its exact tensor-product Bernstein expansion as
\[
 p(t)=\sum_{0\leq\alpha_i\leq n_i}b_\alpha
 \prod_{i=1}^d\binom{n_i}{\alpha_i}
 t_i^{\alpha_i}(1-t_i)^{n_i-\alpha_i}.
\]
Then $\min_\alpha b_\alpha\leq p(t)\leq\max_\alpha b_\alpha$ throughout
the closed box. In particular, all coefficients nonnegative prove $p\geq0$,
and a strictly positive minimum proves $p>0$.
\end{lemma}
\begin{proof}
Every displayed basis function is nonnegative on the cube, and their sum is
one by the binomial theorem in each coordinate. Thus $p(t)$ is a convex
combination of the coefficients. The affine coordinate change is exact and
maps the closed box onto the cube.
\end{proof}

If $p(X)=\sum_j a_jX^j$ and $n$ bounds its coordinate degrees, its
complete coefficient vector is specified by the exact finite formula
\begin{equation}\label{cv:bernstein-coefficients}
 b_i=\sum_{j\le i}a_j\prod_r
             \frac{\binom{i_r}{j_r}}{\binom{n_r}{j_r}},
 \qquad 0\le i_r\le n_r.
\end{equation}
Indeed, $X^j$ has Bernstein coefficients
$\prod_r\binom{i_r}{j_r}/\binom{n_r}{j_r}$, by the binomial theorem,
and linearity gives the formula. A polynomial certificate is therefore
specified by its exact residual, substitution, positive clearing multiplier,
degree vector and coefficient bound, with any canceled faces stated.
For a univariate interval $[a,b]$, put $h=(a+bt)/(1+t)$.
If $p(a+(b-a)x)=\sum_{i=0}^nb_i\binom ni x^i(1-x)^{n-i}$, then
\[
 (1+t)^np\left(\frac{a+bt}{1+t}\right)
 =\sum_{i=0}^n\binom ni b_it^i.
\]
Thus the homographic positive-coefficient test is precisely the same
criterion, and its leading coefficient supplies the endpoint $h=b$.
For midpoint subdivision, in one coordinate set
$b_i^{(0)}=b_i$ and
\begin{equation}\label{cv:bernstein-subdivision}
 b_i^{(r+1)}=\tfrac12(b_i^{(r)}+b_{i+1}^{(r)}).
\end{equation}
The left child coefficients are $b_0^{(r)}$ and the right child
coefficients are $b_{n-r}^{(r)}$ in reverse order, $0\le r\le n$.
Substitution of $x/2$ or $(1+x)/2$ in the Bernstein expansion and the
binomial theorem prove these identities. Applying the rule in one
coordinate at a time gives exact tensor subdivision, including all
closed child faces.

A checker must reconstruct the polynomial identity and every coefficient over
$\mathbb Q$, rather than merely read a saved positive minimum. Clearing a
rational denominator requires a proof of its sign on the entire box; possible
zeros must be retained as separate cases. Nonnegative coefficients with some
zero entries do not prove strict positivity at every boundary point. Strictness
must come from a positive coefficient bound, an independently proved positive
factor, or the actual-source hypotheses used in the particular application.

\begin{theorem}[Exact payment certificate]\label{cv:capacity}
For qualified physical payment rows $Ac\ge b$, $c\ge0$, any rational
$\lambda\ge0$ with $A^{\mathsf T}\lambda\le\mathbf1$ proves
$\mathbf1^{\mathsf T}c\ge\lambda^{\mathsf T}b$, with the strictness and
exact equality conditions of Theorem~\ref{as:capacity}.
\end{theorem}
\begin{proof}
Apply \eqref{as:slack-identity} with $a_{ri}=A_{ri}$.
\end{proof}
The rational vector must be checked exactly, regardless of its discovery
method. Budget equality and exhaustive alternatives are governed by
Corollary~\ref{as:budget-exclusion} and Theorem~\ref{as:branches}; a feasible
primal certifies only its scalar system, not a physical minimizer.

\subsection{Fixed local polynomial calculations}
\label{cv:local-exact}

The electronic packet \texttt{local\_exact} supplies fixed rational and
symbolic calculations for the local estimates. It also contains auxiliary
comparisons; only the calculations expressly invoked in the article are
needed for its proof. Each calculation invoked in the proof starts from
explicit formulas printed there, reconstructs the
polynomial or rational expression, and checks every coefficient in the stated
positive basis. Polynomial identities are checked exactly after clearing
nonzero denominators. There is no fitting of sampled values or search over
unprinted parameters.

The opposite-centre bound $63/500$ is now proved by the monotone
compensation and endpoint argument in Theorem~\ref{fee:opposite126}.
The fixed polynomial comparisons still used by other center profiles
retain the explicit positive-denominator and Bernstein conditions above.
\space
\section{Finite checks and their manuscript claims}
\label{reader:algorithms}\label{reader:appendix-start}\label{new:review}\label{ten:verification-scope}
The entries below are task identifiers in \path{TASKS.json}.
The archive's \path{verify.py} runs each task with the fixed inputs
and its own acceptance conditions. The analytic comparisons in the
cited statements give the normalization, positive denominators, and
branches on which the finite results are used.

\begin{samepage}
\subsection{Local and short-chain checks}\label{reader:local}\label{reader:local:algebra}\label{reader:local:cover}\label{reader:local:m109}
\begin{center}\small
\begin{tabular}{@{}>{\raggedright\arraybackslash}p{.36\textwidth}>{\raggedright\arraybackslash}p{.56\textwidth}@{}}
\toprule
Comparison & Task entries and checked obligation\\
\midrule
Exact local algebra and returned-HIGH absorption
 & \path{local-algebra}. Exact identities, positive-denominator polynomial comparisons and their Bernstein coefficients.
The separate \path{M109} check gives the two-return large-width sign.\\
Three-arc \(49/200\) payment (Theorem~\ref{f245:theorem})
 & \path{F245-positive}, \path{F245-negative}, \path{F245-outer}. Positive residual closed cover, negative branch and outer inequalities.\\
Forward pair and weak five-arc charge
 & \path{local-forward-pair}, \path{local-proper-five}. Complete fixed covers for the comparisons of Theorems~\ref{nx:forward-pair} and~\ref{wk:five600}.\\
First contact for the returned union
 & \path{H2-01} through \path{H2-21}, \path{H3-01} through \path{H3-05}.
Check the closed pieces and scalar tests of Section~\ref{six:app:computation}.
Appendix~\ref{six:app:reduction-proofs} and its analytic supplement give
the physical qualifications for these tests.\\
R07, R15, R17
 & \path{local-r07}, \path{local-r15}, \path{local-r17}. The closed residual systems in Lemmas~\ref{new:r07finite}, \ref{new:r15finite} and~\ref{new:r17center}, including the other finite comparisons of their proofs. R17 retains its 64-significand-bit \texttt{long double} contract.\\
\bottomrule
\end{tabular}
\end{center}
\end{samepage}

\medskip\noindent
\textit{Six arcs (\path{SHORT6-ALL}).}\label{short:source-appendix}\label{short:finite-inputs}\label{reader:local:taylor}\label{sh:arithmetic-six}
The checked derivative in Lemma~\ref{sh:six02-finite} holds \(a,t\)
fixed. The small-\(r\) polynomial, two analytic tails, compact
closed boxes and derivative cube have the distinct roles given in
its proof. On the compact cover the permitted tests are
\(\mathscr J>K_{\rm sh}\), \(\mathscr J_b>0\),
\(\mathscr J_b+8(\mathscr J-K_{\rm sh})>0\), or containment in the separately
verified derivative cube. The task recomputes the tests with outward
integer enclosures at scale \(2^{-80}\), including the continuously
extended logarithmic quotients.

\medskip\noindent
\textit{Reverse seven arcs (\path{reverse7-algebra} and the nine
 \path{reverse7-} family tasks).}\label{sh:arithmetic-reverse}
Proposition~\ref{sh:reverse-seven} gives the actual group-state
inclusion and the lower function \(L\). The finite cost and derivative
predicates use the original compact domains and collar. The seven
analytic tail estimates and the positive terms in the octagon collar
are in its proof. The interval checks use directed integer arithmetic
with first and second derivatives enclosed on complete boxes.

\medskip\noindent
\textit{Forward seven arcs (\path{forward7}).}\label{reader:local:dp}\label{sh:arithmetic-forward}
The complete state graph and all sixteen journals are evaluated.
The two stages use the state inequalities of Lemmas~\ref{sh:state-filter}
and~\ref{sh:state-monotone} and the signed frame chart of the seven-arc
proof. The first stage gives the state restriction, including \(q_i<2\);
the second uses the recovered state and excludes the specified cost
counterexample. The two stages retain their distinct budgets.
This task uses the archive's x86-64 arithmetic guards.

\subsection{The three eight-arc covers}\label{reader:eight}
The tasks \path{R01}, \path{R02} and \path{R03} start from the
original domains of Lemma~\ref{eight:lem:end-cover} and
\eqref{mix:r02-root}--\eqref{mix:r03-root}, respectively. For R01 the
middle grid is included as well. Every change of chart, conservative
contraction and closed split respects the original candidate
inclusion. The terminal predicates for R01 include stationarity and
local-minimum tests for its moving-state cost. For R02 and R03 the
predicates concern closure, state, costs and their rigorous value
bounds. The denominators and guard conditions in the source formulas
are retained throughout. Theorems~\ref{eight:thm:main},
\ref{mix:thm:r02} and~\ref{mix:thm:r03} use the complete covers of
Propositions~\ref{eight:prop:finite} and~\ref{mix:finite-covers}.
The original root constructors, partition records, kernel assignments
and ordered terminal inputs are supplied by the three tasks. R01's
outward arithmetic and endpoint bounds are given in
Section~\ref{eight:app:arithmetic}; the other two covers retain the
Bernstein, Hermite and whole-box Taylor tests of
Section~\ref{mix:verification}.

\subsection{The four-dimensional R04 cover}\label{reader:r04}
The task \path{R04-full} replays the $31$ original runs from the single
root $[0,20]\times[0,1]^2\times[0,20]$ with $A\le E$. It evaluates the
interval predicates and checks the closed subdivisions and transfers
between runs. Proposition~\ref{r04:finite-cover} explains how this
covers the candidate box~\eqref{r04:root-box} and how
Proposition~\ref{r04:tail-bound} discharges the tail boxes.
The local deletion radii are the binary64 values of the input strings
\texttt{.0009} and \texttt{.00109}, both strictly below the radius
$11/10000$ of Proposition~\ref{r04:local-bound}. The task \path{R04-guards}
checks the supporting finite tests.

The archive also records an independent calculation in $93$ segments
whose $16$ roots are obtained by bisecting each coordinate of
$Q_{04}$. The complete replay inputs supplied in the archive are
those of the $31$-run cover used in Proposition~\ref{r04:finite-cover}.

\subsection{The R06 necessary rows and their finite check}
\label{reader:r06}\label{reader:r06:finite}

Read the original binary partition from its single root
\eqref{r06:root-box}.  A split replaces one coordinate interval $[l,h]$
by $[l,(l+h)/2]$ and $[(l+h)/2,h]$, both closed.  Require an exact interior
midpoint, both children, the complete stream, and every original
identifier once.  A historical terminal byte specifies a checking
obligation; it is not evidence for its truth.  The early necessary
comparisons, the original-value contractions, and the lifted
Gauss--Seidel contractions are those specified by the frozen sources.
Each acts only by an enclosure or intersection that retains the same
qualified point.  The complete enclosing $c$ interval is covered by
eight closed pieces when the fibre construction is used.  A surviving
piece retains its full fourteen-dimensional record.

For a Taylor step use the twelve smooth equalities
\eqref{r06:twelve-rows} and the single inequality $C_{06}-J_8\le0$.
Their first ten coordinates are positive; the four mean coordinates
are independent ambient variables linked by the four displayed actual
mean equations.  In particular, differentiate
$\widehat F=G/(1+s^2)$ on the whole positive ambient rectangle.  Do not
substitute $H^2-RL$ away from the frame or differentiate a
source-clipped value formula.

Let $x_*$ belong to the current closed rectangle $X$, and put
$d_j=\max\{x_{*,j}-X_j^-,X_j^+-x_{*,j}\}$.  Whole-box Hessian intervals
$\mathbb H_{jk}$ give, for each of the thirteen scalar functions $T$, the
Taylor remainder enclosure
\[
 \mathcal R_T=\frac12\sum_j\mathbb H_{jj}[0,d_j^2]
       +\sum_{j<k}\mathbb H_{jk}[-d_jd_k,d_jd_k].
\]
The integral Taylor formula on the entire segment in $X$ justifies this
enclosure.  Its mixed coefficient is one.  The function value and
gradient at $x_*$ and the Hessian on all of $X$ are recomputed; a
proposed Hessian is not an input to acceptance.

Write $x_j-x_{*,j}=d_jv_j$, $|v_j|\le1$, taking $v_j=0$ when $d_j=0$.
If $[f^-,f^+]$ encloses the value at $x_*$ plus $\mathcal R_T$, and
$[L_j,U_j]$ encloses the corresponding first derivative at $x_*$, set
\[
 a_j=(L_j+U_j)d_j/2,\qquad E=\sum_j(U_j-L_j)d_j/2.
\]
An equality yields $a\cdot v\le E-f^-$ and
$-a\cdot v\le E+f^+$.  The cost inequality yields only the first row.
Thus there are exactly $25$ rows, in the order
$f_0^+,f_0^-,\ldots,f_{11}^+,f_{11}^-,(C_{06}-J_8)^+$.
Each row is divided by its exact positive maximum absolute
coefficient or bound, using one when all are zero.

Read proposed weights as exact nonnegative rationals, each an integer
multiple of $2^{-48}$.  For their combined inequality
$A\cdot v\le b$, accept an empty box only when
\[
 -\sum_j|A_j|-b>0.
\]
Otherwise a nonzero coefficient may give a necessary coordinate bound;
intersect it with the previous box, round outward, and recompute the
next nonlinear enclosure.  A zero margin supplies no exclusion.
The same acceptance rules apply after a prescribed further subdivision,
with both complete closed children retained.  Unknown predicates,
unsafe arithmetic, missing children, or exhausted resources leave the
original unresolved.

The arithmetic profile uses GMP integer endpoints on scale $2^{-256}$,
integer-directed rational operations and square roots, and directed
MPFR logarithm endpoints at $512$ bits.  The support and chain readers
use exact rational arithmetic.  The actual-mean derivatives and
Hessians retain their proved series tails and integral enclosures.
Each numerical profile and its exact operation order are bound to the
corresponding original records.  An experimental faster kernel was
not used for the completed review.

\medskip
Proposition~\ref{r06:source-root} and Appendix~\ref{app:r06-source}
give the inclusion from the strict source into the root. The
\path{R06} task reconstructs all original endpoints from the
single root and recomputes the nonlinear enclosures, including on
previously accepted branches. Its exact rational checking binds
the support inequalities to the newly enclosed values and requires
complete closed continuation at every terminal. The original
\(806875\) leaves each carry their own complete obligation.
The conclusion is Theorem~\ref{r06:finite-cover}.
\label{reader:r06:evidence}

\subsection{Execution and input identity}\label{reader:compact}
The archive's \path{README.md}, \path{PAPER-PREMISES.md} and
\path{PACKAGE-MANIFEST.json} specify the execution modes, written
premises and file identities. Integrity and partition checks establish
the intended inputs and topology; evaluation recomputes the mathematical
terminal bounds. The archive binds every adopted specialized program
to its original arithmetic profile and exact inputs.
\label{reader:appendix-end}
\space
\section{The source implications for R05}\label{verify:r05-source-inputs}
The source collection for Proposition~\ref{nine:source-reduction} has three mathematical tasks. It identifies the actual physical graph and its necessary values, proves the bounded-domain comparisons, and carries every surviving candidate through closed covers to the fixed original-cell table. We record these interfaces here. The complete original mathematical texts and finite inputs are supplied in the accompanying mathematical verification material; its catalogue identifies their versions and exact source locations.

\subsection{The physical graph and the actual means}
The native reciprocals in \eqref{nine:native} belong to one actual trajectory with the genuine label twist. The original mean--frame argument identifies the three uniform reciprocal-shear integrals in \eqref{nine:three-means}, the complete signed frame equations, and the physical duration of each arc. The exterior owners $E$ and $Z$ remain distinct. The full-frame root-cone argument solves the exterior frame rows and selects the original calibrated inner rectangle; it does not infer physical variations merely from an injective parameter representation.

The local calibrated competitor theorem supplies the variation argument at a strict positive candidate. Its selected calibration root and positive central root depend smoothly on five parameters. The undivided exterior rows give the original positive exterior coordinates. All strict state and duration inequalities persist in a sufficiently small parameter neighbourhood. Three entries, including a nonzero central off-diagonal entry, and determinant one recover the entire signed frame equation. The reconstructed full frame paths vary continuously with fixed endpoints, so their real winding is unchanged. Lemma~\ref{fw:balanced-reconstruction} gives genuine competitors for the original cost, with the initial state allowed to vary subject to its prescribed twist.

Ordinary stationarity is applied only after this family has been constructed. One differentiates the full once-counted cost \eqref{nine:cost}, uses the integral identities for the actual means, and then imposes the two cap equalities at the candidate. This yields \eqref{nine:central-rows} with the stated $c$ and $f$ owners and the full signed contrast in $Q_c$. Nearby competitors are not required to satisfy the candidate's cap, density or cost filters. At a zero duration, an arbitrary auxiliary tuple does not acquire stationarity from continuity of its formulas; only the qualified limits retained by the source theorems are used.

\subsection{The bounded domain and the exact coefficient inputs}
The domain proof retains the $29$ mathematical sources numbered $00$--$28$ and the subsequent shape comparison. They concern the same actual source, with its compact state bounds, proper-source density, selected calibration, positive returns, two true caps, signed determinant and physical budget. In the fixed package, under \path{mathematical_inputs/r05-reduction/components/recovered/proofs/}, Section~7 of \path{r05_joint_residual_domain_research_20261111/sources/03_r05_shape_cap_floor_independent_audit_20261111.md} gives the final lower shape bound $49/125$ and the imbalance bound $327/1000$. The upper shape bound $2341/4000$ is the high-shape comparison in source~13. The actual own-pair and return comparisons give the scale and exterior-coordinate bounds of \eqref{nine:B0}. These restrictions enclose the minimizing candidates.

Two whole-domain coefficient arguments enter the central comparisons. To state their common function, use $\mathcal F,v_a,v_b$ from \eqref{nine:central-functions} and $g_1=\partial_1\phi$ for the potential \eqref{nine:potential}. Set
\[
 q=\frac{(y+c)(y+f)}{1+y^2},\qquad
 E_{\log}(q)=\frac{q-q^{-1}}2-\log q,\qquad
 \mathcal L=g_1(y,c)+g_1(y,f)-2yE_{\log}(q).
\]
Here $q\ge1$ on the scalar domains below. Since
$E_{\log}(1)=0$ and $E_{\log}'(q)=(q-1)^2/(2q^2)$, the complete logarithmic remainder is nonnegative, and is strictly positive for $q>1$. Define
\[
 g(y)=\sqrt{1+y^2}-y,\qquad c_0(y)=\frac{1-y^2}{2y}.
\]
We first verify the scalar domain at the actual source. The proper native pairs $[5,8]$ with phase $+1$ and $[6,9]$ without a phase are both $1201$. At their middle switches Lemma~\ref{fin:forward-force} bounds the third reciprocals
\[
 z=\frac{1-yc}{y+c}<\min\{y,c\},\qquad
 d=\frac{1-yf}{y+f}<\min\{y,f\}.
\]
They are positive, so $c,f<1/y$. The comparisons with $y$ give $c,f>c_0(y)$. The comparisons with $c,f$, together with the reciprocal identity, give $c^2+2yc>1$ and $f^2+2yf>1$, so $c,f>g(y)$.

The two-node defect comparison retained in D1 (Section 6 of the two-node source proof), applied to the unsquared sum of the actual jets, gives
$KD<29(a^2+b^2)/20-47/1250$. The inner mean is positive on its positive-width interval, so $K>3/2$. Since $F=D/(h_ah_b)$,
\[
 F<\frac{8}{3}\frac{29(a^2+b^2)/20-47/1250}{(1+a^2)(1+b^2)}
 \le\frac{26744}{30625}.
\]
The quotient increases in either anchor square on $[0,1/6]^2$, since the derivative numerator in the first square has sign $(29/20)(1-b^2)+47/1250>0$, and the other is symmetric. If $y\le9/20$, lower $c,f$ to $c_0(y)$ in $\mathcal F$. Indeed $\mathcal F=2(1+y^2)(q^2-1)/[(1+c^2)(1+f^2)]$, and its $c$ derivative has the sign of
\[
 \frac{(c+y)(f+y)^2}{(1+y^2)^2}(1-yc)+c>0
\]
throughout the comparison segment. The derivative in $f$ is symmetric. At the lower corner, with $Y=y^2$, this gives
\[
 F=\mathcal F>\mathcal F(c_0(y),c_0(y),y)
 =T_0(Y):=\frac{2(1-3Y)(1+5Y)}{(1+Y)^3}.
\]
Here $T_0'(Y)=2(-1-34Y+15Y^2)/(1+Y)^4<0$ for $0<Y\le81/400$, and
$T_0(81/400)=101108000/111284641>26744/30625$, a contradiction. Thus all the required scalar bounds hold at the same source:
\begin{equation}\label{r05:central-coefficient-domain}
 y>9/20,\qquad \max\{g(y),c_0(y)\}<c,f<1/y.
\end{equation}

On $y\ge9/20$ and $\max\{g(y),c_0(y)\}\le c,f\le1/y$, the first coefficient proof gives
\[
 \mathcal L\ge (v_a)_y+(v_b)_y+\frac{\mathcal F_y}{3}.
\]
On $y>0$ and $g(y)\le c,f\le1/y$, the second gives
\[
 \mathcal L\ge\frac34\bigl((v_a)_y+(v_b)_y\bigr)
       +\left(\frac12+\frac7{64}g(y)^2\mathcal F^2\right)\mathcal F_y.
\]
Their actual application uses the intersection of the two stated domains. Both estimate the same $\mathcal L$, so the remainder $2yE_{\log}(q)$ is restored once. It is not an additional physical payment.

For each inequality, the original rational difference is multiplied by its proved positive denominator and transformed to a closed polynomial comparison domain. The first proof divides at $y=1/\sqrt3$; the second uses $y=(1-k^2)/(2k)$, $0<k<1$, and a scalar exchange of $c,f$ to cover both ordered halves. This exchange is a polynomial symmetry, not a reflection of the physical trajectory. The complete exact inputs are
\begin{center}\small
\begin{tabular}{@{}llcr@{}}\toprule
Comparison&Scalar task&Tensor degrees&Bernstein coefficients\\\midrule
First, $g$ branch&\path{r05-coercivity-g}&$(24,6,6)$&$1\,225$\\
First, $c_0$ branch&\path{r05-coercivity-c0}&$(19,6,6)$&$980$\\
Second, ordered cube&\path{r05-sharp-reserve}&$(39,10,20)$&$9\,240$\\\bottomrule
\end{tabular}
\end{center}
The source proofs retain the full rational numerators, substitutions, all coefficient slots including zeros, and the inverse coefficient identities. Every coefficient in these three tensors is nonnegative. Since the Bernstein basis is also nonnegative on the closed comparison cubes, these exact coefficients prove the whole-box signs. The positive clearing factors are used only on their original domains; auxiliary polynomial endpoints do not license division by a vanishing factor. These are finite polynomial proofs of the stated scalar inequalities, not a cover of all R05 trajectories.

\subsection{From the source domain to the original cells}
\label{reader:r05}\label{reader:r05:reduction}\label{reader:r05:kernel}\label{reader:r05:cell-algorithm}
Together with the actual source comparisons in
Proposition~\ref{nine:source-reduction}, the following tasks give
the domain inclusion and the exclusion of every complete original
cell. A record is identified by its seed, continuation address,
all endpoints and original input. The retained records in each
row are precisely the input of the next applicable row.

\begin{center}\small
\begin{tabular}{@{}p{.23\textwidth}p{.69\textwidth}@{}}
\toprule Task & Input and verified conclusion\\\midrule
\path{R05-coefficients}, \path{R05-domain}
 & The three coefficient tensors above and the necessary bounds for
 the five coordinates in \(B_{05}\).\\
\path{R05-boundary}
 & The complete closed return thickness \([0,1/1000]\) with the
 undivided frame equations, and the same-side owner comparisons.
 The actual return widths exceed \(1/1000\).\\
\path{R05-root}
 & The closed root \(B_{05}\). The retained \(12208\) original
 terminal records are the continuation seeds.\\
\path{R05-continuation}
 & The closed forest from exactly those seeds. Its \(58274\)
 retained records are the overlay input.\\
\path{R05-half-sector}
 & The whole sector \(F\ge1/2\), including equality. The cost
 comparison counts the inner reserve once and uses both endpoints
 in its concave area parameter.\\
\path{R05-overlay}
 & The same continuation records with the half-\(F\) test; the
 surviving \(46383\) records retain their full original endpoints.\\
\path{R05-final}
 & All \(46383\) original cells. For each cell the necessary values
 are enclosed on the same candidate and contradicted by a full-cell
 test or a complete closed internal cover.\\
\bottomrule
\end{tabular}
\end{center}

In the overlay the conservative collar is \([1/20,1/2]\).
Its endpoints are kept in the enlarged comparison. The half-sector
cover uses outward integers at scale \(2^{-160}\) and its
\(96\)-term logarithm bound. The final original-cell tests use the
necessary equations of Section~\ref{nine:system}, the contractions
of Section~\ref{nine:separation}, and the \(2^{-128}\) arithmetic
and \(80\)-term logarithm bound of Section~\ref{nine:finite-cover}.
Every closed split retains both children, so the inclusion and
the final exclusion follow by induction using the actual roots and
terminal inequalities. The same candidate passes through all the
stages. Proposition~\ref{nine:source-reduction} supplies its
source inclusion, and Lemma~\ref{nine:exhaustion} uses the final
exclusion.
\space
\label{revision:appendix-end-11}
\end{proofappendix}
\clearpage%

\space
\end{document}